\documentclass[a4paper]{article}
\usepackage{amsmath,amsfonts,amssymb}
\usepackage{graphicx, epsfig}

\newtheorem{theorem}{Theorem}[section]
\newtheorem{Pseudo-Theorem}{Pseudo-Theorem}[section]

\newtheorem{lemma}[theorem]{Lemma}
\newtheorem{cor}[theorem]{Corollary}
\newtheorem{prop}[theorem]{Proposition}
\newtheorem{defn}[theorem]{Definition}
\newtheorem{ques}[theorem]{Question}
\newtheorem{conj}[theorem]{Conjecture}

\newtheorem{rem}[theorem]{Remark}

\newtheorem{prob}[theorem]{Problem}
\newtheorem{Numerology}[theorem]{Numerology}

\newtheorem{Scholium}[theorem]{Scholium} 
\newtheorem{Notation}[theorem]{Notation}

\newtheorem{quota}[theorem]{Quote}

\newtheorem{principle}[theorem]{Principle}

\newenvironment{proof}[1][Proof]{\textbf{#1.} }
{\hfill\rule{0.5em}{0.5em}\medskip}
\newenvironment{proof*}[1][Proof]{\textbf{#1.} }{}

\def\epsilon{\varepsilon} 
\def\loccit{{\it loc.\,cit.}}
\def\RR{{\Bbb R}} 
\def\CC{{\Bbb C}} 
\def\PP{{\Bbb P}} 

\def\ZZ{{\Bbb Z}} 

\def\disc{{\frak D}} 

\def\mH{{\vert mH \vert}} 

\def\la{\langle} 
\def\ra{\rangle}
\def\vc{\sqcup}

\begin{document}


\title{Ahlfors circle maps and
total reality:\\
from Riemann to Rohlin}


\author{Alexandre Gabard}
\maketitle

\vskip-15pt

\newbox\quotation
\setbox\quotation\vtop{\hsize
7.8cm \noindent

\footnotesize {\it Nihil est in infinito quod non prius fuerit
in finito.}

\noindent Andr\'e Bloch 1926
{\rm \cite{Bloch_1926_EM}, \cite{Bloch_1926_Mem}}.




}

\hfill{\hbox{\copy\quotation}}

\medskip
\newbox\abstract
\setbox\abstract\vtop{\hsize 12.2cm \noindent


\noindent\textsc{Abstract.} This is a prejudiced survey on the
Ahlfors (extremal) function and the weaker {\it circle maps}
(Garabedian-Schiffer's translation of ``Kreisabbildung''), i.e.
those (branched) maps effecting the conformal representation upon
the disc of a {\it compact bordered Riemann surface\/}. The theory
in question has some well-known intersection with real algebraic
geometry, especially Klein's ortho\-symmetric curves via the
paradigm of {\it total reality}. This leads to a gallery of
pictures quite pleasant to visit of which we have attempted to
trace the simplest representatives. This drifted us toward some
electrodynamic motions along real circuits of dividing curves
perhaps reminiscent of Kepler's planetary motions along ellipses.
The ultimate origin of circle maps is of course to be traced back
to Riemann's Thesis 1851 as well as his 1857 Nachlass. Apart from
an abrupt claim by Teichm\"uller 1941 that everything is to be
found in Klein (what we failed to assess on printed evidence), the
pivotal contribution belongs to Ahlfors 1950 supplying an
existence-proof of circle maps, as well as an analysis of an
allied function-theoretic extremal problem. Works by Yamada
1978--2001, Gouma 1998 and Coppens 2011 suggest sharper degree
controls than available in Ahlfors' era. Accordingly, our partisan
belief is that much remains to be clarified regarding the
foundation and optimal control of Ahlfors circle maps. The game of
sharp estimation may look  narrow-minded
``Absch\"atzungsmathematik'' alike, yet the philosophical outcome
is as usual to contemplate how conformal  and algebraic geometry
are fighting together for the soul of Riemann surfaces. A second
part explores the connection with Hilbert's 16th as envisioned by
Rohlin 1978.

 }

\centerline{\hbox{\copy\abstract}}

{\small \tableofcontents}

\section{Introduction}\label{sec1}


{\it Preliminary Warning.} [13.11.12]---Despite its exorbitant
size, the actual mathematical content of the present text is very
limited. It focuses primarily on the Ahlfors map. Neither does our
work
have the pretence of being  the logical sum of all knowledge
accumulated in the past, nor
will it give an accurate picture of real developments taking shape
contemporaneously. Our intention was rather more
to delineate a reasonably clear-cut perception of the  early
branches of the theory as to understand objectively the basic
truths making Ahlfors theorem possible. Failing  systematically,
our pretence converted to that of throwing enough obscurantism
on the whole theory as to motivate others to shed  fresh lights
over the edifice. Even the
primary contribution to the field (that of Ahlfors 1950
\cite{Ahlfors_1950}) has not yet been fully assimilated by the
writer (compare optionally Section \ref{Ahlfors-proof:sec} for our
fragmentary comprehension). We strongly encourage mathematicians
having a complete mental picture of Ahlfors proof to publish yet
another account helping to clarify the original one. We hope
during the next months (or years) to be gradually able to improve
the overall organization of this text, in case our understanding
of classical results sharpens. All of our ramblings starts
essentially in the big-bang of Riemann's Thesis. It looks almost a
triviality alike to expect that subsequent developments
will involve a deeper interpenetration between the
conformal and algebro-geometric viewpoints. One oft encounters in
the field problems requiring serious combinatorial skills or
geometric intuition. For instance how does the moduli space of
bordered surfaces stratifies along gonalities;
Sec.\,\ref{sec:profile-histogram} guesses some scenarios via
primitive methods.  Riemann surfaces or the allied projective
realizations offer an ornithological  paradise requiring
patience and observational skills from the investigator. This is
especially stringent when the complexes are traded against the
real number field, and inside this universe of $3g-3$ real
dimensions one encounters with probability $1/3$ the so-called
real orthosymmetric curves of Felix Klein (1876--1882)
subsumed to the paradigm of {\it total reality\/}. This little
third is actually all what our topic of the Ahlfors map is about.
Last but not least,
experimental studies point to a large armada of potential
counter-examples menacing the improved bound $r+p$ announced in
Gabard 2006 \cite{Gabard_2006}. It seems safe to declare as an
open problem to either corroborate this bound (via other more
analytic or algebraic treatments) or to
reject it.

\smallskip
{\small

[11.04.13] {\it Glossary of synonyms.} (all coined by Klein
between 1876 and 1882).

$\bullet$ Real algebraic curve=symmetric Riemann surface in the
sense of Klein, i.e. with an anti-holomorphic involution, usually
induced by the Galois symmetry of Tartaglia 1535.

$\bullet$ Type~I=orthosymmetric=dividing=separating, when the real
locus of  a real algebraic curve (equivalently a symmetric Riemann
surface) separates the complexification.

$\bullet$ Type~II=diasymmetric=nondividing, when the contrary
occurs, for instance when there is no real points.

}

\smallskip

[19.03.13] The main addition to
the present edition (v.2) of our text, is an
essay to connect Ahlfors' theory with Hilbert's 16th problem (at
least the part thereof pertaining to the topology of real plane
algebraic curves). Again our dancing queen is the paradigm of
total reality (Riemann, Schottky, Klein, Bieberbach,
Teichm\"uller, Ahlfors, Alling-Greenleaf, Geyer-Martens, etc.),
but now  reoriented as a
missile
against Hilbert's 16th
(quite akin to an asteroid menacing peaceful life on  planet
Earth). This trend is not new having been much foreseen in
Rohlin's seminal work 1978 \cite{Rohlin_1978} effecting a
Verschmelzung\footnote{=Fusion in Klein's prose when viewing all
his work (and that of Sophus Lie) as being merely a Galois-Riemann
synthesis.} between the conceptions of Klein and Hilbert (when it
comes to real geometry).
Rohlin never refers back to Ahlfors' work (which he probably
ignored?), yet the connection is very vivid through Rohlin's
(conjectural) philosophy that {\it schemes\/} of type~I (not in
Grothendieck's highbrow sense, but merely Rohlin's
synonym for a {\it distribution of ovals\/} \`a la
Zeuthen-Harnack-Hilbert) are necessarily {\it maximal\/} in the
hierarchy of all schemes of some fixed degree. Through the lines
of Rohlin's text
transpires
the
intuition that what detects (pure) orthosymmetry of schemes is a
vertiginous phenomenon of {\it total reality\/} positing
 existence of  adjoint pencils cutting only real points on the
given curve. The
byproduct
is that total reality should act prohibitively upon all
schemes enlarging those totally flashed by a pencil, which are
so-to-speak already B\'ezout-saturated. Hence total reality should
contribute to Hilbert's problem (isotopic classification of
curves), though this method
can hardly be said to have been systematically exploited as
yet (apart of course in
the prophetical allusions in Rohlin 1978).
%

\smallskip {\small

{\it Added in proof} [11.04.13].---Part of this recalcitrance
may be imputed to the fact that 4D-topology \`a la Rohlin
(1951--72) and the resulting (trinity of) congruences modulo 8
(suspected by Gudkov 1969) seemed to imply all what
seemed desirable to know,  e.g. a type~I criterion for
$(M-2)$-schemes with $\chi\equiv k^2+4 \pmod 8$. (Standard
notation: $M=g+1$=Harnack's bound, $\chi$=Euler characteristic of
the Ragsdale orientable membrane bounding the ovals from
``inside'', and $k=m/2$=semi-degree of an even order curve.)
However we should probably return to the geometric substance to
gain more namely
maximality of such schemes. Further, there is (conjecturally)
another source of total reality (``hence'' maximality) coming from
the operation of satellites (amounting merely to replicate the
curve within a tube neighborhood of it up to a certain
multiplicity $2,3, etc$). This promises to offer ``new''
obstructions that were perhaps slightly denigrated/overlooked by
Rohlin in 1978; compare his prose {\it ``However, all the schemes
that we have so far succeeded in coping with by means of these
devices {\rm [=total reality]} are covered by Theorem~3.4{\rm
[=Kharlamov's congruence (unpublished in 1978)]} and 3.5{\rm
[=extremal properties of the strong Arnold inequalities
(Zvonilov-Wilson)].}''} To be specific we suspect that the scheme
of degree 10, which occurs as 2nd satellite of Harnack's quintic
(of Gudkov symbol like $(1,6\times \frac{1}{1})$, i.e. $6$ nests
of depth 2 enveloped in a larger oval,
Fig.\,\ref{satellite-of-Harnack's-quintic:fig}) is maximal, i.e.
cannot be enlarged algebraically. Likewise in degree 12 we suspect
that the 2nd satellite of any one  of
both Rohlin's schemes (symbols $\frac{6}{1}2$ and $\frac{2}{1}6$)
are maximal. Similarly, the 2nd satellites of any of the three
$M$-schemes of degree 6 (Harnack's $\frac{1}{1}9$, Hilbert's
$\frac{9}{1}1$ and the glamorous Gudkov type $\frac{5}{1}5$)
should be maximal in degree 12. Of course we conjecture generally
a stability of total reality under higher satellites (hence of
type~I and maximal schemes too), but our purpose here is to give
the lowest degree examples in the hope that someone (presumably a
patchworker) can refute our stability prediction.

}
\smallskip

This grand vision of Rohlin could benefit from the Klein-Ahlfors
theory (which in our opinion has been much neglected in the
tradition of the German-US-Italian-Russian school of real geometry
involving such
{\it pointures\/}
as Hilbert 1891, Rohn 1888--1913, Ragsdale 1906, Brusotti
1910--50, Petrovskii 1933/38, Gudkov 1954--69, Arnold 1971, Rohlin
1972--78, Kharlamov, Viro, Marin, Fiedler, Shustin, Itenberg,
etc.). Alternatively it could be reassessed through purely
synthetical procedures that Rohlin himself envisioned
(probably as consequences of deep
4D-topology, notably the type-I-forcing Kharlamov-Marin congruence
modulo 8). Alas Rohlin's synthetical proof even on the simplest
prototype of sextics has never been published (and seems now to be
lost forever), but was recently (partially)
resuscitated in a tour de force of S\'everine Le~Touz\'e 2013
\cite{Fiedler-Le-Touzé_2013-Totally-real-pencils-Cubics}. So
the epoch seems ripe to dream about a big Ahlfors-Rohlin
Verschmelzung
with direct repercussions upon Hilbert's 16th byway of
prohibitions. It is already breathtaking (for a beginner) to
contemplate how the Rohlin-Le~Touz\'e total reality phenomenon for
sextics explains nearly all prohibitions observed in Gudkov's
census (1969) of sextics curves (cf.
Fig.\,\ref{Gudkov-TableTop:fig}), which after all supplies
(nothing less than) the {\it complete\/} solution to Hilbert's
problem in its original formulation (degree $m=6$). Those aspects
are addressed in Sec.\,\ref{Klein-Rohlin-conj:sec}\&ff. For
convenience we wrote a general overview in
Sec.\,\ref{Hilbert's16th-PartII:sec} pointing to some open
questions which
looks semi-urgent to settle in order to build a more solid theory.
All this second part, devoted to Hilbert's 16th, could not have
been written without the constant support and information
generously shared by the leading experts (Viro, Marin, Kharlamov,
Shustin, Orevkov, Le~Touz\'e, Fiedler), whose instructive e-mails
are reproduced in Sec.\,\ref{e-mail-Viro:sec}. Needless to say we
have not yet assimilated all their wisdoms and advices, but have
reproduced faithfully their messages in the hope that other
amateurs of the field can also beneficially profit from their
invaluable expertise.

[11.04.13] It seems (now) a firm conviction that a big piece of
Hilbert's 16th puzzle still remains to be fixed. This should be a
fairly simple matter of assembly between the conceptions of
Riemann-Schottky-Klein-Bieberbach-Teichm\"uller-Ahlfors about
total reality and the theory of
Hilbert-Rohn-Petrovskii-Gudkov-Arnold-Rohlin aiming to predict the
distribution of ovals
traced by algebraic curves (a God-given video game). Precisely,
the isotopic classification of real plane
curves should be
regulated by a sole paradigm (total reality) itself
piloted by the geometry of the canonical series (adjoint curves of
order $m-3$) assigned to visit $(M-3)$ basepoints randomly
selected among the most profound ovals of $(M-2)$-curves (alias
the {\it extended Rohlin-Le~Touz\'e phenomenon\/}, which in degree
4 boils down to the total reality of the G\"urtelkurve quartic
with 2 nested ovals). This looks special but wait a moment.

It is conjectured (\ref{primitive-manifestation-of-tot-real:conj})
that {\it any\/} {\it primitive\/} manifestation of the phenomenon
of total reality on a plane curve
 is of this sort
(i.e. a Rohlin-Le~Touz\'e ``adjunction'' for $(M-2)$-curves
subsumed to the eightfold periodicity $\chi\equiv_8 k^2+4$ of
Rohlin-Kharlamov-Marin), except when it comes to the trivial case
of $M$-curves, where total reality is nearly completely settled by
an extension of (another) Le~Touz\'e's scholium
(\ref{Le-Touzé-extended-in-odd-degree:scholium}). The maximality
of $M$-schemes being so trivial (Harnack-Klein inequality of
1876), this latter case looks a sterile syllogism not worth paying
attention at,
but this is probably not so via satellites.


If {\it not primitive\/}, this is to say that  the scheme is just
derived
as a {\it satellite\/} replicating the curve  up to a certain
multiplicity within its tube neighborhood. For instance satellites
of a single oval of degree 2
reproduce the infinite series of deep-nests total under a pencil
of lines [by the way highly reminiscent to the ``rondelles'' of a
certain artist known as Markus Schneider-Zeitler, Jura Suisse].
Such deep-nests are B\'ezout-saturated hence extremal shapes in
the Hilbert-Gudkov hierarchy.
More generally, the magic formula reads $A+B=R m c^2$, i.e.
Ahlfors plus B\'ezout implies Rohlin's maximality conjecture
(any scheme of {\it type~I\/} kills all its enlargements). At this
stage the architecture of higher Gudkov's pyramids ($m\ge 7$) is
completely
predestined by (the felicity of)
%
%
 pure orthosymmetry \`a la Felix
Klein, and (optionally) Rohlin's theorem on the signature of spin
$4$-manifolds governing the (Gudkov)
$8$-fold periodicity via differential-topology. Paraphrasing,
Hilbert's 16th is virtually solved in {\it all\/} degrees, at
least in its qualitative shape (prohibitions). It remains then of
course to programme a (patchworking) machine doing all the
constructions. This should be merely a matter of passive
contemplation, requiring immortality and much patience from the
investigator. It is
evident that all this programme is not a novel idea, but
much---not to say completely---anticipated by Rohlin 1978, safe
that he does not seem to have been consciously aware of the
Riemann-Ahlfors theory (nor perhaps the conjectural stability
under satellites), but seemed rather
adventurous enough  to rediscover it {\it ab ovo\/} through purely
synthetical processes, without
any intrusion of  analysis or transcendental gadgets, like Abelian
integrals. Our
messy text is just
an invitation to inspect more exactly how the
whole process will sediment itself
within the next
decades, probably via massive usage of Brill-Noether (i.e.,
Riemann for dummies). Then, the night of ignorance {\it (les
t\'en\`ebres de l'ignorance)}
allied to Hilbert's combinatorial mess about the distribution of
16 ovals (recall that $M_7=11+5=16$) should  be completely
dissipated. Whether conversely, all this (ancient) geometry can
acknowledge
in feedback some impact upon
4D-differential-topology, e.g. the question of smooth structures
on $S^4$ or $\CC P^2$ (so-called {\it smooth Poincar\'e
conjecture\/}) is merely speculation of longstanding (reminding
such names as Arnold, Maxwell, Milnor, Kuiper, Massey, Marin,
Akbulut, Donaldson, Taubes, Finashin, Wang, Seiberg, Witten,
etc.). But this is another story.

[22.03.13] Lastly, we adopted (not deliberately but because we
were not clever enough to proceed differently) Arnold's philosophy
of the mushroom (compare Arnold 2004 \cite{Arnold_2004}). That is,
not just presenting overwhelming theorems (arid as they are) but
the slow organical eclosion of truths through mistakes,
conjectures, historical meanders, etc. The drawback is an
intolerable  inflation in size, also partly caused by the
abundance of pictures, which in our opinion form the true core of
any mathematical
truth\footnote{Best example thereof, read Borsuk's article ca.
1936 where a contractible compactum  lacking the fixed-point
property is presented. If you have just the boring (unreadable)
formulas  of Borsuk you understand nothing, but if  you know the
picture that the space in question is a crumpled-cube spiraling
twice around itself as pictured by Bing, you
start to believe why the fact holds true.}. We expect in the
future to reorganize the material \`a la Bourbaki as to offer a
cleaner view of what
happens (after
distillation, the
factual content
should be compressible to ca. 20 pages). Especially crucial is a
rectification due to Fiedler of an erroneous theorem of mine that
would have proved one-half of the (still open) Ragsdale conjecture
for $M$-curves via the Thom conjecture (Sec.\,\ref{Thom:sec}).

\centerline{$\star$$\star$$\star$} \medskip

This is a prejudiced survey on the  Ahlfors (extremal) function
and (improvising terminology) the weaker {\it circle maps},
effecting the conformal representation  upon the disc of an
arbitrary differential-geometric membrane, alias {\it compact
bordered Riemann surface\/}. Our jargon,
 borrowed from Garabedian-Schiffer 1950
\cite{Garabedian-Schiffer_1950},  translates essentially the term
{\it Kreisabbildung} used e.g., by Koebe 1915 \cite{Koebe_1915}
and Bieberbach 1914 \cite{Bieberbach_1914}.

Exciting works by Yamada 1978--2001 \cite{Yamada_1978},
\cite{Yamada_2001}, Gouma 1998 \cite{Gouma_1998} and Coppens 2011
\cite{Coppens_2011} suggest that fewer sheets than required in
Ahlfors' era is expectable, for a clever placement
of the basepoint(s) required to pose the extremal problem.
E.g., is Coppens' (absolute) gonality of a membrane always
sustained by an Ahlfors function? We also started tabulating a
list of known applications in the hope of guessing future ones.
Some applications (e.g. Fraser-Schoen's recent one to Steklov
eigenvalues \cite{Fraser-Schoen_2011}) do not require the full
punch of Ahlfors' extremals, raising the hope that  the improved
control $r+p$ on the degree of circle maps (predicted in Gabard's
Thesis 2004/06 \cite{Gabard_2006}  for surfaces of genus $p$ with
$r$ contours) could imply some `automatic' upgrades (e.g. in the
corona problem with bounds, as studied by Hara-Nakai 1985
\cite{Hara-Nakai_1985}).

As to the foundation of the Ahlfors mapping theory itself, the
issue that the naive qualitative approach (used in Gabard 2004/06
\cite{Gabard_2006}) affords a bound, $r+p$, quantitatively
stronger than Ahlfors' original $r+2p$ is somewhat surprising. It
results a certain psychological tension between topological and
analytical methods, which  hopefully is just a superficial and
temporary state of affairs destined to disappear after renewed
examination of Ahlfors' argument. The latter seems indeed to leave
some free man{\oe}uvring room,  in its ultimate convex geometry
portion (cf. Sec.\,\ref{Ahlfors-proof:sec} for some strategy).

It is our partisan belief that much remains to be clarified both
historically and logically in the theory of the Ahlfors map.
Albeit sembling a retrograde attitude, it is probably not since
Ahlfors bound $r+2p$ certainly fails sharpness, at least for low
values of the invariants $(r,p)$. (Consider for instance the
topological type of Klein's G\"urtelkurve; i.e. $(r,p)=(2,1)$
where a projective realization (of the Schottky double) as a plane
quartic with 2 nested ovals prompts  existence of a total map of
degree $3$ via projection from the inner oval. This beats by one
unit Ahlfors' bound $r+2p=4$.)

Apart from an abrupt claim by Teichm\"uller 1941
\cite{Teichmueller_1941}, that everything (safe bounds) is to be
found in Klein (what the writer was unable to certify from printed
evidence), it is fair to admit that the bulk of the theory
crystallized right after World War II. Several workers like
Ahlfors 1948/50 \cite{Ahlfors_1950}, Matildi 1945/48
\cite{Matildi_1945/48}, Andreotti 1950 \cite{Andreotti_1950},
Heins 1950 \cite{Heins_1950} (perhaps even Courant 1939/40
\cite{Courant_1939}, not to mention Grunsky 1937--40--41--42
\cite{Grunsky_1937}, one of the most brilliant protagonist albeit
his work looks confined to the genus $0$ case) offered quite
overlapping conclusions. It seems fair however to give full credit
to Ahlfors for having first expressed the story in the most
clear-cut fashion. Quite shamefully, I confess that Ahlfors
argument  still escapes me slightly. A non-negligible amount of
literature is devoted to reproving Ahlfors' theorem: Heins
1950/75/85 \cite{Heins_1950} \cite{Heins_1975}
\cite{Heins_1985-Extreme-normalized-LIKE-AHLF}, Garabedian 1950
\cite{Garabedian_1950}, Kuramochi 1952 \cite{Kuramochi_1952}, Read
1958 \cite{Read_1958_Acta} (student of Ahlfors), Mizumoto 1960
\cite{Mizumoto_1960} (topological methods), Royden 1962
\cite{Royden_1962} (Hahn-Banach like Read), Forelli 1979
\cite{Forelli_1979} (extreme points and Poisson integral),
Jenkins-Suita 1979 \cite{Jenkins-Suita_1979} (Pick-Nevanlinna
viewpoint), just to name those authors addressing the positive
genus  case ($p>0$).


Another promising route is Meis' work 1960 \cite{Meis_1960}
validating Riemann's (semi)intuition of the $[\frac{g+3}{2}]$
gonality of closed genus $g$ surfaces via some
Teich\-m\"uller-theoretic
background. It is likely
 that Meis' approach is transmutable to the
bordered setting, reassessing thereby Ahlfors' result (probably
even with the sharp bound $r+p$ in case the latter is reliable).
To put it briefly, it seems that the Gr\"otzsch-Teichm\"uller
mode-of-thinking (of the {\it m\"oglichst konform} mapping) has
not yet fully penetrated the paradigm of the Ahlfors circle map,
more generally that of branched coverings, except of course in
Meis' memoir (alas notoriously difficult to access). Dually, it
also seems desirable to reprove the Riemann-Meis bound via
topological methods (e.g. that used in Gabard 2006
\cite{Gabard_2006}, which perhaps is nothing else than Riemann's
parallelogram method). Poincar\'e's ``Analysis Situs'' (1895
\cite{Poincare_1895-Analysis-Situs}) invented ``homology'' (modulo
the Riemann-Betti=Brioschi [sic!] heritage) with precisely
function theory (Abelian functions) as one of the key motivation
(beside celestial mechanics and the like). This, jointly with the
subsequent work of Brouwer, gives the basic
conceptual
framework
for implementing such topological methods.

{\it User guide.}---This draft is a preliminary version, so avoid
printing it for environmental reasons. A list of hopefully
clear-cut questions is given in Sec.\,\ref{sec:question}. This is
intended to challenge investigators. Several synoptic diagrams
scattered as figures through the text should permit a quick
optical scan of the whole content. More specifically, those
includes:

$\bullet$ an {\it exhaustive\/} list (Fig.\,\ref{Map:fig}) of {\it
all\/} articles supplying (or claiming to supply) a proof of
Ahlfors theorem (existence of circle maps),

$\bullet$ a list of keywords (Fig.\,\ref{Keyword:fig}) tabulating
concepts traditionally related to the Ahlfors map,

$\bullet$ a comprehensive  map (Fig.\,\ref{Geneal:fig}) of authors
involved in the theory (at least those cited in the bibliography).

This essay, as already said, contains no original insights,
instead a series of attempts to contemplate the theory from
different angles. A commented bibliography (of ca. 900 entries)
tries to brush a panorama  of trends related to  the Ahlfors map.
This includes topics like Riemann surfaces, algebraic curves,
conformal mapping, potential theory, Green's functions,
Dirichlet's principle, Riemann mapping theorem, Kreisnormierung,
parallel slit-maps, Bieberbach's least-area map interpretation of
the Riemann map, Bergman and Szeg\"o kernels, minimal surfaces,
Plateau's problem, spectral theory, analytic capacity, removable
singularities, corona problem, operator theory, Gromov's filling
conjecture, etc.). We have not attempted to reach any overwhelming
mathematical density, but rather tried to dilute through
historico-philosophical anecdotes.

There is some interplay between Ahlfors maps and  total reality of
Klein's orthosymmetric curves which gives rise to the gallery of
pictures mentioned in the abstract. For a tourist view, browse the
string of figures starting from Fig.\,\ref{Pencil:fig} up to
Fig.\,\ref{Fcubic3:fig}. For ``do-it-yourself'' purposes, it is
probably more valuable to describe the general recipe used to
manufacture such pictures. Take any configuration of simple
objects like lines and conics, and smooth it in an
orientation-preserving sense to get a dividing curve (one is free
to keep certain nodes unsmoothed). (Rohlin's eminent student
Thomas Fiedler (1981 \cite{Fiedler_1981}) ensures for us that the
smoothed curve is dividing, alias orthosymmetric in the sense of
Klein.) According to Ahlfors theorem there must be a totally real
pencil of auxiliary curves cutting only real points on the given
curve (plus maybe some imaginary conjugate basepoints). Geometric
intuition usually tells us where to locate such a total pencil,
roughly by assigning basepoints among the {\it deepest\/} ovals
(in the sense of D. Hilbert's 16th Problem). Albeit this is just a
Plato cavern style extrinsic manifestation of Ahlfors theorem, the
possibility of finding {\it always\/} such a total pencil reveals
strikingly (in our opinion) some of the depth of Ahlfors theorem.
(Incidentally it is not to be excluded that a deep understanding
of extrinsic algebraic geometry (say \`a la Brill-Noether) could
reprove the full Ahlfors theorem from within the Plato cavern.) In
philosophical terms, {\it real orthosymmetric curves behave on the
reals as if
they were complex
varieties\/}: all intersections prompted by B\'ezout are visible
over the reals. This phenomenon is what we (and others, e.g.
Geyer-Martens 1977 \cite{Geyer-Martens_1977}) call the paradigm of
{\it total reality}. It seems evident that a global study of such
pencils bears some close connection with Poincar\'e index theory,
foliations \`a la Poincar\'e-Kneser-Ehresmann-Reeb, etc., and that
both experimentally and  theoretically much remains to be explored
along the way. In particular we failed to make such totally real
pictures for an $M$-quintic
(Sec.\,\ref{sec:Total-reality-Harnack-max-case}). This could be a
challenging problem of computer visualization.

As to our speculation about a mechanical interpretation of the
Klein-Ahlfors theory of real orthosymmetric curves (and the allied
totally real maps) in terms of gravitational systems, see
Sec.\,\ref{sec:gravitation}. This posits a
broad extension of Kepler's planetary motions around ellipses,
enabling virtually all algebraic curves (and not just conic
sections) to arise as the trajectories of a perfectly stable and
periodic motion.  Of course if such a grandiose connection between
Klein-Ahlfors and Kepler-Newton-Coulomb-Poincar\'e is not
verifiable, this may just be interpreted as a metaphoric language
describing the dynamics of totally real morphisms prompted by
Ahlfors theorem. In fact rather than mere gravitation, it is
really a ``{\it dynamique de l'\'el\'ectron}'' which seems to be
involved; for a toy example on the G\"urtelkurve compare
Fig.\,\ref{FGuert:fig}. The resulting
metaphysics is quite akin to Lord Kelvin's speculations about the
ultimate constitution (and stability) of matter (via the vortex
atom geometrized by a knot), except that in our story the dancing
queen is rather a (naked) bordered Riemann surface.

\subsection{Trying to wet some appetite
out of the blue}

A long time ago (ICM 1908 Rome), Poincar\'e argued that in
mathematics we need a strong principle of economy of thoughts by
conceptualizing such notions as  `uniform convergence' as if the
sole naming
process would spare us repeating long intricate arguments. On the
other hand,
Felix Klein, asserted boldly ``die Franzosen unhistorisch wie Sie
sind'' (exercise recover the source) and liked
the
motto ``Zur\"uck zur Natur, sie bleibt die gr\"o{\ss}te
Lehrmeisterin''. Beside all those psychological tips of the
masters of geometrization, we can safely agree with both of them
that
science
requires---as a
matter of conciliating the principle of economy
with that of historical continuity (of course not so
structurally incompatible as
neo-expressionism seems to assess)
---a certain amount of
respectfulness about wisdoms accumulated during the past. This
explains our ca. 900 references (albeit the explosion was mainly
caused by my lack of internet connection occasioning a manual
references chasing).

In these
notes we propose a
(poorly guided) tour of some geometric function theory (GFT). The
field is an old fashioned one, lying quite dormant with its old
mysteries and legends (e.g., Koebe's Kreisnormierungsprinzip, the
exact determination of the Bloch constant in
quasi-stagnation since Ahlfors-Grunsky 1937
\cite{Ahlfors-Grunsky_1937}, etc.). Function-theory seems a
volcano alike awaiting
anxiously the next explosive eruption, whose
pyroclastic rejections turned out
to act (in the past at least) as a powerful fertilizer over
neighbouring
areas (like Riemannian and algebraic geometry, spectral theory,
etc.). Actually Koebe
had a more picturesque description, when proclaiming (im September
1921, Jena, Jahresversammlung der DMV\footnote{Source=H. Cremer,
Erinnerungen an Paul Koebe, Jahresber. DMV, 1968, p.\,160.
(Mitteilung von Heinrich Behnke).}): ``Es gibt viele Gebiete in
der Mathematik, wo man sich durch Entdecken neuer Ergebnisse
verdient machen kann. Es sind meistens lange und steile
Gebirgsh\"ange f\"ur meckernde Ziegen. Die Funktionentheorie ist
aber mit einem saftigen Marschland zu vergleichen, besonders
geeignet f\"ur dickes
Rindvieh!''

The field itself (GFT) seems to be a strange cocktail of
qualitative-flexible versus quantitative tricks, or as Gauss puts
it {\it geometria situs} versus {\it geometria magnitudinis}.
If
topological
methods look a priori
quite foreign to
the
discipline, it was probably Riemann who first
revealed:

$\bullet$ the
reactivity
of the underlying
topological substratum
(anticipated
maybe by Abel 1826 \cite{Abel_1826}, who first introduced the {\it
genus} (under a different name and  the transcendant
disguise
of differentials of the first
kind). [The word {\it Geschlecht\/}
is first coined in Clebsch 1865 \cite[p.\,43]{Clebsch_1865};
and the allied {\it Geschlechtsverkehr}\footnote{For more
historical details on the theory of quasiconformal mappings
compare Ahlfors 1984 \cite{Ahlfors_1984-The-Joy} or Lehto 1998
\cite{Lehto_1998}. [02.10.12] Alas we were not as yet able to show
any deep connection between the theory of  Ahlfors circle maps and
that of quasiconformal maps, yet it is not unlikely that such a
connection is worth studying, more in Section \ref{sec:question}.
} must have originated about the same period]

$\bullet$ the
amazing
plasticity (inherited from potential-theoretic considerations) of
2D-conformal mappings,
leaving out moduli spaces of finite dimensionality
after conformal evaporation of
all
metrical incarnation of a given surface. [Gromov wrote in 1999
\cite{Gromov_1999}: {\it Shall we ever reach spaces beyond
Riemann's imagination?}]


Our
text will soon be biased
toward a single obsession, the so-called {\it Ahlfors
function}, which
is one (among several other possible)
generalisation of the {\it Riemann mapping theorem} (RMT) to
configurations of higher topological
structure than the disc. Such configurations (compact bordered
surfaces) are topologically determined by the number $r$ of
boundary contours and the genus $p$ (number of handles) (see
Fig.\,1a), as is well-known since the days of M\"obius 1860/63
\cite{Moebius_1863} and Jordan 1866 \cite{Jordan_1866} (and of
course very implicit in all of Riemann's work).

\begin{figure}[h]
\hskip-35pt \penalty0
    \epsfig{figure=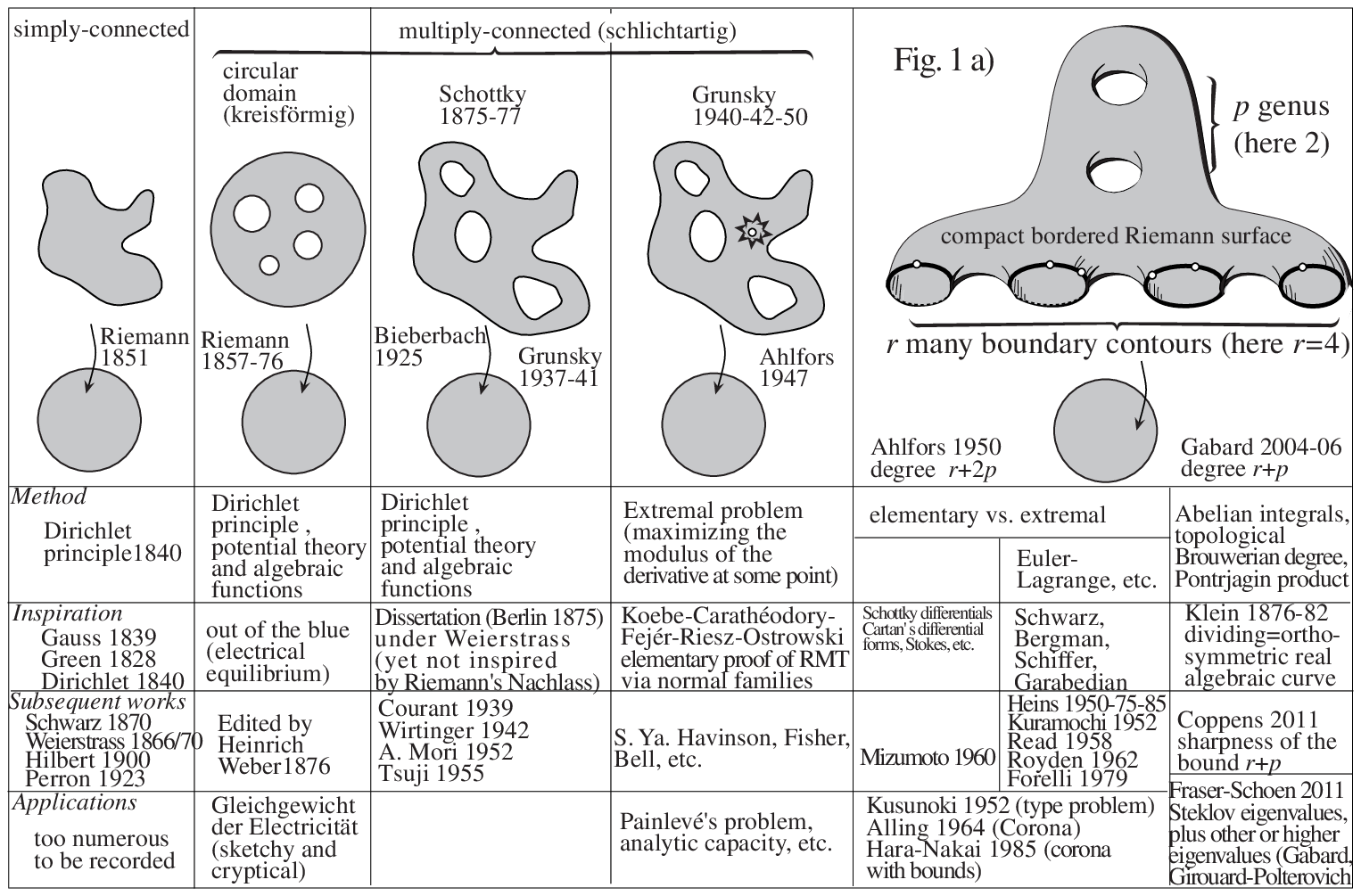,width=146mm}
    \vskip-5pt\penalty0

\caption{\label{RSGBA:fig} Schematic evolution of some mapping
theorems: from Riemann to Ahlfors transiting via
Schottky-Bieberbach-Grunsky}
\end{figure}

The possibility of mapping any bordered surface to the disc
conformally was pioneered by such
towering
figures as:

$\bullet$ Riemann 1857/76 \cite{Riemann_1857_Nachlass}
(manuscript not published during his lifetime), in which
circular domains (hence $p=0$) of
finite connectivity are
mapped upon the disc. This
fragment
was edited by H. Weber and appeared in print only in
1876 in the first edition of Riemann's Werke.
%
The date of 1857 follows some oral tradition
(Schwarz--Schottky), compare Bieberbach 1925
(Quote~\ref{quote:Bieberbach-1925} below), but conflicts
slightly with Summer 1858 as
estimated by Klein (cf. Quote
\ref{Klein-1923:quote:Riemann-1858}). [11.08.12] To pinpoint more
about the exact date, should we recall that Riemann himself
reports in the introduction of
``Theorie
der Abel'schen Functionen'' 1857 \cite[p.\,116]{Riemann_1857}
his involvement with the topic of conformal mapping of
multi-connected ``surfaces'' (Fl\"achen) right after his Thesis
(Fall 1851--Begin 1852), but  was then sidetracked to another
subject ({\it ward aber dann durch einen andern Gegenstand von
dieser Untersuchung abgezogen\/}).

$\bullet$ Schottky 1875--77 \cite{Schottky_1877} (=Dissertation
under Weierstrass, Berlin, 1875), where a similar mapping is
obtained for general real analytic contours. At first sight, it is
natural to speculate that Schottky knew about Riemann's Nachlass,
but Schottky himself describes his trajectory as independent
(cf. Quote~\ref{quote:Schottky-1882}). Apparently, it was
Weierstrass' special pupil, namely H.\,A. Schwarz who made
Schottky aware of this connection, as reported in Bieberbach 1925
\cite{Bieberbach_1925}, compare Quote~\ref{quote:Bieberbach-1925}.
Albeit independent of Riemann's, Schottky's work was likewise
physically motivated as emphasized by Klein 1923
\cite[p.\,579]{Klein-Werke-III_1923}=Quote~\ref{Klein-1923:quote:Riemann-1858}
below, or via Schottky's own recollections
(1882)=Quote~\ref{quote:Schottky-1882}.

$\bullet$ Bieberbach 1925 \cite{Bieberbach_1925}, found some
elementary arguments (or just modernization) of the same
Riemann--Schottky result, while emphasizing the trivial fact that
the degree bound is optimum (apparently Schottky gave no bound),

$\bullet$ Grunsky 1937--41
\cite{Grunsky_1937,Grunsky_1941_KA}, 1940--42--49
\cite{Grunsky_1940, Grunsky_1942, Grunsky_1950}, who
in a first series of papers rederived Bieberbach's result and then
switched to an extremal interpretation of the mapping problem.
This terrible quantitative/competitive weapon
(with historical precedents to be soon
discussed) culminated, finally, in:

$\bullet$ Ahlfors 1947 \cite{Ahlfors_1947}, but
it remained until Ahlfors 1950 \cite{Ahlfors_1950},
to prove a generalization capable of including positive genera
($p>0$),
superseding thereby quite dramatically the planarity
(Schlichtartigkeit) where
 all previous
efforts were
perpetuated.
%
(We shall attempt to
ponder this absolute
originality of Ahlfors, by comparing with others writers (e.g.,
Courant), but only with limited success due to my moderate
competence with minimal surfaces and Plateau.)

For an overall picture of the roots plus some ramifications of
Ahlfors, the reader may glance at the following map
(Fig.\,\ref{Map2:fig}) showing some of the links we are going to
explore in this survey. We have opted for a Riemann surface style
depiction of this
histogram so as to give a quick-view of the varied {\it troncs
vivaces} (in A. Denjoy's prose when alluding to history of
mathematics).
Such trunks or handles are attached whenever some philosophical
dependence (citation) is detected. Alas, it resulted a
prolix accumulation of links creating a somewhat chaotical
picture. For sharper pictures of the  ``Riemann galaxy'', we
recommend
 Neuenschwander 1981 \cite{Neuenschwander_1981}, Gray 1994
\cite{Gray_1994} and Remmert 1998 \cite{Remmert_1998}.

\begin{figure}[h]
\hskip-0pt\penalty0 \centering
    \epsfig{figure=Map2.eps,width=109mm}
\vskip-5pt\penalty0
  \caption{\label{Map2:fig}
A free-style depiction of the Ahlfors-map theory as a Riemann
surface} \vskip-5pt\penalty0
\end{figure}

{\small {\it Caveat.}---The own contribution of the writer (Gabard
2006 \cite{Gabard_2006}) predicting an improved control $r+p$ upon
Ahlfors' degree $r+2p$ is
enormously exaggerated, especially if it turns out to be false.
Other
distortions
only reflects the writer's poor understanding of this
tentacular topic.
For a more extensive compilation of authors involved in the
theory, cf. Fig.\,\ref{Geneal:fig}. If you are not cited on it,
please send me an e-mail.

}

As already said,
our central hero
will be Ahlfors, especially his paper of 1950 \cite{Ahlfors_1950}.
%
In
retrospect, it is
not quite impossible that Riemann himself (or disciples like
Schwarz, Schottky, Klein, Hurwitz, Koebe, Hilbert, Gr\"otzsch,
Teichm\"uller, etc., or also Bieberbach, Grunsky, Wirtinger,
Courant, while not forgetting in Italy, Cecioni, Matildi 1945/48
\cite{Matildi_1945/48}, Andreotti 1950 \cite{Andreotti_1950})
could have succeeded in proving such a version. Such speculations
look not purely science-fictional especially in view of
%
Ahlfors' elementary argument in
\cite[pp.\,124--126]{Ahlfors_1950}, which involves
primarily only  classical tricks (no deep extremal problem), like
annihilating all the periods to ensure single-valuedness of the
conjugate potential, and  basic potential functions arising from
the Green-Gauss-Dirichlet era. All these tricks are standard since
Riemann's days (cf. e.g. Riemann 1857
\cite[p.\,122]{Riemann_1857}, ``{\it so bestimmen da{\ss} die
Periodicit\"atsmoduln s\"ammtlich $0$ werden.}''). Remember also,
despite sembling  dubious historical revisionism, that
Teichm\"uller 1941
\cite{Teichmueller_1941}(=Quote~\ref{quote:Teichmueller-1941})
seems to have possessed a clear-cut conception of the result at
least without precise bound, while ascribing the assertion even
back to Klein.

However it took ca. 91 years---say from Riemann's 1857/58 Nachlass
up to the 1948 Harvard lecture held on the topic by Ahlfors, cf.
Nehari's Quote~\ref{Nehari-1950:quote} of 1950---until somebody
puts it on the paper and it turned out to be no less an authority
than Lars Valerian Ahlfors\footnote{Prose borrowed by Louis de
Branges.}.

It is true that Ahlfors
moved in considerably deeper waters by
solving as well a certain
{\it extremal problem\/}. This extremal viewpoint is more punchy,
yet arguably the corresponding extremals (so-called Ahlfors
functions) are only circle maps of a special character. We gain in
punch but loose in
flexibility.  The extremal functions do not substitute to---nor
are substituted by---circle maps. Deciding which viewpoint is more
useful is another question, probably premature to answer except
for guessing a complementary nature depending on the problem at
hand. Incidentally in Ahlfors paper (1950 \cite{Ahlfors_1950}),
existence of circle maps is required as a preliminary step toward
posing (non-nihilistically) the extremal problem. Ahlfors'
extremal problem stemmed surely not out of the blue, but
was patterned along a tradition, whose first steps should probably
be located in the following works. (We acknowledge guidance by
Remmert's book 1991 \cite[p.\,160--2,
 p.\,170--2]{Remmert_1991}, to which we refer for sharper
historical details.)

$\bullet$ Koebe's elementary proof 1907, 1909, 1912
\cite{Koebe_1912}, 1915 \cite{Koebe_1915}  of the (RMT); ({\it
Quadratwurzeloperationen}, {\it Schmiegungsverfahren}, etc.)

$\bullet$ Carath\'eodory 1912 \cite{Caratheodory_1912}: similar
iterative methods and convergence of his sequence via Montel's
theorem. This revitalized Koebe's interest (cf. again Remmert's
description \cite[p.\,160, p.\,172]{Remmert_1991}); in
Carath\'eodory 1914 \cite{Caratheodory_1914} full details of the
method were given in the Schwarz-Festschrift;

$\bullet$ Fej\'er and F. Riesz 1922 obtain the Riemann mapping
via an extremal problem for the derivative (published in
Rad\'o 1922/23 \cite{Rado_1922-3}).  Montel's normal families
are also used, plus a tedious derivative computation
eradicated in:

$\bullet$  Carath\'eodory 1928 \cite{Caratheodory_1928} and
Ostrowski 1929 \cite{Ostrowski_1929}, where (independently)
ultimate simplifications are provided.

Carath\'eodory wrote about these developments:

\begin{quota}[Carath\'eodory 1928 {\rm \cite[p.\,300]{Caratheodory_1928}}]\label{quote:Caratheodory-1928}

{\small \rm

Nachdem die Unzul\"anglichkeit des urspr\"ung\-lichen {\it
Riemann}schen Beweises erkannt worden war, bildeten f\"ur viele
Jahr\-zehnte die wundersch\"onen, aber sehr umst\"andlichen
Beweismethoden, die {\it H.\,A. Schwarz} entwickelt hatte, den
einzigen Zugang zu diesem Satze. Seit etwa zwanzig Jahren sind
dann in schneller Folge eine gro{\ss}e Reihe von neuen k\"urzeren
und besseren Beweisen [von ihm selbst und von Koebe (Remmert's
addition); in the original Lindel\"of 1916 is also quoted]
vorgeschlagen worden; es war aber den ungarischen Mathematikern
{\it L.~Fej\'er} und {\it F. Riesz} vorbehalten, auf den
Grundgedanken von {\it Riemann} zur\"uckzukehren und die L\"osung
des Problems der konformen Abbildung wieder mit der L\"osung eines
Variationsproblems zu verbinden. Sie w\"ahlten aber nicht ein
Variations\-problems, das, wie das {\it Dirichlet}sche Prinzip,
au{\ss}erordentlich schwer zu behandeln ist, sondern ein solches,
von dem die Existenz einer L\"osung feststeht. Auf diese Weise
entstand ein Beweis, der nur wenige Zeilen lang ist, und der auch
sofort in allen neueren Lehrb\"uchern aufgenommen worden ist.
[Footnote 2: Siehe {\it L. Bieberbach}, Lehrbuch der
Funktionentheorie, Bd.\,2 S.\,5.] Mein Zweck ist nun zu zeigen,
da{\ss} man durch eine geringe Modifikation in der Wahl des
Variationsproblems den {\it Fej\'er-Riesz\/}chen Beweis noch
wesentlich vereinfachen kann.

}

\end{quota}


Let us quote thrice Ahlfors in this connection (the second of
which occurred while celebrating the centennial of Riemann's
Thesis, 1851):

\begin{quota}[Ahlfors 1961 {\rm \cite[p.\,3]{Ahlfors_1961}}]
\label{Ahlfors-1961}

{\small \rm In complex function theory, as in many other branches
of analysis, one of the most powerful classical methods has been
to formulate, solve, and analyze extremal problems. This remains
the most valuable tool even today, and constitutes a direct link
with the classical tradition.

}

\end{quota}


\begin{quota}[Ahlfors 1953
{\rm \cite[p.\,500]{Ahlfors_1953}}]\label{Ahlfors-1953}

{\small \rm Very important progress has also been made in the
use of variational methods. I have frequently mentioned
extremal problems in conformal mapping, and I believe their
importance cannot be overestimated. It is evident that
extremal  mappings must be the cornerstone in any theory that
tries to classify conformal mappings according to invariant
properties.

}

\end{quota}

\begin{quota}[Ahlfors 1958
{\rm \cite[p.\,3]{Ahlfors_1958}}]\label{Ahlfors-1958}

{\small \rm Es ist mir zugefallen, eine \"Ubersicht \"uber die
Extremalprobleme in der Funktionentheorie zu geben. Seit der
Formulierung des Dirichletschen Prinzips ist es klar gewesen,
dass die Cauchy-Riemannschen Gleichungen nichts anderes sind
als die Eulerschen Gleichungen  eines Variationsproblems, und
in diesem Sinne ist alle Funktionentheorie mit
Extremaleigenschaften verbunden. Aber es ist nicht immer von
vornherein klar, wie diese Probleme gestellt werden sollen,
damit sie in wesentlicher Weise die tiefen Eigenschaften der
analytischen Funktionen abspiegeln. Es gibt nat\"urlich
unz\"ahlige Maximaleigenschaften, etwa in der konformen
Abbildung, die ganz nahe an der Oberfl\"ache liegen. Von da
aus soll man zu schwierigeren Problemen aufsteigen. Das
geschieht nicht etwa so, dass man ein beliebiges, wenn auch
verlockendes, Extremalproblem ins Auge fasst und es zu l\"osen
versucht. Im Gegenteil, die Entwicklung ist so vor sich
gegangen, dass man die Aufgaben stellt, die man l\"osen kann.
Dadurch ist ein reiches Erfahrungsmaterial entstanden, und die
Aufgabe des heutigen Funktionentheoretikers besteht darin,
dieses Material zu klassifizieren und dadurch weiter zu
entwickeln.

[\dots, and on page 7, of the same philosophical paper]

Carath\'eodory sagte einmal, dass er immer wieder zur
Funktionentheorie zur\"uck\-kehrt, weil man gerade dort die
verschiedensten und verbl\"uffendsten Methoden verwenden kann. Das
ist sicher wahr, und eben deshalb ist die Funktionentheorie kein
eng spezialisierter Zweig der Mathematik. Im Gegenteil, die
Funktionentheorie scheint fast wie ein Miniaturbild der gesamten
Mathematik, denn es gibt kaum eine Methode in der Geometrie, der
Algebra und der Topologie, die nicht fr\"uher oder sp\"ater in der
Funktionentheorie wichtige Anwendung findet. [\dots]

}

\end{quota}

Such
wisdoms
cultivating the
extremal philosophy---in particular as a growing mode
for
conformal mappings---presumably capture the deepest
telluric part of the
mushroom, out of which
everything derives effortlessly. Alas, our survey is  far from
this ideal conception. In fact, we would be quite challenged if we
were demanded to list a single application
of Ahlfors' extremal property,
except of
course in the planar case where one can easily mention all the
activities centering around Painlev\'e's problem.

\subsection{Applications}

The writer's interest in the topic
was recently revived by
the article of Fraser-Schoen 2011 \cite{Fraser-Schoen_2011},
where the
Ahlfors function
received a clear-cut interaction with spectral theory (Steklov
eigenvalue) with a view toward minimal surfaces.

At a more remote period of time, in the early 1950's, when
classification theory of open Riemann surfaces was a hot topic
(especially in the Finnish and Japanese schools), Kusunoki 1952
\cite{Kusunoki_1952} proposed an application to the type problem,
in the analytic sense of Nevanlinna's Nullrand (null boundary). A
(somewhat misleading but frequently used) synonym is {\it
parabolic type\/} (not to be confused with the geometric sense of
uniformization theory). This (analytic) sense of parabolicity is
the one related to the transience of the Brownian motion
(Kakutani, etc.)

In view of
the extremal r\^ole played by the (round) hemisphere
as a vibrating membranes (compare Hersch 1970
\cite{Hersch_1970}, and less relevantly Gabard 2011
\cite{Gabard_2011}),  the author speculatively expected---yet
failed dramatically to establish (Summer 2011)---the
following:

\begin{conj} {\rm (Gabard, April 2011, ca. 300 pages of sterile
hand-written notes, unpublished)} There is a mysterious connection
between the  Ahlfors function and the
(still
open) {\it filling area conjecture} (FAC) of {\rm Gromov 1983
\cite{Gromov_1983}}, whose genus zero case follows from the Thesis
of {\rm Pu 1952 \cite{Pu_1952}}, under Loewner 1949. More
precisely, the filling area conjecture is true for all
genus $p\ge 0$, and the proof will employ an Ahlfors map, at least
as one of the ingredients [others being Schwarz's inequality, and
group theoretical tricks \`a la Hurwitz--Haar--Loewner like in the
$p=0$ case]. The basic link is of course that conformal maps
supply isothermic coordinates, yielding a way to compute areas via
the infinitesimal calculus (of Newton--Leibniz, etc.).
\end{conj}

The best available result on FAC is still the hyperelliptic case
handled by Bangert-Croke-Ivanov-Katz 2004 \cite{Bangert_2004},
implying
the full conjecture for $p=1$ (as in this case the double is of
genus $g=2$, hence automatically hyperelliptic). Remember the
formulation of the FAC problem: among all compact bordered
(orientable?) Riemannian surfaces bounding the circle without
shortening its intrinsic distance, the round hemisphere has the
least possible area.

The above ``Ahlfors$\Rightarrow$Gromov'' conjecture flashed my
attention, after completing the note (Gabard 2011
\cite{Gabard_2011}) in view of the striking analogy between the
isoperimetric r\^ole of the hemisphere both acoustically (spectral
theory, like  in Hersch 1970 \cite{Hersch_1970}) and
geometrically in the L\"owner-Pu-Gromov isosystolic
($\approx$filling) problem. Of course this analogy is already
explicit in Gromov 1983 \cite{Gromov_1983}, where Hersch 1970
(\loccit) is cited. Incidentally, Gromov's account also let play
to Jenkins, Ahlfors' student and Gr\"otzsch's admirator, a
predominant logical r\^ole via the notion of ``extremal length''.
After more immature thinking (August 2012), it seems safer to
formulate a relaxed version of the conjecture where the impulse
does not necessarily come from the Ahlfors map but from some more
ancestral source like the Green's function (or the allied
Gauss-Riemann isothermic coordinates). Also the (Lorenz-)Weyl's
asymptotic law enabling to ``hear'' the area of a drum from
high-vibratory modes could be involved as well in FAC.
When Marcel Berger describes Gromov's systolic exploits (1983
\loccit), he insinuates  (surely with right) of them
as lying at a much higher level of sophistication than
2D-conformal geometry (\`a la Gauss-Riemann, etc.). This acts as
an optimism killer against anything like the above conjecture.
 Of course
 our conjecture or its relaxed
variant ``Conformal$\approx$Isothermic$\Rightarrow$Gromov'' is far
from prophetical, but only the expectation that the traditional
methods (conformal theory and uniformization) which
settled low-genus cases (Loewner 1949, Pu 1952 \cite{Pu_1952})
will extend soon or later to $p\ge 2$. Yet, who knows?

Remember that even
Marcel Berger,
once validated (or at least quoted) an erroneous proof (ca. 1998)
of the 2D-case of the filling conjecture in question. Compare his
brilliant ``Panoramic view'' (2002
\cite{Berger_2002-A-Panoramic-view-of-RG}), or rather his likewise
excellent survey in JDMV (1998 \cite[p.\,147]{Berger_1998-JDMV}):
``{\it The simplest filling volume, namely that for the circle
$S^1$, was only obtained in ([N.] Katz, 1998).}'', where the
reference is given as (cf. p.\,196) ``{\it Katz, N. (1998).
Filling volume of the circle.\/}'' This work has apparently never
been published and probably turned out to contain a gap. This
reference is still quoted in the ``Panoramic view'' (2002
\cite[p.\,790]{Berger_2002-A-Panoramic-view-of-RG}) modulo a
puzzling shift of authorship from Neil N. Katz to Mikhail G. Katz:
``Entry [794]={M.\,G. Katz}, Filling volume of the circle, to
appear, 1998.'' In the text of ``A Panoramic view\dots'' this
reference is apparently not cited, and at any rate on p.\,367 we
read ``{\it Today there is not a single manifold whose filling
volume is known, not even the circle (for which Gromov conjectures
the value [is] $2\pi$).\/}''

Of course, probably no better guide than Ahlfors himself for
listing applications of his method would have been desired. Alas
it seems that the latter was suddenly sidetracked in the
stratosphere of Teichm\"uller theory in the early 1950's, leaving
  the Ahlfors map topic in some standby ``in absentia'' status.
An exception is the later paper Ahlfors 1958 \cite{Ahlfors_1958},
where Ahlfors discusses again extremal problems, though in a more
philosophical way. Also the work of his student Read 1958
\cite{Read_1958_Acta} is described, which supplies another
existence-proof of circle maps  via a more  abstract viewpoint
(Hahn-Banach) inspired by other works like Macintyre-Rogosinski
1950 \cite{Macintyre-Rogosinski_1950}, Rogosinski-Shapiro 1953
\cite{Rogosinski-Shapiro_1953}, Rudin, etc. This Teichm\"uller
shift  in Ahlfors activities seems to coincide with the 100 years
celebration of Riemann's Thesis (in 1951), where L.~Bers cames up
with his list of urgent questions about Riemann surfaces.
As a partial consolation, Grunsky worked out a brilliant book
(1978 \cite{Grunsky_1978}) where much of the historical continuity
is supplied.

{\it Quoting some first-hand sources.}--- We shall have to
reproduce several quotations from primary sources as an attempt to
observe the mutual influences among the variety of viewpoints. It
resulted some inflation in size, but hopefully excusable as the
information of some relevance to our topic is otherwise dispatched
through a vast amount of literature. Those are given in the
self-explanatory format {\bf Quote (Author, year)}.

\smallskip
{\it Broad-lines organization of the sequel.}---We shall
essentially touch the following aspects (all in reference to the
Ahlfors mapping):

(1) Origins, background: prehistory of Ahlfors
(Sec.\,\ref{Sec:Prehistory-Ahlfors}); potential precursors
(Sec.\,\ref{Sec:Precusors});

(2) How the writer came across this topic? (via Klein); cf.
Sections~\ref{Sec:Klein} and
\ref{Sec:Biased-recollections-of-Gabard};

(3) Potential theory vs. extremal problems (both from the same
variational soup);

(4) Applications (Sec.\,\ref{Sec:Applications-of-the-Ahlf-map}):
equilibrium of electricity Riemann 1857, Pain\-lev\'e's problem,
type problem, Carath\'eodory metric, corona problem, quadrature
domains, spectral theory (Steklov or Dirichlet-Neumann);

(5) Open problems fictionally related to the Ahlfors function
(Sec.\,\ref{Sec:Virtual-applications-Ahlf-map});

(6) (Partial) assimilation of Ahlfors or other works (logical
reconstruction); via Green in Sec.\,\ref{Green:sec} and via
Ahlfors in Sec.\,\ref{Ahlfors-proof:sec};

(7) Sharpening Ahlfors work (for circle maps not necessarily
subjected to the extremal problem).

Roughly speaking our text splits as follows. A first half is
devoted to historical aspects, while a second half (initiated by
Sec.\,\ref{Sec:Starting-from-zero} titled ``Starting from zero
knowledge'') is more ``logical'', or rather liberal and futurist.
This second part tries to explore  what sort of mathematics lies
beyond Ahlfors theorem. Of course it is hard going beyond Ahlfors
without having digested his own work, and consequently much energy
is spent to the original account. His result affords considerable
information, especially the realizability of all gonalities lying
above Ahlfors bound $r+2p$. (The {\it gonality\/} $\gamma$ is
 the least degree of a circle map tolerated by the given
bordered surface.) Classically, some (episodic) penetrations
beyond Ahlfors occurred by Garabedian, Heins, Royden, etc., and
more recently in the spectacular progresses made by Yamada, Gouma
on the extremal function. In the dual direction (of circle maps),
Coppens' work on the gonality is likewise penetrating deep
behind the line fixed by Ahlfors, and raises several questions of
primary importance. This includes that of describing how the
moduli space of bordered surfaces (with fixed topological type
$(r,p)$) stratifies along gonalities. Calculating dimensions of
the varied strata is a first step toward quantifying by how much
and how frequently one can expect to improve Ahlfors bound. We
obtain so the {\it gonality profile\/}, that is, the function
assigning to each gonality $\gamma$ (in the Coppens range $r\le
\gamma \le r+p$, or outside it in case Gabard is wrong) the
dimension of the moduli strata with prescribed gonality $\le
\gamma$ (Section \ref{sec:profile-histogram}). Describing this
gonality profile appears to me a challenging (but hopefully
reasonably accessible) problem. Another ``futurist'' problem is
the one of describing the list of all degrees of circle maps
tolerated by a given surface. This we call the {\it gonality
sequence}. It is full above Ahlfors bound $r+2p$, but what can be
said below? These are perhaps two typical kind of problems hinting
at what sort of games we may encounter  ``beyond Ahlfors''.

\subsection{Bibliographic and keywords chart}

The following chart (Fig.\,\ref{Map:fig}) focuses on the
tabulation of several articles where an existence-proof of Ahlfors
circle maps is given. Such items are marked by full black circular
symbols with eventual decorations. Applications are marked by
triangular symbols. All entries of the picture (e.g. ``Ahlfors
1950'') can unambiguously be located in the bibliography at the
end of the paper. One counts essentially ca. 13 papers addressing
the existential question of circle maps.

Those includes: Ahlfors 1950 \cite{Ahlfors_1950}, Garabedian 1950
\cite{Garabedian_1950}, Heins 1950 \cite{Heins_1950}, 1975
\cite{Heins_1975}, 1985
\cite{Heins_1985-Extreme-normalized-LIKE-AHLF} and in the same
spirit Forelli 1979 \cite{Forelli_1979}. Another trend is Nehari
1950 \cite{Nehari_1950} and Tietz 1955 \cite{Tietz_1955} (alas
those works are a bit confusing, Tietz criticizes Nehari and is in
turn attacked subsequently in K\"oditz-Timmann 1975
\cite{Koeditz-Timmann_1975}). The latter work (KT1975) actually
offers an alternative
existence-proof  without
degree control.  In Japan we have Kuramochi 1952
\cite{Kuramochi_1952} and Mizumoto 1960 \cite{Mizumoto_1960}. (One
should probably add several works of Kusunoki from the early
1950's, but those are often confusing with subsequent errata,
etc.) Another mouvance is the usage of Hahn-Banach in the papers
Read 1958 \cite{Read_1958_Acta} and Royden 1962
\cite{Royden_1962}. Finally there is a work by the writer, Gabard
2006 \cite{Gabard_2006}, which even
claim a better control $r+p$ upon the degree of circle maps. Of
course this work should still be better understood and its result
should be either disproved or consolidated by alternative
techniques.

To this obvious list one can add some more telluric flows or
possible forerunners:

$\bullet$ Teichm\"uller's claim (1941 \cite{Teichmueller_1941})
that everything is already in Klein.

$\bullet$ Courant's works starting say with Courant 1939
\cite{Courant_1939} where a Plateau-style approach \`a la Douglas
is asserted to reproduce the Bieberbach-Grunsky ``schlichtartig''
case of Ahlfors.

$\bullet$ Italian workers: Matildi 1945/48 \cite{Matildi_1945/48}
and Andreotti 1950 \cite{Andreotti_1950}.

\begin{figure}[h]
\hskip-25pt\penalty0
    \epsfig{figure=Map.eps,width=142mm}
\vskip-5pt\penalty0
  \caption{\label{Map:fig}
Synoptic chart of articles with an existence-proof of Ahlfors map
(in various forms). Full black-colored circles include the
positive genus case ($p\ge0$).} \vskip-5pt\penalty0
\end{figure}

\begin{figure}[h]
\hskip-25pt\penalty0
    \epsfig{figure=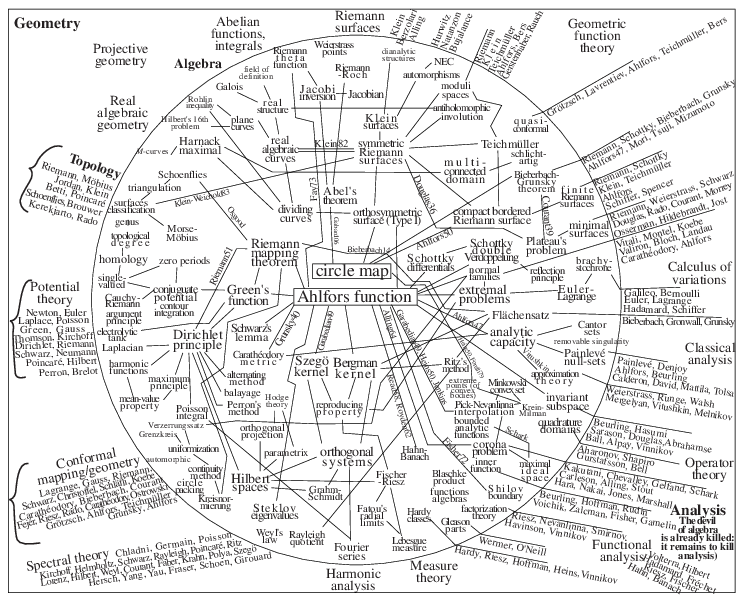,width=142mm}
\vskip-5pt\penalty0
  \caption{\label{Keyword:fig}
  Some of the keywords gravitating around the Ahlfors function}
\vskip-5pt\penalty0
\end{figure}

\smallskip
{\it Picture of Keywords.}---Let us now put Ahlfors 1950
\cite{Ahlfors_1950} at the center of the universe, while trying to
describe the portion of the cosmos visible from this perspective.
Picturing in the non-Euclidean crystal, we obtain  something like
the following chart of keywords (Fig.\,\ref{Keyword:fig}): a
nebulosity of sidereal dusts gravitating in the immediate
conceptual vicinity of the Ahlfors map.

\subsection{Mathematical questions}
\label{sec:question}

In this section we collect
questions raised by our text.
Most of our questions are of the retrograde sort ``Can we reprove
Ahlfors via \dots'', yet striving toward a perfect
crystallography, where each result of the theory is certified by
all methods ever imagined (compare optionally the kaleidoscopic
Fig.\,\ref{Kaleidoscope:fig} much below).

$\bullet$ {\bf Klein $\Rightarrow$ Ahlfors?} [reported 04.11.12]
Is it possible to reprove existence of Ahlfors circle maps via
Klein's R\"uckkehrschnitttheorem (RST) (cf. Klein 1882
\cite{Klein_1882_Ruckkehrschnitt} or Klein 1923 Ges. Math. Abh.
III \cite[p.\,622--626]{Klein-Werke-III_1923})? This paradigm RST
may be
conceived as a positive genus case of the Kreisnormierung (of
Koebe, but implicit in the Latin version of Schottky's Thesis, cf.
Klein's Quote~\ref{Klein-1923:quote:Riemann-1858}). Further recall
that Riemann (1857 \cite{Riemann_1857_Nachlass}) was able to
produce circle maps for domains bounded by circles, and by analogy
it seems plausible that Klein's RST implies (modulo some work \`a
la Riemann) the Ahlfors circle map. Of course Klein himself may
not have been able to prove rigorously his RST, but the result was
completed via some Brouwer-Koebe techniques ca. 1911/12
\cite{Klein-Brouwer-Koebe_1912}. (For a few more details about
this strategy, cf. Sec.\,\ref{sec:Ruckkehrschnittthm}.) [18.11.12]
An allied historical question is whether Teichm\"uller's
accreditation to Klein (1941 \cite{Teichmueller_1941}) of circle
maps is based on the same stratagem (RST) as we are just
suggesting.

$\bullet$ {\bf Witt or Geyer $\Rightarrow$ Ahlfors?} Can we
reprove the theorem of Ahlfors via a purely algebraic method (say
Abel, or Riemann-Roch) as Witt 1934 \cite{Witt_1934}, Geyer
1964/67 \cite{Geyer_1964-67} or Martens 1978 \cite{Martens_1978}
succeeded to do for the Witt mapping (of 1934)? For more on this,
cf. Sec.\,\ref{sec:Witt}.


$\bullet$ {\bf Plateau $\Rightarrow$ Ahlfors?} Can we reprove the
theorem of Ahlfors via the method based on the Plateau problem (as
Courant 1939 \cite{Courant_1939} did for the
Riemann-Schottky-Bieberbach-Grunsky theorem, i.e. the
schlichtartig case $p=0$ of Ahlfors). (See
Sections~\ref{sec:Courant} and \ref{sec:Douglas} for historical
precedents (i.e., Douglas 1931 \cite{Douglas_1931-Solution}),
and precise references about contemporary workers attacking
related questions (Jost, Hildebrandt, von der Mosel). A closely
related historical
 question
is whether the works of Courant do not already contain
(more-or-less explicitly) an existence-proof of Ahlfors circle
maps.

$\bullet$ {\bf Bergman $\Rightarrow$ Ahlfors?} Idem via the method
of the Bergman kernel function. This seems implicit in the
literature (say especially by Bell, e.g. Bell 2002
\cite{Bell_2002}, the great specialist of the technique), but to
the writer's knowledge no
pedestrian account is available to the mathematical public (in the
positive genus case). Compare Sec.\,\ref{sec:Bergman} for some
links to the literature. Of course behind Bergman 1922
\cite{Bergman_1922} one finds Bieberbach's characterization (1914
\cite{Bieberbach_1914}) of the Riemann map via an extremal problem
involving least area. This problems should be in some duality with
Ahlfors extremal problem, more about this soon.

$\bullet$ {\bf Behnke-Stein $\Rightarrow$ Ahlfors?} [reported
05.11.12] The article (of K\"oditz-Timmann 1975 \cite[Satz 3,
p.\,159]{Koeditz-Timmann_1975}) seems to contain a qualitative
version of Ahlfors' theorem based upon an ``Approximationssatzes
von Behnke u. Stein'', yet without any bound on the degree. Can
one improve the argument to get a quantitative control? As to
Behnke-Stein 1947/49 \cite{Behnke-Stein_1947/49} (the famous paper
going back to 1943), it contains the result that any open Riemann
surface (arbitrary connectivity and genus) admits a non-constant
analytic function. Is it possible conversely to deduce this
theorem from Ahlfors theorem by  exhaustion while pasting together
various circle maps defined over a system of expanding compact
subregions?

$\bullet$ {\bf Other techniques?} Koebe's iteration, circle
packings (cf. Rodin-Sullivan 1987 \cite{Rodin-Sullivan_1987}),
Ricci flow, etc. Virtually any technique involved in the proof of
the RMT (=Riemann mapping theorem) or the allied uniformization is
susceptible
to reprove the Ahlfors circle map.

$\bullet$ {\bf Does Ahlfors
imply Ahlfors?} [02.09.12] This repetition is intentional and
intended to emphasize that the writer was not able to digest
Ahlfors argument in full details (compare Sections~\ref{Green:sec}
and \ref{Ahlfors-proof:sec}). If one remembers the proof of
Koebe's Kreisnormierung (say as implemented in Grunsky 1978
\cite{Grunsky_1978} or Golusin 1952/57 \cite{Golusin_1952/57}),
then upon making abstraction of Koebe's proof by iterative
methods, it may be noticed that ultimately the proof depends on a
topological principle (namely Brouwer's invariance of domain). In
comparison, Ahlfors' proof of a circle map (1950
\cite{Ahlfors_1950}) makes no use of any topological principle,
reducing rather to considerations of convex geometry (cf. Ahlfors
1950 \cite{Ahlfors_1950}). Should one deduce that the Ahlfors
function lies somewhat less deep than Koebe's Kreisnormierung? If
not then maybe Ahlfors' argument lacks a global topological
character, and perhaps its validity needs to be reevaluated. (Of
course this is only a superficial objection arising from my own
frustration in not being able to catch the substance of Ahlfors
text.)

$\bullet$ {\bf Does Brill-Noether ($+$ Harnack's trick)
imply Ahlfors?} [26.10.12] Upon using projective models of Riemann
surfaces, especially birational models in the plane, it is common
practice to understand the geometry on a curve via auxiliary
pencils living on the ambient plane. Of particular importance are
the so-called adjoint series passing through the singularities of
the model which have the distinctive feature of cutting economical
series of points on the curve. Such pencils are thus involved in
the description of low-degree pencils living on the (abstract)
smooth curve, hence morphisms to the line. Adapting this
methodology to orthosymmetric curves one can evidently hope to
reprove Ahlfors theorem, provided one is able to ensure total
reality of the corresponding morphism. Details look quite
formidable to implement. If such a proof exists it will probably
be
a happy hour for its discoverer. For more vague ideas about this
strategy, see Sec.\,\ref{sec:Brill-Noether-approach-to-Ahlfors}.

$\bullet$ {\bf Does Ahlfors
imply Gabard?} [09.09.12] Upon using Ahlfors' original argument in
\cite{Ahlfors_1950} for the existence of a circle map of degree
$r+2p$, it seems evident that one could append to Ahlfors argument
a sharper geometric lemma which could produce a better control
than Ahlfors'. Ideally one would like to recover Gabard's bound
$r+p$. For some evidence of why this should be possible compare
Sec.\,\ref{Red's-function:sec}.

$\bullet$ {\bf Gabard true? If, yes analytifiable?} [June 2012] Is
the bound $r+p$ predicted by the writer on the degree of a circle
map true? And if yes is it accessible to more conventional
analytical methods? Remember that the derivation in Gabard use
some topological methods combined with the classical Abel theorem.

$\bullet$ {\bf Gonality profile.} [June 2012] Can we compute the
dimension of the moduli spaces of membranes having fixed gonality
$\gamma\le r+p$. (The {\it gonality\/} is the least degree of a
circle map from the given bordered surface.) The similar question
in the case of complex curves is well-known and easily predicted
by a simple Riemann-Hurwitz count (but established rigorously much
later). Slightly more on this in
Sec.\,\ref{sec:profile-histogram}.

$\bullet$ {\bf Ahlfors extremals as economic as Gabard?} [March
2012] Can the degree of the  Ahlfors extremal function be made as
economical as $r+p$, the circle map degree predicted by the
writer, for a suitable location of the two points required to pose
the extremal problem (resp. of a single point when considering the
derivative maximizing variant of the problem)?

$\bullet$ {\bf Ahlfors extremals as super-economic as Coppens?}
[March 2012] Same question for the sharper {\it (separating)
gonality} introduced by Coppens 2011 \cite{Coppens_2011}, that is,
the minimum sheet-number required to concretize the bordered
Riemann surface as a (holomorphic) branched cover of the disc.

$\bullet$ {\bf Topology$\Rightarrow$Riemann-Meis complex
gonality?} [21.06.12] Can the topological  method (irrigation)
used in Gabard 2006 \cite{Gabard_2006} be adapted to prove that
any complex curve of genus $g$ is $\le [\frac{g+3}{2}]$-gonal,
meaning that there is always a morphism to ${\Bbb P}^1$ of degree
$\le$ than the specified bound. (Perhaps this is already answered
in the lectures of Gunning 1972 \cite{Gunning_1972}, who uses
Mattuck's topological description of the symmetric powers of the
curve).

Conversely, there is a dual problem:

$\bullet$ {\bf
Gr\"otzsch-Teichm\"uller-Meis$\Rightarrow$Ahlfors-Gabard
separating gonality?} [16 June 2012] According to secondary
sources (e.g. Kleiman-Laksov 1974 \cite{Kleiman-Laksov_1974}),
Meis' proof (1960 \cite{Meis_1960}) of the complex gonality $\le
[\frac{g+3}{2}]$ of genus-$g$ curves, is
eminently Teichm\"uller-theoretic. By analogy, it should therefore
be possible to prove the $(r+p)$-gonality of membranes (cf. Gabard
2006 \cite{Gabard_2006}) by using the same (Teichm\"uller-style)
method as Meis. This would incidentally give an ``analytic'' proof
(or if you prefer, a ``geometria magnitudinis'' proof of Gabard
2006 \cite{Gabard_2006}). Notice the fighting interplay between
topology and analysis (or geometry) since Teichm\"uller amounts
essentially to the ``m\"oglichst konform'' map of Gr\"otzsch.

$\bullet$ [05 June 2012] Ozawa 1950 \cite{Ozawa_1950} presents
a genuine extension of the Schwarz lemma to multiply-connected
domain. Can we do the same job for a membrane of positive
genus?

$\bullet$ {\bf Ahlfors$\Rightarrow$Gromov?} [Mai 2011] Does
Ahlfors (or perhaps the non-orientable variant of Witt 1934
\cite{Witt_1934}) implies Gromov's filling area conjecture? Any
solution to this puzzling problem is rewarded by 50 Euros by
Mikhail Katz (cf. his home web-page). Perhaps, some other
ingredients than Ahlfors are required. We (already) loosely
suggested, Weyl's asymptotic law (acoustic proof) or perhaps a
sort of duality between ``Ahlfors'' extremal problem and that of
Bieberbach 1914 \cite{Bieberbach_1914} (more widely known for its
connection to Bergman). Added [02.09.12], maybe it is enough to
consider the isothermic coordinate generated by a single Green's
function (or a dipole avatar) instead of an Ahlfors function.

$\bullet$ {\bf Gromov non orientable} (Easier?) [June 2011] Is the
Gromov filling conjecture also true (and meaningful) for
non-orientable membranes? Can it be generalized to several
contours (desideratum J. Huisman 2011, oral e-mail communication).

We may also drift to related problems like KNP
(Kreisnormierungsprinzip). This asserts that any domain (or planar
Riemann surface) is conformally diffeomorphic to a domain bounded
by circles (we suppose finite connectivity for simplicity).

$\bullet$ {\bf Extremal problem$\Rightarrow$KNP?} Inspired by
the paper Schiffer-Hawley 1962 \cite{Schiffer-Hawley_1962},
where (Koebe's) Kreisnormierung (in finite connectivity) is
derived from a minimum problem of the Dirichlet type, one may
wonder if a suitable variant of Ahlfors extremal function may
not be used to reprove
the Kreisnormierung. More about this is Sec.\,\ref{sec:KNP}
(related to works by Gr\"otzsch, and others.).

$\bullet$ {\bf Bieberbach's (least area) minimum problem.}
Bieberbach 1914 \cite{Bieberbach_1914} considers in a
simply-connected domain $B$ the problem of minimizing the integral
$\int\!\! \int_B \vert f(z)\vert^2 d\omega$ amongst analytic
functions $f\colon B \to {\Bbb C}$ normed by $f'(t)=1$ at some
fixed point $t\in B$ of the domain. He shows that the minimum
gives the Riemann map. (It is well-known that this problem
constitutes the origin of the Bergman kernel theory, cf. besides
Bergman's original paper of 1922 \cite{Bergman_1922}, e.g.
Behnke's BAMS review of Bergman's 1950 book \cite{Bergman_1950}.)
The naive question is what sort of maps are obtained when this
problem is formulated on a multiply connected domain? Do we obtain
a circle map? And if yes, does this $\beta$-function coincides
with the Ahlfors map? Can the problem be extended to Riemann
surfaces?
 More on this is discussed in Sec.\,\ref{Sec:Bieberbach-Bergman}.
Of course this is closely allied to the Bergman kernel, and was
treated by several authors, cf. e.g. Garabedian-Schiffer 1950
\cite{Garabedian-Schiffer_1950}. However as far as the writer
browsed the literature, the qualitative feature of this
$\beta$-map appear to have not been explicitly described. In fact
it seems that ultimately the answer is a bit disappointing in the
sense that the least-area map may lack single-valuedness. This is
well-explained in papers by Maschler (1956--59, e.g.
\cite{Maschler_1956}), and was probably known earlier by Bergman,
Schiffer, etc.

$\bullet$ {\bf Heins' proof?} [28.06.12] Heins 1950
\cite{Heins_1950} proposes another   existence-proof of circle
maps \`a la Ahlfors, by using some theory of Martin and concepts
from convex geometry (minimal harmonic functions and extreme
points of convex bodies).
Unfortunately, he does not
keep a quantitative control upon the degree of the map so
obtained. However, on p.\,571 Heins introduces the number $m$ (of
loops generating the fundamental group), which is easily estimated
as $2p+(r-1)$ for a surface of genus $p$ with $r$ contours. [E.g.,
imagining contours as punctures, the first perforation liberates a
free group of rank $2p$ (twice the genus), and each additional
perforation creates a new generator.] Since this must be augmented
by one (cf. Heins' lemma on p.\,568, i.e. essentially the issue
that each point of a convex body in Euclidean $m$-space is
expressible as a barycentric  sum of $m+1$ extreme points of the
body spanning an $m$-simplex) it seems probable that Heins' proof
reproduces the bound $r+2p$ of Ahlfors. More about this in
Sec.\,\ref{sec:Heins}. (Actually, Heins' convex geometry argument
looks quite akin to the one used ``subconsciously'' by Ahlfors
1950 \cite{Ahlfors_1950}.)

$\bullet$ [22.10.12] {\bf The gonality sequence.}
An emerging
question of some interest
is that of calculating for a given bordered surface $F$ (of type
say $(r,p)$) the list of all integers arising as degrees of a
circle map defined on the given surface. We call this invariant
the {\it gonality sequence\/} of $F$. As a noteworthy issue
Ahlfors upper bound $r+2p$ is  always effectively realized, in
sharp contrast to Gabard's one $r+p$ which can fail to be. For
some messy and premature thoughts on this problem cf.
Sec.\,\ref{sec:gonality-sequence}. Of course the problem looks a
bit insignificant combinatorics, yet studying it properly seems to
require both experimental contemplation of concrete Riemann
surfaces and sharp theoretical analysis of the existence-proofs
available presently. Asking fine quantitative questions should aid
clarifying the qualitative existence theorems.

$\bullet$ [03.11.12] {\bf Generalized Keplerian motions via
Klein-Ahlfors?} It is well known that the motion of  a single
planet around a star
describes an orbit which is a certain algebraic curve, namely an
ellipse (other conics do occur for cold comets escaping at
infinity without periodicity). To visualize  Ahlfors circle maps
on  real plane algebraic curves of dividing type (Klein's
orthosymmetry), one can contemplate totally real pencils of curves
sweeping out the given curve along totally real collections of
points. The prototypical example is the G\"urtelkurve (quartic
with two nested ovals) swept out by a pencil of lines whose center
of perspective is located inside the deepest oval. All such lines
cut the quartic in 4 {\it real\/} points (cf.
Fig.\,\ref{FGuert:fig}). This paradigm of total reality
is the exact algebro-geometric pendant of Ahlfors theorem, and
suggests looking at real dividing curves as orbits of planetary
systems with dynamics governed by a total pencil. For instance the
G\"urtelkurve could occur as the orbit of a system of 4 electrons
gravitating around a proton with
electric repulsive forces explaining the special shape of the
G\"urtelkurve (cf. again Fig.\,\ref{FGuert:fig} below). In
Sec.\,\ref{sec:electrodynamics} we explore the (overambitious?)
idea positing that the real locus of {\it any\/} real
orthosymmetric curve (in the Euclid plane or space) arises as the
orbital structure of an electrodynamical system obeying
Newton-Coulomb law's of attraction/repulsion via a dynamics
controlled by an Ahlfors circle map (incarnated  by a totally real
pencil). This gives quite an exciting interpretation affording
plenty of periodic motions to the $n$ body problem. This idea
probably requires to be better analyzed. Even if physically
irrelevant, one can (by Ahlfors) trace for any orthosymmetric real
curve (in the plane) a totally real pencil generating usually
quite intriguing figures, especially when members of the pencil
are varied through the full color spectrum to create some rainbow
effect. Depictions of such totally real rainbows are given in
Sec.\,\ref{sec:electrodynamics}, but we failed drastically to make
serious pictures for Harnack-maximal curves. This represents
perhaps a certain challenge for computer graphics?

$\bullet$ [27.12.12] {\bf Green-Riemann imply Schoenflies?} This
question is quite outside our main track of 2D-conformal geometry,
 belonging really to highbrow unsettled differential
topology. Remember first that the Riemann mapping theorem (and the
closely allied Green's function measuring the proliferation of a
bacteria in a nutritive
medium) is essentially the best approach toward the (topological)
Schoenflies theorem (ca. 1906) stating that any plane Jordan curve
bounds a disc. Compare the contributions of Osgood  and
Carath\'eodory 1912 \cite{Caratheodory_1912}, plus the recent
discussion in Siebenmann 2005 \cite{Siebenmann_2005}. When it
comes to high-dimensional versions of Schoenflies, we know in the
smooth category by combination of the topological version of
Mazur-Brown with Smale's $h$-cobordism theorem giving uniqueness
of the smooth structure on high-dimensional balls ($\dim\ge 5$).
However remind that presently differential-topologic methods
failed to prove the (so-called) smooth Schoenflies conjecture in
dimension $4$, that any smoothly embedded $S^3$ in ${\Bbb R}^4$
bounds a $4$-ball with its usual smoothness structure. It is
tempting to wonder if the classical tools of potential theory
(especially the Green's function) are able to reprove at least the
high-dimensional cases of Smooth-Schoenflies, and if so, if it is
able to crack the residual remaining exceptional case resisting
all efforts of topologists so far. More details and references in
Sec.\,\ref{Schoenflies:sec}. The intuition behind all this is that
the bacteria expand from any given interior point as concentric
circles (resp. spheres) in the infinitesimally small but soon
realize where there is more free vital room for expanding more
quickly in those directions (cf. Fig.\,\ref{Green:fig}). In
particular all bacteria reach the boundary spheroid
simultaneously. Mathematically this is formalized by considering
the Green's function $G(z,t)$, where $z\in {\Bbb R}^n$ and $t$ is
an interior point of the bounded component of $\Sigma$ the
$S^{n-1}$ embedded in ${\Bbb R}^n$, defined as $\log \vert z -t
\vert - u$ where $u$ is the unique harmonic function with boundary
values given by $\log \vert z -t \vert$ on the boundary $\Sigma$.
Studying the Green's lines,that is the trajectory orthogonal to
the levels of $G(z,t)$ should (possibly) enable one to establish
the required diffeomorphism between the (sealed) interior of
$\Sigma$ with the ball $B^{n}$ endowed with its usual smooth
structure. (It is unknown if $B^4$ supports an exotic differential
structure but that this another question a priori much harder to
decide.)

\subsection{Some vague answers}

This section tried to report question which looks exciting, and to
which I tried some premature answer. It requires to be polished
drastically and reorganized seriously. Hence it is probably safer
to skip, but maybe  readers fluent with techniques like Ahlfors
extremals, Teichm\"uller extremal quasi-conformal maps, Plateau's
problem, etc. may find useful to clarify our vague ideas.

$\bullet$ {\bf Quantum fluctuations of Ahlfors' degree} [20.09.12]
The following problem is somewhat ill-posed, yet it is just an
attempt to excite the imagination. Suppose given a compact
bordered Riemann surface $F$ with $r\ge 1$ contours and of genus
$p\ge 0$. For each interior point $a\in \rm{int} (F)$ there is a
uniquely defined analytic Ahlfors function $f_a$ solving the
extremal problem of making the derivative $f'(a)$ as large as it
can be, while keeping this magnitude positive real and the range
inside the unit disc. This extremal function is uniquely defined
and independent of the local uniformizer used to compute the
derivative. It is known by Ahlfors 1950 that each $f_a$ is a
circle map of degree somewhere in the range from $r$ to $r+2p$,
that is a (surjective) branched cover of the disc. According to
Coppens 2011 \cite{Coppens_2011} the generic bordered surface has
gonality  $r+p$ so that one can  considerably squeeze the Ahlfors
range to the interval $r+p$ to $r+2p$. One would like to
understand in geometric term (if possible?) what phenomena is
responsible of the fluctuation of the Ahlfors degree. Of course,
if $p=0$ there is no fluctuation just because of the Ahlfors
squeezing: i.e. $\deg f_a$ is constant when the center of
expansion $a$ is dragged throughout the surface. However if $p>0$,
it is likely that some jump must occur albeit I know no argument.
Gabard 2006 only showed that there is a circle map of degree $\le
r+p$, but a priori there is no reason forcing such low degree maps
to be realized as Ahlfors maps. Following Coppens we may define
the gonality $\gamma$ of $F$ as the least degree of a circle map
on $F$. By Gabard (2006 \cite{Gabard_2006}) $\gamma\le r+p$ (and
trivially $r\le \gamma$). Coppens tell us that all intermediate
values of $\gamma$ are realized (modulo the trivial exception that
when $r=1$ and $p>0$, $\gamma=1$ cannot be realized). This
gonality invariant infers a sharpened variability for the Ahlfors
degrees, namely $r\le \gamma\le \deg f_a \le r+2p$, where $\gamma
\le r+p$. A priori all intermediate values could be visited
(between $\gamma$ and $r+2p$). However this scenario is
incompatible with the case of hyperelliptic membranes studied in
Yamada and Gouma, where the effective Ahlfors degrees are either
maximal $r+2p$ or minimal (i.e. $2$). Those examples still
indicate that despite a sparse repartition the degree distribution
is in some sense extremal, occupying the maximum space at
disposition. Is this a general behavior? This is the maximum
oscillation (Schwankung) conjecture (MOC). If true, then Coppens
gonality would always be sustained by an Ahlfors map and also
Ahlfors upper bound $r+2p$ would be sharp for any surface,
whatsoever its differential-geometric granularity. MOC displays
the most naive scenario for the fluctuation of Ahlfors degree, and
it would be a little miracle if it is correct. If not, then what
can be said? A very naive idea idea would be that there is a sort
of  conservation law like in the Gauss-Bonnet theorem: whatsoever
you bend the surface the Curvatura integra keeps constant. (Of
course this holds for a closed surface but not for a bordered one,
unless the geodesic curvature of the boundaries is controlled,
e.g. by making it null.) The vague idea would be that if we think
of the Ahlfors degree $\deg f_a$ as a sort of discrete curvature
$\delta(a)$ assigned to the point $a$ then maybe $\int_F \delta
(a) d\omega$ keeps a constant value (independent of the conformal
structure). If so then at least in the cases where there is a
hyperelliptic model (i.e. $r=1$ or $2$) one could conclude that
the Ahlfors degree are somehow balanced. Yet recalling
Yamada-Gouma's investigations it seems that  the maximum degree
$r+2p$ occurs very sporadically for the center $a$ located on the
finitely many Weierstrass points of the membrane, hence high
values have little weight. So in the hyperelliptic case (with few
contours $r=1$ or $2$) the Ahlfors degree are constantly very low
$2$ with exceptional jump taking place on a finite set of points.
Maybe this suggests a low energy scenario valid in general: given
any (finite) bordered surface $F$ the Ahlfors degree is always
equal to the gonality safe for some jump occurring on a finite set
of points. Of course this must be perhaps
refined suitably by saying that there is a stratification
(decomposition) in pieces, where the lowest  degree (i.e. the
gonality) is always nonempty and containing the contours, and then
as we penetrate more deeply inside the surface the degree may
increase (eventually always reaching the extremum value $r+2p$?).

$\bullet$ {\bf Quasiconformal doodlings} [02.10.12] As is well
known, Teichm\"uller 1939 \cite{Teichmueller_1939} exploited the
flexibility of quasiconformal maps to put Riemann's intuition of
the moduli of conformal classes of differential-geometric surfaces
(Riemannian surfaces) on a sound footing. The idea is both soft
and flexible, yet with the devil of capitalism (geometria
magnitudinis) cached just behind for one counts the distortion
effected upon infinitesimal circles into ellipses. Using
Gr\"otzsch idea of the m\"oglischt konform map relating two
configurations produces an extremal map relating both
configurations, and  this least distortion gives  the
Teichm\"uller metric (a first step to endow the moduli ``set'' of
a genuine space structure). Maybe this methodology is also
fruitful in the theory of the (Ahlfors) circle maps. The first
desideratum is to show existence of circle maps, and then the game
refines in finding best possible bounds (over the degree of such
maps).

The framework is as follows (aping again
Gr\"otzsch-Teichm\"uller): given a finite bordered surface (and
maybe also a mapping degree $d\ge r$) we look at all
quasiconformal map (not necessarily schlicht), i.e. (full)
branched cover of the disc (with the same topological feature as
circle maps of taking the boundary to the boundary and the
interior to the interior). Following Gr\"otzsch's idea we may look
at the ``m\"oglischt konform'' map, i.e. the most conformal
quasiconformal map in the family (hoping eventually to find a
beloved  conformal one). Measuring distortion (largest
eccentricity of the ellipses images of infinitesimal circles) one
gets a numerical invariant $\varepsilon(F, d) \ge 0$, namely the
infimum of the dilation among the class of all (differentiable)
maps from the bordered surface $F$ to the disc. This invariant
$\varepsilon(F, d)$ vanishes precisely when $F$ admits a
(conformal) circle-map of degree $d$. Hence it vanishes if $d\ge
r+2p$ by Ahlfors 1950 \cite[p.\,124--126]{Ahlfors_1950}, and even
as soon as $d\ge r+p$ if one believes in Gabard 2006
\cite{Gabard_2006}, where as usual $p$ is the genus and $r$ the
number of contours of $F$. However we are rather interested to use
the Gr\"otzsch-Teichm\"uller theory to rederive an independent
existence-proof. Of course in contrast with the classical setting
of Teichm\"uller's approach to the moduli problem, where one
considers exclusively schlicht(=injective) maps, we tolerate now
multivalent mappings, but this should not be an insurmountable
obstacle.

Our intuition is that it is not just a matter of measuring that is
required, but one must somehow explore the pretzel underlying the
surface to get an existence proof. Yet the flexible-quantitative
viewpoint of measuring eccentricity probably gives an interesting
numerical invariant which is now not a metric (Teichm\"uller
metric), but rather a (potential) function on the moduli space. In
fact we assign  to a given (bordered) surface $F$ a series of
number $\varepsilon(F,d)$ for $r\le d\le r+p$ (larger values of
$d$ give $0$ by Gabard 2006 \cite{Gabard_2006}), which is probably
decreasing (after eventually modifying the original problem by
permitting all maps of degree $\le d$ instead of those having
degree exactly $d$). So we get attached to $F$ a series of
dilations $\varepsilon(F,r)\ge\varepsilon(F,r+1)\ge \dots \ge
\varepsilon(F,r+p)=0$. Of course the sequence can crash to zero
before the $r+p$ bound and indeed do so as soon as Coppens'
gonality $\gamma$ is reached (that is, the least degree of a
circle map for the fixed $F$). [Of course in the exact degree $d$
variant of the problem  one can imagine more romantic behaviors
with oscillation down to zero and then becoming positive again
(touch-and-go phenomenology).] Those $p$ invariants would refine
Coppens gonality in a continuous fashion, yet  fails to be
``moduli'' since there are $3g-3$ of them (Riemann-Klein) where
$g$ is the genus of the double (that is $2p+(r-1)$), hence giving
a total of $3g-3=3(2p+(r-1))-3=6p+3r-3$ free parameters which
exceeds of course our $p$ parameters.

But coming back to the basic existence problem, one can get
started by observing that any topological type of membrane admits
a circle map. One trick is to use symmetric membranes (cf.
Chamb\'ery section \ref{sec:Chambery} below). This amounts to
imagine a membrane in $3$-space symmetric under rotation by 180
degree so that the quotient as genus zero (cf.
Fig.\,\ref{Chambery:fig} below). Once the handles are killed one
is reduced to the simple (planar) case of Ahlfors due to
Bieberbach-Grunsky (and largely anticipated by Riemann, Schottky
(no bound by Schottky?), and Enriques-Chisini (via Riemann-Roch
and a continuity argument, cf. e.g. Gabard 2006
\cite[Sec.\,4]{Gabard_2006}). The degree of the resulting map is
easily computed (and of degree essentially equal to $(r/2) \cdot
2=r$ the minimum possible value, for the rotation identifies pairs
of contours and gyrate all handles over themselves, cf. again
Fig.\,\ref{Chambery:fig}, below). Thinking in the moduli space $M$
we have shown that the set $C$ of all circle-mappable surfaces is
nonempty, and using the connectedness (of $M$) it would suffice to
show that $C$ is {\it clopen\/} (i.e., closed and open). Checking
openness, certainly requires enlarging the mapping degree to
larger values. Now given an arbitrary bordered surface $F$ we can
quasiconformally map it to our symmetric model $S$ and then
compose with the circle-map. The dilatation is then controlled in
term of the Teichm\"uller distance from $F$ to $S$, giving an
upper bound over the eccentricity invariant $\varepsilon$ (for the
appropriate degree). Of course this is still miles away from
reproving even Ahlfors but maybe the idea is worth pursuing.

In fact what is truly interesting is that we get for each $d$ a
numerical function $\varepsilon_d$ (defined as
$\varepsilon_d(F):=\varepsilon(F,d)$) on the moduli space
$M_{p,r}$ of membranes of genus $p$ with $r$ contours, that
vanishes precisely when $F$ has gonality $\le d$. Of course this
sequence of functions is monotone decreasing when the index
increases, and $\epsilon_d\equiv 0$ is identically zero (for $d\ge
r+p$). According to Coppens result each of these functions (let us
call them the Teichm\"uller potentials) vanishes somewhere. It is
then perhaps interesting to look at the gradient flow $\varphi_d$
(w.r.t. Teichm\"uller metric) of these functions $\varepsilon_d$
affording a dynamical system (=flow) in which each bordered
surface evolves in time to a sort of best possible surfaces for
the prescribed gonality. (Morally each surface tries to improve
its gonality along the trajectory of steepest descent.) If the
global dynamics is simplest (say each trajectory finishes its life
on a surface of gonality $d$) it is therefore reasonable to expect
that the whole Teichm\"uller space is retracted by deformation to
a sort of spine consisting of surfaces having the prescribed
gonality $d$. Maybe one can deduce that the global topology of
this spine is that of a cell (like the full Teichm\"uller space).
Further it seems probable that the flows preserve the
stratification by the gonality of $M_{p,r}$ since if $F$ has
gonality say $d$ then its future $F_t$ has lower gonality. [The
situation looks  analog to some  works of Ren\'e Thom (isotopy
lemma, vector fields preserving a stratification, and ``fonction
tapissante'' as it arise in the Thom-Mather problem of the
stability of polynomial mappings??]

[03.10.12] Of course the above can be adapted to the case of
closed (non-bordered) surfaces of genus say $g$, by replacing the
target disc by the (Riemann) sphere. Likewise we define
Teichm\"uller potentials $\varepsilon_d$, measuring the dilatation
of the ``m\"oglichst konform'' map of a fixed degree $d$ from the
surface $F$ to $S^2$, and ideally one can imagine that the theory
is able to reprove the famous (Riemann-Brill-Noether) bound
$[\frac{g+3}{2}]$ first proved by Meis 1960 \cite{Meis_1960}.
Hence all what we are trying to do is surely already well-known
(alas I was never able to find a copy of Meis' work, which is
Teichm\"uller-theoretic according to other sources). Hence if Meis
theory is just a sort of Teichm\"uller theory for branched covers
of the sphere, with the ultimate miracle  that Teichm\"uller not
only affords a solution to Riemann's moduli problem but also to
the gonality question. A priori Meis' theory  should adapt to the
bordered setting and arguably lead to another proof of the Ahlfors
map, and optimistically  with the sharp bound predicted in Gabard
2006 \cite{Gabard_2006}. Sharpness of the bound is due to Coppens
2011 \cite{Coppens_2011}. Recall that, Teichm\"uller himself was
close to this (bordered) topic in the article Teichm\"uller 1941
\cite{Teichmueller_1941}, yet the details (as well as exact
bounds) are probably missing.

$\bullet$ {\bf Ahlfors inflation/injection and generalized
 Ahlfors maps taking values outside the disc (alias, circle)} [09.10.12] The theory of the
Ahlfors function is primarily based upon the paradigm of
maximizing the derivative (its modulus) within the family of maps
with range confined to a (compact) container namely the unit disc.
So it is primarily an inflation/injection (or pressurization)
procedure (by opposition to the dual deflation/suction approach of
Bieberbach-Bergman amounting to minimize the area among maps
normed by $f'(z_0)=1$). Ahlfors 1950 \cite{Ahlfors_1950} showed
that if the source object is any compact bordered Riemann surface
and the target the unit disc then the Ahlfors (inflating) map
turns out to be a circle map, i.e. a full covering of the unit
circle taking
 boundary to boundary. This behavior is not
surprising since maximizing the distortion (scaling factor) at a
given basepoint forces the whole surface to be maximally stretched
over the target, like an elastic skin pushed to its ultimate limit
(in the Hollywoodian context of aesthetical surgery). The
existence of Ahlfors maps relies on a Montel normal family
argument,  in substance  inherited from the compactness of the
disc. This suggests replacing the target disc by any compact
bordered Riemann surface. We formulate then the following extremal
problem:

Given two finite bordered Riemann surfaces $F$ and $G$ and a
given point $a \in F$ and $b\in G$, we look inside the family
of all analytic maps $f\colon F \to G$ taking $a$ to $b$ at
the map maximizing the modulus of the derivative $f'(a)$
computed w.r.t. local parameters introduced at $a$ and $b$.

By analogy with the Ahlfors et ali theory, we expect that the
extremal function exist (compactness of the receptacle $G$), is
unique (this is either less evident or false for in the classical
case $G=\Delta$ the argument relied heavily on the Schwarz lemma
for the disc, so that our only hope in favor of uniqueness is that
what actually counts is the universal covering). Arguably, even if
lacking uniqueness extremals could still be interesting. Finally
it is reasonable to expect that extremals are not oversensitive to
the choice of local uniformizers. So we can speak of the map
$f_{a,b}$ of extremum dilatation at $a,b$. Finally we are
interested about knowing if the extremals are total maps in the
sense of taking boundary to the boundary, as do the classical
Ahlfors map in the circle/disc-valued case. Before proceeding to
examples let us perhaps observe that in the special case where $F$
is given as a subsurface of $G$ and both points $a=b$ coincide,
then the (complex) tangent space are readily identified so that
$f'(a)$ has an intrinsic meaning as scaling factor of this complex
line. Another special case of interest is when $G$ is  a plane
subregion, in which case the tangent bundle is trivialized so that
one can consider a relaxed form of the problem without the
constraint $f(a)=b$, in which no point $b$ is given but the sole
extremalization of $f'(a)$ will actually dictate where $a$ has to
be mapped.

Albeit all we are saying looks a bit messy and unnatural (?), it
should be noted that the whole game can be drastically simplified
by just looking at avatars of circle maps, that is given two
finite Riemann surfaces $F$ and $G$ when does there exist a total
map (taking boundary to boundary) from the first to the second.
(Of course this question is quite standard yet probably hard to
answer precisely, cf. Landau-Osserman 1960
\cite{Landau-Osserman_1960}, and Bedford 1984
\cite{Bedford_1984}.) As we shall soon explain a vague answer is
readily supplied by ``algebraic geometry'', namely when the target
$G$ is not the disc, and if $F$ has general moduli then in general
there in not a single total map from $F$ to $G$. The moral is that
circle maps enjoy a certain privilege due to their unconditional
existence (by Ahlfors 1950 precisely).

A basic obstruction arises from the Riemann-Hurwitz formula.
Indeed given $f\colon F \to G$ a total map, it has no ramification
along the boundary and is a full covering surface (cf. e.g.
Landau-Osserman 1960 \cite[p.\,266,
Lemma~3.1]{Landau-Osserman_1960}). Denoting by $d\ge 1$ the degree
of the map, we have $\chi(F)=d \chi (G)-b$, where $b\ge 0$ counts
the branch points. When $d=1$, there is no branching and the
topological types must agree. Another constraint says roughly that
a total map can only simplifies the topology, precisely $\chi(F)=d
\chi (G)-b\le d \chi (G)\le \chi(G)$, when $\chi(G)\le 0$.

\begin{lemma} If $G$ is not the disc then the existence of a
bordered map $f\colon F \to G$ implies that the Euler
characteristic satisfies $\chi(F) \le \chi(G)$. (Of course the
conclusion persists when $G$ is the disc for it maximizes the
Euler characteristic among bordered surfaces.)
\end{lemma}

Another simple constraint comes from the fact that a total map
$f\colon F \to G$ induces a covering of the boundary $\partial
f\colon \partial F \to \partial G$. Hence if $G$ has $r'$ contours
then $F$ has at most $d \cdot r'$ contours, i.e. $r\le  d \cdot
r'$ where $r$ is the number of contour of $F$. On the other hand
as $\partial f$ is onto, the surjection induced by $\partial f$ on
the $\pi_0$ (=the arc-wise connected component functor from TOP to
SET) implies that $r\ge r'$.

Then there is a little zoology of cases to study.

(Z1) Let us first suppose that the {\it source\/} is just the
disc, then who is the (``Ahlfors'') extremal map? So we assume
$F=\Delta$ and $G$ any bordered surface marked at $a=0$, $b\in G$
respectively. By uniformization (Koebe-Poincar\'e 1907) we know
that the universal cover of the interior of any finite bordered
surface is the disc. Now the extremal map $f_{a,b}\colon \Delta
\to G$ (maximizing the distortion) may be lifted to the universal
cover as say $F \colon \Delta \to \Delta$. Now by the Schwarz-Pick
principle of hyperbolic contraction for analytic maps, the latter
map contracts the hyperbolic metric implying the universal
projection to effect a greater dilatation than the presumed
extremal $f_{a,b}$. It follows that $F$ must be the identity (up
to rotation) and the extremal function get identified to the
universal cover. (Actually, works by Carath\'eodory and Grunsky
actually manage to prove uniformization via the (Ahlfors) extremal
problem, whereas we assumed it.)

(Z2) Now consider the situation were both source and target have
complicated topology. For instance the source is any bordered
surface and the target an annulus. One may expect to get analogues
of circle maps, i.e. {\it total maps\/} taking boundary to
boundary (sometimes known as proper maps). (Such maps are called
{\it boundary preserving\/} in Jenkins-Suita 1988
\cite{Jenkins-Suita_1988}, cf. also Landau-Osserman 1960
\cite[p.\,265]{Landau-Osserman_1960} who speak of maps ``which
takes the boundary into the boundary'', while ascribing to Rad\'o
1922 \cite{Rado_1922-Z-Theorie-mehr} the basic result that such
maps are full coverings taking each value of  the image surface a
constant number of times). Unfortunately, there is severe
obstructions to boundary preservation of such (generalized)
Ahlfors maps. One way to argue is via algebraic geometry and the
Jacobians. It is indeed classic that a generic closed Riemann
surface tolerates only nonconstant maps to the sphere (ruling out
the trivial identity map or automorphisms available incidentally
only for surfaces with specialized moduli). Assuming the Ahlfors
map of $F$ to an annulus to be total, its symmetric extension  to
the Schottky-Klein double is a map from a closed surface to the
torus, which for general moduli cannot exist at all! Of course all
this requires better proofs, but is fairly well-known and
classical (cf. e.g. Griffiths-Harris 1980
\cite{Griffiths-Harris_1980}, who argue as follows (p.\,236--237):
``{\it A general curve $C$ of genus $g\ge 2$ cannot be expressed
as a multiple cover of any curve $C'$ of genus $g'\ge 1$.} This is
readily seen from a count of parameters: the curve $C'$ will
depend on $3g'-3$ parameters, and the $m$-sheeted covering $C\to
C'$ depends on $b$ parameters, where [$\chi (C) = m \chi(C')-b$,
that is]
$$
b=2g-2-m(2g'-2)
$$
is the number of branch
points of the cover. Thus if $m\ge 2$, $C$ will depend on
$$
b+(3g'-3)=b+\frac{3}{2}(2g'-2)
=2g-2-\underbrace{\bigl(m-\frac{3}{2}\bigr)}_{\ge
1/2}\underbrace{(2g'-2)}_{\ge 0}\le 2g-2<3g-3
$$
parameters, and so cannot be general.'' (Another argument is
given in the exercises of Arbarello-Cornalba-Griffiths-Harris
1985 \cite[p.\,367,
Ex.\,C-6]{Arbarello-Cornalba-Griffiths-Harris_1985-BOOK},
which of course we were not able to solve!)

(Z3) Finally one can imagine a bordered surface embedded in a
slightly larger one (say of the same topological type). Then the
inclusion map is permissible in the extremal problem, so the
extremal map will have distortion $\ge 1$ at some basepoint, and
naively should expand the small surface into the larger one.
However by the argument of (Z2) in general it is unlikely that the
extremal will be total, and also a priori it not even  clear that
a true expansion can occur (try to lift the map to the universal
cover a get maybe a conflict with the Schwarz-Pick principle of
contraction??) But of course this looks dubious for when the
subsurface is a disc expansion is possible.

\begin{figure}[h]
\centering
    \epsfig{figure=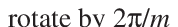,width=112mm}
\vskip-5pt\penalty0
  \caption{\label{Cyclo:fig}
  A cyclotomic Riemann surface}
\vskip-5pt\penalty0
\end{figure}

$\bullet$ {\bf Cyclotomic Riemann surfaces}
[09.10.12] (but similar examples in Cham\-b\'ery Talk ca. 20
December 2004) At this stage we can do perhaps the following sort
of  experiment. As is well-known (Riemann-Prym-Klein 1882
\cite{Klein_1882}) a Riemann surface structure can also be defined
in the most simplest way to visualize, namely as
differential-geometric surface in $3$-space with metric (hence
conformal structure) inherited by the Pythagorean/Euclidean line
element. Consider a hemisphere in Euclidean 3-space surmounted by
$m$ handles cyclotomically distributed as on
Fig.\,\ref{Cyclo:fig}, joining themselves above the north pole.

Ignoring the south hemisphere, we obtain so a bordered surface $F$
with one contour ($r=1$) and of genus $p=m-1$. (Notice here the
standard psychological aberration that the genus is one less than
the ``handles'', for the first handle is not yet coupled to
another one to create a real handle!) On rotating by angle
$\frac{2\pi}{m}$ the configuration $F$ upon itself we obtain a map
from $F$ to the disc (hence a circle map), because the fundamental
domain of the rotation is glued over itself to give a disc. The
circle map so obtained has degree $m=p+1$. This matches with the
general bound $r+p$ predicted in Gabard 2006 \cite{Gabard_2006}.
Let us now assume that a bubbling, i.e. an Euclid-Riemannian
deformation of the metric takes place at one of the handle (yet
not on the remaining ones) then the rotational symmetry is killed
and it becomes much
 nontrivial that a circle map of same degree is still
persistent.  This experiment seems to
damage the truth of Gabard 2006 \cite{Gabard_2006} (but hopefully
is not?) A naive parade would be to use the (Riemann-Schwarz)
uniqueness of the conformal structure on the closed 2-cell to
resorb the cancerigenic bubbling. Yet this looks cavalier (for we
are not living in the soft smooth $C^{\infty}$ category) and this
would not  settle the case of less localized cancerigenic
degenerations not supported over a disc, but along a subregion
having itself moduli. Then one cannot repair easily the
deformation by a simple surgical lifting. At this stage we see
that the result of Gabard 2006 \cite{Gabard_2006}, if true at all,
looks quite formidable
for it should resist all those plastic deformations within the
flexibility of conformal maps.

It could be interesting (by adapting Yamada-Gouma) to study
the degree of the Ahlfors map of such cyclotomic Riemann
surfaces, especially when the basepoint is situated on the 3
fixed points of the rotation.

$\bullet$ {\bf Special triangulations} [10.10.12] Given a circle
map of a bordered surface $F$, one can post-compose it with the
map taking conformally the disc to an equilateral triangle (in the
Euclid plane ${\Bbb C}$). (Recall that this can be done for any
three point prescribed along the boundary). Upon subdividing the
triangle in a mesh of equilateral triangles, and lifting via the
conformal map we generate certain triangulations of $F$ which are
almost equilateral. In fact if the mesh size is chosen so that all
ramification points lye in the interior of the tiny triangles then
the inverse image of such ramified triangles will be small
hexagons. Try to study the differential geometry and specialize to
Gromov's Filling conjecture, or try to find a link with
Belyi-Grothendieck (a Riemann surface is defined of
$\overline{{\Bbb Q}}$ iff it admits an equilateral triangulation).

Another special triangulation of the disc is the hyperbolic
tessellation depicted on the front cover of Grunsky's Collected
papers (by equilateral triangles with angles $\pi/6$). [This
tessellation is supposed via the Ahlfors-Grunsky conjecture (1937
\cite{Ahlfors-Grunsky_1937}) to play an extremal r\^ole in the
Bloch schlicht radius of maps $\Delta \to {\Bbb C}$ for it
dominates the densest circle packing of the Euclidean plane.] Try
to understand if it is useful (or aesthetical) to lift this
tessellation to the bordered surface via a circle map.

$\bullet$ {\bf Plateau heuristics $\Rightarrow$ Ahlfors maps?}
[17.10.12] Soap film experiments of the Belgian physicist have a
certain existential convincing power, albeit the rigorous
mathematical existence proof (Douglas/Rad\'o ca. 1930/31) required
circa 30 years more delay than the allied Dirichlet principle
(Hilbert 1900) itself interpretable at the equilibrium temperature
distribution in a heat-conducting plate with assigned boundary
values. Now Douglas 1931 \cite{Douglas_1931-Solution} observed
that the Riemann mapping theorem (RMT) may be derived by
specializing Plateau's problem to the case where the contour
degenerate to the plane, and Courant pushed the remark further so
as to include the Riemann-Bieberbach-Grunsky theorem (=planar case
of the Ahlfors map). On the other hand Douglas 1936
\cite{Douglas_1936-Some-new-results} envisaged the so-called
Plateau-Douglas problem (PDP, or just PP) for membranes of higher
topological structure. It should thus follow (either logically or
intuitively) a physico-chemical existence proof of the Ahlfors
map.

Let $F$ a finite bordered Riemann surface of genus $p$ with $r\ge
1$ contours, and suppose also given a fixed circle in the plane
interpreted as the prescribed wire frame of PP. More generally one
can imagine a collection of $r$ contours to be given, and we look
at the special case where all these coincide with the unit circle.
Now cultivating the right intuition about PP it should be possible
to deduce the existence of Ahlfors maps perhaps even with the
degree control $r+p$ of Gabard. In fact it should even be possible
to study wide extensions where not all frames are coincident with
the unit circle. Can one take any frame prescription (e.g.
disjoint round circles)? For instance take two unit circles with
centers lying distance ten apart ($\vert z \vert=1$ and $\vert
z-10 \vert=1$). Suppose the membrane to have the topological type
of an annulus ($r=2$ and $p=0$). Then the minimal surface is
something like a flat catenoid, where the inside of each circle is
covered once and a certain tube connecting both circles is covered
twice by the map. Yet notice that the apparent contour (where the
map is folded) of such a film violates the local behavior of
holomorphic maps. As we just saw the folding obstruction makes
unlikely to span contours consisting of disjoint maximal circles.
(Circles being ordered by inclusion of their interior in the
plane.) In contrast a nested configuration of circles should cause
no trouble to holomorphy. Thus it should be possible to render
Ahlfors intuitively obvious via soap film experiments. Of course
this was essentially done in Courant's book (1950
\cite{Courant_1950}), yet the exact juncture with Ahlfors result
probably deserves some extra working. Of course the real challenge
would be to investigate if Plateau-style approaches are
susceptible to vindicate the degree bound $r+p$ advanced by Gabard
2006 \cite{Gabard_2006}.

Another idea is to imagine a Plateau problem with ``wind'' blowing
through 3-space in some prescribed way (along a given vector
field). For instance a soap film spanning a planar disc at rest
could deform under a perpendicular wind into say a hemispherical
membrane. Try  to connect this with Gromov's filling conjecture,
yet unlikely due to the embedded nature of Plateau. Another more
plausible connection would be with Gottschalk's conjecture on
flows  in 3-space (no vector fields in 3-space having only dense
trajectories). This is probably one of the most alienating open
problem in the qualitative theory of dynamical systems.

$\bullet$ {\bf (Ahlfors) circle maps of minimal degree} [19.10.12]
Given a finite bordered Riemann surface $F$ of genus $p$ with $r$
contours, there is always (by Ahlfors) a circle map. The set of
(positive) integers being well-ordered there is a circle map of
minimal degree. Call perhaps such maps {\it minimal circle maps.}
We may ask to which extent such a map is unique and if not can we
describe the  ``moduli space'' of such maps. Of course in the most
trivial case where $p=0$ and $r=1$ (topologically a disc) the
Riemann map is essentially unique ignoring automorphisms of the
disc. Likewise uniqueness holds for surfaces with hyperelliptic
double provided the latter is not Harnack-maximal. Such
hyperelliptic membranes have $r=1$ or $r=2$ and the hyperelliptic
involution induces a totally real morphism of degree $2$. Our
uniqueness assertion follows of course from the uniqueness for
complex curves of the hyperelliptic involution when $g\ge 2$ and
thus  holds in our context provided $p\ge 1$ (recall that
$g=(r-1)+2p$). When $p=0$ and $r=2$ uniqueness fails, for then the
double has genus one and may be concretized as a smooth plane
cubic with two circuits: one being a genuine ``oval'' bounding a
disc in ${\Bbb P}^2({\Bbb R})$, the other being termed a
pseudo-line. Projecting from any point located on the oval gives a
totally real morphism of degree $2$, and correspondingly a circle
map when restricted to the semi-Riemann surface. Another example
is the G\"urtelkurve, i.e. any smooth quartic with two nested
ovals. Then the minimal degree of a circle map (for the half of
the curve) is 3 (argue with the complex gonality of smooth plane
curves), and such maps arise by projecting the curve from a real
point located on the innermost oval of the nest. Hence there
$\infty^1$ circle maps of minimum degree, those being
parameterized by a circle $S^1$. Of course the problem does not
depend only on the topology: the half of the G\"urtelkurve belongs
to the topological type $r=2$ and $p=1$, which contains also
hyperelliptic representatives, those being circle mappable in a
unique fashion via a map of degree 2.

When $F$ is planar ($p=0$) then the double is Harnack-maximal and
either the argument of Enriques-Chisini or that of
Bieberbach-Grunsky shows that any divisor with one point on each
oval moves in a linear system which is totally real (cf. e.g.
Gabard 2006 \cite{Gabard_2006}). So we have now essentially a
torus of dimension $r$ ($r$=number of contours) of circle maps of
minimum degree. A details description is not so evident for such a
divisor $D$ moves in a linear system of dimension $\dim \vert D
\vert \ge \deg D-g$ (Riemann's inequality, a direct consequence of
Abel), where $g=r-1$ is the genus of the double. Thus $\dim \vert
D \vert \ge r-(r-1)=1$ so that $D$ does not necessarily determines
unambiguously a totally real pencil. Despite this difficulty it
seems reasonable to assert that the set of circle maps for a
planar membrane is a torus perhaps of dimension only $r-1$ for one
has to unite divisors lying in the same pencil. Extrapolating such
examples, we may wonder about structural properties of the set of
(minimal) circle maps. Is it always compact? Always a manifold?
Perhaps even always a torus. Is it always connected? Of course
there are various way to formulate the question and there probably
basic experiments giving quick answers to the naive connectedness
assumption. Another question is to understand how the global
degree $d$ of the circle map splits (partitioned) into the
bordered degrees of the restriction to each contours. For instance
in the case of the G\"urtelkurve, albeit both ovals are perfectly
equivalent from the  viewpoint of analysis situs, it seems that on
the Riemann surface the one corresponding to the inner oval can be
mapped with degree 1 whereas the other is less ``economic''
requiring a wrapping of degree 2. Of course it would be nice to
understand this in some intrinsic fashion? But how? (Perhaps via
the uniformizing hyperbolic metric and the length of the
corresponding ovals???)

Let us try a naive approach to the connectedness problem (by
actually trying to corrupt it). Consider in the ``abstract quadric
surface'' ${\Bbb P}^1\times {\Bbb P}^1$ a configuration of
bidegree say $(3,3)$. We have chosen both degrees equal so that
both projections have the same degree. Imagine 3 lines in each
ruling and smooth out the corresponding line arrangement to create
a smooth curve $C_{3,3}$ (cf. Fig.\,\ref{Cyclo:fig}b). Actually we
have performed  sense-preserving smoothings (cf. again the figure)
so that the resulting curve is dividing (Fiedler type argument
\cite{Fiedler_1981}). Contemplating the figure we count $r=3$
``ovals''. Both projections on the factors are totally real
morphisms of degree 3 (the minimum possible degree in view of the
trivial lower bound $\deg f \ge r$). However it seems unlikely
that one can continuously deform one map into the other (while
keeping its degree minimum). Hence this may give some evidence
that the space of minimal circle maps (for the corresponding
bordered surface, namely the half of the orthosymmetric Riemann
surface underlying our dividing curve $C_{3,3}$) is not connected.
However our argument is quite sloppy,  having equally well applied
to bidegree $(2,2)$ in which case the corresponding curve is
Harnack-maximal [recall that $g=(a-1)(b-1)$ for bidegree $(a,b)$],
hence subsumed to the connectivity principle. Of course it is
probable that some basic complex algebraic geometry (gonality of
complex curves) suffices to complete the above argument. Is it
true that a smooth curve of bidegree $(n,n)$ is $n$-gonal in only
two fashions (provided $n\ge 3$) via the natural projections? Of
course the assertion is false for $n=2$, for then $g=1$ and a
smooth plane cubic model creates $\infty^{1}$ pencils of degree 2.
For another plane example seeming to violate the connectivity
principle of minimal maps see Fig.\,\ref{Coppens:fig}(code 313).

\subsection{Some historical
puzzles}

%

\indent $\bullet$ The most
scorching question is whether Klein really anticipate Ahlfors as
suggested in Teichm\"uller 1941 \cite{Teichmueller_1941}? (Compare
Sec.\,\ref{sec:Teichmueller}, especially
Quote~\ref{quote:Teichmueller-1941}.) Of course the question bears
not only historical interest, but has some didactic  importance
for if a Kleinian argument ever existed (and not just in
Teichm\"uller's imagination), it is quite likely to be
 more geometric than Ahlfors' (decidedly analytic) account.
As we already said, it is possible that the Klein-Teichm\"uller
proof rest upon the R\"uckkehrschnitt intuition of Klein. Even in
case Klein himself never anticipated the Ahlfors circle maps, one
may wonder from where Teichm\"uller derived it? In turn one may
wonder if Ahlfors took inspiration by Teichm\"uller 1941
\cite{Teichmueller_1941}? Of course Ahlfors himself never quoted
this Teichm\"uller work, except in Ahlfors-Sario 1960
\cite{Ahlfors-Sario_1960}, where also all the Italian workers are
carefully listed (especially Matildi 1945/48
\cite{Matildi_1945/48} and Andreotti 1950 \cite{Andreotti_1950}).

$\bullet$ Does Courant's paper of 1939 \cite{Courant_1939} (and
the somewhat earlier announcement of 1938 \cite{Courant_1938},
plus the later book treatment of 1950 \cite{Courant_1950}) presage
(modulo a suitable interpretation) any anticipation over the
circle map result of Ahlfors 1950 \cite{Ahlfors_1950}? (For more,
compare Sec.\,\ref{sec:Courant}.)

\section{The province of Felix Klein}\label{Sec:Klein}

\subsection{Felice Ronga and Felix Klein's influence}

In fact the writer himself came across (a weak version of) the
Ahlfors function topic from a somewhat different angle,  namely
via Klein's theory of {\it real algebraic curves} (spanning over
the period 1876--92). For Klein this was probably just a baby case
of his paradigm of the Galois-Riemann Verschmelzung (Erlanger
Program 1873, friendship with Sophus Lie, Ikosaheder and its
relation to quintic in one variable, etc.). Yet, real curves
surely deserved special (Kleinian) attention as it provided a
panoramic view
(through the algebro-geometric crystal) of the just emerging
topological classification of surfaces (M\"obius 1863
\cite{Moebius_1863}, Jordan 1866 \cite{Jordan_1866}, etc.).
This offered also a bordered (even possibly non-orientable)
avatars of Riemann surfaces, as shown in the somewhat
grandiloquent title chosen by Klein ``{\it \"Uber eine neue
Art der Riemannschen Fl\"achen\/}'' (=title of 1874
\cite{Klein_1874}, 1876 \cite{Klein_1876}) . Those works of
Klein were probably not
extremely influential (and still today represent only a
marginal subbranch of the giant tree planted by Riemann).

Klein himself lamented at several places  his work not having
found the quick impact he expected from.
In 1892 \cite[p.\,171]{Klein-Werke-II_1922} (ten years after
his systematic theory presented in 1882 \cite{Klein_1882}), he
writes: ``{\it Inzwischen hat noch niemand, so viel ich
wei{\ss}, die hier gegebene Fragestellung seither
aufgegriffen, {\rm [\dots]}}''. About the same period in his
lectures of 1891/92
\cite[p.\,132]{Klein_1892_Vorlesung-Goettingen}, he wrote:
{\it Was ich bislang von diesen Theoremen publicirt habe (so
die Einteilung der symmetrischen Fl\"achen in meiner Schrift
von 1881), hat nur wenig Anklang gefunden. Ich meine aber,
da{\ss} das nicht am Gegenstande der Untersuchung liegt, der
mir viel mehr das gr\"o{\ss}te Intere{\ss}e zu verdienen
scheint, sondern an der knappen Form, mit der ich meine
Resultate darstellte. }

Of course this impact was first limited to his direct circle of
students, where we count Harnack 1876 \cite{Harnack_1876},
Weichold 1883 \cite{Weichold_1883} and Hurwitz 1883
\cite{Hurwitz_1883} (also a student of Weierstrass). Klein was
also very proud that his results on real moduli supplied a natural
answer to questions addressed (but not solved) at the end of
Riemann Thesis. Klein insists twice on this issue in 1882
\cite{Klein_1882}=\cite[p.\,572, \S 24]{Klein-Werke-III_1923} and
in his subsequent lectures 1891/92
\cite[p.\,154]{Klein_1891--92_Vorlesung-Goettingen}, where he
writes: ``{\it Mit dieser Abz\"ahlung ist implicite die
entsprechende Frage f\"ur \underline{berandete Fl\"achen}
beantwortet, was darum ein gewi{\ss}es Intere{\ss}e hat, weil
diese Frage von Riemann in seiner Di{\ss}ertation aufgeworfen,
aber nicht zu Ende discutirt wird. Riemann denkt nat\"urlich nur
an berandete \underline{einfache} Fl\"achen (nicht an
Doppelfl\"achen; dem deren Existenz wurde erst zehn Jahre sp\"ater
von Moebius bemerkt und wohl erst in meiner Schrift f\"ur
funktionentheoretische Zwecke herangezogen).}''

From the very beginning 1876 \cite{Klein_1876}=\cite[\S
7,\,p.\,154]{Klein-Werke-II_1922}, Klein
noticed that real curves are subjected to the dichotomy of
being dividing or not, where the former case amounts to a
separation of the complex locus through its real part
(consisting of {\it ovals}, a jargon immediately suggesting
Hilbert's 16th problem,
%
%
%
yet used much earlier, e.g. by Zeuthen 1874 \cite{Zeuthen_1874}).
Zeuthen's work seems to have much inspired Klein's investigation
on real curves, starting circa 1876, just two years later (cf.,
e.g. Klein 1892 \cite[p.\,171]{Klein-Werke-II_1922}: ``{\it Ich
hatte 1876 den Ausgangspunkt unmittelbar von den Kurven genommen.
Das war bei $p=3$ m\"oglich, wo ich zahlreiche geometrische
Vorarbeiten, insbesondere diejenigen des Herrn Zeuthen [\dots]
(1874), benutzen konnte.}'')

Perhaps, the more tenacious followers of Klein's viewpoint came
somewhat later and the real demographic explosion of the subject
took place much later, say perhaps in the 1970's. Here is a little
chronology:

$\bullet$ del Pezzo 1892 \cite{del-Pezzo_1892}, where Klein's
trick of assigning the unique real point of an imaginary tangents
is taken as the starting point of a  study of curves of low genus.

$\bullet$ Berzolari 1906 \cite{Berzolari_1906}, who in an
encyclopedia article surveyed in few pages Klein's achievements
and virtually coined the term ``Klein surfaces'' (Kleinsche
Fl\"achen) as a way to designate possibly non-orientable and
eventually bordered avatar of Riemann surfaces. To say the least,
this terminology was dormant during several decades until
Alling-Greenleaf managed in 1969 \cite{Alling-Greenleaf_1969} a
resurrection of Berzolari's coinage, and since then the
nomenclature gained in popularity.

$\bullet$ Koebe 1907 \cite{Koebe_1907_UrAK} who studied
uniformization of real algebraic curves taking advantage of
Klein's
distinction orthosymmetric vs. diasymmetric.

$\bullet$ Severi 1921
\cite[p.\,230--6]{Severi_1921-Vorlesungen-u-alg.-Geom-BUCH}, who
devotes some few  pages of his book to Klein's theory of real
curves, [Note: there Severi writes down the same formula as one
used by Courant in his approach to conformal circle maps,
ascribing it to Cauchy].

$\bullet$ Comessatti 1924-25 \cite{Comessatti_1924/26} in Italy
(full of admiration for Klein), who pushed
the philosophy up to include a study of real abelian
varieties,
rational varieties, etc. (For this ramification we refer to the
remarkable survey by Ciliberto-Pedrini 1996
\cite{Ciliberto-Pedrini_1996}.)

$\bullet$ several works of Cecioni  in the late 1920's
(\cite{Cecioni_1929}, \cite{Cecioni_1933}, \cite{Cecioni_1935}),
and his students (Li Chiavi 1932 \cite{Stella-li-Chiavi_1932})
makes direct allusion to Klein's works.

$\bullet$ In France, the work of Klein found a little echo in
some passages of the book by Appell-Goursat whose second tome
(1930) was apparently mostly
written by Fatou. There, Klein's orthosymmetry occurs at several
places \cite[p.\,326--332 and
p.\,513--521]{Appell-Goursat-Fatou_1930}.

$\bullet$ Witt 1934 \cite{Witt_1934}, where a general existence
theorem for {\it invisible\/} real algebraic curves (those with
empty real locus like, e.g. $x^2+y^2=-1$) was established. This
will be discussed later (Sec.\,\ref{sec:Witt}), and is somehow
quite akin to the Ahlfors function. Witt's work makes explicit
mention of Klein, and was subsequently elaborated by Geyer 1964/67
\cite{Geyer_1964-67}, who arranged a purely algebraic
interpretation of Weichold's work. His pupil G.~Martens, managed
(1978 \cite{Martens_1978}) to determine the lowest possible degree
of the Witt mapping;

$\bullet$ (Jesse) Douglas 1936--39 makes  a systematic use of
Klein's symmetric surfaces in his study of Plateau's problem for
configuration of higher topological structure. (We shall have to
come back to this in Sec.\,\ref{sec:Douglas}.)

$\bullet$ A marked influence of Klein upon Teichm\"uller
1939 \cite{Teichmueller_1939}, 1941 \cite{Teichmueller_1941}. We
shall try to explore this connection in greater detail later
(Section \ref{sec:Teichmueller}).

Then different events occurred at a rather rapid
pace with several schools penetrating into Klein's reality
paradigm through different angles:

$\bullet$ Ahlfors 1950 \cite{Ahlfors_1950},  who never
quotes Klein. Probably with Lindel\"of and Nevanlinna as teachers
one
is more inclined toward  hard analysis \`a la Schwarz, than
innocent looking geometry \`a la Klein. Of course Ahlfors quotes
instead Schottky, as typified by the terminology Schottky
differential, etc. used in Ahlfors 1950 (\loccit). It may then
appear as a little surprise that Ahlfors' result affords a purely
algebraic (in term of real function fields) characterization of
Klein's orthosymmetric curves. However to my knowledge, this
connection---as trivial as it is---was never emphasized in print
until much later, namely in Alling-Greenleaf 1969
\cite{Alling-Greenleaf_1969}.

$\bullet$ Schiffer-Spencer's book 1954
\cite{Schiffer-Spencer_1954} (outgrowing from Princeton lectures
held during the academic year 1949--50) where the book is started
by recalling how Klein assimilated the full Riemannian concept
after a 1874 discussion with Prym revealing him the ultimate
secret of Riemann's function theory developed over arbitrarily
curved surfaces not necessarily spread over the plane. The
original source reads as follows:

\begin{quota}[Klein 1882 {\rm \cite{Klein_1882}}]\label{quote:Klein-Prym}
{\rm Ich wei{\ss} nicht, ob ich je zu einer in sich
abgesch\-lossenen Gesamtaufassung gekommen w\"are, h\"atte mir
nich Herr Prym vor l\"ang\-eren Jahren (1874) bei gelegentlicher
Unterredung eine Mitteilung gemacht, die immer wesentlicher f\"ur
mich geworden ist, je l\"anger ich \"uber den Gegenstand
nachgedacht habe. Er erz\"ahlte mir,} {\it da{\ss} die
Riemannschen Fl\"achen urspr\"unglich durchaus nicht notwendig
mehrbl\"attrige Fl\"achen \"uber der Ebene sind, da{\ss} man viel
mehr auf beliebig gegebenen krummen Fl\"achen ganz ebenso komplexe
Funktionen des Ortes studieren kann, wie auf den Fl\"achen \"uber
der Ebene. }
\end{quota}

From circa 1970 upwards, the study of so-called Klein surfaces
(jargon of Berzolari \cite{Berzolari_1906}) {\it per se} enjoyed a
rather exponential rate of growth as if the simple naming of them
was a stimulus for a big expansion of the topic. After two decades
an impressive
body of knowledge has been accumulated (cf. e.g. the rich
bibliography compiled in Natanzon's survey 1990
\cite{Natanzon_1990/90}). Those developments can be roughly ranged
into 3 main axes:

$\bullet$ {\it Foundational aspects.}---Alling-Greenleaf 1971
\cite{Alling-Greenleaf_1971}, and also in Romania with the
numerous contribution of Andreian Cazacu (1986--88
\cite{Andreian-Cazacu_1986}, \cite{Andreian-Cazacu_1988-Interior})
about the structure of morphism between them (interior influence
of Stoilow).


$\bullet$ {\it Symmetry, automorphisms and NEC(=non-Euclidean
crystallography).}---This is especially active in the Spanish
school but started somewhat earlier with Singerman 1971--88 (5
items), May 1975--88 (9 items), Bujalance 1981--89 (29 items)
Costa, etc.

$\bullet$ {\it Moduli spaces  of Klein surfaces.} This starts of
course in Klein 1882 \cite{Klein_1882}, to reach a certain climax
in Teichm\"uller 1939 \cite{Teichmueller_1939} and the
Ahlfors-Bers school, Earle 1971, Sepp\"ala 1978--89 (6 items on
Teichm\"uller and real moduli), Silhol 1982--89 (Abelian varieties
and Comessatti), Costa, Huisman 1998+, etc.

All those works contributed to feel virtually as comfortable with
real curves as with their complex grand sisters. We just mention
one result of Sepp\"al\"a 1990 (revisited by
Buser-Sepp\"al\"a-Silhol 1995 \cite{Buser-Seppala-Silhol_1995} and
Costa-Izquierdo 2002 \cite{Costa-Izquierdo_2002}) to the effect
that the moduli space of real curves is connected. (This sounds
almost like a provocation to anybody familiar with the
bio-diversity of topological types of symmetric surfaces listed by
Klein). Of course the trick, here, is that those authors regard
this moduli space projected down in that of complex curves (by
forgetting the real involution). In other words we may deform the
structure until new anti-conformal symmetries appears and switch
from one to the other. Hence the subject is sometimes hard to
grasp (due to varying jargon) and more seriously is full of real
mysteries allied to the real difficulty of the subject.

$\bullet$ {\it Geometry of real curves\/}. Here much of the
impulse---very much in  Klein's tradition---came through the paper
of Gross-Harris 1981 \cite{Gross-Harris_1981}. In this or related
direction, we may cite authors like Natanzon, Ballico, Coppens, G.
Martens, Huisman, Monnier, etc. This area proved very active since
the 2000's up to quite recently and a remarkable variety of
difficult question are addressed giving the field arguably some
maturity soon comparable to the complex hegemony.

Of course, another line of thought is the interest aroused by
Hilbert's 16th problem (on the mutual disposition of circuits of
real algebraic varieties esp. curves) especially among the early
German, Italian and then mostly the Russian annexion of the
subject. This captured and probably contributed to mask Klein's
more intrinsic viewpoint for a while. This axis includes the
following workers (precise references listed in Gudkov 1974
\cite{Gudkov_1974/74}):

$\bullet$ Hilbert 1891--1900--09, Rohn 1886--1911--11--13; (it
is interesting to note that Hilbert's first 1891 paper on the
subject is quite synchronized with Klein's lectures of
1891/92, which conjecturally may have stimulated Hilbert's
interest, yet not a single allusion to Klein in this paper,
and recall also that Hilbert was still in K\"onigsberg at that
time).

$\bullet$ Brusotti
1910--13--14--14--15--16--16--16--16--17--21--28--38/39--40--44/45--46--50/51--52--55--55
(characterized by ``{\it la piccola variazione\/}'', i.e. the
method of small perturbation permitting to construct real
algebraic curves with controlled  topology. The writer is indebted
to
Felice Ronga for this method, which of course has some historical
antecedents older that Brusotti. In  Klein 1873, footnote 2 in
\cite[p.\,11]{Klein-Werke-II_1922} the principle is traced
 at least back to  Pl\"ucker 1839 \cite{Plücker_1839}. However Brusotti 1921
\cite{Brusotti_1921} may have been the first---modulo its reliance
over work of Severi---to notice that the Riemann-Roch theorem
admits as extrinsic traduction the possibility of smoothing
independently the nodes of a plane curve. The main issue (as
transmitted by Felice) is that the nodes a plane curve with nodal
singularities impose independent conditions on curves of the same
degree. Hence when  the curve is being imagined  as a point in the
(projective) space of all curves, it sits on the discriminant
hypersurface (parameterizing all singular curves) and nearby our
nodal curve we see several transverse smooth branches
crossing transversally. (In French or Italian, there are better
synonyms like ``falde analytiche'' or ``nappe''.) The net effect
of transversality is that one can leave at will certain strata,
while staying on others. This implies the independency of
smoothing crossings, and thereby a rigorous foundation to the
small perturbation method. (The resulting graphical flexibility of
algebraic curves is a pleasant way to create Riemann surfaces, and
we shall exploit it later in this text as a way to explore degrees
of Ahlfors circle maps.)

$\bullet$ Comessatti (more in the spirit of Klein)
1924--25--27/28--31--32--33, etc.

$\bullet$ Petrovskii 1938--49 (\cite{Petrowsky_1938}), etc. many
joint with Oleinik (real algebraic (hyper)surfaces and Betti
numbers).

$\bullet$ Gudkov 1954--54--62--62--62--65--66--69--69--73 (those
works include in particular the spectacular discovery of a sextic
whose oval configuration was expected to be impossible by
Hilbert).

$\bullet$ Arnold 1971--73.

$\bullet$ Rohlin 1972--72--73.

$\bullet$ Finally the long awaited (?) reunification of forces
(call it maybe the  Klein-Hilbert Verschmelzung) came in the work
of Rohlin (himself apparently inspired by Arnold). Surprisingly,
Rohlin took notice of Klein's work quite late, ca. 1978 (compare
Rohlin 1978 \cite{Rohlin_1978}).

$\bullet$ Then real algebraic geometry exploded through the work
of Kharlamov, Viro, Fiedler, Nikulin 1979, Orevkov, Finashin, etc.
and in the west Risler, Marin, and many others gave a  new golden
age to a discipline reaching a certain popularity.

Sometimes the real theory seems only to adapt over ${\Bbb R}$
whatever has been achieved over ${\Bbb C}$, yielding usually a
kaleidoscopic fragmentation of truths into a real zoology. Thus
for instance the Castelnuovo-Enriques classification of
(algebraic) surfaces can be pushed through reality: K3
(Nikulin-Kharlamov), Abelian surfaces (Comessatti-Silhol),
elliptic surfaces, etc. The topic is then strongly allied to deep
methods in differential topology, Galois cohomology, symplectic
geometry, Gromov-Witten, enumerative problems, tropical geometry,
etc. The present  number of active workers is so impressive and
the recent connections so amazing (Okounkov, etc.) that we prefer
to stop here our impressionist touristic overview of real
algebraic geometry.

\subsection{Digression about Hilbert's 16th problem (Klein 1922,
Rohlin 1974, Kharlamov-Viro ca. 1975, Marin 1979, Gross-Harris
1981) }

The connection between  Klein's theory (especially the ortho- and
diasymmetric dichotomy) with Hilbert's 16th problem (plane curves
in the projective plane ${\Bbb P}^2$) were
profoundly investigated by the Russian school in the early 1970's
especially Arnold, Rohlin, Viro, Kharlamov, etc. Klein himself
always dreamed of such a relationship , without really being able
to formulate its precise shape. Here is a quote which Klein added
(ca. 1922) to his Werke (cf. \cite[p.\,155,
footnote]{Klein-Werke-II_1922}):

\begin{quota}[Klein 1922]\label{Klein-1922-immer-vorsgeschwebt:quote}
{\small \rm Es hat mir immer vorgeschwebt, dass man durch
Fortsetzung der Betrachtungen des Textes Genaueres \"uber die
Gestalten der reellen ebenen Kurven beliebigen Grades erfahren
k\"onne, nicht nur, was die Zahl ihrer Z\"uge, sondern auch,
was deren gegenseitige Lage angeht. Ich gebe diese Hoffnung
auch noch nicht auf, aber ich muss leider sagen, dass die
Realit\"atstheoreme \"uber Kurven beliebigen Geschlechtes
(welche ich aus der allgemeinen Theorie der Riemannschen
Fl\"achen, speziell der ``symmetrischen'' Riemannschen
Fl\"achen ableite) hierf\"ur nicht ausreichen, sondern nur
erst einen Rahmen f\"ur die zu untersuchenden M\"oglichkeiten
abgeben. In der Tat sind ja die doppelpunktslosen ebenen
Kurven $n$-ten Grades f\"ur $n>4$ keineswegs die allgemeinen
Repr\"asentanten ihres Geschlechtes, sondern wie man leicht
nachrechnet, durch $(n-2)(n-4)$ Bedingungen partikularisiert.
Da man \"uber die Natur dieser Bedingungen zun\"achst wenig
weiss, kann man noch nicht von vornherein sagen, dass alle die
Arten reeller Kurven, die man gem\"ass meinen sp\"ateren
Untersuchungen f\"ur $p={ n-1 \cdot n-2 \over 2}$ findet,
bereits im Gebiete besagter ebener Kurven $n$-ter Ordnung
vertreten sein m\"u{\ss}ten, auch nicht, da{\ss} ihnen immer
nur {\it eine} Art ebener Kurven entspr\"ache. \quad K.

}
%
\end{quota}

It took several decades until the experimentally
obvious conjecture (possibly anticipated by Klein, though he
left no trace in print) that  dividing curves in the plane
have at least as many ovals as the half value of its
degree
found place in a
paper of Gross-Harris 1981 \cite[p.\,177, {\it
Note\/}]{Gross-Harris_1981}. In fact, in a
paper by Alexis Marin 1979/81 \cite{Marin_1979} this is stated as
a corollary of a Rohlin formula (1978 \cite{Rohlin_1978}),
involving intersection of homology classes deduced from the halves
of the dividing curve capped off by the interiors of ovals in
${\Bbb P}^2({\Bbb R})$). In the case of $M$-curves (=the Russian
synonym of Harnack-maximal coined by Petrowskii 1938
\cite{Petrowsky_1938}), this technique occurred earlier in Rohlin
1974/75 \cite{Rohlin_1974/75}. Moral: the  tool missing to Klein
was intersection theory of homology classes developed by
Poincar\'e, Lefschetz, etc.
In the little note Gabard 2000 \cite{Gabard_2000} it is
verified that this Rohlin-Marin obstruction ($r\ge \frac{m}{2}$)
is the only one, settling thereby completely the
Klein-Gross-Harris question.
This (simple)
fact was known to Rohlin's students Kharlamov and Viro which were
familiar with this result as early as the middle 1970's (as they
both kindly informed me by e-mail). Of course the crucial ideas
are due to Rohlin.

\subsection{A long unnoticed tunnel between Klein
and Ahlfors (Alling-Greenleaf 1969, Geyer-Martens 1977)}

More importantly, for our present purpose is to keep the abstract
viewpoint of Klein (by opposition to the embedded Hilbert's 16th
problem), and to make the following observation.

\begin{theorem} {\rm (Klein?, Teichm\"uller 1941?, Ahlfors 1948/50,
Matildi 1945/48?, Andreotti 1950?, who else?)} Dividing curves are
precisely those admitting a real morphism (i.e., defined over the
ground field ${\Bbb R}$)
to the projective line ${\Bbb P}^1$ such that all fibers over real
points consist entirely of real points.
\end{theorem}

The non-trivial implication of this fact follows precisely from
Ahlfors 1950 \cite{Ahlfors_1950} (but is made very explicit only
in Alling-Greenleaf 1969 \cite{Alling-Greenleaf_1969}, see also
Geyer-Martens 1977 \cite{Geyer-Martens_1977}). To my actual
knowledge there is no record in print of this fact prior to
Ahlfors' intervention, modulo the cryptical allusion in
Teichm\"uller 1941 \cite{Teichmueller_1941} that the result was
implicit in Klein's works. Another related works are those of
Matildi 1948 \cite{Matildi_1945/48} and Andreotti 1950
\cite{Andreotti_1950}.
As we shall recall later Ahlfors' result was exposed at Harward as
early as 1948 (cf. Nehari 1950 \cite{Nehari_1950} reproduced here
as Quote~\ref{Nehari-1950:quote}).

It is however picturesque to notice that an analog result stating
that a real curve without real points
maps through a real morphism
upon the empty curve $x^2+y^2=-1$ (or projectively
$x_0^2+x_1^2+x_2^2=0$) was established as
long ago as
Witt 1934 \cite{Witt_1934}.
Perhaps both problems are
of comparable
difficulty, and the method employed by Witt---namely Abelian
integrals---turns out
to be likewise relevant to the Ahlfors context (i.e. dividing
curves).
Hence in our opinion, there were no technological obstruction to
Ahlfors result being discovered much earlier, say by Witt in the
1930's, or by Bieberbach in 1925 \cite{Bieberbach_1925}, or by
Klein in the 1876--80's, or even by Riemann in the late 1850's
(especially in view of his {\it Parallelogramm methode/Figuren},
cf. e.g., Haupt 1920 \cite{Haupt_1920}), and
ultimately why not by Abel himself? (Of course all these peoples
were probably involved with more urgent tasks, like some {\it
fl\"uchtigen Versuche} about the Riemann hypothesis, or regarding
Klein the {\it Grenzkreistheorem} (in his health taking contest
with Poincar\'e), which later became known as the uniformization
theorem. The list of competent workers coming also very close to
the paradigm ultimately discovered by Ahlfors could easily be
elongated: especially Schwarz, Hurwitz (esp. in 1891
\cite{Hurwitz_1891-Uber-Riemannsche-Flachen}), Koebe, Courant
(esp. in 1939 \cite{Courant_1939}, 1940 \cite{Courant_1940-Acta}
or 1950 \cite{Courant_1950}).

As to the interesting result of Witt 1934 (on invisible real
curves), we will try to discuss it later in more details
(Sec.\,\ref{sec:Witt}).

\subsection{Motivation (better upper bounds exist)}

Even though Ahlfors' result is approaching 65 years
(a
venerable age for retirement) the basic result
looks still grandiose, and mysterious enough if one wonders about
the exact distribution of Ahlfors' degrees (as suggested in
Yamada-Gouma's penetrating study (1978--1998--2001), discussed in
Sec.\,\ref{Yamada-Gouma:subsec}).

The writer published a paper in 2006 \cite{Gabard_2006} where a
circle map with fewer sheets (viz. $\le r+p$) than that proposed
by Ahlfors (namely $\le r+2p$) is exhibited. This quantitative
improvement is the motivation for much of this survey, and will
hopefully excuses the bewildering variety of topics addressed. An
obvious game is to renegotiate known application of the Ahlfors'
mapping involving a controlled degree in the hope to upgrade the
bound. As tactically simple as it may look, we were not very
successful in this game as it often already requires analytical
skills beyond the competence of the writer. Yet we shall mostly
content to list some articles where some upgrade could be expected
(e.g., Hara-Nakai's quantitative version of the corona with bounds
\cite{Hara-Nakai_1985} looks to be a challenging place to test).
Of course for this {\it bound upgrading procedure}  to work it
requires that the application in question does not use the full
strength of the Ahlfors function, but only its qualitative
property of being a circle map.
A concrete instance were this was accomplished is Fraser-Schoen's
paper 2011 \cite{Fraser-Schoen_2011}.

Alternatively we can dream of certain high powered applications
requiring the full extremal power of the Ahlfors mapping. In this
case it is  known a priori (Yamada-Gouma) that we cannot lower the
degree of the Ahlfors function, except possibly for very
particular choices of base-points.

So the main philosophical issue is roughly the following point:

Is the Ahlfors extremal property truly required in applications,
or just the arcane residue of those attempt to
salvage the Dirichlet principle via extremal methods. Put
differently, is the extremal problem just an artefact of the proof
or something really worth exploiting in practice?

\subsection{Full coverings  versus
Ahlfors' extremals}

To avoid any confusion, one must from the scratch relativize
strongly the importance of the recent contribution on the $r+p$
bound (Gabard 2006 \cite{Gabard_2006}) for several reasons.

First the result is quite recent and probably not sufficiently
verified as yet. In later sections when looking at explicit curves
from the experimental viewpoint it seems that there is a large
armada of potential counterexamples flying at high altitudes
(flying fortresses).

Next, Ahlfors' upper bound $r+2p$ is known to be sharp within the
realm of the extremal problem it solves. Indeed, Yamada 1978
\cite{Yamada_1978} has a rather simple argument showing that the
Ahlfors function centered at the Weierstrass points of a
hyperelliptic membrane has degree precisely $r+2p$ (and not less).
Maybe it is an open question whether a similar sharpness holds for
all membranes.

Hence, one must keep in mind a
subtle distinction between Ahlfors' deep extremal problem
(involving hard analysis
via the paradigm of extremality) and the writer's soft version
(\cite{Gabard_2006}) which leads to a sharper bound but is based
only upon (soft) topological methods, i.e. the Brouwerian degree
and the allied criterion of surjectivity. To put it briefly, we
must distinguish Ahlfors' extremal function from the mere {\it
circle map}, defined as follows (nomenclature borrowed from
Garabedian-Schiffer 1950
\cite[p.\,182]{Garabedian-Schiffer_1950}):

\begin{defn}
A circle map is an  analytic function from a
compact
bordered Riemann surface to the disc, expressing the former as
a (generally branched)
cover
of the disc,
say $f\colon W \to D=\{\vert z\vert \le 1\}$.
Each interior points maps to an interior points of the disc
(otherwise there is a problem as infinitesimally the mapping
is a power map $z\mapsto z^n$, $n\ge 1$). Thus, the restricted
covering $\partial W \to
\partial D=S^1$ is unramified, whereupon
it follows that $r\le \deg (f)$ (i.e. the number of contours is a
trivial lower bound for the degree of a circle map).
\end{defn}

Varied synonyms (or closely allied designations) are used
throughout the literature (here is a little sampling with citation
of the relevant sources):

$\bullet$ $n$ fach ausgebreitete Fl\"ache, $n$ fach bedeckende
Fl\"ache (Riemann 1857--Weber 1876
\cite[p.\,473]{Riemann_1857_Nachlass});

$\bullet$ Schottky 1877 no clear cut terminology, and re-reading
it (25.06.12) in details I realize that the statement about
existence of circle maps is in fact not really proved (thus much
of the written is somewhat biased), note that Bieberbach somewhat
wrongly ascribe the result as well to Schottky, but that remains
to be elucidated... In contrast, Grunsky never (?) credits
Schottky, but rather Bieberbach 1925 \cite{Bieberbach_1925};

$\bullet$ mehrfach bedeckte Kreisscheibe, $n$-bl\"attrige
Kreisscheibe (Bieberbach 1925 \cite[p.\,6]{Bieberbach_1925});

$\bullet$ mehrbl\"attrige Kreise, $n$-bl\"attrige Kreisscheibe
(Grunsky 1937 \cite[p.\,40]{Grunsky_1937});

$\bullet$ ein endlichvielbl\"attriges Fl\"achenst\"uck \"uber der
oberen $z$-Halbebene mit endlich vielen Windungspunkten, das durch
Spiegelung an der reellen Achse eine symmetrische geschlossene
Riemannsche Fl\"ache ergibt (Teichm\"uller 1941
\cite{Teichmueller_1941});

$\bullet$ cerchio multiplo (Matildi 1945/48
\cite[p.\,82]{Matildi_1945/48}, a student of Cecioni);

$\bullet$ full covering surface of the unit circle
(Ahlfors 1950 \cite[p.\,124, p.\,132]{Ahlfors_1950});

$\bullet$ $(2g+m)$-sheeted unbounded covering surface of the unit
disc (Encyclopedic Dictionary of Mathematics 1968/87
\cite[p.\,1367]{EDM_1968/87});

$\bullet$ unbounded finitely sheeted covering surfaces of the unit
disk (Nakai 1983 \cite[p.\,164]{Nakai_1983});

$\bullet$ Schottky functions (Garabedian-Schiffer 1949
\cite[p.\,214]{Garabedian-Schiffer_1949}, K\"uhnau 1967
\cite[p.\,96]{Kuehnau_1967}, and earlier (yet without this
appellation) in several works of Picard, e.g. Picard 1913
\cite{Picard_1913} and Cecioni, e.g. Cecioni 1935
\cite{Cecioni_1935});

$\bullet$ $p$-times covered unit-circle (Bergman 1950
\cite[p.\,87, line 5]{Bergman_1950});

$\bullet$ $n$-times covered circle, multiply-covered circle
(Nehari 1950 \cite[p.\,256, resp. p.\,267]{Nehari_1950}, Stanton
1971 \cite[p.\,289 and 293]{Stanton_1971} Aharonov-Shapiro 1976
\cite[p.\,60]{Aharonov-Shapiro_1976});

$\bullet$ Ahlfors mapping (Nehari 1950 \cite[p.\,256,
p.\,267]{Nehari_1950},  Stanton 1971 \cite[p.\,289 and
293]{Stanton_1971};

$\bullet$ Ahlfors function (Aharonov-Shapiro 1976
\cite[p.\,60]{Aharonov-Shapiro_1976});

$\bullet$ Ahlfors map (Alling 1966~\cite[p.\,345--6]{Alling_1966},
Stout 1967 \cite[p.\,274]{Stout_1967-Interpolation}, and then in
many papers by Bell);

$\bullet$ Ahlfors type function (Yakubovich 2006
\cite[p.\,31]{Yakubovich_2006});

$\bullet$ Einheitsfunktionen (Carath\'eodory 1950
\cite[vol.\,II,
p.\,12]{Caratheodory_1950_Buch_Funktionentheorie}), translated
as:

$\bullet$ {\it unitary} function in Heins 1965
\cite[p.\,130]{Heins_1965}, a jargon also adhered to by Fay 1973
\cite[p.\,108, 111, etc.]{Fay_1973};

$\bullet$ unimodular function (Douglas-Rudin 1969
\cite{Douglas-Rudin_1969}, Fisher 1969
\cite{Fisher_1969-BAMS-convex-combination-unimodular-fct}, Gamelin
1973 \cite{Gamelin_1973-BAMS}, Lund~1974 \cite{Lund_1974});

$\bullet$ many-sheeted disc (A. Mori 1951 \cite{Mori_1951});

$\bullet$ multi-sheeted circle (Havinson 1953
\cite{Havinson_1953});

$\bullet$ finitely sheeted disks (Hara-Nakai 1985
\cite{Hara-Nakai_1985});

$\bullet$ Vollkreisabbildung
(Meschkowski 1951 \cite[p.\,121]{Meschkowski_1951});

$\bullet$ (volle) $n$-bl\"attrige (Einheits)Kreisscheibe
(Golusin 1957 \cite[p.\,240, 412]{Golusin_1952/57}, as translated
by Grunsky or Pirl);

$\bullet$ interior mappings (Sto\"{\i}low, Beurling);

$\bullet$ inner functions (Beurling 1949 \cite{Beurling_1949},
Hoffman 1962 \cite{Hoffman_1962} (esp. p.\,74, where Beurling is
credited of the coinage), Rudin 1969 \cite{Rudin_1969}, Stout 1972
\cite[p.\,343]{Stout_1972}, \v{C}erne-Forstneri\v{c} 2002
\cite[p.\,686]{Cerne-Forstneric_2002}).
This concept usually refers to analytic functions with modulus
a.e. equal to one along the boundary, but some writers
corrupted this sense to mean a circle map, cf. Stout 1966/67
\cite{Stout_1966/67} which is followed by Fedorov 1990/91
\cite[p.\,271]{Fedorov_1991}.

$\bullet$ boundary preserving maps (Jenkins-Suita 1984
\cite{Jenkins-Suita_1988}); maps taking the boundary into the
boundary (Landau-Osserman 1960 \cite{Landau-Osserman_1960}).

$\bullet$ complete covering surfaces (cf. Ahlfors-Sario 1960
\cite[p.\,41--42, \S\,21A]{Ahlfors-Sario_1960}), i.e. one such
that any point in the range has a neighborhood whose inverse
image consists only of compact components; complete Klein
coverings (Andreian Cazacu 2002 \cite{Andreian-Cazacu_2002})
(a direct extension of the former concept shown to be
equivalent in the case of finite coverings to the next
conception of Sto\"{\i}low).

$\bullet$ total Riemann coverings (Sto\"{\i}low 1938
\cite{Stoilow_1938-Lecons}), i.e. one such that any sequence
tending to the boundary has an image tending to the boundary.

$\bullet$ unlimited covering surfaces (Nakai 1988
\cite{Nakai_1988}, EDM=Japanese encyclopedia 1968/87
\cite{EDM_1968/87}, Minda 1979
\cite{Minda_1979-hyperbolic-metric-and-coverings})

$\bullet$ proper (holomorphic) maps (onto the unit disc)
(e.g., Bedford 1984 \cite[p.\,159]{Bedford_1984}, Bell 1999
\cite[p.\,329]{Bell_1999-Ahlfors-maps}, \v{C}erne-Flores 2007
\cite{Cerne-Flores_2007}, Fraser-Schoen 2011
\cite{Fraser-Schoen_2011}).

$\bullet$ distinguished map (Jurchescu 1961
\cite{Jurchescu_1961-A-maximal})

$\bullet$ Myrberg surface over the unit disc (Stanton 1975
\cite[p.\,559, \S\,2]{Stanton_1975} uses this terminology for
a Riemann surface $W$ admitting an analytic function $z\colon
W \to \Delta$  realizing $W$ as an $n$-sheeted, branched, full
covering surface of the unit disc $\Delta$).

As no ramification appears along the boundary, explains the
naming:

$\bullet$ Randschlicht mapping (K\"oditz-Timann 1975
\cite{Koeditz-Timmann_1975}).

\medskip

In fact the writer came across this concept through real algebraic
geometry where I used (2006 \cite{Gabard_2006}) the term {\it
saturated}, whereas Coppens 2011 \cite{Coppens_2011} proposes the
term {\it separating} morphism. In the same context, Geyer-Martens
1977 \cite{Geyer-Martens_1977} coined:

$\bullet$ ``total reell Morphismus''=totally real morphism/map, to
which we shall adhere as it seems to be the most convenient
terminology, especially when abridged just as ``total maps'',
which is quite in agreement with Sto\"{\i}low's jargon.

We shall attempt to reserve the designation {\it Ahlfors
maps/functions} for those solving the extremal problem formulated
in Ahlfors 1950 \cite{Ahlfors_1950}. The latter are known (since
Ahlfors 1950 \cite{Ahlfors_1950}) to be circle maps, but the
converse is wrong. Indeed, circle maps may have arbitrarily
large degrees (post-compose with a power map $z\mapsto z^n$ for
some large integer $n$), whereas Ahlfors maps have degrees $\le
r+2p$ (in view of the deep result in Ahlfors 1950
\cite{Ahlfors_1950}).

Are circle maps of degree compatible with Ahlfors' bound always
realizable via an Ahlfors map? The answer seems to be in the
negative, at least if  attention is restricted to infinitesimal
Ahlfors maps. This follows from Gouma's restriction (1998
\cite{Gouma_1998}) in the hyperelliptic case. Indeed consider a
2-gonal membrane, then post-composing with $z\mapsto z^n$ we get
circle maps of degrees ranging through all multiples $2n$, whereas
only 2 and $r+2p$ are realized as degrees of Ahlfors maps, by a
result of Gouma 1998 \cite{Gouma_1998}. Note  that Gouma restricts
to ponctual Ahlfors maps and our claim is only firmly established
in this context.

A somewhat deeper question is whether any (or at least one) circle
map of smallest degree arises via an Ahlfors map. We were not able
to settle this question, but in a tour de force Yamada 2001
\cite{Yamada_2001}
proved this in the hyperelliptic case. It
amounts to know if the Ahlfors map is flexible enough to capture a
circle map of the lowest possible degree (alias the gonality). Let
us optimistically pose the conjecture, amounting to say that we
can essentially take out the best of the two worlds:

\begin{conj}\label{gonality:conj}
Any (or at least one) conformal mapping realizing the gonality
arises as an Ahlfors extremal function $f_{a,b}$ (perhaps for
coalescing two points yielding then the Ahlfors map $f_a$
maximizing the modulus of the derivative at $a$).
\end{conj}

Recent work by Marc Coppens 2011 \cite{Coppens_2011} supplies a
sharp understanding of the gonality $\gamma$ as spreading through
all
permissible
values $r\le \gamma \le r+p$ when the membrane is varied
through its moduli space.

Paraphrased differently the conjecture wonders if a suitable
Ahlfors map always realizes the gonality. As yet we lack evidence,
but the vague feeling that Ahlfors' method is the best possible
(being distilled by the paradigm of extremality) inclines one to
believe that its economy should be God given. In
contradistinction, it may be argued that Ahlfors maps depend on so
few parameters (essentially one or two points on the surface),
that they are perhaps not flexible enough to explore the full room
of all circle maps.
Such simple minded question
exemplifies that the old subject of the Ahlfors' map still
deserves better understanding.
A fine understanding of the Ahlfors map would truly be worth
studying if we had some clear-cut applications in mind (taking
full advantage of the extremal property of the map). In practice,
one is often content with the weaker notion of circle maps, but in
the long run it is likely that more demanding applications
requires the full punch of the Ahlfors map.

\subsection{Sorting out applications: finite vs. infinite/compact vs. open}

As to applications (of the Ahlfors map), there are several
ramifications, which
---at the risk of oversimplification---may be ranked in two
headings ({\it in finito} vs. {\it in infinito}). By this we have
in mind essentially the sharp opposition between compact and
non-compact Riemann surfaces. The later were intensively
approached by several schools (mostly Finnish, Japanese and US),
but the theory is certainly less complete than for compact
surfaces, which from our viewpoint already represent a serious
challenge. Furthermore it is evident that there is essentially one
and only one road leading from the finite to the infinite namely
the exhaustion process affording a cytoplasmic expansion of a
compact bordered Riemann surface in some ambient open surface. Now
let us enumerate such applications.

(A) {\bf Lifting truth from the disc via conformal
transplantation.} A reliable philosophy is roughly that a
result known to hold good in the disc is lifted via the
Ahlfors map to configurations of higher topological type. This
is the strategy used by Alling 1964 \cite{Alling_1964} to
transplant the corona of Carleson 1962 \cite{Carleson_1962} to
Riemann surfaces. (The corona theorem amounts to say that the
Riemann surface is dense in the maximal ideal space of its
algebra of bounded analytic functions.)

In spectral theory this method
(systematically utilized by Poly\'a-Szeg\"o) is known as
``{\it conformal transplantation}''. Subsequent elaborations
arose through the work of Hersch 1970 \cite{Hersch_1970} and
Yang-Yau 1980 \cite{Yang-Yau_1980} (where branched covering
are admitted, thereby diversifying widely the topology).

Recently Fraser-Schoen 2011 \cite{Fraser-Schoen_2011} applied
the Ahlfors
mapping to spectral theory (Steklov eigenvalues). (This
inspired a
note of the writer \cite{Gabard_2011} extending Hersch 1970's
study of Dirichlet and Neumann eigenvalues on spherical
membranes to arbitrary (compact) bordered surfaces.) Another
spectacular work is
due to Girouard-Polterovich 2012 \cite{Girouard-Polterovich_2012}
where Fraser-Schoen's work is extended to higher eigenvalues.

(B) {\bf Exhaustion and infinite avatars.} Another philosophy
(Nevanlinna, Ahlfors, etc.) is to exploit the fact that (infinite,
i.e. open Riemann surface) may be exhausted by compact subregions
(reminding somehow the finitistic slogan of Andr\'e Bloch, {\it
``Nihil est in infinito...''}) offering thereby a wide range of
application of compact bordered Riemann surfaces to the more
mysterious realm of open Riemann surfaces. This ramifies quickly
to the so-called  classification theory of Riemann surface
(Nevanlinna 1941 \cite{Nevanlinna_1941}, Ahlfors 46, Sario 46--49,
Parreau 1951 \cite{Parreau_1951}, Royden 1952 \cite{Royden_1952},
etc.) much completed by the Japanese school (T\^oki 1951, A. Mori,
Kuramochi, Kuroda, etc.). Several books attempt to give a coherent
account of this big classification theory, e.g. Ahlfors-Sario 1960
\cite{Ahlfors-Sario_1960}, Sario-Nakai 1970
\cite{Sario-Nakai_1970}, where the guiding principle (due to Sario
1946) is to classify surfaces according to the force of their
ideal boundary.

In another infinite direction, S.\,Ya. Havinson 1961/64
\cite{Havinson_1961/64} was the first (with Carleson 1967
\cite{Carleson_1967-book}) to extend the theory of the Ahlfors
function to domains of infinite connectivity , and was followed by
S. Fisher 1969 \cite{Fisher_1969}, which propose some
simplifications.

The Slovenian school of complex geometry (\v{C}erne,
Forstneri\v{c}, Globevnik, etc.) are also employing the Ahlfors
function, often in connection with the open problem (Narasimhan,
Bell, Gromov, etc.) of deciding if any open Riemann surface embeds
properly in ${\Bbb C}^2$. In
Forstneri\v{c}-Wold 2009 \cite{Forstneric-Wold_2009} reduced the
full problem to a
finitary  question as to whether each compact bordered Riemann
surface embeds holomorphically in the plane ${\Bbb C}^2$. (Maybe
this is achievable by a suitable of Ahlfors functions, or more
sophisticated variant thereof like (?) in the broader
Pick-Nevanlinna context). As suggested by those authors, it is
maybe enough to embed one representant in each topological type
(this is possible, compare \v{C}erne-Forstneri\v{c} 2002
\cite[Theorem 1.1]{Cerne-Forstneric_2002}) and try to use a
continuity argument through the Teichm\"uller (moduli) space.

\section{Biased recollections of the writer}\label{Sec:Biased-recollections-of-Gabard}

\subsection{Klein's viewpoint: real curves as symmetric
Riemann surfaces (as yet another
instance of the Galois-Riemann Verschmelzung)}

If the writer is allowed to recollect his own memories about his
involvement with this circles of ideas, it started as follows.
Maybe a natural point of departure is the (basic) algebraic
geometry of curves. While reading Shafarevich's Basic algebraic
geometry (ca. 1998) one encounters  some nice drawings of the real
locus of a plane cubic into its complex locus materialized by a
torus (as we know since time immemorial: Euler?, Abel?, Jacobi,
Riemann, etc.). A torus of revolution reflected across a plane
cutting the torus along two circles yields a plausible
visualization of the embedding of $C({\Bbb R})$ into $C({\Bbb C})$
(even with the symmetry induced by the complex conjugation).

Of course there are also real cubic curves whose real loci possess
only one component. How to visualize the corresponding embedding?
Lee Rudolph quickly helped us by just realizing that the Galois
action (complex conjugation) acts over the torus $S^1\times S^1$
just by exchanging the two factors $(x,y)\mapsto (y,x)$ fixing
thereby the diagonal (circle) $\{(x,x)\}\approx S^1$.

More generally  how to picture out the topology of a real curve?
The first observation is that the complex locus $C({\Bbb C})$ is
acted upon by complex conjugation $\sigma$ relative to some
ambient projective space ${\Bbb P}^n({\Bbb C})$ (where after all
the concrete curve is embedded). Therefore to each real curve $C$
is assigned a {\it symmetric surface} $(C({\Bbb C}),
\sigma)=(X,\sigma)$ consisting of a pretzel $X$ together with an
orientation reversing involution $\sigma\colon X\to X$. (For
aesthetical reasons all of our algebraic curves are
projective and non-singular, prompting thereby compactness of the
allied Riemann surfaces.) With the
invaluable assistance of (overqualified scholars) Claude Weber and
Michel Kervaire, I learned how to classify such objects, according
to the invariants $(g,r,a)$ where $g$ is the genus of $X$, $r$ the
number of ``ovals'' (fixed under $\sigma$), and $a$ is the
invariant counting mod 2 the number of components of $X-{\rm
Fix}(\sigma)$. In other words $a=0$ corresponds to the separating
(or dividing) case where ${\rm Fix} (\sigma)$ disconnects $X$,
whereas $a=1$ means that the fixed locus does not induce a
morcellation of the surface.

I soon realized thanks to the paper Gross-Harris 1981
\cite{Gross-Harris_1981}, that all this material was a well-known
game for Felix Klein, who was essentially the first to classify
symmetric surfaces taking advantage of the just established
classification of compact bordered surfaces (M\"obius 1863
\cite{Moebius_1863}, Jordan 1866 \cite{Jordan_1866}, etc.). The
key trick is of course the yoga assigning to $(X,\sigma)$ its
quotient $X/\sigma=:Y$ by the involution, and moving upward again
via the orientation covering supplied by local orientations. If
the point lies on the boundary then there is no duplication of the
point by local orientations (alias ``indicatrix'' in older
literature).

\begin{theorem} {\rm (Klein 1876
\cite{Klein_1876}=\cite[p.\,154]{Klein-Werke-II_1922}, explicit in
Klein 1882 \cite{Klein_1882}, Weichold 1883 \cite{Weichold_1883})}
There is one-to-one correspondence between symmetric surfaces and
compact bordered surfaces. Moreover the correspondence extends to
the realm of conformal geometry, i.e. Riemann surfaces or Klein
surfaces, if you prefer.
\end{theorem}

\begin{rem} {\rm Modernized treatments of
this Klein correspondence---say compatible with Weyl--Rad\'o's
(1913/1925) abstract conception of the Riemann surface---are
plenty, compare, e.g. Teichm\"uller 1939 \cite[p.\,99--101, Die
Verdoppelung, \S 92,
93]{Teichmueller_1939}=\cite{Teichmueller_1982}, Schiffer-Spencer
1954 \cite[p.\,29--30, \S 2.2]{Schiffer-Spencer_1954},
Alling-Greenleaf 1971 \cite{Alling-Greenleaf_1971}. }

\end{rem}

Via this dictionary, it is plain that the dividing case
corresponds precisely to the orientable case. [As a matter of
terminology, Klein used (since Wintersemester 1881/82) the jargon
{\it orthosymmetrisch} versus {\it diasymmetrisch} corresponding
to the dividing respectively nondividing case. For instance
Weichold 1883 \cite[p.\,322]{Weichold_1883} writes:

\begin{quota}[Weichold 1883]\label{Weichold-1883:quote}
{\small \rm

Was ferner die symmetrischen Riemann'schen Fl\"a\-chen anbelangt,
deren Betrachtung die Grundlage der folgenden Untersuchung bildet,
so sind auch diese schon mehrfach behandelt worden, wenn auch zum
Theil unter ganz anderen Gesichtspunkten. Es hat sich n\"amlich
Herr Professor Klein in den B\"anden VII und X der Mathem. Annalen
in den Aufs\"atzen mit dem Titel: ``\"Uber eine neue Art von
Riemann'schen Fl\"achen'' mit diesen Fl\"achen eingehender
besch\"aftigt und daselbst auch schon die Hauptunterscheidung
derselben in orthosymmetrische und diasymmetrische Fl\"achen
aufgestellt. Diese Bezeichnung findet sich allerdings noch in
keiner Publication angewendet; sie wurde zuerst in einem in
Wintersemester 1881/82 von Herrn Professor Klein abgehaltenen
Seminar eingef\"uhrt, in welchem derselbe auch die weiter unten
erw\"ahnte weitergehende Classification mittheite und bei welchem
auch der Verfasser die unmittelbare Anregung f\"ur die vorliegende
Arbeit empfing.

}
\end{quota}

Perhaps it is worth tracking down further Klein's motivation for
this ``savant'' terminology; for this we
supply the following extract:

\begin{quota}[Klein 1923 {\rm
\cite[p.\,624]{Klein-Werke-III_1923}}]\label{Klein-1923:quote}

{\small \rm

Die Benennungen ``diasymmetrisch'' und ``orthosymmetrisch'' f\"ur
die beiden Klassen symmetrischer Fl\"achen wurden sp\"ater von mir
gerade wegen der im Text ber\"uhrten Verh\"altnisse eingef\"uhrt;
siehe Bd. 2 dieser Ausgabe, S.\,172. Vgl. auch Fu{\ss}note $^{
58)}$ auf S.\,565/566 im vorliegenden Bande. \quad K.

}
\end{quota}

So this brings us at other places, the first cross-reference leads
us to the following quote (whereas Fu{\ss}note $^{ 58)}$ is merely
a text written by Vermeil, not really worth reproducing here):

\begin{quota}[Klein 1892 {\rm
\cite{Klein_1892_Realitaet}=\cite[p.\,172]{Klein-Werke-II_1922}}]
\label{Klein-1892/22:quote}

{\small \rm

Reelle algebraische Kurven ergeben {\it symmetrische} Riemannsche
Fl\"achen und k\"onnen umgekehrt allgemein g\"ultig von letzteren
aus defieniert werden, das ist der hier fundamentale Satz, den ich
in \S 21 meiner Schrift entwickelte. Ich bezeichne dabei eine
Riemannsche Fl\"ache als symmetrisch, wenn sie durch eine konforme
Abbildung zweiter Art von der Periode 2 in sich \"ubergef\"uhrt
wird (i.e. durch eine konforme Abbildung, welche die Winkel
umlegt). Die symmetrischen Riemannschen Fl\"achen eines gegebenen
$p$ zerfallen, wie ich ebendort angab und Herr Weichold a.\,a.\,O.
eingehender ausgef\"uhrt hat, nach der Zahl und Art ihrer
``Symmetrielinien'' in $[\frac{3p+4}{2}]$ Arten. Wir haben
erstlich $[\frac{p+2}{2}]$ Arten {\it orthosymmetrischer}
Fl\"achen bez. mit $p+1, p-1, p-3, \dots$ Symmetrielinien; das
sind solche symmetrische Fl\"achen, welche l\"angs ihrer
Symmetrielinien zerschnitten, in zwei (zueinander symmetrische)
H\"alften zerfallen; --- das einfachste (zu $p=0$ geh\"orige)
Beispiel ist eine Kugel, welche durch ``orthogonale'' Projektion
auf sich selbst bezogen ist~---. Wir haben ferner $(p+1)$ Arten
{\it diasymmetrischer} Fl\"achen bzw. mit $p, p-1, \dots, 1, 0$
Symmetrielinien; das sind Fl\"achen, die l\"angs ihrer
Symmetrielinien zerschnitten gleichwohl noch ein
zusammenh\"angendes Ganzes vorstellen;~---~man vergleiche bei
$p=0$, die durch eine ``diametrale'' Projektion auf sich selbst
bezogene Kugel.~---

 }
\end{quota}

Hence to summarize this explanation of Klein, the fundamental
dichotomy seems to be motivated by the basic case of genus $0$
(the sphere), which may be acted upon in two fashions by a
sense-reversing involution (orthogonal vs. diametral). This basic
motivation is even more emphasized in Klein's lectures, worth
reproducing (despite its very elementary character):

\begin{quota}[Klein 1891/92 {\rm \cite[p.\,138--9]
{Klein_1892_Vorlesung-Goettingen}}]\label{quote:Klein-1891/92-ortho/dia}

{\small \rm

Wir beginnen damit, anzugeben, auf wieviel verschiedene Weisen
eine Kugel mit sich selbst symmetrisch sein kann (d.\,h. durch
eine $\Sigma$ von der Periode $2$ in sich selbst \"ubergehen
kann). Das ist offenbar auf $2$ wesentlich verschiedene Arten
m\"oglich: das eine Mal bezieht man die Kugel auf sich selbst
durch eine Centralprojection, deren Centrum au{\ss}erhalb liegt:

\begin{figure}[h]
\centering
    \epsfig{figure=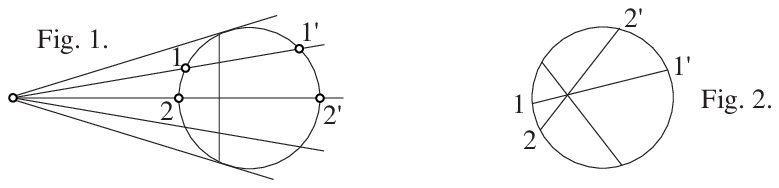,width=92mm}
\end{figure}

\noindent ($1,1';2,2'; \dots$ sind entsprechende Puncte), das
zweite mal durch eine Centralprojection, deren Centrum sich
innererhalb der Kugel befindet.

Im ersten Falle giebt es auf der Kugel eine sogenannte
\underline{Symmetrielinie}, deren Puncte bei der Umformung
s\"ammtlich festbleiben, das ist der Schnitt der Kugel mit der
Polarebene des Projectionscentrums; im $2^{\rm ten}$ Falle giebt
es eine solche Symmetrielinie nicht. Wir haben damit dasjenige
Unterscheidungsmerkmal, nach welchem wir sogleich die
symmetrischen Fl\"achen einteilen: nach \underline{der Zahl und
Art der Symmetrielinien.} Erw\"ahnen wir da gleich die
Terminologie, welche ich anl\"a{\ss}lich der Figuren 1 und 2 in
Vorschlag gebracht habe. Figur 1 kann insbesondere so gezeichnet
werden, da{\ss} das Projectionscentrum unendlich weit liegt. Die
Polarebene wird dann eine Dia\-metralebene und die zugeh\"orige
Centralprojection eine orthogonale Projection. Ich sage
dementschprechend \"uberhaupt von der Figur 1, die Kugel sei bei
der selben \underline{orthosymmetrisch} auf sich selbst bezogen.
Die bei Figur 2 vorliegende Beziehung aber nenne ich
\underline{diasymmetrisch}, insofern bei ihr das
Projectionscentrum, insbesondere in den Mittelpunkt der Kugel
r\"ucken kann, worauf je zwei diametrale Puncte der Kugel
zusammengeordnet erscheinen. Diese Benennungen
``orthosymmetrisch'' u. ``diasymmetrisch'' \"ubertrage ich dann
demn\"achst in noch zu erkl\"arender Weise auf die Fl\"achen eines
beliebigen $p$.

}

\end{quota}

Ahlfors result precisely affords a deeper function-theoretical
propagation of this Kleinian paradigm: orthosymmetric surfaces are
precisely those mapping in totally real way to the orthosymmetric
sphere!


The Russian school (Gudkov, Rohlin, Kharlamov, Viro, etc.) uses
the (less imaginative) nomenclature Type I versus Type II, whose
labelling is pure convention vintage; yet still a heritage from
Klein's initial nomenclature of 1876
\cite{Klein_1876}=\cite[p.\,154]{Klein-Werke-II_1922} reproduced
in the follwing:

\begin{quota}[Klein 1876]\label{Klein-1876:quote}
{\small \rm

Andererseits ergibt sich f\"ur die Kurven, deren Z\"uge\-zahl
$C>0$, $C<p+1$ eine bemerkenswerte Einteilung in zwei Arten.

{\it Die Kurven der ersten Art haben die Eigenschaft, da{\ss} ihre
Riemannsche Fl\"ache, l\"angs der $C$ Z\"uge zerschnitten,
zerf\"allt: bei den Kurven der zweiten Art findet ein solches
Zerfallen nicht statt.}

}
\end{quota}

Rohlin 1978 \cite[p.\,90]{Rohlin_1978} refers explicitly to Klein
as follows:

\begin{quota}[Rohlin 1978]\label{Rohlin:quote}

{\small \rm

Following Klein (see [4], p.\,154), we say that $\alpha$ belongs
to type~I if $A$ splits ${\Bbb C}A$ and to type~II if $A$ does not
split ${\Bbb C}A$. For example, $M$-curves obviously belong to
type~I.

}
\end{quota}

It is worth recalling that Rohlin made a
surprisingly late discovery of Klein's work as shown by the
following extract:

\begin{quota}[Rohlin 1978 {\rm \cite[p.\,85]{Rohlin_1978}}]\label{Rohlin2:quote}

{\small \rm

As I learned recently, more than a hundred years ago, the problem
of this article occupied Klein, who succeeded in coping with
curves of degree $m\le 4$ (see [4], p.\,155). I do not know
whether there are publications that extend Klein's investigations.

}
\end{quota}

It is concomitant  to speculate that the infamous {\it Klein
bottle} (={\it Kleinsche Fl\"ache} which traversed the Atlantic as
a ``Flasche'') probably originated during Klein's study of real
curves. It just amounts to have a real curve of genus one without
real points, whose complex locus will be a torus (of revolution)
acted upon by a diametral involution $(x,y,z)\mapsto (-x,-y,-z)$.

\subsection{Criterion for Klein's orthosymmetry=Type I,
in Russian (Klein 1876--82;
Rohlin 1978, Fiedler 1978 vs. Alling-Greenleaf 1969, Geyer-Martens
1977)}

Klein's dichotomy for symmetric surfaces prompts for criterion
detecting the dividing character of a real curve.

The writer knows of essentially two methods: the first
being genetic and the other qualifiable of synthetic. Despite
their simplicity those criterions where overlooked by Klein, who
relied upon more complicated arguments (cf. the following optional
remark).

\begin{rem} {\rm
Besides, there are several other original methods due to Klein.
One involves the dual curve, and more specifically a
representation assigning to each imaginary point of the curve the
real line passing through it and its conjugate. When the points
becomes real the limiting position of this secant becomes the
tangent. In this way Klein manages to visualize the complex locus
of a plane curve living in the 4D-space ${\Bbb P}^2({\Bbb C})$
onto a the 2D real projective plane as a multiple cover, and to
guess the type of the curve. Beautiful pictures are to be found in
vol. II of his Ges. math. Abhandl. \cite{Klein-Werke-II_1922}).
Another brilliant argument of Klein involves a degeneration to the
hyperelliptic case.}
\end{rem}

\begin{figure}[h]
\centering
    \epsfig{figure=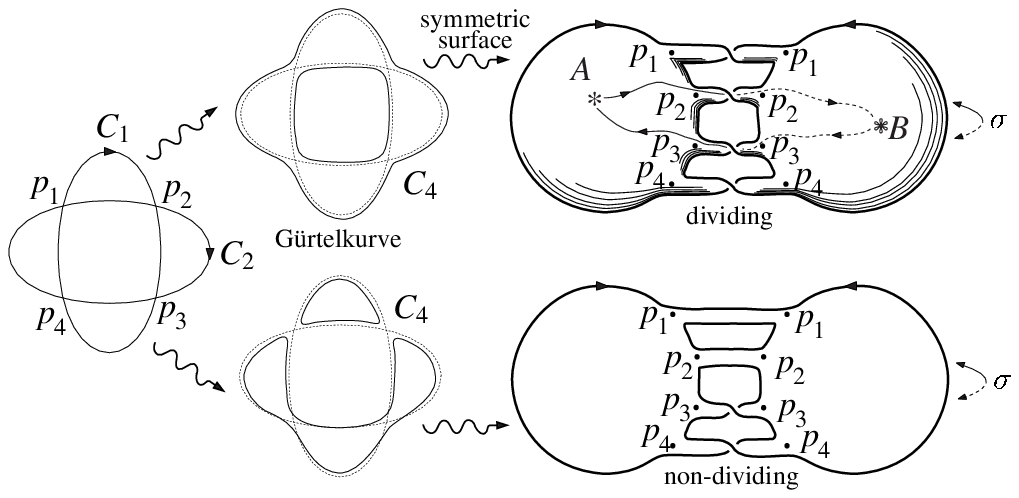,width=102mm}
\vskip-5pt\penalty0
  \caption{\label{Guertel-genetic:fig}
  Dividing character of the G\"urtelkurve via surgery,
  and nondividing character of a curve with unnested ovals}
\vskip-5pt\penalty0
\end{figure}

{\bf Genetic method.} This
is essentially a {\it surgery} (if we may borrow the jargon of
Thom, Milnor, etc.), and applies primarily to curves gained by
small perturbation of two curves whose type is known. Maybe it is
best explained on a specific example. Consider the {\it
G\"urtelkurve} as a small deformation of two conics having two
nested ovals. ({\it G\"urtel} means ``belt'', a nomenclature
coined by Klein in 1876
\cite{Klein_1876_Verlauf}=\cite[p.\,111]{Klein-Werke-II_1922},
presumably as a translation of the term ``{\it quartique
annulaire\/}'' used by Zeuthen in 1874 \cite[p.\,417+Tafel
I.,\,Fig.\,1]{Zeuthen_1874}.) Each conic corresponds to an
equatorial sphere, and
each smoothing amounts attaching a handle. During the process one
can keep track of the two real braids to make a global drawing of
the surface (compare right part of
Figure~\ref{Guertel-genetic:fig}). Some contemplation of the
drawing
shows that when all smoothings are dictated by  orientations then
the resulting curve is dividing. Thus the G\"urtelkurve is
dividing. Indeed in this case all handles contains twisted braids
and thus when travelling in the imaginary locus, say starting from
position $A$ in the north (top) hemisphere of the left sphere and
moving to the right sphere via an handle we reach position $B$ in
the south hemisphere of the right sphere. Coming back to the left
sphere, the twisting forces a return to the north hemisphere. We
are thus never able to visit the south hemisphere of the left
sphere.

\begin{figure}[h]
\centering
    \epsfig{figure=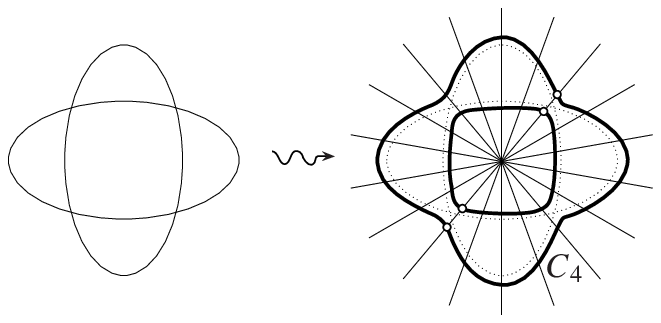,width=62mm}
\vskip-5pt\penalty0
  \caption{\label{Guertel-saturated:fig}
  G\"urtelkurve is dividing via a
  totally real
  morphism to the line.}
\vskip-5pt\penalty0
\end{figure}

{\bf Synthetic method.} Another way to see the dividing character
of the G\"urtelkurve involves looking at the pencil of lines
through a point lying deepest inside the two nested ovals
(Figure~\ref{Guertel-saturated:fig}). Since each real line of this
pencil cuts the quartic $C_4$ along a totally real collection of
points, this induces a map between the imaginary loci $C_4({\Bbb
C})-C_4({\Bbb R})\to {\Bbb P}^1({\Bbb C})-{\Bbb P}^1({\Bbb R})$.
It follows that $C_4$ is dividing since ${\Bbb P}^1$ is. (Just use
the fact that the continuous image of a connected set is
connected.)

More generally, this argument gives the following
criterion (which
quite curiously seems to have escaped Felix Klein's attention, cf.
e.g. his lectures notes 1891--92
\cite[p.\,168--69]{Klein_1892_Vorlesung-Goettingen}, where in our
opinion Klein draws the orthosymmetric character of the
G\"urtelkurve from more complicated arguments than those just
given):

\begin{lemma}\label{saturated:lemma} If a real curve permits a
morphism to the line whose fibers over real points are
exclusively real, then the curve is dividing.
\end{lemma}

Conversely, one may wonder if any dividing curve is expressible as
such a totally real cover of the line. I clearly remember having
asked this question at several experts (ca. 1999),
yet without receiving clear-cut answers, and so decided to embark
on a self-study of this
question. Being a slow and superficial worker, I needed circa 2
years of work until getting an answer, which turned to be
positive:

\begin{prop}\label{saturated:Gabard} {\rm (Gabard 2001, first
published in 2004)} Any dividing real curve admits a totally real
morphism to the line. Moreover the degree of such a morphism can
always be chosen $\le g+1$, where $g$ is the genus of the curve.
\end{prop}

Having completed this work, I started some
detective work, and via papers of Geyer-Martens 1977
\cite{Geyer-Martens_1977} and Alling-Greenleaf 1969
\cite{Alling-Greenleaf_1969} (probably located via the
bibliography of a survey by Natanzon 1990 \cite{Natanzon_1990/90})
realized that L.\,V. Ahlfors already proved this result in 1950
(and even exposed his results at Harvard in 1948 as reported in
Nehari 1950 \cite{Nehari_1950}). This was a great deception, or
rather more my first contact with the (glamorous) L.\,V. Ahlfors.

{\bf Very anecdotic details:} However as Ahlfors' result was not
fairly well-known (among the real algebraic geometry community) I
received a nice invitation to expose this re-discovery in a
RAAG-conference at Rennes in 2001. It was a great pleasure to meet
for the first time great specialists like Johannes Huisman,
Natanzon, Finashin, Viro, etc. My original proof involved an
argument with incompressible fluids and Abel's theorem to prove
(\ref{saturated:Gabard}). Some one week after the talk (or maybe
even during the week of that conference yet preceding my talk), I
confusedly realized that my argument was probably vicious, and
reworked it completely to find a topological parade, amounting to
the paragraphs 5,6 of Gabard 2006 \cite{Gabard_2006}. This
argument looked more tangible and I was again invited to Rennes in
2001--2002 (by J. Huisman) to present it at a specialized seminar.
At this
stage I started to believe that one could improve the bound $g+1$
into $\frac{r+g+1}{2}$, which is the mean value of the number of
ovals $r$ and the so-called {\it Harnack bound} $r\le g+1$. (In
the abstract setting is truly a remark of Klein directly reducible
to Riemann's definition of the genus as the maximal number of
retro-sections practicable on the pretzel without disconnecting,
compare Klein 1876 \cite[\S 7]{Klein_1876}.) I needed some weeks
(or months?) to establish this sharper version which gave
a relative progress over Ahlfors.

\begin{theorem}\label{saturated_new:Gabard}
{\rm (Gabard 2002, published 2004, 2006 {\rm\cite{Gabard_2006}})}
Any dividing real curve admits a totally real morphism to the line
${\Bbb P}^1$ of degree
 $\le \frac{r+g+1}{2}$, where $g$
is the genus of the curve and $r$ the number of ``ovals''
(=reellen Z\"uge).
\end{theorem}

Using the Schottky(-Klein) double of a compact bordered Riemann
surface (whose genus is visually seen to be $g=(r-1)+2p$) this can
be translated as

\begin{theorem}\label{saturated_new_bordered:Gabard}
Any compact bordered Riemann surface with $r$ contours of genus
$p$ is conformally representable as full covering of the disc of
degree
 $\le r+p$.
\end{theorem}


\section{Dirichlet's principle
(\"Uberzeugungskraft vs. mathematical
comedy)}\label{Sec:Dirichlet}


This section (with parenthetical title derived from jokes by
Hilbert 1905 \cite{Hilbert_1905} and Monna 1975
\cite{Monna_1975} resp.) recalls the early
vicissitudes of a principle supported by strong physical
evidence (as early as Green 1828 \cite{Green_1828} in print),
which Riemann
placed as the grounding for the edification of  the theory of
conformal mappings (and the allied Abelian integrals). This
section can be skipped without any further ado, but it fixes
the
context out of which emerged (simpler?) variational problems
 more suited to pure function-theoretical purposes. However,
Dirichlet's principle (after Hilbert's resurrection) pursued
his life (especially in the fingers of Courant) and merged
again to our main topic of the Ahlfors mapping (at least in
the schlichtartig situation handled by
Riemann-Schottky-Bieberbach-Grunsky). Of course, this
``Dirichlet'' line of thought is very active today, e.g., by
Hildebrandt and his collaborators. In short, Dirichlet's
principle flourished above any expectation by Riemann, was
``killed'' by Weierstrass, but resurrected by Hilbert, yet
re-marginalized by
extremal
methods
(Fej\'er-Riesz, Carath\'eodory, Ostrowski, Grunsky, up to
Ahlfors) and
re-flourished by Douglas and Courant as a (reliable)
instrument for the existence of conformal mappings.

\subsection{Chronology (Green 1828, Gauss 1839, Dirichlet ca. 1840,
Thomson 1847, Kirchhoff 1850, Riemann 1851--57, Weierstrass
1859/70, etc.)}

Apart from
a early contribution of Gauss 1825 \cite{Gauss_1825} about
local isothermic
parameters
(conformal mappings in the small),
the ``global'' theory of
such mappings
emerged
from Riemann's Thesis 1851 \cite{Riemann_1851} and his subsequent
work 1857 \cite{Riemann_1857} on abelian functions. A landmark is
the  {\it Riemann mapping theorem} (RMT) (cf.
 Riemann 1851 \cite{Riemann_1851}, and Riemann 1857 \cite{Riemann_1857}), derived from the so-called {\it
Dirichlet principle}. This was apparently formulated by Dirichlet
as long ago as the early 1840's (lectures in Berlin, attended by
Riemann in 1847/49). (The G\"ottingen 1856/57 version of those
were published by Grube in 1876 as \cite{Dirichlet_1840-1876}.)
Independent formulations (or utilizations) of this
principle are due to Gauss 1839 \cite{Gauss_1839}, Thomson
1847 \cite{Thomson_1847} (popularizing the long neglected work
of Green 1828 \cite{Green_1828}) and Kirchhoff 1850
\cite{Kirchhoff_1850}.
It is known that Riemann knew all those works (when exactly in
another question) from a manuscript estimated 1855/60
reproduced below (source=Neuenschwander 1981
\cite[p.\,225]{Neuenschwander_1981}). Riemann does not cite
Thomson and Kirchhoff in 1857 \cite{Riemann_1857}.

\begin{quota}[Riemann 1855/60]

{\small\rm

Mit dem Namen des Dirichlet'schen Princip's habe ich eine
Mehode bezeichnet, um nachzuweisen, da{\ss} eine Function
durch eine partielle Differentialgleichung und geeignete
lineare Grenzbedingungen v\"ollig bestimmt ist,
d.\,h.\,da{\ss} die Aufgabe, eine Function diesen Bedingungen
gem\"a{\ss} zu bestimmen, eine L\"osung und zwar nur eine
einzige L\"osung zul\"a{\ss}t. Es ist diese Methode von
William Thomson in seiner Note Sur une \'equation aux
diff\'erences (Liouville.\,T.\,12.\,p.\,493.) und von
Kirchhoff in seiner Abhandlung \"uber die Schwingungen einer
elastischen Scheibe angewandt worden, nachdem Gau{\ss} schon
vorher eine Aufgabe, welche als ein specieller Fall dieser
Aufgabe betrachtet werden kann, \"ahnlich behandelt hatte
(Allgemeine Lehrs\"atze.\,Art.\,29--34.) Ich habe diese
Methode nach Dirichlet benannt, da ich von Hrn Professor
Dirichlet erfahren hatte, da{\ss} er sich dieser Methode schon
$\langle$seit dem Anfang der vierziger Jahre (wenn ich nicht
irre) [Bl.\,66r]$\rangle$ in seinen Vorlesungen bedient habe.

}

\end{quota}

There is also a letter of Riemann dated 30.\,Sept.\,1852
(cf.\,Neuenschwander 1981 \cite{Neuenschwander_1981-lettres}),
where it is reported that Dirichlet supplied some references
to Riemann. Here is the relevant extract, out of which we may
speculate that Riemann learned the ref. to Thomson and
Kirchhoff  at this occasion (through Dirichlet).

\begin{quota}[Riemann 1852, 30. Sept.]
{\small\rm    Am Freitag Morgen, um in meinem Berichte
fortzufahren, suchte Dirichlet mich in meinem Zimmer auf. Ich
hatte ihn bei meiner Arbeit um Rath gefragt und er gab mir nun
die dazu n\"othigen Notizen so vollst\"andig, da{\ss} mir
dadurch die Sache sehr erleichtert ist. Ich h\"atte nach
manchen Dingen auf der Bibliothek sonst lange suchen k\"onnen.
D.[irichlet] war \"uberhaupt \"au{\ss}erst nett theilte mir
mit, womit er sich in den letzten Jahren besch\"aftigt hatte,
ging meine Dissertation mit mir durch; und so hoffe ich,
da{\ss} er mich auch sp\"ater nicht vergessen und mir seine
Theilnahme schenken wird. }

\end{quota}

 As we know the
principle was disrupted by the (non-fatal) Weierstrass'
critique 1870 \cite{Weierstrass_1870},
but
resuscitated by Hilbert in 1900-1 \cite{Hilbert_1900}
\cite{Hilbert_1901/04} \cite{Hilbert_1905}, after partial results
by Neumann 1870, 1878 \cite{Neumann_1878}, and 1884
\cite{Neumann_1884} Schwarz 1869/70
\cite{Schwarz_1869-70_Zur-Theorie-der-Abbildung}, 1870
\cite{Schwarz_1870}, 1872 \cite{Schwarz_1872} ({\it alternierendes
Verfahren}) and Poincar\'e for fairly general boundary contours.

Dirichlet's principle (as Riemann christened it in 1857
\cite{Riemann_1857-DP}) amounts to solve the first boundary
value problem for the Laplacian $\Delta u=0$ by minimizing the
Dirichlet integral
$$\int\!\!\int \Big\{ \bigl(\frac{\partial u}{\partial x}\bigr)^2+\bigl(\frac{\partial
u}{\partial y}\bigr)^2\Big\} dx dy.
$$
As a such the paradigm of {\it extremality} entered the arena
of geometric function theory since its earliest day, and
governed much of the subsequent developments.

Other noteworthy hot spots in this realm are:

$\bullet$ {\it The Bieberbach conjecture} (1916
\cite{Bieberbach_1916-BC}) $\vert a_n \vert \le n$ on the
coefficients of schlicht (=univalent=injective) functions from
the disc $\Delta=\{ \vert z \vert <1 \}$ to the (finite) plane
${\Bbb C}$ with Koebe's function
$k(z)=\frac{z}{(1-z)^2}=z+2z^2+3z^3+\dots$ as unique
extremals among those satisfying the normalization
$f(0)=0,f'(0)=1$. Completely solved
by de Branges 1984.

$\bullet$ Gr\"otzsch-Teichm\"uller extremal quasi-conformal
mappings (1928--1939 \cite{Teichmueller_1939}), i.e. the
search of the ``m\"oglichst konform'' mapping relating two
configurations. This gave a sound
footing to Riemann's liberal study of the moduli spaces (1857
\cite{Riemann_1857}), and
paved the way to the modern theory of deformation of complex
structures (Kodaira-Spencer).

\subsection{Early
suspicions about the Dirichlet principle (Weierstra{\ss}
1859/70, Schwarz 1869, Prym 1871, Hadamard 1906)}

Weierstrass seems to have been the first
to express doubts about the Dirichlet principle,
pivotal to Riemann's theory. Weierstrass lectured on his critique
in 1870, and this appeared in print as late as 1894 in his Werke.
However it is known that a meeting between Riemann and Weierstrass
took place in Berlin, 1859, where this issue was discussed.
Klein reports upon
Riemann's reaction at several places:

\begin{quota}[Klein 1926
{\rm \cite[p.\,264]{Klein_1926-Vorlesungen-über-die-Entwicklung}}]

{\small \rm

Er [Riemann] erkannte die Berechtigung und Richtigkeit der
Weierstra{\ss}chen Kritik zwar voll an; sagte aber, wie mir
Weierstra{\ss} bei Gelegenheit erz\"ahlte: ``er habe das
Dirichletsche Prinzip nur als ein bequemes Hilfsmittel
herangeholt, das gerade zur Hand war---seine Existenztheoreme
seien trotzdem richtig.'' Weierstra{\ss} hat sich dieser
Meinung wohl angeschlossen. Er veranla{\ss}te n\"amlich seinen
Sch\"uler H.\,A. Schwarz, sich eingehend mit den Riemannschen
Existenzs\"atzen zu befassen und andere Beweise daf\"ur zu
suchen, was durchaus gelang.

}

\end{quota}

\begin{quota}[Klein 1923 {\rm
\cite[p.\,492, footnote~8]{Klein-Werke-III_1923}}] {\small \rm

Ich erinnere mich, da{\ss} Weierstrass mir bei Gelegenheit
erz\"ahlte, Riemann habe auf die Gewinnung seiner
Existenzs\"atze durch das ``Dirichletsche Prinzip'' keinerlei
entscheidenden Wert gelegt. Daher habe ihm auch seine
(Weierstrass') Kritik des ``Dirichletschen Prinzips'' keinen
besonderen Eindruck gemacht. Jedenfalls ergab sich die
Aufgabe, die Existenzs\"atze auf andere Art zu beweisen. Diese
d\"urfte dann Weierstrass seinem Spezialsch\"uler Schwarz
\"ubertragen haben, bei dem er die erforderliche Verbindung
geometrisch-anschaulichen Denkens mit der F\"ahigkeit,
analytische Konvergenzbeweise zu f\"uhren, bemerkt hatte.

}

\end{quota}

A more detailed chronology is roughly as follows (cf.
Elstrodt-Ullrich 1999 \cite[p.\,285--6]{Elstrodt-Ullrich_1999}):

$\bullet$ In the late 1850s Weierstrass notices some gap in
the Dirichlet principle (DP), and
presents his objection
to Riemann in 1859, who is not tremendously affected claiming
that his existence theorems keep however their truths.

$\bullet$ Thieme 1862, who met Riemann and requested from him
some elucidations about his theory of Abelian functions, and
the conversation turned to the foundation of the Dirichlet's
principle. This is materialized by a letter of Thieme to
Dedekind of 1878 (reproduced in Elstrodt-Ullrich 1999
\cite[p.\,270--1]{Elstrodt-Ullrich_1999}, or as
Quote~\ref{quote:Thieme} below)

$\bullet$ Kronecker 1864, in a discussion with Casorati,
also
exposes some criticism of the (DP). This is materialized by
notes taken by Casorati, and published by Neuenschwander 1978

$\bullet$ Schwarz 1869
\cite[p.\,120]{Schwarz_1869-Ueber-einige-Abbildungsaufgaben}
expresses for the first time in print doubts about (DP)
(compare Quote~\ref{quote:Schwarz-1869} below).

$\bullet$ Heine February 1870 \cite[p.\,360]{Heine_1870} also
puts in print the reserves
expressed by Weierstrass and Kronecker, specifically their
objections to the assumption that a minimum must exist.

$\bullet$ Weierstrass July 1870 \cite{Weierstrass_1870}
presents a variational problem where the minimum is not
attained. This note, however, appeared in print only in 1895
in the second volume of Weierstrass's Werke
\cite{Weierstrass_1870}.

$\bullet$ Prym 1871 \cite[p.\,361--4]{Prym_1871} gives the first
(published) counterexample to the (DP) (as formulated, e.g. in
Grube's text 1876 \cite{Dirichlet_1840-1876} based upon
Dirichlet's lectures). Prym gives a continuous function on the
boundary of the unit disc such that the Dirichlet integral for the
associated harmonic solution to the Dirichlet problem is infinite.
However Prym expressly emphasizes that Riemann never stated such a
naive version of DP corrupted by Prym's example. In fact Prym's
example seems rather to attack a vacillating attempt by Weber 1871
\cite{Weber_1870} to rescue the Dirichlet principle.


\begin{quota}[Schwarz 1869 {\rm \cite[p.\,120]{Schwarz_1869-Ueber-einige-Abbildungsaufgaben}}]
\label{quote:Schwarz-1869}

{\small \rm

Dass es stets m\"oglich ist, die einfach zusammenh\"angende
Fl\"ache, welche von einer aus St\"ucken analytischer Curven
bestehenden einfachen Linie begrenzt ist, auf die Fl\"ache
eines Kreises zusammenhangend und in den kleinsten Theilen
\"ahnlich abzubilden, hat {\it Riemann} mit Zuh\"ulfenahme des
sogenannten {\it Dirichlet\/}schen Principes zu beweisen
gesucht.

Da gegen die Zul\"assigkeit dieses Principes bei einem
Existenzbeweise hinsichtlich der Strenge gegr\"undete
Einwendungen geltend gemacht worden sind, war es
w\"unschens\-werth, ein Beweisverfahren zu besitzen, gegen
welches die bez\"uglich des {\it Dirichlet\/}schen Principes
geltend gemachten Bedenken nicht erhoben werden konnten.

}

\end{quota}

\begin{quota}[Thieme 1878\,{\rm
:letter to Dedekind}] \label{quote:Thieme}

{\small \rm

Vielleicht werden Sie sich meiner noch erinnern, als ich mich
im Sommer 1862 in G\"ottingen aufhielt um bei Riemann
Aufkl\"arung \"uber seine Theorie der Abel'schen Funct. zu
erbitten. Ich traf Sie damals in der Krone, wo wir beide
abgestiegen waren, und das Gespr\"ach kam auf die, meiner
damaligen Meinung nach (was seitdem vielseitig anerkannt),
nicht ganz stichhaltige Begr\"undung des Dirichlet'schen
Princips, welches in der Riemann'sche Theorie fundamental ist.

}

\end{quota}

\section{Philosophical remarks}\label{Sec:Philosophical}

\subsection{Flexibility of 2D-conformal maps}

Maybe one way to enlarge slightly the discussion at the
philosophical level is to observe some unifying plasticity in
conformal maps. The underlying principle is roughly as
follows:

\begin{principle}[Conformal Plasticity (CP)]
If there is no topological obstruction to a mapping problem,
then a conformal mapping exist.
\end{principle}

This idea is very close to Koebe's allgemeines
Uniformisierungsprinzip in Koebe 1908 \cite{Koebe_1908_UbaK3},
which is stated as follow {\it Jedes Problem der im Sinne der
Analysis situs eine L\"osung hat kann auch funktiontheoretisch
verwirklicht werden.}

Of course this is not quite true
in view of say Riemann's moduli for closed (or non closed) Riemann
surfaces. However seminal instances where it works are the (RMT),
the uniformization theorem (UNI) [any simply-connected Riemann
surface is biholomorphic to the sphere, the plane or the disc],
and the more general Koebe schlicht theorem to the effect that a
schlichtartig Riemann surface is schlicht. Here the topological
condition of ``Schlichtartigkeit'' (i.e. any Jordan curve divides)
implies the stronger conformal embeddablility in the Riemann
sphere.

\begin{theorem} {\rm (Koebe 1908 \cite{Koebe_1908_UbaK3}, 1910 \cite{Koebe_1910_UAK2})}
Any dichotomic\footnote{This nomenclature is used by Hajek,
compare some arXiv preprints of the author joint with Gauld.}
(=schlichtartig) Riemann surface (i.e. one divided by any Jordan
curve)
embeds conformally in the Riemann sphere.
\end{theorem}


Since simply-connected implies dichotomic this implies the
(UNI) via (RMT).

An even stronger assertion is Bochner 1928 that any Riemann
surface of finite connectivity embeds in a closed Riemann
surface.

 Ahlfors' result about circle maps likewise illustrates the above
principle (CP), especially if we interpret it in Klein's theory of
symmetric surfaces (compare Lemma~\ref{saturated:lemma}).

\subsection{Free-hand pictures of some Riemann-Ahlfors maps}

It would be nice if some general methodology for picturing such
mappings could be developed. Let us try a naive look for domains
(Riemann surface are harder but not hopeless).

Maybe first a
comment by Poritsky 1949--52
\cite[p.\,21]{Poritsky_1949-52:Book}:

\begin{quota}[Poritsky 1949]

{\small \rm

From the above it is clear that analytical methods, at least
as developed thus far, have only limited power in solving the
complicated field problems arising in electric machines.
Electrical engineers have resorted extensively to the use of
``flux plotting'' or {\it free-hand drawing} of the flux lines
and equipotentials. As is well known, these curves, when drawn
for constant equal increments $\Delta \varphi=\Delta \psi$,
form a curvilinear set of {\it small squares}. A certain
aptitude, somewhat between mechanical drawing ability and
artistic drawing, is required for successful flux plotting,
and with practice people possessing such aptitude can learn to
draw flux plots for a great variety of cases with relative
ease.

} \end{quota}

The picture below (Fig.\,\ref{Riemann:fig}) is supposed to depict
the pullback of the radial-concentric bi-foliation of the disc via
a conformal representation of this 4-ply connected circle domain
by a Riemann map to the disc. (Usually the term ``Riemann map'' is
reserved for the simply-connected case, but recall that Riemann
was the first to prove the existence of such maps, cf. Riemann
1857/76 Nachlass \cite{Riemann_1857_Nachlass}). Physically one may
try to interpret it at the galvanic current generated by 4
batteries (electric charge) situated on a conducting plate.
Whenever the potential generated by two charge enter in conflicts
some saddle type singularity is generated (those can be counted
via an Euler characteristic argument \`a la Riemann-Hurwitz). In
the present case there is 6 saddles. In general $\chi D=d \chi
(\Delta) -\deg (R)$, and $\chi D=2-r$ (holed sphere) and the
degree $d$ is $d=r$, hence there is $\deg(R)=2r-2$ ramification
points. (This was of course well-known to Riemann, compare his
Nachlass, or Quote~\ref{quote:Riemann}). Dashed lines are
equipotentials.

\begin{figure}[h]
\centering
    \epsfig{figure=Rieman5.eps,width=122mm}
\vskip-5pt\penalty0
  \caption{\label{Riemann:fig}
  Attempting to plot a Riemann map by free-hand drawing}
\vskip-5pt\penalty0
\end{figure}

This sort of picture as mystical as it is (the reader confesses to
have had some trouble to generate it without grasping completely
the possible physico-chemical interpretation) gives the impression
of grasping slightly Riemann's title to his Nachlass
(Gleichgewicht der Electricit\"at), i.e. equilibrium potential of
 electricity. Our figure is pure
free-hand drawing without much scientific understanding. Thus it
would be nice if the computer can do better pictures, maybe via
the Bergman kernel (an eminently computable object, compare e.g.
Bell papers). In particular albeit it looks physically obvious, it
is not clear if the charge may be placed arbitrarily. (For
instance it is not clear why the corresponding divisor should be
linearly equivalent to its conjugate, compare Lemme~5.2 in Gabard
2006 \cite{Gabard_2006}.)

[Some related references: Henrici Computational conformal map,
Gaier, Konstruktive methoden in Konformen Abbildungen, etc... Or
maybe Crowdy via the Klein's prime]

Extracting some global understanding in the non-schlicht case of
such isothermic coordinates may be of some relevance to Gromov's
filling conjecture.

Further for less contours we may do similar pictures, and we then
obtain the following figures (Fig.\,\ref{Annulus:fig}).
The fact that the boundary contours are circles is not crucial
(but convenient for simple depiction). First we draw the
electrical forces in the case of an annulus. Then we made two
pictures for triply-connected with symmetrically disposed battery
(electrical charge). Geometrically those are supposed to be  the
pull-back of the origin under the Riemann-Ahlfors map. Finally we
would like to make a similar picture in the case where the charge
distribution is not symmetric. Then the picturing becomes very
difficult. Already in the symmetric cases it is hard to be
convinced that what we are doing is really serious. There is a
sort of subconscious algorithm to make such pictures: (1) first
draw the thick black lines where the particles enter in collision,
(2) draw at angle $\pi/4$ the dual saddle at those collision
point, and then the filling by thin lines is essentially a matter
of artistic feeling. Of course it is not always easy to arrange
such that all lines meet perpendicularly, but experience gives
some sort of algorithm to do this. Of course it is quite
convenient to do such pictures on a computer rather than on the
paper, as one can adjust trajectories by successive
approximations. As we used a software Adobe Illustrator; with
B\'ezier curves, thus the mathematical faithfulness of all this
picture is highly questionable, but we hope that the picture are
still of some qualitative value to help visualize such mappings,
and to feel some sort of physical interpretation. One guess is
that it amounts to have some positive electric charge at the
marked point plus perhaps a distribution of such charges on the
border. Then each positive particle is rejected by the charge and
the border. Thus the particle move faster when there is much
free-room in the plate. Alternatively one may have a biological
interpretation where the source are bacteroides and the black line
show the progression of the growing population which is faster in
those direction where there is much vital room (imagine
herbivores). Then the saddle amounts to junctions between various
ethnical population, and at time one the full universe is
explored. In this interpretation the proliferation of species can
be slowed down either by proximity to the border (limitation of
resources), or by vicinity of a competing population.

\begin{figure}[h]
\centering
    \epsfig{figure=Annulus.eps,width=122mm}
\vskip-5pt\penalty0
  \caption{\label{Annulus:fig}
  Another Riemann map by free-hand plotting}
\vskip-5pt\penalty0
\end{figure}

\subsection{Hard problems and the hyperelliptic claustrophobia}

Another unifying theme when it comes to hard problems
regarding Riemann surfaces is the following constat:

Several problems are fully settled in the hyperelliptic case,
but horribly complicated otherwise.

This is a paradigm well known since time immemorial. Probably one
of the first problem were it came acute was Jacobi's inversion
problem occupying Jacobi, then  G\"opel and Rosenhain
(hyperelliptic case) and only Weierstrass and above all Riemann
1857 \cite{Riemann_1857} could handle the general case.
(Weierstrass never managed to put in print his own approach
probably due to the extreme difficulty to follow an arithmetized
path.) Another place is Klein's trick of degeneresence to the
hyperelliptic configurations (cf. Klein 1892
\cite{Klein_1892_Realitaet}).

In some more contemporary problems we already addressed briefly
this hyperelliptic barrier also delineate the current frontier of
knowledge regarding:

(1) The Gromov filling area conjecture (compare the work by
Bangert et al. \cite{Bangert_2004} where the conjecture is
established in the hyperelliptic case, hence in particular for
membranes of genus $p=1$).

(2) The Forstneri\v{c}-Wold conjecture
\cite{Forstneric-Wold_2009} that compact bordered Riemann
surface embeds in ${\Bbb C}^2$ (this is also known in the
hyperelliptic case).

(3) The exact determination of Ahlfors degrees \`a la
Yamada-Gouma (this is also settled in the hyperelliptic case,
but not much seems to be known beyond those configurations).

\subsection{Topological methods}


We started our Introduction by claiming that topological
methods have some relevance to the field of function theory,
Riemann surfaces, and the allied fields.

The experience of the writer in this realm is rather modest and
essentially reduces to his lucky stroke in Gabard 2006
\cite{Gabard_2006}, about lowering the degree of a circle map upon
the prediction made by Ahlfors for his extremal function.

Such topological methods are quite common in function theory
(Riemann, Klein, Poincar\'e, Brouwer, Koebe, etc.) albeit
occupying a marginal place in comparison to potential-theoretic
consideration or the allied quantitative extremum problems. Let us
list some contributions using qualitative
 topological methods in the realm of classical function theory:

(1) The most famous (and probably important) example is the
continuity method of Klein-Poincar\'e related to the
uniformization problem. (Prior to this we may detect earlier trace
of the continuity method, as one learns by reading Koebe 1912
\cite{Koebe_1912_BdKm}, in the work of Schwarz-Christoffel and
Schl\"afli.)

(2) The intuitions of Klein-Poincar\'e were put on a firm footing
by Brouwer 1912 \cite{Brouwer_1912_Modulmannig},
\cite{Brouwer_1912_top-Schwierig}, using invariance of domain
which he was the first able to prove via combinatorial topology.

(3) Closer to our main topic, we cite Garabedian 1949
\cite{Garabedian_1949} who also relies heavily on
combinatorial topology to select appropriately certain auxiliary
parameters.

(4) Mizumoto 1960 \cite{Mizumoto_1960}, who reproves the existence
of (Ahlfors-type) circle map of degree $r+2p$ (i.e. like Ahlfors
bound).

(5) Gabard 2006 \cite{Gabard_2006}, where the degree is
lowered to $r+p$ (also via topological methods).

Roughly, the philosophy is that Riemann surfaces are volatile
objects when fluctuating through their moduli spaces, so that
practically nothing is observable outside the inherent topological
substratum which turns out to behave rather stably say w.r.t.  the
Abel-Jacobi mapping. At least this is philosophical substance of
the proof in Gabard 2006 \cite{Gabard_2006}.

\subsection{Bordered Riemann surfaces and real algebraic curves}

It may seem at first that bordered surfaces are a bit borderline
deserving less  respectableness than the temple of closed Riemann
surfaces. Likewise real algebraic geometry always appears like a
provincial subdiscipline of pure complex algebraic geometry of the
best stock.

Perhaps, less is true. The rehabilitation of reality within
algebraic geometry is in good portion, especially regarding the
connection with bordered (possibly) non-orientable surfaces, the
credit of Felix Klein (especially in 1882 \cite{Klein_1882}).
Moreover, independently of the algebro-geometric analytic
correspondence (\`a la Riemann, etc.) there is another simple
reason for which bordered (Riemann) surfaces took gradually
more-and-more importance during the 20th century, especially under
the fingers of the Finnish and Japanese schools. This pivotal
r\^ole of compact bordered objects results indeed merely from
their intervention as building
elements of general open surfaces. The latter being always
exhaustible through such compact elements as follows easily from
Rad\'o's triangulability theorem of 1925 \cite{Rado_1925}. Once
again this illustrates basically the philosophy ``Nihil est in
infinito\dots''.

Remind also that the device of exhaustion by finite (=compact)
surface is somewhat older than those schools. It may have first
occurred in Poincar\'e 1883 \cite{Poincare_1883} (where analytic
curves or open Riemann surfaces are first taken seriously and
subsumed to the uniformization paradigm) and then Koebe 1907
\cite{Koebe_1907_UbaK1} in same context. To caricature a bit
Koebe's proof, it amounts to use the RMT for compact discs (in a
version cooked by Schwarz) and expand in the large. The exhaustion
device is again used  in Nevanlinna 1941 \cite{Nevanlinna_1941},
where via exhaustions one constructs the corresponding so-called
harmonic measure solving the Dirichlet problem for boundary values
$0$ and $1$ on the initial resp. expanding contours of the
exhaustion $F_n$, yielding the ``Nullrand'' dichotomy according to
whether the $\omega_n$ flatten to $0$ or converge to a positive
function.  Ahlfors 1950 \cite{Ahlfors_1950} also uses (or planned
to use) a similar technique for other problem. This was enough to
launch the big classification programme of open Riemann surfaces.

\subsection{Lebesgue versus Riemann}

[11.10.12] This paragraph is free-style philosophical lucubration
coming to me right after reading the fantastic paper Forelli 1978
\cite{Forelli_1979}. From a
narrow minded viewpoint (the writer having
zero
measure
theoretic knowledge) it  seems that modernism, especially along
the ``capitalistic'' line of thought involving measure theory,
albeit initially quite concomitant with the (older complex)
function theory, ultimately may have drifted a vast body of the
vital fluid in a somewhat
arid valley. (For a somewhat related diagnostic cf. Morse-Heins
1947 \cite{Morse-Heins_1947}.)

Let us be more specific. Circa 1898  the way was paved toward
measure theory starting from function theoretic preoccupations
(not to mention the earlier ``Cantorism'' starting from Fourier
series). We have of course in mind E. Borel 1898
\cite{Borel_1898}, and then the stream along Lebesgue 1902
\cite{Lebesgue_1902}, Fatou 1906 \cite{Fatou_1906}, the old
brother F. Riesz 1907 (Fischer-Riesz effecting an Hilbert-Lebesgue
unification, etc.). All those grandiose efforts/achievements may
have polluted the pureness of (Riemann's) geometric conceptions by
charging the theory with complicated pathological paradigms not
truly inherent to its geometric substance (at least in its
finitistic aspects, which are not completely elucidated yet, e.g.
Gromov's filling conjecture). Of course the
antagonism we are speaking about goes back to older generations,
e.g. already acute in the Hermite vs. Jordan opposition, who were
resp. anti- and pro-Lebesgue\footnote{There is a letter form
Jordan to Lebesgue saying roughly: ``Pers\'everez dans vos
recherches math\'ematiques, vous allez y \'eprouver beaucoup de
plaisirs, mais il va vous falloir apprendre \`a y gouter seul, car
en g\'en\'eral les
g\'eom\`etres ne se lisent m\^eme pas entre eux-m\^emes.'' (quoted
by pure memory, hence highly unreliable).}).

This
tension is also felt when it comes to prove existence of circle
maps, where say proofs like Ahlfors' 1950 \cite{Mizumoto_1960},
Mitzumoto's 1960 \cite{Mizumoto_1960}, and many others (maybe even
Gabard's 2006 \cite{Gabard_2006}) proceeds along essentially
classical lines, often emphasizing the soft topological category
(very implicit by Riemann-Klein-Poincar\'e-Brouwer) instead of
measure theory (again Borel-Lebesgue-Fatou-Riesz). Of course
initially topology also arose from  capitalism over the real line,
namely the notion of metric (distance function). Yet ultimately
the theory (be it axiomatically
Bolzano-Cantor-Hilbert-Fr\'echet-Riesz-Hausdorff-Weyl or through
educated intuition Riemann-Klein-Poincar\'e-Brouwer-Thurston)
reached some higher romantic stratosphere producing some lovely
science essentially the most
remote from capitalistic preoccupation we were able to produce.
Alas or fortunately, Grisha Perelman (and precursors
Thurston/Yau-Hamilton) showed us that the likewise pleasant
Riemannian geometry (albeit slightly more quantitative) turned to
have some important topological impact (typically over
Poincar\'e's conjecture).

In a survey article by Lebesgue (ca. 1927, easy to locate), a
rather primitive mercantile metaphor is appealed upon to argue
that his theory of integration supersedes Riemann's. Lebesgue
argues that when a huge amount of money (delivered as a chaotic
mixture of pieces and bills) requires enumeration, his theory
amounts to count things properly by first enumerating what has
highest value and then paying attention to the more negligible
money pieces. This procedure is tantamount to subdividing rather
the range of the function as do Lebesgue instead of its domain as
did Cauchy or Riemann. The bulk of the US production (Rudin,
Gamelin, Forelli and many others) in the 1950-1970's is much
influenced by measure flavored analysis, and the art-form
continues to prosper with deep paradigms allied to Painlev\'e's
problem (fully solved in Tolsa 2003 \cite{Tolsa_2003}).

In contrast, some older workers, e.g. Koebe (cf. Gray's 1994 paper
\cite{Gray_1994}) as well as Lindel\"of (cf. Ahlfors' 1984
\cite{Ahlfors_1984-The-Joy}) (and probably more recent ones) were
never full partisans of Lebesgue's integral. Of course the latter
theory added a mass of grandiose contributions, yet in some
finitary problems like the one at hand (Ahlfors circle maps) its
significance can probably be marginalized, or completely
eliminated. So measure theory exists, but does it really capture
the quintessence of the problematic we are interested in, which is
more likely to be first of a {\it qualitative\/} nature (coarse
existence theory). Arguably, the next evolution step is the {\it
quantitative\/} phase (e.g. Ahlfors extremal problem, which is
essentially solved modulo fluctuating incertitudes about degree
variations of such maps). Finally any theory should culminate in
the {\it algorithmic\/} era, that is claustrophobic (computer
ripe) era. Remind that Riemann precisely disliked Jacobi's
approach, finding it too algorithmic and not conceptual enough
(according to some forgotten source, try maybe Klein's history
\cite{Klein_1926-Vorlesungen-über-die-Entwicklung}: ``Jacobi war
ihm zu algorithmisch.'' [quoted by memory]). At such a stage it is
safer to let computers do the work, but of course it remains to
find the algorithms. Will the machine not quickly be more fluent
in this game as well? (Compare the little green men survey by
David Ruelle in Bull. Amer. Math. Soc. ca . 1986, who tabulated on
the imminence of machines cracking theorems with more ease than we
are able to do. Hopefully so, since the goal of any science
(indeed any living being) is to reach immortality.

So if measure theory and general open (=non-compact) Riemann
surfaces inclines much to Lebesgue (and the like), it seems
evident that still much work must be clarified at the more basic
(combinatorial) geometric level of
simpler objects,
e.g. super classical algebraic geometry should be cultivated again
to penetrate more deeply in a variety of problems still unsolved.

\section{Prehistory of Ahlfors}\label{Sec:Prehistory-Ahlfors}

This section attempts a fairly exhaustive
tabulation of works antedating Ahlfors 1950
\cite{Ahlfors_1950}, bearing more-or-less direct connection to
it. In some critical cases, some of those may also be
considered as (vague?) anticipations of the Ahlfors mapping by
other ``pretenders''. In chronological order, we shall discuss
contributions of Riemann 1857--58--76, Schottky 1875--77,
Klein ca. 1876--82--92, Koebe ca. 1907, Bieberbach 1925,
Grunsky 1937--41--49, Courant 1937--39--50, Teichm\"uller
1941.

Our history is
 not intended to be a smoothly readable account inclining to
passive somnolence, but rather
one inviting to further active
searches to clarify several puzzling aspects, where in our opinion
historical continuity is violently lacking. Historical turbulences
arise mostly from several links hard to track down due to poor
cross-referencing (especially in the case of Teichm\"uller 1941
\cite{Teichmueller_1941}, who seems to credit Klein for a sort of
qualitative version of the Ahlfors circle map, yet without bound
upon the degree). In contrast, the first steps, i.e. the
affiliation Riemann-Schottky-Bieberbach-Grunsky is well documented
(but confined to planar surfaces, hence inferior
to Ahlfors' work). Courant's contribution is more in the trend
Dirichlet-Riemann-Plateau-Hilbert, but ultimately a bit
sketchy when it comes to
compare with Ahlfors.

Regarding Koebe, he was quite influenced by Klein's orthosymmetry
(which bears a direct connection to  Ahlfors' conformal circle map
via the algebro-geometric viewpoint), but was more involved with
uniformization (in particular of real curves) and the {\it
Kreisnormierungsprinzip} (rooted back in Schottky, if not
Riemann). Koebe's work concentrates more upon conformal
diffeomorphisms than branched covers. Perhaps an exception
concerns his later works ca. 1910 influenced by Hilbert, where he
comes to investigate more closely non-schlichtartig surfaces.
However in the overall we could not find (as yet) in the
torrential
series of Koebe's papers a clear-cut anticipation of Ahlfors'
result. (Relevant works of Koebe will in fact rather be
surveyed in the next section.)

To summarize
we have located essentially 3 potential forerunners of the Ahlfors
circle map:

(1) Klein, through a citation (or rather allusion) of
Teichm\"uller in 1941 (supplied without precise reference!) and to
which we were not able to supply sound footing (despite long
searches through Klein's collected papers, plus his harder-to-find
G\"ottingen lectures in 1891--92
\cite{Klein_1891--92_Vorlesung-Goettingen},
\cite{Klein_1892_Vorlesung-Goettingen}). In case  no trace is to
be found in Klein's work, it is conceivable that Teichm\"uller
distorted somehow his memory about Klein, in which case
Teichm\"uller should be regarded as the genuine forerunner. It may
be imagined that a  micro-tunnel (=logical wormhole) links Klein
to Ahlfors, and this may have existed in Teichm\"uller's brain
(but as far as I know no proofs are to be found in print).

(2) Courant who makes a vague claim that the result of
Riemann-Schottky-Bieberbach-Grunsky extend to configurations
of higher genus. If Courant's claim is
correct, it would be of extreme interest to
present the details, especially if it is possible to write
down the bound arising from Courant's argument (inspired from
Plateau's problem).

(3) Matildi 1945/48 \cite{Matildi_1945/48} and Andreotti 1950
\cite{Andreotti_1950}.

\subsection{Tracing back the early history (Riemann 1857,
Schwarz, Schottky 1875--77, H. Weber 1876, Bieberbach 1925,
Grunsky 1937--50, Wirtinger
1942)}

From Grunsky's papers (1937 \cite{Grunsky_1937}, 1941
\cite{Grunsky_1941_KA}, both cited in Ahlfors 1950's paper) one
can trace down the early history of Ahlfors theorem back to the
very origin (i.e. Riemann) as follows. Grunsky was Bieberbach's
student. The latter proved a version of this theorem (yet without
the extremal interpretation) for planar (schlichtartig membrane,
i.e. $p=0$) in Bieberbach 1925 \cite{Bieberbach_1925}. In this
paper, one detects an early influence of Schottky's
Dissertation (Berlin 1875, under Weierstrass) published 1877
\cite{Schottky_1877}, as well as a Nachlass of Riemann
estimated of 1857 (which was published in his Werke ca. 1876).
Riemann apparently only handles the case of a {\it
Kreisbereich} (circular domain), yet it seems that Heinrich
Weber---who
edited this Riemann's Nachlass---may have considerably
amputated the original manuscript.
(Of course it would be a first class Leistung if
some specialist of Riemann's work would undertake the
difficult project of
producing a more all-inclusive account.)
Let us reproduce the introduction of Bieberbach 1925
\cite{Bieberbach_1925}:

\begin{quota}[Bieberbach 1925]\label{quote:Bieberbach-1925}
{\small \rm Es handelt sich in dieser Arbeit um die Abbildung
eines mehrfach zusammen\-h\"angenden schlichten Bereiches auf
eine mehrfach bedeckte Kreisscheibe. Insbesondere stelle ich
mir die Aufgabe, zu beweisen, da{\ss} ein $n$-fach
zusammenh\"angender Bereich stets auf eine $n$-bl\"attrige
Kreisscheibe abgebildet werden kann. Die erste im Druck
erschienene Arbeit, die sich mit diesen Fragen besch\"aftigt,
ist die Dissertation von Schottky (Berlin 1875), die im 83.
Bande des Crelleschen Journal abgedruckt ist. Die Frage nach
der kleinstm\"oglichen Bl\"atterzahl ist dort nicht behandelt,
aber die Analogie und die Beziehung zur Theorie der
algebraischen Funktionen und ihrer Integrale liegt den
Betrachtungen zugrunde, und auch die Beziehung zur Theorie der
linearen Differentialgleichungen 2. Ordnung kommt zum
Vorschein. Wie mir Herr Schottky erz\"ahlte, machte bald
darauf H.\,A. Schwarz darauf aufmerksam, da{\ss} sich Riemann
im Sommer 1857 bereits mit der eingangs erw\"ahnten Frage
besch\"aftigte. In der von H. Weber bearbeiteten Darstellung
dieses Teils des Riemannschen Nachlasses findet sich freilich
keine volle Erledigung der Frage. Ich finde, da{\ss} auch
nicht alle Gedanken des Riemannschen Manuskriptes zur
Verwendung kamen. (Vrgl. Riemanns Werke 2.\,Auflage
S.\,440--444) Riemann kn\"upft bei seinen \"Uberlegungen an
die Theorie der linearen Differentialgleichungen an. Die
Theorie der algebraischen Funktionen wird nach der Weberschen
Darstellung zur L\"osung des Abbildungsproblems nicht
herangezogen. Dagegen scheinen mir die Riemannschen Notizen zu
lehren, da{\ss} Riemann auch einen \"uber die Theorie der
Abelschen Integralen f\"uhrenden Weg unabh\"angig von dem bei
Weber dargestellten erwogen hat. Welcher von beiden Wegen der
fr\"uhere ist, vermag ich nicht zu entscheiden.

}
\end{quota}

Hence the tension between Abelian integrals and potential theory
seems to have always been surrounded by a little ring of
mysteriousness, even in the passage  of Bieberbach 1925's article
just quoted. Furthermore, after Grunsky completed in 1941 his
series of papers on the question, it looked desirable to Wirtinger
to publish 1942 \cite{Wirtinger_1942} his own interpretation of
Riemann's Nachlass which he probably knew since ca. 1899 during
his duties as publisher of the second edition of Riemann's Werke.

\begin{quota}[Wirtinger 1942]\label{Wirtinger:quote}
{\small \rm

Die Abhandlung des Hrn. Helmut Grunsky, welche in diesen
Berichten, Jahrgang 1941, Nr. 11, unter dem Titel ``\"Uber die
konforme Abbildung mehrfach zusammenh\"angender Bereiche auf
mehrbl\"attrige Kreise II'' erschienen ist, bringt mir
\"Uberlegungen wieder gegenw\"artig, welche unmittelbar an die
klassische Dissertation von F. Schottky (Berlin 1875)
anschlie{\ss}en, welche noch vor dem Bekanntwerden des
Riemannschen Fragmentes \"uber das Gleichgewicht der
Elektrizit\"at auf Zylindern von kreisf\"ormigem Querschnitt
(1876) erschienen ist. Zusammen mit dem dort entwickelten
Symmetrieprinzip reicht die Theorie der algebraischen
Funktionen vollkommen aus, um zu beweisen, da{\ss} ein von
$p+1$ Randkurven, welche v\"ollig getrennt verlaufen und von
denen keine sich auf einen Punkt reduziert, begrenzter Bereich
sich konform auf die $p+1$fach \"uberdeckte Halbebene der
Variabeln $z=x+iy, y\ge 0$ abbilden l\"a{\ss}t, wobei noch auf
jeder Linie der dem Punkte $z=\infty$ entsprechende beliebig
vorgegeben werden kann.

}

\end{quota}

In the above quote, Bieberbach also mentions that H.\,A. Schwarz
was well acquainted with this Riemann's Nachlass. In this
connection, it can be reminded that the whole trend connected to
the so-called {\it Schwarz lemma} involving Schwarz 1869--70
\cite[p.\,109]{Schwarz_1869-70_Zur-Theorie-der-Abbildung},
Carath\'eodory 1907 \cite{Caratheodory_1907}, 1912
\cite{Caratheodory_1912} (where the coinage ``Schwarz lemma'' is
first used), Pick 1916, Ahlfors 1938, with some intermediate steps
due to E.~Schmidt ca. 1906  (as acknowledged in Carath\'eodory
1907 \cite{Caratheodory_1907}) is well known to have been another
inspiring source for Ahlfors' extremal problem.

To be even more mystical, Carath\'eodory mentions---in his 1936
laudation to Ahlfors' reception of the (chocolate) Fields
medal (ICM 1936)---a certain ``{\it \"Olfleckmethode} of Schwarz,
which seems to be related to all this. This intriguing
terminology, probably refers to the common ``\"Olfleck''
experiment consists of taking any ``oil'' region in a water
recipient while exciting it slightly or even strongly with a thin
instrument, yet preferably without causing a rupture of its
connectedness. Observationally, one can then contemplate with
which determination and structural stability the possibly highly
distorted ``\"Olfleck'' restores to the round circle-shape even if
there are thin necks in the initial position. This seems to be one
of the most beautiful way to visualize the Riemann mapping theorem
in nature. Mathematically this \"Ofleck experiment bears perhaps
more analogy to the normal curvature flow (Huisken, etc.), than
the levels of the Riemann mapping function. One can wonder if
there is an identity between the curvature flow and RMT.

\subsection{Schottky 1875--77}

All sources indicate that Schottky discovered the circle mapping
for multiply-connected domains independently of Riemann's
Nachlass. Compare the next 3 quotes of Schottky
(\ref{quote:Schottky-1882}) and Klein
(\ref{Klein-1923:quote:Schottky}),
(\ref{Klein-1923:quote:Riemann-1858}).
In 1882, during the hot Klein-Poincar\'e
``competition'' on automorphic functions vs. Fuchsian functions,
Schottky's Thesis came again to
the
forefront, with Klein asking its writer for some
precision
about
its genesis. Besides, Schottky rectified some (historically)
inaccurate statement made by Klein.
It resulted a letter published 1882 in Math. Annalen
\cite{Schottky_1882_Brief}, which we reproduce in part:

\begin{quota}[Schottky 1882 {\rm \cite{Schottky_1882_Brief}}]
\label{quote:Schottky-1882}

{\small \rm Dass \"ubrigens Riemann bereits die mit dieser Figur
in Zusammenhang stehenden Functionen und ihre
Differentialgleichungen entdeckt hat, wird durch die Stelle pag.
413--416 seiner gesammelten Werke bewiesen. [$\bigstar$\,Footnote:
Gleichgewicht der Electricit\"at auf Cylindern mit kreisf\"ormigem
Querschnitt und parallelen Axen.---Herr Weber f\"ugt als
Herausgeber diesem Aufsatze die Bemerkung zu: ``Von dieser und den
folgenden Abhandlungen [des Riemann'schen Nachlasses] liegen
ausgef\"uhrte Manuscripte von Riemann nicht vor. Sie sind aus
Bl\"attern zusammengestellt, welche ausser wenigen Andeutungen nur
Formeln enthalten.''\,$\bigstar$] Indess m\"ochte ich betonen,
dass meine Dissertation ein Jahr vor der Publication von Riemann's
Nachlass erschienen ist. Auch erfuhr ich von Letzterem
erst\footnote{Apparently via H.\,A. Schwarz, compare Bieberbach
1925 \cite{Bieberbach_1925}.}, als meine Arbeit bereits in ihrer
zweiten Fassung zum Druck \"ubergehen war. Aber ich bin
gl\"ucklich, mit Ihnen die Priorit\"at der Entdeckung Riemann's
constatiren zu k\"onnen. \dots

Sie haben in freundlicher Weise den Wunsch ge\"aussert,
Genaueres \"uber die Pr\"a\-missen meiner damaligen Arbeit zu
erfahren. Die Anregung zum selbst\"andigen Eindringen in die
Potentialtheorie verdanke ich Herrn Helmholtz. Das in der
Arbeit behandelte Problem, der urspr\"unglichen Auffassung
nach der Potentialtheorie geh\"orig, und wesentliche
Anschauungen meiner Arbeit sind aus
mathematisch-physikalischen Autoren gesch\"opft. Ich nenne
neben den Vorlesungen und Schriften von Herrn Helm\-holtz
insbesondere ein mir g\"utig von Herrn O.\,E. Meyer geliehenes
Heft noch nicht publicirter Vorlesungen von Herrn F. Neumann,
dann ferner ein Buch \"uber Elektrostatik von Herrn
K\"otteritzsch, etc. Die Durchf\"uhrung der so gewonnenen
Ideen wurde mir sodann wesentlich erleichtert durch Herrn
Weierstrass' Vorlesungen \"uber Abel'sche Functionen, sowie
besonders durch die von Herrn Schwarz publicirten
Untersuchungen \"uber das Abbildungsproblem einfach
zusammenh\"angender Fl\"achen. Mit R\"ucksicht auf die
letzteren wurde auf den Rath meines hochverehrten Lehrers,
Herrn Weierstrass, der urspr\"unglich \"uberreichte Entwurf
der Arbeit so abgeh\"andert, dass sich dieselbe in beiden
ver\"offentlichten Fassungen an die Untersuchungen von Herrn
Schwarz anschliesst. \dots
\hfill \hbox{Breslau, im Mai 1882.}

}

\end{quota}

This work of Schottky enjoyed
early and great recognition
among colleagues, and still today is frequently cited. The
reasons of this success are multiple, but I cannot resist to
quote first
 Le Vavasseur 1902 \cite{Le-Vavasseur_1902} [since
in Geneva there is a prominent artist bearing a similar name],
himself quoting Picard:

\begin{quota}[Le Vavasseur 1902]

{\small \rm  Dans le Tome II de son {\it Trait\'e d'Analyse}, page
285, M. \'Emile Picard \'ecrit: ``Deux aires planes $A$ et $A_1$,
limit\'ees chacune par un m\^eme nombre de contours, ne peuvent
pas, en g\'en\'eral, \^etre repr\'esent\'ees d'une mani\`ere
conforme l'une sur l'autre. L'\'etude approfondie de ce probl\`eme
a \'et\'e faite par M.~Schottky dans un beau et important
M\'emoire.''

Plus loin m\^eme Tome, page 497, en note, M. \'Emile Picard
\'ecrit encore: ``Nous avons d\'ej\`a eu l'occasion de citer
le beau travail de M. Schottky; c'est un M\'emoire fondamental
\`a plus d'un titre.''

}

\end{quota}

The enthusiasm for Schottky's work diffused from France to Italy,
cf. especially Cecioni 1908~\cite{Cecioni_1908}, who may be
credited for the first rigorous proof of the parallel slit map.
The reason of Schottky's popularity is the quite amazing novelty
of his work in  prolongation of Riemann' ideas---but so in
retrospect only for Schottky was not directly influenced by
Riemann. The methods range from potential-theoretic to algebraic
functions, flourishing into an
breathtaking variety of results. Beside the circle map for
multiply-connected domains, it contains both what later will be
known as the Kreisnormierungsprinzip (KNP), plus the parallel-slit
mappings (PSM).
---{\it Warning.} In fact I am not sure that it contains KNP, but
could easily have on the basis of  a naive parameter count. Also
it is never clear if material was amputated from the first 1875
edition of Schottky's Thesis. According to Klein's quote
(\ref{Klein-1923:quote:Riemann-1858}), it seems however that the
first Latin edition (1875) of Schottky's Thesis
contains the statement of the Kreisnormierung, yet ``{\it nur auf
Grund einer Konstantenz\"ahlung\/}''.---
%
%
At any rate, it contains (explicitly or in embryo) virtually all
of the varied canonical conformal maps which will be re-studied by
Koebe during the period 1904--1930, trying even to extend the
results to infinite connectivity. As is notorious, this ramifies
to deep waters still not completely elucidated today, cf.
He-Schramm 1993 \cite{He-Schramm_1993}, which is still
the best result reached so far on the Kreisnormierung problem.
Schottky's Thesis also contains the idea of symmetric reproduction
of such a domain, where Klein identifies one of the first
instance of automorphic functions. The name ``Schottky
uniformization'' is still of widespread usage today (e.g. Bers,
Maskit, etc.). The influence of Schottky's work is also apparent
in the jargon ``Schottky differentials'' widely used in several of
Ahlfors' papers, especially Ahlfors 1950 \cite{Ahlfors_1950}.
(From the algebro-geometric viewpoint this probably just amounts
to a real differential.) Last but not least, the Schwarz principle
of symmetry (1869
\cite{Schwarz_1869-Ueber-einige-Abbildungsaufgaben}) [which
afterwards Klein liked to identify in Riemann's Nachla{\ss}
\cite{Riemann_1857_Nachlass} already, as testimonies the many
brackets added in his collected papers, e.g. Klein 1923
\cite[p.\,631, line 3]{Klein-Werke-III_1923}] enables one to form
the so-called {\it Schottky double}\/. All this appears first in
this single work of Schottky.

The admiration for Schottky's Thesis propagated long through the
ages, e.g.:

\begin{quota}[Garabedian-Schiffer 1949
{\rm \cite[p.\,187, p.\,214]{Garabedian-Schiffer_1949}}]

{\small \rm

An understanding of all identities between domain functions may be
obtained by sustained application of Schottky's theory of
multiply-connected domains [15](=1877). Schottky proved that there
is a close relation between the mapping theory of these domains
and the theory of closed Riemann surfaces; the identities among
domain functions have their complete analogue in the theory of
Abelian integrals and might be proved by means of the latter.
[\dots, and on p.\,214]

{\bf Schottky functions and related classes.} Schottky
[15](=1877) was the first to consider the family $\frak R$ of
all functions which are single-valued and meromorphic in $D$
[a multiply-connected domain] and have real boundary values on
$C$ [the full contour of $D$]. He developed an interesting
theory of conformal mapping of multiply-connected domains from
the properties of this family and established by means of it
the relation of this theory with the theory of closed Riemann
surfaces. It is evident that functions $f(z)\in \frak R$ are
very useful in the method of contour integration.

}

\end{quota}

Schottky's Thesis originated in the ambiguous context of physical
intuition vs. Weierstra{\ss}ian rigor. It is notorious that the
ultimate redaction was a hard gestation process subjected to
incessant revisions demanded by Weierstra{\ss}. As we know (from
Schottky himself (\ref{quote:Schottky-1882}), plus the next two
quotes by Klein) the first impulse was physically motivated
(Helmholtz, F. Neumann, the father of C. Neumann, etc.), and then
only lectures of Weierstra{\ss} and papers of Schwarz came to
influence the mathematical treatment.
%
For Klein this excessive Weierstrassization is regarded from a
sceptical angle (cf. again the next two quotes).

It is
a delicate question to wonder about the rigor reached in Schottky,
despite its ultimate foundation over technology of Schwarz as a
substitute to the Dirichlet principle. To ponder its ultimate
rigor, it suffices to say that all  of Schottky's results where
subsequently revisited, by the following workers:

$\bullet$ Koebe  for KNP in a
(torrential) series of paper spread from 1906 to 1922.

$\bullet$ Cecioni 1908 \cite{Cecioni_1908} for the PSM
(=parallel-slit mapping);  the latter even mark (discretely) the
superiority of his proof by emphasizing that Schottky's argument
relies on a parameter count,
whereas he proposes to prove PSM ``{\it direttamente\/''} (cf.
\loccit\,p.\,1). This technical ``gap'' was of course known to
Klein, cf. the next Quote~\ref{Klein-1923:quote:Riemann-1858}.
Also in Salvemini 1930 \cite[p.\,3]{Salvemini_1930} (a student of
Cecioni) the critique is made more explicit: ``{\it Questo
risultato [=PSM] era stato enunciato dallo Schottky in base ad un
computo di parametri, computo che non \`e poi esauriente.\/}''

$\bullet$ Bieberbach 1925 \cite{Bieberbach_1925} for the circle
mapping problem.

None of those writers attacks frontally the standards of rigor in
Schottky (as based upon the  complicated but solid foundations
laid by Schwarz). Still, the technical complications was seen as a
need to find simpler derivations of the geometrical results. After
sufficiently time
elapsed, the subsequent generation tends to ascribe the
(rigorous) proof of Schottky's result to this second
wave of
workers. E.g., Grunsky 1978 \cite{Grunsky_1978} ascribes
Schottky's circle maps to Bieberbach 1925 \cite{Bieberbach_1925}
(cf. Quote~\ref{quote:Grunsky-1978}), and Bieberbach
1968 \cite{Bieberbach_1968-Das-Werk-Paul-Koebes} credits Koebe for
the proof of KNP (in finite connectivity).
All these redistributions are done without specific objections
upon
the original arguments Schottky's. This is a usual loose process
relegating methodologies just due to their cumbersomeness, as a
sufficient reason for lack of rigor.


In contrast, even more contemporary workers still credits Schottky
for the first proof of the Kreisnormierung result (cf. e.g.,
Schiffer-Hawley 1962 \cite[p.\,183]{Schiffer-Hawley_1962}). So it
is a subtle socio-cultural game to pinpoint  precisely about which
writer furnished the first acceptable proof.

\subsection{Klein's comments about Riemann-Schottky}

In the third volume of his collected papers Klein makes
several
comments about Riemann
and
Schottky Thesis. He insists first on the physical motivations of
Schottky, which were progressively
``censured'' under Weierstrass' influence.

\begin{quota}[Klein 1923
{\rm \cite[p.\,573]{Klein-Werke-III_1923}}]
\label{Klein-1923:quote:Schottky}

{\small \rm Ich greife gern noch einmal auf die wie\-der\-holt
genannte Arbeit Schottkys in Crelle Journal, Bd. 83 (1877)
zur\"uck, zumal ich weiter unten (S. 578/579) ohnehin
ausf\"urlicher auf sie zur\"uckkommen mu{\ss}. Die gro{\ss}e
\"Ahnlich\-keit der auf einen besonderen Fall bez\"uglichen
Schottkyschen Untersuchungen mit den allgemeinen meiner
Schrift war mir von vornherein aufgefallen. Ich schrieb also
damals an Herrn Schottky und fragte ihn nach der Enstehung
seiner Ideen. Hierauf antwortete her mir in einem Briefe von
Mai 1882 (der in Bd. 20 der Math. Annalen abgedruckt wurde),
da{\ss} er in der Tat urspr\"unglich auch von der Bertrachtung
der Str\"omungen einer inkompressiblen Fl\"ussigkeit
ausgegangen sei und diesen physikalischen Ausgangs\-punkt nur
auf Rat von Weierstrass bei der Drucklegung durch die
Bezugnahme auf Schwarz' Untersuchungen \"uber konforme
Abbildung ersetz habe.

}
\end{quota}

Then Klein recollects some more details in the following passage.
This contains an anecdotic conflict (with Bieberbach 1925
\cite{Bieberbach_1925}=Quote~\ref{quote:Bieberbach-1925}) about
the estimated date of Riemann's Nachlass. More interestingly,
Klein expresses the view that Schottky's theorem (to the effect
that a multiply-connected domain is conformal to a circle domain)
may be seen as the planar case of Klein's {\it
R\"uckkehrschnitttheorem}, which in turn seems to be one of the
weapon that Klein used in his
early
strategy toward uniformization (an approach not successfully
completed until Brouwer-Koebe ca. 1911
\cite{Klein-Brouwer-Koebe_1912}).

\begin{quota}[Klein 1923
{\rm \cite[p.\,578--579]{Klein-Werke-III_1923}}]
\label{Klein-1923:quote:Riemann-1858}

{\small \rm

\"Ubrigens hat Riemann ja auch die andere Art automorpher
Funktionen, die enstehen, indem man an einen von Vollkreisen
begrenzten Bereich der Ebene an diesen Kreisen fortgesetzt
symmetrisch reproduziert (Siehe das von H. Weber bearbeitete
Fragment XXV in der ersten (1876 erschienenen) bzw. XXVI in
der zweiten (1892 erschienenen) Auflage der Ges. math. Werke
von Riemann.) Die Pr\"ufung der Originalbl\"atter hat ergeben,
da{\ss} Webers Mitteilungen den Vorbereitungen zu einer im
Sommer 1858 gehaltenen Vorlesung entnommen sind. Und zwar geht
 Riemann dabei zun\"achst von der Aufgabe aus, f\"ur ein von
 mehreren Kugeln gebildetes Konduktorsystem das Gleichgewicht
 elektrostatischer Ladungen zu bestimmen. Hierf\"ur war die
 Benutzung des Symmetrieprinzipes in den Arbeiten von W.
 Thompson vorgebildet, die als Briefe an Liouville in dessen
 Journal von 1845 an erschienen. Also auch hier sind die
 mathematischen Entwicklungen aus physikalischen Anregungen
 erwachsen.

Auf dieselben Funktionen ist dann unabh\"angig in seiner
Berliner Dissertation 1875 Herr Schottky gekommen. Von seinem
physikalischen Ausgangspunkte ist schon oben auf S.\,573, die
Rede gewesen. Im \"ubrigen sind die Schicksale der
Schottkyschen Arbeit, wie sie sich nach pers\"onlicher
Mitteilung des Verfassers ergeben, so merkw\"urdig, da{\ss}
ich gern die Gelegenheit ergreife, sie hier mitzuteilen. Es
erfolgten nach einander drei verschiedene Redaktionen:

a) Eine lateinische Fassung, die nicht publiziert ist, sondern
nur der Philosophischen Fakult\"at in Berlin vorgelegen hat,

b) Eine deutsche Bearbeitung, welche 1875 in Berlin als
Dissertation gedruckt wurde,

c) Die umgearbeitete Darstellung in Crelles Journal, Bd. 83
(1877).

\noindent Bei Niederschrift von a) hat der Verfasser noch
keine F\"uhlung mit Weierstrass gehabt, daf\"ur aber ganz
seiner freien Ideenbildung folgen k\"onnen. Aus dem Gutachten,
da{\ss} Weierstrass \"uber a) seinerzeit f\"ur die Fakult\"at
abgegeben hat und von dem ich durch die Freundlichkeit von
Herrn Schottky eine Abschrift vor Augen habe, scheint mit
Gewi{\ss}heit hervorzugehen, da{\ss} Schottky hier, freilich
nur auf Grund einer Konstantenz\"ahlung, das
``R\"uckkehrschnitttheorem'' f\"ur den besonderen, von ihm
betrachteten Fall ausgeschprochen hat, d. h. die
M\"oglichkeit, einen von $p+1$ regul\"aren Randkurven
begrenzten eben Bereich auf einen von $p+1$ Vollkreisen
begrenzten Bereich konform abzubilden (also das
R\"uckkehrschnitttheorem f\"ur den obersten orthosymmetrischen
Fall, wie ich mich ausdr\"ucke).

Die Redaktion b) ist dann durch eine erste F\"uhlungnahme mit
Weierstrass bedingt. Bei der umfassenden Beherrschung
ausgedehnter Teile der Mathematik und seiner stark
ausgepr\"agten Pers\"onlichkeit, die sich zu bestimmten
Beweisg\"angen durchgearbeitet hatte, \"ubte Weierstrass auf
j\"ungere Forscher je nachdem einen au{\ss}erordentlich
f\"ordernden, oder auch, wo ihm die Gedankeng\"ange fremdartig
waren, einen hemmenden Einflu{\ss}. [\dots]. Schottky scheint
\"ahnliche Erfahrungen gemacht zu haben, so da{\ss} er in b)
sich blo{\ss} auf die Konstantenz\"ahlung beschr\"ankt, ohne
ihre Tragweite f\"ur das Fundamentaltheorem anzudeuten
[\dots]. Die physikalische Ideenbildung aber, von der doch der
Autor ausgegangen war, wird g\"anzlich ausgeschaltet und durch
Zitate auf die das Existenzproblem der konformen Abbildungen
betreffenden Arbeiten von Schwarz ersetzt.

In c) endlich ist auch noch besagte Konstantenz\"ahlung
weggeblieben. [[[Fu{\ss}note: Dagegen hat Schottky in c)
(S.\,330 daselbst), wiederum auf Grund blo{\ss}er
Konstantenz\"ahlung, den Satz ausgeschprochen, da{\ss} sich
jedes ebene, von $p+1$ Randkurven begrenzte, Gebiet umkehrbar
eindeutig konform auf die Vollebene mit Ausnahme von $p+1$
geradlinigen, zur $x$-Achse parallelen Strecken abbilden
l\"a{\ss}t. Bereiche der letzteren Art spielen in der modernen
Literatur unter dem Namen {\it Schlitzbereiche} bekanntlich
eine wichtige Rolle.]]] Statt dessen finden sich wertvolle,
vorher nicht publizierte, Angaben \"uber die verschiedenen
Normalformen, die Weierstrass bei den Gebilden $p>2$
unterschied; [\dots]

}
\end{quota}

Incidentally this {\it R\"uckkehrschnitttheorem}, may have
some connection with the Ahlfors function albeit probably no
direct link is evident, there is still some striking analogy
developed in the next section.

\subsection{A historical puzzle: why Klein missed the Ahlfors
circle mapping?}

[27.04.12] After reading quite closely the above comments of
Klein, plus having a vague idea of the content of Schottky's
Dissertation one is
puzzled
by how close Klein
might have been (ca. 1882) to anticipate by circa 70 years the
circle map of Ahlfors (1948--1950). Here is our reasoning.

First, Schottky's Thesis (and in cryptical form already Riemann's
Nachlass) contains two striking results:

$\bullet$ the {\it circle map} (CM) of a (compact)
multiply-connected domain to the disc, and beside

$\bullet$ what later came to be known (in Koebe's era, cf. e.g.,
Koebe 1922 \cite{Koebe_1922}) as the {\it Kreisnormierungsprinzip}
(KNP) to the effect that any such domain is conformally equivalent
to a circular domain. (Recall from Klein's quote
(\ref{Klein-1923:quote:Riemann-1858}) that this occurs only in the
original Latin version of Schottky's Thesis.)

Both results are natural extensions of RMT (=Riemann mapping
theorem) either by allowing branched coverings or just by using
faithful conformal diffeomorphisms (but then of course the target
depends upon moduli).

Now loosely speaking one may consider both results (CM and KNP) as
lying at the same order of difficulty (at least both are to be
found in Schottky's Thesis).

Next, Klein points out (cf. right above
Quote~\ref{Klein-1923:quote:Riemann-1858}) that he was able in
1881--82 to prove an  extension of (KNP) to positive genus
$p>0$, which he calls (apparently with Fricke's
assistance---cf. Klein 1923 \cite[p.\,623,
footnote~4]{Klein-Werke-III_1923}) the {\it
R\"uckkehrschnitttheorem} (RST). Klein was very proud of this
result (cf. especially Klein~1923
\cite[p.\,584]{Klein-Werke-III_1923}, where this discovery is
dated from September 1881 (Borkum)), comparing it (as a
psychological experience) to Poincar\'e's discovery of his
general {\it fonctions fuchsiennes}.

\begin{figure}[h]
\centering
    \epsfig{figure=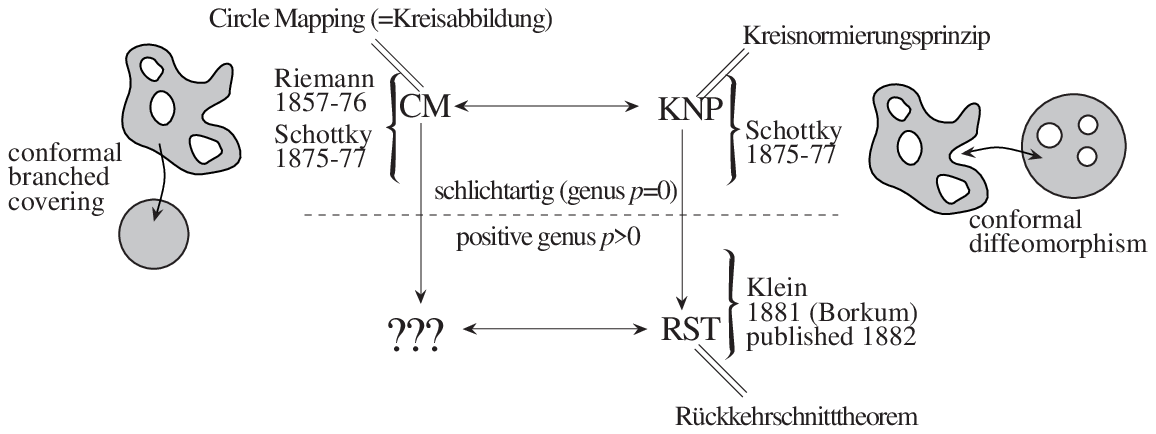,width=126mm}
    \vskip-5pt\penalty0

\caption{\label{KNP-RST:fig} How Klein could miss the Ahlfors
map?}
\end{figure}

Thus, there is an obvious commutative diagram
(Fig.\,\ref{KNP-RST:fig}), and whatsoever the actual meaning of
Klein's (RST) should be, there is only a single natural candidate
to fill in the diagram at the (triple) question-marks, namely the
Ahlfors circle map. This accentuates once more why Klein may have
been a serious candidate to anticipate the Ahlfors circle map, at
least without extremal interpretation.

Furthermore, in view of say Gabard 2006 \cite{Gabard_2006}, the
Ahlfors mapping amounts essentially to Jacobi's inversion problem
in the real case, and here again this was one of Klein's major
preoccupation (cf. e.g., Weichold's Thesis 1883
\cite{Weichold_1883}, Hurwitz's work 1883 \cite{Hurwitz_1883},
plus other sources, e.g. Klein 1892 \cite{Klein_1892_Realitaet}).

Of course, it would be an excellent
project to try getting acquainted with Klein's techniques so as to
inspect if they
lead to another elementary existence-proof of
Ahlfors maps. Again it should be recalled that even if Klein
himself was never able to complete his programme some helping hand
from Brouwer-Koebe ultimately vindicated all of Klein's
intuitions.

\subsection{R\"uckkehrschnittteorem (Klein 1881--82)}\label{sec:Ruckkehrschnittthm}

Klein found this theorem in 1881, and published it 1882 in
\cite{Klein_1882_Ruckkehrschnitt}. From the start the paper
confesses to use some irregular methods.

What does Klein in this paper? First he takes a closed Riemann
surface of genus $p>1$ (w.l.o.g) and traces on it $p$ disjoint
R\"uckkehrschnitten (retrosections) and asserts that the cutted
Riemann surface may be mapped to a $2p$-ply connected domain on
the sphere, whose corresponding boundaries $A_i'$, $A_i''$ are
related by  a linear substitution. He then uses these $p$
substitutions
to reproduce the conformal mapping ad infinitum (a trick already
present in Riemann's Nachlass 1857 \cite{Riemann_1857_Nachlass}).
He notes that the construction depends on the right number of free
constants $3p$ compatible with Riemann's moduli $3p-3$, thus
yielding a sound evidence for some sort of uniformization. Of
course it is not yet the standard uniformization as the reproduced
domain filling more and more the sphere still avoids an infinite
set (a Cantor set). In fact this construction gives an unramified
infinite cover of the given closed Riemann surface by a subregion
of the sphere (which is however not simply-connected).

Then he applies a similar method to the case of symmetric Riemann
surfaces by using a symmetric system of retrosections while
showing that the above construction may be done equivariantly. For
instance in the simpler to visualize dividing case the above
dissection process leads to a similar symmetric domain, symmetric
with respect to the orthogonal symmetry of the sphere (whence the
name {\it orthosymmetric}). In the nondividing case the structural
symmetry is rather the diametral one (antipodal map), whence the
name {\it diasymmetric}.

Logically it seems that Klein's method depends on Schottky's
inasmuch as first doing the retrosections one is reduced to
the schlichtartig case which turns out to be schlicht. (This
result was extended by Koebe in 1908 \cite{Koebe_1908_UbaK3}
%
to schlichtartig surfaces of infinite connectivity:
schlichtartig implies schlicht!)

Clearly something remains to be understood on this RST, and our
guess that it is sufficiently strong to imply Ahlfors theorem is
quite disputable.
At any rate Klein seems to have had a
clear-cut conception of how his dichotomy ortho- vs. diasymmetric
is reflected into the Riemann sphere  with its two real structures
(equatorial symmetry vs. antipody). However the issue that
dividing curves are precisely those mapping to the equatorial
sphere in a
totally real fashion may have escaped his attention and does not
seem to be logically reducible to his RST. Yet since RST is
supposed to be the positive genus case of KNP (cf. Klein's quote
(\ref{Klein-1923:quote:Riemann-1858})) it may be expected that one
first establishes KNP and from here one deduces a circle map, much
like Riemann was able to do for the zero genus case in his
Nachlass \cite{Riemann_1857_Nachlass}. This suggests yet another
strategy to approach Ahlfors theorem.

[04.11.12] A more naive idea could be to start from a bordered
surface of type $(r,p)$, and make $p$ retrosections to get it
planar (but with $r+2p$ contours). Then there is on the dissected
surface a circle map of degree $r+2p$. Of course the map is a
priori not assuming the same values on both ridges of the
retrosections and even if we can arrange this, we would like those
points to get mapped in the interior and not the boundary of the
circle.

\subsection{Grunsky's bibliographical notes (Grunsky 1978)}

Let us  now reproduce Grunsky's historical comments (in his
monumental book 1978 \cite[p.\,198]{Grunsky_1978}) about circle
maps. (Brackets are ours additions. We added author's names in
front of the bracket-references to improve readability, plus page
numbers, and finally inserted the symbol $\bigstar$ when
disagreeing with Grunsky's cross-references.)

\begin{quota}[Grunsky 1978]\label{quote:Grunsky-1978}
{\small \rm Theorem 4.1.1. goes back to Riemann 1857/58/76
\cite{Riemann_1857_Nachlass}, who gave some hints for the proof if
[the domain] $D$ is bounded by circles. The first proof is due to
Bieberbach 1925 \cite{Bieberbach_1925}, who used the
Schottky-double and deep results in the theory of algebraic
functions. Elementary proofs were given by Grunsky 1937
\cite{Grunsky_1937}, 1941 \cite{Grunsky_1941_KA}; for 4.1.3. [a
sort of auxiliary lemma in linear algebra] see Furtw\"angler 1936,
Bourgin 1939. Related proofs in Akira Mori 1951 \cite{Mori_1951},
Komatu 1953 \cite{Komatu_1953} (containing generalizations), Tsuji
1956 \cite{Tsuji_1956}; cf. Golusin 1952/57
\cite{Golusin_1952/57}, Tsuji 1959
\cite{Tsuji_1959-BOOK/Chelsea1975}. A proof based on the method
for Plateau's problem: Courant 1937 \cite{Courant_1937}
[$\bigstar$ in Gabard's opinion this paper does not reprove the
circle mapping, but rather the mapping to a Kreisbereich, due to
Schottky--Koebe, cf. p.\,709 and p.\,717, of {\it loc.\,cit.}],
generalized in Courant 1939 \cite{Courant_1939}; cf. Courant 1950
\cite{Courant_1950} [especially p.\,183--187]. Another proof,
using, like Bieberbach 1925 \cite{Bieberbach_1925}, the Schottky
double in Wirtinger 1942 \cite{Wirtinger_1942}; cf. also
Rodin-Sario 1968 \cite{Rodin-Sario_1968} [where ???]. Triply
connected domains: Limaye 1973 \cite{Limaye_1973}. Representation
of the mapping function (Ahlfors function, see 4.3.) by an
orthonormal system in Meschkowski 1952 \cite{Meschkowski_1952}, by
the Bergman kernel in Nehari 1950 \cite{Nehari_1950}. Proofs using
extremal properties in papers quoted in 4.3. and 6. [More about
this below (Quote~\ref{quote:Grunsky-1978-B}).] An extension to
certain domains of infinite connectivity in R\"oding 1975
\cite{Roeding_1975}. A more general type of image domain for
doubly connected domains in Bieberbach 1957
\cite{Bieberbach_1957}.

Some generalizations, based on ideas used in the
aforementioned papers, mainly concerning Riemann surfaces in
Nehari 1950 \cite{Nehari_1950}, Tietz 1955 \cite{Tietz_1955}
(cf. K\"oditz-Timmann \cite{Koeditz-Timmann_1975}), Mizumoto
1960 \cite{Mizumoto_1960}, Timmann 1969 (Diss., Hannover)
\cite{Timmann_1969}, R\"oding 1972 (Diss., W\"urzburg)
\cite{Roeding_1972}, R\"oding 1977 \cite{Roeding_1977_mero}.
Cf. Ahlfors-Sario 1960 \cite{Ahlfors-Sario_1960} [$\bigstar$
where?], Carath\'eodory 1950
\cite{Caratheodory_1950_Buch_Funktionentheorie} [$\bigstar$
where?], Sario-Oikawa 1969 \cite{Sario-Oikawa_1969}
[$\bigstar$ where?].

}

\end{quota}

{\bf Comments (Gabard, Mai 2012):} Alas, regarding the three last
books no pagination is supplied by Grunsky, and as far as I
browsed through them, I failed to locate any place where  Ahlfors'
circle mapping is established anew.

Now we reproduce Grunsky 1978 \cite[p.\,199]{Grunsky_1978}:


\begin{quota}[Grunsky 1978]\label{quote:Grunsky-1978-B}
{\small \rm Theorem 4.3.1., a generalization of Schwarz' lemma
to multiply connected domains, is a special case of a more
general theorem (individual bounds on each boundary component,
prescribed zeros) proved by Grunsky in 1942
\cite{Grunsky_1942} (save for uniqueness, see Grunsky 1950
\cite{Grunsky_1950}). Cf. Herv\'e 1951 \cite{Herve_1951}.
Another proof of 4.3.1. was given by Ahlfors in 1947
\cite{Ahlfors_1947}, completed in Ahlfors 1950
\cite{Ahlfors_1950} (cf. Golusin 1952/57
\cite{Golusin_1952/57}) and the extremal function is called
the ``Ahlfors function'', a term frequently used in the
broader sense of any function mapping [the domain] $D$ [in a]
$(1,n)$ onto $U$ [the unit disc]; the result was carried on to
characterization of the additional zeros of the extremal
function. The method used by Ahlfors, Euler-Lagrange
multipliers (also pointed out in Grunsky 1946
\cite{Grunsky_1946} and applied in Grunsky 1940
\cite{Grunsky_1940}) is likewise a basis for our \S 6. -- For
further proofs of our theorem see Nehari 1951
\cite{Nehari_1951} and Nehari 1952 \cite[pp.\,378
ff.]{Nehari_1952-BOOK}, and some of the papers quoted for
theorem 4.6.4. -- Ahlfors function in a ring domain Kubo 1952
\cite{Kubo_1952}. -- Applications of the Ahlfors function in
Alenicyn 1956 \cite{Alenicyn_1959}, 1961 \cite{Alenicyn_1961},
(cf. Mitjuk 1965 \cite{Mitjuk_1965}).

}
\end{quota}

\subsection{Italian school: Cecioni 1908,
Stella Li Chiavi 1932, Matildi 1945/48, Andreotti
1950}\label{sec:Italian-school}

Of course in the overall Grunsky's comments and references are
essentially sharp (especially a deep knowledge of
Russian/Ukrainian works). Maybe only some contribution of the
Italian school are ignored. (Those  are however meticulously
listed in the book Ahlfors-Sario 1960 \cite{Ahlfors-Sario_1960}.)

For instance the simple continuity argument in the Harnack-maximal
case based upon Riemann-Roch (without Roch) gives a simple proof
in this case (compare e.g., Huisman 2001 \cite{Huisman_2001} or
Gabard 2006 \cite{Gabard_2006}). This simple argument goes back to
Enriques-Chisini seminal book 1915/18
\cite{Enriques-Chisini_1915-1918}, and may have been implicit in
Riemann's original manuscript (not published), compare
Bieberbach's quote (\ref{quote:Bieberbach-1925}).

Further, closely allied work is to be found in works of Cecioni
1908 \cite{Cecioni_1908}, and his students: Salvemini 1930
\cite{Salvemini_1930},  Stella Li Chiavi 1932
\cite{Stella-li-Chiavi_1932}, etc.

Those works certainly deserve  closer studying, but they do not
seem to anticipate Ahlfors circle map. One notable exception is
the article Matildi 1945/48 \cite{Matildi_1945/48} (discovered by
the writer as late as [13.07.12]), where existence of
circle maps for
surfaces
bounded by a single contour seems to be established. (This Italian
work was known to Ahlfors (or Sario?) at least as late as 1960,
again  being quoted in Ahlfors-Sario 1960
\cite{Ahlfors-Sario_1960}.) Of course it would be interesting to
see if Matildi's method adapts to more
contours, while trying to make (more) explicit the degree bound
obtained by him.  Andreotti 1950 \cite{Andreotti_1950} seems to go
precisely in this sense by including several contours (alas the
writer's Italian declined fast enough to have failed understanding
properly  Andreotti's achievements).

\section{Is there any precursor to Ahlfors 1950?}\label{Sec:Precusors}

\subsection{What about Teichm\"uller 1941?}\label{sec:Teichmueller}

One can wonder about the content of Teichm\"uller's Werke.
Does it
overlap with the Ahlfors function?
While reading the long memoir of Teich\-m\"uller 1939
\cite{Teichmueller_1939}
it
transpires to anybody familiar with Klein's work how strong
the
latter's influence is;
in particular Teich\-m\"uller gives a  thorough account of the
(now) so-called {\it Klein surfaces} (and their moduli). Of course
such results were anticipated by Klein (at least at the heuristic
level). Hence, it seems quite natural to wonder if Teichm\"uller
anticipated
the existence
of Ahlfors function
(for orientable membranes). Here is a
report of those portions of Teichm\"uller's
works which
looks closest to this goal, but it should still be debated how
much of the Ahlfors circle maps was anticipated by Teichm\"uller.

The most relevant passage in Teichm\"uller's writings seems to be
the following extract of Teichm\"uller 1941
\cite{Teichmueller_1941} (reedited in
\cite[p.\,554--5]{Teichmueller_1982}):

\begin{quota}[Teichm\"uller 1941]\label{quote:Teichmueller-1941}

{\small \rm

Wir besch\"aftigen uns nur mit {\bf orientierten end\-lichen
Riemannschen Mannigfaltigkeiten.} Diese k\"onnen als Gebiete
auf geschlos\-senen orientierten Riemannschen Fl\"achen
erkl\"art werden, die von endlich vielen gesch\-los\-senen,
st\"uckweise analytischen Kurven begrenzt werden. Sie sind
entweder geschlos\-sen, also selbst geschlossene orientierte
Riemannsche Fl\"achen, die man sich endlichviel\-bl\"attrig
\"uber eine $z$-Kugel ausgebreitet vorstellen darf, oder
berandet. Im letzteren Falle, kann man sie nach Klein durch
konforme Abbildung auf folgende Normalform bringen: ein
endlichvielbl\"attriges Fl\"achenst\"uck \"uber der oberen
$z$-Halbebene mit endlich vielen Windungspunkten, das durch
Spiegelung an der reellen Achse eine symmetrische geschlossene
Riemannsche Fl\"ache ergibt; [\dots]

(So l\"a{\ss}t sich z.\,B. jedes Ringgebiet, d.\,h. jede
schlichtartige endliche Riemannsche Mannigfaltigkeit mit zwei
Randkurven, konform auf eine zweibl\"attrige \"Uberlagerung
der oberen Halbebene mit zwei Verzweigungspunkte abbilden.)

}

\end{quota}

Unfortunately, no precise cross-reference  to Klein is
given
and one needs to browse
Klein's works (with the option of some G\"ottingen Lectures
Note 1891/92 \cite{Klein_1891--92_Vorlesung-Goettingen},
\cite{Klein_1892_Vorlesung-Goettingen} not reproduced in
Klein's collected papers).
This absence of precise location is quite annoying. A charitable
excuse is the  World War II context in which the paper was
written: ``{\it Weil mir nur eine beschr\"ankte Urlaubzeit zur
Verf\"ugung steht, kann ich vieles nicht begr\"unden, sondern nur
behaupten.\/}'' (compare {\it loc.\,cit.}
\cite[p.\,554]{Teichmueller_1982} 2nd parag.)

\subsection{Detective work: Browsing Klein through the
claim of Teichm\"uller}

Regarding Teichm\"uller's cryptical allusion to Klein (as
discussed
 in the previous section) we have the following
candidates in Klein's works (none of which at the present stage of
our historical search truly corroborates Teichm\"uller's
crediting):

(1) Klein 1882
\cite[p.\,75]{Klein_1882}=\cite[p.\,567]{Klein-Werke-III_1923}
where one reads:

\begin{quota}[Klein 1882] {\small \it

Man hat also eine komplexe Funktion des Ortes, welche in
symmetrisch gelegenen Punkten geiche reelle, aber
entgegengesestzt gleich imagin\"are Werte aufweist.

}

\end{quota}

This looks quite close to the desired assignment, yet in
reality only corresponds to the existence of a real morphism
on any real curve; equivalently the existence for any (closed)
symmetric Riemann surface of an equivariant holomorphic map to
the sphere acted upon by the (usual) complex conjugation
fixing an equator. Hence, in our opinion, this passage of
Klein cannot be regarded as a genuine forerunner of the
Ahlfors circle mapping.


(2) Another place where Klein comes quite close to Teichm\"uller's
assertion occurs in the same 1882 booklet ``{\it \"Uber Riemanns
Theorie \dots}'', where Klein computes the
moduli of real algebraic curves---equivalently symmetric
Riemann surfaces (cf.
\cite[p.\,568--9]{Klein-Werke-III_1923}):

\begin{quota}[Klein 1882] \label{quote:Klein-1882-moduli}

{\small \rm

Indem wir uns jetzt zu den {\it symmetrischen} Fl\"achen
wenden, haben wir noch eine kleine Zwischenbetrachtung zu
machen. Zun\"achst ist ersichtlich, da{\ss} zwei solche
Fl\"achen nur dann ``symmetrisch'' aufeinander bezogen werden
k\"onnen, wenn sie neben dem gleichen $p$ dieselbe Zahl $\pi$
der \"Ubergangskurven [=real ``ovals''] darbieten und
\"uberdies beide entweder der ersten Art oder der zweiten Art
angeh\"oren. [This is the dichotomy ortho- vs diasymmetric.]
Im \"ubrigen wiederhole man speziell f\"ur die symmetrischen
Fl\"achen die Abz\"ahlungen des \S 13 betreffs der Zahl der in
eindeutigen Funktionen enthaltenen Konstanten unter der
Bedigung, da{\ss} nur solche Funktionen in Betracht gezogen
werden, welche an symmetrischen Stellen konjugiert imagin\"are
Werte aufweisen. Hiermit kombiniere man sodann nach dem Muster
des \S 19 die Zahl solcher \"uber der $Z$-Ebene konstuierbarer
mehrbl\"attrigen Fl\"achen, welche in bezug auf die Achse der
reellen Zahlen symmetrisch sind. [\dots] Die Sache ist dann so
einfach, da{\ss} ich sie nicht speziell durchzuf\"uhren
brauche. Der Unterschied ist nur, da{\ss} die in Betracht
kommenden, fr\"uher unbeschr\"ankten Konstanten nunmehr
gezwungen sind, entweder {\it einzeln reell} oder {\it
paarweise konjugiert komplex} zu sein. Infolgedessen
reduzieren sich alle Willk\"urlichkeiten auf die H\"alfte. Wir
m\"ogen folgenderma{\ss}en sagen:

{\it Zur Abbildbarkeit zweier symmetrischer Fl\"achen $p>1$
aufeinander ist neben der \"Ubereinstimmung in den Attributen
das Bestehen von $(3p-3)$ Gleichungen zwischen den reellen
Konstanten der Fl\"ache erforderlich.}

}

\end{quota}

If this passage sounds a bit sketchy to the reader, we may refer
to Klein's subsequent lecture notes of 1892
\cite[p.\,151--4]{Klein_1892_Vorlesung-Goettingen}, where full
details are given.

The basic idea of this (Riemann-style) moduli count is to
represent a given curve of genus $g$ as an $m$-sheeted cover of
the line. If $m$ is
large enough
(so as to avoid exceptional cases of Riemann-Roch's theorem),  a
group $g_m$ of $m$ points will move in a linear system of
dimension $m-g$. To specify a map to ${\Bbb P}^1$ we may
send the divisor $g_m=:D$ to $0$, say, and another $D'$
(linearly equivalent to the former) to $\infty$, leaving the
possibility of a scaling factor. Thus the function depends on
$2m-g+1$ constants. On the other hand by Riemann-Hurwitz such
maps have $2m+2g-2$ branch points. Hence considering the
totality of such covers modulo those yielding the same curve
leaves $2m+2g-2-(2m-g+1)=3g-3$ essential constants.
(cf. also Griffiths-Harris 1978
\cite[p.\,256]{Griffiths-Harris_1978/94}).

[15.12.12] It is legitimate to wonder if this method (\`a la
Riemann-Klein) is powerful enough to compute the gonality profile
(cf. Definition~\ref{gonality-profile:def}).

Klein adapts this counting argument to the real case (again for
full details we recommend Klein 1892
\cite[p.\,151--4]{Klein_1892_Vorlesung-Goettingen}). Doing so we
may hope that he anticipates the Ahlfors mapping when the
construction is particularized to the orthosymmetric case.

Since a
totally real morphism lacks real ramification, we must prescribe
imaginary conjugate branch points. However this necessary
condition is not sufficient as
shown by a quartic smoothing a visible conic plus an invisible one
like $x^2+y^2=-1$ (alternatively consider the Fermat curve
$x^4+y^4=1$ projected from the inside of the unique oval). In this
case the projection from the interior of the oval yields a real
map without real ramification, but
 not
totally real.

We see no obvious link from Klein's equivariant branched covers to
the stronger assertion that fibres over real points consists only
of real points, and consequently one of the orthosymmetric halves
maps conformally to the upper half-plane (as Teichm\"uller
credits to Klein).


Of course it is not impossible that a suitable complement to
Klein's method yields something like an Ahlfors mapping. By a
continuity argument in Gabard 2006
\cite[Lemme~5.2]{Gabard_2006}, it would be enough to chose
$g_m=:D$ as an {\it unilateral} divisor, i.e. one
supported entirely by one half of the curve. Then we would be
finished if the symmetric divisor $D^{\sigma}$ is linearly
equivalent to $D$. But this condition is far from automatic and
involves probably some lucky choice in the position of the initial
divisor $D$.

Alternatively, one may try to specify the ramification and work
out the L\"uroth-Clebsch sort of argument to construct explicitly
the finitely many conformal type of Riemann surfaces lying above
the prescribed ramification. But the writer failed to draw any
serious conclusion.

\subsection{More is less: Teichm\"uller again (1939)}

For those
not overwhelmed by
German prose, the following passage also bears some
resemblances to the Ahlfors function:

\begin{quota}[Teichm\"uller 1939
{\rm\cite[p.\,103]{Teichmueller_1939}}]

{\small \rm Falls $\frak M$ eine orientierte und berandete
Mannigfaltigkeit ist, braucht man $f$ nur auf $\frak M$ zu
kennen, um $f$ auf $\frak F$ berechnen zu k\"onnen. [The
latter is of course the doubled surface.] $f$ mu{\ss} dann auf
den Randkurven von $\frak M$, die ja zu sich selbst punktweise
konjugiert sind, reelle Werte haben. Umgekehrt ist eine
Funktion der Fl\"ache, die in unendlich vielen Randpunkten von
$\frak M$ reell ist, eine Funktion von $\frak M$, denn sie
stimmt mit der konjugierten in unendlich vielen Punkten
\"uberein und ist darum gleich ihrer konjugierten Funktion.
Ja, wir k\"onnen die Funktionen $f$ von $\frak M$ sogar ganz
auf $\frak M$ charakterisieren:

Die Funktionen der orientierten berandeten endlichen
Riemannschen Mannigfaltig\-keit $\frak M$ sind genau die
Funktionen $f$, die in $\frak M$ bis auf Pole regul\"ar
analytisch sind und die am Rande von $\frak M$ reell werden.
D.\,h. die Punkte, wo die Funktion Werte eines abgeschlossenen
Kreises der oberen oder der unteren Halbebene annimmt, sollen
eine kompakte Menge im Innern von $\frak M$ bilden. In der Tat
lassen sich diese Funktionen durch Spiegelung zu Funktionen
von $\frak F$ machen, insbesondere sind sie auf den Randkurven
von $\frak M$ stetig.}

\end{quota}

In this passage we note that just adding the single word
``nur'' in the third line of the 2nd parag. to read ``die nur
am Rande von $\frak M$ reell werden'' would essentially lead
to an anticipation of Ahlfors 1950.

However taken literally this assertion of Teichm\"uller is
weaker than Ahlfors' and indeed the previous Quote
\ref{quote:Teichmueller-1941} is perhaps just a logical
distortion (through hasty writing!) of the above more precise
(but logically weaker) formulation. Under this hypothesis then
we agree perfectly with Teichm\"uller 1941 (cf. again
Quote~\ref{quote:Teichmueller-1941}) that this reality
behaviour of functions was known to Klein.

The crucial distinction is between functions real on the
boundary and those which are real only on the boundary. Now a
priori a real function may be real on an interior point of the
membrane, in which case the range (of the function) will not
be contained in one of the half-plane, but overlap with both
of them. In contrast a stronger reality behaviour arises when
fibres of real points excludes imaginary conjugate points, in
which case the range is contained in one of the half-plane,
which is the context of Ahlfors' circle mapping.

\subsection{Courant 1937, 1939, 1950}\label{sec:Courant}

In the paper Courant 1939 \cite{Courant_1939}, one detects another
approach to the existence of circle maps via the methods of
Plateau's problem (at least so is claimed by Grunsky 1978, cf.
Quote~\ref{quote:Grunsky-1978}). We cite some portion of Courant's
introduction:

\begin{quota}[Courant 1939]\label{quote:Courant-1939}

{\small \rm

The theory of Plateau's and Douglas' problem furnishes
powerful tools for obtaining theorems on conformal mapping.
Douglas emphasized (1931) that Riemann's mapping theorem is a
consequence of his solution of Plateau's problem; then he
treated doubly connected domains and in a recent paper (1939)
multiply connected domains. With a different method I gave in
a paper on Plateau's problem (1937) a proof of the theorem
that every $k$-fold connected domain can be mapped conformally
on a plane domain bounded by $k$ circles. The same method can
be applied to the proof of the parallel-slit theorem and, as
will be shown in the thesis of Bella Manel, to mapping
theorems for various other types of plane normal domains. It
is the purpose of the present paper first to give a
simplification of the method by utilizing an integral
introduced by Riemann in his doctoral thesis, and secondly, to
prove a mapping theorem of a different character referring to
normal domains which are Riemann surfaces with several sheets.
[\dots]

For the case $p=0$, the theorem was stated by Riemann,
according to oral tradition. [See Bieberbach 1925, where a
proof is indicated; and Grunsky 1937, where another proof is
given.]

}

\end{quota}

It should still be
elucidated if this work by Courant (officially overlapping
with Bieberbach-Grunsky) may also be connected to the Ahlfors
circle mapping. This is still not completely clear to the
writer.

The topic is
addressed again in Courant's book of 1950, e.g., as follows:

\begin{quota}[{\rm Courant 1950  \cite[p.\,183,
Thm\,5.3]{Courant_1950}}]

{\small

{\it Theorem 5.3:} Every plane {\rm [footnote 12: As said
before, in view of the general result of Chapter II the
assumption that $G$ is a plane domain is not an essential
restriction.]} $k$-fold connected domain $G$ having no
isolated boundary points can be mapped conformally onto a
Riemann surface $B$ consisting of $k$ identical disks, e.g.
interiors of unit circles, connected by branch points {\rm
[footnote 13: The conformality of the mapping is of course
interrupted at the branch points.]} of total multiplicity
$2k-2$. [\dots]

}

\end{quota}

This somewhat loose footnote 12 of Courant may advance him as a
forerunner of the Ahlfors circle map.
Courant does not specify
the degree
derived by his method, but reading him literally one recovers
(quite strikingly!) Ahlfors' bound $r+2p$ (compare the
following numerology):

\begin{Numerology} {\rm
The connectivity $k$ of a membrane of genus $p$ with $r$
contours is equal to $r+2p$ (each handles contributes 2 units
to the connectivity). [Alternatively, we may interpret the
connectivity $k$ as $b_1+1$, where $b_1$ is the first Betti
number. The Euler characteristic is $\chi=2-2p-r$, but also
expressible as $\chi=1-b_1$ (since $b_2=0$). Back to the
connectivity, we find
$k=b_1+1=(1-\chi)+1=2-\chi=2-(2-2p-r)=2p+r$, as desired.]

Adopting Courant's branching multiplicity $b:=2k-2$,
we
compute the corresponding degree $d$. By Riemann-Hurwitz
$\chi= d \cdot \chi (D^2)-b$, hence
$d=\chi+b=(2-2p-r)+(2k-2)=2k-2p-r=2(r+2p)-2p-r=r+2p$.
{\it q.e.d.}}
\end{Numerology}

This is pure numerology, without much control of the underlying
 geometry. More insight is suggested by Courant's subsequent statement in
{\it loc.\,cit.} \cite[p.\,183--4, Thm\,5.3]{Courant_1950}, which
we reproduce:

\begin{quota}[Courant 1950]

{\small \rm

Moreover, an arbitrarily fixed point $F_{\nu}$ on each
boundary circle $\beta_{\nu}$ can be made to correspond to a
fixed boundary point $P_\nu$ on the boundary continuum
$\gamma_\nu$ of $G$, and the position of one simple branch
point in $B$ may be prescribed. The class $\frak N$ of these
domains depends on $3k-6$ real parameters: the $2k-3$ freely
variable branch points represent $4k-6$ coordinates, while
fixing the points $F_\nu$ reduces the number of parameters by
$k$.

}

\end{quota}

Extending this reasoning to (non-planar) membranes, we derive
again Ahlfors' bound, as follows:

\begin{Numerology} {\rm

We assume the membrane $F_{r,p}$ (of genus $p$ with $r$ contours)
conformally mapped as a $d$-sheeted cover of the disc $D^2$ with
$b$ branch points. As usual the Riemann-Hurwitz relation reads
$\chi=d\cdot \chi(D^2)-b$. From the $b$ branch-points, one of them
can be normalized to a definite position (through a conformal
automorphism of the disc). Now the fibre over a boundary point of
the disc gives $d$ points on $\partial F$. Those $d$ boundary
points can be thought of as having a prescribed image. Thus the
mapping itself is fully determined by $2(b-1)-d$ real constants.

On the other hand, we know since Klein 1882 (cf. our
Quote~\ref{quote:Klein-1882-moduli}) that $F_{r,p}$ has $3g-3$
real moduli where $g$ is the genus of the double $2F$, i.e.
$g=2p+(r-1)$. Positing the Ansatz that the family of $d$-sheeted
covering surfaces has enough free-parameters to fill the full
moduli space leads to the inequation $2(b-1)-d \ge 3g-3$. But
$b=d-\chi$ and $2\chi=\chi(2F)=2-2g$. Hence
$2(d-\chi-1)-d=d+(2g-2)-2\ge 3g-3$, i.e. $d\ge g+1$, which is
Ahlfors' bound $r+2p$.

}

\end{Numerology}

Of course, this happy numerology (noticed by the writer the
[20.05.12]) is
no substitute to a serious proof of the Ahlfors circle map.
However Courant formulates a variational problem \`a la
Plateau-Douglas (or Dirichlet-Riemann-Hilbert)
affording existence of a circle map (presumably with the same
bound as predicted by Ahlfors as prompted by our heuristic count).
Unfortunately, in Courant's book the presentation is
not directly adapted to the case of general membranes of positive
genus ($p>0$), making the reading somewhat hard to digest.
Hopefully someone will manage in the future to present a
self-contained account based upon Courant's method. (This project
involves some hard analysis and will be deferred to a subsequent
technical section. ABORTED: I had not the time/force to adapt
Courant's text to
higher
genera as suggested by his sloppy
footnote 12.)
Of course in view of Carath\'eodory's philosophy (cf.
Quote~\ref{quote:Caratheodory-1928}) one may wonder which of
Courant's vs. Ahlfors approach enjoys methodological superiority?
Further remind that Ahlfors (1950 \cite[p.\,125--6]{Ahlfors_1950})
has also an elementary argument for
circle mapping involving
no extremal problem.

Another puzzling feature of the above numerology is that it gives
the impression that any $r+2p$ points prescribed on the contours
may be mapped to a fixed point of the circle. Whether this is
really true deserves to be investigated.

Trying to read Courant's book 1950 \cite{Courant_1950} with the
focus of the Ahlfors circle map is not an easy task (in our
opinion). We may then hope that reading the original 1939 article
\cite{Courant_1939} is easier due to its more restricted content.
Let us write down its main statement:

\begin{quota}[Courant 1939 {\rm \cite[p.\,814]{Courant_1939}}]
\label{quote:Courant_1939:statement}

{\small \rm

We consider a Riemann surface on a $u,v$-plane consisting of
the interior of $k$ unit circles which are connected in branch
points of total multiplicity $2k-2$; to this surface we affix
$p\ge 0$ full planes with two branch points each. Thus we
define a class of domains $B$ with the boundary $b$ on the
plane of $w=u+iv$.

Now our theorem is: Each $k$-fold connected domain $G$ in the
$x,y$-plane with the boundary curves $g_1,g_2,\dots,g_k$
[\dots] can be mapped conformally on a domain $B$ of our class
for any fixed $p$.

In this mapping the branch points on the full planes and one
more branch point may be arbitrarily prescribed and, moreover,
on each boundary circle $b_\nu$ of $B$ a fixed point may be
made to correspond to a fixed point of $g_\nu$.

}

\end{quota}

Personally, I find this statement hard-to-read for several
reasons, I shall list subsequently. Moreover it is not clear
if suitably interpreted, it really implies the Ahlfors circle
mapping.

How to interpret this statement of Courant? Here are some
critics probably due to the writer's incompetence (rigid
brain)! On the one hand, we have $B$, which moves in a class
of domains. Perhaps those are Riemann surfaces? For instance
the operation of affixing $p$ full planes may give a surface
of genus $p$, at least this is what is suggested by a latter
publication of Courant 1940 \cite{Courant_1940-Acta}, whose
relevant portion we quote again for definiteness:

\begin{quota}[Courant 1940 {\rm \cite[p.\,67]{Courant_1940-Acta}}]
{\small \rm

On the basis of the previous results, the proof of the
characteristic relation $\varphi(w)=0$ for the solution of the
variational problem becomes very simple, if the underlying
class of domains $B$ is chosen not as a domain in the plane
but as a Riemann surface all of whose boundary lines are unit
circles. This class is defined as follows:

We consider for the case of genus zero a  $k$-fold connected
domain $B$ formed by the discs of $k$ unit circles which are
connected in branch points of the total multiplicity $2k-2$.
For higher genus $p$, we obtain domains $B$ by affixing to the
$k$-fold circular disc $p$ full planes each in $4$ branch
points [footnote 2: Each such full plane represents a
``handle'' and increases the genus by $1$.].

}

\end{quota}

Well, but then the domain $B$ of
Quote~\ref{quote:Courant_1939:statement} would have genus $p$.
Then how is it possible for him to get mapped conformally (in
a one-to-one fashion?) to the domain $G$, which seems to be
planar since its connectivity is equal to the number of
boundaries! Perhaps $G$ should be assumed to be $(k+2p)$-fold
connected (or put more briefly $G$ should have genus $p$ and
$k$ contours)?

If so then Courant gives  a (conformal) one-to-one(?) map
(=diffeomorphism) $G\to B$ onto a ``normal'' domain $B$. To
make a link with Ahlfors, it would be desirable to know if $B$
maps to the disc even after the affixing of the $p$ full
planes. (Incidentally, this operation is somewhat poorly
defined, but perhaps better exposed in other publications, cf.
e.g., Courant's book 1950 \cite[p.\,80 and ff.]{Courant_1950}
or Courant 1949/52 \cite{Courant_1949-52:Book}.)

Hence the crucial point would be to know if $B$ is a
many-sheeted cover of the disc, and if yes: how many sheets
are required? Very naively $k+p$ could suffice, in which case
Courant would not only compete with Ahlfors 1950
\cite{Ahlfors_1950}, but also with Gabard 2006 \dots (NB: This
$(k+p)$-sheeted-ness occurs again in Courant 1940
\cite[p.\,78]{Courant_1940-Acta}, and it would be of interest
to decide if this constitutes an anticipation of Gabard 2006.)

If we push our misunderstanding of Courant to its ultimate
limit, we may have the impression that what he do, is an
attempt to mix the parallel-slit mapping he learned from
Hilbert 1909 \cite{Hilbert_1909}, with the
Riemann-Schottky-Bieberbach-Grunsky theorem, but that the
resulting surgery/transplantation does not lead to any really
viable
creature.

Of course, probably much of our misunderstanding is caused not
merely from the difficult mathematics but also from a shift in
language (plus perhaps some inaccuracies due to the torrential
number of publications?), yet we may still hope that either an
appropriate reading (or reorganization) of Courant's thoughts may
lead to an anticipation of the Ahlfors circle map. Hence, we
encourage strongly any reader able to take the defense of Courant
to publish an account in this direction.

Finally, we cite another papers of Courant about conformal maps,
which could be of some relevance:

$\bullet$ Courant 1937 \cite{Courant_1937}, especially p.\,682,
footnote~7, where we read: ``{\it If we assume the possibility of
a conformal mapping on the unit circle for all surfaces admitted
to competition [\dots]}''. This could have some connection with
Ahlfors circle maps, but probably does not. Later on, this article
contains some conformal mapping theorems, which are only announced
without proof. Perhaps, those could be of some relevance.
Especially Fig.\,11, p.\,722, seems to be close to Klein's
R\"uckkehrschnitt-Theorem, and could eventually leads to a proof
of Ahlfors? This paper also relates the ideas of J. Douglas about
minimal surfaces (especially his extended version of the Plateau
problem for surfaces of higher topological structures, where
Douglas uses systematically Klein's symmetric surfaces). One may
therefore wonder if Ahlfors' circle maps may somehow find
application in this grandiose theory of minimal surfaces \`a la
Plateau-Douglas-Rad\'o-Courant, etc. As far as the writer knows no
direct connection is presently available in print, despite the
probable vicinity of both topics.

$\bullet$ Courant 1938 \cite{Courant_1938}, especially p.\,522
``{\it Every plane $k$-fold connected domain can be mapped
conformally to a $k$-fold unit circle\/}''. Hence the result---we
are mostly interested in---occurs here already in 1938. In
contrast to the 1939 version \cite{Courant_1939}, here neither
Riemann, nor Bieberbach 1925 \cite{Bieberbach_1925},
not even Grunsky 1937 \cite{Grunsky_1937} are cited.
Did Courant rediscovered the result independently?

$\bullet$ Finally we quote, Courant 1919 \cite{Courant_1919},
where (under some influence of Hilbert 1909, and Koebe 1909)
conformal mappings to ``normal domains'' are discussed for
non-schlichtartig surfaces (of finite genus). This is also
re-discussed in Courant's book of 1950 \cite{Courant_1950}.

Last but not least, it is perhaps relevant to remind that some
doubts where expressed by Tromba 1983 \cite{Tromba_1983-PREPRINT}
about the validity of Courant's argumentation regarding higher
genus cases of the Plateau-Douglas problem (compare also Jost 1985
\cite{Jost_1985}). It is not clear to the writer if Tromba's
objections compromise seriously the validity of Courant's
assertions (regarding higher genus conformal maps re-derived via
the method of Plateau). This could be a another obstacle toward
completing a Courant-style approach to the Ahlfors map.

\subsection{Douglas 1931--36--39}\label{sec:Douglas}

Having discussed (very coarsely) Courant, it would be unfair to
neglect J. Douglas. His resolution of Plateau's problem interacts
strongly with conformal
mapping,
with the distinctive attitude (partially successful) of not
getting subordinated to the latter. As already pointed out (in
Courant's Quote~\ref{quote:Courant-1939}), Douglas re-derived the
(RMT) as the 2D-case of Plateau (cf. Douglas 1931
\cite[p.\,268]{Douglas_1931-Solution}). Subsequently, Douglas
extended his Plateau solution to configurations of higher
topological structure (cf. Douglas 1936
\cite{Douglas_1936-Some-new-results}, 1939
\cite{Douglas_1939-min-surf}, 1939
\cite{Douglas_1939-The-most-general}).
Thus, it is nearly natural to ask if Douglas (himself, or at least
his methods) may anticipate/recover the Ahlfors circle map?
Ironically, Douglas' work relied on Koebe's, and interestingly
took a systematic advantage of (Klein's) symmetric Riemann
surfaces (e.g., orthosymmetry). Without entering the details
 of all those exciting connections, we just refer to the
cited original works, plus the account of Gray-Micallef 2008
\cite{Gray-Micallef_2008}, of which we quote some extracts:

\begin{quota}[Gray-Micallef 2008 {\rm
\cite[p.\,298, \S 4.3; p.\,299, \S 4.5]{Gray-Micallef_2008}}]

{\small \rm

An unexpec\-ted bonus of Douglas's method is a proof of the
Riemann-Carath\'eodory-Osgood Theorem, which follows simply by
taking $n=2$. [\dots] Douglas was rightly proud that his
solution not only did not require any theorems from conformal
mapping but that some such theorems could, in fact, be proved
using his method.

However, Douglas did have to use Koebe's theorem in order to
establish that his solution had least area among discs
spanning $\Gamma$. He had hoped to fix this blemish, but he
never succeeded. That had to wait for contributions from
Morrey [1948] and, more recently, from Hildebrandt and von der
Mosel [1999]. [\dots]

Even before working out all the details for the disc case, Douglas
was considering the Plateau problem for surfaces of higher
connectivity and higher genus. [\dots] As early as 26 October
1929, Douglas announced that his methods could be extended to
surfaces of arbitrary genus, orientable or not, with arbitrarily
many boundary curves in a space of any dimension. He may well have
had a programme at this early stage, but it is doubtful that he
had complete proofs. Even when he did publish details in [3](=1939
\cite{Douglas_1939-min-surf}), the arguments are so cumbersome as
to be unconvincing. One should remember that Teichm\"uller theory
was still being worked out at that time and that the description
of a Riemann surface as a branched cover of the sphere is not
ideally suited for the calculation of the dependence of the
$A$-functional on the conformal moduli of the surface. Courant's
treatment in [7](=1940 \cite{Courant_1940-Acta}) was more
transparent but still awkward. The proper context in which to
study minimal surfaces of higher connectivity and higher genus had
to wait until the works of Sacks-Uhlenbeck [19](=1981), Schoen-Yau
[20](=1979), Jost [11](=1985) and Tombi-Tromba [21](=1988).
[\dots]

}

\end{quota}

Finally, we mention the recent work of Hildebrandt-von der
Mosel 2009 \cite{Hildebrandt-von-der-Mosel_2009}, plus the
survey Hildebrandt 2011 \cite{Hildebrandt_2011}. Here we
learn, that Morrey 1966 \cite{Morrey_1966} was the first to
re-prove Koebe's KNP (=Kreisnormierungsprinzip) via Plateau,
modulo a gap  fixed by Jost 1985 \cite{Jost_1985}. The
ultimate exposition of 2009 (of loc.\,cit.
\cite{Hildebrandt-von-der-Mosel_2009}) is intended to be
``{\it possibly simpler and more direct\/}'' (loc.\,cit.,
2009, p.\,137) and ``{\it are complete analogs of the approach
of Douglas and Courant\/}'' (loc.\,cit., 2011, p.\,77).

As an agenda curiosity,  the ``Plateau-ization'' of conformal
mapping theorems does occur along diabolic
chronological regularity. From Riemann 1851 \cite{Riemann_1851} to
Douglas 1931 \cite{Douglas_1931-Solution}, gives an elapsing
period of 80 years. For circle maps, we have from Riemann 1858 to
Courant 1939($-1$) also 8 decades, and from Koebe 1904
(announcement of KNP, in his Thesis talk, yet without convergence
proof until 1907/08) to Jost 1984 \cite{Jost_1985} gives the same
interval of time. Thus Ahlfors 1950 \cite{Ahlfors_1950} can safely
wait up to 2030, before getting reproved via the method of
Plateau?

Again, from our focused viewpoint, the critical question is
whether within the   problem of Plateau (\`a la
Douglas--Rad\'o--Courant, etc.)
germinates an alternative proof of the Ahlfors mapping. As far
as we know, the paper closest to this goal his Courant 1939
\cite{Courant_1939}. Yet, we cannot readily claim that it
includes the result of Ahlfors 1950.


\subsection{Cecioni and his students, esp.
Matildi 1945/48, and Andreotti 1950}

Among several interesting works of Cecioni and his students (cf.
Sec.\,\ref{sec:Italian-school}) we point out especially the
article by Matildi 1948 \cite{Matildi_1945/48} (discovered by the
writer as late as [13.07.12]). In it the existence of an
(Ahlfors-type) circle map in the special case of surfaces with a
single contour seems to be established via classical
potential-theoretic tricks, plus at the end some algebraic
geometry. This work was known to Ahlfors (or Sario?) at least as
late as 1960, being quoted in Ahlfors-Sario 1960
\cite{Ahlfors-Sario_1960}. It would be interesting to see if
Matildi's method adapts to an arbitrary number of contours, and
also try to make (more!) explicit the degree bound obtained by
him. In that case Matildi should be considered as a serious
forerunner of Ahlfors 1950 \cite{Ahlfors_1950}, at least at the
qualitative level (no extremal problem). Perhaps, il professore
Cecioni himself has---and may have---several works (some of which
we could not consult as yet) coming quite close to the circle
mapping thematic \`a la Ahlfors.

The idea that Matildi's argument should extend easily to the case
of several contours looks an accessible exercise. Andreotti 1950
\cite{Andreotti_1950} seems to go precisely in this sense.

\subsection{A global picture (the kaleidoscope)}

The place occupied by RMT (Riemann mapping theorem) is quite
pivotal in conformal mapping with an organical explosion of
results around it, like:

\noindent $\bullet$ KNP=Kreisnormierungsprinzip (implicit in
Riemann 1857/8, Schottky 1875 (Latin version of his Thesis, cf.
Klein's Quote~\ref{Klein-1923:quote:Riemann-1858}), in full by
Koebe 1905-10-20).

\noindent$\bullet$ RS=Riemann-Schottky mapping of a
multiply-connected domain to the disc. (This is also known as the
Bieberbach-Grunsky theorem, so RS$\approx$BG, if you want.)

\noindent$\bullet$ AM=the Ahlfors mapping (of a compact bordered
surface to the disc).

\noindent$\bullet$ GKN=generalized Kreisnormierung
(of a compact bordered surface to a circular domain inside a
closed Riemann surface of constant curvature having the same
genus $p$): apart from some anticipation for $p=1$
Strebel 1987 \cite{Strebel_1987} and Jost (unpublished), the full
result is due to Haas 1984 \cite{Haas_1984} (existence), and
Maskit 1989 \cite{Maskit_1989} (uniqueness). For an approach via
circle packings, compare also He 1990 \cite{He_1990} and
He-Schramm 1993 \cite{He-Schramm_1993}.

\noindent$\bullet$ RST=R\"uckkehrschnitttheorem of Klein 1882
\cite{Klein_1882_Ruckkehrschnitt} is yet another form of
generalized Kreisnormierung to positive genera, and for simplicity
we identify it loosely to GKN. The first (rigorous) proof of RST
is to be found in Koebe 1910 UAK2 \cite{Koebe_1910_UAK2}, see also
Bers 1975 \cite{Bers_1975} for a modern account via quasiconformal
deformations.

Of course (at least modulo some sloppiness) we have universal
implications (just by specializing the topological structure) like
$$
\textrm{AM}\Rightarrow \textrm{RS} \Rightarrow \textrm{RMT}
\Leftarrow \textrm{KNP} \Leftarrow \textrm{GKN}.
$$
Besides it is desirable that GKN or RST$\Rightarrow$AM, at least
this would resolve our big historical puzzle about
Klein-Teichm\"uller as anticipating Ahlfors. This desideratum is a
bit cavalier, yet akin to the implication KNP$\Rightarrow$RS,
which is cryptical since Riemann's Nachlass (1857
\cite{Riemann_1857_Nachlass}).

On the other hand there is
a large panoply of methods including:

\noindent$\bullet$ algebraic functions (Abel 1826, Jacobi 1832,
Riemann 1857, etc.),

\noindent$\bullet$ potentials (Dirichlet ca. 1840, Green 1828,
Gauss 1839, Thomson 1848, etc.),

\noindent$\bullet$ iterative methods (Koebe, Carath\'eodory
1905--12),

\noindent$\bullet$ extremal problems (Fej\'er-Riesz 1922,
Carath\'eodory 1928, Ostrowski 1929, etc.),

\noindent$\bullet$ orthogonal systems (Bergman kernel 1922,
Szeg\"o 1921)

\noindent$\bullet$ Plateau-Douglas functionals (Plateau 1849,
Douglas 1930, Courant 1939 via Dirichlet resurrected),

\noindent$\bullet$ circle packings (originally in Koebe 1936,
rediscovered
by Andreev and Thurston 1985 with convergence proof by
Rodin-Sullivan 1986),

\noindent$\bullet$ Ricci flow (Hamilton 1988 \cite{Hamilton_1988},
which specialized to 2D enables one to recover the uniformization
theorem); idem via Liouville's equation (desideratum Schwarz,
followed by Picard 1890--93, Poincar\'e 1899, Bieberbach 1916
\cite{Bieberbach_1916-Delta-u-und-die-automorphen-Funkt}, etc.,
cf. Mazzeo-Taylor 2002 \cite{Mazzeo-Taylor_2002} for a modern
account), and also e.g., Zhang {\it et al.} 2012
\cite{Zhang-et-al_2012}, where a mixed Ricci flow/Koebe's
iteration is advocated.

\begin{figure}[h]
\hskip-55pt \penalty0
    \epsfig{figure=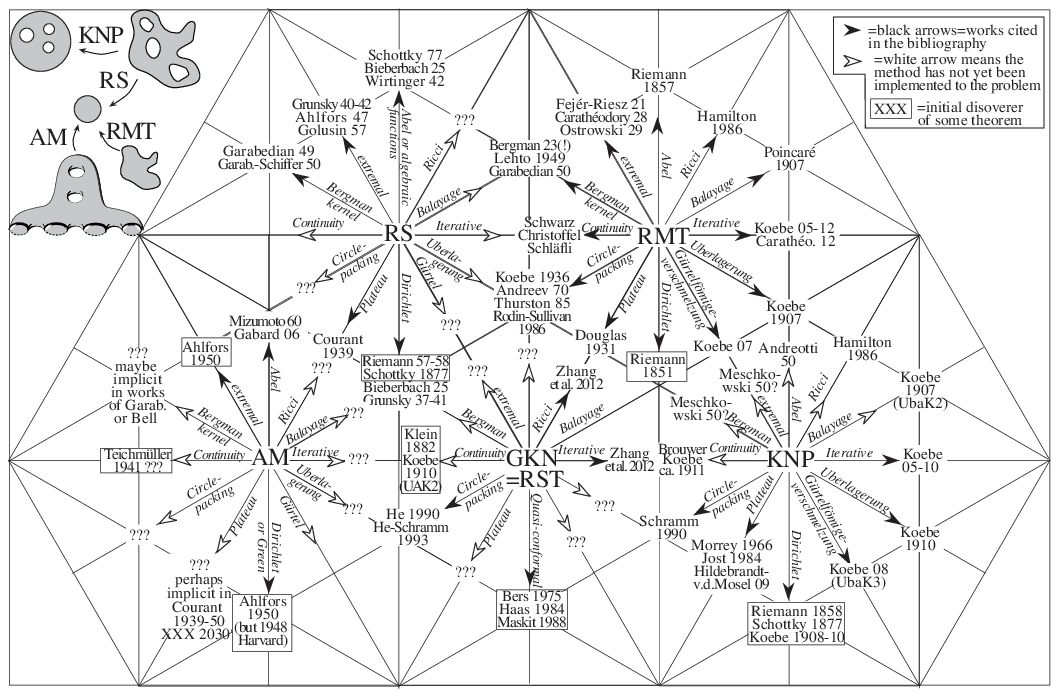,width=162mm}

\caption{\label{Kaleidoscope:fig} The kaleidoscope: several
mapping theorems and some methods used to prove them}
\end{figure}

Blending all these results with  all those methods accessing them,
we get the  kaleidoscope depicted below
(Fig.\,\ref{Kaleidoscope:fig}) attempting to classify a body of
results
in a (more-or-less) systematic fashion. Black arrows stress out
methods effective in solving a certain mapping problem, whose
extremity points to the source (listed in our bibliography).
Starting around RMT,  arrows  are propagated by translation to
other locations (e.g., RS, or KNP). Arrows turn to white colored,
if the corresponding method has not yet been applied to solve the
relevant mapping problem.  Of course several methods (like the
balayage of Poincar\'e 1907, or some of Koebe's method may be
slightly
outdated having few living practitioners). In contrast, Koebe's
iteration method is still quite popular due
to its
computational efficiency (see e.g., Zhang et al. 2012
\cite{Zhang-et-al_2012}), and presumably theoretically fruitful as
 well (\loccit{ }where it is used in conjunction with the
 Ricci flow).

Of course the picture is hard to make completely reliable, yet it
may aid feeling the power (or popularity) of some methods (e.g,
the extremal problem method seems to apply quite universally,
except presently to GKN). On the other hand some recent methods
like circle packings look very powerful, and may not  have as yet
explored their full range of applicability (e.g. regarding AM). As
discussed in the previous section, we do not know if the Plateau
method could crack the AM. Another powerful method is that of the
Bergman kernel, which probably also leads to a derivation of the
AM. When reading papers of the golden period (1948--1950, Bergman,
Schiffer, Garabedian, etc.) this seems to be almost folklore, as
well as in some papers of Bell (e.g. 2002 \cite{Bell_2002}).
While spending some time reading precisely what is put on the
paper, the writer  rather developed the feeling that the positive
genus case is never handled in full details. (As a general
lamentation, it is an easy challenge to cite about 20 papers where
results proved in the planar case are followed by the
apocryphal  allusion that the proof works through {\it mutatis
mutandis\/} without
planarity proviso.)

\section{Digression on Bieberbach and Bergman}
\label{Sec:Bieberbach-Bergman}

\subsection{The Bergman kernel}\label{sec:Bergman}

Among the variety of methods mentioned in the previous section,
one especially popular is the Bergman kernel function. This
emerges in Bergman's Thesis 1921/22 \cite{Bergman_1922}. The point
of departure is an area
minimal problem going back to Bieberbach 1914
\cite{Bieberbach_1914} capturing some salient geometric feature of
the Riemann mapping.
Interestingly, Bergman 1922 (\loccit{}
\cite[p.\,245]{Bergman_1922}) confesses to be not able to reprove
the RMT with this method:

\begin{quota}[Bergman 1922] {\small \rm

In dem betrachteten Spezialfall (Minimalabbildung durch
analytische Funktion) ist die erhaltene Minimalfunktion die
Kreisabbildungsfunktion. Wie oben gezeigt, kann man die
Existenz der ersteren unabh\"angig von dem Hauptsatze der
Funktionentheorie beweisen; es besteht somit die
M\"oglichkeit, den Hauptsatz auf diesem Wege von neuem zu
beweisen, was mir aber bis jetzt nicht gelungen ist.

}

\end{quota}

A similar lamentation is expressed by Bochner 1922
\cite[p.\,184]{Bochner_1922}:

\begin{quota}[Bochner 1922] {\small \rm

Aus der M\"oglichkeit der Kreisuniformisierung eines einfach
zusammenh\"angenden Bereiches folgt aber, wie Bieberbach
bemerkt hat (l.\,c.), da{\ss} die Minimalabbildung mit eben
der Kreisabbildung identisch ist, indes ist es mir nicht
gelungen, aus der Minimalabbildung der Kreisuniformisierung
aufs neue herzuleiten. }

\end{quota}

In a similar vein, some 3 decades later one among the prominent
aficionados of the method wrote (source=Math.-Reviews for Lehto's
Thesis 1949 \cite{Lehto_1949}):

\begin{quota}[Nehari 1950] {\small \rm

Despite its great intrinsic elegance and its adaptability for
numerical computations, the theory of complex orthonormal
functions (centering about the concept of the Bergman kernel
function) had the drawback of being a mere representation
theory; the fundamental existence theorems had to be borrowed
from other fields. In $\S 4$ the author fills this gap in one
important instance by giving an existence proof for the
parallel-slit mappings (in the case of simply-connected
domains this is identical with the Riemann mapping theorem
[provided the slit is extended to $\infty$]) within the
framework of the orthonormal function theory.

}

\end{quota}

So  somewhere in between 1922--1949 some technological turning
point must have occurred amplifying  dramatically the power of the
Bergman kernel method. When and how did this occurred exactly?
Probably through the Bergman--Schiffer collaboration in the 40's,
plus some fresh blood like Garabedian or Lehto. In several
subsequent publications of Garabedian and Schiffer, it is
emphasized that parallel-slit mappings are easier than circle maps
(cf. Quotes~\ref{quote:Garabedian-Schiffer_1950} and
\ref{quote:Garabedian_1949-52}). However the Ahlfors circle
mapping seems accessible to the Bergman-Szeg\"o orthogonal system
method as
suggested in Garabedian-Schiffer 1950
\cite{Garabedian-Schiffer_1950}, where only the planar case is
handled in detail. Often in literature, yet not in the just cited
paper, it is sloppily insinuated that a method implemented in the
planar case extends to Riemann surfaces. A typical
specimen is the earlier paper by Garabedian 1950
\cite{Garabedian_1950} claiming another proof of the full Ahlfors
theorem by deploying a broad spectrum of techniques (yet not
readily reducible to the Bergman kernel) ranging from
Teichm\"uller 1939 \cite{Teichmueller_1939-Dreikreisesatzes},
Grunsky 1941--42 \cite{Grunsky_1940}, \cite{Grunsky_1942}, Ahlfors
1947 \cite{Ahlfors_1947} and Schiffer's inner variations.

Inspecting back the Bergman method itself, it is not hard to
understand why it is most readily implementable in the planar
case. It seems indeed to require a sort global ambient coordinate
system. Let us look at the beautiful original paper Bergman 1922
\cite[p.\,240]{Bergman_1922}. Here the key idea is a
characterization of the Riemann function $w\colon B \to \Delta$
(of a [simply-connected] domain $B$ to the disc) as the one whose
range $w(B)$ has smallest possible area amongst
all functions $f\colon B \to {\Bbb C}$
constrained by $f'(0)=1$ (and $f(0)=0$ after harmless translation
so that $0\in B$).
The area swept out by $f$ is calculated by the integral
$$
\int \int_B f'(z) \overline{f'(z)} d\omega,
$$
where $d \omega$ is the surface element in the $B$-plane, while
the integrand  $\vert f'(z)\vert^2$ measures  the distortion
effected at $z$.
Following Bieberbach 1914 \cite{Bieberbach_1914} (who in turn
seems inspired by Ritz 1908 \cite{Ritz_1908}), Bergman plugs in
place of $f(z)$ a polynomial
(recall the  finitistic motto of Bloch ``{\it Nihil est
infinito \dots }''):
$$
w_n(z)=a_0+a_1z+\dots+a_n z^n,
$$
with coefficients determined as to minimize the above integral
under the constraint $w_n'(0)=1$. There is always a unique such
polynomial, which is computed by usual methods
(finite extremum-problem). The limiting function $\lim_{n\to
\infty} w_n$ gives---again Bieberbach is cited---the required
mapping. The method is so simple and elegant that it is hardly
conceived why it fails to reprove the RMT  (which Bergman and
others call the {\it Hauptsatz der Theorie der konformen
Abbildung} (\loccit{ }p.\,240)). The reason is however a quite
simple vicious circle, namely that the above (Bieberbach)
``areal'' characterization of the Riemann function logically rests
on RMT. Hence the minimum function (of Bergman) is eminently
computable,
but the resulting power series may not have a priori the required
geometrical property of univalence and the right disc-range. I
guessed the latter property follows from Bieberbach 1914
\cite{Bieberbach_1914}, hence the real problem is univalence.
%
However on [13.06.12], after reading Bergman 1947
\cite[p.\,32]{Bergman_1947}, the opposite looks true: namely
univalence is easy but the disc-range issue is not. There are
mentioned two contributions, one by Bergman 1932
\cite{Bergman_1932} and also Schiffer 1938
\cite{Schiffer_1938-CRAS-domaines-minima} where the
desideratum (of reproving RMT) is established for starlike
domains. So almost as importantly, this source (of 1947)
points out that Bergman's dream of 1922 (new proof of RMT via
the area extremum problem) was not borne out until 1947, and
therefore seems really to be credited to the newer generation
like Garabedian and Lehto.

Generally speaking, extremum problems are often solvable (even
uniquely soluble), but it is another
piece of careful analysis to control precisely the geometric
behavior of solutions, e.g. in the hope to re-crack RMT. Of
course, the problem was ultimately solved, cf. e.g., Garabedian
1950 \cite{Garabedian_1950} or the already mentioned Thesis of
Lehto 1949 \cite{Lehto_1949}, which are the first completed
Bergman-style approaches to RMT.

The point for re-exposing the hearth of the method is to emphasize
the r\^ole of polynomials generated by $z^n$ as a preferred system
of global functions on the domain $B$ out of which an ideal object
is processed through an extremum
procedure handled {\it in finito}. How can one adapt this on a
Riemann surface where no global parameters are supplied a priori?
This is a little puzzle to the writer [06.06.12], but the masters
(Bergman, Garabedian, Bell, etc.) often claim the method to
suit the broader context with minor changes. Compare, e.g., the
following sources:

$\bullet$ Bergman 1950 \cite[p.\,24, Remark]{Bergman_1950}
justifies in this book extensibility to Riemann surfaces by
referring to results of Sario 1949--50.

$\bullet$ Garabedian 1950 \cite[p.\,361]{Garabedian_1950}, where
one reads ``{\it For the sake of a simple presentation of results
we have merely stated the theorem for the case of schlicht domains
of finite connectivity. However the theorem is true with only one
change if $D$ is a Riemann surface {\rm [\dots]}.   The reader
will easily verify that the proof which we shall give of the
theorem carries over with minor changes to the more general
situation.}'' If not pure bluff, it is sad that Garabedian did not
write down the details at that time. If we believe in the unity of
mathematics especially the algebro-geometric curves and analytic
Riemann surfaces at the compact level, then the existential aspect
of circle maps is frankly more trivial in the ``schlicht'' and
even ``schlichtartig'' case, compare e.g., the argument in Gabard
2006 \cite{Gabard_2006}  (reproduced below as
Lemma~\ref{Enriques-Chisini:lemma}), which in substance is the one
of Bieberbach 1925 \cite{Bieberbach_1925}, Wirtinger 1942
\cite{Wirtinger_1942}, but perhaps slightly streamlined by the
mere usage of algebro-geometric language.

\subsection{Minimizing the integral
vs. maximizing the derivative (suction vs. injection), i.e.
Bieberbach 1914-Bergman 1921/22 vs. Fej\'er-Riesz 1922, etc.}


Trying to avoid the vicissitudes of life concomitant with the
Dirichlet principle, the early 1920's imagined two
methods of attack to the RMT via extremum problems. Given $B \ni
a$ a simply-connected domain in the complex plane ${\Bbb C}$,
which is not the plane and therefore can easily be assumed to be
bounded via a suitable transformation, RMT amounts to find a
conformal map to the disc. The following (animalistic) acronyms
are derived by contracting the contributors' names:

\medskip \noindent$\bullet$ (BIBER)=(Bieberbach 1914
\cite{Bieberbach_1914} and Bergman[n] 1922
\cite{Bergman_1922}). [Biber=German for ``beaver''
(=``castor'' in French).]

\noindent$\bigstar$ {\it Amongst analytic functions $f\colon B
\to {\Bbb C}$
normed by $f(a)=0$ and $f'(a)=1$ minimize the integral $\int
\int_B \vert f'(z)\vert^2 d\omega$, where $d\omega$ is the surface
element of the Euclidean metric. }

\medskip
\noindent$\bullet$ (FROG)=(Fej\'er-Riesz 1922, Carath\'eodory 1928
\cite{Caratheodory_1928}$\leftrightarrow$Ostrowski 1929
\cite{Ostrowski_1929}, and Grunsky 1940 \cite{Grunsky_1940},
Ahlfors 1947 \cite{Ahlfors_1947} in the multiply-connected
context)

\noindent$\bigstar$ { \it Amongst analytic functions $f\colon
B \to \Delta$(= unit disc)
normed by $f(a)=0$ maximize the modulus $ \vert f'(a)\vert$.}
\medskip

As remembered in the previous section, the problem BIBER was not
prompt in supplying an autonomous proof of RMT, while succeeding
only in the late 1940's (Garabedian's or Lehto's Thesis). Further
this succeeded perhaps only under the  proviso of smooth boundary
(Jordan curve), cf. e.g. Garabedian-Schiffer 1950
\cite[p.\,164]{Garabedian-Schiffer_1950}: ``{\it The most serious
drawback in our method is, perhaps, that we must make assumptions
upon the smoothness of the boundary of the domain we consider, so
that the general case is reached only after a topological
approximation argument is given.\/}''.

In contrast FROG
met earlier success (cf. e.g., Ostrowski 1928/29
\cite{Ostrowski_1929} and Carath\'eodory 1928
\cite{Caratheodory_1928}) streamlining previous work of
Fej\'er-Riesz 1922 (published in Rad\'o 1923
\cite{Rado_1923-Uber-konf-Abb}).

For extensions  to multiple-connectivity, or even Riemann
surfaces, we have the following contributions:

$\bullet$ FROG leads to the works of Grunsky 1940--42
\cite{Grunsky_1940,Grunsky_1942} (schlicht domains of finite
connectivity) and Ahlfors 1950 \cite{Ahlfors_1950} (non-planar
compact bordered Riemann surfaces), where the derivative $f'(a)$
is computed w.r.t. any local chart. In fact Ahlfors rather
considers the  variant where given two points $a,b$ the modulus of
$f(b)$ has to be maximized amongst functions with $f(a)=0$.

$\bullet$ BIBER is somewhat harder to formulate on a Riemann
surface $F$ (taking the r\^ole of the domain $B$) as the
magnitudes involved in the problem require something more than the
Riemann surface structure. A Riemannian metric would make the
problem meaningful, but which metric to choose? Of course there is
the canonical conformal metric given by uniformization of the
doubled membrane $F$. Of course we deviate slightly from a
self-contained proof of RMT or Ahlfors (=AMT), but this is maybe
not a dramatic concession.

Thus, even in its basic formulation, some ideas are required
to
set a perfect analogue of the problem BIBER for a (bordered)
surface. If this can be done, it is
likely (or desirable) that the extremal function (whose existence
and uniqueness is derived by Hilbert's spaces arguments) is a
circle map, i.e. effects a conformal representation over the disc.
(This is a priori not the unit disc, but renormalize so.) In the
simply-connected case, both extremals of BIBER and FROG (denoted
$\beta$ and $\alpha$ respectively) yield the one and the same
object, namely the Riemann mapping $B\to\Delta$ (again after a
harmless scaling of $\beta$, cf. Bergman 1950
\cite[p.\,24]{Bergman_1950} for its exact value in terms of the
Bergman kernel). Hence, it is plausible that the least area map
for the surface $F$ coincides  with the Ahlfors function. So this
would be a sort of conformal identity, perhaps of some practical
significance.

Of course, the primary interest would be to reobtain (via
BIBER) a novel proof of Ahlfors 1950 \cite{Ahlfors_1950}.
(This game may be already implicit in several works, as those
of Bergman and Garabedian itemized in the previous section,
but no pedestrian redaction is available in our opinion.)
Yet, the real novelty
would be the resulting
``binocular view'' of the one and same object (i.e., the Ahlfors
extremal) through two different angles, yielding a sharper
perception of it. Perhaps, this gives sharper
differential-geometric insights about the Ahlfors map of a
membrane, and
may have some implications toward
Gromov's
FAC(=filling area conjecture).
Remember our naive conviction that this problem FAC should succumb
just under the powerful methods of 2D-conformal geometry.

\subsection{Bergman kernel on Riemann surfaces}

[13.06.12] Consulting other sources (e.g. Weill 1962
\cite{Weill_1962}), it seems that  the theory of the Bergman
kernel can be developed over any Riemann surface. The idea is to
use the Hilbert space structure on the space of analytic
differentials. A complete exposition is e.g., Ahlfors-Sario 1960
\cite[p.\,302]{Ahlfors-Sario_1960}.
Whether or not this leads to another proof of Ahlfors circle maps
is another question.

[15.06.12] Other references for the Bergman kernel on Riemann
surfaces include Nagura 1951 \cite{Nagura_1951}, and Nehari 1950
\cite{Nehari_1950} where the Ahlfors function is expressed in
terms of the Bergman kernel.

[25.06.12] In fact the key observation is probably the {\it
conformal invariance\/} of the integral involved in the minimum
problem BIBER (of the previous section). Thus it may be hoped that
this problem leads to
an independent treatment of the Ahlfors mapping, treated from a
Hilbert space [of ``areally'' ({\it a\'erolaire})
square-integrable holomorphic functions] viewpoint. This would
give some culmination to the
device  of Bieberbach 1914 \cite{Bieberbach_1914}.

So having in mind the possibility of extending the BIBER minimum
area problem of the previous section to compact bordered Riemann
surfaces (which looks reasonable in view of the conformal
invariance of this area functional) we would like to reprove the
existence of a circle map (\`a la Ahlfors 1950
\cite{Ahlfors_1950}).

Relevant literature on this problem (but from our naive viewpoint
not completely satisfactory) includes in chronological order:

$\bullet$ Bieberbach 1914 \cite{Bieberbach_1914} (simply-connected
schlicht case)

$\bullet$ Bergman 1950 \cite[p.\,24]{Bergman_1950}, where the fact
that the range of the minimizing function  is a circle is
considered as well-known (with reference to Bieberbach's Lehrbuch
(1945 edition) \cite{Bieberbach_1945-Lehrbuch}). Later in
Bergman's book 1950 \cite[p.\,87]{Bergman_1950} the circle map
$B\to \Delta$ is recovered through the function $F(z,
\zeta)=\frac{\hat K(z, \bar \zeta)}{\hat L(z,  \zeta)}$ defined on
p.\,86, but it is not clear if this function solves the least area
problem. (Perhaps the connection is easy to do.)

$\bullet$ Garabedian-Schiffer 1950
\cite[p.\,166--7]{Garabedian-Schiffer_1950} where the BIBER
problem is again formulated, but somehow only in the purpose of
showing existence of the reproducing kernel function, in the optic
of re-deriving the PSM (parallel-slit maps). In particular one may
wonder if it possible to show by a direct analysis if the minimum
function is a circle map. Circle maps are reobtained later in the
paper (p.\,182) however through a
different procedure.

$\bullet$ Nehari's book 1952 \cite{Nehari_1952-BOOK} where the
BIBER minimum problem appears on p.\,362 (for multiply-connected
domain only) and its relation to the Bergman kernel is made
explicit in the subsequent pages (esp. p.\,368-9). However I do
not think that the issue about the circle mapping property of the
minimum function of BIBER is handled. Later in the book (p.\,378)
the Ahlfors extremal function is treated, yet a priori there is no
clear-cut identity between the Bieberbach and Ahlfors extremal
function. Nehari's book borrows a lot of ideas from other writers
without referring to them, thus it is an easy task to observe
strong overlap with the previous literature (e.g. Bergman 1950 and
Garabedian-Schiffer 1950).

\subsection{$\beta$ and $\alpha$ problems}\label{sec:beta-and-alpha-problems}

[27.06.12] As already discussed in Section 7.8, there are
essentially two problems BIBER and FROG amounting respectively to
minimize an integral and to maximize a derivative. We may
rebaptize them respectively the $\beta$-problem (for
Bieberbach-Bergman) and the $\alpha$-problem for Ahlfors (albeit
this should truly be Fej\'er-Riesz 1922, for historical
sharpness).

For simplicity we restrict to domains,
though the ultimate dream is to concoct
didactic expositions pertaining to Riemann surfaces.

Regarding the $\beta$-problem (of minimizing the areal integral)
it has a direct Hilbert-space interpretation (recall the
affiliation
Dirichlet-Riemann-Hilbert-(Schmidt)-Ritz-Bieberbach-Bergman), as
finding the vector of minimal length on the hyperplane defined by
the prescription $f'(t)=1$, where $t$ is some fixed point
(previously denoted $a$). Such minimization traduces into
orthogonality to this hyperplane, yielding the so-called {\it
reproducing property\/} while permitting to identify the
$\beta$-extremal with the Bergman kernel (function). For a
detailed execution, cf. e.g. Garabedian-Schiffer 1950
\cite[p.\,166--7]{Garabedian-Schiffer_1950} (henceforth abridged
GS50).

Likewise the $\alpha$-problem received ultimately a similar
treatment through Garabedian's Thesis 1949 \cite{Garabedian_1949}
(recast in the just cited Garabedian-Schiffer article), but the
treatment is somewhat more involved appealing to the Szeg\"o
kernel instead, characterized via an orthogonalization taking
place along the boundary of the domain (hence in substance the
idea of length rather than area). It follows in particular an
explicit formula for the derivative of the Ahlfors function $\vert
f'(t) \vert=2\pi k(t,t)$ in term of the Szeg\"o kernel.
(Garabedian's work is such a tour de force that it was represented
in virtually all major texts of that period, e.g. Bergman 1950
\cite{Bergman_1950} and Nehari 1952 \cite{Nehari_1950}, plus also
the paper GS50.)

Can we understand better the connection between both extremal
problems? Our naive question is whether the $\beta$-map is a
circle map. Remember that Bieberbach 1914 \cite{Bieberbach_1914}
has an argument in the case where the domain is simply-connected
(via his first Fl\"achensatz saying that a map from the disc with
normalized derivative expands the area of the disc unless it is
the identity). Combining this with the Riemann mapping, Bieberbach
argues that the $\beta$-map must be disc-ranged, for otherwise we
could deflate the area by post-composing with the Riemann map,
hence violating the minimum property.

Alas, it seems that this argument is hard (impossible?) to extend
to the multiply-connected case. Thus it is puzzling to wonder if
the $\beta$-map is a circle map. If it is the case, then we could
inject the $\beta$-solution into the Ahlfors problem and compare
them. In view of the explicit formula of Garabedian we can even
try a direct comparison of the respective derivatives at $t$ and
hope to find an equality in which case by uniqueness we would have
$\beta=\alpha$ (modulo scaling), i.e. a perfect coincidence
between the Bieberbach and Ahlfors functions.

Of course, ideally everything should be done geometrically from
the extremal problem, without duelling with hard analysis. Recall
that each problem has its allied reproducing kernel, which serves
to express its solution. In particular we may hope to derive the
circle mapping property of the $\beta$-function from the property
of its allied (Bergman) kernel (cf. GS50, p.\,167). And if not, we
may hope to connect the $\beta$ to the $\alpha$-map through a
somewhat accidental identity between their  kernels functions. As
far as the writer knows this is not explicitly made, and perhaps
wrong.

Let us emphasize a naive duality between the $\alpha$- and
$\beta$-problem. The first amounts to a maximal pressurization
(inflation) within a limited container (the unit disc), whereas
$\beta$ is a free vacuum deflation
leading ineluctably toward a big-crunch to a point (constant map)
if there were not the initial explosion sustained by the
derivative normalization $f'(t)=1$. Hence it is not so surprising
that the Ahlfors map is a circle-map but the same issue for the
Bieberbach least-area map
seems more like an isoperimetric miracle.

We learned from Gaier's 1978 survey \cite[p.\,34--35, \S
C]{Gaier_1978-JDMV} the following piece of information. Gaier's
article contains a proof of a striking fact due to Gr\"otzsch 1931
(see also Gaier 1977 \cite{Gaier_1977-Roumaine}, where the precise
ref. is identified as Gr\"otzsch 1931 \cite{Groetzsch_1931}) that
a map (non-unique!) minimizing the area integral $\int\int \vert
f'(z) \vert^2 d \omega$ (\`a la Bieberbach 1914
\cite{Bieberbach_1914}--Bergman[n] 1922 \cite{Bergman_1922}, but
extended to the multiply-connected setting) under the schlichtness
proviso (and the normalization $f(a)=0, f'(a)=1$) maps the domain
upon a {\it circular slit disc\/} (with concentric circular slits
centered at the origin). According to Gaier, Gr\"otzsch's paper
contains no details outside the indication of using his {\it
Fl\"achenstreifenmethode} (striptease method). Gaier's proof
combines Carleman's isoperimetric property of rings (relating the
modulus to the area enclosed by the inner contour) with Bieberbach
1914 \cite{Bieberbach_1914} (first area theorem) to the effect
that a schlicht normalized map from the disc inflates area, unless
it is
the identity. A natural question [13.07.12] is what happens if we
relax schlichtness of the map? Do we recover an Ahlfors circle
map?

As a  historical curiosity, Gaier 1977 \cite{Gaier_1977-Roumaine}
remarks that the above least-area problem for schlicht functions
was reposed as a research problem as late as 1976 in the Durham
meeting by Aharonov (compare for the exact ref. the Math. Review
by Burbea of Gaier 1977 \cite{Gaier_1977-Roumaine}). It is
apparently K\"uhnau (Gr\"otzsch's eminent student) who pointed
Gr\"otzsch's priority in the reference just cited (Gr\"otzsch 1931
\cite{Groetzsch_1931}). It should be remembered that several
treatments existed in print (prior to Aharonov's question), e.g.
the one in Sario-Oikawa's book of 1969 \cite{Sario-Oikawa_1969}
(see pages as in MR of Gaier 1977 \cite{Gaier_1977-Roumaine}),
which is inspired from Reich-Warschawski 1960
\cite{Reich-Warschawski_1960}. All these treatments are quite
involved, and Gaier 1977 \cite{Gaier_1977-Roumaine} claims to
simplify them.

A paper related to Gaier's and to this circle of ideas---i.e.
Bieberbach's area minimization, yet, alas not exactly furnishing
our naive desideratum---is Alenicyn 1981/82
\cite{Alenicyn_1981/82}: this gives the exact reference to the
relevant work of Carleman 1918 \cite{Carleman_1918} as well as to
that of Vo Dang Thao 1976 \cite{Vo-Dang-Thao_1976} (the latter
being however slightly criticized for  mistakenly assuming the
schlichtness of some function).

Philosophically such Bieberbach-type area minimization problem
amounts to a deflation as opposed to the inflation of Ahlfors-type
problem maximizing  the derivative. According to popular wisdom,
both viewpoints could coincide since a semi-empty bottle is the
same as a half-filled one. (This reminds the story of Ahlfors'
whiskey bottle used as a defense-weapon against an aggressor.)

[17.07.12] We can also switch completely of extremal problem by
looking at an Ahlfors (for short $\alpha$-type) extremal
(inflationist) problem of maximizing the derivative among schlicht
functions. Given $D$ a multiply-connected domain  marked
interiorly at the point $a\in D$, find among all schlicht
functions $f\colon D \to {\Bbb C}$ bounded-by-one $\vert f\vert
\le 1$ the one maximizing the modulus of the derivative $f'(a)$.
It is reasonable to guess that ``the'' (unique?) extremal map will
take $D$ upon the full circle with circular slits (schlichtness
being only fulfilled on the interior). It seems that this behavior
is the one described in Meschkowski 1953 \cite{Meschkowski_1953}
(basing his analysis upon a distortion result of Rengel 1932
\cite{Rengel_1932-33}), and see
also the treatment by Reich-Warschawski 1960
\cite{Reich-Warschawski_1960}. Added 27.07.12: Compare also Nehari
1953 \cite[p.\,264--5]{Nehari_1953-Inequalities}, where another
treatment of this problem is given, and credits is given to
Gr\"otzsch 1928 \cite{Groetzsch_1928} and Grunsky 1932
\cite{Grunsky_1932}.

{\bf Optional digression:} Asking schlichtness up to the boundary,
we get maybe the Kreisnormierung of Koebe? This would be
interesting since as pointed out in one of Meschkowski's paper
cited in the bibliography (locate where exactly!?, but anecdotic
because cf. also Schiffer-Hawley 1962 \cite{Schiffer-Hawley_1962},
Hejhal 1974 \cite{Hejhal_1974}, etc.) there was in the 1950's no
clear-cut extremal problem leading to the Kreisnormierung (even in
finite connectivity). Maybe the situation changed slightly after
several works of Schiffer (and his collaborator Hawley) where some
Fredholm eigenvalues came into the dance (compare several refs.
cited below in the period 1959--1963).

\begin{figure}[h]
\centering
    \epsfig{figure=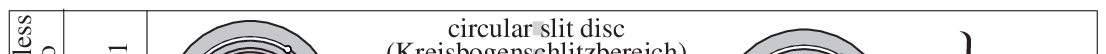,width=102mm}

\caption{\label{alpha:fig} Two prominent types of extremal
problems: on the left maximizing the derivative $\vert f'(a)
\vert$ (inflationist) and on the right minimizing the integral
$\int\!\!\int_D\vert f'(z) \vert^2 d \omega$ (deflationist). On
the top-part schlichtness of the mapping is imposed, whereas on
the bottom part no-schlichtness is imposed allowing all analytic
functions in the competition. It is
tempting to conjecture [extrapolating the symmetry of the
top-part] (and perhaps already known) that the minimum
$\beta$-function is a circle map. (Works cited to be found in the
bibliography)}
\end{figure}

At this stage combining the analysis of Gaier 1978
\cite{Gaier_1978-JDMV} for the $\beta$-problem and that of
Meschkowski/Reich-Warschawski for the $\alpha$-problem (refs. as
in the penultimate paragraph) we see a perfect duality between the
behavior of the extremal {\it schlicht\/} functions (at least
qualitatively since both mappings carry the domain upon the same
canonical region of a circular slit disc). Maybe one can even
identify both functions (after harmless scaling). Those works
raise some hope that
the schlicht-relaxed $\beta$-problem (area minimization \`a la
Bieberbach) produces again the Ahlfors map (or at least enjoys the
same property of being a circle map). As far as we know
[20.07.12], there is no such published account corroborating this
intuition. This would be highly desirable to complete the symmetry
of the picture below (Fig.\,\ref{alpha:fig}) summarizing our
discussion.

[22.07.12] On reading Alenicyn 1981/81
\cite[p.\,202]{Alenicyn_1981/81}, 1981/82 \cite{Alenicyn_1981/82},
where one is referred for the least-area problem back to Nehari's
book of 1952 \cite{Nehari_1952-BOOK}, especially pp.\,340 (one can
safely add p.\,341) and p.\,362. Nehari's pages\,340--341 are
perhaps not so relevant as it is merely a set of exercises. What
is truly relevant is page 362, where the least area problem is
posed and partially analyzed. In fact, this least area problem is
handled earlier (with somewhat sharper information) in
Garabedian-Schiffer 1949 \cite[p.\,201]{Garabedian-Schiffer_1949}
where the solution is represented as $M(z,a) \,
M'(a,a)^{-1}=:M^{\ast}(z,a)$, where $M(z,a)=[A(z,a)-B(z,a)]/2$ is
a combination of $A,B$ the two canonical parallel slit maps of the
domain $B$ upon horizontal (resp. vertical) slit domains taking
$a$ to $\infty$ as a simple pole with residue $+1$ (compare
\loccit\,p.\,200).

[26.07.12] In fact this solution is already announced in Grunsky's
Thesis 1932 \cite[p.\,140]{Grunsky_1932}! As to the geometry of
this map $M^{\ast}$, Garabedian-Schiffer (\loccit\,p.\,201) add
the fact that it is at most $n$-valent ($n$ being the number of
contours of the domain, equivalently, its connectivity). (This
information is not to be found in Nehari 1952
\cite{Nehari_1952-BOOK}.) Alas, Garabedian-Schiffer (1949 \loccit)
never seem to assert that the least-area map $M^{\ast}(z,a)$ is a
circle map. On p.\,217, they show that any unitary function $E$
(=unit-circle map) may be expressed as a linear combination of the
least-area maps $M(z, n_{\nu})$ centered at the $N$ zeros
$n_{\nu}$ (`Nullstellen') of $E$ (assumed to be simple), compare
Eq. (131) and (131'). Finally, on p.\,219 it is observed that the
area of any such $E$, mapping the domain $D$ upon the unit-circle
covered $N$ times, is exactly $N\cdot \pi$ (since area as to be
counted with multiplicities). Of course, if our conjecture about
the circle-mapping nature of least-area maps (there is one for
each center $a$) is correct, then we could sharpen
Garabedian-Schiffer's assertion about the ``at most $n$-valency''
into an exact $n$-valency of those maps.

[27.07.12] It could be the case, that our conjecture about the
circle mapping nature of the least area map is settled in Lehto's
Thesis 1949 \cite{Lehto_1949} (see especially p.\,41).

[29.07.12] However on consulting M. Maschler 1959
\cite{Maschler_1959} (esp. p.\,173) it seems to be  asserted that
the range of the least area maps are unknown for domains of
connectivity higher than $2$.

[26.07.12] To our grand surprise, we notice that the least-area
problem is handled in full generality (i.e., for compact bordered
Riemann surfaces) in Schiffer-Spencer 1954
\cite[p.\,135]{Schiffer-Spencer_1954}. However again (as in
Garabedian-Schiffer 1949 \cite{Garabedian-Schiffer_1949}) it is
not shown that the resulting extremal function is a circle map.

At this stage we see that there is a wide variety of extremal
problems, and as a rough rule we may split the most common of them
into the $\alpha$- and $\beta$-type (for Ahlfors and Bieberbach
resp.) Each problem is hard to analyze precisely but there is a
large body of wisdoms accumulated about them by the masters
(Koebe, Carath\'eodory, Bieberbach, Gr\"otzsch, Grunsky, Ahlfors,
Schiffer, Garabedian, Golusin, etc.) Optionally by a nebulous
bottle principle there may be a certain duality (even possibly an
identity) between $\alpha$- and $\beta$-solutions. At least so is
the case in the simplest simply-connected setting according to
Bieberbach 1914 \cite{Bieberbach_1914}, and apparently in the
multi-connected setting we have at least coincidence of the range
when considering the restricted schlicht problems. We may also
speculate that a careful analysis of a suitable extremal problem
may lead to a solution of the Gromov filling area conjecture.

Finally we  mention a related  extremal problem treated in
Schiffer 1938 \cite{Schiffer_1938}, namely that of minimizing the
maximum modulus in the family of schlicht functions $f\colon B \to
{\Bbb C}$ normalized by $f(a)=0$ and $f'(a)=1$. The (or rather
any) extremal is shown to map (conformally) the Bereich $B$ upon a
circular slit disc.

\subsection{Least area problem vs. least momentum}

[03.08.12] The menagerie of extremal problems leading to the
Riemann mapping can still be further enlarged. Each extremal
problem exploits the ordered structure of the real line
via some real-valued functional. One may incidentally get some
feeling of regression about this massive usage of real numbers in
complex geometry problems, but this is common and respectable
practice since Dirichlet's principle. Regarding the problem of
circle maps {\it per se\/} it is not perfectly clear what is the
{\it ideally suited\/} extremal problem (if any beside that of
maximizing the derivative)? What is somehow missing is an extremal
principle selecting the best extremal problem! The competitive
nature of such extremal problems
fascinated generations but requires strong classification
aptitudes in view of the difficulty of each problem and the
diversity of them.

First the {\it least-area problem\/} consists in minimizing the
area of the range of an analytic function counted
by multiplicity. This is measured by the functional $A[f]=\int
\int \vert f' \vert^2 d\omega$ (which seems much allied to the
Dirichlet integral).  (To extend the problem to Riemann surfaces
one just needs to take notice of the conformal invariance of this
integral upon conformal change of metrics.) To avoid the
minimizers collapsing to the (uninteresting) constant functions,
one imposes the side condition $f'(t)=1$ at some inner point $t$
of the domain $B$. The least-area map (which exists uniquely  by
Hilbert space theory) effects when $B$ is simply-connected nothing
but than the Riemann mapping (Bieberbach 1914
\cite{Bieberbach_1914}). This viewpoint was widely pursued
especially by Bergman, yielding in particular the concept of
Minimalbereich. See for instance Bergman 1922 \cite{Bergman_1922},
Bergman 1929 \cite{Bergman_1929-Hermite} where the concept seems
to emerge, yet no precise definition. As noted in Maschler's
papers e.g. 1959 \cite{Maschler_1959} it seems that the nature of
those minimal-domains was not completely elucidated in the late
1950's. However, Maschler---extending a result of Schiffer 1938
\cite{Schiffer_1938-CRAS-domaines-minima}---observes that such
minimal domains satisfy the mean property. Therefore on applying
the result of Davis (as quoted in Aharonov-Shapiro 1976
\cite{Aharonov-Shapiro_1976}) characterizing the circle as the
unique domain with a one-point quadrature identity (i.e. such that
the mean value property holds for all harmonic functions) one may
hope to infer our desideratum that the least area map has a range
which is a disc.

Another problem is that of the ``least momentum'' where one
minimizes instead the integral $\int\int_B \vert f(z) \vert^2
d\omega$ (notice the suppression of the derivative) and again to
avoid the trivial solution $f=0$ we impose $f'(t)=1$ at some point
$t\in B$ of the domain. Another possible normalization is to ask
$f(t)=1$, like in Fuchs 1945 \cite{Fuchs_1945/46}. Here again it
seems reasonable to expect circularity of the range of the minimum
mapping. The intuition being that the inertia-momentum of a
rotating body gets minimized for a circular body (granting some
atomical resistance avoiding a complete gravitational collapse of
matter).

[07.08.12] Of course all those problems are super-classical, yet
we still find it hard to delineate the relevant clear-cut results
among the super-massive literature. Our naive intuition would be
that such least-area (or momentum) map are closely allied to
circle maps. However it is not sure that this is the pure truth
for non-simply-connected domains (and a fortiori for bordered
surfaces). As we already said the relevant sources includes for
the area problem:

$\bullet$ Grunsky 1932 \cite[p.\,140]{Grunsky_1932}, alas no
details, some more details in Garabedian-Schiffer 1949
\cite{Garabedian-Schiffer_1949} (but no assertion of circularity)
only the Grunsky formula expressing the least-area map as
combination of the two slit-maps.

$\bullet$ for the least momentum see many works of Bergman
starting from his Thesis 1922 \cite{Bergman_1922}.

Perhaps it should be observed that the least-momentum problem is
perhaps somewhat less easily extensible to Riemann surfaces in
view of the lack of conformal invariance of its functional.

Finally, we can mention Walsh's 1935 survey (M\'emorial
\cite{Walsh_1935}) where all such problems are united under a
generalized form where more points $z_1, \dots z_n$ are prescribed
in the domain joint with some prescribed values $\gamma_1, \dots
\gamma_n$ and one is required to find the map minimizing the
functional under the interpolating condition $f(z_i)=\gamma_i$.
Alas, in Walsh's survey attention is confined to the
simply-connected case and the multi-connected variants where at
that time not systematically understood.

\subsection{A digression about Nehari's paper of 1955}\label{Nehari-digression:sec}

In Nehari 1955 \cite{Nehari_1955}, the author presents a nice
application of Bieberbach's 1925 \cite{Bieberbach_1925} existence
theorem of a circle map for an $n$-ply connected domain upon the
disc of degree $n$. Precisely Nehari deduces a bound on the number
of linearly independent solutions to a certain extremal problem
(akin to those treated by Szeg\"o 1921 \cite{Szego_1921}). It
seems plausible that this Nehari argument is sufficiently
universal to extend directly to the more general setting of
compact bordered Riemann surfaces (membranes for short) upon
invoking Ahlfors 1950 \cite{Ahlfors_1950} instead of Bieberbach
1925 \cite{Bieberbach_1925}. As the argument uses only the circle
mapping nature of the Ahlfors map, we may even appeal to Gabard
2006 \cite{Gabard_2006} to obtain a sharper bound. In reality what
is truly relevant is the absolute invariant of the (separating)
gonality \`a la Coppens 2011 \cite{Coppens_2011}. Let us try to
explore this connection, albeit some details require to be better
worked out in order to really understand this technique of Nehari.

We try first to go quickly to the hearth of Nehari's ideas.
The starting point is the following extremal problem
formulated for $D$ a compact domain bounded by $n$ analytic
curves (for simplicity) forming its complete boundary contour
$C$. Further in the interior of $D$ a set $C_1$ consisting of
a finite number of rectifiable Jordan arc and/or curves is
given. [Warning: in his paper \cite[p.\,29]{Nehari_1955}
Nehari writes ``$C_1$ will stand for a subset of $C$'', which
in our opinion is just a misprint! $C$ should be $D$!? Of
course, our domain $D$ differs from Nehari's as ours includes
the contours.] Let also $L^2=L^2(D)$ be the (Hilbert) space of
analytic functions on $D$ with finite integral $\int_C \vert
f(z) \vert^2 ds<\infty$ where $ds$ is the (Euclidean) length
element.

\medskip
{\bf Problem (P)}. Find the functions $f\in L^2$ minimizing
the norm $\int_C \vert f(z) \vert^2 ds $ under the constraint
$\int_{C_1} \vert f(z) \vert^2 ds=1$.
\medskip

This problem suggests looking at the functional
$$
J(f)=\frac{\int_C \vert f(z) \vert^2 ds}{\int_{C_1} \vert f(z)
\vert^2 ds}
$$
whose minimizers are (up to scaling) the solution of problem
(P).

Next Nehari sets up a certain integral equation whose eigenspace
attached to the lowest eigenvalue parametrize the extremals of
(P). We skip the details, but the key issue is just the linearity
of the set of solutions to Problem (P). With this at hand, we
can plunge directly to the core of Nehari's argument, namely the:

\begin{prop} {\rm (Nehari 1955 \cite[p.\,36]{Nehari_1955})}
Assuming (as above) the domain $D$ of connectivity $n$
(=number of contours), problem {\rm (P)} admits at most $n$
linearly independent solutions.
\end{prop}

\begin{proof}  Nehari's argument splits in 4 short steps:

{\bf Step 1 (Bieberbach 1925)} According to the latter
(\cite{Bieberbach_1925}) there is a circle map $f\colon D \to
\overline{\Delta}=\{\vert z\vert \le 1\} $ of degree $n$. This
means that $\vert f(z) \vert=1$ exactly on the contours (i.e.
$f^{-1}(\partial\overline{\Delta}=S^1)=C$) and upon changing the
origin to an unramified place we may assume that $f$ has exactly
$n$ zeroes, say $z_1,\dots,z_n$.

{\bf Step 2 (Nehari's trick in linear algebra)} Assume by
contradiction that (P) has $n+1$ linearly independent
solutions $f_i$ ($i=1,\dots,n+1$). We consider the linear map
$$
{\Bbb C}^{n+1} \to L^2 \to {\Bbb C}^n\,,
$$
where the first arrow maps $(A_1,\dots, A_{n+1})\mapsto
\sum_{i=1}^{n+1} A_i f_i$ and the second is the evaluation $
\varphi\mapsto (\varphi(z_1),\dots, \varphi(z_n))$ at the
zeroes of the (Bieberbach) function $f$. For dimensionality
reasons, there is a non-zero vector $(A_i)$ in the kernel
which creates the function $f_0:=\sum_{i=1}^{n+1}A_i f_i$
vanishing at all $z_i$, yet without being identically $0$ (the
$f_i$ being linearly independent).

{\bf Step 3 (Nehari factorizes)} The function $g$ defined by $g
\cdot f=f_0$ is regular in $D$ (since writing $g=f_0/f$ we see
that the zeroes of $f$ are cancelled out by those of $f_0$ which
by construction englobe those of $f$). Now using the property of
the circle map $f$ we find the following strict inequality
$$
J(f_0)=\frac{\int_C \vert f_0(z) \vert^2 ds}{\int_{C_1} \vert
f_0(z) \vert^2 ds}=\frac{\int_C \vert g(z)\vert^2 \overbrace{
\vert f(z) \vert^2}^{=1} ds}{\int_{C_1} \vert g(z) \vert^2
\underbrace{\vert f(z) \vert^2}_{< 1} ds}>\frac{\int_C \vert
g(z) \vert^2 ds}{\int_{C_1} \vert g(z) \vert^2 ds}=J(g)\,.
$$
(Moreover reading backwards the numerators we see that the
norm of $g$ equals that of $f_0$ so that $g\in L^2$.) The just
obtained inequation $J(g)<J(f_0)$ shows that $f_0$ fails to
solve (P).

{\bf Step 4 (Using the linear structure)} However the $f_i$
($i=1,\dots,n+1$) solve (P), hence by virtue of the linear
structure of the extremals to (P) [which Nehari derives from
an interpretation as the eigenspace attached to the lowest
eigenvalue, but which perhaps may be derived more directly] it
follows that $f_0$ solves also (P) [after scaling
appropriately], violating the conclusion of Step~3.
\end{proof}

Albeit our presentation is not completely polished (and
Nehari's maybe not perfectly organized for the beginner), we
see that the basic trick looks sufficiently universal, as to
extend to the following context.

Instead of the finitely-connected domain $D$, we consider $F$ a
compact bordered orientable Riemannian surface of genus $p$ and
with $r$ contours. Now $ds$ denotes the induced length element
attached to the (Riemannian) metric. As above, we specify a subset
$C_1$ of the interior of $F$ consisting of a finite ``drawing'' of
Jordan arcs and curves (perhaps they do not even need to be
pairwise disjoint). Then we set up the extremal problem (P) in
this context, and they above proof seems to work mutatis mutandis,
except for trading Bieberbach 1925 \cite{Bieberbach_1925} by
Ahlfors 1950 \cite{Ahlfors_1950} or Gabard 2006
\cite{Gabard_2006}. Precisely, we may consider a circle map
$f\colon F \to \overline{\Delta}$ of least possible degree, say
$\gamma$. By Gabard 2006 \cite{Gabard_2006} we know that
$\gamma\le r+p$. So we arrive at the following statement:

\begin{prop}
Let $F$ be a membrane of genus $p$ with $r$ contours. Assume that
$F$ has the gonality $\gamma$, i.e. the least degree of a circle
map to the disc. (We know $\gamma\le r+p$) Then the extremal
problem {\rm (P)} admits at most $\gamma$ linearly independent
solutions.
\end{prop}

\section{Ahlfors' extremal problem}

\subsection{Ahlfors extremal problem (Grunsky 1940--42, Ahlfors 1947--50)}

Ahlfors' method involves solving the following extremal
problem:

\begin{theorem}\label{Ahlfors-extremal:thm} {\rm (Ahlfors 1950
 \cite{Ahlfors_1950})} Given any compact bordered Riemann
surface (membrane for short) and two interior points $a,b$,
find among all (analytic)
functions bounded-by-one
taking $a$ to $0$ the
one maximizing the modulus $\vert f(b) \vert$.

Such a function exists (normal families argument \`a la
Vitali-Montel) and is unique up to a rotation (=multilication by
an unimodular complex number $\omega=e^{i \theta}$). Hence it is
unambiguously defined by the points $a,b$ if $f(b)$ is required to
be positive real, and we denote $f_{a,b}$ the corresponding
function.

Furthermore Ahlfors' extremal function $f_{a,b}$ concretizes
the given surface as a full-covering of the disc $\Delta$, of
degree
\begin{equation}
r\le \deg f_{a,b} \le r+2p, \label{Ahlfors:pinch}
\end{equation}
where $r$
is the number of contours and $p$ the genus (of the given
membrane).
\end{theorem}

It is nowadays quite customary---following (another) Russian
school (Golusin, S.\,Ya. Havinson, etc.)---to call the
extremal an {\/\it Ahlfors function}, albeit even Ahlfors
seems
to have been rather
embarrassed by this probably unearned distinction (cf. his
comments in Collected Papers
\cite[p.\,438]{Ahlfors_1982_Coll_papers}). The same idea occurred
somewhat earlier in works of Grunsky 1940--42 \cite{Grunsky_1940},
\cite{Grunsky_1942}, yet the latter confined attention to plane
domains (as did Ahlfors 1947 \cite{Ahlfors_1947}).
Being close colleagues---as materialized by their
joint note (Ahlfors-Grunsky 1937 \cite{Ahlfors-Grunsky_1937})
about the best conjectural value for the {\it Bloch constant}
(still open up to present days)---it is puzzling that both were
not very aware of overlapping studies (admittedly imputable to the
difficult World War II  context).

\subsection{Semi-fictional reconstruction of Ahlfors' background
(Fej\'er-Riesz 1922, Cara\-th\'eodory 1928, Ostrowski 1929)}


Where does Ahlfors' extremal problem come from? This is surely a
non-trivial question  yet let us attempt to give some elements of
answers. The narrative is made more plausible by looking a bit
around while  trying to keep track of the historical continuity.
We shall  thus use several
indirect sources, especially Remmert.

As notorious, the Dirichlet principle suffered ill-foundations
during a long period of about 40 years (1860-1900). This was
beneficial to Schwarz-Christoffel who developed some constructive
methods for the RMT for polygons. Another trend involves directly
rescuing the Dirichlet principle via the ``alternierendes
Verfahren'' of Schwarz and the parallel work of C. Neumann. This
influenced Picard's {\it m\'ethodes des approximations
successives}, as well as Poincar\'e's balayage.

Then came Hilbert's breakthrough. Yet, alternative methods
circumventing the intricacies of potential theory seemed worth
attention. As reported in Remmert 1991 \cite{Remmert_1991}, one
can ascribe to Fej\'er-Riesz ca. 1921 (published by Rad\'o 1923
\cite{Rado_1923-Uber-konf-Abb}) the first purely complex variables
(potential-theoretic free) proof of the RMT by using the extremal
problem of making the modulus of the derivative as large as it can
be. Several technical simplifications were then obtained by
Carath\'eodory 1928 \cite{Caratheodory_1928} and Ostrowski 1929
\cite{Ostrowski_1929} (independently). This leads in principle to
the most elementary proof of the RMT. Extending this idea to
multiply-connected domains (say first of finite connectivity)
leads directly to the extremal problem considered by Grunsky
1940--42 \cite{Grunsky_1940}, \cite{Grunsky_1942}, and Ahlfors
1947 \cite{Ahlfors_1947}, and Ahlfors 1950 \cite{Ahlfors_1950}
when extended to Riemann surfaces.

In fact prior to Fej\'er-Riesz, it is fair to refer to Koebe's
(and Carath\'eodory's) elementary proofs of the RMT, also via an
extremal problem or at least iterative methods (compare e.g.
Garabedian-Schiffer 1950~\cite{Garabedian-Schiffer_1950}).

As a matter of digression, it can be recalled that this extremal
viewpoint leads as well to a proof of the uniformization theorem
(without potential theory). Compare Carath\'eodory 1950
\cite{Caratheodory_1950}, plus several papers by Grunsky (easily
located in his collected papers).

\subsection{Extremal problems and pure function-theoretic
proofs of the RMT (Koebe, Carath\'eodory, Bieberbach)}

The previous section is a bit caricatural and the real history is
marvellously detailed in Gray 1994 \cite{Gray_1994}. Let us
summarize the chronology of this period, in the center of which
there is probably one of the main inspiring force toward the
Ahlfors extremal function (namely the {\it Schwarz lemma} as
Carath\'eodory christened it in 1912).

$\bullet$ Painlev\'e 1891 \cite{Painleve_1891}: boundary
behavior of the Riemann mapping for a contour having an
everywhere continuously varying tangent.

$\bullet$ Harnack 1887 \cite{Harnack_1887} provides a
satisfactory proof for solving a suitable version of
Dirichlet's principle, and states what has become known as
{\it Harnack's theorem} on monotone limits
of harmonic functions.

$\bullet$ Osgood 1900 \cite{Osgood_1900} applies Harnack's theorem
to draw the existence of a Green's function for any
simply-connected plane domain thereby resolving the Riemann
mapping theorem (RMT). This dependance is eliminated in Koebe 1908
and Carath\'eodory 1912 (cf.\,items below), where  Schwarz's lemma
is substituted.

$\bullet$ Poincar\'e 1907 \cite{Poincare_1907} (and independently
Koebe 1907, cf. below) proves uniformization  (rigourously). For
this Poincar\'e combines his {\it m\'ethode de balayage} (of 1890
\cite{Poincare_1890}) and simplifies it using Harnack's theorem.
From the Green's function he deduces the conformal map of a
Riemann surface (\`a la Weierstrass) to the disc, and uses earlier
works of Osgood.

$\bullet$ Koebe 1907 also proves uniformization (UNI). In Koebe
1907c \cite{Koebe_1907_UbaK2} he compares his method to
Poincar\'e's. Like Poincar\'e he had relied on Schwarz's method,
but unlike him made a much more modest use of Harnack's theorem.
Koebe also insists upon his avoiding of the use of modular
functions.

$\bullet$
Koebe 1908 \cite{Koebe_1908_UbaK3}
supplies another proof (of UNI) avoiding completely Harnack's
theorem. [Subsequently Koebe interacted widely with Fricke's
attempt to modernize the original {\it continuity method} of
Klein-Poincar\'e, and showed how this could be rigorized
overlapping thereby with simultaneous work by Brouwer. This
interaction with Brouwer seems to have ended quite
contentiously.]

$\bullet$ Koebe 1909 \cite{Koebe_1909_UAK1}, 1910
\cite{Koebe_1910_UAK2} proof of his {\it Verzerrungssatz\/}
(distortion theorem).
From it he derives, the first elementary proof of the
(RMT) appealing to a long list series of name going back via
Arzel\`a and Montel 1907 \cite{Montel_1907} to Ascoli 1883.

$\bullet$ Carath\'eodory 1912 \cite[p.\,109]{Caratheodory_1912}
notes that {\it Schwarz's lemma} (which he was the first to call
by this name, and which he locates in Schwarz's Ges. Abh., vol. 2,
p.\,109) can act as a substitute to Harnack's theorem (upon which
Osgood 1900 relied heavily). [Interrupting the present narrative
this will have to play a major r\^ole in Ahlfors' extremal
problem.] Using the Schwarz's lemma and Montel's theorem,
Carath\'eodory obtains the Riemann mapping using an exhaustion of
the domain $G$ by subdomains $(G_n)$ each mapped via $f_n$ to the
disc and studied under which condition on $G_n$ the $f_n$
converges to a function $f$ giving the Riemann mapping (again
without potential theory).


$\bullet$ Carath\'eodory 1913a \cite{Caratheodory_1913a} proves
Osgood's conjecture, that the Riemann map extends to a
homeomorphism of the boundary iff the boundary is a Jordan curve.
In Carath\'eodory's opinion this achievement is mostly a byproduct
of Lebesgue's far-reaching theory of integration (1902
\cite{Lebesgue_1902}),
and the consequences drawn from it by Fatou 1906
\cite{Fatou_1906}. This reliance upon Lebesgue-Fatou was soon
disputed by Koebe 1913 (cf. item below).

$\bullet$ Carath\'eodory 1913b \cite{Caratheodory_1913b} discusses
the boundary behavior when the boundary curve is not a Jordan
curve. This paper is oft regarded as inaugurating the concept of
{\it prime ends} (although earlier origins are in the work of
Osgood, and related ideas in Study-Blaschke 1912
\cite{Study-Blaschke_1912}).

$\bullet$ Koebe 1913 \cite{Koebe_1913} disputes the need for
Lebesgue's theory in Carath\'eodory's treatment, showing how
to generalize a theorem of Schwarz to the same effect. A
similar result is claimed independently by Osgood-Taylor 1913
\cite{Osgood-Taylor_1913}.

$\bullet$ Bieberbach 1913 \cite{Bieberbach_1913} wrote a short
paper disputing the (in his opinion) excessive
Carath\'eodory's reliance on Schwarz's lemma, proposing to use
only Montel's theorem. The next year Bieberbach 1914
\cite{Bieberbach_1914} invokes another
extremum principle (area minimization of the  range of the mapping
suitably normalized) to simplify Carath\'eodory's work. This
freed the theory from any reliance upon Montel's theorem (but
uses
instead ideas of Ritz).

$\bullet$ Back to Koebe, in 1912 \cite{Koebe_1912} could not
resist after the stimulus aroused by Carath\'eodory's work to go
back to some old idea of his own ({\it Quadratwurzeloperationen})
to create his {\it Schmiegungsverfahren} (squeezing methods) for
solving the Riemann mapping by the iterated taking of square
roots. This presentation was entirely elementary.

$\bullet$ Carath\'eodory 1914 \cite{Caratheodory_1914}
incorporated all these criticisms in his paper for the Schwarz
Festschrift, which was to remain his final account until the newer
methods of Perron were introduced. [Here we may have also
mentioned the argument of Fej\'er-Riesz 1921.]

$\bullet$ Bieberbach 1915 in his pocket book G\"oschen
\cite[p.\,95]{Bieberbach_1915} also proposes to deal entirely
within pure function theory, while rejecting the
potential-theoretic approach (despite Hilbert's work). This
actually presents a version of Koebe's Schmiegungsverfahren and
concludes to the Riemann mapping theorem via Koebe's
Verzerrungssatz (seen as a preferred alternative over Schwarz's
lemma).

\subsection{Interlude: Das Werk Paul Koebes}

In this section we digress slightly from our main path to look
closer at the monumental works of Koebe.
A useful guide is Bieberbach's overview of Koebe's work in 1968
\cite{Bieberbach_1968-Das-Werk-Paul-Koebes}. The main point of
overlap of Koebe with our main theme (Ahlfors) lies in the
Riemann-Schottky mapping (albeit for Koebe the mapping to a
Kreisbereich is given full attention neglecting thereby the circle
mapping). Of course, the other main
aspect of Koebe's
life is the uniformization theorem of
(Klein-Poincar\'e-Schwarz).

Again some chronology:

$\bullet$ Riemann 1857--58 \cite{Riemann_1857_Nachlass} and
Schottky 1877 \cite{Schottky_1877} (maybe only in the 1875 Latin
version?) proved that any $n$-ply connected domain maps
conformally to a Kreisbereich (circular domain). [Bieberbach and
indeed Koebe 1910 \cite{Koebe_1910_JDMV} ascribe this to Riemann,
albeit we are not sure to be in
total agreement with this assertion.]

$\bullet$ In Bieberbach's opinion the above Riemann-Schottky
Kreisbereich-mapping is first rigourously proved by Koebe in a
series of four papers written in 1906, 1907, 1910, 1920 (which
we attempt to summarize in more details):

(1)
Koebe 1906 \cite{Koebe_1906_JDMV}: this starts with a rigidity
result for two Kreisbereiche
as being conformal to each other only through linear
transformations. The proof uses potential theory (and the
Cauchy integral).
%
%
It follows that $(\varrho+1)$-ply-connected Kreisbereiche
depend upon $3\varrho -3$ essential constants when $\varrho
\ge 2$,  the same quantity as predicted by Schottky for
general multiply-connected domains of the same connectivity.
This yields some evidence for the possibility of mapping those
to a Kreisbereich. Actually Koebe (p.\,150) reminds that
the Kreisbereich mapping is (essentially) solved by Schottky
and by Poincar\'e
(referring loosely to the first volumes of Acta). [In the next
paper Koebe
adopts a more critical position, and does not take this as
granted.] Next, he claims the result extends to schlichtartig
surfaces. His argument amounts to fill  the Riemann surface by
discs, to get a closed surface of genus 0, and appeal to
Schwarz
1870 \cite{Schwarz_1870} to map this to a sphere. Next, Koebe
proposes to relax the schlichtartig character
to formulate a similar result for positive genus. Again one
fills the surface by discs to
gain a closed surface of genus $p$. This can be mapped as a
ramified cover of the sphere of degree $p+1$ (as well-known since
Riemann, but for Koebe being Schwarz's pupil Riemann is taboo and
an ad hoc [somewhat sketchy] argument is supplied). At any rate
the result is that any compact bordered Riemann surface of genus
$p$ is conformally embeddable in a closed Riemann surface of genus
$p$, hence representable as a $(p+1)$-sheeted cover of the sphere.
Although this result concerns like Ahlfors 1950
\cite{Ahlfors_1950} compact bordered surfaces, it seems that this
Koebe mapping lies not so deep as the image of the contours of the
map are poorly controlled, in particular they need not coincide.

(2)
Koebe 1907 \cite{Koebe_1907_UrAK}: this starts by quoting again
his rigidity result of the previous paper. Then more critically
Koebe notices that the mapping of a planar $(\varrho+1)$-ply
connected domain upon a Kreisbereich of the same connectivity is
not so easily established, making abstraction of the
Klein-Poincar\'e {\it Kontinuit\"atsmethode} ({\it m\'ethode de
continuit\'e}) not yet effective in 1907. (This had to wait until
the work of Brouwer and Koebe ca. 1911
\cite{Klein-Brouwer-Koebe_1912},
and Koebe 1912 \cite{Koebe_1912_BdKm}.) The rigidity result
affords an essentially unique solution of the mapping problem.
Then Koebe proceeds to show that a
Kreisbereich mapping exists for triply-connected domains
($\varrho=2$), and generally if the domain is symmetric under
complex conjugation provided the real axis cuts all contours. For
triply connected domains, he takes the Schottky double, which
conformally maps to a closed Riemann surface of genus $2$ (via
massive quotations to Schwarz, Ges. Abh. II, S.\,133--143,
S.\,144--171, S.\,175--210). As any curve of genus 2 this is
hyperelliptic (canonical mapping via holomorphic 1-forms). As to
the more general case, the problem involves cutting the domain
along the real axis, yielding a simply-connected region. This is
mapped conformally to the upper half-plane, and symmetrically
reproduced. Then Green's function is constructed via Harnack's
theorem (quotation to Harnack 1887 \cite{Harnack_1887}, Poincar\'e
1883 \cite{Poincare_1883}, Osgood 1900 \cite{Osgood_1900} and
Johansson 1905).

(3)
Koebe 1910 \cite{Koebe_1910_JDMV}: the paper starts again with
the objective to solve the {\it Problem der konformen
Abbildung eines $(p+1)$-fach zusammenh\"angenden Bereiches auf
einen von $p+1$ Vollkreisen begrenzten Bereich} (which he
proposes to call {\it Kreisbereich} for short). Koebe
recalls that the problem was first
addressed by Schottky 1877 in his {\it Doktordissertation},
and earlier in Riemann's Nachlass. He reminds from his first
work [item (1)] that {\it je zwei Kreisbereiche aufeinander
nur durch lineare Funktionen konform abgebildet werden
k\"onnen}. Then he repeats the two special cases he was able
to solve previously, and now proposes to tackle the general
case via two different methods (of his own): {\it
\"Uberlagerungsfl\"ache} and {\it iterierendes Verfahren} [cf.
items (A) and (B) below]. He proudly
emphasizes that both methods have a larger applicability than
to the present
Kreisbereich problem, since their combination, allowed him to
settle the whole series of classical mapping problems
of Klein and Poincar\'e (1881--84) in their pioneering works on
automorphic functions, and the allied uniformization. Hilbert's
22th Problem (1900) is mentioned for reposing the uniformization
question especially in connection to Poincar\'e 1883's paper
\cite{Poincare_1883}. Schwarz is again (justly) regarded as the
father of the method  {\it der \"Uberlagerungsfl\"ache}, which
plays a key r\^ole in the newer developments in the automorphic
theory,
as exemplified through the
work of Poincar\'e 1907
\cite{Poincare_1907} himself and Hilbert 1909
\cite{Hilbert_1909}. After these general remarks Koebe
proceeds to prove the general Kreisbereich mapping. [As warned
in Bieberbach's report, the present paper of Koebe does not
contain full details, yet some lovely geometric ideas worth
sketching. Complete details appear in the last contribution item
(4), but then it is easy to get lost in  technicalities.]

(A) Koebe assumes the contours of the domain $B$ to be analytic
curves. Via some abstract {\it Spiegelungsprozesses} (ascribed to
Schwarz) he constructs via symmetric reproduction of $B$ a
schlichtartig Riemann surface $B^{(\infty)}$. (One must imagine
$B$ glued with replicas thought of as the back-side of the
domain.) Then he can apply his {\it allgemeines Abbildungsprinzip}
to the effect that schlichartig implies schlicht (first
established in Koebe 1908 \cite{Koebe_1908_UbaK3}, with
subsequent approaches by Hilbert 1909 \cite{Hilbert_1909} and in
Courant's Thesis 1910/12 \cite{Courant_1912}).
The new schlicht domain $B^{(\infty)'}$ is tesselated by replicas
of the conformal copy $B'$ of $B$. Hence $B'$ admits a complete
infinite system of symmetric reproduction. This is enough (for
Koebe) to characterize a Kreisbereich. (Here we may agree with
Bieberbach's diagnostic that Koebe's exposition is sketchy, but
details were supplied later in Koebe 1920 \cite{Koebe_1920}.)

\begin{figure}[h]
\centering
    \epsfig{figure=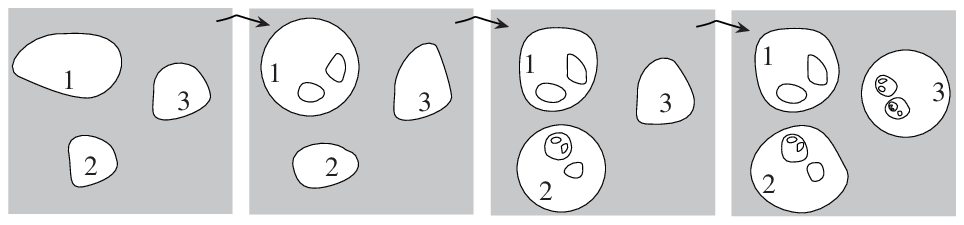,width=122mm}

\caption{\label{KoebeiV:fig} Koebe's {\it iterierendes Verfahren}:
successive circularization of the contours of a multiply-connected
domain via the Riemann mapping and a magical convergence to a
Kreisbereich (first established by Koebe on the basis of his {\it
Verzerrungssatz\/} in 1908 \cite{Koebe_1908_UbaK3})}
\end{figure}

(B) Then is exposed the promised {\it iterierendes Verfahren}.
This is a beautiful device based upon successive applications of
the RMT to circularize a specific contour and then reflecting by a
{\it Spiegelgung} (inversion by reciprocal radii) the domain
across this circularized contour. Koebe draws nice pictures (like
below Fig.\,\ref{KoebeiV:fig})
suggesting  that this iteration scheme produces domains with {\it
sukzessive Steigerung der Spiegelungsf\"ahigkeit des Bereichs}
whereupon it is made plausible that when repeated ad infinitum the
resulting domain has an infinite aptitude of symmetric
reproduction, hence must be a Kreisbereich. The convergence proof
uses his {\it Verzerrungssatz} (distortion theorem).

(4)
Koebe 1920 \cite{Koebe_1920}, where
full details are supplied.

\medskip
$\bullet$ In parallel, Koebe concentrates his efforts on the
uniformization problem starting with Koebe 1907
\cite{Koebe_1907_UrAK} devoted to the uniformization of real
algebraic curves, yet the real technological breakthrough
occurs in the next paper.

\begin{figure}[h]
\hskip-35pt \penalty0
\epsfig{figure=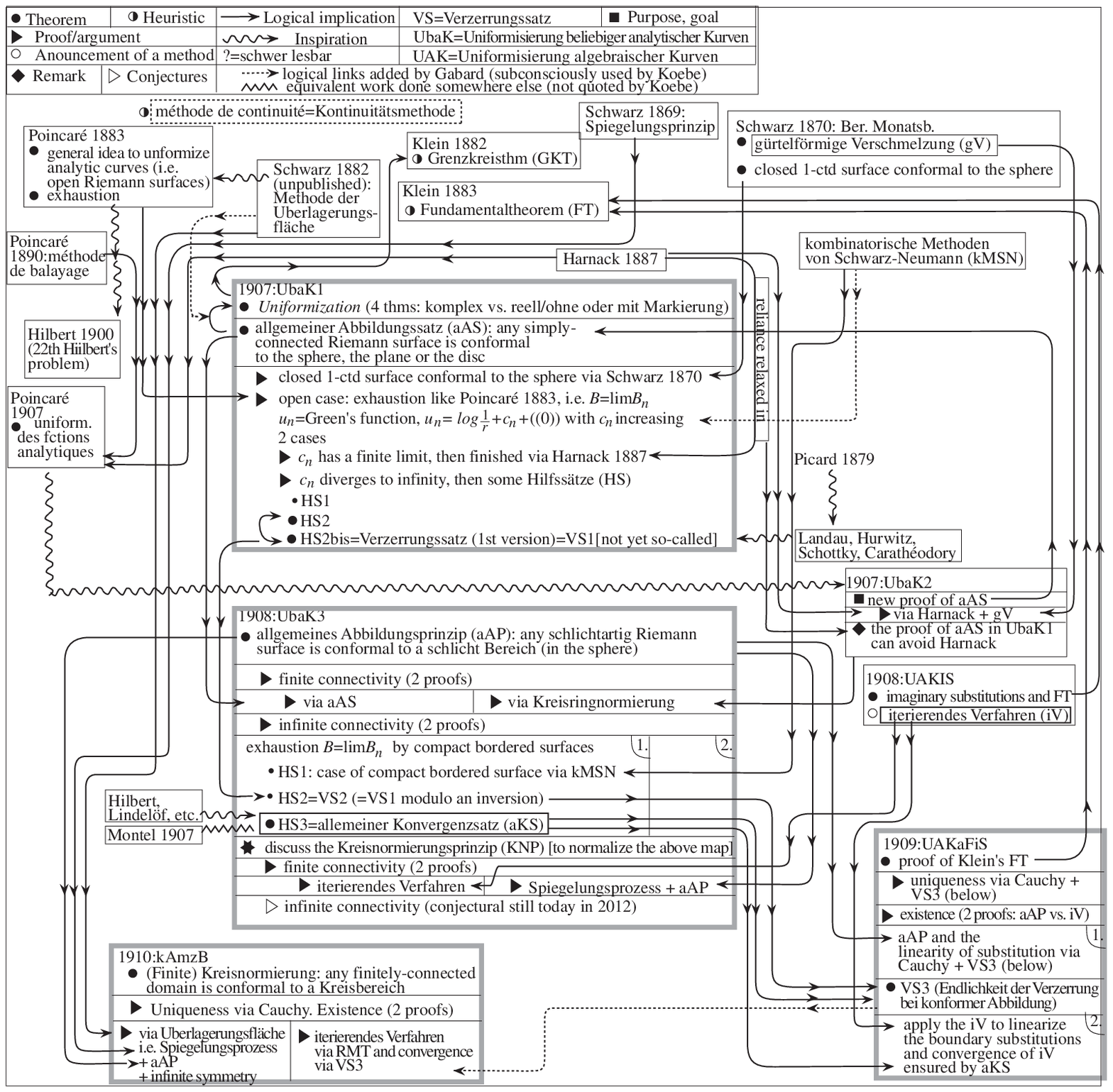,width=152mm}
%
\caption{\label{KoebeMa:fig} Logical dependance between
Koebe's early theorems}
\end{figure}

$\bullet$ Koebe 1907 \cite{Koebe_1907_UbaK1} discovers a first
version of his {\it Verzerrungssatz} (VZS), which turns out to be
relevant both to the Riemann-Schottky Kreisbereich-mapping, as to
uniformization. As forerunners of the (VZS) Bieberbach mentions
the works of Landau, Schottky related to Picard's theorem (1879
\cite{Picard_1879}). This Koebe's paper also contains (what later
came to be known) as the {\it Viertelsatz} to the effect that the
range of any schlicht function on the unit disc normalized by
$f(0)=0$ and $\vert f'(0) \vert=1$ contains a disc of some
universal positive radius $\varrho$. The sharp value $\varrho=1/4$
is conjectured, but only established by Bieberbach 1915
\cite{Bieberbach_1915}. Armed with this
Verzerrungssatz (yet without the precise bound) Koebe manages to
prove uniformization. This represents a generalization of the RMT
to simply-connected Riemann surfaces. Bieberbach recalls that
according to oral tradition the trick of the universal covering
surface is due to H.\,A. Schwarz (ca. 11. April 1882, as carefully
reported in Klein's Werke \cite[p.\,584]{Klein-Werke-III_1923}).

$\bullet$ Simultaneously and independently Poincar\'e 1907
\cite{Poincare_1907} also proves the uniformization theorem
via his {\it m\'ethode de balayage}.

$\bullet$ Koebe 1907 \cite{Koebe_1907_UbaK2} inspects
Poincar\'e's proof and proposes a variant using Harnack's
theorem (in potential theory) circumventing thereby the
Viertelsatz, as well as Poincar\'e's balayage.

$\bullet$ The new ingredient (Verzerrungssatz of Koebe) turned out
to act  usefully in other uniformization problems envisioned by
Klein (e.g., the {\it R\"uckkehrschnitt\-theorem}, etc.) In
Koebe's formulation this resulted to the  conformal mapping of a
schlichtartig Riemann surface to a schlicht domain of the Riemann
sphere. This result appears in Koebe 1908 \cite{Koebe_1908_UbaK3}.
Its proof uses beside the Verzerrungssatz a general convergence
theorem (\`a la Montel-Vitali), which Koebe discovered
independently [according to Bieberbach].

$\bullet$ Koebe 1909 \cite{Koebe_1909_UAK1} gives a sharper
version of the {\it Verzerrungssatz} and applications to Klein's
general uniformization problem (via groups of linear
transformations).

$\bullet$ Hilbert 1909 \cite{Hilbert_1909}, using a variant of the
Dirichlet principle, gives another method for the schlicht mapping
of a schlichtartig surface (to the sphere), via a so-called
parallel-slit mapping [extending the Schottky-Cecioni result to
infinite connectivity].

$\bullet$ In response Koebe 1909 \cite{Koebe_1909_UbaK4}, 1910
\cite{Koebe_1910_Hilbert} and independently Courant 1910/12
\cite{Courant_1912} proves anew the above Hilbert's Ansatz about
parallel-slit mappings.

$\bullet$ Already Schottky 1877 \cite{Schottky_1877} tried [in
Bieberbach's opinion] to prove the [Riemannian] theorem that every
$n$-ply connected planar domain conformal-maps bijectively to a
parallel Schlitzbereich. Hilbert's new method proves this for
arbitrary schlichartig Riemann surfaces. Koebe in the
aforementioned two works, sharpens Hilbert's theorem by noticing
that the range of the mapping fill the full plane save a set of
measure zero. At this occasion Koebe also formulates his {\it
Kreisnormierungsprinzip}, still open today, despite the
spectacular progress by He-Schramm 1993 \cite{He-Schramm_1993}.

$\bullet$ Bieberbach emphasizes that the {\it iterierendes
Verfahren\/} may really have first emerged through the
Kreisbereich mapping problem. [This conflicts slightly with
Koebe's claim that he employed it earlier for uniformization.] At
any rate Bieberbach writes ``{\it Solche iterierenden Verfahren
entwickelt Koebe \"uber Jahr\-zehnte hin immer weiter, bis alle
Uniformisierungsprobleme algebraischer Gebilde dem iterierenden
Verfahren zug\"anglich werden.}''

$\bullet$ The proof of the (RMT) via repeated {\it
Quadratwurzelabbildungen}  itself constitutes an iterative method,
which Koebe calls the {\it Schmiegungsverfahren}. Credit for this
discovery is to be shared with Carath\'eodory.

$\bullet$ A rigorous foundation to the {\it Kontinuit\"atsmethode}
of Klein-Poincar\'e is paid by Koebe much attention in a
torrential series of paper starting with 1912
\cite{Koebe_1912_BdKm}, 1912 \cite{Koebe_1912_BdKm2}, 1914
\cite{Koebe_1914_UAK4}, etc. Those works overlaps (and then may
supplement) the works of Brouwer on the invariance of domain (and
dimension), and its application to Riemann surfaces. The resulting
priority question is very intricate. Even Klein in 1923
\cite[p.\,734]{Klein-Werke-III_1923} writes: {\it Die
entscheidende Wendung trat aber erst 1911/12 durch das Einsetzen
der Untersuchungen von Brouwer und Koebe ein. (Ich halte um so
mehr an der alphabetischen Reihenfolge fest, als die gegenseitige
Beziehung der beiden Forscher nicht ganz gekl\"art ist.)} Soon
afterwards Klein also cites footnote 2) in Brouwer 1919
\cite{Brouwer_1919}, where Brouwer seems to revendicate some
priority over Koebe, while reporting some falsification of his own
(G\"ott. Nachr.) article via a citation to Koebe added after
proof-reading.

\subsection{Koebe and
his relation to Klein or Ahlfors}

In the overall Koebe's monumental work is quite intricate with
deep influences by methods of Schwarz (ca. 1870), results of
Schottky (1875/77), visions of Klein and Poincar\'e (early 80's),
supplemented by methods of his own. The following chart
(Fig.\,\ref{KoebeMap:fig}) gives an \"Uberblick maybe helping
navigation through Koebe's works and the logical links between his
results.

\begin{figure}[h]
\hskip-85pt \penalty0
    \epsfig{figure=KoebeMap.eps,width=182mm}

\caption{\label{KoebeMap:fig} Logical dependance between
Koebe's early theorems}
\end{figure}

From our Ahlfors' biased viewpoint  several points are worth
noticing:

(1) Koebe frequently refers to Klein's orthosymmetry for real
algebraic curves. In view of the close connection between
orthosymmetry and the Ahlfors circle mapping, it is tempting to
wonder if Koebe was ever close to discover the Ahlfors circle
mapping. Of course Koebe's focus seems to have been more attracted
by the uniformization problem (in particular for real algebraic
curves), cf. Koebe 1907 \cite{Koebe_1907_UrAK}. However Klein's
orthosymmetry appears in many subsequent papers (e.g., 1919
\cite[p.\,29, p.\,35]{Koebe_1919:47}), and we would not bet that
one day someone discovers in Koebe some anticipation of the
Ahlfors map (as it occurred say with the circles packing of
Andreev--Thurston). If not directly, it could
via the R\"uckkehrschnitttheorem of Klein (cf.
Sec.\,\ref{sec:Ruckkehrschnittthm}), which Koebe was the first to
prove seriously (cf. Koebe 1910 UAK2 \cite{Koebe_1910_UAK2}).
Hence schematically, there might exist a
(harsh style) path like:
$$
\textrm{Koebe}\Rightarrow\textrm{Klein} \Rightarrow
\textrm{Teichm\"uller} \Rightarrow \textrm{Ahlfors}.
$$

(2) Koebe also notices at several places (e.g., 1907 UbaK1
\cite[p.\,199]{Koebe_1907_UbaK1}) that the orthosymmetry concept
for real algebraic curves extend to analytic real curves. One can
then wonder if there is likewise a function theoretical
characterization of orthosymmetry in terms of (totally real)
mapping to the sphere. This would amount to say that any bordered
surface is expressible as a total cover of the disc (taking
boundary to boundary). Of course this might be a bit
fantasist, but perhaps deserves to be analyzed more carefully.
(Maybe this fails already for planar domains, cf. a work of Heins
ca. 1954.)

\subsection{Ahlfors' background (Bergman 1941, Schiffer 1946,
Schottky differentials)}

Let us quote the introduction of Ahlfors 1950 \cite{Ahlfors_1950}:


\begin{quota}[Ahlfors 1950]\label{Ahlfors-1950:quote}

{\small \rm In the handling of the extremal problems we are in
close contact with the methods of Bergman
[1941]=\cite{Bergman_1941} and Schiffer
[1946]=\cite{Schiffer_1946}, which they have developed for plane
regions. A convenient tool for applying these methods to regions
on Riemann surfaces is found in the class of Schottky
differentials, and it was the recognition that Bergman's
kernel-functions are in fact Schottky differentials that led us to
undertake this study.

The second part of the paper (\S\S\,4--5) deals with an
extremal problem that we have previously solved for plane
regions. There are great simplification over my original proof
for which I am partly indebted to my student P. Garabedian. An
interesting point is that the extremal functions are again
defined by means of Schottky differentials.}

\end{quota}


As a complement, we may reproduce a passage of Ahlfors'
comments in his collected papers
\cite[p.\,438]{Ahlfors_1982_Coll_papers}:


\begin{quota}[Ahlfors 1982]\label{Ahlfors-1982:quote}

{\small \rm The purpose of [36](=Ahlfors 1950 \cite{Ahlfors_1950})
was to study open Riemann surfaces by solving extremal problems on
compact subregions and passing to the limit as the subregions
expand. The paper emphasizes the use of harmonic and analytic
differentials in the language of differential forms. It is closely
related to [35](=Ahlfors-Beurling 1950
\cite{Ahlfors-Beurling_1950}), but differs in two respects: (1) It
deals with Riemann surfaces rather than plane regions and (2) the
differentials play a greater r\^ole than the functions.

I regard [36] as one of my major papers. It was partly inspired by
R. Nevanlinna, who together with P.\,J. Myrberg had initiated the
classification theory of open Riemann surfaces, and partly by M.
Schiffer (1943) and S. Bergman (1950), with whose work I had
become acquainted shortly after the war\footnote{Notice that
Ahlfors never quote Teichm\"uller 1981 \cite{Teichmueller_1941},
except in Ahlfors-Sario 1960 \cite{Ahlfors-Sario_1960}, where also
all the Italian works of Matildi 1945/48 \cite{Matildi_1945/48}
and Andreotti 1950 \cite{Andreotti_1950} are cited.}. The paper
also paved the way for my book on Riemann surfaces with L. Sario
[1960], but it is probably more readable because of its more
restricted contents.

I would also like to acknowledge that when writing this paper
I made important use of an observation of P. Garabedian to the
effect that the relevant extremal problems occur in pairs
connected by a sort of duality. This is of course a classical
phenomenon, but in the present connections it was sometimes
not obvious how to formulate the dual problem.

}

\end{quota}


\subsection{The allied infinitesimal form of the extremal
problem}

The input required to pose Ahlfors' extremal problem
(Theorem~\ref{Ahlfors-extremal:thm}), is a membrane with two
interior marked points, denoted $a,b$. When the point $b$
converges to the point $a$ (becoming infinitely close to it), we
may think of a unique point of multiplicity two. This limiting
process mutates the extremum problem into:

\begin{prob}
Let $a$ be a single point in the membrane $W$. Among all functions
$f$ analytic on $W$ with $\vert f \vert \le 1$ on $W$ it is
required to find the one which makes the modulus of the derivative
$f'(a)$ to a maximum. Here the derivative is computed w.r.t. any
holomorphic chart. Its maximum value  has no intrinsic meaning,
yet the extremal function exists and is uniquely defined (up to a
rotation) and denoted by $f_{a,a}=f_a$.
\end{prob}

It seems to be folklore that such functions are also circle maps
subjected to the same Ahlfors bound $\deg f_{a}\le r+2p$.
Presumably a continuity argument reduces to the case of (bipolar)
functions $f_{a,b}$, or maybe adapt the whole argument in  Ahlfors
1950 \cite{Ahlfors_1950}. At any rate, the result is taken for
granted in Yamada 1978 \cite{Yamada_1978}, Gouma 1998
\cite{Gouma_1998}. This can maybe deduced as a special case of
Jenkins-Suita 1979 \cite{Jenkins-Suita_1979}.

\subsection{Higher extremal problems=HEP$\approx$High
energy physics, alias Pick-Nevanlinna interpolation}

What happens if we take more than two points? For instance three
points $a,b,c$? Should we then maximize the area of the simplex
spanned by the image points? If yes for which metric on the disc
(Euclid vs. hyperbolic)? How does the problem reformulate when the
3 points coalesce at the subatomic level into a point affected by
a multiplicity 3. Does the problem amount then to maximize the
modulus of the first two derivatives?

Maybe this brings us in the realm of Pick-Nevanlinna
interpolation, a theory initially developed in the disc. Compare
e.g. Garabedian 1949 \cite{Garabedian_1949}, Heins 1950
\cite{Heins_1950}, Jenkins-Suita 1979 \cite{Jenkins-Suita_1979}.

Perhaps for any (effective) divisor $D=d_1 p_1+\dots +d_n p_n$
interior to the membrane there is an extremal problem denoted
$EP(D)$. Then how much of Ahlfors' theory extends: existence,
uniqueness and qualitative circle mapping nature of the function,
and estimates over the degree of the extremals. In the classic
theory where $\deg(D)=2$ we have $\deg f_{a,b}=r+2p$. Maybe in
general denoting by $f_D$ the extremal function allied to the
divisor $D$ we find $\deg f_{D}\le r+ \deg(D) p$. Compare
Jenkins-Suita 1979 \cite{Jenkins-Suita_1979} for more serious
answers. If we could find a divisor of degree one then this would
recover Gabard's bound $r+p$. Maybe not a divisor is required but
an ordered collections of points, as in Ahlfors' original problem
where $a$ seems to have a preferred r\^ole over $b$, getting
mapped to zero.

Such higher extremal problems depending upon a higher number of
free parameters are probably more flexible in the sense that if
the original $\deg(D)=2$ case of Ahlfors fails to realize the
gonality, then maybe higher versions succeed. Perhaps there is
even a universal quantum limit of such problem $EP\infty$ for a
divisor of infinite degree, leading thereby to a branched (yet
Randschlicht) version of the Bieberbach coefficient problem. This
is to mean a version of the Ahlfors map where all derivatives are
simultaneously maximized as a large convey? One can speculate
about the existence of such an universal extremal problem whose
solution would be a branched avatar (non schlicht) of the Koebe
extremal function (involved in the Bieberbach-de Branges theorem).
This would be for the given bordered surface the best circle
mapping and arguably it ought to realize the gonality.
[05.11.12]~In the classic Bieberbach problem involving the disc
the coefficients of schlicht power series are estimated by $\vert
a_n\vert\le n$ when $f'(0)=1$. If we replace the disc by a finite
bordered surface $F$ we could expect that all maps $F\to \Delta\to
{\Bbb C}$ factorizing as a circle map (of minimal degree) followed
by a schlicht map also admit universal estimates upon the
coefficients w.r.t. to a chart. Perhaps the upper bounds sequence
involved in Bieberbach-de Branges (regularly spaced integers $n$)
has to be replaced by certain spectral eigenvalues of $F$
conceived as a vibrating membranes. So the problem is the
following. Given a bordered surface $F$ marked interiorly at some
point $a$. We look at all analytic maps $F\to {\Bbb C}$ with
$f(a)=0$ and $f'(a)=1$ w.r.t. some chart. We develop $f$ in power
series and expect some universal estimates on the coefficients at
least when $f$ factorizes as a circle map of minimal degree
followed by a schlicht map. The dream would be that there is a
unique extremal function maximizing simultaneously all
coefficients and this would be essentially the best possible
Ahlfors map post-composed with the Koebe function.

Of course it may happen that all this generality is not necessary
in case the basic Ahlfors map $f_{a,b}$ is already the most
ergonomic object, in the sense of realizing the gonality for
suitable centers $a,b$.

A more orthodox way to formulate higher versions of Ahlfors'
extremal problem involves the theory of Pick-Nevanlinna
interpolation. Cf. for instance Jenkins-Suita 1979
\cite{Jenkins-Suita_1979}. The original theory being formulated in
the disc $\Delta$, one may hope to lift things via an Ahlfors map
but this probably leads nowhere. Genuine avatars of Ahlfors
extremal problem are formulated by prescribing Taylor section
(jets) at a given collection of points. Compare again
Jenkins-Suita 1979 \cite{Jenkins-Suita_1979}, building upon a
paper of Heins 1975 \cite{Heins_1975}. In this extended context
all features of the Ahlfors map persist: existence of an extremal
(via normal families), uniqueness of the solution (Heins 1975),
finite sheeted covering of the disc, and upper bound over the
mapping degree. Again a crucial question is whether such problems
always achieve the gonality.

\section{Ahlfors' proof}\label{Ahlfors:sec}

[January 2012] This section is a superficial glimpse into Ahlfors'
original resolution of his extremal problem emphasizing that
Ahlfors requires first the qualitative existence of a circle map.
A more detailed analysis will be attempted later
(Sec.\,\ref{Ahlfors-proof:sec}).

\subsection{Soft part of Ahlfors 1950: circle maps with $\le
r+2p=g+1$ sheets}\label{sec:Ahlfors-soft}

When writing the paper Gabard 2006 \cite{Gabard_2006} (and a
fortiori in my Thesis 2004 \cite{Gabard_2004}), I was very
ignorant about the depth of Ahlfors' paper (and the massive
literature around it). To be honest I am still today quite
ignorant having only a very fragmentary understanding of Ahlfors
arguments. I take this opportunity, to rectify
the arrogant claim (in
\loccit\,\cite{Gabard_2006}) to the effect that
a simplified proof of Ahlfors' theorem is proposed.
Of course, my paper only recovers the weaker assertion about
existence of  circle maps (in contradistinction to the deeper
extremal problem analyzed by Ahlfors).

Furthermore even in the weaker circle maps context, I only
realized recently [January 2012]
that a much shorter portion of Ahlfors' paper achieves this
goal (cf. Ahlfors 1950 \cite[p.\,124--126]{Ahlfors_1950}),
even with the $r+2p=g+1$ bound on the degree. We reproduce the
relevant extracts (p.\,124 and then p.\,126):

\begin{quota}[Ahlfors 1950]\label{Ahlfors-1950-circle-map:quote}

{\small \rm [p.\,124] It must first be proved that the class
of functions with $F(a)=0$ and $\vert F \vert=1$ on $C$ [=the
boundary contours] is not empty. In other words, we must show
that $\overline{W}$ can be mapped onto a full covering surface
of the unit circle.

[p.\,126] The function [\dots] maps $\overline{W}$ onto a
covering surface of the unit circle[=disk], and a standard
argument[=just number conservation] shows that every point is
covered exactly $P+1$ times. [$P$ is the genus of the double
in Ahlfors' notation]

}
\end{quota}

Thus, we have the following historical:

\begin{conj}
As early as Spring 1948, Ahlfors had an existence-proof of circle
maps of degree $\le r+2p$.
\end{conj}

This conjecture is supported by the remarks made in Nehari 1950
\cite{Nehari_1950} (cf. our Quote~\ref{Nehari-1950:quote}).
In contrast, the issue that the same upper bound $r+2p$ holds true
for Ahlfors extremals may have required Garabedian formulation of
the dual extremal problem for differentials. This is somehow in
line with Jenkins-Suita 1979 \cite{Jenkins-Suita_1979}, who speak
of the {\it Garabedian bound\/} following a coinage of Heins 1975
\cite[p.\,4]{Heins_1975}.

At any rate, it seems first crucial to understand the easy part of
Ahlfors' argument (existence of a circle map of degree $\le
r+2p$). Even here we failed as yet.

{\bf Anecdote (skip!)} Ahlfors' argument bears some vague
resemblance with the argument exposed by myself in the
RAAG-conference of 2001. Here the game was that (in view of
Riemann without Roch) any group of $g+1$ points on the curve
moves. The orthosymmetric curve in question is of course the
Schottky double of the given bordered surface. If such points are
chosen on the real locus we are forced in the non-Harnack-maximal
case ($r<g+1$, $r\equiv g+1 \pmod 2$) to select two points on the
same oval (pigeon hole principle). All the subtlety is to ensure
that those points will circulate along the complex orientation (as
the border of one half) without doing collision repulsing them in
the imaginary locus, and thereby violating total reality. Using
Abel's theorem plus some incompressible fluid argument I tried to
argue that this is always possible for a clever choice of (totally
real) divisor. However the argument was slightly vicious, and it
would require me too hard work to repair it. If I have enough
energy I should try to write down this argument, while trying to
analyze it properly.

In Ahlfors' paper (1950 \cite{Ahlfors_1950}), one starts with a
circle map of degree $\le r+2p$, and by a miraculous intervention
of Garabedian the same bound turns out to be valid for all Ahlfors
extremals. Let us refer to this vague principle as the
Ahlfors-Garabedian divination (AGD). (Vagueness only alludes to my
own poor understanding of their methods.)

Now in view of Gabard 2006 \cite{Gabard_2006}, as well as the
deeper investigation in Coppens 2011 \cite{Coppens_2011}, we know
that circle maps of lower degrees $\le r+p$ exist. Thus, granting
the AGD-divination, we may expect to find Ahlfors extremals of
correspondingly  low degrees. Of course this amounts to take the
best from two different worlds, and is extremely far from a
serious argument.
%
%
Hence a thorough study of Ahlfors 1950 \cite{Ahlfors_1950} perhaps
suitably adapted (and augmented by other tricks) could lead to a
confirmation of the naive Conjecture~\ref{gonality:conj}. Of
course this is pure speculation, and arguably the emphasis could
be a study circle maps {\it per se\/} without getting obnubilated
by Ahlfors extremal problem.

\subsection{Ahlfors hard extremal problem}

We have nothing to add for the moment, suffices to say that
Lagrange multipliers play a crucial r\^ole (as in earlier work of
Grunsky). Yet it would be nice to summarize the idea (and the
logical structure):

(1) {\bf Existence of extremals.} Ahlfors first needs the
existence of a circle map so as to arrange a nonempty set of
competing functions (giving some ground under the foots to get
started). Of course a function bounded-by-one would have been
sufficient to get started, but Ahlfors achieves much more. Of
course the normal families argument alone cannot supplant this
preliminary study.

In papers subsequent to Ahlfors', namely  Read 1958
\cite{Read_1958_Acta} and Royden 1962 \cite{Royden_1962} existence
is derived via more abstract functional-analysis (Hahn-Banach).
More on this in Sec.\,\ref{Read-Royden:sec}.

Other treatments Heins 1950 \cite{Heins_1950} appeals to Martin's
theory and elementary convexity consideration, which expressed in
more highbrow setting essentially amounts to Krein-Milman
existence of extreme points in convex bodies (cf. esp. Forelli
1979 \cite{Forelli_1979} and the discussion in Heins 1985
\cite{Heins_1985-Extreme-normalized-LIKE-AHLF}).

(2) {\bf Uniqueness of the extremal.} Looks easy (essentially like
when defining something by a universal property in category
theory). Universal properties of category theory are essentially
akin to extremal problems in geometry. This is not completely true
for some natural extremum problems admits several solutions). In
the case at hand uniqueness is essentially a version of Schwarz's
lemma.

\section{Other accounts of Ahlfors' extremal problem}

Ahlfors' paper of 1950 \cite{Ahlfors_1950} aroused quick interest
among the conformal mappers community (Nehari 1950, Heins 1950,
Garabedian 1949--50, Schiffer, etc.). Numerous papers seems to
reprove Ahlfors' theorem along (better?) routes (e.g., Read 1958,
slightly optimized in Royden 1962). The latter article seems to be
among the most popular revision of Ahlfors 1950
\cite{Ahlfors_1950}, with identic results but possible
simplifications in the proof. The present section tries to review
those (second generation) contributions while providing link to
subsequent critiques (e.g., Nehari 1950 is criticized by Tietz
1955, who in turn is attacked by K\"oditz-Timmann 1975).

\subsection{Garabedian 1949, 1950}

Garabedian qualifies himself as a hard-worker, who could absorb
simultaneously the influence of three giants: Ahlfors, Bergman and
Schiffer. As a result, he seems to have exerted a notable
influence over the final shape of Ahlfors 1950
\cite{Ahlfors_1950}, and is even apparently able to reprove the
full result of Ahlfors 1950 \cite{Ahlfors_1950} in the paper
Garabedian 1950 \cite[p.\,361]{Garabedian_1950}. (A little
Riemann-Hurwitz computation is required to convince that
Garabedian reobtains exactly the same degree $r+2p$ as Ahlfors.)
The proof deploys a rich mixture of techniques (Teichm\"uller,
Grunsky, Ahlfors, plus the variational method of Schiffer).

Another point worth noticing is the following issue oft
emphasized by Garabedian
\cite[p.\,182]{Garabedian-Schiffer_1950}:

\begin{quota}[Garabedian-Schiffer 1950]
\label{quote:Garabedian-Schiffer_1950}

{\small \rm Thus our procedure leads to the existence of the
circle mapping $F(z)$ which is associated with Schwarz's
lemma. It is to be noted that the existence of this function
lies somewhat deeper than the existence of the slit mappings
$\varphi(w)$ and $\psi(w)$ in multiply-connected domains, and
therefore it is not too surprising that the present section is
more difficult that the preceding ones. Of course, for $n=1$,
$F(z)$ is just the function found in the elementary Koebe
proof of the Riemann mapping theorem.

}

\end{quota}

Garabedian alone repeats a similar comment in Garabedian 1949
\cite[p.\,207]{Garabedian_1949-52:Book}:

\begin{quota}[Garabedian 1949]\label{quote:Garabedian_1949-52}
{\small \rm

The conformal mappings which we obtain here are closely
related to the generalization of Schwarz's lemma to multiply
connected domains in sharp form [1, 7] [=resp. Ahlfors 1947
\cite{Ahlfors_1947}, and Garabedian 1949, Duke Math. J.], and
their existence lies somewhat deeper than that of the more
standard canonical maps in a multiply connected region.

}
\end{quota}

\subsection{Nehari 1950, Tietz 1955, K\"oditz-Timann 1975}

Regarding the first two mentioned papers (Nehari 1950
\cite{Nehari_1950}, Tietz 1955 \cite{Tietz_1955}),  I suggested in
Gabard 2006 \cite[p.\,946]{Gabard_2006}, that those papers  may
have conjectured the improved control $r+p$ on the degree of
circle maps. (When discovering the $r+p$ bound ca. 2001/02, I was
not influenced by those papers which I located only later in 2005
while polishing the ultimate shape of Gabard 2006
\cite{Gabard_2006}.)

Nehari 1950 \cite{Nehari_1950} does not seem to give a new proof
of circle maps (Ahlfors' theorem), but inspired by it proposes to
describe canonical slit maps (incidentally those for which
Garabedian seems to have a lesser esteem, cf.
Quotes~\ref{quote:Garabedian-Schiffer_1950} and
\ref{quote:Garabedian_1949-52}). Nehari also shows how to express
the Ahlfors function in terms of the Bergman kernel function. (If
I understand well the situation, this is just a representation
theory yet not an alternative existence-proof.) Nehari's paper
shows that Ahlfors was in possession of the degree $r+2p$ as early
as Spring 1948, at least for a circle map. It is a delicate
question if the same bound for extremal maps requires Garabedian's
remark about the dual extremal problem. Heins' paper 1975
\cite{Heins_1975} using the term ``Garabedian's bound'' may
suggest a positive answer.
The reader is not well placed to guess the answer, but remember
that the (published) proof in Ahlfors 1950 \cite{Ahlfors_1950}
requires (and acknowledges) Garabedian's dual problem. Let us
quote the crucial extract of  Nehari:

\begin{quota}[Nehari 1950]\label{Nehari-1950:quote}

{\small \rm It was recently shown by Ahlfors [1](=L. Ahlfors,
Material presented in a colloquium lecture at Harvard University
in Spring 1948.) that the well known canonical conformal mapping
of a schlicht domain of connectivity $n$ onto an $n$-times covered
circle [5,7] (=Bieberbach 1925 \cite{Bieberbach_1925}, Grunsky
1937--41 \cite{Grunsky_1937}, \cite{Grunsky_1941_KA}) can be
generalized, in the case of an open Riemann surface, in the
following manner: an open Riemann surface of genus $g$ which is
bounded by $n$ closed curves can be mapped conformally onto a
multiply-covered circle, the number of coverings not exceeding
$n+2g$.

}

\end{quota}

Soon afterwards, Tietz 1955 \cite[p.\,49]{Tietz_1955}
criticizes (slightly) some of Nehari's asserted results:

\begin{quota}[Tietz 1955]\label{Tietz-1955:quote}

{\small \rm Bei der Herleitung seiner Schlitztheoreme kommt
Herr Nehari ebenfalls auf diese Frage; sein Beweis f\"ur die
genannte Vermutung ist jedoch unhaltbar.

Nimmt man jedoch diese Neharische Behauptung als richtig an,
so hie{\ss}e das, da{\ss} $R$ aus $p+r$, und damit $R^2$ aus
$2p+r=G+1$ Bl\"attern best\"unde; \dots

}

\end{quota}

This seems to be a forerunner of the bound $r+p$ (by commutativity
of addition!), at the conjectural level at least.
[Parenthetically, I do not understand Tietz's claim about the
sheet number of the double $R^2$. I believe that the degree keeps
the same value $p+r$, as one has to double the map not just the
space.]
Finally, Tietz concludes his paper \cite[p.\,49]{Tietz_1955}
as follows:

\medskip
{\small Die selben \"Uberlegungen, die zu unserem
Abbildungssatz f\"uhrten, erm\"oglichen auch einen neuen
Existenzbeweis f\"ur die Ahlforsche Normalform, wiederum
jedoch ohne eine Schranke f\"ur die Anzahl der ben\"otigten
Bl\"atter zu ergeben.

}
\medskip

So
Tietz does not seem to be able reprove a result as
strong as the one of Ahlfors 1950. In fact, the situation looks
even worse, since even Tietz's weak version is questioned in the
paper by K\"oditz-Tillmann 1975
\cite[p.\,157]{Koeditz-Timmann_1975}, as shown in the following
extract (parenthetical reference are ours addition):

\medskip
{\small Derartige randschlichte Abbildungen wurden von Tietz in
[4] (=Tietz, 1957, Faber-Theorie \dots) ben\"otigt, um die
Faber-Theorie auf nicht kompakte Riemannsche Fl\"achen zu
\"ubertragen. Sein in [3] (=Tietz 1955 \cite{Tietz_1955})
angegebener Beweis der Existenz solcher Funktionen ist jedoch
l\"uckenhaft. [\dots]

}
\medskip
The extract is followed by a specific objection (not reproduced
here). The article (of K\"oditz-Timmann 1975 \cite[Satz 3,
p.\,159]{Koeditz-Timmann_1975}) seems however to contain a proof
of Ahlfors' theorem based upon an ``Approximationssatzes von
Behnke u. Stein'', yet without any bound on the degree.

A propos Behnke-Stein 1947/49 \cite{Behnke-Stein_1947/49} (the
famous paper going back to 1943), it contains the result that any
open Riemann surface (arbitrary connectivity and genus) admits a
non-constant analytic function.

\begin{ques}
Can one deduce this Behnke-Stein theorem  by  agglomerating
Ahlfors extremals (or weaker circle maps) relative to compact
subregions of an adequate exhaustion?
\end{ques}

In this connection let us remember the paper by Nishino 1982
\cite{Nishino_1982}, where Ahlfors is applied to prove  existence
of (non-constant) analytic functions on certain complex surfaces
(four real dimensions). Since this Nishino paper employs Ahlfors
bound $r+2p$, it would be nice to understand it thoroughly to see
if some better constant leads to some sharpened result. (Alas it
seems that a subsequent paper of Nihino ca. 1983
 proves a stronger result (pertaining to arbitrary complex
dimension) while eradicating apparently any logical dependence
upon Ahlfors 1950. The MR-reviewer, M. Herv\'e, seems to have been
a bit overwhelmed by the work.)

\subsection{Heins 1950}\label{sec:Heins}

Heins being one of the most prolific and pleasant-to-read writers
of the U.S. school (student of Walsh), it is not surprising to
find several first classes contributions regarding our special
Ahlfors map topic. Specifically, the paper Heins 1950
\cite{Heins_1950} reproves Ahlfors' result in presumably its full
strength (this even without quoting Ahlfors 1950
\cite{Ahlfors_1950} but the closely allied work Garabedian 1949
\cite{Garabedian_1949}). Remember that Ahlfors' result was exposed
at the Harvard seminar in Spring 1948 (cf. Nehari's Quote
\ref{Nehari-1950:quote}), and must have widely circulated since
then. Taking a closer look to Heins' paper, it is at first sight
not completely evident that a bound on the degree derives from his
method but is quite  likely to do since  his quantity $m$ (number
of generators of the fundamental group, cf. p.\,571) is easily
recognized to be $2p+(r-1)$, where $p$ is the genus and $r$ the
number of contours. Thus one certainly recovers exactly Ahlfors'
result with its bound. In some sense, Heins' paper goes even
deeper than Ahlfors by treating Pick-Nevanlinna interpolation.

Several subsequent works in Heins' spirit (overlapping with
Ahlfors theorem) are worth mentioning: Heins 1975
\cite{Heins_1975}, Forelli 1979 \cite{Forelli_1979} and Heins 1985
\cite{Heins_1985-Extreme-normalized-LIKE-AHLF}.

\subsection{Kuramochi 1952}

The paper Kuramochi 1952 \cite{Kuramochi_1952} also seems to
recover Ahlfors bound for circle maps using the extremal problem.
This is maybe the sort of technical paper with too much {\it
notatio\/} and not enough {\it notio\/}? (This is a joke of
Hellmuth Kneser, compare his paper in JDMV.) Kuramochi's work
seems to be inspired by Ahlfors 1950 \cite{Ahlfors_1950} and by a
1951 paper by Nehari (confined to the planar case). Nehari offers
a positive review (in MathReviews):

\begin{quota}[Nehari 1953]\label{Nehari-Kuramochi:quote}

{\small \rm Generalizing a method developed by the reviewer
for the case of plane domains [Amer. J. Math. 73 (1951),
78--106], the author discusses extremal problems for bounded
analytic functions on open Riemann surfaces of positive genus.
The procedure is illustrated by a detailed treatment of the
case corresponding to the classical Schwarz lemma which had
previously been discussed, by different methods, by L.\,V.
Ahlfors [1950]. A complete characterization of the extremal
function is obtained and Ahlfors' positive differential is
constructed.}

\end{quota}

\subsection{Read 1958, Royden 1962 (via
Hahn-Banach)}\label{Read-Royden:sec}

We start with:

$\bullet$ Royden 1962 \cite{Royden_1962}, where existence of a
solution to Ahlfors' extremal problem is achieved via conjunction
of Hahn-Banach with Riesz's representation theorem (circumventing
thereby both Euler-Lagrange and normal families). Exploiting the
duality pointed out by Garabedian (pair of extremal problems with
a dualizing Schottky differential, i.e. one extending to the
double), the control on the degree is achieved by the usual index
formula $\deg(\vartheta)=2g-2$ (Poincar\'e 1881--85, but already
in Riemann in the holomorphic case at hand). Ahlfors'  upper bound
$\deg f_{a,b}\le r+2p$ follows.
Royden's paper is therefore quite remarkable for supplying
alternatives arguments. It seems to have been inspired mostly by:

$\bullet$ Read 1958 (two papers \cite{Read_1958_Fenn},
\cite{Read_1958_Acta}). Read is also a student of Ahlfors (as one
may learn in Ahlfors 1958 \cite{Ahlfors_1958}) and already relies
on Hahn-Banach to prove existence of an Ahlfors function (but, as
Royden observes, does not take care of making the
argument with the Schottky differential so as to bound the
degree). The technique employed (by Read) to prove extremals is to
relate the dual extremal problems (\`a la Garabedian-Ahlfors,
1949--1950) to conjugate extremum problems of the Lebesgue classes
$L_p$ and $L_q$, where $p^{-1}+q^{-1}=1$, where one maximizes an
$L_p$-norm versus vs. minimizing an $L_q$-norm. Such problems
classically reduce to Hahn-Banach. For this reduction of
Garabedian-Ahlfors to Hahn-Banach, Read employs a converse to
Cauchy's theorem (itself an application of Stokes) due to Rudin
1955 \cite{Rudin_1955-class-Hp} in the planar case. Methods of
Rogosinski--Shapiro 1953 \cite{Rogosinski-Shapiro_1953} are
another ingredient to the proof.

To summarize, the Read-Royden approach via Hahn-Banach (functional
analysis, coinage of Hadamard) effects a little drift from the
traditional Euler-Lagrange variational approach (used in Grunsky
1940--46 (\cite{Grunsky_1940}, \cite{Grunsky_1946}), Ahlfors 1947
\cite{Ahlfors_1947}, 1950 \cite{Ahlfors_1950}). As conceptually
brilliant as it is, this new method does not lead to an improved
degree control. The reason is quite simple, namely Ahlfors' bound
$r+2p$ is sharp within the extremal problem it solves
(contribution of Yamada 1978 \cite{Yamada_1978} in the
hyperelliptic case).

The game naturally splits in existence of extremals (either via
Montel's normal families or via Hahn-Banach) and then to analyze
its geometric properties. Ahlfors' 1950 treatment (apparently
influenced by Garabedian's dual extremal problem) supplies the
trick to bound the degree via a Schottky differential, and
Royden's argument looks, in this second geometric step,  virtually
osculant to Ahlfors' original.

Remember yet  that Ahlfors' original proof---presented in Spring
1948 at Harvard as reported in Nehari 1950 \cite[p.\,258,
footnote]{Nehari_1950}), and perhaps nearly similar with pages
124--126 of the published paper 1950 \cite{Ahlfors_1950}---manages
without Garabedian's influence to supply existence of circle maps
of degree bounded by $r+2p$.

\section{Existence of
(inextremal) circle maps}

This section focuses on existence of circle maps on membranes
(=finite bordered Riemann surfaces) without appeal to the extremal
problem. In fact those are logically required (at least in
Ahlfors' account but not in Royden's 1962 \cite{Royden_1962}) as a
qualitative preparation to the analysis of the quantitative
problem.

\subsection{Ahlfors 1948/50, Garabedian 1949}

[09.06.12] We mean the papers Ahlfors 1950 \cite{Ahlfors_1950} and
Garabedian 1949 \cite{Garabedian_1949}. The additional 1948 date
is intended to reflect
that Ahlfors lectured
on this material somewhat earlier, as shown by Nehari's
Quote~\ref{Nehari-1950:quote}. Those writers address the deeper
extremal problem $\max \vert f'(a)\vert$ amongst functions
bounded-by-one $\vert f\vert\le 1$, however it seems that they are
well-acquainted with topological methods (e.g., Garabedian 1949
\cite{Garabedian_1949} cites Alexandroff-Hopf's classical 1935
treatise ``{\it Topologie}'').

Such a qualitative topological inspection
seems a logical prerequisite to their treatments of the
quantitative extremal problem.
Prior to posing any extremal problem, it is vital to
ensure non-emptiness of the set of permissible competing
functions.
Perhaps the following trivial remarks are
worth doing. For domains bounded by $r$ Jordan curves, we have
clearly some function bounded-by-one (take the identity map
suitably scaled to shrink the domain inside the unit-disc). For a
general compact bordered surface, it is less obvious that such
functions exist at all. Of course one can take the Schottky double
to apply Riemann's existence theorem (of a morphism to ${\Bbb
P}^1$) and look at the image of the (compact) half. However the
latter can still cover the full Riemann sphere, which is annoying
for our purpose. [05.11.12] Using Klein's work one can certainly
find an equivariant map from the double to the sphere acted upon
by orthosymmetry (standard complex conjugation), yet it may still
be the case that the full sphere is covered by the half of the
double. [As a simple example we may take a conic $C_2$ with real
points and project it from a real point $p$ outside of the unique
oval. The corresponding map $C_2\to {\Bbb P}^1$ is equivariant and
surjective when restricted to one half of the complex locus of the
conic $C_2$. Indeed given a point of ${\Bbb P}^1$ is tantamount to
give a line $L$ through the center of perspective  $p$. This line
$L$ cuts $C_2$ in two points (except for the two real tangents).
If $L$ is a real line cutting the real locus $C_2({\Bbb R})$ we
can take as antecedent a point  on the border of the half Riemann
surface. If $L$ does not cut $C_2({\Bbb R})$ its intersection with
$C_2$ is a conjugate pairs of points one of them lying in the
fixed half of $C_2$. Finally if $L$ is an imaginary line then its
intersection with $C_2$ consists of two points distributed in both
halves of $C_2$. Indeed $L$ can by continuity be degenerated to a
real line $L_0$ missing the real locus of $C_2$ (recall that the
pencil of line is just an equatorial sphere with equator
corresponding to real lines) and since during the process no
points of $L\cap C_2$ became real it follows that both $L$ and
$L_0$ have the same distributional pattern when intersected with
the conic.]

Hence in general some preparatory qualitative ``topological''
investigation is required to see that the extremal problem is
non-vacuously posed. Remember that Ahlfors directly attacks the
existence of a circle map, where it may have been sufficient to
prove existence of a function bounded-by-one. His argument is in
part topological inasmuch as it involves annihilating the periods
of the conjugate differential of a suitable harmonic function, but
also contains a great deal of non-trivial analysis, plus basic
principles of convex geometry. We shall  try later to penetrate in
more details in Ahlfors proof.

Regarding Garabedian 1949 \cite{Garabedian_1949}, topological
considerations also plays a vital r\^ole in conjunction with
Abelian integrals, etc. We refer the reader to the original paper.
In retrospect, it may just be too sad
that this brilliant work was not directly written in the broader
context of Riemann surfaces.

\subsection{Mizumoto 1960}

This is the paper Mizumoto 1960 \cite{Mizumoto_1960} (which I
discovered only in March 2012), yet it looks quite original making
use of a topological argument involving (Brouwer's) topological
degree of a continuous mapping. So it is spiritually close to
Gabard 2006 \cite{Gabard_2006}. However Mizumoto \cite[p.\,63, Thm
1, with $N$ defined on p.\,58]{Mizumoto_1960} only recovers the
old bound of Ahlfors $r+2p$.

\subsection{Gabard 2004--2006}

The proof published in the writer's Thesis 2004 \cite{Gabard_2004}
is essentially the same as the one in 2006 \cite{Gabard_2006}
(modulo slight modifications suggested by the referee, presumably
J. Huisman). In fact J. Huisman already on the 2004 version
supplied some corrections about naive little mistakes that I made
(esp. a wrong statement of Abel's theorem forgetting to ask both
divisors to be of the same degree).
Of course it is to be hoped that the new bound $\le r+p$ will stay
correct in the long run. In case the result is true, it would be
desirable if alternative more conventional analytic methods are
able to reprove this bound $r+p$. Recent results of Coppens 2011
\cite{Coppens_2011} show the bound $r+p$ to be best-possible, at
least for generic curves in the moduli space. Coppens' work
actually supplies a much sharper understanding of all intermediate
gonalities (compare Sec.\,\ref{Coppens:subsec} for more).

There is a little historical inaccuracy in Gabard 2006
\cite{Gabard_2006}. When writing the paper, I did not realized
properly that Ahlfors has also a quite elementary proof of the
existence of circle maps of degree $r+2p$. Alas, I still do not
completely understand Ahlfors' proof yet it is clear-cut that its
elementary part does not use the extremal problem! Accordingly,
the sentence in Gabard 2006, p.\,946 reading as follows is quite
inaccurate: ``{\it {\rm [\dots]} un r\'esultat \'equivalent fut
d\'emontr\'e par L.\,V. Ahlfors en 1950, qui d\'eduit d'un
probl\`eme d'extr\'emalisation la possibilit\'e de repr\'esenter
toute surface de Riemann \`a bord compacte comme rev\^etement
holomorphe (ramifi\'e) du disque.}

\section{Related results
}

Some closely allied problems involves {\it
Parallelschlitzabbildung} (parallel-slit mapping), the
relationship with the Bergman kernel, etc. Although a bit outside
our main theme of the Ahlfors map, the methods employed are quite
similar and therefore a thorough knowledge of those proximate
mapping problems can only reinforce the general understanding. In
fact it is not to be excluded that the Kreisnormierung, or its
positive genus case avatar, known as Klein's
R\"uckkehrschnitttheorem, is logically stronger than the Ahlfors
mapping (but this is for the moment just a naked speculation).

\subsection{Parallel slit mappings (Schottky 1877,
Cecioni 1908, Hilbert 1909)}

Those mappings (abridged PSM) involve several tentacles using
varied technologies tabulated as follows:

$\bullet$ (Classical) Schottky 1877 \cite{Schottky_1877}, Cecioni
1908 \cite{Cecioni_1908} (via methods of Schwarz, and Picard).
Classically Schottky's argument is criticized (by e.g. Klein,
Cecioni, Salvemini, etc.) for depending only upon a constant count
not fully sufficient to establish the mapping existence (this
critique appears e.g., in Cecioni \loccit) It is likely that
subsequent rigorous continuity methods as developed by Brouwer
upon topological ground can easily supplement Schottky's heuristic
argument (browse through Koebe's works, etc.)]

$\bullet$ (Dirichlet resurrected) Hilbert 1909
\cite{Hilbert_1909}, Koebe 1910 \cite{Koebe_1910_Hilbert},
Courant 1910/12 \cite{Courant_1912} (those writers extend the
PSM to domains of infinite connectivity)

$\bullet$ (Extremal problem \`a la FROG
Fej\'er-Riesz-Rad\'o-(Carath\'eodory)-Ostrowski-[Grunsky]) de
Possel 1931 \cite{de-Possel_1931}, 1932 \cite{de-Possel_1932},
Gr\"otzsch, Rengel 1932/33 \cite{Rengel_1932-33}, 1934
\cite{Rengel_1934},

$\bullet$ (Bergman kernel) Nehari 1949 \cite{Nehari_1949},
Lehto 1949 \cite{Lehto_1949}, Meschkowski 1951
\cite{Meschkowski_1951}, etc.

A philosophical curiosity is that PSM is somewhat easier
(according to specialists, cf. e.g. Garabedian's
Quote~\ref{quote:Garabedian-Schiffer_1950} and Hejhal 1974
\cite{Hejhal_1974}) than the Kreisnormierung (KNP) (cf. next
section), and this already in finite connectivity (cf. e.g. the
very subtle approach to KNP imagined in Schiffer-Hawley 1962
\cite{Schiffer-Hawley_1962}). One may wonder about this sharp
discrepancy of difficulty, since it is easily conceivable that for
such canonical regions (bounded by elementary curves of the most
elementary stock) one could easily pass from one normal-form to
the other through explicit maps (at least in finite connectivity).
[Of course I do not claim that this is an easy game for me, but I
suspect so for people like Schwarz-Christoffel or Schl\"afli it
could be accessible. Of course there is maybe a difficulty in
choosing the ``accessory parameters'' but this should be
pulverizable through modern topological arguments \`a la Brouwer?]
Another striking asymmetry of the theory is that PSM hold true in
infinite connectivity (since Hilbert 1909 \cite{Hilbert_1909} and
the subsequent work of Koebe 1910 \cite{Koebe_1910_Hilbert}),
whereas  KNP is still wide open in infinite connectivity. (A very
naive guess would be to deduce KNP$\infty$ from PSM$\infty$
through a continuity method for infinite dimensional manifolds.
Maybe this suggests  using Leray-Schauder theory as an
infinite-dimensional avatar of the Brouwer degree?

Regarding PSM, lucid remarks are made in Burckel 1979
\cite[p.\,357--8]{Burckel_1979}, namely:

$\bullet$ the result in infinite connectivity is due to
Hilbert, Koebe, Gr\"otzsch, Rengel, de Possel (as we just said
also),

$\bullet$ excellent book expositions are credited to Bieberbach
193?/67 \cite{Bieberbach_1967-BUCH-Einfuehrung-in-die-konf-Abb},
Golusin 1952/57 \cite{Golusin_1952/57} and Nehari 1952
\cite{Nehari_1952-BOOK},

$\bullet$ de Possel's proof in 1931 \cite{de-Possel_1931} (and
the allied work by Rengel and Gr\"otzsch) via an extremum
problem is recognized as reminiscent of Fej\'er-Riesz's proof
of RMT. However at one point of the proof RMT is invoked.
Later de Possel 1939 \cite{de-Possel_1939} found a (short)
constructive way around this (see also Garabedian 1976
\cite{Garabedian_1976}).

$\bullet$ for an approach to PSM, and the other canonical regions
(radial or circular slits), via the Dirichlet principle see
Ahlfors 1966 \cite{Ahlfors_1966-BOOK-Cplx-Anal}.

\subsection{Kreisnormierungsprinzip (Riemann 1857, Schottky 1875/77,
Koebe 1906-08-10-20-22, Denneberg 1932, Gr\"o\-tzsch 1935,
Mesch\-kowski 1951--52, Strebel 1951--53, Bers 1961, Sibner
1965--68, Morrey 1966, Haas 1984, He-Schramm 1993)}\label{sec:KNP}

This (cavalier?) principle (abridged KNP) starts with the fact
that a multiply-connected domain of finite connectivity
maps conformally  to a circular domain.

This was already implicit in Riemann's Nachlass 1857/58/76
\cite{Riemann_1857_Nachlass} (according to Bieberbach 1968
\cite[p.\,148--9]{Bieberbach_1968-Das-Werk-Paul-Koebes} who
apparently saw a copy of Riemann's original manuscript, cf. our
Quote~\ref{quote:Bieberbach-1925} reproducing Bieberbach 1925
\cite{Bieberbach_1925}; cf. also Koebe 1910
\cite[p.\,339]{Koebe_1910_JDMV}: ``{\it Den Hauptgegenstand dieser
und des gegenw\"artigen Vortrages bildet das Problem der konformen
Abbildung eines $(p+1)$-fach zusammenh\"angenden Bereiches auf
einen von $p+1$ Vollkreisen begrenzten Bereich, ein Problem,
welches in der Literatur zuerst bei Schottky (Dissertation, Berlin
1875, umgearbeitet erschienen in Crelle 1877) in seiner bekannten
Doktordissertation auftritt, jedoch fr\"uher bereits von Riemann
in Betracht gezogen worden ist, wie aus seiner nachgelassenen
Schriften hervorgeht.}'').

The statement resurfaced more explicitly in Schottky's Thesis
1875/77 \cite{Schottky_1877} (at least in the Latin 1875 version).
The latter's argument rests again only on a naive parameter count
of moduli. Indeed, a circular domain with $r$ contours depends
upon $3r$ free parameters to describe centers and radii of those
$r$ circles, while removing the 6 (real) parameters involved in
the automorphism group of the (Riemann) sphere, we get $3r-6$
essential constants.
Invoking the (Schottky) double of the domain, whose genus is
$g=r-1$, this number agrees with Riemann's count  of $3g-3$ moduli
(where of course attention is restricted to ``real'' moduli). This
adumbrates why circular domains are flexible enough to conformally
represent any domain.
Such naive counting arguments usually turn into rigorous proofs by
appealing to some topological principles (like Brouwer's
invariance of the domain) vindicating the so-called continuity
method. This sort of game occupy several of Koebe's papers, who
probably arranged this already; see also Grunsky 1978
\cite{Grunsky_1978} for an implementation of KNP in 12 pages
(p.\,114--126).

Koebe devoted several papers to the KNP question in 1906
\cite{Koebe_1906_JDMV}, 1907 \cite{Koebe_1907_JDMV}, 1910
\cite{Koebe_1910_JDMV} (\"Uberlagerunsfl\"ache and iteration
method), 1920 \cite{Koebe_1920}.

As early as 1908 \cite{Koebe_1908_UbaK3}, Koebe advanced
conjecturally the validity of this principle for domains of
infinite connectivity: an issue still undecided today (2012), but
corroborated in He-Schramm 1993 \cite{He-Schramm_1993} in the case
of countably many boundary components (via the method of circle
packings).
%
Most of the contributions (listed in our subtitle) are carefully
referenced in He-Schramm's paper just cited.

Other proofs of the basic
(finitary) KNP result are obtained by:

$\bullet$ Courant 1950 \cite{Courant_1950} (via a Plateau-style
approach) [Micro-Warning:  Hilde\-brandt-von der Mosel 2009
\cite{Hildebrandt-von-der-Mosel_2009} and also Hildebrandt 2011
\cite{Hildebrandt_2011} credit rather Morrey 1966
\cite{Morrey_1966} for the first rigorous proof, modulo yet
another gap filled by Jost 1985 \cite{Jost_1985}].

$\bullet$ Schiffer and Hawley in several papers: Schiffer 1959
\cite{Schiffer_1959} (via the Fredholm determinant) and
Schiffer-Hawley 1962 \cite{Schiffer-Hawley_1962} (via an extremal
problem of the Dirichlet type).

It is common folklore that the Kreisnormierung, like the
uniformization and even the Ahlfors circle map belong to a
somewhat deeper class of problems than the parallel-slit mapping
succumbing quickly to elementary techniques of potential theory.
(Compare Garabedian-Schiffer's
Quotes~\ref{quote:Garabedian-Schiffer_1950} and
\ref{quote:Garabedian_1949-52}, and also Hejhal 1974
\cite[p.\,19]{Hejhal_1974} who makes similar remarks, for instance
``{\it We remark that the Koebe [circular] mapping is similar to
the universal covering map, in that neither an explicit formula
nor an explicit differential equation is known for it.\/}'')

Such higher stock problems make it challenging to ask whether
KNP(finite) could not be handled via an extremal problem \`a la
Ahlfors, or to be  historically sharper in the spirit of
FROG=Fej\'er-Riesz-(Carath\'eodory)-Ostrowski-Grunsky.

[$\bigstar$ Warning the sequel looks attractive yet erroneous, cf.
the next paragraph for a rectification $\heartsuit$] Maybe the
relevant extremal problem (under educated guess)
is to maximize the modulus of the derivative at a fixed point $a$
of the domain amongst functions bounded-by-one (in modulus) while
imposing
schlichtness to the mappings (otherwise we recover Ahlfors'
many-sheeted discs). Intuitively, this maximum pressurization
exerted at the point $a$ ensures surjectivity of the mapping while
filling most of the container in which the function in constrained
by the condition $\vert f \vert \le 1$, yet roundness of the
residual set of the image looks less intuitive. Speculating
further,   this ``Ahlfors-schlicht'' extremum problem could crack
the fully general KNP in infinite connectivity (KNP($\infty$)).
However, it suffices to remind that several complications are
reported for the usual Ahlfors function in infinite connectivity
(existence easy and uniqueness due to Havinson 1961/64
\cite{Havinson_1961/64}, Carleson 1960/67
\cite{Carleson_1967-book}, see also Fisher 1969
\cite{Fisher_1969}) by subsequent investigators like R\"oding 1977
\cite{Roeding_1977_Ahlfors}, Minda 1981
\cite{Minda_1981-image-Ahlfors-fct}, Yamada 1983--92
\cite{Yamada_1983-rmk-image-Ahlfors-fct}
\cite{Yamada_1992-Ahlfors-fct-on-Denjoy},  where the Ahlfors
extremal function ceases to be a circle map and start to omit
values). It is therefore quite overoptimistic to hope an
Ahlfors-type (=FROG) strategy toward the prestigious
KNP($\infty$).

[05.11.12] $\heartsuit$ {\it Correction.}---The beginning of the
previous paragraph is quite erroneous since  the analogue of the
Ahlfors map under the schlichtness proviso (=injectivity) is known
to take a multi-connected domain not on a Kreisbereich but on a
circular-slit disc. This result is due to Gr\"otzsch 1928
\cite{Groetzsch_1928}, Grunsky 1932 \cite{Grunsky_1932}, Nehari
1953 \cite[p.\,264--5]{Nehari_1953-Inequalities} (another proof
while crediting the just two cited works by Gr\"otzsch and
Grunsky), Meschkowski 1953 \cite{Meschkowski_1953} and finally
Reich-Warschawski 1960 \cite{Reich-Warschawski_1960}. (Those
references were already listed in
Sec.\,\ref{sec:beta-and-alpha-problems}.) It is yet to be observed
that such circular-slit-disc ranged maps fail schlichtness up to
the boundary, and one can
legitimately speculate about a suitable extremal problem akin to
Ahlfors' establishing KNP (in finite connectivity at least).

\subsection{R\"oding 1977 (still not read)}

The paper R\"oding 1977 \cite{Roeding_1977_mero} is perhaps quite
dangerous (for Gabard 2006 \cite{Gabard_2006}), yet I could not
procure a copy as yet.

\subsection{Behavior of the Ahlfors function
in domains of infinite connectivity}

There is a series of works studying the behavior of the Ahlfors
function for domains of infinite connectivity. Traditionally those
works look more confined to the domain case.

The basic existence and uniqueness result are addressed by
Havinson 1961/64 \cite{Havinson_1961/64}, Carleson
\cite{Carleson_1967-book} with simplifications in Fisher 1969
\cite{Fisher_1969}. In contrast to the finite case, the image of
the Ahlfors function does not necessarily fill the full unit
circle (=disc). We just list some main contributions:

R\"oding 1977 \cite{Roeding_1977_Ahlfors} (2 points are omitted),
Minda 1981 \cite{Minda_1981-image-Ahlfors-fct} (fairly general
discrete subset of omitted values), Yamada 1983
\cite{Yamada_1983-rmk-image-Ahlfors-fct} (omission of a fairly
general set of logarithmic capacity zero), Yamada 1992
\cite{Yamada_1992-Ahlfors-fct-on-Denjoy} (characterization of
omitted point-sets of the Ahlfors function in  case of Denjoy
domains).

\section{The quest of best-possible bounds}

The writers's own contribution $r+p$ seems, at first glance, a
dramatic improvement upon Ahlfors' upper bound $r+2p$ (at least so
sounded the diagnostic of the generous Zentralblatt reviewer of my
article, namely Bujalance). In the long run it may be that
Ahlfors' extremals are always as good for suitable choices of
points $a,b$, but only meagre evidence is presently available.

\subsection{Distribution of Ahlfors' degrees (Yamada 1978--2001,
Gouma 1998)}\label{Yamada-Gouma:subsec}

The papers by Yamada 1978 \cite{Yamada_1978}, 2001
\cite{Yamada_2001} and Gouma 1998 \cite{Gouma_1998} address
the delicate question about the exact values realized as degrees
of Ahlfors functions.

Ahlfors' pinching $r\le \deg(f_{a,b})\le r+2p$ collapses for
planar surfaces ($p=0$) to an equality, and the question is
trivially settled in this case.

Yamada and Gouma rather consider the infinitesimal form of the
problem where just a single interior point $a$ is prescribed while
maximizing $\vert f'(a) \vert$. They obtain spectacular complete
results for membranes having a hyperelliptic double ({\it
hyperelliptic membranes\/}), yet without being planar ($p=0$) in
which case we are in the trivial range already discussed.

\begin{theorem}
For a hyperelliptic membrane, the followings
hold true:

{\rm (1) (Yamada 1978)} The ponctual Ahlfors function $f_a$ has
degree $g+1$ at the fixed points of the hyperelliptic involution
(so-called Weierstrass points).

{\rm (2) (Gouma 1998)} The degree of $f_a$ can only assume values
$2$ or $g+1$.

{\rm (3) (Yamada 2001)} The case of degree  $2$ is always realized
at suitable points.
\end{theorem}

[05.11.12] Gouma's result shows
large discrepancy between degrees taken by Ahlfors extremals and
those of general circle maps. Of course the latter are more
flexible with a specimen of degree 2 (just quotient by the
hyperelliptic involution), whence circle maps exist in all even
degrees (post-compose with a power map $z\mapsto z^k$).

Those works promise a grandiose link between Ahlfors and the
classic tradition of Weierstrass points,  which probably also
regulate  the degree of Ahlfors maps for  general
(non-hyperelliptic) surfaces.

\subsection{Separating gonality (Coppens 2011)}\label{Coppens:subsec}

In another direction of dramatic depth, Marc Coppens 2011
\cite{Coppens_2011} is able to show sharpness of the bound $r+p$
claimed in Gabard 2006 \cite{Gabard_2006}. Actually, Coppens
establishes the more spectacular realizability of all intermediate
values for the gonality.
Even if Coppens' result looks at first sight subsumed to that of
Gabard, it is in reality  logically independent, so that  a
possible misfortune of Gabard's result should not necessarily
affect the truth of Coppens' one.
To be more specific, we introduce the following definition:

\begin{defn}\label{gonality:def}
The gonality (denoted $\gamma$) of a membrane
(i.e. a compact bordered Riemann surface) is  the least degree of
a full (or total) covering map to the disc.
\end{defn}

[05.11.12] A full (or total) covering map can be defined just as
non-constant analytic map taking boundary to boundary. Then it
makes good sense to Schottky-double the map and classic theory
ensures the local power-map $z\mapsto z^k$ character of analytic
functions, whence the branched cover nature of the map, in
particular its surjectivity (via a clopen argument). The jargon
``total'' is borrowed from Stoilow 1938 \cite{Stoilow_1938-Lecons}
and quite compatible with the ``total reality'' jargon (of
Geyer-Martens 1977 \cite{Geyer-Martens_1977}) incarnating the
algebro-geometric pendant of Ahlfors circle maps.

It is easy to show that a total map lacks ramification along the
boundary. (Possible argument: Else it behaves locally like
$z\mapsto z^2$ near a boundary uniformizer, but then the
half-space is wrapped to a full domain expanding outside the
permissible range of the map.)

In particular such a total map induces a usual (unramified) cover
of the circle $\partial W \to
\partial D=S^1$, whereupon the trivial lower bound $r
\le \gamma$ follows, where $r$ is the number of boundary contours
of the membrane $W$. On the other hand Gabard's main result in
2006 \cite{Gabard_2006} asserts the upper bound $\gamma \le r+p$,
where $p$ is the genus of $W$. Coppens's striking result states:

\begin{theorem} {\rm (Coppens 2011)}
Practically, all intermediate values of the gonality
compatible with the pinching $r\le\gamma\le r+p$ are realized
as the gonality of a suitable membrane of topological type
$(r,p)$. More accurately, there is a single trivial exception
when $r=1$ and $p>0$, in which case the value $\gamma=1$ must
be excluded.
\end{theorem}

Taking $\gamma=r+p$ supplies  sharpness of Gabard's upper bound.
On the other hand,  Coppens' theorem tightens
%
considerably Ahlfors' squeezing
$$
r\le \deg f_{a,b} \le r+2p
$$
into
$$
r\le\gamma\le \deg f_{a,b}\le r+2p,
$$
yielding a notable contribution to Yamada-Gouma's general question
on the distribution of Ahlfors degrees (cf. previous section). Of
course the contraction becomes most stringent when the gonality
$\gamma$ attains its maximum value (i.e., $\gamma =r+p$ if
thrusting Gabard), as it does for generic membranes in the moduli
space ${\cal M}_{r,p}$ (parameterizing isomorphism classes of
bordered Riemann surfaces).

Of course the moduli space stratifies through the gonalities.
Imitating Riemann's original count in our context it should be
possible to predict dimensions of the varied strata. Such a deeper
investigation  looks  desirable to complement the theory of
Ahlfors circle maps. More on this in
Sec.\,\ref{sec:gonality-sequence}.

[05.11.12] We are presently not aware of any total-bordered avatar
of the simple Riemann-type counting argument, so efficient for
closed surfaces in predicting correctly the gonality
$[\frac{g+3}{2}]$ as well as the dimensions of moduli strata of
lower gonalities. It is suspected that this asymmetry is inherent
to the boundary behavior of total maps which causes certain
difficulties. Of course the difficulty is somewhat akin to the
intricacies arising when doing real instead of complex algebraic
geometry. Yet the problem is certainly not insurmountable.

\subsection{Naive question: Ahlfors degree vs. the gonality}

All information mentioned so far is summarized in the string of
estimates:
\begin{equation}
r\le \gamma \le \begin{Bmatrix} \le\deg f_{a,b} \\ \le r+p
\end{Bmatrix} \le r+2p.
\end{equation}
An obvious question is whether inequality $\gamma \le \deg
f_{a,b}$ is best-possible:

\begin{ques} Is Ahlfors extremal problem flexible enough that
each membrane has an Ahlfors map $f_{a,b}$ of degree as low as the
gonality $\gamma$ for suitable centers?
\end{ques}

Yamada's deep result (2001 \cite{Yamada_2001}) positively answers
the case of hyperelliptic membranes (those which are 2-gonal
$\gamma=2$).

\subsection{Other sources (Fay 1973, \v{C}erne-Forstneri\v{c} 2002)}


In Fay 1973 \cite[p.\,116]{Fay_1973}, one reads the following
assertion:

\begin{quota}[Fay 1973]
\label{Cerne-Forstneric-2002:quote}

 {\small \rm It has been proved in [3, p.\,126](=Ahlfors 1950
\cite{Ahlfors_1950}) that there are always unitary functions with
exactly $g+1$ zeroes {\it all\/} in $R$; and when $R$ is a planar
domain, it is shown in Prop.\,6.16 that $S_{0,\dots,0}\cap
\Sigma_a$ is empty for $a\in R$ and that the unitary functions
holomorphic on $R$ with $g+1$ zeroes are parametrized by the torus
$S_0$.''

}
\end{quota}

A similar comment is to be found in \v{C}erne-Forstneri\v{c} 2002
\cite[p.\,686]{Cerne-Forstneric_2002}

\begin{quota}[\v{C}erne-Forstneri\v{c} 2002]
\label{Cerne-Forstneric-2002:quote}

{\small \rm It is proved in Ahlfors 1950
\cite[pp.\,124--126]{Ahlfors_1950} that on every bordered
Riemann surface of genus $p$ with $r$ boundary components
there is an inner function with multiplicity $2p+r$ (although
the so-called Ahlfors functions may have smaller
multiplicity).

}
\end{quota}


Actually, it seems that Ahlfors' proof shows even the slightly
stronger fact that each integer $d\ge r+2p$ do arise as the degree
of a circle map.

On page 684, Rudin 1969 \cite{Rudin_1969} is quoted. Also on page
693 we find an interesting stability of inner functions of degrees
$\ge r+2p-1$.

\section{Applications of the Ahlfors mapping}\label{Sec:Applications-of-the-Ahlf-map}

This section lists some of the known applications of Ahlfors maps.
Those applications either require the extremal property or merely
conformality and the essentially topological feature of circle
maps.

\subsection{Gleichgewicht der Electricit\"at (Riemann 1857)}

This source (Riemann 1857/58/76 \cite{Riemann_1857_Nachlass}) is
the very
origin of all our story. Alas the physical applications Riemann
had in mind were apparently only
partially reproduced in H. Weber's reconstruction of the
original manuscript.
Can someone  imagine what Riemann had exactly in mind (eventually
on the basis of the original manuscript, which must still be
dormant somewhere in G\"ottingen)?

Here are some well-known remarks concerning this posthumous
fragment; compare the ``original'' (as edited by H. Weber and
reproduced in part below) as well as the remarks in Bieberbach
1925 \cite[p.\,9, \S 7]{Bieberbach_1925}.

Interestingly, Riemann starts with the first boundary value
problem for plane domains, and actually uses the conformal circle
map to solve it, whereas the reverse engineering may look more
natural in view of his Dirichlet principle philosophy. Strikingly,
Riemann anticipates both the Schwarz symmetry/reflection principle
(Schwarz 1869
\cite[p.\,106]{Schwarz_1869-Ueber-einige-Abbildungsaufgaben}) as
well as the Schottky double (Schottky 1875--77
\cite{Schottky_1877}). Typical to Riemann,  an equality sign is
virtually put between potential theory and algebraic functions:
the Green theorem is used and
Abelian integrals (of the third species) and their periods
(Periodicit\"atsmoduln)
enter the scene. [05.11.12] Recall also that Bieberbach 1968
\cite{Bieberbach_1968-Das-Werk-Paul-Koebes} asserts that Riemann's
work also contains a trace of the Kreisnormierung, and so does
earlier Koebe 1910 \cite{Koebe_1910_JDMV}. Besides, Bieberbach
1925 \cite{Bieberbach_1925} (cf.
Quote~\ref{quote:Bieberbach-1925}) gives full credit to Riemann
for the proof of circle maps in the planar case (both via
potential theory and algebraic functions) emphasizing that Weber's
account is not completely faithful of the original manuscript. In
contrast when based only on
Weber's account, reviewers of Riemann's work tend to be more
minimalist. E.g., Grunsky 1978 \cite[p.\,198]{Grunsky_1978}
writes: ``{\it Theorem~4.1.1. [i.e. full covers of the disc for
multi-connected domains] goes back to Riemann, {\rm [423]}, who
gave some hints for the proof when $D$ is bounded by circles. The
first proof is due to Bieberbach {\rm [88](=1925
\cite{Bieberbach_1925})}, who used the Schottky-double and deep
results in the theory of algebraic functions. Elementary proofs
were given by Grunsky {\rm [195](=1937--41); [\dots]} \/}''

\noindent\begin{quota}[Riemann
1857/58-1876]\label{quote:Riemann}
\leavevmode  \linebreak


 {\small
{\center{\sc Gleichgewicht der Electricit\"at auf Cylindern
mit kreisf\"ormigem\\ Querschnitt und parallelen Axen.

Conforme Abbildung von durch Kreise begrenzten Figuren.}

} \medskip \rm {\footnotesize[Footnote (Weber): Von dieser und
den folgenden Abhandlungen liegen ausgef\"uhrte Manu\-scripte
von Riemann nicht vor. Sie sind aus Bl\"attern
zusammengestellt, welche ausser wenigen Andeutungen nur
Formeln enthalten.

Der zweite Theil der Ueberschrift bezeichnet wohl besser die
allgemeine Bedeutung des Fragmentes, als die in der ersten
Auflage allein genannte specielle Anwendung. \qquad W.]

}

\medskip
Das Problem, die Vertheilung der statischen Electricit\"at
oder der Temperatur im station\"aren Zustand in unendlichen
cylindrischen Leitern mit parallelen Erzeugenden zu bestimmen,
vorausgesetzt, dass im ersteren Fall die vertheilenden
Kr\"afte, im letzteren die Temperaturen der Oberf\"achen
constant sind l\"angs geraden Linien, die zu den Erzeugenden
parallel sind, ist gel\"ost, so bald eine L\"osung der
folgenden mathematischen Aufgabe gefunden ist:

In einer ebenen, zusammenh\"angenden, einfach ausgebreiteten,
aber von beliebigen Curven begrenzten Fl\"ache $S$ eine
Funktion $u$ der rechtwinkligen Coordinaten $x,y$ so zu
bestimmen, dass sie im Innern der Fl\"ache $S$ der
Differentialgleichung gen\"ugt:
$$
\frac{\partial^2 u}{\partial x^2} +\frac{\partial^2
u}{\partial y^2}=0
$$
und an den Grenzen beliebige vorgeschriebene Werthe annimmt. [So
this is `just' the first boundary value problem, alias Dirichlet
problem.]

Diese Aufgabe l\"asst sich zun\"achst auf eine einfachere
zur\"uckf\"uhren:

Man bestimme eine Function $\zeta=\xi + \eta i$ des complexen
Arguments $z=x+iy$, welche an s\"ammtlichen Grenzcurven von
$S$ nur reell ist, in je einem Punkt einer jeden dieser
Grenzcurven unendlich von der ersten Ordnung wird, \"ubrigens
aber in der ganzen Fl\"ache $S$ endlich und stetig bleibt. Es
l\"asst sich von dieser Function leicht zeigen, dass sie jeden
beliebigen reellen Werth auf jeder der Grenzcurven ein und nur
einmal annimmt, und dass sie im Innern der Fl\"ache $S$ jeden
complexen Werth mit positiv imagin\"arem Theil $n$\,mal
annimmt, wenn $n$ die Anzahl der Grenzcurven von $S$ ist,
vorausgesetzt, dass bei einem positiven Umgang um eine der
Grenzcurven $\zeta$ von $-\infty$ bis $+\infty$ geht. Durch
diese Function erh\"alt man auf der obern H\"alfte der Ebene,
welche die complexe Variable $\zeta$ repr\"asentirt, eine
$n$\,fach ausgebreitete Fl\"ache $T$, welche ein conformes
Abbild der Fl\"ache $S$ liefert, und welche durch die Linien
begrenzt ist, die in den $n$ Bl\"attern mit der reellen Axe
zusammenfallen. Da die Fl\"ache $S$ und $T$
gleich\footnote{Gabard micro-comment: Here the last edition of
Riemann's Werke contains a little misprint $F$ instead of the
obvious $T$, not present e.g., in the French translation of
Riemann by Laugel, Paris 1898.} vielfach zusammenh\"angend
sein m\"ussen, n\"amlich $n$-fach, so hat $T$ in seinem Innern
$2n-2$ einfache Verzweigungspunkte (vgl. Theorie der Abelschen
Functionen, Art.\,7, S.\,113) und unsere Aufgabe ist
zur\"uckgef\"uhrt auf die folgende:

Eine wie $T$ verzweigte Function des complexen Arguments
$\zeta$ zu finden, deren reeller Theil $u$ im Innern von $T$
stetig ist und an den $n$ Begrenzungslinien beliebige
vorgeschriebene Werthe hat.

Kennt man nun eine wie $T$ verzweigte Function
$\widetilde{\omega}=h+ig$ von $\zeta$, welche in einem
beliebigen Punkt $\varepsilon$ im Innern von $T$ logarithmisch
unendlich ist, deren imagin\"arer Theil $ig$ ausser in
$\varepsilon$ in $T$ stetig ist und an der Grenze von $T$
verschwindet, so hat man nach dem Greenschen Satze (Grundlagen
f\"ur eine allgemeine Theorie der Functionen einer
ver\"anderlichen complexen Gr\"osse Art.\,10. S.\,18\,f.):
$$
u_{\varepsilon}=-\frac{1}{2\pi} \int u \frac{\partial
g}{\partial \eta} d \xi\,,
$$
wo die Integration \"uber die $n$ Begrenzungslinien von $T$
erstreckt ist.

Die Function $g$ aber l\"asst sich auf folgende Art bestimmen.
Man setze die Fl\"ache $T$ \"uber die ganze Ebene $\zeta$
fort, indem man auf der unteren H\"alfte (wo $\zeta$ einen
negativ imagin\"aren Theil besitzt) das Spiegelbild der oberen
H\"alfte hinzuf\"ugt. Dadurch erh\"alt man eine die ganze
Ebene $\zeta$ $n$\,fach bedeckende Fl\"ache, welche $4n-4$
einfache Verzweigungspunkte besitzt und welche sonach zu einer
Klasse algebraischer Functionen geh\"ort, f\"ur welche die
Zahl $p=n-1$ ist. (Theorie der Abel'schen Functionen Art.\,7
und 12, S.\,113, 119.)

Die Function $ig$ ist nun der imagin\"are Theil eines
Integrals dritter Gattung, dessen Unstetigkeitspunkte in dem
Punkt $\varepsilon$ und in dem dazu conjugirten $\varepsilon'$
liegen, und dessen Periodicit\"atsmoduln s\"ammtlich reell
sind. Eine solche Function ist bis auf eine additive Constante
v\"ollig bestimmt und unsere Aufgabe ist somit gel\"ost,
sobald es gelungen ist, die Function $\zeta$ von $z$ zu
finden.

Wir werden diese letztere Aufgabe unter der Voraussetzung
weiter behandeln, dass die Begrenzung von $S$ aus $n$ Kreisen
gebildet ist. Es k\"onnen dabei entweder s\"ammtliche Kreise
ausser einander liegen, so dass sich die Fl\"ache $S$ ins
Unendliche erstreckt, oder es kann ein Kreis alle \"ubrigen
einschliessen, wobei $S$ endlich bleibt. Der eine Fall kann
durch Abbildung mittelst reciproker Radien leicht auf den
andern zur\"uckgef\"uhrt werden.

Ist die Function $\zeta$ von $z$ in $S$ bestimmt, so l\"asst
sich dieselbe \"uber die Begrenzung von $S$ stetig fortsetzen,
dadurch dass man zu jedem Punkt von $S$ in Bezug auf jeden der
Grenzkreise den harmonischen Pol nimmt und in diesem der
Function $\zeta$ den conjugirt imagin\"aren Werth ertheilt.
Dadurch wird das Gebiet $S$ f\"ur die Function $\zeta$
erweitert, seine Begrenzung besteht aber wieder aus Kreisen,
mit denen man ebenso verfahren kann, und diese Operation
l\"asst sich ins Unendliche fortsetzen, wodurch das Gebiet der
Function $\zeta $ mehr und mehr \"uber die ganze $z$-Ebene
ausgedehnt wird.
[\dots]

}
\end{quota}

This last paragraph is the one where Klein identifies (by Riemann)
early examples of ``automorphic functions'' (compare
Quote~\ref{Klein-1923:quote:Riemann-1858}).

\subsection{Painlev\'e's problem
(Painlev\'e 1888, Denjoy 1909, Besicovitch, Ahlfors 1947,
Ahlfors-Beurling 1950, Vitushkin, Melnikov, Garnett, Marshall,
Jones, Tolsa 2003)}

This connection is first explored in Ahlfors 1947
\cite{Ahlfors_1947}. The point of departure is usually identified
(modulo notorious sloppiness on finding the modern formulation) in
Painlev\'e's Thesis 1888 \cite{Painleve_1888} concerned with
generalizations of Riemann's removable singularity theorem: {\it
when do all bounded analytic functions defined in the
vicinity of a compactum extend across the compactum?} Riemann's
theorem settles  removability of singletons.

A necessary and sufficient condition for removability is the
vanishing of  a certain  numerical invariant directly attached to
the Ahlfors function, the so-called {\it analytic capacity}. This
is nothing but the maximum possible  distortion $\vert
f'(\infty)\vert$ measured at infinity among all analytic functions
defined on the complement of the compactum and bounded-by-one
there. This  characterization (due to Ahlfors 1947
\cite{Ahlfors_1947}) is not regarded as a satisfactory
answer to Painlev\'e problem requiring a purely geometric
(quasi-optical) recognition procedure of removable sets. If the
compact set lies on a rectifiable curve of the plane, removability
is tantamount to zero length (Denjoy's conjecture 1909
\cite{Denjoy_1909-Painleve/Sur-les-fct-anal-unif-a-sing-discontinues},
initially a theorem which turned out to be ``gapped'', but
confirmed via Calder\'on 1977 \cite{Calderon_1977} in Marshall
\cite{Marshall_1978?}). In the general case, Vitushkin 1967
\cite{Vitushkin_1967} proposed a characterization via ``invisible
sets'' (due to Besicovitch in the 1930's), i.e. those sets having
orthogonal projections of zero Lebesgue measure along almost all
directions. Verdera 2004 \cite{Verdera_2004} explains brilliantly
a metaphor with ghost objects virtually impossible to photography.
Alas, Vitushkin's expectation turned out to be not entirely
correct, cf. Jones-Murai 1988 \cite{Jones-Murai_1988} for a
counterexample, yet it gives already an approximate idea of the
whole problem. For instance, the prototypical example is the {\it
one-quarter Cantor set\/}: a unit-square subdivided in $4\times
4=16$ congruent subsquares whose only  4 extreme
``corner-squares'' are kept, with this operation iterated ad
infinitum (Fig.\,\ref{Cantor:fig}). The resulting Cantor set turns
out to be removable (Garnett 1970 \cite{Garnett_1970}), but has
positive (Hausdorff) length since its projection on the line of
slope $1/2$ fills a whole interval (again Fig.\,\ref{Cantor:fig}).
Here the 1/2-slope photography of the set is Lebesgue massive
(hence ``visible''), yet most other projections give sets of zero
measures, in accordance with the removability of the set. Of
course Garnett argues
differently  using in particular the classic analytic theory of
Ahlfors-Garabedian.

\begin{figure}[h]
\centering
    \epsfig{figure=Cantor.eps,width=62mm}

\caption{\label{Cantor:fig} One-quarter Cantor set and its
projection upon the line of slope $1/2$ (a full interval)}
\end{figure}

The Painlev\'e problem
engaged
 many investigators (Painlev\'e 1888 \cite{Painleve_1888},
Denjoy 1909
\cite{Denjoy_1909-Painleve/Sur-les-fct-anal-unif-a-sing-discontinues},
Urysohn, Besicovitch 1930's, Ahlfors 1947 \cite{Ahlfors_1947},
Ahlfors-Beurling 1950 \cite{Ahlfors-Beurling_1950}, Vitushkin 1958
\cite{Vitushkin_1958}, 1967 \cite{Vitushkin_1967}, Garnett 1970
\cite{Garnett_1970}, Melnikov 1967 \cite{Melnikov_1967}, 1995
\cite{Melnikov_1995}, Calder\'on 1977 \cite{Calderon_1977}, 1978
\cite{Calderon_1978-ICM}, G. David 1998 \cite{David_1998}, and
many others up to its ultimate solution in Tolsa 2003
\cite{Tolsa_2003}. This tour de force blends a vast array of
technologies (Melnikov's Menger curvature, stopping processes \`a
la Carleson already involved in the corona, etc.)

As a naive question how much of this theory extends to Riemann
surfaces, using say Ahlfors 1950 \cite{Ahlfors_1950} instead of
Ahlfors 1947 \cite{Ahlfors_1947}. [06.11.12] Of course it may be
argued that most compactums of interest are phagocytable in a
chart or a schlichtartig region hence planar via Koebe's theorem.
However it may seem that non-planar compactums exist as well on
Riemann surfaces? What is the simplest example if any? Of course I
certainly miss(ed) something trivial. A naive example is to take
the $1/4$-Cantor set and project it down to the torus ${\Bbb
R}^2/{\Bbb Z}^2$, but of course the latter set may be planarized
again (via suitable retrosections).

\subsection{Type problem (Kusunoki 1952)}

In a 1952 paper \cite{Kusunoki_1952}, Kusunoki found a clever
application of the Ahlfors function to the type of open Riemann
surfaces. Beware that the type of open Riemann surfaces is here
understood in the analytic sense due the Finnish school (Myrberg
1933 \cite{Myrberg_1933}, Nevanlinna 1941 \cite{Nevanlinna_1941})
of having a Nullrand.

More precisely, Nevanlinna 1941 (\loccit) introduced a notion of
surfaces with {\it null-boundary} (Nullrand). This amounts to
exhaust the open Riemann surface by compact subregions $F_n$,
while solving via Dirichlet (rescued by Schwarz, Hilbert, etc.)
the boundary problem $\omega_n$ equal $0$ on $\Gamma_0=\partial
F_0$ and equal to $1$ on $\Gamma_n=\partial F_n$. As the
subregions $F_n$ expand to infinity, two scenarios are possible:

$\bullet$ either the $\omega_n$ converges to $0$, or

$\bullet$ the sequence $\omega_n$ converges to a positive harmonic
function, $\omega$.

\begin{defn}
In the first case, the open Riemann surface $F$ is said to
have null-boundary, and in the second case to have positive
boundary.
\end{defn}

Null-boundary is equivalent to having no Green's function, or a
recurrent Brownian motion. More relevant to Kusunoki's work is
Nevanlinna's equivalent formulation in terms of the convergence to
$0$ of the Dirichlet integral $d_n=D[\omega_n]$.

Now, Kusunoki proves the following estimate (yielding a
null-boundary criterion in case the right hand-side explodes
to infinity):

\begin{lemma} {\rm (Kusunoki, 1952)}
$$
\frac{1}{2\pi \lambda_n} \log{\frac{1}{\bar{r}_n}} \le
\frac{1}{d_n},
$$
where $\lambda_n\le r_n+2p_n$ is the degree of an Ahlfors
function $f_n\colon F_n\to D$ and $\bar{r}_n$ the maximum
value of $f_n$ achieved on the Anfangsbereich $F_0$ of the
exhaustion.
\end{lemma}

Kusunoki's argument
does not seem to use in any fundamental way the extremal property
of the Ahlfors function. Thus perhaps any circle map (of possibly
lower degree, e.g. $\le r_n+p_n$ via Gabard 2006
\cite{Gabard_2006}) accomplishes the job as well. This option is
also corroborated by the fact that Kusunoki also appeals to
Bieberbach 1925 \cite{Bieberbach_1925} where no extremal property
is put in the forefront. Accordingly there is some hope to derive
a sharper Kusunoki's estimate. Alas the magnitudes $r_n$ change
during the process so the net bonus is hard to quantify.

\subsection{Carath\'eodory metric (Carath\'eodory 1926,
Grunsky 1940, etc.)}

Cf. for instance Grunsky 1940 \cite[p.\,232, \S 3]{Grunsky_1940},
Burbea 1977 \cite{Burbea_1977-Caratheodory}.

\subsection{Corona
(Carleson 1962, Alling 1964, Stout 1964, Hara-Nakai 1985) }

In Alling 1964 \cite{Alling_1964}, the explicit degree bound
$r+2p$ of the Ahlfors map is {\it not\/} employed. In fact any
``innocent''  circle map (of finite degree and not necessarily
solving Ahlfors' extremal problem) suffices to transplant the
truth of Carleson' corona theorem (1962 \cite{Carleson_1962}) from
the disc to any finite bordered Riemann surface. Assuming that
Ahlfors circle mapping theorem is really involved to prove, or
speculating on a very apocalyptic earthquake destroying
simultaneously all the ca. 13 proofs presently available, it is
still true that the Alling/Stout extension of the corona persists
all such crashes. Recall that K\"oditz-Timmann 1975
\cite{Koeditz-Timmann_1975} prove existence of a circle map (via a
Behnke-Stein approximation theorem) without any control on the
mapping degree. This weak form of Ahlfors is enough to complete
Alling's proof.

In contrast,  Hara-Nakai 1985 \cite{Hara-Nakai_1985} exploit fully
Ahlfors bound $r+2p$ for the finer {\it corona problem with
bound\/}. The obvious problem is whether one can produce better
corona bounds using circle maps of lowered degrees (e.g. those in
Gabard 2006 \cite{Gabard_2006}). What probably plagues the game is
that even in the disc case sharp estimation of the best corona
constant is still an open difficult matter. Cf. e.g. Treil 2002
\cite{Treil_2002}, where the best upper estimate of Uchiyama 1980
(Preprint) is supplemented by a lower bound improving one of
Tolokonnikov 1981.

Literature includes:

$\bullet$ For the disc: Carleson 1962 \cite{Carleson_1962},
H\"ormander,  Gamelin 1980 \cite{Gamelin_1980-Wolff's-proof}
(Wolff's proof), Garnett's book 1981 \cite{Garnett_1981-BOOK},
etc.

$\bullet$ For bordered surfaces: Alling 1964 \cite{Alling_1964},
Hara-Nakai 1985 \cite{Hara-Nakai_1985}, Oh 2008 \cite{Oh_2008}.

\subsection{Quadrature domains (Aharonov-Shapiro 1976,
Sakai 1982, Gustaffson 1983, Bell 2004, Yakubovich 2006) }

This is another discipline bearing deep connections with the
Ahlfors function. For instance Aharonov-Shapiro 1976
\cite{Aharonov-Shapiro_1976} prove that Ahlfors maps associated to
quadrature domains are algebraic. Combining this with works by
Gustafsson 1983 \cite{Gustafsson_1983-Quadrature}, Bell 2005
\cite{Bell_2005-Quadrature-domains} arrives at the striking
conclusion: ``{It is proved that quadrature domains are ubiquitous
in a very strong sense in the realm of smoothly bounded multiply
connected domains in the plane. In fact they are so dense that one
might as well assume that any given smooth domain one is dealing
with is a quadrature domain, and this allows access to a host of
strong conditions on the classical kernel functions associated to
the domain.}''

Compare also Yakubovich 2006 \cite{Yakubovich_2006}, and the
references therein.

\subsection{Wilson's optical recognition of dividing  curves
(Gabard 2004)}

[30.12.12] Another highbrow (yet poorly explored) application of
Ahlfors theorem was sketched in Gabard's Thesis (2004
\cite[p.\,7]{Gabard_2004}). This was an answer to Wilson's
question (1978 \cite[p.\,67]{Wilson_1978}) on deciding the
dividing character of a plane curve by sole inspection of its real
locus. Here again Ahlfors theorem affords an answer:  a real curve
is dividing iff it admits a total pencil (with possibly imaginary
conjugate basepoints). Yet it must be admitted that the answer,
albeit perfectly geometric, has probably little algorithmic value
unless complemented by further insights.

\subsection{Steklov eigenvalues (Fraser-Schoen 2010,
Girouard-Polterovich 2012)}

Compare the paper by Fraser-Schoen 2010/11
\cite{Fraser-Schoen_2011} where, for the first time, the Ahlfors
map is applied to spectral theory (the first Steklov eingenvalue).
Of course the basic trick of conformal transplantation is akin to
the closed case (Yang-Yau 1980 \cite{Yang-Yau_1980}), yet in the
bordered case it seems that the Ahlfors map respects precisely
what should be, when it comes to take care of the Neumann boundary
condition. In this respect the
 Fraser-Schoen contribution looks extremely original.

Building upon a paper of Payne-Polya-Schiffer,
Girouard-Polterovich 2012 \cite{Girouard-Polterovich_2012} are
able to extend the (Fraser-Schoen) estimate to higher eigenvalues.

\subsection{Other (Dirichlet-Neumann) eigenvalues (Gabard 2011)}

Inspired by Fraser-Schoen exciting paper, I also tried the game
with the modest arXiv note Gabard 2011 \cite{Gabard_2011}, where
the second inequality of Hersch 1970 \cite{Hersch_1970} is adapted
to  configurations of higher topological structure. Note that the
other two remaining inequalities of Hersch are probably likewise
extensible (involving the quadrant and octant of a sphere).

\subsection{Klein's intuition (Klein 1876, Marin 1988,
Viro 2013, Gabard 2013)}

Another little application of the Ahlfors map can be given to
Klein's intuition that a orthosymmetric (i.e. dividing or of
type~I) curve in the plane cannot acquire a solitary double point
by progressive variation of its coefficients. This goes back to
Klein 1876, and was probably justified by several workers though
in a somewhat different shape from this original statement (e.g.
Marin 1988 \cite{Marin_1988}, and based upon him Viro 1986/86
\cite{Viro_1986/86-Progress}). For a clear-cut arguement using the
deep Ahlfors theorem, cf. our
Lemma~\ref{Klein-via-Ahlfors(Viro-Gabard):lem} below, which was
essentially suggested to me by Viro (though in the modern
Marin-Viro formulation differing somewhat from Klein's original
assertion). Our lemma below
(\ref{Klein-via-Ahlfors(Viro-Gabard):lem}) is however exactly
Klein's assertion, though proved by the device of Ahlfors maps.
Even if Teichm\"uller 1941 \cite{Teichmueller_1941} should be
right by ascribing to Klein the existence of Ahlfors circle maps,
it is quite unlikely that Klein disposed of this as early as 1876
(the critical range being rather ca. 1882 right before the
psychological collapse of Klein due to overwork). Accordingly our
(Viro inspired) proof of Klein's assertion via Ahlfors might be a
bit too eclectic, yet it is quite hygienical while requiring
little topological concentration.

\section{Eclectic applications of the Ahlfors
map}\label{Sec:Virtual-applications-Ahlf-map}

Those are only
oneiric applications of Ahlfors maps, i.e. topics bearing only
vague analogy to our main focus.

\subsection{Filling area conjecture (Loewner 1949, Pu 1952,
Gromov 1983)}

\begin{figure}[h]
\centering
    \epsfig{figure=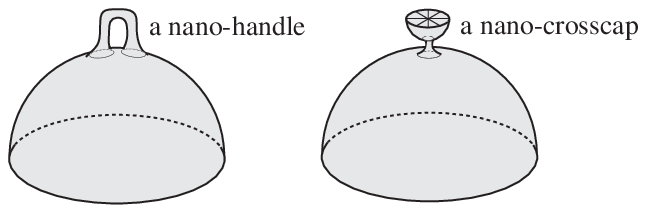,width=72mm}

\caption{\label{Wineglass:fig} Filling area conjecture for a
M\"obius band?}
\end{figure}

This was already discussed in the Introduction. One may wonder
whether the FAC is also meaningful (and true) for non-orientable
membranes. It seems so, imagine, e.g, a hemisphere surmounted by a
microscopic cross-cap over a ``glass-of-wine shaped'' protuberance
at the north pole (Fig.\,\ref{Wineglass:fig}). This membranes
satisfy FAC, for it effects no shortening of the intrinsic
distance of the circle while having an area slightly larger than
that of the hemisphere. One (possibly more accessible) question is
whether the FAC holds true for membrane having the topological
structure of a M\"obius band (equivalently a disc with a single
cross-cap). This case in view of simplicity of the topological
structure is perhaps already known, or at least accessible via the
traditional methods of Loewner-Pu, etc. (Alas, I am not aware of a
specific reference.)

Another option is to generalize Gromov's problem to membranes
filling several contours (as suggested by J. Huisman ca. Sept.
2011). Arguably, a disjoint union of hemispheres is the best
filling, at least when the contours are completely insulated (at
infinite distance). Perhaps specifying some finite
distance-functions $\rho_{i,j}$ between each pair of circles one
can expect a least-area connected filling (without shortenings),
but I have presently no clear view on how to pose properly such
generalized problems.

\subsection{Open Riemann surfaces embed in ${\Bbb C}^2$
(Narasimhan, Gromov, Slovenian school,
etc.)}\label{Open-RS-embed-in-C2:sec}

The Slovenian school of complex geometry (\v{C}erne,
Forstneri\v{c}, Globevnik, etc.) are also frequently employing the
Ahlfors map. One among the most notorious elusive open problem
(Narasimhan,
Gromov, Forstneri\v{c}, Wold, etc.) is:

\begin{conj}\label{Narasimhan:conj} Any open Riemann surface embeds properly in ${\Bbb C}^2$
(equivalently such that the image is a closed set).
\end{conj}

In Forstneri\v{c}-Wold 2009 \cite{Forstneric-Wold_2009}, the full
problem (or at least the case of interiors of finite Riemann
surfaces) is reduced to the following
finitary version, seemingly much more accessible, yet apparently
still out of reach:

\begin{conj} {\rm (FW2009)} \label{FW2009:conj}
 Each compact bordered Riemann surface $F$ embeds
holomorphically in the plane ${\Bbb C}^2$.
\end{conj}

[06.11.12 (based on ideas of ca. Sept. 2011)] Such an embedding is
possible whenever the corresponding real curve $C$ (namely the
Schottky double of $F$) admits a {\it totally real pencil of
lines}. This is for instance the case for Klein's G\"urtelkurve
(any real plane quartic with 2 nested ovals).
Fig.\,\ref{Coppens:fig} below provides plenty of other baby
examples (alas most of them being only immersed). Indeed in this
situation (total pencil of lines) the corresponding projection is
totally real and the allied morphism $C\to {\Bbb P}^1$ induces a
continuous map between the imaginary loci, i.e. $C({\Bbb
C})-C({\Bbb R})\to {\Bbb P}^1({\Bbb C})-{\Bbb P}^1({\Bbb R})$. It
follows that an imaginary line of the pencil cuts the curve $C$
{\it unilaterally\/} (i.e. only along one half of the
orthosymmetric Riemann surface). Removing such an imaginary line
from ${\Bbb P}^2({\Bbb C})$ leaves a replica of ${\Bbb C}^2$
containing entirely the original bordered surface $F$. This simple
method fails miles-away from the full Forstneri\v{c}-Wold
desideratum. Indeed Ahlfors theorem (1950 \cite{Ahlfors_1950})
only implies existence of a totally real pencil but a priori
involving auxiliary curves of order higher than one. On the other
hand when starting from the abstract bordered surface (and its
double) we may have first a projective model in ${\Bbb P}^3$,
which projected down to the plane ${\Bbb P}^2$ may develop
singularities. Hence the model in question is only immersed in
general. Our naive approach only helps grasping the notorious
difficulty of the question, yet still permits to settle a limited
collection of special cases. Actually the method, requiring a
totally real pencil of lines, applies only
 to real dividing smooth curves of order $m=2k$ having a deep nest
of profundity $k$ (that is, higher order avatars of the
G\"urtelkurve).

Another classical idea was to exhibit the required embedding
$F\hookrightarrow {\Bbb C}^2$ via a suitable pair of Ahlfors
circle maps (not necessarily extremals). This works in special
cases, e.g. hyperelliptic configurations; see
\v{C}erne-Forstneri\v{c} 2002 \cite{Cerne-Forstneric_2002}, and
also the related paper Rudin 1969 \cite{Rudin_1969}. Maybe  more
sophisticated variants of Ahlfors maps arising in the broader
Pick-Nevanlinna context could do the job, but this looks extremely
delicate.

Another natural strategy is to embed one representant in each
topological type (this is actually possible by
\v{C}erne-Forstneri\v{c} 2002 \cite[Theorem
1.1]{Cerne-Forstneric_2002}), while trying to use a continuity
argument inside Teichm\"uller (moduli) space as suggested in
Forstneri\v{c}-Wold 2009 \cite{Forstneric-Wold_2009}.

\subsection{Naive approaches to the  Forstneri\v{c}-Wold question}

This section tries (unsuccessfully) to connect some highbrow
geometry on the isometric resp. conformal embedding problem with
the FW-desideratum of the previous section. Available are some
rather formidable weapons cooked resp. by Nash-Kuiper-Gromov and
Teichm\"uller-Garsia-R\"uedy-Ko, which alas lack some rigid
analytic character upon the image model as to assess anything like
the FW-conjecture. Of course the conformal embedding technique is
most likely to
pierce the hearth of the FW-problem, yet the merely smooth
character of the conformal model hinders realizability as a
holomorphic curve. However it is still conceivable that a suitable
tour de force, somewhat akin to Garsia's (1962/63
\cite{Garsia_1962/63-algebraic-surfaces}) conformal realizability
as a real algebraic surface in $E^3$, is able to unlock the secret
of the FW-problem.

[06.11.12]
Another little puzzle is whether there is a connection with
(Gromov's and probably others) question as to whether {\it any
Riemannian surface embeds isometrically in Euclidean $4$-space
$E^4$}. (Compare Gromov 1999 \cite{Gromov_1999} delightful
preprint ``Spaces and questions'', note yet the article
Gromov-Rohlin 1969 \cite{Gromov-Rohlin_1969} where the (real)
projective plane with its round ``elliptic'' geometry is shown to
lack such an embedding.) Thus orientability is required. Via a
bordered version, we can probably embed our Riemann surface $F$
(equipped with a conformal Riemannian metric) in $E^4$
isometrically hence conformally. Via hasty thinking, the
FW-desideratum~(\ref{FW2009:conj}) follows, but alas it does not
due to the lacking rigid analytic nature of the image-model.

[16.11.12] A more rigid constant curvature version of Gromov's
isometric embedding conjecture would be actually sufficient:

\begin{conj}\label{Space-form-embedding:conj} {\rm (Space-forms embedding)} Any orientable bordered Riemannian surface of
constant Gaussian curvature $K\equiv -1$ (and totally geodesic
boundary) isometrically embeds in $E^4$.
\end{conj}

When combined with the uniformization theorem, one should be able
to deduce the FW2009 conjecture (but again this is illusory unless
one is able to ensure complex analyticity of the image).

\begin{prop}
The space-form embedding {\rm (\ref{Space-form-embedding:conj})}
implies the {\rm FW2009} conjecture (and perhaps the full proper
embedding problem {\rm (\ref{Narasimhan:conj})} via an exhaustion
trick).
\end{prop}

\begin{proof} Given the bordered Riemann surface, we take its
double $2F$, which is acted upon by a canonic involution $\sigma$.
On this closed Riemann surface, there is by the uniformization
theorem (Poincar\'e--Koebe 1907 \cite{Poincare_1907},
\cite{Koebe_1907_UbaK1}) a conformal hyperbolic metric (whenever
$\chi (F)=\frac{1}{2} \chi(2F)<0$) and  the involution $\sigma$
becomes isometric. (This equivariant uniformization is due to
Koebe 1907 \cite{Koebe_1907_UrAK}, cf. also Jost 1985
\cite{Jost_1985} for another approach via Plateau.) It follows
that the boundary of $\partial F$ are geodesics. On applying the
space-form embedding (\ref{Space-form-embedding:conj}) to $F$ we
get the FW2009 desideratum.
\end{proof}

Naively it seems that bordered hyperbolic space-forms already
embed isometrically in $E^3$, cf. Fig.\,\ref{SpaceForms:fig} for
some qualitative pictures. Of course finding a hyperbolic model
for a membrane of type $(r,p)=(1,1)$ is more tricky to visualize.

\begin{figure}[h]
\centering
    \epsfig{figure=SpaceForms.eps,width=122mm}

\caption{\label{SpaceForms:fig} Some hyperbolic space-forms in
3-space}
\end{figure}

However on tessellating the hyperbolic pants one would (under
suitable junctures) get probably trouble with the Cebyshev-Hilbert
obstruction to realizing the hyperbolic geometry in 3-space. So
maybe one must still accept variable (negative) curvature.

More flexible and suited to the problem at hand is the theory of
Teichm\"uller-Loewner-Garsia-R\"uedy realizing in the vicinity of
any smoothly embedded closed surface in $E^3$ any conformal type
of Riemann surface having the same topology via normal
deformations. In particular:

\begin{theorem}\label{Garsia's-thm}
{\rm (Garsia 1961 \cite{Garsia_1961})} Any closed Riemann surface
embeds conformally in Euclidean $3$-space $E^3$. (The image model
can also be made real-algebraic by techniques \`a la Nash, etc.,
cf. {\rm Garsia 1962/63
\cite{Garsia_1962/63-algebraic-surfaces}}.)
\end{theorem}

(R\"uedy 1971 \cite{Ruedy_1971} extended the result to open
Riemann surfaces, and may also have contributed to the embedded
version of Garsia, if the latter only showed an immersion, at
least so is claimed in Ko 1993 \cite{Ko_1993-finite-type}, yet not
so in the papers by R\"uedy.)

A direct bordered version of Garsia's result is the following:

\begin{prop}
Any compact bordered Riemann surface embeds conformally in $E^3$.
\end{prop}

\begin{proof}
Let $F$ be the given surface. Take its (Schottky) double, to get
the closed Riemann surface $2F$. By Garsia's theorem
(\ref{Garsia's-thm}) $2F$ is conformally diffeomorphic to a
classical surface in $E^3$. We conclude by restricting this
embedding to the original half of the orthosymmetric Riemann
surface $2F$.
\end{proof}

The  Garsia-R\"uedy theorem climaxes
the Riemann-Prym-Klein conception of the Riemann surface seen as a
classic (Euclid-Gauss) differential-geometric curved surface in
$3$-space (compare the introduction of Klein 1882
\cite{Klein_1882}, equivalently Quote~\ref{quote:Klein-Prym}).

Next there is a series of papers by Ko starting with his Thesis in
1989 where the Garsia-R\"uedy conformal embedding  is extended by
trading ambient 3-space by an arbitrary preassigned Riemannian
manifold of dim $\ge 3$. Specifically, he obtains the following
results:

\begin{theorem} {\rm (Ko 1989, 1991, 1999, 2001)}
Given any ambient orientable Riemannian manifold $\frak M$ of
dimension $\ge 3$, then any Riemann surface $F$ embeds conformally
in $\frak M$ provided:

{\rm (1)} $F$ is compact(=closed) {\rm (Ko's Thesis 1989
\cite{Ko_1989-compact} reissued as Ko 2001 \cite{Ko_2001})}.

{\rm (2)} $F$ has finite topological type, i.e. $\pi_1$ is of
finite generation or equivalently homeomorphic to a finitely many
punctured closed surface {\rm (Ko 1993
\cite{Ko_1993-finite-type})};

{\rm (3)} {nothing!}, i.e. $F$ is a completely arbitrary open
surface {\rm (Ko 1999 \cite{Ko_1999-open})}.
\end{theorem}

Specializing (1) to $\frak M=E^4$ seems to approach the
desideratum of FW2009. However there is a serious plague, for when
applied to a closed surface we get a conformal embedding in ${\Bbb
R}^4={\Bbb C}^2$, while complex-analyticity of the image is
inhibited by the lack of non-constant bounded analytic functions
on closed Riemann surfaces. Again the whole point is that the
conformal model (of the Garsia-R\"uedy-Ko=GRK theory) are only
smooth $C^{\infty}$-surfaces and not holomorphic curves in ${\Bbb
C}^2$.  While it is impossible to ensure complex-analyticity in
the closed case, there is no evident obstruction in the (compact)
bordered realm, except that the maximum modulus of any linear
projection on a complex line must take its maximum modulus on the
boundary.

[27.12.12] Assuming $F$ holomorphically embedded in ${\Bbb C}^2$,
we get a family of holomorphic maps $\pi_t\colon F\to {\Bbb C}_t$
parameterized by the Riemann sphere of all (complex) lines through
the origin. Since $F$ is compact bordered each such holomorphic
mapping has compact range with maximum modulus reached on a
boundary point. Of course the image is  a priori not a disc (in
which case we would have a circle map), but some more complicated
shadow of the Riemann surface $F$. Concentrating much on this
geometry it may be hoped that for some $F$ (alas not all remind
the half of the G\"urtelkurve) some obstruction is detected and FW
is false. Alternatively effecting a linear projection on the
Riemann sphere of all lines through a point $p=(x,y)\in{\Bbb C}^2$
gives the family of projections $\lambda_p\colon F \to {\Bbb
P}^1_p$, where the latter symbol is the pencil of lines through
$p$. Using translation in ${\Bbb C}^2$ all such pencils identify
to ${\Bbb P}^1_0$, where $0=(0,0)$ is the origin, and we get
$\Lambda_p\colon F \to {\Bbb P}^1_0={\Bbb P}^1({\Bbb C})$ a
holomorphic family parameterized by $p\in{\Bbb C}^2$. Again some
obstruction could occur here, but it is hard to capture. Last one
could look at the tangent line assignment yielding a sort of Gauss
mapping $F\to {\Bbb P}^1({\Bbb C})$ and explore its geometry in
the hope of detecting some fine geometric obstruction. All this is
merely canary singing without tangible grounding.

In the other optimistic direction it  can be hoped that a much
boosted version of the Garsia-R\"uedy-Ko theory proves the
FW-conjecture.

[19.12.12] First we know (from \v{C}erne-Forstneri\v{c} 2002
\cite{Cerne-Forstneric_2002}) that any topological type of finite
bordered surface contains a representative holomorphically
embedded in ${\Bbb C}^2$. Applying the high-dimensional version of
Garsia (due to Ko 1989 \cite{Ko_1989-compact}, plus subsequent
articles) we can realize all Riemann surfaces within a normal
tubular neighborhood via an (infimal) normal variation. This is
akin to a cellulite bubbling, alas destroying a priori the
holomorphic character of the initial model. However it is not to
be excluded that better controlled vibrations of the pudding
permit to explore the full moduli space. This would assess the
full Forstneri\v{c}-Wold conjecture. Of course what we are saying
here is nothing new that was not already said in FW2009, and one
requires serious new idea to make progresses.

A naive idea would be to take a holomorphic tube around the
bordered surface $F\subset {\Bbb C^2}$, i.e. a neighborhood $N$ of
$F$ together with a framing, i.e. a biholomorphic trivialization
$N \to F\times \Delta$ where $\Delta$ is the unit disc. (I hope
that a such exists but I am not sure!) One may also suppose with a
compact tube involving the closed disc $\overline{\Delta}$, so
$t\colon \overline N \to F\times \overline \Delta$. Via this
trivialization any holomorphic normal variation amounts to a
circle map, provided the amplitude of the variation is maximum
along the boundary.
Conversely given a circle map $f\colon F \to \overline \Delta$ we
construct a holomorphic deformation by considering the image under
the map
$$
F \buildrel{id\times f}\over{\longrightarrow} F \times \overline
\Delta \buildrel{t^{-1}}\over{\to} \overline  N.
$$
However the image curve is biholomorphic to the original $F$ and
our variation is trivial. Perhaps one should use quasi-conformal
avatars of circle maps (QCM for short) to perform a genuine
variation of the complex structure. Such maps clearly exist, and
we get so perhaps a tangible strategy toward Forstneri\v{c}-Wold
2009. Note however that the normal variation effected by a QCM
destroys the analytic character of the image. On the other hand it
is not essential to work with circle maps to get normal
deformations so perhaps there is some freedom to be gained here.
As another vague idea on how to construct the required
trivialization of the normal bundle (or thickening $N$) one could
imagine that it is fixed in the smooth category and then sliced by
the normal 2-planes orthogonal to $F$. Each slice would be
essentially a simply-connected domain and one would construct the
trivialization by a version of the Riemann mapping theorem with
parameters. This looks dubious but maybe leads somewhere? It is
likely now that the trivialization $t$ is not holomorphic globally
but only so in restriction to each slice (=fibre of the normal
bundle). This is good for varying moduli, but of course disrupting
the holomorphic character of the deformation.

Let us try to summarize a naive strategy toward FW2009:

{\it Step 1}.---Fix a semi-holomorphic trivialization $t$ of the
normal neighborhood of $F$ in ${\Bbb C}^2$. (Seems accessible via
a Riemann mapping theorem with parameters).

{\it Step 2}.---Among all normal deformations (cellulite bubbling)
described by $C^\infty$ maps $f$ to the disc $\overline \Delta$
those inducing ``rigid'' analytic curves under the above displayed
map.

{\it Step 3}.---By an avatar of the Teichm\"uller-Garsia-R\"uedy
technique try to calculate the moduli of the resulting
deformations. Ko's theorem ensures that all moduli are realized
via soft $C^{\infty}$ deformations, but what happens when we
restrict to the class of rigid perturbations?

\subsection{Pick-Nevanlinna interpolation}

Compare the paper Jenkins-Suita 1979
\cite{Jenkins-Suita_1979}.

\subsection{Klein-Rohlin maximality conjecture(s) (Gabard 2013)}

[11.01.13] The first paper were the notion of dividing curves
appeared (namely Klein 1876 \cite{Klein_1876}) is concluded by
some cryptical allusions which Klein might have derived from
experimental data or by a theoretical argument involving his deep
geometric intuition of Riemann surfaces.

Those intuitions were nearly forgotten for ca. 102 years until
Rohlin picked them up again in his seminal paper 1978
\cite{Rohlin_1978} enriching the complete solution ca. 1969 by
Gudkov's of Hilbert's 16th problem for sextics by the data of
complex characteristics, i.e. Klein's types I/II (erster und
zweiter Art). This Rohlin achievement was in good part made
possible by Arnold 1971 \cite{Arnold_1971/72} breakthrough of
filling the half of an $M$-curves (or more generally) one of
type~I by discs bounding the ovals. Once this is in place Klein's
assertion or intuition found quite spectacular evidence and were
somehow distorted by Rohlin in a related but somewhat more
Hilbertian and grandiose conjecture: ``a real scheme is of type~I
iff it is maximal''. One of the application was apparently
destroyed in Shustin 1985/85
\cite{Shustin_1985/85-ctrexpls-to-a-conj-of-Rohlin}, yet the
direct half ``type~I$\Rightarrow$maximal'' could still be true.

A somewhat elusive desideratum would be that this reputed
difficult conjecture of Rohlin follows from Ahlfors theorem. More
about this in Sec.\,\ref{Klein-Rohlin-conj:sec} below.

\section{Starting from zero knowledge}\label{Sec:Starting-from-zero}

As yet the text was mostly historiographical, but from now on our
intention is to elevate to the higher sphere of complete
mathematical arguments. (Of course the title of this section is
borrowed from a joke of academician V.\,I. Arnold.)

\subsection{The Harnack-maximal case (Enriques-Chisini 1915,
Bieberbach 1925, Wirtinger
1942, Huisman 2001)}

The theorem of Ahlfors (existence of circle maps) is easier in the
planar case (and due in this case to
Riemann-Schottky-Bieberbach-Grunsky, etc.). Using the
corresponding Schottky double which is a real curve (of
Harnack-maximal type), the assertion follows quite immediately
from Riemann-Roch (Riemann's inequality) via a simple continuity
argument. This argument is implicit in Enriques-Chisini 1915
\cite{Enriques-Chisini_1915-1918} (perhaps even in Riemann 1857/58
manuscript \cite{Riemann_1857}), and was then rediscovered by many
peoples including Bieberbach 1925 \cite{Bieberbach_1925},
Wirt\-inger 1942 \cite{Wirtinger_1942}, Johannes Huisman 1999/01
\cite{Huisman_2001}, and myself Gabard 2006 \cite{Gabard_2006}.
The nomenclature {\it Bieberbach-Grunsky theorem} used say by much
of the Japanese school (e.g. A. Mori 1951 \cite{Mori_1951}, etc.)
is thus slightly
in jeopardy.


\begin{lemma}\label{Enriques-Chisini:lemma} {\rm
(Riemann 1857/8, Schottky 1875/77 ($\pm$), Enriques-Chisini 1915
($\pm$), Bieberbach 1925, Wirtinger 1942, etc.)} Any planar
bordered Riemann surface with $r$ contours has a circle map of
degree $r$. Moreover the fibre over a boundary point may be
prescribed as any collection of points having one point on each
contour.
\end{lemma}

\begin{proof} Double the surface to get a closed one of genus
$g=r-1$. On the corresponding Harnack-maximal curve (i.e.
$r=g+1$), pick one point $p_i$ on each oval to get a divisor $D_0$
of degree $g+1$. Riemann's inequality states $\dim \vert D \vert
\ge d-g$, where $\vert D \vert$ is the complete linear system
spanned by the divisor $D$ and  $d$ is its degree. (This is
Riemann-Roch without Roch, and follows easily from Abel's
theorem.) It follows that the divisor $D_0=\sum_{i=1}^{g+1} p_i $,
moves in its linear equivalence class. We may thus choose in the
linear system $\vert D_0 \vert$  a line (classically denoted)
$g^1_d$, consisting of groups of $d=g+1$ points. Subtracting
eventual basepoints, this $g^1_\delta$ ($\delta\le d$) induces a
totally real morphism to ${\Bbb P}^1$, since by continuity the
points $p_i$ cannot escape their respective ovals. Indeed looking
at Fig.\,\ref{Enriques:fig} while imagining one point evading the
real locus $C({\Bbb R})$ another one must instantaneously jump to
locate himself symmetrically w.r.t. the (Galois-Klein) symmetry
$\sigma$ induced by complex conjugation.
Since a totally real morphism has degree $\ge r$, the final degree
$\delta$ must be $g+1=r$.

\begin{figure}[h] \hskip-35pt \penalty0
\centering
    \epsfig{figure=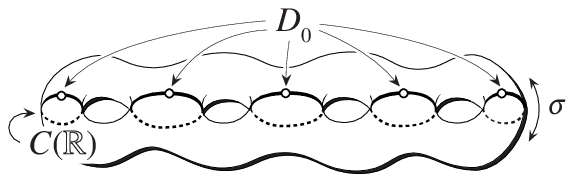,width=66mm}

\caption{\label{Enriques:fig} Totally real morphism in the
Harnack-maximal case}
\end{figure}
\end{proof}

\subsection{Gabard's argument: circle maps of
$\deg\le r+p$}\label{sec:Sketch-of-Gabard}

The basic principle used in Gabard 2006 \cite{Gabard_2006} to
prove existence of circle maps is some topological stability of
the embedding of a closed Riemann surface into its Jacobian via
the Abel map, which is quite insensitive to variations of the
complex structure. This is
how we derived universal
existence theorem valid for all Riemann surfaces with upper
control on the degree of such maps.

We suspect that the same method (suitably adapted to closed
surfaces) enables one to recover the Riemann-Meis bound
$[\frac{g+3}{2}]$ for the minimal sheet number concretizing a
genus $g$ curve as a branched cover of the line ${\Bbb P}^1({\Bbb
C})$ (cf. Riemann 1857 \cite{Riemann_1857} and Meis 1960
\cite{Meis_1960}). Yet we failed presently to write down the
details.


[22.10.12] Let us sketch rapidly the argument in Gabard 2006
\cite{Gabard_2006}, to which we refer for more details.

\begin{theorem}\label{Gabard:thm} Any bordered surface $\overline{F}=\overline{F}_{r,p}$
of type $(r,p)$ supports a circle map of degree $\le r+p$.
\end{theorem}

\begin{proof}  Using the Schottky double $C=2 \overline{F}$, it is
enough locating an unilateral divisor $D$ (i.e. one
supported in the interior denoted $F$)  linearly equivalent to its
conjugate $D^{\sigma}$. By a simple continuity argument the pencil
spanned by the pair $D, D^{\sigma}$ is totally real, hence induces
a circle map; compare Lemme 5.2 in Gabard 2006 \cite{Gabard_2006}
which we reproduce now:

\begin{lemma}
The morphism induced by a pencil spanned by an unilateral pair of
linearly equivalent divisors $D, D^{\sigma}$ is totally real.
\end{lemma}

\begin{proof} Consider $g^1_d$ the linear series spanned by $D,
D^{\sigma}$. It is readily verified to be real. To check total
reality imagine $D$ degenerating toward a point $D_0$ on the
equator $g^1_d({\Bbb R})$ of the pencil $g^1_d$ (cf.
Fig.\,\ref{Continuity-Gab:fig}).

\begin{figure}[h]
\centering
    \epsfig{figure=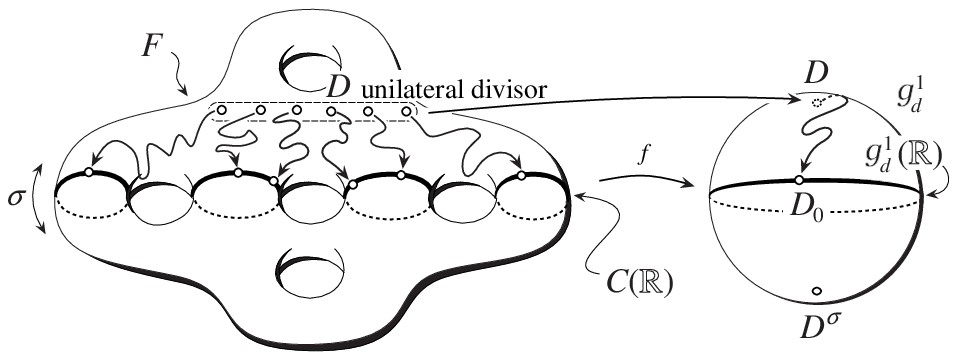,width=90mm}
    \vskip-5pt\penalty0

\caption{\label{Continuity-Gab:fig} The continuity argument
ensuring total reality of a pencil spanned by an unilateral pair
of linearly equivalent divisors (courtesy of Gabard 2004/06)}
\end{figure}

As long as $D$ stays imaginary it cannot acquire a real point
(else as the morphism induced by $g^1_d$ is real it would have a
real image). Therefore $D$ is so-to-speak magnetically confined to
the original half, hence itself unilateral. Yet when $D$ becomes
real it corresponds to a symmetric divisor (invariant under the
involution $\sigma$), which must be the limit of unilateral
divisors. The only possibility is for $D_0$ to be totally real.
Since in a sphere, any point of the equator is accessible from the
north pole, it follows that $D_0$ is always totally real. This
completes the proof.
\end{proof}

The
task is thus reduced to exhibit an unilateral divisor such that
$D\sim D^{\sigma}$ (linear equivalence on the curve $C$). Using
Abel's map $\alpha\colon C \to J$ to the Jacobian (variety) this
amounts to say that $\alpha(D)$ is a real point of the Jacobian.
Looking in the quotient $J/J({\Bbb R})$ this amounts to express
zero as a sum of unilateral points. Taking any point $x_d$ in $F$,
we search points $x_i\in F$ so that
$$
x_1+\dots +x_{d-1}=-x_d.
$$
To solve this equation we use a principle of topological
irrigation (subsumable to Brouwer's theory of the mapping degree),
but whose essence lies in the periodic behavior of the Abel map.
Specifically, we know that $\alpha$ induces an isomorphism on the
first homology. In a similar way (cf.
Fig.\,\ref{Orthosym-basis:fig}), the $r-1$ semi-cycles $\beta_1^+,
\dots, \beta_{r-1}^+$ (linking one contour to the others) and the
$2p$ cycles $\widetilde{\alpha}_1, \dots,
\widetilde{\alpha}_p,\widetilde{\beta}_1, \dots,
\widetilde{\beta}_p$ winding around the $p$ handles form a basis
of the first homology of  the quotient $T^g:=J/J({\Bbb R})$, a
$g$-dimensional (real) torus. Note that the extremity of the
semi-cycles $\beta_i^+$ are pasted together when passing to the
quotient.

\begin{figure}[h]
\centering
    \epsfig{figure=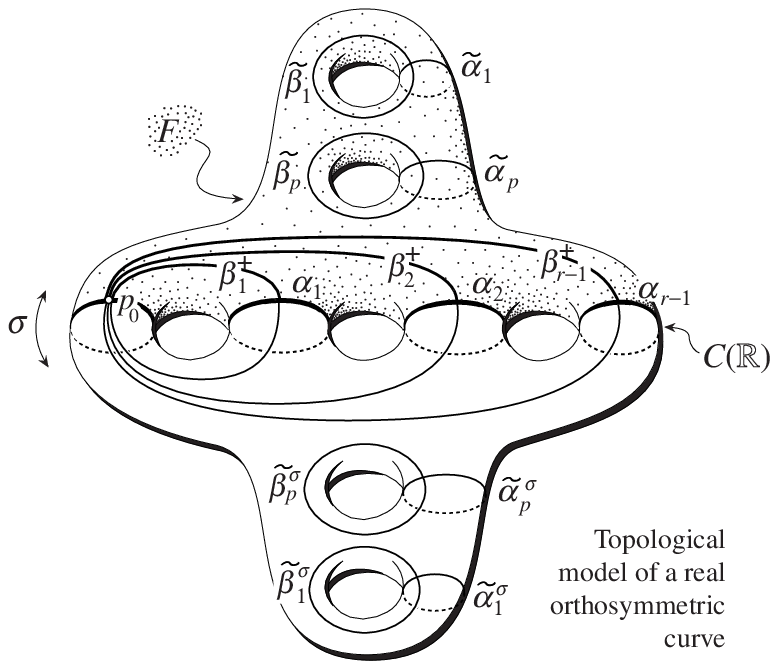,width=70mm}
    \vskip-5pt\penalty0

\caption{\label{Orthosym-basis:fig} Tracing a unilateral
collection of cycles irrigating the Jacobi torus moded out by its
real locus (courtesy of Gabard 2004/06)}
\end{figure}

The irrigation principle says that if we have $g$ cycles
representing a basis of the 1st-homology of a $g$-dimensional
torus $T^g$ then any point of the torus is expressible as the sum
of at most $g$ points situated on the given cycles. Applying this,
we can solve the above equation for $d-1\le (r-1)+2p$, i.e. $d\le
r+2p$ recovering Ahlfors bound. Since we are presently unable to
reprove Ahlfors theorem via his original argument, let us state
this as an independent theorem. (Note at the didactic level that
our proof merely use Abel 1826 \cite{Abel_1826}, and perhaps some
Riemann in as much as we use that a curve of genus $g$ supports
$g$ many holomorphic $1$-forms involved in the definition of the
Abel map, yet nothing like say Green 1828 \cite{Green_1828} which
is pivotal in Ahlfors' implementation, although this is only
stressed subconsciously.)

\begin{theorem} {\rm (Ahlfors 1950, via Gabard's method)}
\label{Ahlfors-via-Gabard:thm} Any bordered surface
$\overline{F}=\overline{F}_{r,p}$ of type $(r,p)$ supports a
circle map of degree $\le r+2p$.
\end{theorem}

Now our points $x_i$ are situated on curves traced in advance
around the handles. This constraint is not inherent to our
problem, where only unilaterality is required. Thus the points
enjoy more freedom and this is how we discovered (ca. 2002) the
possibility of improving Ahlfors. More formally, we can imagine
instead of the two cycles
$\widetilde{\alpha}_i,\widetilde{\beta}_i$ winding around a handle
a 2-cycle $\widetilde{\alpha}_i \star\widetilde{\beta}_i$ having
the shape of a 2-torus. The latter torus is not traced on our
surface $F$, but a vanishing cycle operation  makes the torus
visible. This torus is interpreted as a cycle with stronger
irrigating power. Summarizing, we have in the quotient $T^g$ the
$(r-1)$ semi-cycles and $p$ many 2-tori of stronger irrigating
power. An (evident) variant of the irrigation principle gives
solubility of the above equation for $d-1\le (r-1)+p$, i.e. $d\le
r+p$ (Gabard's bound).
\end{proof}

{\it Warning.}---[06.11.12]~Presenting the full details in some
less intuitive manner occupies the last 7 pages of Gabard 2006
(\loccit). It is hoped that the $r+p$ result is correct, but it
should not be excluded that something wrong happened (or at least
that the proof is not convincing enough). Thus more investigations
require to be made to assess or disprove the above theorem. Of
course the first part, where only Ahlfors' bound $r+2p$ is
recovered, seems less subjected to ``corrosion'', because the
irrigating cycles are readily traced on the bordered surface
(without appeal to vanishing cycles, homologies, etc.).

\subsection{Assigning zeroes and the gonality sequence}
\label{sec:gonality-sequence}

[22.10.12] Here we explore some little new ideas inspired by the
irrigating method discussed in the previous section. Alas, details
are a bit messy (mostly due to severe degradations of the little I
knew about algebraic curves). Most propositions of this section
suffer the plague of hypothetical character. We hope that, despite
vagueness of conclusions, the thematic addressed is worth
clarifying.
A general question of some interest is that of calculating for a
given bordered surface the list of all integers arising as degrees
of circle maps tolerated by the given surface. We call this
invariant the {\it gonality sequence\/}. Another noteworthy issue
is that apparently Ahlfors' upper bound $r+2p$ is always
effectively realized, in sharp contrast to Gabard's one $r+p$
which can fail to be.

\medskip

In the above argument (proof of (\ref{Gabard:thm})) we may replace
the point $x_d\in F$ by a collection of $k$ points say $z_1,
\dots, z_k\in F$. By the irrigation principle it is still possible
to solve the following equation in the group $T^g=J/J({\Bbb R})$
$$
x_1+\dots +x_{d-1}=-(z_1+ \dots+ z_k)
$$
for $d-1\le (r-1)+p$. Alas, if the divisor $z_1+ \dots+ z_k$ is
linearly equivalent to its conjugate the right hand side vanishes
in $T^g$, and  all $x_i$ could lye on the boundary of the
semi-cycle (violating the unilaterality condition). However, in
this case there is a circle map of degree $k$ exactly given by the
divisor $D=\sum_{i=1}^k z_i$. Thus, we can still conclude the
following:

\begin{prop} (Circle maps with assigned zeroes)
Given any collection $z_i$ of $k$ points in  a bordered surface
$\overline{F}$ of type $(r,p)$ there is a circle map of degree
$\le (r-1)+p+k$ vanishing on the assigned points $z_i$.
\end{prop}

\begin{proof}
It must just be observed that the pencil through $D,D^{\sigma}$,
where $D=x_1+\dots +x_{d-1}+(z_1+ \dots+ z_k)$ is basepoint free
due to the unilaterality of this divisor.  (This holds true even
if some of the $x_i$ or $z_i$ come to coincide.)
\end{proof}

It seems even that there exists circles maps of any degree $d\ge
r+p$, but I am not sure about this point. Checking the truth of
this requires the assertion that any point in the torus is
expressible as the sum of the exact number of cycles available in
the irrigating system. At first glance, this looks untrue in the
trivial irrigating system for the flat 2-torus ${\Bbb R}^2/{\Bbb
Z}^2$ consisting of  the 2 factors. Yet  the origin may be
redundantly expressed as sum of two points. Idem for a point on
the vertical axis, there is an expression as that point plus the
origin. So maybe it works. The general (hypothetical) statement
would be:

\begin{lemma} (Hypothetical lemma=Sharp irrigation principle)
Given cycles $\gamma_1, \dots, \gamma_k$ of  dimensions (say  one
and two, yet this is certainly not essential) in a $g$-torus $T^g$
such that their Pontrjagin product $\gamma_1 \star \dots \star
\gamma_k$ represents the fundamental class of $T^g$. Then any
point of $T^g$ is expressible as the sum of $k$ points $x_i$, one
situated on each $\gamma_i$. (Some $x_i$ may coincide.)
\end{lemma}

Granting this we seem to get a sharpener version of the previous
proposition.

\begin{prop} (Very hypothetical!!!) \label{hypothetical:prop}
Given any collection $z_i$ of $k$ points in  a bordered surface
$F$ of type $(r,p)$ there is a circle map of degree exactly
$(r-1)+p+k$ vanishing on the assigned points $z_i$. In particular
there exists circles maps of any degree $d\ge r+p$.
\end{prop}

In fact the real problem is that our irrigating system involves
the $r$ semi-cycles on $F$ (which  close up into $J/J({\Bbb R})$).
If the sum involves points located on the boundary of those
semi-cycles, then those points must be discarded to ensure
unilaterality of the divisor. Thus our method gives only an upper
bound on the degree of the final map, but never an exact control.

Basic examples show that special Riemann surfaces may well admit
circle maps of degree $d<r+p$ (cf. e.g. Fig.\,\ref{Chambery:fig}).

\begin{defn} The {\it gonality\/} $\gamma=\gamma(F)$
of a compact bordered Riemann surface $F$ is the least possible
degree of a circle map tolerated by $F$.
\end{defn}

Evidently $r\le \gamma\le r+p$ (the second estimate being Gabard's
claim). One can ask if each value $d\ge \gamma$ above the gonality
occurs as the degree of a circle map. Alas, the above irrigation
technique fails close to imply this. Our guess is that the
response is in the negative, that is, there may be  ``gaps'' in
the sequence of all circle mapping degrees.

Thus to detect a gap it is natural to look among ``special''
surfaces of small gonality in comparison to its generic value
$r+p$. A rapid glance at the combinatorics of
Fig.\,\ref{Coppens:fig} (below) helps us identify the simplest
such example as a hyperelliptic surface with $(r,p)=(2,1)$. Then
$\gamma=2<r+p=3$. Borrowing an idea of Klein, we can think of the
corresponding real curve as a doubled conic. This occurs actually
via the so-called canonical mapping (of algebraic geometry) which
fails injectiveness for hyperelliptic curves. (Note: we switch
constantly from bordered surfaces to real dividing curves,
committing oft slight abuses of language.) Klein regards this
doubled conic as a degeneration of the general G\"urtelkurve (with
two nested ovals) when both of them come to coalesce. This
projective model of the hyperelliptic surface suggests that when
projected from the doubled curve it has degree 2, but if the
center of projection moves in the inside of the conic then the
projection acquires degree 4 suddenly, without visiting the value
3. However substituting to the bordered surface this double conic
is a bit fraudulent, e.g. because the latter is reducible and
correspond rather to a disconnected Riemann surface. Also the
doubled conic looks 2-gonal in $\infty^1$ ways whereas the
original surface is uniquely 2-gonal. Thus another more reliable
argument requires to be given. (This must probably be akin to the
lemma proving uniqueness of the hyperelliptic involution.)

($\bigstar$) If I remember well there is a lemma saying that any
basepoint free pencil $g^1_d$ on a hyperelliptic curve is composed
with the hyperelliptic involution $g^1_2$. In more concrete words,
any morphism to the line factors through the hyperelliptic
projection, and so has even degree. If this is correct,
Prop.\,\ref{hypothetical:prop} is corrupted since the gonality
sequence is exactly the set of even integers $2{\Bbb N}$. This
remark would equally apply to any hyperelliptic membrane with
$r=1$ or $2$, $p$ arbitrary.

However this conclusion conflicts with the
\v{C}erne-Forstneri\v{c} claim (cf. 2002
\cite{Cerne-Forstneric_2002}) that Ahlfors proved  any surface to
exhibit a circle map of degree $r+2p$ exactly (take $r$ odd equal
to 1). Of course the mistake is mine and to be found in the parag.
($\bigstar$) right above, as shown by the following example.
\medskip

\begin{defn} (Convention) Below and in the sequel,
we shall often say just total morphism instead of totally real
morphism.
\end{defn}

\begin{figure}[h]
\centering
    \epsfig{figure=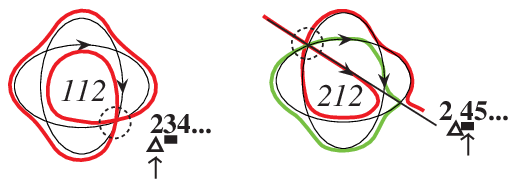,width=90mm}
    \vskip-5pt\penalty0

\caption{\label{F112:fig} Towards detecting gaps in the gonality
sequence}
\end{figure}

{\bf Example 1.} Consider a quartic with one node (so of genus
$g=2$). This is hyperelliptic (alias 2-gonal) when projected from
the node. However the curve also admits maps to the line of degree
3 (projection from a smooth point). Manufacturing a real picture
gives the picture nicknamed 112 on Fig.\,\ref{F112:fig} deduced
via sense-preserving smoothings of both ellipses (ensuring the
dividing character of the resulting curve by Fiedler 1981
\cite{Fiedler_1981}). The dashed circle indicates the node left
unsmoothed. To avoid any mysticism, our nicknaming coding consists
in writing the 3 invariants $r,p, \gamma$ as the string
$rp\gamma$.

\noindent Picture 112 shows total morphisms (i.e. with totally
real fibers over real points) of degree 2 (projection from the
node), of degree 3 (projection from the inner loop) and of degree
4 (projection from inside the inner loop).
One would like to know if $5$ is also the degree of a circle map,
etc. We believe the answer to be positive: Ahlfors exhibits circle
maps of degree exactly $r+2p$, and of all higher values too (as
follows from his convexity argument). Hence  the gonality sequence
seems to be $\gamma=2=r+p, 3=r+2p, \dots$, where the dots mean all
subsequent gonalities do occur after Ahlfors bound. This being a
bit notationally messy,  we introduce a graphical punching-card
system on the figure, where the gonality sequence 234\dots is
decorated by a triangle for $r+p$, an underlining of Ahlfors'
bound $r+2p$ (after which the gonality sequence is full), and the
arrow indicating the least position from whereon the sequence is
full. The given example does {\it not} confirm our initial guess
about gaps in the gonality sequence, so let us examine another
example.

\medskip

{\bf Example 2.}
Consider a hyperelliptic model of type $(r,p)=(2,1)$. Then the
genus $g$ of the double is $g=(r-1)+2p=1+2=3$. This prompts
looking at plane smooth quartics having the right genus $3$, but
alas the wrong gonality 3 (not 2). Thus we move to quintics
(``virtual'' smooth genus 6) and to lower down to $g=3$ we
introduce one triple point (counting like 3 double points since
perturbing slightly 3 coincident/concurrent lines creates 3
ordinary nodes). This gives the correct gonality $5-3=2$. Doing a
real picture one may draw picture 212 on Fig.\,\ref{F112:fig}.
(Keep in mind   the orientation-consistent smoothing ensuring the
dividing=orthosymmetric character of the curve). It has $r=2$, and
$p=\frac{g-(r-1)}{2}=1$. Notice total maps of degrees 2
(projection from the ``tri-node''=triple point), degree 4
(projection from the inner circuit) and degree 5 (projections from
the inside of this inner circuit). Yet we missed degree 3. Over
the complexes such a curve is not 3-gonal (because it is 2-gonal
from the tri-node and 4-gonal when projected from a smooth point).
Consequently, the allied bordered surface has circle maps of
degrees $2,4,5$ but not $3$, which is missing. Hence this example
probably corrupts our naive Prop.\,\ref{hypothetical:prop}. Also,
Gabard's bound $r+p$  needs {\it not\/} to be exactly the degree
of a circle map. Further this example shows  the gonality sequence
to be gapped in general.
\medskip

Now one general question is to wonder what can be said about the
following invariant.

\begin{defn}\label{def:gonality-sequence} The
{\it gonality sequence\/} $\Lambda=\Lambda(F)$ consists of the
ordered list $\gamma<\gamma_1<\gamma_2<\dots$ of all
integers occurring as degrees of circle maps tolerated by
the given bordered Riemann surface $F$.
\end{defn}

Fragmentary information includes the following facts, gathered as
a theorem. (To nuance reliability of  the varied constituents we
assign them some percentages of truth likelihood, with frankly
Schopenhauerian scepticism!)

\begin{theorem}
For any bordered Riemann surface with topological invariant
$(r,p)$ (viz. number of contours $r$ and genus $p$) and gonality
$\gamma$ (i.e. the least degree of a circle map), the following
estimates hold good ({\it en principe\footnote{Joke of Ivan
Babenko, yet irritating the western auditors coming down from the
``alpage''.
}}):

{\rm [100\%] $\bullet$  (T) (Trivial)} $r\le\gamma$.

{\rm [99\%] $\bullet$  (KTA)} $\Lambda$ is nonempty or
equivalently $\gamma<\infty$ is finite {\rm (Ahlfors 1950
\cite{Ahlfors_1950}, or Teichm\"uller 1941
\cite{Teichmueller_1941} crediting Klein for the result; cf. also
K\"oditz-Timmann 1975 \cite{Koeditz-Timmann_1975} for a proof via
Behnke-Stein)}.

{\rm [100\%] $\bullet$ (Semigroup property)} the set $\Lambda$ is
``multiplicative'', i.e. whenever it contains an element $\lambda
\in \Lambda$ it contains all integral multiples $k \lambda$. (This
follows by composing the corresponding circle map by a power map
$z\mapsto z^k$ from the disc to itself.) In particular  {\rm
(KTA)} implies that $\Lambda$ is always infinite.

{\rm [98\%] $\bullet$  (A50) (Ahlfors 1950)} $\gamma\le r+2p$.

{\rm [75\%] $\bullet$  (G06) (Gabard 2006)} $\gamma\le r+p$.

{\rm [79\%] $\bullet$  (C11) (Coppens 2011)} $\gamma $ takes all
intermediate values $r \le\gamma\le r+p$ (if $r=1$ the lower bound
$r$ must be modified as $2$, excepted when $p=0$).

{\rm [97\%] $\bullet$ (AFCF) (Ahlfors 1950
\cite[p.\,126]{Ahlfors_1950}, adhered to in Fay 1973
\cite[p.\,116]{Fay_1973} and \v{C}erne-Forstneri\v{c} 2002
\cite{Cerne-Forstneric_2002})} Ahlfors bound $r+2p\in \Lambda$
always belongs the gonality sequence; and
so do all higher values.
\end{theorem}

\begin{proof}
The last assertion follows from Ahlfors proof (1950
\cite[pp.\,124--126]{Ahlfors_1950}) where the origin is expressed
as convex sum of points lying on a collection of circuits in
${\Bbb R}^g$. This is always feasible for $g+1=r+2p$ points, and a
fortiori for more points. We shall try to digest Ahlfors argument
in subsequent sections.
\end{proof}

In contrast to (AFCF), Example 2 (=212 on Fig.\,\ref{F112:fig})
above shows (or at least indicates strongly) that Gabard's bound
$r+p$ is not necessarily in the gonality sequence.


Further evaluations of the gonality sequence are tabulated on
Fig.\,\ref{Coppens:fig} as bold fonts. As before, the underlined
number is Ahlfors (universal) bound $r+2p$, after which all
gonalities are realized. The position pointed onto, by a triangle,
is Gabard's bound $r+p$. The little arrow is a pointer indicating
the lowest integer after which the gonality sequence is full.

At an early stage of the tabulation, it seemed realist to advance
the following.

\begin{conj}\label{conj:full-above-Gabard} (Naive,
destroyed by Example $4$) Strictly above $r+p$ each gonality
occurs.
\end{conj}

This is pure guessing, but if true it would considerably lower
Ahlfors' universal lower bound $r+2p$ for ``fullness''. The next
example still supports the guess, but the next Example 4
ought to violate it.

\begin{figure}[h]
\centering
    \epsfig{figure=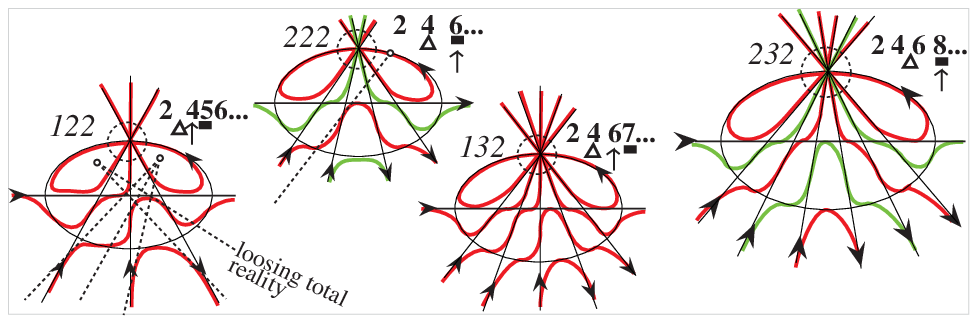,width=122mm}
    \vskip-10pt\penalty0

\caption{\label{F122:fig} Some hyperelliptic curves and their
gonality sequences. Those constructs are primarily intended to
disprove the guess that above Gabard's bound $r+p$, all gonalities
arise. Yet, the real outcome is rather that for hyperelliptic
curves one might be able to be completely explicit about the
gonality sequence.}
\end{figure}

\medskip
{\bf Example 3.} Consider, within the topological type
$(r,p)=(1,2)$ where $g=(r-1)+2p=4$,  a hyperelliptic model
($\gamma=2$). Looking at quintics (virtual genus 6) requires 2
nodes to correct the genus, but then the (complex) gonality is
still 3 (and not 2 as we would like). The trick is (like in
Example~2) to increase further the degree to permit a high order
singularity lowering drastically the gonality. So we move to
sextics (virtual genus 10) with a 4-node (counting for 6 ordinary
nodes) decreasing correctly the genus to 4. As initial
configuration we consider 3 coincident lines plus a conic through
the coincidence and another  line (in general position). An
appropriate smoothing generates picture 122 on
Fig.\,\ref{F122:fig} with $r=1$ (all real circuits being connected
through $\infty$). The gonality sequence seems to be $2, 5, 6,
\dots$. However 4 must be added to the list (being a multiple of
2). Hence the true sequence is $2,4,5,6, \dots$. Gabard's bound is
$r+p=3$, and strictly above it all values are realized (Ahlfors
bound is $r+2p=5$).

[23.10.12]---{\it Vague philosophy.}
An interesting feature of this example (122 on
Fig.\,\ref{F122:fig}) is that when gonality is very low in
comparison to topological complexity, the Riemann surface, after
having dispensed much energy to reach such a low gonality, seems
falling into some dormant state without creating many new
gonalities (missing the value 3). Perhaps this is a general
phenomenon prompted by a principle of energy conservation.

[28.12.12]---{\it Warning.} On looking carefully at picture~122
above (Fig.\,\ref{F122:fig}) it is seen that as the center of
perspective is dragged from the $4$-node toward the two red loops
or even their insides we may loose total reality for some lines of
the pencil become tangent to the circuit somewhere (cf. dashed
lines on picture~122) so that a suitable perturbation let
disappear two intersections in the imaginary locus. Accordingly,
it is not even evident form the picture that degrees $5, 6$ occur
as degrees of total maps.

\medskip
{\bf Example 4.} We now consider, in the topological type
$(r,p)=(2,2)$ for which $g=(r-1)+2p=5$, again a hyperelliptic
model. Looking at sextics with smooth genus $10$, we must use a
correction by 5 (alas not a triangular number as those involved in
multiple points). Thus we move to septics (order $m=7$) of smooth
genus $\tilde{g}=\frac{(m-1)(m-2)}{2}=15$, and a 5-node (counting
for $1+2+3+4=10$ ordinary nodes) effects the desired correction
upon the genus. Smoothing a suitable configuration gives
picture~222 on Fig.\,\ref{F122:fig} with $r=2$ (two real circuits
red and green colored). The gonality sequence includes the values
$2,4,6=r+2p, \dots$. Six being Ahlfors bound the sequence is full
from there on. When projected from the 5-node the degree is 2.
Dragging the center of perspective along one of the two red loops
gives total maps of degrees $7-1=6$ ({\it warning:} this is not
even true, cf. again the dashed line on the picture!). The value 4
is not visible on the projective model, yet arises by the
semigroup property. Studying the gonality over the complexes, it
seems evident that 3 and 5 are not even complex gonalities, and we
should be able to conclude that $2,4,6,\dots$ is the exact
gonality sequence. (Here the ``dots'' refer again to the issue
that all higher values belong to the gonality list, according to
Ahlfors.) But then our conjecture \ref{conj:full-above-Gabard} is
violated (as $5$ does not belong the list). Incidentally, this
example shows  sharpness of Ahlfors bound $r+2p$ as the place from
where the sequence is full.

\medskip
{\bf Example 5/6.} Those examples can be iterated for higher
values of the invariant $(r,p)$ while staying in the hyperelliptic
realm. The arithmetical issue is the possibility to compensate the
genus by a high-order singularity. We obtain for $(r,p)=(1,3)$,
hence $g=6$, the figure 132 on Fig.\,\ref{F122:fig}, an octic
(smooth genus 21) with a 6-node (counting for $\frac{6\cdot
5}{2}=15$ ordinary nodes) hence lowering down the genus to $6$.
The gonality sequence is $2,4,6,7=r+2p, \dots$. Similarly, for
$(r,p)=(2,3)$, $g=7$. Browsing through increasing degrees the
genus are $10, 15, 21, 28, \dots$, whereas the nodes give the list
$1,3,6,10,15,21, \dots$. The right pair is thus $28-21=7$. So we
take a 9-tic (smooth genus $\tilde
g=\frac{(9-1)(8-1)}{2}=\frac{56}{2}=28$) with a 7-node. We
construct easily picture 232 on Fig.\,\ref{F122:fig}, a curve
whose gonality sequence is $2,4,6,8=r+2p, \dots$. (Note that in
this case Ahlfors bound is sharp for the fullness of the sequence,
but it was not in the previous example. It may again be observed
that in the first example the $r+p$ bound occurs as a gonality,
but it does not in the second example.)

The real outcome of these constructions is that for (certain,
all?) hyperelliptic curves we can be totally explicit about the
gonality sequence. Iterating ad infinitum we have:

\begin{prop}
For any topological type $(r,p)$ there is a surface of
hyperelliptic type $(r,p)$ (with $r=1$ or $2$) whose gonality
sequence $\Lambda$ is known explicitly. Namely,

$\bullet$ if $r=1$, then $\Lambda=\{2,4,6, \dots, 2p, r+2p, \dots
\}$, where the first ``dots'' runs through even values and the
second means fullness after Ahlfors bound $r+2p$.

$\bullet$ if $r=2$, then $\Lambda=\{2,4,6, \dots, r+2p, \dots \}$,
where the first ``dots'' runs through even values and the second
means fullness after Ahlfors bound $r+2p$.
\end{prop}

The natural question is of course to know if this spectrum
distribution is specific to our models or generally valid for all
hyperelliptic surfaces. (This looks likely, we think, maybe just
by counting moduli.) Of course it is evident by the semigroup
property that the gonality sequence contains the value listed, and
is full after Ahlfors bound, yet the assertion that it reduces to
this requires some argument.

\subsection{A conjecture about fullness}

[23.10.12] At this stage the situation is admittedly a bit messy.
We try to clarify it by
bringing into the picture
the {\it fullness invariant\/} $\varphi$, that is the least
integer from whereon the gonality sequence is full. (On the
pictures discussed this is nothing but the little arrow used
previously.) We have the string of inequalities:
\begin{equation}
r\le \gamma \le \begin{Bmatrix}  \buildrel{{\rm Ga}}\over{\le} r+p \le \\
\hskip4pt \le \hskip8pt \varphi \hskip8pt \buildrel{{\rm
Ah}}\over{\le}
\end{Bmatrix} \le r+2p.
\end{equation}
Is any comparison possible between $r+p$ and $\varphi$? On example
212 of Fig.\,\ref{Coppens:fig} $r+p=3$
beats the fullness $\varphi=4$. Many examples on
Fig.\,\ref{Coppens:fig} do satisfy $r+p\le \varphi$, but there is
also several counter
indicators, e.g. pictures 313, 414 or 223.

The following is a trivial consequence of inequation $\gamma\le
\varphi$:

\begin{lemma}
Fullness below Gabard's bound (i.e. $\varphi < r+p$) implies
low-gonality (i.e. $\gamma < r+p$).
\end{lemma}

The converse fails, see pictures 212 or 222.

On the pictures of Fig.\,\ref{Coppens:fig} the fullness $\varphi$
is indicated by a little upward arrow. Examining examples on this
figure it seems that when the surface has generic gonality (i.e.
$\gamma=r+p$) then its fullness coincides with the gonality (i.e.
$\varphi=\gamma$). It would be interesting to know if  a general
theorem hides behind this experimental observation.

\begin{conj} (Pressing up and down: fullness conjecture)
\label{conj:fullness} If $\gamma=r+p$, then $\varphi=\gamma$. In
other words if $\gamma$ achieves maximum value (granting the truth
of Gabard's bound!) then $\varphi$ collapses to its minimum value
(namely $\gamma$). In particular the gonality sequence of a
generic surface would be
perfectly explicit, as being full  from $r+p$. This would also
show that generically Ahlfors bound $r+2p$ for fullness can be
drastically lowered.
\end{conj}

It seems plausible that an adaptation of Gabard 2006 could prove
this conjecture. (Ahlfors' original proof can also be useful.) The
idea would be that in the irrigation method the equation
$x_1+\dots +x_{d-1}=-x_{d}$ which is soluble for $d\le (r-1)+p$
points is, by the assumption made on $\gamma$, not soluble for
fewer points. One would then like to extend this ``exact
solubility'' to the equation $x_1+\dots+
x_{d-1}=-(z_1+\dots+z_k)$, where the $z_i$ is a collection of
points assigned in $F$. A vague idea is then that if some $x_i$
(or their lifts to $\overline{F}$ belong to the border) then upon
dragging the $z_i$ we may hope to displace them to avoid this
circumstance (incompatible with unilaterality). This would
construct an unilateral divisor of any assigned degree $\ge r+p$,
producing in turn circle maps of all such degrees. The conjecture
would follow.

\begin{figure}[h]
\centering
    \epsfig{figure=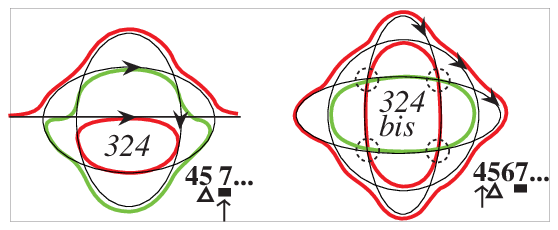,width=70mm}
    \vskip-5pt\penalty0

\caption{\label{F324:fig} Two curves with same $(r,p,\gamma)$ but
diverging gonality sequences}
\end{figure}

Another basic phenomenon is that even when two surfaces have the
same invariants $(r,p)$ and the same gonality $\gamma$ their
gonality sequences may differ. (See for such an example both
pictures 324 on Fig.\,\ref{F324:fig}). Interestingly the left
figure 324 is 4-gonal in $\infty^1$ ways (projection from the
inner oval), whereas its companion 324bis, is 4-gonal in only 4
ways (projection from the nodes). Again some  conservation law
seems involved for all the energy absorbed by the many pencils of
degree 4 living on the first model seems to
provoke the
missing of pencils of degree 6. It could be imagined that the
right curve is totally real when swept out by the pencil of conics
through the 4 nodes of degree $2\cdot 6 -4 \cdot 2=12-8=4$, but it
fails to be total for the circular conics through the 4 nodes
certainly misses the outer oval.

To investigate the fullness conjecture (\ref{conj:fullness})
further, we test curves of higher topological structure.

\begin{figure}[h]
\centering
    \epsfig{figure=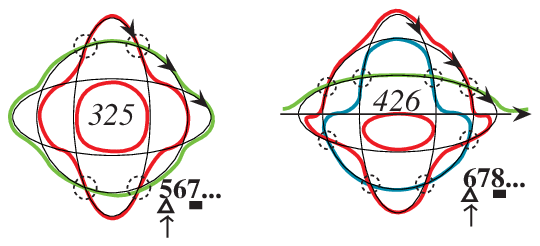,width=90mm}
    \vskip-5pt\penalty0

\caption{\label{Generic:fig} Testing the fullness conjecture}
\end{figure}

$\bullet$ For $(r,p)=(3,2)$, we seek a surface with maximum
gonality $\gamma=r+p=5$. If we imagine this gonality arising via
linear projection it is natural to look at a sextic having a deep
nest. The virtual genus is then $10$, but we want genus
$g=(r-1)+2p=2+4=6$. Thus we introduce 4 nodal singularities. To
keep the gonality maximum those nodes must not be accessible from
the inner oval, and consequently we distribute the dashed circles
(indicating unsmoothed nodes) in the ``periphery''. We thus obtain
curve nicknamed 325 (on Fig.\,\ref{Generic:fig}). It has
$\gamma=5$ and the gonality sequence is $5,6,7=r+2p,\dots$. In
fact $\gamma$ could be $<5$ via some nonlinear pencil harder to
visualize. A pencil of conics with 4 basepoints matching with the
$4$ nodes creates a series of degree $2\cdot 6-4\cdot 2=12-8=4$,
more economical than our $5$. However looking at picture 325, the
special conic consisting of two horizontal lines fails to
intersect the inner oval. Thus this pencil is not total, and we
safely conclude that $\gamma=5$, exactly. In particular, the
fullness conjecture (\ref{conj:fullness}) is verified on this
example.

$\bullet$ For $(r,p)=(4,2)$, we seek a surface with maximum
gonality $\gamma=r+p=6$. Imagine again this gonality arising via
linear projection, we consider a septic with a deep nest. The
virtual genus is then $15$, but we want $g=(r-1)+2p=3+4=7$, hence
we conserve 8 nodal singularities. We obtain so the curve labelled
426 on Fig.\,\ref{Generic:fig}. It has $\gamma=6$ and the gonality
sequence is $6,7,8=r+2p,\dots$. The fullness conjecture
(\ref{conj:fullness}) seems verified on this example. Warning
[25.10.12]: now the claim $\gamma=6$ is possibly an optical
illusion, for a pencil of cubics with basepoints assigned on the
nodes has degree $3\cdot 7-8\cdot 2=21-16=5<r+p$. If the latter is
total then $\gamma\le 5$, violating our claim $\gamma=6$. Of
course tracing pencil of cubics is not an easy game. Experience
tell us that total pencils arise when basepoints are deeply rooted
inside the deepest ovals. In the case at hand (curve 426), this
feature is not fulfilled. The 8 basepoints of the cubics pencil
lye outside the inner oval, yet, it could be that the 9th
basepoint falls (by a lucky stroke) inside this oval. In fact, it
is enough to observe that the cubic, consisting of the ellipse
through the 6 points lying highest on figure 426, plus the line
through the remaining 2 points (lying lowest on the same figure),
fails to cut the inner oval. This gives evidence that our pencil
of cubics is not total. We conclude $\gamma=6$, exactly. In fact
one must check that the 8 assigned basepoints impose independent
conditions on cubics, and so our pencil is forced to contain the
special reducible cubic just described. Independence is checked by
the usual stratification method, where one imposes more and more
conditions while verifying that each extra condition drops
 dimension by checking that the corresponding inclusion is
strict. The method seems to apply to our situation, and we
conclude $\gamma=6$ (with reasonable self-confidence). Of course
another detail that must be taken care is our somewhat tacit
supposition that the gonality (or the gonality sequence) do not
depend tremendously on the choice of smoothing. The classical
method  of small perturbation (Brusotti, etc.) asserts existence
of a curve effecting the assigned smoothings, but there is an
infinitude of choices for the coefficients. A priori the fine
gonality invariants are sensitive to the choice effected. Remember
that Brusotti's method relies on the fact that the initial curve
(thought of as a point in the discriminant hypersurface) has a
neighborhood consisting of several ``falde analytice'', i.e.  a
divisor with normal crossings each branch of which corresponding
to preserving a certain node. This explains the liberal way to
smooth away nodes of our initial configurations. Yet more
hazardous is the claim of a smoothing conserving the exact
location of all nodes. This remark hinders slightly the previous
argument made on figure 426.

$\bullet$ We next test the  invariant $(r,p)=(5,2)$, and within it
seek again a representative of maximum gonality $\gamma=r+p=7$.
Using the same device as above, we are inclined to look at an
octic (order $m=8$) with an interior oval kept protected from
intrusion of singularities. The smooth genus is then $\tilde
g=\frac{7\cdot 6}{2}=21$, but need be lowered down to
$g=(r-1)+2p=8$. We thus consider a distribution of 13 nodes
distant from the inner oval to produce the curve nicknamed 337 (on
Fig.\,\ref{Coppens:fig}, see also Fig.\,\ref{Gabard:fig} for a
larger depiction). This curve has $r=3$ (not 5 as desired!). This
means that I am a bad experimentalist, but the curve 337 is worth
looking at closer. Since $g=8$ by construction, and $r=3$ we have
$p=3$ (recall $p=\frac{g-(r-1)}{2}$). When projected from a point
on the inner oval the curve is 7-gonal. This degree is {\it
larger\/} than Gabard's bound $r+p=6$! The example seems to
violate Gabard's bound $r+p$.

{\it Summary.}---While testing the fullness conjecture, we rather
arrived to a counterexample to Gabard's estimate $\gamma\le r+p$.
We thus switch slightly of game, but try to keep in mind the
fullness problem for later.

\subsection{Potential counterexamples to Gabard 2006 ($\gamma\le r+p$)}

[24.10.12] The curve just discussed (337) seems a potential
violation of the theorem $\gamma\le r+p$ asserted in Gabard 2006
\cite{Gabard_2006}. Can we solve this
paradoxical situation? Either Gabard's bound $r+p$ is false or
something wrong happened. A possible explanation is that we were
too cavalier when claiming $\gamma=7$; in fact the total lines
pencil on curve 337 just shows  $\gamma\le 7$. A priori there
might be optical illusion about evaluating gonality. For instance
sweeping our octic with 13 nodes by a pencil of cubics with 9
basepoints located on the nodes gives a linear series of degree
$3\cdot 8-9\cdot 2=24-18=6$, rescuing the $r+p=6$ bound. Of
course, it is another story to  convince that such a map can be
chosen total. Thus 337 still represents a severe aggression
against $\gamma\le r+p$.

Similar  counterexamples (be they illusory or real) can be
manufactured in lower topological complexity. Starting with a
configuration of 3 conics, we conserve the deep nest, but keep the
maximum number of singularities in the periphery so as to lower
the genus as much as possible. Keeping 7 nodes unsmoothed, but
smooth away all others crossings in a sense-preserving fashion (to
ensure the dividing character of the curve), we obtain  the curve
215 with $r=2$ (on Fig.\,\ref{Gabard:fig}). Its genus is
$g=10-7=3$. Thus the genus of the half (complex locus split by the
real one) is $p=\frac{g-(r-1)}{2}=\frac{3-1}{2}=1$. Projecting
from the interior oval gives $\infty^{1}$ total maps of degree
$5$, and the hasty guess is that $\gamma=5$. Since $r+p=3$,
Gabard's bound $\gamma\le r+p$ looks again corrupted.

\begin{figure}[h]
\hskip-25pt \penalty0
    \epsfig{figure=Gabard.eps,width=132mm}

\caption{\label{Gabard:fig} (Pseudo?)-counterexamples to Gabard's
 bound $r+p$}
\end{figure}

However the curve at hand (215) having $g=3$ (and being dividing),
elementary knowledge of Klein's theory prompts that the canonical
map $C\to {\Bbb P}^{g-1}$, here ${\Bbb P}^2$, will exhibit the
curve as a ``G\"urtelkurve'', i.e. a quartic with two nested
ovals. Then the gonality is reevaluated as $\gamma=3$, and
Gabard's bound is
vindicated again (by the rating agency!).

Another way to argue, would be to take a pencil of cubics with 7
basepoints assigned on the 7 nodes and another basepoint on the
curve. The degree is then $3\cdot 6- 7\cdot 2- 1\cdot
1=18-14-1=3$. The bound  $r+p$ looks rescued again. Yet some hard
work is required to check total reality of a suitable pencil.
Perhaps there is some conceptual argument, else one really
requires tracing carefully the pencil after an educated guess of
where to place the extra assigned basepoint.

It is even possible to construct a quintic with ``visual''
gonality exceeding $r+p$. The cooking recipe is the same as above.
Start from a configuration of 2 conics and one line, keep the
inner oval while maximizing the number of peripheral
singularities. It results picture 214 on Fig.\,\ref{Gabard:fig}.
We see  $r=2$ real circuits. The genus is $g=6-3=3$ (3 nodes must
be subtracted), and thus $p=1$. The naive gonality seems to be
$4$, exceeding (hence violating) Gabard's bound $r+p=3$. Again to
resolve the paradox one can either argue via the canonical map
carrying the curve to a G\"urtelkurve, or find a total linear
series of lower degree. Here this would involve a pencil of conics
through the 3 nodes plus one assigned basepoint inside the deep
oval. The resulting series has degree $2\cdot 5 - 3\cdot 2- 1\cdot
1=10-6-1=3$, in agreement with the $r+p$ bound.

A drawback of figure 214 is that the 3 remaining nodes are nearly
collinear, rendering nearly impossible the depiction of the conics
pencil. (In reality the 3 nodes are not aligned, else the line
through them cuts the quintic in 6 points.) It is convenient to
consider rather a related quintic $214bis$ on
Fig.\,\ref{Pencil:fig}, where the line has penetrated the inner
oval (yet without destroying it). All invariants $r,p$ (as well as
the naive gonality) keep the same value as on the previous example
214. The new curve makes it easier to trace a total series cut out
by a pencil of conics, where the extra basepoint has been chosen
most symmetrically. Each member of it has beside the 4 assigned
basepoints (counting for $3\cdot 2+ 1\cdot 1=7$ intersections) 3
moving points which are permanently all real, as follows (only?)
through patient inspection of the picture.

\begin{figure}[h]
\hskip-45pt \penalty0
    \epsfig{figure=Pencil.eps,width=157mm}
    \vskip-35pt\penalty0

\caption{\label{Pencil:fig} Constructing a total pencil of conics
on a quintic}
\end{figure}

\noindent Fig.\,\ref{Pencil:fig} attempts to show various members
of the conics pencil, as well as the 3 mobile points of the
series. Those sections are depicted by the same letters e.g. 1,1,1
corresponds to the section by the ellipse invariant under symmetry
about the vertical axis. Ultimately, the whole figure has to be
extended by this symmetry, but this is better done mentally for
not surcharging the figure. The dynamics (circulation) is quite
tricky to understand, but  the motion looks much accelerated
(hence hard-to-follow) when the (red) curve crosses the basepoints
of the pencil. This is a bit if the particle motion would be much
accelerated by a gravitational black hole. Once the picture is
carefully analyzed, it is evident that all 3 mobile points stay
permanently real. Thus our quintic curve has gonality $\gamma\le
3$. (Gabard's bound is rescued on this simple example.)

Perhaps similar miracles (via high-order pencils
hard-to-visualize) produce for all other pseudo-counterexamples to
$\gamma \le r+p$. Yet this probably requires considerable work
even just for the previous curves (of Fig.\,\ref{Gabard:fig}). In
full generality,  some divine act of faith is required to imbue
with chimeric respect the last vestiges of truth imputable to
Gabard's result. Note that pencil of cubics are required for the
examples (Fig.\,\ref{Gabard:fig}): even our sextic with 7 nodes
achieves, via conics, only gonality $2\cdot 6- 4\cdot 2=4$, not as
economic as the $r+p$ bound.

{\it Summary.}---Two scenarios are possible: either Gabard's bound
$\gamma\le r+p$ is false (which is not quite improbable as its
proof is intricate and there is an infinite menagerie of potential
counterexamples), or it is true in which case it might be proved
extrinsically by a highbrow extension of the last example
described (Fig.\,\ref{Pencil:fig}). This brings us to the next
section, which albeit not very tangible in our fingers is perhaps
technically implementable (at least at the level of Ahlfors bound
$\gamma \le r+2p$).

\subsection{Brill-Noether-type (extrinsic) approach
to Ahlfors via total
reality}\label{sec:Brill-Noether-approach-to-Ahlfors}

[26.10.12]
Neutralizing all virtual counterexamples (of the previous section)
to $\gamma\le r+p$ amounts  a sort of high-powered Brill-Noether
theory for totally real pencils able to reprove Ahlfors theorem
$\gamma \le r+2p$ (and then optionally to corroborate Gabard's
$\gamma\le r+p$) in a purely synthetic way. This section touches
superficially this grandiose programme, we are quite unable to
complete.

Let us be more explicit. Any smooth projective curve (or, what is
the same, closed Riemann surface) embeds in ${\Bbb P}^3$. A
generic projection will realize the curve (like a knot projection)
as a nodal model in the plane ${\Bbb P}^2$ having at worst
ordinary double points. Specializing to real (orthosymmetric)
curves we get a model in the plane, on which one can hope to first
prove existence of a total pencil while evaluating the least
degree of such a pencil. This should amount considering adjoint
curves passing through the nodes so as to lower most the degree.
The procedure would be as follows.

Let $F$ be a bordered Riemann surface of invariant $(r,p)$. We
consider its Schottky double $C=2F$, interpreted as a real
orthosymmetric curve of genus $g=(r-1)+2p$ with $r$ real circuits.
Using a generic immersion in the plane gives a model $\Gamma_m$ of
the curve $C_g$ of order $m$ having $r$ real circuits, and a
certain number of nodes $\delta$. For simplicity assume the nodes
to be simple, though the more general situation must perhaps not
be excluded. We have of course $g=\frac{(m-1)(m-2)}{2}-\delta$.
Let $\Delta$ be the divisor of double points of $\Gamma_m$. (Those
can occur in conjugate pairs under complex conjugation.) Consider
in the complete linear system $\vert kH\vert:=\vert {\cal O}_{\Bbb
P^2}(k)\vert$ of all curves of degree $k$,  a linear pencil $L$ of
curves  passing through the nodes $\Delta$ of $\Gamma$ (adjunction
condition). The resulting series has degree $\le k\cdot m- 2\cdot
\delta$. In fact a better control must be possible. First $k$ has
to be chosen large enough so that the adjunction condition is
possible at all. Since $\dim \vert k H\vert=\binom{k+2}{2}-1$, the
integer $k$ is chosen as the least integer such that this
dimension exceeds $\delta$. Then we may have some excess
permitting to assign other (simple) basepoints.

Let us be even more explicit (we work first over the complex, for
simplicity). So assume given $C$ a curve of genus $g$. We look
first at the canonical embedding $\varphi\colon C \to {\Bbb
P}^{g-1}$. The image curve has degree $2g-2$. We manufacture a
plane model via successive projections from points chosen on the
curve. This
lowers the degree by one unit after each projection. We arrive
ultimately at a nodal model $\Gamma_m\in {\Bbb P}^2$ of degree
$m=(2g-2)-[(g-1)-2]=g+1$. Experimental study or an inspired guess
suggests considering adjoint curves of degree $k=m-3$. This value
is  calibrated so that our $k$-tics have enough free parameters to
visit all $\delta$ nodes of $\Gamma$. Indeed
\begin{align*}
\dim \vert k H\vert
=\binom{k+2}{2}-1=\frac{(k+2)(k+1)}{2}-1&=\frac{(m-1)(m-2)}{2}-1\cr
&\ge\frac{(m-1)(m-2)}{2}-g=\delta.
\end{align*}
We look at all curves of degree $k$ going through the nodes
$\Delta$ of $\Gamma$. Denote  $\frak d=\vert kH(-\Delta)\vert$ the
corresponding linear system, and let $\varepsilon$ be its
dimension. Obviously
$$
\varepsilon \ge \dim\vert kH\vert-\delta.
$$
(In fact since nodes of an $m$-tic impose independent conditions
upon adjoint curves of degree $m-3$ this is an equality. But we do
not this deep fact essentially equivalent to Riemann-Roch.) Both
displayed formulas show that $\varepsilon\ge 0$, and we may thus
impose to our $k$-tics to pass through $\varepsilon-1$ extra
points, while still moving inside a linear system of dimension
$\ge 1$ (a so-called pencil). This gives a pencil $L\subset \frak
d$ of degree
\begin{align*}
&\le k\cdot m- 2\delta -1 \cdot (\varepsilon-1)\cr
&\le k\cdot m- 2\delta-\Bigl[\binom{k+2}{2}-1-\delta\Bigr]+1\cr
&= k\cdot m- \delta-\binom{k+2}{2}+2\cr
&= k\cdot m+ g-\binom{m-1}{2}-\binom{k+2}{2}+2 \quad \textrm{[now
recall $m=g+1$]}\cr
&= (m-3) m+ (m-1)-2\binom{m-1}{2}+2\cr
&= (m-3) m+ (m-1)-(m-1)(m-2)+2\cr
&= (m-3) m+ (m-1)\underbrace{[1-(m-2)]}_{-(m-3)}+2\cr
&= (m-3) [m-(m-1)]+2\cr
&= m-1=g.
\end{align*}
This proves that any curve of genus $g$ admits a pencil of degree
$\le g$, which made basepoint-free induces a map of, eventually,
lower degree. (Of course our assertion fails when $g=1$, but true
otherwise granting some knowledge.) This ``degree $g$'' bound is a
bit sharper than the usual degree $g+1$ prompted by
Riemann(-Roch)'s inequality, but much weaker than the Riemann-Meis
bound $[\frac{g+3}{2}]$ for the complex gonality.
A natural wish is obtaining the Riemann-Meis bound via the above
strategy, hoping that special configurations of $\varepsilon-1$
points on the curve impose less conditions than expected, leaving
some free room for additional constraints lowering further the
degree. This is essentially what Riemann was able to do (at least
heuristically) via transcendental methods, and (exactly) what
Brill-Noether's theory is about at the pure algebro-geometric
level. Recall, yet, that both works apparently fail satisfying
modern standards, cf. e.g. Kleiman-Laksov 1972
\cite{Kleiman-Laksov_1972} and H.\,H. Martens 1967
\cite{Martens_Henrik_1967}, where the problem was not yet solved
apart via Meis' analytic (Teichm\"uller-style) approach.

At this stage, starts the difficulties. The  big programme would
be to adapt the above trick to  real orthosymmetric curves, in
order to tackle Ahlfors theorem. The latter prompts the bound
$g+1$ rather, but this little
discrepancy should not
discourage us. So in some vague sense a ``real'' Brill-Noether
theory is required, combining probably also principles  occurring
in Harnack's proof (1876 \cite{Harnack_1876}) of the after him
named inequality.

From the real locus $\Gamma({\Bbb R})$ one shall identify deep
nests, and it is favorable to choose them as the extra basepoints
to ensure total reality of the pencil we are trying to construct.
Then there is also a foliation on the projective plane induced by
the members of the pencil. Inside each oval, the foliation must
exhibit singularities (otherwise total reality is violated). In
fact total reality imposes  the foliation to be transverse to the
real circuits. Hence if there is no singularity we would have a
foliation of the disc which is impossible. Perhaps Poincar\'e's
index formula is also required. To be brief there is some little
hope that a very careful analysis of the geometry establishes
existence of a total pencil of degree $g+1=r+2p$, recovering so
Ahlfors result.

This would be pure geometry (or the allied devil of algebra)
without intrusion of either potential theory, neither
transcendental Abelian integrals, nor even topological principles.
Perhaps only elementary topological tricks are required to ensure
total reality by gaining extra intersections via a continuity
argument akin to Harnack's. This offers maybe another approach to
Ahlfors, yet it requires some deep patience. It looks perhaps
somewhat cavernous as (extrinsic) plane curves with singularities
are just a ``Plato cavern''-style shadow of the full Riemannian
universe.

If this dream of a synthetic proof of Ahlfors theorem is possible,
then it would be nice (if possible) to boost the method at the
deeper level of special groups of points to gain the sharper
Riemann-Brill-Noether-Meis sharp control upon the gonality, whose
real orthosymmetric pendant is expected to be the $r+p$ bound (of
Gabard).

Last, I know (only through cross-citations) the work of Chaudary
1995 \cite{Chaudary_1995-Thesis} where a real Brill-Noether theory
is developed. This probably helps clarifying the above ideas.

{\it Philosophical remark.}---Everybody
experimented
difficulties when playing with extrinsic models of Riemann
surfaces. A typical instance occurs with Harnack's inequality
$r\le g+1$, whose extrinsic proof (Harnack 1876
\cite{Harnack_1876}) is pretty more intricate than Klein's
intrinsic version (same year 1876 \cite{Klein_1876}) based on
Riemann's  conception of the genus. By analogy, one can predict
that any synthetic programme toward Ahlfors will ineluctably share
some unpleasant features of Harnack's proof. The substance of the
latter is a spontaneous creation of additional intersection points
forced by topological reasons, leading to an excess violating
B\'ezout.
Arguments similar to Harnack's
might be required to ensure total reality of a well chosen pencil.
Instead of being obnubilated by real loci (of both the curve and
the plane), it is sometimes fruitful to move in the ``complex
domain'' to understand better reality. A typical example is
Lemme~5.2. (in Gabard 2006 \cite{Gabard_2006}) about unilateral
divisors linearly equivalent to their conjugates. This was one of
the key in my approach to Ahlfors maps. Perhaps this lemma is also
relevant to the problem at hand ensuring total reality quite
automatically. In the series of adjoint curves $\frak d$, one then
imposes passing not through deeply nested ovals, but rather
through imaginary points all located on the same half. The
difficulty is of course showing existence of such a curve
intersecting the fixed one only along one half (unilaterality
condition), except eventually for some assigned basepoints (either
real or imaginary conjugate).

\subsection{Extrinsic significance of Ahlfors theorem}

[07.11.12] Another (less retrograde) desideratum is to explicit
the extrinsic significance of Ahlfors theorem for real algebraic
(immersed) plane curves. We touched this already in the Slovenian
section \ref{Open-RS-embed-in-C2:sec} but now a sharper idea is
explored. The point is delicate to make precise and  already quite
implicit in my Thesis (2004 \cite{Gabard_2004}, especially p.\,7
second ``bullet'') plus of course in Rohlin 1978
\cite{Rohlin_1978} (albeit the latter may never have been aware of
Ahlfors theorem). Today I discovered a certain complement which is
perhaps
worth presenting.

First Ahlfors theorem traduces in the following.

\begin{lemma}
Any real orthosymmetric (=dividing) algebraic curve admits a
totally real morphism to the line.
\end{lemma}

\begin{proof} The half of the dividing curve is a bordered
surface. By Ahlfors 1950, the latter tolerates a circle map, which
Schottky-doubled gives the required total map. For another proof
cf. e.g. the first half of Gabard 2006 \cite{Gabard_2006}.
\end{proof}

This pertains to abstract curves (equivalently Riemann surfaces)
but it acquires some extra flavor when the curve becomes concrete.
Of course the ontological problem of concreteness is that there
are plenty of ways for an abstract object to become concrete. Thus
concreteness is oft the opposite extreme of canonicalness.
Arguably, there is perhaps still a preferred ``Plato cavern''
namely the projective plane which can be used as an ambient space
where to trace all Riemann surfaces provided we accept nodal
singularities. Concretely this is done via generic projections
from a higher projective space (${\Bbb P}^3$ actually suffices to
embed any abstract curve), and then projecting down to the plane
${\Bbb P}^2$ gives a nodal model.
All this being pure synthetic geometry it transpierces matters
regarding fields of definition (A. Weil's jargon) and so adapts to
the reality setting. As yet just trivialities, but now we aim
interpreting synthetically the (non-trivial) Ahlfors theorem.

Starting from a real dividing curve in some projective space,
suitable projections exhibit a birational model, $C$, in the plane
as a nodal curve. Existence of a total morphism traduces into
 that of a total pencil, i.e. one all of whose member cut only
real points on the curve $C$, at least as soon as they are mobile.
 A priori basepoints may include
conjugate pairs of points. (A simple example arises when we look
at the pencil of circles through 2 points. Recall that circles
always pass to the so-called cyclic points at $\infty$, but
 this is just an affine conception).

In
extrinsic terms, Ahlfors theorem takes essentially the following
form.

\begin{theorem} {\rm (IAS=Immersed Ahlfors via Kurvenscharen)}
Given a dividing (real algebraic) curve $C$ immersed nodally in
the plane ${\Bbb P}^2$. There is a totally real pencil of
(auxiliary) curves  of some order $k$, all of whose members cut on
$C$
solely  real points plus eventually imaginary conjugate pairs of
basepoints.
\end{theorem}

\begin{proof}
This reduces to the basic theorem that any abstract morphism of
algebraic geometry admits a concrete description in terms of
ambient linear systems when the abstract object is projectively
concretized. In substance this is just the spirit of Riemann
(algebraic curves=Riemann surfaces) but extended to the realm of
morphisms. So the required theorem is just basic algebraic
geometry but I forgot all the foundations. Historically add to
Riemann, certainly Cayley-Bachach, Brill-Noether, (Klein?), all
the Italians, and finally Weil, Grothendieck,  plus of course many
others.
\end{proof}

Now the new observation
[07.11.12] is that we may always assume $k=1$ (in the theorem IAS)
up to changing of birational nodal model. The idea is that we may
first reembed the curve $C$ via the complete linear system of all
curves of degree $k$ (alias Veronese embedding) in some higher
space ${\Bbb P}^N$, where $N=\binom{k+2}{2}-1$. Then the image
curve $C'$ is (totally) swept out by a pencil of hyperplanes
corresponding to the original total pencil $L$ of $k$-tics in the
plane ($k$-tics=curves of degree $k$). If we project from the base
locus of the hyperplane pencil which is a linear variety of
codimension 2 we arrive down again in ${\Bbb P}^2$, but now with a
new model
total under a pencil of lines. It seems to me
that this trick works and we get the:

\begin{theorem} {\rm (IAP=Immersed Ahlfors via lines pencils)}
Given an abstract dividing (real algebraic) curve, there is always
a  nodal(ly immersed) model in the plane ${\Bbb P}^2$ which is
total under a pencil of lines.
\end{theorem}

This permits to remove one of the obstruction in our discussion of
the Forstneri\v{c}-Wold problem (already touched in
Sec.\,\ref{Open-RS-embed-in-C2:sec}). We now deduce the stronger
assertion:

\begin{cor}
Any finite bordered Riemann surface immerses in ${\Bbb C}^2$.
\end{cor}

\begin{proof} Let $F$ be the bordered surface, and $C:=2F$ be
its Schottky double which is real orthosymmetric. By the theorem
(IAP) we find a nodal model in the plane ${\Bbb P}^2$ total under
a pencil of lines. The pencil being real its unique basepoint $p$
is forced to be real. Since the allied morphism (projection) is
total the fibre of an imaginary point is an unilateral divisor,
i.e. confined to one half of the curve. This means that all
imaginary lines through the basepoint cuts unilaterally the curve.
It suffices thus to remove (from ${\Bbb P}^2({\Bbb C})$) an
imaginary line through $p$ to obtain an immersed replica of $F$ in
${\Bbb C^2}$. Note that if $p$ lies on the (nodal) curve then only
the open half (interior of $F$) is so embedded, but we can
probably arrange this by displacing slightly the center of
perspective $p$ outside the curve while conserving total reality.
The net bonus is  that the whole bordered surface (boundary
included) is in ${\Bbb C}^2$.
\end{proof}

Of course this is still millions of lightyears away from
Forstneri\v{c}-Wold
postulated embedding (for all finite bordered surfaces), yet
represents already a nice application of Ahlfors. Of course the
corollary is also the special (finitary) case of the famous
Gunning-Narasimhan theorem (1967 \cite{Gunning-Narasimhan_1967}),
immersing any open Riemann surface in ${\Bbb C}^2$. Maybe their
immersions are proper also, whereas ours are not. Maybe the
Fatou-Bieberbach trick arranges this issue always, cf. e.g.
Forstneri\v{c}-Wold 2009 \cite{Forstneric-Wold_2009}. Anyway using
the quantitative form of Ahlfors (not used as yet) one can go
perhaps further, maybe saying things on the degree of the model.
Note also that the viewpoint of nodal model of orthosymmetric
curves affords another numerical invariant, namely:

\begin{defn} {\rm (quite implicit in
Matildi 1945/48 \cite{Matildi_1945/48})} Given an abstract
dividing real curve $C$. The least order $\delta$ of a nodal
birational model of $C$ is termed (by us) the nodality of the
curve $C$. Via Schottky-doubling this invariant also makes sense
for finite bordered Riemann surfaces.
\end{defn}

Projecting down to ${\Bbb P}^2$ the canonical model in ${\Bbb
P}^{g-1}$ of a curve of genus $g$, we get a nodal model of degree
$g+1=(2g-2)-[(g-1)-2]$ (each projection from a point on the curve
decreases the degree by one unit). Hence $\delta \le g+1$.

If the theorem (IAP) is correct, one could also try to define the
linear gonality of a bordered surface (or the allied
orthosymmetric double) as the least degree of a nodal plane model
totally real under a pencil of lines. This gives perhaps yet
another invariant $\lambda$, which seems to satisfy $\gamma\le
\lambda+1$.

Another  dream of longstanding (Gabard's Thesis 2004
\cite{Gabard_2004}) is whether Ahlfors' theorem implies Rohlin's
inequality $r\ge m/2$ for a smooth dividing curve of order $m$. If
such a curve $C=C_m$ is total under a pencil of lines, then
sweeping out the curve by the pencil gives collections of $m$ real
points. When rotating the line around the basepoint, those $m$
points never enter in collision (else smoothness is violated), nor
do they disappear in the imaginary locus (else total reality is
violated). After a 180 degree rotation already, the line returns
to its initial position while the group of $m$ points recover its
initial position giving raise to a monodromy permutation. Total
reality forces each circuit of the curve $C$ to be transverse to
the foliation underlying the pencil of lines. It follows that the
monodromy transformation is an involution (order 2) and we deduce:

\begin{lemma} {\rm (Rohlin essentially)}
Let $C_m$ be a smooth real curve of order $m$ totally real under a
pencil of line. Then the real locus $C_m({\Bbb R})$ consists of a
deep nest of $m/2$ ovals when $m$ is even, and if $m$ is odd there
is as usual one pseudoline and ovals distribute in a nest of depth
$(m-1)/2$.
\end{lemma}

In particular Rohlin's inequality $r\ge m/2$ follows in this
special case where total reality is given by a pencil of lines.
The general case of Rohlin still appeals to some formidable work,
but perhaps may be derived via a linear pencil on a nodal model.
Alas we are unable to complete this project.

Let us however try to be more explicit. Given a smooth dividing
$C_m$. Let $L$ be a total pencil of $k$-tics given by Ahlfors
(theorem IAS). Then one can either try to study directly the
corresponding foliation appealing to Poincar\'e's index formula,
and hope to mimic the above argument. Alternatively one can try to
use the reembedding trick, where we use another model total under
a pencil of lines. Now on the new nodal model of degree say
$\lambda$, we apply the same sweeping procedure. We see on one
initial line $L_0$ (assumed generic, i.e. avoiding the nodes)
$\lambda$ points all real. When rotating by a half-twist the line
we see groups of $\lambda$ points which now may cross themselves,
but one can still assign a monodromy permutation. Naively any
point finishes its trajectory on the other side of the basepoint
(alas this makes no sense since a projective real line is a circle
not disconnected by a puncture). The number of real circuits $r$
of the curve $C$ is equal to the number of cycles of the monodromy
permutation, but a priori the latter number can be very low since
crossings are permitted. (Imagine e.g. a spiral which after
 several growing revolution times closes up to form a single
circuit.) Note that  we do not yet exploited the smoothness
hypothesis of the original model $C_m$. A naive way to exploit
this is via the complex gonality $\gamma_{\Bbb C}$. We have indeed
$m-1=\gamma_{\Bbb C}\le \gamma$ ($C$ being smooth). On the other
hand $\gamma\le \lambda-1$. Hence $m\le \lambda$. This is
interesting yet certainly not enough to conclude Rohlin's
inequality. So we give up the question for the moment.

\subsection{The gonality spectrum}

An idea perhaps worth exploring is to enrich the gonality sequence
(Definition~\ref{def:gonality-sequence}) into what could be called
the   {\it gonality spectrum\/}. This would just be the former
weighted by the dimension of the space of all circle maps having
the prescribed degree.

As we already observed earlier (hyperelliptic examples) it seems
that when a surface has a very low gonality then it ``somnolates''
without creating new gonalities. Thus more generally, the
intuition behind this spectrum invariant would be a conservation
law somewhat akin to Gauss-Bonnet: whatever the Riemannian
incarnation of a topological surface the curvatura integra keeps
constant value equal to the Euler characteristic ($ \int_F K
d\omega= 2\pi \chi(F)$).

Of course experiments requires to be made (using e.g. the
specimens on Fig.\,\ref{Coppens:fig}). Alas I had not presently
the time to do serious investigations about this spectrum. It
seems also expectable that from a certain range on, the spectrum
is independent from the conformal structure. (At least so is the
case for the gonality sequence which is always full after $r+2p$.)

Of course some convention is required, probably consider only maps
up to automorphisms of the disc.

Example the only example where the spectrum is very easy to
describe is the disc: in this case the $\gamma$-sequence is full
starting from 1, and there is essentially only one map of degree
one (the Riemann map). Given any unilateral group $D$ of $d$
points in the disc, thought of as the north hemisphere of the
Riemann sphere the pencil through $D$ and its complex conjugate
$D^{\sigma}$ induces a totally real map. (cf. Lemme 5.2 in Gabard
2006 \cite{Gabard_2006}). Conversely, given the map its fibre over
0 gives an unilateral divisor, which up to a range automorphism
may be assumed to contain 0. Normalizing by a rotation there are
thus the map depends upon $2d-3$ real constants. (Make this more
precise\dots). Such maps are (in the complex function literature)
often called finite Blaschke products.

Once the setting is well understood, this gonality spectrum
encodes valuable information upon all circle maps. Of course one
perhaps still want to know more; e.g. to understand the incidence
relation among the varied maps, especially how high-degree maps
may degenerate to lower degree ones. Fig.\,\ref{Coppens:fig} shows
some interesting examples. Considering e.g. picture 313 we see
that both maps of degree 3 are limit of maps of degree 4 (actually
can be connected by such), and both of them are also limits of
maps of degree 5.

Looking at picture 112 (again on Fig.\,\ref{Coppens:fig}) we see
that the unique (total) map of degree 2 is also the limit of maps
of degrees 3 and 4. The gonality sequence $2,3,4,\dots$ can be
enriched by weighting by dimensions to get $2_0,3_1, 4_2, \dots$.
Beware that probably there are other maps of degree 4 than those
visible on the picture as linear projection, namely the unique
2-gonal map post-composed by circle maps of degree 2 from the disc
to itself. Our guess is that such Blaschke maps may degenerate to
their originator (the hyperelliptic projection) but not to maps of
degree 3.

\subsection{More lowbrow counterexamples to $\gamma\le r+p$}

[27.10.12] We now pursue the project of multiplying and
diminishing further the order of virtual counterexamples to
Gabard's estimate $\gamma\le r+p$ (cf. Fig.\,\ref{Gabard:fig} and
Fig.\,\ref{Pencil:fig}). There we found curves (via an uniform
recipe) seemingly violating the gonality upper bound $r+p$. The
simplest example had order 5, but it is easy to get examples of
order 4. The game is again to depict total pencils vindicating
Gabard's bound. Albeit very modest corroboration of the bound, we
found instructive to visualize the corresponding total pencils.

First remind the general recipe: to manufacture an (at least
virtual) counterexample to $\gamma\le r+p$, we leave
tranquil the inner oval but maximize the number of singularities,
so as to lower the genus $g=(r-1)+2p$, and hence $(r,p)$. Having
left quiet the inner oval the virtual gonality via linear
projection is one less than the degree, but $r+p$ may go lower
down this value.

\begin{figure}[h]
\centering
    \epsfig{figure=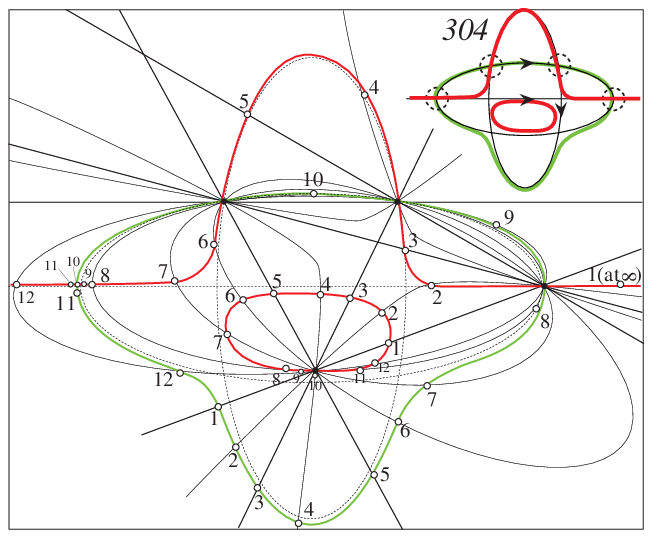,width=122mm}

\caption{\label{F304:fig} Tracing a totally real pencil of conics
on an orthosymmetric quintic, whose underlying bordered Riemann
surface has genus zero. This illustrates the Ahlfors map (rather
its Bieberbach-Grunsky special case) at the extrinsic level. The
 circulation of real points is (violently) accelerated  when the cutting conics
nearly osculate the cutted quintic.}
\end{figure}

We first consider a configuration of order 5 consisting of 2
conics plus one line, see picture 304 below
(Fig.\,\ref{F304:fig}). Smoothing it as dictated by orientations
while keeping unsmoothed the dashed circles gives a curve with
$r=3$ real circuits of genus $g=6-4=2$. Hence
$p=\frac{g-(r-1)}{2}=0$. The virtual gonality
is $\gamma^\ast=4$ (projection from the inner oval). This seems to
violate $\gamma\le r+p=3$. Looking at the pencil of conics through
the 4 nodes gives a series of degree $2\cdot 5-4\cdot 2=10-8=2$.
This violates the trivial bound $r\le \gamma$, but of course this
pencil is not total: e.g. the conic consisting of the 2 horizontal
(or better oblique) lines misses the inner oval. Assigning instead
one of the 4 basepoints on the inner oval gives a pencil of degree
$3$, which is claimed to be total.

\noindent Totality of the morphism requires examining (patiently)
that each conic of the pencil cuts only real points on the quintic
$C_5$. This is depicted on the large part of Fig.\,\ref{F304:fig},
where each triad of moving points of the series are labelled by
triples $1,1,1$, then $2,2,2$, etc. Let us start from the conic
consisting of the oblique line through $1,1$, plus the horizontal
line. The latter cut the red pseudoline at infinity. This pair of
lines deforms to a hyperbola cutting the triad $2,2,2$. This
hyperbola is in turn pinched toward a pair of lines cutting the
group $3,3,3$, etc, up to $7,7,7$. From here on, things becomes
harder to visualize. (Alas our picture is not optimally designed.)
The conic of the pencil now becomes very close to the primitive
conic involved in the generation of the quintic $C_5$ via small
perturbation. The net effect is that points on the green branch
nearly ``osculated'' by the primitive ellipse are (violently)
accelerated (like in CERN's particles accelerator). At this stage
it is quite delicate to make a consistent picture, but total
reality seems to work: all particles stay real during the motion
without disappearing as ghost in the imaginary locus (as conjugate
pairs of  points under Galois).

We promised a similar example of degree 4; this will be pictured
later (Sec.\,\ref{sec:degree-four}), being now sidetracked to
another topic which looks more exciting.

\section{Some crazy ideas about gravitation
and unification of forces}
\label{sec:gravitation}

\subsection{From gravitation to electrodynamics}\label{sec:electrodynamics}

Now we arrive at the following crazy interpretation (discovered
the 27.10.12 at ca. 13h58). It would be nice if there is some
relation of the Ahlfors maps with periodic solutions to the $n$
body problem in gravitation (celestial mechanics). The 4
basepoints of Fig.\,\ref{F304:fig} may be thought of as
supermassive black-holes, so massive that there is no interaction
between them (imagine purely static objects lying in different
sheets of the multiverse). Dually, the moving points of the linear
system are imagined as massless microparticles  (electrons, or
better photons). There is also no gravitational interaction
between them. Thus the sole interactions reigning are those
between black holes and photons. It is also imagined that a photon
can traverse a black-hole (without captivation).

As a wild speculation,  the trajectories described by the 3
photons on Fig.\,\ref{F304:fig} may satisfy exactly Newton's law
of gravitation. In particular the full trajectory would be the
real locus of an algebraic curve! This would of course be a wide
extension of Kepler's law (on the r\^ole of conic sections in the
simplest case of one sun and one planet).

If this is true we see a deep connection between Klein's
orthosymmetric curves, Ahlfors maps of conformal geometry and the
totally real circulations positing periodic stable motions along
circuits of an orthosymmetric curve. Exaggerating a bit this
should explain the ultimate constitution of matter (and its
relative stability) not via knots (as Lord Kelvin desired via
Helmholtz vortices) but via bordered Riemann surfaces (probably
quite ubiquitous already in the so-called string theory).

Note that our basic experiment (with Fig.\,\ref{F304:fig}) is---as
far as  speed of motion is concerned---quite in line with this
interpretation.

\begin{figure}[h]
\centering
    \epsfig{figure=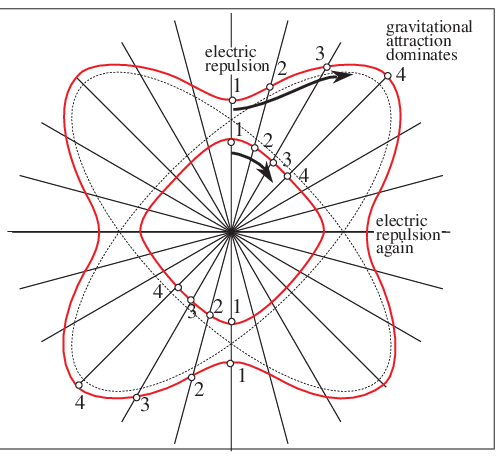,width=72mm}
    \vskip-10pt\penalty0

\caption{\label{FGuert:fig} An electro-gravitational toy model
with Klein's G\"urtelkurve. Four particles of electronic type
(electrons) are gravitating around a single star (better a
proton). Unusual shapes of trajectories are explained by electric
repulsion.}
\end{figure}

Let us look at one of the simplest example of orthosymmetric
curve, namely the (Zeuthen-Klein) G\"urtelkurve (aka {\it courbe
annulaire}). This is a quartic with two nested ovals arising by
smoothing two transverse ellipses having 4 intersections. The
picture is given below (Fig.\,\ref{FGuert:fig}). One can
convincingly argue that the shapes of trajectories (especially the
outer oval) are unlikely to be gravitational orbits. It seems that
some hidden force repulses the particles (labelled 1 on the
figure).
Invoking some other (electric) force effecting repulsion between
particles, then the  trajectories of the G\"urtelkurve look again
physically tolerable. Thus the ``physical'' model should include
two types of interactions: gravitational and electromagnetic.

Of course one can drag the position of the sun while still having
a totally real pencil. This gives the next figure
(Fig.\,\ref{FGuert2:fig}). Note that we did not changed the curve,
yet it is still plausible that for suitable initial conditions
(velocity vectors) the orbits of our 4 bodies follows exactly the
same quartic curve.

\begin{figure}[h]
\centering
    \epsfig{figure=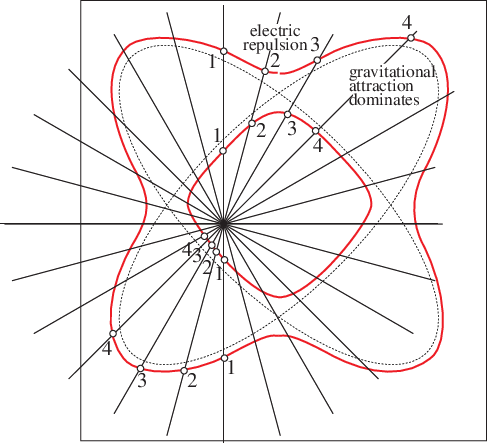,width=72mm}
    \vskip-5pt\penalty0

\caption{\label{FGuert2:fig} Moving the sun.}
\vskip-10pt\penalty0
\end{figure}

We arrive at the following metatheorem [14h57]:

\begin{theorem} {\rm (Kepler generalized?)} \label{metatheorem:thm}
Given any orthosymmetric real (algebraic) curve embedded (or
immersed) in the Euclid plane ${\Bbb R^2}$ and a totally real
pencil (existence ensured by Ahlfors theorem). There exists
initial conditions (velocity vectors) such that the trajectories
of  particles
obeying
the inverse square law of
Newtonian  attraction resp. Coulombian repulsion
match exactly
the real circuits (``ovals'') of the given real algebraic
curve. Further the dynamics (speed of
motion) is dictated by
the pencil. In particular there is plenty of periodic solutions to
the $n$-body problem, essentially one for each such curve.
\end{theorem}

How to prove this? Philosophically,
algebraicity
might be not so surprising: recall Laplace's potential-theoretic
interpretation of Newton, and from Laplace there is just one step
to Riemann, hence to Klein. The miracle should be essentially akin
to Riemann's existence theorem prompting  any closed Riemann
surface (an a priori completely fluid object) to rigidify
canonically as an algebraic curve. Even if true the metatheorem is
quite modest because in practice (meteorites, apocalyptic black
holes scenarios, etc.) one is given the initial conditions and the
goal is to predict the future evolution of the system. Here in
contrast, we know in advance the trajectories (hence the destiny)
while claiming existence of initial conditions compatible with the
orbital structure. Generally, integrating the differential
equations governing some motion, we meet a highly complex
dynamical system subjected to the paradigm of chaotic determinism
\`a la Poincar\'e. Note that a Euclidean model of the projective
curve is required  to give sense to Newton's inverse square law.

Several questions naturally occurs assuming the truth of the
metatheorem. The theorem affords plenty of periodic motions.
Essentially we obtain as many periodic motions as there are real
orthosymmetric curves. Even more than that, one requires an
Ahlfors circle map (equivalently a totally real morphism \`a la
Klein-Teichm\"uller). A first naive question is: do this recipe
exhausts all periodic motions? Certainly not, try Euler and
Lagrange's periodic motions. Roughly all algebraic motions are
periodic, but the converse has no chance to be true.

Observationally, Fig.\,\ref{FGuert2:fig} looks anomalous because
the series 1,2,3,4 closest to the sun looks much slowed down,
whereas we are accustomed (Kepler) to rapid motions near a massive
star. One requires perhaps a third type of interaction, say the
{\it strong interaction\/}, to explain this. Namely both particles
the
proton and the electron are of a dualistic nature, hence they tend
to ``love'' themselves like
partners staying close together over a long period of time. This
third force would have the net effect of diminishing the real
speed by a factor proportional to the (squared?) distance
separating the bodies. What is then the fourth force, alias {\it
weak interaction\/} in contemporary physics? Maybe none is
required in our model? Perhaps dually, particles of the same
nature (namely electrons) dislike themselves like competitors and
the {\it weak force\/} just produces  some acceleration of the
motion when they are in close vicinity. Visually this behavior is
perhaps observed near the groups labelled 2,3 on the top part of
Fig.\,\ref{FGuert2:fig}.

We have now a model with 4 fundamental forces. One must of course
still define time. This would, on our example, just be the angular
parameter of the pencil. Presumably the metatheorem should take
into account these two extra forces, becoming somewhat
sophisticated, yet probably still completely deterministic and
hopefully reasonably easy to integrate. The miracle would be that
it admits
dividing (=orthosymmetric) real curves as periodic orbits.
Of course to relativize, one can do similar games with real
diasymmetric curves, but then there is no total reality prompted
by Ahlfors theorem and particles sometimes disappear in the
imaginary locus. We leave to the reader's imagination appropriate
physical interpretations (ghost particles, anti-matter, etc.)

Perhaps there is a more elementary way to explain slowness of the
motion near the star (without appealing to the exotic forces at
the subatomic level). Recall Kepler's law in the elliptic case,
that identic sectorial areas are swept out during
the same amount of time. This suggests that the time parameter is
not the angular parameter but the areal one. Of course one gets
other troubles since the distant
electron is supposed to move synchronously
with the one
closer to the
proton (cf. Fig.\,\ref{FGuert:fig}).

\subsection{Some little objections}

[28.10.12] Another objection to our metatheorem
(\ref{metatheorem:thm}) is the following one. Assume the given
orthosymmetric curve to be of the simplest stock, namely a line
swept out by a total pencil of lines. Then one must assume that
there is no forces between the two bodies to explain the
rectilinear motion.

A more serious objection arises when $C$ is an ellipse swept out
by a total pencil of lines through the middle of both foci. If all
(four) fundamental forces involved satisfy the inverse square law,
then so does the resulting force. Hence all interactions reduce to
a single one which is attractive (to get an elliptic trajectory).
However according to Kepler the orbit must be an ellipse with the
sun located at one of the foci. Hence our geometric model where
the basepoint of the total pencil lies at the center of the
ellipse is not physically relevant.

This example suggests that the metatheorem requires corrections.
Maybe one is only given in advance the orthosymmetric curve but
not the total pencil, while the metatheorem states existence of a
pencil physically observable. For an ellipse we would only be
allowed to take pencil of lines through one of both foci; if the
ellipse degenerates to a circle only the center would be
permissible. [30.12.12] At this stage it might be relevant to
remind that there is a vast theory of foci for high-order
algebraic curves, due it seems to Pl\"ucker first and then Siebeck
1864 \cite{Siebeck_1864}, etc. cf. e.g. Casas-Alvero 2013
\cite{Casas-Alvero_2013}.

Of course this Kepler obstruction should not preclude physical
systems obeying more complicated interactions laws with say
several fundamental forces, maybe not all subsumed to the inverse
square law. Such could validate exotic orbital structures, e.g.
 an ellipse with a sun at its center, as physically
reasonable.

\begin{figure}[h]
\centering
    \epsfig{figure=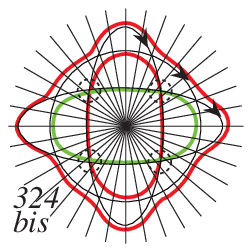,width=42mm}
    \vskip-10pt\penalty0

\caption{\label{F324bis:fig} A totally real pencil of lines
without repulsion (rebounding instead)}
\end{figure}

[28.12.12] Another possible objection comes from the following
curve Fig.\,\ref{F324bis:fig}.
This possesses a total
pencil of lines, yet some particles do not repulse, rather
crossing themselves
unsensitive
of each other. On relabelling the particles one can posit a
repulsion acting rather as a
bounce of billiard balls (elastic shock). Another
explanation could involve some quantumchromodynamics like
assigning spins to the electrons neutralizing interaction between
some of them.

\subsection{Back to insignificant geometry}

\begin{figure}[h]
\centering
    \epsfig{figure=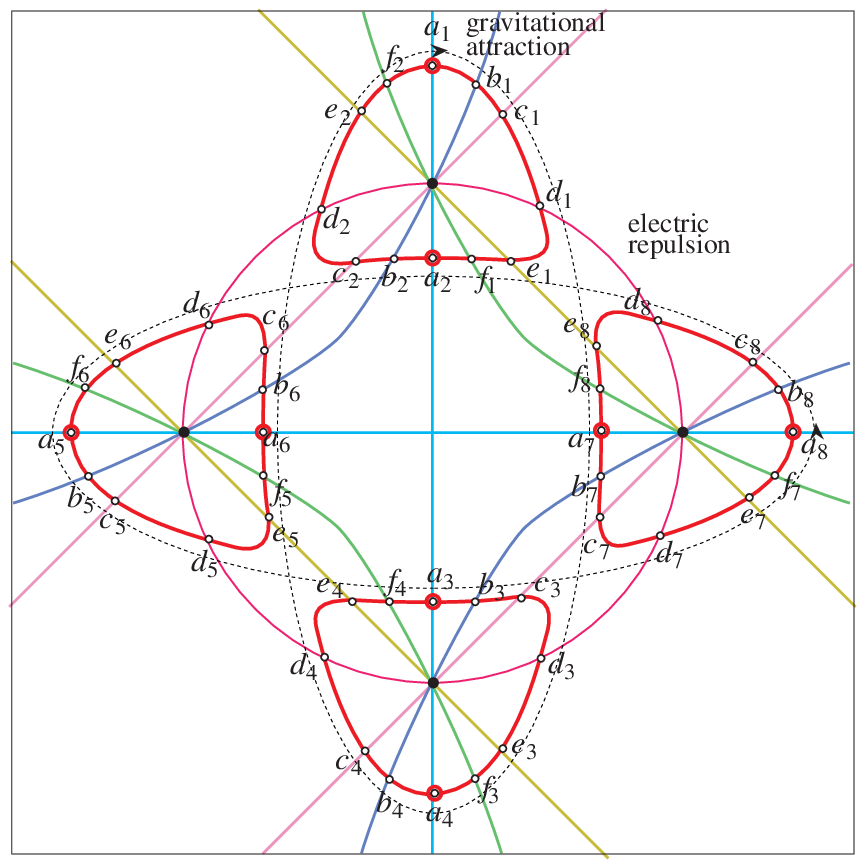,width=82mm}
    \vskip-10pt\penalty0

\caption{\label{F4oval:fig} Harnack-maximal quartic with a total
pencil of conics traced through 4 basepoints located inside the
ovals. Physically, this is a magneto-gravitational system with 4
stars and 8 electrons. The initial position is labelled
$a_1,a_2,\dots, a_8$ whose forward orbit consists of $b_i, c_i,
d_i,\dots$ (increasing alphabetic order).}
\end{figure}

Let us now leave such complex modelling question, to contemplate
more complicated systems arising from other curves than the
G\"urtelkurve, especially some of higher order. First staying of
order 4 there is, dual to the G\"urtelkurve, the curve arising by
reversing orientation of one of the ellipses (cf. arrows on
Fig.\,\ref{F4oval:fig}). This gives a quartic with 4 ovals when
smoothing compatibly with the prescribed orientations. A total
pencil arises from all conics through 4 basepoints distributed
inside the ovals (Fig.\,\ref{F4oval:fig}).

Initially the point $a_1$ animated by a suitable horizontal
velocity vector is mostly subjected to the attraction of the
nearby star (=upper basepoints of the conics pencil). If $a_1$ and
this star were to be alone in the universe, $a_1$'s orbit  would
be close to the dashed ellipse of ``vertical eccentricity'',
provided the upper star coincides with the focus of this ellipse.
Yet in reality, as the body $a_1$ arrives near position $d_1$ and
meanwhile body $a_8$ reached position $d_8$, electric repulsion is
becoming predominant
causing a (finally violent)
deviation from the elliptic trajectory.

Instead of appealing to gravitation one can just imagine the
basepoints (alias ``stars'' previously) as positively charged
protons, the whole system  reducing to an electrodynamical one
obeying only Coulomb's law of attraction resp. repulsion. The
fixed protons would however not repulse, maintaining their fixed
positions due to  some nuclear cohesion (strong/weak forces).

\begin{figure}[h]
\centering
    \epsfig{figure=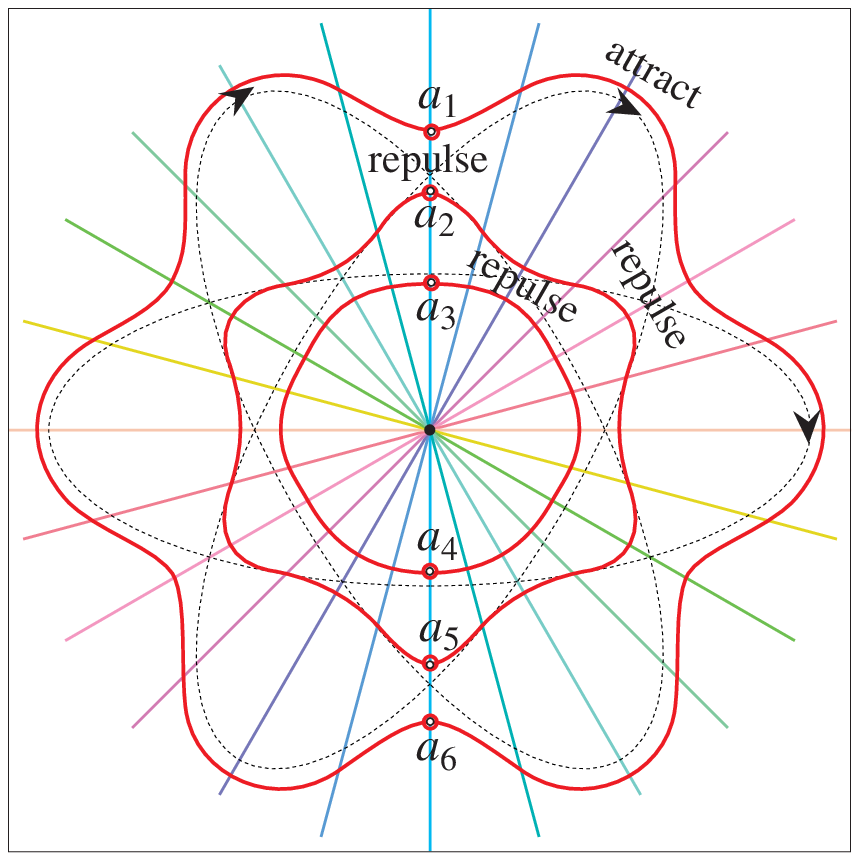,width=62mm}
    \vskip-10pt\penalty0

\caption{\label{Fsextic:fig} A sextic akin to the G\"urtelkurve
(nest of depth 3).}
\end{figure}

It is easy to produce examples of higher topological complexity
via curves of higher orders. Instead of starting with two
ellipses, take three of them and smooth the configuration in a
sense-preserving way to get Fig.\,\ref{Fsextic:fig}.

\begin{figure}[h]
\centering
    \epsfig{figure=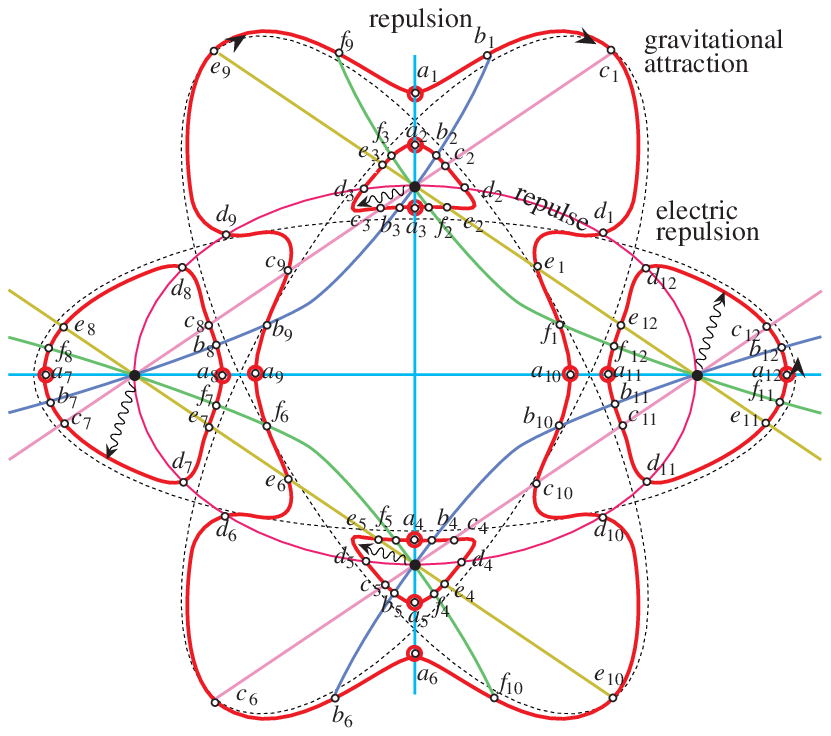,width=82mm}
    \vskip-10pt\penalty0

\caption{\label{Fsext2:fig} A sextic with a total pencil of conics
through 4 deep nests. Physico\-chemically, an atomic nucleus
consisting of 4 protons is gravitated around by 12 electrons
dancing around in a fairly complicated way. The long orbit is
circulated by 4 electrons (all others by 2). Is there some
relation between Poincar\'e indices of the foliation and the
number of moving points?}
\end{figure}

Reversing orientation of one of the ellipses (say that with
horizontal major axis) gives the more interesting
Fig.\,\ref{Fsext2:fig} requiring a pencil of conics to exhibit
total reality.

Again we use the same labelling as before, namely the first
(cyan=pale blue colored) conic consisting of the vertical and
horizontal lines cuts on the sextic $C_6$ the group of points
labelled $a_1,\dots, a_{12}$, all of them being real. Moving
clockwise from the top, a subsequent conic (blue colored) cuts the
series denoted $b_i$, etc. One checks easily
all conics of the pencil to cut only real points on the $C_6$.
Looking at the corresponding dynamical process, we note that $a_1$
is first repulsed against $a_2$ being rejected  as far as $c_1$,
then attraction of the 4 protons (mostly the North and East one)
track back the orbit to position $d_1$ where a repulsion against
$d_{12}$ takes place, deflecting again the orbit along the way of
the North proton, but then vicinity of $d_2$ causes another
repulsion towards $e_1$ and $f_1$, which is finally gently
repulsed by $a_{11}$. etc. The sequel of the story reproduces
symmetrically.

It is now fairly evident how to construct similar dynamical
systems of ever increasing complexity. It may be observed that the
totally real map induced by the pencil gives a circle map of
degree 12. Now the topological invariants are $r=5$ and $g=10$.
Hence the half-genus is $p=\frac{g-(r-1)}{2}=3$. Hence this map
has degree exceeding Ahlfors bound $r+2p=11(=g+1)$. However a
parietal degeneration of the 4 basepoints against the ovals
immediately enclosing them (cf. ``squigarrows'' on
Fig.\,\ref{Fsext2:fig}) exhibits a total map of degree $2\cdot
6-4\cdot 1=12-4=8$. This is actually in accordance with the $r+p$
bound predicted in Gabard 2006 \cite{Gabard_2006}.

It is tempting to consider the (mildly singular) foliation induced
by the pencil (of conics). It seems clear from the picture that
there is a relation between the sum of Poincar\'e indices extended
to the interior of an oval and the number of points circulating on
the oval. Observe also that the foliation is transverse to  the
boundary of the disc bounding the oval. This property is general
and follows at once from the fact that totally real maps lack real
ramification points. Using Ahlfors total reality paradigm combined
maybe with Poincar\'e's index formula we suspect that some old
(and perhaps new?) information on the topology of real plane
(dividing) curves can be re-derived. In particular we suspect that
it must be possible to recover Rohlin's inequality. This states
$r\ge m/2$, i.e. any smooth dividing plane curve of order $m$ has
at least so many circuits as the half value of its order. This is
a fantastic project, but we leave it aside for now. [vague details
p.\,32 of hand-notes].

[08.11.12] Another highbrow (yet poorly explored)
application of Ahlfors theorem was sketched in Gabard's Thesis
(2004 \cite[p.\,7]{Gabard_2004}). This was an answer to Wilson's
question (1978 \cite[p.\,67]{Wilson_1978}) on deciding the
dividing character of a plane curve by sole inspection of its real
locus. Here again Ahlfors theorem affords an answer:  a real curve
is dividing iff it admits a total pencil (with possibly imaginary
conjugate basepoints). Yet it must be admitted that the answer,
albeit perfectly geometric, has probably little algorithmic value
unless complemented by further insights. Of course another
question is to decide the dividing character from the sole data of
a ternary form (homogeneous polynomial in 3 variables with real
coefficients). The simplest case of Wilson's question is that of a
deep nest, i.e. a smooth curve $C_m$ of say even degree $m=2k$
with a completely nested collection of $k$ ovals. Then linear
projection from a point on the deepest oval is total of degree
$m-1$. Since the complex gonality is also $m-1$, we deduce that
the gonality $\gamma$ is also $m-1$. On the other hand the
topological invariants are $r=k$ and $g=\frac{(m-1)(m-2)}{2}$.
Hence in this case Ahlfors bound $r+2p=g+1$ is strongly beaten by
the gonality $\gamma=m-1<\!\!<g+1=[1+2+3+\dots+(m-2)]+1$. Gabard's
bound $r+p$ is also much greater than the exact $\gamma=m-1$;
indeed $r+p$ is nothing but the mean value of $r$ and $g+1$ and in
the case at hand the former is $m/2$ but the latter is quadratic
in $m$.

\begin{figure}[h]
\centering
    \epsfig{figure=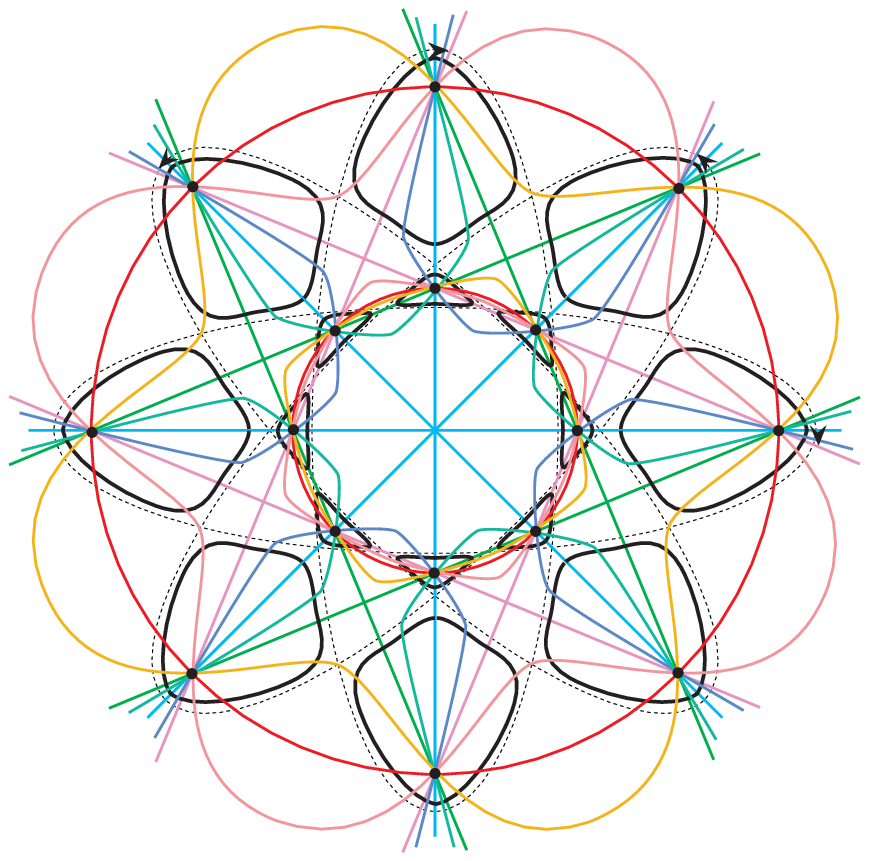,width=72mm}
    \vskip-10pt\penalty0

\caption{\label{Foctic:fig} An octic (reminiscent of gothical
rosewindow)  swept out by a total pencil of quartics with 16 real
basepoints distributed in the innermost ovals}
\end{figure}

[29.10.12] We consider next an octic (Fig.\,\ref{Foctic:fig})
arising from a sense-preserving perturbation of 4 ellipses rotated
by 45 degrees. Of course if all ellipses are oriented clockwise we
get a nest of depth 4 and accordingly a total pencil of lines
through the innermost oval. Here instead, we reverse some
orientations to create 16 ovals and no nesting (cf. black curve on
Fig.\,\ref{Foctic:fig}). The theorem of Ahlfors predicts existence
of a total pencil. The general principle is to impose basepoints
inside the deepest ovals, hence the desired pencil must have
degree 4. At this stage depiction can be a fairly difficult
artform (reminiscent of gothical ``rosaces''= rosewindows). Our
trick was to use a ground ellipse of pretty large eccentricity so
that oblique line of (angular) slope different from $\pi/4$ (the
green and lilac colored ones) also passes through the deep nests.
Of course this trick is not supposed to affect the generality of
the method (i.e. Ahlfors theorem) but just intended to simplify
the artwork!

\begin{figure}[h]
\centering
    \epsfig{figure=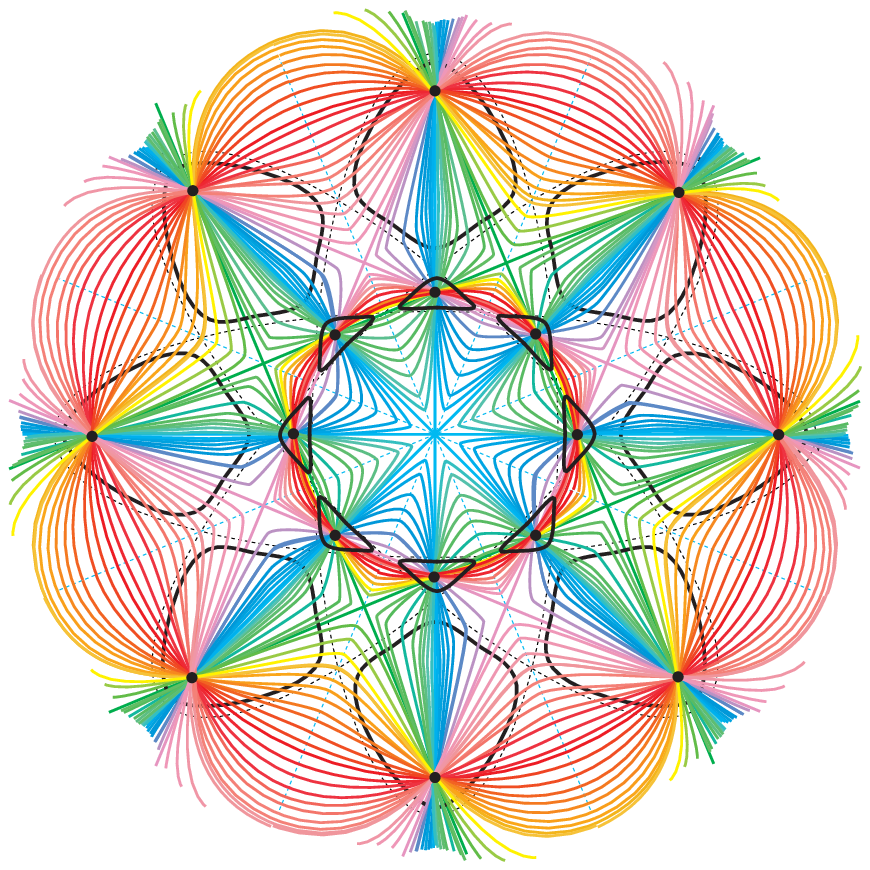,width=72mm}
    \vskip-10pt\penalty0

\caption{\label{Foctic2:fig} The rainbow effect}
\end{figure}

As to the arithmetics, recall that (plane) quartics depends upon
$\binom{4+2}{2}-1$ parameters (coefficients counting), hence one
is free to assign 13 basepoints. On the other hand, our dividing
octic $C_8$ has genus $g=\frac{(m-1)(m-2)}{2}=\frac{7\cdot
6}{2}=21$ and $r=16$ ovals, thus the genus of the half (semi
Riemann surface) is $p=\frac{g-(r-1)}{2}=3$. Imagine now that
among all 16 basepoints of the pencil 13 moves against the ovals,
then a series of (reduced) degree $4\cdot 8- 13\cdot 1=32-13=19$
is obtained. This matches with the $r+p$ bound on the degree of
circle maps predicted in Gabard 2006 \cite{Gabard_2006}.
Geometrically it is pleasant to observe that certain members of
the pencil are G\"urtelkurven (see the lilac-colored curve). Those
are not connected. Hence total reality of a pencil is not
necessarily allied to connectedness of the auxiliary curves. For
the fun of depiction, one can increase the number of curves of the
pencil while sweeping out more and more of the full color
spectrum, creating a sort of rainbow effect (cf.
Fig.\,\ref{Foctic2:fig}).

At this stage one gets the impression that the theory (or rather
the pictures) works only for highly symmetric patterns. However
the strength of Ahlfors result lies in its universal validness for
all curves regardless of symmetry. This imbues some suitable
respect plus a certain feeling of vertigo about the whole Ahlfors
result.

\begin{figure}[h]
\centering
    \epsfig{figure=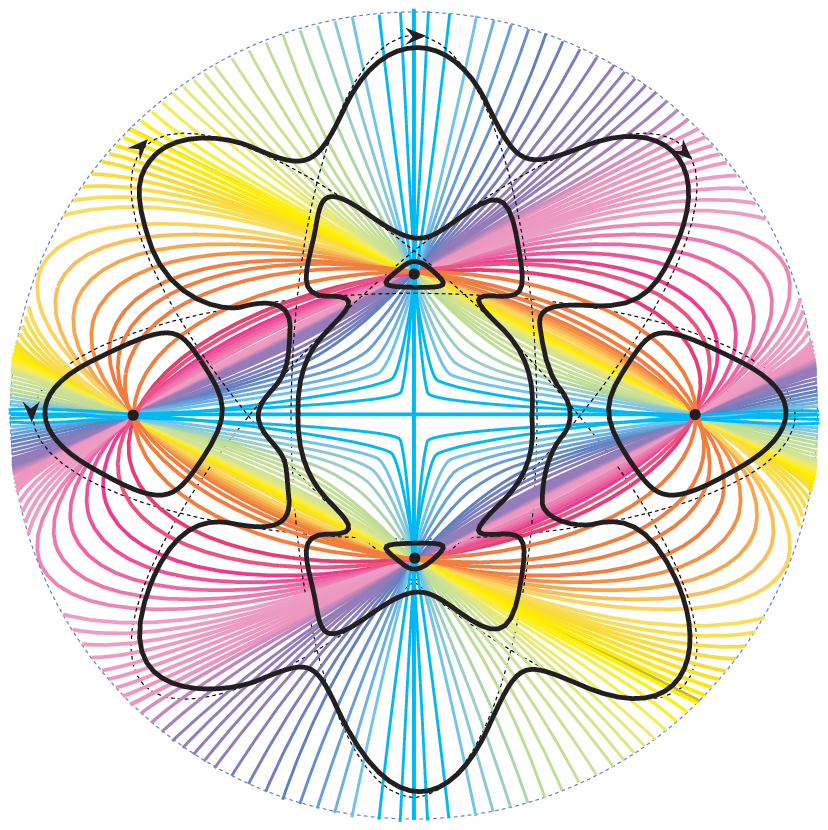,width=72mm}
    \vskip-10pt\penalty0

\caption{\label{Foctic3:fig} Another octic with a total pencil of
conics}
\end{figure}

Of course there is another possible orthosymmetric smoothing of
our configuration of 4 ellipses. This is given by reversing one of
the orientations of the ellipses, and we obtain the black-traced
curve on Fig.\,\ref{Foctic3:fig}. This times there is only 4 deep
nests and a pencil of conics suffices to exhibit total reality.

As to arithmetic matters, this octic has still $g=21$ but now only
$r=6$ ovals. Hence the semi-genus $p=\frac{g-(r-1)}{2}=8$.
Dragging the 4 basepoints against the deep ovals gives a total map
of degree $2\cdot 8-4\cdot 1=16-4=12$. This is more economical
that the $r+p$ bound, here equal to $14$.

\begin{figure}[h]
\centering
    \epsfig{figure=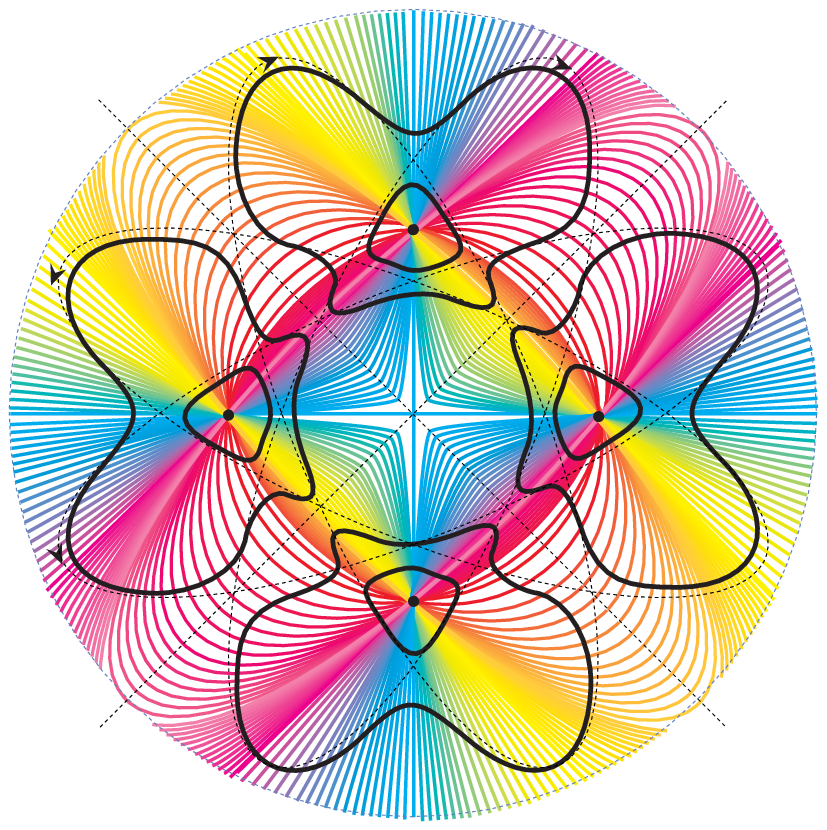,width=72mm}
    \vskip-10pt\penalty0

\caption{\label{Foctic4:fig} Yet another octic with a total pencil
of conics}
\end{figure}

Finally there is yet another  smoothing of our 4 ellipses
producing Fig.\,\ref{Foctic4:fig} with 4 nests of depth 2. A
pencil of conics suffices to show total reality.

Regarding the topological invariants we have $r=8$, hence $p=7$.
As before there is a total map of degree 12 (via parietal
degeneration), which is better that Gabard's bound $r+p=15$.
Naively this relative improvement over the previous example (in
comparison to the $r+p$ bound) could be explainable by the higher
symmetry of the new curve probably reflecting a further
particularization of the ``moduli''. (Recall that if one believes
Gabard 2006 \cite{Gabard_2006} and especially Coppens 2011
\cite{Coppens_2011} a bordered surface of type $(r,p)$ has
generically gonality $\gamma =r+p$.)

Note yet that our total pencil of conics persists for any octic
with 4 nests of depth 2, hence the symmetry of the pattern can be
greatly damaged by large deformation of the coefficients without
affecting the (estimated) gonality. So we certainly have the:

\begin{prop} Any octic curve with $4$ nests of depth $2$ has gonality
$\gamma\le 12$ (and presumably not lower, yet this remains to be
elucidated).
\end{prop}

\begin{figure}[h]
\centering
    \epsfig{figure=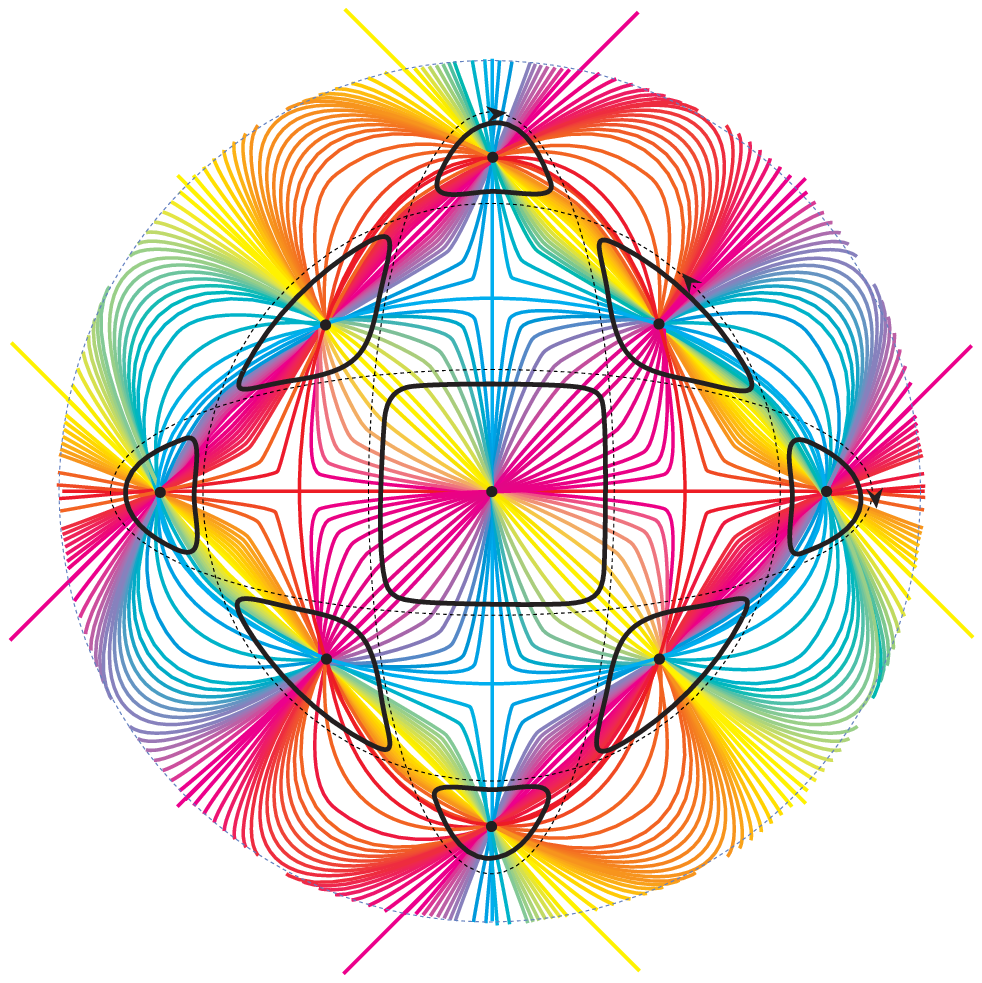,width=82mm}
    \vskip-10pt\penalty0

\caption{\label{Fcubic:fig} A sextic with a total pencil of
cubics}
\end{figure}

Having clearly exhausted the smoothing options of our 4 ellipses,
one is somehow disappointed that pencils of cubics were  not yet
required. Looking on p.\,7 of my Thesis \cite{Gabard_2004} I
rediscover a simple such example involving only a sextic. Let me
reproduce this with the rainbow technology. We start now from a
configuration of 3 ellipses one of which is a circle and get
Fig.\,\ref{Fcubic:fig}.

The sextic has $g=10$ and $r=9$ (hence pre-maximal amond dividing
curves), and thus $p=1$. Cubics depends on
$\binom{3+2}{2}-1=10-1=9$ parameters, hence 8 basepoints may be
freely assigned. Pushing them along ovals gives a total map of
degree $3\cdot 6-8\cdot 1=10$. This matches with Gabard's bound
$r+p$, hence the curve should be considered has having general
moduli. Of course if the smoothing is done very symmetrically and
if moreover we play with the radius of the initial circle, we can
perhaps arrange that all 9 basepoints lands on the sextic curve in
which case the gonality would lover to 9 the minimum value (recall
$r\le \gamma$).

Starting from the above sextic, one can perform a large
deformation of the coefficients staying inside the space of all
smooth sextic curves. The real locus picture may then undergo
drastic change of shape yet its topological type keeps unaltered
and so in particular the orthosymmetric character. It is not clear
anymore that our simple minded pencil of cubics (spanned by 2
pairs of 3 lines) suffices to exhibit total reality. This amounts
essentially to the claim that for any 8 basepoints distributed
among the ovals then the ninth basepoint luckily falls into the
remaining one. This luckiness phenomenon becomes even more
hazardous when it comes to vindicate Gabard's bound by a
synthetical procedure. The latter seems equivalent to the claim
that given any such curve (orthosymmetric with 9 non-nested ovals
it is always possible to choose 8 points one on each oval) so that
the pencil through them creates an extra basepoint inside the
remaining oval. This lucky-stroke phenomenon should perhaps be
further explored either as an application of the $r+p$ bound or as
a way to disprove it.

[08.11.12] Let us fail to be more specific as follows. Remember
first that a real sextic curve with 9 unnested ovals needs not to
be dividing, cf. e.g. Gabard's Thesis 2004 \cite{Gabard_2004}
p.\,8, but this is of course well-known since at least  the
Rohlin-Fiedler era, e.g. Rohlin 1978 \cite{Rohlin_1978}. Second it
is not even clear a priori that the conditions ``dividing plus 9
unnested ovals'' specifies a unique
rigid-isotopy type of curves, i.e. a unique
chamber in the space of all smooth sextics. This is a projective
space of dimension $\binom{6+2}{2}-1=28-1=27$ parcelled into
chambers by the discriminant hypersurface of degree $3(m-1)^2=3
\cdot 5^2=75$. (Inserted [24.01.13]: However this is true by a
deep result of Nikulin 1979 \cite{Nikulin_1979/80}.)

\begin{conj} (very hypothetical!!)
Any dividing sextic with $9$ unnested ovals  admits a total pencil
of cubics with $8$ basepoints on the sextic and the $9$th
basepoint inside the remaining oval.
\end{conj}

\begin{proof} (pseudo-proof!) Since the curve  is dividing we know by Ahlfors that
there is a total pencil. We have very poor control on the degree
of the curves of the  pencil. We only know Ahlfors bound
$r+2p=g+1=11$, Gabard's one  $r+p=10$ and the complex gonality
$\gamma_{\Bbb C}=5$ which is completely useless. Stronger
information comes from the trivial bound $r\le \gamma$. So the
gonality $\gamma$ is fairly well squeezed as $9=r\le \gamma \le
r+p=10$. A priori a least degree total map could be given by a
pencil of quartics. Then the degree could be as low as $4\cdot
6-16=24-16=8$; for quintics as low as $5 \cdot 6-25=5$; for
sextics as low as $6\cdot 6-36=0$; septics $7\cdot 6-49=-7$;
$k$-tics $k\cdot 6- k^2$ highly negative! Hence we have virtually
no control on the degree of (members of)  a total pencil, despite
the bounds on the degree of the abstract total map. Let us thus
shamefully postulate that the pencil in question can be chosen
among cubics. For foliated reasons it is clear that the nine
basepoints (elliptic points or ``foyers'' of Poincar\'e index
$+1$) must be surjectively distributed among the 9 ovals. Indeed
the total pencil is transverse to the real circuits and the disc
bounding an oval cannot be foliated transversely (Euler-Poincar\'e
obstruction). Hence we have the:

\begin{lemma}
All basepoints of a total cubics pencil on a smooth sextic with
$9$ unnested ovals  are real, distinct, and
surjectively(=equitably) distributed between the $9$ ovals (either
in their insides or their periphery).
\end{lemma}

Applying the parietal degeneration trick we can take any 8 of the
basepoints and drag them to the ovals. During the process we get
new pencils (of possibly jumping dimension?) while the 9th
basepoint could a priori escape its enclosing oval. The difficulty
looks so insurmountable that we have to abort the project.
\end{proof}

In fact the following principle is worth noticing. It gives a
basic lower bound on the degree of total pencils, yet as we saw
the real difficulty is rather upper bounds! As a matter of
annoying nomenclature crash, note that the degree of the pencil is
not that of the allied map but that of its constituting curves, so
we should perhaps rather speak of the order of a (total) pencil.

\begin{lemma}
{\rm (Poincar\'e-style lower bound on the order of total
pencils)}\label{Poincare-lower-bound} Given
a (smooth) (dividing) plane curve with a total pencil of $k$-tics
with $D$ many deepest ovals (i.e. the minimal elements of the
nesting ordered structure). Then $D\le k^2$ or $k\ge \sqrt{D}$.
\end{lemma}

\begin{proof} Each deep oval must enclose at least one singularity
of the foliation. Remember that the latter is transverse to the
curve by total reality. Poincar\'e's index formula (1882/85) says
that the sum of all indices equates the Euler characteristic.
Applied to the disc bounding a deepest oval this forces the latter
to encloses at least one singularity of index +1. Warning: one
must explain why the disc could not be foliated by say two
singularity of index $1/2$, so-called thorn singularities. The
pencil has at most $k^2$ singularities of the foyer type
(index=+1) materialized by the basepoints. Thus $D\le k^2$. Indeed
for each deepest oval chose one foyer inside it. We get a map from
 the set of deepest oval to that  of basepoints, which is
injective since the deepest ovals are disjoint at least for a
smooth curve. Try to clarify if smoothness is really required as a
hypothesis!
\end{proof}

[30.10.12] Let us look at another intriguing example. Start again
with 2 ellipses invariant under rotation by 90 degrees, and add a
concentric circle as the dashed one on Fig.\,\ref{Fcubic:fig}, but
shrink its radius slightly beyond the critical radius where the
circle passes through the 4 intersections of the 2 ellipses.
Smoothing this configuration along our choice of arrows gives
Fig.\,\ref{FcubicA:fig}: a sextic with $r=9$ ovals one of them
enclosing all others.

\begin{figure}[h]
\centering
    \epsfig{figure=FcubicA.eps,width=92mm}
    \vskip-10pt\penalty0

\caption{\label{FcubicA:fig} The archipelago: a sextic arising by
smoothing 2 ellipses $E_1,E_2$ plus a circle of radius $R$ pinched
between the distance of $E_i$ to the origin and that of $E_1\cap
E_2$ to the origin.}
\end{figure}

The picture has the annoying property that ovals are pretty small,
 challenging a bit the visual perception of homo habilis. Since
the curve is dividing, Ahlfors theorem predicts the existence of a
total map. It is evident that no pencil of lines, nor of conics,
is total. (This is either optically clear or deduced from
Poincar\'e's bound $k\ge \sqrt{D}=\sqrt{8}=2.828\dots$, i.e.
Lemma~\ref{Poincare-lower-bound}). The 8 deep ovals prompts
seeking among pencil of cubics. Of course we may just assign 8
basepoints inside those deep ovals and hope for total reality. Yet
to manufacture a concrete picture it is natural to assign
basepoints in the most symmetric way. Once this is done one try to
identify special singular curves passing through the 8 points. We
find 4 degenerate cubics consisting of a line plus a conic (cf.
colored curves on Fig.\,\ref{FcubicA:fig}). Once those are
detected it is an easy matter to interpolate between them (by
continuity) to trace a qualitative picture of the pencil
(Fig.\,\ref{Fcubic3:fig}).

\begin{figure}[h]
\centering
    \epsfig{figure=Fcubic3.eps,width=92mm}
    \vskip-10pt\penalty0

\caption{\label{Fcubic3:fig} The archipelago sextic swept out by a
total pencil of cubics with 8 basepoints in the deep ovals and one
extra basepoint at the origin}
\end{figure}

This archipelago sextic $C_6$ has $g=10$ as usual, and $r=9$, thus
$p=1$. The total pencil can be lowered to degree $3\cdot 6- 8\cdot
1=10$, as predicted by the $r+p$ bound.

[08.11.12] Again several  questions poses themselves naturally.
(The sequel uses some jargon of Rohlin 1978 \cite{Rohlin_1978},
for instance the {\it real scheme\/} of a smooth plane real  curve
is the isotopy class of the embedding of its real locus in the
real projective plane):

(1) Is any sextic $C_6$ belonging to the real scheme of the
archipelago
 (i.e. 8 unnested ovals altogether surrounded by an outer oval) of
dividing type? The answer is probably known to Rohlin and his
students, especially  if there is a nondividing counterexample?
[Update 24.01.13: yes there is one and this was well-known at
least since Rohlin 1978 \cite{Rohlin_1978}, yet his article was
far from explicit when it comes to constructions. Personally I
understood this point only after reading Marin 1979
\cite{Marin_1979}, compare our Fig.\,\ref{GudHilbMarin:fig} much
below (virtually copied from Marin), which is  a (clever) variant
of Hilbert's method of vibration.]

Rohlin distinguishes real schemes as definite or indefinite
depending on whether all its representatives belong to the same
type or not, w.r.t. Klein's dichotomy (ortho- vs. diasymmetric).
(cf. Rohlin 1978 \cite{Rohlin_1978})

(2) What is the exact gonalities occurring in this archipelago
scheme (of course restricting attention to dividing models in case
the scheme is indefinite)? If we believe in Gabard's bound
$\gamma\le r+p$, we have $9=r\le \gamma \le r+p=10$.

Perhaps answers are to be searched along the following direction.
Maybe it is true that for any 8 basepoints (injectively)
distributed in the 8 deepest ovals the corresponding cubics pencil
is total. On counting intersections, we get roughly $8\cdot 2=16$
many coming from the 8 deep ovals and the outer oval should also
contributes for 2 intersections. This is at least evident if the
real part of the cubics are connected since the real circuit of
each such cubic has to go at ``infinity'' (in the sense of moving
outside the outer oval, for otherwise it would be contractible
inside the bounding disc of the latter, whereas we know the cubic
circuit to realize an ``odd'' nontrivial class in the fundamental
group $\pi_1({\Bbb R}P^2)$ or just the allied homology). On the
other hand, the cubic circuit must also visit the 8 assigned
basepoint inside the outer oval, and so is forced to intercepts
the latter. We arrive at a total of 18 real intersections, the
maximum permissible by B\'ezout ($3\cdot 6=18$). Total reality
would follow.

I remind vaguely of a standard result claiming that for a generic
collection of 8 points there is a pencil of rational (hence
connected) cubic interpolating them. (Cf. e.g. Kharlamov-Degtyarev
survey ca. 2002). Now if all this is true, the archipelago scheme
is dividing, and any such curve admits plenty of total cubics
pencil of degree $3\cdot 6- 8\cdot 1=10$ (essentially one for each
selection of 8 points on the deep ovals). It seems however hard to
lower the gonality $\gamma$ up to the absolute minimum $r=9$, but
I know no argument.

\subsection{Total reality in the Harnack-maximal case}
\label{sec:Total-reality-Harnack-max-case}

[08.03.13] Much of this section is by now much illuminated by
Le~Touz\'e's observation in Le~Touz\'e 2013
\cite{Fiedler-Le-Touzé_2013-Totally-real-pencils-Cubics}, where it
is remarked that a very simple prescription of basepoints ensure
total reality of a pencil of cubics on an $M$-quintic.

[ca. 31.10.12] Quite paradoxically it is much harder to depict
total pencils on Harnack-maximal curves, alias $M$-curves (in
Russia since Petrovskii 1938 \cite{Petrowsky_1938}, cf. Gudkov
1974 \cite[p.\,18]{Gudkov_1974/74}), especially when the order is
$m\ge 5$. (For lower orders $m\le 4$ everything is essentially
trivial: since $m=4$ just requires a pencil of conics passing
through the 4 ovals of the quartic (with $g=3$).) Recall indeed
that Ahlfors theorem is much easier in the planar case $p=0$,
where it goes back to Bieberbach-Grunsky, if not earlier.
Logically the argument simplifies much via Riemann-Roch and the
absence of collision, cf. e.g. Gabard 2006
\cite[Prop.\,4.1]{Gabard_2006} or
Lemma~\ref{Enriques-Chisini:lemma} above in this text.

Shamefully, the following section climaxes the poor level of
organization of the present text. Of course the game is quite
outside the main stream of our subject (Ahlfors theorem), yet we
think that some
 phenomena require to be clarified. In particular we
were not able to make any reliable picture of a total pencil on a
Harnack-maximal (smooth) plane curve of order $m\ge 5$. After some
three days of pictorial tergiversation  we found a sort of weak
obstruction to manufacturing such pictures involving a basic type
of pencil spanned by two special cubics. This obstruction is
described at the end of the section, which otherwise reduces to a
messy gallery of failing attempts of the desired easy depiction!
Yet the abstract theorem of Bieberbach-Grunsky implies the
existence of total pencil but they probably involve
delicate-to-visualize pencil of cubics (in the quintic case). We
would like to challenge gifted amateurs to picture them
appropriately.

Let us first recall the construction of such $M$-curves due to
Harnack (in the variant of Hilbert). We start with degree 5.
Consider as primitive configuration an ellipse $E_2$ plus a line
$L_1$. Take further 3 parallel lines $l_1,l_2,l_3$. There is some
psychological difficulties to know if we should first smooth $E_2
\cup L_1$ and then perturb along $l_1\cup l_2\cup l_3$ or if we
can directly perturb $E_2 \cup L_1$ without taking care of
smoothing. Let us adopt  the shorter route (actually so do
Hilbert) by putting directly $C_3=(E_2\cup L_1)+\varepsilon
\vartheta_3$. This cubic (in black thick stroke) oscillates across
the ellipse $E_2$ meeting it in the maximum number of 6 real
points. Next smoothing their union (=product) $C_3 \cup E_2$ we
get the (red-colored) quintic $C_5$ realizing the maximum number
$r=7(=g+1)$ of ovals (one of them being in fact a pseudoline i.e.
a Jordan curve in ${\Bbb R}P^2$ not bounding a disc).

\begin{figure}[h]
\centering
    \epsfig{figure=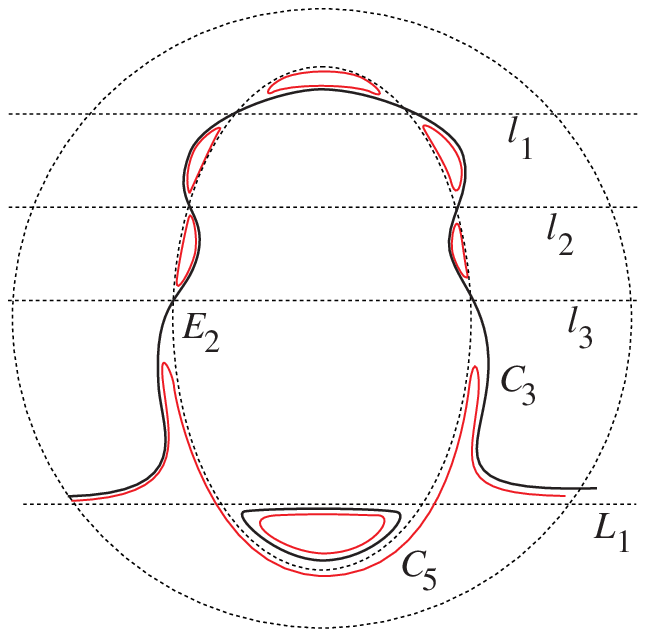,width=62mm}
    \vskip-10pt\penalty0

\caption{\label{Harna0:fig} Harnack-Hilbert oscillation trick
creating a quintic with $r=7$ circuits}
\end{figure}

\begin{figure}[h]
\centering
    \epsfig{figure=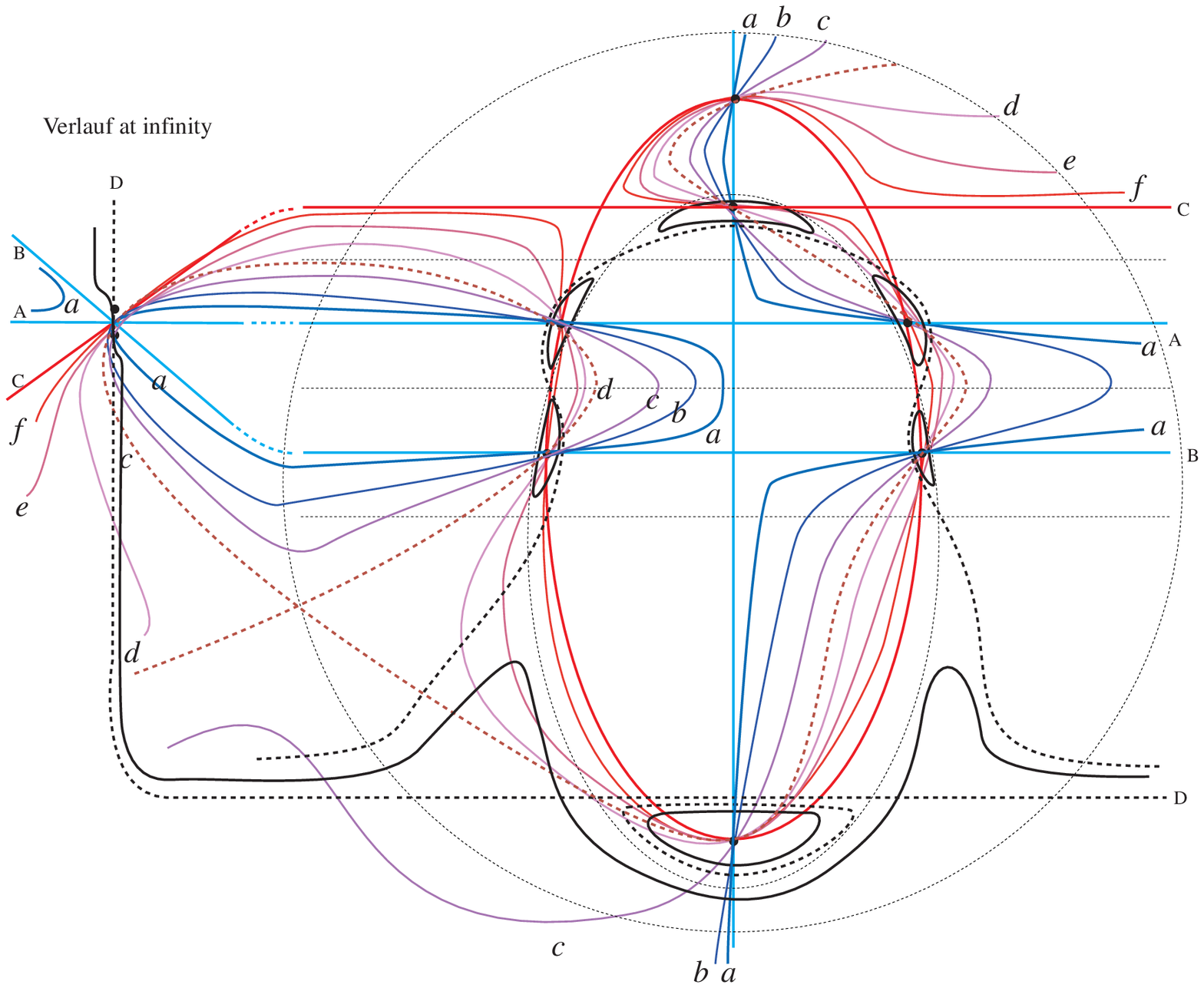,width=112mm}
    \vskip-10pt\penalty0

\caption{\label{Harna1:fig} Trying to construct a total pencil on
our $M$-quintic}
\end{figure}

[31.10.12] Now the (perpetual) game is to find a total pencil on
this dividing curve $C_5$ (recall that Harnack-maximal curves are
always dividing). As usual the recipe is to distribute imposed
basepoints $p_1, \dots, p_6$ in the deepest ovals. Those are fixed
once for all and marked by black points on Fig.\,\ref{Harna1:fig}.
Since there are 6 ovals, pencil of lines or conics are not
flexible enough to reveal the total reality of our $C_5$. We thus
have to look among pencils of cubics. In view of the (vertical)
symmetry of the curve $C_5$ it is natural to seek a symmetric
pencil. We shall define them by specifying two of its members. A
first vertically symmetric cubic through the 6 basepoints  is the
union of the 3 cyan-colored lines. This special (cyan) cubic $C_3$
cuts our quintic $C_5$ twice along each oval and once on the
pseudoline, hence in $12+1=13$ points. Those are at finite
distance but looking at infinity both horizontal cyan lines cuts
the pseudoline branch of $C_5$ in two extra points, yielding a
total of 15 point, the maximum possible (all of them being real).
Beside, we consider another vertically symmetric cubic, namely the
red-colored cubic $R_3$ consisting of the red ellipse through 5
points $p_i$ plus the red horizontal line (denoted $C$) through
the remaining $p_i$. We can now consider the corresponding pencil
spanned by the cyan and red cubics (equation $\lambda C_3+ \mu
R_3=0$). Unfortunately, the red cubic cuts $C_5$ along $2\cdot
6=12$ points on the ovals and only once at infinity. Indeed the
pseudoline branch of $C_5$ is asymptotic to the line $D$ which in
transverse to the red line $C$. Hence the intersection $R_3\cap
C_5$ is not totally real. Of course this defect does not prevent
us from tracing the corresponding (non-total) pencil.

{\it Note.}---A pencil of cubic may be defined by assigning 8
basepoints. By letting degenerate those against the 6 ovals (or
the pseudoline) we get a series of degree $3\cdot 5 - 8\cdot
1=15-8=7$ as predicted by Bieberbach-Grunsky (cf. e.g.
Lemma~\ref{Enriques-Chisini:lemma}). But it is far from evident to
ensure total reality. Of course a coarse calculation would
stipulate that the 6 ovals contributes for $2\cdot 6=12$ many
intersections and imposing 2 extra basepoints on the pseudoline
gives 2 additional intersection, totalizing 14 many hence the last
man surviving is forced to be real as well. This argument
certainly holds good if we know that all cubics of the pencil are
connected but a priori a cubic may well have an oval which could
be nested in one of the tiny ovals of our sextic. If so is the
case then this one cubic's oval only visits one of the 8
basepoints, without spontaneous creation of intersection on one
oval of the quintic $C_5$. Maybe this scenario is quite improbable
but I missed some argument.

A modest  improvement over our previous attempt is to take a
red-colored cubic satisfying total reality. This is given by
changing the red-colored ellipse by taking the one passing through
the 5 ``highest'' (relatively to our figure
Fig.\,\ref{Harna2:fig}) black-colored basepoints. Symmetry forces
us then to take an additional red-colored line passing through the
``lowest'' basepoint. We obtain the following
Fig.\,\ref{Harna2:fig}. Alas it is not evident that total reality
is satisfied.

\begin{figure}[h]
\centering
    \epsfig{figure=Harna2.eps,width=112mm}
    \vskip-10pt\penalty0

\caption{\label{Harna2:fig} Trying to construct a total pencil on
our $M$-quintic}
\end{figure}

A third option is to change the cyan configuration of 3 lines and
we get the following Fig.\,\ref{Harna3:fig}, which alas again
seems to fail total reality.

\begin{figure}[h]
\centering
    \epsfig{figure=Harna3.eps,width=82mm}
    \vskip-15pt\penalty0

\caption{\label{Harna3:fig} Failing again to construct a total
pencil on our $M$-quintic}
\end{figure}

[01.11.12] Of course we would like ultimately to extend the game
to sextic. Let us first reproduce a picture in Hilbert 1909
\cite{Hilbert_1909-Ueber-die-Gestalt-sextic}. The idea is again
that a union of two ellipses is vibrated into a quartic $C_4$
oscillating across one of the ellipse $E_2$ (which is a circle on
Fig.\,\ref{Hilb1:fig}, left), and next $E_2\cup C_4$ is smoothed
to a sextic with 11 ovals (compare Fig.\,\ref{Hilb1:fig}, right).

\begin{figure}[h]
\centering
    \epsfig{figure=Hilb1.eps,width=102mm}
    \vskip-10pt\penalty0

\caption{\label{Hilb1:fig} Hilbert's picture of an $M$-sextic:
just vibrate and smooth!}
\end{figure}

Again the challenge would be to trace a total pencil of curves on
this $C_6$. We have 10 deep ovals, thus  pencils of cubics look
overwhelmed already with their only 8 assignable basepoints (and
maximally 9 of them). Quartics have $\binom{4+2}{2}-1=15-1=14$
free parameters hence we can impose 13 basepoints. Choosing them
in the deep ovals and doing a parietal degeneration gives a series
of degree $4\cdot 6- 13 \cdot 1=24-13=11$. This matches with the
Bieberbach-Grunsky bound, however it is far from evident that
total reality is ensured.

In general if $C_m$ is a Harnack-maximal curve of order $m$, the
previous examples (with $m=5,6$) suggest to consider auxiliary
curves of degree $m-2$ forming a space of dimension
$\binom{(m-2)+2}{2}-1=\binom{m}{2}-1$ and thus assigning
$\binom{m}{2}-2$ basepoints will define a pencil. By parietal
degeneration the resulting series has degree
$(m-2)m-[\binom{m}{2}-2]$, and this is easily calculated as being
equal to
\begin{align*}
(m-2)m-[\binom{m}{2}-2]&=(m-2)m-\frac{m(m-1)}{2}+2\cr
&=\frac{1}{2}[2(m-2)m-m(m-1)+2]+1=\frac{1}{2}[m^2-3m+2]+1\cr
&=\frac{(m-1)(m-2)}{2}+1=g+1, \end{align*} where $g$ is the genus.
This again agrees with the Bieberbach-Grunsky theorem, but of
course does not reprove it, be it  just for the simple reason that
smooth plane curves have specialized moduli among all curve sof
the same genus. Still it would be exciting  to manufacture
tangible pictures of such total pencils in the planar case.

Now let us try again to do better pictures of the $M$-quintic. Any
such $M$-quintic has 6 ovals and one pseudoline. By B\'ezout no
three ovals can be aligned (otherwise 6 intersection with a line).
Thus the six ovals are somehow distributed along a configuration
resembling  a hexagon. This raises some hope to draw reasonable
pencil of cubics spanned by two configurations of 3 lines
according to one of the following patterns (left of
Fig.\,\ref{Hilb2:fig}). This suggested to draw another model whose
6 ovals are nearly situated like a regular hexagon. Using
cyclotomy, we get quickly the right part of Fig.\,\ref{Hilb2:fig}.

\begin{figure}[h]
\centering
    \epsfig{figure=Hilb2.eps,width=82mm}
    \vskip-40pt\penalty0

\caption{\label{Hilb2:fig} An $M$-quintic with hexagonal
distribution of ovals}
\end{figure}

A little piece of comment on the last Fig.\,\ref{Hilb2:fig}: of
course we started with a circle divided primarily in 6 equal
parts, and have chosen the 3 horizontal lines as passing through
the cyclotomic points. Those three lines are those used for the
Harnack-Hilbert vibration trick, and the rest of the picture
should be self-explanatory. Alas the bottom portion is quite
difficult to observe. Yet a clear-cut portrait of Lars Valerian
clearly emerges: the bottom oval is the mouth, then just above two
big eyes ``with an air of determination'', as well as some hairs
emanating from the beret. In fact the portrait looks more like an
alien, but the resemblance with Lars is much more flagrant when
the circle is depicted as a ``vertically  oblong'' ellipse. [I
apologize for adding some extra prose as otherwise the figures
desynchronize from the text.]

Now we consider the following pencil spanned by the cyan and red
collections of lines (Fig.\,\ref{Hilb3:fig}). Alas it fails to be
totally real, for it contains the green cubic cutting only 13 real
points on the quintic $C_5$. Of course the advantage of our pencil
is that it is simple to draw, yet its disadvantage is that it has
only 6 among the 8 assignable points located on the quintic.
Somehow one should try to conciliate both properties.

\begin{figure}[h]
\centering
    \epsfig{figure=Hilb3.eps,width=82mm}
    \vskip-25pt\penalty0

\caption{\label{Hilb3:fig} Total reality fails again!}
\end{figure}

Testing the other configuration (of 2 pairs of 3 lines through the
hexagon) one gets Fig.\,\ref{Hilb4:fig}. The situation is not much
improved. Now the 3 additional basepoints (intersection of pairs
of parallel lines) are ejected at infinity but are not lying on
the (black-colored) quintic curve $C_5$ whose pseudoline is
asymptotic to the horizontal line. The corresponding pencil of
cubics (spanned by the cyan and red colored lines) is probably not
total, for it should contain a nearly circular ellipse through the
hexagon plus the line at infinity, and the aggregated
corresponding cubic seems to cut the $C_5$ only along 12, plus one
at infinity, so a total of only 13 real points!?

\vskip-0pt
\begin{figure}[h]
\centering
    \epsfig{figure=Hilb4.eps,width=82mm}
    \vskip-10pt\penalty0

\caption{\label{Hilb4:fig} Total reality fails again (for
radioactive configuration of lines)!}
\end{figure}

One can also make the following picture Fig.\,\ref{Hilb5:fig},
where the  3 additional basepoints are marked by circles, one of
them lying, alas, quite outside the range of the picture. A
possible, yet delicate, desideratum would be to distort the
configuration (pair of 3 lines arrangements) so that 2 of those
circled basepoints lands on the quintic $C_5$. Then we would get a
good candidate for an easy to depict total pencil of cubics on our
quintic. Evidently this desideratum is probably impossible to
arrange (a so-called ``Irrweg'').

\vskip-0pt
\begin{figure}[h]
\centering
    \epsfig{figure=Hilb5.eps,width=82mm}
    \vskip-10pt\penalty0

\caption{\label{Hilb5:fig} Another attempt!}
\end{figure}

Maybe another arrangement worth looking at is the following
Fig.\,\ref{Hilb6:fig}. Now among the 3 extra basepoints at least
one (that one corresponding to the intersection of both horizontal
lines) is located on the quintic $C_5$ (at infinity). Hence 13
points are ensured to be real for all members of the pencil. It is
easily checked that both fundamental curves of the pencil (cyan
and red cubics) cut the $C_5$ in a totally real fashion (15 real
points). For symmetry reasons (along the axe at 120 degrees) the
nearly circular ellipse through the 6 points at finite distance
plus the line at angle 120 degrees belongs to the pencil, but alas
its intersection with the $C_5$ it hard to understand. Note by the
way that the hexagonal configuration of 6 points is slightly
perturbed thus there is no perfectly well defined such ellipse. At
this stage the whole exercise is akin to a dolorous acupuncture
session.  Note that our symmetry deduced member of the pencil has
the wrong behavior through the basepoints at infinity, hence the
right curve belonging to the pencil includes rather the line at
infinity (or at least a slight perturbation thereof). Thus we
count 12 intersections with the oval coming from the nearly
circular circuit, and just one intersection at infinity. This
underscored total of 13 seems to indicate that this pencil again
fails total reality.

\begin{figure}[h]
\centering
    \epsfig{figure=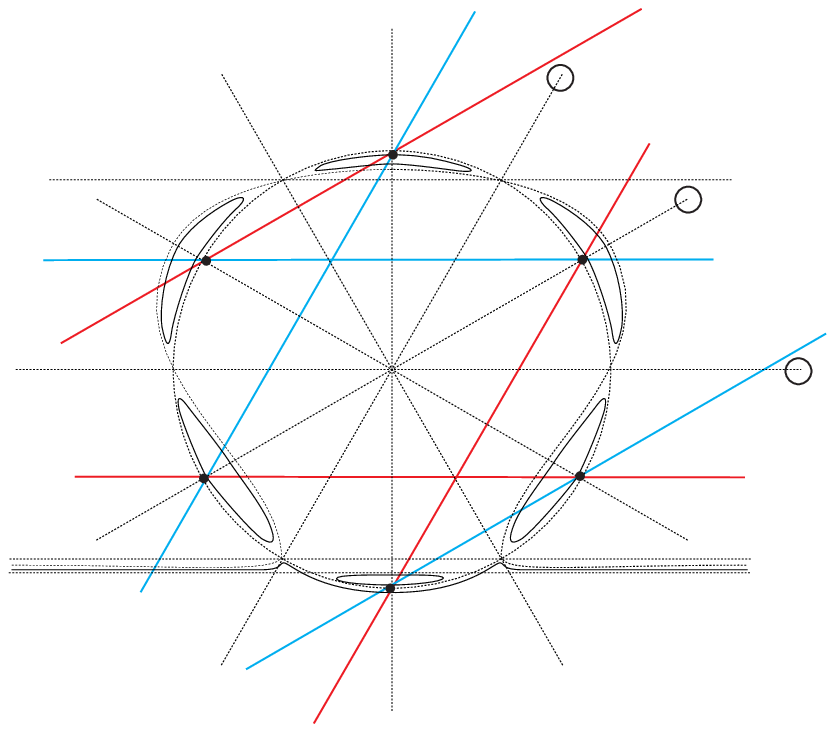,width=82mm}
    \vskip-10pt\penalty0

\caption{\label{Hilb6:fig} Another figure raising some hope but
soon failing again!}
\end{figure}

Albeit our exposition is not from the best stock, we hope at least
to have demonstrated that the synthetic construction of total
pencils on $M$-curves is not an easy matter. Of course it is not
improbable that I missed something fairly easy!

{\bf Isoperimetric digression.} During the session I wondered if
the following problem makes sense. One of the notorious difficulty
when trying to do real  pictures of algebraic curves  is that some
ovals tend to be microscopic (especially for Harnack-maximal
curves). Is there some optimal curve best suited for depiction?
Admittedly the problem makes sense only for Euclidean affine
models as opposed to projective curves (which could be pictured on
the sphere up to a double cover). One could for instance ask  the
curves to enclose maximum area for a given length of the circuits.
(Of course this makes sense only for curves of even degrees,
except if we neglect the pseudoline.) This would be a sort of
isoperimetric problem for curves competing among algebraic ones
(of some   fixed degree). Of course for degree two the
isoperimetric   solution is the circle. What about degree 4? A
candidate is perhaps the Fermat curve $x^4+y^4=1$ whose real
picture is somewhere between a circle $x^2+y^2=1$ and a square
$x^{\infty}+y^{\infty}=1$. Of course one could argue that the
optimal quartic is just a circle counted by multiplicity 2, but
then the length of the circuit has to be counted twice. We have no
certitude that our problem is well posed, nor that it is truly
interesting. The naive scenario would be that the optimum is
always the Fermat curves of higher even orders, yet what about
$M$-curves? Maybe we need to restrict the problem to them, and ask
for the best Euclidean realization of an $M$-curve? So for
instance what is the best $M$-quartic? The best $M$-quintic? Does
it looks like Ahlfors' portrait (on Fig.\,\ref{Hilb2:fig})?

Let us a last time return to our main problem of tracing a totally
real pencil for an $M$-quintic. Remember once more that
theoretical existence is ensured by the baby case
(Bieberbach-Grunsky)  of Ahlfors theorem on circle maps. Our dream
would be that for such a quintic there is  a simple-to-draw pencil
generated by  2 configurations of 3 lines. Psychologically it is
helpful to reverse the viewpoint. Instead of starting from the
quintic $M$-curve $C_5$ and trying hard to depict the pencil, we
shall start from the pencil and try to construct a  curve tailored
to it.

\begin{figure}[h]
\centering
    \epsfig{figure=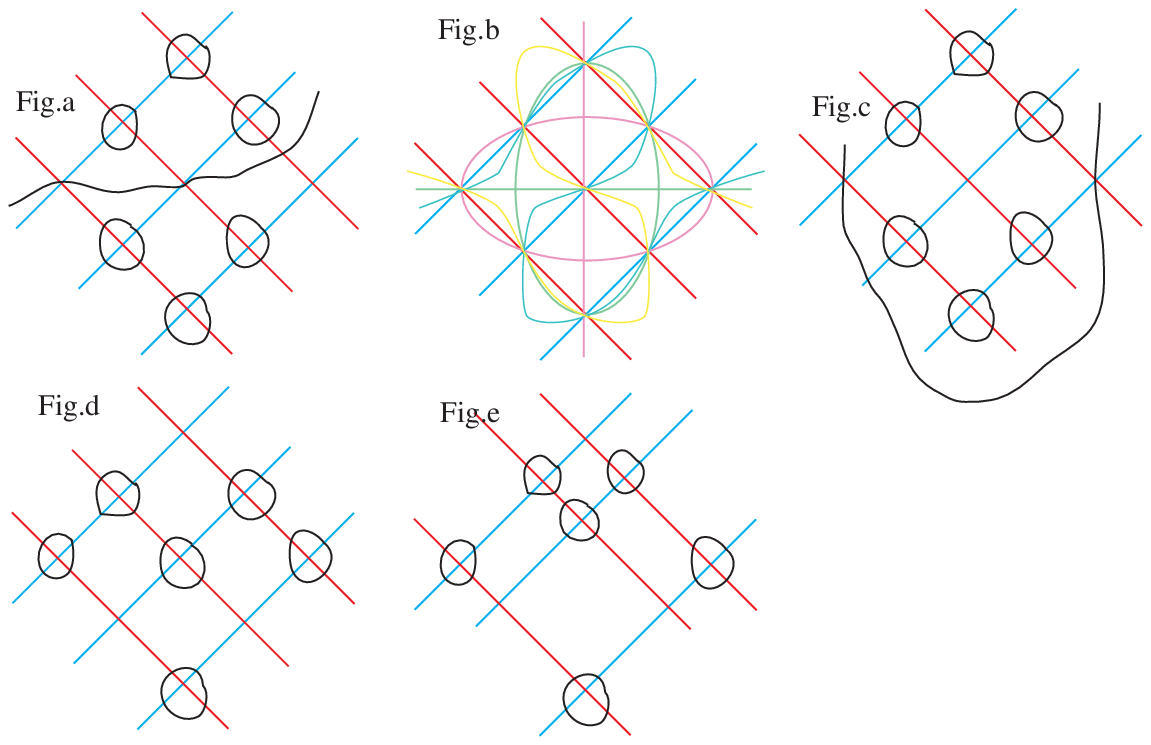,width=122mm}
    \vskip-10pt\penalty0

\caption{\label{Hilb6b:fig} Trying to find an obstruction}
\end{figure}

So we consider the pencil generated by 2 systems of parallel lines
(colored cyan and red) with 9 basepoints (multicolored
intersections) and  try to build around this perfectly explicit
pencil (cf. the previous Fig.\,\ref{Fcubic:fig}
Fig.\,\ref{Hilb6b:fig}b below) a quintic having the following
schematic picture (Fig.\,\ref{Hilb6b:fig}a). This is to mean that
each of the 6 ovals encloses one of the 9 basepoints, with the
B\'ezout restriction that no aligned triad are enclosed (else 6
intersections in $C_5$ with a line) and further the pseudoline
passes through 2 other basepoints. If such a ``real scheme''
(Rohlin's jargon) exists then each curve of the pencil will cut on
the $C_5$ a total of 15 real points. Indeed the 6 ovals contribute
each for twice (now Fig.\,\ref{Hilb6b:fig}b ensures connectedness
of all cubics forming the pencil!) and the pseudoline for 2, hence
a total of 14 and the last one is forced to be real as well (for
algebraic ``Galois theoretic'' reasons).

So exhibiting this scheme would complete our goal. Note the
absence of B\'ezout-type obstruction to the posited real scheme
(Fig.\,\ref{Hilb6b:fig}a). Yet maybe there is deeper topological
obstructions involving say the foliation underlying the pencil. In
fact the argument is more modest. The two basepoints connected by
the pseudoline are separated by the green ellipse. So the arc
joining them (choose one!) is forced to have an extra intersection
with the green ellipse (on Fig.\,\ref{Hilb6b:fig}b). Topology
forces the creation of a second intersection (intuitively the
pseudoline once trapped in the green ellipse has to escape it).
Thus we arrive at a total of $12$ (6 ovals), plus the $2$ assigned
basepoints on the pseudo-line and plus the 2 extra-points just
created. This gives 16 intersections between $C_5$ and the green
cubic (enough to overwhelm B\'ezout). This prohibits the desired
scheme.

Another (a priori) tangible real scheme is the one depicted on
Fig.\,\ref{Hilb6b:fig}c. Then it seems that arguing with the lilac
conic we may repeat something like the previous argument. More
precisely, if the pseudoline never penetrates inside the lilac
ellipse $L_2$ then it has to be tangent to it at the 2 assigned
basepoints but this gives already 4 extra-points which added to
the 12 God-given produce an excess $16>15$! Thus we may assume the
pseudoline $P$ to penetrate in the lilac ellipse (total of 13
intersection). Then several cases may occur. If $P$ tries to evade
from the lilac ellipse $L_2$ then we have $14$ intersections, yet
it must still pass to the second basepoint and (being now outside
the $L_2$) this creates at least 2 intersections (counted by
multiplicity). So eventually the pseudoline $P$ is forced to reach
the other basepoint while staying inside the lilac $L_2$, and
hence to cut the lilac axis of this ellipse. The latter axis being
contained in the inside of the green ellipse, we get again  4
extra  intersections with the green cubic (beside the 12 arising
de facto from the ovals); too much for B\'ezout.

All this (if correct?, and suitably simplified!) should prove the
following:

\begin{prop} It is impossible to sweep out in a totally real
fashion an $M$-quintic via a basic pencil of cubics spanned by two
arrangements of parallel lines.
\end{prop}

If true and suitably generalized to other configurations (see
$\bigstar$ right below) this explains perhaps why we had so much
trouble to make an appropriate depiction of the desired pencil.
Again totally real pencils exist in abstracto hence in concreto,
yet are probably of a somewhat more elaborated vintage.

[02.11.12] $\bigstar$ For instance it should be noticed that there
is another possible scheme (distribution of 6 ovals) satisfying
the ``no-three-in-line'' condition prompted by B\'ezout. This is
depicted on Fig.\,\ref{Hilb6b:fig}d which is admissible provided
the horizontal diagonal is not aligned. Hence the real picture
looks rather like Fig.\,\ref{Hilb6b:fig}e. Of course it would be
too cavalier to claim that the previous obstruction to the case at
hand as the ellipses were destroyed during the process.

We leave the problem in this very unsatisfactory state of affairs,
but let us perhaps try to motivate why the explicit depiction
project could be fruitful!

From  the viewpoint of gravitational systems (cf. the previous
Sec.\,\ref{sec:electrodynamics}) the interest of $M$-curves is
that they express in some sense the most complex orbital structure
permissible for a given genus (at least the maximum number of real
circuits). Hence if  Metatheorem~\ref{metatheorem:thm} is reliable
such $M$-curves should display some remarkable motions. The
intricacy of the trajectories is already suggested by Hilbert's
$M$-sextic on Fig.\,\ref{Hilb1:fig}. However until the total
pencil (of Bieberbach-Grunsky-Ahlfors) is not made explicit the
dynamics of the electrons is imbued by mystery and darkness.
Remind from Bieberbach-Grunsky
(=Lemma~\ref{Enriques-Chisini:lemma}) that Hilbert's $M$-sextic is
not only static object but one  animated by a circulation (total
pencil) having one electron on each oval. We can from the static
picture  vaguely try to guess where repulsions take places and
arrive at something like Fig.\,\ref{Hilb7:fig}.
\begin{figure}[h]
\hskip63pt \penalty0
    \epsfig{figure=Hilb7.eps,width=92mm}
    \vskip-10pt\penalty0

\caption{\label{Hilb7:fig} Trying to guess the dynamics from the
static locus}
\end{figure}

On Fig.\,\ref{Hilb7:fig}, italics numbers enumerate ovals while
roman numbers indicates positions at various times $1,2,3,4$. Note
that our Harnack-maximal curve being dividing, it has a complex
orientation (as the border of one half). This orientation agrees
with that inherited from the smoothing. Further it has to be
respected by the circulation due to the holomorphic character of
the (Bieberbach-Grunsky) circle map. Having this is mind it is
straightforward to make the picture above (Fig.\,\ref{Hilb7:fig})
using the rule that whenever a repulsion is observed then
electrons must be in close vicinity and thus any pair of points
minimizing the distance between two neighboring ovals must be
synchronized, hence labelled by the same time unit. In contrast
when two close ovals do not repulse them (like ovals {\it 1} and
{\it 10}) then they must be anti-synchronized in the sense that
both particles do not visit the contiguity zone at the same
moment.
 For instance there is also a repulsion between
electrons on ovals {\it 1} and {\it 11} at time 1. So far so good.
However on completing the picture one sees between ovals {\it 6}
and {\it 11} some anomalous (asynchronic) repulsion. Maybe one can
explain this via distant repulsion involving other particles of
the system (especially the electron on oval {\it 10}).

All this is very informal and  saliently illustrates the sort of
obscurantism
caused by a lack of explicit knowledge of the total pencil. This
perhaps motivates once more to complete the programme of the
present section (construction of total pencils in Harnack-maximal
cases). Ultimately one could dream of a computer program showing
in real time the  circulation of electrons prompted by the
Bieberbach-Grunsky Kreisabbildung(en) along an Hilbert $M$-sextic.

Let us finally observe that there are other $M$-sextics
(Harnack's, Hilbert's and even Gudkov's). Basically the one we
depicted (Hilbert's) is
gained by smoothing the configuration $E_2\cdot C_4=0$ consisting
of an ellipse $E_2$ (circle on the picture) and an $M$-quartic
$C_4$ one of whose oval oscillates across the ellipse $E_2$. It
may  be noticed that the oscillating oval lies mostly inside the
ellipse (cf. the left-top part of Fig.\,\ref{Hilb8:fig}). [This
schematic---yet B\'ezout compatible---style of depiction is
borrowed from Gudkov 1974 \cite[p.\,20]{Gudkov_1974/74}.] One can
reverse this situation, by putting the vibrating oval outside the
ground ellipse to get another $M$-sextic (cf. the right-top part
of Fig.\,\ref{Hilb8:fig}). A concrete construction this is
achieved on the bottom part of Fig.\,\ref{Hilb8:fig}).

\begin{figure}[h]
\centering
    \epsfig{figure=Hilb8.eps,width=112mm}
    \vskip-10pt\penalty0

\caption{\label{Hilb8:fig} Two constructions of $M$-sextics
(Harnack and Hilbert) in Walt-Disney mode of depiction borrowed by
Gudkov 1974 \cite[p.\,20]{Gudkov_1974/74}.}
\end{figure}

This curve has one ``big'' oval enclosing nine ``small'' ovals and
the other lies outside. Of course if our metatheorem
(\ref{metatheorem:thm}) is plausible then it is challenging to
interpret the dynamics especially the orbit along this long oval
enclosing all others but one. Of course this would essentially
boils down to visualize a total pencil for this $C_6$.

Finally let us make a little remark. We see that there must a deep
reaching connection between Ahlfors theory of circle maps and the
extrinsic geometry of real dividing curves, the link being given
by the notion of total pencil. Another basic application of total
pencils could arise in curve plotting problems. Assume given an
algebraic equation $f(x,y)=0$  and a machine supposed to make a
plot of the real locus. Suppose e.g. that the polynomial has
degree 5, defines a smooth curve and that we have already traced
within reasonable accuracy 2 ovals and a pseudoline and finally
that both ovals are nested. Then the theory of total maps (but in
fact B\'ezout suffices) ensures that the real locus has already
been exhausted and we may stop the ``root finding'' algorithm. Of
course the story becomes even more grandiose on appealing to
Newton-Cayley iteration method and the allied fractals appearing
as attracting basins. Likewise if an octic has 4 nests of depth 2
its real locus has already been exhausted (compare
Fig.\,\ref{Foctic4:fig}). Indeed in that case the pencil of conics
through the 4 deeply nested ovals imposed to pass through another
hypothetical point would create an excess of $8\cdot 2+1=17>16$
intersection points.

\subsection{A baby pseudo-counterexample in degree 4}\label{sec:degree-four}

\begin{figure}[h]
\centering
    \epsfig{figure=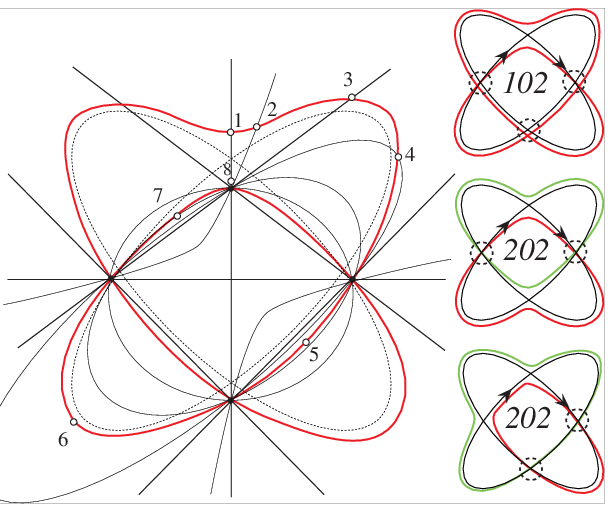,width=122mm}

\caption{\label{F102:fig} Tracing a totally real pencil of conics
on an orthosymmetric quartic}
\end{figure}

We now give an example in degree 4. The recipe is always is the
same and we get the example 102 below (Fig.\,\ref{F102:fig}). It
has $g=0$, $r=1$, thus $p=0$. At first glance the visual gonality
as measured via a pencil of lines is $\gamma^{\ast}=2$ (projection
from one of the nodes). This seems of course to violate Gabard's
bound $\gamma\le r+p$. However using a pencil of conics passing
through the 3 nodes plus the point (labelled $8$ on the figure)
gives a total pencil of the right degree. Of course the example is
a paroxysm of triviality, yet it is still a nice case to visualize
the fairly complex dynamics of total pencils. The forward
semi-orbit of the series is depicted by points $1,2,\dots, 8$
after which the motion reproduces symmetrically.

Another example arises when we keep less singularities unsmoothed.
We obtain so a linkage of ``heartsuits'' (cf. the middle picture
202 on Fig.\,\ref{F102:fig}). Now $r=2$, $g=3-2=1$, and so $p=0$.
A linear projection from one of the 2 nodes suffices to exhibit
total reality, and so the gonality is $\gamma=2$. One can still
trace pencils of conics through the nodes plus 2 extra points on
the curve to get series of degree $2\cdot 4- 2\cdot 2-2\cdot
1=8-4-2=2$. Those gives more maps realizing the gonality. Of
course one can also materialize such a curve as a smooth plane
cubic, in which case we also see $\infty^1$ total pencils induced
by linear projection from the unique oval. (Projecting from the
pseudo-line, the oval of the cubic has some ``apparent contour''
and total reality fails.) One can also get the bottom picture 202
on Fig.\,\ref{F102:fig}, which has the same invariants.

\subsection{Low-degree circle maps in all topological types
by Harnack-maximal reduction}\label{sec:Chambery}

[Source=Gabard 2005, Chamb\'ery talk (unpublished as yet)] Once
Ahlfors theorem is known in the simple Harnack-maximal case  (cf.
Lemma~\ref{Enriques-Chisini:lemma}) one can easily exhibit in any
topological type some very special surfaces (in Euclid's 3-space)
admitting a circle map to the disc having very low degree. Of
course this is far remote from reassessing the full Ahlfors
theorem, yet it is an interesting construction, which perhaps
could lead to a general proof when combined with some
Teichm\"uller theory. But this is only a vague project we shall
not be able to pursue further.

\begin{figure}[h]
\centering
    \epsfig{figure=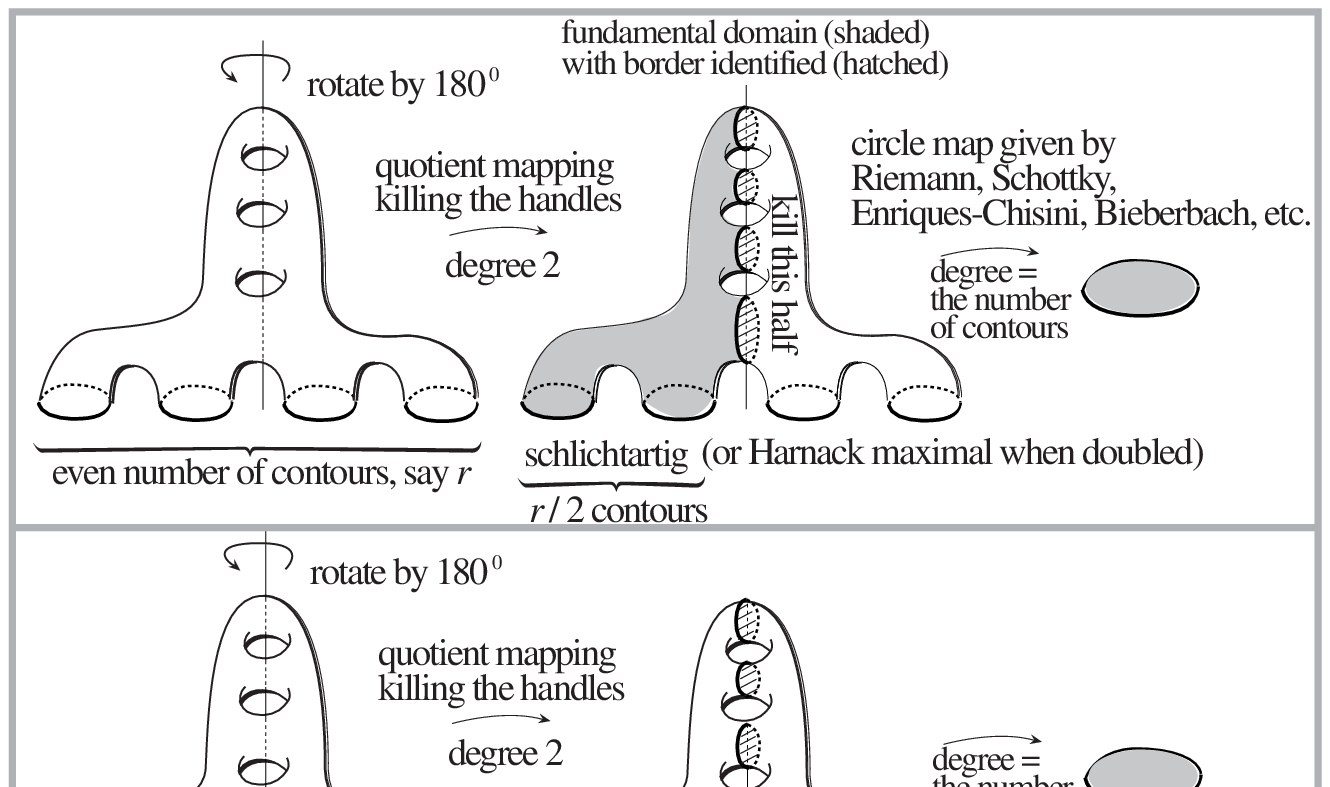,width=106mm}

\caption{\label{Chambery:fig} Low-degree circle maps for pretzels
with 2-fold rotational symmetry}
\end{figure}

Let us start with a membrane in Euclidean 3-space (endowed with
the conformal structure induced by the Euclidean metric). Suppose
the surface invariant under a symmetry of order two (cf.
Fig.\,\ref{Chambery:fig}). The
key
feature of this figure is that the axis of rotation ``perforates''
each ``hole'' of the pretzel.
Hence, when taking the quotient all handles are
killed, and we get a proper(=total) morphism to a schlichtartig
configuration (i.e. of genus $p=0$). This in turn admits a circle
map of degree equal to the number of contours (by the
Bieberbach-Grunsky theorem=Lemma~\ref{Enriques-Chisini:lemma}).

The composed mapping gives a circle map of degree $2\cdot
\frac{r}{2}=r$ when $r$ is even, and of degree $2\cdot
\frac{r+1}{2}=r+1$ if $r$ is odd. (Compare again
Fig.\,\ref{Chambery:fig}.)

Of course this has little weight in comparison to the general
theorem of Ahlfors (1950 \cite{Ahlfors_1950}), yet it is  a simple
example showing that the degree of circle maps can be fairly lower
than the degrees $r+2p$ or even $r+p$.

\subsection{Experimental evidence for Coppens' gonality}


[24.03.12]/[19.10.12] In this section we discuss Coppens result
(2011 \cite{Coppens_2011}) on the realizability of all gonalities
compatible with the $r+p$ bound (Gabard 2006 \cite{Gabard_2006})
on the degree of an Ahlfors circle map. Our superficial approach
will not recover Coppens full result, yet is worth presenting for
it enhances the depth of Coppens' result. Looking at explicit
projective models of Riemann surfaces always makes Riemann-type
existence theorems (like Ahlfors maps) look quite formidable
jewels (not to say miracles) when looked at experimentally through
the Plato cavern of extrinsic algebraic geometry. The game is also
pleasant because sometimes one gets the impression that Gabard's
bound $r+p$ looks blatantly violated. Also interesting is the
issue that such basic experimental studies (akin to the CERN
particles collider at a modest scale) are quite useful for
understanding the failure of connectivity of the space of minimal
circle maps (those of lowest possible degree). Further experiments
should contribute to add some valuable insights over Ahlfors'
theory. (A. Einstein puts it as follows: ``Any knowledge of the
world starts and ends with experiments.'')

Coppens' result is the following. To stay closer to Ahlfors'
viewpoint, we paraphrase it in the language of {\it compact
bordered Riemann surfaces\/} (abridged membranes) instead of that
of real dividing curves. Albeit most of our examples are derived
via
algebraic geometry, we will never have to write down any (boring)
equation due to the graphical flexibility of plane curves \`a la
Brusotti/Klein/Pl\"ucker
(reverse historical order).
So we are drifted to a sort of synthetic geometry.

\begin{theorem} {\rm (Coppens 2011 \cite{Coppens_2011})}
Given any two integers $r\ge 1$ and $p\ge 0$, and any integer
$\gamma$ satisfying $\max\{2, r \}\le \gamma \le r+p$, there is
membrane $F_{r,p}$ with $r$ contours of genus $p$ whose gonality
is the assigned value $\gamma$.
\end{theorem}

Recall that the gonality of the membrane is understood as the
least degree of a circle map from the given membrane (to the
disc).

$\bullet$ For $(r,p)=(1,0)$, the statement becomes vacuous, but of
course we can alter the range of permissible values as $r\le
\gamma\le r+p$.

$\bullet$ When $p=0$, $\gamma$ can take only the value $r$ and the
latter is realized via the Bieberbach-Grunsky theorem
(Lemma~\ref{Enriques-Chisini:lemma}).

$\bullet$ For $(r,p)=(1,1)$, the double has genus $g=(r-1)+2p=2$
hence is hyperelliptic. This actually proves the existence of a
circle map of degree 2 ($=r+p$) in accordance with Gabard's bound
$r+p$. Coppens's realizability theorem is trivially verified in
this case for $\gamma$ can only assume value 2.

$\bullet$ For $(r,p)=(2,1)$, the range of $\gamma$ is $2\le \gamma
\le 3$. The value $\gamma=2$ is realized by a hyperelliptic model.
The value $\gamma=3$ is obtained by considering a smooth quartic
$C_4$ with two nested ovals while projecting it from a point on
the innermost oval. This gives a {\it totally real\/} morphism of
degree $3$. Total reality means that fibers above real points
consists entirely of real points. We use also the abridged jargon
{\it total map\/} which is quite in line with terminology used by
Sto\"{\i}low 1938 \cite{Stoilow_1938-Lecons} or Ahlfors-Sario 1960
\cite{Ahlfors-Sario_1960}, who use ``complete coverings''.

\begin{figure}[h]
\hskip-75pt \penalty0
    \epsfig{figure=Coppens.eps,width=172mm}

\caption{\label{Coppens:fig} Tabulation of bordered surfaces with
assigned gonality: for each value $(r,p)$ the array of permissible
gonalities is depicted as a dashed line imagined as lying over the
grid. Italicized integers indicate the corresponding gonality.}
\end{figure}

$\bullet$ For $(r,p)=(3,1)$, the genus of the double is
$g=(r-1)+2p=2+2=4$. This is not the genus $g=\frac{(m-1)(m-2)}{2}$
of a smooth plane curve of order $m$
 which belongs to the list $0,1,3,6,10,
\dots$ of triangular numbers, yet suggests looking at a quintic
$C_5$ with two nodes. We thus consider a configuration of two
conics plus a line and smooth it out in a orientation preserving
way (so as to ensure the dividing character of the curve by a
result of Fiedler 1981 \cite{Fiedler_1981}). We obtain so the
curve depicted on Fig.\,\ref{Coppens:fig} bearing the nickname
313. This actually encodes the value of the invariant $(r,p,
\gamma)$ written as the string $rp\gamma$, yet a priori the
gonality $\gamma$ is not known and its value must be justified. On
that figure 313 the dashed circles indicate those crossings that
were {\it not\/} smoothed. The half of this curve is a bordered
surface of type $(r,p)=(3,1)$, since $p=\frac{g-(r-1)}{2}$. It
remains to evaluate its gonality. The idea is always to look at
the curve from the innermost oval. In the case at hand, we project
the curve from one of the two nodes to get a total morphism of
degree $5-2=3$. Since $r=3$ is a lower bound on the gonality
$\gamma$, it follows that $\gamma=3$, exactly. Note that this
example seems to {\it answer in the negative our question about
the connectivity of the space parameterizing minimal circle maps}.
Further one can drag one point to the other while travelling only
through {\it total\/} maps of degrees 4 (namely projections from
points located in the intersection of the interiors of the blue
resp. red ovals). [09.11.12]---{\it Warning.} Remember that a
similar picture (Fig.\,\ref{F102:fig}, right-middle part) gave an
example where the curve looked 2-gonal in only 2 ways, but another
model of the curve (as a plane cubic) prompted the same gonality
in $\infty^1$ fashions. So some deeper argument is required either
to assess (or disprove) the italicized assertion.

Next, still for the same topological invariants $(r,p)=(3,1)$, we
would like to find a membrane of gonality $\gamma=4$. This may be
obtained from the same initial arrangement while moving the
location of the dashed circles (of inert crossings) to get picture
labelled 314 on Fig.\,\ref{Coppens:fig}. The corresponding quintic
projected from a point situated on the inner (blue-colored) oval
has $\gamma\le 4$. Over the complexes, this quintic has gonality
$\gamma_{\Bbb C}=3$ (projection from one of the nodes) and this is
the only way for the curve to be trigonal. Yet over the real
picture (our 314) none of these (trigonal) projections is total
(since the inner oval has an apparent contour, i.e. some tangent
to it passes through the node). It follows that $\gamma=4$,
exactly.

$\bullet$ Let us next
examine $(r,p)=(4,1)$. Then $g=(r-1)+2p=5$, so we look at quintics
with one node. To create as many ovals, it proves convenient to
reverse the orientation of one of the conics. We obtain so the
figure coded 415. After noting that $r=4$, we project the curve
via a pencil of conics assigned to pass through 4 points chosen in
the innermost ovals (asterisks on the figure). Letting those 4
points degenerate against the ovals while exploiting the
possibility of pushing one of them toward the node (so as to lower
by 2 units the degree) we find $\gamma\le 2\cdot 5- 3\cdot 1
-1\cdot 2=10-3-2=5$. Over the complexes, the curve at hand
(uninodal quintic) is trigonal only when seen from its unique node
and 4-gonal only when projected from a smooth point. Inspection of
the figure shows that none of these maps is total. It follows that
$\gamma=5$ exactly.

It remains to find an example with $\gamma=4$. For this we just
drag below the dashed circle (cf. label 414 on
Fig.\,\ref{Coppens:fig}), do the prescribed smoothing (always in
the orientation consistent way). The resulting curve has $r=4$ (as
it should). The novel feature is that the node is now  accessible
from 2 basepoints of the pencil of conics assigned in the deep
ovals. This permits a lowering of the degree to $\gamma\le 2\cdot
5- 2\cdot 1 -2\cdot 2=10-2-4=4$. Remarking that the unique
morphism of lower degree 3 (linear projection from the node) is
not total we deduce that $\gamma=4$ exactly. The other morphisms
of degrees 4 (namely projections from real points on the curve)
obviously fails to be total, thus we infer that the curve (or the
allied membrane) is uniquely minimal (i.e. there is a unique
circle map of minimum degree).

Before embarking on larger values of the invariants $(r,p)$, we
make a general remark, related to the previous
Sec.\,\ref{sec:Chambery}. There  a suitable membrane in 3-space
invariant under rotation by $\pi=180^{0}$ with a totally vertical
array of handles (cf. Fig.\,\ref{Chambery:fig}) showed the
following:

\begin{lemma} {\rm (Barbecue/Bratwurst principle)}\label{Barbecue:lem}
$\bullet$ If $r$ is even, there is for any value of $p$ a membrane
of type $(r,p)$ admitting a circle map of degree $r$ (the minimum
possible value), whose gonality is therefore $\gamma=r$ exactly.

$\bullet$ If $r$ is odd ($p$ arbitrary), there is a membrane of
type $(r,p)$ admitting a circle map of degree $r+1$, whose
gonality $\gamma$ is therefore $r\le \gamma\le r+1$. (Alas, the
exact value remains a bit undetermined!)
\end{lemma}

This lemma fills quickly several positions of our
Fig.\,\ref{Coppens:fig}, namely those marked by a square. In the
special case $r=1$ (belonging to the indefinite odd case), we can
get rid off the annoying indetermination, because as soon as $p\ge
1$ the minimal value $r$ of the range $r\le \gamma\le r+1$ cannot
be attained. Corresponding invariants are reported by rhombuses
(squares rotated by $\pi/4$) on Fig.\,\ref{Coppens:fig}.

$\bullet$ Next we study $(r,p)=(5,1)$. Then $g=(r-1)+2p=6$,
prompting to look at smooth quintics (without nodes). Consider the
curve denoted 516 on Fig.\,\ref{Coppens:fig}, which has $r=5$.
When projected via a pencil of conics through the assigned 4
basepoints (depicted by asterisks on the figure) and letting them
degenerate toward the ovals gives a total map of degree $2\cdot 5-
4\cdot 1 =10-4=6$. Hence $\gamma \le 6$. Morphisms of lower
degrees exist in degree 4 (linear projection from a point situated
on the curve), and degree 5 (projection from points outside the
curve). Clearly none of these maps is total, so that $\gamma=6$
exactly. Of course the minimal degree maps considered are plenty
(no uniqueness), yet their parameter space is connected.

Next we require a specimen with $\gamma=5$. It seems evident that
we have exhausted the patience of quintics (at least for the given
arrangement), hence let us move to sextics of genus 10 (when
non-singular). To get the right genus $g=6$, we have to conserve 4
nodes. Starting from a configuration of 3 conics suitably oriented
and smoothed we obtain the figure denoted 515 with $r=5$ (still on
Fig.\,\ref{Coppens:fig}). Using a pencil of conics with 4 assigned
basepoints (asterisks on the figure) gives a (probably total) map
of degree $2\cdot 6-1\cdot 2-3\cdot 1=12-2-3=7$. This agrees with
Ahlfors bound $r+2p$, but seems to challenge Gabard's bound
$r+p=6$. Maybe a pencil of cubics is required instead.
Such a cubics pencil has 9 basepoints but only 8 of them may be
assigned. Hence creating some 4 new basepoints (denoted by bold
letters {\bf 1,2,3,4} on the figure) and letting them degenerate
to the ovals or better the nodes (when some are accessible) gives
a map of degree $3\cdot 6-3\cdot 1- 5\cdot 2=18-3-10=5$, rescuing
Gabard's $r+p=6$ and also giving the desired gonality $\gamma=5$.
Admittedly this example is quite complex and perhaps not the best
suited to illustrate Coppens' gonality result. Its interest is
still that it seems to corrupt Gabard's bound $r+p$, and the
latter can  only be rescued by appealing to fairly sophisticated
pencils. Of course it could be the case a priori that our curve
(515) admits a pencil of conics of lower degree than 7, but under
the totality condition basepoints must be distributed in the deep
ovals by a Poincar\'e index argument (cf.
Lemma~\ref{Poincare-lower-bound}). This impedes a lowering of the
degree via a more massive degeneration of the base-locus to the
nodes of picture 515 on Fig.\,\ref{Coppens:fig}. Admittedly the
predicted total pencil of cubics ought to be described more
carefully.

{\it Summary of the situation.}---Of course one should still work
out the higher values of $r$ while keeping $p=1$. As you notice
our method is far from systematic. (All the difficulties
encountered so far already enhance the power of Coppens' result.)

$\bullet$ Then one must also handle higher values of $p$, starting
with $(r,p)=(1,2)$. The case $\gamma=2$ is easy (via the barbecue
construction, Lemma~\ref{Barbecue:lem}). For $\gamma=3$ we can
imagine a surface with 3-fold rotational symmetry (cf. picture
123X on Fig.\,\ref{Coppens:fig}). For it $\gamma\le 3$, but how to
show equality? Alternatively, one may consider an algebraic model.
Since $g=(r-1)+2p=4$,  we look among quintics with 2 nodes. A
suitable smoothing gives figure named 123, with $r=1$ (one
circuit). Linear projection for the ``inner'' node gives a total
map of degree $1\cdot 5-1\cdot 2=5-2=3$, so $\gamma\le 3$. But the
complex gonality of such a quintic is $\gamma_{\Bbb C}=3$. Since
$\gamma_{\Bbb C}\le \gamma$ it follows that $\gamma=3$ exactly.

$\bullet$ Let us next explore $(r,p)=(2,2)$. Then $g=5$. A surface
with $\gamma=2$ is easily found (barbecue rotational symmetry). To
realize the other gonalities we look among quintics with one node.
We first meet figure 223, which has a total morphism of degree 3
(projection from the node). Hence $\gamma\le 3$ which is in fact
an equality, since 3 is also the complex gonality of an uninodal
quintic. To get a curve with $\gamma=4$ we just drag the
unsmoothed singularity to get figure 224. Projection from the node
is not total anymore, but a total map arises when projecting from
the (green) oval giving rise to degree $5-1=4$. Since such a
quintic is uniquely trigonal (via projection from the unique node,
which which failed to be total), we infer that $\gamma=4$,
exactly. Coppens's theorem is verified for this topological type.

{\bf Premature conclusion/State of the art.} It is clear that one
can continue the game to tackle higher and higher values of the
invariants. Instead of looking solely in ${\Bbb P}^2$ it is also
pleasant to trace curves in ${\Bbb P}^1\times{\Bbb P}^1$, albeit
${\Bbb P}^2$ is a universal receptacle (any Riemann surfaces
nodally immerses in the projective plane). However it is clear
that our naive approach is quite time consuming and as yet we did
not deciphered a combinatorial pattern permitting to boost the
speed of the procedure to the level of an inductive process.
(Curves or Riemann surfaces of higher topological structures are
like {\it homo sapiens\/}, the result of a long, intricate
morphogenesis.) Coppens  proved the full result in one stroke by
somehow penetrating the genetic code governing the evolution of
all species.

\subsection{Minimal sheet number of a genus $g$ curve as a
cover of the line}\label{Minimal-sheet:sec}

It is classical (since Riemann 1857 \cite[\S 5,
p.\,122--123]{Riemann_1857}) that a general curve of genus $g$ is
expressible as a branched cover of the sphere ${\Bbb P}^1$ of
degree the least integer $\ge \frac{g}{2}+1$ (equivalently of
degree $[\frac{g+3}{2}]$). [Indeed if $g$ is even $g=2k$ the first
value is $k+1$ and $[\frac{g+3}{2}]=[(2k+3)/2]=(2k+2)/2=k+1$; if
$g=2k+1$ is odd then the first value is $ g/2+1$ which rounded
from above gives $(2k+2)/2+1=k+2$, and
$[\frac{g+3}{2}]=[(2k+4)/2]=k+2$.]

Riemann's truly remarkable argument (involving Abelian integrals)
is beautifully cryptical (I should still study it properly). It is
not clear (to me) if it includes the stronger assertion that {\it
any} curve of genus $g$ admits a sphere-map of degree $\le
[\frac{g+3}{2}]$. At any rate, all modern specialists agree that
the first acceptable proof of this pi\`ece de r\'esistance is
Meis' account (1960 \cite{Meis_1960}). (Meanwhile the
algebro-geometric community devised several alternative
approaches.)

Another allied (but different?) argument is the one to be
found in Klein's lectures 1892
\cite[p.\,98--99]{Klein_1891--92_Vorlesung-Goettingen},
cf. also Griffiths-Harris 1978
\cite[p.\,261]{Griffiths-Harris_1978/94}.

The latter's argument works as follows. Assume there is a
$d$-sheeted map $C_g \to {\Bbb P}^1\approx S^2$ of a genus $g$
surface to the Riemann sphere. Then Euler characteristics are
related by $\chi(C_g)= d \chi(S^2)-b$, where $b$ is the number of
branch points. This gives $b$ ramified positions, whose locations
determine the overlying Riemann surface up to finitely many
ambiguities. So the $d$-sheeted surface depends upon $b-3$
essential parameters (after substraction of the 3 arising from the
linear transformations on ${\Bbb P}^1$). This quantity has to be
$\ge 3g-3$ the number of moduli of genus $g$ curves. This implies
$b\ge 3g$, i.e. $2d-\chi(C_g)\ge 3g$, or $2d\ge 3g+(2-2g)=g+2$.
{\it q.e.d.}

So far as we know,
a similar computation as never been written down for the case of a
bordered Riemann surface expressed as a $d$-sheeted cover of the
disc (i.e., the context of  Ahlfors circle maps). The reason is
probably quite mysterious, yet also quite simply that the naive
parameter count seems to lead nowhere.

Let us attempt the naive computation. Suppose $F_{r,p}\to D^2$ to
be a membrane of genus $p$ with $r$ contours expressed as a
$d$-sheeted cover of the disc. Euler characteristics are related
by $\chi(F)=d\chi(D^2)-b$, where $b$ is the number of branch
points.
The group of conformal automorphisms of the disc as (real)
dimension $3$. Hence our $d$-sheeted surface depends upon $2b-3$
real constants, whereas the membrane $F$ itself depends on $3g-3$
real constants, where $g$ is the genus of the double (cf. Klein
1882 \cite{Klein_1882}). The Ansatz $2b-3\ge 3g-3$ gives $2b\ge
3g$, and since $b=d-\chi$ and $g=(r-1)+2p$, this gives $2d \ge
3g+2\chi=3[(r-1)+2p]+2(2-2p-r)=r+2p+1$, equivalently $d\ge
(r+1)/2+p$.


This
beats the value $r+p$ predicted in Gabard 2006 \cite{Gabard_2006},
but
looks blatantly overoptimistic. For instance taking $p=0$, gives
the degree $\frac{r+1}{2}$ violating  the absolute lower bound $r$
on the degree of a circle map. The provisory conclusion is that
the naive parameters count  leads nowhere in the bordered case.
Does somebody know an explanation?

[21.10.12] A crude attempt of explanation is that the above count
merely uses the Euler characteristic which is  a complete
topological invariant only for closed surfaces, but not for
bordered ones. (Since $\chi(F_{r,p})=2-2p-r$, trading one handle
against two contours leaves $\chi$ unchanged.) Of course the above
counting uses also $g$ (the genus of the double $2F$) but the
latter is also uniquely defined by $\chi(F)$, via the relation
$2-2g=\chi(2F)=2\chi(F)$. Thus it is maybe not so surprising that
Riemann(-Hurwitz)'s count predicts correctly the gonality of
closed Riemann surfaces but fails seriously to do so in the
bordered case. It could be challenging to find a moduli count
existence-proof of Ahlfors circle maps supplemented probably by an
adequate continuity method. For an (unsuccessful) attempt cf.
Sec.\,\ref{Hurwitz-type}.

\begin{figure}[h]
    \epsfig{figure=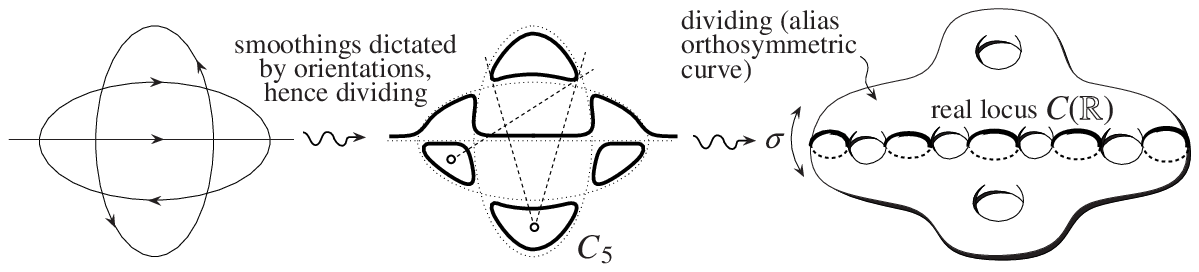,width=122mm}

\caption{\label{quintic:fig} A quintic with gonality $\gamma=6$?}
\end{figure}

A very naive (numerological) parade is to introduce a new bound
$\nu:=\max\{p+ \frac{r+1}{2}, r \}$ between the one found above
and $r$ the absolute minimum of a
total morphism.
However a simple example probably shows this to be overoptimistic
as well. Consider the plane quintic $C_5$ derived via a sense
preserving smoothing of the depicted configurations of $2$ conics
and a line (cf. Fig.\,\ref{quintic:fig}).

Its genus is $g=\frac{(m-1)(m-2)}{2}=\frac{4\cdot3}{2}=6$, and we
see $r=5$ real circuits. The relation $g=(r-1)+2p$ gives $p=1$
(genus of the half). Hence the new bound is $p+
\frac{r+1}{2}=1+3=4$, but $r=5$ so $\nu=\max=5$. However
 the membrane (corresponding to one half of the dividing
curve $C_5$) cannot be represented with 5 sheets over the disc.
Indeed a morphism of degree 5 from $C_5$ to the line ${\Bbb P}^1$
can \underline{{\it only}} arise through linear projection of the
quintic $C_5$ from a point not on the curve (else degree $4$), but
no such projection is totally real (compare central part of
Fig.\,\ref{quintic:fig}, or argue via the Poincar\'e index, cf.
Lemma~\ref{Poincare-lower-bound}).

[09.11.12] {\it WARNING about the underlined ``only''.}---This
argument looks at first sight quite convincing, yet it appears to
be insufficient, and possibly the  assertion itself on the
gonality $\gamma(C_5)=6$ is erroneous. First, a
total morphism of degree $r+p=5+1=6$ (as predicted in Gabard 2006
\cite{Gabard_2006}) should exist. This is corroborated  by taking
a pencil of conics through 4 points inside the 4 ovals of the
above depicted $C_5$ (cf. Fig.\,\ref{FGuerN:fig}, left part) and
letting them degenerate against the ovals,  giving a total map of
degree $2 \cdot 5-4=6$. This tell us only $\gamma \le 6$. A
priori, it could be the case that higher order pencils access the
low degree $5$, and with some good-fortune do it in a totally real
way. In that case the gonality lowers down to $\gamma=5$ (the
minimum permissible as $r=5$). Let us quickly discuss how this
could happen, at least over the complexes. A priori pencils of
cubics may have degrees as low as $3\cdot 5-3^2=6$ (hence not
violating the previous token); quartics as low as $4\cdot 5-
4^2=4$, but quartics have dimension $\binom{4+2}{2}-1=14$ so that
in reality only 13 basepoints may be assigned freely, hence the
right value is $4\cdot 5- 13=7$; for quintics this is as low as
$5\cdot 5- 5^2=0$ (yet all values $<4$ violates already the
complex gonality of a smooth quintic, cf. e.g. Arbarello et al.
1985 \cite[p.\,56,
Exercise~18]{Arbarello-Cornalba-Griffiths-Harris_1985-BOOK}). In
fact the dimension of quintics is $\binom{5+2}{2}-1=20$ and thus
the minimum degree is $5\cdot 5- 19=6$. For sextics the degree is
as low as $6\cdot 5- 6^2=-6$, but since the sextics dimension is
$\binom{6+2}{2}-1=27$, the real minimum degree is $30-26=4$ (and
this beats linear projections from outside the curve). Recall
incidentally that this is the value of the universal Riemann-Meis
bound $[\frac{g+3}{2}]=[9/2]=4$, which was already attained by
linear projections from the curve but nobody will exclude a priori
a second return.
Actually all 26 assigned basepoints fails to impose independent
conditions on sextics, because our quintic $C_5$ aggregated to any
line is a sextic meeting the requirement and varying among
$\infty^2$ parameters (and not just the expected $\infty^1$
pencil). Thus we seem to fail  getting a genuine pencil, but
contrast this with the just remembered Riemann-Meis gonality. The
situation is quite more tricky than initially expected. Another
torpedo against the  naive belief that a smooth $C_5$ has only
$\infty^1$ series of type $g^1_5$ is the existence theorem of
Brill-Noether-Kempf-Kleiman-Laksov theory (cf. e.g. Arbarello et
al.
\cite[p.\,206]{Arbarello-Cornalba-Griffiths-Harris_1985-BOOK}).
The latter states the following.

\begin{theorem} Let $C$ be a (complex) curve of genus $g$.
Every component of the variety $G^r_d$ parameterizing all linear
series $g^r_d$ of dimension $r$ and degree $d$ has dimension at
least equal to the so-called {\it Brill-Noether number\/} $\rho$,
symbolically:
$$
\dim_{\ast} G^r_d\ge \rho:=g-(r+1)(g-d+r).
$$
In particular when the latter number $\rho$ is $\ge 0$ the variety
$G^r_d$ is nonempty.
\end{theorem}

In the case at hand it follows that $\dim_{\ast} G^1_5 \ge
6-(1+1)(6-5+2)=6-2\cdot2=6-4=2$. Hence there are other pencils of
degree 5 than those readily visualized on the projective
realization! This shows how vicious the Plato cavern is! Of course
our appeal to the above general theorem, is a violation against
the principle of do-it-yourself-ness, since low genus cases are in
best treated by hand (cf. Arbarello et cie
\cite[p.\,209--211]{Arbarello-Cornalba-Griffiths-Harris_1985-BOOK}
for a possible treatment, alas not perfectly self-contained).

The following summarizes the swampy situation (while trying to
extend the generality):

\begin{lemma} {\rm (To be clarified with percentages of truth)}

$\bullet$ [100 \%] Any smooth real quintics $C_5$ with $r=5$
(hence $4$ ovals and one pseudoline) is unnested (otherwise the
line through the nest plus another oval gives 6 intersections,
corrupting B\'ezout).

$\bullet$ [80 \%] Furthermore taking a pencil of conics through
the $4$ nests gives a total pencil (why exactly? clear on the
Fig.\,ref{FGuerN:fig}(left part) but why in general?).

$\bullet$ [79 \%] Assuming  the previous point,  the gonality is
$\gamma\le 2\cdot5-4\cdot 1= 6$ (in accordance with Gabard's bound
$r+p$, but it is preferable to mistrust this!).

$\bullet$ [100 \%=0 \%] Alas it is not clear a priori that pencils
of orders $\ge 6$ do not induce total pencils of possibly lower
degree $=5$. (Recall that $r=5$ is an absolute lower bound for
total maps!)
\end{lemma}

[10.11.12] In the light of the Kempf-Kleiman-Laksov existence
theorem of special divisors (ESD) in the case of complex curves
one may wonder about its relativization in the Ahlfors context of
total maps. The point is of course that for $g^1_d$'s the
existence theorem (ESD) boils down to the Riemann-Meis bound
$\gamma_{\Bbb C}\le [\frac{g+3}{2}]$ for the gonality of complex
curves. (Plug $d\ge g/2+1$ in the Brill-Noether number $\rho$ and
notice its non-negativity.) Since Ahlfors 1950 $\gamma\le r+2p$ or
maybe Gabard 2006 $\gamma \le r+p$ is to be considered as the
genuine bordered (or orthosymmetric) avatar of the Riemann-Meis
theorem one can dream of an orthosymmetric(=dividing) version of
the whole special of divisor theory. It is not clear how to extend
total reality for higher series $g^r_d$ which are not pencils
$g^1_d$. Of course one can ask that all real members are totally
real but this seems too restrictive. Is there any example at all?
Perhaps not for simple dimension reason. For $g^2_d$'s this would
amount to a plane model of the curve cut by all real lines in real
points only. This looks overambitious by just perturbing a tangent
at a non-inflection point outside the sense of curvature.

At any rate the theory surely works  for pencils and the bonus is
that we have a certain variety akin to $G^1_d$ parameterizing all
total pencils of degree $d$ on a given dividing curve. How to
denote it? I never understood for what the ``$g$ or $G$'' of resp.
$g^1_d$ or $G^1_d$ is standing? (Candidates: groups of points,
Gerade, Gebilde, Grassmann, ?) Improvising notation,  we define
$T^1_d$ the variety of total linear series of degree $d$ on a
given dividing curve. We dream about repeating all the
phenomenology of the classic theory, cf. e.g. p.\,203 of ACGH 1985
\cite{Arbarello-Cornalba-Griffiths-Harris_1985-BOOK}:

``{\it A genus $g$ curve depends on $3g-3$ parameters, describing
the so-called moduli. Our goal is to describe how the projective
realizations of a curve vary with its moduli, and what it means to
say that a curve is general or special. Accordingly, we would like
to know, what linear series can we expect to find on a general
curve and what the subvarieties of the moduli space corresponding
to curves possessing a series of specified type look like. [\dots]
A natural question is, how can we tell one curve from another by
looking at these configurations [$G^r_d$], or more precisely, what
do these look like in general, and how---and where---can they
degenerate?}''

For our ``totality'' varieties $T^1_d$ of total pencils we would
gather them into a ``telescope'' $T^1:=\cup_{d=1}^{\infty} T^1_d$
naturally embedded in $C^{(\infty)}$, the infinite symmetric power
of the (dividing) curve $C$. We have the degree function
$\deg\colon C^{(\infty)}\to {\Bbb N}$, and the image of $T^1$ is
nothing but than the gonality sequence $\Lambda$ (Definition
\ref{def:gonality-sequence}), whose least member is the
(separating) gonality $\gamma$ (of Coppens). One would like to
understand how total pencils may degenerate to lower degrees w.r.t
the natural topology induced by $C^{(\infty)}$. We probably get a
sort of telescope with high strata attached to lower dimensional
ones (like in a CW-complex) and the game would be to understand
the geometry or combinatoric of this tower. Understanding all this
is arguably the most refined form of Ahlfors theorem one could
desire. One would then like to know not only the gonality spectrum
telling one the dimension of each strata $T^1_d$, but also know
how they can degenerate to lower strata. Degeneration could still
be encoded combinatorially in a  simplicial-complex
$\Lambda^{\ast}$ with vertices $\Lambda$ (gonality sequence). Two
vertices $d_1<d_2\in \Lambda$ are linked by an edge if a total
$g^1_{d_2}$ can degenerate to a $g^1_{d_1}$. More generally
$d_1<d_2<\dots<d_{k+1}\in \Lambda$ form a $k$-simplex whenever
each integer of the sequence admits a representative $g^1_d$
degenerating to  its immediate predecessor, hence to all
predecessors.

Working out this explicitly looks tedious already for simple
example. For the G\"urtelkurve (any smooth quartic $C_4$ with 2
nested ovals) the variety $T^1_3$ is a circle and $T^1_4$ is a
2-cell attached to the former in a natural way. Of course when a
total $g^1_4$ degenerates to a total $g^1_3$ it acquires a
basepoint, which as to be deleted (particle destruction). Total
$g^1_d$ will ultimately be denoted as $t^1_d$'s. In view of the
Brill-Noether theorem (ESD) the variety $G^1_4$ has dimension $\ge
\rho=3-2(3-4+1)=3$ and so we have a priori more than the
$\infty^2$ evident total pencils $t^1_4$ arising via projection
from the inner oval. For instance pencils of conics may have
degree as low as $2\cdot 4- 4 \cdot 1=8-4=4$. Can they be total? I
would have guessed not, but it seems that they can. Compare
Fig.\,\ref{FGuerN:fig} below.

\begin{figure}[h]
\centering
    \epsfig{figure=FGuerN.eps,width=122mm}

\caption{\label{FGuerN:fig} A total pencil on the G\"urtelkurve
cut out by conics}
\end{figure}

It would be desirable if  some continuity principle can ensure
total reality, e.g. if the 4  basepoints are distributed both
inside and outside the nested resp. unnested oval. Then like a
salesman traveller, the conic has to visit all 4 basepoints and
thus  creates at least $8$ real intersections. Our picture would
just be the limiting position of such a bipartite pencil, and the
variety $T^1_4$ would be $\infty^4$, a much larger dimension than
initially expected. Further if 3, among the 4 basepoints, become
collinear then it may be argued that the conics pencil specializes
to one of lines (after removing the static line). All this remains
to be better analyzed.

\subsection{Heuristic moduli count to justify Ahlfors
or Gabard (Huisman 2001)}

It is still plausible that one may gain some evidence in favor of
the Ahlfors circle map (either with Ahlfors $r+2p$ or preferably
the improved Gabard's bound $r+p$) by arguing via a moduli count.
(The reader reminds to have discussed orally this option with
Natanzon and Huisman in Rennes in  Summer 2001, resp. December
2001.) I do not know if
it is possible to supply a better count than the unrealistic one
of the previous section.

[14.10.12] In fact at a time when I only conjectured the bound
$r+p$, Huisman (December 2001 or 2002?) reacted
instantaneously with a parameter count giving some evidence to
the conjecture. Let me reproduce this faithfully from hand
written notes.

We adopt the viewpoint of dividing real (algebraic) curves. So let
$C$ be a such with $r\ge 1$ ovals and of genus $g$. I mentioned to
Johannes Huisman the intuition that there is a totally real
morphism $C\to {\Bbb P}^1$ (i.e. inverse image of real locus
contained in the real locus) whose degree is the barycenter of $r$
and $g+1$, that is $\frac{r+(g+1)}{2}$. (The heuristic reason
behind this 2001 intuition are given in Gabard 2006
\cite{Gabard_2006}, and in its most  primitive form in the
previous Section\ref{sec:Sketch-of-Gabard}.) ``Let us count
parameters!''. Thus spoke Huisman, like Zarathustra.

First the Riemann-Hurwitz relation written for the Euler
characteristic is $\chi(C)= d \cdot \chi({\Bbb P}^1)-b $, where
$d$ is the degree and $b$ the number of branch points (with
multiplicity). Now we count real moduli. The ramification divisor
of any totally real morphism actually lies in the imaginary locus
of the sphere (not on the equator), but is of course symmetric
w.r.t. the involution. Hence we may imagine  the $b/2$ branch
points prescribed only in the north hemisphere, thus depending  on
$2\cdot (b/2)=b$ real constants. The curve itself depends on $b-3$
moduli (subtract the dimension of the automorphism group of ${\Bbb
P}^1$ defined over ${\Bbb R}$), that is
\begin{align*}
b-3= d \cdot \chi({\Bbb P}^1)-\chi(C)&=\frac{r+g+1}{2} \cdot
2-\chi(C)-3\cr
&=(r+g+1)-(2-2g)-3=3g-4+r\ge 3g-3.
\end{align*}
This prompts enough free parameters to sweep out the full moduli
space. Of course this does not reprove the existence of circle
maps of the prescribed degree, yet give some evidence to the
assertion.

[15.10.12] A notable defect of this Huisman count is that it is a
posteriori, giving no hint why the degree value should be given by
our Ansatz. It is thus preferable to make the same computation in
a more organical way. As above the curve $C$ depends on $b-3$ real
moduli, and we demand $b-3\ge 3g-3$. This gives $d\cdot
\chi(S^2)-\chi(C)\ge 3g$, i.e. $2d\ge 3g+(2-2g)=g+2$, or $d\ge
g/2+1$.

Two remarks are in order. The above is exactly the same heuristic
calculation as the that (going back to Riemann) to be found in
Griffiths-Harris for the complex gonality of a curves, and which
we remembered before. (The least integer $d\ge (g+2)/2$ is
$[\frac{g+3}{2}]$, obvious for $g$ even and also obvious when $g$
is odd.) Hence in substance this modification of Huisman's count
truly just assert that Gabard's bound is compatible with the
gonality of the underlying complex curve, yet does not predict the
bound $(r+g+1)/2$.

Perhaps there is a better way to count, compare the section
devoted to Courant (Sec.\,\ref{sec:Courant}).

\subsection{Other application of the irrigation method
(Riemann 1857,  Brill-Noether 1874, Klein, etc.)}

The method used in Gabard 2006 \cite{Gabard_2006} is primarily
based upon an irrigation principle in a torus, which in turn
is logically reducible to the surjectivity criterion via the
(Brouwer) topological degree of a mapping to a manifold.

Via this method we obtained (in \loccit) the existence of an
(Ahlfors) circle map of degree $\le r+p$. As pointed out there,
the method also supplies a purely topological proof of Jacobi
inversion theorem, to the effect that the Abel-Jacobi mapping from
the symmetric powers $C^{(d)}$ of a complex curve to its Jacobian
is surjective as soon as dimension permits (that is for $d\ge g$).

Of course the complex (or closed) avatar of the Ahlfors
mapping is just the mapping of a closed genus $g$ surface as a
branched cover of the sphere.
In this situation it is classically known since Riemann 1857
\cite{Riemann_1857} and Brill-Noether 1874
\cite{Brill-Noether_1874}  (but disputed by the modern writers)
that the most economical sheet number required is
$[\frac{g+3}{2}]$.

Contributions on this
problem is vast (and according to the modern consensus first {\it
rigorously} proved in Meis 1960 \cite{Meis_1960} for linear series
of dimension one, whereas some classic references includes the
more case of arbitrary dimensional series, esp. Brill-Noether and
Severi)


$\bullet$ Riemann 1857 (Theorie der Abel'schen Functionen)
\cite[\S 4]{Riemann_1857},

$\bullet$ Brill-Noether 1874 \cite{Brill-Noether_1874} (working
with plane curves with singularities, so a pure algebraization of
Riemann's theory if one does not fell claustrophobic in the Plato
cavern.)

$\bullet$ Klein's lectures of 1891
\cite[p.\,99]{Klein_1891--92_Vorlesung-Goettingen} (based on
Abelian integrals and Riemann-Roch, essentially akin to Riemann's
original derivation)

$\bullet$ Hensel-Landsberg  1902 \cite[Lecture
31]{Hensel-Landsberg_1902} (probably quite similar to
Brill-Noether or inspired by Dedekind-Weber)

$\bullet$ Severi 1921 \cite[Anhang
G]{Severi_1921-Vorlesungen-u-alg.-Geom-BUCH}

\noindent Then the modern era begins with:

$\bullet$ Meis  1960 \cite{Meis_1960} (Teichm\"uller theoretic)
[alas, this monograph is notoriously difficult to obtain]

$\bullet$ H.\,H. Martens  1967 \cite{Martens_Henrik_1967} (no
proof, but a remarkable study of the geometry assuming
non-emptiness)

$\bullet$ Kempf 1971 \cite{Kempf_1971} the first existence
proof (simultaneous with the next contributors) of special
divisors in general case (linear series of arbitrary
dimension, extending thereby the pencil case first established
by Meis 1960)

$\bullet$ Kleiman-Laksov 1972--74 \cite{Kleiman-Laksov_1972}
\cite{Kleiman-Laksov_1974} (using resp. Schubert calculus, plus
Poincar\'e's formula  and resp. singularity theory \`a la Thom,
Porteous)

$\bullet$ Gunning 1972 \cite{Gunning_1972} using MacDonald
computation of the homology of the symmetric power of the
curve

$\bullet$ Griffiths-Harris 1978
\cite[p.\,261]{Griffiths-Harris_1978/94}, where the heuristic
count \`a la Riemann-Klein is reproduced; and latter a rigorous
argument (p.\,358) is supplied (along the line of Kempf's Thesis
ca. 1970).

In view of the interest aroused by this Riemann-Meis bound, and
the apparent difficulty to prove it (appealing to a variety of
ingenious devices), it seems reasonable to wonder if there is not
a much simpler argument based upon the same ``irrigation method''
as the one used by the writer in relation with the Ahlfors map.
This would merely use simple homology theory and the allied
surjectivity criterion in term of the Brouwer degree.
Heuristically, this amounts to see the genus $g$ pretzel inside
its Jacobian and let it homologically degenerate over a bouquet
(wedge) of $g$ $2$-tori irrigating the Jacobian. Thus it seems
evident that with roughly $g/2$ points we may find a pair of
(effective) divisor of that degree collapsing to the same point of
the Jacobian. This pair of disjoint divisors
 serves to define the desired morphism to ${\Bbb P}^1$.
The writer as yet did not found the energy to write down the
details, but is
quite confident that the strategy is worth paying attention. Of
course it could be the case that this merely boils down (up to
phraseological details) to the already implemented attack of
Gunning 1972 \cite{Gunning_1972}. (Shamefully, I did not yet had
the time to consult this properly.)
%
%
%
%
%
%
Of course ``irrigation'' would not establish the sharpness of
Meis' bound (which is another question), but could predict its
value as universal upper-bound upon the gonality.

\subsection{Another application: Complex manifolds
homeomorphic to tori}

This section deviates from the mainbody of the text, but serves to
illustrate another spinoff of the irrigation method. The writer
wondered about the following naive question (ca. 2001/2?). Assume
given a complex (analytic) manifold (arbitrary dimension), and
suppose also the underlying  manifold to be homeomorphic to a
torus. {\it Must such a manifold be biholomorphic to a complex
torus, i.e. ${\Bbb C}^n$ modulo a lattice?} The answer is easy in
dimension one (Abel essentially). In general the answer is
negative, by virtue of a construction of Blanchard (Thesis ca.
1955) closely allied to  the Penrose twistor. Basically there is
over $S^2$ a certain bundle parametrizing quaternionic structures,
and taking a fiber product with an elliptic curve yields on the
torus $T^6$ (of 6 real dimensions) a complex structure which turns
out to be not K\"ahler. This answers negatively the question when
the complex dimension is 3. (For more details cf. also work by
Sommese (ca. 1978), etc.)

All this is rather exotic complex geometry, but one may wonder if
the assertion becomes true under the  K\"ahler assumption. Then
Hodge theory applies, and we dispose of a bona fide analog of the
Abel mapping (sometimes called the Albanese mapping). The latter
is also a map to a complex torus (called Albanese) and using the
irrigation principle it is easy to show that $\alpha$ induces an
isomorphism on the top-dimensional homology. First, it induces an
isomorphism on the $H_1$, but the latter elevates up to the
top-dimension since tori have a total homology $H_{\ast}$ modelled
upon the exterior algebra over the $H_1$. By the Brouwer degree
argument (irrigation intuitively), it follows that $\alpha$ is
surjective. Then one can show that it is injective as well (I have
forgotten the exact argument, but essentially if Albanese collapse
a submanifold then like by Abel it collapses linear varieties
which are simply-connected projective spaces, hence liftable to
the universal cover of the Albanese torus).

\begin{lemma}
Any torus shaped K\"ahler manifold is biholomorphic to its
Albanese torus.
\end{lemma}

Of course this is surely well-known, but we just wanted to
remember this as another high dimensional---but baby---application
of the irrigation principle. Further Kodaira's classification of
(complex analytic) surfaces plus a deformation argument of
Andreotti-Grauert (which I learned from R. Narasimhan) implied
also a positive answer to the basic question in (complex)
dimension 2. But I take refuge in my failing memory, and to not
remember the exact details. Thus in principle, Blanchard's
3-dimensional counterexample is sharp.

\subsection{Invisible real curves (Witt 1934, Geyer 1964,
Martens 1978)}\label{sec:Witt}

Ahlfors' theorem bears some analogy with Witt's theorem (1934
\cite{Witt_1934}) stating that a (smooth) real curve without real
points admits a morphism (defined over the reals ${\Bbb R}$) to
the invisible real line (materialized by the conic
$x_0^2+x_1^2+x_2^2=0$). The analogy is again that when there is no
topological obstruction, then a geometric mapping exists.

Subsequent works along Witt's direction are due to:

$\bullet$ Geyer 1964/67 \cite{Geyer_1964-67} (alternative proof of
Weichold, and Witt via Galois cohomology and Hilbert's Satz 90);

$\bullet$ Martens 1978 \cite{Martens_1978}, where the precise
bound on the degree of the Witt mapping has been determined.

Philosophically, it seems challenging to examine if such strongly
algebraic techniques (Riemann-Roch algebraized \`a la
Hensel--Landsberg 1902, Artin, etc.) are susceptible to crack as
well the
Ahlfors mapping?
Geyer, Martens or others are perhaps able to
address this challenge?
(So far as we know,
no such account exist in print.)

Martens' statement (quantitative version of Witt) is the
following.

\begin{theorem}[{\rm Martens 1978 \cite{Martens_1978}}]
Given a closed non-orientable Klein surface with algebraic
genus $g$ (i.e. the genus of the orientable double
cover\footnote{This is not explicitly specified in the paper,
but is the (common) jargon in Klein surface theory, probably
due to Alling-Greenleaf 1971 \cite{Alling-Greenleaf_1971}.})
there is a morphism to the projective plane of degree $\le
g+1$. Moreover this is the best we can hope for, i.e. for each
$g$ there is a Klein surface not expressible with fewer
sheets.
%
\end{theorem}

Perhaps the first portion of the statement is already in Witt 1934
\cite{Witt_1934}.
Of course this can---via the Schottky-Klein
Verdoppelung---also be stated in term of symmetric Riemann
surfaces (equivalently real algebraic curves) as follows:

\begin{theorem}[Martens 1978 \cite{Martens_1978}]
Given a symmetric Riemann surface of genus $g$ without fixed
point, there is an equivariant conformal mapping to the
diasymmetric sphere of degree $g+1$. Moreover the bound is
sharp.
\end{theorem}

This formulation of Martens's result also appears in Ross 1997
\cite[p.\,3097]{Ross_1997}, who supplies additional comments which
are quite in accordance with our own sentiments, especially the
issue that the short argument by Li-Yau 1982
\cite[p.\,272]{Li-Yau_1982} does not appear as very convincing.
Moreover Ross supplies some attractive differential geometric
applications of this Witt-Martens mapping theorem, e.g. to the
effect that the totally geodesic ${\Bbb R}P^2$ is the only stable
minimal surface in ${\Bbb R}P^3$.

\subsection{The three mapping theorems (Riemann 1857, Ahlfors 1950,
Witt 1934)}

From the conformal viewpoint we have thus three basic mapping
theorems enabling a  gravitational collapse of all compact
surfaces to their simplest representatives (the sphere, the disc
or the projective plane) depending on whether the original surface
is:

$\bullet$ closed orientable (Riemann 1857 \cite{Riemann_1857});

$\bullet$ compact bordered orientable (Ahlfors 1950
\cite{Ahlfors_1950});

$\bullet$ closed non-orientable  (Witt 1934 \cite{Witt_1934}).

None of those results tells what to do with a compact bordered
non-orientable surface (whose simplest specimen is the M\"obius
band/strip). The latter does not carry positive curvature, which
implies finiteness of the fundamental group for complete metrics
(else punctured sphere). Alternatively the orientable double cover
of M\"obius is the torus, which has already moduli. Hence it is
quite clear that the above three theorems form an exhaustive list
of truths positing  a fundamental trichotomy (of definitive
crystallized shape).
The motto ``Alle guten Dinge sind drei'', is quite ubiquitous in
life and mathematics!
It is reasonable to expect that each of those mappings will pursue
to find valuable applications in the future, yet much work remain
to be done as to the stratification of the moduli space induced by
the degree of such representations, etc.

For each of these 3 concretization problems one is interested in
the exact determination of the lowest possible sheet number. In
principle the answer is already known as follows:

\begin{theorem}
For all $3$ types of conformal mapping to elementary surfaces of
positive Euler characteristics $\chi >0$ (including $\chi (S^2)=2,
\chi (\Bbb R P^2)=1, \chi (\Delta=D^2)=1$) the sharp universal
bound on the degree of such representation is known. More
precisely,

$\bullet$ $[\frac{g+3}{2}]$ always concretizes closed genus $g$
surfaces expressed as cover of the sphere {\rm (Riemann, Meis 1960
\cite{Meis_1960})}, and the bound is sharp {\rm (again  Meis 1960
\cite{Meis_1960})}.

$\bullet$ $g+1$ always concretizes  non-orientable closed surface
of algebraic genus $g$ (i.e. genus of the orientation double
cover) expresses as cover of the projective plane {\rm (Witt 1934
\cite{Witt_1934})}, and the bound is sharp {\rm (Martens 1978
\cite{Martens_1978})}.

$\bullet$ $r+p$ always concretizes bordered orientable surfaces
with $r$ contours and  $p$ handles as (full or total) cover of the
disc {\rm (Gabard 2006 \cite{Gabard_2006})}, and the bound is
sharp {\rm (Coppens 2011 \cite{Coppens_2011})}.
\end{theorem}

Adhering to Klein's viewpoint of symmetric surfaces, one can
always interpret such objects as real curves of some genus $g$
(the first class is an exception except if one tolerates
disconnected surfaces). In the third bordered case $g=(r-1)+2p$.
The $r+p$ bound can be rewritten as $\frac{r+(g+1)}{2}$. If $r$ is
lowest, i.e. $r=1$, this is statistically equal to $g/2$, as so is
the first Riemann-Meis bound. In contrast the Witt-Martens bound
looks much higher. Of course if $r=g+1$ is highest
(Harnack-maximality) then $r+p=r+0=g+1$, agreing with
Witt-Martens's bound. In the overall it may be argued that both
Martens' and Gabard's bound are fairly less economical that
Riemann-Meis', and that this is due to the equivariance or even
total reality of the corresponding maps. On the other hand Ahlfors
bound $r+2p=g+1$ looks much more compatible with Martens' and if
one is sceptical about Gabard's version one could imagine  that
Ahlfors is asymptotically sharp for  large values of the
invariants. This scenario remains hypothetically possible in case
we are unable to reassess through other mean Gabard's $r+p$ or
able to disprove its validity.

 The following
tabulation summarizes the key contributions:

(1) Riemann 1857: any (or at least the general) closed
Riemann(ian) surface maps conformally to the sphere with $\le
[\frac{g+3}{2}]$ sheets, where $g$ is the genus. It is not
clear-cut
 if Riemann showed sharpness of the bound.

Related works includes (in chronological order):

$\bullet$ Brill-Noether 1874 \cite{Brill-Noether_1874};

$\bullet$ Klein 1891
\cite[p.\,99]{Klein_1891--92_Vorlesung-Goettingen};

$\bullet$ Severi 1921
\cite{Severi_1921-Vorlesungen-u-alg.-Geom-BUCH};

$\bullet$ B. Segre 1928 \cite{Segre_1928};

$\bullet$ Meis 1960 \cite{Meis_1960};

$\bullet$ Kempf 1971 \cite{Kempf_1971} and Kleiman-Laksov 1972--74
\cite{Kleiman-Laksov_1972} \cite{Kleiman-Laksov_1974};

$\bullet$ Gunning 1972 \cite{Gunning_1972};

$\bullet$ Griffiths-Harris 1978
\cite[p.\,261]{Griffiths-Harris_1978/94};

$\bullet$ Arbarello-Cornalba 1981 \cite{Arbarello-Cornalba_1981}.

This sharp bound $[\frac{g+3}{2}]$ as applied to spectral theory
is observed in El Soufi-Ilias 1983/84
\cite{El-Soufi-Ilias_1983/84} (Yang-Yau 1980 \cite{Yang-Yau_1980}
contented themselves with the weaker value $g+1$.) An interesting
aspect of the Italian works is that they not only focus on the
gonality upper bound, but also compute the dimensions of the lower
dimensional strata for a prescribed gonality. Of course, the
answer is the expected one (as easily predicted by
Riemann-Hurwitz). [The above Italian works, especially Segre has
however a little objection to the simplicity of the exercise.] We
point out this is issue as it could be interesting to make a
similar count for the Ahlfors circle map (bordered case). This
topic will be briefly addressed in the next
Sec.\,\ref{sec:profile-histogram}.

(2) Ahlfors 1950 \cite{Ahlfors_1950}: any compact bordered Riemann
surface maps conformally to the disc with $\le r+2p$ sheets (where
as usual $r$ is the number of boundary contours and $p$ the
genus). This bound is not sharp (at least for low values of the
invariants $(r,p)$, e.g. for the G\"urtelkurve type
$(r,p)=(2,1)$). Modulo a mistake by the writer (in Gabard 2006
\cite{Gabard_2006}), Ahlfors bound can be improved as $\le r+p$.
The latter is in turn sharp according to Coppens 2011
\cite{Coppens_2011}.

(3) Witt 1934 \cite{Witt_1934}: any closed non-orientable surface
maps conformally to the projective plane ${\Bbb R}P^2$. Witt does
not specify a bound (?), or maybe he does but sharpness was
obtained by Martens 1978 \cite{Martens_1978}.

Witt's result received, arguably, only sporadic spectral
applications, except in the article  Li-Yau 1982
\cite{Li-Yau_1982}, which however does not quote Witt, but whose
authors were apparently able to reprove the result by their own
[compare their argument on p.\,272].
(As already mentioned, Ross 1997 \cite{Ross_1997} does not seem to
be convinced by the Li-Yau argument.)

Of course all this ``diaporama'' is the direct heritage of Riemann
(plus maybe indirectly some Abel!) the first result being often
called Riemann's existence theorem. The 2 avatars of Ahlfors and
Witt are akin to the absolute case of Riemann, via the trick of
the Schottky-Klein double (or Verdoppelung as Teichm\"uller calls
it) but then some equivariance or total reality is required,
acting as a sort of boundary condition explaining probably why
those versions took longer to emerge. Of course such
equivariance/or boundary behaviors just
hide
a reality condition encoded in the field of definition of the
allied Riemann surfaces. All this is best summarized
diagrammatically:

\begin{figure}[h]
\centering
    \epsfig{figure=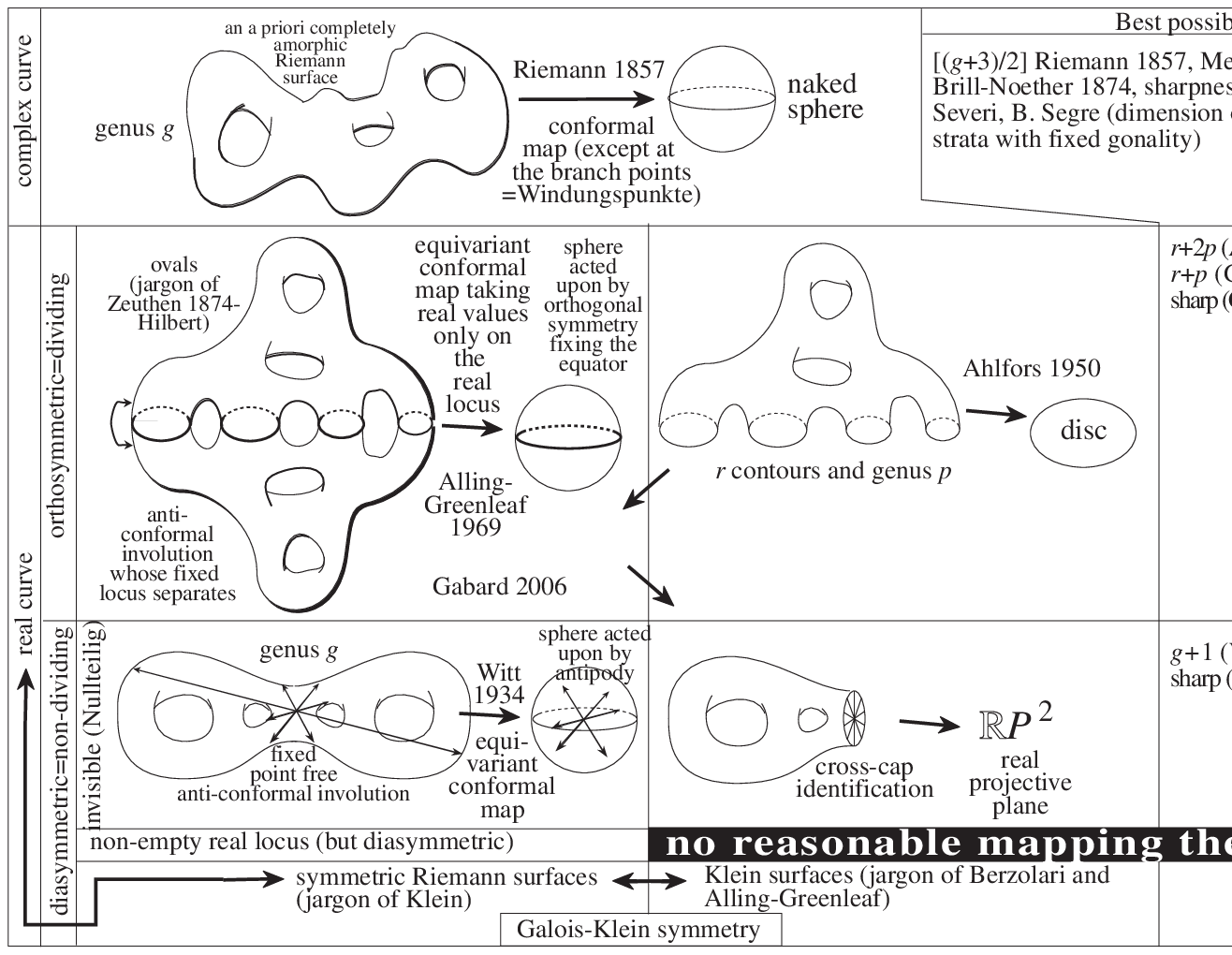,width=126mm}

\caption{\label{RAW:fig} The 3 types of conformal maps to the
simplest configurations ($\chi>0$)}
\end{figure}

\subsection{The gonality profile, moduli strata
and the Ahlfors space}\label{sec:profile-histogram}

[10.11.12] Heuristically (cf. e.g. Segre 1928 \cite{Segre_1928})
one can count the dimension of closed Riemann surfaces expressible
as coverings of degree $d$ of the sphere as follows. By
Riemann-Hurwitz $\chi(C_g)=d \chi(S^2)-b$. Hence there are $b-3$
free parameters, that is
$$
b-3=(2d-\chi)-3=2d-(2-2g)-3=2d+2g-5.
$$
In particular the strata of given gonality $\gamma=d$ increases
its dimension along a very simple arithmetic progression (as a
function of $d$) until the full moduli space is exhausted for $d$
the least integer $\ge g/2+1$ (Riemann-Meis bound). The smallest
strata is the {\it hyperelliptic locus\/} ($d=2$) of dimension
$4+2g-5=2g-1=(2g+2)-3$, in accordance with the $2g+2$ ramification
points visible as invariant points of an half twist acting upon a
purely vertical pretzel in 3-space. I do not know if such a
regularity occurs for bordered surfaces. Coppens's theorem states
another form of regularity, namely full realizability of all
intermediated gonalites, but it does not pertain to the dimensions
of the corresponding moduli strata.

On behalf of Coppens's theorem the situation could be as follows.
For a given topological type $(r,p)$, Coppens tells us that all
intermediate $r\le \gamma \le r+p$ are realized. So we have $p+1$
possible gonalities, the largest of which $\gamma=r+p$ fills the
full moduli space of real dimension $3g-3$ (Klein's count
conjointly with Gabard's bound). As usual $g=(r-1)+2p$, so
$3g-3=6p+3r-6$. If we knew the number of moduli of the minimal
strata $\gamma=r$ we could try a linear interpolation as a
possible scenario for the dimensions increments of the gonality
strata. Naively our rotationally invariant picture
(Fig.\,\ref{Chambery:fig}) could act as a bordered substitute to
the hyperelliptic closed case (at least for $r$ even). If so is
the case can we count its moduli? Everything would be determined
by the quotient planar surface with $r/2=r'$ contours. This planar
surface (whose double has genus $g'=r'-1$) depends on $3g'-3$
moduli. This expressed in terms of $r$, gives he following
$3g'-3=3r'-6=3/2 \cdot r-6$. This a candidate for the dimension of
the lowest strata. Looking for a progression in $p$ steps toward
the maximum value, we consider the difference $[6p+3r-6]-[3/2
\cdot r-6]=6p+3/2\cdot r=1/2[12p+3r]$, which is not easily divided
by $p$. $\bullet$ In fact we have looked at the quotient but
barely omitted the branched locus. Taking this into account we get
rather a dependance on $3g'-3+2(2p+2)$ (real) moduli for the
lowest strata. Expressing this in terms of $(r,p)$, gives $3/2
\cdot r+4p-2$. Hence the difference of the top and lowest strata
would be $2p+3/2\cdot r-1$, which is alas still not nicely
divisible by $p$. $\bullet$ Another idea is just to use maps from
$F_{r,p}$ to the disc of minimum degree $r$. Then we have
$\chi(F)=r \chi(\Delta)-b$. Hence there are $2b-3$ free real
parameters. Expressed in terms of $(r,p)$, this is
$2b-3=2(r-\chi)-3=2(r-(2-2p-r))-3=4r+4p-7$. Hence the difference
between the top dimensional and the lowest dimensional strata is
$\delta=(6p+3r-6)-(4r+4p-7)=2p-r+1$, which is not even positive in
general. It looks again dubious to divide this in $p$ equal parts
as suggests Coppens result. Again this just confirms what we
already noticed (earlier in the text) that the Riemann-Hurwitz
count looks seriously jeopardized in the bordered case, at least
as long as we apply it so naively as we do.

 One can reverse the game: instead of speculating on the size of
the lowest strata we can speculate on the increment as being by 2
real units (like in the complex case) and draw the dimension
$\lambda$ of the lowest strata. This would give $\lambda=\dim
({\rm top \;strata})-p\cdot 2=(6p+3r-6)-2p=4p+3r-6$. Testing this
on the type of the G\"urtelkurve $(r,p)=(2,1)$ gives
$\lambda=4+6-6=4$, whereas the hyperelliptic model depends on
$2g+2-3=2\cdot3+2-3=5$ real parameters. Hence the later has
codimension 1 in the full moduli of the G\"urtekurve type, which
as dimension $3g-3=3\cdot 3-3=9-3=6$.
 This motivates modifying the
increment to one of only 1 unit. This leads to  the following
Ansatz: $\lambda=5p+3r-6$. This gives for $(r,p)=(2,1)$,
$\lambda=5+6-6=5$ the correct number. But if we look at the type
$(r,p)=(2,2)$ we get $\lambda =10+6-6=10$; but on the other hand
the hyperelliptic models have $2g+2-3=2\cdot 5+2-3=9$ moduli
conflicting the new Ansatz for $\lambda$.

Of course the real scenario about the increments might be pretty
more complicated than the linear progression observed in the
complex case (corresponding to closed Riemann surface).

Another more neutral way to look at the question is as follows.
Given is $(r,p)$ a pair of integers. Allied to this there is a
moduli space ${\cal M}_{r,p}$ of all bordered (Riemann) surfaces
of type $(r,p)$. Its dimension is $3g-3$ (Klein 1882
\cite{Klein_1882}), where $g$ is the genus of the double. We
imagine the range of all possible gonalities $r\le\gamma\le r+p$
as a horizontal array of entries above each of which is reported
the dimension of the moduli space of curve having gonality $\le
\gamma$. This is depicted as a vertical bar. At first, only the
top dimension attached to $\gamma=r+p$ is known as $3g-3$. By
Coppens we know that there will be $p$ descents of this highest
bar to the lower gonalities between $r$ and $r+p$. Pause at this
stage to notice that assigned to the sole data $(r,p)$ there is
assigned unambiguously such a histogram of gonalities (cf.
Fig.\,\ref{Histogr:fig}).

\begin{figure}[h]
\centering
    \epsfig{figure=Histogr.eps,width=56mm}

\caption{\label{Histogr:fig} Histogram encoding the dimensions of
each gonality strata (alias the gonality profile)}
\end{figure}

One special case in which we can hope to be more explicit
regarding the lowest strata is when $r$ equals 1 or 2. In this
case we know that the moduli space contains hyperelliptic
membranes. Assuming $p$ large enough ($p\ge 1$) the lowest
gonality is $\gamma=2$. It is tautological that the hyperelliptic
locus has this gonality, and conversely. So we control explicitly
the dimension of the lowest strata. We find $(2g+2)-3$ real
constants. Thus the dimension difference $\delta$ of the top and
lowest strata is $\delta=3g-3-[(2g+2)-3]=g-2$. This rewritten in
terms of $(r,p)$ is also $g-2=(r-1)+2p-2=r+2p-3$.

$\bullet$ If $r=1$, this gives $\delta=2p-2=2(p-1)$. Positing
linearity of the increment, this ought to be divided in $p-1$
equals parts (since $r=1$ itself is not a gonality when $p\ge 1$),
and we get exactly a progression by 2 units. Hence under the
Ansatz of linearity the histogram would be completely known.

$\bullet$ If $r=2$, this gives $\delta=2p-1$. Assuming linearity
of the increment, this ought to be divided in $p$ equals parts,
and we get something like a progression by 2 units. However the
non-divisibility implies that in this case it is impossible to
have a linear progression of the histogram. Hence some jumps must
occurs.

So in these cases there is some hope to be completely explicit
about the histogram attached to $(r,p)$. It would essentially
suffices to decide where occur some irregular jumps.

Let us formalize a bit. Given a pair of integers $(r,p)$, we have
a moduli space ${\cal M}:={\cal M}_{r,p}$ of all bordered Riemann
surfaces of type $(r,p)$. (To allege notation with omit the
indices $(r,p)$, as the topology is fixed once for all.) Its
dimension is invariably $3g-3$, where $g=(r-1)+2p$ is the genus of
the double.

\begin{defn}\label{gonality-profile:def}
{\rm Inside the full moduli space ${\cal M}:={\cal M}_{r,p}$,
consider the sublocus $M_d$ of all surfaces with gonality
$\gamma\le d$, and let $\mu_d=\dim M_d$ be its dimension. The
histogram we were speaking about is essentially the function
$d\mapsto \mu_d$, which we call the {\it gonality profile}.}
\end{defn}

It is evidently monotone but a priori not strictly.
Misinterpreting Coppens's result one would guess strict monotony,
but Coppens states only that each gonality is exactly realized,
hence in symbols that $M_d-M_{d-1}$ is non-void (at least for $d$
in the range $[r,r+p]$). Thus a priori it could be the case that
when incrementing the parameter $d$ we get new surfaces but their
variety is not of larger dimension. Of course this scenario may
look a bit unlikely due to the algebro-geometric character of the
whole topic, but I do not know an argument. The domain of our
function $d\mapsto \mu_d$ is the set of all integers but the
interesting range is $[r, r+p]$ at least taking Gabard for
granted. The latter amounts to say that $\mu_{r+p}=3g-3$.

Now if $r=1$ or $2$, then the moduli space ${\cal M}={\cal
M}_{r,p}$ contains hyperelliptic representatives, and the latter
exhaust the locus $M_2$. We calculate easily $\mu_2=(2g+2)-3=2g-1$
and deduced the difference
$\delta=\mu_{r+p}-\mu_2=(3g-3)-(2g-1)=g-2$. From here we inferred
that:

$\bullet$ when $r=1$ (and $p\ge 1$) then $r$ itself is not a
gonality and so there is really only $p-1$ descents. Since
$\delta=g-2=r+2p-3=2p-2=2(p-1)$, we can divide (without rest) this
by the number of $p-1$ descents, to get a statistical increment of
2 units. If one believes in the linearity regularity then the
histogram would be completely known in that case.

$\bullet$ when $r=2$ then $r$ is  a gonality, and we have exactly
$p$ admissible descents along the range $[r,r+p]$. Now
$\delta=g-2=r+2p-3=2p-1$, which is not divisible by $p$. We infer
an obstruction to the scenario of linearly evolving histogram. (In
other words the function is not linear on the segment $[r, r+p]$.)
Perhaps it is just doing a gentle seesaw at some early place?

At this stage we may have exhausted all what can be said on
trivial arithmetical grounds. Going further probably requires some
geometric impetus, like looking at explicit models (extending the
hyperelliptic case). So one needs probably to describe large
families of $d$-gonal surfaces for $d\ge 2$. If a general result
describing the gonality profile $d\mapsto \mu_d$ looks out of
reach, one can start examining low values of $(r,p)$ to explore
the situation.

[11.11.12] {\it Examples.}---$\bullet$ E.g. for $(r,p)=(2,1)$
(thus $g=3$) (the G\"urtelkurve type) then the profile is
completely known, namely $\mu_2=5$ (hyperelliptic locus of
dimension $(2g+2)-3$) and $\mu_3=6$ (equal to $3g-3$). Of course
Gabard's $\gamma\le r+p$ follows in this case via the canonical
embedding realizing the curve as a G\"urtelkurve in ${\Bbb P}^2$.

$\bullet$ For $(r,p)=(2,2)$ (thus $g=5$), we have again the
hyperelliptic locus giving $\mu_2=(2g+2)-3=9$. The top locus
$M_{r+p}=M_4$ has dimension $\mu_4=\mu_{r+p}=3g-3=15-3=12$ (Gabard
is used but maybe there is an argument by hand). What about
$\mu_3$? To seek an answer we refer back to the table of
Fig.\,\ref{Coppens:fig}, where we traced a picture (label 223) of
an uninodal quintic with gonality $\gamma=3$. Quintics depends on
$\binom{5+2}{2}-1=\frac{7\cdot 6}{2}-1=20$ parameters, but modulo
the collinearity group $PGL(3)={\rm Aut}({\Bbb P}^2)$ of $3^2-1=8$
dimensions, this boils down to 12 effective parameters. Of course
the uninodal quintic we consider is really compelled
to live on the smaller discriminant hypersurface of dimension $19$
and so our curve 223 truly depends on only 11 essential
parameters. Assuming that a full neighborhood of curve $223$
consists of curves keeping the same gonality $\gamma=3$ suggests
therefore the value $\mu_3=11$ (at least as a lower bound).
Observe that the picture 223 is total under a pencil of lines, and
it seems reasonable to expect that when the curve is
 slightly perturbed along the discriminant hypersurface,
total reality of the pencil persists on the ground of some
topological {\it stability}. Remember e.g., that total reality
amounts to the transversality of the foliation (induced by the
pencil) along the curve, and  transversality is the mother of any
topological stability (Thom-style philosophy). Note of course that
our curve (being uninodal) represents actually a smooth point of
the discriminant and so we safely dispose of the required
parameters of deformation. This is perhaps worth saying if one
remembers certain plane cubics (or even conics) as examples of
real algebraic varieties having an isolated real point. Maybe the
above stability argument adapts to situations where there are
several nodes via  Brusotti's theorem describing the infinitesimal
structure of the discriminant near a multi-nodal curve (with say
$\delta$ nodes) as an union of smooth branches crossing
transversally (normal crossing).

$\bullet$ For $(r,p)=(3,2)$ (thus $g=6$), we have no hyperelliptic
locus. The top locus $M_{r+p}=M_5$ has dimension
$\mu_5=\mu_{r+p}=3g-3=18-3=15$ (Gabard is used but maybe there is
an argument by hand). What about $\mu_3$ and $\mu_4$? We look
again back to Fig.\,\ref{Coppens:fig}, where we find curve 324.
This is merely a smooth quintic with 2 nested ovals hence with
gonality $\gamma=4$. Remember that smooth plane $m$-tics have in
general complex gonality $(m-1)$. As quintics depends on 12
essential parameters, the above stability argument shows that the
strata $M_4$ contains the locus of all such quintics, and we infer
$\mu_4\ge 12$. Is this an equality? How to estimate $\mu_3$? Due
to time limitation, we have to leave all this (in our opinion)
exciting topic at a fragmentary stage. Perhaps a last word, if we
use picture 324bis (still on Fig.\,\ref{Coppens:fig}), which is a
sextic with 4 nodes also having $\gamma\le 4$ (projection from the
node), then we get a model depending on
$\binom{6+2}{2}-1=\frac{8\cdot 7}{2}-1=27$ constants, minus the
$8$ coming from $PGL(3)$ gives $19$, of which must be subtracted 4
units (using  Brusotti's normal crossings description). The final
result is 15. Repeating the above stability argument implies that
$\mu_4\ge 15$. This is a much stronger lower estimate, which in
fact must be an equality since we have already attained the
dimension of the full moduli space. Hence we conclude $\mu_4=15$;
strikingly as big as $\mu_5$! This answer is quite intriguing in
case it is correct at all? It would show that the gonality profile
does not need to be strictly increasing! [11.12.12] Alas all of
our counts are completely erroneous in view of some basic examples
shown in the next section. Of course the mistake is that not all
birational (conformal) equivalence giving rise to Riemann's moduli
space need to be induced by a collineation as an automorphism of
the ambient plane ${\Bbb P}^2$.

[11.11.12] Finally, it is perhaps fruitful to keep a view on the
space of all (total) circle maps. This is the {\it Ahlfors space}
(improvised jargon) quite akin to so-called Hurwitz spaces. All
what we were concerned with in this subsection is arguably just a
shadow of this larger space dominating the moduli space ${\cal
M}_{r,p}$. Precisely,  the {\it circle maps (or Ahlfors) space\/}
$C_{r,p}$ consists for a fixed pair $(r,p)$ (number of contours
and handles resp.) of all circle maps $f\colon F \to \Delta$ on a
``variable'' bordered Riemann surface of specified topological
type $(r,p)$. Forgetting the circle map $f$ induces a natural map
$C_{r,p}\to {\cal M}_{r,p}$ to the moduli space. Of course one
must consider the space $C$ modulo the equivalence relation of a
conformal diffeomorphism commuting with the maps to the disc. The
strata $M_d$ of all surfaces of gonality $\gamma\le d$ appear then
as the projections of the fibres of the degree function on
$C_{r,p}$. The fibre of the map $C\to M$ (indices omitted) is the
space of all total maps on a fixed bordered Riemann surface $F$.

\subsection{Correcting the previous section}

[11.12.12] There are many counterexamples to our naive moduli
count via plane nodal models. For instance considering curves of
$g=2$, and using the projective realization as a quartic with one
node, we get the dependence upon $\dim \vert 4H
\vert=\binom{4+2}{2}-1=\frac{6\cdot 5 }{2}-1=14$ parameters, of
which must be subtracted one unit to be on the discriminant (due
to the presence of the node) and finally one has to remove the $8$
dimensions of $PGL(3)$. The end result is $14-1-8=5$, which
exceeds by $2$ units Riemann's $3g-3=3$. Of course this excess is
due to the fact that we moded out only by (linear) automorphisms
of the plane whereas there might be more mysterious birational
equivalence relating to configurations of our family of uninodal
quartics.

This demonstrates that the estimate we got in the previous section
are completely erroneous and unreliable, and one must find some
completely new ideas (i.e. old stuff \`a la Riemann) to really
penetrate the intrinsic nature of the problem. For the moment I
have no idea on how to attack the problem of describing the size
(=dimensions) $\mu_d$ of the varied gonality strata
$$
M_d=\{ F\in {\cal M}_{r,p} : \gamma (F)\le d \}.
$$
Let us try anew to attack this problem of describing the gonality
profile $d\mapsto \mu_d$ for each pair $(r,p)$.

[12.12.12] First complete information is obtained in the easy case
of planar membranes ($p=0$) as a consequence of the
Bieberbach-Grunsky theorem (Lemma~\ref{Enriques-Chisini:lemma}).

\begin{lemma}
When $p=0$ the gonality profile is a skyscraper concentrated at
the single place $d=r=r+p=r+2p$, i.e. $\mu_r=3g-3$, where
$g=(r-1)+2p=r-1$ is the genus of the double.
\end{lemma}

\begin{proof}
This follow at once from the trivial lower bound $r\le \gamma$ on
the degree of circle maps (or the allied gonality $\gamma$), and
the Bieberbach-Grunsky theorem
(Lemma~\ref{Enriques-Chisini:lemma}). (Notice that neither Ahlfors
1950 \cite{Ahlfors_1950} ($\gamma \le r+2p$), nor Gabard 2006
\cite{Gabard_2006} ($\gamma \le r+p$) is required.)
\end{proof}

After that let us examine the cases with $p=1$. We start with:

$\bullet$ $(r,p)=(1,1)$: then we have $1=r\le \gamma\le r+p=2$
(using Gabard). However since the genus of the double is
$g=(r-1)+2p=2$, the curve is hyperelliptic and we may avoid
Gabard. The lower bound $r=1$ cannot be realized, since $p\neq 0$.
We deduce:

\begin{lemma}
For $(r,p)=(1,1)$, the gonality profile is a skyscraper
concentrated at $d=2$, i.e. $\mu_d=3g-3=3$ for $d=2$ and $\mu_d=0$
elsewhere $(d\neq 2)$.
\end{lemma}

$\bullet$ $(r,p)=(2,1)$ (with $g=(r-1)+2p=3$): then we have
$2=r\le \gamma\le r+p=3$ (using Gabard). However one can dispense
Gabard by using the canonical embedding taking the double of the
bordered Riemann surface to a G\"urtelkurve $C_4 \subset {\Bbb
P}^2$, i.e. a quartic with 2 nested ovals. This proves $\gamma \le
3$ (via projection from the inner oval). So $\mu_3=3g-3=6$. Of
course we have also a hyperelliptic locus, whose dimension is
$(2g+2)-3=5$, so $\mu_2=5$. This proves the:

\begin{lemma}
For $(r,p)=(2,1)$, the gonality profile is a ``twin tower''
concentrated at two places $d=2,3$, and $\mu_3=3g-3=6$ whereas
$\mu_2=5$ (all other $\mu_d$ are zero).
\end{lemma}

$\bullet$ $(r,p)=(3,1)$ (with $g=(r-1)+2p=4$): then we have
$3=r\le \gamma\le r+p=4$ (using Gabard). Without using Gabard, one
can look at the canonical model in ${\Bbb P}^{g-1}={\Bbb P}^{3}$
of degree $2g-2=6$. This is probably a complete intersection of a
cubic surface with a quadric, weighting bidegree $(a,b)=(3,3)$ on
the latter, hence of genus $g=(a-1)(b-1)=2\cdot 2 = 4 $ (the
expected value). One can then draw a picture by smoothing two
pairs of 3 lines distributed in each ruling. When the lines are
oriented in the most trivial way (each inducing the same integral
homology class on the torus ${\Bbb P}^1({\Bbb R})\times {\Bbb
P}^1({\Bbb R})$) we get a total map of degree $3$ by projection on
the factors of ${\Bbb P}^1\times {\Bbb P}^1$ (do a picture).
Taking (somewhat cavalier) Gabard for granted we get
$\mu_4=3g-3=9$. How to estimate $\mu_3$? Let us try several
strategies:

(1) {\it Extrinsic plane projective realizations}.---The naive
idea is to look at Fig.\,\ref{Coppens:fig} (picture 313). This is
a quintic with 2 nodes and $\gamma\le 3$ (hence equal to $3$ by
the trivial lower bound $r\le \gamma$). If we count the (naive)
moduli of such a curve we obtain: $\dim \vert 5H
\vert=\binom{5+2}{2}-1=\frac{7\cdot 6}{2}-1=20$, of which must be
subtracted $2$ for the two nodes, and $8=\dim PGL(3)$ to get $10$.
This exceeds by one unit the full moduli space $3g-3$, and so we
get an alienating count. Of course as already said the reason is
that we only took into account linear collineations (ambient
automorphisms) whereas one should mod out by all inherent
isomorphisms of the family of curves. One way to remedy the
situation would be to look at Cremona transformations (birational
transformations of the plane), but it is not even evident that
this would give the right answer on abstract moduli. Another idea
is to look at higher order plane models with $\delta$ many nodes
as to adjust the genus to $g=4$. For instance sextics with 6 nodes
$C^6_6$, septics with 11 nodes $C_7^{11}$, octics with 17 nodes
$C_8^{17}$, etc. However the same calculation shows that such
family of curves depends on $\dim C^6_6=\frac{8\cdot
7}{2}-1-\delta-8=13$, $\dim C^{11}_7=\frac{9\cdot
8}{2}-1-\delta-8=16$, $\dim C^{17}_8=\frac{10\cdot
9}{2}-1-\delta-8=19$, etc. It seems that there is perpetually an
increment by 3 units, and never get something realistic via naive
counting.

(2) {\it Extrinsic projective realization as a branched cover of
the line (or the disc), i.e. circle maps via
Riemann-Hurwitz.}---We fix as an Ansatz $\gamma=3$ (inside our
fixed topological type $(r,p)=(3,1)$), and by the Riemann-Hurwitz
relation applied to a circle map $F \to \overline{\Delta}$ of
degree $3$, we find $\chi (F)= 3 \chi (\overline \Delta)-b$, so
$b=3-\chi (F_{3,1})=3-(-3)=6$ many branch points. Moving those 6
points arbitrarily in the disc, and quotienting by automorphisms
of the disc we arrive at $2\cdot 6 - 3=9$ (real) moduli. This
looks again anomalous for we receive the same answer as for the
full $3g-3$ moduli. (We already experimented this failure of
Riemann-Hurwitz in the bordered setting, and we are in depressive
mode.) The mystery is perhaps that we do not enjoy complete
freedom in moving branch points in the bordered setting, but I
lack any understanding of which sort of geometric restrictions
have to be taken into account.

(3) {\it Intrinsic count \`a la Nielsen-Fenchel.}---Another
possible strategy, is to adapt the Nielsen-Fenchel count of moduli
via a decomposition in pants. Remember that this works at the
perfection to predict the dimension of the full moduli space (cf.
e.g. our Sec.\,\ref{Nielsen-Fenchel:sec} below). The idea would be
that if we prescribe a lower gonality then an appropriate
decomposition in pants (somehow calibrated on the circle map)
should predict the moduli dimension of the restricted class too.

Alas for $\gamma=3$, I do not really see how to proceed, but let
us first experiment the method on the simpler hyperelliptic case.

Consider e.g. a membrane with $(r,p)=(2,p)$, $p$ arbitrary, of
gonality $\gamma=2$ (hyperelliptic case). On drawing the
configuration, and decomposing it into pants invariant under the
hyperelliptic involution (visualized as a half-twist rotation) we
obtain Fig.\,\ref{Pants-hyper:fig} (left part).

\begin{figure}[h]
\penalty0 \centering
    \epsfig{figure=Pants-hyper.eps,width=82mm}
    \vskip-5pt\penalty0

\caption{\label{Pants-hyper:fig} Pants dissection  applied to (the
baby) hyperelliptic membranes}
\end{figure}

Introducing on the surface its uniformizing metric of constant
curvature $-1$ (alias hyperbolic metric), we count moduli as the
lengths of loops bounding pants affected by certain twist
parameters. We get (reading contributions from the top to the
bottom of the figure):
$$
2+2 \cdot 2 p-1=4p+1
$$
free parameters. Indeed the first term ($2$) arises from the top
loop (its length plus its twisting aptitude). Next we see $p$
shaded pants whose contours exhaust all junctures of the pants
decomposition. However all bottom parts of the shaded pants are
permuted via the hyperelliptic involution (half-turn rotation),
hence of the same length. So each shaded pants really contributes
for $2$ lengths each susceptible of a twist, whence the second
term ($2 \cdot 2 p$). As to the last term ($-1$), notice that the
very bottom contours of the surface have no gluing companion (to
be twisted with), so one unit must be subtracted. The announced
count follows.

On the other hand, such hyperelliptic curves depend (via a count
\`a la Riemann-Hurwitz) on $(2g+2)-3$ parameters where
$g=(r-1)+2p=2p+1$. Hence on $(2g+2)-3=[2(2p+1)+2]-3=4p+1$, in
accordance with the result as calculated via the pants method.

A similar count works for hyperelliptic membranes with
$(r,p)=(1,p)$, cf. right part of Fig.\,\ref{Pants-hyper:fig}. In
that case we obtain
$$
2+2\cdot 2(p-1)+1=4p-1
$$
moduli, and on the other hand $(2g+2)-3=[2(2p)+2]-3=4p-1$
parameters. Both counts are again in accordance.

Of course we could even dream that the pants dissection method
(cf. Fig.\,\ref{Pants:fig}, right part) affords yet another full
proof of Ahlfors circle maps, but this looks a bit tricky to
implement. Perhaps even more ambitious one could hope that pants
dissection affords a proof of the Forstneri\v{c}-Wold 2009
\cite{Forstneric-Wold_2009} desideratum that each finite bordered
surface embeds holomorphically in ${\Bbb C}^2$. (Notice that this
is much stronger (viz. complex analytic) than the conformal
embedding in $E^4$ prompted by the Garsia-R\"uedy-Ko theory as
implemented in Ko 1999 and 2001 \cite{Ko_2001}. At least the
latter shows that there is no conformal obstruction to the Gromov
conjecture/question (1999 \cite{Gromov_1999}) that any Riemannian
surface should isometrically embed in $E^4$.)

Coming back to our problem of calculating $\mu_3$ for membranes of
type $(r,p)=(3,1)$ we severely lack any reliable technique of
calculation.

\section{Existence of Ahlfors maps via the Green's
function (and the allied Dirichlet
principle)}\label{Green:sec}

All what follows is extremely classical, yet the writer
confesses to have
assimilated (the first steps of the argument) as late as the
[04.08.12]! First it is well-known that the solubility of the
Dirichlet problem (say on a bordered Riemann surface) is
tantamount to the existence of the Green's function $G(z,t)$
with pole at $t$, for each $t$. (Actually, we primarily need
that the former implies the latter.) This
``Dirichlet-to-Green'' mechanism will be recalled below
along with the definition and some geometric (biochemical)
intuition about the Green's function. The latter has also strong
electrostatic or hydrodynamic connotations. The definition of the
Green's function is somewhat easier in the case of plane domains,
and its extension to bordered surface---while still laying in the
range of Dirichlet---implicates
some
conceptual difficulties.

The {\it Green's function\/} $G(z,t)$ with pole at $t$ (a fixed
interior point) is a completely canonical function characterized
by the properties: it is harmonic off $t$, vanishes along the
boundary and its germ near has the singular behavior prescribed by
the function $\log\vert z-t\vert$ in any local uniformizer $z$. It
will be verified that $G(z,t)$ is negative on the interior of the
bordered surface (consequence of Gauss' mean value property of
harmonic function and the resulting maximum principle). Then we
shall try to approach the existence of the Ahlfors function by
duplicating the Green-type proof of the Riemann mapping theorem
(simply-connected case), which just amount to write down the magic
formula $f(z)=e^{G(z,t)+i G^{\ast}(z,t)}$, where $G^{\ast}$ is the
conjugate potential. Note that $G(z,t)\le 0$ ensures $\vert f(z)
\vert=e^{G(z,t)}\le 1$ with equality precisely along the boundary.
The main difficulty about extending this ``Green-to-Riemann''
trick to the multiply-connected setting is to arrange
single-valuedness of the conjugate potential $G^{\ast}$. This
amounts to kill all periods of the $1$-form $dG^{\ast}$ from which
$G^{\ast}$ arises through
line-integration. To achieve this one is invited to introduce
enough free parameters in the problem by considering a
superposition of various Green's functions $\sum_i \lambda_i G(z,
t_i)$ for several poles $t_i$ sufficiently abundant so as to
enable the  killing of all periods (via linear algebra). Since a
planar domain with $r$ contours has $r-1$ essential cycles (up to
homology) and attaching $p$ handles creates 2 new essential
cycles, we need annihilating $(r-1)+2p$ periods.
Taking one more pole (raising the total number to $r+2p$) supplies
enough parameters for linear algebra to ensure  existence of a
non-trivial solution in the kernel of the period mapping. This
prompts
(almost) the existence of an Ahlfors circle map of degree $r+2p$
(as predicted in Ahlfors 1950 \cite{Ahlfors_1950}). Alas, a
serious technical difficulty occurs, namely  ensuring the
positivity of all $\lambda_i$. Ignoring this issue,  any $r+2p$
points (in the interior) could be the zeroes of a circle map.
Presently, we lack a complete existence of an Ahlfors map through
this procedure. Of course it would be even more challenging to
arrive at Gabard's bound (mapping degree $\le r+p$) through  this
classical strategy (\`a la Green, Riemann, Grunsky, Ahlfors,
Kuramochi, etc.). In Riemann the trick of annihilating periods
appears of course very explicitly in the following jargon: ``{\it
so bestimmen da{\ss} die Periodicit\"atsmoduln s\"ammtlich $0$
werden.}'' (cf.\, e.g. Riemann 1857 \cite[p.\,122]{Riemann_1857}).
The core of Heins' argument 1950 \cite{Heins_1950} is also exactly
in this spirit and Heins seems able to
complete the program via consideration of convex geometry. Our
intention is first to recall the basic procedure, and we hope to
be able later to settle the positivity problem. A priori it is not
evident that the latter condition is always achieved for an
arbitrary selection of poles $t_i$ of Green's functions (which
will mutate into zeroes of the ``Riemann-Ahlfors map'' $f$ after
exponentiation).

[25.08.12] {\it Corrigendum.}---The above linear superposition
$\sum_i \lambda_i G(z,t_i)$ on Green's functions is maybe somewhat
too continuous in nature.
This may be seen by exponentiating and looking at the local
behavior of $f$. Near some $t_i$, $G(z,t)\sim \lambda_i G(z,
t_i) \sim \lambda_i \log \vert z \vert$ so that $\vert
f(z)\vert \sim \exp(\lambda_i \log\vert z\vert)=\vert z
\vert^{\lambda_i}$ so that $f$ has not the character of a
holomorphic function when $\lambda_i$ is not integral.

Another way to argue in the same sense is suggested by Ahlfors
1950 \cite[p.\,126--7, \S4.3]{Ahlfors_1950}. Assume that $f(z)$ is
a circle-map $f\colon F \to \overline \Delta$ with zeros at
$t_1,\dots, t_d$ (counted with multiplicities), then upon
post-composing with the  function $\log\vert z \vert$ (harmonic
off the origin) we get the function $\log \vert f (z) \vert$
harmonic on $F$ save at the $t_i$  where it has logarithmic poles.
Therefore this function must coincide with superposition
$G:=\sum_{i=1}^d G(z,t_i)$ of Green's potentials. Indeed, the
difference $\log \vert f (z) \vert-G$ is throughout $F$ harmonic
(cancellation of singularities) and vanishes along the border
$\partial F$, hence is identically zero. [NB: the above remark is
to be found in Ahlfors (\loccit), who (in our opinion) fails to
insist on the assumption that $f$ is a circle-map (i.e. $\vert f
\vert =1$ along the border), which is crucial to ensure that $\log
\vert f (z) \vert$ vanishes along the border $\partial F$.]

So given a circle-map $f$ with $d$ zeros $t_i$ we have the
formula
$$
\log \vert f (z) \vert=\textstyle\sum_{i=1}^d G(z,t_i).
$$
Conversely, given points $t_i$, we may consider the right-hand
side of the previous equation
\begin{equation}\label{Green-super:eq}
G:=\textstyle\sum_{i=1}^d G(z,t_i)
\end{equation}
and the following
formula will define a circle-map
$$
f(z)=e^{G+iG^{\ast}}
$$
provided $dG^{\ast}$ (the conjugate differential of $G$) has
all its periods integral-multiples of $2\pi$. (It follows
incidentally, that a circle-map is uniquely determined up to a
rotation by the geographic location of its zeros. This can
also be seen algebro-geometrically, by considering the
Schottky double, where the divisor of zeros $D$ becomes
linearly equivalent to its symmetric conjugate $D^{\sigma}$,
spanning together a pencil $g^1_{d}$ defining a
total morphism to ${\Bbb P}^1$ of degree $d$, cf. Lemme~5.2 in
Gabard 2006 \cite{Gabard_2006}.

The desired  integrality of periods resembles a {\it Diophantine
condition} (at least is qualified as a such by Ahlfors 1947
\cite[p.\,1]{Ahlfors_1947}), emphasizing from the outset the
relative difficulty of the problem. All of our  freedom relies on
dragging the
points $t_i$ through the surface $F$ hoping that for a lucky
constellation the $1$-form $dG^{\ast}$ acquires simultaneous
integrality of all its periods along $\gamma_1, \dots, \gamma_g$
the $g:=(r-1)+2p$ many essential $1$-cycles traced on $F$ (cf.
Fig.\,\ref{Green:fig}e).

As a personal trouble, $dG^{\ast}$ seems to have singularities
where $G$ does, but maybe they disappear. Bypassing this point,
Ahlfors' Diophantine problem (1947) looks well-posed and one may
hope a direct attack upon arranging integrality of all periods.
(Ahlfors 1950 \cite{Ahlfors_1950} (p.\,127) first reformulates the
condition in term of Schottky differentials and then switches
quickly to the extremal problem, so does not seem to  attack
directly the Diophantine question. In fact, its elementary proof
on p.\,124--126 follows a somewhat different route by constructing
a half-space map involving avatars of Green's function with poles
situated along the boundary. We shall come back to this
subsequently.)

{\bf Trying a direct attack.} Assuming the problem well-posed,
we can consider a period mapping
$$
\wp\colon R^d\longrightarrow {\Bbb R}^g \longrightarrow ({\Bbb
R} / 2\pi {\Bbb Z})^g=:T^g,
$$
where $R={\rm int} (F)$ is the interior of the bordered surface
$F$, and the first map takes the periods along the fixed basis of
the first homology $\gamma_i$ of the $1$-form $dG^{\ast}$
corresponding to the points $(t_1, \dots, t_d)\in R^d$ via formula
\eqref{Green-super:eq}. The second map is just the natural
quotient map.

Now one may hope to apply the usual surjectivity criterion for
a continuous map to a closed manifold (here $R^d\to T^g$)
saying that if the representation induced on the
top-dimensional homology of the target-manifold is non-zero
then the mapping is surjective. For definiteness we recall its
statement and short proof.

\begin{lemma} Let $f\colon X \to T$ be a continuous map
from a (topological) space $X$ to a (target) manifold $T$ of
dimension $n$, say. It is assumed that $T$ is closed (i.e.
compact borderless). It is also essential to assume that $T$
is a Hausdorff manifold. If the induced homomorphism $H_n(f)$
is non-zero, then $f$ is onto.
\end{lemma}

\begin{proof} One considers the map induced on the homology $H_n$
of dimension $n$ equal to that of the manifold $T$. If $f$
fails to be surjective,  it factors through the punctured
manifold $X \to T-\{t\}$ for some point $t$. Now it is  a
simple fact that the top-dimensional homology of a (Hausdorff)
manifold vanishes, so in particular $H_n(T-\{t\})$ is trivial.
By functoriality it follows that $H_n(f)=0$, violating our
assumption.
\end{proof}

In particular $0=(0,\dots, 0)\in T^g$ would be the image of
some $(t_1, \dots , t_d)\in R^d$ and the corresponding
potential $G$ given by \eqref{Green-super:eq} would have a
conjugate differential $dG^{\ast}$ meeting the Diophantine
requirement.

This strategy requires a good understanding of the mapping $\wp$
perhaps in the sense that when one pole $t_i$ is dragged along the
cycle $\gamma_j$ then the image winds once around the
corresponding factor of the torus $T^g$. Choosing $d=g$ and in the
K\"unneth factor of $H_g(R^d)$ the element $\gamma_1 \otimes \dots
\otimes \gamma_g$ which has  the correct weight $g$ so as to be an
element of $H_g(R^g)$ whose image would be the fundamental class
of the torus $T^g$. This would establish the surjectivity of $\wp$
for $d=g$. Alas, this is a bit too optimistic in the planar case
($p=0$). So our argument must be foiled at some place.
The reasonable result to be  expected is $d=g+1$ (like Ahlfors
1950 \cite{Ahlfors_1950}) and boosting the method upon choosing
$\gamma_1 \otimes\dots \otimes \gamma_{r-1} \otimes (\alpha_1
\star \beta_1)\otimes \dots \otimes(\alpha_p \star \beta_p)$ where
the $\alpha_i, \beta_i$ are the cycles winding around the handles
(cf. Fig.\,\ref{Green:fig}e) one may expect to achieve $d=r+p$ as
predicted in Gabard 2006 \cite{Gabard_2006}.

\subsection{Digression on Dirichlet (optional)}

The Dirichlet solution may be interpreted as the permanent
equilibrium state of temperature in a heat-flow conducting medium.
Arguably (physico-chemical intuition?), this phenomenology is
completely insensitive to the topology. Hence Dirichlet's problem
is always soluble whatever the topological complexity of the
bordered manifold is. One only requires a Riemannian metric to
give a good sense to the (Beltrami) Laplacian (or the allied mean
value property). Hence any metric bordered smooth manifold, say
compact to stay in the reasonable realm of finiteness is suitable
to pose and solve the first boundary value problem. [Remember
maybe that there is vast jungle of non-metric manifolds, those of
Cantor 1883 and Pr\"ufer 1922 being the most prominent examples,
but the latter do not enter the scene of function theory at least
in complex dimension 1.]
%
%
%
%
Hence Dirichlet makes sense also on non-orientable manifolds, but
the case of immediate interest is that of compact bordered Riemann
surfaces
({\it ipso facto\/} orientable). Solid existence proofs were
primarily devised  by H.\,A. Schwarz, alternating method (ca.
1870), etc. with many subsequent extensions, e.g. Nevanlinna 1939
\cite{Nevanlinna_1939}, several works of Ahlfors, H. Weyl 1940
\cite{Weyl_1940} (method of orthogonal projection), not to mention
Neumann, Poincar\'e, Korn-Lichtenstein, etc., cf. e.g. Neumann
1900 \cite{Neumann_1900}). Another source is Hilbert-Courant's
book cited e.g. for this purposes in Royden's Thesis 1950/52
\cite{Royden_1952}. [For those inclined toward modern
expressionism, there is surely a concept of ``Dirichlet space''
(Brelot, Beurling, Deny, etc.) which should englobe any bordered
Riemannian manifold and much more.]

In the appropriate Hilbert space, minimizing the Dirichlet
integral amounts to minimize the length of a vector lying on a
certain hypersurface $M$ corresponding to the boundary data
$f\colon\partial F\to {\Bbb R}$. A priori this hypersurface could
spiral around the origin impeding existence of a minimum or be
bumpy enough as to violate uniqueness. But one rather imagine it
to be a linear manifold implying a unique minimum of the distance
function (norm). Of course the hypersurface in question
(corresponding to a certain boundary prescription) is readily
shown to have linear character, as subtracting any member of it,
its translate through the origin identifies with the set of
functions vanishing along the boundary. The latter is vectorial,
being the kernel of a linear mapping (restriction to the
boundary). Dirichlet principle looks thus immediately imputable to
an Euclid-Hilbertian conception of space, yet with  difficulty
concentrating on the existence question of a member (=point) in
this hypersurface $M$ (i.e., of a function matching the boundary
prescription having with finite Dirichlet integral). As we know
Hilbert's solution primarily involved the compactness paradigm,
formalized as a such some few years later by Fr\'echet. The naive
minimization procedure is not fairly evident, and indeed plagued
by the counterexample of Hadamard 1906 \cite{Hadamard_1906}, and
the earlier one of Prym 1871 \cite{Prym_1871}.
Prym (1871 \loccit) describes a continuous function on the
boundary of the unit disc such that the Dirichlet integral for the
associated harmonic extension of the boundary function is
infinite. [The latter harmonic extension is known to exist
independently of the Dirichlet principle, e.g. on the ground of
Poisson's formula which solves Dirichlet in the disc-case.] Later
Hadamard (1906 \loccit) gave a similar example where any
(continuous) function matching the boundary data has infinite
Dirichlet integral. (Perhaps, any Prym data is also explosive in
the sense of Hadamard?)
The moral is quite
subtle to grasp: roughly the Dirichlet principle fails but not the
Dirichlet problem which is always uniquely soluble! Hilbert's
solution (ca. 1900 \cite{Hilbert_1900}, \cite{Hilbert_1901/04})
under special hypotheses (involving only the space and not the
boundary data?!) is certainly sufficient for the purpose at hand.
Hilbert's hypothesis where weakened in subsequent works by B. Levi
1906 \cite{Beppo-Levi_1906}, Fubini 1907 \cite{Fubini_1907},
Lebesgue 1907 \cite{Lebesgue_1907}, compare also the
historiography in Zaremba 1910 \cite{Zaremba_1910}). For practical
purposes (e.g. for the construction of the Green's function) one
can probably restrict attention to reasonable boundary data, as
those arising via geometric construction (e.g., the logarithmic
charge allied to the construction of the Green's function of a
plane smoothly bounded domain). Possibly, for tame boundary data
the original Dirichlet principle remains an efficient tool for a
direct variational treatment of the boundary value problem.

Alternatively, of Dirichlet-Riemann-Hilbert one may use the
classical but cumbersome alternating method of Schwarz (or
Neumann's variant) to solve the Dirichlet problem. To summarize we
need the result:

\begin{theorem} {\rm (Dirichlet, Riemann,
Schwarz 1870, Hilbert 1900, etc.)} Given a compact bordered
Riemann surface $F$, and a continuous boundary function $f\colon
\partial F \to {\Bbb R}$. There is a unique harmonic function
$u\colon F \to {\Bbb R}$ extending $f$.
\end{theorem}

\begin{proof} First (rigourously) obtained in
Schwarz 1870 \cite{Schwarz_1870-alternirendes-Verfahren} via the
alternating method. Variation of this technique Picard's method of
successive approximation (cf. Picard, Zaremba, Korn,
Lichtenstein). Another variant of proof is Hilbert's resurrection
of the Dirichlet principle (direct variational method). Reference
in book form cf. Hilbert-Courant.
Another more modern trend is to use Perron's method which affords
great simplification. Compare for instance Ahlfors-Sario 1960
\cite[p.\,138--141, esp. 11G]{Ahlfors-Sario_1960} for an execution
of Perron's method (joint with Harnack's principle) in the context
of abstract Riemann surfaces.
\end{proof}

\subsection{From Green to  Riemann}\label{sec:From-Green-to-Riemann}

In term of the Green function for a simply-connected domain
one may  write down the Riemann map as
$$
f(z)=e^{G+iG^{\ast}}\,,
$$
where $G^{\ast}$ is the conjugate potential (satisfying the
Cauchy-Riemann equations). [This is basically the second proof
given by Riemann in 1857 \cite{Riemann_1857-DP}, and see also e.g.
Picard 1915 \cite{Picard_1915}.] That $f$ is a circle map follows
from $G\le 0$ with vanishing precisely on the boundary, and the
fact that $G^{\ast}$ is single-valued since the domain is
simply-connected. Details are supplied during the next Steps,
where we examine the more delicate multiply-connected domains or
even general compact bordered Riemann surfaces.

\medskip
{\it Step 3 (Memento about the conjugate potential)} The conjugate
$G^{\ast}$ potential is defined by the desideratum that
$G+iG^{\ast}$ is holomorphic, i.e. ${\Bbb C}$-linearizable in the
small. This gives the Cauchy-Riemann equations
$$
\frac{\partial G}{ \partial x}=\frac{\partial G^{\ast}}{
\partial y}, \qquad
\frac{\partial G^{\ast}}{ \partial x}=-\frac{\partial G}{
\partial y}.
$$
Writing formally $G^\ast$ as the integral of its differential,
gives
\begin{align}
G^{\ast}= \int dG^{\ast}=\int (\frac{\partial G^{\ast}}{
\partial x} dx+\frac{\partial G^{\ast}}{ \partial y}dy)
=\int (-\frac{\partial G}{
\partial y} dx+\frac{\partial G}{ \partial x}dy),
\end{align}
whose integrand (a $1$-form) coincides actually with the $dG$
twisted by multiplication by $i$ on the tangent bundle. Therefore
$dG^{\ast}$ is a genuine $1$-form canonically attached to the
function $G$. ({\it Warning}.---The symbol $G^\ast$ (taken alone)
as no intrinsic meaning at least as a single-valued function
unless $dG^{\ast}$ is period free.)

\subsection{The Green's function}

But what is the Green's function at all about? It is a sort of
logarithmic potential attached to an electric charge placed at
$t$. It is easier to define in the case of a plane domain bounded
by smooth curves. The case of ultimate interest (compact bordered
Riemann surfaces) will be discussed later. Given a domain
$B\subset {\Bbb C}$ (smoothly bounded) marked at an (interior)
point $t$ one considers the function $\log \vert z-t\vert$ which
induces (by restriction) a charge (temperature)  on the boundary
$\partial B=C$ and one solves the Dirichlet(=first boundary-value)
problem for this data. It results an (everywhere regular) harmonic
function $u=u(z,t)$, which subtracted from the original
logarithmic potential gives the (so-called) {\it Green's function
with pole at $t$}
\begin{equation}
( \log\vert z-t \vert )- u(z)=:G(z,t). \label{Green:eq}
\end{equation}
By construction, it vanishes along the contour $\partial B$
and possesses a logarithmic singularity near the point $t$.
%
This is a canonical function attached to the sole data of $B$ and
a certain interior point $t$. Note that $G(z,t)$ tends to
$-\infty$ as $z$ approaches $t$, so one may think of the Green's
function as a black hole centered at $t$ with a  vertiginous sink
plunging into deep darkness.

One can interpret this Green's function as some electric potential
(Galvanic current) on a conducting plate. If one prefers a
biological metaphor one can visualize $G(z,t)$ as the
proliferation of bacteroides originating from $t$ while expanding
through the medium $B$ driven by
an
apparent global knowledge of the shape of the universe. To be more
concrete, the expansion is more rapid where more free resources
are available. In particular all bacteria reach synchronously
 the boundary having consumed all
resources of the nutritive substratum in what looks to be the
most
equitable way. Compare the pictures in the simply-connected case
(Fig.\,\ref{Green:fig}a) and then for a multi-connected region
(Fig.\,\ref{Green:fig}b). Trying to imagine the same proliferation
occurring on a bordered surface realized say in Euclidean
$3$-space we get something like Fig.\,\ref{Green:fig}d.

\begin{figure}[h]
\penalty0 \centering
    \epsfig{figure=Green.eps,width=122mm}
    \vskip-5pt\penalty0

\caption{\label{Green:fig} The levels of Green's functions of two
planar domains with pole at $t$, and an attempt to draw Green on a
surface of genus $2$ (an aggressive bulldog?)}
\end{figure}

Now it is clear that the above formula
$f(z)=f(z,t)=e^{G(z,t)+iG^{\ast}(z,t)}$ supplies the Riemann map
with $f(t)=0$ and $f(z)\in S^1$ (unit circle) whenever $z$ lies on
the boundary, where $G$ vanishes. Of course the map is only
defined up to rotation, coming from an arbitrary additive constant
in $G^{\ast}$. [Compare for instance Riemann 1857
\cite{Riemann_1857-DP}, Picard 1915 \cite{Picard_1915}, etc.]

If one tries to adapt this
proof  to multi-connected domains one meets the notorious
difficulty that the conjugate potential $G^{\ast}$ is not
single-valued, a priori. So the
efforts focus on eliminating the periods of its differential
$dG^{\ast}$ by choosing appropriately some accessory parameters.
[This universally known device goes back at least to Riemann 1857
\cite[p.\,122]{Riemann_1857} Schottky 1877 \cite{Schottky_1877},
see also Picard 1913 \cite{Picard_1913}, Koebe 1922, Julia 1932
\cite{Julia_1932}, Grunsky 1937 \cite{Grunsky_1937}, etc.]

Using this idea we may concoct a circle map $B\to
\overline{\Delta}$. [cf. Grunsky 1937 \cite{Grunsky_1937} or
Grunsky 1978 \cite{Grunsky_1978} and also Ahlfors]. The natural
trick is probably to take several
poles $t_i$ (say $d$ many). Those will ultimately become the
zeroes of the circle map we are looking for as $e^{-\infty}=0$.
One now form the combination of the corresponding Green's
functions
$$
G(z):=\textstyle\sum_{i=1}^d \lambda_i G(z,t_i) \quad
(\lambda_i\in {\Bbb R}).
$$
This gives a (finite) constellation of black holes
scattered through the domain $B$ and we shall try to choose the
constants $\lambda_i$ so that $dG^\ast$ has no period. Since the
combination $G$ vanishes on the contour $\partial B$ (being a
superposition of Green's functions) the allied function $f(z)=f(z;
t_i, \lambda_i):=e^{G(z)+iG^{\ast}(z)}$ will map $\partial B$ onto
$S^1$. To arrange it as a circle map $f\colon B \to
\overline{\Delta} $ requires the basic remarks of the next
section, plus the more delicate issue of being able to choose
positive $\lambda_i>0$.

\subsection{Quasi-negativity of Green}

The following property of the Green's function is basic, yet
important.

\begin{lemma} Each Green function $G_t(z):=G(z,t)$ is
quasi-negative (i.e. $\le 0$ throughout the domain and
strictly $<0$ in its interior).
\end{lemma}

\begin{proof} From its definition \eqref{Green:eq} it is clear
that $G_t(z)\to -\infty$ as $t$ approaches the pole $t$. Thus
choosing a very large negative (real) constant $C<0$ the
corresponding level line $L_C$ of Green $G_t^{-1}(C)$ will be a
nearly circular (Jordan) curve enclosing the pole $t$ in its
interior. Further it looks evident that for $C<0$ large enough (in
absolute value) this Jordan curve bounds a (topological) disc in
the domain. (One could uses the general Schoenflies theorem
requiring just to check that $L_C$ is null-homotopic in the domain
$D$.) Next it is  intuitive (but need to be arithmetized) that
within this sufficiently small disc-shaped domain (i.e. the inside
of $L_C$ for $C<0$ sufficiently large) the Green function $G_t$ is
negative (indeed $\le C$).

Cutting away from the domain $D$ the interior of $L_C$ we obtain
an excised domain $D^{\ast}$ with one more contour. On this new
domain, the Green's function $G_t$ solves  Dirichlet (first
boundary-value) problem for the data $0$ on all contours but $C<0$
on the newly created contour $L_C$. We now conclude via the next
lemma.
\end{proof}

\begin{lemma} \label{Depressive:lem} {\rm (Depressiveness of Dirichlet,
or rather the allied harmonic functions)} Let $F$ be a compact
bordered Riemann surface. If the (continuous) boundary data
function $f\colon \partial F \to {\Bbb R}_{\le 0}$ is
non-positive, then so is its Dirichlet solution $u:=u(f)$, i.e.
$u\le 0$ throughout~$F$.
\end{lemma}

\begin{proof} If not then $u(z_0)>0$ (positive) at some
interior point $z_0$ of the surface $F$. By compactness $u$
achieves its maximum, which is positive. Since $f\le 0$ the
latter would not be achieved on the boundary violating the
maximum principle (compare the next lemma).
\end{proof}

\begin{lemma} {\rm (Maximum principle)}
Any harmonic function $u$ on a compact bordered surface $F$
achieves its maximum on the boundary $\partial F$. In fact, if
the maximum is achieved at some interior point then the
function $u$ is constant.
\end{lemma}

\begin{proof} Assume $z_0$ to be an interior point realizing
the maximum $M$ of the harmonic function $u$ defined on $F$. We
trace a little (metric) circle about $z_0$ of sufficiently small
radius as to lye entirely inside $F$ (together with its interior
disc $D$). Harmonicity may be characterized via the {\it
mean-value property} (Gauss, it seems):
\begin{equation}
\int \int_D u(z) d\omega = area(D) \cdot u(z_0).
\label{mean-value-prop:eq}
\end{equation}
As $u(z)\le M$, we get $M\cdot area (D)\ge \int \int_D u(z)
d\omega = area(D) \cdot u(z_0)$. Since $M= u(z_0)$, both extreme
members coincide and so does the last inequality.
This forces constancy on
the little disc $D$ ($u$ being continuous).

It follows by `propagation' that $u$ is globally constant.
(Alternatively use general topology: the set of points where $u$
achieves its maximum is both nonempty (compactness), closed and
open.) Indeed choosing a path from $z_0$ to any point $z\in F$
covered  by a chain of little discs $D_1, \dots, D_k$, each $D_i$
centered on the border of the previous one $D_{i-1}$, one argues
that two successive discs have enough overlap to ensure constancy
over the next disc.
\end{proof}

{\small {\it Micro-Warning} [11.08.12] {\rm There is a Garabedian
paper 1951 (A PDE..., p.\,486) were it is asserted that the
Green's function of a convex clamped plate need not be of one
sign; but of course this is not relevant to our matter were we use
the usual the Laplacian $\Delta$ and not the bi-Laplacian
$\Delta^2$ corresponding to clamped plated, instead of vibrating
membranes. This is the seminal work of Garabedian (but others were
also involved) were the famous Hadamard conjecture on the
bi-Laplacian was disproved.}

}

\subsection{Killing the periods}

The previous section  ensures that any superposition of Green's
functions $G:=\sum_i \lambda_i G(z,t_i)$ will be likewise
quasi-negative provided all $\lambda_i$ are positive. In this
circumstance the function $f=e^{G+iG^{\ast}}=e^G \cdot e^{i
G^{\ast}}$ (whose modulus is $e^G$) is a unit-circle map ($\vert f
\vert \le 1$), because the real exponential takes nonpositive
values $(-\infty, 0]$ to $(0,1]$. It is consistent by continuity
to
send the $t_i$ on $0$.

If $r$ is the connectivity of the domain $B$ (number of its
contours) then there are homologically $r-1$ non-trivial loops
$\gamma_1, \dots \gamma_{r-1}$ running around the $r-1$ holes
in our domain (cf. Fig.\,\ref{Green:fig}c illustrating the
case $r=3$).
We consider the linear period mapping
\begin{align}\label{period-mapping:eq}
{\Bbb R}^{d} &\longrightarrow {\Bbb R}^{r-1} \cr
(\lambda_1, \dots \lambda_d) &\mapsto
(\textstyle\int_{\gamma_1} dG^{\ast}, \dots,
\int_{\gamma_{r-1}} dG^{\ast})
\end{align}
%
%
By linear algebra
if $d$ is large enough (precisely already for $d=r$) we have
enough free constants so as to find non-trivial $\lambda_i$
extincting all periods. [Heuristically the electric poles of the
multi-battery in the electrolytic tank (nomenclature as in e.g.
Courant 1950/52 (Conformal book)) are affected by suitable charges
so as to generate an ``ideal'' potential with single-valued
conjugate.] Exponentiating gives $f=e^{G+iG^{\ast}}$ a circle map
with $d=r$ zeroes, provided one is able to ensure all $\lambda_i
>0$.

Without taking care of this last proviso, one  may reach too
hastily the impression that we have complete freedom in
prescribing the location of the $d=r$ poles (of the Green's
functions, which convert ultimately to zeroes of the related
circle map).
The linear-algebra argument gives only a real-line inside the
kernel of the (linear) period-mapping \eqref{period-mapping:eq},
but a priori this line could miss the ``octant''  ${\Bbb
R}_{>0}^{d}$ consisting of totally positive coordinates. In fact
upon letting vanish some of the $\lambda_i$ what is only required
is a non-trivial penetration of this line $\ell$ into the closed
octant $\overline{O}={\Bbb R}_{\ge 0}^{d}$, i.e. the intersection
$\ell \cap \overline{O}$ should not reduce to the origin. A true
penetration of this line in the interior of $O$, or a degenerate
one where the line meet along one of its face would be enough to
complete the existence-proof. The latter case amounts to extinct
some Green's ``batteries'' by assigning a vanishing coefficient
$\lambda_i=0$. The net effect would be degree lowering of the
circle map $f$. Beware, that for planar domains (which correspond
to Harnack-maximal Schottky doubles) no such lowering of the
degree is possible for simple topological reasons ($r\le \gamma$).
However the  described theoretical
eventuality  may well happen in the non-planar case to be soon
discussed. Understanding how and why to arrange degenerate
penetrations could well offer a strategy toward improving Ahlfors
$r+2p$ bound.

\subsection{Extra difficulties in the surface case}

It is obvious that the above method via Green's functions adapts
to bordered Riemann surface $F=F_{r,p}$ of (positive) genus $p$
with $r$ contours (Rand). Remember however that at this stage we
did not offered a complete treatment of the planar case ($p=0$).

First note a conceptual difficulty regarding Green's function,
which, in the plane case of a domain $B\subset {\Bbb C}$, is
constructed via $\log\vert z-t \vert$ appealing to a global
coordinate system.
In the abstract bordered setting, there is no such ambient medium.
One could try to work with a (conformal) Riemann metric and the
allied logarithmic distribution $ \log \varrho (z,t), $ where
$\varrho$ is the intrinsic distance (defined as usual as the
infimum of lengths of rectifiable pathes joining  two given
points). Note however that this construction specialized to the
domain case does not duplicate the former, since the intrinsic
distance $\varrho(z,t)$ does not coincide with the extrinsic one
$\vert z-t\vert$, unless the domain $B\subset {\Bbb C}$ is
starlike about $t$.

Bypassing this difficulty [which will be resolved later], we first
 note that each handle creates two $1$-cycles yielding a total of
$(r-1)+2p$ many essential loops (compare
Fig.\,\ref{Green:fig}\,e). Thus introducing $d:=r+2p$ poles $t_i$
we dispose of enough free parameters  to arrange (via linear
algebra) the vanishing of all periods of the conjugate
differential  $dG^{\ast}$ of the  potential
$G=\sum_{i=1}^{d}\lambda_i G_{t_i}$. This explains quite clearly
why Ahlfors discovered (about 1948) the upper-bound $r+2p$ for the
degree of a circle map. Of course there is still the subtlety of
explaining why it is possible to choose all $\lambda_i>0$ at least
for a clever choice of the poles $t_i$.

{\small All this is probably when suitably interpreted the
quintessence of the Ahlfors mapping (of degree $r+2p$). Again the
writer does not mask his happiness after having understood this
point (as late as the 04.08.12).
Now it is evident to reconstruct (even if somewhat fictionally)
what must have happened in Ahlfors' brain (at least as early as
1948, and presumably much earlier, yet no record in print). With
this piece of information and, on the other hand, being well-aware
of the modern purely function-theoretic proofs of RMT (\`a la
Koebe-Carath\'eodory, Fej\'er-Riesz 1922 (published by Rad\'o
1923), Carath\'eodory 1928 and Ostrowski 1929) it must have seemed
highly desirable (or trendy) to reinterpret the above (somewhat
heuristic but fruitful potential theory) in terms of a
function-theoretic extremal problem. This leads e.g. to the
problem we discussed at length of maximizing either the modulus of
the derivative at some inner point $t=a$, or to maximize the
distance of two points $a,b$ where the first maps to $0$ and the
second is repulsed at maximum distance from the origin. In both
case the competing  functions are analytic and bounded-by-one in
modulus $\vert f \vert \le 1$.
So we get the Ahlfors function $f_a$ or $f_{a,b}$. It seems
obvious that all those Ahlfors functions are included in the above
trick \`a la Green-Riemann (GR), and  thus subsumed to an
electrolytic interpretation. Yet the exact dependance and location
of the corresponding logarithmic poles of Green's $G$ (becoming
the zero of Riemann's $f$, after exponentiation) must be a
transcendentally sublime  business. Also the corresponding degree
of the Ahlfors function is another mystery.

}

It is  conceivable that less than the $r+2p$ generically required
poles suffices in case the linear period mapping ${\Bbb R}^d \to
{\Bbb R}^{(r-1)+2p}$  along  fundamental loops has a degenerate
image permitting to economize
some poles $t_i$. The task is reduced to find the lowest $d$ such
that the kernel of the period map is non-trivial and contains a
non-zero element
 all of whose coordinates are $\ge 0$. Remember, that Gabard 2006
\cite{Gabard_2006} showed---using another method, based on a
topological argument of irrigation
(Riemann-Betti-Jordan-Poincar\'e's homologies, and Brouwer's
degree plus some basic Pontrjagin theory in the Jacobian torus as
a very special commutative Lie group---that there is a circle map
of degree $\le r+p$ (i.e. with one unit economized for each
handle). Assuming that any circle map is allied to a Green-Riemann
map there would be a fewer number namely $d\le r+p$ of batteries
required to generate this mapping. Of course, the first part of
the assertion looks evident: given a degree $d$ circle map $f$
with zeroes at $t_i$, then $\log \vert f(z) \vert$ coincides with
$\sum_{i=1}^d G(z,t_i)$. This is Ahlfors formula following from
the fact that both functions vanishes on the border and have the
same singularities.

{\it Philosophy.} [08.08.12] Modulo elusive details, it is fair to
resume the situation by saying that the Ahlfors circle maps (if
not all existence theorems of function theory) derives form the
Dirichlet principle (or the allied Green's functions). [This was
of course best incarnated by Riemann, 1851 and 1857, where in
bonus  the whole algebraic geometry of curves was subsumed to this
principle!] Conversely one could hope that the Ahlfors function
could be used to lift the Dirichlet solubility on the disc (via
Poisson integral formula) to an arbitrary bordered surface.
However it seems obvious that there is no way to descend the
boundary function to the disc since the Ahlfors branched covering
is multi-valent. We arrive at the conclusion that the true
mushroom is the Dirichlet principle, while Ahlfors function being
just one tentacle of the mushroom. Of course, the only paradigm
susceptible of competing with Dirichlet are the function-theoretic
extremal problems \`a la
Koebe-Carath\'eodory-Fej\'er-Riesz-Bieberbach-Ostrowski, etc. For
plane domains the Kreisnormierung (instead of the Ahlfors map) may
be used as normal domains where the Dirichlet problem is easier to
solve. This is akin to Poisson's formula for the round disc case
of Dirichlet, and quite implicit in Riemann's Nachlass 1857
\cite{Riemann_1857_Nachlass} (cf. also Bieberbach 1925
\cite{Bieberbach_1925}). A similar reduction of Dirichlet for
bordered surfaces occurs is also likely on the ground of Klein's
R\"uckkehrschnitttheorem (cf. Section
\ref{sec:Ruckkehrschnittthm}), supposed to be an extension of the
Kreisnormierung.

Regarding the detailed execution of the removal of the period as
to construct an Ahlfors-type mapping one should compare also the
paper of Heins 1950 \cite{Heins_1950}, Kuramochi 1952
\cite{Kuramochi_1952} and (albeit confined to planar domains) the
paper by D. Khavinson 1984 \cite{Khavinson-Dimitri_1984}, whose
argument is considered by its author akin to the arguments of
Grunsky.

\subsection{The Green's function of a compact bordered Riemann
surface=CBRS}

[14.08.12] This section examines the  issue that the Green's
function $G(z,t)$ with pole at $t$ is a canonically defined
function in the generality of a CBRS. This is super-classical, cf.
e.g., the treatises Ahlfors-Sario 1960 \cite{Ahlfors-Sario_1960}
or Schiffer-Spencer 1954 \cite{Schiffer-Spencer_1954}. It is to be
expected to find older treatments by Riemann, Schwarz, Klein,
Koebe, etc. Several accounts by Nevanlinna proceed via Schwarz's
alternating method, a viewpoint which looked most convenient to
adhere with.

As already noticed, the case of a plane domain $B\subset {\Bbb C}$
(bounded by smooth curves) it is easy to define Green's function
$G(z,t)$ via the (logarithmic) potential $\log \vert z-t \vert$
from which we subtract the Dirichlet solution matching the
logarithmic potential restricted to the boundary $\partial B$.
Alas, for a CBRS $F$ one lacks an ambient space like ${\Bbb C}$
permitting an analogous construction.

Of course, $\log \vert z-t \vert$ bears some significance only
locally within a uniformizer chart about $t$. Taking another local
chart, one may argue that in the small the expression will mutate
into $\log \vert \alpha (z-t) \vert$ for some $\alpha\in{\Bbb
C}^{\ast}$ incarnating the derivative of the transition between
the two charts. Thus the log-potential w.r.t. the new chart is
$\log\vert \alpha \vert + \log \vert z-t \vert$, hence equal to
the old one modulo an additive constant. Presumably some
philosophical argument can corroborate the vague feeling that the
asymptotic of the logarithmic pole is unaffected by such additive
constant. [Added in proof: compare Pfluger 1957 \cite[p.\,110,
28.3]{Pfluger_1957} for an accurate formulation,  or Farkas-Kra
1980/1992 \cite[p.\,182, Remark]{Farkas-Kra_1980/1992}.] It seems
then meaningful to set:

\begin{defn} The Green's function of a CBRS $F$ with pole at
 $t$ (an interior point of $F$) is the unique harmonic function on $F$
save $t$ with singularity $\log\vert z-t \vert$ near $t$ which
vanishes continuously on the boundary $\partial F$.
\end{defn}

Compare  (modulo a different sign convention) Ahlfors-Sario 1960
\cite[p.\,158, 4B]{Ahlfors-Sario_1960}. Uniqueness is considered
as evident there. Indeed, a chart change affect the logarithmic
potential by an additive constant and harmonic functions are quite
rigid (being determined by their values on any open disc). Hence
knowledge of the function on any punctured chart about $t$ via
$\log \vert z-t \vert$ determines it uniquely.
The delicate point is  existence.  Choose around $t$ a nice
analytic Jordan curve $J$ and via RMT construct a holomorphic
chart taking $D$ (the ``sealed'' interior of $J$, i.e. $J$
included) to the unit disc $\overline{\Delta}$. Consider
$\log\vert z\vert$ in the unit circle and transplant to $D\subset
F$ and then after adding an additive constant we try to solve a
Dirichlet-Neumann problem on $F-D$ piecing together smoothly the
logarithmic piece with the Dirichlet-Neumann solution. In this
procedure  the Green's function looks highly non-unique depending
on  the ``ovaloidness'' of the Jordan curve $J$ chosen. In fact
$J$ cannot be chosen at will but must somehow be a level-line of
Green (still undefined). Infinitesimally $J$ should be a perfect
circle, and this is perhaps the key to put the naive pasting
argument on a sound basis via a convergence procedure.
(Infinitesimal circles are well-defined on Riemann surfaces via
the conformal structure.)
Existence and uniqueness look then plausible, but involve a
considerable sophistication over the plane-case where the Green's
function reduced  straightforwardly to the Dirichlet problem.

Let us paraphrase the above more formally.  Take any chart
$\varphi\colon U \to \Delta$ about the ``pole'' point $t$ (sending
$t$ to the origin $0\in {\Bbb C}$), write down $\log\vert z\vert$
in that chart and shrink gradually attention to the (round) disc
$\Delta_{\varepsilon}$ of radius $\varepsilon$. Let
$D_{\varepsilon}$ be $\varphi^{-1}(\Delta_{\varepsilon})$. For
each (positive) value of $\varepsilon$ one can solve the Dirichlet
problem in $F-{\rm int} D_{\varepsilon}$ with boundary value $0$
on $\partial F$ and $\log \varepsilon$ on $\partial
D_{\varepsilon}$. Denote by $u_\varepsilon $ the corresponding
solution. By construction $u_\varepsilon$ pasts continuously with
the $\varphi$-pullback of the log-potential (i.e. $(\log \vert
z\vert) \circ \varphi$). Of course this glued function is a
Frankenstein creature lacking a smooth juncture. For instance, if
$\varepsilon=1$ then $u_{\varepsilon}$ is identically zero,
whereas in $D_1$ we have the logarithmic ``trumpet'' with
derivative $1$ along the normal direction. However as
$\varepsilon$ decreases from $1$ to $0$, $u_\varepsilon$ becomes
$\le 0$ (having prescribed the negative value $\log \varepsilon$
on $\partial D_{\varepsilon}$) and the dependence of
$u_{\varepsilon}$ is perhaps monotonic. So it seems arguable
(Harnack?)  that while $\varepsilon \to 0$ (say via dyadic numbers
$\varepsilon_n=1/2^n$) the $u_n$ converges to a harmonic function
on $F-t$ which is the desired Green's function $G(z,t)$. It seems
evident (since $\partial D_n$ becomes more and more circular in
$F$ as $n$ grows to infinity) that the limit is harmonic and
independent of the gadgets used along the way (chart $\varphi$,
dyadic sequence $\varepsilon$).

This vaguely explains existence and uniqueness can maybe be
derived by a similar trick (combined with a ``leapfrog''
argument). Try to locate a reference along this naive line: maybe
Schwarz?, Klein? Koebe? Weyl? Pfluger? and otherwise try
Ahlfors-Sario \cite{Ahlfors-Sario_1960}, Sario-Oikawa
\cite{Sario-Oikawa_1969}. (Sometimes Sario's formalism of the
normal/principal operator  is a bit awkward to digest.) For
treatments of the Green's function on a CBRS cf. Schiffer-Spencer
1954 \cite[p.\,33, and 93--94]{Schiffer-Spencer_1954}. See also
Sario-Oikawa 1969 \cite[p.\,49--50]{Sario-Oikawa_1969}. We
summarize the discussion by the

\begin{theorem}
Given a CBRS $F$ and an interior point $t$, there is a
uniquely defined Green's function $G(z,t)$ with pole $t$ which
is characterized by the following conditions: it is harmonic
on $F-t$, vanishes (continuously) on the boundary $\partial F$
and it has the prescribed singularity $\log \vert z-t \vert$
near $t$.
\end{theorem}

\begin{proof} For complete details, compare several sources:

$\bullet$ first Ahlfors-Sario 1960
\cite[p.\,158]{Ahlfors-Sario_1960} and Sario-Oikawa 1969
\cite[p.\,50]{Sario-Oikawa_1969} (both via Perron's method,
and Sario's formalism of the normal operator).

$\bullet$ Then also Pfluger 1957 \cite[p.\,110, end of \S
28.2, as well as p.\,110, \S28.3 and last 3 lines of
p.\,111]{Pfluger_1957}

$\bullet$ Schiffer-Spencer 1954 \cite[\S4.2,
p.\,93--94]{Schiffer-Spencer_1954}

$\bullet$ Nevanlinna 1953
\cite{Nevanlinna_1953-Uniformisierung} via Schwarz's
alternating method (SAM). We detail this argument in the next
section.
\end{proof}

\subsection{Schwarz's alternating method to construct the
Green's function of a compact bordered surface (Nevanlinna's
account)}

[15.08.12] As promised, in this section we attempt to understand
Nevanlinna's exposition of the existence of Green's functions on
compact bordered surfaces. All pagination given refers to
Nevanlinna 1953 \cite{Nevanlinna_1953-Uniformisierung}, the book
``Uniformisierung''. Nevanlinna follows Schwarz's alternating
method (SAM) quite closely. The argument is a bit tedious but
quite elementary. It uses merely the following facts on a bordered
surface $F$:

(1) if $f\le g$ on $\partial F$ then the associated Dirichlet
solution $u(f)\le u(g)$ (compare Lemma~\ref{Depressive:lem}).

(2) maybe Harnack's theorem is required?

In Nevanlinna's book the relevant information is scattered at two
places (at least) so we attempted to compactify the presentation
for our own understanding.

First Nevanlinna introduces a concept of ``Kreisbereich''. Alas
the jargon is not very fortunate being already
consecrated by Koebe in a different context, so let us speak
rather of a ``celluloid'' (or a ``Kreisgebilde''). This is
[cf.\,p.\,142] a connected finite union of (closed) discs in a
Riemann surface (whose images by a (parametric) chart are round
discs in ${\Bbb C}$). On each such disc the first boundary value
problem (abridged DP=Dirichlet problem) is soluble via Poisson's
formula. Assuming the contours of each pair of discs to have
finite intersection, SAM enables one to solve DP on the union,
hence on any
 celluloid.

So for instance it is clear that any CBRS is a celluloid. (A
formal proof certainly requires Rad\'o's triangulation theorem
1925 \cite{Rado_1925}.) To absorb the boundary in one stroke one
could add annular regions where the  DP is also soluble by an
explicit recipe, sometimes ascribed to Villat 1912
\cite{Villat_1912}.

The Green's function will be obtained by specializing the
following technical lemma [cf.\,p.\,148] [due to Schwarz and
probably related to what Koebe's calls the ``g\"urtelf\"ormige
Verschmelzung''($\approx$belt-shaped fusion)]. Intuitively, the
lemma amounts to construct a harmonic function $u$ with prescribed
boundary values and with prescribed singularity $u_0$ near a point
$t$, or rather on a ring enclosing the pole $t$.
%
(At first, it is not perfectly transparent how to deduce the
Green's function from the lemma, but we shall try elucidate this
issue later.)

\begin{lemma} Let $F$ be a compact bordered Riemann surface
and $t\in {\rm int} F$ an interior point. Let $U$ be a
neighborhood of $t$ mapped to the unit-disc $\Delta=\{ \vert z
\vert <1\}$ via a chart. In $U$, let $K$ be the ring corresponding
to $r_1\le \vert z \vert \le r<1$. {\rm [It seems that $r_1=0$ is
permissible and needed for the application to the Green's
function.]} Let further $X$ be a celluloid containing the external
contour  of the ring $K$ $(\vert z \vert=r)$ in its interior, as
well as the boundary $\partial F$ but missing the small disc in
$U$ corresponding to $\vert z\vert \le r_1$. Set finally $A:=X
\cup K$ {\rm (compare  Fig.\,\ref{Nevanl:fig}}).

Then given  $u_0 \in H (K)$ (harmonic on the ring $K$) and
$f\colon \partial F \to {\Bbb R}$ continuous, there is a
unique $u \in H(A)$ such that $u_{\vert \partial F}= f$ and
with $u-u_0$ extending harmonically to $B$, the disc
corresponding to $\vert z\vert\le r$.
\end{lemma}

\begin{figure}[h]
\centering
    \epsfig{figure=Nevanl.eps,width=92mm}
\vskip-5pt\penalty0
  \caption{\label{Nevanl:fig}
  Schwarz's alternating method (Nevanlinna's implementation) to
  construct the Green's function of a compact bordered
  surface
  with pole at $t$}
\vskip-5pt\penalty0
\end{figure}

Before detailing the proof, let us see how this helps defining
$G(z,t)$ the Green's function. [The difficulty is just so trivial
to be not completely explicit in Nevanlinna [p.\,198--199, \S2,
Art.\,6.4].] First we impose $f\equiv 0$. Then we choose the
singularity function $u_0=\log\vert z\vert$ which has to be
defined on $K$, hence we shrink $r_1$ to $0$ via a sequence of
dyadic radii $r_n=1/2^n$. On applying the lemma we get a sequence
of solution $u=u_n$ defined on $A_n$ a sequence of expanding
subsurfaces (the outsides of the shrinking discs $\vert z\vert<
1/2^n$). Now observe that $u_n$ for a large $n$ solves the problem
of the lemma for all smaller values of $n$: just take the
restriction (and use uniqueness). Consequently all the $u_n$ form
a telescopic system of functions (each restricting to all its
predecessors) defined on larger and larger compact subregions
$A_n$ ultimately expanding to the punctured surface $F-t$. The
very constant (indeed completely monotone) limit of those $u_n$
gives the desired Green's function $G(z,t)$.

It is harmonic on $F$ save $t$, vanishes on the boundary and
$G(z,t)-\log\vert z \vert$ extends harmonically through $t$ (on a
little neighborhood). It remains to check that those 3 properties
defines $G(z,t)$ unambiguously. This is again the same sort of
argument. Assume there were two Green's solutions $G_1,G_2$, then
$G_i-u_0=:h_i$ harmonic on some neighborhoods $V_i$ of $t$. So
$G_1-G_2=(h_1+u_0)-(h_2+u_0)=h_1-h_2$ which is harmonic on the
intersection $V_1 \cap V_2$. Hence the difference $G_1-G_2$ is
harmonic throughout $F$, but with vanishing boundary value on
$\partial F$. Consequently it must be identically zero (by the
uniqueness part of Dirichlet) which follows from the maximum
principle.

\begin{proof} This is a matter of implementing Schwarz's
alternating method [see p.\,148--150]
and we follow exactly Nevanlinna's text (annotating our copy by
the symbol $\bigstar$ to indicate the sole cosmetic difference).

$\bullet$ {\it Uniqueness} Assuming the existence of two functions
$u_1, u_2$ solving the problem, their difference $u_1-u_2$ will be
harmonic on $A$, and 0 on $\gamma:=\partial F$. But each
difference $u_i-u_0=:h_i\in H(B)$ extends harmonically across $B$
($i=1,2$). Hence on $B$, $u_1-u_2=(h_1+u_0)-(h_2+u_0)=h_1-h_2 \in
H(B)$, and therefore $u_1-u_2 \in H(A\cup B)$, and $A\cup B$ is
all of $F$. It follows (Dirichlet's uniqueness) that $u_1-u_2$
vanishes identically.

$\bullet$ {\it Optional remark.} It is clear that the case
$f\equiv 0$ is typical, since the general case just requires
adding the Dirichlet solution for the data $f$. [This explains why
I had the impression to find many misprints!]

$\bullet$ {\it Existence (after Schwarz)} First it is observed
that DP is solvable on both $A$ and $B$ ($B$ is just a ball  and
$A$ is a celluloid, yet of the general type involving a ring). Of
course $A$ is also a compact bordered surface and therefore one is
ensured of Dirichlet solvability, thereby bypassing the concept of
a celluloid, and accordingly one can shorten slightly the
statement of Nevanlinna's lemma, with the direct bonus that one
can make abstraction of all the little discs drawn on the picture.
[CAUTION: here it is perhaps NOT permissible to take $r_1=0$]

We denote by $\alpha$ and $\beta$ the internal resp. external
contour of $K$, and let $\gamma:=\partial F$. Set first
$v_0\equiv 0$. We define inductively sequences $u_n\in H(A)$
and $v_n\in H(B)$ by their boundary values ($n\ge 1$)
\begin{equation}\label{SAM:eq}
u_n=\begin{cases} v_{n-1}+u_0 \quad &\textrm{on } \alpha, \cr
0 [\bigstar or f] \quad &\textrm{on } \gamma,
\end{cases}
\end{equation}
and
$$
v_n=u_n-u_0 \quad \textrm{on } \beta.
$$
[This Ansatz comes a bit out of the blue, but notice that passing
to the limit both definitions leads to the identity $u-u_0=v$
holding on $\alpha \cup \beta$ which is the full contour of the
ring $K$, so that anticipating harmonicity this will hold
throughout $K$, and $v$ will afford the required extension of
$u-u_0$ (only defined on $A\cap K=K$) to the disc $B$ containing
the ring $K$. Of course, it is also crucial to notice that both
sequences $u_n,v_n$ are ``interlocked'' or ``leapfrogged''
requiring an alternating progression of one term to go one step
further with the other.]

The successive differences are given by
\begin{equation}\label{success-diff:eq}
u_{n+1}-u_n=\begin{cases} v_{n}-v_{n-1} \quad &\textrm{on }
\alpha \cr 0  \quad &\textrm{on } \gamma
\end{cases}
\end{equation}
and
$$
v_{n+1}-v_n=u_{n+1}-u_{n} \quad \textrm{on } \beta.
$$
Let us write
$$
M_n:=\max_{\beta} \vert u_{n}-u_{n-1}\vert=\max_{\beta} \vert
v_{n}-v_{n-1}\vert,
$$
then by the maximum- and minimum-principle $ \vert
v_{n}-v_{n-1}\vert\le M_n$ in $B$, and so in particular on
$\alpha$. Hence  by \eqref{success-diff:eq}, $ \vert
u_{n+1}-u_{n}\vert\le M_n$ on $\alpha$. Further, the
difference $u_{n+1}-u_{n}$ vanishes on $\gamma$ (cf.
\eqref{success-diff:eq}), and so it is bounded on the boundary
of $A$ (and therefore throughout $A$) by the potential $M_n
\cdot \omega$, where $\omega$ is the harmonic function
vanishing along $\gamma$ and equal to $1$ on $\alpha$. Hence
\begin{equation}\label{diff:eq}
\vert u_{n+1}-u_n \vert\le M_n \cdot \omega\quad  \textrm{in }
A.
\end{equation}
In the interior of $A$, one has $0<\omega < 1$. If $q$ is the
maximum of $\omega $ on $\beta$, then $0<q<1$. Further on
$\beta$ we have
$$
\vert u_{n+1}-u_n \vert\le q \cdot M_n,
$$
and also (by definition of $M_n$)
$$
M_{n+1}\le q \cdot M_n.
$$
By induction, it follows that
$$
M_{n+1}\le q^{n} \cdot M_{1},
$$
and recalling again the definition of $M_n$ we get (first on
$\beta$ and thus on $B$)
$$
\vert v_{n+1}-v_n \vert \le M_{n+1}\le q^{n} \cdot M_{1}.
$$
When particularized to $\alpha$, this implies in view of
\eqref{success-diff:eq}
$$
\vert u_{n+1}-u_n \vert \le q^{n-1} \cdot M_{1} \quad \textrm{in }
\alpha,
$$
and by the maximum principle this extends to $A$ (recall that
$\partial A= \alpha \cup \gamma$ and the function
$u_{n+1}-u_n$ vanishes on $\gamma$). Consequently, both series
$\sum_{n}(u_{n+1}-u_{n})$ and $\sum_{n}(v_{n+1}-v_{n})$
converges uniformly on $A$ resp. $B$.

The limiting functions $u$ and $v$ of $u_n$ resp. $v_n$ are
therefore harmonic on $A$ resp. $B$, and taking the limit in
the definition of $u_n$ (see \eqref{SAM:eq}) we see that $u$
vanishes on $\gamma$ [$\bigstar$ equals $f$ on $\gamma$].

We show finally that $u-u_0=v$ on $B$ [$\bigstar$ $K$
probably?]. Indeed, taking the limit in the first line of
\eqref{SAM:eq} gives $u=v+u_0$ on $\alpha$, and the definition
of $v_n$ pushed to its limit gives $v=u-u_0$ on $\beta$.
Therefore the same identity $u-u_0=v$ holds on both contours
of the ring $K$, and consequently its validity propagates
throughout $K$.

Finally, as $v$ is harmonic on $B$ we are happy to
conclude that $u$ fulfills all of our requirements: namely $u\in
H(A)$, $u=f$ on $\gamma=\partial F$ and $u-u_0$ defined on $A \cap
K=K$ coincide there with $v$ defined on the larger set $B\supset
K$, yielding the asserted harmonic extension.
\end{proof}

[NAIVE AND WRONG---see rather the argument given above]
Finally, [compare p.\,198--199] one obtains the Green's
function $G(z,t)$ by taking $u_0=\log \vert z \vert$, $f\equiv
0$ and $r_1=0$ [Caution: this point is not made explicit in
Nevanlinna]. For this choice of $r_1$, note that $A=F-t$. The
lemma supplies a unique $u\in H (F-t)$ such that $u_{\vert
\partial F}=0$ and so that $u-\log\vert z\vert=:h$ is harmonic
on $F$. The function $u$ is the desired Green's function
$G(z,t)$.

\section{Little green's men dreams (extraterrestrial
applications of Green's)}

The following three subsections are optional reading containing
more questions than answers. The reader interested primarily in
the Ahlfors map should preferably skip them.

\subsection{From Green to Gromov? (directly bypassing
Riemann and L\"owner)}\label{sec:Green-to-Gromov}

To mention once more a
deep frustration (the Gromov filling conjecture) it looks not
completely crazy to hope that a careful examination of the Green's
function and the allied isothermic coordinates could prompt a
solution of this problem. We tried quickly the [14.08.12] but
failed dramatically as usual
(along with circa 10 attempts of essentially the same vein).
Roughly the idea would be to look at the streamlines of Green and
its equipotentials, and remove every trajectory ending to the
(finitely many) critical points of Green while attempting to
estimate area via this (isothermic) parametrization. Of course,
Schwarz's inequality enters into the game but I only arrived  at
weak estimates like $\pi$ or $\pi/2$ (in place of $2\pi$!) upon
doing highly fallacious calculus.

\subsection{Schoenflies via Green?}\label{Schoenflies:sec}

A notorious topological paradigm is the so-called Schoenflies
theorem to the effect that a reasonably embedded sphere $S^{n-1}$
in ${\Bbb R}^n$ bounds a topological ball $B^n$. (There is a large
debate (cf.\,e.g.\,Siebenmann 2005 \cite{Siebenmann_2005}) about
who (and more broadly speaking which community) proved  first the
case $n=2$. In the topological-combinatorial realm there is a
contribution of Schoenflies reaching full maturity ca. 1906, and
somewhat earlier there is the contribution of Osgood which may
have reached full stability with Carath\'eodory 1912
\cite{Caratheodory_1912}. Of course the statement (for $n=2$ and
maybe even $n=3$) was largely anticipated heuristically by other
workers, e.g. Moebius 1863 \cite{Moebius_1863}. Schoenflies's
theorem was extended to higher dimensions by J.\,W. Alexander
($n=3$ ca. 1922), B. Mazur and M. Brown (all $n$ ca. 1960) for any
locally flat (e.g. smooth) hypersphere in ${\Bbb R}^n$. From Thom
or Smale's $h$-cobordism theorem (early 1960's) it is inferred
that the closed ball $B^n$ carries a unique smooth structure when
$n\neq 4$ (the case $n=4$ being still largely unsettled). It
follows that the interior of the smoothly embedded sphere is a
ball differentiably. Another unsolved problem of longstanding
is the truth of the same conclusion for $n=4$ (the so-called {\it
smooth Schoenflies} in dimension $4$, SS4,
see e.g. papers by Scharlemann).
Naive physical (or bacteriological, cf. Fig.\,\ref{Green:fig})
intuition about the Green's function makes hard to visualize why
there should be any anomaly for $n=4$, yet
nobody ever succeeded to prove or disprove SS4.
This belongs to the charming mysteries of low-dimensional
differential topology at the critical dimension $n=4$. One may
speculate about a naive approach to SS4 through the ca. 200 years
older potential theory (of Laplace, Poisson, Green, Gauss,
Dirichlet and Riemann's era). Alas, there is few records in print
of analysts feeling confident enough about the explorative
aptitudes of the Green's function (compare Fig.\,\ref{Green:fig})
to claim the required diffeomorphism with $B^4$. Of course in the
very small vicinity of the pole $t$ the levels of $G(z,t)$ (now
$z\in {\Bbb R}^n$)
look  alike round spheres, and by the synchronization principle
stating that each bacteria reaches the boundary at the same moment
it may look immediate how to write down the diffeomorphism. Can
somebody explain why this Green's strategy fails to establish SS4.
Less ambitiously can somebody reprove SSn (for $n \neq 4$) via the
Green's function. If yes with some little chance his/her proof
will possibly include the case $n=4$.

\subsection{Green, Schoenflies, Bergman and Lu Qi-Keng}

[06.08.12] As discussed in the previous section, a dream would be
to show SS4 (smooth Schoenflies conjecture) via the Green's
function in 4D-space ${\Bbb R}^4$. On reading an article by Boas
1996 (PAMS), where Suita-Yamada 1976 \cite{Suita-Yamada_1976} is
cited we see a potential connection between both problems.

The problem of Lu Qi-Keng asks for domains where the Bergman
kernel is zero-free (so-called Lu Qi-Keng=LQK-domains). Since
Schiffer 1946 \cite{Schiffer_1946}, there is an identity
connecting the Bergman kernel to the Green's function. It seems
that the zeros of Bergman corresponds to the critical points of
Green. Of course the latter is forced to have critical points as
soon as the topology is complicated (not a disc). Suita-Yamada's
result that the Bergman kernel necessarily exhibits zeroes for
membranes which are not discs looks nearly obvious. Hence
LQK-bordered surfaces are precisely those topologically equivalent
to the disc.

Now Boas in 1986 found a counterexample showing that no
topological characterization of LQK-domains holds in higher
dimensions: there exists in ${\Bbb C}^2$ a bounded, strongly
pseudoconvex, contractible domain with $C^{\infty}$ regular
boundary whose Bergman kernel does have zeroes. [Addendum
[18.09.12]: in fact upon reading Boas original paper (1986), Boas'
domain is diffeomorphic to the ball $B^4$.]

{\bf Optimistic scenario (Green implies Schoenflies)} It would
be interesting to know what the topology of Boas' hypersurface
$S=\partial \Omega$ is. In view of
Poincar\'e-Alexander-Lefschetz duality $S$ must be a homology
sphere, if I don't mistake. Now upon speculating that SS4 is
true (by naive geometric intuition), and even more that it is
provable via the streamlines of Green's function, and granting
a persistence of Schiffer's Green-Bergman identity (in the
realm of two complex variables), it may seem that Boas's
counterexample must have an ``exotic'' boundary (not
diffeomorphic to $S^3$). [Of course, not so in view of the
just given Addendum.]

{\bf Pessimistic scenario (Green does not implies Schoenflies).}
The other way around, assuming that Boas' boundary is the
3-sphere, there would be critical points of the Green's function
$G(z,t)$ and Boas's example may foil any naive attempt to reduces
SS4 to the streamlines of the Green's function. But even so maybe
the Green-Bergman identity of Schiffer is specific to one complex
variable, leaving some light hope that there is a
potential-theoretic proof of the differential-topology puzzle of
SS4.

So a bold conjecture (somewhat against Boas' philosophy that
there is no topological characterisation of LQK-domains) would
be that any domain in ${\Bbb C}^2$ bounded by  a smoothly
embedded 3-sphere is a LQK-domain (i.e. its Bergman function
is zero-free). [This is wrong in view of Boas 1984 (addendum
just mentioned)]. However it could be true that the Green's
functions $G(z,t)$ for any center $t$ located in the inside of
$\Sigma$ is critical point free, whereupon an elementary
integration of its gradient flow should establish a
diffeomorphism of the inside the spheroid with the ball $B^4$
with its usual differential structure. (Recall that it is yet
another puzzle of low-dimensional topology, whether the
$4$-ball has a unique smooth structure! All others balls
(maybe except the five-dimensional one) do enjoy uniqueness by
virtue of Smale's $h$-cobordism theorem.) Note that the
Bergman kernel is defined without reference to a basepoint
whereas Green's function requires a basepoint (its pole).

\subsection{Arithmetics vs. Geometry
(Belyi-Grothendieck vs. Ahlfors)} \label{sec:Belyi-Grothendieck}

\def\Qbar{\overline{\Bbb Q}}

[10.08.12] Closed Riemann surfaces are subsumed to the
(alienating) theorem of Belyi-Grothendieck, that {\it a surface is
defined over $\Qbar$ iff it admits a morphism to the line ${\Bbb
P}^1$ ramified at only $3$ points\/} (so-called {\it Belyi map}).
Another characterization (due to Shabat-Voevodsky 1989
\cite{Shabat-Voevodsky_1989/89}) is the possibility to triangulate
the surface by equilateral triangles (with or without respect to
the hyperbolic uniformizing metric). Basically this follows as one
may sent homographically the 3 points to the vertices of the
regular tetrahedron inscribed in the sphere. (Compare Belyi
1979/80 \cite{Beyli_1979/80}, Grothendieck 1984
\cite{Grothendieck_1984/1997-esquisse-d'un-programme} ``Esquisse
d'un progamme'', Shabat-Voevodsky 1989
\cite{Shabat-Voevodsky_1989/89}, Bost 1989/92/95
\cite{Bost_1989/92/95} (p.\,99--102), Colin de Verdi\`ere-Marin,
etc.)

Is there an analog of this result for bordered surfaces in the
context of Ahlfors (circle) mapping to the disc, and if so what is
its precise shape? In the Riemann sphere any 3 points are
transmutable through a Moebius rigid motion. The analog statement
in the disc involves either one boundary point plus one interior
point or 3 boundary points. Those are of course just the
(heminegligent) hemispherical trace of real triads on the
equatorial sphere corresponding to ${\Bbb P}^1$ with its standard
real structure. (Remember that there is an exotic twisted real
structure projectively realized by the invisible conic
$x_0^2+x_1^2+x_2^2=0$.) This lack of canonical choice of a real
triad on ${\Bbb P}^1$ could plague slightly an appropriate
bordered version of Belyi-Grothendieck. [12.11.12] More seriously
the ubiquity of real points in both those triads of the disc looks
incompatible with Ahlfors maps lacking real ramification (when
Schottky doubled to the realm of Klein's orthosymmetric curves).
Of course since bordered surfaces are in bijective correspondence
with real orthosymmetric curves, one may expect first an answer
along the line: {\it a real orthosymmetric curve is defined over
$\Qbar\cap {\Bbb R}\supset {\Bbb Q}$ iff it admits a totally real
map ramified solely at $3$ real points or at one real point and
$2$ imaginary conjugate points}. Remember yet that total reality
 means that the inverse image of the real line is the real
locus of the (orthosymmetric) curve, and since such maps lack real
ramification our naive real version of Belyi-Grothendieck looks
foiled. There seems to be a structural incompatibility between
Belyi-Grothendieck and Klein-Ahlfors. Of course our desideratum of
a simultaneous realization of Belyi-Grothendieck and
arithmetization of Ahlfors may well just be a nihilist folly. By
an ``arithmetization of the Ahlfors map'' we just mean something
in much the  same way as Belyi-Grothendieck arithmetizes Riemann's
existence theorem (any closed Riemann surface admits a morphism to
the sphere ${\Bbb P}^1({\Bbb C})$). Possibly, one should be
content with a reality version of Belyi-Grothendieck without
bringing Ahlfors' total reality into the picture. Then we have
something like {\it a real curve is defined over $\Qbar$ iff it
admits a real morphism to the line ramified above only one of the
two real triads, i.e. $0,1, \infty$ or $0,\pm i$.} A priori this
statement tolerates  both types of real curves (ortho- and
diasymmetric) and thus be more liberal than Ahlfors theorem (which
tolerate only orthosymmetric curves). Adhering instead to the
geometric interpretation of Belyi-Grothendieck (due to
Shabat-Voevodsky 1989/89 \cite{Shabat-Voevodsky_1989/89}) in terms
of equilateral triangulations might be more appealing. For
instance one can imagine an orthosymmetric real curve with an
equilateral triangulation invariant under (complex) conjugation. A
such would according to BG be defined over $\Qbar$. It is clear
that such a triangulation would contain the real circuits as
subcomplex of the triangulation. In particular what is the
significance of the corresponding vertices, e.g. as rational
points of the curve. Also the tetrahedron plays some r\^ole in
Belyi-Grothendieck-Shabat-Voevodsky and what are the r\^ole of the
other Platonic solids? In particular the octahedron looks
particularly well suited for getting pull-backed by the  Ahlfors
map? etc. [14.11.12] Of course invariant equilateral
triangulability is not reserved to orthosymmetric patterns, as
shown e.g. by the sphere acted upon by the antipodal map endowed
with a Platonic triangulation invariant under the involution
(octahedron and icosahedron). One can also consider in genus $1$ a
rhombic lattice in ${\Bbb C}$ leading to a diasymmetric (non
dividing) curve with $r=1$ real circuit. When the lattice is
equilateral say spanned by 1 and $\omega$ a cubic root of $-1$,
 we have an obvious invariant equilateral triangulation by 8
triangles (with vertices at $0,
1/2,1,\omega/2,\omega/2+1/2,\omega$ and their conjugates).

[10.08.12] Back to the closed case, we know (Mordell-Faltings ca.
1981) that when the genus is $g\ge 2$ then the curve has finitely
many rational points in any number field (finite extension of
$\Bbb Q$). Of course this fails if we raise up to the full $\Qbar$
(as slicing a plane model by rational lines gives infinitely many
$\Qbar$-points on the curve). One can dream on a  connection
between the ``canonical'' equilateral triangulation (ET) and the
finitely many rational points evaluated in the various number
fields.

Of course given an ET of an arithmetic (Riemann) surface we can
imagine a subdivision into another ET. Given a Euclidean
equilateral  triangle it is obvious how to subdivide it in 4
smaller equilateral triangles (bisecting the edges). Is there an
equivalent subdivision for hyperbolic equilateral triangles? (I
cannot see one...) Thus maybe there is some rigidity. At any rate
among all ET of an arithmetic $\Qbar$-surface there is one
involving the least number of triangles. This gives an integer
invariant for any Riemann surface defined over $\Qbar$. Can this
value be related to the finitely many rational points when $g\ge
2$?

By Gauss(-Bonnet) [$\alpha+\beta+\gamma=\pi+\int_T K dA$]
which reduces to $3\alpha=\pi-area$ for an equi-triangle in
constant negative curvature equal to $-1$ we see a direct
relation between the area and its angle of an equi-triangle.

[11.12.12] For a more lucid Real Belyi theory than our vague
ideas, compare the account in K\"ock-Singerman 2006
\cite{Koeck-Singerman_2006}, where however the Ahlfors maps does
not seem not enter the arena.

\section{Ahlfors' proof}\label{Ahlfors-proof:sec}

[27.08.12] This section is our modest attempt to examine and
understand Ahlfors' existence proof  of a circle map (of degree
$\le r+2p$). Alas we failed this basic goal, but it is perhaps of
some interest to discuss the original text while trying to capture
some mental pictures (made real) which may have circulated in
Ahlfors' vision. More objectively we also try to identify if
Ahlfors argument can be boosted to reassess the prediction of maps
with smaller controlled degree $\le r+p$ (Gabard 2006
\cite{Gabard_2006}). We emphasize once more that Gabard's result
is potentially false, but even if so, it is evident that for low
values of the invariants $(r,p)$ Ahlfors bound $r+2p$ fails
sharpness. Near its completion, Ahlfors proof takes a geometric
``tournure'' (convex geometry) where there seems to be some free
room suitable for improvements. We tried to imagine some
(topological) strategy which could possibly sharpen Ahlfors result
along his method (at least for low invariants). This is, apart
from didactic interest, the only original idea of the present
section.

In the original paper Ahlfors 1950
\cite[p.\,124--126]{Ahlfors_1950}, the existence proof, we are
interested in, occupies only a short 2 pages argument which looks
essentially self-contained albeit not quite easy to digest. I
would (personally) be extremely grateful if someone understanding
Ahlfors proof could publish a more pedestrian account than
Ahlfors', explaining it in full details. Some of the background
required is dispatched earlier in the text (esp. p.\,103--105 in
\loccit), hence trying some rearrangement could improve
readability.
We were personally not able to follow all the (boring)
computations or formulas required by Ahlfors.

{\small
Alas, big masters tend to give only cryptical output of
boring computations.
Ahlfors is further typical for his annoying (arrogant?) style ``it
is clear that'', etc. and one often suffers a lot just to fill
some details. Of course, nothing is clear in mathematics
especially when it comes to follow  mechanical computations. Maybe
the presence of those just reveals a lack of conceptual grasp over
the underlying geometry.
Trying to be more optimistic and less severe due to frustration,
it would be nice---I repeat myself intentionally---if somebody
could take the defense of Ahlfors by presenting an argument as
close as possible to the original (meaning perhaps just
eradication of misprints, if any?) which further would be
completely mechanical, i.e. where each
identity is decorated by the appropriate tag referring to the
formula under application.

}

Of course, Ahlfors' proof seems to involve nothing more than the
formalism of differential forms (\`a la Cartan, de Rham, etc.,
which he learned from A.~Weil's visit in Scandinavia during World
War II), plus Stokes' formula (already a nightmare to prove, at
least for Bourbaki) and the allied integration-by-part formula
(consequence of Leibniz's rule). We were personally unable to
produce a perfectly pedestrian (accessible to anybody, in
particular myself!) exposition of Ahlfors' account, lacking both
intelligence and patience to make his text perfectly intelligible.
The writer probably read this Ahlfors' argument several times in
diagonal (since ca. 2001/02), but never completely understood the
details. My motivation for looking at it more closely became more
acute, after realizing (August 2012) that it is not completely
trivial to complete the Green's function strategy to the problem
(cf. previous Section \ref{Green:sec}). It should be noted that
Ahlfors' argument does not employ exactly  the Green's function,
but a close relative cousin with pole located on the boundary
instead of the interior. As a matter of joking we refer to it as
the {\it Red's function\/}, and as far as we know there is no
(standardized) terminology to refer to this object! Accordingly,
Ahlfors rather constructs an half-plane map instead of a circle
map. Of course both moneys are ultimately convertible, yet both
geometrically and analytically this implies a little alteration of
the viewpoint. One may then may get a bit confused about wondering
on the optimal strategy.

Finally, remember that several workers in Japan or the US seem to
have found necessary to rework Ahlfors' proof in a more
do-it-yourself fashion. Several other authors, having to cite
Ahlfors work, often cross-cited those alternative  proofs, like
those produced  by Heins 1950 \cite{Heins_1950} or Royden 1962
\cite{Royden_1962} (cf. e.g. Stout 1972 \cite{Stout_1972} or
Gamelin 1973 \cite[p.\,3]{Gamelin_1973-Extremal-I}, who both cite
Royden for the piece of
work originally due to Ahlfors). For a more complete list of
``dissident'' authors
drifting from Ahlfors' account as the optimal source compare
Sec.\,\ref{dissident:sec}. The latter tabulation is supposed to
illustrate that I may not be isolated in having missed the full
joy of complete satisfaction with Ahlfors' output. Yet, personally
we still would like to believe that Ahlfors account is superior in
geometric quintessence to all of what followed, but only regret to
have missed some crucial details.  As far as we know, nobody ever
raised a fatal objection against Ahlfors' proof. (Personally, I
only criticize a lack of details in the execution, plus a matter
of organization\footnote{Of course this can hardly be taken
seriously, in view of the messy nature of the present text!} and
finally a lack of geometric visualization.) It may also be
speculated that the argument published by Ahlfors 1950
\cite{Ahlfors_1950} (and reproduced below) is not the way Ahlfors
originally discovered the statement (as early as 1948, cf. Nehari
1950~\cite{Nehari_1950}), which looks more intuitive when
approached from the Green's function viewpoint, or just bare
Riemann-Roch theorem (yet with dangerous probability of collision,
cf. the remark in Gabard 2006 \cite[p.\,949]{Gabard_2006}). In the
sequel we shall  attempt to conciliate Ahlfors' analytic treatment
with the  geometric intuition behind it.

The goal is (as usual) to prove:

\begin{theorem} {\rm (Ahlfors 1950 \cite[p.\,124--126]{Ahlfors_1950})}
Let $\overline{W}$ be a compact bordered Riemann surface of genus
$p$ with $r\ge 1$ contours. Then there exists a circle map
$f\colon \overline{W}\to \overline{\Delta}$ of degree $\le
r+2p=g+1$, where $g:=(r-1)+2p$ can be either interpreted as the
genus of the (Schottky) double or as the number of essential
$1$-cycles on $F$ considered up to homologies (the so-called Betti
number).
\end{theorem}

\subsection{The core of Ahlfors' argument}

For the proof Ahlfors uses the concept of a Schottky
differentials. Those are differentials on the bordered surface
which extends to the Schottky double. The following subclass plays
a special r\^ole:
$$
S_r=\textrm{ the space of analytic Schottky differentials
which are real along $C=\partial \overline{W}$}.
$$

\begin{lemma}\label{bipole:lemma} Given $g+1$ distinct points $z_j$ on the contour
$C=\partial \overline{W}$ and corresponding reals $A_j \in
{\Bbb R}$, it is possible to construct an analytic
differential $\theta_0$ which is real on\footnote{Perhaps it
would be more corrected to say ``along'' here. Compare in this
respect Ahlfors, p.\,108, the text just preceding footnote~3)}
$C$ and whose only singularities are double poles at the $z_j$
with singular parts:
$$
A_j \frac{dz}{(z-z_j)^2},
$$
where the local variable $z$ at $z_j$ is chosen so as to map
$C$ onto the real-axis ${\Bbb R}$ and inner points of $W$ into
the upper half-plane.

Further such a differential $\theta_0$ is uniquely determined
up to a differential $\theta\in S_r$, and for a proper choice
of the latter we can make vanish the periods and half-periods
of the imaginary-part $\Im \theta_0$.
\end{lemma}

Ahlfors prefers to construct instead of a circle map a upper
half-plane mapping $F\colon \overline W \to \overline H=\{ \Im z
\ge 0\}$ which will ultimately arise through the equation
$\theta_0=dF$, after arranging exactness of $\theta_0$ for a
suitable location of the $z_j$ and some $A_j\ge 0$.

Once this is achieved we may write $\theta_0=dF$ for some
analytic function $F$ on $\overline W$. The latter is uniquely
defined modulo an additive constant and can be chosen real on
$C=\partial \overline W$, except at the $z_j$ where $\Im F$
becomes positively infinite. The maximum principle ensures
$\Im F >0$ on the whole interior $W$, and therefore $F$ is the
desired half-plane mapping of degree $\le r+2p$.

This is the bare strategy of the argument, but it is time to
adventure into the details.

A first ingredient is the fact (compare the second Corollary
on p.\,109):

\begin{lemma}\label{g-dimensional:lemma}
The real vector space $S_r$ (of Schottky differentials real
along the border) has real dimension $g$.
\end{lemma}

This looks rather plausible upon thinking with the Schottky
double and explains the second (uniqueness) clause of the
above lemma. Notice indeed that there is $(r-1)$ half-periods
corresponding to pathes on the bordered surface $\overline W$
joining a fixed contour $C_1$ to the remaining ones $C_2,
\dots, C_r$ and $2p$ full periods arising by winding around
the $p$ handles.

To arrange exactness of $\theta_0$, Ahlfors employs the inner
product $(\theta_0, \theta)$ and a corresponding criterion for
exactness in terms of orthogonality to the space $S_r$ (cf.
Lemma~\ref{orthogonality:lemma} below). (The reader can skip
the proof of the next two lemmas to move directly to the core
of the argument which in our opinion is
Lemma~\ref{clever-choice:lemma}.)

Before attacking the proof we first recall the pertinent
definitions. The {\it inner product} of two differentials on a
Riemann surface is defined by:
$$
(\omega_1, \omega_2)=\int_W \omega_1
\overline{\omega_2}^{\ast},
$$
where the star denotes the {\it conjugate differential} and
the bar is the {\it complex conjugate} (compare Ahlfors,
p.\,103). (Locally if $\omega=a\, dx+b\,dy$ then
$\omega^{\ast}=-b\, dx+a\,dy$ and $\overline\omega=\bar a\,
dx+\bar b\, dy$)

Further we need probably Stokes
$$
\int_W d\omega=\int_C \omega,
$$
which combined with Leibniz
$$
d(f \omega)=df \cdot \omega + f d \omega.
$$
gives the so-called integration by parts formula
$$
\int_W (df \cdot \omega + f d \omega)\buildrel{\rm
\tiny{Leibniz}}\over=\int_W d(f\omega) \buildrel{\rm
\tiny{Stokes}}\over{=}\int_C f\omega,
$$
which can be rewritten as
$$
\int_W df \cdot \omega =\int_C f\omega- \int_W f d \omega,
$$
which is hopefully the exact form used (subconsciously) in the
sequel.

Further he requires an expression of this inner product in
term of local variables. Namely the following:

\begin{lemma}\label{loc-formula:lemma}
If $\theta=\alpha dz$ near $z_j$, then we have the following
formula for the inner product
\begin{equation}\label{loc-formula:eq} (\theta_0,
\theta)=-\pi\textstyle\sum_{j=1}^{g+1}A_j \alpha(z_j),
\end{equation}
where $\theta_0$ is the differential of
Lemma~\ref{bipole:lemma}.
\end{lemma}

\begin{proof} As in the first lemma, once we have arranged vanishing
of the period and the half-period of the imaginary part $\Im
\theta_0$ we may write something like
$$
\theta_0-\overline{\theta_0}=i \, dG,
$$
where $G$ vanishes on $C$ except at the $z_j$. Then
brute-force computation gives
\begin{equation}\label{inner-product:eq}
(\theta_0,
\theta)\buildrel{?}\over{=}(\theta_0-\overline{\theta_0},
\theta)=(i \, dG, \theta)=\dots=-\int_C G \bar \theta,
\end{equation}
where the ``dots'' indicates steps left un-detailed by Ahlfors. Of
course one should first apply the definition of the inner product
and then use  integration-by-part, as we just recalled, while
noticing that the second term vanish involving the differential of
an analytic function. [Alas, the writer had not the energy to
complete the detailed computation.]

Now writing $\theta=\alpha dz$ near $z_j$, Ahlfors claims the
following local expression for $G$
$$
G \sim i\, A_j (\frac{1}{z-z_j}-\frac{1}{\bar z-z_j}),
$$
whereupon he claims that the singularity at $z_j$ contributes
the amount $-\pi A_j \alpha(z_j)$ to the last integral of
\eqref{inner-product:eq}. The announced formula should follow
easily.
\end{proof}

\begin{lemma}\label{orthogonality:lemma}
$\theta_0$ is exact iff $(\theta_0, \theta)=0$ for all $\theta
\in S_r$.
\end{lemma}

\begin{proof} A priori we could expect to save forces by proving
only sufficiency (i.e. the implication $[\Leftarrow]$), but alas
Ahlfors' proof requires the direct sense as well, plus the
previous lemma involving the rather (unappealing) computation in
local coordinate. Enough philosophy and lamentation, and let us
follow along Ahlfors' exposition.

$[\Rightarrow]$ Write $\theta_0=dF$. Then Ahlfors write
cryptically
$$
(\theta_0,\theta)\buildrel{?}\over{=}(\theta_0,\theta+\bar
\theta)=\int_W dF \overline{\cdots}=i\int_C F(\bar \theta-
\theta)=\pi \textstyle\sum_{j=1}^{g+1 } A_j \alpha(z_j),
$$
and comparison with Equation~\eqref{loc-formula:eq} shows that
$(\theta_0, \theta)=0$, as required.

$[\Leftarrow]$ Conversely, suppose $(\theta_0, \theta)=0$ for
all $\theta\in S_r$, and let $\varphi$ be the analytic
Schottky differential making $\theta_0-\varphi$ exact. Then by
the former implication\footnote{Here our argument shorten
slightly the prose of Ahlfors, hopefully without loosing in
precision?!} $( \theta_0-\varphi, \theta)=0$ and so $(\varphi,
\theta)=0$ for all $\theta \in S_r$. This implies $\varphi=0$,
and we conclude that $\theta_0$ is exact.
\end{proof}

Combining both those lemmas, the exactness of $\theta_0$ is
reduced to the following (tricky) lemma, involving a mixture of
convex geometry and Stokes formula (which Ahlfors calls the {\it
fundamental formula\/} probably due its anticipation by Green or
Gauss and others).

\begin{lemma}\label{clever-choice:lemma}
It is possible to choose the $z_j$ and the $A_j\ge 0$ so that
\begin{equation}\label{Ahlfors_sum-Aj:eq}
\textstyle\sum_{j=1}^{g+1}A_j \alpha(z_j)=0
\end{equation}
for all $\theta \in S_r$ locally
expressed as $\theta=\alpha dz$.
\end{lemma}

\begin{proof} Let $\theta_i\in S_r$ ($i=1, \dots ,g$) be a basis of the
$g$-dimensional space $S_r$ (cf.
Lemma~\ref{g-dimensional:lemma}).  Locally we can write
$\theta_i= \alpha_i dz$ near $z_j$.
Equation~\eqref{Ahlfors_sum-Aj:eq} can be satisfied with
$A_j\ge 0$ iff the simplex with vertices
$$
(\alpha_1(z_j), \dots, \alpha_g(z_j))\in {\Bbb R}^g
\qquad\textrm{ for } j=1, \dots, {g+1}
$$
contains the origin $0\in {\Bbb R}^g$.

If this condition is not full-filled for any choice of the
$z_j$, the convex-hull of the set of points
$$
K:=\{(\alpha_1(t), \dots, \alpha_g(t)) : \textrm{ for } t\in C
\}
$$
would fail to contain $0$. (One can think of this set as a
sort of link (in the sense of knot theory) traced in ${\Bbb
R}^g$ with $r$ components. However the latter is not perfectly
canonical since the $\alpha_i(t)$ depends on the local chart.

{\it Expressing some naive doubts.} So here Ahlfors argument looks
a bit fragile (or at least sketchy) as one probably requires to
fix a finite system of holomorphic charts covering the full
contour of the bordered surface). [We do not have a specific
objection, yet it should be noted that the whole Ahlfors theory
even that of the refined extremal problem depends on the
non-emptiness of the class of bounded functions, hence upon the
present argument! In principle even if there should be a
global crash of Ahlfors' proof here, then the theorem should
conserves its validity in view of several subsequent treatments
hopefully logically more reliable, we cite:

$\bullet$ Kuramochi 1952 \cite{Kuramochi_1952} (alas quite
unreadable?),

$\bullet$ Mizumoto 1960 \cite{Mizumoto_1960},

$\bullet$ Royden 1962 \cite{Royden_1962} (alas a bit
functional-analytic, whereas the statement sentimentally
belongs to pure geometric function theory), and maybe

$\bullet$ Gabard 2006 \cite{Gabard_2006} (hopefully reliable, at
least it first part not improving Ahlfors' $r+2p$).

However it is likely that the set $K$ can be defined according
to the totality of possible $\alpha_i(t)$ arising through a
fixed system of permissible charts covering the contour $C$.

Now a (Euclidean) set of ${\Bbb R^g}$ whose convex-hull
misses the origin is contained in a closed half-space [maybe even
an open half-space?]. Thus there exists scalars $a_1, \dots, a_g
\in {\Bbb R}$ (not all zero) so that
$$
\textstyle\sum_{i=1}^g a_i \alpha_i(t)\ge 0 \quad \textrm {for
all } t\in C.
$$
(Geometrically, this is to be interpreted as the scalar
product with the vector $(a_1, \dots, a_g) \in {\Bbb R}^g$
orthogonal to the hyperplane whose half contains the set $K$.)
Hence the corresponding differential $\theta= \sum_{i=1}^g a_i
\theta_i$ is $\ge 0$ along $C$. [Maybe strict???] However this
violates the fact that $\int_C \theta=0$, as prompted by
Stokes' formula
$$
\int_{C=\partial \overline W} \theta=\int_{\overline W} d
\theta,
$$
and the fact that $\theta$ belongs to $S_r$, hence analytic,
and thus closed, i.e. $d\theta =0$.
\end{proof}

\subsection{Geometric interpretation as dipoles}

[28.08.12] Let $F$ be a membrane (=compact bordered Riemann
surface), then Ahlfors constructed (cf. previous subsection) a
half-plane map $F\to \overline H:=\{\Im z\ge 0 \}$ to the closed
upper-half plane. We get a circle map after post-composing with
the natural conformal map to the unit-disc $ \overline H \to
\overline{\Delta}$. Under such a map, the horizontal lines
transforms to a pencil of circles tangent to the boundary and
vertical lines mutate to arc of circles orthogonal to the
boundary. (cf. Fig.\,\ref{Dipole:fig}a). One recognizes
essentially the so-called Hawaiian earrings (cf.
Fig.\,\ref{Dipole:fig}b).

Given a circle map, one can pull-back the isothermic
(=right-angled) Hawaiian bi-foliation to obtain a graphical
representation of the circle map.

\begin{figure}[h]
\centering
    \epsfig{figure=Bipole.eps,width=122mm}
\vskip-5pt\penalty0
  \caption{\label{Dipole:fig}
  Pictures of dipoles: another attempt to visualize
  Ahlfors circle maps}
\vskip-5pt\penalty0
\end{figure}

Starting with with a (doubly-connected) ring, one obtains
Fig.\,\ref{Dipole:fig}c or Fig.\,\ref{Dipole:fig}d. Going to
higher connectivity on gets for instance Fig.\,\ref{Dipole:fig}e.
The Bieberbach-Grunsky theorem (or just Riemann-Roch, cf. e.g.
Lemma~\ref{Enriques-Chisini:lemma}) tell us that we can prescribe
 a point on each contour and there is a circle
map taking all those points to the same image in the
unit-circle $S^1=\{ \vert z \vert=1\}$. Hence, we enjoy
complete freedom in picturing the isothermic bi-foliation of
circle maps, at least in the planar case. This situation is to
be contrasted with the situation for the zeros, where some
hidden symmetry requires to be fulfilled (compare e.g. Gabard
2006, where we have the condition $D\sim D^{\sigma}$ of linear
equivalence of the divisor with its conjugate, an also Fedorov
1991 \cite{Fedorov_1991} who speaks of an opaque condition
that must be satisfied to prescribe the zeros).

Of course, the contemplation (and manufacture) of such pictures
raises more questions than clarifying the perception of Ahlfors'
theorem. One can hope some guidance via physical intuition (if one
feels comfortable with the mineral world) or appeal again to the
metaphor about  proliferation of bacteria in some nutritive
medium. We do not repeat the long discourse we made already for
Green (cf. Sec.\,\ref{Green:sec}, esp. Fig.\,\ref{Green:fig})
where one had radial expansions emanating from an inner point.
Presently, the bacteria are rather located on the boundary,
whereupon their local expansion is more of the Hawaiian type, or
if you prefer look alike the Doppler effect at the critical speed
of sound. The dipole of our title would occur upon considering the
symmetric Schottky double of the membrane. This new
Hawaiian/Doppler mode of expansion can again be explained via
lacking nutritive resources caused by the boundary where the world
stops.

On Fig.\,\ref{Dipole:fig}f, we have attempted to picture the
pull-back of the Hawaiian foliation under a circle-map of degree
$r+p=1+1=2$ (for the value $r+p$ predicted by Gabard). This
picture looks anomalous for the following reason. Letting grow the
population, there is a first junction of the 2 populations right
``under'' the handle, then there is 2 self-junction at 2 points
aside the handle. From now on the bacteria starts invading the
handle from both ``sides'' and will actually merge on the core
circle of it. This is problematic since ultimately the expansion
should finish along the boundary contours (by definition of a
circle-map). It easy to manufacture  a picture where no such
anomaly occurs (cf. e.g. Fig.\,\ref{Dipole:fig}g which admittedly
requires some little effort of concentration to contemplate its
morphogenesis). Of course, similar pictures can be made by
prescribing less boundary points than the degree of circle-maps
predicted by Ahlfors $r+2p$ or $r+p$, e.g. by a choosing a single
dipole, cf. Fig.\,\ref{Dipole:fig}h and Fig.\,\ref{Dipole:fig}i.
However those patterns cannot correspond to circle-map due to
obvious topological obstructions: first the degree of a circle-map
must be $\ge r$ impeding Fig.\,\ref{Dipole:fig}h to be allied to a
circle-map. As to Fig.\,\ref{Dipole:fig}i the degree would be one,
implying the circle-map to be unramified and covering theory (of
the simply-connected disc) implies the membrane $F$ to be the
disc, violating its genus $1$ nature.

We stop this graphical discussion at this primitive stage, yet it
is to be hoped that a deeper study of such figures could lead to
some theoretical results complementing the understanding of the
Ahlfors maps. Perhaps such (dipole) isothermic drawings are of
some relevance to Gromov's filling conjecture, as we already
suggested in the case of Green's function
(Sec.\,\ref{sec:Green-to-Gromov}).

[29.08.12] In fact there is a another more convincing obstruction
impeding Fig.\,\ref{Dipole:fig}f to represent a circle-map. This
consists in identifying the counter-images of the growing Hawaiian
circles past the critical levels while checking if they contribute
to the correct numerical multiplicity permissible with the degree
of the branched covering. To be concrete we enumerate a series of
typical
smooth levels on Fig.\,\ref{Dipole:fig}f. The first one denoted
$1$ consists of $2$ little circles. Past the first critical level,
we see the curve $2$ with $1$ component. After the next critical
level, we pick a curve $3$, which has $3$ components. This is too
much for our mapping to be of degree 2. This proves that
Fig.\,\ref{Dipole:fig}f do not correspond to a circle-map.

In contrast repeating the same counting exercise for
Fig.\,\ref{Dipole:fig}g, no such excess occurs. The level $1$ has
2 components, level 2 (chosen after the first critical level) has
one component, level 3 has 2 components and finally level 4 has 1
component. Thus the picture looks topologically coherent, but it
is evident that it is far from metrically realist. Naively
speaking we were forced to distort the propagation so has to have
a virtually planar mode of depiction for the levels.


\subsection{Trying to recover Ahlfors from the Red's
function}\label{Red's-function:sec}

[29.08.12] Let $F$ denote a finite (=compact) bordered Riemann
surface of genus $p$ and with $r$ contours.
From the previous section, it seems evident that there is some
canonical function akin to the Green's function yet with pole
pushed to the boundary (dipole singularity when doubled). Call
them perhaps the {\it Red's function\/} as an {\it ad hoc} acronym
honoring writers like Riemann, Schwarz, Klein, Koebe, Ahlfors,
etc. Such a Red's function denoted $R(z,t)=R_t(z)$ with (di)pole
at $t\in \partial F$ (a boundary-point) is defined by the property
of being harmonic, null along $\partial F$ save at $t$ where it
becomes positively infinite according to a specific local
singularity (maybe like ${\rm Re}(1/z^2)$). [18.10.12] As a more
intrinsic definition one can define $R_t$ as the unique positive
harmonic function vanishing continuously along $\partial F-\{
t\}$. The function then looks unique up to scalar multiple. Note
however that Heins (in e.g. Heins 1985 \cite[p.\,241, right after
Thm\,3.1]{Heins_1985-Extreme-normalized-LIKE-AHLF}) defines the
function $u_{\zeta}$ our $R_t$ by adding the requirement of
minimality (in the sense of Martin 1941 \cite{Martin_1941}). A
positive function $u$ is minimal if whenever there is a smaller
function $0<v<u$, $v$ is a constant multiple of $u$.

The sudden explosion of $R_t$ at just one boundary point looks at
first almost paradoxical, but see again our previous
Fig.\,\ref{Dipole:fig}b-h-i) for a depiction of their levels and
one can of course imagine such a function just as a ``borderline''
degeneration of the usual Green's function. Now one can attempt to
construct a half-plane-map (HP-map, for short), by considering a
superposition $R(z):=\sum_{i=1}^{d} R(z, t_i)$ of such Red's
functions $R(z,t_i)$ for several points $t_i$ on the border. The
formula
$$
\varphi:=R+i R^{\ast},
$$
where $R^{\ast}$ is the conjugate function would then define the
HP-map provided the conjugate potential is single-valued in other
word that the conjugate differential of $R$, $(d{R})^{\ast}$ is
period free. Since $F$ has $(r-1)+2p=:g$ essential cycles
(homologically independent), a parameter count suggests that if
$d=g+1$ there is enough freedom to annihilate all the $g$ periods
of $d{R}^{\ast}$.

Maybe this approach (which presumably differs not very much from
Ahlfors') has some technical advantage over the Green's technique
(presented in Sec.\,\ref{Green:sec}). First it seems that the
dipole singularity has some linear character contrasting with the
arithmetical rigidity of the logarithmic singularity. Thus it is
permissible to form a more general linear combination
$$
R(z):=\textstyle\sum_{i=1}^{d} \lambda_i R(z, t_i),
$$
with some reals $\lambda_i$ which must however be $\ge 0$.
Hence killing the periods essentially reduces to  linear
algebra.
Another advantage over the Green's approach stems from the
fact that in the interior we meet no singularity thus the
period mapping looks less dubious.

As usual we write down the period mapping by integrating the
$1$-form $dR^{\ast}$ along the $g$ many 1-cycles $\gamma_1, \dots,
\gamma_g$ and obtain for each fixed $t_1, \dots t_{g+1}\in
\partial F$ a linear map ${\Bbb R}^{g+1} \to {\Bbb R}^g$.
Thus there is some non-zero vector in the kernel, and the
corresponding $(\lambda_i)$ would solve the problem, provided
one is able to check that they can be chosen $\ge 0$. This is
non-trivial and a priori it is not evident (and nobody ever
asserted) that this can be done for any choice of the
$(g+1)$-tuple $t_i$.

So it is just here that the difficulty starts, and  that some
idea is required to complete the proof.

[04.09.12] Due to a lack of creativity/energy, I was blocked
here for a couple of days. So let me make a list of writers
who seem to have grasped the geometric quintessence of
Ahlfors' argument:

$\bullet$ Gamelin-Voichick 1968
\cite[p.\,926]{Gamelin-Voichick_1968}: ``According to [1,
\S\,4.2](=Ahlfors 1950 \cite{Ahlfors_1950}), there exist $r+1$
($r=g$ in our notation) points $w_1, \dots, w_{r+1}$ on $bR$
such that if $B_j$ is the period vector of the singular
function $T_j$ corresponding to a unit point mass at $w_j$,
then $B_1, \dots, B_{r+1}$ are the vertices of a simplex in
${\Bbb R}^r$ which contains $0$ as an interior point.''
[10.09.12]

$\bullet$ Fisher 1973 \cite[p.\,1187/88]{Fisher_1973}: ``By a
theorem of Ahlfors [A1; \S4.2] there is a set of $r+1$ points
$p_j$ in $\Gamma$ such that if $v_j$ is the period vector of a
unit mass at $p_j$, then $v_0,\dots, v_r$ form the vertices of a
simplex in ${\Bbb R}^{r}$ which contains the origin as an interior
point.'' [this looks alike verbatim copy of the previous source,
yet reinforce confidence in the viewpoint]

[07.09.12] In fact some little hope to complete the argument is
raised by borrowing ideas of convex geometry used by Ahlfors, yet
in our context which is perhaps not so reliable (albeit it seems
to match with the Gamelin-Voichick twist of Ahlfors). Alas, we
failed to recover Ahlfors statement, but we see obvious room for
improving upon Ahlfors by using essentially his method of proof
augmented by some further geometric tricks. Ideally one would like
to recover the bound predicted in Gabard 2006 \cite{Gabard_2006}
by using an argument very close to Ahlfors'. Let us now be more
concrete.

Again we fix some $d$ points $t_1, \dots t_d$ on the boundary
$\partial F$, with at least one point one each contour $C_i$
(forming the boundary $\partial F$). For any point $t\in
\partial F$ the function $R_t(z):=R(z,t)$ is uniquely defined
once a chart around $t$ is specified (otherwise it is unique only
up to a positive scaling factor). Let us assume $R_t$ fixed once
for all with a continuous dependence over the parameter $t$. (Alas
the writer has no clear-cut justification of this possibility.
[09.09.12] Maybe use a boundary uniformizer for an annular tubular
neighborhood of each contour, cf. e.g. Hasumi 1966
\cite[p.\,241]{Hasumi_1966}, also Gamelin-Voichick 1968
\cite[p.\,926]{Gamelin-Voichick_1968}. [18.10.12] Of course since
$R_t$ is unique up to scalar multiple, we are somehow choosing a
section of a ray-bundle and even if after winding once around an
oval of $\partial F$ the $R_t$ should not return to its initial
position $R_{t_0}$, it seems easy to apply a sort of ``closing
lemma'' so that $R_t$ comes back to the original choice.)

We now introduce $\Pi(t)$ the period of $(dR_t)^{\ast}$ along the
fixed representatives $\gamma_1, \dots, \gamma_g$ of the first
homology, that is,
$$
\Pi(t)=(\textstyle\int_{\gamma_1} (dR_t)^{\ast}, \dots,
\int_{\gamma_g}(dR_t)^{\ast}) \in {\Bbb R}^g.
$$
We seek $R$ of the form $R=\sum_{i=1}^{d} \lambda_i R_{t_i}$ with
$\lambda_i>0$ such that the conjugate differential $(dR)^{\ast}$
is period-free. Period-freeness amounts to say that {\it the
simplex of ${\Bbb R}^g$ spanned by the $\Pi(t_1), \dots, \Pi(t_d)$
contains the origin in its interior}\footnote{In the combinatorial
sense, by opposition to the topological sense.}. Then positive
masses $\lambda_i$ can be assigned to the $\Pi(t_i)$ so that the
origin occurs as barycenter of this masses distribution.

The italicized condition is equivalent to saying that the
convex-hull of the set $X:=\Pi(\partial F)$ contains the origin
(say then that the set $X$ is balanced). Balancedness paraphrases
also into the condition that the set is not contained in a
half-space delimited by a hyperplane through the origin.

Ahlfors derives his result from the following simple lemma applied
to $X=\Pi(\partial F)$.

\begin{lemma} Let $X$ be a subset of some number space ${\Bbb
R}^g$. Any point in the convex-hull of $X$ is the barycenter
(=convex combination involving positive coefficients) of at most
$g+1$ points of $X$.
\end{lemma}

\begin{figure}[h]
\centering
    \epsfig{figure=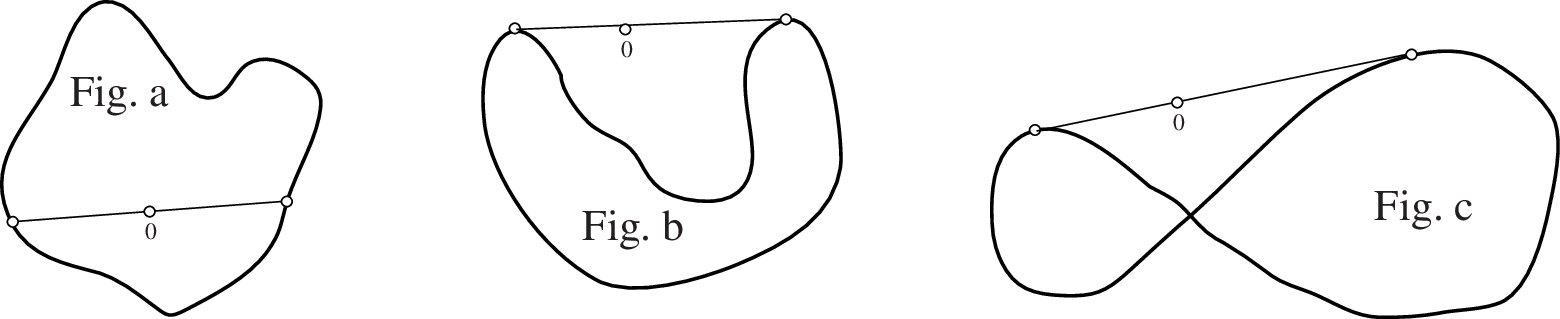,width=122mm}
\vskip-5pt\penalty0
  \caption{\label{Convex:fig}
  Improving upon Ahlfors by using Ahlfors}
\vskip-5pt\penalty0
\end{figure}

Of course the lemma is sharp in general: consider $X\subset {\Bbb
R}^2$ a set of 3 points in general position (not collinear) then
any point chosen in the interior of the convex-hull of $X$ (a
simplex) requires all 3 points in a barycentric combination.
However if $X$ is a more continuous shape like a topological
circle in ${\Bbb R}^2$ it is clear that 2 points situated on $X$
will suffice (cf. Fig.\,\ref{Convex:fig}a). Indeed, imagine first
that $X$ is a Jordan curve and that the point lies in its
interior. Any line through the point intercepts the Jordan curve
in at least 2 points which can be used for a convex combination of
the given point. If the point is not in the interior, one can meet
an ``U-shaped'' Jordan curve where the point is situated near the
top of the ``U'' (Fig.\,\ref{Convex:fig}b), yet still expressible
as the barycenter of 2 points on the top of the ``U''.  This
already raises some hope upon improving Ahlfors, and
optimistically a careful inspection could recover the $r+p$ bound
of Gabard 2006 \cite{Gabard_2006}.

Let us summarize the situation. The lemma shows is that if the
convex-hull of $\Pi(\partial F)$ contains the origin $0$, then one
can certainly find  $g+1$ points $t_i$ (eventually fewer) and
corresponding $\lambda_i>0$ such that $R=\sum \lambda_i R_{t_i}$
has a period-free conjugate differential. This implies the
existence of a half-plane map (via $f=R+iR^{\ast}$) of degree $\le
g+1=r+2p$, recovering therefore Ahlfors' result of 1950.

Thus the problem splits in two parts:

$\bullet$ Step (1): explain why the convex-hull of $\Pi(\partial
F)$ contains the origin $0$ (implying Ahlfors' $r+2p$ bound);
(Ahlfors is able to do this, yet hopefully the ambient context of
his argument can be slightly simplified to our present setting
which is  closer say to Heins' accounts in 1950 \cite{Heins_1950}
or 1985 \cite{Heins_1985-Extreme-normalized-LIKE-AHLF})

$\bullet$ Step (2): try to lower Ahlfors degree $r+2p$ by taking
advantage of the fact that $X=\Pi(\partial F)$ is not an arbitrary
set but the continuous image of $r$ circles; (ideally try to
recover the $r+p$ upper bound predicted in Gabard 2006
\cite{Gabard_2006}, or at least partial improvements of Ahlfors
bound $r+2p$ for low values of the invariants $(r,p)$).

As to the first point (1), we notice that if it is violated then
the set $\Pi(\partial F)$ is contained in a half-space of ${\Bbb
R}^g$. Thus there is a non-zero vector $a=(a_1,\dots, a_g)\in
{\Bbb R}^{g}$  such that the scalar product $(a, \Pi(t))>0$ for
all $t\in \partial F$. This means
$$
\sum_{i=1}^g a_i \int_{\gamma_i}(dR_t)^{\ast}>0 \textrm{ for
all } t\in \partial F.
$$
Alas, the writer failed to find a reason why this should be a
contradiction. (In Ahlfors' presentation Stokes' theorem plays a
crucial r\^ole.)

Even if the present geometric strategy (cooked by the writer via
slow assimilation of the very classical strategy of annihilating
periods) should be impossible to complete, nothing forbids to
switch again to the original treatment of Ahlfors, and apply our
Step (2), whose tangibleness relies on Fig.\,\ref{Convex:fig}. The
essential point is that ultimately the geometric setting is
invariably the one and same problem of convex geometry, whether we
start from Ahlfors ``analytic'' approach or from our more
geometric reformulation via the Red's functions.

Let us be more explicit. We have a map $\Pi\colon
\partial F=:C \to {\Bbb R}^g$. (``$C$''
for contours, like in Ahlfors notation.) In Ahlfors' paper this
occurs as the map $C\ni t \mapsto (\alpha_1(t), \dots ,
\alpha_g(t))$ cf. p.\,125 of his article. (From the
algebro-geometric viewpoint this must probably be the vectorial
lift of the so-called {\it canonical map\/} $\varphi\colon C \to
{\Bbb P}^{g-1}$ (usually ascribed to Noether or Klein) allied to
the canonical series $\vert K \vert$ living over the curve $C$,
obtained by doubling the bordered Riemann surface.)

We try to address the second issue (2). The setting is a map
$\Pi\colon C \to {\Bbb R}^g$ whose image is balanced (i.e. the
convex-hull of the image contains the origin, or equivalently the
set $\Pi(\partial F)$ is not contained in any open half-space of
${\Bbb R}^g$ delimited by a hyperplane through the origin). The
whole problem is then reduced to the following geometric question.

\begin{prob} \label{problem:Ahlfors-circuit} Given two integers $r\ge 1$ and $p\ge 0$.
Let $g:=(r-1)+2p$, and suppose given in the corresponding
Euclidean space ${\Bbb R}^g$ a collection of $r$ (possibly
singular) circles $C_1, \dots, C_r$. It is assumed that the union
of all these circles is balanced. Find the minimum cardinality of
a group of $d$ points with at least one point on each $C_i$
spanning a simplex containing the origin.
\end{prob}

The previous lemma solves the problem for degree $d=g+1=r+2p$
(recovering Ahlfors' result). To do better we start from such a
group and try to move the vertices, while taking care that the
simplex still contains $0$. From the $r+2p$ points, we  imagine
$r$ many as essentially fixed and the other coupled in $p$ many
pairs. The initial simplex is top-dimensional matching the
dimension $g$ of the ambient number-space ${\Bbb R}^g$. Moving
vertices, it looks reasonable that we may coalesce two points of
the $g$-simplex to get a $(g-1)$-simplex still containing $0$.
This presupposes both coalescing points being located on the same
circuit $C_i$ (try to argue with the pigeon hole principle). After
$p$ such
collisions (one for each pair) we reach the degree $r+p$ predicted
by Gabard 2006~\cite{Gabard_2006}.

Alas this ``piano mover'' argument is not easy to believe, nor to
prove. Perhaps a less naive  variant involving an adequate trick
(most probably of a topological nature akin say to the Borsuk-Ulam
proof of the ham-sandwich theorem) could recover the $r+p$ bound.
Less optimistically, it may happen that the above problem is not
always soluble with $d\le r+p$, but only for circuits $C_i$
arising from  bordered Riemann surfaces via the period map recipe.

At any rate, we see the prominent r\^ole of convex geometry in the
question of the least possible degree of the Ahlfors function. In
principle there is a canonically defined set $\Pi(C)\subset {\Bbb
R}^g$ (we shall call the {\it Ahlfors figure}) whose spanning
simplices going through the origin affords a complete
understanding (in theory at least) of the minimal degree of a
circle map concretizing the given bordered surface $F$.

[11.09.12] Perhaps one can solve the above problem
(\ref{problem:Ahlfors-circuit}) for $d=r+p$ by an inductive
procedure. Let us sketch an attempt that fails (reasonably close
to the goal). Recall that given two integers $(r,p)$ and a
balanced configuration of $r$ circles $C_i$ in ${\Bbb R}^g$, where
$g:=(r-1)+2p$. We would like to show that the origin is the
barycenter of at most $r+p$ points with at least one on each
$C_i$. Of course the assertion holds true when $p=0$, because we
know (by the lemma) that $d\le g+1=r+2p=r$ and on the other hand
we have the trivial lower-bound $r\le d$ imposed by the fact that
each circle supports at least one point. It follows that
$d=r=r+p$, and the claim is vindicated.

Thus one can try an induction reducing to the ``planar case''
$p=0$. This can be done in several ways via the moves
$(r,p)\mapsto (r,p-1)$, or $(r,p)\mapsto (r+1,p-1)$ or finally
$(r,p)\mapsto (r+2,p-1)$. The latter of which has the advantage
that the new value of $g$, denoted $g'$ stays invariant. Now given
a geometric configuration of type $(r,p)$ in the number-space
${\Bbb R}^g$ we construct one of type $(r+2, p-1)$ in the same
${\Bbb R}^g$, maybe naively just by duplicating two of the circles
(i.e., assigning them a multiplicity). This new configuration is
still balanced, so by induction hypothesis the origin is
expressible as the barycenter of $r'+p'=(r+2)+(p-1)=r+p+1$ points
located on the $C_i$. Alas, this exceeds by one unit the desired
$r+p$.

\begin{figure}[h]
\centering
    \epsfig{figure=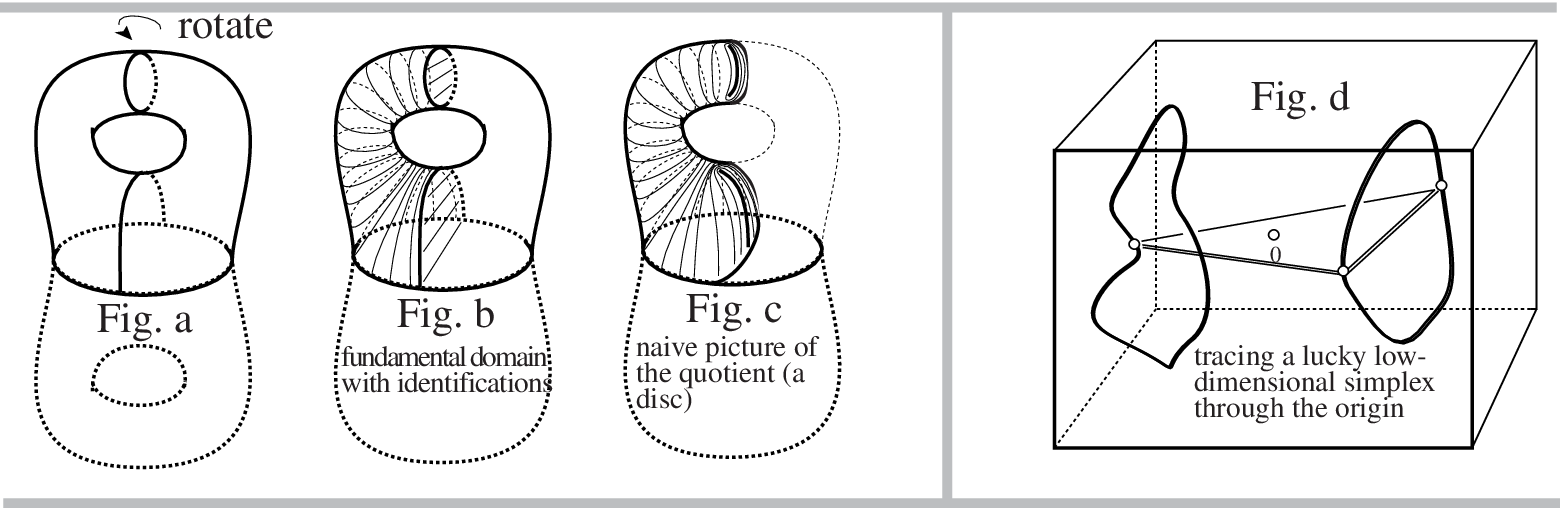,width=122mm}
\vskip-5pt\penalty0
  \caption{\label{Convex2:fig}
  Low-dimensional improvements of Ahlfors' convexity argument}
\vskip-5pt\penalty0
\end{figure}

[18.10.12] Low-dimensional examples may help to give some weak
evidence toward  solving Problem~\ref{problem:Ahlfors-circuit}
with Gabard's bound $d=r+p$. Let us discuss this aspect. If we
take $(r,p)=(1,1)$, then $g=2$. So geometrically we have one
circuit in the plane ${\Bbb R}^2$. In this situation our
Fig.\,~\ref{Convex:fig} prompts solubility of the problem with
$d=2$. Note the agreement with Gabard's bound $r+p$. This proves
the (modest) theorem that {\it a bordered surface with one contour
and of genus one always admits a circle map of degree $2$\/},
whereas Ahlfors only predicts degree $r+2p=1+2=3$. Another
evidence comes from the well-known hyperellipticity of genus 2
curves. Indeed the double of such a membrane having genus 2, it is
hyperelliptic and can therefore be visualized in 3-space as
something like Fig.\,\ref{Convex2:fig}a. Doing a rotation of angle
$\pi$ we find the required circle map of degree 2 (look at the
Figs.\,\ref{Convex2:fig}b and \ref{Convex2:fig}c).

Let us next examine the case $(r,p)=(2,1)$, then
$g=(r-1)+2p=1+2=3$. So we have 2 circuits in space ${\Bbb R}^3$
(like in knot or link theory). Since the set of circuits is
balanced, we have something like Fig.\,\ref{Convex2:fig}d
(assuming no knotting for simplicity). Balancing amounts
picturesquely to say that if you dispose of a 180 degrees angular
vision  (like any respectable homo sapiens) you will never be able
from the origin to contemplate the full link. Paraphrased
differently, whatever the direction you choose to focus your
vision the link will always move in your back. It seems plausible
that, instead of the 4 points prompted by Ahlfors' top-dimensional
$3$-simplex, 3 points actually suffices to span a $2$-simplex
passing through the origin (see again Fig.\,\ref{Convex2:fig}d).
Justifying this intuition could again corroborate the $r+p$ bound
(at least for low invariants). Of course the genus (of the double)
being now $g=3$ there is no hyperelliptic reduction, yet appealing
to the canonical map $C_g\to {\Bbb P}^{g-1}$ (an embedding
precisely when the curve is not hyperelliptic) our curve is
concretized as a plane quartic (the canonical divisor $K$ having
degree $2g-2$). Some basic knowledge of Klein's theory then
prompts that our orthosymmetric real quartic with $r=2$
 must consist of two nested ovals. Projecting from a real
point on the inner oval gives a totally real morphism of degree
$4-1=3$, in accordance again with the $r+p$ bound.

All these little
experiments raise the hope that Ahlfors original approach suitably
sharpened by a  geometric lemma about balanced collections of
circuits in ${\Bbb R}^g$ should enable some improvements, and
eventually confirm the prediction of the $r+p$ bound. However we
confess that the required positive solution to
Problem~\ref{problem:Ahlfors-circuit} with $d=r+p$ looks difficult
to obtain and perhaps only true for special circuits arising
through period maps. It is quite hard to connect Ahlfors method
with the one in Gabard 2006 \cite{Gabard_2006} in which Abel's map
was exploited more systematically. Since both maps,  $\Pi$ an
Abel, involve periods, a natural guess is that {\it Ahlfors'
figure\/}, that is the set $\Pi(\partial F)\subset {\Bbb R}^g$, is
closely related to the Abel map or at least the so-called
(Noether-Klein) canonical map $C\to {\Bbb P}^{g-1}$ which is just
the Gauss map of the Abel map: each tangent to the curve seen in
its Jacobian is reported to the origin via translation in the
Jacobi torus. If so interpretable, it is perhaps no surprise that
Ahlfors approach is cumbersome because one is working in the Plato
cavern where the essence (embedded-ness) of things is lost.

Still, the Ahlfors figure is perhaps useful for other questions.
For instance if we take a top-dimensional spanning simplex with
$g+1=r+2p$ vertices containing $0$ in its interior,  it is clear
that we may perturb slightly the vertices  keeping the origin
inside the  simplex. This shows a sort of topological stability of
Ahlfors maps having degrees $r+2p$. (This phenomenon is not new,
compare \v{C}erne-Forstneri\v{c} 2002
\cite{Cerne-Forstneric_2002}.) The same stability cannot be
expected with the more economical $r+p$ bound, for a slight
perturbation of our hypothetical simplex will generally miss the
origin. Ahlfors' figure also shows existence of circle maps for
each degree $\ge r+2p$. For those of degrees $>r+2p$ there is a
menagerie of convex combinations expressing $0$ and accordingly
plenty of circle maps having the same fibre above a boundary
point. Such results look not easily accessed via Gabard's method
(in Gabard 2006 \cite{Gabard_2006}).

Trying to make the last ``menagerie'' point more accurate could
lead to interesting result. For simplicity imagine ${\Bbb R}^g$ as
the plane ${\Bbb R}^2$ and in it a $2$-simplex spanning the
origin. If we have more than $(g+1)$ points, say $g+2=4$ then we
may interpret the convex-hull of those 4 points as the shadow
(projection) of a $3$-simplex living in ${\Bbb R}^3$. Hence above
the origin there is a segment in this higher $3$-simplex each
elements of which is a convex sum of the 4 vertices. Hence we get
$\infty^1$ circle maps having the same 4 points as prescribed
value. This requires of course to be better presented but should
be straightforward application of Ahlfors method.

\subsection{Strip mappings (Nehari, Kuramochi)}

[31.08.12] As we saw instead of a circle map, Ahlfors 1950
\cite{Ahlfors_1950}  prefers to construct a half-plane map.
Ultimately this amounts to the same except that the disc instead
of being decorated by the polar coordinates it is by the Hawaiian
dipole (Fig.\,\ref{Dipole:fig}a). A third option is to envisage
(as Nehari and Kuramochi 1952 \cite{Kuramochi_1952}) a strip
mapping to the strip $S:=\{z: -1\le\Re(z)\le 1\}$. When
rectangular coordinate on the strip are transplanted to the disc
we obtain a dipole looking like a mitosis. This yields yet another
isothermic system on the disc.

To synthesize, the disc can be decorated by 3 types of
isothermic coordinates (systems):

(1) the monopole attached to an inner point of the disc, which
when the pole is the center is just the foliation by
concentric circles plus the orthogonal rays. We may from here
drag the pole away from the center to get other isothermic
systems best interpreted as the geodesic expansion w.r.t. to
the hyperbolic metric on the disc. Upon letting degenerate the
pole to the boundary circle we get:

(2) the dipole depicted on Fig.\,\ref{Dipole:fig}a and finally
upon disintegrating this source of multiplicity 2 into two
separate elements of multiplicity one we get:

(3) a genuine dipole which ultimately can be the mitosis about
antipodal points of the circle.

In principle to each of these geometric decoration of the disc
corresponds an existence-proof of the Ahlfors function differing
so-to-speak just in the ``cosmetic details''.

Finally,  each isothermic system suggests an angle of attack
to Gromov's filling conjecture. Eventually, it seems plausible
that the totality of those isothermic systems could be
exploited collectively upon using an averaging process
(somehow reminiscent to L\"owner-Pu's trick).

\section{Hurwitz type proof of Ahlfors maps?}\label{Hurwitz-type}

[21.10.12] This section wonders about an elementary
existence-proof of circle maps via a continuity method reinforcing
some naive moduli count. As we noted (in
Sec.\,\ref{Minimal-sheet:sec}) the disaster with bordered surfaces
is that their gonality is not prompted by a naive moduli count,
and thus the project looks from the scratch a bit hazardous.
However it is not impossible that we missed something crucial.

The general philosophy  would be not to fix a surface and try hard
to find a map, but rather to look at all possible maps and lift
the complex structure of the disc while hoping that if the degree
is large enough there are sufficiently many free parameters to
paint the full moduli space. Hence any Riemann surface would be
expressible as a branched cover of the disc of some controlled
degree. (Natanzon suggested to me this strategy during an oral
conversation at the Rennes conference 2001, and I came again to
this idea by reading Natanzon et al. 2001
\cite{Natanzon-Shapiro-Vainshtein_2001/XX}.)

The basic idea may be formalized as follows. We fix a topological
type $(r,p)$ encoding the number of contours and the genus. We
introduce the (Hurwitz) space
$$
H_{r,p}^d:=\textrm{ set of all circle maps from surfaces of type }
(r,p) \textrm{  having degree } \le d.
$$
An element of this natural set (hence a space!) is a branched
cover of which we may keep in mind only the ``total space''. This
gives a map
$$
\tau\colon H_{r,p}^d\to {\cal M}_{r,p},
$$
to the moduli space of bordered surfaces of type $(r,p)$. We want
to show that this mapping is surjective for $d$ sufficiently large
(but controlled \`a la Ahlfors). First, we know (since Klein
essentially) that ${\cal M}_{r,p}$ is connected. Thus it would be
enough to find a suitable $d$ so that the $\tau$-image is closed,
open and nonempty.

As $(r,p)$ is fixed we may omit it from the notation. Of course
$H^d:=H_{r,p}^d$ is empty when $d<r$. The example of rotational
surfaces (cf. Fig.\,\ref{Chambery:fig}) shows that $H^d$ is
non-void for $d=r$ or $d=r+1$ when $r$ is even resp. odd.

It seems also trivial (since we have defined $H^d$ by the
condition $\deg(f)\le d$) that the image $\tau(H^d)$ is closed for
any $d$. Intuitively a map can degenerate to a map of lower
degree, but will never degenerate to one of higher topological
complexity. Observationally, this is well seen on the example of
the G\"urtelkurve (plane quartic with two nested ovals): when
projected from a point in the interior of the oval we get a total
map of degree 4, which can degenerate to one of degree 3 if the
center of projection is specialized toward the inner oval.
However, if we take a sequence of maps of degree 3 given by such
projections the limit will be a similar projection (the oval being
closed) and we never reach a map of degree 4. Of course an
abstract explanation requires be given (perhaps just by
compactness of $H^d$).

The hard part is to show that $\tau$ is open for some large $d$.

Naively one could hope to do this via Brouwer's invariance of the
domain requiring something like $\tau=\tau_d$ being \'etale for a
suitable $d$.

Another idea is perhaps to factorize $\tau$ by taking the fibre of
the circle map $f\colon F \to \Delta$ ($\Delta$=closed disc,
here!) over the origin $0$ of the disc to get a surface marked by
a group of $d$ points. The nice feature is that $(F,f^{-1}(0))$
permits one to recover uniquely (up to rotation) the map $f$ (cf.
Lemma~5.2 about unilateral divisors in Gabard 2006
\cite{Gabard_2006}). Taking instead the fibre over the real unit
$1\in \Delta$ gives a surface marked by a group of $d$ distinct
along the boundary. Taking simultaneously the fibre over $0$ and
$1$ gives a surface marked by $d$ points on both the interior and
the border.

So we have 3 natural spaces of marked surfaces living above the
moduli space ${\cal M}={\cal M}_{r,p}$, namely $I^d$ (interior
marking); $B^d$ (bordered marking); and $M^d$ (mixed marking).
Forgetting the markings gives varied arrows descending to $M$. The
map $\tau$ factorizes through all these marked moduli space.

An idea could be to show that the lift of $\tau$ (which is an
embedding especially when we factor through the mixed marking) is
sufficiently horizontal w.r.t. to the fiber bundle projection
afforded by the forgetful map.
Alas, this is not very evident and should of course hold for some
special value of $d$.

Another route to explore is to make a L\"uroth-Clebsch/Hurwitz
type analysis of trying to understand from ramification and
monodromy how one reconstruct the Riemann surface.

\section{Miscellaneous}

\subsection{Moduli counts via dissection in pants
(Klein, Fricke, Nielsen, Fenchel,
etc.)}\label{Nielsen-Fenchel:sec}

[25.11.12] This section presents a well-known argument to count
moduli of Riemann surfaces, which applies both to the closed and
bordered cases. The argument uses a decomposition in pants and the
hyperbolic metric, so differs somewhat from the original arguments
of Riemann 1857 \cite{Riemann_1857} and Klein 1882
\cite{Klein_1882}, respectively.

\begin{theorem}\label{Nielsen-Fenchel:thm}
{\rm (Riemann 1857 \cite{Riemann_1857},
Klein 1882 \cite{Klein_1882}, Teichm\"uller 1939
\cite{Teichmueller_1939})} The closed genus $g$ surface $F_g$
depends on $3g-3$ complex moduli or $6g-6$ real moduli, while
compact bordered Riemann surfaces $F_{r,p}$ with $r$ contours and
$p$ handles depend upon $3g-3$ real moduli, where $g=(r-1)+2p$ is
the genus of the double.
\end{theorem}

\begin{proof}
First consider the closed case. Introduce on $F_g$ a uniformizing
metric of curvature $K\equiv -1$ and choose a decomposition in
{\it pants\/} (alias {\it trinion\/} by M\"obius 1860/63
\cite{Moebius_1863}). Each pant is a bordered surface with $3$
contours and of genus $p=0$. The conformal structure is
unambiguously defined by the lengths of the contours, plus some
twisting parameters (rotation like) permissible at the junctures
of pants. Looking at the left-hand side of Fig.\,\ref{Pants:fig},
we count $(g-2)$ shaded pants each contributing for 3 lengths, and
one must add one loop on the top and two on the bottom part of the
figure. We arrive at a total of
$$
1+3(g-2)+2=3g-3
$$
many loops. Since each such loop is a juncture we add as many
twisting parameters to  get finally the dependence upon
$$
2(3g-3)=6g-6
$$
real moduli.
\begin{figure}[h]
\vskip-0.0cm\penalty0
    \hskip-0.0cm\penalty0\epsfig{figure=Pants.eps,width=122mm}
  \caption{\label{Pants:fig}%
  Dissection in pants to count moduli}
\vskip-5pt\penalty0
\end{figure}

In the bordered setting, we proceed similarly by looking at a
pants decomposition of the bordered surface $F_{r,p}$ as depicted
on the right-hand side of Fig.\,\ref{Pants:fig}. Counting from the
top to the bottom we get
$$
1+3(p-1)+(r-1)+r=3p+2r-3
$$
many loops. Each of these loops is twistable by a parameter,
except the $r$ boundary loops which have no companion loops. So we
get $3p+r-3$ additional parameters, hence a total of
$$
6p+3r-6
$$
real moduli. On the other hand the genus of the double of
$F_{r,p}$ is $g=(r-1)+2p$, so that the above quantity is nothing
but $3g-3$. This completes the proof.

Of course a more algebro-geometric count do as well the job while
using the reality paradigm of the Galois-Riemann Verschmelzung.
More concretely inside the complex moduli space one defines an
antiholomorphic involution, and the moduli of ``real surfaces''
appears as the real (=fixed-point) locus of that involution so has
half dimensionality. Such an argument has the advantage of
encompassing
directly the diasymmetric case, which leads to
non-orientable Klein surfaces. For more details, cf. Klein 1882
\cite{Klein_1882}, Teichm\"uller 1939 \cite{Teichmueller_1939},
Earle 1971 \cite{Earle_1971-On-the-moduli}, Sepp\"ala 1978
\cite{Seppala_1978-Teich-spaces-of-Klein-surfaces}, Huisman, etc.
\end{proof}

\section{Part II: Hilbert's 16th}
\label{Hilbert's16th-PartII:sec}

\subsection{General overview}

[26.03.13] As announced in the introduction, we enter now in the
second part to our text dealing with Hilbert's 16th problem. The
switch from Ahlfors to Hilbert's 16th flashed us when reading in
more details Rohlin's work of 1978 (cf. the next
Sec.\,\ref{Klein-Rohlin-conj:sec}). Since our assimilation of the
material evolved in slow organical mode (with several mistakes of
ours), it seems worth summarizing which waters were
investigated and what seems to be urgent open problems in the
field. This section should thus replace the reading of all the
sequel which needs severe reorganization at several places.
Further what we understood is still miles away of the fine
jewellery reached by Russian scholars in this field, but to defend
our messy text we also feel that the philosophy \`a la Ahlfors or
Rohlin has not yet been fully exploited, nor elucidated.

First, Hilbert's 16th includes the topological classification of
real algebraic (smooth=non-singular) curves in $\RR P^2$. In its
original formulation the critical degree was $m=6$ (sextics),
where Hilbert's intuition produced both the best (the Ansatz that
an $M$-curve\footnote{Since Petrovskii 1938 \cite{Petrowsky_1938}
it is customary to call an $M$-curve, any curve realizing
Harnack's bound $r\le g+1$ of 1876.} cannot have all its $11$
ovals unnested) as well as a misconception that persisted through
several decades, until being refuted through Gudkov's seminal 1969
construction of the curve $\frac{5}{1}5$ (5 ovals enveloped in a
larger oval, plus 5 ovals outside). This was a big surprise as
Hilbert expected that $M$-curves appear only along the scheme
$\frac{1}{1}9$ discovered by Harnack 1876, and the one constructed
by himself $\frac{9}{1}1$ in Hilbert 1891 (compare the top-row of
Fig.\,\ref{Gudkov-Table3:fig}). The Gudkov symbol $\frac{x}{1}y$
encodes a distribution of ovals where $x$ ovals are directly
nested in one oval, while $y$ unnested ovals are lying outside
(compare again Fig.\,\ref{Gudkov-Table3:fig} if necessary).
Petrovskii's own scepticism about the unexpected twist of Gudkov's
solution, launched the work of Arnold 1971, and Rohlin 1972, where
Gudkov hypothesis  $\chi \equiv_8 k^2$ went verified through
revolutionary insights on the ``complexification''. Here $\chi$
always denotes the Euler characteristic of Ragsdale's orientable
membrane bounding the curve from ``inside'', while $k=m/2$ is the
semi-degree of a curve of even degree $m=2k$. The modern era of
real algebraic geometry was launched, characterized by deep
interconnections with 4D-differential topology (Rohlin's early
work on spin 4-manifolds, etc.)

What has this topic to do at all with Ahlfors maps? To say the
least very few factual links have been tied up presently, but we
can dream of a big connection. The sequel is our attempt to
enhance the r\^ole which Ahlfors theory could play in Hilbert's
16th. We should warn the reader
that our viewpoint is much partisan (biased by what produced such
masters as Ahlfors and Rohlin) and it may well be the case that
the real mathematical terrain is not as plastic and smooth  as the
expressed in the next lines. First, we should stress that there is
no anachronism in expecting such a connection with Hilbert since
(modulo technical details) the quintessence of Ahlfors theory
truly goes back to the Riemann-Schottky-Klein era (resp.
1851/57--1875/77--1876/82), which is much prior to Hilbert
(1862--1943),
and a fortiori to Hilbert's geometrical period ca. 1891---when he
left Algebra, Invariant theory, Number theory---to move in the
softer realms of geometry, calculus of variations, or
``functional'' analysis, especially Dirichlet, Fredholm, etc.

To be honest, our connection was already envisioned by Rohlin 1978
\cite{Rohlin_1978} (apparently completely unaware of Ahlfors work,
as we were ourselves  ca. 2001 when rediscovering the result
independently), but who also used implicitly what we call {\it
total reality\/} as a tool detecting the dividing character of
curves. More strikingly, in a genius stroke without any
antecedents, Rohlin asserts a phenomenon of total reality for
certain $(M-2)$-sextics
explaining a posteriori (nearly) all prohibitions of Gudkov's
table of periodic elements (=Fig.\,\ref{Gudkov-Table3:fig}). The
latter table affords nothing less than the complete solution to
Hilbert's 16th by way of a curious pyramidal structure encoding
all possible distributions of ovals realized by algebraic curves
of degree 6 with real coefficients. Rohlin's (unproved)
synthetical assertion stayed dormant for more than 3 decades until
Le~Touz\'e 2013
\cite{Fiedler-Le-Touzé_2013-Totally-real-pencils-Cubics} recently
managed to establish a slightly weaker form thereof. We can thus
now feel confident in expecting that Ahlfors theory will have to
play some major r\^ole  in the future destiny of Hilbert's 16th,
i.e. for curves of degree $m\ge 8$. (The case $m=8$ looks nearly
settled if one is expert enough in the field and willing to
sacrifice a long period of his time to assemble many bits of
knowledge scattered through the literature.) Hilbert's problem
(like any existential puzzle) splits naturally into constructions
versus prohibitions. Now the r\^ole of Ahlfors could be as
follows. If one has a distribution of ovals (\`a la Hilbert) such
that all curves representing it are
dividing(=type~I=orthosymmetric) in the sense of Klein (what
Rohlin calls {\it a scheme of type~I}\/), then it seems a
reasonable folly to expect the phenomenon of {\it total
reality\/}, namely  existence of a pencil of ``adjoint'' curves
cutting only real points on the given curve. At least Ahlfors
theorem implies
no conformal obstruction to the scenario.

Incidentally, it should be no surprise
that both Ahlfors 1950 and Rohlin's maximality claim (1978) refers
back to a common denominator, namely works of Felix Klein. In
Ahlfors' case this is indirect since reference is more readily
confessed to Schottky's results somewhat prior to Klein's (but
also more schlichtartig than Klein's). Via Teichm\"uller 1941
\cite{Teichmueller_1941} some return to Klein is implicit though
poorly cross-referenced. In Rohlin's case the analogy with Klein
is inherent though disputed in Viro's survey 1986
\cite{Viro_1986/86-Progress} via Marin's assessment of Klein's
assertion that curves of type~I cannot gain an oval by crossing
the discriminant. Apart from those details it is evident that
Klein (and before him Riemann) gave the impulse for all what
followed, and the fusion awaited upon is probably merely a matter
of reunifying the original conception of Riemann-Klein before it
diverged into pure conformal geometry (Schwarz, Schottky, Klein,
Koebe, Bieberbach, Gr\"otzsch, Ahlfors, Grunsky, Teichm\"uller,
Ahlfors again) versus plane curves in Hilbert's 16th (Harnack,
Hilbert, Ragsdale, Rohn, Brusotti, Petrovskii, Gudkov, Arnold,
Rohlin). One can wonder how much knowledge went lost just through
older generations passing away and how much time consuming it will
be for us to revive old wisdoms that are probably the key to most
of our naive questions.

Our 1st fundamental problem is to decide if Ahlfors theorem
particularized to the setting of plane curves implies  existence
of such a total pencil. I personally always thought this being a
triviality (see optionally Gabard 2004 \cite[p.\,7]{Gabard_2004}),
but recently Marin warned me that life might not be so easy (cf.
letter in Sec.\,\ref{e-mail-Viro:sec}). Alas, meanwhile I forgot
nearly all the little I ever knew about the foundations of
algebraic geometry, so that what I thought to be trivial is now
floating in some suspense (``ombre propice'' as would say Thom).

Even if not true (or rather implementable), synthetical procedures
\`a la Rohlin-Le~Touz\'e (=RLT) could redeliver the phenomenon of
total reality {\it ab ovo\/} (independently of Ahlfors conformal
geometry). This seems to require a vertiginous and lengthy
verification process climbing ad infinitum. In more gently slope,
one can expect a gradual
propagation of RLT from degree 6 to 8, and so on, that could be
relevant to detect new prohibitions in Hilbert's 16th.

Why so optimistic? As exemplified by the case of sextics $m=6$
(much influential upon Gudkov, Arnold, Rohlin, etc. and as
demonstrated by those smart guys fairly typical of the general
case $m=2k$) it is likely that a scheme of type~I is totally real
under a suitable pencil (or viceversa) and this should in turn
imply the scheme being maximal in the hierarchy of all schemes.
This produces prohibitions in Hilbert's 16th, which a priori could
be new, and governed by an uniform paradigm.

So a 2nd fundamental problem is to decide if Rohlin's maximality
conjecture (RMC) positing that ``type~I implies maximal'' is true.
At first sight, it seems that a positive solution to the 1st
problem implies this as a byproduct, but  there seems to be severe
obstacles in completing the programme. For explicitness, it is
worth sketching the (naive, uncomplete) argument. Given a scheme
of type~I, there is by Ahlfors a total pencil, which cuts only
real points on the curves. Hence the curve is already saturated,
and cannot be enlarged by adding an additional oval without
violating B\'ezout. The difficulty however is that the enlargement
is not a priori involving the same (or even a nearby curve)
augmented by some other ovals, but can be a priori
very distant of the original curve. ({\it Added in proof}
[13.03.13], for a loose strategy using isotopies, cf.
Sec.\,\ref{RMC-via-Mangler:sec}].)

However a 3rd route is that whenever we encounter a synthetic
phenomenon of
total reality \`a la Rohlin-Le~Touz\'e looks (akin to a
concretization of Ahlfors abstract theorem within the Plato cavern
of Hilbert's 16th  involving only plane curves), then it seems
evident (via B\'ezout-saturation) that Rohlin's maximality
conjecture will hold true for this specific scheme. Again the
proof is not easy to formalize, but it is perhaps realist to
expect a positive solution in the case of curves of degree 6, and
hopefully somewhat higher as to produce new truths.

This brings us to the 4th problem. How to extend Rohlin's total
reality claim to high-degree curves $m> 6$. Is there any algorithm
telling one where to assign basepoints in order to assure total
reality of the corresponding pencil? One could dream that this can
be done from the sole knowledge of Rohlin's complex orientations,
cf. optionally Sec.\,\ref{Galton-brett:sec}.

For plane $M$-curves of degree $m$, we prove below a basic
Theorem~\ref{total-reality-of-plane-M-curves:thm} stipulating a
total pencil swept by curves of order $(m-2)$. This merely
traduces  the so-called Bieberbach-Grunsky theorem, which (apart
from phraseological details like Dirichlet's illness) truly
belongs to Riemann 1857, Schottky 1875--77, Enriques-Chisini 1915,
and only then Bieberbach 1925, Grunsky 1937, Wirtinger 1942
\cite{Wirtinger_1942}, etc. A structural asymmetry appears: while
$M$-curves are crudely-put reputed hardest-to-construct within
Hilbert's 16th, their conformal geometry is most trivial, due to
the planar=schlichtartig character of the half of the curve. Total
reality is simplest to ensure in the $M$-setting, just because it
is like having one train on each track, hence no risk of
collision. Precisely, the trick is just to choose one point on
each oval getting so a group of $g+1$ points which moves (by
Riemann-Roch or Abel), and total reality is automatically granted
(cf. Lemma~\ref{Enriques-Chisini:lemma} for more details). Making
this abstract argument concrete proves the theorem.

Can we extend this to  non-maximal curves? The risk is then an
overpopulation of $g+1$ points scattered
on $r<g+1$ ovals, hence 2 of them are forced to cohabit on the
same oval (pigeonhole principle due to Dirichlet apparently)
exposing us to a possible collision jeopardizing total reality! So
what is demanded is controlling a dextrogyration of points when
moving along linear equivalence. Ahlfors 1950 \cite{Ahlfors_1950}
or Gabard 2006 \cite{Gabard_2006} affords basic tricks to achieve
dextrogyration in the abstract setting. Can we transplant them
directly inside the Plato cavern of plane curves, as we  just
managed  to do for $M$-curves? As yet we never succeeded, but this
should not discourage
more serious attempts.

If we think more concretely \`a la Rohlin-Le~Touz\'e (or perhaps
to go back earlier in history \`a la Brill-Noether), numerological
reasons make evident that $(M-2)$-curves of type~I (or even
schemes of type~I) and degree $m$ will have their total reality
exhibited by a pencil of curves of degree $(m-3)$. Evidence is
given later in this text, but readily follows by analogy with
Rohlin-Le~Touz\'e and a simple constants count (cf.
Remark~\ref{M-2-curve-degree-like-Gabard:rem}, which is just the
end product of numerological coincidences observed for $m=6,8,10,
\dots$). Remind perhaps at this stage old Italian works
(recognized as possible competitors to Ahlfors 1950), like Matildi
1945/48 \cite{Matildi_1945/48}, Andreotti 1950
\cite{Andreotti_1950}. Those could
already anticipate our present desideratum. Since already Rohlin's
proof (which is lost) and that of Le~Touz\'e (2013
\cite{Fiedler-Le-Touzé_2013-Totally-real-pencils-Cubics})  is
quite delicate (or the sequel of this text which contains ca. 30
pages of unsuccessful attempts to prove Rohlin's original claim),
we are not claiming that total reality will be easy to prove in
full generality but perhaps for degree $8, 10$ this remains
manageable (at least within the next 4 decades). By experience we
are accustomed in the field to slow progresses (remind Hilbert,
Gudkov, etc.), and it is quite unlikely (but not impossible) that
a new Abel or Riemann will crack the full puzzle in a single
stroke.

Some part of our text tries to take the census of all such schemes
of type~I in degree 8. Alas our optical faculties tend to be much
more limited than those of aliens,  like insects with 8 eyes
looking at their preys via  pencil of cubics, chameleons with
mobile ocular systems, or any sort of creature with 19 eyes (when
it comes to look at the world through a pencil of quintics,
\dots), and generally, $M-3$ basepoints (i.e. 19 when $m=8$, $34$
for $m=10$). Accordingly, we are presently (and probably for the
rest of our life) confined to deduce total reality not from
optical skills but via boring arithmetics, namely the subliminal
(Rohlin?)-Kharlamov-Marin congruence $\chi\equiv k^2+4 \pmod 8$,
which forces type~I under this ``shifted'' Gudkov-style congruence
mod 8.

This is a crucial weapon (whose proof we have not yet studied in
full details). This harpoon detects for us a menagerie of schemes
of type~I, all possibly subsumed to total reality, and
conjecturally (via RMC=Rohlin's maximality conjecture) acting
prohibitively upon all schemes(=distribution of ovals) pretending
to enlarge the given one. Such schemes crystallize therefore
B\'ezout-extremal (or saturated) shapes of Hilbert's 16th, which
as would say Klein cannot develop further without exploding the
latent degree.

Call an RKM-scheme any scheme satisfying the RKM-congruence
$\chi\equiv_8 k^2+4$. It is not clear to the writer,
and
the experts were a bit silent
on this aspect as yet, if conversely any $(M-2)$-scheme of type~I
is forced to respect the RKM-congruence. This deserves perhaps to
be clarified at the occasion. [{\it Added in proof} [13.04.13] An
answer is probably implicit in Rohlin 1978, Art.\,3.5, on p.\,93
(extremal property of Zvonilov-Wilson).]

With some sloppiness, we arrive at some big picture along the
following philosophy (in our opinion fairly implicit in Rohlin
1978). Any scheme of type~I is detected:

$\bullet$ either trivially because it is an $M$-scheme whose total
reality is exhibited \`a la Bieberbach-Grunsky (yet no direct
impact upon Hilbert's 16th by virtue  of Harnack's bound (1876),
or more simply its intrinsic variant $r\le g+1$ due to Klein 1876
proved via retrosections \`a la Riemann), [but some indirect
impact by using satellites!! (13.04.13)]

$\bullet$ or it is an $(M-2)$-scheme verifying the RKM-congruence,
in which case total reality is flashed by a pencil of
$(m-3)$-tics.

$\bullet$ or finally it arises as  ``satellite'' of a scheme of
lower degree dividing the given degree.

The idea of satellites arises simply by noting (or expecting) that
total reality propagates when the curve is doubled,  tripled and
so on, by replicating several copies of the curve  within a
tube-neighborhood of it (\ref{satellite-total-reality:sec}). For
the conic with a single oval (unifolium) this just leads to the
series of deep nests total under a pencil of lines, while for a
quartic with 4 ovals (quadrifolium in the jargon of Zeuthen 1874
\cite{Zeuthen_1874} who inspired much Klein 1876) this leads to
the series of curves in degrees multiple of $4$,  totally real
under a pencil of conics. The case of degree 8 is explicitly
mentioned in Rohlin 1978, being just the double of the
quadrifolium. We expect that satellites do extend to curves of odd
degrees (\ref{Satellite-odd-degree:sec}), yielding some
interesting prohibitions on schemes of degree 10 when applied to
an $M$-quintic doubled. Likewise doubling the Rohlin-Le~Touz\'e
sextic $\frac{6}{1}2$ (or its mirror $\frac{2}{1}6$) gives a
scheme of degree 12 which should be maximal (hence killing all
extensions of it).

The general philosophy is now clear. Total reality (basically due
to  Ahlfors 1950, though Teichm\"uller 1941 ascribes it to Klein
directly) acts as an upper-bound on the complexity of Hilbert's
16th problem, by killing all distributions of ovals adventuring
above one totally flashed by a pencil. In substance everything
boils down to a phenomenon of B\'ezout-saturation, with in the
background of the scene an extension of the Riemann mapping
theorem to surfaces of higher topological structures (so-called
Ahlfors maps).

This looks a fundamental truth (or philosophy?) since it seems
robust, and implementable when the flashing is as explicit as
Rohlin-Le~Touz\'e's as opposed to the abstract nonsense of
Ahlfors. If skeptical, the just predicted maximality of satellites
in degree 10 and 12 should be tested against highbrow methods of
constructions of the modern era (Viro-Itenberg). If the
Ahlfors-Rohlin philosophy resists the shock against this
structural test,
then some experimental evidence is gained that the Ahlfors-Rohlin
Verschmelzung is a deep reality governing a substantial part of
Hilbert's 16th at the universal scale (all degrees). If not, then
the whole story of the 16th problem could be even more chaotical
and unruly than it presently is, i.e. just a combinatorial mess
only worth deserving the attention of computing machines. Of
course the latter are quite likely to show us hidden patterns of
symmetries, maximality, etc. that were not yet appreciated due to
a lack of experimental data.

More pragmatically, it must feasible to inspect if in degree 8,
the Ahlfors-Rohlin scenario of total reality and the allied
extremal principle of saturated schemes is compatible with factual
data, and optimistically even able to preclude schemes that were
not yet prohibited. Alas, we are not expert enough in the field to
tell an answer, but peoples like Viro, Fiedler, Korchagin,
Orevkov, Le~Touz\'e, must already have a clear-cut vision along
this idea. In degree 6 it is clear that the saturation principle
of Rohlin is entirely covered (or re-explain) by the congruences
mod 8 due to Gudkov and Gudkov-Krakhnov-Kharlamov, but it is not
clear to me if the same subordination holds true in degree 8
(maybe in general). On the other hand if the RKM-congruence fails
to detect a type~I scheme, then there could be some sporadic
phenomenon of total reality explaining it, and this would be a new
source of saturation (perhaps not covered by the congruences
mod~8).

This is the main-body of our quest, but during the trip we
went sidetracked to other  connected topics. Here are some aspects
perhaps worth putting in evidence centering around the theme of
rigid-isotopy, and the allied contraction principles where the
end-point of the path is permitted to touch the discriminant
(parameterizing all singular curves).

Total reality takes its simplest incarnation for the deep nests
swept out by a total pencil of lines. A theorem by Nuij 1968
\cite{Nuij_1968} (later revisited by Dubrovin 1983
\cite{Dubrovin_1983/85}) states that such deep nests are rigid in
the sense than one can pass from any 2 curves representing it by a
continuous deformation of the coefficients without encountering
any singular curve during the deformation. Such large deformation
pertains to what is called {\it rigid-isotopies\/}, which actually
refines Hilbert's 16th problem. This topic always
attracted geometers even prior to Hilbert's era, e.g. Schl\"afli
(apparently known for having the most massive human brain ever
weighted with ca. 1.936 kg for only 157 cm of body height), or
Zeuthen and Klein adding several
contributions regarding curves and surfaces of low orders (quartic
and cubics resp.). For quartic curves Klein established (1876
\cite{Klein_1876_Verlauf}) that the rigid-isotopy class is  fully
determined by the real scheme already.

It seems natural to ask if Nuij's rigidity result (for deep nests)
has equivalents whenever total reality holds true. Alas this fails
by the Marin-Fiedler locking technique which refutes this Ansatz
for $M$-curves of degree 7 (cf. Fig.\,\ref{Marin:fig}, for a
hopefully lucid exposition of Marin's trick). Despite this
disruption of the naive scenario, it seems to us
 likely that rigidity holds true for satellites of the
 quadrifolium. We confess however to have not yet studied
Nuij's proof, nor do we know (a fortiori) if his proof extends
mutatis mutandis.

Though rigid-isotopy  merely involves the $\pi_0$ (=nullest
homotopy group measuring the arcwise-connected components) of the
space of curves excised along the discriminant, very little is
known on such problems. A naive conjecture of us---based
essentially on the failure of the Marin-Fiedler locking device,
plus the fact (subsequent to Rohlin's formula) that curves of
type~I have $r\ge m/2$ ovals (also valid if $m$ is
odd)---postulates that curves with few ovals are necessarily {\it
rigid}, i.e., are unambiguously determined up to a large
deformation by their sole real schemes. Precisely this could hold
true as soon as the curve has strictly less ovals than $DEEP+2$,
where $DEEP:=\Delta:=[(m+1)/2]$ is the number of circuits of the
deep nest of degree $m$. The intuition is simply that by Rohlin's
formula (\ref{Rohlin-formula:thm}) this is the first dividing
scheme encountered (as $r$ the number of ovals increases), and two
units above this ($r=\Delta+2$) it is a simple matter to exhibit a
scheme of indefinite type (Rohlin's jargon to say that there is
curves of both types I.\,vs\,II realizing a prescribed
configuration of ovals).

This conjecture (called LARS, for
low-altitude-rigidity-conjecture) is merely a cavalier extension
of:

(1) Nuij's theorem of 1968, which is not specific to curves (but
valid for  algebraic hypersurfaces, where there is an evident
notion of deep nest via concentric spheres, plus an eventual
pseudo-plane).

(2) Nikulin's rigid classification of sextics in 1979/80
\cite{Nikulin_1979/80} implying the case $m=6$ of our LARS, and
telling much more namely the fact that the real scheme (as
tabulated on Gudkov's table) enhanced by the data of
Klein-Rohlin's types affords a complete system of invariants under
rigid-isotopy. Hence for $m=6$, Nikulin is stronger than LARS as
it prompts rigidities at all altitudes. However as soon as $m\ge
7$ this is foiled (cf. again Marin's
example=Fig.\,\ref{Marin:fig}).

(3) A unofficial conjecture of Rohlin (reported in a Viro letter
in Sec.\,\ref{e-mail-Viro:sec}) that curves of odd degrees with a
single (pseudoline) component are rigid-isotopic, cf. also Viro
2008 \cite{Viro_2008-From-the-16th-Hilb-to-tropical}. By analogy,
curves of even degrees with a single oval could  be
rigid-isotopic. Those questions are settled in degrees $m\le 6$
($m=4$ Klein 1876 \cite{Klein_1876_Verlauf}, $m=5$ Kharlamov 1981
\cite{Kharlamov_1981/81}, $m=6$ Nikulin 1979/80
\cite{Nikulin_1979/80}), but  still resist in degrees $\ge 7$.
Hence our conjecture LARS appears very presumptuous, and it may be
a more reasonable challenge trying to disprove it. Alas the
locking method of Fiedler-Marin looks (as far as we experimented
in the sequel) quite impuissant to destroy LARS.

What techniques could be used to prove LARS or more modest
rigidity conjectures? Our naive idea is that geometric flows
(amounting to look at orthogonal trajectories of suitable
functionals like calculating the length or area of ovals) could
prove this and related results of rigid-isotopies. This would give
some intrusion of differential-geometric methods in problems of
rigid-isotopies, a priori of a purely algebraic nature. Presently
we were never able to complete any serious proof along this way,
but our text contains ca. 20 pages of (dubious) trials along such
lines.
Viro's survey 2008 \cite{Viro_2008-From-the-16th-Hilb-to-tropical}
also contains a brief desideratum to know more about geometric
properties of curves, and this could evidently pertain to
rigid-isotopies, in a way perhaps reminiscent of the \"Olfleck of
H.\,A. Schwarz (where the Riemann mapping theorem is visualized by
an oil-flake restoring to the circular shape), or the eclectic
Ricci flow of Yau-Hamilton-Perelman, where a similar phenomenology
appears in the abstract Riemannian setting (convergence to the
round metric, with the well-known bonus about Poincar\'e's
conjecture).

Affiliated to those rigid-isotopy questions there is a conjecture
of Itenberg-Viro (cf. Viro's preface of the volume containing
Itenberg 1994 \cite{Itenberg_1994}) to the effect that some empty
oval of any curve can always be shrunk toward a solitary node.
[{\it Added in proof} [13.04.13].---Similar (but more vague) ideas
are actually ubiquitous in Klein, e.g. 1892
\cite{Klein_1892_Realitaet}.] This is still wide open, but
Itenberg's article just cited establishes the case $m=6$. Again
one could hope that the flow minimizing the length of an oval
could achieve such a contraction. Inspired by this conjecture we
advanced a strengthened version CCC(=collective contraction
conjecture) saying that all empty ovals can simultaneously
contract  toward solitary nodes. (This is like a perfect landing
in Flight-Simulator v.18.5 with aircraft Antonov~72 having its 94
wheels touching the ground simultaneously!) If this (unlikely)
miracle is true,  one gets e.g. a 2-seconds proof of Hilbert's
Ansatz of the non-existence of an $M$-sextic without nesting by
reduction to B\'ezout. Of course this is also more hygienically
derived from Rohlin's formula (or Arnold's congruence mod 4),
which involves softer homological intersection theory \`a la
Poincar\'e-Lefschetz, etc. We were not as yet able to disprove our
strong CCC-version of Itenberg-Viro. Its real impact being still
obscure we did not pursued this issue in any serious fashion.

A philosophical consequence, of large deformations is that they
should (like total reality) act prohibitory, whereas small
perturbations are classically  exerting their swings at the
constructive level (Harnack 1876, Hilbert 1891, Brusotti 1914/21,
Gudkov 1969/72, Viro 1980, Itenberg 1993, etc.) Of course a
clear-cut realm of  where to corrupt CCC could be a dividing curve
without nesting, for those could after {\it strangulation\/} be
split into two complex-conjugate halves intersecting in as many
points as there were ovals initially. Alas either Thom
(\ref{Thom-Ragsdale:thm}) or better Rohlin's formula
(\ref{Rohlin-formula:thm}) forces such a curve to have $\chi=r\le
k^2$ resp. exactly $r=k^2$ ovals, hence we fail to corrupt
Monsieur \'Etienne B\'ezout. It looks so quite challenging to kill
CCC, albeit its truth looks very fragile, as it incarnates an
extreme flexibility of algebraic objects reputed ``rigid'' in the
large. Yet it should be remembered (though at some more local viz.
regional scale) that Brusotti's theorem gives via
Riemann-Roch-Brill-Noether-Severi a remarkable flexibility of
algebraic curves (independence of the smoothing of nodes). So one
should not be surprised at last, if sometimes algebraic curves
appear more plastic than expected a priori. However, as we shall
soon discuss, Shustin disproved (in degree 8) a flexibility
conjecture of Klein (1876) that nondividing curves can always
acquire a solitary node through continuous variation of the
coefficients ({\it champagne bubble phenomenon\/}). In slight
contrast, building over the previously cited works of Nikulin
1979, Itenberg 1994, and the whole diagrammatic of the
Gudkov-Rohlin table (1969--78), we think that Klein's intuition of
champagne bubbling is correct in degree 6 (cf.
Prop.\,\ref{Klein-vache-deg-6:prop}). The philosophical impact of
Shustin's disproof of Klein (though his aim was refuting a related
assertion of Rohlin) is that we cannot expect to have solely
topological obstructions regulating  large algebraic deformations.

A last theme involves the impact of Thom's conjecture (meanwhile
Kron\-heimer-Mrowka's theorem 1994) upon Hilbert's 16th
(Sec.\,\ref{Thom:sec}). A classical trick (called the Arnold
surface) is to fill Klein's half of a dividing curve $C_m$ (of
degree $m=2k$) by the real Ragsdale orientable membrane bounding
the curve from inside. This gives a homology class of half-degree
$k=m/2$, smoothly represented (after rounding corners, if
necessary). This object looks ideally suited to an application of
Thom's genus estimate. Taking for granted orientability of the
Arnold surface, we found an erroneous estimate $\chi\le k^2$ (for
all dividing plane curves of degree $2k$ (cf.
(\ref{Thom-Ragsdale:thm})). Albeit wrong in general (as Fiedler
kindly pointed out to us) it holds in special cases when all
(primitive) pairs of ovals are positive in the sense of Rohlin,
i.e. when complex and real orientations match together. Real vs.
complex orientations may even disagree yet along pieces not
connected by the Ragsdale membrane (cf.
Lemma~\ref{Arnold-surface-orientable-iff-oddly-charged:lem}). If
optimistic Thom or even Rohlin's formula gives a way to attack the
(still open) Ragsdale's conjecture for $M$-curves, which amounts
to $\vert \chi \vert \le k^2$
(\ref{Thom-implies-one-half-of-Ragsdale:lem}). Is the Arnold
surface (=Klein's half glued with the Ragsdale membrane) of an
$M$-curve always orientable? If yes, then the proof of our
(erroneous) Theorem~\ref{Thom-Ragsdale:thm} implies Ragsdale's
conjecture via Thom's estimate on the genus. Unfortunately,
Arnold's surface is nonorientable already for Hilbert's
$M$-sextic, cf.
Lemma~\ref{disproof-orientability-Arnold-M-curve:lem}.

Maybe the theorem \`a la Bieberbach-Grunsky specialized to plane
$M$-curves (i.e. our
Theorem~\ref{total-reality-of-plane-M-curves:thm}) could give (via
dextrogyration\footnote{This concept is not really meaningful for
$M$-curves.}) enough control on complex orientations of $M$-curves
as to imply Ragsdale, either via Thom or directly via Rohlin's
formula $2(\pi-\eta)=r-k^2$ (where $\pi:=\Pi^+$, $\eta:=\Pi^-$ to
abridge notation). (This amounts then to check that $\pi-\eta\le
n$, the number of negative ovals.) This admittedly looks naive,
but we cannot exclude such a coarse strategy for the moment.

More modestly, it may be noted that filling Klein's half with
Ragsdale's membrane of an $M$-sextic without nesting reduces
Hilbert's nesting Ansatz to the ``baby'' case of Thom for homology
classes of degree 3, acquitted by Kervaire-Milnor 1961
\cite{Kervaire-Milnor_1961} building upon Rohlin's early work 1951
on spin $4$-manifolds.

As said, our erroneous estimate $\chi\le k^2$ was corrected by
Fiedler in a series of letters where  he learned us the Petrovskii
estimates on $\chi=p-n$, and  Arnold's strong avatars thereof
involving hyperbolic ovals. This is again  closely connected to
the Ragsdale conjecture, which is still a {\it pi\`ece de
resistance} in the case of $M$-curves. Moreover though our
estimate $\chi\le k^2$ was erroneously founded it turned out to be
quite difficult to find an explicit counterexample. At least we
failed via classical methods (cf.
Figs.\,\ref{HilbGab1:fig}--\ref{HilbGab4:fig}), which rather
inclined to think that $\chi\le k^2$ was sharp if true at all.
Namely using Hilbert's construction we find an infinite series of
$M$-curves or $(M-2)$-curves such that $\chi=k^2$ exactly, but
failed to beat $k^2$. We presume this is exactly the sort of
experiments that led Ragsdale to her conjecture. However the story
does not finish here, and the big surprise arrives now.

It is notorious that a marvellous construction of Viro-Itenberg
(patchwork and $T$-construction, cf. Fig.\,\ref{Itenberg:fig})
killed the Ragsdale conjecture in degree 10 (even in its relaxed
shape of Petrovskii), yet leaving intact  the $M$-curves case. The
Itenberg-Viro construction supplied us with the apparently
simplest counterexample
 to our erroneous
estimate $\chi \le k^2$ (for type~I curves of degree
$2k$).
It produces namely an $(M-2)$-curves of degree 10 with $\chi=29
\nleqslant k^2=25$, hence necessarily dividing by the
RKM-congruence $\chi\equiv_8 k^2+4$. This was the fatal stroke
(coup de gr\^ace) against our estimate $\chi \le k^2$, which is
quite robust as it seems  incorruptible via Harnack-Hilbert and
challenging to refute in the $M$-case.

Last but not least, there is a disproof due to Shustin 1985
\cite{Shustin_1985/85-ctrexpls-to-a-conj-of-Rohlin} of one-side of
Rohlin's maximality conjecture namely ``type~I$\Leftarrow$
maximal''. This disproof is not so dramatic for Rohlin's
prohibitive programme which uses rather the converse (still
hypothetical) implication ``type~I$\buildrel{?}\over\Rightarrow$
maximal''.  Shustin's note looks historically pivotal as it kills
the second part of Klein's intuition, pertaining to large
deformations of nondividing curves as always admitting the
apparition of a champagne bubble created by crossing the
discriminant through a solitary node. Due to its extreme concision
we had  first not understood Shustin's argument (and {\it
unduly\/}
mistrusted his result for a while). Finally, we
understood its logic, but confess to have not yet assimilated all
the results required to complete its proof. It suffices to say
that Shustin's work exploits Viro's construction on the one side,
and also advanced B\'ezout-style obstructions due to Fiedler and
extended by Viro. Some details perhaps assisting beginners to
grasp the structure of Shustin's proof are to be found in
Sec.\,\ref{Shustin-understood:sec}.

This a brief summary of the territories we managed to explore in
ca. 3 months of investigation. Besides our text may have some
didactic value on the following aspects.

(1).---We give a self-contained account of Gudkov solution to
Hilbert's 16th problem in degree $m=6$, by exposing the original
constructions of  Harnack, Hilbert and Gudkov. Those issues are
well-known and described in Gudkov's seminal survey (and at
several other places like A'Campo's Bourbaki survey 1979
\cite{A'Campo_1979}, etc.). Yet not all species are always
accompanied by decent pictures requiring sometimes clever twists
of Harnack's construction (oft messy to implement if one wants to
realize a type given in advance). So we had long hours of trials
with computer-assisted depictions. This can hopefully be of some
use to some nonspecialist readers. Our intention was to reproduce
all (including the infructuous) trials, but that generated
``microfilm'' pictures often too heavy for the purpose of
arXivation.
By the way our microfilm though still readable in pdf-format at
600 dpi resolution will still be hard to contemplate on the
screen. [{\it Added in proof} [13.04.13] This technical problem
was settled by shrinking the size of pictures in the Adobe
software, permitting so to economize much memory space, yet
without altering the optical size of pictures.]

(2).---We give also full details (and a graphical
view=Fig.\,\ref{Gudkov-Table3:fig}) of Rohlin's enhancement of
Gudkov's census of sextics by adding the complex topological
characteristics of Klein (that were much
neglected during the era of Hilbert, Ragsdale, Rohn, Petrovskii,
Gudkov) up to the Arnold-Rohlin revival of the complexification
(which turned to be the conceptual key to explain Gudkov's
experimental phenomenology). This is merely a simple exercise yet
that can be quite time-consuming if one starts from
zero-knowledge. Of course an excellent account of this, differing
form ours only in the minor details, is already given in the
masterpiece Marin 1979 \cite{Marin_1979}. (Our account differs
just in using more primitive configurations of 3 ellipses.)

(3).---We give in Sec.\,\ref{Prohibitions:sec} a reasonably
exhaustive list of classical obstructions, especially a (nearly
complete) proof of Rohlin's formula (\ref{Rohlin-formula:thm}). In
the original source (Rohlin 1974 \cite{Rohlin_1974/75}) this is
not presented in its full generality (only $M$-curves), though the
adaptation to general dividing curves is very minor. This Rohlin's
formula looks extremely fundamental as it appears as the most
universal obstruction that can be derived by nearly abstract
nonsense (i.e. using very little from the assumption of
algebraicity), yet still affording strikingly  precise information
while staying completely elementary. For instance it covers
Hilbert's Ansatz of nesting, and extends it to all degrees $m\ge
6$. It also formally implies the  Arnold congruence mod 4, which
is a weak form of Gudkov hypothesis for $M$-curves, yet an
extension thereof to arbitrary dividing curves. Then there is a
series of avatars of the Gudkov congruence mod 8, that truly
requires more advanced topological tools, essentially in the
spirit of Rohlin 1952 \cite{Rohlin_1952-4-manifolds}. Those more
advanced results are not proved in our text, and we hope to be
able to offer a lucid view on them in the future.
%
%
Hence, to assimilate the marvellous congruences due to Rohlin,
Gudkov-Krakhnov/Kharlamov, Kharlamov-Marin, etc., our reader is
invited to consult the original sources (Rohlin 1972
\cite{Rohlin_1972/72-Proof-of-a-conj-of-Gudkov} (with a gap but
essentially correct and repaired by Marin-Guillou, e.g. 1986
\cite{Guillou-Marin_1986} or Marin 1979 \cite{Marin_1979}), and
also Kharlamov-Viro 1988/91 \cite{Kharlamov-Viro_1988/91}, giving
a synthesized view).

\subsection{Challenging vs. less challenging open problems}
\label{Challenging-open-prob:sec}

[28.03.13] This section summarizes what looks to us major open
questions in the field investigated. The reader is warned that our
list is a mixture of hard Soviet conjectures of longstanding with
newcomers (due to myself), therefore probably much easier to
settle down when not ill-posed. To distinguish among them the
symbol $\bigstar$ marks venerable Russian conjectures, while our
more modest variants are marked by ``$\bullet$''.

$\bigstar$ (R6) Can somebody reconstruct Rohlin's lost proof that
the $(M-2)$-sextics with schemes $\frac{6}{1}2$ or $\frac{2}{1}6$
are totally real under a pencil of cubics assigned to pass through
8 points distributed on (or inside) the empty ovals. This is
nearly solved in Le~Touz\'e 2013
\cite{Fiedler-Le-Touzé_2013-Totally-real-pencils-Cubics}, but she
uses the RKM-congruence (\ref{Kharlamov-Marin-cong:thm}) to infer
a priori the curve being dividing. It could be more natural to
draw dividingness from total reality via a purely synthetical
procedure {\it a priori}. At any rate the conjunction of the
RKM-congruence ($\chi\equiv_8 k^2+4$) with Le~Touz\'e's result
implies that Rohlin's assertion is true. Hence, it should be
already ``safe terrain'' to explore. If much more pessimistic
Rohlin's claim is wrong and then either RKM, or Le~Touz\'e is
false, which is very unlikely.

$\bigstar$ (RLT6$\to$RMC6) Can someone complete the proof that the
Rohlin-Le~Touz\'e phenomenon of total reality (RLT6) prevents all
sextic schemes extending those described in the previous problem
(R6), so as to infer nearly all obstructions of Hilbert's 16th via
the paradigm of total reality (TR) and the allied phenomenon of
B\'ezout-saturation. Cf. the diagrammatic of the Gudkov table
(Fig.\,\ref{Gudkov-Table3:fig}) to appreciate this issue in degree
$m=6$. Of course this problem can be considered as very implicit
in Rohlin 1978 \cite{Rohlin_1978}, but in our opinion not solved
there.

$\bullet$ (RLT$m>6$) How does the Rohlin-Le~Touz\'e phenomenon
described in (R6) above extend to higher degrees $m>6$? Cf.
Sec.\,\ref{total-(M-2)-schemes:sec} for a germ of answer.

$\bullet$ (A50$\to$R78) How valuable is the abstract theory of
Ahlfors to assess Rohlin's vision of total reality? Cf. e.g.
Sec.\,\ref{Esquisse-dun-prog-deja-esquiss:sec} for some scenarios.
In particular is it true that any dividing plane curve admits a
total pencil (i.e. whose real members cut only real points)? If
yes, is it always of degree $\le (m-2)$ when the given curve has
degree $m$? For the case of $M$-curves, cf.
(\ref{total-reality-of-plane-M-curves:thm}) which gives a total
pencil of order $(m-2)$.

$\bullet$ (R78$\to$G13) Is it true as conjectured in our text
(\ref{M-2-curve-degree-like-Gabard:rem}) that any curve belonging
to an $(M-2)$-scheme of type~I and degree $m$ has its total
reality exhibited by a pencil of curves of degree $(m-3)$. Further
what is the exact r\^ole of Riemann, and Brill-Noether adjoint
curves, in this game? Notice still in
(\ref{M-2-curve-degree-like-Gabard:rem}) a strange concomitance
between Rohlin-Le~Touz\'e's role of cubics and Gabard's $r+p$
bound on the gonality. The latter improves Ahlfors by replacing
$g+1$ by the mean-value of Harnack's bound $g+1$ and the number
$r$ of real circuits. All this numerology looks to match too
nicely for this being merely a fortuitous
coincidence. In  particular for quartics, quintics, sextics, etc.
the total reality of $(M-2)$-curves of type~I seems always
exhibited by such a pencil of degree $(m-3)$.

$\bullet$ (RKM$\leftarrow$type~I) The RKM-congruence $\chi\equiv_8
k^2+4 $ detects many $(M-2)$-schemes of type~I, but does it detect
all of them? The answer is yes for $m=6$ (cf. the Gudkov-Rohlin
table Fig.\,\ref{Gudkov-Table3:fig}). Hence $m=8$ is the first
place to look for a counterexample. Assuming there is one, then it
could be that total reality detects type~I schemes at places where
RKM fails. [{\it Added in proof} [13.04.13].---The answer to this
question must be implicit in Rohlin 1978 \cite[p.\,93,
Art.\,3.5]{Rohlin_1978} (and due to Zvonilov-Wilson).]

$\bigstar$ (RMC) Is Rohlin's maximality conjecture true, i.e. all
schemes of type~I are maximal in the hierarchy of all schemes of
fixed degree? Can this be disproved by the Viro-Itenberg
patchwork, as it was possible to refute the converse sense of
Rohlin's conjecture (cf. Shustin's note 1985
\cite{Shustin_1985/85-ctrexpls-to-a-conj-of-Rohlin}).

$\bullet$ (SAT) Are satellites of schemes of type~I still of
type~I? For instance what about the 2nd satellite of the
Rohlin-Le~Touz\'e sextics of point (R6). Can this be disproved via
patchwork? Assuming a positive answer to the first question (even
in a special case) points to potential place where to corrupt RMC.
Personally, we would be much more happy if RMC holds true, as
then, and only then, there is some chance to make a big
Riemann-Hilbert or Ahlfors-Rohlin synthesis.

$\bigstar$ (H8) Complete the solution of Hilbert's 16th in degree
8, and analyze objectively if Rohlin's maximality principle (RMC)
has some things to say in this realm, as it did in degree $m=6$.
In particular does  Rohlin's maximality conjecture (RMC) still
persists in degree $8$. (Some hints are given in Orevkov's letter
in  Sec.\,\ref{Orevkov:sec}.)

$\bullet$ (LARS)  Can someone disprove our low-altitude rigidity
speculation (LARS) positing that a curve with less ovals than 2
units above the deep nest is entirely determined up to large
deformations by its real scheme. Cf. (\ref{LARS:conj}).

$\bullet$ (URS) [02.04.13] The unnested rigidity speculation (URS)
is akin to LARS, and posits that any unnested curve is rigid
provided the number of ovals is not the square of the semi-degree
($r\neq k^2$, and assume $m$ even). Motivation comes from Rohlin's
formula (which forces such curves being of type~II), and the case
$m=6$ which follows from Nikulin. Another (weak) evidence comes
from the fact that the locking technique of Fiedler-Marin seems to
have little grip on such schemes as there is no way to choose a
canonical triangle (moving frame).

$\bigstar$ (OOPS)=(One oval postulation).---In particular what
about the much more modest (than LARS or URS, yet still wide open)
rigidity conjecture for curves having only one component. Are such
unifolium curves rigid as conjectured by Rohlin, Viro, etc. How
useful are geometric flows to do this? Cf.
(\ref{OOPS:one-oval-rigid-isotopic:conj}). Actually Viro ascribes
to Rohlin (cf. e-mail in Sec.\,\ref{e-mail-Viro:sec}) the rigidity
of curves of odd degree with a unique component, but the case of
even degrees looks likewise open. Further it seems evident (at
least for Viro, cf. the same letter) that OOPS is implied by CC,
i.e. the Itenberg-Viro contraction conjecture for empty ovals. By
analogy it seems evident that our CCC (cf. right below) implies
URS. Sketch of proof: contract all ovals simultaneously (which are
all supposed empty) as to reach the connected empty chamber, and
do this twice. Of course when $r=k^2$ the real scheme can be of
both types, and this case has to be ruled out (or optimistically
the type is the sole obstruction to rigid isotopy).

$\bullet$ (CCC)=(Collective contraction conjecture).---Can someone
disprove our strong version CCC (of the Itenberg-Viro contraction
principle for empty ovals) positing a simultaneous and collective
contraction of all empty ovals toward solitary nodes. Cf.
(\ref{CCC:conj}).

$\bigstar$ (CC) [30.03.13] There is a (still open) conjecture of
Itenberg-Viro (cf. Itenberg 1994 \cite{Itenberg_1994} and Viro's
preface of the same volume) positing that {\it any} empty oval of
a real plane curve can be contracted to a point (solitary node).
This is true in degree 6 as proved by Itenberg (\loccit), and
weaker than our (CCC) above (cf. (\ref{CCviaCCC-Brusotti:lem})).

$\bigstar$ (CC vs. TR) A noteworthy consequence of CC is that all
obstructions in degree 6 derived (clumsily) via total reality (TR)
are likewise derived by this contraction principle of Itenberg
(CC6) (modulo knowledge of the RKM-congruence and Klein's Thesis
which is fairly easy to prove since Marin 1988 \cite{Marin_1988}).
The problem is first to decide which method ``total reality''
versus ``contraction'' is more easily implemented in degrees $\ge
8$, while trying to make a comparative study of the prohibitions
resulting from both procedures. In particular one may wonder if
the Itenberg-Viro conjecture implies (formally or not)  Rohlin's
maximality conjecture. Sketch of proof: Take any scheme of type~I,
and a curve enlarging it. Contract an empty oval so as to recover
the initial scheme (note here an obvious difficulty, namely the
additional oval of the extended scheme is not necessarily an empty
one!), and conclude via Klein's Thesis (a curve of type~I cannot
champagne-bubble).

$\bullet$ (Refuting CC via Shustin?) [31.03.13] By the proof of
Prop.~\ref{Klein-vache-deg-6:prop}, we see that the Itenberg
contraction principle combined with the diagrammatic of the
Gudkov-Rohlin table (of all typed-schemes) implies {\it
Klein-vache\/} (KV), i.e. the possibility for diasymmetric curve
to acquire a solitary node and then a new oval ({\it comme surgit
du n\'eant\/}). Now as Klein-vache is disproved in degree 8, it
seems that it is just a matter of waiting
completion of
Hilbert's 16th problem in degree 8, until the Itenberg-Viro
contraction conjecture get refuted. This is merely a crude
scenario but of course one needs to keep track of a massive
diagrammatic to get an extension in degree 8 of Rohlin's theorem
(\ref{Rohlin-type:thm}) classifying all sextics according to their
types.

$\bullet$ (GR8)=(Gudkov-Rohlin census in degree 8).---Assume
someone has completed Hilbert's 16th in degree 8 (i.e. isotopy
classification of real schemes), how difficult will it be to
complete the corresponding Rohlin table enhancing schemes by their
types I or II. Assume this information available, does it follow
(by analogy with our proof of Prop.\,\ref{Klein-vache-deg-6:prop})
that under the contraction principle (CC), Klein-vache holds true
in degree 8? If yes, then Shustin 1985 would refute the
Itenberg-Viro contraction conjecture in degree 8 (CC8).

$\bigstar$ (RAG)=(Ragsdale).---While our erroneous Thom-style
estimate $\chi\le k^2$ (cf. \ref{Thom-Ragsdale:thm}) is disproved
by the Itenberg-Viro $(M-2)$-curve (Fig.\,\ref{Itenberg:fig}), is
this estimate still true for $M$-curves? This amounts to one-half
of Ragsdale conjecture $\vert \chi\vert \le k^2$ (still open in
the $M$-context). A priori a ``random'' computer-assisted search
along the Itenberg-Viro method could detect an $M$-curve refuting
Ragsdale. How difficult is it to program a machine adventuring
blindly and by brute force in such a random quest? In
contradistinction, how difficult is it to write down a proof of
Ragsdale's conjecture in case it should be true. Could it be that
a clever use of Thom, or Rohlin's formula and even some knowledge
of complex orientations derived maybe from our synthetic version
(\ref{total-reality-of-plane-M-curves:thm}) of the
Bieberbach-Grunsky theorem (planar case of Ahlfors) assesses the
full puzzle. If feasible this would be a spectacular application
of conformal geometry to the Hilbert-Ragsdale-Petrovskii 16th
problem, boiling down in quintessence  to Riemann's Nachlass 1857
\cite{Riemann_1857_Nachlass}. Of course we do not claim this to be
an easy project.

$\bullet$ In degree 6, it may be observed that among the trinity
of congruences mod 8 (due to GR, GKK, RKM, where G=Gudkov,
R=Rohlin, K=Krakhnov or Kharlamov (twice), M=Marin, cf.
(\ref{Gudkov-hypothesis:thm}),
(\ref{Gudkov-Krakhnov-Kharlamov-cong:thm}),
(\ref{Kharlamov-Marin-cong:thm})), the latter, i.e. RKM,  implies
the 2 formers, when combined with Rohlin's maximality principle
(RMC). Is this subsuming a general feature due to trivial
geographical/arithmetical reasons? If yes can we condense, i.e.
proceed to an unification of forces by reducing nearly all
prohibitions of Hilbert's 16th to the phenomenon of total reality.

$\bullet$ Can we write down an explicit ternary form with integral
coefficients $F\in \ZZ[x_0,x_1,x_2]$ whose real locus is Gudkov's
curve $\frac{5}{1}5$, and estimate the smallest size of the
coefficients involved? As discussed in
Sec.\,\ref{Diophantine-and-proba:sec}, can we compute the natural
masses (w.r.t. Lebesgue measure on the space of coefficients) of
each of the 64 chambers (past the discriminant) of smooth sextic
curves given by the census of Gudkov-Rohlin-Nikulin (i.e.
Fig.\,\ref{Gudkov-Table3:fig}).

$\bullet$ [31.03.13] A more modest but fundamental problem is to
publish (in the West side of Ural) an avatar in degree 8 of the
Gudkov-Rohlin table (Fig.\,\ref{Gudkov-Table3:fig}). For a partial
depiction of just the simplest planar face of this 4D-pyramid, cf.
Fig.\,\ref{Degree8:fig}. I presume that one can by mean of an
Atlas consisting of ca. 20 pages dress a list of all
combinatorially possible schemes after  taking into account the
obvious B\'ezout-style obstructions (B\'ezout, Zeuthen, Hilbert's
bounds on the depth of nest, Gudkov, plus the total reality
obstructions allied to the deep nest and doubled quadrifolium,
etc.). Once this atlas of all octics is made available it should
be a trivial matter to appreciate:

---how far/close we are to solve Hilbert's 16th in degree 8
(soft-isotopy);

---how the paradigm of total reality (resp. the contraction
principle) explain the prohibitions, and finally,

---whether the contraction principle (CC) implies Klein-vache, in
which case CC would be disproved by Shustin's refutation of
Klein-vache in degree~8.

$\bullet$ (CG6)=(Contiguity graph for $m=6$) [01.04.13] Can we
describe all the contiguity relation realizable via algebraic
Morse surgeries on the Gudkov-Rohlin table of periodic elements
(in degree 6). To be more specific, is some result along our
Conjecture~\ref{eversion-and-other surgeries:conj} true. The proof
of this could be merely a matter of adapting the work by Nikulin,
and Itenberg, yet it seems quite challenging to decide precisely
which eversions are realized algebro-geometrically.

$\bullet$ (KV7) Klein-vache (KV) was disproved in degree 8 by
Shustin 1985 via a conjunction of Viro's method and advanced
B\'ezout obstructions due to Viro (and Fiedler). On the other
hand, we prove below (\ref{Klein-vache-deg-6:prop}) that KV is
true in degree 6. So one may wonder about the case $m=7$, where to
my knowledge KV is undecided.

$\bullet$ (II/II) Is the ``toutou'' conjecture true? This posits
that any scheme of type~II and even degree $2k$ augmented by a
pseudoline to a scheme of degree $2k+1$ is of type~II too. Cf.
(\ref{toutou:conj}) for some surgical motivation (\`a la Fiedler)
and inspiration coming from reading Gross-Harris 1981, who were
unable to settle the case of quintics with 2 unnested ovals (and a
pseudoline of course). Perhaps a general solution of this problem
merely follows from a conjunction of Rohlin's and Mishachev's
formulae. If not, then one could use a large deformation
principle.

$\bullet$ (Klein's bipolarity conjecture). Is it possible for two
real plane curves to have distinct distributions of ovals, yet
conformally equivalent underlying symmetric Riemann surfaces
(under an equivariant diffeomorphism). This can be paraphrased in
the algebro-geometric language as the quest of two real planes
curves with distinct distribution of ovals, but bi-rationally
equivalent over $\RR$ as abstract curves. For more see
(\ref{Klein-bipolarity:conj}), where it is explained that the
first place where to look for this (hypothetical but likely)
phenomenon is degree $m=6$. It would be interesting to see if this
question due to Klein 1922 (safe misunderstanding on my side) can
be settled via Cremona transformations not inducing
diffeomorphisms of $\RR P^2$.

$\bullet$ (RIG/SAT) Is rigidity stable under satellites?  This is
a wild speculation based on Nuij's rigidity of the deep nest
caricatured as  reducible via satellite to the rigidity of the
conic (known since time immemorial). Likewise the more highbrow
rigidity result of Klein 1876 for quartics could induce rigidity
of all satellites of the 6 possible quartic schemes, in particular
of the quadrifolium (whose satellites are totally real under  a
pencil of conics). Further Nikulin's rigidity result for sextics
could imply also a vast array of rigidity results in degrees $6k$,
by satellitosis of all schemes of degree 6 which are not of
indefinite type (and which are explicitly known $64-2\cdot 8=48$
types by the Gudkov-Rohlin table=Fig.\,\ref{Gudkov-Table3:fig}).
Perhaps all this stability of rigidity under satellites has to be
combined with total reality, in which case the analogy with Nuij's
rigidity is still deeper. In that case we would only take
satellites of Rohlin's $(M-2)$-schemes of degree $6$, cf.
(\ref{satellite-Rohlin-(6)-schemes-rigid:conj}), and for quartics
only the quadrifolium (and the deep nest) would be permissible.

$\bullet$ (ANTI-GAB) It seems that the case of $(M-4)$-sextics of
type~I offers a possible corruption of Gabard's bound $r+p$.
Compare~Scholium~\ref{(M-4)-sextics-corrupt-Gabard:scholium}. If
not, this is  at least a {\it pi\`ece de r\'esistance\/} against
the principle that any abstract pencil is concrete, and therefore
Ahlfors abstract theorem is unlikely to apply without friction in
Hilbert's 16th problem. In other words Riemann's canary feels
claustrophobic in the Plato cavern of Brill-Noether-Hilbert.

$\bullet$ (LETOUZE-SCH) Inspired by a Scholium of Le~Touz\'e 2013
(\ref{LeTouze-quintic:scholie}), we extended her result to all
$M$-curves of odd degrees, cf.
Theorem~\ref{Le-Touzé-extended-in-odd-degree:scholium}. It seems
of interest to extend her method to even degrees as well. We had
just the time to treat the case of degree 6, cf.
Lemma~\ref{Le-Touzé-scholium-deg-6:lem} which uses imaginary
basepoints yet without affecting total reality. It could be
challenging to see if this method of total pencil (becoming more
and more explicit) could be used to reprove the deep prohibitions
for $M$-curves due to Hilbert-Rohn-Gudkov-Rohlin.

$\bullet$ (RMC via Mangler 1939 and Ahlfors 1950, maybe implicit
in Rohlin 1978).----Rohlin's maximality conjecture looks nearly
implied by Ahlfors, safe for the difficulty that the enlargement
of the type~I scheme is a priori very distant from the enlarged
curve realizing the orthosymmetric scheme. Using triviality of the
mapping class group of $\RR P^2$ (Mangler 1939, probably a student
of H. Kneser?), one can try to isotope the distant enlargement to
make it identic with the original curve. The latter being swept
out by a total pencil (Ahlfors 1950, plus epsilon!), one could get
a corruption of the homological version of B\'ezout (i.e.
intersection theory \`a la Poincar\'e, Lefschetz, etc.) Of course
one requires a procedure to extend the (Mangler) isotopy to $\CC
P^2$, and one may object that our sketch of proof equally well
applies to curves of type~I whose scheme is however of indefinite
type (but non-maximal). So there is perhaps some obstruction to
extend Mangler's isotopy as to preserve positivity of
intersection-indices. Understanding this obstruction, and
supposing one able to show its vanishing in case of a scheme of
type~I, could
procure a proof of the elusive RMC. It seems very likely that
Rohlin thought about this strategy, but never wrote something
down. Perhaps experts like Marin can complete this game? Cf.
Sec.\,\ref{RMC-via-Mangler:sec} for slightly more details.


\medskip

To keep some slight control on all these conjectures, see
Fig.\,\ref{CCvsCCC:fig} showing how they interact and their
validity range.

\begin{figure}[h]
\centering
    \epsfig{figure=CCvsCCC.eps,width=122mm}
\vskip-5pt\penalty0
  \caption{\label{CCvsCCC:fig}%
  A zoo of standard conjectures (Klein 1876, Rohlin 1978,
  Itenberg-Viro 1994, Le~Touz\'e 2013, Gabard 2013)} \vskip-5pt\penalty0
\end{figure}

\smallskip

[30.03.13] Let us conclude with a historical remark. It should
always be remembered, and amazing to rediscover everyday, that
``everything'' in this topic goes back to Klein. Himself expected
that the type of the symmetric Riemann surface (underlying a real
curve acted upon by complex conjugation) has some interplay with
Hilbert's problem on the distribution of ovals. Compare a footnote
added ca. 1922 in his Ges.\,Math.\,Abhdl., reproduced  as
Quote~\ref{Klein-1922-immer-vorsgeschwebt:quote}, but of which we
now reproduce the most prophetical side:

\smallskip
{\small

Es hat mir immer vorgeschwebt, dass man durch Fortsetzung der
Betrachtungen des Textes Genaueres \"uber die Gestalten der
reellen ebenen Kurven beliebigen Grades erfahren k\"onne, nicht
nur, was die Zahl ihrer Z\"uge, sondern auch, was deren
gegenseitige Lage angeht. Ich gebe diese Hoffnung auch noch nicht
auf, aber ich muss leider sagen, dass die Realit\"atstheoreme
\"uber Kurven beliebigen Geschlechtes (welche ich aus der
allgemeinen Theorie der Riemannschen Fl\"achen, speziell der
``symmetrischen'' Riemannschen Fl\"achen ableite) hierf\"ur nicht
ausreichen, sondern nur erst einen Rahmen f\"ur die zu
untersuchenden M\"oglichkeiten abgeben.

}

\smallskip

This is worth translating (in the poor English of the writer):

{\small

It always puzzled me, to infer more about shapes of real plane
curves of arbitrary degrees by pursuing  considerations of the
text, not only regarding the number of circuits, but also their
mutual dispositions. I do not abort this hope, but must alas
confess, that the reality theorems on curves of arbitrary genus
(which I deduce from the general theory of Riemann surfaces,
specially that of symmetric Riemann surfaces), are not sufficient
for this purpose, affording instead merely a
framework for the menagerie of possibilities to be investigated.

}

It is striking to notice how this Kleinian prose remains very much
actual, reflecting best our own frustration to make the
Ahlfors-Rohlin Verschmelzung, we are dreaming about, a true
reality.
It shows also how much Klein
would have appreciated the developments made possible in the
1970's by Gudkov, Arnold, and especially Rohlin, etc., and perhaps
even more, something like anticipating  the vision of total
reality by Rohlin.

[05.04.13] The last sentence of this same footnote, reads:

\smallskip

{\small

Da man \"uber die Natur dieser Bedingungen zun\"achst wenig weiss,
kann man noch nicht von vornherein sagen, dass alle die Arten
reeller Kurven, die man gem\"ass meinen sp\"ateren Untersuchungen
f\"ur $p={ n-1 \cdot n-2 \over 2}$ findet, bereits im Gebiete
besagter ebener Kurven $n$-ter Ordnung vertreten sein
m\"u{\ss}ten, auch nicht, da{\ss} ihnen immer nur {\it eine} Art
ebener Kurven entspr\"ache. \quad K.

}

\smallskip

Here, one realizes that Klein anticipated the simple phenomenon of
what Rohlin calls schemes of indefinite type, i.e. that the real
scheme alone (i.e., distribution of ovals) does not need to
determine the type (i.e. dividingness or not). Klein also
emphasizes the issue that not all topologically permissible
symmetric Riemann surfaces have to appear in the plane. In both
cases the first examples appear in degree 5, and then massively in
degree 6. For instance a quintic with only one pseudoline cannot
be of type~I , albeit since its genus is even (namely~6) the
corresponding Riemann surface exists. (Compare
(\ref{Klein-Marin-quintic:lem}) which is based on Klein-Marin, or
Gross-Harris argument via theta-characteristics discussed at the
same place that was probably known to Klein in 1892.

Finally, it is also puzzling to see that Klein 1892 anticipated
somewhat the contraction conjecture of Itenberg-Viro, cf.
historical note right after
(\ref{Itenberg-Viro-contraction:conj}).

Albeit Klein missed some basic modern tricks (like Rohlin's
formula, or Fiedler surgical smoothing law), he also mastered
perfectly the Riemannian theory (conformal maps, the allied circle
maps and total reality as credited by Teichm\"uller 1941,
theta-characteristics, allied deep enumerative problems of
bitangents to quartics \`a la Pl\"ucker-Zeuthen). Further he
appealed to contraction principles, as well as his own singular
geometric method to represent  complex loci as multiple cover of
the projective plane as to infer the ``complexified'' topology of
real curves (in the 1874--76 articles ``\"Uber eine neue Art der
Riemannschen Fl\"achen''). Hence Klein's legacy on the topic is
massive (ca. 300 pages, if one counts the G\"ottingen lectures
\cite{Klein_1892_Vorlesung-Goettingen}), and much remains to be
learned from it.

\section{The Klein-Rohlin conjecture on real schemes of type~I}
\label{Klein-Rohlin-conj:sec}

[01.01.13] A fascinating question is raised by (the master) V.\,A.
Rohlin in his 1978 \cite[p.\,95]{Rohlin_1978} survey
looping back directly to a (prophetic) allusion of Klein. Remember
that Rohlin was
fluent with German language, being
involved during World War II as translator on the front-line, cf.
Guillou-Marin's book 1986 \cite[p.\,ix]{Guillou-Marin_1986}:
``{\it En 1941, quand l'Allemagne attaqua l'U.R.S.S., Rohlin
rejoignit le corps des volontaires du Peuple (unit\'es militaires
non entra\^{\i}n\'ees). Son unit\'e fut encercl\'e et Rohlin fait
prisonnier par les allemands. Ensuite il r\'eussit \`a
s'\'echapper, \`a rejoindre l'arm\'ee sovi\'etique et finit la
guerre comme traducteur militaire (Rohlin parlait couramment
l'allemand). Imm\'ediatement apr\`es la guerre Rohlin fut
emprisonn\'e par la s\'ecurit\'e de l'arm\'ee (comme ce fut le cas
pour de nombreux anciens prisonniers de guerre) mais fut
lib\'er\'e \`a la fin de l'ann\'ee 1945.}''

\begin{quota}[Rohlin 1978] {\small {\bf 3.9 A conjecture
about real schemes of type~I}. {\rm A study of the available
factual material suggests that possibly a real scheme belongs to
type~I iff it is {\it maximal}, that is, it is not part of a
larger real scheme of the same degree. This conjecture is true for
$m\le 6$, and there is much to be said in its
favour\footnote{[28.03.13] I would personally be much interested,
if someone can guess more explicitly what Rohlin had in mind at
this place!} for $m>6$. There is an allusion to it in Klein: see
[4], p.\,155 (=Klein 1922=Ges. Math. Abh. II
\cite{Klein-Werke-II_1922}).}

}
\end{quota}

The passage Rohlin had in mind is unambiguously identified as the
following (going back actually to Klein 1876 \cite{Klein_1876}),
which is worth reproducing albeit it is first quite hard to
interpret (cf. also Viro 1986/86
\cite[p.\,67--68]{Viro_1986/86-Progress} or Degtyarev-Kharlamov
2000 \cite[p.\,785]{Degtyarev-Kharlamov_2000}, and especially
Marin 1988 \cite{Marin_1988}, clarifying earlier work partially
incorrect of Cheponkus 1976 \cite{Cheponkus_1976}):

\begin{quota}[Klein 1876]\label{Klein_1876-niemals-isolierte:quote}
{\small {\rm Die Kurven der selben Art zeigen
eine gro{\ss}e Reihe gemeinsamer Eigenschaften. Z.~B. kann bei den
Kurven der ersten Art durch allm\"ahl\-iches \"Andern der
Konstanten niemals eine isolierte reelle Doppeltangente neu
enstehen, um dann einen $(C+1)$-ten Kurvenzug zu liefern;
w\"ahrend die Kurven der zweiten Art in dieser Richtung nicht
beschr\"ankt sind. Die Kurven der zweiten Art sind sozusagen noch
entwicklungsf\"ahig, w\"ahrend es die Kurven der ersten Art nicht
sind. Doch soll hier auf diese Verh\"ahltnisse noch nicht n\"aher
eingegangen werden.}

}
\end{quota}

[17.01.13] It is essential to note that Klein's quote contains two
very distinct parts. The first part on which Klein is affirmative
may be translated as the assertion that a dividing (=type~I or
orthosymmetric) plane curve cannot acquire a new oval by
transgressing the discriminant at a solitary node (with imaginary
conjugate tangents like the germ $x^2+y^2=0$). With the strong
word ``kann \dots niemals'' (=never never!!), Klein emphasizes his
complete self-confidence about the truth of his assertion. Alas no
proof (as far as I know) were ever given by him, even in his
G\"ottingen lectures 1892 \cite{Klein_1892_Vorlesung-Goettingen}.
The first proof had to wait 112 years until Marin 1988
\cite{Marin_1988} write down a two-lines argument (of a somewhat
stronger assertion). A recent e-mail exchange with Viro suggested
that Klein's Ansatz may easily be deduced from the Ahlfors map
(cf. Lemma~\ref{Klein-via-Ahlfors(Viro-Gabard):lem}).

The second part of Klein's text, starting with ``w\"ahrend die
Kurven \dots'',  is pretty subtle to interpret and definitively
less categoric. It is suggested that curves of type~II are in
contrast susceptible of acquiring new ovals springing {\it ex
nihilo} from a solitary double point like a champagne bubble. The
vague wording ``sozusagen noch entwicklungsf\"ahig''
(=``so-to-speak still developable'') emphasizes that Klein did not
saw any nondividing curve champagne-bubbling, but merely that he
found no (topological) obstruction to such an eventuality. As we
shall see, this second clause which we shall call ``Klein-vache''
is refutable (in degree $8$) via the disproof of one half of
Rohlin's conjecture by Shustin 1985 \cite{Shustin_1985}. Alas, we
have not yet completely digested Shustin's work, which relies on
deep B\'ezout-style obstructions due to Fiedler-Viro. In the
positive sense, we proved via a
 cocktail of Russian results that ``Klein-vache'' holds true
in degree~6 (cf. Prop.\,\ref{Klein-vache-deg-6:prop}).

[01.01.13] Rohlin's conjectural criterium looks a
pearl of observational skills. What does it mean, or rather how
practical is it  if true at all? One should make a list of all
real schemes (i.e. isotopy classes) of curves of a given order.
Then assuming one competent and patient enough to have tabulated
the exhaustive list one could detect the dividing types by
inspecting maximum elements in the lattice ordered by inclusion.
({\it Insertion} [28.03.13].---This is not really what happens in
practice, especially since Shustin's disproof, and it seems more
likely that the residual half of Rohlin's conjecture acts by means
of prohibitions, that are anyway required to dress a table of all
schemes.)

Let us work out low-order examples to gain some experimental
evidence Rohlin is referring to.

First in degree $1$ there is just the line, which is of type~I
(=dividing). Then in degree 2, there is two isotopy classes
represented either by the circle $x^2+y^2=+1$ and the invisible
conic $x^2+y^2=-1$ (empty real locus). (This follows e.g. from
Sylvester's law of inertia, alias diagonalization of quadratic
forms, also to be found earlier by Jacobi, and presumably many
others? and on the case of 2 variables this can safely goes back
to ancient Greeks, Euclid, etc.) In order 3 we have cubics
(extensively studied by Newton and Pl\"ucker), but up to isotopy
the story becomes much simpler and we have two isotopy classes
differentiated merely by the number of real circuits $r=1,2$; the
latter being dividing while the other is not. This follows readily
from the abstract Klein-Weichold classification of symmetric
surfaces (the latter being merely a mirror image of the
M\"obius-Jordan classification of abstract topological surfaces).

\subsection{The little zoo of all quartics
(Pl\"ucker 1839, Zeuthen 1874, Klein 1876, Rohlin 1978)}

\begin{figure}[h]
\centering
    \epsfig{figure=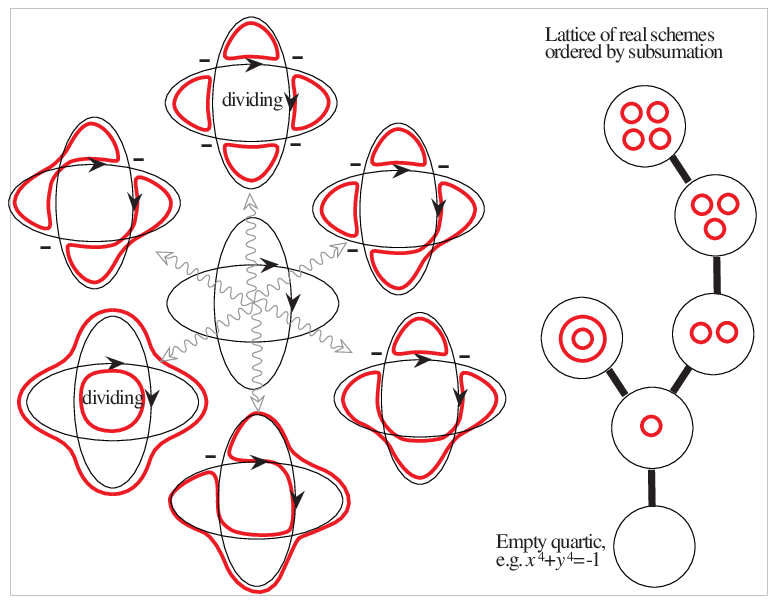,width=122mm}
\vskip-5pt\penalty0
  \caption{\label{KleinRohlin-quartic:fig}%
  The zoo of all quartics up to (rigid) isotopy}
\vskip-5pt\penalty0
\end{figure}

[01.01.13] Next it comes to quartics (order 4). Here the number of
real circuits $r$ fails to classify isotopy classes for there
exist quartics with 2 ovals being either nested or not. The first
basic thing-to-do (going back apparently to Pl\"ucker 1839
\cite{Plücker_1839}) is exploring varied examples by smoothing a
pair of conics with 4 intersections (cf.
Fig.\,\ref{KleinRohlin-quartic:fig}).

Among all those curves only those marked by the attribute
``dividing'' are dividing as they result from sense-preserving
smoothings. One can also remember Klein's congruence $r\equiv g+1
\pmod 2$ in the dividing case to detect nondividing curves, e.g.
those with $r=3,1$. Also under Harnack-maximality, i.e. $r=g+1$
(here $4$) then dividingness is automatic (by Riemann-Klein 1876
\cite{Klein_1876} or the more tedious synthetic argument of
Harnack 1876 \cite{Harnack_1876}). So $r=2$ is the only ambiguous
value. Here, as shown on the figure, the 2 ovals can either be
nested or not. In the first case the curve is dividing (due to
total reality under a pencil of lines through the innermost oval),
while  the unnested case is always nondividing (as Klein knew as
early as 1876).

\begin{lemma}\label{Klein-unnested-quartic-nondividing:lem} {\rm (Klein 1876)} All quartics
with $r=2$ unnested ovals (and not just the two traced above) are
nondividing.
\end{lemma}

\begin{proof} This is already a nontrivial result. We sketch several
proofs:

(1) How Klein derived the result? Maybe as follows. Klein knew as
early as 1876 \cite{Klein_1876_Verlauf} (basing himself on deep
works of Schl\"afli and Zeuthen on cubic surfaces and their
apparent contours which are quartics) that quartics are rigidly
classified by the real scheme. This is to mean that any two
quartic curves having the same distribution of ovals can be
continuously deformed through a large deformation of the
coefficients without ever meeting a singular curve. Hence Klein
had only to check the nondividing character of a specific quartic
to get that of all curves with 2 unnested ovals. Klein used a
special device of representation of the curve as a branched cover
of the projective plane by assigning to each point (of the
complexification) the unique real point of the tangent and so
could see the curve. Nowadays the surgical recipe of Fiedler looks
also best suited to do this. For an elementary graphical proof
compare our Fig.\,\ref{Guertel-genetic:fig} earlier in this text.

(2) Perhaps one way to argue could involve Ahlfors theorem, yet
some nontrivial details deserve being worked out. [03.01.13]
Assuming Gabard there is a total map of degree $r+p=3$ if dividing
(and not less via the complex gonality), yet since the ovals are
unnested it cannot be induced by a pencil of lines. So the
auxiliary curves are of order at least two. Assume first the order
to be two, so we have a pencil of conics. Since the degree of the
morphism is 3, we have $2\cdot 4-3=5$ basepoints on the $C_4$, but
a pencil of conics has only 4 basepoints by B\'ezout. (Note that
Ahlfors bound $r+2p=g+1=4$ would not be strong enough for this
purpose!) If the auxiliary curve are of order $3$, then we must
have $3\cdot 4-3=9$ basepoints on the $C_4$. No basic corruption
is detected?

(3) Another argument via theta-characteristics is implicit in
Klein 1892 \cite{Klein_1892_Realitaet} (see also his G\"ottingen
lectures 1891/92 \cite{Klein_1892_Vorlesung-Goettingen}) and
appears in modernized form in Gross-Harris 1981
\cite{Gross-Harris_1981}.

(4) Another more elementary (and purely topological) proof follows
from Rohlin's formula 1974--78 (valid for dividing curves), cf.
Sec.\,\ref{Rohlin-formula:sec}. This formula reads
$2(\pi-\eta)=r-k^2$, where $\pi$, $\eta$ are the number of
positive, resp. negative pairs of ovals. This distinction appears
by comparing the orientation induced as boundary of the half of
the Riemann surface underlying the curve, with that of the
 annuli bounding a pair of nested ovals. In our case there is no
 nesting hence $\pi=\eta=0$, and so $r=k^2=4$ violating the
 assumption $r=2$.

(5) A related proof involves Arnold congruence 1971
\cite{Arnold_1971/72} for $M$-curves of degree $2k$ (with an
obvious extension to dividing curves in Wilson 1978
\cite{Wilson_1978}). This reads $\chi:=p-n \equiv k^2 \pmod 4$ and
suffices. Here $p, n$ are notation coined in Petrovskii 1938
\cite[p.\,190]{Petrowsky_1938} for positive and negative ovals,
also interpretable as the number of even and odd ovals.  An oval
is said to be {\it even} if it is lying within an even number of
consecutive ovals. For the case at hand (2 unnested ovals), both
are even (being subsumed to zero ovals), hence $p-n=2-0\equiv
k^2=4 \pmod 4$ is violated, and the nondividing character of the
curve follows. The difference $p-n$ is readily interpreted as the
Euler characteristic $\chi$ of the Ragsdale membrane bounding
(orientably) the curve from inside. In the case at hand, the
Ragsdale membrane is the disjoint union of 2 discs, whence
obviously $\chi=2$.
\end{proof}


Now using the theorem of B\'ezout, it is clear (cf. Zeuthen 1874)
that our picture above (Fig.\,\ref{KleinRohlin-quartic:fig})
exhaust all possible shapes traced by quartics. For instance a
such cannot have 2 ovals nested in a third one, etc. So it is a
simple matter to convince that we have listed all real schemes of
quartics (with all of them safe the empty curve $x^4+y^4=-1$
arising through small perturbation of 2 ellipses). Of course a
priori a quartic could have 5 ovals but this was precluded by
Harnack 1876 \cite{Harnack_1876}, and of course already by Zeuthen
1874 \cite[p.\,411]{Zeuthen_1874} using a prototype of Harnack's
device. Indeed if a quartic had 5 ovals (or more) the conic
through them would cut it in $5\cdot 2=10>8=2\cdot 4$ overwhelming
B\'ezout. At any rate Klein's argument of 1876 via the underlying
Riemann surface gives the general Harnack bound $r\le g+1$  in
some more intrinsic fashion.

Even stronger is the following  result (due to Klein 1876
\cite{Klein_1876_Verlauf}, though his proof makes a d\'etour
through surfaces and it could be interesting to find a more direct
argument). Perhaps the transition through cubic surfaces is
necessary as it rationalize the irrationality of quartics curves,
though Klein in 1876 seems to have add a direct argument staying
in the realm of curves, but he did not exposed details.

\begin{lemma} The real scheme is a
complete invariant for rigid-isotopy classes of quartics.
\end{lemma}

(Rigid isotopy refers to the morcellation of the space of all
curves of some fixed order effected by the discriminant
hypersurface parametrizing singular curves.) Modulo such knowledge
one can draw the lattice of all real schemes (right part of
Fig.\,\ref{KleinRohlin-quartic:fig}) on which the (Klein-)Rohlin
intuition is verified: a real scheme is dividing (or of type~I)
iff it is maximal.

[28.03.13] In fact it is tempting to make a baby Gudkov table in
degree 4 (inspired by the case of degree 6, cf.
Fig.\,\ref{Gudkov-Table3:fig}) as to visualize the situation. Here
the Gudkov symbol $\frac{x}{1}y$ is merely a symbolical way to
mean that $x$ ovals are nested in a big oval (the denominator
$1$), while $y$ ovals are lying outside. It is noteworthy that
Rohlin's maximality principle is fully validated here and
prohibits  all the schemes lying above the configuration
$\frac{1}{1}$ of the nest of depth 2, which is already
B\'ezout-saturated. It is also pleasant to notice the presence
already of the highbrow Gudkov-Rohlin sawtooth (dashed on the
picture) so typical of the solution of Hilbert's 16th in degree 6
(cf. again Fig.\,\ref{Gudkov-Table3:fig}). This extends to all
degrees by the congruences of Gudkov-Rohlin $\chi\equiv_8 k^2=4$
for $M$-curves (\ref{Gudkov-hypothesis:thm}),
Gudkov-Krakhnov-Kharlamov for $(M-1)$-curve and Kharlamov-Marin
for $(M-2)$-curves of type~I). This sawtooth, which looks like a
piecewise linear sine-curve, forces the scheme below its
depressions, to be of type~I. Of course it tends to pass unnoticed
here ($m=4$) as it is such a trivial consequence of B\'ezout with
lines.

\begin{figure}[h]
\centering
    \epsfig{figure=Gudkov-Table-quartic.eps,width=122mm}
\vskip-5pt\penalty0
  \caption{\label{Gudkov-Table-quartic:fig}%
  A mini-Gudkov table for quartics}
\vskip-5pt\penalty0
\end{figure}

As already announced we conjecture in general that the whole
sawtooth can be explained by invoking the phenomenon of total
reality for $(M-2)$-curves via adjoint curves of order $(m-3)$. At
this stage a comparison with degree 6
(Fig.\,\ref{Gudkov-Table3:fig}) makes it puzzling to wonder if all
schemes lying below the sawtooth are always realized, yielding a
sort of denseness below the sawtooth (alias Gudkov line). The
answer is no in degree 8 (cf. Fig.\,\ref{Degree8:fig}), where
Petrovskii's estimate of 1933/38
(\ref{Petrovskii's-inequalities:thm}) starts to act prohibitively.
Yet taking this and other conjectural hypothesis of Ragsdale into
account, one can still wonder about the question of denseness of
schemes, namely the issue as to whether prohibition are
essentially confined to the high level of the pyramid (i.e. above
$(M-2)$-curves), or if in contrast there is some sort of porism
(or lacunae) killing schemes at low altitudes. As we shall see
latter if the conjectural maximality principle of Rohlin as well
as our stability of type~I under satellites holds true, it is
likely that for high degree $m$ (especially when the integer
$m=2k$ as a rich factorization into primes) then there will be a
myriad of cone-like {\it no man's land\/} zone
where schemes are killed because they extend a
B\'ezout-Ahlfors-Rohlin saturated scheme subsumed to the paradigm
of total reality.

Of course the situation of low degrees $m=4,6$ may give the wrong
impression that the whole paradigm of obstruction by the
saturation allied to total reality are already explained by the
trinity of Russian congruences mod 8 (of all the workers already
cited starting with Gudkov-Arnold-Rohlin). Yet in reality this is
not even true for low degrees because under the disguise of
B\'ezout it is already total reality which assures the planar
character of the lowest Gudkov tables $m\le 6$. Otherwise we had
to consider a menagerie of other schemes with more nested
structures. So this gives some intuition a priori that Rohlin's
maximality principle (in our opinion much allied to Ahlfors) will
not be subsumed to the trinity of congruence mod 8, albeit the
lowest of it pertaining to $(M-2)$-curves may act as vivid
generator of total reality phenomena.

\subsection{Quintics (Klein 1892?,
Rohlin-Mishachev 1976, Fiedler 78, Marin 79, Gross-Harris 1981)}
\label{quintic-table-Klein-Gudkov:sec}

[04.04.13] This survey has some repugnancy against curves of odd
order for reasons hard-to-explain, perhaps allied to the
cumbersomeness of the avatar of Rohlin's formula (due to
Mishachev). However the theory especially our main focus of total
pencils works as well in this case. Let us take a small look at
the ``Gudkov-Rohlin'' table in degree $m=5$. We recommend however
to skip this section on first reading as our understanding is
lacunary (in part because we do not discuss Mishachev's formula,
or because we do not entered into the Klein-Gross-Harris theory of
real theta characteristics). Yet, the case of quintics and more
generally curves of odd degrees (especially those of the shape
$m=5+4n$, else Klein's congruence suffices) offer a pleasant
application of the Klein-Marin principle, when it comes to check
that the scheme with only one circuit is of type~II (see
Lemmas~\ref{Klein-Marin-quintic:lem} and
\ref{Klein-Marin-odd-degree:lem}).

\begin{figure}[h]
\centering
    \epsfig{figure=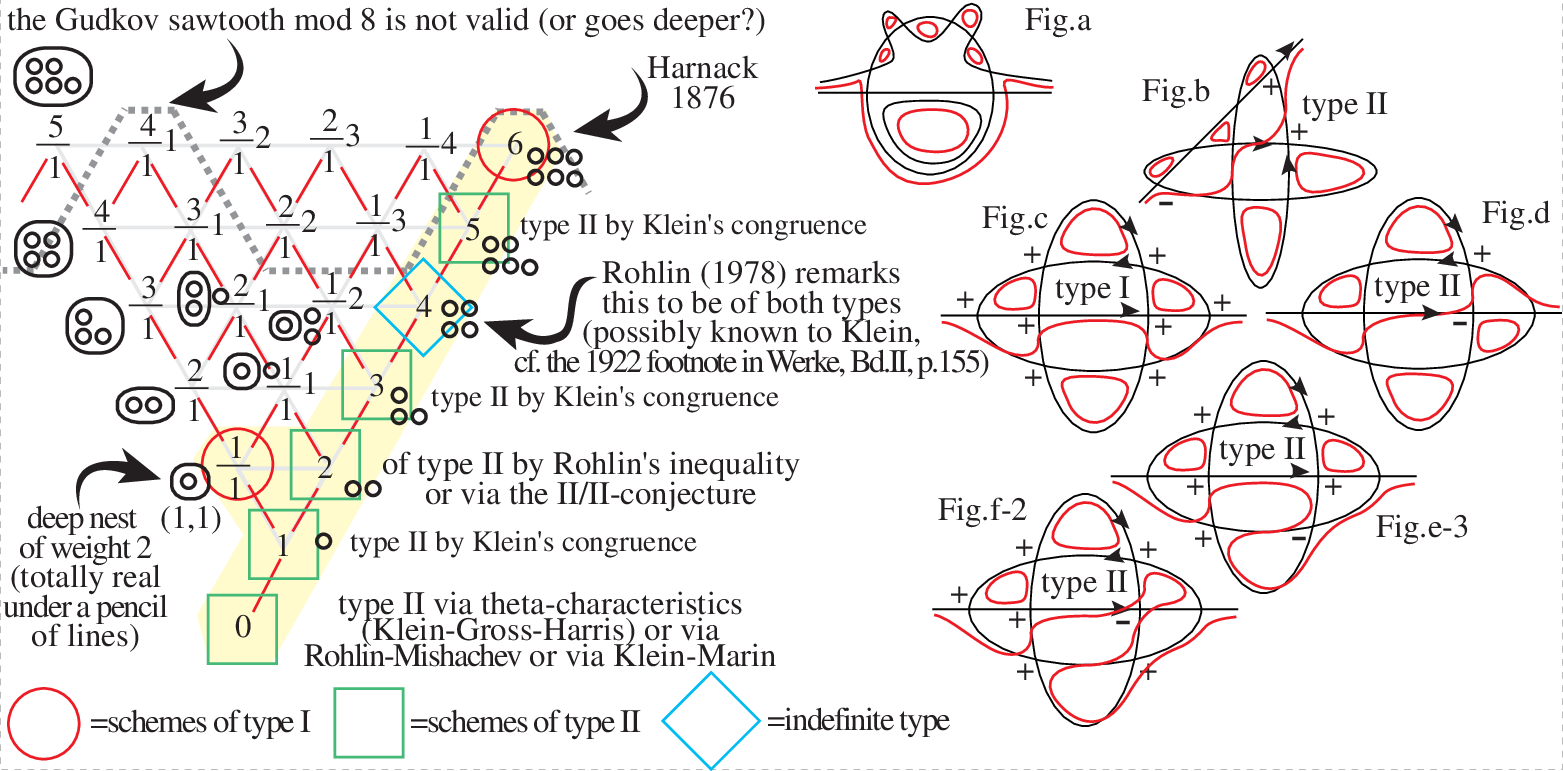,width=122mm}
    \vskip-10pt\penalty0

\caption{\label{Gudkov-Table-quintic:fig} A Gudkov table for
quintics}
\end{figure}

First, in degree $m=5$, Harnack's bound is $M=g+1=7$, since
$g=\frac{(m-1)(m-2)}{2}=(4\cdot 3)/2=6$. In odd degrees there is
always a unique pseudoline (M\"obius 18XX \cite{Moebius_18XX}, von
Staudt 18XX \cite{von-Staudt_18XX}, Zeuthen 1874
\cite{Zeuthen_1874}, Harnack 1876 \cite{Harnack_1876}, Hilbert
1891 \cite{Hilbert_1891_U-die-rellen-Zuege}, etc.), and we may
omit it from the Gudkov symbols. Hence on the table below
(Fig.\,\ref{Gudkov-Table-quintic:fig}) we suppress the pseudoline
from the real scheme depiction. As usual one of the most
noteworthy configuration is the deep nest, here
$\frac{1}{1}=(1,1)$ which is total under a pencil of lines, hence
of type~I.

{\footnotesize

{\it Bibliographical puzzle}.---This argument looks to us much
more elementary than the one of Gross-Harris 1981
\cite[p.\,175]{Gross-Harris_1981} via theta-characteristics which
 goes back probably to Klein 1892. Despite its extreme
 elementariness, it seems also to have escaped Klein's attention.
 In fact Klein uses it in 1892 (p.\,177 of Ges. Math. Abh, II) but
 only after contracting the empty oval of the G\"urtelkurve. We
 cite the relevant passage:

\smallskip

Sollen wir diese geometrischen Verh\"altnisse durch Beispiele
belegen, so nehmen wir viel\-leicht zun\"achst den Fall der
G\"urtelkurve $p=3$. Hier hat es ersichtlich keine Schwierigkeit,
das innere Oval auf einen Punkt zusammenzuziehen. Von diesem aus
projizieren wir jetzt die Kurve auf eine gerade Linie. Die Gerade
wird dann nach ihrer ganzen Erstreckung von den Bildpunkten
doppelt \"uberdeckt, so zwar, da{\ss} dabei kein reeller
``Scheitel'' auftritt\footnote{This shows that Klein anticipated
the phenomenon of total reality.}. Das entspricht in der Tat dem
orthosymmetrischen Falle $\lambda=1$ des Geschlechtes $p=2$.

\smallskip

Later, the argument of total reality seems to have escaped the
attention of another great master, namely Alexis Marin 1979
\cite{Marin_1979}, compare especially on p.\,56 his complicated
argument for ``N'existe pas'' in the bottom-right angle of the
tabulation, as well as the question p.\,59: ``Est-ce qu'une courbe
ayant cette disposition s\'epare sa complexifi\'ee.''

}

 This total reality (or
B\'ezout-saturation) kills all schemes enlarging it (cf. unframed
white-colore schemes on Fig.\,\ref{Gudkov-Table-quintic:fig}).
Apart from this obstruction there is essentially no other. First
we can construct an $M$-curve necessarily unnested (by
B\'ezout-saturation) with symbol $6$ (again the pseudoline $J$ is
omitted). This is constructed via a Hilbert-method on
Fig.\,\ref{Gudkov-Table-quintic:fig}a (for a less schematic
picture cf. Fig.\,\ref{Harna0:fig}). We presume that nobody knew
 existence of such a curve prior to Harnack 1876 (but this is
just a historical challenge, perhaps Pl\"ucker, Zeuthen, Klein
before but not sure).

Next Fig.\,b. shows a quintic with 5 ovals (and one pseudoline),
which arises by slight perturbation of 2 ellipses plus a line.
This was probably known to Pl\"ucker 1839 \cite{Plücker_1839}, or
earlier. The curve constructed is of type~II by Fiedler's signs
law, or just by Klein's congruence $r\equiv_2 g+1$. The latter
forces actually all the schemes $5,3,1$ being of type~II.

Below we have the scheme $4$. This admits realizations of both
types (as knew Rohlin 1978, and perhaps Klein?), as shown by
Fig.\,c and Fig.\,d (using Fiedler's smoothing law). The curve of
Fig.\,c has actually a total pencil of conics assigned to pass
through the 4 ovals (which was depicted earlier in this text).
Again it is interesting to note that Gross-Harris (p.\,176--177)
used a somewhat more synthetic and complicated argument than just
Fiedler's law, to show existence of curves in both types I/II.

{\it Historical note}.---The above phenomenon is the first
instance of where the type of a curve is not determined by the
distribution of ovals. It admits as a simple consequence the fact
that there exists obstructions to rigid-isotopy lying beyond the
real scheme (remark due to Rohlin 1978). It is clear however that
Klein knew (or at least suspected) this basic phenomenon, compare
his footnote of 1922:

\smallskip
{\footnotesize

Da man \"uber die Natur dieser Bedingungen zun\"achst wenig weiss,
kann man noch nicht von vornherein sagen, dass alle die Arten
reeller Kurven, die man gem\"ass meinen sp\"ateren Untersuchungen
f\"ur $p={ n-1 \cdot n-2 \over 2}$ findet, bereits im Gebiete
besagter ebener Kurven $n$-ter Ordnung vertreten sein
m\"u{\ss}ten, auch nicht, da{\ss} ihnen immer nur {\it eine} Art
ebener Kurven entspr\"ache. \quad K.

}


[Inserted 05.04.13].---In fact this can be interpreted either \`a
la Rohlin, by saying that a real scheme can have realizations in
both types (I/II=ortho- or diasymmetric). Somewhat more crazy
would be the following interpretation.

\begin{conj} \label{Klein-bipolarity:conj}
{\rm (Kleinian bipolarity---Klein 1922, Gabard 2013)}.---An
abstract symmetric Riemann surface (SRS) can admit plane
realizations with distinct distributions of ovals.
\end{conj}

Actually I do not know if this phenomenon of ``bipolarity'' can
occur. Of course it does trivially occur, with a line and a conic
both representing the Riemann sphere with its standard real
structure (equatorial involution). In degree 3, it cannot occur
since the real scheme determines the type, and likewise in degree
4 (cf. previous section). In degree 5, predestination of the type
by the real scheme is not true any more, yet the combinatorics of
Fig.\,\ref{Gudkov-Table-quintic:fig} is sufficiently simple as to
preclude bipolarity. Indeed if the SRS is fixed, hence in
particular the number $r$ of real circuits, the only height at
which there are several distributions of ovals is $r=3$, where we
have the nested ($\frac{1}{1}$) and the unnested ($2$) schemes.
Yet both of them are differentiated by the type. Hence the first
place to look for bipolarity is degree 6. Here we have (see the
Gudkov-Rohlin table Fig.\,\ref{Gudkov-Table3:fig}) a myriad of
sextics having the same underlying (topological) symmetric
surface. It is unclear if they can be conformally equivalent while
exhibiting different distributions of ovals. It could be imagined
that a Hilbert sextic is sometimes conformally diffeomorphic to
one of Harnack, or even Gudkov. This question looks a bit
artificial or puzzling, yet has perhaps of some importance if one
likes to link with the abstract theory of Ahlfors taking into
account only the abstract Riemann surface. To settle the
bipolarity question we can look at the natural map from the
hyperspace of smooth plane curves to the real moduli space of SRS,
i.e. $\mH-\disc \to M_g$. Perhaps then two chambers may have
overlapping images. As noted by the old Felix  Klein (aged 73 at
the moment of his 1922 footnote) plane curves have specialized
moduli. Hence the images in question are fairly small subloci of
the moduli space, but this does not prevent  overlap. Another
approach is to use Cremona transformations of the plane defined
over $\RR$ which do not induce diffeomorphisms of the plane $\RR
P^2$ (this remembers works by Ronga-Vust ca. 2002, or their
student J.~Blanc). By this procedure we can perhaps alter the
distribution of ovals, yet without distorting the conformal
structure, as the curve and its image are in birational
equivalence. Can this vague idea be implemented? Otherwise the
approach can be the Teichm\"uller-theory of the map described
above from concrete plane curves to the moduli space of abstract
real(=symmetric) Riemann surfaces, while trying to study exactly
the coincidences of this mapping. One could try to determine
exactly which among the 64 chambers of sextics (Nikulin's theorem
(\ref{Nikulin:thm})) are in bipolarity, i.e. contains conformal
replicas of the same symmetric Riemann surface. This defines an
additional graph structure on the Gudkov-Rohlin table where edges
are traced whenever two vertices(=chambers) contains curves
abstractly isomorphic over $\RR$. Of course the edges have to
preserve the height $r$ on the Gudkov pyramid, as well as the
types. Those (topological) obstructions to bipolarity could be the
sole ones, in which case the graph in question would show plenty
of edges. In particular it restricts to the complete graph on each
levels at height $r$ not congruent to $g+1 \pmod 2$ where
diasymmetry reigns ubiquitously (Klein's congruence). Few other
levels have pure types too, e.g. $M$-curves (type~I only), and the
level $r=3$ (types~II only) (via Rohlin's formula), see again
Fig.\,\ref{Gudkov-Table3:fig}. At those levels the bipolar graph
could be  complete too.

After this digression, we return to the basic classification of
quintics. The scheme $3$ is of type~II by Klein's congruence (as
we already noted), and exists as shown by Fig.\,e.

Next the scheme 2 exists in type~II as shown by Fig.\,f. It is
however more tricky to prove that the scheme $2$ is of type~II.
This follows either from the avatar in odd degree of Rohlin's
formula (i.e. Mishachev's formula). Perhaps there is a more
elementary argument, say by using a pencil of lines while trying
to permute 2 imaginary points during a sweeping. Also Gross-Harris
have probably an argument via theta-characteristics, but alas
those authors confess being not able to prove this compare p.\,175
of Gross-Harris 1981 \cite{Gross-Harris_1981} where we read: ``In
the non-nested case, we suspect that $a(X)$ is always 1 [i.e.
type~II, or nondividing] but have no proof''. Was Felix Klein
(1892 paper \cite{Klein_1892_Realitaet} or lectures
\cite{Klein_1891--92_Vorlesung-Goettingen}) able to tackle this
case?

A crude principle that do this work is the postulation that
whenever we add a pseudoline to a scheme of type~II, it remains of
type II. Recall that we know since Klein that the quartic scheme
$2$ is of type~II, cf.
Lemma~\ref{Klein-unnested-quartic-nondividing:lem}. Some evidence
comes from surgeries on the Riemann surface while noticing that
diasymmetry is a dominating character in the genetical sense. This
is implicit in Fiedler's law of smoothing and really a simple
matter of visualizing the corresponding Riemann surfaces. So let
us posit the:

\begin{conj}\label{toutou:conj} {\rm (Gabard 2013,
but probably standard by Rohlin-Fiedler, if not
erroneous)}.---When a scheme of even degree $2k$ is of type~II,
then the same scheme of degree $2k+1$ augmented by a pseudoline is
of type~II too. (Il y a trop de toutous dans la langue anglaise,
mon ostie!)
\end{conj}

Proving this could again involve a large deformation principle (as
discussed in the sequel) like minimizing the length of the
pseudoline as to make it a line. There will then be a
strangulation of the Riemann surface and we are reduced to
Fiedler's genetic law. Perhaps there is an elementary proof of the
conjecture based on a conjunction of  Rohlin's and Mishachev's
formulae. Note also that the conjecture holds true for the empty
scheme by Lemma~\ref{Klein-Marin-odd-degree:lem} below.

Another idea to show that the quintic (unnested) scheme $2$ is of
type~II could be to use Klein's 1876 remark that a curve of type~I
cannot acquire a solitary node, cf. below for a proof
(\ref{Klein-Marin:lem}) essentially along the lines of Marin 1988.
Then we are reduced to showing that any quintic with scheme 2
(again we omit the pseudoline $J$) can indeed acquire a new oval.
This looks a priori hard, but in view of the diagrammatic of the
table Fig.\,\ref{Gudkov-Table-quintic:fig} (mostly prompted by
B\'ezout) we could just make a deformation along a pencil spanned
by the curve plus a curve with more ovals (e.g. Harnack's or just
$5$ of Fig.\,b). The difficulty however is that the deformation is
not forced to raise immediately the number of ovals, as it may
first lower down the number of ovals. Incidentally if this
argument via Klein-Marin would have worked it would also have
prohibited the type~I realization of the scheme 4.

Next we have the scheme $1$ forced to be of type~I, by Klein's
congruence, and easily constructed (e.g. by a slight alteration of
the picture Fig.\,\ref{Gudkov-Table-quintic:fig}f above).

Finally, the scheme $0$ poses again a little problem, but can also
be shown to be of type~II. This follows either from
Rohlin-Mishachev, or via theta-characteristics. In this case
Gross-Harris 1981 \cite[p.\,175]{Gross-Harris_1981} were able to
conclude type~II via theta-characteristics (see their proof of
Prop.\,7.1, p.\,173, which contains some minor misprints, namely
``Prop.\,4.1'' should be ``Prop.\,5.1'', and ``$h^0({\bf
a})=(d^2-1)/2$'' should be ``$h^0({\bf a})=(d^2-1)/8$'').

A somewhat more conceptual argument can be based on Klein's Thesis
(as Viro calls it) of 1876 to the effect that {\it a curve of
type~I cannot gain an oval (at least when crossing a solitary
node)\/}. This was perhaps historically the first known proof,
albeit Klein did not mentioned this consequence explicitly in
print (1876 paper, nor later). On writing down the proof below, we
realized that one needs the stronger version due to Marin 1988 of
Klein's Thesis relaxing the parenthetical proviso above. Hence our
claim of historical priority is somewhat sloppy, but in substance
Klein could have anticipated it.

\begin{lemma}\label{Klein-Marin-quintic:lem}
{\rm ($\approx$Klein 1876, 1892, Rohlin-Mishashev ca. 1974--76,
Gross-Harris 1981, Marin 1988, Gabard 2013 trying to assembly all
this today)}.---Any quintic with only one pseudo-line is
necessarily of type~II (i.e. nondividing or diasymmetric).
\end{lemma}

\begin{proof}
Take such a curve $C_5$ (with only a pseudoline) and any auxiliary
(smooth) curve with at least one oval (and so $r\ge 2$). Pass a
line through both curves (in the hyperspace of curves) and perturb
it slightly to ensure transversality w.r.t. the discriminant.
Since the initial curve $C_5$ has the least possible number of
real circuit (namely one), the first contact (along one of the 2
possible pathes inside the pencil) with the discriminant will be a
``Morse'' surgery (jargon Thom-Milnor) {\it forced to increase the
number of ovals\/} ($\bigstar$).

This last (italicized) assertion ($\bigstar$) requires perhaps
more
substantiation. Let us admit it to conclude quickly. If the new
oval raises from a solitary-node then Klein's Thesis of 1876 (alas
left unproven by the great geometer) suffices to conclude. If not,
e.g. if the pseudoline self-collides with itself as to split off a
new oval (Fig.\,\ref{Eversionpseudo:fig}a), then Marin's version
of Klein completes the proof, cf. (\ref{Klein-Marin:lem}) or Marin
1988 \cite{Marin_1988}.

To justify better ($\bigstar$) we should check that all Morse
surgeries of a pseudoline forces an augmentation of the number of
circuits. In the case of an oval this is not true due to
``eversions'' (cf. Sec.\,\ref{Eversion:sec} especially
Fig.\,\ref{Eversion:fig}), whence our extreme prudence. However
doing naive experimental pictures deforming a pseudoline, it seems
impossible to evert a pseudoline
(Fig.\,\ref{Eversionpseudo:fig}b). It remains of course to find a
theoretical explanation.
\end{proof}

\begin{figure}[h]
\centering
    \epsfig{figure=Eversionpseudo.eps,width=122mm}
\vskip-5pt\penalty0
  \caption{\label{Eversionpseudo:fig}%
  Trying to evert a pseudoline}
\vskip-5pt\penalty0
\end{figure}

It is clear that the above lemma extends to all other odd degrees:

\begin{lemma}\label{Klein-Marin-odd-degree:lem}
Any curve with a unique real circuit is of type~II, safe if it is
a line or a conic (degree $m=1,2$).
\end{lemma}

{\it Remark}.---The argument of (Klein-)Gross-Harris only works
under the (extraneous and stringent) assumption  $m \equiv 5 \pmod
8$ (cf. their Prop.\,7.1, p.\,173).

\smallskip

\begin{proof}
The case of odd degrees follows by the same method using the
Klein-Marin theorem. The case $m=3$ is of course more elementary
and can be reduced to the uniformization of elliptic curves e.g.
\`a la Weierstrass via the doubly-periodic $\wp$-function defined
on a rhombic lattice. Sorry, it suffices actually to use  Klein's
congruence, or to remember---if you do not want to sell your soul
to the devil of arithmetics---that a symmetric torus with one
fixed circuit is forced to be $S^1\times S^1$ acted upon by
exchange of both factors (while fixing the diagonal circle).

Actually Klein's congruence $r\equiv_2 g+1$ settles the lemma
whenever $m=3+4n$. Indeed then
$g=\frac{(m-1)(m-2)}{2}=\frac{(2+4n)(1+4n)}{2}=(1+2n)(1+4n)=1+6n+8n^2$,
which is odd, and so Klein's congruence (forced by type~I) is
corrupted, whence type~II.

For the other cases $m=1+4n$, Klein's congruence tells nothing and
one make appeal to the Klein-Marin argument instead.

For even degrees, one can again treat half of the cases via
Klein's congruence, namely when $m=4, 8, 12, \dots$, i.e. $m=4n$
as then $g=\frac{(4n-1)(4n-2)}{2}=(4n-1)(2n-1)$ which is odd, and
so Klein congruence is violated for $r=1$. For the other cases
$m=2+4n$, the congruence tells nothing. However we can still
conclude type II (of course provided $m\ge 4$), either via
Rohlin's formula (\ref{Rohlin-formula:thm}) or maybe a variant of
the Klein-Marin argument. However now the configuration with one
circuit has not the  minimal number  of circuits and so we may
first descend to the empty chamber and the Klein-Marin method
looks impuissant. Of course perhaps some extra trick can ensure
that we can increase the number of component immediately yet I do
not see any obvious argument.
\end{proof}

\subsection{The impressive landscape of all sextics
(Harnack 1876, Hilbert 1891/00/09, Rohn 1911/13, Petrowskii
1933/38, Gudkov 1948/54/69, Arnold 1971, Rohlin 1972/74/78)}

[31.12.12] Perhaps the  Klein-Rohlin conjecture follows from
Ahlfors theorem interpreted in terms of total reality. Intuitively
having a total pencil, no real circuit can be added without
corrupting B\'ezout (more on this in
Sec.\,\ref{Rohlin-via-Ahlfors}). Yet perhaps this is too naive as
shown by an example of order 6 to be found in Gabard's Thesis 2004
\cite[p.\,8]{Gabard_2004} (I should acknowledge Kalla-Klein 2012
\cite{Kalla-Klein_2012-Computation-cite-Gabard} for reminding me
that my Thesis contained this example).

\begin{figure}[h]
\centering
    \epsfig{figure=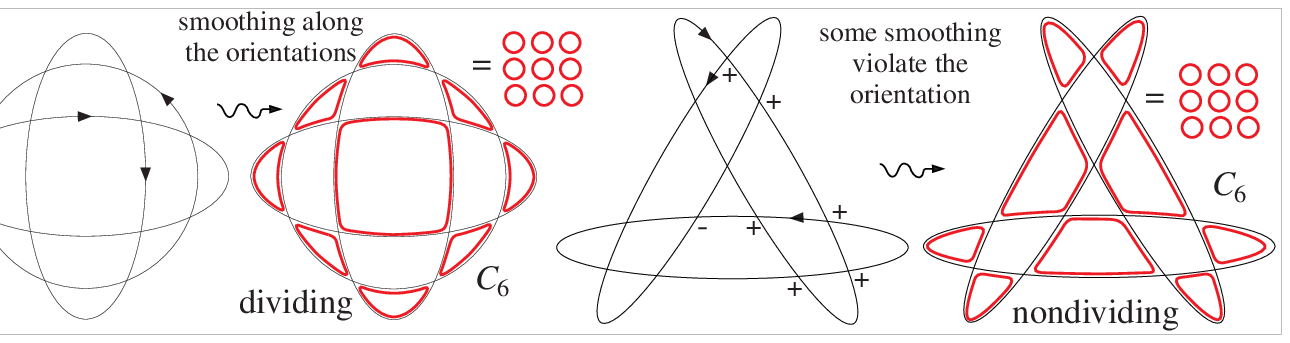,width=122mm}
\vskip-5pt\penalty0
  \caption{\label{KleinRo-sextic:fig}%
  A counterexample to the Klein-Rohlin conjecture? (Not at all!)}
\vskip-5pt\penalty0
\end{figure}

This shows that the real scheme alone fails to determine the
dividing character (alias type~I=erster Art in Klein 1876
\cite{Klein_1876}). At first I thought this corrupts Rohlin's
assertion that his conjecture is true in degree 6. Of course
Rohlin 1978 \cite{Rohlin_1978} knew very much this phenomenon,
which he calls ``real schemes of indefinite type'' (on p.\,94 of
\loccit), i.e. real schemes admitting representatives of both
types (dividing or not). Hence it was first puzzling to wonder why
he made such a basic mistake, or more likely why we first failed
to interpret correctly his simple message?

\begin{figure}[h]
\centering
    \epsfig{figure=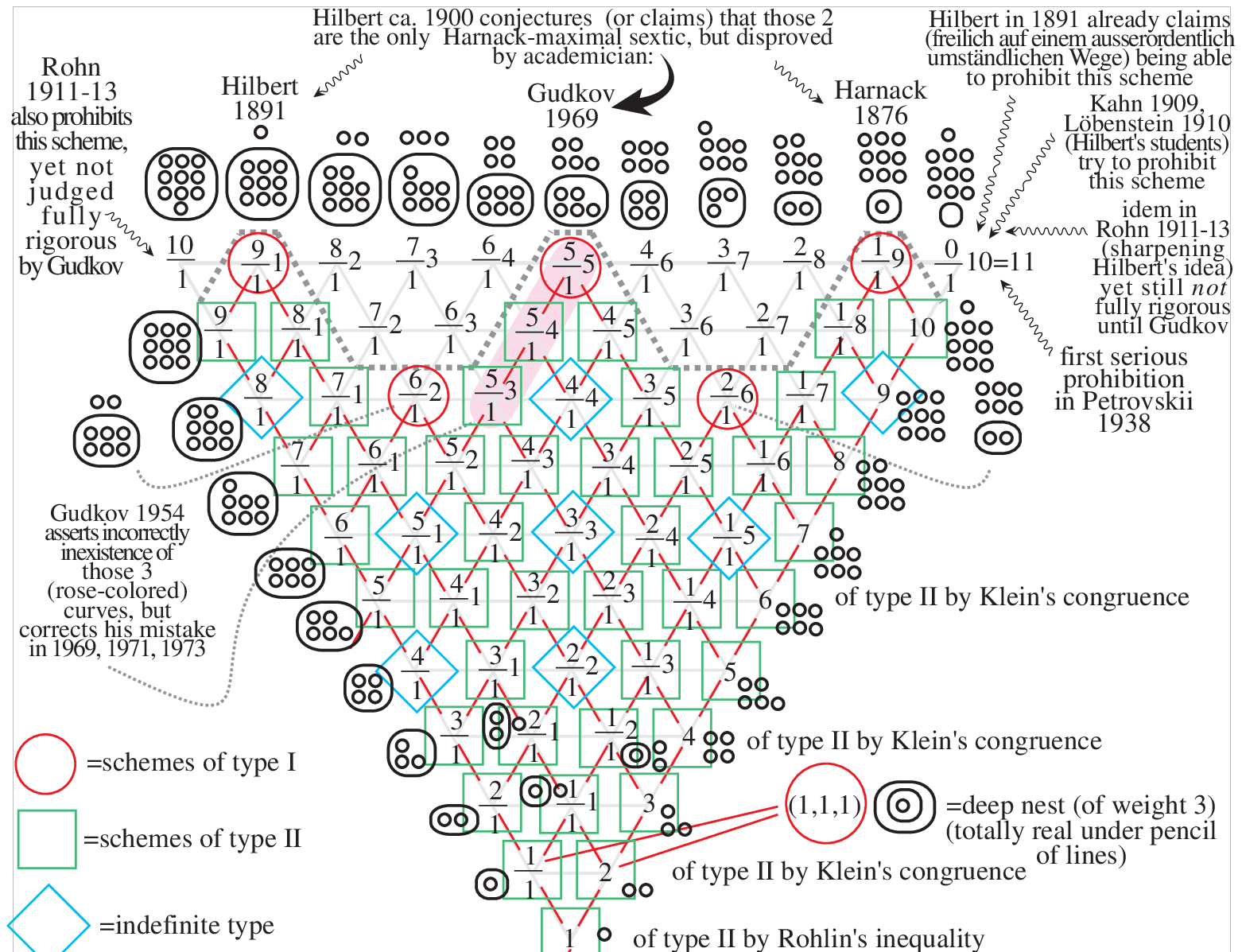,width=122mm}
\vskip-5pt\penalty0
\caption{\label{Gudkov-Table3:fig}%
Gudkov's table of all sextic real schemes (enhanced by
Klein-Rohlin data of complex characteristics): those below the
broken line are realized (Harnack, Hilbert, Gudkov), while those
above are not (Hilbert, Rohn, Petrovskii, Gudkov). Red-circles
mark the maximal schemes which turn out to be those of type~I
(giving some support to the Klein-Rohlin conjecture 1876-1978 that
the hierarchy is dominated by pure orthosymmetry)}
\vskip-5pt\penalty0
\end{figure}

Understanding the Klein-Rohlin conjecture requires some more
mature thinking. One should list all schemes dominating the ``nine
unnested ovals'' scheme. On p.\,95 Rohlin 1978 \cite{Rohlin_1978}
refers to the census (tabulation) set up by Gudkov 1974
\cite[p.\,40]{Gudkov_1974/74} listing all the logically possible
(real schemes of) sextic curves (taking into account B\'ezout for
lines). This is worth reproducing as
Fig.\,\ref{Gudkov-Table3:fig}. Recall that in Gudkov's symbolism,
$\frac{x}{1} y$ denotes the scheme consisting of $x$ ovals
enclosed by one ``big'' oval, while there is $y$ ovals living
outside. This gives a total of $1+2+3+\dots+11=\frac{12 \cdot
11}{2}=6\cdot 11=66$ logically possible curves (counting inside
the ``triangle''), to which must be added the empty real scheme
(denoted $0$) and the deep nest of depth $3$ (denoted $(1,1,1)$ or
$\frac{1}{{\frac{1}{1}}}$). We get so the 68 schemes ({\it ann\'ee
\'erotique}) mentioned by Gudkov (p.\,40).

Not all of those 68 schemes are actually realized.
%
%
If they would this would roughly mean that all obstructions are
B\'ezout-like prompted by tracing a single line. However the plane
is swept out by a myriad of other curves. Quite eclectically, the
architecture of  Hilbert-Gudkov's  table of elements is a bit like
a pharaohs pyramid (turn Fig.\,\ref{Gudkov-Table3:fig} upsidedown)
and those are known to have sanctuary galleries forming tunnels.
The first to have spotted this porousness of the pyramid is no
less an authority than Hilbert 1891
\cite{Hilbert_1891_U-die-rellen-Zuege}, albeit it took several
decades until his work got consolidated (especially by Rohn 1911,
Petrovskii 1933/38) and pushed forward to its ultimate perfection
(thanks to the efforts of Gudkov). Soon afterwards, Arnold and
Rohlin offered quite dramatic simplifications based on pure
topology,  and extensions of the prohibitions to all degrees. We
suspect this Hilbert-Gudkov pyramid to have some overlap---not yet
much elucidated except perhaps for allusions in Rohlin 1978
\cite[p.\,94]{Rohlin_1978}---with some more ancient force, namely
Abel-Riemann-Klein-Teichm\"uller-Ahlfors and their circle maps.
The latter are of course specific to dividing curves concomitant
with the paradigm of total reality. In the case at hand (real
plane smooth sextics), total reality is exhibited according to a
Rohlin's claim 1978 (not yet fully understood by the writer) via
total pencil of cubics for certain schemes, e.g. $\frac{6}{1}2$
and its mirror $\frac{2}{1}6$. More trivial is the nest of depth
$3$, totally real under a pencil of lines. This can be interpreted
as a prohibition of diasymmetry for those schemes. Likewise the
porous portion of Gudkov's pyramid (=white cases above the broken
line on Fig.\,\ref{Gudkov-Table3:fig}) can be prohibited (and this
is how I understand vaguely the Hilbert-Rohn method) by pure
synthetic geometry. Paraphrasing, not  merely linear B\'ezout
obstructions do exist, but also those via the menagerie of all
other curves grooving nonlinearly the plane. More than that, not
just static curves but dynamical collections of such (e.g.,
pencils) have to be considered. It is charming to note a strong
parallel between Hilbert's and Rohlin's claims that pure geometry
is able to prohibit schemes, especially as both look
insufficiently justified, but intuitively plausible. How much
Rohlin's synthetic proof of the type~I of the schemes
$\frac{6}{1}2$ and its mirror has in common with Hilbert-Rohn's
method?

A last word of caution for pyramids builders: one of the first
ever constructed in Ancient Egypt had a somewhat pathetic destiny.
Once arriving near the $2/3$ of the planned final size, fissures
started to appear menacing the whole foundations to crack under
pressure. The only reasonable option
left to the
engineers was to diminish the slope for the last third as to lower
pressure. It is not known if this sufficed to ensure immortality
of the Pharaoh. Thus, it should be no surprise that the most
telluric  part of the pyramid (near the funerary chamber  where
the pressure is highest) is the most secrete part of the edifice.
This needed to wait the contribution of Gudkov 1969
\cite{Gudkov-Utkin_1969/78} who exhibited the most elusive schemes
$\frac{5}{1}5$ of Fig.\,\ref{Gudkov-Table3:fig} (cf. also
Fig.\,\ref{GudkovCampo-5-15:fig} for the explicit construction).
At this stage Hilbert's 16th problem was completely solved (at
least for sextics which is arguably the official context of
Hilbert's question).

Now let us be more formal. As explained by Gudkov (1974 \loccit),
Kahn 1909 \cite{Kahn_1909} and L\"obenstein 1910
\cite{Löbenstein_1910} published dissertations under Hilbert's
direction---(if I understood well the story, vgl. Hilbert 1909
\cite{Hilbert_1909-Ueber-die-Gestalt-sextic}, both were feminine
candidates)---attempting to prohibit sextics with 11 unnested
ovals. (Challenge: try to prove this via Ahlfors 1950 or rather
via Bieberbach-Grunsky (1925/1937). Philosophically, this would
just, as it should, put Little Hilbert in the baskets of Big
Riemann!) This follows also from Rohlin's formula of 1974--78, cf.
Sec.\,\ref{Rohlin-formula:sec}, or from Arnold's congruence of
1971. Soon later Rohn 1911--1913 \cite{Rohn_1913} devoted two
articles attempting by the same method to exclude sextics of type
$\frac{10}{1}$ or $11$, making a big contribution to the
development of Hilbert's idea. The resulting prohibition method
was christened by Gudkov (1974 \loccit) the {\it Hilbert-Rohn
method}. In Gudkov's view, even Rohn's proof is not perfectly
sound due to some messy combinatorics impeding Rohn to take care
of all logically possible cases. Gudkov then mentions several more
Western attempts, by Wright 1907 (same idea as Hilbert, but not
rigorous prohibition of type $11$). In Donald 1927, the same
non-rigorous attempt is repeated (apparently without knowledge of
Kahn, L\"obenstein or Rohn's work). Hilton 1936 devoted a paper
criticizing Donald's work.

The next step is essentially Gudkov's work (yet do not miss what
did Petrovskii 1933/38 though its impact upon the case of sextics
is nearly covered by Hilbert-Rohn). Ultimately Gudkov was able to
prohibit in 1969 and probably much earlier (Gudkov-Utkin 1969
\cite{Gudkov-Utkin_1969/78}) all schemes above the broken line of
Fig.\,\ref{Gudkov-Table3:fig}. This breakthrough originated in
1948 when Andronov suggested (to Gudkov) applying the concept of
{\it roughness} (also known later as {\it structural stability} in
the West since Lefschetz, and adhered to by Thom, etc.) to the
topology of real algebraic surfaces. Petrovskii's advice (1950)
suggested focusing rather on the case of sextic curves. Combining
those novel Russian ideas with the Hilbert-Rohn method, enabled
Gudkov in 1954 \cite{Gudkov_1954} to get solid prohibitive proofs
above the critical line, and even beyond (sic! cf. lilac schemes
$\frac{5}{1}y$, $3\le y \le 5$ on Fig.\,\ref{Gudkov-Table3:fig})
but that turned out to be too massive amputation for the pyramid
to support its own structural mass. Nowadays there are  simpler
proofs from the Arnold-Rohlin era (early 1970's) or via Rohlin's
complex orientation formula 1974--78 which prohibit only a portion
of those (namely those {\it not\/} lying on the continuation of
the lattice by blue rhombs and red circles on
Fig.\,\ref{Gudkov-Table3:fig}). However Rohlin's proof (1972/72
\cite{Rohlin_1972/72-Proof-of-a-conj-of-Gudkov}) of the Gudkov
hypothesis
inhibits all $M$-schemes above the broken line, but the price to
pay is highbrow differential topology \`a la Rohlin from the early
1950's (cf. Sec.\,\ref{Gudkov-hypothesis:sec}). Related work by
Gudkov-Krakhnov/Kharlamov prohibits all the four $(M-1)$-schemes
above the broken line.

What happens under the critical line? Short-cutting a century of
efforts, the answer is rapid: all of them are realized. In fact
all specimens (except the 3 lilac-colored ones) are easily
construct by (slight variants) of Harnack's and Hilbert's method.
At least this is what we read in Gudkov's survey 1974
\cite{Gudkov_1974/74}, yet the cases of $\frac{4}{1}5$,
$\frac{3}{1}5$ are a bit tricky (but see our
Fig.\,\ref{HarnaGudkov4-15:fig} and
Fig.\,\ref{HarnaGudkov3-15XXL:fig} resp.). The three remaining
ones $\frac{5}{1}5$, $\frac{5}{1}4$, $\frac{5}{1}3$ needed to wait
until Gudkov's trick (1969--1973) of using some Cremona
transformations. Beware yet that historically, the very first
argument of Gudkov's Thesis 1969 was a pure existence proof along
the line of Hilbert-Rohn's method (ca. 20 pages long and extremely
hard-to-follow according to Russian experts, cf. Polotovskii 1996
\cite{Polotovskii_1996-D-A-Gudkov}, Viro, etc.), and was not
constructive at all. Remind also that Gudkov himself at some early
stage, in 1954, asserted incorrectly inexistence of those 3
difficult birds, quite in line with Hilbert's intuition at the
Paris Congress of 1900. As like to emphasize Arnold, it seems that
Petrovskii himself was at first very skeptical about the twist
taken by Gudkov's solution.

\subsection{Harnack's and Hilbert's constructions}

\begin{figure}[h]
    \hskip-2.3cm\penalty0
    \epsfig{figure=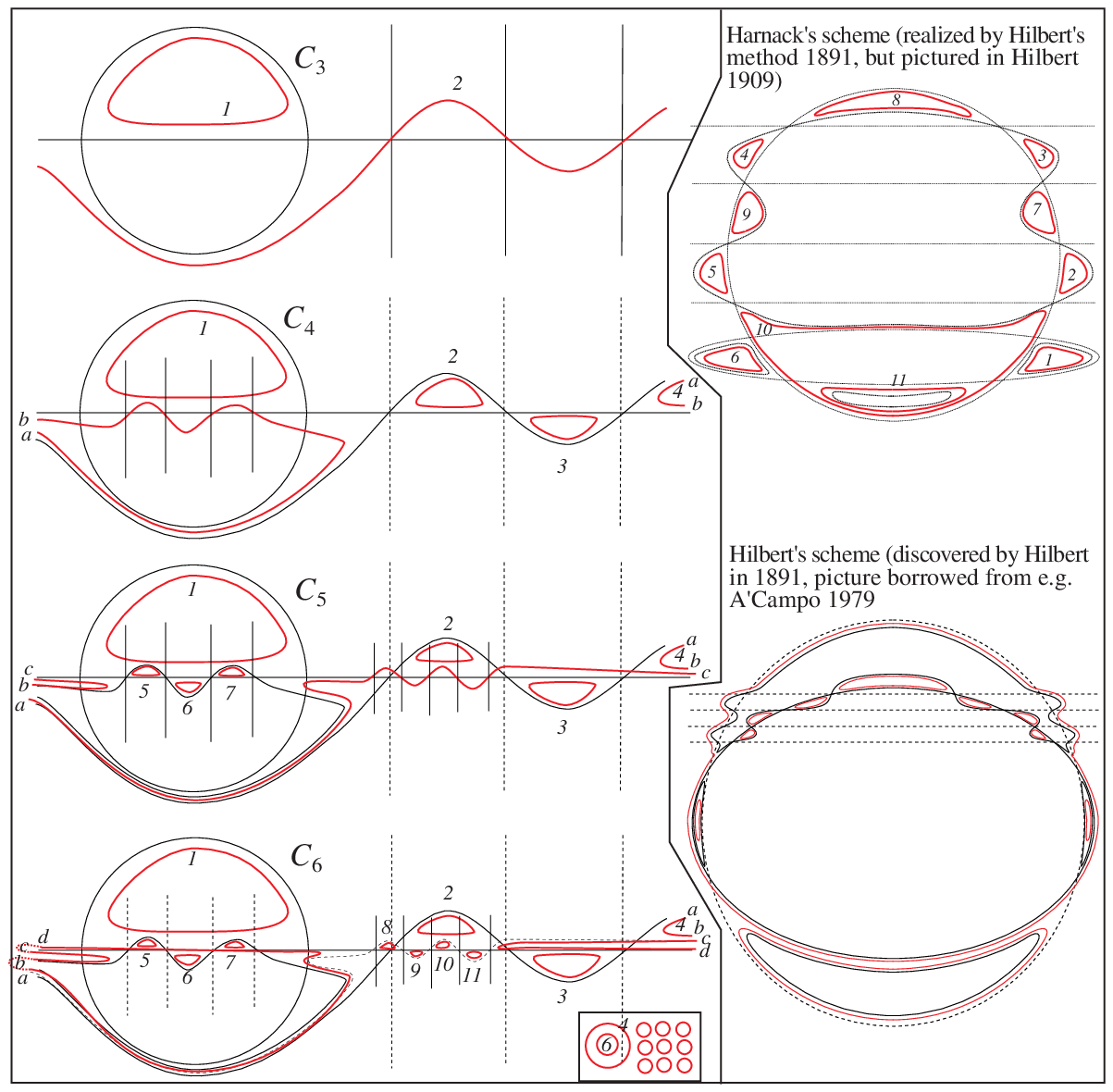,width=165mm}
\vskip-5pt\penalty0
\caption{\label{Harnack-original:fig}%
Left: Harnack's original construction (after Viro 2008); Right:
Hilbert's method, yielding the (new) scheme $\frac{9}{1}1$}
\vskip-5pt\penalty0
\end{figure}

First Fig.\,\ref{Harnack-original:fig}(left) recalls Harnack's
method of construction (trying to keep reasonably close to the
original 1876 \cite[p.\,195]{Harnack_1876}, but making more
explicit pictures (assisted by Viro 2008
\cite[p.\,188]{Viro_2008-From-the-16th-Hilb-to-tropical}). The
top-right of Fig.\,\ref{Harnack-original:fig} reproduces Hilbert's
more
expeditious way to realize this scheme. The bottom-right is the
new scheme discovered by Hilbert in 1891
\cite{Hilbert_1891_U-die-rellen-Zuege} (yet no pictures until
Hilbert 1909 \cite{Hilbert_1909-Ueber-die-Gestalt-sextic}, who
traces only the top-right picture, whose scheme is Harnack's). If
The bottom nice picture is borrowed from A'Campo 1979
\cite[p.\,08--09]{A'Campo_1979}.

Once all the knowledge synthesized in Gudkov's table
(Fig.\,\ref{Gudkov-Table3:fig}) is understood (or admitted for
short) one  starts appreciating Rohlin's maximality claim. Indeed
having listed all real schemes (there remains $68-12=56$ many
below the broken line) it is an easy matter to spot maximal
elements in the lattice. We find 6 types red-circled on
Fig.\,\ref{Gudkov-Table3:fig} corresponding to

$\bullet$ the three $M$-schemes $\frac{9}{1}1$, $\frac{5}{1}5$,
$\frac{1}{1}9$ of Hilbert, Gudkov, Harnack respectively,

$\bullet$ plus two $(M-2)$-schemes namely $\frac{6}{1}2$ and
$\frac{2}{1}6$,

$\bullet$ and finally one $3$-schemes $(1,1,1)$ corresponding to
the deep nest.

Rohlin's assertion is that those (distinguished) 6 schemes are
precisely those which are definite of type~I (i.e. universally
orthosymmetric). Of course the assertion is trivial for the 3
possible $M$-schemes (since Harnack 1876 or via Klein's 1876
intrinsic proof of Harnack's inequality via the topology of
Riemann surfaces). The deep nest of weight 3 is likewise trivially
of type~I, for it is enough to sweep it out by a total pencil of
lines. It remains thus to analyze the two $(M-2)$-schemes with
$r=9$, i.e. $\frac{6}{1}2$ and $\frac{2}{1}6$ by showing that they
are definite of type~I. At first one can imagine to prove this via
pencil of conics (or maybe cubics pencils?). A complete argument
must be given in Rohlin. {\it Insertion.}---[28.03.13] In fact
Rohlin claimed a proof which is now lost via pencils of cubics, so
that there is strictly speaking presently still only one known
proof which involves a congruence modulo 8 due to
Rohlin-Kharlamov-Marin (\ref{RKM-congruence-reformulated:thm}).
Le~Touz\'e 2013
\cite{Fiedler-Le-Touzé_2013-Totally-real-pencils-Cubics} was able
to validate the total reality assertion of Rohlin, yet only after
supposing the curve dividing. It should however not be impossible
that methods of Le~Touz\'e suitably modified establish the full
Rohlin claim. This seems to be an urgent problem to deal with. We
personally tried a lot but failed dramatically. This sort of
problem seems to require extreme cleverness.

Rohlin 1978 \cite{Rohlin_1978} states the following theorem
summarizing all those efforts (up to Gudkov plus his own input
reconciliating with Klein's viewpoints):

\begin{theorem}\label{Rohlin-type:thm} {\rm (Rohlin 1978)}
The $56$
possible real schemes for sextics (Harnack, Hilbert, Rohn, Gudkov)
split as follows according to Klein's types:

$\bullet$ There are $6$ schemes of type~I (red-circles on Gudkov's
table=Fig.\,\ref{Gudkov-Table3:fig});

$\bullet$ There are $42$ schemes of type~II (green-squares on
Gudkov's table=Fig.\,\ref{Gudkov-Table3:fig});

$\bullet$ There are $8$ schemes of indefinite(=mixed) type
(blue-rhombs on Gudkov's table=Fig.\,\ref{Gudkov-Table3:fig}).
\end{theorem}

Note that the resulting distribution of types to be nearly
symmetric (on Fig.\,\ref{Gudkov-Table3:fig}), modulo some anomaly
at the place $5$ (five unnested ovals). A similar remark is made
in Fiedler 1981 \cite[p.\,13]{Fiedler_1981}: ``{\it Bemerkung.
Eventuell ist die Tabelle nicht vollst\"andig. Aber es ist schon
ersichtlich, da{\ss} die Tabelle der zerteilenden Kurven im
Unterschied zur Tabelle aller existierender Typen von
singularit\"atfreien Kurven sechster Ordnung (vgl. {\rm
[1](=Gudkov 1974 \cite{Gudkov_1974/74})}\/) nicht symmetrisch
ist.}''. Indeed the asymmetry we (and Fiedler) notice (but of
course implicit in Rohlin's survey) is the (unique) symmetry
breaking occurring between the schemes $\frac{4}{1}$ and the
scheme ${5}$. The latter turns out to be of type~II, as it cannot
satisfy Rohlin's formula (\ref{Rohlin-formula:thm}), whereas the
former scheme is easily seen to be indefinite (cf.
Fig.\,\ref{R4-1:fig}). Note that Arnold's congruence $5-0=p-n=k^2
\pmod 4=3^2=9=5 \pmod 4$ is not fine enough to detect this break
of symmetry.

Let us try to understand this spectacular statement of Rohlin
(\ref{Rohlin-type:thm}).

(1) From the easy Klein congruence $r\equiv g+1 \pmod 2$ if
type~I, we draw that all schemes with an even number $r$ of ovals
belong to type~II (this explains all the green-squares at heights
$r=0,2,4,6,8,10$, cf. Fig.\,\ref{Gudkov-Table3:fig}).

(2) As already explained all $M$-schemes (here $r=11$) are
trivially of type~I (Harnack's inequality or Klein's argument of
1876 \cite{Klein_1876}). For another reason the scheme $(1,1,1)$
(deep nest of profundity 3) is easily shown to be of type~I (total
pencil of lines).

(3) For similar reasons (but deeper) is the assertion that the 2
circled schemes with $r=9$ belongs to type~I. This is truly the
work of Rohlin, albeit philosophically akin to
Klein-Teichm\"uller-Ahlfors' total reality. {\it
Insertion}---[28.03.13] As just said this is still unproven
synthetically, and the only proof available involves deep
differential topology (Rohlin-Kharlamov-Marin).

(4) Appurtenance to the indefinite type is usually easy requiring
merely exhibiting two curves, one in each type. So for instance
the scheme $9$ (consisting of $9$ outer ovals without nesting) is
indefinite (cf. our Fig.\,\ref{KleinRo-sextic:fig}). Of course
here the basic theoretic tool is Fiedler's observation that the
type is governed by the smoothing effected in the
Pl\"ucker-Klein-Brusotti method of small perturbation. Full
details are worked out in Sec.\,\ref{indefinite-types:sec}.

(5) Another piece of information (now purely Rohlinian) is
Rohlin's inequality $r\ge m/2$ for a smooth plane curve of degree
$m$. (Remind this to follow for Rohlin's formula, in turn derived
by a intersection theory argument of halves of the dividing curve
capped off by real ovals and brought into general position by
perturbing via a vector field normal to the real locus).
Conceptually this involves Poincar\'e homology theory, and the
allied intersection theory (e.g. by Lefschetz, etc.). From
Rohlin's inequality, one deduces that the scheme with $r=1$ is of
type~II.

(6) Using the stronger Rohlin formula, one must be able to treat
all schemes with $r=3$ to belong to type~II (except the deep nest)
and likewise assess type~II for all other schemes. {\it
Insertion\/} [28.03.13].---Yes this is essentially true. More
precisely Rohlin's formula admits the Arnold congruence as
corollary (cf. \ref{Rohlin-implies-Arnold:lem}), and the latter
$\chi\equiv_4 k^2=9\equiv 1$, forces a curve of type~I to live on
the grid formed by blue rhombs (and red-circles) of
Fig.\,\ref{Gudkov-Table3:fig}. So Rohlin's assertion is evident
safe for the scheme $5$, $\frac{1}{1}1$ and $1$. But all those
cases are prohibited by Rohlin's formula $2(\pi-\eta)=r-k^2$.
Indeed in case of no-nesting Rohlin's formula reduces to
$r=k^2=9$, hence rules out the schemes $1$ and $5$. For
$\frac{1}{1}1$, we have only one pair so $\pi+\eta=1$, while
Rohlin's formula says $\pi-\eta=-3$, whence $2\pi=-2$, which is
impossible as $\pi$ is a cardinal (namely the number of positive
pairs).

At this stage the proof of Rohlin's theorem
(\ref{Rohlin-type:thm}) is complete.

\subsection{What can be proved via Ahlfors?}

[17.01.13] From the viewpoint of our survey, it is of some
interest to decide which results of the theory aroused from
Hilbert's 16th problem (Hilbert-Rohn, etc. and the Russian school
Gudkov-Arnold-Rohlin, just to quote the 3 supermassive black
holes) can be (re)proved via the Ahlfors map.

As we said already all this section was actually motivated by the
guess that Ahlfors could be used to prove the still unsettled
Rohlin maximality conjecture (at least what remains thereof post
Shustin 1985 \cite{Shustin_1985/85-ctrexpls-to-a-conj-of-Rohlin}).
However as yet we failed to complete this grandiose project.

Another more didactic aspect (yet perhaps not to be neglected as a
first step toward subsequent progresses) would be to see if
Ahlfors implies the (Gudkov-)Arnold congruence mod 4:
$\chi=p-n=k^2 \pmod 4$ for dividing curves of degree $2k$.

{\it Insertion} [02.04.13].---This game looks quite artificial
since Arnold's congruence is, e.g., a fairly trivial consequence
of Rohlin's formula, cf. (\ref{Rohlin-implies-Arnold:lem}).

The method would be to examine the Ahlfors foliation, i.e. that
induced by the total pencil of curves while trying to apply
Poincar\'e(-Bendixson-Kneser-Hamburger) index formula for
foliations. Recall that $\chi$ is the Euler characteristic of the
``Ragsdale-Petrovskii'' (orientable) membrane of ${\Bbb R}P^2$
bounding the ovals. Of course it looks hard for Ahlfors to beat
the elegance of Arnold's argument based on intersection theory and
the divisibility by 8 of the signature of spin manifold (as
prompted by the algebra of integral quadratic symmetric form).
However, as just mentioned, there is an even simpler proof of
Arnold based on Rohlin's formula.

Of course there is a myriad of sub-Arnoldian truths that could be
treated via the Ahlfors foliation, e.g. Hilbert-Rohn prohibition
of an $M$-sextic without nesting, or the type~II of a quartic with
two unnested ovals (all these assertions being implied by Arnold's
congruence).

Another game out of reach to Arnold, but proved via Rohlin's
formula is the prohibition of the sextic scheme $5_I$ of five
unnested ovals in the type~I case. This could perhaps also be
proved via the Ahlfors map.

The only point which we managed (presently) to prove via Ahlfors
is the (easy sense) of Klein's Ansatz that a dividing curve cannot
gain an oval while crossing the discriminant through a solitary
node (with a complex conjugate pair of tangents). Compare for this
Lemma~\ref{Klein-via-Ahlfors(Viro-Gabard):lem} suggested by a
letter of Viro. However, we always use as a premiss the issue that
for a plane curve the abstract total map of Ahlfors extends to the
ambient projective plane. We should acknowledge a letter of Marin
(cf. Sec.\,\ref{e-mail-Viro:sec}) for having made us aware of this
subconscious short cut. We still hope this to be true via basic
algebraic geometry, of which we forgot all the foundations.

Klein's Ansatz (1876) can also be proved without Ahlfors by using
some Picard-Lefschetz and Dehn stuff, or rather just some
``Anschauung'' that might have been folklore as early as 1876.

Here is an argument (cf. also the next
Sec.\,\ref{Klein-Marin:sec}). At the level of the
complexification, one can only explain the apparition of a
solitary node as the strangulation of some vanishing cycle $\beta$
on the Riemann surface. Then we analyze all possibilities. By the
reality of our deformation, the cycle $\beta$ must be invariant
under conjugation $\sigma$, hence either be a real circuit (or
``oval''\footnote{Alas the word ``oval'' is quite ambiguous, as it
is either just a real component of the abstract curve, or
sometimes used in the much more specific sense of a component of a
plane curve which is null-homotopic (equivalently bounds a disc)
in ${\Bbb R}P^2$. Whenever we use the term oval in this abstract
sense, due to a lack of better synonym (in German there is a good
one ``Zug/Z\"uge''), we write it ``entre guillemet'' (=inverted
comma in English, according to my Dictionary).}) (pointwise
invariant under $\sigma$), or an ortho-cycle (2 fixed points under
$\sigma$) or a dia-cycle (no fixed point under $\sigma$). This is
exhaustive via the classification of involutions on the circle,
which via the quotient map and covering theory, reduces to the
classification of one-dimensional manifolds (Hausdorff and
metric).

As we already noted a dia-cycle cannot exist in the orthosymmetric
case (Lemma~\ref{antioval:lem}). For an ``oval'' it can indeed
shrink to a point (hence a solitary node) but then disappear of
course (cf. Fig.\,\ref{Klein-Marin:fig}, right). If we have an
ortho-cycle then two cases are to be distinguished. It can either
cross two distinct ``ovals'', in which case both ovals merges
together after the Dehn twist (cf. Fig.\,\ref{Klein-Marin:fig},
left). The last possibility is an ortho-cycle cutting only one
``oval''. In this case Fig.\,\ref{Orthoovals2:fig} shows that $r$
stays constant, and of course we do not cross a solitary node in
that case, but rather a non-isolated one with 2 real tangents.

This ``proves'' Klein's Ansatz  (modulo some Picard-Lefschetz
theory), and even Marin's stronger assertion that a dividing curve
cannot increase its number component when crossing the
discriminant. Indeed in all 3 cases analyzed, either $r$ drops by
one unity (first two cases), or stays constant. Remark however
that Marin 1988 \cite{Marin_1988} has a more conceptual proof.

What is Picard-Lefschetz theory in our context? Since any crossing
of the discriminant can be interpreted as a smooth arc traversing
the discriminant transversally,  we may (in the small) always
replace this little arc by a linear pencil, and  are reduced to
classical Picard-Lefschetz theory, where in our case we have
holomorphic fibration of the plane by a pencil of curves. The
theory in question tell us the geometric monodromy when winding
around a singular member of the pencil, but also gives the Dehn
twist description of what happens when we (more cavalier) cross
frontally the singularity.

Recall that Picard's thesis (the first work of Picard on another
subject) is dated 1879 \cite{Picard_1879}, while Klein's Ansatz
(no proof but probably Klein had one) is dated 1876 (3 years
younger). So clearly our approach is somewhat historically
contorted. Still, it is not impossible that Klein (and many
others) were aware of the geometry behind our argument (via Dehn
twists, ca. 1910). It is also possible that Klein's argument was
closer to Marin's, albeit the latter result is perhaps slightly
different (and of course stronger).

\subsection{Rohlin's conjecture almost implied by Klein-Marin
(Klein 1876, Marin 1988, Viro 1986)} \label{Klein-Marin:sec}

[04.01.13] This section presents another tactic (pseudo-proof) of
Rohlin's conjecture that probably everybody had in mind
(especially Rohlin and Marin), yet nobody write it down as it
fails blatantly. Although being a ``pot-pourri'' it is worth
presenting as it helps clarifying the relation (or absence
thereof) between Klein's original assertion 1876 \cite{Klein_1876}
as interpreted by modern workers, notably Marin 1988
\cite{Marin_1988} and Viro 1986/86 \cite{Viro_1986/86-Progress}
(the latter being based on a ``private communication'' of the
former).

First there is
a remarkable observation of Klein 1876 that Marin 1988
\cite{Marin_1988} was probably the first to supply with a proof.
To be perfectly accurate we believe that Marin's result is
slightly stronger than Klein's original
asserting only that
a dividing curve cannot acquire a new oval  like a champagne
bubble emanating from a solitary node (compare Klein's Quote
\ref{Klein_1876-niemals-isolierte:quote} especially the
phraseology ``isolierte reelle Doppeltangente''). In Marin 1988
article, full credit is ascribed to Klein, either by over-modesty
or because
Marin overlooked to notice the little nuance between his and
Klein's weaker assertion. (Compare the recent e-mail exchanges in
Sec.\,\ref{e-mail-Viro:sec}.)

\begin{lemma} \label{Klein-Marin:lem} {\rm ($\approx$
Klein 1876, but
in the formulation of Marin 1988)}.---A (plane) dividing curve
cannot increase its number of ovals when crossing a  node
(non-degenerate
double point).
\end{lemma}

For Viro (1986 {\loccit}) the curve does not actually need to be
plane.

\begin{proof}
Perhaps Klein gained evidence from the case of quartics. Imagine a
G\"urtelkurve (quartic with 2 nested ovals), then there cannot be
created a new oval without violating B\'ezout. Hence either both
ovals
amalgamate
or the inner oval evanishes. In both cases the number of ovals
decreases (by one unit). [02.04.13] This is not an exhaustive
discussion, for there can be also an eversion
(Sec.\,\ref{Eversion:sec}), keeping $r=2$ constant.

Bringing into the picture the Riemann surface (of orthosymmetric
type) underlying the dividing curve, then, as the latter traverses
the discriminant (at some smooth point of it)  our curve becomes
uninodal via some vanishing cycle pinching the Riemann surface.

\begin{figure}[h]
\vskip-0.2cm\penalty0
\hskip-1.2cm\penalty0 \epsfig{figure=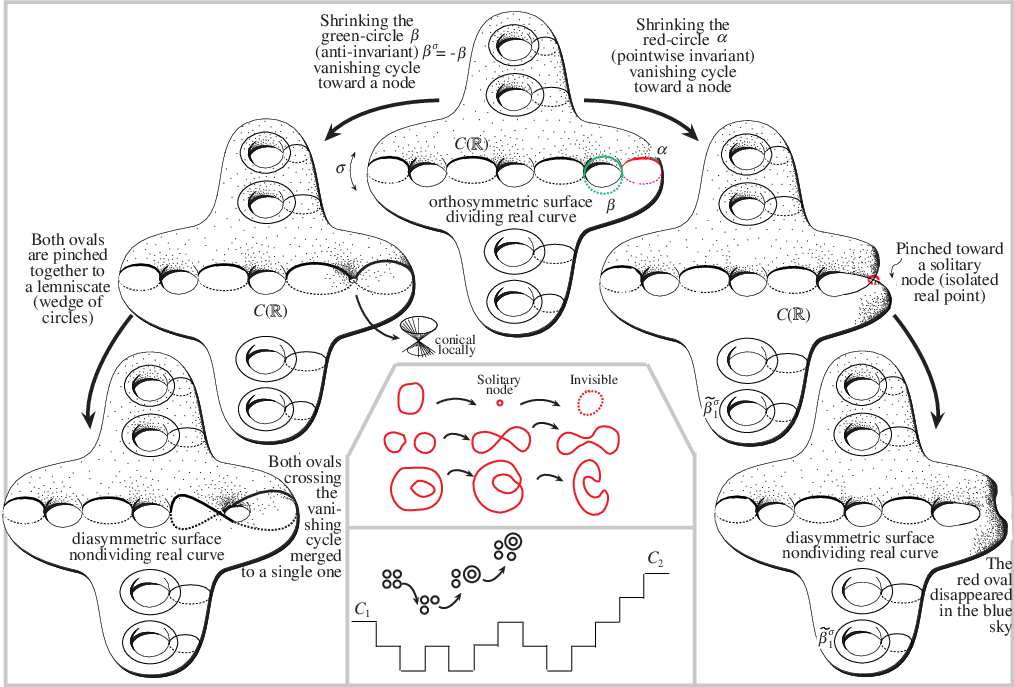,width=147mm}
\vskip-5pt\penalty0
  \caption{\label{Klein-Marin:fig}%
  Two possible real degenerations of an orthosymmetric
  surface toward an uninodal Riemann surface,
  and the blue-sky catastrophe occurring right after having
  crossed the singularity to reach the diasymmetric type with one
  less ovals(interpretable as Dehn twists I presume?)
  For a 3rd option cf. Fig.\,\ref{Orthoovals2:fig}.} \vskip-5pt\penalty0
\end{figure}

This vanishing cycle can actually be an oval  of the Riemann
surface: once shrunk to a point it disappears and one oval
gets lost (cf. Fig.\,\ref{Klein-Marin:fig}, right). This is not
the only possibility as shown by the G\"urtelkurve whose 2 ovals
may coalesce. It is just a little harder to visualize the
corresponding surgery on the Riemann surface. The key is to
imagine an anti-invariant vanishing cycle $\beta$ whose
contraction is depicted on Fig.\,\ref{Klein-Marin:fig}, left. The
two ovals traversed by the cycle $\beta$ have merged together.

{\it Insertion} [02.04.13].---Further there is a 3rd possibility,
of when the ortho-cycle $\beta$ intersects only one oval (cf.
Fig.\,\ref{Orthoovals2:fig}). Then the corresponding Morse surgery
is an eversion keeping $r$ constant but destroying the dividing
character of the curve.
\end{proof}

\begin{proof} (Marin's proof in the dirty fingers of
Gabaredian\footnote{This appellation is now a common joke in
Geneva, based on a mixture of the writer's name with themuch more
eminent Paul Garabedian, the notorious student of Ahlfors,
Schiffer, Bergman, who seems to have played a pivotal r\^ole in
the ultimate shape of Ahlfors theorem, as published in 1950
\cite{Ahlfors_1950}.}) Marin's proof is somewhat different and in
substance as follows (please refer to the French original for
faithfulness). The initial curve is orthosymmetric. Such curves
satisfy the Klein congruence $r\equiv g+1 \pmod 2$. On the other
hand when traversing the discriminant the curve is uninodal and
the real part undergoes a ``Morse surgery'' which alter the number
of ovals by one unit.

{\it Insertion} [02.04.13].---Warning. Possibly $r$ can stay
constant in case of an eversion, cf. Fig.\,\ref{Orthoovals2:fig},
but then the post-critical Riemann surface becomes diasymmetric.

At any rate, the new curve (past the discriminant) is necessarily
diasymmetric, either by Klein's congruence when $r$ moves by one,
or by the Dehn-twist argument of Fig.\,\ref{Orthoovals2:fig} when
$r$ is kept constant.
Marin concludes  by arguing that a path between two conjugate
points avoiding the real locus subsists in all nearby curve. (Alas
I confess to  have not  properly understood this argument which is
presumably much superior to the above via vanishing cycles.)
Compare also Marin's e-mail, where he explained us more details.
(NB: Marin's argument is also repeated in Degtyarev-Kharlamov 2000
\cite[p.\,785, 4.6.8]{Degtyarev-Kharlamov_2000}, and was
apparently always easily digested by  Russian scholars, Viro
included.)
\end{proof}


Armed with the lemma, let us try to attack Rohlin's conjecture.

\smallskip
\begin{proof} (Pseudo-proof of Rohlin via Klein-Marin)
Suppose $S_1$ to be a scheme of type~I which is not maximal, say
embeddable in $S_2$. Take algebraic models $C_i$ of each $S_i$
($i=1,2$). By general position we may assume that the line $L$
through $C_1$ and $C_2$, in the hyperspace of all curves of degree
$m$ (\`a la Cayley, etc.), crosses transversally the discriminant
hypersurfaces (in smooth points of it) at uninodal curves. This
unique node is necessarily real (when we look at real members of
the pencil $L$). Hence when we join $C_1$ to $C_2$ we get real
curves (finitely many of them being singular). Whenever we cross
the discriminant the real locus undergoes a ``Morse surgery''
which is (up to reversing time) is either

$\bullet$ the death of an oval (shrinking to a point)

$\bullet$ the fusion of two unnested ovals

$\bullet$ the fusion of two nested ovals.

Each operation effects a fluctuation of $\pm 1$ on the number $r$
of ovals. (Warning [02.04.13].---This is not even true due to
eversions!) So we can imagine a staircase starting from $C_1$ to
$C_2$ recording the history of the varied fluctuations of $r$
during the transition from $C_1$ to $C_2$ along the pencil
$\lambda C_1 + \mu C_2=0$ (cf. Fig.\,\ref{Klein-Marin:fig},
center-bottom). By Klein-Marin the first staircase is moving
downwards, and as we ultimately reach $C_2$ having more ovals, we
are naively inclined to claim that we shall revisit the same
scheme $S_1$ at some step after which $r$ only increases. This
would be true if a scheme would be completely encoded by its
number $r$ of circuits. In this naive world, we get an
intermediate curve $C_1'$ also representing the scheme $S_1$
(hence dividing since $S_1$ is of type~I) and after which $r$ only
increases. This would violate the Klein-Marin theorem.
\end{proof}

The moral is that Klein-Marin seems to imply, but fails implying,
the Rohlin maximality conjecture (for an explicit objection see
the little diagrammatic of ovals on Fig.\,\ref{Klein-Marin:fig},
bottom). Nonetheless, the Klein-Marin lemma certainly implies the:

\begin{cor}
The chambers
past the discriminant corresponding to orthosymmetric curves are
local maxima of the function $r$ counting the number of real
circuits. Further all chambers adjacent to an orthosymmetric
chamber are diasymmetric.
\end{cor}

\begin{proof}
The last assertion follows directly from Klein's congruence
$r\equiv g+1 \pmod 2$.
{\it Insertion} [02.04.13] This is true when $r$ varies (by one
unit), but if it stays constant one has to invoke the Dehn-twist
argument of Fig.\,\ref{Orthoovals2:fig}.
\end{proof}

Hence orthosymmetric chambers are never contiguous (along a wall
of codimension 1), but a priori they could still have closures
with non-void intersections.

\subsection{Back to degree 6: Rigid isotopy (Nikulin 1979 via K3's,
Torelli of Pyatetsky-Shapiro-Shafarevich 1971)}

[05.01.13] Even if Rohlin's conjecture (type~I $\Rightarrow$
maximal) looks out of reach, it might be easier in degree 6 (we
mean by a theoretical argument independent of Rohlin's census). In
that case granting orthosymmetry of the schemes $\frac{6}{1}2$ and
$\frac{2}{1}6$ one recovers all of Gudkov's obstructions
(prohibition of the semi-hexagons above those schemes, cf.
Fig.\,\ref{Gudkov-Table3:fig}, safe those that were established by
Hilbert and Rohn, namely the schemes $11$ and $\frac{10}{1}$). Of
course this is nothing new, yet methodologically distinct from the
topological arguments \`a la Arnold-Rohlin explaining the Gudkov
hypothesis. So we are asking for a fighting interplay between pure
geometry and topology.

Also in view of the Morse surgery inherent in the Klein-Marin
theorem, one can ask several questions about the contiguity of
chambers in the space of all sextics and correlate this with the
diagrammatic of Gudkov's table (Fig.\,\ref{Gudkov-Table3:fig}).
The first basic  point is that when we cross a wall $r$ fluctuates
by $\pm 1$. Hence we do not have complete freedom to random-walk
on the triangular lattice underlying Gudkov's table (all
horizontal edge cannot be used).

{\it Insertion} [02.04.13].---This is a naive misconception, since
in fact there is also eversion (cf. Sec.\,\ref{Eversion:sec})
keeping $r$ constant!

Define the {\it distance\/} between two chambers as the minimum
number of walls needed to be crossed to join them (by a path
transverse to the discriminant). Another ``distance''
is defined by restricting
to pathes along (linear) pencils of curves.

A great miracle (specific to order 6)
is the following result due to joint efforts of Kharlamov and
Nikulin 1979/80 \cite{Nikulin_1979/80}:

\begin{theorem} {\rm (Nikulin 1979)}\label{Nikulin:thm}
The real scheme enhanced by the (Klein-Rohlin) type affords a
complete invariant of the rigid-isotopy class of sextics. Thus via
Rohlin's classification (Theorem~\ref{Rohlin-type:thm}) there is
$6+42+2\cdot 8=56+8=64=2^{8}$ ``typed'' schemes and so many
rigid-isotopy classes. (This number being a power of 2 is perhaps
just good fortune? probably because for quartics the  number of
chambers is $6$.)
\end{theorem}

\begin{proof} This is yet another
``tour de force''. It uses (but strangely does note cite!) the
topological classification of Rohlin 1978 \cite{Rohlin_1978}
(making already a fusion between Klein 1876 and Hilbert
1891/1900's 16th problem as solved by Gudkov 1969
\cite{Gudkov-Utkin_1969/78}). But that is not all! It also
combines this with the complex geometry of K3 surfaces
(Kummer-K\"ahler-Kodaira as coined by Weil), especially the
contribution of Pyatetsky-Shapiro--Shafarevich 1971/71
\cite{Pyatetsky-Shapiro-Shafarevich_1971/71} on the global Torelli
theorem, as well as the surjectivity of the period mapping
(Kulikov 1977 \cite{Kulikov_1977}). Rohlin, and especially
Kharlamov's r\^ole in this proof seems to have been quite pivotal
(and acknowledged as a such).
\end{proof}

\subsection{The Gudkov-Rohlin-Nikulin
pyramid and the contiguity graph}

{\it Warning} [02.04.13].---This section is a miscellany of
mistakes about the combinatorial structure of the hyperspace of
all sextics. We kept our text intact in its original shape (modulo
Insertions and corrections) since we think that it is more
important to avoid the basic mistake than to reach the ultimate
verity of what is quite likely to become a combinatorial mess if
pushed to its ultimate perfection. We still encourage the
indulgent reader to follow our output as it may contain
interesting problems. In particular is it possible for a pencil of
sextics to visit only a single chamber? This could be the case if
there is a chamber contiguous to itself. You move in your chamber
and tries to get out of it by traversing a wall, but alas fall
again trapped in the same room as you started where.

[05.01.13] Looking (once again) at Gudkov's table
Fig.\,\ref{Gudkov-Table3:fig} (while duplicating all the
blue-rhombs) we get a complete picture of all chambers past the
discriminant (i.e. rigid-isotopy classes of smooth curves). One
would like to understand their contiguity relation to enhance this
set into the {\it contiguity graph}.

\begin{figure}[h]
\hskip-1.2cm\penalty0
\epsfig{figure=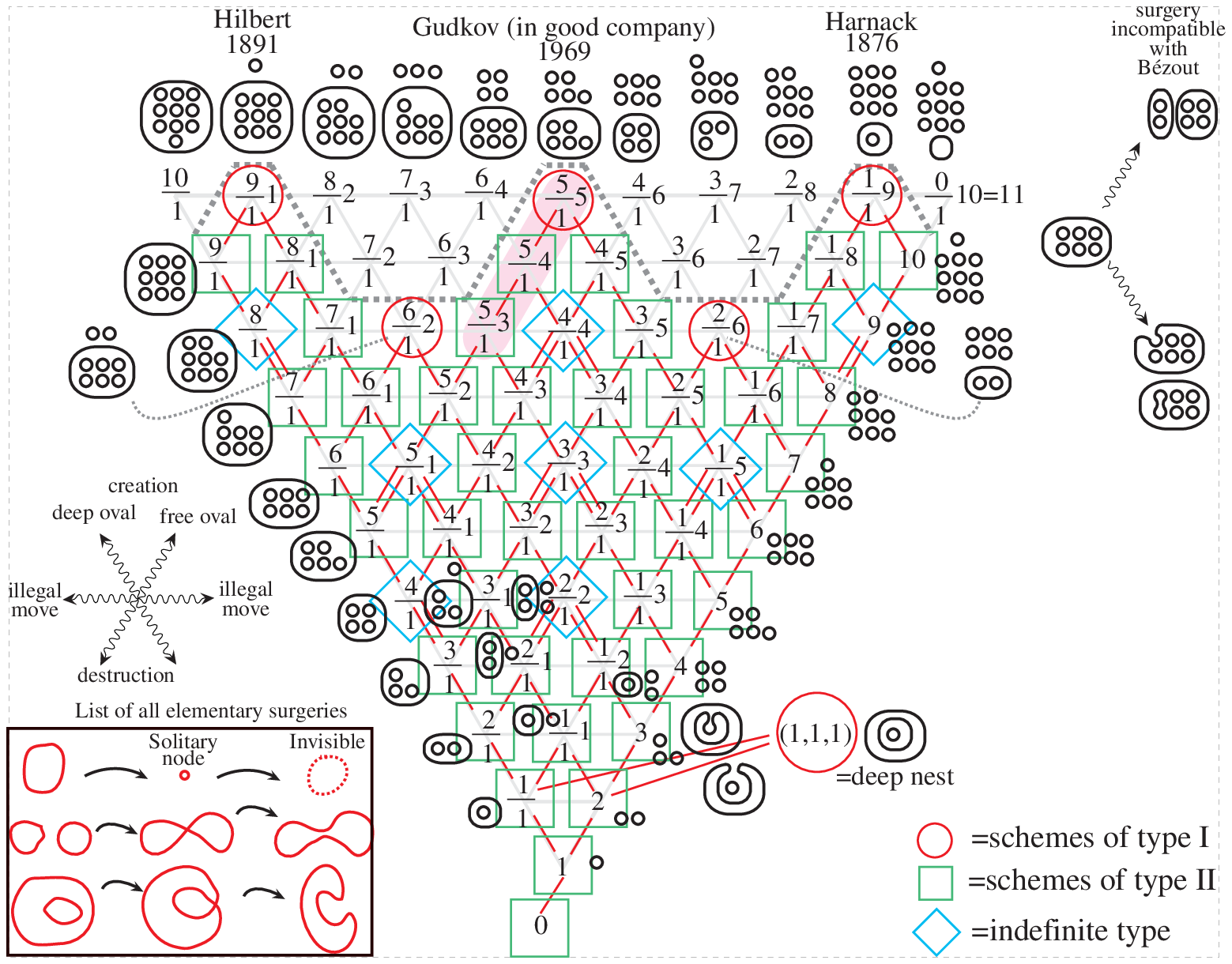,width=147mm} \vskip-5pt\penalty0
  \caption{\label{Gudkov-contig:fig}%
  The Gudkov-Rohlin table seen as the
  chambers of all rigid-isotopy classes via Nikulin's theorem
  (1979/80 \cite{Nikulin_1979/80})} \vskip-5pt\penalty0
\end{figure}

Basically we have 6 moves prompted by the equilateral lattice
underlying Gudkov's table. But as all Morse surgeries amounts to
the creation or destruction of an oval we can rule out the two
horizontal moves (keeping $r$ unchanged and corresponding resp. to
the evasion or encapsulation of an oval).
Next certain Morse surgeries are of course incompatible with
B\'ezout (e.g. that depicted on the top-right of
Fig.\,\ref{Gudkov-contig:fig}). A little moment thought shows that
all admissible Morse surgeries correspond to one of the $4$ legal
moves.

A further obstruction comes from the Klein-Marin theorem
(\ref{Klein-Marin:lem})
impeding the $2$ creationist moves (going up $r\mapsto r+1$) as
soon as the chamber is of type~I.

Naively one is tempted to say that this is a complete list of
legal moves. If so is the case then we would have a complete
description of the contiguity graph (whose edges are depicted by
red strokes on Fig.\,\ref{Gudkov-contig:fig}).
As we shall see later (eversions) this answers is quite unlikely
to be the definitive answer.

(Exercise: count the number of edges of this graph: counting edges
going to the North-West: we have $2\cdot 9+6+3\cdot 5+2+3\dot
1=18+6+15+2+3=44$. This must be doubled by symmetry to $88$. Then
adding those edges going to the $8$ indefinite types adds $2\cdot
8-3=13$ edges. Next the empty scheme $0$ gives one edge and the
deep nest for 2. In definitive, $88+13+3=104=2\cdot 52=2^2\cdot
26=2^3 \cdot 13$ edges. This is the number of contiguity zones of
the discriminant hypersurface.)

It seems also that the most remote pair of vertices are the
schemes $\frac{9}{1}$ and $10$ lying at distance 18 apart. On the
other hand the discriminant has degree $\delta=3(m-1)^2$. (This
can be proved via a Euler characteristic count in the fibration
induced by a pencil of $m$-tics after blowing up the basepoints,
and is also implied by the so-called Zeuthen-Segre formula in the
algebro-geometric community, which is merely the avatar of
Riemann-Hurwitz in one more dimension.) For sextics $m=6$ this
gives $\delta=3 \cdot 5^2=75$. It follows that the linear distance
between two chambers (as measured inside a linear pencil) is at
most $[75/2]=[37.5]=37$ (the temperature of the human body).
Probably this bound is far from sharp (except of course if there
is a line hitting the discriminant 75 times on the reals), and one
could try to find a least upper bound.

\begin{defn}
{\rm Given two chambers (=rigid-isotopy classes) define their
distance $\delta$ as the minimum number of wall-crossings
separating them. This is also the combinatorial distance in the
contiguity graph.

$\bullet$ Define also their \'ecart $\varepsilon$ as the minimum
number of wall-crossings in a generic pencil (transverse to the
discriminant) through two curves belonging to the given chambers.}
\end{defn}

\begin{lemma}
We have always $\delta \le \varepsilon$; and when $m=6$, $\delta
\le 18$ and $\varepsilon \le 37$.
\end{lemma}

Since the degree of the discriminant is 75 ($\delta=3(m-1)^2$ is
odd whenever $m$ is even) any pencil of sextics intersects the
discriminant (in a real point) and so the curves undergo at least
one Morse surgery, and assuming genericity there is an odd number
of such surgeries. However the structure of the graph only permits
loops of even length (the girth=systole of the graph is 4). This
is almost a contradiction in mathematics. How to resolve it?

{\it Insertion} [21.01.13].---Just look at eversions, cf.
Sec.\,\ref{Eversion:sec}. Further using eversions it is likely
that $\delta$ is much smaller than  above, and we predict rather
something like  $\delta\le 11$ (cf. \ref{eversion-and-other
surgeries:conj} and the semi-conjectural
Prop.\,\ref{Erdos-number-of-sextics=11:prop}).

Then there is a host of combinatorial-geometric questions arising.
E.g. is there a pencil cutting 75 times the discriminant (a sort
of total reality of B\'ezout). If not what is the maximal number
$\mu$ of real intersections a line can have with the discriminant?
(Since the distance between the extreme $(M-1)$-schemes is 18,
taking the line joining them we get a pencil with $18+18=36$
(aller-retour, no one way ticket!) real intersections, to which
one can safely add one unit due to oddness of the degree, so
$\mu
\ge  37$ the temperature of the human body.)

{\it Insertion} [01.04.13].---Alas this argument is foiled as it
does not take into account eversions. With eversions the maximal
distance seems to be 11 (between $M$-curves and the empty one),
and so arguing as above gives only $\mu\ge 22+1=23$.

What is the least (resp. maximum) number of chambers visited by a
pencil of sextics? Denote them $\alpha$ resp. $\omega$. Naively a
pencil could stay entirely inside a chamber, but this is precluded
by B\'ezout as $\deg \disc=75$ is odd (so $\alpha\ge 2$). A priori
among the 64 chambers all could be visited since the discriminant
has degree 75. Same question for the length of the loop induced in
the contiguity graph by a (generic) pencil. (A priori this length
can be as long as $75$, but not longer.) Is this loop
 always
non-contractible (in the contiguity graph)? Can it be embedded
(i.e. visits only once each chambers it visits)? Is there a pencil
visiting only diasymmetric chambers? (It is evident from Klein's
congruence and surgeries affecting $r$ by $\pm 1$ that a pencil
cannot visit only orthosymmetric chambers.) What is the maximum
number of orthosymmetric chambers visitable by a pencil? (Of
course at least two, take the line spanned by 2 points in two
ortho-chambers, but probably some lucky
Stonehenge
alignment exists.)

What is the minimum height of a pencil? The height being just the
invariant $r$, number of real circuits. An interesting result is
the theorem of Cheponkus-Marin (cf. Marin 1988
\cite[p.\,192]{Marin_1988}):

\begin{theorem} {\rm (Cheponkus-Marin 1988)}
In any generic (linear) pencil of curves of even degree $m> 2$
there is a curve having at most $M-3$ components (\/$r\le M-3$\/).
\end{theorem}

Looking at Gudkov's table (Fig.\,\ref{Gudkov-Table3:fig}) this
looks almost evident, but is not. Since $M-3$ is the  highest line
full of squares, this Marin result implies that a pencil cannot
confine its visits in one of the 3 regions lying above this line.
One can define the depth $d$ of a pencil as the lowest value of
$r$. So Marin's result implies $d\le M-3$. Is this sharp at least
for $m=6$?

Define three chambers as aligned if there is a line hitting them
simultaneously. In view of Gudkov's table enhanced by the
Klein-Marin theorem one sees that there are $3+1=4$ special
chambers which are contiguous to a single chamber (vertices of
valency 1 in the contiguity graph), namely those of type~I with
schemes $\frac{8}{1}$, $9$, $\frac{4}{1}$ as well as the empty
scheme $0$ ($\bigstar$).  It follows that the triad consisting of
any of those four, plus its unique neighbor and any chamber are
aligned. It would be interesting to find a triad of chamber which
are {\it not} aligned.

{\it Insertion} [02.04.14] Again the assertion right before the
($\bigstar$) above, looks foiled due to eversions. It looks more
realist to expect that the empty scheme is the unique chamber
contiguous to a single chamber. It could be interesting to
describe the chambers adjacent to only 2 chambers. By the theory
of eversion (developed latter), those includes the three
$M$-schemes, and 2 Rohlin maximal $(M-2)$-schemes, plus apparently
the 3 orthosymmetric chambers corresponding to symbols on the
``boundary'' of the pyramid, namely $\frac{8}{1}$, $9$,
$\frac{4}{1}$. However it should not include a scheme like $10$,
which by eversion is potentially related to $\frac{9}{1}$ (cf.
Fig.\,\ref{Gudkov-eversion:fig}). The above list could be
exhaustive, but beware that the median schemes in type~I, like
$\frac{4}{1}4$, $\frac{3}{1}3$, $\frac{2}{1}2$ have also only two
connections except for being potentially related to themselves
under eversion. Yet later we shall see that eversion necessarily
destroy the orthosymmetry, so that those schemes are eversively
related to their type~II twins lying below the sheet of paper. So
those schemes have really valency 3. This raises however the
question if a chamber can be contiguous to itself. By the
diagrammatic of all Morse surgeries (eversion included) a
necessary condition is that the chamber lies on the median line of
the Gudkov table. By what as been said (orthosymmetry destroyed by
eversions), the sole candidate for self-contiguity are
$\frac{1}{1}1$ and $1$. Under this phenomenon of self-contiguity
it could be the case that a pencil of sextics stays entirely
within such a chamber safe for a quick perforation of the
discriminant (forced by by B\'ezout), yet bringing us directly
back to the same chamber. In that case the invariant $\alpha$
discussed above could be as low as $1$.

[13.01.13] Another aspect of Nikulin's isotopic classification of
sextics via the Rohlin pyramid is that it affords a broad
generalization of the Nuij 1968 \cite{Nuij_1968} and Dubrovin
1983/85 \cite{Dubrovin_1983/85} theorem stating that the deep nest
schemes represents a unique rigid-isotopy class of curves.
(Actually Nuij's theorem holds in arbitrary dimension.) By
Nikulin's theorem this uniqueness determination by the real scheme
holds true more generally for all sextic schemes which are not
hermaphrodite (i.e. of indefinite type).

All this just amount to the connectivity of the chambers, yet one
may wish to know more on their individual topology (in the large).
One obvious tool is the monodromy representation
$$
\pi_1(\textrm{ some chamber} ) \to {\frak S}(\textrm{ovals})
$$
 acting upon the
ovals by permutation while following a loop inside some fixed
chamber. Now for a scheme having both inner and outer ovals there
is an obvious constraint preventing the permutation to  shuffle
inner ovals with outer ovals. In other words for a scheme of type
$\frac{k}{1}\ell$ the range of monodromy would lye inside ${\frak
S}_k\times {\frak S}_{\ell}$. A (naive?) conjecture would be that
this are the sole restrictions on the monodromy (i.e. the
restricted morphism is epimorphic). If so is the case then all
chambers are not simply-connected, except perhaps the 5 ones
corresponding to the non-permutable schemes, i.e.  $0$, $1$,
$\frac{1}{1}$, $\frac{1}{1}1$ and $(1,1,1)$. For those schemes the
monodromy of ovals is a trivial representation, and so there is no
obstruction for those chambers to be simply-connected. Can one of
those chambers even be contractible? A natural tactic is to ask if
it can be starlike, in the sense of having a special viewpoint
(curve) inhabiting the chamber so that each curve of the same
chamber is accessible by the half-circle of the line joining the
base curve to the ``variable'' one. Some obvious candidate are the
Fermat equations $x^6+y^6=-1$ for chamber $0$ and  $x^6+y^6=+1$
for chamber $1$, yet it is not clear at all if those are
``visibility curve''.

{\it Insertion} [02.04.13].---Much sharper and complete results of
the monodromy of sextics are due to Itenberg 1994
\cite{Itenberg_199X-monodromy-deg-6} extending results of
Kharlamov. We shall come back this this latter.

\subsection{Weak reformulation \`a la Marin-Viro of the
Klein-Rohlin maximality conjecture}

[13.01.13] After some discussions with Marin (12--13 Jan. 2013,
cf. Sec.\,\ref{e-mail-Viro:sec}), the following issue came quite
clear. First let us contemplate once more the Gudkov-Rohlin
pyramid as depicted as the contiguity figure
\ref{Gudkov-contig:fig}. On it we imagine the blue rhombs schemes
doubled with a ``$\Lambda$'' shaped pair of edges raising to the
orthosymmetric chambers (provided not on the periphery of the
pyramid), whereas the diasymmetric chamber have generically a
$X$-shaped quadruplets of edges in the contiguity graph.

We can consider the POSET of all real schemes enhanced by the
type~I/II of Klein. This is essentially what did Rohlin 1978, safe
that instead of declaring indefinite those ``hermaphrodite''
schemes tolerating both type of representatives (type~I and II) we
duplicate those schemes to see them as independent elements. This
amounts considering all Gudkov's symbols decorated by signs $\pm$
telling the ortho/dia-symmetry, and of course only those realized
algebro-geometrically. This is a well-defined finite set of
$64=2^8$ elements. How to define an ordered structure to make it
into a POSET? Answer just as the picture
Fig.\,\ref{Gudkov-contig:fig} suggests, namely a type~I scheme
(alias ortho-scheme) has two legs going down (some leg may be
amputated if the scheme is peripheral), whereas dia-schemes have
two legs (going down) and two arm (going up), except if it lies in
the periphery. For instance the scheme $0$ is maximally amputated
having one arm but no legs.

Of course this order structure looks somewhat ad hoc yet quite in
line with the remarks of Klein and the theorem of Marin 1988,
which is a stronger variation thereof (apparently Marin did not
noticed that his statement looks stronger than Klein's original
statement, compare our discussion in Sec.\,\ref{e-mail-Viro:sec}).
Let me call the purified pyramid this poset. Paraphrasing Rohlin's
census (diagrammatically encoded in Fig.\,\ref{Gudkov-contig:fig})
we plainly have:

\begin{prop}
The maximal elements of the purified pyramid of sextic ortho- and
dia-schemes are exactly the orthosymmetric ones.
\end{prop}

It seems evident that this statement extends trivially to all
degrees as a mere paraphrase of Marin's  theorem (1988
\cite{Marin_1988}). Yet is non trivial to make a picture even for
degree $8$. Yet this is worth trying to depict at the occasion
(cf. Fig.\,\ref{Degree8:fig}).

\subsection{An eclectic proof of Rohlin's conjecture via
Ahlfors}\label{Rohlin-via-Ahlfors}

[04.01.13] Let us summarize the situation. Rohlin in 1978
\cite{Rohlin_1978} advanced (taking some indirect inspiration by
Klein 1876) the bold conjecture that
$$
\textrm{a scheme is of type~I iff it is maximal,}
$$
 in the hierarchy of all
real schemes of some fixed degree.

One sense of the conjecture turned
wrong,  in degree 8 by a conjunction of Polotovskii 1981
\cite{Polotovskii_1981} and Shustin 1985/85
\cite{Shustin_1985/85-ctrexpls-to-a-conj-of-Rohlin} works (see
also the remarks in Viro's survey 1986/86
\cite[p.\,67--68]{Viro_1986/86-Progress}). Namely Shustin showed
existence of a maximal scheme in degree 8  of type~II. So the
``$\Leftarrow$'' implication of Rohlin's conjecture is disrupted.
It remains the hope that the ``$\Rightarrow$'' implication is
correct (still open in 2013):

\begin{conj} {\rm (Rohlin's maximality conjecture---post Shustin)}
Fix any integer $m\ge 1$, and consider only schemes of that degree
$m$. If a real scheme is of type~I, then it is maximal in the
lattice (POSET) of all real schemes.
\end{conj}

This is perhaps a trivial consequence of Ahlfors theorem:

\begin{Pseudo-Theorem} {\rm (Gabard 31.12.12 and 04.01.13)} If a real
scheme is of type~I, then it is maximal (among all schemes of the
same degree).
\end{Pseudo-Theorem}

\begin{proof}
Suppose the given scheme, say $S_1$, to be of type~I. By
contradiction assume it non-maximal so that it embeds in some
larger scheme $S_2$ as a strict subset. But our schemes are
algebraically realized by real algebraic curves say $C_1$ and
$C_2$ (defined over ${\Bbb R}$)
so that the inclusion $C_i({\Bbb R})\subset {\Bbb P}^2({\Bbb R})$
belongs to the respective isotopy classes of $S_i$ ($i=1,2$).
Since $S_1$ is of type~I, $C_1$ is dividing, and thus there is by
Ahlfors 1950 \cite{Ahlfors_1950} a total pencil $\pi$ of auxiliary
curves all of whose real members cut only  real points on $C_1$
(at least as soon as they are mobile). Now $C_2$ has at least one
extra real circuit over $C_1$ (which in fact must be an oval, as
curves of odd order necessarily have a pseudoline). Naively, one
would like to choose any point $p_0$ on $C_2({\Bbb R})$ and let
pass through it a curve of the pencil $\pi$, say $\Gamma_0\ni p_0$
while arguing that this curve has supernumerary intersection with
$C_2$, violating thereby B\'ezout. This works (effortlessly) if we
could assume $C_1({\Bbb R})\subset C_2({\Bbb R})$, but this
corrupts  rigidity of algebraic curves.

We see that Ahlfors nearly implies Rohlin, but some gigantic gap
requires to be filled. Obviously the problem has to be embedded in
some more flexible medium so as to bridge the gap between
algebraic rigidity and softness of isotopy classes \`a la
Hilbert-Rohlin.
\end{proof}

[10.01.13] The little flash  on how to complete the argument came
to me ca. [05h20] in the morning after some too early
waking up. It is as follows. Suppose our curve $C_1$ to be of
type~I. By Ahlfors there is a total pencil of curves. If the
scheme of $C_1$ is not maximal it can be enlarged, so there is a
curve $C_2$ with larger scheme. But $C_1$ is transverse to the
foliation induced by the total pencil, and transversality is a
robust feature (structural stability \`a la Thom, etc.)
Accordingly a small perturbation of $C_1$ towards $C_2$ is still
maximally cut by the curves of the pencil. Propagating this so
forth we see (assuming genericity of the pencil) that the first
Morse surgery decreases the number of ovals. (This was anticipated
by Viro yesterday [09.01.13] (cf. Sec.\,\ref{e-mail-Viro:sec}),
and goes back to Klein 1876.) However one would  more, namely that
$C_2$ cannot have more ovals than $C_1$. To be more precise one
should compare the pencil $L$ spanned by $C_1, C_2$ to the total
pencil for $C_1$, while understanding perhaps the filmography of
the deformation in reference to this foliation. Thinking of the
latter as locally vertical,  the first Morse surgery is like the
hyperbola $(x-y)(x+y)=x^2-y^2=\varepsilon<0$ transverse to the
vertical foliation while degenerating to the pair of lines of
slope $\pm 1$ and then becoming another hyperbola
$x^2-y^2=\varepsilon>0$ no longer transverse to the vertical
foliation. (In fact this is only the scenario of when the node is
not a solitary one.) The effect of crossing the first critical
level is that some member of the pencil loose their total reality.
Yet not all of the total reality is lost. In fact merely an
interval of ``imaginariness'' is inserted.

Pushing the analysis in the large (several Morse surgeries), while
also treating the other case one may hope to ensure that when
$C_2$ is reached still some totally real curve persists in the
pencil.

In fact Ahlfors theorem (only) implies Klein's thesis (cf. Klein's
Quote~\ref{Klein_1876-niemals-isolierte:quote}):

\begin{lemma} \label{Klein-via-Ahlfors(Viro-Gabard):lem}{\rm (Klein 1876, in a presentation of Gabard
inspired by Viro, while using Ahlfors)} When a dividing curve
crosses the discriminant it cannot acquires a solitary double
point. In particular all (discriminantal) walls bounding an
orthosymmetric chamber correspond to nodes with real tangents.
\end{lemma}

\begin{proof} (This was
anticipated the [10.01.13] by Oleg Viro (e-mail communication in
Sec.\,\ref{e-mail-Viro:sec}), but I only understand it now after
ca. 10 hours of delay and some sleep in between.)

Imagine the curve moving, with suddenly, a solitary double point
appearing in the real locus, like an ufo(=unidentified flying
object) raising into
the blue sky. If not believing in extraterrestrial flying saucers,
think of this as the ex-nihilo creation of a champagne bubble from
a point inflating slightly to a little oval. The initial curve
$C_{-1}$ is dividing, hence admits by Ahlfors 1950
\cite{Ahlfors_1950} a total pencil of curves\footnote{[30.03.13]
Strictly speaking I do not know how to proof this, but hope this
to be a triviality of basic algebraic geometry.}. The later
induces a (mildly singular) foliation of the real projective
plane, which we call the {\it Ahlfors (or total) foliation}. It
can be assumed transverse to the given curve $C_{-1}$. A priori
the (Ahlfors) foliation may be singular along the curve, yet upon
dragging away the center(s) of perspectives it should always be
possible to avoid this (compare upper row of Fig.\,\ref{Tube:fig}
for an implementation on the G\"urtelkurve). W.l.o.g. suppose the
curve (already) close to the discriminant and the deformation
$C_t$ ($t\in[-1,1]$) to be  a small one traversing that
hypersurface.
By continuity, transversality to the Ahlfors foliation persists
after the critical level, and so a nascent champagne bubble would
violate B\'ezout. Indeed, passing a curve of the total pencil
through a point inner to  the newly created oval (or on that oval)
gives an excessive intersection with the post-critical curve
$C_{+1}$.

Formalizing requires to fix a tubular neighborhood of the initial
curve while noticing that its product structure is the trace of
the total foliation. If the degree $m$ is odd then there is one
pseudoline whose tubular neighborhood is a twisted bundle
(M\"obius band).

\begin{figure}[h]
\epsfig{figure=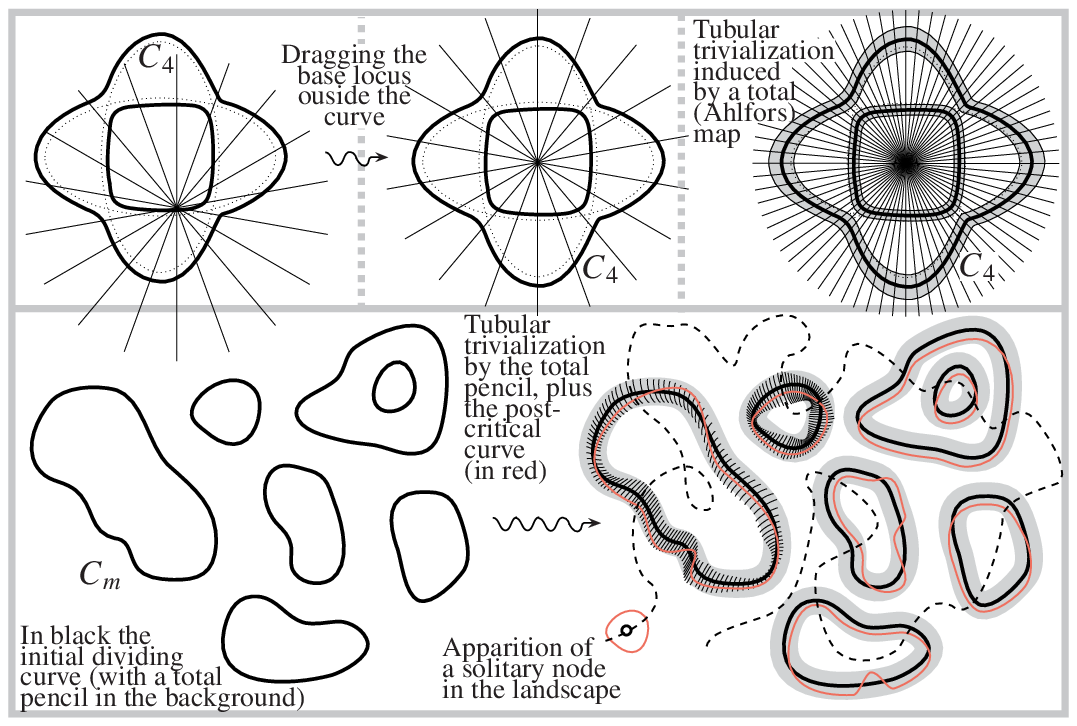,width=122mm} \vskip-5pt\penalty0
  \caption{\label{Tube:fig}%
Tubular neighborhood framed by an (Ahlfors) total pencil. The
post-critical curve (red-colored) is slaloming across the first
within the tube while staying transverse to the Ahlfors foliation.
The curve of the total pencil through the solitary node (dashed
line) has as many (nearby) intersections with the red curve as
with the black one, plus some
extra intersections with the newly created oval surrounding the
solitary node.} \vskip-5pt\penalty0
\end{figure}

As long as the discriminant is avoided, a slight continuous
variation of the coefficients engenders a small perturbation of
the real locus (continuity lemma for rigid-isotopies). When the
discriminant is crossed at a solitary node (ordinary double point
with 2 imaginary conjugate tangents),
 the real locus
acquires
 a (single) champagne
bubble, while the rest of the curve is isotoped within the
prescribed tube (compare Fig.\,\ref{Tube:fig} bottom). It can be
assumed (for simplicity, but not vital) that the deformed curves
$(C_t)_{t\in[-1,+1]}$
stay transverse to the foliation. Inside of the tube, the
intersection of the curve $\Gamma$ of the total pencil through the
solitary node, with the post-critical curve $C_{+1}$ is in natural
bijection with that of the pre-critical curve $C_{-1}$ . Yet the
former intersection $C_{-1} \cap \Gamma$ is totally real, whereas
the second contains two additional points when $\Gamma$ cuts the
new created oval of $C_{+1}$ (better
ask $\Gamma$ to pass through a point of this new oval). B\'ezout
is violated, and our (Viro inspired) proof of Klein's assertion is
complete.
\end{proof}

Klein's thesis (as discussed in Marin 1988 \cite{Marin_1988} or
Viro 1986) is the stronger assertion that a dividing curve cannot
see its number of circuits increase when crossing the
discriminant. This probably also follows from Ahlfors after some
suitable thinking, yet perhaps is less close to Klein's original
statement. It is only now that Klein's allusion (Quote
\ref{Klein_1876-niemals-isolierte:quote}) appears to me quite
transparent (yet via the powerful Ahlfors theorem).

{\it Insertion} [30.03.13].---It is unlikely that this was the
original proof of Klein (despite the fact that Teichm\"uller 1941
ascribes to Klein the theorem usually ascribed to Ahlfors).
Klein's original reasoning (alas unpublished in details, but only
claimed in Klein 1876 \cite{Klein_1876}) might have been rather
purely topological, essentially like Marin's (though the latter's
statement is somewhat stronger). More on this will be discussed
below, especially in (\ref{Klein-Marin:lem}).

[11.01.13] One may wonder if, conversely, a nondividing curve can
always acquire a solitary node and so a new oval. This is also
implicit in Klein 1876 intuition, and probably true up to degree
$6$ (cf. Gudkov's table=Fig.\,\ref{Gudkov-Table3:fig} for some
evidence, while a rigorous proof probably rests on Nikulin's
rigid-isotopy classification via Rohlin's enhanced Gudkov table by
complex characteristics).
(For a verification of Klein's intuition in degree $6$, see
Prop.\,\ref{Klein-vache-deg-6:prop}, and for a disproof in degree
8, cf. Shustin 1985 and our accompanying comments in
Sec.\,\ref{Shustin-understood:sec}.)

Gudkov's table shows however that the location for the apparition
of a bubble cannot be chosen in advance. Indeed starting say from
the scheme $10$ we could by bubbling create the scheme $11$
violating the Hilbert-Rohn-Petrovskii-Gudkov theorem, that such a
scheme is not realized algebraically. Another more obvious
argument is just to take any scheme (not on the ``visible faces''
of Gudkov's pyramid, equivalently, such that $\frac{1}{1}1$ is a
subscheme) and create a bubble inside the outer oval so that the
new real scheme contains the subscheme $\frac{1}{1}\frac{1}{1}$
consisting of 2 nests of depth 2, which violates B\'ezout.

\section{Prohibitions}
\label{Prohibitions:sec}

[28.03.13] From now on, we do not follow historical order, but
rather logical necessity. Admittedly there is no universal measure
of simplicity as it depends much on the background of the
investigator. From a radical viewpoint, the unique measure of
simpleness could be the natural historical time-arrow. Yet
sometimes big surprises arise. Arguments extremely powerful and
strikingly simple (nearly stemming from nowhere) tend to
trivialize much of the past efforts. Such an example is Rohlin's
formula discussed below, which bears some antecedents only by
Arnold, plus the topological heritage of Riemann, Betti,
Poincar\'e, Lefschetz, Weyl, Pontryagin, etc (homological
intersection theory).

The source of prohibitions in Hilbert's 16th problem, are
multiple. Albeit we are not expert in the field let us brush a
brief historical sketch.

First there are evident restrictions coming from B\'ezout. Those
were used by Zeuthen 1874 \cite{Zeuthen_1874}, and exploited in
full in Hilbert 1891 (boring bounds on the depth of ovals). A
major prohibition (directly affiliated to Zeuthen) is the Harnack
inequality $r\le g+1=\frac{(m-1)(m-2)}{2}$ of 1876, in no way
specific to plane curves. Klein 1876, then aged 27, (but already
Harnack's teacher) gave a more intrinsic justification boiling
down to a basic fact on the topology of surfaces directly
imputable to Riemann's definition of the genus (or rather its
allied connectivity). Recall that the jargon of the genus is due
to Clebsch.  Harnack's inequality is something very robust, as it
extends to all dimensions via Smith theory, as was noticed by Thom
and Milnor, yielding something like $b_\ast(\RR X)\le b_\ast (\CC
X)$, for $b_\ast$ the total Betti number. We shall not need this
as we confine attention to curves (where enough work remains to be
done).

As discussed above (\ref{Klein-unnested-quartic-nondividing:lem}),
Klein 1876 also used large deformations (rigid-isotopies), to
prove e.g. that a quartic with 2 unnested ovals is nondividing.
Later he also exploited theta-characteristics. The first method is
unlikely to extend to curves of higher orders (despite Nikulin's
rigid classification), while the second has been poorly explored
further since Klein 1892, and Gross-Harris 1981
\cite{Gross-Harris_1981}, and does not seem able to compete
seriously with information distilled by Rohlin's formula. Perhaps
those old Jacobi-Riemann-Klein methods deserve to be revived. As
to Nikulin, it seems at first that it will tell nothing being
rather built upon the Gudkov-Rohlin classification by types (cf.
Fig.\,\ref{Gudkov-Table3:fig}). However in the fingers of Itenberg
1994 \cite{Itenberg_1994} (contraction theorem of empty ovals), we
can expect (at least if this strengthens to our
CCC=(\ref{CCC:conj})) to rederive via strangulation the
diasymmetry(=type~II) of the schemes $1$ and $5$ of degree 6
(gaining so some analogy with Klein's rigid-isotopy argument for
the bifolium $2$ in degree 4). Yet, the difficulties are so great
that this looks quite artificial as compared to the topological
straightforwardness of Rohlin's formula. By the way this would
miss  the scheme $\frac{1}{1}1$. Further this seems much limited
to degree 6 as we lack precise information on rigid-isotopy in
high-degrees. Then the connection with K3-surfaces is lost, and so
the tool making Nikulin's theorem possible (deep transcendental
algebraic geometry, global Torelli theorem, etc.).

After Zeuthen-Harnack-Klein, came Hilbert's 1891
intuition\footnote{Unpublished, but see Rohn 1911
\cite{Rohn_1911}, 1913 \cite{Rohn_1913}, yet not judged complete
by Gudkov 1974.} that an $M$-sextic is forced to nest. This has no
antecedents (as far as I know), yet it could be challenging to
reprove it via conformal geometry (i.e. the
Riemann-Schottky-Bieberbach-Grunsky theorem). This has never been
implemented and is probably a hard game, if feasible at all.

After Hilbert came several things like Ragsdale 1906, and Rohn
1911 who consolidated Hilbert's method. This involves a deep
analysis of the stratification of the space of curves and the
usage of pencils. In Gudkov's fingers, this produced an exhaustive
list of prohibitions in degree 6. Perhaps an extension of this
method also implies (or rather converges) with
Rohlin-Le~Touz\'e's phenomenon of total reality. At least the
diagrammatic of the Gudkov table (Fig.\,\ref{Gudkov-Table3:fig})
strongly suggests this. In degree 6, all the information gained
via Hilbert-Rohn is recovered for topological reasons \`a la
Gudkov-Arnold-Rohlin (and extended to all degrees). It is not
clear to me if this subsuming of HR to GAR is specific to degree 6
or a general feature. Probably not if I remember well a seminal
talk by Orevkov (Geneva, ca. 2011) where Hilbert-Rohn was still
much on the appetizer. After all it is unlikely that deep
geometrical methods get completely phagocytozed by topological
ones.

Enriques-Chisini 1915 \cite{Enriques-Chisini_1915-1918} gave a
proof of Harnack's inequality based on Riemann-Roch and a
continuity argument (compare our
Lemma~\ref{Enriques-Chisini:lemma}). This is much akin to the
phenomenon of total reality, and need to be extended to less
trivial cases. Recall that from the viewpoint of total reality,
$M$-curves constitute the trivial case. This desideratum is the
main motivation of the present text yet we  still have very few
factual things to present.

The next great step is Petrovskii 1933/38, who seems to be the
first to find universal obstructions (valid in all degrees). This
is based on Euler-Jacobi-Kronecker's interpolation formula plus
some Morse theory.

Then there is Gudkov breakthrough (apparition of congruences mod 8
as opposed to mere estimates), and the theorists Arnold, Rohlin,
etc. validating them via 4D-topology or Atiyah-Singer. In this
move we have the trinity of congruences modulo 8 for $M$, $(M-1)$
and $(M-2)$-curves due to GR, GKK, RKM, respectively. Here
G=Gudkov, R=Rohlin, 1st K=Krakhnov, 2nd and 3rd K=Kharlamov, while
M=Marin. The importance of those can hardly be underestimated.
First, the conjunction of GR and GKK explains all prohibitions in
degree 6 on real schemes (i.e. Hilbert's 16th), while the 3rd GKK
(forcing orthosymmetry(=type~I) of schemes with $\chi\equiv_8
k^2+4$) seems even to imply (via the hypothetical Rohlin
maximality conjecture=RMC) the conjunction of GR+GKK. Even without
the elusive RMC, it can be that explicit instances of total
reality (e.g., Rohlin-Le~Touz\'e's) imply in low-degrees (say
$m=6,8$) the truth of RMC in special situations. This looks after
all plausible, since totality involves a geometrization of the
type~I topological condition by a stronger geometric property
(total pencil). Here and in the sequel, we shall often abridge
``total reality'' by ``totality''.

Then appears Rohlin's formula 1974--78. This is very strong and
completely elementary. In degree 6, it rules out all schemes above
the broken-line of Gudkov's table (Fig.\,\ref{Gudkov-Table3:fig})
safe 6, namely the 2 triangles involving the symbols
$\frac{7}{1}3$, $\frac{7}{1}2$, $\frac{6}{1}3$ and its mirror
$\frac{3}{1}7$, $\frac{2}{1}7$, $\frac{3}{1}6$.
Rohlin's formula is very powerful, yet somewhat too elementary to
grasp the full mystery. It need therefore to be complemented by
more advanced weapons like the Gudkov congruence (GR), and GKK, or
by Rohlin's maximality principle allied to total reality.

In 1978, we have Rohlin's maximality principle (RMC), still
conjectural and not yet fully exploited in our opinion. This could
loop-back to conformal geometry \`a la Riemann, Schwarz, Schottky,
Klein, Koebe, Bieberbach, Grunsky, Teichm\"uller, Ahlfors. As said
above, if RMC looks impossible to implement in universal
generality it could be verifiable in special cases by using
totality as a geometric strengthening of the (topological)
type~I-condition. For instance  Rohlin-Le~Touz\'e's totality
should suffice (either with or without RKM) to kill all expansions
of the 2 orthosymmetric $(M-2)$-schemes of degree 6. This would
unify all prohibitions in degree 6 safe the schemes $11$ and
$\frac{10}{1}$ (easily ruled out via Rohlin's formula).

Ca. 1978--80, we have advanced B\'ezout-style obstructions \`a la
Fiedler-Viro (\ref{Viro-Fiedler-prohibition:thm}) that really
pertains to curves of degree $8$. Those plays a pivotal r\^ole in
Shustin's disproof of Klein's champagne bubbling principle for
nondividing curves, as well as the disproof of the reverse
implication of Rohlin's maximality principle. More generally those
look indispensable in the higher cases $m=7,8$ of Hilbert's 16th.

We have also the locking trick of Marin-Fiedler (also founded on
B\'ezout for lines) that provides obstruction to rigid-isotopy on
$M$-curves of degree $\ge 7$. Here the idea is that if we have a
triangle (3 lines) which is B\'ezout-saturated and canonically
attached to a scheme (typically a disc with 3 holes), then during
a rigid-isotopy ovals cannot traverse this moving frame. Hence the
distribution of ovals past such a fundamental triangle is an
invariant of the rigid-isotopy class. Of course this method is not
a method of prohibition of schemes, but prohibits the existence of
pathes in the hyperspace of smooth curves.

Finally, we have probably a r\^ole of Thom's conjecture on
genus-bound (verified since Kronheimer-Mrowka 1994
\cite{Kronheimer-Mrowka_1994}), yet whose role is not so clear-cut
as initially expected. The simple case of Thom, due to
Kervaire-Milnor 1961 \cite{Kervaire-Milnor_1961}, may be used to
settle Hilbert's nesting ``theorem'' for $M$-sextics. In general
the role of Thom, is perhaps marginalized by Rohlin's formula and
other strong results, yet seems to give new information in the
work of Mikhalkin 1994 \cite{Mikhalkin_1994-adjunction-Thom} when
it comes to split curves (communication of Fiedler, not yet
digested by the writer=Gabard).

This is a brief overview of nearly all what exists. In contrast
one can ask for more conciseness when it comes to explain all the
prohibitions  of Hilbert's problem (in degree 6) to a classroom.
As often repeated, nearly everything could reduce to the
(Klein-Ahlfors-)Rohlin-Le~Touz\'e's phenomenon of total reality.
Remind that 2 technical points are still obscure, but
philosophically trivial. The first is a complete proof of Rohlin's
claim (preferably without employing the RKM congruence mod 8). The
second is to verify that Rohlin-Le~Touz\'e's total reality is
strong enough to imply maximality of the two Rohlin's
$(M-2)$-schemes. Assuming this settled, we still miss the
prohibition of Hilbert's unnested scheme $11$, and  Rohn's
maximally nested scheme $\frac{10}{1}$. This is paradoxical
inasmuch as those 2 guys were historically the first ruled out by
the Hilbert-Rohn method. The 1st scheme $10$ can be killed by the
Kervaire-Milnor 1961 \cite{Kervaire-Milnor_1961} elementary case
of Thom's conjecture in degree $k=3$, but the second
$\frac{10}{1}$ fails to succumb under Thom. However both of them
are killed by Rohlin's formula. Hence a good cocktail (for the
classroom or the economical reader) is to mix total reality with
Rohlin's formula. This reduces all prohibitions in degree 6 to
only 2 paradigms. As far as we know, apart form the Hilbert-Rohn
method (as developed by D.\,A. Gudkov) there is no universal force
unifying all prohibitions in a single one (even in degree 6). A
substitute to Thom-Kervaire-Milnor is to use Petrovskii 1933/38
(\ref{Petrovskii's-inequalities:thm}). This rules out $11$ but not
Rohn's scheme $\frac{10}{1}$. The latter is not even killed by the
strong Petrovskii inequality of Arnold
(\ref{Strong-Petrovskii-Arnold-ineq:thm}), i.e., $n-p^-\le
\frac{3}{2}k(k-1)=\frac{3}{2}3\cdot 2=9$, where $p^{-}=1$ is the
number of hyperbolic positive ovals, so $n\le 10$ while Rohn's
scheme has $n=10$.

\subsection{Obstructions via Rohlin's formula
(Rohlin 1974, 1978)}
\label{Rohlin-formula:sec}

[03.01.13] We
 repeat the proof of the following pivotal
result (whose proof puzzled me a lot as I was young, and still
imbues some suitable respect\footnote{Prose borrowed by Jack
Milnor, when he speaks about non-metric manifold, cf. his preprint
on foliated bundles.} when getting older). Crudely put, Rohlin's
formula is nothing less than the most universal obstruction that
one may derive by abstract non-sense (i.e. using virtually nothing
from the algebraicity assumption).

\begin{theorem} {\rm (Rohlin 1974--78)} \label{Rohlin-formula:thm}
For any (real, smooth, algebraic, plane) dividing curve of even
order $m=2k$ (odd orders were treated by Rohlin's student
Mishachev), the following equation holds:
\begin{equation}
2(\Pi^+ -\Pi^-)=r-k^2, \label{Rohlin-formula:eq}
\end{equation}
where $r$ is the number of ovals, while $\Pi^{\pm}$ are the number
of positive (resp. negative) pairs of nested ovals. Each pair of
nested ovals  bounds a ring=annulus in ${\Bbb R}P^2$, and upon
comparing with the complex orientation (as the border of the
semi-Riemann surface) one defines a positive pair when both
orientations (complex vs. real) agree, and a negative pair when
they disagree (cf. Fig.\,\ref{Rohlinsformula:fig}).
\end{theorem}

\begin{figure}[h]
\vskip-0.3cm\penalty0 \centering
    \epsfig{figure=Rohlinsformula.eps,width=92mm}
\vskip-5pt\penalty0
\caption{\label{Rohlinsformula:fig}%
Rohlin's positive and negative pairs} \vskip-5pt\penalty0
\end{figure}

\begin{proof}
The idea involves computing the self-intersection of the half of
the dividing curve after
capping off  by discs bounding the ovals in the real projective
plane, or rather the intersection with the conjugate capped off
membrane. This argument seems
inspired by a similar device used by Arnold in 1971
\cite{Arnold_1971/72}, but now slightly more punch is acquired.
The proof is very elementary using merely intersection of homology
classes (available since the days of Poincar\'e, Lefschetz, Hopf,
Pontryagin, etc.) and Poincar\'e's index formula (available since
Gauss?, Kronecker 1868, Poincar\'e 1885), plus some basic trick
about the ``Lagrangian property'' of real parts of algebraic
varieties (jargon used in Degtyarev-Kharlamov 2000
\cite{Degtyarev-Kharlamov_2000}). One can then nearly wonder why
such a formula escaped Felix Klein's attention, but this is of
course just historical slowness of the revelation of
brutal combinatorial truths.

[06.03.13] The detailed proof is given in Rohlin 1974/75
\cite[p.\,332]{Rohlin_1974/75}, but  stated there only for
$M$-curves. The adaptation to the general case requires only minor
notational changes, even  simplifying a bit Rohlin's original. For
convenience let us thus copy Rohlin's prose (while adapting it to
the broader context, brackets are our additions):

Denote by $B_C\subset \RR P^2$ the bounding disc for the oval $C$.
Complete the half $C^+$ of the dividing curve to a closed surface
$\Sigma$ by adding nonintersecting copies of the disc $B_C$. Let
$T$ be the closed surface obtained from the other half $C^{-}$ by
the same procedure, and let $\varphi \colon \Sigma \to \CC P^2$
and $\psi \colon T \to \CC P^2$ be mappings fixed\footnote{I.e.
identity.} on $C^+$ and $C^-$ and superimposing copies of $B_C$
onto these discs. Further, let $\xi$, $\eta$ be elements of the
(integral) homology group $H_2(\CC P^2)$ determined by the
mappings $\varphi$ and $\psi$ and the natural orientations on the
pretzels $\Sigma$ and $T$ (i.e. the orientation obtained from
$C^+$ and $C^-$). We shall establish
Eq.\,\eqref{Rohlin-formula:eq} by computing the intersection index
$\xi \eta$ by two procedures [geometrically and algebraically].

1.---The first procedure is based on the fact that $\xi \eta$ can
be interpreted as the algebraic number of points in the
intersection of the oriented singular pretzels $\varphi\colon
\Sigma \to \CC P^2 $ and $\psi \colon T \to \CC P^2 $. This number
cannot be determined directly, since the intersection consists of
wholes disks, and we begin by applying a deformation to $\varphi$,
making the intersection more regular[=nearly transverse]. Let $u$
be some tangent vector field on $\RR P^2$ with a finite number of
zeros, not having zeros on $A:=C_m(\RR)$ and normal to $A$ on $A$.
Since the field $iu$ is normal to $\RR P^2$ in $\CC P^2$ and
normal to $\CC A:=C_m(\CC)$ on $A$, it can be normally extended to
some field $v$ on $\RR P^2 \cup C^+$ (the latter, of course, will
have zeros inside of $C^+$); let $\gamma\colon \RR P^2 \cup C^+
\to \CC P^2$ be a geodesic translation defined by the field
$\delta v$, where $\delta$ is a sufficiently small positive
number, and $\varphi'\colon \Sigma\to \CC P^2$ be the mapping
defined by the formula $\varphi'(x)=\gamma(\varphi(x))$. For
$\varphi'$ the algebraic number of points of intersection with
$\psi$ is determined directly and can be found in the following
way. Since the sum index of the singularities of $u$ in each of
the disks $B_C$ is equal to 1 and multiplication by $i$
anti-isomorphically maps the tangent bundle of $\RR P^2$ onto its
normal bundle in $\CC P^2$, the sum index of $v$ on each of the
disks $B_C$ is equal to $-1$. Consequently, the contribution added
by the pair of disks $B_C$ and $B_C'$ to the algebraic number
[$\xi \eta$] of intersection points that we are interested in is
equal to:

$\bullet$ $+1$ if $C=C'$;

$\bullet$ $+2$ if the pair $C,C'$ is negative, and equal to

$\bullet$ $-2$ if the pair $C,C'$ is positive.

This number itself is thus equal to $r-2(\Pi^+-\Pi^-)$. Since
$\varphi'$ is homotopic to $\varphi $, the index $\xi \eta$ is
also like that and thus
$$
\xi \eta=r-2(\Pi^+-\Pi^-).
$$

2.---The second procedure reduces to two remarks. First the class
of $\xi + \eta$ is realized by the surface $\CC A$ and therefore
coincides with $2k \alpha$, where $\alpha$ is the natural
generator of the group $H_2(\CC P^2)$. Second, since the
homomorphism $conj_\ast\colon H_2(\CC P^2)\to H_2(\CC P^2)$
represents multiplication by $-1$ [as it flips the orientation of
the generator interpreted as the fundamental class of a line
defined over $\RR$] and takes $\xi$ to $-\eta$, we have $\xi=
\eta$. From these remarks it follows that $\xi=k\alpha$,
$\eta=k\alpha$ and $\xi \eta=k^2$.

Comparing the last equations obtained along each procedure, we
obtain the announced formula \eqref{Rohlin-formula:eq}.
\end{proof}

We list some consequences. First a (promised) remark about
quartics:

\begin{cor} {\rm (Klein 1876, Rohlin 1978)}
Any quartic with $2$ unnested ovals is nondividing.
\end{cor}

\begin{proof}
Since there is no nesting there in no pairs of ovals and the
left-side of Rohlin's formula \eqref{Rohlin-formula:eq} vanishes,
while the right-side is equal to $r-k^2=2-2^2=2-4=-2$.
\end{proof}

\begin{cor}
The sextic scheme $5$ (five unnested ovals) is of type~II. More
generally the sextic scheme $r$ (\/$0\le r\le 11$ excepted $r=9$)
is of type~II (actually $11$ is not realized by Hilbert, Kahn
1909, L\"obenstein 1910, Rohn 1911--13, Petrovskii 1938, Gudkov,
but a more limpid proof follows from Rohlin's formula).
\end{cor}

\begin{proof} (due to Rohlin 1978 \cite{Rohlin_1978}, also in Fiedler
1981 \cite[p.\,13]{Fiedler_1981}). Since there is no nesting
$\Pi^{\pm}$ are both zero, while the left-side $r-k^2=r-3^2$ of
Rohlin's formula vanishes only for $r=9$.
\end{proof}

\begin{cor} \label{Rohlin's-inequality:cor} {\rm (Rohlin's inequality)}
A dividing plane curve of (even) order $m$ has at least $r\ge m/2$
ovals. Further if equality $r=m/2$ holds (and the curve is
dividing) then its real scheme must be a deep nest (i.e. $m/2$
ovals each pair of them being nested).
\end{cor}

\begin{proof} (explicit in Marin 1979 \cite{Marin_1979}, or
Gabard 2000 \cite[p.\,148]{Gabard_2000}, but due to Rohlin). Let
$\Pi=\Pi^+ +\Pi^-$ be the total number of nested pairs of ovals.
We have
$$
\Pi\le \textstyle\binom{r}{2},
$$
(binomial coefficient counting the number of pair of a finite set
of size $r$). Equality occurs only for a deep nest! Rohlin's
formula gives:
$$
r=k^2+2(\Pi^+-\Pi^-)\ge k^2-2\Pi^-\ge k^2-2\Pi\ge k^2-2
\textstyle\binom{r}{2}=k^2-r(r-1),
$$
whence (looking at the extremities) $r^2\ge k^2$, i.e. $r\ge k$.
If an equality each intermediate estimates crunch to equality, in
particular the estimate $\Pi\le \binom{r}{2}$, which is fulfilled
only for a deep nest.
\end{proof}

The sequel  studies Rohlin's consequence in degree 6. This is a
bit pedestrian, and can be omitted as we gave a somewhat more
conceptual explanation before, by noticing that Rohlin implies
Arnold, etc.

Assume again no-nesting  ($\Pi=0$). Then Rohlin's formula gives
$0=2(\Pi^+-\Pi^-)=r-k^2=r-9$, it follows $r=9+ 0= 9$ (in
accordance with Fig.\,\ref{Gudkov-Table3:fig}). It is quite
remarkable to notice that this gives an instant proof of Hilbert's
conjecture (and semi-theorem of his students Kahn-L\"obenstein and
Rohn, etc.) to the effect that there is no $M$-curve with 11
unnested ovals.

Next we assume that there is one pair of nested ovals ($\Pi=1$).
Then Rohlin's formula gives $\pm 2=2(\Pi^+-\Pi^-)=r-k^2=r-9$, it
follows $r=9\pm 2= 11, 7$ (in accordance with
Fig.\,\ref{Gudkov-Table3:fig}).

Next suppose $2$ nested pairs ($\Pi=2$). Hence $\{ 4, 0,-4 \} \ni
2(\Pi^+ -\Pi^-)=r-k^2=r-9$, it follows $r=9+ \{ 4,0,-4 \}= 13,
9,5$ (in accordance with Fig.\,\ref{Gudkov-Table3:fig}).

For 3 nested pairs, $\{ 6, 2,-2, -6 \} \ni 2(\Pi^+
-\Pi^-)=r-k^2=r-9$, it follows $r=9+ \{ 6, 2,-2, -6 \}= 15, 11,
7,3$ (in accordance with Fig.\,\ref{Gudkov-Table3:fig}).

For 4 nested pairs, $\{ 8, 4, 0, -4,-8 \} \ni 2(\Pi^+
-\Pi^-)=r-k^2=r-9$, it follows $r=9+ \{ 8, 4, 0, -4,-8 \}= 17, 13,
9,5,1$ (in accordance with Fig.\,\ref{Gudkov-Table3:fig}).

For 5 nested pairs, $\{ 10, 6, 2, -2,-6, -10 \} \ni 2(\Pi^+
-\Pi^-)=r-k^2=r-9$, it follows $r=9+ \{ 10, 6, 2, -2,-6, -10 \}=
19, 15, 11, 7, 3, -1$ (in accordance with
Fig.\,\ref{Gudkov-Table3:fig}).

Etc., at this stage it is clear how to link the arithmetics of
Rohlin's formula to the geometry of Gudkov's table enhanced by
Rohlin's data, and we have proven Rohlin's claim:

\begin{cor}
All green-squared schemes on Gudkov's
table=Fig.\,\ref{Gudkov-Table3:fig} are of type~II.
\end{cor}

One noteworthy feature of the diagrammatic is that Rohlin's
formula gives a tiling by squares rooted on our blue-rhombs plus
the red circles of Fig.\,\ref{Gudkov-Table3:fig}. Hence all the
schemes not situated on this grid are necessarily of type~II. In
particular since $M$-schemes are forced to type~I it follows from
Rohlin that all $M$-schemes not on the square grid are prohibited
as real schemes (yielding a significant contribution to Hilbert's
16th problem). Explicitly we have:

\begin{cor}
All the $M$-schemes outside the grid are not realized
(algebraically), that is $\frac{10}{1}$, $\frac{8}{1}2$,
$\frac{6}{1}4$, $\frac{4}{1}6$, $\frac{2}{1}8$, $11$.
\end{cor}

However Rohlin's formula alone fails to prohibit the schemes
$\frac{7}{1}3$ and $\frac{3}{1}7$ (which are situated on the
grid). Those are however prohibited either by the
Hilbert-Rohn-Gudkov method, or by the Gudkov hypothesis proved by
Rohlin 1972/72 (as detailed in the next
Sec.\,\ref{Gudkov-hypothesis:sec}).

This last corollary helps the beginner to catch the substance of
the following remark by Degtyarev-Kharlamov 2000
\cite[p.\,736]{Degtyarev-Kharlamov_2000}: {\it ``Another
fundamental result difficult to overestimate is Rokhlin's formula
for complex orientations. The notion of complex orientation of a
dividing real curve (see below), as well as Rokhlin's formula and
its proof, seem incredibly transparent at first sight. The formula
settles, for example, two of Hilbert's conjectures on 11 ovals of
plane sextics, which Hilbert himself tried to prove in a very
sophisticated way and then included in his famous problem list (as
the sixteenth problem).''}

To remind it seems that Hilbert conjectured (wrongly) that only
the schemes $\frac{9}{1}1$ and $\frac{1}{1}9$ do exist among
$M$-schemes. This turned out to be wrong when Gudkov exhibited the
scheme $\frac{5}{1}5$. So Rohlin's formula settles actually six
(!) of Hilbert's conjectures (if taken as individual prohibition).
Presumably what Degtyarev-Kharlamov had in mind were the extreme
schemes $\frac{10}{1}$ and $11$ (eleven unnested ovals). The
philosophical outcome of this spectacular Rohlin formula is how
much information can be derived by basic topological methods,
basically emanating from the Riemann-Betti-Poincar\'e tradition,
yet to which workers like Klein or Hilbert were not enough
familiar with. Of course a first class topologist like Rohlin was
needed to reveal this truth.

\subsection{Gudkov hypothesis (Gudkov 1969, Arnold 1971,
Rohlin 1972, etc.)}\label{Gudkov-hypothesis:sec}

[07.01.13] For $M$-curves, the congruence $\chi\equiv_8 k^2$ was
conjectured by Gudkov on the basis of experimental data gathered
along his Hilbert-Rohn approach for sextics,
and of course by looking as well to higher degrees via the
Harnack-Hilbert construction. Figs.\,\ref{HilbGab1:fig},
\ref{HilbGab2:fig}, \ref{HilbGab4:fig} below illustrate  with
which metronomic precision the Hilbert construction always produce
$M$-curves respecting the congruence $\chi\equiv_8 k^2$. As
pointed somewhere in Viro's writings, nothing could thus have
impeded Miss Ragsdale  to detect this congruence in 1906 already.
Yet it is the full-credit of Gudkov to have spotted this
regularity. Once Arnold knew this, it was just a matter of
hard-work toward elaborating the right strategy of proof, and some
extra-skills of Rohlin turned to be indispensable.

So the full proof belongs  to Rohlin 1972/72
\cite{Rohlin_1972/72-Proof-of-a-conj-of-Gudkov} (alas contain a
little bug),  boosting ideas initiated by Arnold 1971
\cite{Arnold_1971/72} (who got the weaker  congruence mod 4).
Rohlin's proof extract his punch not just from algebra
(divisibility by 8 of an even integral unimodular quadratic form)
but from the deeper divisibility by 16 coming from his own old
``grand cru''
 of 1952 (Rohlin 1952 \cite{Rohlin_1952-4-manifolds}) on the
signature of spin smooth $4$-manifolds. It is notorious that
Rohlin's proof (1972/72 \loccit) contains a mistake that was
repaired by Guillou-Marin 1977 \cite{Guillou-Marin_1977} (compare
e.g. Degtyarev-Kharlamov 2000
\cite[p.\,736]{Degtyarev-Kharlamov_2000} and also Wilson 1978
\cite{Wilson_1978}, who seems to have noticed (the same?) gap).

\begin{Notation}\label{Ragsdale-Petrovskii:notatio} {\rm (Ragsdale 1906 \cite{Ragsdale_1906},
Petrovskii 1938 \cite{Petrowsky_1938}) Given a plane curve of even
order (or more generally a real scheme of ovals), it is customary
to denote by $p$  the number of {\it even ovals\/} (those included
in an even number of ovals) and by $n$ the number of {\it odd
ovals} (defined analogously). The difference $p-n$ can always be
interpreted as the Euler characteristic $\chi$ of the orientable
membrane of $\RR P^2$ bounding the curve. The notation $p,n$ are
Petrovskii's, intended to stand for positive and negative ovals.}
\end{Notation}

\begin{theorem}\label{Gudkov-hypothesis:thm} {\rm
(Gudkov hypothesis/conjecture =Rohlin's theorem of 1972, modulo a
correction by Guillou-Marin 1977), and another proof in Rohlin
1974} A plane $M$-curve of degree $m=2k$ satisfies the
Gudkov-Rohlin congruence:
$$
\chi=p-n\equiv k^2 \pmod 8.
$$
\end{theorem}

\begin{proof}
The technique is akin to the subsequent Rohlin's complex
orientation formula of 1974--78, namely fill the halves of the
orthosymmetric curve to a closed membrane and calculate the
resulting intersection. However here the proof use (an extension
of) the seminal Rohlin theorem 1952 \cite{Rohlin_1952-4-manifolds}
on the divisibility by 16 of the signature of spin 4-manifolds.
(At this stage there is a huge constellation of coincidence around
Hilbert's heritage: the 16th problem and as well his student, H.
Weyl, whose ``Analisis situs combinatorio'' of 1922  is the first
place where the signature of $4k$-manifolds is defined). So  a
breathtaking connection between differential topology and the more
rigid algebraic geometry is accomplished in the Arnold-Rohlin era.

Logically this proof is a bit tricky to implement for one is
warned of some mistakes in Rohlin's initial paper. Hence pivotal
is the Guillou-Marin extension of Rohlin's signature formula, for
a full exposition cf. Guillou-Marin 1986
\cite{Guillou-Marin_1986}. Once this is understood its application
to Gudkov's hypothesis is exposed in A'Campo 1979
\cite{A'Campo_1979} (following a presentation due to Marin).
\end{proof}

\begin{cor}
Among all sextic $M$-schemes only those of Harnack, Hilbert and
Gudkov exist.
\end{cor}

\begin{figure}[h]
\centering
    \epsfig{figure=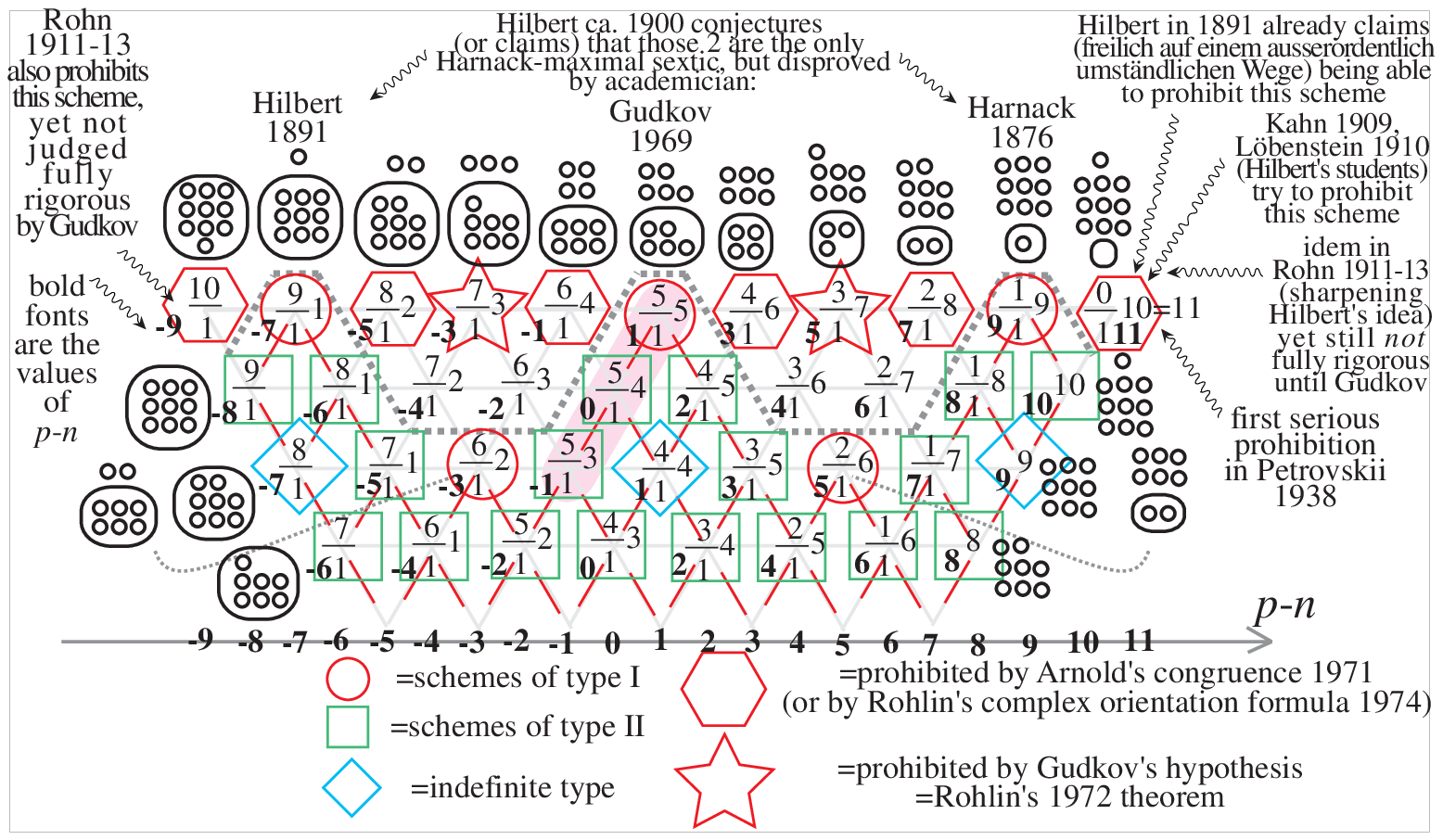,width=122mm}
\vskip-5pt\penalty0
\caption{\label{Gudkov-TableTop:fig}%
The top of Gudkov's table: reporting obstructions of his
followers. Bold-faced number shows the value of $\chi=p-n$.
Arnold's congruence implies that all schemes lying on the  2-by-2
square grid extending downwards the red hexagons are necessarily
of type~II.} \vskip-5pt\penalty0
\end{figure}

\begin{proof} Originally the proof was achieved by Gudkov 1954
\cite{Gudkov_1954} via the Hilbert-Rohn(-Gudkov) method, i.e.
supplemented by the concept of roughness coming from the
Andronov-Pontrjagin theory of dynamical (structural stability).
However here we derive it rather from the above theorem (Rohlin
1972). For schemes of degree 6, written in Gudkov's notation
$\frac{k}1 \ell$, we obviously have $p=\ell+1$ and $n=k$. Some
boring computation is required to check that this prohibit all
$M$-schemes above the broken line in Gudkov's table. Indeed:

$\bullet$ for\vadjust{\vskip2pt} $\frac{10}{1}$,
$p-n=1-10=-9=-1\neq +1=k^2 \pmod 8$.

$\bullet$ for\vadjust{\vskip2pt} $\frac{9}{1} 1$, $p-n=2-9=-7=
+1=k^2 \pmod 8$ (no obstruction),

$\bullet$ for\vadjust{\vskip2pt} $\frac{8}{1} 2$, $p-n=3-8=-5\neq
+1=k^2 \pmod 8$,

$\bullet$ for\vadjust{\vskip2pt} $\frac{7}{1} 3$, $p-n=4-7=-3\neq
+1=k^2 \pmod 8$, etc. (progression by 2 units),

\noindent so the rest is better done mentally on looking at
Gudkov's table (Fig.\,\ref{Gudkov-Table3:fig}) of which we
reproduce the top portion below (Fig.\,\ref{Gudkov-TableTop:fig}).
\end{proof}

Interestingly, prohibiting sextic $M$-schemes is much easier (no
deep differential topology \`a la Rohlin) for the schemes not
situated at the two centers of the semi-hexagon of Gudkov's table,
i.e. $\frac{k}{1} \ell$ with $k$ even, whereas the ``hexagonal''
schemes $\frac{7}{1}3$ and its mirror $\frac{3}{1}7$ are much
harder to disprove (at least in the modern Arnold-Rohlin theory).
For sextics, one may wonder what is more elementary:
Hilbert-Rohn-Gudkov or Rohlin 1952--1972.

Remember that Arnold 1971 \cite{Arnold_1971/72} proved  the weaker
congruence modulo 4 of Gudkov's hypothesis:

\begin{theorem} {\rm (Arnold 1971, Wilson 1978)}
For $M$-curves of degree $2k$ (or more generally dividing curves,
cf. {\rm Wilson 1978 \cite[p.\,67--69]{Wilson_1978}}), we have
$$
\chi=p-n\equiv k^2 \pmod 4.
$$
\end{theorem}

This theorem of Arnold is more elementary than Gudkov-Rohlin,
while prohibiting exactly the same $M$-schemes as those excluded
by Rohlin's
formula (i.e. fails to exclude
the $\frac{7}{1}3$ and its mirror $\frac{3}{1}7$). (Of course this
is not so surprising as Rohlin owed some inspiration from Arnold.)
Note also that Arnold's congruence forces all schemes on the
square-grid (extending the red-hexagons where $p-n\equiv -1 \pmod
4$) to be of type~II as do Rohlin's formula. The latter is however
a bit stronger for ascribing type~II to the schemes $5$,
$\frac{1}{1}1$, and $1$.
[06.03.13] In fact, as suggested in Degtyarev-Kharlamov 2000
\cite[p.\,737]{Degtyarev-Kharlamov_2000}:

\begin{lemma} \label{Rohlin-implies-Arnold:lem}  Rohlin's
formula {\rm (\ref{Rohlin-formula:thm})}  implies
straightforwardly the (extended) Arnold congruence $\chi\equiv k^2
\pmod 4$ (for dividing curves).
\end{lemma}

\begin{figure}[h]
\centering
    \epsfig{figure=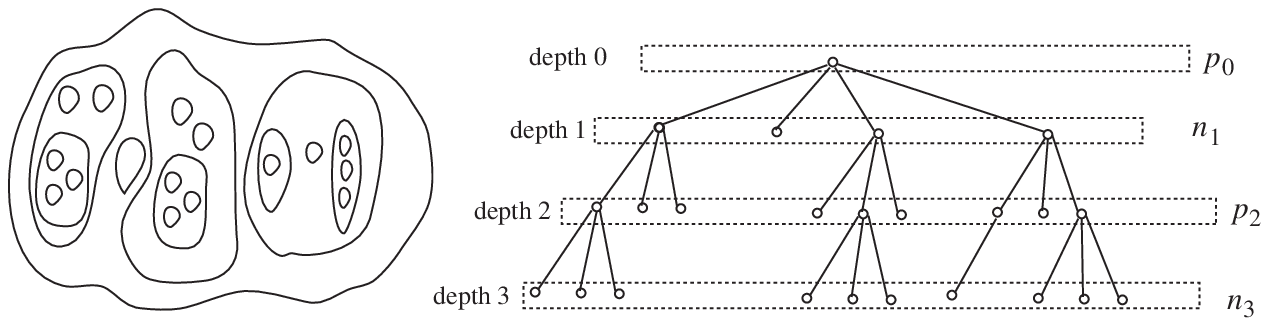,width=102mm}
\vskip-5pt\penalty0
\caption{\label{Stalin3:fig}%
The ``Hilbert tree'' of a plane curve encoding the distribution of
ovals} \vskip-5pt\penalty0
\end{figure}

\begin{proof}
This involves some abstract nonsense (yet pleasant) combinatorics.
Using the usual notation of Petrovskii (cf.
\ref{Ragsdale-Petrovskii:notatio}), we have $\chi=p-n$ (Euler
characteristic of the Ragsdale membrane), $r=p+n$ (total number of
ovals split into $p$ even ones and $n$ odd ones), and Rohlin's
formula $2(\Pi^+-\Pi^-)=r-k^2$. Assembling this gives
\begin{align}\label{Rohlin-to-Arnold:eq}
\chi=p-n=(p+n)-2n&=r-2n\cr
&=[2(\Pi^+-\Pi^-)+k^2]-2n\cr
&=k^2+2(\Pi^+-\Pi^- -n).
\end{align}
It remains to check that the ``corrector term'' $(\Pi^+-\Pi^- -n)$
is even. Modulo 2 we have, $ \Pi^+-\Pi^- \equiv_2 \Pi^+ +\Pi^- =
\Pi$. Hence we can ultimately ignore Rohlin's complex
orientations. The following lemma concludes the proof, via the
usual construction (like on Fig.\,\ref{Stalin3:fig}) assigning to
a plane curve its nested hierarchy of ovals ordered by inclusion
of their insides (i.e. the unique bounding disc of the oval
afforded by ``Schoenflies theorem'' applied in $\RR P^2$). Recall
that a Jordan curve on any surface is null-homotopic iff it bounds
a disc. (Cf. e.g.  Reinhold Baer's proof ca. 1927, Thesis under H.
Kneser, reproduced in Gabard-Gauld 2010
\cite{Gabard-Gauld_2010-Jordan-and-Schoenflies}.)
\end{proof}

\begin{lemma}\label{Stalin:lemma}
Given a finite tree with a directed structure upward so that the
tree really looks like the roots of a tree (or better a mushroom
in Arnold's metaphor). Formally we have a finite POSET where each
element admits at most one superior, i.e. an element larger than
it and minimal with this property (like in capitalistic or
feodal hierarchies). Then  the number $\Pi$ of pairs $x<y$ and the
number $n$ of vertices lying at odd depths are congruent modulo
$2$:
$$
\Pi\equiv n \pmod 2.
$$
\end{lemma}

\begin{proof} By additivity we may assume the tree connected. Then
there is a unique maximal element  in the hierarchy (Stalin), and
we can draw from him all his subordinated elements as a ``tree''
growing downwards (cf. Fig.\,\ref{Stalin:fig}) with several
elements lying at different depths(=combinatorial distance to
Stalin).

\begin{figure}[h]
\centering
    \epsfig{figure=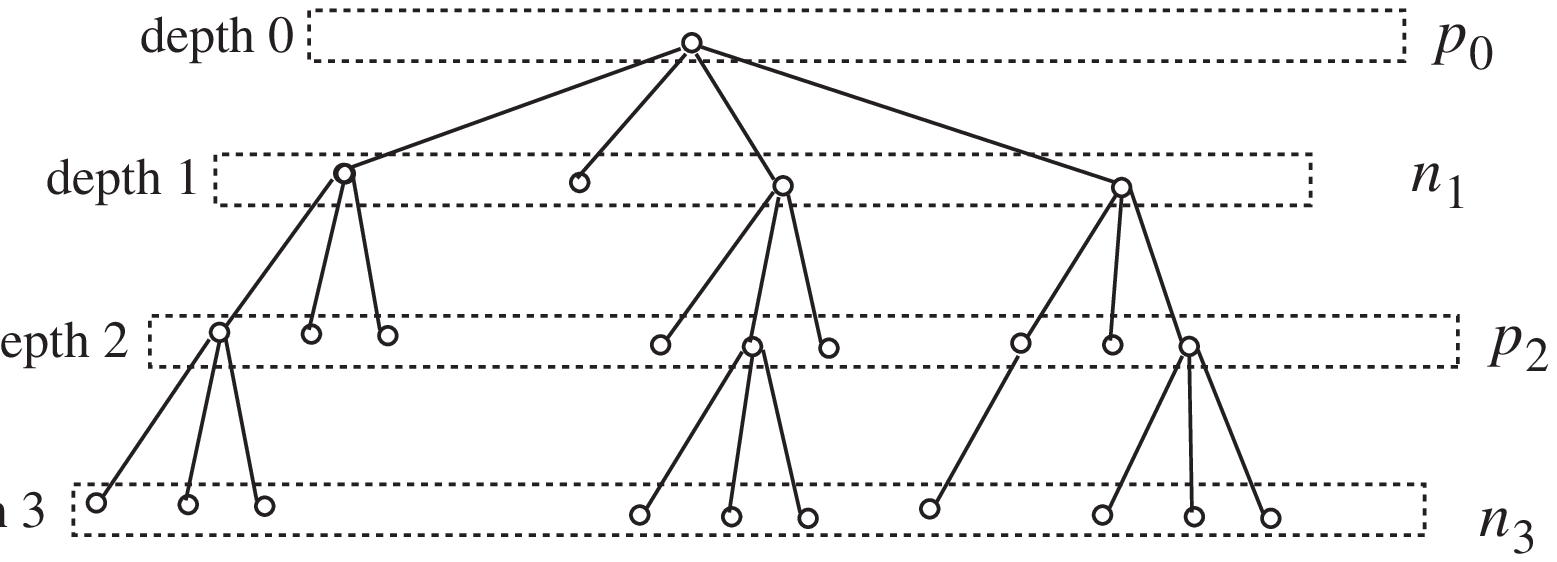,width=72mm}
\vskip-5pt\penalty0
\caption{\label{Stalin:fig}%
Count of injective pairs subsuming  Arnold 1971 to Rohlin 1974-78}
\vskip-5pt\penalty0
\end{figure}

Let $p_0, p_2, p_4, \dots$ be the number of elements at even
depths $0, 2, 4, \dots$, and $n_1, n_3, n_5, \dots$ be those at
odd depths $1,3,5, \dots$ respectively. To count $\Pi$ the number
of subordinations of the hierarchy, we range them by order of
importance (proximity to Stalin).  Since an element has as many
superiors as its depth, this gives
$$
\Pi=n_1+2p_2+3n_3+4p_4+5n_5+\dots\equiv_2 n_1+n_3+n_5+\dots=n.
$$
This enumeration clearly exhausts all possible hierarchical pairs,
and the proof is complete.
\end{proof}

\subsection{Gudkov-Rohlin congruence via Rohlin's formula?}
\label{Gudkov-hyp-via-Rohlin's-formula:sec}

[11.03.13] In fact the programme of this section looks extremely
dubious, just by virtue of the diagrammatic of the Gudkov table in
degree 6 (=Fig.\,\ref{Gudkov-Table3:fig}). Indeed for a scheme
like $\frac{7}{1}3$ the Rohlin equation is trivially soluble (cf.
e.g. Theorem~\ref{no-chance-to-reduce-Gudkov-to-Rohlin:thm}).
Therefore there is little chance to reduce Gudkov hypothesis to
Rohlin's formula and the sole signs-law on the Rohlin tree, unless
one is able to infer sharper information on complex orientations
from geometrical considerations, maybe via total pencils that are
fairly easy to construct (cf.
Theorem~\ref{total-reality-of-plane-M-curves:thm}) but it is
probably another matter to visualize their dynamics. Hence we
recommend to skip reading this section.

[08.03.13] In view of the previous reduction of Arnold's mod 4
congruence to Rohlin's formula, an evident idea  is to get better
control on the residue modulo 4 of the term $(\Pi^+- \Pi^- -n)$
occurring in Equation~\eqref{Rohlin-to-Arnold:eq} to draw sharper
congruences (than Arnold's). In the above proof we ignored
completely (and could do so) the sign of Rohlin's pairs, yet there
is an evident composition law when the Hilbert tree of the scheme
is decorated by signs (dictated by Rohlin's pairs with complex
orientations, see Fig.\,\ref{Rohlinsformula:fig}). It seems
however unlikely that one can boost the method up to include a
proof of Gudkov hypothesis based on Rohlin's formula. (If
feasible, this would have certainly been mentioned in the
Degtyarev-Kharlamov survey \cite{Degtyarev-Kharlamov_2000}).

If optimistic one could use a total pencil (like in
Theorem~\ref{total-reality-of-plane-M-curves:thm}) as  to control
 complex orientations. This could give an elementary proof of
Gudkov hypothesis via basic algebraic geometry instead of highbrow
topology (like Rohlin 1952, or Marin 1977--79). Let us look if
Rohlin's formula implies the Gudkov-Rohlin congruence. As above,
we start from the Rohlin-to-Arnold equation
\eqref{Rohlin-to-Arnold:eq}
$$
\chi=k^2+2(\Pi^+-\Pi^- -n),
$$
and try now to control the residue modulo $4$ of the corrector
term $(\Pi^+-\Pi^- -n)$ under the assumption that the curve
$C_{m=2k}$ is an $M$-curve. (En passant, it seems that this
corrector term is always $\le 0$ by Thom's conjecture, cf.
Theorem~\ref{Thom-Ragsdale:thm}. [29.03.13] Warning: this is
false!) If one is able to show that this corrector term is $0$
modulo 4 then Gudkov hypothesis follows.

As in the previous reduction of Arnold-to-Rohlin, we consider the
hierarchy of the scheme (alias tree or mushroom), but now one
takes into account complex orientations. The latter induce a
distribution of signs on all injective pairs of the tree. (An
injective pair is any hierarchical pair $x>y$ like in Hegel's
dialectic ``du ma\^{\i}tre et de l'esclave'' where $x$ is not
necessarily the direct superior of $y$.)

\begin{lemma} \label{Signs-law:lem} {\rm (Signs-law)}.---In the Rohlin tree with (injective) pairs decorated with signs
$\pm 1=\sigma_{x,y}$ we have given two (composable) pairs $x<y$,
and $y<z$ the following ``twisted'' signs-law:
$$
\sigma_{x,z}=(-1)\sigma_{x,y} \cdot \sigma_{y,z}.
$$
This looks a priori exotic, as it amounts to say that $+\times
+=-$, $-\times-=-$, i.e. consanguinity is bad, while mixing the
genes is good, i.e.  $+\times -=+$, $-\times+=+$.
\end{lemma}

\begin{proof}
This exotic signs law  is justified by looking at
Fig.\,\ref{Signs-law-dyad:fig}:

\begin{figure}[h]
\centering
    \epsfig{figure=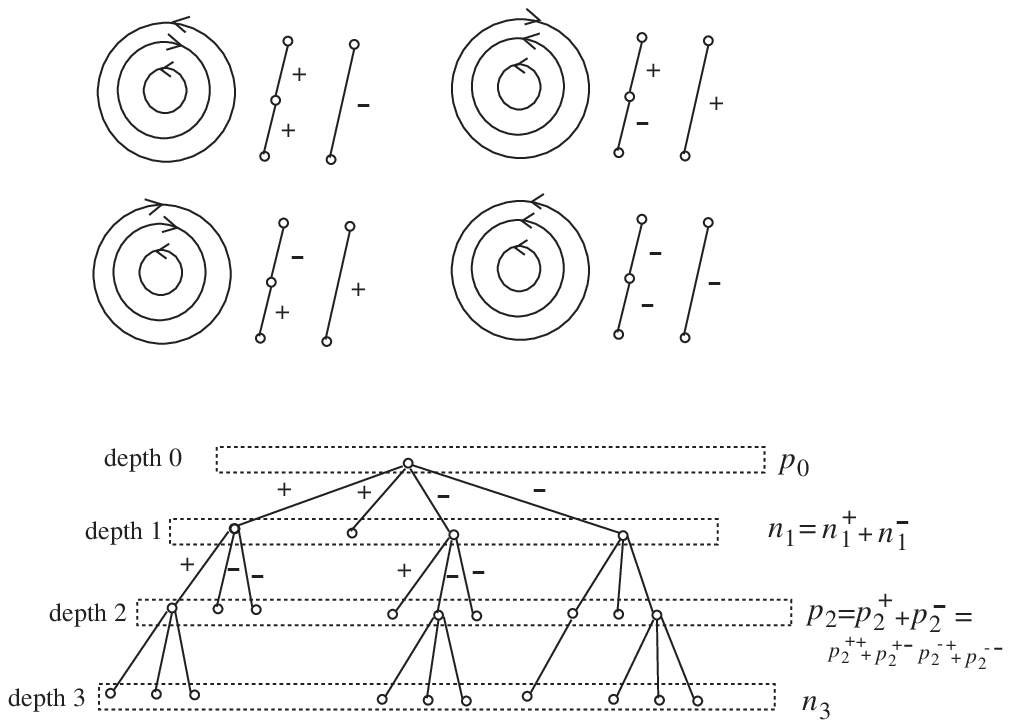,width=102mm}
\vskip-5pt\penalty0
\caption{\label{Signs-law-dyad:fig}%
Signs-law for composable edges of the Rohlin tree}
\vskip-5pt\penalty0
\end{figure}
\end{proof}

Of course the boring sign $(-1)$ in the lemma could be avoided if
we flipped the convention in Rohlin's definition, but we are too
conservative to risk such a modification. Actually, Rohlin's
convention is perfectly sound, cf. again
Fig.\,\ref{Rohlinsformula:fig}.

Next, it may be observed that the difference $\Delta
\Pi=\Pi^+-\Pi^-$ computed locally in reference to a pair of
consecutive edges is $\sigma_{x,y}+\sigma_{y,z}+\sigma_{x,z}=
\sigma_{x,y}+\sigma_{y,z}+(-1)\sigma_{x,y} \cdot \sigma_{y,z}$
which is always either $+1$ or $-3$ (compare again
Fig.\,\ref{Signs-law-dyad:fig}), hence always $+1$ modulo 4.

Globally on the whole Rohlin tree, we have the formula
$$
\Delta \Pi:=\Pi^+-\Pi^-= \sum_{all\; edges (x<y)} \sigma_{x,y}.
$$
One could hope via the signs-law to evaluate this modulo 4, and
all should boil down to $n$ modulo $4$ under the assumption of an
$M$-curve. Along each triad $x<y<z$ the contribution is $+1$
modulo 4. However it becomes soon a combinatorial mess, and
 one cannot hope this in full generality, as
otherwise the Gudkov-Rohlin congruence mod 8 would hold for all
dividing curves and not merely for $M$-curves. This violates
experimental knowledge,  e.g. the Gudkov-Rohlin table in degree 6
(Fig.\,\ref{Gudkov-Table3:fig}).

A naive idea is to write down a cumbersome formula evaluating
$\Delta \Pi$, but alas this still does not use the $M$-curve
assumption. Maybe this is a first necessary step unless one has
some better idea.

Denote as before $p_0,n_1,p_2,n_3,p_4,\dots$ the number of ovals
at depths $0,1,2,3,4,\dots$ respectively. Using Rohlin's signs we
define an oval at depth $\ge 1$ as positive if the edge
immediately above it is a positive pair and as negative otherwise.
Accordingly we get splittings:
\begin{align}\label{splitting-rel:eq}
n_1&=n_1^+ + n_1^-\cr
p_2&=p_2^+ + p_2^- =p_2^{++}+p_2^{+-}+p_2^{-+}+p_2^{--} \cr
n_3&=n_3^+ + n_3^-
=n_3^{++}+n_3^{+-}+n_3^{-+}+n_3^{--}=n_3^{+++}+n_3^{++-}+etc.,
\end{align}
where $n_1^+$ is the number of oval at depth 1 which are positive,
$p_2^{++}$ is the number of ovals at depth 2 such that the 2 edges
right above it are positive, while $p_2^{+-}$ is the number of
ovals at depth 2 surmounted by 2 edges of signs $+$ and $-$ (in
this order when moving up), and so on. Once all this notation is
introduced we can write down a cumbersome formula for $\Delta \Pi$
enumerating all edges (=injective pairs) weighted by their signs
according to the depth of their starting-point:
\begin{align*}
\Delta \Pi =& n_1^+ - n_1^-\cr
&+p_2^+ -p_2^- + (-p_2^{++}+p_2^{+-}+p_2^{-+}-p_2^{--})\cr
&+n_3^+ -n_3^- + (-n_3^{++}+n_3^{+-}+n_3^{-+}-n_3^{--})\cr
&+(+n_3^{+++}-n_3^{++-}-n_3^{+-+}+n_3^{+--}
-n_3^{-++}+n_3^{-+-}+n_3^{--+}-n_3^{---})\cr &+etc.
\end{align*}
Using the splitting relations above \eqref{splitting-rel:eq} this
can be somewhat condensed  as
\begin{align*}
\Delta \Pi =& n_1^+ - n_1^-\cr
&+ 2p_2^{+-}-2p_2^{--}\cr
&+(+n_3^{+++}-n_3^{++-}+n_3^{+-+}-n_3^{+--}
-n_3^{-++}+n_3^{-+-}-n_3^{--+}-3n_3^{---})\cr
&+etc.
\end{align*}
alas some intelligence is required to decipher the hidden
structure. Even if properly done we still require to put into
action the $M$-curve assumption. Since each non maximal vertices
of the Rohlin tree defines a unique edge above it we have the
relation $r=p_0+ number\; of\; edges$. This is only a weak grip.

All this mess is just given as to motivate someone to arrange a
combinatorial proof of the Gudkov hypothesis $\chi\equiv k^2 \pmod
8$ on the basis of Rohlin's formula alone. This seems  a serious
combinatorial challenge. Since all  classical proofs---(i.e., the
first erroneous one of Rohlin 1972, the latter one by Rohlin 1974
via Atiyah-Singer, plus the Marin-Guillou rescue of Rohlin's
original misproof, yet still via an extension of Rohlin's deep
result on signatures of spin 4-manifolds)---use some deep results
it is quite unlikely that our naive programme can be completed.
Still someone gifted in combinatorics with a clever idea on how to
exploit the $M$-curve assumption (Harnack maximality) can perhaps
crack the problem in a very elementary fashion. A vague suggestion
is to exploit the total reality result for $M$-curves given in
Theorem~\ref{total-reality-of-plane-M-curves:thm} prompting
perhaps some information on complex orientations via the usual
dextrogyration argument.

Another project along this reductionism to Rohlin's formula would
be to attack our naive converse conjecture to the RKM-congruence,
cf. Conjecture~\ref{RKM-converse:conj}. But this is merely a naive
conjecture quite unlikely to hold true.

\subsection{A
trinity of congruences: Gudkov-Krakhnov-Khar\-lamov and
(Rohlin)-Kharlamov-Marin}

[07.01.13] To prohibit $(M-1)$-schemes (above Gudkov's broken
line) one can use (beside the Hilbert-Rohn-Gudkov method) the
following analogue congruence (paralleling that of Gudkov-Rohlin)
due to Kharlamov 1973/73 \cite{Kharlamov_1973/73} and
independently Gudkov-Krakhnov 1973/73
\cite{Gudkov-Krakhnov_1973/73}:

\begin{theorem} {\rm (Kharlamov 1973, Gudkov-Krakhnov 1973)}
\label{Gudkov-Krakhnov-Kharlamov-cong:thm} A plane $(M-1)$-curve
of degree $m=2k$ satisfies the congruence
$$
\chi=p-n\equiv k^2\pm 1 \pmod 8.
$$
\end{theorem}

\begin{proof}
Several proof are available:

$\bullet$ The original sources just referred to.

$\bullet$ Since $(M-1)$-curves are not dividing the technique is
different from the capping-off trick \`a la Arnold-Rohlin. However
Marin is able to get an unified proof (\`a la Rohlin) by using the
Guillou-Marin extension of Rohlin's signature formula. For an
exposition cf. A'Campo 1979 \cite{A'Campo_1979}.
%
\end{proof}

[11.01.13] We have also the following remarkable congruence (due
independently  to Kharlamov and Marin (first reported in print in
Rohlin 1978 \cite[3.4]{Rohlin_1978} and the first detailed proof
is given in Marin 1979/80 \cite{Marin_1979}):

\begin{theorem} {\rm (Kharlamov 197?, Marin 1979/80, first reported
in print in Rohlin 1978)} \label{Kharlamov-Marin-cong:thm}
 A plane $(M-2)$-curve of degree $m=2k$
and type~II satisfies the congruence
$$
\chi=p-n\equiv k^2 \textrm{ or } k^2\pm 2 \pmod 8.
$$
This can be paraphrased by saying that an $(M-2)$-curve with
$\chi\equiv k^2+4 \pmod 8$ is necessarily of type~I.
\end{theorem}

\begin{proof} Compare Rohlin 1978 \cite[3.4]{Rohlin_1978} or Marin
1979/80 \cite{Marin_1979}, or also Kharlamov-Viro 1988/91
\cite{Kharlamov-Viro_1988/91} for an unified account (and various
approaches). For the paraphrase either look at the Gudkov table in
degree 6, or more seriously do some boring arithmetics, cf.
(\ref{RKM-congruence-reformulated:thm}), which we reproduce
quickly. The paraphrase follows from the fact that {\it an
$(M-2)$-curve of order $m=2k$ verifies universally $\chi \equiv
k^2 \pmod 2$.} This is easy to prove using the relation
$\chi=p-n$, $r=p+n=M-2$
$$
\chi=p-n=(p+n)-2n\equiv_2 p+n = r=M-2,
$$
while by Harnack's bound and the genus formula
$g=\frac{(m-1)(m-2)}{2}$ we have
$$
M=g+1=\textstyle\frac{(2k-1)(2k-2)}{2}+1=(2k-1)(k-1)+1=2k^2-3k+2,
$$
whence
$$
\chi\equiv_2 M-2=2k^2-3k \equiv_2 k \equiv_2 k^2.
$$
\end{proof}

Specializing to sextics ($m=6$, so $k=3$) implies the following
(compare Fig.\,\ref{Gudkov-TableTop:fig}):

\begin{cor} {\rm (Rohlin 1978)}
The two real sextic schemes $\frac{6}{1}2$ and $\frac{2}{1}6$ are
of type~I.
\end{cor}

It seems that this result had no classical counterpart \`a la
Hilbert-Rohn-Gudkov prior to the Arnold-Rohlin revolution.
Nonetheless Rohlin 1978 \cite{Rohlin_1978} mentions the
possibility of a synthetic argument involving pencil of cubics.
More on this in  Sec.\,\ref{total-(M-2)-schemes:sec}.

[10.03.13] In fact this synthetic argument of Rohlin is now lost
but was partially reconstructed by Le~Touz\'e 2013
\cite{Fiedler-Le-Touzé_2013-Totally-real-pencils-Cubics}. This
issue of Rohlin-Le~Touz\'e should have a strong interaction with
Ahlfors theorem, while affording the first nontrivial phenomenon
of total reality. We will have the occasion to dwell more on this
later in this text.

[10.03.13] Another little remark (valid in degree 6 but perhaps
more universally) is the:

\begin{lemma} Once we know Arnold's congruence and the
Gudkov-Krakhnov-Kharlamov(=GKK) congruence then  Gudkov hypothesis
follows formally.
\end{lemma}

\begin{proof} Indeed contemplating  Gudkov's table
(Fig.\,\ref{Gudkov-Table3:fig}), the Arnold congruence prohibits
all while-colored $M$-schemes safe those at the center of the
semi-hexagons (i.e. $\frac{7}{1}3$ and it mirror $\frac{3}{1}7$).
Further the GKK-congruence prohibits all white-colored
$(M-1)$-schemes on that same Gudkov table. Hence if one of the two
schemes $\frac{7}{1}3$ (or its mirror) existed, it would appear as
an isolated island in the ocean. Yet, by transversality (\`a la
Bertini-Morse-Whitney-Sard-de~Rham-Thom\footnote{Thom learned Sard
from de Rham, cf. the 1954 Commentarii article.}) a generic pencil
of curves through an (hypothetical) algebraic representant and any
other curve with less oval, e.g. the anti-Fermat (invisible) curve
with zero oval (equation $x_0^6+x_1^6+x_2^6=0$) would produce a
combinatorial path on the Gudkov table. This violates isolation.
(NB: even an eversion (Sec.\,\ref{Eversion:sec}) can only take the
scheme to its mirror, and so elementary Morse surgeries
necessarily create an elementary path on the Gudkov table).
\end{proof}

[10.03.13] A more radical intuition of Rohlin 1978 (now partially
justified by Le~Touz\'e 2013
\cite{Fiedler-Le-Touzé_2013-Totally-real-pencils-Cubics}) is that
owing to their total reality the $(M-2)$-schemes $\frac{6}{1}2$
(and its mirror) are maximal. This explains all the prohibitions
materialized by the  (white) semi-hexagons on the Gudkov Table
(Fig.\,\ref{Gudkov-Table3:fig}), safe the 2 schemes $11$ and
$\frac{10}{1}$ that were prohibited since the Hilbert-Rohn era (at
least modulo some German sloppiness, made perfectly rigorous by
Academician D.\,A. Gudkov). Nowadays prohibiting them is a trivial
consequence of either Arnold's congruence or Rohlin's formula
(\ref{Rohlin-formula:thm}). So at least in degree 6, we see that
the phenomenon of total reality acts as a strong unifying
principle for classical prohibitions. Rohlin probably had the
intuition that this phenomenon perpetuates in higher degrees. More
along this vertiginous idea (potentially allied to Ahlfors
theorem) will be discussed in
Sec.\,\ref{Esquisse-dun-prog-deja-esquiss:sec}.

\subsection{Total reality of the two maximal sextic $(M-2)$-schemes
(Rohlin 1978, Le~Touz\'e 2013)} \label{total-(M-2)-schemes:sec}

[03.01.13] This is akin to Ahlfors' theorem, yet somewhat
different and actually the hard part of the game. Rohlin 1978
\cite[p.\,94]{Rohlin_1978} writes the following cryptical note:

\begin{quota}\label{Rohlin1978-total-reality:quote}
{\rm (Rohlin 1978)}
{\rm \small {\sc Note on the method.} After a suitable
modification, these arguments can be applied to some other
schemes. For example, when we replace real lines by real curves of
degree 2 we can establish that a real scheme of degree~8
consisting of 4 nests of depth 2 lying outside one another belongs
to type~I\footnote{[29.03.13] This is just a very special case of
a more general satellite principle, cf.
Sec.\,\ref{satellite-total-reality:sec}.}, and when we apply it to
curves of degree 3, we can establish (in a rather complicated way)
that the schemes $\frac{6}{1}2$ and $\frac{2}{1} 6$ of degree 6
[considered below in \S 3.8]\footnote{Omit this bracketing for it
is just to refer to Gudkov's notation.} belong to type~I. However,
all the schemes that we have so far succeeded in coping with by
means of these devices are covered by Theorem~3.4 and 3.5.

}
\end{quota}

What is crucial here is the parenthetical comment ``(in a rather
complicated way)''. This is highly  reminiscent of some Hilbertian
prose ``{\it freilich auf einem au{\ss}erordentlich
umst\"andlichen Wege\/}'', cf. Hilbert 1891 (p.\,418, in Ges.
Abh., Bd.\,II)): ``{\it Diesen Fall $n=6$ habe ich einer weiteren
eingehenden Untersuchung unterworfen, wobei ich\,---\,freilich auf
einem au{\ss}erordentlich umst\"andlichen Wege\,--- fand, da{\ss}
die elf Z\"uge einer Kurve 6-ter Ordnung keinesfalls s\"amtlich
au{\ss}erhalb un voneinander getrennt verlaufen k\"onnen. Dieses
Resultat erscheint mir deshalb von Interesse, weil er zeigt,
da{\ss} f\"ur Kurven mit  der Maximalzahl von Z\"ugen der
topologisch einfachste Fall nicht immer m\"oglich ist.\/}'' Of
course both problems are slightly different but perhaps there is
some common difficulty in both  games, while it is not impossible
that Rohlin's made a direct winking at Hilbert's prose.

So Rohlin claims being able to prove the following (on the basis
of pure geometry):

\begin{theorem} {\rm (Rohlin 1978, no published proof)}
The two real sextic $(M-2)$-schemes  $\frac{6}{1}2$ (\/$6$ ovals
encapsulated in one oval and $2$ outsides) and $\frac{2}{1} 6$
(\/$2$ ovals encapsulated in one oval and $6$ outsides) are of
type~I, i.e. any smooth real curve realizing one of those schemes
is necessarily orthosymmetric (=dividing) in the sense of Klein.
\end{theorem}

On the basis of Rohlin's Quote (right above) one guesses that the
proof involves looking at a pencil of cubics through 8 points
inside the deep ovals while checking total reality of the
resulting morphism to the line. (As usual it results a circle map
\`a la Ahlfors, which is of degree $3\cdot 6-8=10$ after
degenerating the basepoints on the ovals. This quantity coincides
with Gabard's bound
$r+p=\frac{r+(g+1)}{2}=\frac{(g-1)+(g+1)}{2}=g=10$.) Of course, it
is quite sad that Rohlin did not found the place to write down the
details.

Naively the proof could be as follows. Take a cubic in the pencil
based at some  8 points inside the $2+6=8$ deep ovals
(equivalently those containing no ovals). If the cubic is
connected then it visits all 8 points. Counting intersection we
have $2\cdot 8=16$ intersections coming from the deep ovals, plus
two intersections coming from traversing the enclosing oval. This
gives 18 the maximum permitted by B\'ezout, whence the desired
total reality.

This looks simple, but this by no mean a complete argument. What
to do if the cubic is not connected? One could of course try to
arrange a pencil of connected cubics. Recall that the discriminant
parametrizing singular cubics has degree $3(k-1)^2=3 \cdot 2^2=12$
of even degree. Thus there is no ``Galois'' obstruction to finding
a line in the space of cubics $\vert 3 H \vert \approx {\Bbb P}^9$
missing the discriminant.

\begin{figure}[h]
\centering
    \epsfig{figure=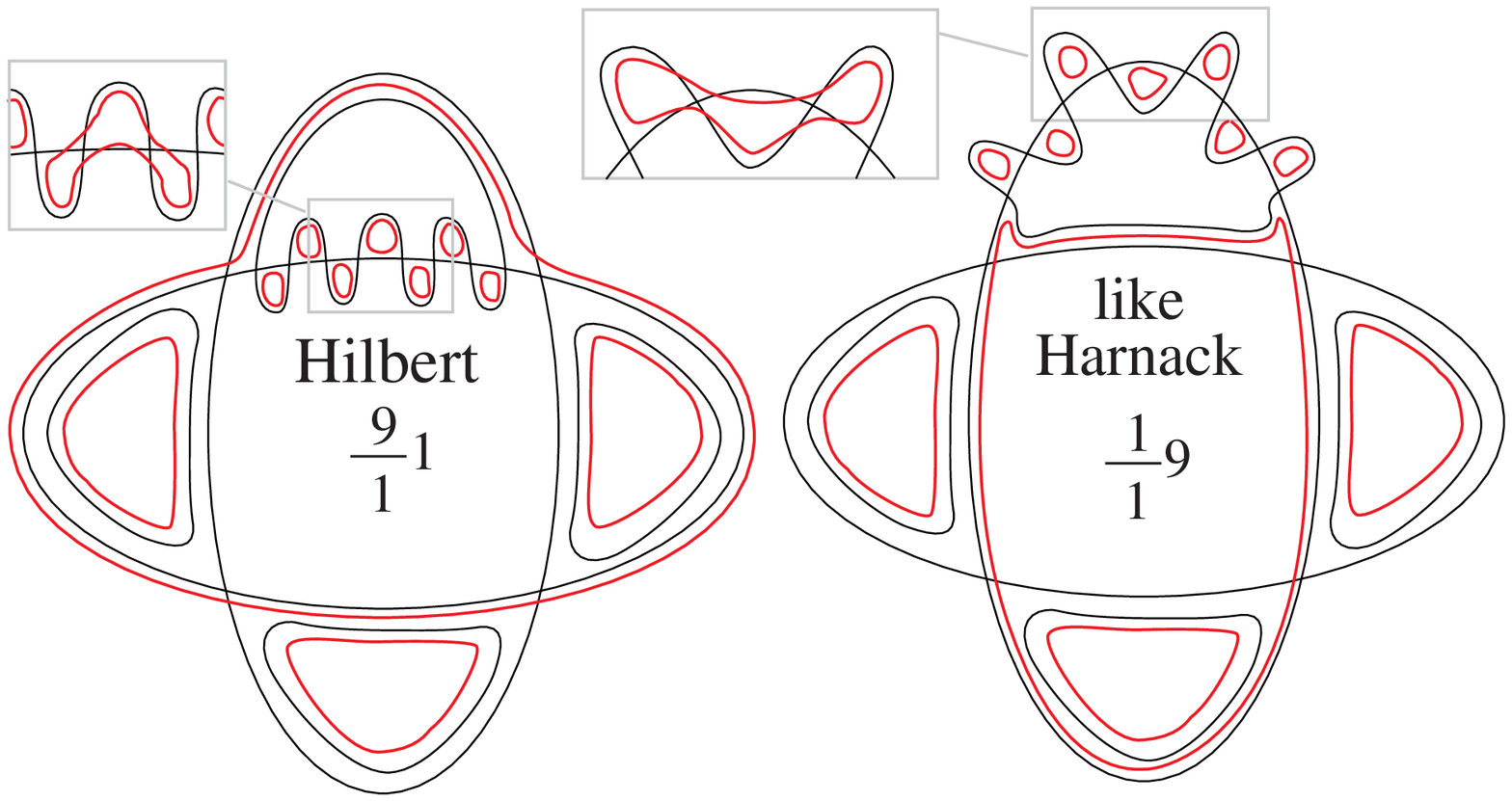,width=72mm}
\vskip-5pt\penalty0
\caption{\label{TypeII:fig}%
Killing ovals in the $M$-sextics of Hilbert and ``Harnack''}
\vskip-5pt\penalty0
\end{figure}

Another objection is that our simpleminded proof  equally well
applies to all other $(M-2)$-schemes excepted $9$. Indeed this is
clear for all of them since the ovals are split in two packets by
the enclosing oval. In the case of $\frac{8}{1}$ the enclosing
oval is also necessarily cut twice, since $C_3({\Bbb R})$ is not
null-homotopic. This would imply  that all $(M-2)$-schemes safe
the unnested one ($9$) are of type~I. This is however too radical
and incompatible with experience (or with theory, e.g. Arnold's
congruence). For instance it is easy to alter Hilbert's method to
get the scheme $\frac{7}{1}1$ in a nondividing way as switched
some signs of smoothing (cf. Fig.\,\ref{TypeII:fig}).

[08.01.13] Maybe a fruitful idea is to look at special pencils of
cubics  spanned by 2 reducible (split) cubics well-understood. For
instance one could try with reducible cubics aggregating a conic
and a line. (This looks too rigid as a conic can pass only through
5 points while the 3 remaining one are not necessarily aligned.)
Another idea is to look for rational cubics (uni\-nodal).

Inside any cubics pencil (e.g. one based on the 8 deep ovals), we
have by the theorem of Cheponkus-Marin 1988 (cf. Marin 1988
\cite[p.\,192]{Marin_1988}) a curve with at most
$M-3=(g+1)-3=(1+1)-3=-1$ components. This looks exotic and wrong
for cubics for it would lack real circuits violating thereby
Galois-B\'ezout!
Shame on me, I forgot that Cheponkus-Marin assume the pencil to be
even degree $>2$!

However it is a well-known
fact that a generic pencil of cubics has a curve with $r=1$, which
sounds quite likely (though not forced by the even degree 12 of
the discriminant). Even accepting this we get only one curve of
the pencil which is connected, then of course a family of such,
but we are by no mean ensured that all will be so.

The case of the G\"urtelkurve, or that of the sextic having a nest
of depth 3,  inclines to believe that the total reality of
Rohlin's sextic schemes enjoy some structural stability in the
sense that it is enough to assign the 8 basepoints inside the deep
nests to ensured total reality. So we are  looking for something
quite rare yet reasonably robust.

Suppose we have some connected cubics in our pencils.
When moving in the pencil it may splits in two components. This
can occur (assuming genericity, i.e. transversality to the
discriminant) either through the birth of an oval after crossing a
solitary node or by a self-coalescence of the pseudoline of the
cubic crossing in this case a non-solitary node (with two real
branches).

$\bullet$ In the first case the newly created oval cannot
contribute to additional intersections.

Here is our argument requiring some hypothesis. First our generic
pencil will exhibit at most 12 uninodal curves (either solitary or
non-solitary nodes). By general position it may be assumed that
those 12 curves (as well as the allied 12 nodes) are not located
on the given sextic $C_6$. When crossing the discriminant, the
solitary node will inflate into a little surrounding oval (or at
least nearby ovals). Thus by continuity this little oval do not
interact with the $C_6$. In fact what is important is that our 12
nodes are not on the 8 basepoints. Hence our just born oval do not
contribute visiting the base locus, which is therefore entirely
visited by the residual pseudoline. This is connected and so total
reality is preserved (even after crossing the wall), at least at
some instantaneous future right after the crossing. Then the
little bubble (oval) can inflate, yet as the number of real
intersections is already maximum via the pseudoline, the oval
cannot cross the $C_6$. Its motion is in some sense confined to
its complement (of the $C_6$). Then two scenarios are possible,
either the oval deflates again to some solitary points or it
merges with the pseudoline. In both cases we come back to a curve
with one circuit and total reality is ensured for topological
reasons. This story has to be repeated perhaps 12 times but we
seem finished, modulo analyzing the other case.

$\bullet$ In the second case (real normal crossing) our pseudoline
of the moving cubic self-collides with itself and then splits in
two branches. Of course it may then result a loss of total
reality. Imagine the crossing (non solitary node) to be located on
the $C_6$, then at the critical time there are locally two
intersections and soon afterwards these may disappear loosing two
intersections. Yet assuming by general position that our nodes are
never located on the $C_6$ (after eventually perturbing the 8
basepoints) we still have right after the critical time two real
intersections, and total reality is conserved in the short run.
Now our curve is decomposed in two branches, and accordingly so
are the basepoints. If the pseudoline (at the post-critical level)
contains a mixture of points both inside and outside the nonempty
oval of $C_6$ we are happy for two extra intersections are created
while entering in and evading out this separating oval. If not,
then the oval of $C_3$ visits the 6 inner ovals of $C_6$ and the
pseudoline the 2 outer oval of $C_6$. (The other situation is also
possible but then total reality is obvious for the pseudoline must
evade the surrounding oval.) Now our oval in the long run may
loose two intersections. However as
the intersection
$C_3\cap C_6$ was totally real, this can only occur by a
retraction of a tongue slipping
inside the separating oval, and then the 6 inner ovals of $C_6$
are trapped inside an oval of a cubic. The latter is reasonably
rigid and convex. This oval still has $2\cdot 6=12$ intersection
with the 6 inner ovals. In particular it cannot shrink to a point.
It cannot also evolve to an ellipse, as otherwise the 6 inner
basepoints would be co-elliptic, which can be avoided by general
position. OF COURSE the proof is still not finished and some idea
need to be discovered. Perhaps let pass a conic through 5 of the
inner basepoints, while noting that a 6th intersection is created
by Galois-B\'ezout, etc... Of course all this needs much more
substantiation!

[08.01.13] The end of Rohlin's Quote (above) shows that there must
be alternative non-synthetic but crudely speaking topological
proofs of the theorem. This is indeed what was discussed in the
previous section.


\section{Constructions}

[29.03.13] Construction of algebraic curves seems a syllogism
since they are nearly God-given. Perhaps the word contemplation
looks more appropriate, but clumsy. Despite existence of Gods, the
art of tracing of algebraic (plane) curves goes back to time
immemorial Newton, Pl\"ucker 1839 \cite{Plücker_1839}, Zeuthen
1874 \cite{Zeuthen_1874}, with the modern era usually identified
by Harnack 1876 \cite{Harnack_1876}, Klein 1876, Hilbert 1891
\cite{Hilbert_1891_U-die-rellen-Zuege}). The game is especially
interesting over $\RR$, else nearly everything follows from
Riemann.

Much of the elementary aspects can be treated by the primitive
method of small perturbation, which nearly gives a good picture of
what happens in degree 6. This is how worked Pl\"ucker, Klein,
Harnack, Hilbert, Ragsdale, Brusotti, etc. However already in
degree 6 the classical method starts showing some limitation.
Albeit the Gudkov curve can be distilled by small perturbation, it
requires an extra twist by means of Cremona transformations (at
first difficult to visualize). It took the community ca. 8 decades
(including such masters as Hilbert, Ragsdale, Brusotti,
Petrovskii, Gudkov first not an exception)  until to discover the
fairly trivial picture traced by Gudkov ca. 1972
(Fig.\,\ref{GudkovCampo-5-15:fig}) exhibiting a curve with
topology $\frac{5}{1}5$.

Ca. 1980 Viro  described how to dissipates more complicated
singularities, allowing experts to create more funny curves
refuting most of the conjectures erected along the primitive
method. For instance also Gudkov's sextic appears fairly trivially
when one knows how to smooth a triplets of ellipses tangent at 2
points, compare Fig.\,\ref{Viro3-15:fig}c. A variant of Viro's
patchwork
 due to Itenberg (called the $T$-construction) is
purely combinatorial and permitted to disprove severely the
Ragsdale conjecture (cf. Fig.\,\ref{Itenberg:fig}), as well as our
naive Thom estimate $\chi\le k^2$ for dividing curves.

\subsection{Constructing the two maximal $(M-2)$-schemes}
\label{const-total-(M-2)-schemes:sec}

\begin{figure}[h]
\centering
    \epsfig{figure=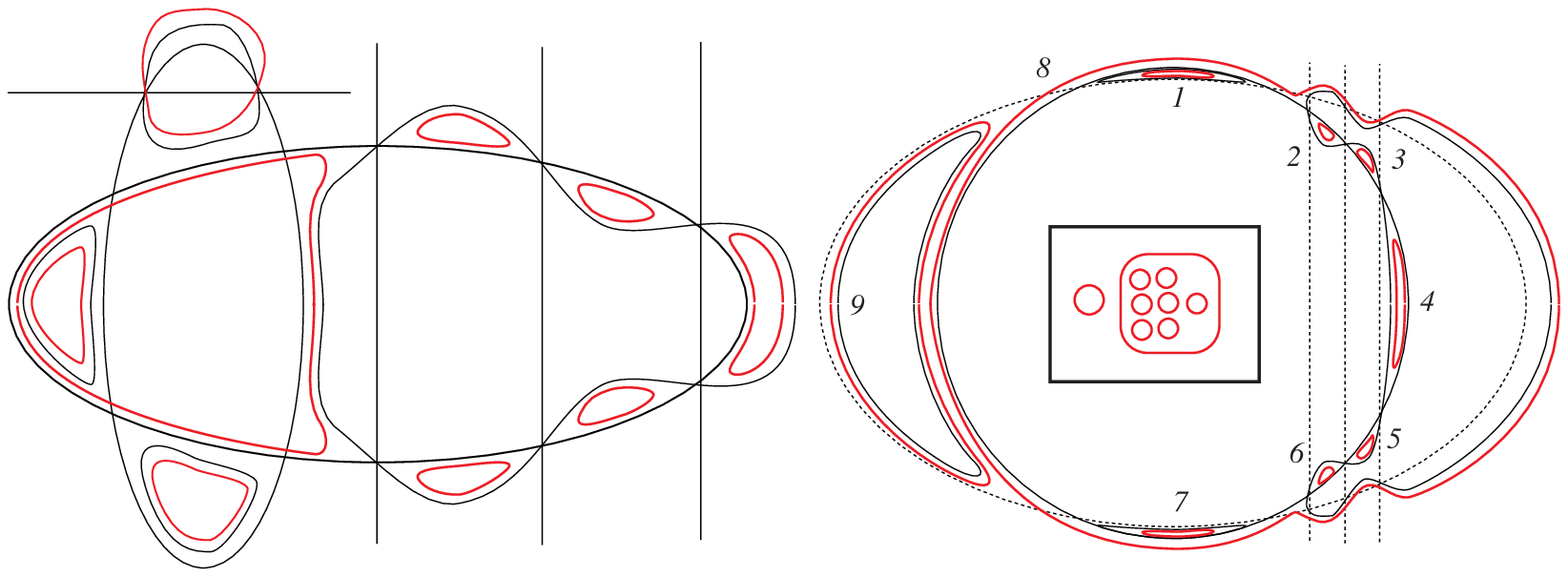,width=122mm}
\vskip-5pt\penalty0
\caption{\label{Gudkov-Hilbert-modified:fig}%
A picture from Gudkov 1974, but not convinced by his text.}
\vskip-5pt\penalty0
\end{figure}

[05.01.13] As to the existence of the two maximal $(M-2)$-schemes
(namely $\frac{6}{1}2$ and $\frac{2}{1}6$), they can be
constructed (as observed in Gudkov 1974/74
\cite[p.\,42]{Gudkov_1974/74}) by a slight modification of
Hilbert's method. Let us reproduce his figure (Fig.\,4, p.\,16).
This gives (after smoothing) the left-side of
Fig.\,\ref{Gudkov-Hilbert-modified:fig}. Alas, this is not the
desired scheme. Is there a mistake in Gudkov at this place?
Apparently not as it seems approved in A'Campo 1979
\cite{A'Campo_1979} (alas no detail).

\noindent The right-part of
Fig.\,\ref{Gudkov-Hilbert-modified:fig} is just a variant inspired
from Hilbert's configuration. This has again the wrong real
scheme. Since Gudkov does not seem to give exactly what he claims,
we must rely on some do-it-yourself endeavor. A naive idea gives
Fig.\,\ref{GudHilb2:fig}, failing again to have the correct
scheme.

\begin{figure}[h]
\centering
    \epsfig{figure=GudHilb2.eps,width=62mm}
\vskip-5pt\penalty0
\caption{\label{GudHilb2:fig}%
A variant of Hilbert's method} \vskip-5pt\penalty0
\end{figure}

\begin{figure}[h]
    \hskip-1.0cm\penalty0
    \epsfig{figure=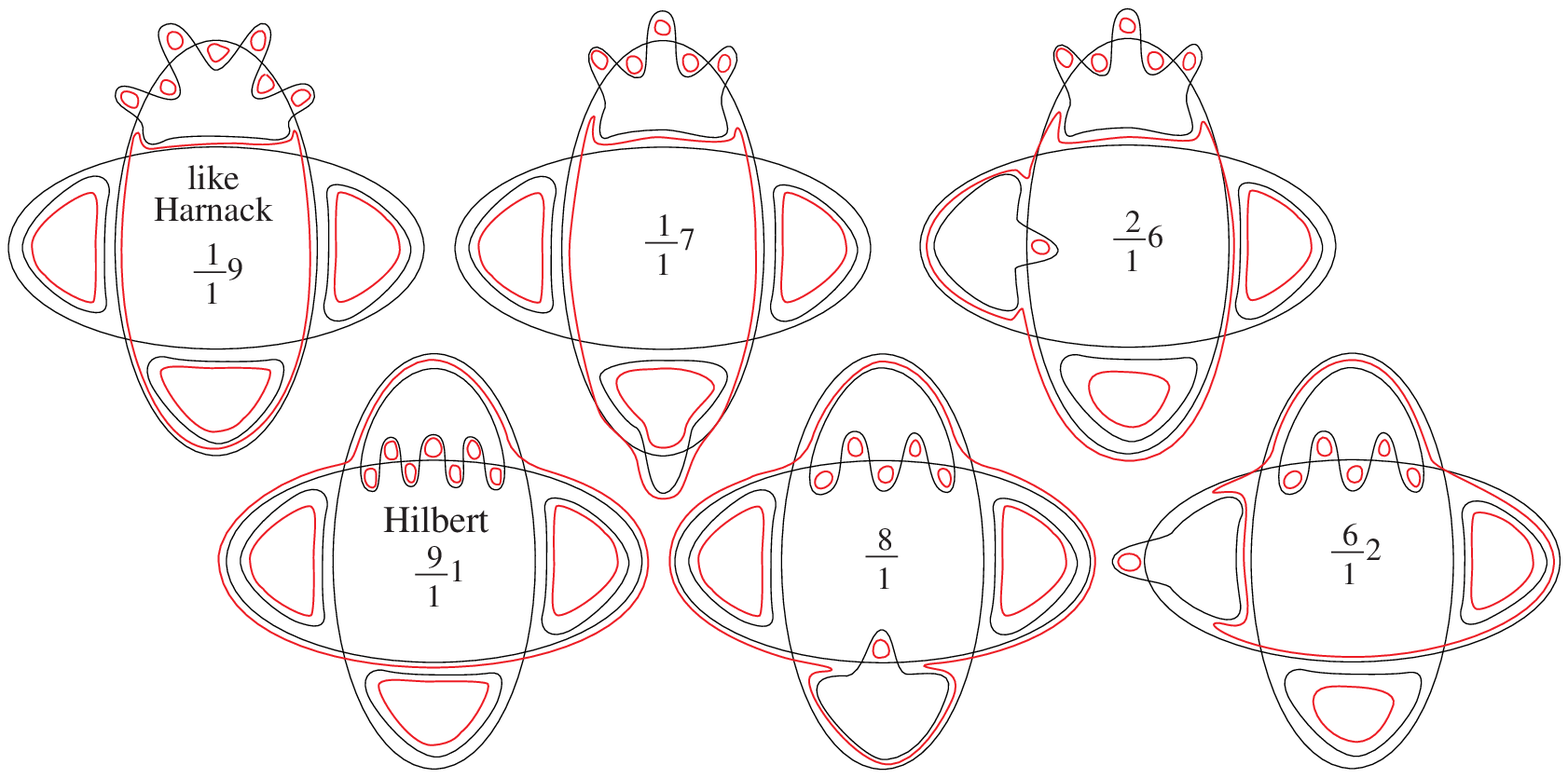,width=142mm}
\vskip-5pt\penalty0
\caption{\label{GudHilb3:fig}%
Hilbert's method and some variants with several oscillating ovals}
\vskip-5pt\penalty0
\end{figure}

Let us now work more systematically. The key is first to make
Walt-Disney pictures of Hilbert's method \`a la Gudkov. This
involves a simplified  art-form, far from geometrically realist,
but topologically faithful and more malleable.
This produces the following pictures (Fig.\,\ref{GudHilb3:fig}).
The trick of Hilbert's method is to let  oscillate an oval across
an ellipse while smoothing their union (cf. left part of
Fig.\,\ref{GudHilb3:fig}). Such oscillations are B\'ezout
compatible: each oscillating quartic intersects 8 times the
ellipse. A posteriori it is a simple matter (Hilbert's method) to
realize such oscillation by rigid algebraic curves suitably
perturbed by lines arrangements. Yet it is valuable  first
exploring the softer topological figures as to understand which
oscillation is able producing a prescribed topology (e.g.
$\frac{6}{1}2$ and its companion $\frac{2}{1} 6$).

 The variant of Hilbert's construction  involves letting
oscillate various ovals across the ground ellipse. On the
middle row
of Fig.\,\ref{GudHilb3:fig} we let oscillate thrice one oval and
once the opposite oval with respect to some ground ellipse.
Smoothing gives some $(M-2)$-curves not realizing the desired
schemes. Choosing instead a triple oscillation of the upper oval
of the quartic combined with a simple oscillation of the nearby
oval gives the desired schemes (right row of
Fig.\,\ref{GudHilb3:fig}) with either 6 outer ovals (top) or 6
nested ovals (bottom). Gudkov was right albeit his discourse was
not in perfect adequation with his picturing.

It remains  to geometrize such oscillations \`a la Hilbert. This
is an easy matter, except that realist pictures require judicious
scalings to make things visible. The first mode of vibration
leading to $\frac{2}{1}6$ is geometrized on
Fig.\,\ref{GudHilb5:fig} below.

\begin{figure}[h]
\centering
    \epsfig{figure=GudHilb5.eps,width=122mm}
\vskip-5pt\penalty0
\caption{\label{GudHilb5:fig}%
Several oscillating ovals for $\frac{2}{1}6$ compare with
top-right of Fig.\,\ref{GudHilb3:fig}} \vskip-5pt\penalty0
\end{figure}

\noindent Some few comments on this figure: first one has the two
blue ellipses forming a quartic $C_4$. Next one has  4 dashed
lines (another quartic). Perturbing slightly the former along the
other (within the pencil spanned by both) gives another $C_4$
traced in black. This has the virtue of oscillating across the
circular ellipse (in blue). Finally, smoothing their union gives
the sextic in red realizing the desired real scheme $\frac{2}{1}
6$ (i.e. 2 ovals captured in one and 6 outside).

The other scheme $\frac{6}{1}2$ is obtained similarly via the
following system of oscillations (Fig.\,\ref{GudHilb6-12:fig})
geometrizing the bottom-right part of Fig.\,\ref{GudHilb3:fig}:

\begin{figure}[h]
\centering
    \epsfig{figure=GudHilb6-12.eps,width=92mm}
\vskip-5pt\penalty0
\caption{\label{GudHilb6-12:fig}%
Multi-oscillation for $\frac{6}{1}2$: compare with bottom-right
part of Fig.\,\ref{GudHilb3:fig}} \vskip-5pt\penalty0
\end{figure}

Using either the schematic pictures or the more geometric one it
is an easy matter to see that both curves just traced are
dividing. This follows as usual (Fiedler's law) by checking that
all smoothings are compatible with complex orientations (cf.
Fig.\,\ref{GudHilbdividing:fig} below).

\begin{figure}[h]
\centering
    \epsfig{figure=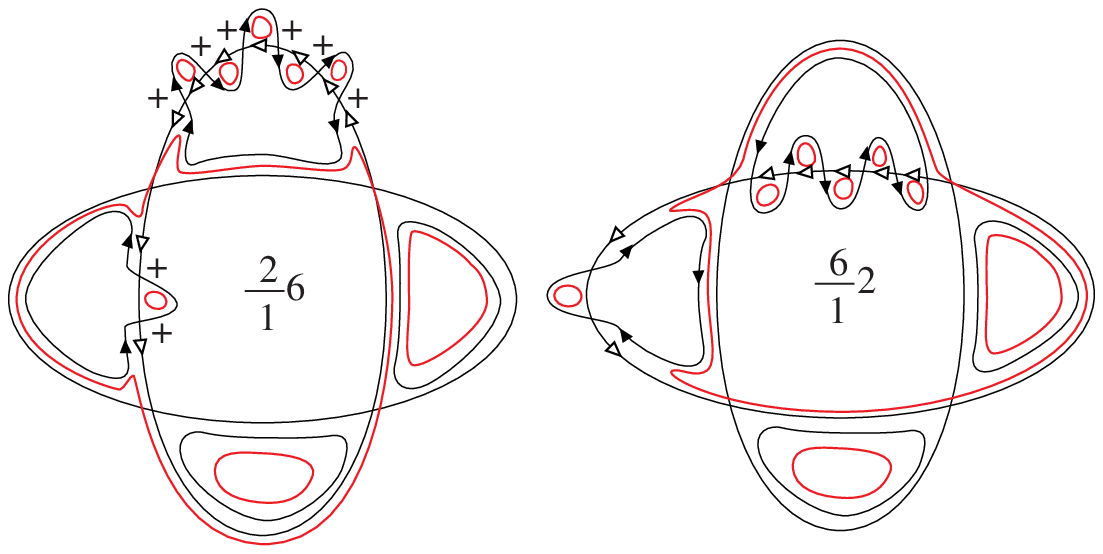,width=102mm}
\vskip-5pt\penalty0
\caption{\label{GudHilbdividing:fig}%
Checking the dividing character of Rohlin's curves via Fiedler's
law} \vskip-5pt\penalty0
\end{figure}

Hence according to Ahlfors theorem there must be a total pencil of
curves. Actually Rohlin claims much more that any sextic realizing
those schemes is totally real under a pencil of cubics, but his
argument has never been published\footnote{[29.03.13] A
nearby corroboration of Rohlin's claim is now available in
Le~Touz\'e 2013
\cite{Fiedler-Le-Touzé_2013-Totally-real-pencils-Cubics}.}. A bit
like Hilbert in
1900, Rohlin 1978 says that his proof is
too cumbersome to be written down. As we know
Hilbert's 2nd assertion that there is only two $M$-sextics  was
refuted by Gudkov some 7 decades latter, so it is not impossible
that Rohlin's claim is fallacious as well. Of course it can also
be the case that Rohlin's claim on the type~I of the maximal
$(M-2)$-schemes is correct, but that total reality via a pencil of
cubics is erroneous. However as we noted cubics leads to a
mapping-degree of $3 \cdot 6 -8=10$, in adequation with Gabard's
bound $r+p$ on the degree of circle maps.

\begin{figure}[h]
    \hskip-1.2cm\penalty0
    \epsfig{figure=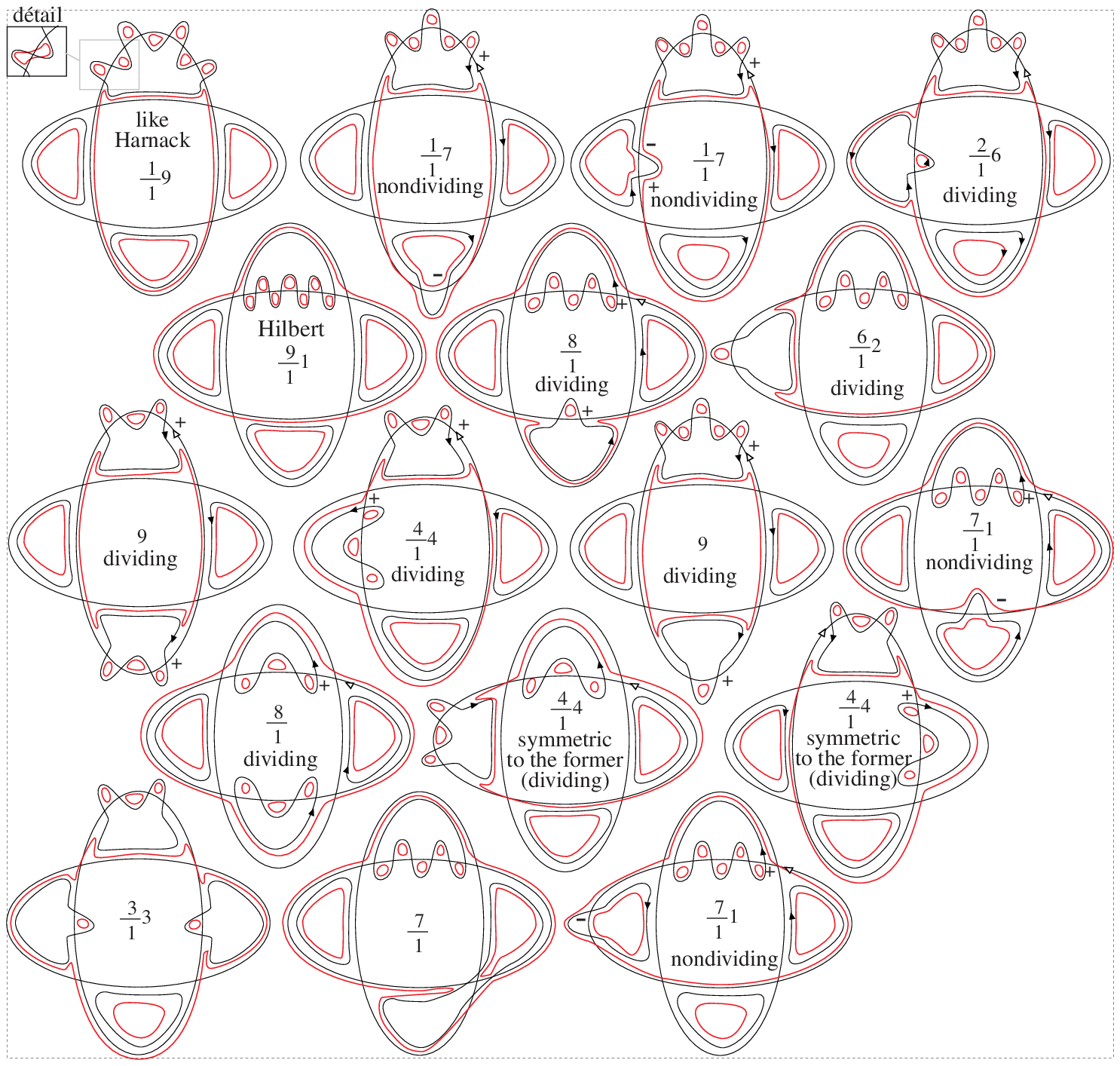,width=152mm}
\vskip-5pt\penalty0
\caption{\label{GudHilb8:fig}%
Flexible pictures with several oscillating ovals (variant of
Hilbert's)} \vskip-5pt\penalty0
\end{figure}

At this stage the naivest thing-to-do is to convince that there is
no trivial counterexample to Rohlin's claim. So we trace more
oscillations to get the following pictures
(Fig.\,\ref{GudHilb8:fig}). Some noteworthy species appear
especially the remarkable scheme $\frac{4}{1}4$, occupying the
central position of Gudkov's table
(Fig.\,\ref{Gudkov-Table3:fig}). Our specimen is dividing and it
looks hard to get the same scheme in the nondividing way (though
Rohlin 1978 asserts its existence). Speculating that this scheme
is of type~I, while admitting the truth of Rohlin's maximality
conjecture, then all 3 sextics schemes dominating $\frac{4}{1}4$
would agonize along a blue sky catastrophe! Gudkov would be wrong
and Hilbert right! Of course this seems a too apocalyptic
scenario, yet up to now our text does not entail this option!

Also difficult to find are the schemes $\frac{5}{1}3$ and its
mirror $\frac{3}{1}5$. Apart those exceptions, Hilbert's method
offers all possible $(M-2)$-schemes.

{\it Insertion} [08.02.13] For a Harnack method realization  of
$\frac{3}{1}5$, see Fig.\,\ref{HarnaGudkov3-15XXL:fig} much below,
while $\frac{5}{1}3$ truly requires the method of Gudkov (cf.
Fig.\,\ref{GudkovCampo-5-15:fig}).

If, via a small perturbation, one merges together $2$ small ovals
on Harnack's or Hilbert's curve (cf. d\'etail on
Fig.\,\ref{GudHilb8:fig}), then one gets the $(M-1)$-schemes
$\frac{1}{1}8$ resp. $\frac{8}{1}1$. The other $(M-1)$-schemes of
Gudkov's table are somewhat harder to exhibit, except of course if
one is aware of a large deformation able to extinct the inner oval
(case of Harnack) or the outer oval (in Hilbert's case). This
contraction of empty ovals is actually possible via Itenberg 1994
\cite{Itenberg_1994}---using the apparatus of Nikulin's  (1979/80
\cite{Nikulin_1979/80}) (rigid-isotopy classification via K3
surfaces)---but of course this is surely not the most economical
argument for our purpose (known to Gudkov 1969 or earlier).

It seems clear that we have exhausted the faculty of Hilbert's
method (and its variation where the vibration is dissipated on
several ovals). Some naive questions: what can be obtained by
perturbing an arrangement of lines (in general position)?

To go further one is helped once more by Gudkov 1974/74
\cite[p.\,42]{Gudkov_1974/74} asserting that the scheme
$\frac{4}{1} 5$ can be
gained by a modification of Harnack's method. Alas, Gudkov makes
no picture but was aware of this at least since 1954.

{\it Insertion} [08.02.13].---For a picture of this cf.
Fig.\,\ref{HarnaGudkov4-15:fig} much below.

\begin{figure}[h]
\centering
    \epsfig{figure=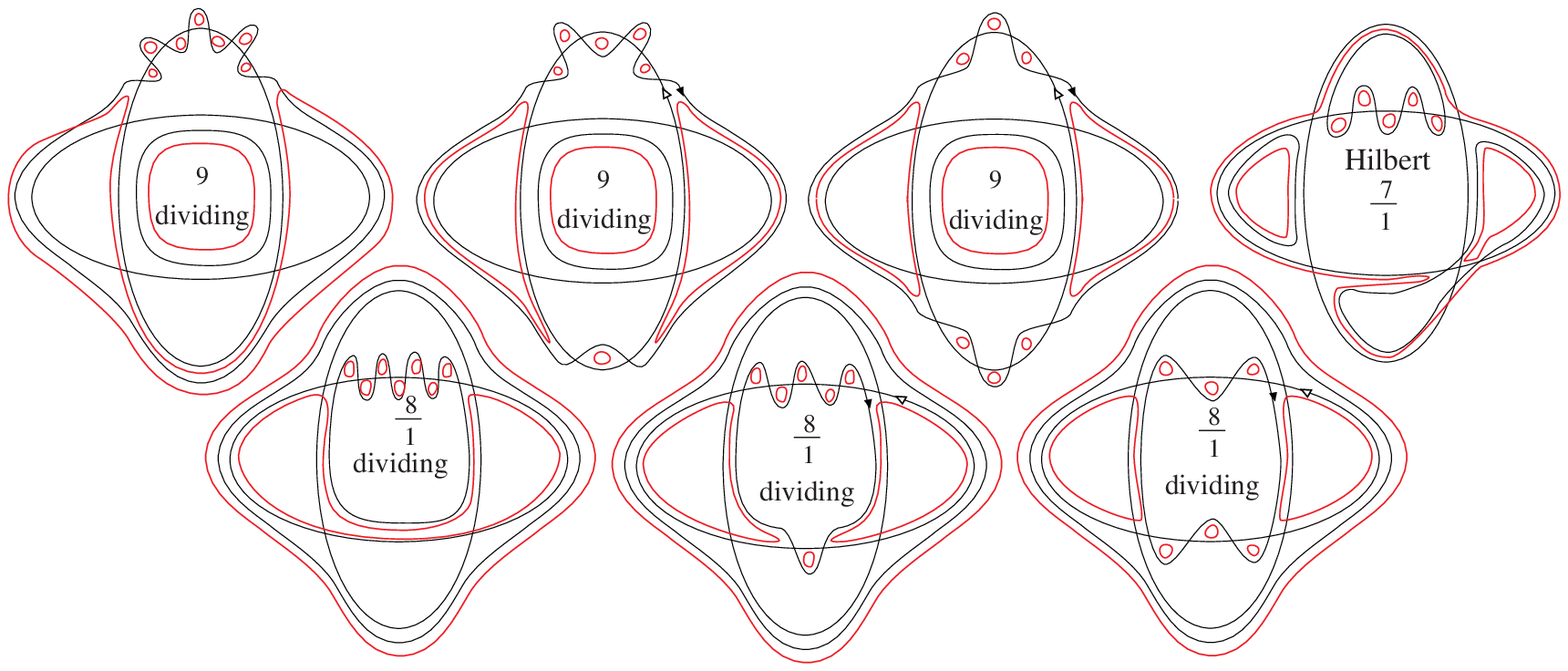,width=122mm}
\vskip-5pt\penalty0
\caption{\label{GudHilb9:fig}%
Hilbert oscillating Klein's G\"urtelkurve} \vskip-5pt\penalty0
\end{figure}

[07.01.13] Of course it is also possible to apply Hilbert's
oscillations to a G\"urtelkurve (or other quartics), cf.
Fig.\,\ref{GudHilb9:fig}. Alas the list of schemes so obtained is
not very exciting (no new species over the previous vibrations).

\subsection{Constructing the indefinite types
(Brusotti, Rohlin 1978, Fiedler 1978, Marin 1979, plus
do-it-yourself)} \label{indefinite-types:sec}

[02.01.13] Recall the definition (Rohlin 1978), a scheme is of
{\it indefinite type\/} if it admits representatives of both
types~I and II, in the sense of Klein 1876, i.e. curves which are
both dividing and not. This section aims to construct all schemes
of indefinite type in degree 6 as to understand in full details
Rohlin's theorem (\ref{Rohlin-type:thm}) enhancing Gudkov's table
by the data of Klein's types. All the strategic information is
tabulated in Fig.\,\ref{Gudkov-Table3:fig}, but each bit of
coloring involves a little fight with the geometrical substratum.
Again the ideas are purely those of Rohlin and his school,
especially Fiedler. In Rohlin's 1978 survey \cite{Rohlin_1978}
detailed constructions are not given.  After completion of this
section, we noted that full details are given in Marin 1979
\cite[p.\,57--58]{Marin_1979}, whose constructions differ slightly
from ours, but  settling one case we failed to detect alone,
namely the type $\frac{8}{1}_{II}$.

$\bullet$ First, consider the scheme $\frac{2}{1}2$.
Smoothing
positively a triad of conics gives the dividing curve on the left
of Fig.\,\ref{R2-12:fig}.

\begin{figure}[h]
\centering
    \epsfig{figure=R2-12.eps,width=122mm}
\vskip-5pt\penalty0
\caption{\label{R2-12:fig}%
Checking indefiniteness of the type $\frac{2}{1}2$ and
$\frac{1}{1}5$} \vskip-5pt\penalty0
\end{figure}

On the other hand, starting from a triangular configuration of
ellipse (center part of Fig.\,\ref{R2-12:fig}) one may by
free-hand drawing (without taking care of orientation) arrange the
real scheme to be the prescribed one. After reporting signs of our
chosen smoothing we find them to be all negatives. At this stage I
thought the curve to be dividing. However, right below one of the
orientation is reversed but the smoothing effected left unchanged.
Now it is of mixed signs, so the curve is in fact nondividing. On
the right-top part of Fig.\,\ref{R2-12:fig} is depicted another
free-hand drawing realizing the given given real scheme. Mixture
of signs implies the nondividing type of this curve.

$\bullet$ On smoothing positively the configuration on the
center-bottom part of Fig.\,\ref{R2-12:fig} we get the (dividing)
curve on the left-bottom of Fig.\,\ref{R2-12:fig} which belongs to
the real scheme $\frac{1}{1}5$. Next starting from the 3 ellipses,
we got the miniature figure on the bottom of Fig.\,\ref{R2-12:fig}
who alas had not the right number of ovals. We thus started anew
for the ``radioactive'' (triangular) triad of ellipses to  find
the right-bottom curve on Fig.\,\ref{R2-12:fig} which has the
correct real scheme and is nondividing (as it involves mixed
signs).

$\bullet$ On smoothing positively  the ``radioactive triad'' of
ellipses for the prescribed orientation gives the dividing curve
on the left of Fig.\,\ref{R4-1:fig}. This belongs to the scheme
$\frac{4}{1}$ (i.e. 4 ovals nested in one big oval and nothing
outside). It is easy to trace the same scheme using as template
the ``atomic triad'' of ellipses, cf. middle-part of
Fig.\,\ref{R4-1:fig}, and checking orientation one finds a mixture
of signs imposing the nondividing character of this curve. Another
option also yielding a nondividing curve is given on the
right-part of Fig.\,\ref{R4-1:fig}

\begin{figure}[h]
\centering
    \epsfig{figure=R4-1.eps,width=112mm}
\vskip-5pt\penalty0
\caption{\label{R4-1:fig}%
Checking indefiniteness of the type $\frac{4}{1}$}
\vskip-5pt\penalty0
\end{figure}

$\bullet$ Consider now a triad of ellipses with two ellipses
invariant under rotation by $90$ degrees, plus one circle pinched
in between. A positive smoothing creates the archipelago sextic on
Fig.\,\ref{R8-1:fig}, which is dividing and of real scheme
$\frac{8}{1}$ (i.e. $8$ ovals captured in a bigger one and nothing
outside). It remains to find a nondividing realization of this
scheme, cf. for this Marin's picture=Fig.\,\ref{GudHilbMarin:fig}
below.

\begin{figure}[h]
\centering
    \epsfig{figure=R8-1.eps,width=102mm}
\vskip-5pt\penalty0
\caption{\label{R8-1:fig}%
Indefiniteness of $\frac{8}{1}$ (but failing); but see Marin's
Fig.\,\ref{GudHilbMarin:fig} below.} \vskip-5pt\penalty0
\end{figure}

\begin{figure}[h]
\centering
    \epsfig{figure=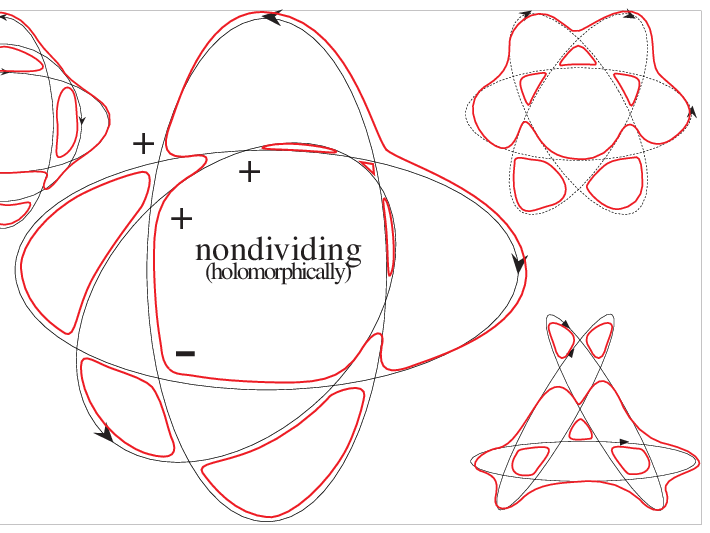,width=122mm}
\vskip-5pt\penalty0
\caption{\label{R3-13:fig}%
Penguin in high pregnancy (of an egg): indefiniteness of the type
$\frac{3}{1}3$} \vskip-5pt\penalty0
\end{figure}

$\bullet$ Dragging down the
archipelago circle gives a configuration of ellipses smoothable
positively to the scheme $\frac{3}{1}3$, cf.
Fig.\,\ref{R3-13:fig}(left). After several infructuous attempts
(depicted as miniatures) one finds the strange triad of conics on
the right-part of Fig.\,\ref{R3-13:fig} which admits a smoothing
belonging to the same real scheme, but which is nondividing due to
mixed signs.

\begin{figure}[h]
\centering
    \epsfig{figure=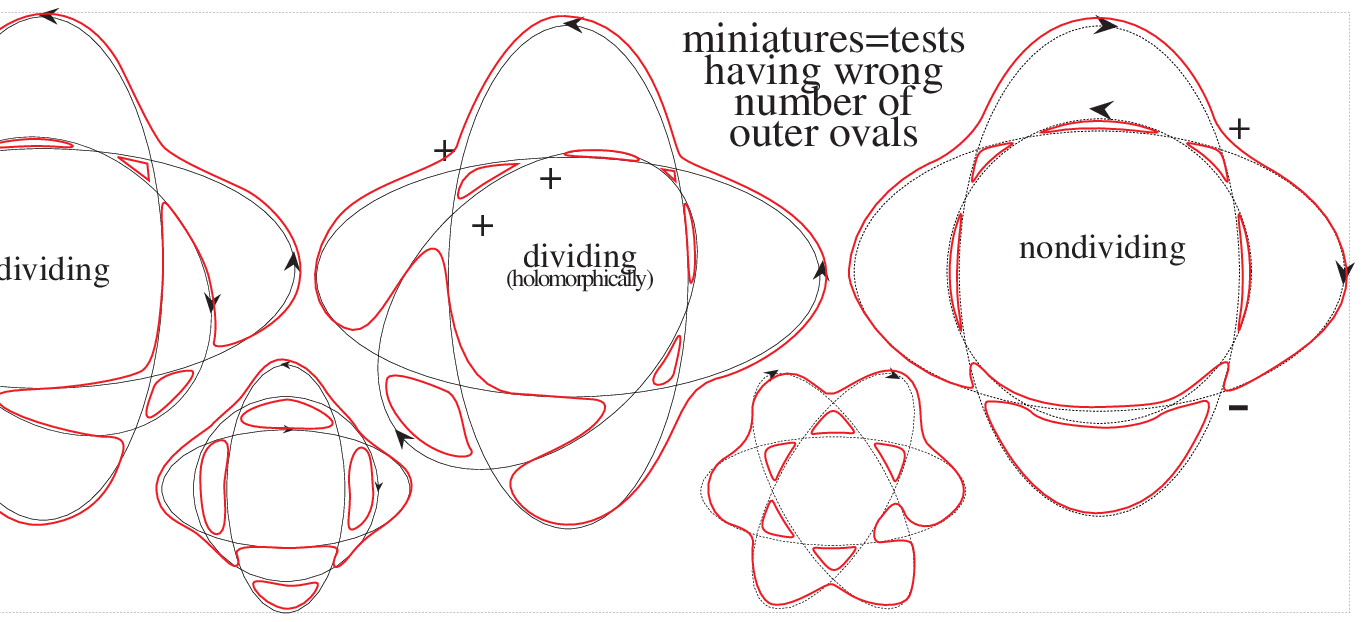,width=122mm}
\vskip-5pt\penalty0
\caption{\label{R5-11:fig}%
Bulldog or Schmetterling=Butterfly: Indefiniteness of the type
$\frac{5}{1}1$} \vskip-5pt\penalty0
\end{figure}

$\bullet$ After some patience and many trials (especially if one
is tired) one finds another configuration of ellipses smoothable
positively to the scheme $\frac{5}{1}1$, cf. left of
Fig.\,\ref{R5-11:fig}. Besides, one finds quickly the right-part
of Fig.\,\ref{R5-11:fig}
belonging to the
same real scheme, yet nondividing due to mixed signs.

{\bf What remains to be constructed?} At this stage we are nearly
finished (compare the list of schemes we explored with those
marked by rhombs on Fig.\,\ref{Gudkov-Table3:fig}). It remains us
to find the scheme $\frac{4}{1}4$. As we did not found it
presently as a perturbation of 3 ellipses, and since this lies
quite near (on Gudkov's table=Fig.\,\ref{Gudkov-Table3:fig}) to
Gudkov's $M$-sextic (notoriously difficult to construct) one is
imbued of some suitable respect. Possibly it is impossible to
exhibit as a deformation of (transverse) 3 ellipses. Notice yet
that the curve $\frac{4}{1}4$ exists as shown by a variant of
Hilbert's method (cf. Fig.\,\ref{GudHilb8:fig}). However presently
this only realizes the scheme in the dividing way, whereas Rohlin
claims this type to be indefinite.

[12.01.13] A somewhat mystical way to solve this question involves
taking a curve lying just above the scheme $\frac{4}{1}4$, while
contracting an empty oval via passage through a solitary node.
(Remember this to be possible by Itenberg 1994
\cite{Itenberg_1994}.) Reading the deformation backward it follows
from Klein's remark (1876)(=Marin's theorem 1988
\cite{Marin_1988}) that the resulting curve has type~II. However
there is probably a more elementary proof by looking at the scheme
$\frac{4}{1}5$, which according to Gudkov (1974
\cite{Gudkov_1974/74}) can be exhibited by a variant of Harnack's
method, while its mirror $\frac{5}{1}4$ is harder to construct
(Gudkov 1954 \cite{Gudkov_1954} even claiming erroneously its
non-existence).

[07.02.13] One elementary way to realize $\frac{4}{1}4$ in type~II
involves a modification of Harnack's method depicted below
(Fig.\,\ref{indef414:fig}). (This is inspired by Gudkov's text,
but alas no picture there). This is the sort of bird hard to
tackle down if one is tired. Beware also that in practicing
Harnack's method one never finds directly what one is seeking (I
found this while searching $\frac{4}{1}5$.)

The trick here is that we leave much room between the  vertical
lines effecting  Harnack's oscillations. So we start with the 3
high vertical lines, and a slight perturbation of the circle union
the horizontal line produces a cubic $C_3$ oscillating thrice
about the horizontal line $L$. The reducible quartic $C_3\cup L$
is then perturbed by a quadruplet of lines, which again is much
stretched so as to effect another intermediate vibration. Then we
have an $M$-quartic $C_4$ oscillating 4 times across $L$. Via the
same trick $C_4\cup L$ is perturbed by a quintuplet of vertical
lines to produce a $C_5$ oscillating 5 times across $L$. Then
using Brusotti's theorem (that German workers used subconsciously
it seems or used ad hoc tricks to complete their perturbations) we
have two ways to smooth $C_5\cup L$ to get a smooth $C_6$. Taking
caring of orientations, the first depicted choice leads a curve of
type~I, whereas the second involves a negative sign and therefore
produces type~II. It is easily checked that both curves belong to
the scheme $\frac{4}{1}4$.

\begin{figure}[h]
\hskip-2.2cm\penalty0 \epsfig{figure=indef414.eps,width=162mm}
\vskip-15pt\penalty0
  \caption{\label{indef414:fig}%
A type~II realization of
$\frac{4}{1}4$ by a variant of
Harnack's method} \vskip-5pt\penalty0
\end{figure}

In a similar way, we do not have yet constructed the type~II
incarnation of the scheme $\frac{8}{1}$. Again in somewhat sloppy
fashion, one could argue by contracting successively two empty
ovals in Hilbert's $M$-curve (scheme $\frac{9}{1}1$), namely the
one outside and one inside the nonempty oval. Granting  such a
deformation through two (successive) solitary nodes, Klein's
remark implies the resulting curve being of type~II, and we are
done. Yet I presume there must be a more elementary construction.

[17.01.13] Indeed one such is sketched in Marin 1979 \cite[p.\,57,
very bottom left of the table]{Marin_1979}. Let us reproduce
Marin's picture as Fig.\,\ref{GudHilbMarin:fig}b:

\begin{figure}[h]
\epsfig{figure=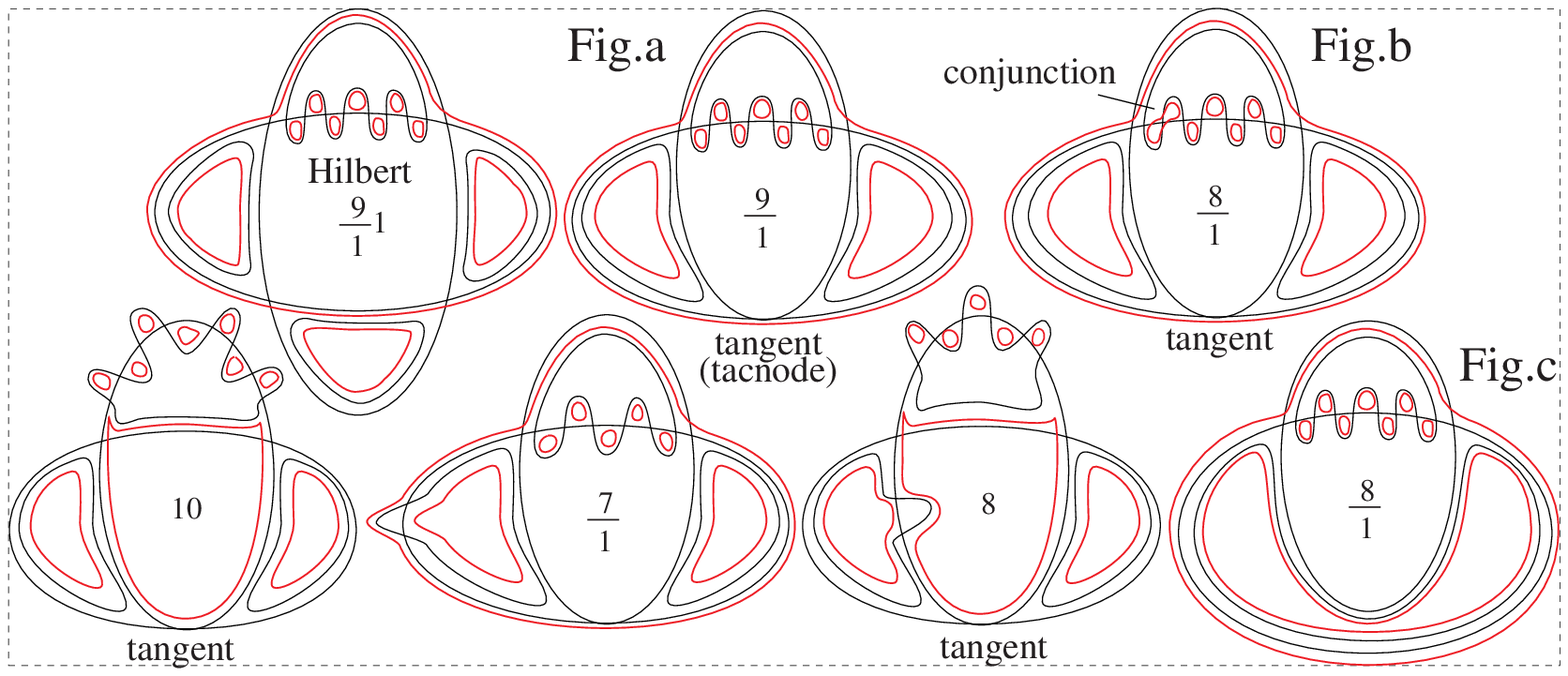,width=122mm} \vskip-5pt\penalty0
  \caption{\label{GudHilbMarin:fig}%
Marin's construction of the scheme $\frac{8}{1}_{II}$ via a
tacnodal configuration of ellipses (as a variant of Hilbert's
method). Compare also Viro's comment in Viro 1989/90
\cite[p.\,1124]{Viro_1989/90-Construction}: ``Degree 6 curves with
the schemes $\langle 10 \rangle$ and $\langle 1\langle 9 \rangle
\rangle$ can easily be constructed by Hilbert's method (the former
can also be constructed by Harnack's method); see
\S\S1.9--1.10.''} \vskip-5pt\penalty0
\end{figure}

Marin's trick here is to start from 2 ellipses  tangent at one
point but transverse elsewhere (Fig.\,a). Perturbing this by a
suitable quadruplet of lines as in Hilbert's method gives a
quartic $C_4$ oscillating as depicted on Fig.\,a and with 3 ovals
only. Hence the $C_4$ is nondividing (Klein's congruence), and so
is a fortiori the resulting sextic $C_6$ (as the nondividing
character is dominant in the genetic sense), which realizes  the
$(M-1)$-scheme $\frac{9}{1}$ (Fig.\,a). A simple conjunction of
two inner ovals yields the $(M-2)$-scheme $\frac{8}{1}$ (Fig.\,b),
we were really interested in (and which is again of type~II for
the same genetical reason).

With this trick we can construct several other curves, depicted on
the second row of Fig.\,\ref{GudHilbMarin:fig}, in particular we
get the scheme $10$ as well as $\frac{8}{1}_{II}$ via a variant of
Marin avoiding tacnodality. The little price to pay is that we
concede two imaginary intersections between the ground ellipses so
that the nondividing character of the $C_4$ (unnested) has to be
derived by some ad hoc argument (e.g. Klein's in
(\ref{Klein-unnested-quartic-nondividing:lem}), or
Arnold's congruence $2=\chi=p-n\equiv k^2 \pmod 4$, or Rohlin's
formula $0=2(\Pi^{+}-\Pi^{-})=r-k^2$ or even B\'ezout modulo the
highbrow contraction conjecture CCC, cf. Sec.\,\ref{CCC:sec}). Of
course there must also be an elementary argument by noticing that
the two imaginary  intersections of both ellipses are
``connecting'' different halves, so that when smoothed as shown
the resulting curve is nondividing. Once this $C_4$ is known to be
nondividing the depicted $C_6$ is likewise by virtue of the
genetical dominance of nondividingness. All this argument looks
tricky but is in reality trivial (think-yourself, and compare
optionally Rohlin 1978, Fiedler 1981 \cite{Fiedler_1981}, Marin
1979 \cite{Marin_1979}, and maybe Gabard 2000 \cite{Gabard_2000}).

{\footnotesize [24.01.13] {\it Intermezzo.} As knowledge  advances
it will perhaps become as difficult to find new truths as to
discover old mistakes. E.g., is Falting's proof of Mordell
correct? Is Freedman's proof \`a la Bing reliable? Is Perelman's
proof of Poincar\'e really eclectic? If not should we retire him
the million. No because because it was never accepted. Finding
mistakes in those venerable implementations will perhaps be as
challenging as claiming new truths? At any rate the game is always
pleasant.

}

\subsection{Gudkov's sextic $\frac{5}{1}5$
(Gudkov 1969, 1973, etc.)} \label{Gudkov:sec}

[24.01.13] Several constructions are available, but first some
historical remarks.

$\bullet$ The very first treatment appears in D.\,A. Gudkov's
Doctor Thesis (1969 \cite{Gudkov_1969-Doctor's-Thesis}) under
Petrovskii and the liberal supervision of Arnold (apparently none
of the supervisors were able to digest the full swing of the
candidate Dmitrii Andreevich). Upon this Polotovskii 1996
\cite{Polotovskii_1996-D-A-Gudkov} comments as follows: ``It is
interesting to remark that the first proof of this fact in
[18](=1969 \cite{Gudkov_1969-Doctor's-Thesis}) was extraordinarily
complicated. It takes up $28$ pages of text, is a ``pure existence
proof'', and was obtained by means of a combination of the
Hilbert-Rohn method with quadratic transformations. Shortly after
D.\,A. Gudkov suggested significantly simpler {\it constructions}
of curves having this scheme, see [19](=1971
\cite{Gudkov_1971-const-new-ser-M-curv}), [21](=1973
\cite{Gudkov_1973-const-curve-deg-6-type-515}), [23](=1974/74
\cite{Gudkov_1974/74}).''

$\bullet$ This complicated proof was published in Gudkov-Utkin
1969/78 \cite{Gudkov-Utkin_1969/78} (English transl. issued in
1978).

$\bullet$ New simpler constructions, are due to Gudkov  and to be
found in Gudkov 1971 \cite{Gudkov_1971-const-new-ser-M-curv}, or
in \cite{Gudkov_1973-const-curve-deg-6-type-515}), reproduced in
his survey Gudkov 1974 \cite{Gudkov_1974/74}.

$\bullet$ This is  also reexposed in A'Campo 1979
\cite{A'Campo_1979}.

$\bullet$ Viro 1989/90 \cite[p.\,1076]{Viro_1989/90-Construction}
also emphasizes Gudkov's initial construction ``was rather
complicated'' (an euphemism as compared to Polotovskii's prose
above).
His second proof reduces ``to the first stage of Brusotti's
construction, i.e., the classical small perturbation of the union
of the curve and the line.'' Yet the whole difficulty is to find a
quintic oscillating 5 times across the line while enveloping 5
ovals in one ``wave oscillation'' while leaving one oval outside
(cf. Viro's figure 12 in \loccit, p.\,1077). According to Viro
(\loccit, p.\,1076): ``It was only in 1971 that Gudkov [11](=1971
\cite{Gudkov_1971-const-new-ser-M-curv}) found an auxiliary curve
of degree 5 that did this.''

Of course since Viro in  the early 1980's,  Gudkov's sextic may
also be exhibited by Viro's patchwork; or as a perturbation of
three ellipses tangent at 2 points like Hawaiian earrings. This
involves yet a deep understanding of how to dissipate such higher
singularities. The interested reader can look at
Fig.\,\ref{Viro3-15:fig}c.

\begin{figure}[h]
\centering
\epsfig{figure=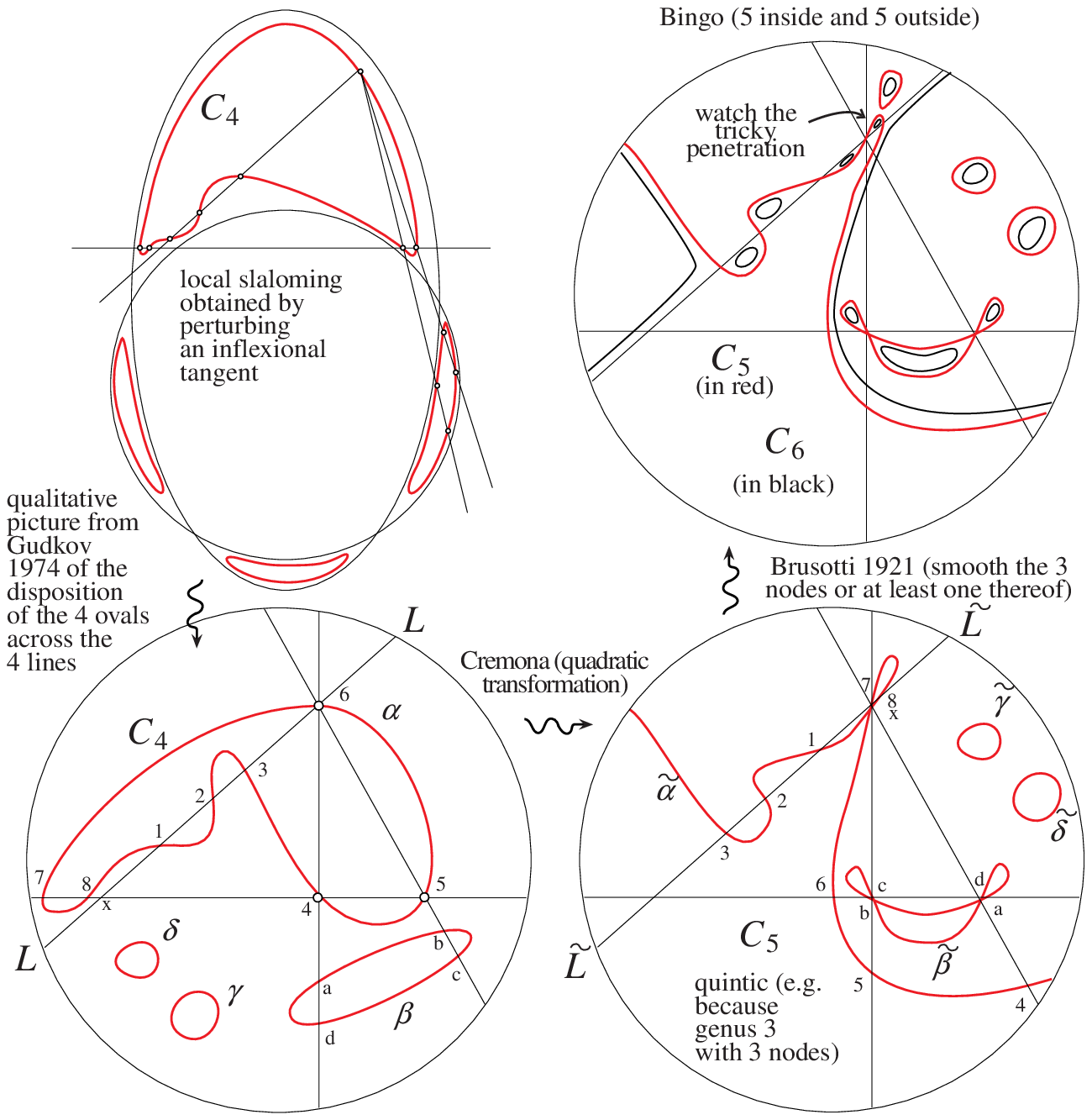,width=122mm} \vskip-5pt\penalty0
  \caption{\label{GudkovCampo-5-15:fig}%
Gudkov's (second) ``elementary'' existence proof of the scheme
$\frac{5}{1}5$: note that nearly all steps look magical (nearly
fallacious)} \vskip-5pt\penalty0
\end{figure}

Now let us describe once more Gudkov's trick (source used Gudkov
1974 \cite[p.\,42--43]{Gudkov_1974/74} and some more d\'etail in
A'Campo 1979 \cite[p.\,12--13]{A'Campo_1979}). This is artwork of
the best stock (cf. Fig.\,\ref{GudkovCampo-5-15:fig}).

$\bullet$ The first step is easy and consists to perturb a line at
one of the 8 flexes of the quartic $C_4$ with $r=4$ slightly so
that it creates 3 nearby intersections. Look at the fourth
intersection, and from here trace two secant intercepting some
other oval as shown, while cutting 2 nearby point on the ``large''
oval such that the line through them cut a little teats on the
large oval. That all this can be achieved is already clever and
explained in detail, first in Gudkov 1971
\cite{Gudkov_1971-const-new-ser-M-curv} or in A'Campo 1979
\cite{A'Campo_1979}.

$\bullet$ The 2nd picture right below is merely a qualitative
redrawing of the first.

$\bullet$ The 3rd picture shows the transformation of the $C_4$
experimented under the Cremona transformation centered at the 3
points $4,5,6$, mutating it into a quintic $C_5$. One way to argue
is via the birational invariance of the genus, keeping the value
$g=3$ constant. Hence as the image curve has 3 nodes (arising as
the intersection of the fundamental triangle through $4,5,6$ with
$C_4$), it must be a quintic. Another way to argue is to remember
the definition of Cremona as the projective (rational) map induced
by the linear system of conics through the 3 basepoints (located
on the large oval). Hence the pullback of a line is a member of
the system, cutting the $C_4$ along 8 points, but 3 of them being
assigned, we find $5$ for the degree of the image of $C_4$.
Likewise the image of the diagonal line $L$ intersecting only the
large oval $\alpha$ is of degree 1, hence a line. To understand
the Cremona-map picture of Gudkov, one must keep in mind that
Cremona contracts any edge of the triangle $4,5,6$ to the opposite
edge of this triangle (and viceversa its explodes each basepoint
to the opposite side of the triangle). The map being actually an
involution (order 2). So the 4 chambers residual to the triangle
are preserved. It is then fairly easy to check that Gudkov's
picture is realist, where tildes are images under Cremona. Life
becomes easier if we number some few points on the $C_4$, while
denoting by the same letters their images under Cremona (omitting
the tilde for simplicity), compare
Fig.\,\ref{GudkovCampo-5-15:fig}. It remains to convince that the
location of $\tilde{\gamma}$, $\tilde{\delta}$ is as depicted by
Gudkov.

The line $L$ is imagined as invariant under Cremona. In fact if we
remove the 3 fundamental lines it remains 4 open triangles
(homeomorphic to a cell ${\Bbb R}^2$) which are preserved. An
involution of the plane has necessarily a fixed point (Brouwer,
Kerekjarto, Smith, etc.) in fact a line or a singleton of
fixed-points depending on whether it reverse or preserve
orientation. The usual formula for Cremona
$$
(x_0, x_1, x_2)\mapsto(x_0 x_1, x_1 x_2, x_0 x_2)
$$
shows that $(1,1,1)$ is fixed, and solving the fixed point
equation  (outside the fundamental triangle whence all $x_i\neq
0$) gives $(1,1,1)=\lambda(x_1, x_2, x_0)$ as unique solution. So
the fundamental triangle splits ${\Bbb R}P^2$ in 4 chambers
preserved under Cremona. How are they permuted? If we normalize
the sign of the first coordinate as positive then we have the
following signs distribution corresponding to the chambers
$$
I=(+,+,+),\; II=(+,+,-),\; III=(+,-,+),\; VI=(+,-,-).
$$
 The first
chamber is preserved by Cremona. The second mutates to the fourth.
The third maps to the second, and the fourth maps to the second.
This looks a bit anomalous for by functoriality we would have
expected that an involution induces an involution on the set of
components (functor $\pi_0$).

Changing the formula to the one written down in Gudkov 1974
\cite[p.\,43]{Gudkov_1974/74} gives
$$
(x_0,x_1,x_2)\mapsto (x_1 x_2, x_0 x_2, x_0 x_1),
$$
and we get
$$
I\mapsto I, II\mapsto(-,-,+)=II, III\mapsto (-,+,-)=III, \textrm{
and } VI\mapsto (+,-,-)=VI,
$$
which is more pleasant. Thus after mutation the ovals $\gamma,
\delta$ stays in the same chamber. Yet this is not enough for if
they would lye like $ \gamma, \delta$ (without tilde) then we
would get the scheme $\frac{3}{1}7$, which is prohibited either by
Hilbert-Rohn-Gudkov or by Rohlin 1972
\cite{Rohlin_1972/72-Proof-of-a-conj-of-Gudkov} proof of Gudkov's
conjecture. (Remind Rohlin's original proof to contain a little
flow, repaired either by Rohlin via Atiyah-Singer or by
Guillou-Marin!). So here we are quite close to adding another
dramatic twist in the Hilbert-Gudkov saga.

However it is more realist that a more thorough examination of the
Cremona map shows the location of $ \tilde\gamma, \tilde\delta$ to
be the one depicted by Gudkov. Indeed the chamber (say $III$)
containing $ \gamma, \delta$ is invariant (like any other).
However the line $L$ is also invariant and divides the chamber
$III$ in two pieces which have to be exchanged by the Cremona
involution. It suffices indeed to use the topological
classification of involutions in the plane ${\Bbb R}^2$ \`a la
Brouwer, to notice that in all cases (orientation reversing or
not) the involution is either a reflection about a line or about a
point (rotation). In both cases the residual components of an
invariant line are exchanged. Hence chamber $III$ splits in two
halves ${III}_{+}$ and ${III}_{-}$, where the former contains
$\gamma, \delta$. Their images have to lye in the other chamber
${III}_{-}$, and Gudkov's depiction is verified. In fact looking
at the image of the point $x$ as mapped to $6=7=8$ shows that
Cremona restricted to $III$ acts as a rotation (having one fixed
points). More algebraically, solving the fixed-point equation
$(x_0,x_1,x_2)=\lambda (x_1x_2, x_0 x_2, x_0x_1)$ shows that
$$
x_0=\lambda x_1 x_2=\lambda^2 x_0 x_2^2
$$
so that $x_2^2=1/\lambda^2$, and likewise---by repeating the
calculation or anticipating it by symmetry---we find
$x_0^2=1/\lambda^2$, and $x_1^2=1/\lambda^2$. Thus up to homothety
we have $(x_0,x_1,x_2)=(1,1,1)$ modulo the 4 possible variations
of signs $(+,+,+)$, $(+,+,-)$, $(+,-,+)$ and $(+,-,-)$. We
conclude that Cremona has exactly 4 fixed points (one in the
barycenter of each chamber). So in particular Cremona is
orientation preserving (within each chamber).

$\bullet$ The fourth picture (of Fig.\,\ref{GudkovCampo-5-15:fig})
contains also a little trick, namely the possibility to smooth the
node (at $7=8$) of the trinodal quintic $C_5$ is such a way that
its pseudoline penetrates slightly inside the line $\widetilde L$.
Once this is done it suffices to smooth \`a la Brusotti $C_5 \cup
\widetilde{ L}$ to obtain the desired Gudkov sextic. And the
miracle is full. Why did it took so long (ca. one century from
Harnack up to Gudkov) to discover this curve? Why Hilbert missed
it? Admittedly the construction is quite tricky, but completely
elementary. Up to our knowledge there is not any further
simplification in this second Gudkov proof, apart perhaps via
Viro's patchwork or dissipation method of higher singularities,
which probably require more highbrow technologies making them
didactically hard to compete with Gudkov's
construction\footnote{[30.03.13] This is not exactly Viro's
opinion, cf. his letter in Sec.\,\ref{e-mail-Viro:sec}.}. As a
last (sentimental) outcome, look how the quintic $C_5$ of
Fig.\,\ref{GudkovCampo-5-15:fig} resemble  a portrait of its happy
discoverer, especially $\tilde \gamma, \tilde \delta$ are like the
eyes, and $\tilde \beta$ the smiling mouth of Gudkov near to crack
the centennial problem.

Finally it is plain from Gudkov's curve to derive curves with less
ovals, especially the $(M-1)$-curve $\frac{5}{1}4$ and the
$(M-2)$-curve $\frac{5}{1}3$ (e.g. by changing the smoothing at
the nodes $2$ and $1$). Those curves were notoriously hard to
construct, and no construction independent of Gudkov's is known.
Using Fiedler's signs-law it is plain that the curve
$\frac{5}{1}3$ so constructed is of type~II, as it should by
virtue of say Arnold's congruence.

If instead we change the smoothing in the inside of the oval along
the smiling mouth $\tilde \beta$ of Gudkov, then we get the
$(M-1)$-scheme $\frac{4}{1}5$, and the $(M-2)$-scheme
$\frac{3}{1}5$. Those were however much easier to construct by a
variant of Harnack's method (as reported in Gudkov 1974); compare
indeed our Fig.\,\ref{HarnaGudkov4-15:fig} and
\ref{HarnaGudkov3-15XXL:fig}.

Finally we note that we may also obtain the $(M-2)$-scheme
$\frac{4}{1}4$ in type~II by smoothing the Gudkov configuration
$C_5\cup \widetilde{L}$. However there  is surely a more
elementary approach via $\frac{4}{1}5$ constructed by a variant of
Harnack; yes indeed compare Fig.\,\ref{indef414:fig}.

\subsection{Diophantine and probabilistic aspects}
\label{Diophantine-and-proba:sec}

[26.01.13] Why did it took so long to discover Gudkov's sextic? Is
it only because it is the most secret part of the pyramid
(Fig.\,\ref{Gudkov-Table3:fig}), or because we have difficulty to
visualize Cremona transformations? Is there some more intrinsic
reason.

One boring algebro-arithmetic game is to think of curves as
ternary forms $F(x_0,x_1,x_2)=\sum_{i,j,k:i+j+k=m} a_{i,j} x_0^i
x_1^j x_2^{k}$ with real coefficients. Up to rounding a bit the
real coefficients randomly we may assume them rational numbers in
${\Bbb Q}$, and this can be done without affecting the topology
nor the rigid-isotopy class. So we find nearby the given curve a
smooth one defined over ${\Bbb Q}$, and we may put all
coefficients in ${\Bbb Z}$ after scaling. As usual we may chase
the common divisor of the equation to get a Diophantine equation
with coefficients primes together $(\gcd (a_{i, j})=1)$. This we
call the reduced equation of the rational curve (in the sense of
Diophante as opposed to having genus $0$). It is unique up to
sign. In particular there is a {\it height\/} defined as the
largest coefficient of the equation.

Then there is a myriad of question. For instance, given an isotopy
type of real curve (or even a rigid-isotopy class) what is the
smallest height of a Diophantine equation realizing this type? To
make this concrete imagine the case of sextics. The Fermat
equation $x_0^m+x_1^m-x_2^m=0$ shows that the corresponding
chamber (unifolium) has always height 1. Similar remark for the
invisible curve $x_0^m+x_1^m+x_2^m=0$ (anti-folium) when $m$ is
even (empty real locus). However it is unknown if the curve with
$r=1$ real branches always corresponds to a unique chamber of the
discriminant (cf. Viro 2008
\cite{Viro_2008-From-the-16th-Hilb-to-tropical}). What is the
height of Gudkov's curve? Can we write down (the best) explicit
equation?

Another question is to look for some fixed  integer $N$ (altitude)
the set of all Diophantine equation $F(x_0,x_1,x_2)\in {\Bbb Z}
[x_0,x_1,x_2]$ of height $\le N$ and consider how they distribute
between the chambers of the discriminant. If $m=6$ there is 64
many chambers by Klein-Rohlin-Kharlamov-Nikulin 1979
\cite{Nikulin_1979/80}) encoded by the chromatic Gudkov table of
Rohlin (Fig.\,\ref{Gudkov-Table3:fig}). Of course some sporadic
equations may land on the discriminant. Now count the
corresponding $65$ (or rather $64$, maybe I added one for the
discriminant but this will tend to zero) frequencies and consider
the corresponding probabilities $p_{i,N}$ (indexed by the Gudkov
symbols $i=\frac{k}{1} \ell$ (plus $(1,1,1)$ deep nest) enhanced
sometimes by the type as on Fig.\,\ref{Gudkov-Table3:fig}). Is the
probability assigned to Gudkov's chamber $\frac{5}1 5$
particularly low, say as compared to Hilbert's or Harnack chamber?
Paraphrasing slightly, how long would it take to a stupid computer
to discover Gudkov's sextic by merely tracing with clever
algorithms the real locus of an explicit Diophantine equation,
while randomly trying one equation after the other.

In contrast one may expect that when $N \to \infty$ there is some
equidistribution, with  all probabilities tending to be equal.
Perhaps some special r\^ole is played by the empty chamber which
is connected by Nikulin 1979 \cite{Nikulin_1979/80}, or better by
the more elementary argument valid in all degrees, cf.
(\ref{empty-chamber-connected-Shustin:lem}). Of course a priori it
is not even clear that the limiting probabilities converge as
$N\to \infty$.

What about the height of Gudkov's chamber, i.e. the least size of
the coefficient of a defining equation. Idem for Harnack and
Hilbert's chambers. Are they lower? Can we estimate the heights
from above using the classical constructions made effective over
${\Bbb Q}$?

Of course all these questions look perhaps a bit unnatural or
somewhat out of reach. Also they depend on the height function
(maximum  coefficient), while there is perhaps other more natural
ways to measure the complexity of an equation, e.g. by the
Pythagorean distance (sum of all spares of the coefficient
$\sum_{i,i} a_{i,j}^2$). This grows like a ball instead of like a
cube, but perhaps the corresponding probabilities are independent
of the exhaustion process? In that case there would be canonical
probabilities and their estimation could be interesting.

All this seems out of reach even when $m=6$, e.g. because we lack
serious algorithms to detect the type from the  equation (compare
e.g. Viro 2008 \cite{Viro_2008-From-the-16th-Hilb-to-tropical}).
Of course the asymptotic probability as $N\to \infty$ of landing
in the discriminant will tend to zero (being a hypersurface of
Lebesgue measure zero). So we should really have a distribution
between $64$ numbers $p_i\in [0,1]$ (some possibly zero? yet
unlikely) weighting the Gudkov-Rohlin pyramid
(Fig.\,\ref{Gudkov-Table3:fig}) by real masses. Are those
probabilities all equal (equidistribution), rational numbers,
etc.? Is the empty chamber much more heavy than the other?

A crude intuition is that when coefficients get larger and larger,
 we get some thermodynamic excitation with all topological
schemes (as complicated as they may be) fairly represented.

Another less arithmetical way to pose the question of the
frequency (e.g. of curves as Gudkov's) is just to put the
natural(?) round elliptic volume element \`a la Riemann-Lebesgue
on the space $\vert mH \vert\approx {\Bbb P}^N$ of all
curves-coefficients dominated by the round (unit) sphere $S^N$.
The latter is calibrated to volume  $2$  as to arrange unity
volume for its quotient ${\Bbb R}P^N$. Each of the 64 chambers
(when $m=6$) has then a (natural) mass, which demands only to be
explicitly determined. It would be again exciting to compare the
mass of Gudkov's $M$-chamber with those of Hilbert's or Harnack's.
Now it is clear that the discriminant has measure zero being a
hypersurface, whereas all other chambers are affected by positive
masses.


How does a random equation (curve) look alike? Letting $p_i$
($i=1,\dots, 64$) be the probabilities assigned to each of the
(Rohlin-Kharlamov-Nikulin) chambers. Those are either all equal
(equidistribution) which is quite unlikely, or some ``curve''
occurs more frequently? From zero-knowledge all what can be said
is that some $p_i\ge 1/64$. What is the largest $p_i$? Maybe the
empty chamber is the most massive?

Of course then there is also refined questions about the
Riemannian geometry of those chambers. Assume for simplicity
equidistribution of masses. Then the whole hotel $\vert mH
\vert-\frak D$ is shared by 64 families having chambers of the
same volume, yet perhaps some are much more comfortable to live
in. Annoying might be chambers highly contorted where there is
little room to plug mobiliary inside. For instance we could look
at the largest Riemannian ball expansible inside a given chamber,
etc.

\subsection{Perturbation of lines (Pl\"ucker 1839, Klein 1873,
Finashin 1996)}\label{Line-perturbation:sec}

[08.04.13] This short section can be skipped. It was written at an
early stage when we had not yet found all schemes asserted by
Gudkov-Rohlin, primarily because we did not mastered sufficiently
the Harnack method. So we attempted to realize schemes by
perturbing lines. In this primitive context it could still be of
interest to understand precisely what schemes are realized. If I
remember well Felice Ronga (ca. 1999) once mentioned this problem
as one challenging his imagination. Perhaps it is worth at the
occasion trying to understand what can be said.

[12.01.13] As yet we missed several schemes whose existence is
asserted in Rohlin 1978 \cite{Rohlin_1978}. This is one motivation
for trying to look at what is obtainable by perturbing an
arrangements of lines. Of course some more ancestral motivation
like the work of Pl\"ucker 1839 \cite{Plücker_1839} as credited
for by Klein 1873 \cite{Klein_1873-Uber-Flächen-dritter-Ordn}
gives also such a motivation. In fact this section was motivated
by a figure in Finashin 1996 \cite[Fig.\,10]{Finashin_1996} which
we shall now reproduce while trying to explore other choices. The
general question could  be which (typed) schemes of the
Gudkov-Rohlin table (Fig.\,\ref{Gudkov-Table3:fig}) can be
realized by perturbing a line arrangement.

\begin{figure}[h]
\centering
\epsfig{figure=Lines.eps,width=122mm} \vskip-5pt\penalty0
  \caption{\label{Lines:fig}%
  Perturbations of 6 lines} \vskip-5pt\penalty0
\end{figure}

\section{Contraction conjectures (Klein 1876, Rohlin 1978,
Shustin 1985, Itenberg 1994, Viro 1994)}

\subsection{``Klein-vache'': Nondividing
implies champagne bubbling? (Klein 1876, disproof Shustin 1985)}


[14.01.13]  As early as 1876 \cite{Klein_1876}, Klein asserted the
firm conviction that curves of type~I cannot gain an oval by
crossing a solitary node. It required ca. 110 years until Marin
1988 took the pain to write down a proof of a somewhat stronger
assertion (cf. Sec.\,\ref{Klein-Marin:sec}).
In the same paper, Klein (1876) speculated about a much more
metaphysical converse
 allowing any nondividing curve to gain an oval after
crossing a solitary node. This was never rigidly asserted by the
cautious Felix Klein, but disproved 99 years later by Shustin 1985
\cite{Shustin_1985/85-ctrexpls-to-a-conj-of-Rohlin}. Personally,
we have not yet assimilated the full details of Shustin argument,
as it uses  much technology, but all experts (Shustin, Viro,
Fiedler, Orevkov, etc.) have validated this disproof.

\begin{conj} {\rm (Klein's  hypothesis of 1876, abridged
``Klein-vache'' in the sequel, disproved in Shustin 1985
\cite{Shustin_1985/85-ctrexpls-to-a-conj-of-Rohlin})}
\label{Klein-1876:conj-noch-entwicklungsfahig} Given any
nondividing plane curve of arbitrary degree $m$, it is possible to
let it cross the discriminant through a solitary node via a path
of curves $(C_t)_{t\in [-1,+1]}$ traversing only once the
discriminant. In other words any diasymmetric chamber bounds a
solitary wall.
\end{conj}

[15.01.13] The conjecture is nearly evident when $m=6$ in view of
Rohlin's enrichment of Gudkov's table by types and the subsequent
rigid-isotopic classification of Nikulin 1979/80
\cite{Nikulin_1979/80} (Theorem~\ref{Nikulin:thm}). With this data
available one gets a
bijection between chambers past the discriminant and Rohlin's
enriched schemes (cf. Fig.\,\ref{Gudkov-contigBIS:fig} below).

\begin{figure}[h]
\centering
\epsfig{figure=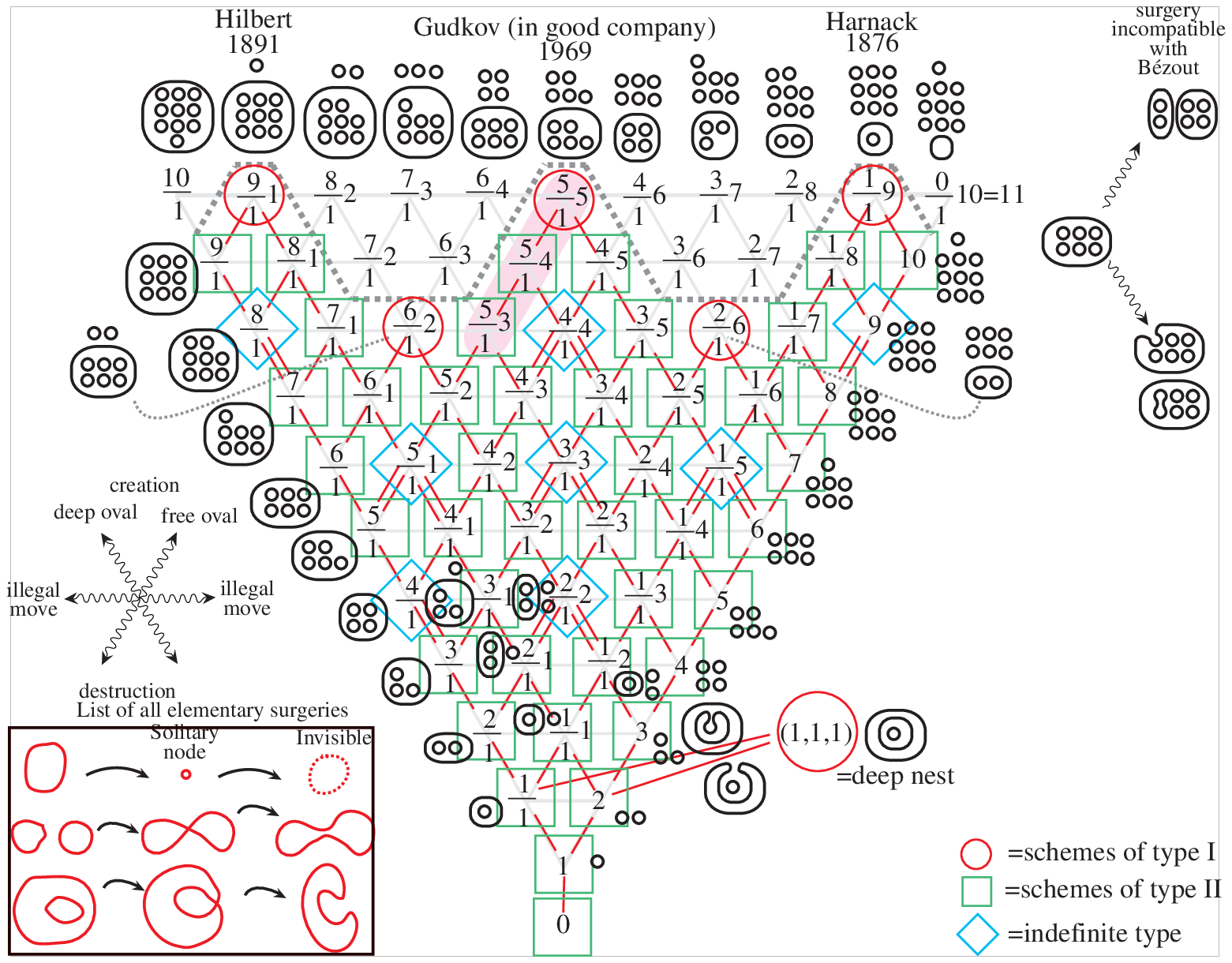,width=122mm}
\vskip-5pt\penalty0
  \caption{\label{Gudkov-contigBIS:fig}%
  Klein-vache in degree 6 via
  Klein-Gudkov-Rohlin-Nikulin-Marin-Itenberg} \vskip-5pt\penalty0
\end{figure}

A moment contemplation of this table shows that all diasymmetric
chambers admit at least one edge moving upwards in the hierarchy
incrementing the number of ovals $r$ by one unit. Of course a
priori such an increment does not necessarily correspond to the
formation of a solitary node (isolated double point) but can also
traduce the subdivision of an oval shrinking to a lemniscate.
Hence more work is required, yet we believe strongly that
``Klein-vache'' holds true for sextics. In fact here is a complete
proof:

\begin{prop}\label{Klein-vache-deg-6:prop}
 {\rm (Gabard 15.01.13, but a trivial
corollary of Rohlin 1978, Nikulin 1979/80, Itenberg 1994
\cite{Itenberg_1994} and Klein-Marin 1876--1988)}.---The
conjecture ``Klein-vache'' {\rm
(\ref{Klein-1876:conj-noch-entwicklungsfahig})\/} holds true for
$m=6$, i.e. any nondividing sextic can acquire a solitary double
point by a rigid-isotopy
crossing only once  the discriminant transversally.
\end{prop}

\begin{proof}
It is first a matter of paying attention to  the combinatorics of
Rohlin's classification into types
(Fig.\,\ref{Gudkov-contigBIS:fig} above). The rest of the proof is
then nearly self-explanatory. In slight contrast to Rohlin 1978
\cite{Rohlin_1978} we forbid the ``hermaphrodite'' {\it indefinite
schemes\/} (allowing projective realizations of both types I/II)
but rather imagine them as two superposed (but distinct) elements,
with the dividing schemes (especially the blue rhombs on
Fig.\,\ref{Gudkov-contigBIS:fig}) levitating slightly above the
sheet of paper of that figure. By Nikulin's theorem 1979/80
\cite{Nikulin_1979/80} those combinatorial symbols (with
levitating twins above the blue-rhombs) are in one-to-one
correspondence with the chambers past the discriminant. Now
imagine on Fig.\,\ref{Gudkov-contigBIS:fig} a sort of random flow
moving downwards along the red-edges of that figure.

Let us be more precise. By a result of Itenberg 1994
\cite[Prop.\,2.1, p.\,196]{Itenberg_1994} (based upon techniques
used by Nikulin (\loccit)) {\it
each empty oval of a sextic can be contracted to a solitary node
before disappearing in the blue sky}. (An oval is said to be {\it
empty\/} if it contains no oval in its interior.)
Pick a curve in each chamber and pick  two contractions (given by
Itenberg) shrinking either an outer oval or an inner oval,
provided both are available on the real scheme. If  only inner or
outer ovals are available,    pick only one contraction. This can
be visualized as a ``random'' vector field moving downward along
the diagrammatic of Fig.\,\ref{Gudkov-contigBIS:fig}.
Each Itenberg contraction necessarily lands in  type~II
(diasymmetric) chambers. Else if landing in an orthosymmetric
(=dividing) chamber, then reading the Morse surgery backwards
corrupts the Klein-Marin theorem (even in its weak original
formulation of Klein 1876 \cite{Klein_1876}, though the latter
gave no proof but see (\ref{Klein-via-Ahlfors(Viro-Gabard):lem}),
or (\ref{Klein-Marin:lem}), or the 1st hand source Marin 1988
\cite{Marin_1988}). Hence our random vector field has its
``trajectories'' ending on the bottom sheet of paper (as we
imagine orthosymmetric chambers levitating somewhat above the
sheet of paper, see again the blue-rhombs on
Fig.\,\ref{Gudkov-contigBIS:fig}). It is plain now that all
diasymmetric chambers (green squares on
Fig.\,\ref{Gudkov-contigBIS:fig}, plus those lying behind the blue
rhombs) do occur as extremities of our vector field encoding the
varied Itenberg contractions chosen. Interpreting this process
backward-in-time proves ``Klein-vache'' in degree $6$. The proof
is complete.
\end{proof}

{\it Insertion} [30.03.13] It should be noted that Itenberg's
contraction theorem affords in degree 6 another proof (independent
of total reality) of Rohlin's maximality principle (in degree 6),
at least if we take for granted the RKM-congruence
(\ref{Kharlamov-Marin-cong:thm}). This prompts another strategy
toward Rohlin's maximality conjecture (independent of total
reality) and perhaps worth exploring further. Of course the hearth
of the problem seems to be the Itenberg-Viro contraction
conjecture for any empty oval
(\ref{Itenberg-Viro-contraction:conj}), but this does not seem to
imply Rohlin's maximality conjecture. In contrast to
``Klein-vache'' the Itenberg-Viro contraction conjecture is still
open and certainly worth investigating further. It is also worth
noting that at the earth of the above proof
(\ref{Klein-vache-deg-6:prop}) we have Itenberg's contraction
theorem. Thus roughly Itenberg implies Klein-vache, yet this is
not the sole ingredient for otherwise in degree 8 Shustin's
disproof of Klein-vache  would refute the contraction principle
(which is still open in degree 8). So the above proof really uses
more than just the contraction principle. In some sense it uses
results by Nikulin but only as a mean to get Itenberg
contractions. What looks more pivotal is the role of the
Gudkov-Rohlin table. One may thus wonder if in degree 8, we can
get sufficient grasp on the Gudkov-Rohlin table as to infer the
logical move from the contraction principle to Klein-vache. If
feasible, then Shustin's disproof (1985) would refute the
Itenberg-Viro contraction conjecture (1994) in degree 8. This
scenario looks a priori quite likely and requires perhaps just
completing the full diagrammatic of Hilbert's 16th  in degree 8,
plus the extra-data of types. (This is perhaps available within
the next decade, if we appreciated correctly the optimism of
experts). Factually, the above proof can be summarized by saying
``Itenberg contraction+Gudkov-Rohlin
diagrammatic$\Rightarrow$Klein-vache'', yet without that it is
crucial to have a bijection between typed-schemes and
rigid-isotopy classes \`a la Nikulin. This correspondence being
disrupted in degree 7 (and so probably 8) by Marin 1979 (cf.
Fig.\,\ref{Marin:fig}). Hence it seems likely that a completion of
the Gudkov-Rohlin table in degree 8, will imply a refutation of
the Itenberg-Viro contraction conjecture.

The above proof of Klein-vache (in degree 6) is quite attractive,
but to be really sublime it should extend to higher orders.
Several obstacles arise. First Itenberg's contraction principle
becomes conjectural for $m>6$ (compare Viro's preface in the same
volume). Next our argument rests on the deep combinatorial
classification of Rohlin 1978 \cite{Rohlin_1978}, plus Nikulin's
rigid-isotopy classification via real schemes enriched by the type
data (I/II). This ceases to be true for orders $m\ge 7$ (Marin
1979/80 \cite{Marin_1979}, Fiedler 1982/83
\cite{Fiedler_1982/83-Pencil}). Thus the above proof looks
jeopardized for higher orders. Of course, if one believes in
Shustin 1985 \cite{Shustin_1985/85-ctrexpls-to-a-conj-of-Rohlin}),
then ``Klein-vache'' is actually false when $m=8$.
Historiographically, it is of course quite improbable that Klein's
(weak) intuition about ``Klein-vache'' was based upon the above
procedure (Torelli for K3's being needed by Nikulin), yet it is
also not completely impossible that a more elementary proof than
the one above exists (cf. optionally
Sec.\,\ref{Klein-vache-proof:sec}). At any rate Klein's power of
prediction via geometric intuition is once more quite amazing.
More modestly, it should be stressed that Klein, interpreted in
the lowbrow fashion, merely asserts that there is no topological
obstacle toward implementing  ``Klein-vache'', yet he is prudent
enough in not claiming this as a theorem (compare again Klein's
original Quote~\ref{Klein_1876-niemals-isolierte:quote} which is
beautifully ambiguous).

[11.01.13] A first natural question is whether
Klein-vache implies the
direct sense of Rohlin's 1978 conjecture (i.e. ``type~I implies
maximal''). In fact Klein-vache shows rather that if a scheme is
not of type~I (so contains a nondividing representative) then it
is non-maximal. Paraphrasing, ``type~I is implied by maximal''.
This is however the part of Rohlin's conjecture that was refuted
by Shustin 1985/85
\cite{Shustin_1985/85-ctrexpls-to-a-conj-of-Rohlin}. So indirectly
Shustin's counterexample also destroys (the hard half of) Klein's
intuition (i.e. Klein-vache).
Shustin's result is somewhat stronger:

\begin{theorem} {\rm (Shustin 1985)} \label{Shustin:thm}
There exists a maximal scheme of degree $8$, which is of type~II.
\end{theorem}

\begin{proof} {\rm  (copied from the source)}
Shustin proves first the following assertion.

\begin{lemma}
There exists $(M-2)$-curves of degree $8$ with the schemes $10 
\sqcup 1 \langle 1 \rangle
\sqcup 1 \langle 2 \rangle
\sqcup 1 \langle 4 \rangle$, and $6 
\sqcup 1 \langle 2 \rangle 
\sqcup 1 \langle 4 \rangle 
\sqcup 1 \langle 5 \rangle$, in the notation of Viro (i.e. the
notation $1 \langle k \rangle$ means one ovals enveloping directly
$k$ empty ovals).
\end{lemma}

\begin{proof} One starts with a certain quintic $C_5$ having
controlled topology with respect to the $3$ axes (constructed in
Polotovskii 1977 \cite{Polotovskii_1977/77}). Then applying a
quadratic transformation gives a singular octic with
``complicated'' singularities. On dissipating such complicated
singularities (Viro's method 1980) one may create the 2 required
schemes.
\end{proof}

{\it End of the proof of Theorem~\ref{Shustin:thm} (compare also
Sec.\,\ref{Shustin-understood:sec} for our slow assimilation of
Shustin's proof)}. Applying a result of Viro 1983
\cite{Viro_1983/84-new-prohibitions}, the $(M-2)$-schemes
constructed above are of type~II. It remains now to check that
they are maximal.

{\it Insertion} [31.03.13].---The Euler-Ragsdale $\chi$ of the
first scheme is $\chi=10+(1-1)+(1-2)+(1-4)=6$, while $k^2=16$.
Hence Arnold's congruence mod 4 (or the allied Rohlin's formula)
suffices to establish type~II of the curve. For the second, $\chi
=6+(1-2)+(1-4)+(1-5)=-2$, and again Arnold/Rohlin suffices to show
type~II.

First Shustin says that the $(M-1)$-schemes obtained from them by
addition of an oval (if they exist) are (always) of type~II,
referring to Rohlin 1978 \cite[point~3.2]{Rohlin_1978}. Needless
to say,  this is actually a trivial consequence of Klein's
congruence (1876) $r\equiv g+1 \pmod 2$. Yet more seriously it
seems to me (Gabard) that we do not need only to know these
schemes being of type~II, but rather that they do not exist at
all!? So in my opinion there may be a trivial misconception here?
In fact we can apply the Gudkov-Krakhnov-Kharlamov congruence
(Theorem~\ref{Gudkov-Krakhnov-Kharlamov-cong:thm}) for
$(M-1)$-curves to all possible enlargements (cf.
Sec.\,\ref{Degree8:sec} for details) yet this fails prohibiting a
specimen. Shustin's argument looks uncomplete at this stage, or
presumably rests on stronger obstructions used subconsciously by
the author!?) ([24.01.13] Compare again
Sec.\,\ref{Shustin-understood:sec} for our assimilation of
Shustin's proof; what is required is a prohibition of Viro.)

Next Shustin argues that the $M$-schemes obtained from the given
ones by the addition of two ovals are forbidden by the extremal
comparison in Rohlin 1978 \cite[point~1.3]{Rohlin_1978}, and Viro
1980 \cite[Theorem~4]{Viro_1980-degree-7-8-and-Ragsdale}.
\end{proof}

{\it Conclusion.}---Beside Polotovskii 1977, Shustin's result
relies massively on Viro's revolutionary technique of construction
via dissipation of complicated singularities (which came to be
known as ``patchworking''). Yet the basic logics of Shustin's
reasoning looks a bit elusive and perhaps flawed. ([24.01.13] Not
all, cf. again Sec.\,\ref{Shustin-understood:sec}.) Hence it is
not clear to me if it really destroys the hard-half of Klein's
intuition (i.e.
Conjecture~\ref{Klein-1876:conj-noch-entwicklungsfahig}).

\smallskip
Let us repeat once more the crucial quote of Klein 1876: {\it
Z.~B. kann bei den Kurven der ersten Art durch allm\"ahl\-iches
\"Andern der Konstanten niemals eine isolierte reelle
Doppeltangente neu enstehen, um dann einen $(C+1)$-ten Kurvenzug
zu liefern; w\"ahrend die Kurven der zweiten Art in dieser
Richtung nicht beschr\"ankt sind. Die Kurven der zweiten Art sind
sozusagen noch entwicklungsf\"ahig, w\"ahrend es die Kurven der
ersten Art nicht sind. Doch soll hier auf diese Verh\"ahltnisse
noch nicht n\"aher eingegangen werden.}

Translated  in English (while adhering to Russian notation and
jargon) gives something like:

{\it For instance, for curves of type~I an isolated solitary node
can never rise
as to produce a new real circuit through
progressive variations of the coefficients; whereas curves of
type~II are not restricted in this way. Curves of type~II are
so-to-speak still developable, while those of the first type are
not.}

\smallskip

This demonstrates that Klein only cautiously  asserted that curves
of type~II are not obstructed to acquire a solitary node, yet not
claiming something so radical as our
Conjecture~\ref{Klein-1876:conj-noch-entwicklungsfahig}, albeit
his second sentence goes closer to suggesting this interpretation.
[24.01.13] At any rate this Ansatz of Klein turns out to be
corrupted by Shustin's article, relying heavily on the new
prohibition detected by Viro (cf. again
Sec.\,\ref{Shustin-understood:sec} for our ultimate assimilation
of this).

\subsection{Degree 8: the Grand pyramid of Gizeh}\label{Degree8:sec}

\begin{figure}[h]
\hskip-3.2cm\penalty0 \epsfig{figure=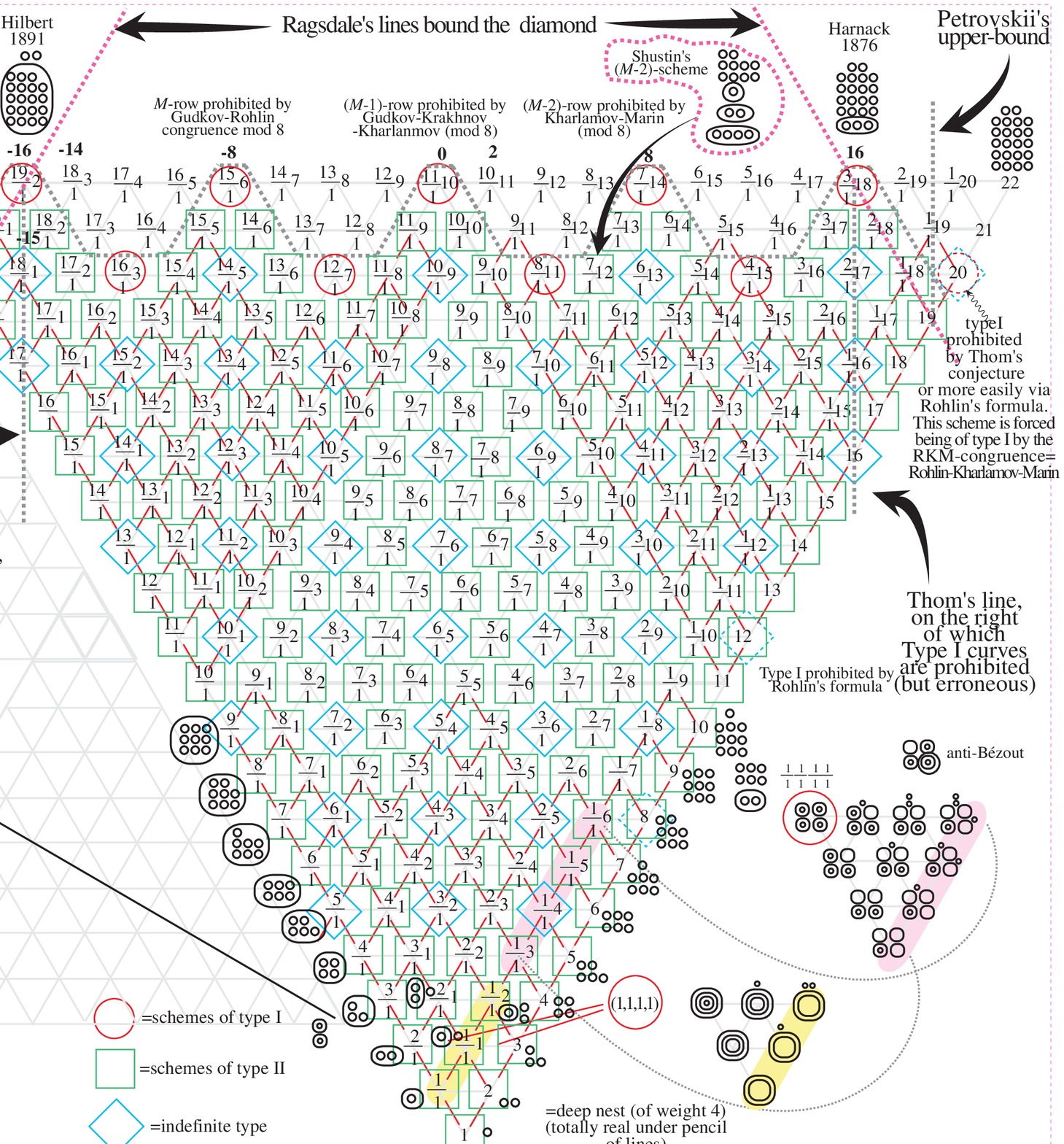,width=182mm}
\vskip-5pt\penalty0
  \caption{\label{Degree8:fig}%
  Schemes of ``all'' octics: a messy pyramid unless one finds a good
  diagrammatic. Warning this is just the basic plate, yet probably
  one needs several other pages (or dimensions)
  to visualize the full pyramid. Puzzle: I do not know if Russian
  workers (Viro, Korchagin, etc.) have all made their own
  road-map catalogue? For
  $\frac{2}{1}17$ via Hilbert's method, cf. Fig.\,\ref{HilbGab3:fig}.} \vskip-5pt\penalty0
\end{figure}

[12.01.13] Can we picture out the Gudkov-Rohlin pyramid in order
8? Since $m=8$ we have $g=\frac{(m-1)(m-2)}{2}=\frac{7\cdot
6}{2}=7\cdot 3=21$. So $M=g+1=22$. It is first quite easy to
extend upwards the Gudkov symbols as to build a larger pyramid
(Fig.\,\ref{Degree8:fig}). Yet this contains only schemes with 1
(or less) nonempty oval. One can easily report the modulo~8
prohibitions coming from Gudkov-Rohlin, etc., as discussed in
Sec.\,\ref{Gudkov-hypothesis:sec}.

{\it Insertion} [02.04.13].---A naive trick is  to think of the
whole pyramid as fibred over the depicted one
(Fig.\,\ref{Degree8:fig}) which shows the range of $(\chi, r)$. So
a basic procedure is to start from the given elementary
configuration with symbol $\frac{x}{1}y$ and to make a menagerie
of transfer of ovals conserving $\chi$. This involves a
Bonsai-cutting art-form of the Hilbert tree. Yet this does not
really solve our puzzle of making a good chart of all possible
schemes in degree 8.

The others schemes having $\ge 2$ nonempty ovals are a bit messy
to report. In particular it seems unrealistic willing to report
all schemes on a single table! Should we try several charts, but
then how to track their interrelations and overlaps? Can we split
in several classes? Let us try to use the number $N$ of nonempty
ovals as a splitting recipe.

$\bullet$ if $N=0$ we have the schemes $\ell$ ($0\le\ell\le 22$),
so 23 schemes.

$\bullet$ if $N=1$ we have the schemes $\frac{k}{1}\ell$.

$\bullet$ if $N=2$ we have $\frac{k}{1}\frac{\ell}{1}m$ (if $m=0$
this is still pictured as the left semi-triangle, yet for larger
$m$'s one must imagine several layers lying above the sheet of
paper). Of course it may be assumed $k\ge \ell$.

$\bullet$ if $N=3$ we have $\frac{k}{1}\frac{\ell}{1}\frac{m}{1}n$

$\bullet$ if $N=4$ we have
$\frac{k}{1}\frac{\ell}{1}\frac{m}{1}\frac{n}{1}o$, but using a
pencil of conics we see that $o=0$, and that $k,\ell, m, n \le 1$
and so we have unique such scheme, namely 4 nests of depth 2.

$\bullet$ Schemes with $N\ge 5$ are prohibited by B\'ezout with
conics.

Okay but all this is a bit overwhelming to depict (except if one
is able to visualize a pyramid in 4D!). Yet we could ask if there
is a reasonable classification of all schemes according to their 3
types (as did Rohlin 1978 for $m=6$). Apart from the obvious
schemes of type~I, and the natural consequences of Arnold-Rohlin,
etc. giving a complete answer looks again  a herculean effort.
Incidentally, it is an open problem as still some few cases are
resisting to the experts of Hilbert's 16th.

Yet we can ask more specific questions like (as did Shustin 1985
\cite{Shustin_1985/85-ctrexpls-to-a-conj-of-Rohlin}) to corrupt
one half of Rohlin's maximality conjecture.

The trick is that under an enlargement of the scheme the number of
nonempty ovals can only increase. So to see what lies above
Shustin $(M-2)$-schemes with $N=3$, it is enough to contemplate
the face of the pyramid with $N=3$, since $N=4$ is nearly empty.

Now writing one of Shustin's scheme in Gudkov's notation gives
$\frac{4}{1}\frac{2}{1}\frac{1}{1}10$.  Note that
$\chi=p-n=(1-4)+(1-2)+(1-1)+(10)=+6=-2 \pmod 8$ so the
Kharlamov-Marin congruence (\ref{Kharlamov-Marin-cong:thm}) says
nothing,
but as observed above the more elementary Arnold congruence forces
type~II.

\smallskip
{\footnotesize

(Elementary B.A.-BA,
hence skip).---To compute the value of $p-n$
(positive minus negative ovals also called even[=pair in French]
and odds) one may use the trick of filling the ovals by an
orientable membrane in ${\Bbb R}P^2$ bounding them in the obvious
way, i.e. we take the interior of all the outer ovals,  then
remove the interior of the subsequent generation of ovals
immediately nested inside, and aggregate again the inside of the
next generation, etc. One has then the psychologically useful
formula $p-n=\chi$, where $\chi$ is the Euler characteristic of
this orientable planar membrane (which M\"obius would call a
reunion of binions, trinions, etc.)

}

\smallskip

\begin{figure}[h]
\centering
\epsfig{figure=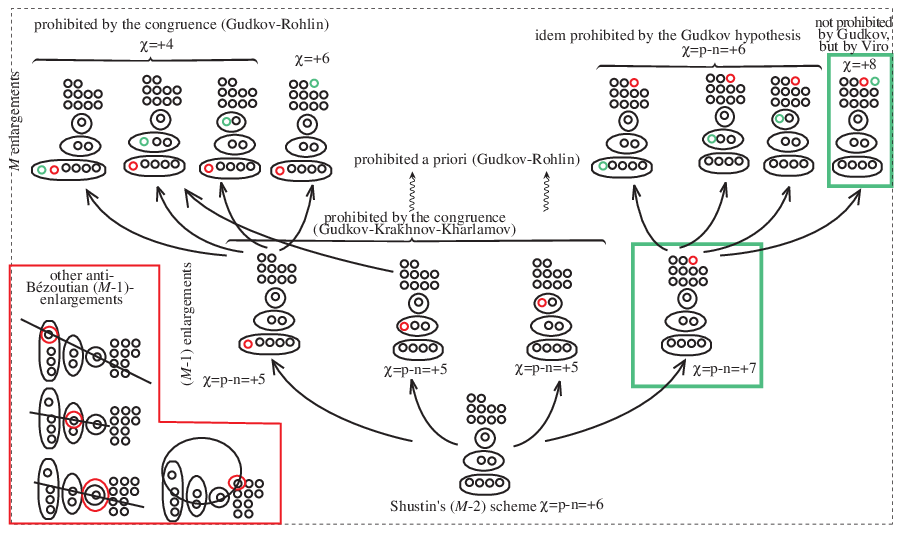,width=122mm} \vskip-5pt\penalty0
  \caption{\label{Shustin:fig}%
  Shustin's scheme
  and its enlargements} \vskip-5pt\penalty0
\end{figure}

Let us now examine the enlargements of Shustin's scheme. First, we
find four $(M-1)$-schemes ruling out those which are not B\'ezout
permissible (cf. Fig.\,\ref{Shustin:fig}). One of them
$\frac{4}{1}\frac{2}{1}\frac{1}{1}11$ (framed on the figure) is
not prohibited by the Gudkov-Krakhnov-Kharlamov-congruence
$\chi=p-n\equiv k^2\pm 1 \pmod 8$
(Theorem~\ref{Gudkov-Krakhnov-Kharlamov-cong:thm}). Whether this
scheme is actually realized is another question. If it is then
Shustin's result (1985
\cite{Shustin_1985/85-ctrexpls-to-a-conj-of-Rohlin}) would be
erroneous. [[24.01.13] No sorry this is a misconception of Gabard,
cf. Sec.\,\ref{Shustin-understood:sec} for a clarification, but
detailed right now for the impatient reader. The point is that if
this $(M-1)$-scheme is realized, then it will be the
counterexample to Rohlin (hence to Klein) for there is nothing
above it by a Viro prohibition (stating that $M$-schemes of degree
8 have an odd content trough out, cf.
(\ref{Viro-Fiedler-prohibition:thm})). Further if it does not
exist then the $(M-2)$-scheme
$\frac{4}{1}\frac{2}{1}\frac{1}{1}10$ will be a counter-example to
Rohlin's conjecture, since it would be maximal but of type~II. So
we get the promised disproof of one-half of Rohlin and of
Klein-vache, without having to know precisely what happens above
Shustin's $(M-2)$-scheme. As a matter of fact it seems, the first
alternative correspond to reality, i.e. the $(M-1)$-scheme
$\frac{4}{1}\frac{2}{1}\frac{1}{1}11$ is realized.]

{\footnotesize

Of course one can play the same game for the other scheme proposed
by Shustin, namely $\frac{5}{1}\frac{4}{1}\frac{2}{1}6$, and we
get the following Fig.\,\ref{Shustin2:fig}. Alas again one larger
$(M-1)$-scheme is not prohibited by Gudkov-Krakhnov-Kharlamov.

\begin{figure}[h]
\centering
\epsfig{figure=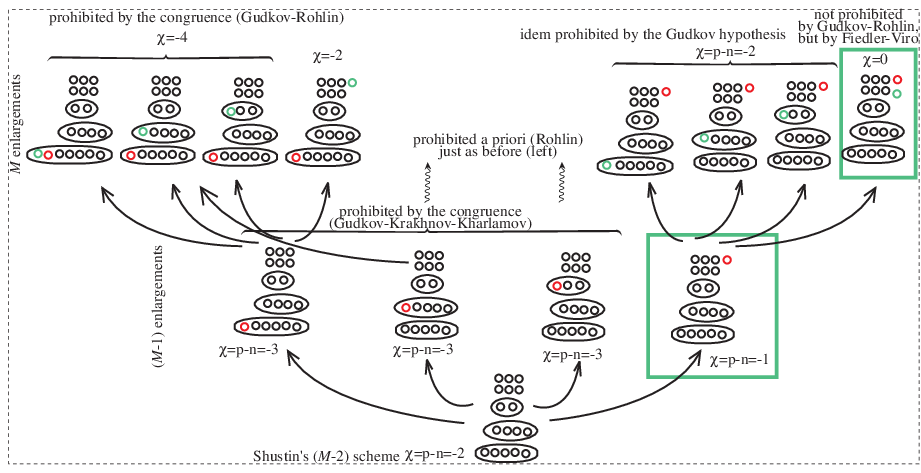,width=122mm} \vskip-5pt\penalty0
  \caption{\label{Shustin2:fig}%
  The other Shustin scheme and its enlargements} \vskip-5pt\penalty0
\end{figure}

Of course it can be the case that the $(M-1)$-scheme lying above
Shustin's scheme (we write in the singular as we fix attention to
one of his scheme) can be prohibited by a stronger prohibition
that Gudkov-Krakhnov-Kharlamov, which Shustin might have used
subconsciously. Yet with what is written down in the article I
could not verify his argument. So perhaps ``Klein-vache'' is still
true (cf.
Conjecture~\ref{Klein-1876:conj-noch-entwicklungsfahig}). In fact
Shustin's paper is the only (published) obstruction I am aware of
against Klein's intuition.

}

{\small

[15.01.13] Idea to explore (but skip as it leads nowhere
[02.04.13]).---Maybe Shustin used subconsciously some inequalities
stronger than the congruences, perhaps those \`a la
Ragsdale-Petrovskii. (Those are discussed in Rohlin 1978
\cite[p.\,87]{Rohlin_1978}.) Let us try the Ragsdale conjecture.
It states
$$
p\le \frac{3}{2} k(k-1)+1.
$$
Rohlin asserts (\loccit) that the first chance of refuting this
inequality is the case $m=10$ (cf. indeed the breakthrough of
Itenberg-Viro of Fig.\,\ref{Itenberg:fig}). So Ragsdale should be
true for $m=8$ (despite having been subsequently disproved by
Itenberg via Viro's patchwork of a variant thereof (again
Fig.\,\ref{Itenberg:fig}). Calculating on the first Shustin's
enlarged $(M-1)$-scheme we find
$p=3+11=14$ versus $\frac{3}{2} k(k-1)+1=\frac{3}{2}
4(4-1)+1=2\cdot 9+1=19$. So Miss Ragsdale is far from violated.

Another idea would be to use a pencil of cubics through some deep
ovals of the enlarged $(M-1)$-schemes. Yet some easy counting
shows that we may force 22 real intersection by taking a connected
cubics through 8 basepoints specified inside the deep ovals, but
this not enough to overwhelm B\'ezout accepting $3\cdot 8=24$
intersections.

}

\subsection{Finally understanding Shustin's argument
(with the help of Orevkov's letter)}\label{Shustin-understood:sec}

[24.01.13] Thanks to a letter of Stepa Orevkov, and a survey of
Viro 1989/90 \cite[p.\,1126]{Viro_1989/90-Construction}, we
learned the presence of new B\'ezout-like prohibitions on
$M$-schemes derived in Viro 1983
\cite{Viro_1983/84-new-prohibitions}. Those have the special
feature of not being of topological origins, but rather
algebro-geometric. (This also  enabled Viro in 1979 to complete
the (soft) isotopic classification of septics solving thereby the
next case of Hilbert's 16th,
compare e.g. Viro 1989
\cite[p.\,1124]{Viro_1989/90-Construction}.)

Here is the relevant Viro's result of which we actually just need
the first clause (transcribed in conservative Gudkov's notation,
having in our opinion a slight advantage of compactness  when it
comes to put the symbols on a pyramid):

\begin{theorem}
{\rm (Viro 1983 \cite{Viro_1983/84-new-prohibitions})}
\label{Viro-Fiedler-prohibition:thm}

$\bullet$ ($M$)---If $
\frac{\alpha}{1}\frac{\beta}{1}\frac{\gamma}{1} \delta$ is the
real scheme of an $M$-curve of degree $8$ with $\alpha, \beta$ and
$\gamma$ nonzero, then $\alpha, \beta$ and $\gamma$ are odd.

$\bullet$ ($M-2$)---If $
\frac{\alpha}{1}\frac{\beta}{1}\frac{\gamma}{1} \delta$ is the
real scheme of an $(M-2)$-curve of degree $8$ with $\alpha, \beta$
and $\gamma$ nonzero and with $\alpha+\beta+\gamma\equiv 0 \pmod
4$, then two of the numbers  $\alpha, \beta, \gamma$ are odd and
one is even.
\end{theorem}

\begin{proof}---[25.01.13] We postpone the proof to a
latter occasion, and merely reproduce now the remark to be found
in Viro 1983 \cite[p.\,416]{Viro_1983/84-new-prohibitions}: ``The
special case of Theorem~2.2.E when $\delta=0$ and $\beta=1$ is due
to Fiedler [11](=Fiedler 1982/83 \cite{Fiedler_1982/83-Pencil}).
Theorem~2.2E was stated as a conjecture by A.\,B. Korchagin in
connection with my results on realization of the real schemes of
$M$-curves of degree $8$. The theorem rules out $40$ real schemes
which are not ruled out by Theorems~2.2.A--2.2.D (of these forty,
four are ruled out by the special case of Theorem~2.2.E which was
proved by Fiedler).''

The proof is completed on p.\,422 of Viro's text. It starts as
follows: Let $C=C_8$ denote a smooth octic with real scheme $\la
\alpha \vc 1\la \beta\ra \vc 1 \la \gamma\ra \vc 1\la \delta \ra
\ra$, where $\beta, \gamma$ and $\delta$ are nonzero.
The crucial result is Theorem~4.2 in Viro, which itself is based
upon Fiedler. So the proof looks too technical to be reproduced
here. A self-contained account encompassing Fiedler and Viro's
article would require several pages, and we postpone this to a
future occasion.

One may wonder if the special case implemented by Fiedler does not
suffice actually to corrupt ``Klein-vache'', i.e. Klein's Ansatz
that nondividing curve can bubble out a new solitary node out of
the blue sky. However Fiedler's result prohibit only the scheme $
\la 1 \la 1\ra \vc 1 \la \alpha \ra \vc 1\la \beta\ra  \ra$ with
even nonzero $\alpha$ and $\beta$ (cf. e.g. Viro 1983
\cite[p.\,420]{Viro_1983/84-new-prohibitions}), and a priori this
is not enough to prohibit the enlargeability of some suitably
chosen $(M-1)$-scheme (compare e.g. the constructions proposed by
Orevkov in the next Sec.\,\ref{Orevkov:sec}).
\end{proof}

The first assertion prohibits the remaining $M$-scheme of
Fig.\,\ref{Shustin:fig}.

As to the second clause of (\ref{Viro-Fiedler-prohibition:thm})
pertaining to $(M-2)$-schemes, I do not know what to do with it.
At this stage I read again Orevkov's letter (cf.
Sec.\,~\ref{e-mail-Viro:sec}), which I have some pain to interpret
properly. Let us reproduce it right below for convenience, while
adding some brackets of mine.

Before completing this reading, I finally understood Shustin's
argument. The point is that whether or not the $(M-1)$-scheme
(framed on Fig.\,\ref{Shustin:fig}) exists do not matter. Indeed
if it does exist (algebraically) then it is of type~II (by Klein's
trivial congruence) and maximal (by Viro's prohibition in the
above theorem), whereas if does not exist then Shustin's
$(M-2)$-scheme is maximal but of type~II, by  construction (or by
Arnold). So in both cases Rohlin's reverse implication
``type~I$\Leftarrow $ maximal'' is foiled.

{\it Insertion} [02.04.13].---Of course it would be  interesting
to know if Shustin's $(M-1)$-scheme enlargement do exist
(algebraically), i.e. the scheme
$\frac{4}{1}\frac{2}{1}\frac{1}{1}11$. If I interpret correctly
the letter below of Orevkov (while removing a little misprint from
it, namely trading the ``11'' for a ``10''), it seems that the
$(M-1)$-scheme written above is realized algebraically.

\subsection{Stepa Orevkov's letter}\label{Orevkov:sec}

We now reproduce Orevkov's letter (brackets=[ ], are our
additions):

$\bullet$$\bullet$$\bullet$ [16.01.13--14h56: Stepa Orevkov]

A small remark:

It is wrong that $11 \cup 1\langle 1\rangle  \cup 1\langle
2\rangle  \cup 1\langle 4\rangle $ is not a part of an
$(M-1)$-scheme. It is. [Not clear how to interpret this? Does it
mean that Shustin's claim is wrong, or simply that this scheme is
an $(M-1)$-scheme. My question  was whether this $(M-1)$-scheme is
realized algebraically, of course. Yet, I admit that my question
was a bit ill posed. In fact I wonder if Orevkov not intended to
write a ``10'' instead of the above eleven.] Moreover, there is no
known example of $(M-2)$-curve of type~II which cannot be obtained
from an $(M-1)$-curve by removing an empty oval. [So Klein looks
still plausible for $(M-2)$-schemes, while Shustin looks wrong. No
sorry, in fact I misunderstood Shustin for a long time, as he does
not claim that the framed $(M-1)$-scheme does not exist.]

In contrary, there are $(M-1)$-curves of degree $8$ (which are
necessarily of type~II) which do not come from any $M$-curve.
These are\footnote{I presume this list is not exhaustive, as
Shustin's scheme above ought to be also realized? If I have well
understood the former ``It is.''.}:

$3\langle 6\rangle $\footnote{This is Viro's notation, and mean 6
ovals enveloped in one, and this thrice. So 21 ovals.}

$4 \cup 1\langle 2\rangle  \cup 2\langle 6\rangle $

$8 \cup 2\langle 2\rangle  \cup 1\langle 6\rangle $

$12 \cup 3\langle 2\rangle $

Construction (inspired by Shustin's construction of $4 \cup
3\langle 5\rangle $ [should locate the reference]):

Consider a tricuspidal quartic $Q_{sing}$ symmetric by a rotation
$R$ by $120$ degree and perturb is[=it] so that each cusp gives an
oval (we assume that this perturbation is very small). Let $Q$ be
the perturbed curve. Two flex points appear on $Q$ near each cusp
of $Q_{sing}$. We chose flex points $p_0, p_1, p_2$ (one flex
point near each cusp) so that $R(p_0)=p_1, R(p_1)=p_2,
R(p_2)=p_0$. We choose homogeneous coordinates $(x_0 : x_1 : x_2)$
so that the line $x_i = 0$ is tangent to $Q$ at $p_i$ $(i =
0,1,2)$.

Let $C$ be the image of $Q$ under the Cremona transformation $(x_0
: x_1 : x_2) \mapsto (x_1x_2 : x_2x_0 : x_0x_1)$. Then $C$ has 3
singular points, each singular point has two irreducible local
branches: a branch with $E6$ and a smooth branch which cuts it
``transversally''. By a perturbation of $C$ we obtain all the four
curves mentioned above.

The fact that these curves cannot be obtained from $M$-curves
immediately follows from the fact that, for any $M$-curve of
degree 8 of the form $b \cup 1\langle a_1\rangle  \cup 1\langle
a_2\rangle  \cup 1\langle a_3\rangle $, all the numbers $a_1$,
$a_2$, $a_3$ are odd\footnote{[24.01.13] The exact reference for
this result is Viro 1983 \cite{Viro_1983/84-new-prohibitions}}.

Best regards Stepa O

This letter helped me much to understand finally Shustin's proof,
and is of course worth studying for its own (especially to make a
picture of it). It gives another counterexample to Rohlin's
maximality conjecture, hence to Klein's Ansatz of champagne
bubbling nondividing curves.

[25.01.13] Now here is an attempt to vizualize Orevkov's example.
As he said we start with a tricuspidal quartic. This is known
since time immemorial (maybe Euler 1745, Steiner 1857, cf. e.g.
Briekorn-Kn\"orrer 1981/86
\cite[p.\,32]{Brieskorn-Knörrer_1981/1986} where it is described
as a hypocycloid, cf. also Lawrence p.\,135, where it is called
the Deltoid). This being given we smooth out the cusps to create
some little ovals. I presume this can be done by hand, otherwise
there is a theorem of Gudkov 1962 \cite{Gudkov_1962} extending to
cusps that of Brusotti 1921. The more difficult task is to
understand what happens under the Cremona transformation. Here I
was much aided by the prototype of Gudkov's example (cf.
Sec.\,\ref{Gudkov:sec}), which is the first place where Cremona
maps were applied to topology of real varieties. Remind that
Orevkov's example, is inspired from Shustin, himself being a
direct student of Gudkov.

So the first steps are fairly easy (say classical for Gudkov's
era), yet it took me some times to trace appropriately the Cremona
transform of the $C_4$. I hope my picture is correct (ask Orevkov
if needed)? It is imagined (I presume) that (like in Gudkov's
construction, cf. again Fig.\,\ref{GudkovCampo-5-15:fig}) the
flecnodal tangent is slightly perturbed to become transverse to
the $C_4$. This implies then the funny behaviour
``forth-back-and-forth'' of the image $C_8$ at the place $1$, say.
So there is an octic as depicted. To trace the  picture it is
useful to keep in mind that the Cremona map takes edges of the
fundamental triangle to the opposite vertices of the triangle,
while preserving the 4 residual component of the triangle.

\begin{figure}[h]
%
\hskip-3.7cm\penalty0 \epsfig{figure=Orevkov2.eps,width=182mm}
\vskip-5pt\penalty0
  \caption{\label{Orevkov2:fig}%
  Attempting to picture Orevkov's variant of Shustin (UNFINISHED)}
\vskip-5pt\penalty0
\end{figure}

Now we arrive at the hard step, namely the dissipation of such
singularities. Here comes the contribution of Viro (if we do not
misunderstood history). As explained say in Viro 1989/90
\cite[p.\,1111--12]{Viro_1989/90-Construction}, there is a myriad
of high order singularity, and one would like to understand their
dissipation. The singularity at hand in our case has 4 branches,
while 3 of them have a second order contact (tangency) at the
singular point. For this specific singularity, he quotes Korchagin
(1988). So we have a table of possible template of dissipation,
which may be locally glued in place of the singularity (exactly as
in Brusotti's method of small perturbation which amounts to the
simplest singularity ``$A_1$'').

Substituting  one of this template, one may hope to find the
schemes announced in Orevkov's letter. Yet, it must be hoped that
I use the right singularity (???, their naming being non-canonical
apparently?), and as yet I failed

Worse if we take $\alpha=5$ and $\beta=1$ and the very first
dissipation of Viro's picture (right of Fig.\,\ref{Orevkov2:fig})
(which is permissible in view of the congruence mod 4), then we
get a curve with $18+3+3=24$ ovals violating frankly Harnack's
bound $M=g+1=(7\cdot 6/2)+1=21+1=22$. I got something wrong!!!

[26.01.13] Another explanation could be that for higher
singularities there is no analog of Brusotti's theorem on the
independence of simplification. The latter is brilliantly
explained in Gudkov 1974 \cite{Gudkov_1974/74}, as reducing
ultimately to Riemann-Roch, but also a theorem of Max Noether, and
even special series. Note that in this passage of Gudkov, he seems
to be not completely up-to-date with the problem of special series
on curves, as was solved by Meis 1960 \cite{Meis_1960} in the
special case of pencils, and by Kempf 1971 \cite{Kempf_1971}, and
Kleiman-Laksov 1972 \cite{Kleiman-Laksov_1972} independently in
the early 1970's.

\subsection{Mistrusting Shustin 1985,
while trying to prove ``Klein-vache'' 1876
(via Garsia-R\"uedy or vanishing cycles \`a la
Poincar\'e-Severi-Lefschetz-Deligne-Mumford,
etc.)}\label{Klein-vache-proof:sec}

[06.03.13] All of our (initial) mistrusting of Shustin's proof is
not really justified anymore, being in part clarified above
(Sec.\,\ref{Shustin-understood:sec}) modulo assimilation of the
Viro-Fiedler advanced  B\'ezout-style prohibition
(\ref{Viro-Fiedler-prohibition:thm}). Hence the sequel has to be
read with suitable discernment, but was not completely censured as
it may contains geometric ideas worth exploring further, and other
issues of independent interest.

[13.01.13] Could it be that Shustin 1985 was wrong, while Klein
1876 is correct!? If so how to prove Klein-vache
(\ref{Klein-1876:conj-noch-entwicklungsfahig})? Of course this
amounts to an amazing topological flexibility of Riemann surfaces
as flying-saucers moving in the Plato cavern of plane projective
geometry, where smooth curves are known to have ``particular''
moduli. More concretely one could imagine that this is always
possible via pure Anschauung, namely the process dual to the
subsequent Fig.\,\ref{Klein-Marin:fig}(left) read in reversed
time. One would take a (globally) invariant cycle (=circle) traced
on the diasymmetric surface which is however acted upon without
fixed point by the symmetry (antipodal map on the circle). Such
circles deserve a name:

\begin{defn}\label{antioval:def}
{\rm An {\it antioval\/} of a symmetric (Riemann) surface is a
topological circle traced on the surface invariant under the
involution and acted upon antipodically by the symmetry.}
\end{defn}


First, note as a trivial topological issue, the following.

\begin{lemma}\label{antioval:lem}
Antiovals only exist on diasymmetric surfaces,
all of them admitting one.
\end{lemma}

\begin{proof}
Assume the surface orthosymmetric (i.e. dividing) and containing
an antioval. By definition an antioval lacks real points,
being acted upon antipodically (by Galois). Take one point of the
antioval and its conjugate (which is distinct) and look at
an arc of the antioval linking $p$ to $p^\sigma$. This arc is in
the imaginary locus, yet connects two conjugate points, violating
the orthosymmetry assumption.

Conversely suppose given a diasymmetric surface $(S, \sigma)$,
hence the quotient $S/\sigma$ is non-orientable. Choose a loop
reversing the indicatrix (local orientation) and avoiding the
boundary of $S/\sigma$. The counter-image of this circle in $S$
gives a circle $C$, since the orientation reversing loop lifts to
an arc via the quotient map which is  a genuine double cover
outside the boundary (alias contour by analysts). Since the
symmetric surface is recovered from the quotient via the
orientation cover (Klein-Weichold yoga), the circle $C$ is the
desired antioval.
\end{proof}

By Klein 1876 (and Riemann),  the number $r$ of ovals (better real
circuits) is bounded by $r\le g+1$ (so-called Harnack bound, under
the supervision of Klein who found a more intrinsic reason). It is
natural asking about a similar bound for antiovals. The antipodal
sphere $S^2$ shows that each great circle is an antioval, whence
an infinity of such. Consider next an antipodal torus of
revolution in 3-space invariant under rotation about the $z$-axis
and acted upon by  central symmetry $(x,y,z)\mapsto (-x,-y,-z)$.
We see 2 evident antiovals by sectioning with the horizontal plane
$z=0$ (Fig.\,\ref{Antioval-torus:fig}a).

\begin{figure}[h]
\centering
\epsfig{figure=Antioval-torus.eps,width=82mm} \vskip-5pt\penalty0
  \caption{\label{Antioval-torus:fig}%
  Cutting the torus by planes through the origin} \vskip-5pt\penalty0
\end{figure}

Varying the slope of this plane  gives an infinitude of antiovals
until we reach the critical tangent plane
(Fig.\,\ref{Antioval-torus:fig}b) after which the 2 ovals are not
nested anymore. The argument certainly generalizes to any genus
upon placing some symmetric pretzels in 3-space. So any
diasymmetric surface without fixed points has infinitely many
antiovals. In fact the above proof (part 2) shows that infinitude
is a general feature without resorting to a Euclidean realization
of Klein's symmetric surfaces.

\begin{proof}[Pseudo-proof of Klein-vache
(\ref{Klein-1876:conj-noch-entwicklungsfahig})]
 After this topological
triviality let us try to attack
Conjecture~\ref{Klein-1876:conj-noch-entwicklungsfahig}, i.e.
``Klein-vache'' (allusion to Lefschetz' vache coined by
Grothendieck? Weil? and used by Deligne, etc.) The idea we try to
exploit (but we are unable to complete the argument) involves
another crazy intuition of Klein
validated by Garsia-R\"uedy building over works by Teichm\"uller,
namely
the fact that any Riemann surface admits a Euclidean realization
in 3-space.

Suppose given a diasymmetric (real) plane curve $C_m$ of
(arbitrary) degree $m$. We may as usual look at the underlying
Riemann surface. According to Klein's intuition (validated by
Garsia and R\"uedy, building over a contribution of Teichm\"uller)
we know that all closed Riemann surfaces admit a conformal model
in Euclidean $3$-space. Let us dream that this adapts as well in
some equivariant form for symmetric surfaces (with an
anticonformal involution).

Note en passant: of course closed non-orientable surfaces do not
embedded in $E^3$ (Euclidean $3$-space), but if bordered they do.
So the only boring diasymmetric surfaces are the invisible ones
(no fixed points) but those luckily enough admit a centrally
symmetric model in $E^3$.

Choose a conformal model of $C_m({\Bbb C})$ in $E^3$ supplied by
Garsia-R\"uedy. By the lemma we know that there is an antioval on
the diasymmetric surface. Shrink the latter to a point via a
(plastical) deformation in $E^3$ akin to
Fig.\,\ref{Klein-Marin:fig}(left) read in reverse sense.

Note at this stage that not all antiovals pinch to a ``connected''
surface (e.g. a pretzel of genus 2 with a belt dividing into two
pieces). So even the topological aspect deserves to be precised,
by looking at ``good'' (i.e. nondividing antiovals). Let us assume
that those always exists.

Next look at our isotopic deformation in $E^3$ to a pinched
pretzel,
generating a one-parameter family of Riemann surfaces.  The
difficulty is to ensure that they stay planar (embeddable in
${\Bbb P}^2$) during the
deformation. This looks a priori
quite implausible, but we have not yet exploited
the full punch of Garsia-R\"uedy. Their result states
that all Riemann surfaces arise in the tubular vicinity of any
classical surface in $E^3$ via a normal deformation of arbitrarily
small amplitude.
Picturesquely, if you have any old woman(=Riemann) surface, but
feel erotically bored by her
due to an acute case of
cellulitis just let vibrate her skin  to get any girl you ever
dreamed about. In the oldest
lady
hides  any beautified young girl with taught epiderm, at least so
conformally! So there is some chance that even if our initial
plastical shrink
deviate outside the realm of plane curves (seen as a stratum $\Pi$
in the moduli space ${\cal M}_g$, where $g=\frac{(m-1)(m-2)}{2}$)
we can still rectify the trajectory so as to stay scotched along
the planarity manifold $\Pi$ (for each time). We get so an
``abstract isotopy'', i.e. a path in $\Pi$ the planarity manifold.

Next we have a canonical map $\vert mH \vert-\disc \to \Pi \subset
{\cal M}_g$ from the space of
smooth curves of order $m$ to the moduli space ($\disc$ being the
discriminant hypersurface). It should be easy to lift our abstract
isotopy to $\vert m H \vert$ while having only the extremity
ending in $\disc$ (necessarily at a solitary node by
construction). Then one continues by letting emerge an oval.
If all this is feasible (taking further better care of the
involution) then  Garsia-R\"uedy implies ``Klein-vache''.
\end{proof}

Perhaps the above strategy requires  to be adapted in $E^4$ to
gain more flexibility and more care about the symmetry. Also if
the Garsia-R\"uedy trick in $E^3$ is  not best suited to the
problem at hand, a more direct approach could be to stay in ${\Bbb
C}P^2$. Recall indeed that Ko 2001 \cite{Ko_2001} has a fairly
general extension of the theorem to any ambient Riemannian
manifold. Alternatively more classic  algebro-geometric methods
(Severi's Anhang F, etc., e.g. as modernized by Harris, etc.) are
perhaps quite likely to imply ``Klein-vache'' if such methods of
degeneration adapt to the reality context (equivariance w.r.t.
Galois which is quite a rule when pure synthetic geometry is
involved). But I must seriously refresh my memory on those works.

\subsection{Another tactic toward ``Klein-vache''
via Itenberg-Viro suggesting a general evanescence principle}

[14.01.13] Is there a relation between ``Klein-vache'' and the
natural Itenberg-Viro conjecture (cf. Itenberg 1994
\cite{Itenberg_1994}, and the preface of that volume by Viro)
positing that:

\begin{conj}\label{Itenberg-Viro-contraction:conj}
{\rm (Itenberg-Viro 1994, abridged CC=contraction
conjecture).}---Any empty oval of a (real, smooth) plane curve can
be contracted to a point (solitary node) via a rigid-isotopy.
\end{conj}

{\small

{\it Historical note} [04.04.13].---In Klein 1892
\cite{Klein_1892_Realitaet} (p.\,176 in the pagination of Ges.
Math.\,Abh. 1922 \cite{Klein-Werke-II_1922}) there is discussed
what he calls the ``Doppelpunktsmethode'' amounting essentially to
contract any symmetry-line of the Riemann surface. This seems to
anticipate the Itenberg-Viro contraction principle. It is not
clear however that Klein ever formulated something as precise as
the above conjecture (specific to plane curves). On p.\,176--177,
Klein's prose extracted from its context sounds a bit
overoptimistic, namely: ``Bei allen anderen F\"allen hat die
Durchf\"uhrung des genannten Prozesses und damit die
Zusammenziehung eines beliebiegen Ovals der Kurve zu einem
isolierten Doppelpunkte keine Schwierigkeit.'' This seems to
trivialize the Itenberg-Viro conjecture but probably does not
because Klein thinks really with abstract Riemann surfaces where
there is much more flexibility than with plane curves. However it
is not impossible that refining Klein's argument/ideas could prove
CC, but it is also quite likely that CC is false.

}

At first sight one may expect a direct logical subsuming of
``Klein-vache'' to ``Itenberg-Viro's contraction conjecture''.
However some moment thought shows that there is no such direct
``rapport de force'', i.e. ``Klein-vache'' is not implied, nor
does it imply, the Itenberg-Viro contraction of empty ovals.
However Prop.~\ref{Klein-vache-deg-6:prop} gives a logical
subordination of Klein-vache to the contraction principle in
presence of additional combinatorial knowledge available in degree
$6$. Via Nikulin's theorem (\ref{Nikulin:thm}) on the
rigid-isotopy classification of sextics it is nearly evident that
Conjecture~\ref{Itenberg-Viro-contraction:conj} holds true for
sextics. This is actually the object of
Itenberg's article just cited.

{\it Insertion} [31.03.13] In view of Viro's isotopy
classification in degree 7, and the philosophy that contraction
plus combinatorial knowledge implies Klein-vache, one can also
wonder if Klein-vache holds true in degree 7. Alas we lack a tool
like K3's in degree 7, and so the situation is somewhat obscure in
degree 7. Possibly, Shustin's disproof of Klein-vache descends
from degree 8 to 7, and then maybe that the contraction principle
is already disrupted in degree 7. Recall, that presently the
contraction principle is wide open in degree 8, yet perhaps
disprovable via Shustin (and a completed classification).

The true relationship between ``Klein-vache'' and Itenberg-Viro
contraction hypothesis (\ref{Itenberg-Viro-contraction:conj})
could be rather an analogy in the principle of proof that one
might naively develop, namely the possibility of shrinking a cycle
invariant under Galois(=complex conjugation). Indeed
``Klein-vache'' amounts essentially to shrink an antioval (cf.
Def.~\ref{antioval:def}), whereas Itenberg-Viro amounts shrinking
an empty oval.

Hence it may be suspected that there is a general
strangulation
principle
specializing to both ``Klein-vache'' and ``Itenberg-Viro''
stipulating the following:

\begin{conj} {\rm (Shrinking principle)}
\label{shrinking-principle-vague:conj} Any (Galois) invariant
cycle(=circle) on a  smooth plane curve of degree $m$ can be
strangulated through a path in the hyperspace of curves
crossing only once the discriminant at a smooth point of the
latter (whenever there is no topological obstruction to do so).
\end{conj}

The parenthetical proviso is required, for one cannot shrink a
nonempty oval without shrinking all its inner ovals, creating
thereby a singularity of higher complexity than nodal.

The proof of (\ref{shrinking-principle-vague:conj}) could be
similar to the eclectic one  sketched for ``Klein-vache'' in the
previous section, i.e. either via Garsia-R\"uedy (hence
Teichm\"uller theoretic) or algebro-geometric via vanishing cycles
\`a la Poincar\'e-Picard-Lefschetz-Severi, etc.

Let us examine the combinatorial possibilities for such a
Galois-cycle. Being (by definition) invariant under complex
conjugation $\sigma$, it can either be:

(1) an oval (pointwise fixed  by the Galois-Klein symmetry
$\sigma$);

(2) an antioval or dia-oval (acted upon
antipodically by $\sigma$);

(3) a pseudo-oval or ortho-oval (acted upon by $\sigma$ with two
fixed points, hence like $(x,y)\mapsto (x,-y)$ on
$S^1=\{x^2+y^2=+1\}$).

Our terminology ortho- and dia-oval  is directly inspired by the
figure of Klein 1892 (reproduced in our
Quote~\ref{quote:Klein-1891/92-ortho/dia}), where given a circle
and a point outside it one considers the involution of $S^1$
exchanging the 2 intersections of each line of the pencil. When
the point lies inside the circle we get a diasymmetry (antipode
like), while if it is outside an orthosymmetry (mirror with 2
fixed points).

Given a real curve (equivalently a symmetric Riemann surface in
the sense of Klein), an oval  exists except in the lowest
diasymmetric case $r=0$ (of Klein's classification). A dia-oval
exists only in the diasymmetric case (Lemma~\ref{antioval:lem}).
An ortho-oval can exist in both the dia- and orthosymmetric cases.
An example of an ortho-oval is traced as the cycle $\beta$ on
Fig.\,\ref{Klein-Marin:fig}.

Specializing the shrinking principle
(\ref{shrinking-principle-vague:conj}) to an oval implies the
Itenberg-Viro contraction hypothesis
(\ref{Itenberg-Viro-contraction:conj}), to a dia-oval implies
``Klein-vache'' (\ref{Klein_1876-niemals-isolierte:quote}).
Finally shrinking an ortho-oval leads to another natural:

\begin{conj} Any two contiguous ovals can coalesce
after crossing
an ordinary node with real tangents.
\end{conj}

{\it Contiguous} means here that both ovals can be joined in
${\Bbb R}P^2$ by an arc having only its extremities on the ovals.
Two contiguous ovals can either be directly nested or unnested yet
unseparated by a larger oval. One should not forget the
possibility of a single oval subdividing himself. The latter
operation is subsumed to no topological obstruction, except that
one might enter in conflict with B\'ezout.

So we may dream of such an unifying principle explaining the
perfect topological flexibility of ``rigid-isotopies'' permitted
to traverse only once the discriminant transversally. In some
sense (to be made precise) our shrinking conjecture asserts that
any Galois-cycle shrinks provided there is no topological
obstruction either in ${\Bbb R}P^2$ nor in the complex locus.

Alas our crude principle does not seem  compatible with the:

\begin{lemma}\label{Finashin-obstruction-to-coalesce-Harnack:lem} {\rm (Admitted, but not understood!,
to whom is it due? Stated in Finashin 1996)} Harnack's (sextic)
scheme $\frac{1}{1}9$ can only degenerate toward the scheme $10$
by contraction of the inner oval, yet not by coalescence of the
two nested ovals.
\end{lemma}

\begin{proof} Cf. e.g. Finashin 1996 \cite[p.\,68, proof of
Thm~6.2]{Finashin_1996}, who alas does not give a  precise
reference for this assertion.
\end{proof}

So here we have a clear-cut example of a Galois-cycle (namely an
ortho-oval) linking the inner oval with the nonempty one of
Harnack's curve (Fig.\,\ref{Finashin:fig}), yet which cannot be
shrunk.

Why is it so? Remember in contrast that the G\"urtelkurve $C_4$
(quartic with 2 nested ovals) can see both its ovals coalesce
(Fig.\,\ref{Finashin:fig}). What is the difference between
Zeuthen-Klein $C_4$ and Harnack's $C_6$? If we take the pain of
tracing the complex orientation (by
Fiedler's algorithm) we get the following pictures. It is seen
that for the G\"urtelkurve $C_4$ the complex orientation (in
red-arrows) disagree from the orientation as the boundary of the
annulus (grey-shaded), while for Harnack's $C_6$ the ${\Bbb
C}$-orientation matches  that as boundary of the ring. Could
positive pairs of ovals be an obstruction to coalescence?

\begin{figure}[h]
\centering
\epsfig{figure=Finashin.eps,width=122mm} \vskip-5pt\penalty0
  \caption{\label{Finashin:fig}%
  Trying to detect an obstruction to coalescing via complex orientations} \vskip-5pt\penalty0
\end{figure}

\begin{figure}[h]
\centering
\epsfig{figure=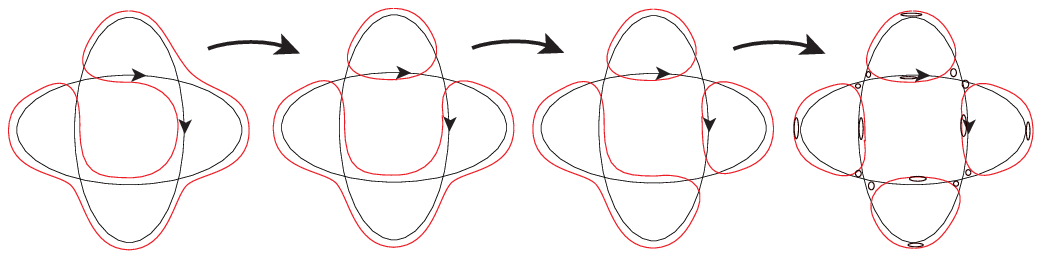,width=122mm} \vskip-5pt\penalty0
  \caption{\label{Subdivide:fig}%
  Subdividing ovals of a quartic} \vskip-5pt\penalty0
\end{figure}

[16.01.13] As we said we may also take an ortho-oval cutting twice
the same oval. Shrinking this would effect a (cellular)
subdivision of the oval. A first example is a hyperbola pinching
to a pair of lines to become another hyperbola. Projectively we
have permanently a conic with a single oval, so there is no naive
minded subdivision like that of a cell in the naive organical
sense. Incidentally 2 ovals for a conic  corrupt either B\'ezout
or Harnack, especially in the formulation of Klein. Likewise the
unique oval of a cubic (if available) cannot be subdivided
(without corrupting either B\'ezout or Harnack-Klein $r\le g+1$).
However an oval of a quartic can sometimes subdivides (cf.
Fig.\,\ref{Subdivide:fig}). (If this figure is realist it is
tempting to create an octic by small perturbation with $16$
unnested oval, yet let us not be sidetracked by this.)

Visualizing the corresponding surgeries on the Riemann surface
must be a pleasant exercise. If the curve is dividing (hence its
Riemann surface orthosymmetric) then  a subdivision is impossible
without corrupting the strong form of Klein's Ansatz proved by
Marin 1988 \cite{Marin_1988}.

\begin{figure}[h]
\centering
\epsfig{figure=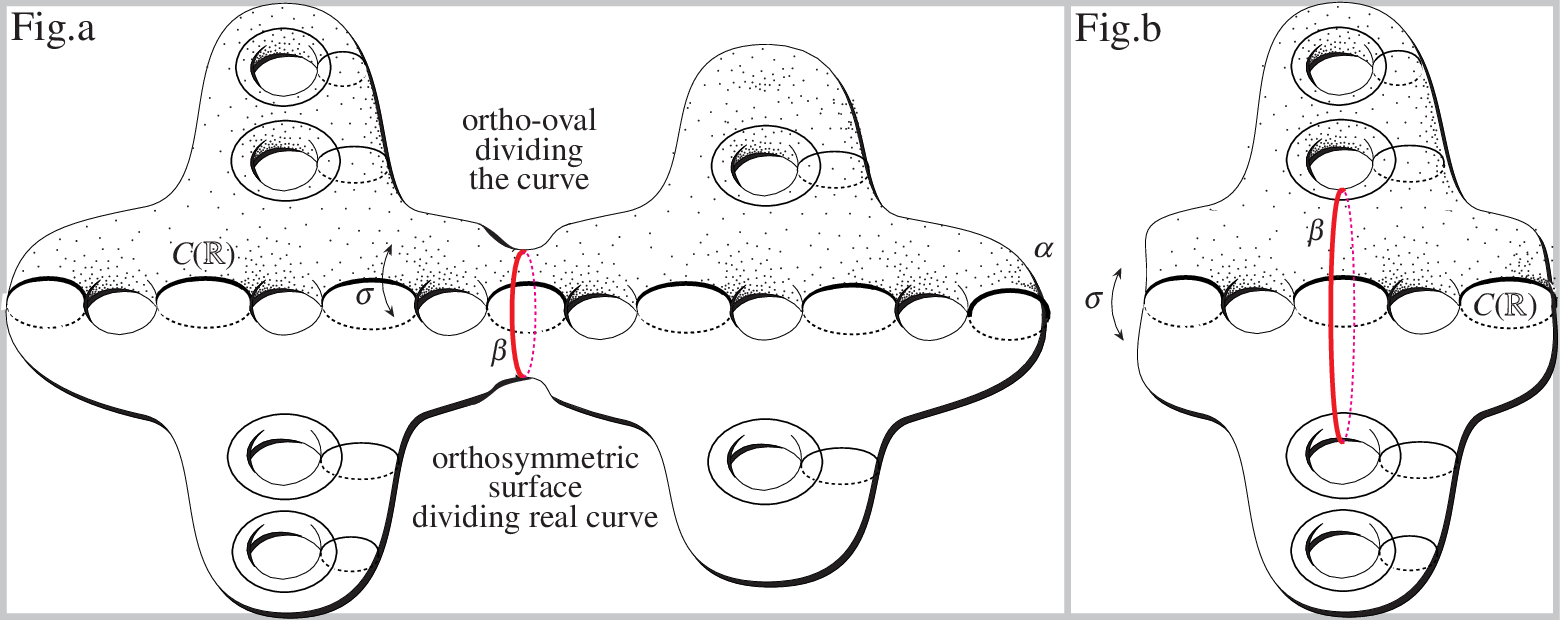,width=122mm} \vskip-5pt\penalty0
  \caption{\label{Orthoovals:fig}%
  Some ortho-ovals cutting only one oval.} \vskip-5pt\penalty0
\end{figure}

It is still tempting to imagine an ortho-oval on a orthosymmetric
curve, especially if it cuts only one oval. A myriad of choices
are possible. Let us depict few of them.
Fig.\,\ref{Orthoovals:fig}a shows an ortho-oval dividing the
surface. Hence when contracted we would get a reducible curve. As
long as our naive picture (be it embedded or abstract) is
respected this is incompatible with B\'ezout (or if you prefer the
intersection theory of ${\Bbb P}^2({\Bbb C})$ whose generator $H$
of the second homology $H_2$ satisfy $H^2=+1$) unless both sides
of the cycle $\beta$ have genus $0$. This proves the following
lemma, whose significance is of course not confined to the case of
real curves.

\begin{lemma} \label{strangulation-impossible:lem}
A dividing cycle on a smooth plane curve $C_m$ of degree $m\ge 3$
cannot be strangulated by a rigid-isotopy crossing only once the
discriminant.
\end{lemma}

\begin{proof}
By contradiction, assume strangulability possible along the given
dividing cycle via a path of curves $(C_t)_{t\in [-1, +1]}$
starting from the given curve, i.e. $C_{-1}=C_m$ and so that only
$C_0$ is singular and uninodal (smooth point of the discriminant).
Denote by $S_t$ the corresponding Riemann surfaces, $S_t=C_t({\Bbb
C})$, where of course $S_0$ is mildly singular. Then the
strangulated surface $S_0$ splits in two (smooth) orientable
surfaces $S_1, S_2$ each porting a fundamental class $\sigma_i$ in
$H_2({\Bbb C}P^2)\approx {\Bbb Z}$ ($i=0,1,2$). Hence we get in
homology $\sigma_0=\sigma_1+\sigma_2$, and so taking respective
degrees $m=m_0=m_1+m_2$, where $m_i=\deg \sigma_i$ (degree in the
homological sense). The intersection $\sigma_1 \cdot \sigma_2$
computes as $m_1 \cdot m_2$, which have to be equal to $1$ (as the
critical curve $C_0$ as just one normal crossing). It follows that
$m_1=m_2=1$, violating the assumption $m\ge 3$.
\end{proof}

The case $m=2$ is entirely different as a conic may degenerate to
a pair of lines.
The interesting  option is to take  a nondividing ortho-cycle
$\beta$, as depicted on
Fig.\,\ref{Orthoovals:fig}b.

Let us now shrink such a nondividing cycle $\beta$
to a point getting something like Fig.\,\ref{Orthoovals2:fig}b,
which is a nodal curve (still irreducible, because its Riemann
surface is connected). After the critical level we could expect to
find Fig.\,\ref{Orthoovals2:fig}d, but this is impossible for the
genus drops by one (remind
that all smooth curves of some fixed degree have the same genus,
since on the complexes the discriminant does not disconnect having
real codimension $2$). In fact as soon as the handle is
strangulated by the vanishing cycle it reappears instantaneously
as depicted on Fig.\,\ref{Orthoovals2:fig}c. On meditating
slightly this occurs like a  twisting, compare the miniature
figures below depicting the pre- and post critical levels near the
singularity. During the twisting one see that the north hemisphere
of the orthosymmetric surface is suddenly connected with the south
hemisphere forcing the diasymmetric nature of the post critical
curve. (All this phenomenology is of course allied to the name of
Picard-Lefschetz and Dehn.)

\begin{figure}[h]
\centering
\epsfig{figure=Orthoovals2.eps,width=122mm} \vskip-5pt\penalty0
  \caption{\label{Orthoovals2:fig}%
  Shrinking a (nondividing) ortho-ovals cutting only one oval.} \vskip-5pt\penalty0
\end{figure}

\subsection{Discovering eversions (Gabard 16.01.13,
but surely in M\"obius, von Staudt, Hilbert, Morosov, Gudkov,
Kharlamov, Finashin, etc.)}\label{Eversion:sec}

[16.01.13]
How can this process (Fig.\,\ref{Orthoovals2:fig}) occur at all if
it is supposed to occur in the plane? In the naive Euclidean plane
${\Bbb R}^2$,  any self-coalescence of a Jordan curve leads to a
subdivision (compare the center of Fig.\,\ref{Klein-Marin:fig}
read backwardly) increasing the number of real circuits. However
during our Riemann surface surgery
(again Fig.\,\ref{Orthoovals2:fig}) the number
of real circuits is kept constant. Hence there seems to be a basic
topological obstruction to our shrinking process, yet some more
mature thinking shows this not to be the case. In reality we live
in the projective plane ${\Bbb R} P^2$, so one oval may well
expand ``to infinity'' to self-coalesce while keeping one
component after having been ``Morse surgered''. For varied
depictions of this phenomenon, see Fig.\,\ref{Eversion:fig} where
as usual ${\Bbb R} P^2$ is depicted as a disc with contour
antipodically identified.

\begin{figure}[h]
\centering
\epsfig{figure=Eversion.eps,width=122mm} \vskip-5pt\penalty0
  \caption{\label{Eversion:fig}%
  Self-coalescing an oval of ${\Bbb R}P^2$ without
  subdividing (increasing the component)} \vskip-5pt\penalty0
\end{figure}

Thus an oval can be Morse surgered without splitting off a new
oval. Let us call this process {\it eversion}\footnote{Our
(non-standard) terminology: V. Kharlamov  explained us (cf.
Sec.\,\ref{e-mail-Viro:sec}) that this phenomenon was quite
crucial a motivation when Morosov suspected some anomaly in
Gudkov's initial solution to Hilbert's 16th problem along the
lines expected by Hilbert. Compare for more Viatcheslav's e-mail
in Sec.\,\ref{e-mail-Viro:sec}, and his terminology ``partner
relationship''. [01.04.13] It would be interesting to know if the
low Gudkov-schemes $\frac{5}{1}4$, and $\frac{5}{1}3$ can be
constructed from their mirrors.} of an oval. Note that after the
eversion all the inside of the oval appears suddenly outside of
it! This basic phenomenon resolves several misconceptions or
paradoxes that foiled for ca. 6 days my understanding of that
theory, especially when it comes to a parity anomaly between the
degree of the discriminant for sextics $3(m-1)^2=75$ and the legal
moves in the Gudkov pyramid (Fig.\,\ref{Gudkov-Table3:fig})
encoding all real schemes combinatorially. If eversions are
overlooked,
the contiguity graph between chambers permits
only closed circuits of even length, whereas by B\'ezout or Galois
a generic pencil of sextics (defined of ${\Bbb R}$) has to cut an
odd number of times the discriminant of degree 75. Nearly a
contradiction in mathematics if eversion would not
exist! Consequently,

\begin{lemma}\label{eversion-deg-6-or-more-forced-by-loops:lem}
Any generic pencil of real curves of even order $m$ contains at
least one eversion. At least this is clear for $m=4,6$, and
hopefully correct in general.
\end{lemma}

\begin{proof} (inserted [01.04.13], but contains a gap!) The discriminant has odd degree $3(m-1)^2$ when $m$
is even.   Our pencil is generic in the sense of being transverse
to the discriminant, hence induces a sequence of Morse surgeries.
Those surgeries (if not ``eversive'') can be visualized on the
Gudkov table (Fig.\,\ref{Gudkov-Table3:fig}) in degree 6 as moves
along the lattice of red-rhombs which permits only closed pathes
of even length, whence the assertion for $m=6$. In general we can
introduce the invariants $(\chi, r)$ and notice that any Morse
surgery which is not an eversion acts as one of the 4
transformations $(\chi, r)\mapsto ( \chi \pm 1,  r \pm 1)$ where
signs can be chosen independently. (Not all of them being possible
as  shown for $m=6$.) Alas this does not seem to be enough to
conclude, because those sole invariants $(\chi, r)$ amounts to a
planar projection of the whole pyramid (which in general is not a
``planar'' object say for $m=8$). So one really needs to
understand the crystallography of higher pyramids which hopefully
still contain merely loops of even length when eversions are
omitted.
Hopefully our lemma is still true (cf. maybe a related argument in
Degtyarev-Kharlamov 2000 \cite{Degtyarev-Kharlamov_2000}).
\end{proof}

Of course an oval belonging to a certain real scheme can be
(topologically) everted iff it is {\it maximal} (i.e. not included
in any larger oval). Let us consider some examples. Suppose the
given scheme to be $\frac{1}{1} 9$, i.e. Harnack's $M$-sextic.
Then there are 10  maximal ovals available. Everting the unique
nonempty oval of Harnack's scheme gives Hilbert's scheme
$\frac{9}{1} 1$ (cf. Fig.\,\ref{Eversion2:fig}), while everting of
of the 9 outer ovals leads to a configuration which is not
B\'ezout-tolerable. The net consequence is that in the hyperspace
of all curves two schemes (better chambers) may be in reality much
closer than they look
far apart on Gudkov's pyramid (Fig.\,\ref{Gudkov-Table3:fig}).
There seems to be some secret passages permitting quick travelling
in the pyramid. Of course whether Harnack's chamber is really
contiguous to that of Hilbert is another story! It would be fine
so for them to sleep
in good company, yet more mature thinking bring us back to the
Riemann surface picture. The ortho-cycle effecting the
strangulation is (since both are $M$-curves) necessarily dividing
hence not strangulable
by Lemma~\ref{strangulation-impossible:lem}. For instance one can
imagine the top orthosymmetric surface with $r=g+1$ as the double
of a planar domain $D$ (with $r$ contours). Suppose given on this
an ortho-cycle $\beta$ meeting twice the same oval. The image of
$\beta$ in one half (our plane domain) is an arc $\beta^{+}$
joining twice the same contour, and some moment thought shows that
$\beta$ divides the surface. Indeed the arc $\beta^+$ can be
completed to a Jordan curve in $D$ by aggregating an arc of the
boundary and we apply Jordan. Better, argue that $\beta^+$ divides
$D$ because we may shrink to a point the contour containing the
extremities of $\beta^{+}$, and then apply Jordan separation. The
separation effected by $\beta^+$ readily implies that by $\beta$.
We have proven:

\begin{lemma}\label{eversion-impossible-for-M-curves:lem}
An $M$-curve (of degree $m\ge 3$) cannot undergo an eversion
(while crossing normally the discriminant). In particular
Harnack's chamber in the hyperspace of sextics is not contiguous
to Hilbert's.
\end{lemma}

\begin{proof} If it could be everted, the corresponding path of curves
would be
materialized by the evanescence (strangulation) of an ortho-cycle
$\beta$ cutting twice the same oval. But $M$-curves correspond to
the top-orthosymmetric case with planar half. Thus by Jordan
separation our cycle $\beta$ divides and therefore cannot be
strangulated (by Lemma~\ref{strangulation-impossible:lem}).
\end{proof}

\begin{figure}[h]
\centering
\epsfig{figure=Eversion2.eps,width=122mm} \vskip-5pt\penalty0
  \caption{\label{Eversion2:fig}%
  From Harnack to Hilbert via a single eversion!} \vskip-5pt\penalty0
\end{figure}

One can imagine more complicated eversions of Harnack's scheme
(cf. bottom of Fig\,\ref{Eversion2:fig}), yet the result is still
the same. Everything depends merely on the oval being everted, for
what was inside becomes outside and conversely. For instance
Gudkov's scheme $\frac{5}{1}5$ turns into itself under eversion of
the nonempty oval, whereas everting a nonempty oval leads to a
scheme enlarging $(1,1,1)$, the deep nest of depth 3, hence
B\'ezout incompatible.

Is the Gudkov chamber self-contiguous to itself via an eversion?
Again this is merely a topological possibility, but it requires a
deeper investigation to see if it is really so. This would imply
Gudkov's chamber to be highly contorted like a banana-shaped, and
it is quite likely that its closure is not simply-connected.
However the lemma above
(\ref{eversion-impossible-for-M-curves:lem}) precludes a
self-contiguity of the Gudkov chamber to itself.

The real option however is that there are two (non-maximal)
chambers past the discriminant related by eversion,
 and actually
we know this phenomenon to exist a priori in view of the degree
argument of
Lemma~\ref{eversion-deg-6-or-more-forced-by-loops:lem}.

\begin{figure}[h]
\centering
\epsfig{figure=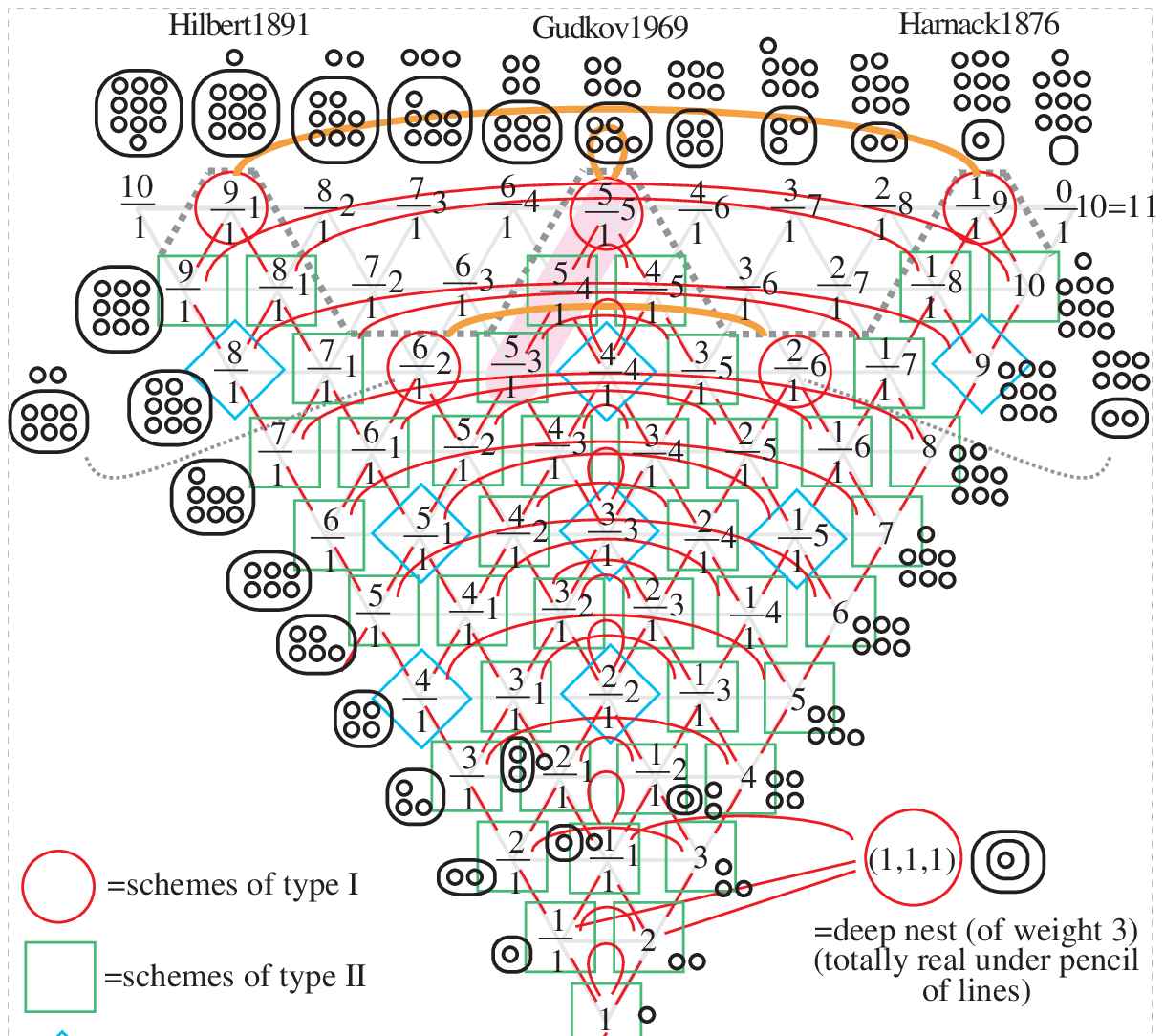,width=122mm}
\vskip-5pt\penalty0
  \caption{\label{Gudkov-eversion:fig}%
  Gudkov's table with eversions
  depicted as red edges or loops} \vskip-5pt\penalty0
\end{figure}

It seems of interest to understand the secrete passages between
Gudkov symbols of Gudkov's pyramid, at least those topologically
permissible under eversion. For sextics we get the following
enhancement of the Gudkov-Rohlin pyramid with curvilinear-edges
amounting to the varied eversion
(Fig.\,\ref{Gudkov-eversion:fig}). Note that we may only evert the
nonempty oval without corrupting B\'ezout (maximality of the deep
nest  $(1,1,1)$). The sole exception arise with the unnested
schemes, plus the scheme $\frac{1}{1} {1}$ whose empty-oval
eversion is precisely the deep nest, whereas the nonempty-oval
eversion flips back the scheme to itself. We get something like
the following messy picture (Fig.\,\ref{Gudkov-eversion:fig})
attempting to keep track of all logically possible eversions.

The question is to decide which among those are effectively
realized algebraically. We already know that those interconnecting
$M$-schemes cannot be realized so due to a topological obstruction
(Lemma~\ref{eversion-impossible-for-M-curves:lem}).

As usual blue-rhombs have to be duplicated according to their
types and really correspond to 2 distinct chambers of the
discriminant.

\begin{lemma}
If a dividing curve (of degree $m\ge 3$) undergoes an eversion
then the post-critical curve is  nondividing.
\end{lemma}

\begin{proof} Fig.\,\ref{Orthoovals2:fig} nearly proves this via the
occurrence of a Dehn-twist in the Picard-Lefschetz transformation.
Again the proviso $m\ge 3$ is evident since Dehn twisting an
equatorial sphere leads to the same equatorial sphere which is
still orthosymmetric.
\end{proof}

Incidentally this gives another proof of the impossibility of
everting $M$-curves
(Lemma~\ref{eversion-impossible-for-M-curves:lem}) since the
latter are necessarily of type~I (by Klein 1876
\cite{Klein_1876}).

{\it Insertion} [01.04.13] Moreover the lemma implies that both
Rohlin's chambers $\frac{6}{1}2$ and $\frac{2}{1}6$ are not
connected by an eversion, since those schemes are of type~I
(either by the RKM-congruence (\ref{Kharlamov-Marin-cong:thm}) or
via Rohlin's (lost) proof of total reality via a pencil of
cubics). Note here that the newly discovered version of Le~Touz\'e
is not strong enough for this purpose (as it uses RKM).

Thus if we imagine the type~I chambers levitating somewhat higher
than the sheet of paper, the eversion starting from a dividing
chamber always moves down to the ground floor of the diagram. Can
we conclude that conversely the diasymmetric type always rises up
to orthosymmetric via eversion?

As already noticed, eversions are impossible for $M$-curves
(except of course if $m=2$, i.e. conics). Thus the 2 top
$M$-curves eversions are actually impracticable. Looking one stage
lower at $(M-1)$-curves we see 3 eversions. Examining the
corresponding Riemann surface, we can imagine something like
Fig.\,\ref{Orthoovals2:fig}, i.e. an orthocycle cutting only one
real circuit while Dehn twisting the handle it is strangulating.
So the real picture is exactly the same as
Fig.\,\ref{Orthoovals2:fig} safe that one real circuit has to be
imagined missing. It
seems plain that the eversion will conserve the diasymmetric
character.
At this stage the problem becomes quite fascinating: for instance
we could via eversion travel from $5_{II}$ to $\frac{4}{1}_{II}$
but not to $\frac{4}{1}_{I}$. Hence $\frac{4}{1}_{I}$ could not be
everted to $5_{II}$. This anti-commutativity looks a bit puzzling,
since any path can be travelled backwardly or forwardly. All this
properly understood could help unravel the mystery allied to the
break of symmetry  prompted by Rohlin's complex orientation
formula, forbidding the scheme $5_{I}$. Despite this and other
intricacies, it seems reasonable to put forward the:

\begin{conj}\label{eversion-and-other surgeries:conj}
All the red-edges of Fig\,.\ref{Gudkov-eversion:fig} (except the
top curvilinear ones linking $M$-curves, and the one linking the
Rohlin's $(M-2)$-schemes $\frac{6}{1}2$ and $\frac{2}{1}6$) can be
realized by an eversion (crossing only once transversally the
discriminant).
\end{conj}

This would give a complete picture of the contiguity graph between
chambers residual to the discriminant via elementary algebraic
Morse surgeries. Nothing forbids that some edges actually
correspond to various Morse surgeries, hence different wall
crossings (e.g. coalescing two inner ovals amounts to coalesce one
inner oval with the nonempty oval).
Hence it is clear that our conjecture is only a crude
approximation, for one might really want to catalogue all walls
between chambers (including self-contiguous one) while describing
the corresponding Morse surgeries. This number of walls is very
likely to be much bigger than the number of edges depicted on
Fig.\,\ref{Gudkov-eversion:fig}.

Modulo little adjustments about the combinatorics, such a
spectacular result is perhaps not completely out of reach, as its
proof could be akin to that implemented in Itenberg 1994
\cite{Itenberg_1994} for contracting empty ovals of sextics. (The
latter implies our conjecture (\ref{eversion-and-other
surgeries:conj}) for all the {\it straight\/} edges.) The rough
philosophy  (at least for sextics) is that when there is no
topological obstruction to shrink, one can strangulate in a
rigid-isotopic way.

However it is quite evident that there must be some extra
obstruction, at least if Finashin's claim
(\ref{Finashin-obstruction-to-coalesce-Harnack:lem}) that the
schematic move $\frac{1}{1}9 \mapsto 10$ (from Harnack's scheme to
the configuration with 10 unnested ovals) can only be accomplished
via contraction of the nested oval but not by its coalescence
with the nonempty oval.
Perhaps Finashin's claim is merely subsumed to a topological
obstruction, which we did not yet understood properly.

\subsection{Strangulation principle (infarctus, etc.)}

[27.03.13] Infarctus=hearth-attack seems to be the generic
mortality cause by geometers too much in love with their topis
(Dirichlet, Gudkov, Rohlin, etc.)

[18.01.13] A message from Viro (16.01.13, cf.
Sec.\,\ref{e-mail-Viro:sec}) suggested
us the following naive remark. Assume an oval (or a priori just a
real circuit) to be contracted to a solitary node via a
rigid-isotopy $C_t \in \vert m H \vert$ (having only one extremity
in the smooth locus of the discriminant parametrized by uninodal
curves). Call such a path a {\it pseudo-isotopy\/} for short. Of
course our contraction supplies a homotopy shrinking the oval to a
point, hence our circuit is null-homotopic in ${\Bbb R}P^2$ and so
forced to be an oval. (By an extension of the theorem of
Schoenflies ca. 1906 a Jordan curve in a surface is null-homotopic
iff it bounds a disc; compare e.g. R. Baer 1928, Epstein 1966, or
Gabard-Gauld 2011 \cite{Gabard-Gauld_2010-Jordan-and-Schoenflies},
etc.)

The presence of this canonical bounding disc gives some evidence
to the Itenberg-Viro contraction conjecture (IVO) of empty ovals,
which supplies some membrane for the strangulation to occur.

{\it Added in proof} [01.04.13] If pessimistic, it may be also the
case that the contraction conjecture is violently false (and
perhaps deducible via Shustin 1985), along the line our
Prop.\,\ref{Klein-vache-deg-6:prop} which shows crudely that
contraction plus classification implies Klein-vache.

Note at this stage that if there were some  dividing plane curve
with only one oval and of degree $m\ge 3$ then strangulating this
oval would be impossible by
Lemma~\ref{strangulation-impossible:lem}. Quite fortunately such
curves do not exist by Rohlin's inequality $r\ge m/2$
(\ref{Rohlin's-inequality:cor}).

So the contraction conjecture for empty ovals (CCEO) suggests
perhaps having some bounding disc for a cycle to be strangulable,
as to refine slightly the statement of the Shrinking principle
(\ref{shrinking-principle-vague:conj}):

\begin{conj} {\rm (Strangulation principle)}
\label{shrinking-principle:conj} Any (Galois) invariant cycle
(=circle) $\beta$ on a  smooth plane curve $C$ of degree $m$ can
be
strangulated (through a path in the hyperspace of curves
crossing only once the discriminant at a smooth point of the
latter), whenever there exist in ${\Bbb C}P^2$ a smooth disc $D$
bounding $\beta$ which is invariant under complex conjugation, and
intersects the
complexification only along $\beta$ (i.e.,  $D\cap C({\Bbb
C})=\beta$). Say in this case that $\beta$ is fillable.
\end{conj}

This is of course a true extension of (CCEO) perhaps susceptible
to imply ``Klein-vache'', i.e. any nondividing curve can acquire a
solitary node by a pseudo-isotopy. Of course any nondividing curve
admits an anti-oval (just lift an orientation reversing loop from
the non-orientable quotient $C(\Bbb C)/\sigma$), but is another
story to find one which is fillable. If so and if the above
conjecture extending Itenberg-Viro's conjecture is right we could
deduce ``Klein-vache'', which is however
disproved in Shustin 1985 \cite{Shustin_1985}.

If there is a filling disc $D$ then as it is invariant under
conj=$\sigma$, we have an involution on the disc (so with a fixed
point by Brouwer). In the case of an anti-oval acted upon by
antipody, it seems that the involution has to act as an antipody
on the whole disc. In general involutions on the disc are of 3
types (identity of order 1, orthosymmetry fixing a diameter, and
antipody fixing the center). (This must  be ex/implicit in work by
Brouwer, Ker\'ekj\'art\'o ca. 1914-1922.)

As argued by Viro's e-mail, CCEO looks more natural that
``Klein-vache'' since the filling is virtually God-given, just
taking the (sealed) inside of the oval. Yet perhaps this has some
analog for an anti-oval in term of differential-geometric
fillings, e.g. minimal surfaces in ${\Bbb C}P^2$ endowed with its
``round'' Fubini-Study metric coming from $S^5\to {\Bbb C}P^2$
(Hopf fibering).

One problem with anti-ovals is that there are plenty of them (not
just finitely many like for ovals), and so one's idea could be to
select some preferred one, maybe as ``the'' systole. In the case
of the diasymmetric sphere this is not enough to ensure
finiteness, but perhaps suffices to single out some natural class
of anti-ovals. Recall that systoles are geodesics, and so are
usual ovals.

A recipe could be as follows: given a nondividing curve $C_m$ of
some degree. Endow it with the natural Fubini-Study metric of
${\Bbb C}P^2$ to get a Riemannian metric on the Riemann surface
$C_m(\Bbb C)$. Since the curve is diasymmetric it contains an
antioval (invariant circle acted upon by antipody by $\sigma$).
Hence by compactness there is also  a such of minimal length, the
so-called systole, not perfectly unique of course, but choose one
such systolic antioval. Consider the latter as a circuit in the
ambient ${\Bbb C}P^2$ and  solve the Plateau problem for that
contour, in its classical setting of soap films diffeomorphic to
the disc. Plateau makes also sense for membranes of higher
topological structure, but ignore them to stay closest to the
Itenberg-Viro conjecture. Plateau is always soluble but the
notorious difficulty is to ensure embeddedness of the solution.
Perhaps this is true in $E^3$ and also in ${\Bbb C}P^2$ due to
some simple-connectivity, or perhaps the special systolic
properties of the boundary data. ([21.01.13] Beware that a minimal
surface has vanishing mean curvature, while the natural
Itenberg-Viro ``reality'' membrane is positively curved. But the
former assertion is specific to $E^3$\dots)

\def\CP2{${\Bbb C} P^2$}

As the given contour is invariant under $\sigma$ (an isometry of
${\Bbb C} P^2$) it is likely that Plateau's solution enjoys a
similar invariance, and we would be essentially finished (modulo
the difficulties enumerated).

At this stage we would have a perfect analog of the bounding disc
of Itenberg-Viro's empty oval, via our Plateau filling of ``the''
systole realizing the anti-oval of shortest length. For the
analogy to be perfect one should ensure that the Plateau film
intersects the Riemann surface $C_m({\Bbb C})$  only along the
contour (systolic anti-oval). This looks either hard or trivial.
For instance recall (from Wirtinger, cf. also Mumford's book
\cite{Mumford_197X-BOOK-complex-alg-var}) that algebraic
subvarieties of ${\Bbb C}P^n$ endowed with Fubini-Study are
(precisely?) minimal surfaces. So there is perhaps some chance to
prove disjointness. (If they intersect interiorly then try to
build a canal surface by surgering a piece of $C_m({\Bbb C})$ to
the Plateau film, trying so to violate  its area minimization
\dots.)

Maybe some interesting twist of Plateau's problem is
that one may be able to reconstruct the whole complex locus via
Plateau if we are given only the real locus of the curve. Of
course as there is now handles (except for $M$-curves) and several
contours this will necessarily involve the so-called
Plateau-Douglas problem permitting membranes of higher topological
structure (than the disc). As hazardous as it is, this claim would
perhaps only work for dividing curves.

Assume all this to work then we have a perfect analogy with
Itenberg-Viro, but it is still not explained why the empty oval or
our anti-ovals are strangulable via a pseudo-isotopy of algebraic
character.

Naively from the given data consisting of minimal film bounding
the cycle $\beta$ we can hope to shrink concentrically the disc
(put via the Riemann mapping in conformal equivalence with the
round disc) to its center. Solving the higher Plateau-Douglas
problem for this shrinking contour gives a minimal surface (which
by the converse of Wirtinger) would be an algebraic curve
realizing the deformation we are looking for. Since this
concentric shrinking respect the symmetry (of the round disc
whatever its type, i.e. identity, antipodal diasymmetry  or
orthosymmetric mirror like $z\mapsto \bar z$), the given smaller
contours are invariant under $\sigma$ in \CP2 and so we get real
curves by solving Plateau-Douglas, and therefore the desired
pseudo-isotopy from $C_m$ to a nodal curve (with a solitary
point).

All this is a bit reminiscent of Riemann's spirit (except for
being of lesser vintage) yet dreaming like a canary   is quite
pleasant. The above strategy (with all its gaps) suggests even
that in the Itenberg-Viro shrinking of an empty oval it could be
arranged that the subsequent curves all have their ovals
progressively shrinking inside the initial one. Whether this
stronger conjecture has some chance to hold true is unknown to me.

The above argument suggests that a real algebraic curve should be
reconstructible from a single oval. This is certainly true via
something like the Nullstellensatz, yet the assertion that all
this algebra can be supplanted by differential geometry \`a la
Plateau looks a bit doubtful for we are not controlling the full
contours of the membrane. Actually the dividing case looks
psychologically more comfortable to a direct appeal of Plateau.

As ``Klein-vache'' is probably false (cf. Shustin 1985
\cite{Shustin_1985/85-ctrexpls-to-a-conj-of-Rohlin} and our
partial discussion in Sec.\,\ref{Shustin-understood:sec}), we
shall from now on concentrate on the Itenberg-Viro contraction
conjecture for empty ovals, which is still non refuted (hence more
likely to be true) and easier technically due to the canonicalness
of the film bounding the vanishing cycle.

So suppose given any empty oval of a smooth plane real algebraic
curve. Goal: strangulate it algebraically via a pseudo-isotopy,
i.e. a rigid-isotopy except for its extremity which is a solitary
nodal curve. To the disc bounding the (marked) empty oval, apply
the Riemann mapping theorem  to take it conformally to the unit
disc $\Delta:=\{z\in {\Bbb C}: \vert z\vert \le 1\}$ via a map
$f\colon D\to\Delta$. This mapping is canonical up to conformal
automorphism of the disc, hence unique once a center and a
boundary point are chosen (variant choose $3$ boundary points).
Consider the pullback $C_\rho:=f^{-1}(\Gamma_\rho)$ of the circles
$\Gamma_\rho: \vert z\vert=\rho$ of radii $\rho$ ($0\le \rho\le
1$). Question: are those still algebraic ovals when $\rho<1$? In
other words are the Riemann levels $C_\rho$ of an algebraic oval
(which is empty) still algebraic curves (at least part thereof)?
The truth of this assertion is of course a necessary condition for
our above strategy of constructing the pseudo-isotopy via Plateau.

If the initial oval is a circle (or even an ellipse), algebraicity
of the Riemann levels looks evident (at least classical I think).
For the ellipse it could involve Schwarz's explicit solution to
the Riemann mapping.

Another idea:  it would we nice if there is some flow effecting
the contraction conjectured by Itenberg-Viro. One idea could be to
take any empty oval, and look at its normal curvature flow \`a la
Huisken. Usually this is presented in $E^2$ but there is surely a
variant in $S^2$ the double cover of ${\Bbb R}P^2$ on which
Fubini-Study induces the round metric (I think). Is it true that
the normal curvature flow preserves algebraicity of ovals? If so,
the flow would shrink one of the empty oval to a point (yet not
necessarily one oval chosen in advance like by Itenberg-Viro), and
perhaps it will shrink all ovals in some succession it is alone
able to decide until all get shrunk. This could prove a weak form
of the contraction conjecture.

\subsection{Toward a naive dynamical treatment of the
Itenberg-Viro conjecture}\label{CC-via-dynamics:sec}

[18.01.13] How large can an oval  be? If we imagine a real
projective curve as traced on the sphere $S^2$ we can wonder what
the area or length of an oval can be. In degree 2 a quadratic cone
can be as large as we please and so the inside area of the oval
can be as close to $2\pi$ as we please, but of course not larger
as its ``twin'' occupy the antipodal area of the sphere. If we
restrict to even degrees then we have the Ragsdale orientable
membrane bounding the ovals (that one with $\chi=p-n$). How large
can its area be when lifted to $S^2$? For quartics and when $r=1$
we can enlarge the Fermat equation $x^4+y^4=+1$ as much as we
want by taking $x^4+y^4=R^4$ which conserve the same shape
(homothety), while covering more and more space of the sphere. So
here the upper bound is again $4\pi$ (the full sphere). (We count
now the full lift to $S^2$). But what about other quartics, e.g.
the G\"urtelkurve. Again we may imagine the latter just as a
perturbation of two concentric circles. While the outer circle
enlarges the inner contracts and so we get the Ragsdale membrane
of nearly full area.  It seems that there is little chance to find
nontrivial upper bound. But what about an $M$-quartic with 4
ovals. How large can the area of its interior be? Again we can
imagine a configuration of 2 transverses ellipses one very large
but of small eccentricity and one fairly small but orthogonal and
smooth it in the usual way to get 4 ovals; one very large and all
3 remaining fairly small. Expanding this at infinity shows that
the Ragsdale membrane can cover nearly all the sphere and the
upper bound is again $4\pi$. Such consideration could extend to
all curves constructed say via Hilbert's method. But how to treat
the general case? At any rate to each rigid-isotopy class of
curves we can consider the supremum of all areas of the
corresponding Ragsdale membrane. Is this always equal to $4\pi$?
Looking at the area of the Ragsdale membrane assigns to each curve
a numerical value in $]0, 4\pi[$. Perhaps it is interesting to
look at the orthogonal trajectories of this Ragsdale function? The
allied gradient flow could provide a dynamical flow shrinking some
empty ovals.

There is plenty of other functionals perhaps better suited to a
dynamical treatment of the Itenberg-Viro conjecture. For instance
instead of Ragsdale area we could look at the {\it empty area\/}
defined as the
cumulated area of all empty ovals. Looking at the descending
orthogonal trajectories of this function is likely to shrink
ovals. Another choice is the function looking at the area of the
smallest oval. Perhaps this has the drawback of lacking smoothness
in case two ovals enter in competition for the infimum? Another
strategy more suited to the Itenberg-Viro contraction problem is
that we are given an empty oval, and during a rigid-isotopy we can
follow him continuously. Hence given a curve with a marked empty
oval we can define in the whole chamber residual to the
discriminant (alias the rigid-isotopy class of the curve) the
functional ascribing the area of the inside of this marked but
moving oval. Of course when dragging the curve around a loop in
its chamber the oval can be
to another oval, so the function looks multivalued. Yet we get it
single-valued on the space of curves with a marked oval. So the
space of curves with a marked oval is actually an $r$ sheeted
cover of the usual space of smooth curves (with $r$ variable on
the different chambers of course).

As long as we keep the marked oval into view there is a way of
steepest descent diminishing maximally the area of the oval. For
this to make good sense we require orthogonal trajectories and so
a metric on the space of all curves. The canonical choice seems to
be the elliptic geometry on $\vert m H \vert \approx {\Bbb
P}^N({\Bbb R})$ the space of coefficients double covered by the
round sphere $S^N$.

Now follow the corresponding trajectory of steepest descent. What
can happen?
By construction our marked oval will decrease in area, but will it
docilely shrink to a point? Here are some evident difficulties
(D.$n$, $n=1,2,\dots$)

(D.1) {\it Wrong attractor (stable equilibrium).} First one can
imagine that our function as a sink trapping us into some
``depression'' like the basin of a lake yet not at zero altitude
(e.g. lake Baikal). Then our motion stops and the goal fails
blatantly, having only reached an algebraic curve realizing a
local minimum of the area yet still positive. Perhaps some clever
argument precludes such depression (e.g. if our functional turned
out to be harmonic by some miracle?)

(D.2) {\it Saddles points (unstable equilibrium).} We may of
course also reach something like a saddle point, where we need
then to choose quite randomly one of the two (or more if not
Morse) way of steepest descent. Generically up to perturbing the
initial curve, we can avoid such accidents.

(D.3) {\it Controlling the limit.} Hence let us assume that the
area shrinks to zero (assuming (D.1) to have been overcome).
Naively one can imagine the oval shrinking to a complicated
dendrite (though B\'ezout unlikely) or to a segment (again
algebro-geometrically improbable). The sole possible limit seems
to be a point. (This seems  an easy task via implicit function
theorem, B\'ezout, etc.)

(D.4) {\it Choosing the right functional (i.e. arranging a
``convex'' or ``harmonic'' landscape).} We have as yet only
considered the area functional yet the length functional looks
nearly as appropriate or better? Or even one could use a mixture
of both like the isoperimetric ratio. Note that nearby a solitary
node the behavior become nearly circular or at least elliptical
like in the local model $x^2+y^2=\varepsilon^2$ or $ax^2+by^2 =
\varepsilon^2$ with $\varepsilon \to 0$.

Among all difficulties the most stringent seems to be (D.1) which
of course as to be settled by playing with (D.4), i.e. choice of
the functional. To settle (D.1) it is enough showing that nearby
all curves there is one of smaller ``energy''. For the area
functional one could imagine an oscillation by perturbing slightly
like in the Harnack-Hilbert method our marked oval by a collection
of $m$ lines. Alas it result a vibration of the oval slaloming
across its initial position so it is not obvious how to decrease
area. So to impede getting blocked by (D.1) we are reduced to the
``local'' problem of finding an appropriate function which can
always be decreased by small (algebraic) perturbations. What about
the total curvature of the oval (or the inverse thereof as to go
to zero for a shrinking circle), etc.
If we work with area functional and if the oval is nearly
circular, we can plug in it a smaller circle and taking this
equation $k$ times (assume $m=2k$ even for simplicity). Perturbing
the curve along this multiple circle  may decrease area. This is
of course just a very special case. It seems  to apply to the case
when the curve is a deep nest in which case the innermost (=empty)
oval as to be ``convex'' (else B\'ezout is violated, compare
Zeuthen 1874 \cite{Zeuthen_1874}).

The whole problem seems reduced to that of finding a good
functional $\varphi$ which has no local minimum. Above it is not
fundamental that the oval is convex, to plug a circle inside it.
This just uses the fact  that the interior of the oval is open,
and tracing a little circle inside the oval while taking its $k$th
multiple gives another curve $k\cdot E_2$, along which  to deform
inside the spanned linear pencil $\lambda C_m +\mu (k\cdot E_2)$.
Alas
nothing ensure the oval to stretch
within its interior.

Another idea is to apply the Riemann mapping theorem and shrink
the radius of the representing circle, while hoping that this new
smaller Riemann level is still an algebraic curve of the {\it
same} degree. If this work we are able to decrease the area
functional.

This reminds me some work of Bell and Aharonov-Shapiro
\cite{Aharonov-Shapiro_1976} to the effect that the Riemann (and
more generally Ahlfors) map of a quadrature domain is algebraic,
and that quadrature domains are dense in the space of all domains
so that virtually any Riemann map is algebraic.

But in our context we have an algebraic contour(=oval) and the
following would simplify life:

\begin{conj}\label{Riemann's-level-algebraic:conj}
Given a nonempty  oval of a real plane algebraic curve of degree
$m$ and suppose the corresponding spherical calotte  conformally
mapped (via Riemann) to the unit disc $\{z: \vert z\vert \le 1
\}$. Then the pullbacks of the smaller circumferences $\vert
z\vert =r $ are still algebraic curves of the same degree
$m$!(???)
\end{conj}

If so then we can decrease area thus solving difficulty (D.1), and
perhaps the whole conjecture of Itenberg-Viro. Of course this
looks a bit optimistic (due to the a priori highly transcendental
nature of the Riemann map), but at least the dynamical strategy
looks quite stimulating. Needless to say we have not proved the
Itenberg-Viro conjecture, but in case it is true, perhaps the
above vague ideas are quite close (at least in broad lines) to its
ultimate technical solution.
We Summarize the discussion as follows:

\begin{lemma}
Let $\vert m H\vert$ be the space of all real algebraic curves of
fixed degree $m\ge 1$ and $\frak D$ be the corresponding
discriminant parametrizing singular curves. The complement
$S_m=\vert mH\vert - \frak D$ is the space of smooth curves.
Suppose given  some real positive-valued smooth functional
$\alpha$ on the space MEO of all curves with a marked empty oval
which has

(H) no stable equilibrium (local minimum) and such that if
$\alpha(C)\to 0$ then the curve $C\in S_m$ tends to a solitary
nodal curve.

Then the trajectory of steepest descent (gradient flow) always
converges toward a curve with a solitary node, after possibly
slight perturbation of the initial data (permissible as we work up
to rigid-isotopy). In particular the Itenberg-Viro contraction
conjecture holds true.
\end{lemma}

Of course this is just the formal aspect of the story (i.e.
imputable to the theory of ordinary differential equations). Yet
the real problem is to find a functional $\alpha $ suiting
hypothesis (H). A candidate is to take $\alpha$ the area of the
marked oval, and then hypothesis (H) could follow from the
optimistic Conjecture~\ref{Riemann's-level-algebraic:conj} on the
algebraicity of Riemann's level.

\subsection{Call for an attack via the
Riemann mapping (yet another Irrweg=aberration?)}
\label{CC-via-Riemann:sec}

\begin{figure}[h]
\centering
\epsfig{figure=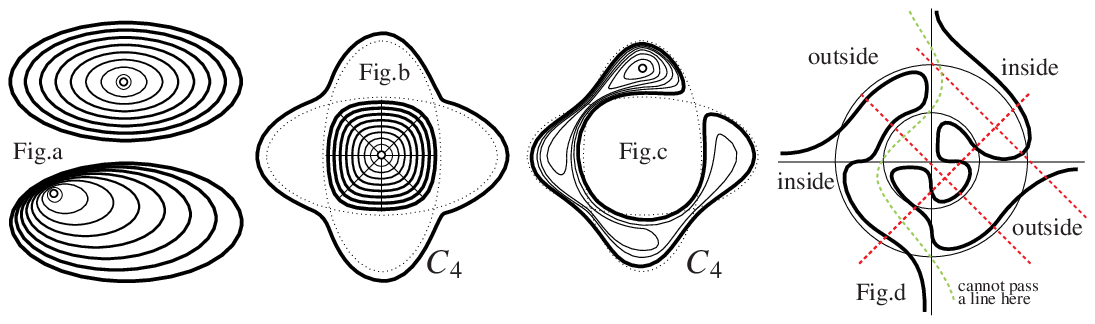,width=122mm}
\vskip-5pt\penalty0
  \caption{\label{ItenbergViroRiem:fig}%
  Some (Green)-Riemann levels of a Riemann map applied to
  the interior of an empty oval of some familiar real algebraic
  curves. Naively it is believed
  that all Riemann levels (pullback of the concentric
  circumferences under the Riemann map) are algebraic ovals of
  a curve of the same degree. On the right is Ronga's sextic
  which cannot be avoided by a line; its outside
  (homeomorphic to  M\"obius) being too contorted so-to-speak.} \vskip-5pt\penalty0
\end{figure}

[18.01.13]
Let us do some experiments. Suppose given some real algebraic
curve and take an empty oval on it. Mark an interior point and
consider the Riemann map $f$ taking the domain $D$ interior of the
oval to the unit disc $\Delta$. Pull-backing polar coordinates on
the disc gives an isothermic system of coordinates on $D$.
Fig.\,\ref{ItenbergViroRiem:fig} gives some qualitative pictures
for an ellipse or with ovals of a cubic or even of some quartics.

Naively the levels $\vert f \vert = \rho$ of the Riemann mapping
look again like curves of the same degree. This is especially
striking for the contorted quartic with one oval
of Fig.\,\ref{ItenbergViroRiem:fig}c. It is like a volcano
spreading its lava on the whole territory available in the island
interior to the oval of this $C_4$.
If algebraicity is true and degree conserved, then
the levels of the Riemann map gives directly the contraction of
the Itenberg-Viro conjecture
(\ref{Itenberg-Viro-contraction:conj}).
Is there such a miracle? Maybe giving to Riemann an algebraic
contour, the Riemann map itself is algebraic of the same degree
and so are all its sublevels $\vert f\vert = \rho= const$.

Let us first consider the more basic converse assertion that if
the Riemann map is algebraic then so is its contour.
So let $f\colon D \to \Delta$(=unit disc) be a Riemann map which
is algebraic, i.e. the power series of this analytic function is
finite, i.e. a polynomial of finite degree $f(z)=c_0+c_1 z+c_2
z^2+ \dots +c_n z^n$. Then as $f^{-1}(\partial
\Delta=S^1)=\partial D$ we see that
$$
f(z)\overline{f(z)}=1
$$
identically on the contour $\partial D=:\Gamma$. Doing the usual
splitting in real and imaginary parts by letting $z=x+iy$ and
$c_n=a_n+ib_n$ one sees that the product $f(z)\overline{f(z)}$
becomes a polynomial $P(x,y)+iQ(x,y)$, where $P(x,y)$ has degree
$n$ [{\it Warning}.---this is disproved in the sequel!] while
$Q(x,y)$ is identically zero. It follows that $P(x,y)-1=0$
identically on the contour $\Gamma$ which is
therefore algebraic. Of course the same holds true for any
sublevels of the Riemann map, i.e. sets $f(z)\overline{f(z)}=\rho$
are also algebraic. This proves the

\begin{lemma}\label{Riemann-level-algebraic}
If the Riemann map $f$ is algebraic then the contour is a real
algebraic curve and so are all sublevels. Further the degree of
all this Riemann levels $f \bar f =\rho$ (\/$0\le \rho \le 1$)
have the same degree as $f$.
\end{lemma}

This is the trivial sense. What about the converse?
Suppose given a contour which is an  oval of some algebraic curve
can we conclude
that the Riemann map is algebraic? This is what Arnold would call
pure Riemannian predestination. We think this to be true as
follows:

\begin{theorem}\label{Riemann-map-algebraic:thm}
Suppose given an (empty) oval $\Gamma$ of a real algebraic curve
$C_m$ of degree $m$. Mark any point $p$ inside this oval and let
$D$ be the sealed interior of the oval. Let $f\colon D \to \Delta$
be the Riemann mapping taking $p$ to the origin $0\in \Delta$
(unique up to rotation). Then the Riemann map is a polynomial of
degree $m$.
\end{theorem}

\begin{proof} (vague sketch) The idea is simply to look at the
function $f \bar f$ which is analytic and vanishes on an algebraic
locus, namely $\Gamma$. So by an appropriate Nullstellensatz
(Arnold's predestination? or Bloch's slogan ``Nihil est in finito
quod non prius fuerit in finito'') it must follow that the
function $f \bar f$ is itself algebraic. The proof is complete.
(???)
\end{proof}

At this stage we would have deduced Itenberg-Viro's contraction
conjecture (as a simple corollary of the Riemann mapping
theorem!!!):

\begin{cor}\label{Itenberg-Viro-via-Riemann:cor}
Any nonempty oval of a plane real algebraic smooth curve (of
degree $m$) can be contracted to a solitary node.
\end{cor}

\begin{proof} Apply Theorem~\ref{Riemann-map-algebraic:thm}
to the nonempty oval under consideration, to obtain a Riemann map
$f$ which is algebraic. By Lemma~\ref{Riemann-level-algebraic} all
the levels $\vert f\vert=\rho$ are algebraic curves of degree $m$,
while shrinking the radius to $\rho=0$ the curve acquires a
solitary node.
\end{proof}

Some objections to the method (or details to be filled):

(DET.1) How to ensure that the nodal curve so obtained has only
this solitary node as sole singularity?

(DET.2) We have worked as if the curve were affine and not
projective. This is usually a harmless nuance via the usual yoga,
(des)homogenization of the equations. More severe is our
supposition that our oval can be put in some affine chart! As I
learned from (the late) Felice Ronga (ca. 1999--02), there is a
sextic with one oval only, such that any line cuts it in at least
two (real) points. In other words no line avoids this sextic. This
is simple to construct by perturbing \`a la Pl\"ucker-Brusotti a
configuration consisting of two concentric circles, plus the two
axes of coordinates. Smoothing this sextic arrangement quite
randomly (slalom as much as you can), one gets easily the required
curve (Fig.\,\ref{ItenbergViroRiem:fig}d).
[Ronga's original picture is to be found as the front cover of his
book ``Analyse r\'eelle post-\'el\'ementaire, 1999 apr\`es
J.\,Christ.'' \cite{Ronga_1999-BOOK}.]

Such a ``Ronga curve'' causes some trouble to our procedure, which
requires putting the oval in an affine chart which we identify to
the complex plane. Perhaps this is not fatal as the Riemann
mapping theorem has some more intrinsic character, namely any
Riemannian membrane topologically equivalent to the disc is
conformal to the unit disc. So maybe one should use this more
general version, and adapt the above affine argument in this more
global setting. Algebraically this would amount to work always
with homogeneous coordinates and think with cones in ${\Bbb R}^3$.
We may then apply the Riemann mapping theorem to the spherical
calotte bounding the oval (there is two of them, but choose one),
and use as cut function homogeneous polynomials. A variant is
perhaps just to pass from the sphere covering ${\Bbb R}P^2$ to the
complex plane via stereographic projection (from a point outside
the oval). This projection is conformal, but does it preserve the
degree of polynomials?

[19.01.13] Here is a more fatal destruction of the above
pseudo-proof of the contraction conjecture via the Riemann mapping
theorem.

The key is a little computation of $f(z)\overline{f(z)}$, where
$$
f(z)=c_0+c_1z+c_2 z^2+\dots +c_n z^n,
$$
$$
\overline{f(z)}=\overline{c_0}+\overline{c_1}\bar z+\overline{c_2}
{\bar z}^2+\dots +\overline{c_n} {\bar z}^n.
$$
When expanding this product it is  useful to write it as a
``diamond'':
\begin{align*}
c_0 &\overline{c_0} \cr
c_0 \overline{c_1} \bar z  + & c_1 \overline{c_0}  z \cr
c_0 \overline{c_2} {\bar z}^2 +c_1 \overline{c_1}&  z \bar z+ c_2
\overline{c_0} {z}^2 \cr
c_1 \overline{c_2} z {\bar z}^2  + & c_2 \overline{c_1}  z^2 \bar
z \cr
c_2 \overline{c_2}& z^2 {\bar z}^2.
\end{align*}
This is the end result when $n=2$, but otherwise it will expand to
larger rhombs. At any rate, the weightiest term is $c_n
\overline{c_n} z^n\bar{z}^n=\vert c_n\vert^2 (x^2+y^2)^n$, and so
we get indeed a polynomial, but one of degree twice as big as
that of $f$.

Now suppose the given contour to be an ellipse (which is not a
circle), then even if the qualitative part of
Theorem~\ref{Riemann-map-algebraic:thm} ought to be true, i.e., if
the Riemann map $f$ of an algebraic contour is algebraic, then the
degree of $f$ cannot be $1$ (for then $f$ is  a similitude which
preserves circles). Hence the degree of the Riemann map $f$ of an
ellipse is at least $2$, hence by the above computation the degree
of $f\bar f$ is at least $4$, and so the sublevels of the Riemann
map are at least quartics, and not ellipses as we initially
imagined (compare e.g. the misleading
Fig.\,\ref{ItenbergViroRiem:fig}a). Hence even granting
algebraicity of the Riemann mapping of an algebraic contour,
the corresponding levels would not be of the
same degree.

Actually the above argument shows that our strategy fails for a
conic, but is it really a disproof in general? Assume the initial
curve $C_m$ to be a quartic, then the degree of the Riemann map
could be $n=2$, and thus the degree of $f\bar f$ is four and so
 the sublevels would stay a deformation within quartics. More
 generally this numerology makes sense for any $C_{2k}$ curve of
 even degree $2k\ge 4$, by taking $n=k$.

A more severe objection is surely, as already apparent by looking
at the highest power $c_n \overline{c_n} z^n\bar{z}^n=\vert
c_n\vert^2 (x^2+y^2)^n$, the fact that $f\bar f-1=0$ is not the
most general curve of degree $m=2n$ (e.g. the monomial $x^{2n-1}y$
is missing). Hence there is no chance the Riemann map of any curve
$C_{m=2k}$ being algebraic of degree $k$, and our strategy is
definitively foiled.

As a modest consolation one can apply the above method to some few
ad hoc curves with equations of the shape $f\bar f-1=0$ where $f$
is some algebraic Riemann map. The corresponding curves (say
``Riemann curves'') could be contracted by the above recipe. This
is of course far from settling the initial desideratum of
Itenberg-Viro.

What can be retained from this attempt? Let us start with a
polynomial $f(z)\in {\Bbb C}[z]$ of degree $n$. This induces a
holomorphic map ${\Bbb C}\to {\Bbb C}$ of degree $n$
(Gauss--D'Alembert fundamental theorem of algebra) which is a
branched covering. Sometimes it turns out that the unit disc is a
trivializing open set for this covering. In the language of
complex analysts (Bloch, Landau, Ahlfors, etc.) this is also what
they would call a schlicht unit disc. Taking one among the $n$
many sheet lying over $\Delta$ gives a simply connected domain
which by $f$ is conformally mapped to the disc, and so we recover
a Riemann map. For all those domains (which are algebraic of
degree $2n$ via the equation $f\bar f-1=0$ but with $n$ unnested
ovals), we may a priori implement our contraction algorithm.
(Overlooking the unnested condition, we could hope that such
Riemann curves are spread in all chambers of the discriminant and
hope a general attack.) Now back to our setting it must be noticed
that during the shrinking $f \bar f=\rho$, with $\rho \to 0$ our
algebraic curve sees all its $n$ ovals being simultaneously
shrunk. Hence in this very favorable setting, the solitary node
condition fails blatantly. Conclusion: our strategy via Riemann
leads nowhere, if it pertains to implement the contraction
conjecture of empty ovals.

\section{CCC: collective contraction conjecture,
as an avatar of Itenberg-Viro (Gabard 2013)}

{\it Insertion} [04.04.13].---As we noticed only today, in
substance it seems that the conjecture posited below bears some
 analogy with Klein 1892 \cite{Klein_1892_Realitaet} (p.\,177
of Ges.\,Math.\,Abhd.) who wrote: ``Wir k\"onnten z.\,B. mehrere
Z\"uge unserer Kurve gleichzeitig in isolierte Doppelpunkte
\"uberf\"uhren''. More generally Klein 1892 discusses at this
place (p.\,176) what he calls the ``Doppelpunktsmethode''
amounting essentially to contract any symmetry-line of the Riemann
surface, and this of course seems to anticipate what we called
before the Itenberg-Viro contraction principle. It is not clear
however that Klein ever formulated something as precise as the
Itenberg-Viro contraction conjecture (specific to plane curves).
On p.\,176--177, Klein writes something which taken out from its
context looks a bit overoptimistic namely: ``Bei allen anderen
F\"allen hat die Durchf\"uhrung des genannten Prozesses und damit
die Zusammenziehung eines beliebiegen Ovals der Kurve zu einem
isolierten Doppelpunkte keine Schwierigkeit.''

\subsection{Failing with Riemann suggests
a variant of Itenberg-Viro, viz. CCC=collective contraction
conjecture: deformation in the large as a method of
prohibition}\label{CCC:sec}

[19.01.13]
The above discussion suggests the following variant of the
contraction conjecture, which has maybe some spontaneous appeal
and independent interest. Suppose given a projective smooth plane
real curve $C_m$. Look at all empty ovals simultaneously. Is it
possible to shrink all of them in one single stroke toward
solitary nodes (via a deformation of smooth curves sole for its
end-point being in the discriminant)?

It seems likely that the contraction conjecture
(\ref{Itenberg-Viro-contraction:conj}) implies this, roughly by
shrinking one oval, and then the second, etc. Of course one then
needs to arrange a bit things so that prior to
extinct  one oval completely, one waits until the second empty
oval becomes ``small'' enough, etc. Finally one synchronizes the
ultimate ``coup de gr\^ace'' to kill all the empty ovals at the
same time (``time'' being just the parameter of the path $[0,1]
\to \vert mH \vert$ in the space of all curves of order $m$).

{\it Insertion} [02.04.13].---It may help reading the sequel to
remarked first that the reverse process, of deducing an individual
(solitary) contraction (\`a la Itenberg-Viro) from our collective
one, is much easier and a trivial consequence of Brusotti, if we
did no mistake, cf. Lemma~\ref{CCviaCCC-Brusotti:lem} below. So
the conjecture posited right below is stronger than the one of
Itenberg-Viro (yet perhaps equivalent, or at least easier to
disprove).

So let us (somewhat
cavalier) formulate the:

\begin{conj} {\rm (Collective contraction conjecture=CCC,
[19.01.13, 22h40])}\label{CCC:conj} Given any smooth real curve
$C_m$ of degree $m$, it is possible to shrink all the empty ovals
simultaneously toward solitary nodes. (Solitary but
synchronized death of all ovals.)
\end{conj}

This is obviously true for $m=2,3$ (being actually equivalent to
the individual contraction principle) since there is at most one
empty oval available.

The case $m=4$ is already more tricky, yet still compatible with
B\'ezout. If $r=4$ ($M$-quartic), we would have a quartic with
four isolated (solitary) nodes. This exists just take an imaginary
conic $C$, and aggregate it with its conjugate $C \cdot C^\sigma$
(this is real but a priori the four intersections need not all be
real). More simply take two transverse conics, look at signs and
arrange a level so that there are 4 isolated points by making a
naive picture of the graph of $E_2 \cdot F_2$. Since the real
scheme encodes completely (in degree $m=4$) the rigid-isotopy
class (Klein 1876, etc.) it follows that CCC holds true in degree
$m=4$. The case $r=2,3$ are treated similarly by looking at the
graph of a special equation and passing a plane tangent to the 3
(or lesser) hills.

Now what about degree 6?
%
Deciding the truth of the above conjecture in degree 6 (CCC6), is
already more tricky. Perhaps this follows from Itenberg's CC6, if
not formally by the method used therein, i.e. Nikulin's theory
with $K3$-surfaces. As said at the start it could be that CC
implies CCC in all generality. In the sequel we assume CCC as
granted and look what can be derived from it.

Let us first suppose that there is an $M$-sextic $C_6$
with 11 unnested ovals (what Hilbert, Rohn, Petrovskii, Gudkov,
Arnold, etc. were fighting hard against). Shrink all of them to a
point according to CCC (\ref{CCC:conj}). Then the Riemann surface
is strangulated along all its oval in two (algebraic) pieces which
are topological spheres. Since the nodes are supposed to be
solitary these two pieces are smooth curves of genus 0
intersecting transversally. Therefore (via the genus formula
$p=\frac{(d-1)(d-2)}{2}$) they are of degree $1$ or $2$, but have
to intersect in 11 points. B\'ezout is overwhelmed! (Alternatively
the genus formula is corrupted, since we have a degeneration of
$C_6$ toward two cubics $C_3$ and its conjugate $C_3^\sigma$!)
This gives a new ``proof'' of Hilbert-Rohn-Petrovskii via CCC. Of
course it is quite tempting to wonder if Hilbert (or subsequent
workers) did not knew about this argument at least as a heuristic
tool.

More generally:

\begin{prop} {\rm (like Hilbert 1891)}
Under axiom CCC, a smooth $M$-curve (of even degree) cannot have
all its ovals unnested unless its degree $m$ is less than four
($m\le 4$). In particular an $M$-sextic cannot have all its $11$
oval unnested (which is Hilbert's original claim as early as
Hilbert 1891 \cite{Hilbert_1891_U-die-rellen-Zuege}.)
\end{prop}

\begin{proof}
Shrinking collectively all the empty ovals of $C_m$ (via CCC)
gives a splitting $C_m \to C_d \cup C_d^\sigma$ in two algebraic
curves of degree $d=m/2$ intersecting transversally in $r$ points.
So by B\'ezout $d^2=r$. Since both strangulated halves have genus
$p=0$ (for we started from an $M$-curve), their common degree $d$
can only be $1$ or $2$. Hence $d\le 2$, and so $r=d^2\le 4$. Since
$r=g+1$ ($M$-curve assumption) it follows $g\le 3$ and so $m\le
4$. (Variant: conclude more directly via $d=m/2$.)
\end{proof}

As a matter of philosophical dilettantism (?), it may be wondered,
whether Hilbert himself used the above argument, at least as a
heuristic tool. To my knowledge there is no record in print along
this sense. Yet, Hilbert, say unlike Poincar\'e was a formalist,
in particular never writing down crazy ideas. Thus, it may be not
be impossible (our subjective speculation) that Hilbert may have
argued along this route. In fact it is probably more realist  that
Hilbert argued along another idea, cf. e.g. the passage of Gudkov
1974 \cite{Gudkov_1974/74}, where Hilbert's method is described as
implemented by his students Kahn and L\"obenstein.

Even more generally:

\begin{prop}
Under CCC, a smooth dividing curve (of even degree $m=2k$) cannot
have all its ovals empty (equivalently unnested) unless:

{\rm (1)} its number of ovals $r$ is a square (\/$1,4,9,\dots$),
and actually the square of its semi-degree $k$\footnote{This
conclusion actually holds true unconditionally as follows from
Rohlin's formula $2(\Pi^{+}-\Pi^{-})=r-k^2$.};

{\rm (2)} some stringent arithmetical conditions (say
predestination or coincidence) are verified, namely all the
displayed formulas in the proof below have to be satisfied.
\end{prop}

\begin{proof}
By CCC shrink all empty ovals of the curve $C_m$. The Riemann
surface $C_m({\Bbb C})$ (``complexification'') has genus
$$
g=(m-1)(m-2)/2.
$$
Once strangulated, it splits in two Riemann surfaces of genus
$$
p=[g-(r-1)]/2
$$
(since $g=(r-1)+2p$ by visualizing the orthosymmetric surface).
Both halves are algebraic smooth curves intersecting in $r$
points. Being actually interchanged by conjugation, they have some
common degree, say $d$, verifying
$$
 p=(d-1)(d-2)/2.
$$
So by B\'ezout (or homological intersection theory) we infer
$$
d^2=r.
$$
Hence $r$ must be  a square.
In fact a sharper argument based on the degeneration $C_{2k}\to
C_k\cup C_k^\sigma$ shows that $r=k^2$ directly by applying
B\'ezout to both halves of the limiting curve of the collective
contraction.
\end{proof}

A this stage ``Eureka'' [23h41] we have already proved  that the
sextic scheme $5$ is necessarily of type~II (as followed first
from Rohlin's complex orientation formula).

Likewise the sextic schemes $\ell$ ($\ell = 9$ excepted)  cannot
admit a type~I incarnation (though this was already implied by
Klein's congruence modulo 2, and Arnold's congruence mod 4), safe
for $\ell =1$
where either Rohlin or our suggestive
geometric argument do instead the job.

Indeed if $r=1$ and $m=6$, then both strangulated parts have genus
$p=5$ (imagine the Riemann surface of genus $g=10$ split by the
one oval), but this is not even the genus of a smooth curve. Hence
strangulation is impossible violating axiom CCC.

Of course  the philosophy behind CCC is quite akin to the filling
trick of Arnold-Rohlin  safe that the closing is God given by some
postulated (but hypothetical) shrinking procedure in the rigid
algebraic category.

In general it remains the boring task of extracting the exact
arithmetical consequences of CCC, while checking if it is really
compatible with factual data. In degree $6$, CCC seems to live in
perfect harmony with Rohlin's enhancement of Gudkov's table
(Fig.\,\ref{Gudkov-Table3:fig}). Since the usual (individual)
contraction conjecture CC holds true in degree $6$ (by Itenberg
1994 \cite{Itenberg_1994}), it is likely that the collective
variant CCC holds good as well.
%
%
Of course all the arithmetical relations are in reality less
stringent that they look at first glance, since they are all
coming from the genus formula which itself may be interpreted as a
surgical process regulated by B\'ezout. (Recall the simple proof
of the genus formula based on the morphogenesis of lines getting
smoothed under surgeries.)

Applying CCC to Harnack's sextic configuration leads nowhere since
the Riemann surface keeps connected.

Another exercise: assume there is a dividing quartic  with two
unnested ovals. Apply CCC to both ovals. Then the Riemann surface
of genus $g=3$ is strangulated in two surfaces of genus $p=1$,
hence of degree $d=3$. But the latter cut themselves in 9 points
and not two. This contradiction reproves (modulo CCC) the
well-known fact (\ref{Klein-unnested-quartic-nondividing:lem}) due
to Klein 1876, Arnold 1971, Rohlin 1972--1978, Wilson 1978, Marin
1979, Gross-Harris 1981, etc., that a quartic with two unnested
ovals is necessarily nondividing. (Variant of the argument: $r=2$
is not a square.)

The principle emerging is that large deformations prompted by
contraction conjectures affords a puissant method of prohibition,
as opposed to the method of small perturbations which is merely a
toolkit for construction.

[20.01.13] At this stage, the method CCC looks quite powerful, at
least as a heuristic tool, reducing to B\'ezout several deep
assertions and results of Klein, Hilbert, Rohlin, etc. However as
yet the method is quite limited to the case where all the ovals
are empty so that the strangulation really implies an algebraic
splitting of the dividing Riemann surface.

Perhaps the method can be extended beyond this proviso. For
instance after shrinking the first generation of all empty ovals
and making them effectively disappear from the real locus, a
second generation of empty ovals appears, which would be
contracted in turn to solitary nodes, etc. Iterating so collective
contractions  would contract all ovals and so
achieve a splitting of the dividing Riemann surface. For the
method to be effective it seems that the solitary nodes of the
first generation ought to resurface as such right after the
contraction of the second generation of empty ovals. All this
looks very dubious, yet some more clever intelligence can perhaps
extract something from this procedure.

Consider a specific example to make the difficulty more concrete.
Consider the sextic scheme $\frac{10}{1}$, and let us call it
(improvising terminology) the Rohn scheme (for Rohn 1911--13
\cite{Rohn_1913}) was the first attempting to disprove its
existence via
a substantiation of Hilbert's method. Let us contract all empty
ovals. The resulting Riemann surface is still connected, and we
lack a splitting suited to an application of B\'ezout. The obvious
idea is to make first an eversion of the nonempty oval so as to
reduce to the unnested scheme $11$. (For the definition of
``eversion'' cf. Sec.\,\ref{Eversion:sec}.) So we need another
highbrow large deformation principle, dual to the contraction
principle stating that any maximal oval can be everted provided
the resulting scheme is not prohibited by B\'ezout.
Alas, eversions are not permitted for $M$-curves by virtue of
Lemma~\ref{eversion-impossible-for-M-curves:lem}, and this
stratagem looks jeopardized.

{\it Insertion} [01.04.13].---After the collective strangulation
of all empty ovals of a Rohn curve of scheme $\frac{10}{1}$, we
would get a Riemann surface of genus 0 and degree 6, but with 10
nodes. Alas this is still permissible! Of course our argument
being merely abstract (i.e. using the abstract topology) would
equally well apply to the veritable $M$-curves, and this explains
that.

Consider next the sextic $M$-scheme $\frac{2}{1}8$. Shrink all
empty ovals to get a curve of genus $p=0$ with one oval (once it
is desingularized). So it is really just a rational curve of the
dividing type (like our orthosymmetric equatorial planet Earth).
By a trivial case of Ahlfors (actually the Riemann mapping
theorem) this has a unique conformal structure. So there is a
total map of degree $1$. This in turn gives a total pencil of
curves, of order say $k$. It seems clear that all solitary nodes
of our (contracted) $C_6$ must be in the base locus of the pencil.
Indeed else the pencil is sweeping out some node, and so the curve
through it has one intersection (counting double) but zero nearby
whence a disappearance in the imaginary locus (violating total
reality). So our pencil must exhibit at least $10$ basepoints.
Naively the disc inside the nonempty oval looks foliated by two
foyers (index$=+1$) violating Poincar\'e's index formula, but this
looks too naive because it would kill as well Hilbert's or
Gudkov's sextics. So we are again confronted to some complicated
foliation argument which we know to be quite difficult to
implement. We could dispense using CCC, by applying instead
directly Bieberbach-Grunsky to the smooth curve $C_6$ (i.e. the
genus zero case of Ahlfors 1950 \cite{Ahlfors_1950}). All this
looks difficult and let us abort here shamefully.

[23.01.13] Let us insert here an optional side remark. It is
tempting to wonder what follows from Thom's conjecture
(Kronheimer-Mrowka theorem meanwhile of 1994) to the effect that
algebraic curves minimize the genus among smooth orientable
surfaces embedded in ${\Bbb C}P^2$ realizing the same fundamental
homology class. If one fills by the half of a hypothetical real
$M$-sextic of real scheme $11$ (eleven unnested ovals) by the
interior of all ovals,
 one obtains a smooth surface of genus $0$ (round the corners) of
 degree $6/2=3$. This violates Thom-Kronheimer-Mrowka, which
 therefore implies again the Hilbert-Rohn-Petrovskii prohibition.

{\it Insertion} [01.04.13].---The degree 3 case of Thom is
really due to Kervaire-Milnor 1961 \cite{Kervaire-Milnor_1961},
based on deep works by Rohlin, ca. 1952.

Some thinking shows however that Thom's conjecture does not imply
 much more, for the filled membrane has then  genus $\ge 1$.

\def\C++{C++}

\subsection{Do iterated contractions (\C++) imply Rohlin's formula?}

[02.04.13] We recommend to skip this section which is neither
exciting nor seriously written.

[22.01.13]
One may wonder if there is not a stronger mode of degeneration
(alias contraction) of a real smooth plane curve $C_m=C_{2k}$ than
CCC yielding  the Rohlin complex orientation formula
$2(\Pi^{+}-\Pi^{-})=r-k^2$ as a corollary of B\'ezout. In the case
of no nesting this is precisely what did the previous section.
Indeed the solitary node degeneration of CCC implies symbolically
$C_{2k}\to C_k \cup C_k^\sigma$ (a topologically dividing curve
divides algebraically!), whence by B\'ezout $r=k^2$. This
coincides with Rohlin's formula since there is no nesting.

Let us call such a hypothetical mode of degeneration C++ (like
Turbo Pascal?). If this exists this would be a geometrization of
Rohlin's formula, in the sense that topology (homological
intersection theory \`a la
Poincar\'e-Lefschetz-Weyl(1923)-Pontrjagin-de Rham (1930), who
else?) would be subsumed to Monsieur \'Etienne B\'ezout (ca. 1768)
alone. This would also posit a wide extension of the Itenberg-Viro
conjecture. It seems evident that such a contraction C++ should
exist. It remains only to be not overwhelmed by the combinatorics.

So suppose given a dividing curve $C_{2k}$ in the plane. The idea
is to contract all its ovals so as to split the curve in two
algebraic pieces exchanged by Galois(=conj), all this being just
caused by a strangulation of the underlying Riemann surface.

Of course we apply first CCC to contract all empty ovals toward
solitary nodes. Then it  appears a second generation of nearly
empty ovals (those which formerly were at height 1 in the tree of
the nesting structure). We may hope to shrink those in turn while
necessarily coalescing together all the solitary nodes inside this
oval. Perhaps we can do this while keeping the tangents distinct
at those solitary nodes gravitationally clumped together. One
continues this big crunch process and once all ovals have been
contracted one get a splitting $C_{2k}=C_k\cup C_k^{\sigma}$.
Counting properly intersections with B\'ezout should give Rohlin's
formula
$$
k^2=r-2(\Pi^+ - \Pi^{-}).
$$

Of course we need to be much more explicit (as if Rohlin would not
have influenced us). Recall that $\Pi^{+}$ is the number of
positive (injective) pairs of ovals that is with complex
orientation matching that of the bounding annulus of ${\Bbb
R}P^2$, and likewise $\Pi^{-}$ being the number of pairs with
disagreeing orientation when induced from the complexes versus the
real bounding annulus. Recall that all pairs are taken into
consideration not just oval succeeding themselves immediately.

On applying first CCC we can shrink all empty ovals to solitary
nodes. Naively one would then like to shrink all the nearly empty
ovals containing only solitary nodes, and so on. So we would have
a degeneration $C_{2k}\to C_k\cup C_k^\sigma$. Computing the
intersection $C_k\cap C_k^\sigma$ with B\'ezout gives $k^2$
algebraically. Geometrically, as each oval is shrunk to a pair of
conjugate lines, and this explains the presence of the term $r$ on
the RHS of Rohlin's formula (as displayed above). Further  the
Riemann surface of the reduced curve can be naively imagined in
3-space as a pair of paraboloid of revolution together with their
orthosymmetric replicas. Each branches intersect its conjugate in
one point, but those contribution where already taken into
account. So it remains to count the intersection of the top small
paraboloid with the bottom large paraboloid, and vice-versa the
large top with the bottom small. So we get 2 additional
intersections, and this explains the term $+2\Pi^{-}$ of Rohlin's
formula. (Alas the term $-2\Pi^+$ looks much harder to explain.)

More clarification is required. A negative pair of ovals can be
shrunk  simultaneously. An example is provided by the
G\"urtelkurve, quartic $C_4$ with two nested ovals. Either via
Fiedler or by Ahlfors it is plain that this $C_4$ has a negative
pair of ovals. Looking at an equation like two concentric circles
and perturbing slightly to get away from the reducible locus (and
the discriminant) we have
$$
(x^2+y^2-\rho^2) (x^2+y^2-R^2)=0,
$$
and if $0\le \rho<R$ then we can shrink $R\to 0$ and obtain the
required multi-contraction. This example obviously extends to deep
nest in any (even) degree, as to shrink negative towers of ovals.

Our guess is in contrast that positive pair of ovals resist
simultaneous shrinking. Probably there is an evident topological
obstruction which I missed to notice as yet.

If so then there is no possibility to reduce Rohlin's formula to
B\'ezout via a super strong contraction principle C++ reducing the
whole curve to a microcosm of solitary nodes. If so, the whole
curve under the action of some gravitational clumping would truly
reduce to a constellation of isolated points with real scheme
condensed at the atomic scale. (Imagine points and then
infinitesimal circle surrounding the first generation of point,
etc.) It seems more likely that there is an obstruction to shrink
everything (algebraically), and so Rohlin's proof is surely the
best one can implement. Yet there could still be some geometry
behind it suggested by the contraction principle.

{\it Insertion} [01.04.13] The latter is true in degree 6, and
actually stronger than Rohlin's formula,  when combined with RKM
(\ref{Kharlamov-Marin-cong:thm}) since it rules out all schemes
lying above the $(M-2)$-schemes of type~I (e.g. $\frac{7}{1}3$),
what Rohlin's formula  is unable to do alone.

\subsection{Failing to reduce Rohlin to B\'ezout suggests again a
dynamical approach}\label{CCCviaDynamics:sec}

[22.01.13] What is this geometry and is it worth paying attention
at? Before trying answering this, note that even if a positive
pair of ovals resists to shrinking it could undergo another type
of Morse surgery, namely coalescence to a figure eight (lemniscate
looking like a ``sweetheart'', i.e. with one branch lying inside
the other). Yet this operation corresponds to the contraction of
an ortho-cycle without disconnecting the Riemann surface so as to
produce an algebraic splitting $C_{2k}\to C_k \cup C_k^\sigma$.

Let us turn to the geometric aspect. Our goal is essentially to
shrink the ovals at least those which are empty (CCC), and then
eventually push further the deformation as to shrink the negative
pairs (memno-technic trick imagine
negative=depressive=shrinkable). The other positive pairs may
offer some resistance (due to a topological obstruction, which we
should still understand better).

To achieve such a shrinking it looks natural to look at the
length-functional  of {\it all\/} ovals (not just the empty ones).
Consider the round metric on the unit sphere $S^2$ lying above
${\Bbb R}P^2$, and measure lengths in this metric. Given $C_m$ a
curve (i.e. a homogeneous ternary form $F_m(x_0,x_1,x_2)$ with
real coefficients up to homothety), look at the set $(F_m=0)\cap
S^2$ which is obviously rectifiable (Lebesgue, Jordan, Riemann,
Gauss, Archimedes, etc.). Denote its length by $\lambda (C_m)$.
This is zero iff $C_m({\Bbb R})$ is empty or contains merely
isolated points. Further there is an obvious way to take into
account the multiplicity of branches; e.g. a conic (degree 2)
consisting of a double line has length not just $2\pi $ but twice
that quantity. This is crucial to ensure continuity of $\lambda$
on the parameter space of all $m$-tics. Since the latter space is
compact (actually an ${\Bbb R}P^N$, $N=\binom{m+2}{2}-1$) the
length functional $\lambda$ reaches a maximum.

How long can an $m$-tics be? By the above compactness argument
there is some universal constant $L_m$ bounding the length of all
curves $C_m$ of some fixed degree $m$:  $\lambda(C_m) \le L_m$,
and the maximum is actually realized (a priori not by smooth
curve). A configuration of $m$ lines produces the lower estimate
$2\pi m \le L_m$, and by Brusotti 1921 \cite{Brusotti_1921} there
are smooth curves $C_m$ of length as close as we please to $2\pi
m$, but slightly longer. Imagine a crossing getting smoothed then
the geodesic of $S^2$ are entailed
by curvilinear arcs which are longer (triangle inequality or
Pythagoras in the small). So there certainly exist longer curves!
But how long can an $m$-tic be? This is   probably very difficult
to answer.

When $m=2$, the above argument via Brusotti still makes sense. If
we imagine quadratic cones in 3-space $E^3$ (say with elliptical
affine cross-section at $x_2=1$), then they may cut strange ovals
on $S^2$ possibly longer that $2\pi \cdot 2$???

When $m$ is odd then there is a pseudoline (or at least a circuit
possibly singular) not null-homotopic in ${\Bbb R}P^2$. Obviously
its length is at least $2\pi$, which is the lower bound of the
functional $\lambda$ when $m=2k+1$ is odd.

The estimation of $L_m$ is surely an attractive problem, but let
us try to be not sidetracked by this. Our goal would be rather to
study the gradient flow of $\lambda$ as a dynamical process
susceptible to implement the collective contraction conjecture
CCC, or more elaborated versions thereof if feasible (e.g. some
iterated contractions like C++).

Assume first $m$ even. Then $\lambda$ vanishes identically on the
chamber corresponding to empty curves, as well as on its adherence
which consists of curves with isolated real points (either
solitary double points or of higher multiplicity necessarily
even). Thus $\lambda $ cannot be analytic, but seems rather being
$C^\infty$. (As usual with distance functions (e.g. $\vert x
\vert$) they sometimes lack even smoothness until tacking their
squares. So perhaps take $\lambda^2$ squared.)

Consider the gradient flow of this functional $\lambda$, while
hoping that the corresponding trajectories materialize the
collective contraction conjecture (CCC).

Usually ovals fails severely to be geodesics on $S^2$, but are
perhaps so when we look them in the Riemann surface $C_m({\Bbb
C})$ endowed with the Fubini-Study(=FS) metric on ${\Bbb C}P^2$.
Is the corresponding length of the ovals the same, in other words
does FS induce the round elliptic metric on $S^2$? It is (always)
tempting to regard ${\Bbb C}P^2$ as the variety of groups of two
points on the Riemann (round) sphere. In this model how to
describe the FS-metric?

Another idea, at least if we restrict to smooth curves, is to
take the uniformizing hyperbolic metric (of
Schwarz-Klein-Poincar\'e-Koebe) on $C_m({\Bbb C})$ with curvature
$K\equiv -1$, when $m\ge 4$ (so $g\ge 3$). Then we get another
measure of length of the ovals, which we shall denote $h$. The
problem here is that this length functional is not a priori
defined on the full hyperspace of curves $\vert m H\vert$.

As soon as we look also in the complex domain, there is a myriad
of other functionals like the systole of the Riemann surface, the
area of one half in the dividing case, etc. We just remark that
from the systolic viewpoint there might by an ortho-cycle of much
shorter length than the real ovals, and which dynamically might be
advantageous  being first contracted.

Of course the technical difficulties look immense, but the problem
involves a mixture of Poincar\'e-Morse versus Hilbert-Petrovskii,
i.e. a synthesis thereof. So the game is certainly worth paying
attention at. What seems called upon is a dynamical study of
algebraic equations governed by motions regulated by (natural)
geometric functionals on the corresponding varieties (zero loci).
In particular find appropriate functionals whose trajectories
converge (generically) to curves with solitary nodes as to
implement the Itenberg-Viro contraction conjecture or its
collective variant CCC.

To shrink all empty ovals simultaneously it seems  not so fruitful
to shrink the shortest oval. More  collective optimization is
asked as if one had to bring fastest to the harbor a convey of
ships each carrying rough materials involved in the manufacture of
some complex end-product (Poly\`a's metaphor for Rayleigh
eigenvalues). Here we are in a similar situation. If all empty
ovals have to dye simultaneously (scenario posited by CCC), it is
important to shrink faster the longer ovals. Perhaps this suggests
looking and $\lambda^2$ the squared length functional penalizing
longer ovals.

It is also tempting to look at the area $\alpha$ (of the interior
of all empty ovals measured on $S^2$), and to play perhaps with
the isoperimetric inequality. For instance the functional
$\lambda^2/\alpha$ looks natural, and is bounded from below in the
small by $(2\pi \rho)^2/(\pi \rho^2)=4\pi$, so it admits a finite
limit when it shrinks. The isoperimetric functional
$\iota=\lambda^2/ \alpha$  intuitively forces ovals to dye in a
round manner, penalizing agonies along eccentric ellipses.


If optimistic, integrating the gradient flow of either $\lambda$
or $\alpha$, length resp. area of the empty ovals directly leads
to a solution of CCC.

The serious obstacle is that there may be a sink, i.e. a local
minimum of the functional preventing convergence to a curve with
solitary nodes arising as contractions of the $r_0$ empty ovals.

\begin{figure}[h]
\centering
\epsfig{figure=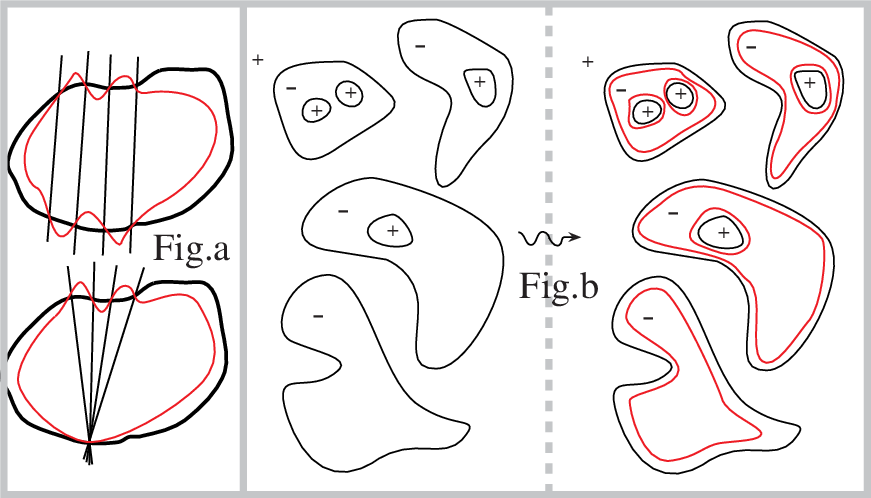,width=122mm} \vskip-5pt\penalty0
  \caption{\label{Vibrate:fig}%
  Vibrating as an attempt to decrease area} \vskip-5pt\penalty0
\end{figure}

One naive idea is to let vibrate the ovals via a configuration of
lines (as in the Harnack-Hilbert method),  hoping to decrease area
through this perturbation. The oval then oscillates inside and
outside itself but on a larger portion it would move inside
himself, hence area decreases (cf. Fig.\,\ref{Vibrate:fig}a). The
similar assertion for the length functional looks even more
fantasist. Hence the area functional looks better suited to the
problem.

{\it Insertion} [01.04.13] One problem is that if we have several
ovals it is not clear that decreasing the area of one will not
enlarge area of the other empty ovals. This problem dissipates
somewhat if we look only at the usual Itenberg-Viro conjecture,
but of course also the latter is subject to doubts, e.g. those
allied with Shustin's disproof of Klein-vache.

A more naive idea is to look at some sublevel of the equation. For
simplicity assume the given curve completely inside some affine
chart (cf. however Ronga's counter-example on
Fig.\,\ref{ItenbergViroRiem:fig}) with equation $f(x,y)=0$, and
w.l.o.g. positive on the outside unbounded (in ${\Bbb R}P^2$
nonorientable) component residual to the curve. Then the sublevel
$f(x,y)=-\varepsilon<0$ (small negative constant) ought to have a
smaller area $\alpha$ (compare Fig.\,\ref{Vibrate:fig}b).

Let us be more precise. Look at all empty ovals of the curve
$C_m$. Then following Ragsdale-Petrovskii, some are positive and
some negative (or, even and odd depending on the parity of the
number of ovals surrounding it). Maximal ovals are even (being
surrounded by zero ovals), those immediately inside them are odd,
etc. Under our sign convention for the equation $f$ (positive
outside) we see that even ovals decrease in area when considering
the sublevel $f=-\varepsilon$, while  odd ovals increase in area
(compare Fig.\,\ref{Vibrate:fig}b). The net bilan is hard to
quantify, but on the situation of the picture where the even empty
oval is much larger than the other empty ovals there is some
chance to decrease the functional $\alpha$.

Is there some chance to deduce a general argument from our naive
picture? Split all empty ovals in even  and odd ones (denoted
resp. 0 and 1 depending on their class modulo 2). Look which of
both collections has more massive total area, i.e. compare
$\alpha_0$ vs. $\alpha_1$. If $\alpha_0>\alpha_1$ then take
$\varepsilon$ positive (and vice-versa if $\alpha_0<\alpha_1$ then
take $\varepsilon$ negative). Of course if unlucky both magnitudes
$\alpha_i$ are equal, in which case we are a bit lost (perhaps
avoidable by genericity, as we work up to isotopy, hence can
always perturb slightly the data).

Let us assume to be in first case $\alpha_0>\alpha_1$ (as on the
picture). By a classical continuity lemma (cf. e.g. Gudkov 1974
\cite{Gudkov_1974/74}) we can preassign a tubular neighborhood of
the curve in which the perturbation will stay confined. Thus we
can probably control from above the area expansion of the odd
ovals. It seems indeed that the tubular expansion is maximum for a
circle as follows from the isoperimetric inequality. However it is
not clear that the deflation of area of the even ovals ought to
supersede this as it may be a very thin penetration (imagine the
big south island as very mountainous hence poorly affected by a
raise of the ocean level).

This leads to the idea of looking at the normal derivative of the
defining (polynomial) function $f(x,y)$ across the sea level
($f=0$). The inflation of area resp. deflation of area of the
empty ovals ought to be
proportional to this normal slope
and the length of those ovals via some explicit formula given by
differential calculus. So we
compare both quantities for even an odd ovals,
and choose the right sign for $\varepsilon$ in order to create a
deflation of area while keeping the degree constant. Of course if
both quantities coincide we are a bit disturbed, but perhaps
avoidable by genericity (two random real numbers are generically
distinct.)

If this trick works there is some hope to show that the area
functional $\alpha$ lacks local minima, and the corresponding
orthogonal trajectories of steepest descent ought to converge
toward curves with solitary nodes.

Note another phenomenon: imagine one empty oval shrinking
prematurely before the others. Soon after this death another oval
(formerly nonempty) may become suddenly empty implying a large
jump of the area functional $\alpha$ which looks therefore
discontinuous. The only reasonable parade against this catastrophe
is that the trajectory of steepest descent will not kill abruptly
the small oval as it is much more profitable to shrink first the
voluminous ovals. Intuitively,  the $\alpha$-flow would promote a
collective contraction. Nonetheless the $\alpha$-functional can be
discontinuous with big jumps across walls of the discriminant, as
caused by the fact that we measure only empty ovals, which in
contrast to all ovals, is subsumed to violent fluctuation.

Have we proved something? Maybe yes if quite sloppy. Let us resume
some of the difficulties:

(D.1) First Ronga's example of a $C_6$ not confined to an affine
chart is presumably not a severe obstacle. For even degrees the
sign of the projective equation $F(x_0,x_1,x_2)$ is always
well-defined so that there is a variation $C_m^{\varepsilon}$ with
disjoint real locus ($C_m\cap C_m^{\varepsilon}=\varnothing$).
(For odd degrees nothing similar can be done so easily, but we
confine  attention to even degrees. That is already hard enough.)

(D.2) We look at the penetration index under a small perturbation
as measured by the normal derivative of the landscape pondered
against the length element of the oval. More precise, calculate
along each point of $C_m$ the normal derivative
$
\frac{\partial f}{\partial n},
$
and then integrate this against the length element $ds$ of some
fixed oval $O$ to get
$$
\int_{O=an oval} \frac{\partial f}{\partial n} ds=: \pi (O).
$$
This real number measures the
rate of area change under a flood (variation of $\varepsilon$). Of
course the normal derivative above
can be interpreted as the gradient of $f$ on the coast line $f=0$.
(If it is big in norm then the slope of the coast is low hence the
territory much affected by floods, while if  small then the coast
slope is
steep and the island has little to fear from inundations).

Call the above quantity $\pi(O)$ the ``piaf'' (protection index
against floods). It is well-defined for any oval (up to sign and
anodyne choices effecting a collective change). It seems to make
also good sense in the projective context.

Next  look at all the empty ovals $O_1, \dots, O_{r_0}$ of some
smooth curve $C_m$, splits them in even and odd, and look which of
both collections have the highest piaf. Depending on this
knowledge, an appropriate choice of $\varepsilon$ create a
variation of the curve with smaller inner area $\alpha$ .

(D.3) What to do exactly when both (even and odd) piafs are equal?
Can we avoid this just by slight perturbation, i.e. is there
always a small perturbation making them different? (Perhaps
Petrovskii thought about such questions\dots)

Assume now that (D.3) can be overcome. Then the functional
$\alpha$ has no stable equilibrium (local minimum) and we
interpret it as a Morse function on the space $\vert mH \vert$ of
all curves. In this generality it may rather look like a Grand
canyon with big ravine when one cross a solitary node due to the
brusque change of empty oval. Two attitudes are possible:

(A.1)---either we localize $\alpha$ to one chamber of the
discriminant\footnote{We often commit an abuse of language, as we
should say one chamber residual to the discriminant. Such an abuse
is harmless like when speaking of the group of a knot, when it is
really that of its complement.} or

(A.2)---we look at the whole space of curves hoping that the
trajectories dictated by the functional never cross such ridges
(=ravine of $\alpha$).

After having introduced a natural metric on $\vert m H \vert$,
e.g. the elliptic geometry available on each projective space, we
look at the trajectories of steepest descents w.r.t. the function
$\alpha$ (area of the nonempty ovals). This is as usual obtained
by integrating the vector field ${\rm grad} \alpha $. By the above
(D.2), for almost every initial condition $C_m$ it decreases
endlessly up to reach level $\alpha=0$, which must necessarily be
a curve with solitary nodes. Could the trajectory starts
oscillating like a $\sin (1/x)$ curve without reasonable
convergence? Looks unlikely due to the algebraic nature of our
problem, but requires perhaps an argument. If the trajectory of
$C_m$ converges to a saddle point (unstable critical point of
$\alpha$) it suffices to perturb slightly the initial condition
$C_m$ (which is
allowable up to small rigid-isotopic perturbation). (In fact it is
likely that such exceptional saddles correspond precisely to
curves having the same even and odd piafs, especially if we have a
rigorous proof that there is no local minimum for $\alpha$, as we
tried to argue in Step~D.2.)

At this stage we believe the proof would be completed (no
additional difficulties) and we would conclude:

\begin{theorem} (Hypothetical!!!)
Given any (non-void) curve $C_m$ (of even degree $m$ for
simplicity), the trajectory of the gradient flow of the
empty-ovals area $\alpha$ generically converges to a curve with
solitary nodes in finite time, while the empty ovals themselves
converges to the solitary nodes. If not then it finishes its
trajectory to an unstable equilibrium and it suffices to perturb
slightly to ensure convergence toward a solitary nodal curve
(soliton for short, as compression of solitary and singleton). In
particular, CCC holds true, i.e. there is a path in the space  of
curves such that all empty ovals contract to solitary nodes.
\end{theorem}

\begin{proof}
That the extinction of  all the empty ovals occurs in finite time
merely follows from the fact that the time parameter of any
gradient flow is just the ``height'' function, here the functional
$\alpha\colon \vert mH \vert \to {\Bbb R}$ but taking value in
$[0, 4\pi]$ (where $4\pi$ is the area of the full sphere or $2\pi$
if you count this area  divided by two).
\end{proof}

Here we have looked at the empty-oval area functional $\alpha$.
What happens if we look the same functional for all ovals. A
priori the functional looks more continuous but be careful with
eversions.

[23.01.13] Metaphor.---Problems of rigid-isotopy (or large
deformations of curves) are like a video game in the sense that
there is a joystick upon which one may act by freewill by varying
the coefficients while there is in reaction a canonical picture
emerging on the screen (the corresponding real locus of the
algebraic curve conceived as an  optical object). In some sense it
is like a flight simulator (you move the ``manche \`a balais'' and
the
aircraft responds accordingly). The contraction
conjecture CCC says that using the full freedom of the joystick
one can always shrink the empty ovals simultaneously. The above
theorem states roughly that  there is some predestination, i.e.
that a very sleepy autopilot or video game player suffices
to land safely the aircraft,
while performing actually a perfect landing (all wheels touch the
ground simultaneously!). Of course in reality the autopilot in
question is very well programmed for its action is governed by a
principle of least action.
The  only little impulse required is when the aircraft arrives at
critical points (global maximum of $\alpha$ or its saddle critical
points), where some jiggling is required to perturb the initial
condition.

\subsection{Some few other applications of CCC}
\label{application-of-CCC:sec}

[23.01.13] What can be deduced from CCC? Quite a lot and alas no
so much, compare the case of $M$-sextics. Applying it to Gudkov's
$C_6$ gives a rational (genus 0) real sextic with 10 solitary
nodes equidistributed as $5$ inside and  $5$ outside the unique
oval.
Of course there is no obstruction given by the genus formula to
such a eventuality.

In general, given any curve $C_m$ with say $r$ ovals, then can be
split as $r_0\le r$ empty ovals, which can be contracted to
solitary nodes. Then applying Brusotti 1921 \cite{Brusotti_1921}
we can let them disappear all, and so appears a new generation of
empty ovals, to which the collective contraction process can be
applied again, and so on. So we can reach the empty curves (or a
pseudoline if $m$ is odd) after some few iterated contraction (as
many as the  height of the oval-graph encoding the nested
structure of the original $C_m$). This we call the height of the
curve $C_m$, denoted $h(C_m)$

This implies (as a crude estimate) that any two curves $C_m, D_m$
can be related by a rigid-isotopy crossing only $h(C_m)+h(D_m)$
times the discriminant. Of course it is not a transverse crossing
in general for our critical curves have several solitary nodes.
However by perturbing slightly we may cross the discriminant
transversally, and each initial crossing through a multi-solitary
curves implies as many intersection with $\frak D$ as there are
nodes.

So counting properly we deduce:

\begin{lemma} {\rm
(modulo CCC and the connectedness of invisible curves=CIC)} Any
two smooth curves $C_m, D_m$ of even degree $m$ can be joined by a
path
 in the hyperspace of curves transverse to the discriminant while
crossing it exactly $r(C_m)+r(D_m)$ times, where $r$ is the number
of  ovals (composing the real locus).
\end{lemma}

\begin{proof}
Applies iteratively  CCC (conjointly with Brusotti) to both
curves, to derive two curves with empty real locus. The latter are
known to form a unique chamber of the discriminant, in other words
to be rigid-isotopic by
Lemma~\ref{empty-chamber-connected-Shustin:lem}.
\end{proof}

A similar assertion holds perhaps true in case of odd degrees,
however it is still unknown whether two curves of the same odd
degree are rigid-isotopic provided their real loci reduce to a
pseudoline (compare  Viro 2008
\cite[p.\,199]{Viro_2008-From-the-16th-Hilb-to-tropical}).

This goes in the sense of showing that the contiguity graph of
chambers residual to the discriminant is a ``small world'', in the
sense that is has high connectivity and much ``consanguinity''.

How good is the above estimate? More precisely the distance
$\delta$ (or Erd\"os number) in the contiguity graph (of chambers)
is majored by $\delta(C,D) \le r(C)+r(D)$. This is fairly good as
compared to the estimate coming from the degree of the
discriminant $3(m-1)^2$, which implies $\delta(C,D)\le
3/2(m-1)^2$, or the integral part thereof, as we may always choose
the one side of the circle  hitting less many times the
discriminant.

For $m=6$, the discriminant estimate gives $\delta\le
[75/2]=[37.5]=37$, while the CCC estimate gives $\delta\le
11+11=22$. We can be  more economical by not going down to the
empty chamber but that having only one oval, which form already a
unique rigid-isotopy class by Nikulin 1979 \cite{Nikulin_1979/80}.

This raises  the following question: we know either by Rohlin's
formula (or less rigourously by CCC) that any curve with one oval
(hence of even degree) is nondividing provided $m=2k\ge 4$. Indeed
if dividing apply CCC to get a splitting $C_{2k}\to C_k\cup
C_k^\sigma$, whence by B\'ezout $C_k \cap C_k^\sigma=k^2=r$, where
$r=1$ whereas $k^2\ge 4$. Thus there is no obstruction to
rigid-isotopy given by the Klein's type between any two curves
having only one oval, and extrapolating (violently) we arrive at
the:

\begin{conj}\label{OOPS:one-oval-rigid-isotopic:conj}
{(OOPS=One oval postulation)} Any two smooth curves having only
one oval are rigid-isotopic. (``Oval'' is interpreted here in the
strong sense of a Jordan curve which is null-homotopic, hence our
curves are of even degree.)
\end{conj}

(Perhaps there is an obstruction \`a la Fiedler-Marin, but
unlikely as it seems to require a splitting of ovals, cf. Marin's
argument exposed below.) Remind also from Viro 2008
\cite{Viro_2008-From-the-16th-Hilb-to-tropical} that replacing
above ``oval'' by pseudoline is still an open problem.
Further Viro in his e-mail (dated 26.01.13 in
Sec.\,\ref{e-mail-Viro:sec}) confirmed me that this is
still an open problem
and goes  back to Rohlin.

{\it Insertion} [02.02.13].---A naive approach to this would be to
assert that any $C_{2k+1}$ may degenerate by a large deformation
to $L_1\cup C_{2k}$ a line $L_1$ union an invisible curve $C_{2k}$
of even degree. Perhaps one can demand that both curves in the
limit are transverse. Call this process a rectification of the
pseudoline. Now given two smooth curves whose real scheme consist
of a unique pseudoline, one may apply twice rectification, and
then isotope both corresponding empty curves. This could imply the
Rohlin-Viro conjecture.

In degree 6 the Erd\"os number of the graph (supremum of $\delta$)
is probably much smaller in view of
Fig.\,\ref{Gudkov-eversion:fig}, i.e. Gudkov's table with
eversions. Of course this figure posits that all logically
possible eversions are realized geometrically (safe those linking
$M$-curves and both $(M-2)$-schemes of type~I). Under this
circumstance the Erd\"os number looks hardly greater than 8, i.e.
$\delta \le 8$ universally. It seems indeed that the maximal
distance is realized by Hilbert vs. Gudkov or Gudkov vs. Harnack.
Naively on the table (without eversions) Hilbert's scheme and
Harnack's looks far apart, but using the eversion $\frac{8}{1}1
\to \frac{1}{1}8$ shows that their real distance is only 3, i.e.
$\delta (Hilbert, Harnack)=3$. This is the answer if we confine
attention to the top of the table, but of course the most distant
chambers are the $M$-curves as separated from the empty scheme
$0$. Those are at distance $11$ apart. So the correct Erd\"os
number is $\delta=11$, with this maximal distance being realized
thrice (empty vs. Hilbert, Gudkov, Harnack respectively). Morally
eversions do not shorten the vertical distance, and we have proven
the following (modulo Conjecture~\ref{eversion-and-other
surgeries:conj}).

\begin{prop} {\rm (Semi-conjectural)}
\label{Erdos-number-of-sextics=11:prop}
The Erd\"os number of the contiguity graph of sextics is actually
equal to the Harnack bound $M=11$.
\end{prop}

It is tempting to wonder if it so in general. Perhaps there is
some little chance to answer this without having to work out the
exact  rigid-isotopy classification in each degree (an
insurmountable task!?). The trick would be that eversions collapse
sufficiently horizontal distances, so as to make only the vertical
chain the only plausible candidate for maximizing $\delta$.

As to the above OOPS conjecture
(\ref{OOPS:one-oval-rigid-isotopic:conj}), one could of course
also imagine a dynamical proof. The whole task reduces to finding
the right functional. Heuristically the obvious attractor ought to
be a circle (possibly multiple). This inclines to look at the
isoperimetric functional looking the length squared divided by the
area of the unique oval.

More pragmatically,

\begin{Pseudo-Theorem} {\rm (but apparently validated by Viro, cf.
comments right-after the proof)}.---It seems clear that CCC
(and of course the weaker formulation CC)
implies the one oval postulate (OOPS).
\end{Pseudo-Theorem}

\begin{proof} (pseudo, of course!)
Indeed take any two curves having only one oval. By CCC (CC
suffices) each of them can be shrunk to a solitary node. By
Brusotti there is slight perturbations to empty curves. The latter
can be linked by a path, by virtue of the connectivity of the
empty chamber (Lemma~\ref{empty-chamber-connected-Shustin:lem}).
Since the empty chamber of the discriminant is a manifold with
corners (like probably any other chamber by the way) which
topologically is a manifold with boundary there is such a path
staying close to the boundary. (Actually the boundary of the empty
chamber has faces consisting of uninodal solitary curves with one
isolated real point.) Pushing this path slightly outside the empty
chamber would give the required isotopy. However there is a
serious difficulty, if our path meets another wall of ${\frak D}$
outside the empty chamber. However by genericity of this path (as
avoiding sets of codimension 2) we may assume this crossing to be
a transverse one of a wall which must keep $r=1$ constant since we
stay in the vicinity of the empty chamber with $r=0$. Note of
course by CCC as applied to quartics for instance that chambers
with higher $r$'s are also contiguous to the empty chamber, yet
the are like cubes hitting the empty chamber imagined as a cube at
some vertices of codimension 2. So we have some wall crossing
keeping $r=1$ constant, which as a Morse surgery must correspond
to an eversion. But this means that there is an eversive wall
falling down to the boundary of the empty chamber like a tripod.
This looks incompatible with the  local structure of algebraic
sets? Maybe there is  a more direct argument in ${\Bbb R}P^2$.
Another obstruction comes also from Brusotti's description of the
discriminant as branches with normal crossing. All this is very
confuse, we confess.
\end{proof}

[30.01.13] As kindly informed by Viro (cf. his e-mail, dated
26.01.13 in Sec.\,\ref{e-mail-Viro:sec}), it seems that the
implication CC$\Rightarrow$OOPS causes no problem. It would be
nice to write down  complete details. Perhaps this simply follows
from the fact that inside the empty locus the discriminant has
real codimension~2 (cf. Lemma
\ref{invisible-discriminant-codim-2:lem}). Hence there is no wall
inside it, hence no wall outside it by the ``implicit function
theorem''.

{\it Insertion.} [01.04.13] If CC$\Rightarrow$OOPS causes no
troubles it would be interesting to extend the method to higher
schemes (having more ovals than one).

\begin{conj}\label{URS:conj} (URS)=(Unnested rigidity
speculation).---The unnested configurations (which are of type~II
provided $r<k^2$ or even $r\neq k^2$, as shown by Rohlin's
formula) are always rigid (i.e. any $2$ curves representing the
unnested scheme are rigid-isotopic provided $r\neq k^2$). When
$r\neq k^2$ it could be that the type is the sole obstacle to
rigid isotopy.
\end{conj}

This holds true in degree 6 by Nikulin's theorem
(\ref{Nikulin:thm}). Further it looks hard to disprove this by the
Fiedler-Marin method as we lack a canonical choice for the
fundamental triangle. Finally, it could be that the same argument
as above shows that CCC implies URS.

{\it Sketch of proof that CCC$\Rightarrow$URS}.---Contract all
ovals simultaneously to land in the connected empty chamber of
curves without real points. So given 2 curves connect them by
respective contractions to the empty locus, and therein by a path
of invisible curves. The hard part is then to push this path by a
small deformation again in the visible locus, and this in such a
way that we never meet the discriminant. This looks feasible as
the (closured) empty locus seems to be a bordered manifold, with
connected boundary (being essentially fibred over $\RR P^2$ via
assignment of the unique solitary node). So the joining-path in
the empty locus can be pushed on the boundary and then further
inside the chamber. Of course it remains the difficulty of
ensuring that we do not meet other nappes of the discriminant. As
above the loose argument is that since the discriminant as
codimension 2 inside the empty locus, it will appear outside along
the same dimension, and not effect any separation. This would
supply the required rigid-isotopy between our pair of unnested
curves. Of course it is essential to assume $r\neq k^2$ since
otherwise there is curves of both types, yet this condition did
not as yet appeared frankly in our argument, which is far from a
serious proof. Of course one could expect that for $r=k^2$ the
type is sole additional obstruction.

\subsection{Looking around (in vain?) for counterexamples to CCC
(=collective contraction conjecture)}

\begin{figure}[h]
\hskip-1.2cm\penalty0 \epsfig{figure=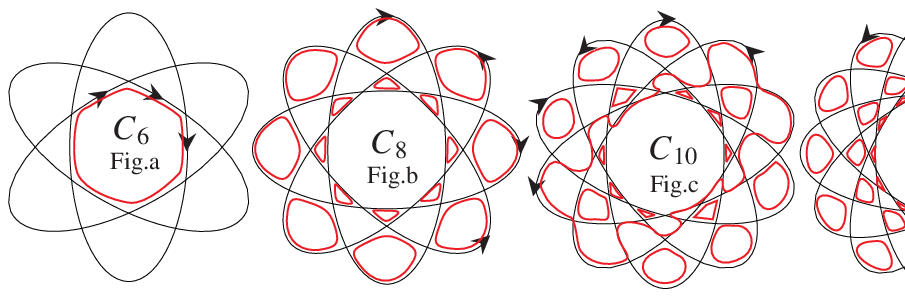,width=147mm}
\vskip-5pt\penalty0
  \caption{\label{CCCRoses:fig}%
  Constructing dividing curves without nesting as roses} \vskip-5pt\penalty0
\end{figure}

[20.01.13] An a priori easier game is to test if CCC
(\ref{CCC:conj}) has really some chance to be true. As usual
experimentation is required. We may first consider a curve $C_8$
arising through perturbation of 4 ellipses rotated by $180/4=45$
degrees (cf. Fig.\ref{CCCRoses:fig}b).

The smoothing being compatible with orientations, this curve $C_8$
is dividing (Fiedler). It has $r=16$ unnested ovals (cf.
Fig.\ref{CCCRoses:fig}b). Assume it shrinkable via CCC, it will
degenerate and decomposes as $C_8 \to C_d \cup C_d^{\sigma}$,
where $C_d$ has genus $p=\frac{g-(r-1)}{2}=\frac{21-(16-1)}{2}=3$,
hence of degree $d=4$ (this can be inferred more directly via the
degree of the degeneration). The relation $d^2=r$ is verified and
there is no numerical obstruction to CCC. Note in contrast that in
the same construction for 3 or 5 ellipses, we cannot arrange all
ovals unnested while smoothing in a sense preserving way. Look and
see! (Figs.\,\ref{CCCRoses:fig}a and c).

The construction of our $C_8$ generalizes whenever the degree
$m=4\ell$ is a multiple of four. Indeed rotate an ellipse by $\pi/
\ell$. Orient the ellipses ``alternatively'' and smooth in a
 sense preserving way. The resulting curve has $r=4 \ell^2$ unnested
ovals (as easily counted by extrapolating the figures $C_8$ and
$C_{12}$ of Fig.\,\ref{CCCRoses:fig}, while noting that the ovals
in $\ell=m/4$ couches containing each $m=4\ell$ ovals). Shrinking
$C_m$ via CCC gives two curves of genus $p=\frac{g-(r-1)}{2}$,
where $g=\frac{(m-1)(m-2)}{2}=(4\ell-1)(2 \ell-1)$. Hence
\begin{align*}
p
=[(4\ell-1)(2 \ell-1)-(4 \ell^2-1)]/2
&=[(2\ell-1) [(4\ell-1)-(2\ell+1) ]]/2 \cr
&=[(2 \ell-1) (2\ell-2)]/2,
\end{align*}
so that the half has degree $2\ell$, and intersects its conjugate
in $4\ell^2$ points, which is precisely $r$. This little numerical
miracle implies
an absence of numerical obstruction to CCC via our rosewindows
constructions. Of course, all the above computation can be
shortcuted by noticing a degeneration $C_{m=4\ell} \to C_d \cup
C_d^{\sigma}$, where $d=2\ell$ necessarily (since the conjugation
$\sigma$ preserves the degree of equations as it just acts upon
the coefficients).

As yet the center of rotation was chosen inside the ground
ellipse. Another series of picture arise when choosing it outside
instead (but sorry this is probably not the right thing to do).
Let us instead manufacture the  dividing $C_6$ with $r=9$ unnested
ovals as on Fig.\ref{CCCRoses2:fig}a. Next we tried to construct a
$C_{10}$ with 25 unnested ovals but failed somewhat. Note that
this is not obstructed by Arnold's congruence $\chi=p-n=k^2 \pmod
4$, as both sides are 25.

\begin{figure}[h]
\hskip-1.2cm\penalty0 \epsfig{figure=CCCRoses2.eps,width=147mm}
\vskip-5pt\penalty0
  \caption{\label{CCCRoses2:fig}%
  Constructing dividing curves $C_{4\ell+2}$ without nesting
  as roses fails, but as ``hyperbolic asters'' works} \vskip-5pt\penalty0
\end{figure}

Assume there is such a $C_{10}$ with real scheme $25$ (Gudkov's
notation) which is dividing. On applying CCC we find a splitting
$C_{10} \to C_5 \cup C_5^{\sigma}$. The genus $g=9\cdot 8/2=9\cdot
4 =36$, hence $p=[g-(r-1)]/2=(36-24)/2=6$ which is indeed the
genus of a quintic. So no obstruction on this side. Note also that
Rohlin's formula $2(\Pi^+ -\Pi^-)=r-k^2$ implies no obstruction,
since $\Pi^{\pm}=0$ (no nesting) and so $r=k^2$ with $k=5$.

After our stupid trials we find ultimately the right ground
configuration of ellipse as Fig.\,\ref{CCCRoses2:fig}y, which is
dividing with $25$ ovals. It is clear now how to extend this in an
infinite series of curves $C_{4 \ell +2=2( \ell +1)}$ with
$(\ell+1)^2$ ovals lying outside each other, i.e. the scheme
$(\ell+1)^2$ in type~I. Fig.\,\ref{CCCRoses2:fig}z gives the case
of a $C_{14}$.

All this pictures are
pleasant, yet they do not help at all to corrupt the conjecture
CCC. Of course the latter has some chance to be true especially if
the Itenberg-Viro conjecture
(\ref{Itenberg-Viro-contraction:conj}) is true.
A dividing curve without nesting has to satisfy $r=k^2$ by
Rohlin's formula, and therefore applying CCC gives a degeneration
$C_{2k}\to C_{k} \cup C_{k}^{\sigma}$ yielding no chance to
corrupt B\'ezout. Also there is no chance to corrupt Klein
$p=[g-(r-1)]/2$. Indeed $g=(2k-1)(2k-2)/2$, hence
\begin{align*}
p
=[(2k-1)( k-1)-(k^2-1)]/2
&=[(k-1) [(2k-1)-(k+1) ]]/2 \cr
&=[(k-1) (k-2)]/2,
\end{align*}
so that the half has degree $k$, as it should.

\subsection{CCC versus Brusotti: large deformations vs.
small perturbations}

[20.01.12] Philosophically, it seems that such contraction
conjectures (Itenberg-Viro (\ref{Itenberg-Viro-contraction:conj})
or our collective version thereof (\ref{CCC:conj}))---if they turn
out to be true by some lucky
stroke---incarnate a sort of large deformation principle
illustrating once more the perfect graphical flexibility of
algebraic curves despite their intrinsic rigidity. In some sense
this is an
avatar in the large of the small perturbation method \`a la
Pl\"ucker-Klein-Harnack-Hilbert-Brusotti-Viro. Hence it seems
indeed (as Viro advocated) being of some primary interest to
establish the contraction conjectures.
Flexibility in the small gives rise to the perturbation method
which is primarily a method of construction, whereas flexibility
in the large (contraction principles) implies as a byproduct
prohibitions (as we superficially experimented at the beginning of
Sec.\,\ref{CCC:sec}
and more convincingly  because Itenberg's contraction theorem for
sextics re-explain all Gudkov-style prohibitions by reduction to
the RKM-congruence). At this stage we feel some big duality: local
versus global and constructions versus prohibitions
(to be or not to be). What would be the net impact of the
contraction principle for Hilbert's 16th problem (in the extended
sense of high degrees)? Somewhat optimistically it would reduce
the whole task (or rather adventure) to a combinatorial video game
best suited for machines. So exaggerating slightly,  the
contraction conjecture seems quite close to reveal the ultimate
secret of the whole problem. Of course some supplementary large
deformation principles ought also to complete the picture, e.g.
certain permissible eversions compatible with B\'ezout, and more
generally the full morphogenesis of all algebraic Morse surgeries.
If all this is available, the video game solving Hilbert's 16th
problem would
show in real time all the possible perestroikas which the real
loci of projective curves of some fixed degree can undergo,  while
dragging at free will the
joystick in the parameter space.

As sketched in the previous Sections~\ref{CC-via-dynamics:sec} and
\ref{CC-via-Riemann:sec}),  in order to prove the contraction
conjectures (CC or CCC), we could either imagine a dynamical proof
via orthogonal trajectories, hence akin to Morse theory, or a
direct intervention of conformal geometry \`a la Riemann (albeit
our implementation failed seriously). Whatsoever the exact details
it is quite likely that the proof of CC or CCC will employ the
calculus of variation in the large over which practically every
deep geometrical theorem is based upon (from the brachystochrone,
to the Dirichlet principle, via the Riemann mapping theorem up the
recent solution of the Poincar\'e conjecture via the Ricci flow.)

Note finally another very modest piece of evidence in favor of
CCC. Remember that for small perturbations \`a la Brusotti, there
is complete freedom to smooth the nodes of a plane curve with
normal crossings (compare e.g. Gudkov 1980
\cite{Gudkov_1980/80-Brusotti}) in the sense that all crossings
may be smoothed away or some may be conserved. By analogy CCC is
just the case where all empty ovals are contracted simultaneously,
while Itenberg-Viro's CC is just the contraction of a single empty
oval. Of course there ought to be the full panoply of intermediate
contractions.

\subsection{CCC implies CC (i.e. Gabard stronger than Itenberg-Viro)}

[21.01.13] Let us now observe that CCC implies CC, just via
Brusotti's theorem (1921 \cite{Brusotti_1921}):

\begin{lemma}\label{CCviaCCC-Brusotti:lem} Suppose given a
collective contraction of all the empty ovals of a smooth real
curve $C_m$, then it is possible to construct all partial
contractions via Brusotti. In particular if CCC holds true then so
does the Itenberg-Viro contraction conjecture {\rm
(\ref{Itenberg-Viro-contraction:conj})}.
\end{lemma}

\begin{proof}
Let $C_m$ be a smooth real curve. W.l.o.g. let us assume it having
some empty ovals, say $r_0\ge 1$ many (otherwise the curve just
reduces to a pseudoline or to empty real locus). By CCC there is a
path $c\colon [0,1] \to \vert mH \vert$ such that $c(0)=C_m$ and
$c(t)\in \frak D$(=the discriminant) only for $t=1$ where $c(1)$
is a nodal curve with solitary nodes only ($r_0$ many). By
Brusotti's theorem the neighborhood of the nodal curve $c(1)$
consists of $r_0$ ``falde analytiche'' (=analytic branches or
better {\it nappes\/}) meeting transversally at $c(1)$. Further
each of those nappes corresponds to the conservation of some node
in the vicinity. In other words the chamber of $C_m$ looks like
manifold-with-corner near the nodal curve $c(1)$, locally
diffeomorphic to ${\Bbb R}^N$ with $r_0$ many distinguished
hyperplanes of coordinates. It is now plain how to construct all
other contractions, in particular all the ones of Itenberg-Viro CC
contracting just a single empty oval (compare
Fig.\,\ref{CCCBrusotti:fig}).

\begin{figure}[h]
\centering
\epsfig{figure=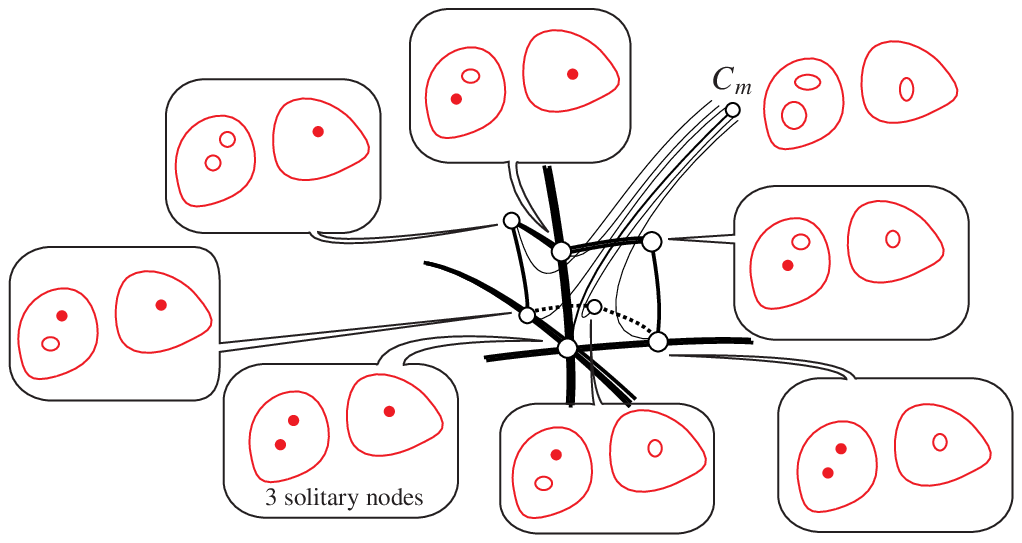,width=122mm} \vskip-5pt\penalty0
  \caption{\label{CCCBrusotti:fig}%
  Construction from the maximal (collective) contraction
  all the partial contractions via Brusotti's theorem on
  the independence of smoothing nodes} \vskip-5pt\penalty0
\end{figure}

\end{proof}

So if optimistic about the truth of CC, it may even seem that CCC
is a fairly reasonable angle of attack. Of course one then needs a
good functional (e.g. the area of all empty ovals).

Alas, it is less evident that CC implies CCC, as we suggested at
the beginning of the investigation. This could be slightly easier
if there is a contraction principle extended  to solitary nodal
curves while keeping the ``solitons'' in place.

Naively, one may dream of a contraction principle effecting a
retraction of a whole chamber (past the discriminant) to its
boundary. Yet  compact bordered manifolds never retract to their
boundaries (as shown by homology mod 2, as we learned from J.-C.
Hausmann). This is evidently no obstacle against the contraction
conjectures, for the retraction may be undefined on  small loci,
as since we work up to isotopy such equilibrium points can be
avoided. More lucidly nobody ever asserted that the contractions
should depend continuously on their initial point(=curve). Imagine
as a very naive picture, the chamber as being a disc with a radial
projection upon the boundary. Then it is undefined on the center
of the disc but this is not a problem for perturbing it slightly
it will get mapped somewhere. The whole analogy with retraction of
bordered manifolds is not extremely pertinent as in general the
chamber will have a boundary consisting not merely of faces
touching the empty chamber. Hence under a retraction a curve close
to the discriminant could first coalescence 2 ovals instead of
shrinking one empty oval.

\begin{figure}[h]
\centering
\epsfig{figure=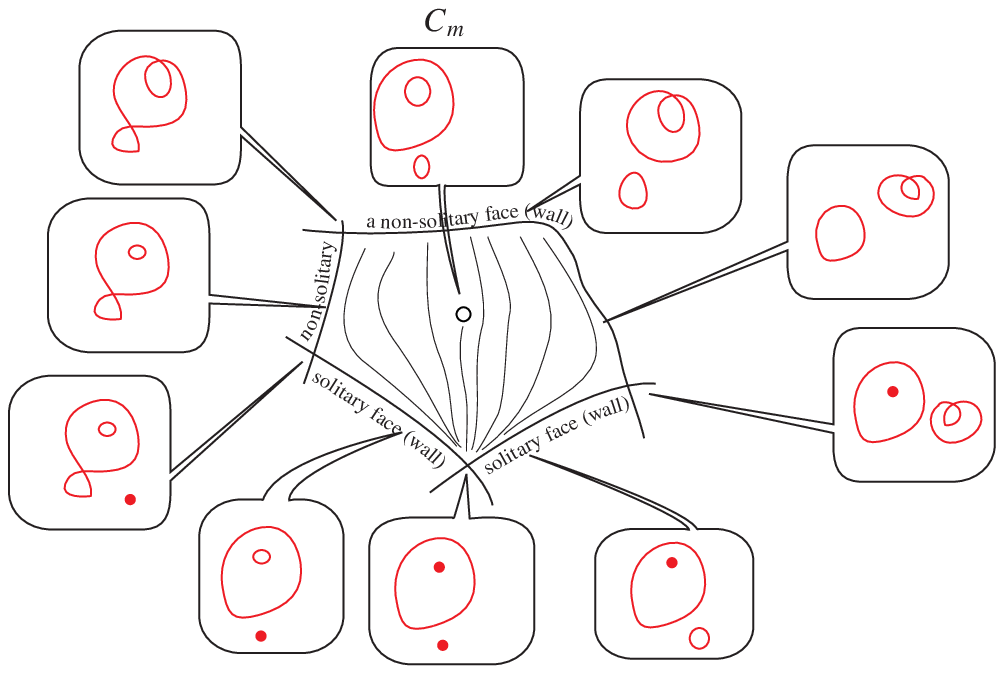,width=122mm} \vskip-5pt\penalty0
  \caption{\label{CCCflow:fig}%
  Schematic portrait of the flow implementing CCC
  (collective contraction conjecture)} \vskip-5pt\penalty0
\end{figure}

To show CC or even CCC we could employ a dynamical system
(continuous flow) spreading nearly all curves toward the boundary
of this chamber at curves having solitary nodes. This could occur
as the orthogonal trajectories of some functional. (Another idea
would be to look at the Green function of the chamber yet the
Green's lines, streamlines of the flow, would often finish at
curves with non-solitary nodes.) Our flow should be strongly
attracted by the multi-solitary nodal curves manifold, whereas all
other walls corresponding to non-isolated singularities have to
repulse the flow. A qualitative picture is given on
Fig.\,\ref{CCCflow:fig}.

This picture proves nothing, safe maybe the absence of topological
obstruction (a priori) to find such a flow (especially if it is
allowed some equilibriums when the chamber has complicated
topology). Again to get a proof it is likely that one should
consider the gradient flow of some real-valued functional on the
chamber. This could be the area or length of all empty ovals (as
measured on the round metric of $S^2$ double covering ${\Bbb
R}P^2$). Note  that CCC implies a little technical simplification
over our previous pseudo-proof of CC, where we were troubled by
marking one oval.

As already discussed, the main obstacle occurs if our function has
some global attracting basin inside the chamber preventing us to
reach the desired multi-solitary nodal curve shrinking all empty
ovals. This is basically the sole difficulty yet it looks quite
insurmountable.

At least for the area or length functional, we saw no obvious way
to produce small variations diminishing the ``energy''. Perhaps
there is some more clever (projective) invariants like degree of
roughness \`a la Gudkov (cf. e.g. Gudkov 1974
\cite{Gudkov_1974/74}). For instance one could look at the largest
(quadratic) cones in ${\Bbb R}^3$ which can be nested inside the
ovals, and the corresponding area intercepted on the unit sphere.
This is another functional measuring the conical area of the empty
ovals. Can we show that this functional is ``good'', i.e. no sink
inside the chamber? Naive idea: trace inside each oval some
maximal ellipse and try to deform the $C_m$ along suitable
multiples of those ellipses.

Now the problem looks reduced to a fantastic game probably
only  soluble by such authorities as Andronov, Leontovich, Gudkov,
etc. mixing the qualitative theory of differential equations with
that of algebraic curves. So it is truly a Poincar\'e-Hilbert
Verschmelzung(=fusion in Klein's prose) which seems demanded to
settle CCC (or its avatar CC).

[03.04.13] Of course, it also safe to say pessimistic and expect
maybe that Shustin's disproof of Klein-vache in degree 8 also
implies a disproof of CC. We do not repeat our vague strategy for
 this, but refer to Sec.\,\ref{Challenging-open-prob:sec}.


\section{Problems of rigid-isotopy}

[09.02.13] Two (real, plane) curves are said to be {\it
rigid-isotopic} if one can pass from one to the other by
continuous deformation of the coefficients of the defining
equations which avoids the discriminant. This involves again the
paradigm of {\it large deformations\/}  like the contraction
conjectures discussed in the previous section.
%
In fact there should be some direct connections between both
topics.

As a rule  very little is known about the phenomenon of rigidity.
All what is trivial is that any topological characteristic
persists during a rigid-isotopy, so for instance the {\it real
scheme\/} (i.e. the isotopy class of $C_m(\RR)\subset \RR P^2$) as
well as Klein's type~I, II measuring the situation of the curve in
its complexification. Such invariance were intuitively clear since
the era of Schl\"afli, Zeuthen, Klein 1876
\cite{Klein_1876_Verlauf}, but requires perhaps Ehresmann's lemma
that a locally trivial fibering over a contractible
(paracompact$\approx$metric) base is globally trivial.
Paracompactness is essential as shown by the simply-connected
(indeed contractible) {\it Pr\"ufer surface\/} (cf. works by
Pr\"ufer 1922, Rad\'o 1925 \cite{Rado_1925}, Calabi-Rosenlicht
1953, Spivak's book on Diff. Geom., Vol.~I, Appendix, Baillif).

Up to degree $\le 4$ the real scheme (or Klein's types I/II)
suffices to
encode the rigid-isotopy class as knew Klein 1876
\cite{Klein_1876_Verlauf}, building over Schl\"afli and Zeuthen's
works. In degree 5, and 6, the same real scheme plus Klein's type
suffices to ensure a rigid-isotopy. This spectacular result is
joint work of Nikulin, with the collaboration of Kharlamov
building over two pillars, namely:

(1) the Gudkov-Rohlin census solving Hilbert's 16th problem (for
sextics) while revitalizing the earlier Riemannian conceptions of
Klein about the complexification and,

(2) the theory of K3 surfaces (Torelli, etc.).

The situation changes drastically from degree 7 upwards as shown
by the Fiedler-Marin trick using a locking triangle which consists
of $3$ B\'ezout-saturated lines, hence which cannot be crossed by
ovals during a rigid-isotopy. Here the basic idea is that if one
can associate to a curve in some canonical manner an auxiliary
curve called the lock then the distribution of ovals past the lock
is rigid-isotopically invariant. A typical example in degree 7 is
the lock consisting of 3 lines through the inner ovals of a curve
having the scheme $\frac{3}{1}\ell J$ for some $\ell$. When $\ell
2$ the fundamental triangle can separate in different ways the
$\ell$ outer ovals, and  curves having different splittings past
the lock will not be rigid-isotopic.

This obstruction to rigid-isotopy requires several ovals, and does
not seem suited to curves with few ovals where the rigidity
problem looks fully open. It is presently very unclear
which sort of scenario is to be expected. For instance it is still
undecided for $m\ge 7$ whether curves having only one real circuit
are always rigid-isotopic. As informed by Viro, this is at least
for odd degrees a (confidential) conjecture due to Rohlin.

More generally we posit the following speculation (not that we
strongly believe in it, but just as a way to confess our
ignorance):

\begin{conj}\label{LARS:conj}
(LARS).---Curves of degree $m$ with less than $DEEP+2$ real
branches, where $\Delta(m)=DEEP(m)=[(m+1)/2]$ is the number of
components of the deep nest of degree $m$, are always
rigid-isotopic provided they have the same real scheme.
\end{conj}

The basic motivation for this conjecture is that below altitude
$r\le \Delta+1$ all schemes are of type~II except the deep nest
which is of type~I (by total reality under a pencil of lines).
This follows from Rohlin's formula (\ref{Rohlin-formula:thm}),
especially its corollary known as Rohlin's inequality
(\ref{Rohlin's-inequality:cor}), as well as Klein's congruence
$r\equiv_2 g+1$ forcing dividing curves to have their numbers of
real circuits $r$ jumping by quanta of 2 units, hence we gain one
type~II level right above the deep nest. Further the deep nest is
known to be rigid by Nuij 1968 \cite{Nuij_1968}.

Apart form those elementary facts, we have little evidence for
this ``low-altitude rigidity speculation'' (LARS), except
suspecting that if the assertion is true it will use a geometric
flow permitting a degeneration to curves of lower orders after
splitting off a line or a conic. For instance given a curve of odd
degree it is tempting to look at the flow shortening the length of
the unique pseudoline. Orthogonal trajectories of this functional
should abut to a curve splitting off a line, which after all is
the shortest pseudoline. This could be the basis of a grand
inductive process reducing the rigidity of curves with few
branches ($r\le \Delta(m)+1$) to that of curves of lower orders.

It seems that much can be explored along this line of geometric
flows, somewhat reminiscent of say M\"obius 1863, Poincar\'e,
Morse, etc., up to Perelman's proof of Poincar\'e's conjecture,
except that in our case the dynamics  lives merely on a
finite-dimensional manifold (the hyperspace of all algebraic
curves of some fixed degree).

Another
source of rigidity
comes from the empty scheme, which is rigid. In fact I started to
doubt about this issue, until Shustin kindly
remembered me the following simple  argument. If we have two
curves with empty real locus ({\it invisible curves\/} for short),
then after choosing equations of the same sign, the linear
deformation $(1-t)P+tQ$ will connect both curves while conserving
the same sign (provided $0\le t\le 1$), hence producing a path of
invisible curves. However it is not a priori (nor a posteriori!?)
evident that our path avoids curves with singularities, which
could occur in imaginary conjugate pair. Hence the complete proof
seems to use the fact that empty curves with singularities form a
locus of codimension 2, since there are two conjugate nodes
generically.

Once rigidity of the empty scheme is known, rigidity of curves
with one oval should follow simply by contracting the one oval and
letting it then disappear. This gives the basic connection with
the former section (which is primarily a remark of Viro). Also as
we remarked earlier, the stronger version CCC of the contraction
conjecture could accomplish stronger rigidity results like URS,
cf. (\ref{URS:conj}).

{\it Insertion} [03.04.13].---Conversely if the one-oval scheme
(unifolium) is rigid then it suffices to contract the Fermat curve
(of even degree) to establish CC, but alas only for this unifolium
chamber. Incidentally the validity of CC even for the unifolium
scheme is not completely obvious, for taking the linear pencil
between such a curve and an empty one leads (after assuming
general position w.r.t. the discriminant, i.e. transversality) to
a sequence of Morse surgeries a priori much more complex than just
the death of the oval.

[10.02.13] Another question of didactic interest is to study the
interplay between Fiedler-Marin locking method and  Nikulin's
rigid classification. For sextic schemes of the form
$\frac{3}{1}\ell$ we have an obvious lock given by the 3 lines
through pairs among the 3 deep inner ovals. Each such line cuts
twice the ovals it visits and twice the (nonempty) surrounding
oval, hence is B\'ezout-saturated. Further this fundamental
triangle is canonically assigned to the configuration in the sense
that the position of 3 points in the insides of the ovals is
parametrized by the 3rd symmetric power of a cell which is a
contractible space. Hence the distribution of the outer ovals past
the deep (Bermudian) triangle is invariant under rigid-isotopy.
This adumbrates a strategy toward corrupting Nikulin's theorem,
but in reality the latter rather implies an invariance of this
Bermudian distribution of outer ovals for all curves having the
same real scheme (and the same type in Klein's sense). It is
therefore of interest to determine this outer distribution past
the fundamental triangle for some specific curves as it will imply
the same for all isotopic curves. This question is elaborated in
Sec.\,\ref{Nikulin-corruption:sec}.

{\it Insertion} [03.04.13].---All this problematic went in
decrepitude after an illuminating message of Le~Touz\'e (cf.
Sec.\,\ref{LeTouze:sec}) yielding a conceptual explanation of why
it is impossible to corrupt Nikulin via Fiedler-Marin. The reason
is a simple chromatic law for conics passing through $5=3+2$
points with 3 of them black-colored (situated or defining a
triangle), while the location of the 2 remaining points (white
colored) past the triangle will determine how the sequence of 5
points distributes on the conic interpolating them. When the 2
white-points belongs to different component of the (black)
triangle the distribution will be dichromatic in the sense that
the 2 white points are not standing nearby, but  separated by
black points ($1$ or $2$ depending on the path chosen on the
topological circle underlying the conic). Applying this lemma on
conics to the  above setting, shows that all ovals of a sextic
enlarging the scheme $\frac{3}{1}$ are necessarily not separated
by the deep triangle, for otherwise we can trace a conic with 4
transitions black-and-white (i.e. inside-vs.-outside of the
nonempty oval) with therefore $5\cdot 2+4=14>12=2\cdot 6$ real
intersections violating B\'ezout.

A last phenomenon is the rigidity of the deep nest established in
Nuij 1968 \cite{Nuij_1968} or Dubrovin 1983
\cite{Dubrovin_1983/85}. This seems connected with Ahlfors total
reality, since the deep nest is totally real under a pencil of
lines. Extrapolating the Nuij-Dubrovin rigidity one can speculate
that curves (or schemes) totally real under other pencils  are
likewise rigid. (A real scheme is {\it rigid\/} if all curves
belonging to it are rigid-isotopic.) For instance the scheme of
degree 8 consisting of 4 nests of depth 2 is totally really under
a pencil of conics and thus could be rigid. Note however that
total reality in the abstract sense which is actually (by Ahlfors
theorem) synonymous to ``type~I'' is not sufficient to ensure
rigidity as exemplified by Marin's construction of two isotopic
$M$-septics, yet not rigid-isotopic (cf. Fig.\,\ref{Marin:fig}
below). Therefore if there is any connection between total reality
and rigidity it must be a more subtle one. This theme is explored
in Sec.\,\ref{Nuij-Dubrovin-extended:sec}, but we lack any serious
result presently.

Problems of rigid-isotopy amount studying the residual components
past the discriminant $\disc$ which is a hypersurface of degree
$3(m-1)^2$ in the hyperspace of all curves of degree $m$. Call
such components,  {\it chambers} ``of'' the discriminant. When
$m\le 6$ virtually everything is known, e.g. there are precisely
64 chambers of sextics by Nikulin's theorem built upon the
Gudkov-Rohlin census (cf. Fig.\,\ref{Gudkov-Table3:fig}). A myriad
of questions occur which are hard to handle systematically. For
instance how many chambers in function of $m$? Is there an
universal upper bound on the number of chambers residual to a
hypersurface in function of its degree and  dimension. This seems
perhaps accessible via a conjunction of
Harnack-Klein-Smith-Thom-Milnor and Jordan-Brouwer separation
(plus Phragmen?). Asymptotic results in this sense were studied by
Kharlamov-Orevkov.

One would like to describe the contiguity graph between chambers
where edges label Morse surgeries while crossing the discriminant
transversally along a principal stratum of codimension 1
(so-called {\it walls\/}). One can also investigate the topology
of the varied chambers. Here one tool is the monodromy
representation encoding  how ovals permute when the curve is
travelled along a loop in the given chamber. This and other issues
is the object of next section, which is probably not extremely
relevant to our main topic of the Ahlfors map, yet pleasant for
its own. It can be left with loss of continuity.

\subsection{The topology of chambers, symmetry, monodromy
and transmutation (Kharlamov 1980, Itenberg 1994)}

[10.02.13] To each real scheme is attached a (Zeuthen)-Hilbert
(multi-)tree (``forest'') with vertices the ovals and with edges
whenever there is a nesting. Since any oval is immediately
enveloped in at most one other oval this forest looks like a
forest of pines (or a mushroom if you prefer). Hence, it is a
directed set branching only downwards. The monodromy acts on this
tree respecting its combinatorial structure. So for instance the
deep nest is a ``naked'' tree having only a trunk but no branches.
The automorphism group of this trivial tree is trivial, and so
must be the monodromy representation. There is no obstruction to
the deep-nest chamber having a simply-connected topology, and  we
can conjecture it to be simply-connected.

What about the empty chamber $E$?  Define the {\it invisible
locus\/} $I$ as the set of all curves with empty real locus. The
empty chamber is $E=I-\disc$. A simple argument (detailed in the
sequel) shows that $I$ is simply-connected and even contractible.
Another simple argument based on Brusotti shows that $\disc \cap
I$ is nonempty (when $m\ge 4$). Further $I\cap \disc$ has real
codimension 2 in the hyperspace of curves $\mH$ (or in $I$) and so
is like a knot (possibly with singularity). In any case, it seems
to follow ($I$ being noncompact) that the fundamental group
$\pi_1(E)$ is non-trivial despite triviality of the monodromy.
Perhaps $E$ is an Eilenberg-MacLane space $K(\pi, 1)$, i.e.
aspheric. Can we compute $\pi_1(E)$ as a function of $m$? What
about $m=6$ or even $m=4$?

Likewise  for the deep chamber $D$ it seems hasty to expect
simple-connectivity from trivialness of the monodromy. Consider
for  instance the G\"urtelkurve $C_4$ quartic (with 2 nested
ovals) and assume it a very symmetric perturbation of 2 transverse
ellipses rotated by $90$ degrees. Assume further the existence of
say a symmetry $\tau$ about the vertical axis. Since the group
$G=PGL(3, \RR)$ is connected we can connected the identity to
$\tau$ by a continuous path $c$. This $c$ induces a loop $\gamma$
in the space of quartics from $C_4$ to itself. As $C_4$ belongs to
the deep chamber $D$ upon which the group $G$ acts, it makes sense
to ask whether $\gamma $ is trivial or not in $\pi_1(D)$. For $c$
we may choose the path in $SO(3)$ given by $180$ degrees gyration
about the vectorial line of $\RR^3\ni (x,y,z)$ parallel to the
axis of symmetry of $C_4$ (viewed in the affine chart $z=1$).
(Warning actually it seems that the line orthogonal to that is
required!)

Actually choosing any path $c$ from $id$ to $\tau$ in $G=PGL$, its
image in $D\subset \mH$ is a loop $\gamma$  likely to be  not
null-homotopic. Alas the map from $G$ to the orbit of $C_4$ is not
really a covering, the argument looks a bit sloppy. A more
convenient way to argue is to consider the double cover of the
deep chamber $D$ by polarized curves, i.e. with a preferred half
of the underlying orthosymmetric  Riemann surface. Polarizing
amounts specifying a complex orientation \`a la Rohlin (by taking
the oriented boundary of the preferred half w.r.t. the canonical
orientation induced by the complex structure). The loop $\gamma$
based at $C_4$ lifts---w.r.t. the polarized cover---to a {\it
non-closed\/} path, since $\tau$
 exchanges both halves of the complexification. Imagine
indeed the symmetric surface underlying the G\"urtelkurve as a
pretzel of genus 3 with 2 ovals acted upon by an involution
($\tau$) with 4 fixed points then it must necessarily be a
rotation by a half-twist about a line in 3-space perforating the
ovals in 4 points.

It is clear that this argument extends to all deep nests and we
obtain the:

\begin{lemma}
For any integer $m$ (odd or even
do not matter) the chamber of the deep nest (alias deep chamber)
is not simply-connected.
\end{lemma}

\begin{proof} Any curve
in the deep chamber $D$ is of type~I (Klein's orthosymmetry) since
there is a totally real pencil of lines. (This is the trivial
sense of Ahlfors theorem so-to-speak.) We consider for each plane
orthosymmetric curve the two possible ways to paint one half of
the curve in black, and call the corresponding painted object a
polarized curve (or Riemann surface). If $O$ is the union of all
orthosymmetric chambers, we have a natural way to topologize the
space $O_2$ of all polarized curves to turn it in a double cover
of $O$, the orthosymmetric locus. In particular we have a double
cover of the deep chamber $D_2\to D$.

Any member of the deep chamber admits a representative $C_m$ with
a mirror involution $\tau$ given as $(x,y)\mapsto (-x,y)$ in
affine coordinates. It suffices indeed to define $C_m$ as a small
perturbation of an union of concentric circles (plus a horizontal
line outside them when $m$ is odd). Either by inspecting the
Riemann surface or just by noticing that $\tau$ reverses
orientation of the ovals (and the pseudoline if $m$ odd) we infer
that $\tau$ takes the polarized curve to its opposite (where the
other half is preferred).

Using connectedness of $G=PGL(3, \RR)$ (more generally $PGL(n,
\RR)$ is connected whenever $n$ is odd because then both
components of $GL(n, \RR)$ given by the sign of the determinant
coalesce together since the identity matrix $I_n$ and its opposite
$-I_n$ are homothetic yet of opposite determinants), we infer
existence of a path $c$ in the Lie group $G$ connecting $id$ to
the symmetry $\tau$. Applying this path to $C_m$ gives a loop
$\gamma$ in the deep chamber $D$ based at $C_m$. Lifting this loop
to the $O_2$ cover continuously amounts tautologically to apply
the path $c$ to the polarized curve, whose end-point $c(1)$ is the
opposite polarization as the one we started with. Hence the lift
of the loop $\gamma$ is not a loop, and covering theory tell us
that $\gamma$ is not null-homotopic.
\end{proof}

If now $O$ denotes a specific orthosymmetric chamber we also have
the double cover $O_2\to O$ (by polarized Riemann surfaces) and it
is likely that the above argument extends to all or at least some
orthosymmetric chambers having a representative with  a mirror
symmetry $\tau$. Each such chamber would not be simply-connected.

Abstractly  an orthosymmetric surface can always be rotated by an
half-twist permuting both halves. However it is not evident that
this can be done in the plane, at least we know about no general
argument. Thus we retract to  examples in degree $m=6$, where due
to the combined efforts of Harnack-Hilbert-Gudkov-Rohlin we know
exactly what happens (cf. the Gudkov
table=Fig.\,\ref{Gudkov-Table3:fig}). In each orthosymmetric
chamber we look for a symmetric representant under an involution
fixing a line.

Consulting this table and gathering  earlier constructions  on a
single plate (Fig.\,\ref{Symmetry:fig} below) gives the following
symmetric realizations of dividing sextics:

$\bullet$ the $M$-schemes of Hilbert and Harnack can both  be
given a symmetric realization as evidenced by the picture below.
In both cases the invariant line intercepts 3 ovals.

$\bullet$ for Gudkov's scheme $\frac{5}{1}5$ the existence of a
symmetry is less obvious in Gudkov's original construction (cf.
Fig.\,\ref{GudkovCampo-5-15:fig}). The situation appears more
pleasant on Viro's construction of the latter (compare
Fig.\,\ref{Viro3-15:fig}c), but alas since we cannot choose
$\alpha=\beta=2$ both in V1 and V2 (please refer to the notation
of that figure) we cannot conclude the existence of a global
symmetry. Should we conjecture that Gudkov is somehow asymmetric?

$\bullet$ $\frac{8}{1}$ admits a symmetric realization as shown by
a variant of Hilbert's method (cf. figure below). Notice also the
model with double (dihedral) symmetry. Again 3 ovals are
intercepted by the (vertical) axis of symmetry.

$\bullet$ $\frac{6}{1}2$ in Hilbert's realization below severely
lacks symmetry. Appealing to Viro's method
(Fig.\,\ref{Viro3-15:fig}c or below) does not aid  (V5 twice looks
promising yet do not confuse the values of $\alpha,\beta$!).

$\bullet$ $\frac{4}{1}4$ along Hilbert's realization again lacks
symmetry, and Viro's method does not seem to help.

$\bullet$ $\frac{2}{1}6$ in Hilbert's realization  is asymmetric.
Via Viro's method this is realized by taking in V1 bottom
$(\alpha, \beta)=(1,1)$ and in V2 top  $(\alpha, \beta)=(0,4)$.
Alas this is highly asymmetric.

$\bullet$ $9$ is symmetric under Hilbert's construction,
or a more elementary (Pl\"ucker-style) deformation of 3 ellipses,
which is even more symmetrical.

$\bullet$ $\frac{5}{1}1$ is symmetric as shown below via a
primitive perturbation of ellipses \`a la Pl\"ucker-Klein (pre
Harnack-Hilbert oscillation trick). Again 3 ovals are intercepted
by the symmetry-axis.

$\bullet$ $\frac{3}{1}3$\vadjust{\vskip2pt} is symmetric  by a
perturbation of ellipses depicted below.

$\bullet$ $\frac{1}{1}5$\vadjust{\vskip2pt} is likewise symmetric
as shown by the depiction below.

$\bullet$ $\frac{4}{1}$\vadjust{\vskip2pt} is highly symmetric as
shown by a perturbation of ellipses below.

$\bullet$ $\frac{2}{1}2$ is symmetric as shown by the perturbation
of ellipses below (2 models).

This is the  exhaustive list of sextics of type~I, modulo the
omission of the deep nest (which is certainly symmetric). When
taking 3 concentric circles one gets the impression of a
continuous Lie group of symmetries, yet any perturbed curve will
be more rigid (recall finiteness of automorphisms due to
Schwarz-Klein-Poincar\'e-Hurwitz and the bound $84(g-1)$).

\begin{figure}[h]
\hskip-1.2cm\penalty0 \epsfig{figure=Symmetry.eps,width=152mm}
\vskip-5pt\penalty0
  \caption{\label{Symmetry:fig}%
  Catalogue of dividing curves with symmetry; ``A''=asymmetric
  realizations.} \vskip-5pt\penalty0
\end{figure}

From this investigation it follows the:

\begin{lemma}
All orthosymmetric chambers of sextics are not simply-connected
except perhaps the $4$ ``antechambers'' nearby Gudkov's scheme
that is $\frac{5}{1}5$, $\frac{6}{1}2$, $\frac{4}{1}4$, and
$\frac{2}{1}6$. In fact all sextic orthosymmetric chambers $O$
have nontrivial polarization covering $O_2\to O$, safe perhaps the
$4$ above
asymmetrical schemes.
\end{lemma}

[11.02.13] Perhaps there is an obstruction for those 4 schemes to
admit a symmetry about a line. In case of a Gudkov curve (of type
$\frac{5}{1}5$), the symmetry has to leave invariant 3 ovals for
the unique nonempty oval has to be preserved while the number of
inner and outer ovals are odd. At this stage (or earlier) it is
pleasant to visualize the Riemann surface in 3-space. Besides the
horizontal orthosymmetry imagine a rotational symmetry under
half-twist (180 degrees) leaving 3 ovals invariant while the 8
remaining one are pairwise exchanged. Of course per se this is  no
obstruction since Hilbert or Harnack have such a symmetry. So the
obstruction is necessarily a subtle one if it exists perhaps say
\`a la Arnold-Rohlin. Another idea is to smooth the Gudkov curve
along the axis of symmetry and hope to get a septic violating
B\'ezout (Fig.\,\ref{Symmetry:fig}d), but looks improbable.

Let us look at Rohlin's formula $2(\Pi^+ -\Pi^-)=r-k^2$ (see
(\ref{Rohlin-formula:thm})). Applying it to a Gudkov curve we have
$r=11$ and $k^2=9$, hence $(\Pi^+ -\Pi^-)=1$. But on Gudkov's
curve we have $5$ injective pairs of ovals, i.e.
$5=\Pi=\Pi^{+}+\Pi^{-}$, and it follows $2\Pi^+=5+1=6$, whence
$\Pi^{+}=3$ and $\Pi^{-}=2$. Can it be inferred that there is no
symmetry? A priori not, since the 5 inner ovals could have their
complex orientations being reversed by the symmetry while one
oval is kept invariant. Doing the same calculation for Hilbert's
curve we find $(\Pi^+ -\Pi^-)=1$ and $9=\Pi=\Pi^{+}+\Pi^{-}$, so
$2\Pi^+=9+1=10$ and $\Pi^{+}=5$ while $\Pi^{-}=4$.

The symmetry could be not a reflection about a line but a rotation
about a point. Yet from the projective viewpoint this seems to be
equivalent. At any rate Gudkov's curve in Viro's realization is
anyway not symmetric under a rotation.

Let us apply Rohlin's formula to the $(M-2)$-schemes which are
potentially asymmetric, e.g., $\frac{6}{1}2$. Then $2(\Pi^+
-\Pi^-)=9-k^2=0$, hence $\Pi^+ -\Pi^-=0$, but $\Pi^+ +\Pi^-=6$ so
that $2\Pi^+=6$, and $\Pi^+=3=\Pi^{-}$. Hence the symmetry cannot
exchange the $6$ inner ovals in pairs without fixing any of them.

This requires some  explanation. Recall we are looking for
holomorphic involutions of some plane dividing curve $C_m$ induced
by an element of $PGL(3, \RR)$ exchanging both halves. Call such
an involution a {\it mutation\/}. If a curve has a mutation then
its chamber $O$ has nontrivial polarized covering $O_2\to O$.

For sextics (except the deep nest and the unnested curves
$1,2,\dots,10$) there is a unique nonempty oval. Distinguished as
a such, this must be preserved by the  mutation $\tau$, which must
 reverse its orientation. If not,  orientation is
preserved and $\tau$ acts as a rotation on this circle. Taking an
invariant tube-neighborhood one deduces that both halves are
preserved as $\tau$ respects orientation of the surface (hence of
this tube), violating the mutating assumption.

Supposed fixed a complex orientation of the dividing curve. The
mutation reverses orientation of the nonempty oval, and also
 the  complex orientations of all other ovals because
$\tau$ preserves orientation but exchanges both halves.
Symbolically, we may see this by writing $\tau(\partial
C^+)=\partial (\tau C^+)=\partial (C^-)=-\partial (C^+)$.

At this stage we are ripe for picturing. Imagine the mutation
given by a symmetry about a line (this is probably no loss of
generality in projective geometry, as the other candidate namely a
rotation about a point fixes the line at infinity).
Consider the following schematic pictures (Fig.\,\ref{Sym2:fig}).
The first (Fig.\,a) is not mutating the orientation, hence
precluded. Fig.\,b is mutating the complex orientation, but
violates Rohlin's formula. In fact the mutation condition in case
where no inner ovals are invariant imposes an even number of
positive injective pairs of ovals and Rohlin's formula cannot be
fulfilled. Hence Rohlin's formula forbids  a mutation without
invariant inner oval. Fig.\,c shows a configuration where both
mutation and Rohlin's formula are satisfied. This is good
schematically but bad theoretically, as it fails obstructing
mutability of a curve of type $\frac{6}{1}2$. We are not much
advanced in our problem. The other figures of Fig.\,\ref{Sym2:fig}
show that for each of the other asymmetric types there is always a
schematic symmetry compatible with both  mutation and Rohlin's
formula. So Rohlin fails to detect any structural asymmetry.
Paraphrasing, asymmetry may just be a defect of our models
(Hilbert and Viro) yet not an intrinsic property of the chamber.

\begin{figure}[h]
\centering
\epsfig{figure=Sym2.eps,width=122mm} \vskip-5pt\penalty0
  \caption{\label{Sym2:fig}%
  No obstruction to symmetry via Rohlin's formula} \vskip-5pt\penalty0
\end{figure}

In fact,   the above Rohlin's formula argument only shows for the
two $(M-2)$-schemes of type~I ($\frac{6}{1}2$ and its mirror),
that if they have a mutation the latter must preserves 2 inner
ovals.

\begin{defn} [12.02.13] $\bullet$ A smooth dividing plane curve
$C_m$ defined over $\RR$ is transmutable if there is a
rigid-isotopy switching its half (called a transmutation).

$\bullet$ A mutation is a linear automorphism $\tau \in
G=PGL(3,\RR)$ of the curve $C_m$ permuting both halves of the
curve, say in this case that the curve is mutable.
\end{defn}

Since the group $PGL(3,\RR)$ is connected any mutation induces
(non canonically) a transmutation. Indeed  choose a path $c$ in
$G$ joining $id$ to $\tau$ and its operation upon $C_m$ defines a
loop in the corresponding chamber of the space of curves (past the
discriminant) which is a transmutation. It is not essential that
$\tau$ has order 2, but then speak of a $2$-mutation.

So any $2$-mutable (dividing) curve is mutable, and in turn
transmutable. The converses looks a priori quite improbable. Are
all dividing plane curves transmutable? or even mutable, or
$2$-mutable after some rigid-isotopy? The question looks of
interest because a non transmutable curve would have a preferred
 half (privileged so-to-speak) which looks a bit against the
flavor of Galois-theory and French revolution
``\'egalit\'e, fraternit\'e, etc.''.

A mutation (like any self map of $\RR P^2$) has a (real) fixed
point (e.g. via Lefschetz fixed point theorem using homology over
${\Bbb Q}$), and any $2$-mutation is a mirror about a line fixing
also a real isolated point, as inferred from linear algebra
(existence of real eigenvalues for an endomorphism of a real
vector space of odd dimension).

Perhaps the above questions  can be handled via C. Segre's
classification of real structures on projective spaces, especially
the fact that the plane ${\Bbb P}^2$ has a unique real structure,
but looks unlikely as the curves are not taken into account.

One way to approach the problem in general would be to look at the
action of $G=PGL(3, \RR)$ on the chamber past  the discriminant
containing a dividing curve. Existence of a mutation in each such
chamber amounts this action being never free, i.e. with nontrivial
isotropy subgroup $G_{C}$ at some suitable curve $C=C_m$. This is
not enough for the automorphism in question needs not permute both
halves. Omitting this difficulty, there would be a free Lie group
operating, hence an induced foliation of the chamber by leaves of
dimension 8 (=$\dim G$). Alas such a foliation also exists when
the action is only locally-free, i.e. discrete isotropy are
allowed. It would be nice to know if all chambers of the
discriminant (not only the orthosymmetric chambers) contains a
curve with (linear) automorphism in themselves. For orthosymmetric
chambers we would further like to know if there is such an
automorphism permuting the halves of the curve (i.e. a mutation).

In a remarkable article extending earlier work by Kharlamov,
Itenberg 199X \cite{Itenberg_199X-monodromy-deg-6} is able to
compute the monodromy groups of each chamber of sextics.
Extracting from his tabulation, only the type~I cases gives the:

\def\ZZ{\Bbb Z}
\def\triv{1}

\begin{lemma} {\rm (Kharlamov, Itenberg)}
The monodromy groups of smooth sextics of dividing type are given
by the following list (where $\triv$ is the trivial group, $S_n$
the symmetric group on $n$ letters, and $D_n$ the dihedral group):

$\bullet$ $\frac{9}{1}1 \rightsquigarrow \ZZ_2 $,
 $\frac{5}{1}5 \rightsquigarrow \triv $,
 $\frac{1}{1}9 \rightsquigarrow S_3 $, \vadjust{\vskip2pt}

$\bullet$ $\frac{8}{1} \rightsquigarrow D_4 $,
 $\frac{6}{1}2
\rightsquigarrow \triv $,
$\frac{4}{1}4 \rightsquigarrow \triv $,
$\frac{2}{1}6 \rightsquigarrow \ZZ_2 $,
$9 \rightsquigarrow S_9 $,\vadjust{\vskip2pt}

$\bullet$ $\frac{5}{1}1 \rightsquigarrow \ZZ_2 $,
$\frac{3}{1}3 \rightsquigarrow \ZZ_2 $
$\frac{1}{1}5 \rightsquigarrow D_5 $,\vadjust{\vskip2pt}

$\bullet$ $\frac{4}{1} \rightsquigarrow S_3 $,
$\frac{2}{1}2 \rightsquigarrow \ZZ_2 \times \ZZ_2 $,
\vadjust{\vskip2pt}

$\bullet$ $(1,1,1) \rightsquigarrow \triv$.
\end{lemma}

It is interesting to compare this result with our picture
Fig.\,\ref{Symmetry:fig}, as sometimes the whole monodromy group
can be realized by rigid projective motions.

Besides, it seems interesting to compare this monodromy of ovals
to the monodromy upon the halves. Albeit the latter viewpoint is
less rich in general (being only a representation on the group
with 2 elements) it is sometimes of complementary nature in
detecting non-triviality of the fundamental group of the fixed
chamber. For the moment, our halves-monodromy is only more
sensitive in the deep-nest case. (Question: does the $\pi_1$ of
the deep chamber reduces to $\ZZ_2$?)

Finally, note that the oval-monodromy is also fairly small for the
4 exceptional schemes, asymmetric in Hilbert's (or Viro's)
realization (again Fig.\,\ref{Symmetry:fig}). Hence both methods
oval-monodromy and half-monodromy (at least via rigid symmetries)
fails to detect nontrivial elements in $\pi_1$ of the
corresponding chamber for the schemes of Gudkov $\frac{5}{1}5$, of
left-Rohlin $\frac{6}{1}2$ and $\frac{4}{1}4$. Can we extrapolate
that those chambers are simply-connected? If yes then those 3
curves are not transmutable, hence not mutable and therefore
structurally asymmetric (i.e. there is no model invariant under a
mirror).

More factually, Itenberg's calculation prevents a curve of type
$\frac{6}{1}2$ to accept a mirror like Fig.\,\ref{Sym2:fig}c.
Indeed otherwise if $\tau$ is such a mirror, it suffices to take a
path in $PGL(3, \RR)$ joining the identity to this $\tau$ to get a
loop in the space of curves with non trivial monodromy. Likewise
no curve of type $\frac{4}{1}4$ can accept a mirror like
Figs.\,\ref{Sym2:fig}d,e, and no curves of the Gudkov type
$\frac{5}{1}5$ can accept a mirror like Figs.\,\ref{Sym2:fig}i.
Hence Itenberg's calculation implies the following answer to one
of our basic question:

\begin{theorem}
None of the $3$ monodromically-trivial dividing curves---as listed
by Itenberg, i.e. Gudkov's, the left Rohlin curve $\frac{6}{1}2$,
and that of type $\frac{4}{1}4$---can support a mirror.
\end{theorem}

\begin{proof}
A mirror is a linear involution of $PGL(3,\RR)$ which fixes a line
(plus a point at $\infty$). By B\'ezout, the fixed line can
intercept at most 3 ovals which are then invariant. Examining the
following Fig.\,\ref{Sym3:fig} shows that there is always some
pair of ovals permuted by the mirror.

\begin{figure}[h]
\centering
\epsfig{figure=Sym3.eps,width=122mm} \vskip-5pt\penalty0
  \caption{\label{Sym3:fig}%
  Mirrors obstructed by Itenberg's calculation of the monodromy} \vskip-5pt\penalty0
\end{figure}

Alas, an oval can be invariant under the mirror without having to
intersect the fixed line, say by running to infinity. Then one can
try to argue with the line at infinity and B\'ezout.

A better way to argue is that the nonempty oval of the sextic has
to be invariant under the mirror, hence at most 2 inner ovals can
be left invariant under the mirror (else B\'ezout corrupted). So
we infer existence of inner ovals permuted under the mirror
$\tau$. Taking a path from $id $ to $\tau $  yields a loop in the
chamber whose monodromy is nontrivial. This contradicts
Itenberg's calculation of the monodromy.
\end{proof}

In particular this shows  existence of curves without a
$2$-mutation, since a $2$-mutation is certainly a mirror. However
it is not clear if our 3 curves lack a mutation. Given a mutation
(i.e. a linear automorphism) permuting the halves it is not a
priori of order 2. However it has finite order (by
Klein-Poincar\'e-Hurwitz finiteness of the automorphism group of
closed Riemann surfaces of genus $\ge 2$). Since $\tau$ permutes
the halves its order must be even say $n=2 e$. Hence $\tau^e$ is
an involution, and a mutation if $e$ is odd. Alas in the other
case, i.e. when $n$ is divisible by $4$, we cannot say much.

Assume given a mutation $\tau$ on a Gudkov curve, then it has to
preserve the unique nonempty oval, and the inner and outer ovals
have to be respected. Further $\tau$ (being a mutation) it has to
reverse the complex orientation. As we may assume the order $n$ of
$\tau$ divisible by $4$, it should follow that $\tau$ permutes the
ovals along a 4-cycle, etc. Via some group theory there is
presumably an obstruction to mutate the Gudkov curve?

Alternatively assume there is a mutation on some Gudkov curve
$C_6$ then it will permute the ovals while preserving the nonempty
one and the inner/outer ovals subdivision. Further it must reverse
the complex orientation. If any  permutation of the ovals is
detected (which is precluded by the Kharlamov-Itenberg's
calculation of the monodromy as being trivial), we are finished.
Assume so that the mutation preserves all ovals. Since it has to
reverse their orientations, the mutation has 2 fixed points on
each of them (an orientation reversing transformation of the
circle has 2 fixed points e.g. via Lefschetz or via covering
theory). But globally our mutation is of finite order (since it
induces an automorphism of the Riemann surface of genus $10\ge
2$), and any element of finite order in $PGL(3, \RR)$ preserves a
line. So we get a line intersecting the $C_6$ in 22 points,
overwhelming B\'ezout. As another variant, once our mutation is
known to have 2 fixed points on each ovals we have $2(g+1)$ of
them, which is the maximum permitted by Lefschetz trace formula.
We would like to conclude that $\tau $ is the hyperelliptic
involution, which cannot exist on a smooth $C_6$ (which is only
$5$-gonal).

Let us clarify this argument with the:

\begin{defn} {\rm A curve is {\it antidromic\/}
if it monodromy group is trivial. A {\it symmetry\/} of a plane
real curve is a linear automorphism of the plane (i.e. an element
of $PGL(3,\RR)$) preserving globally the curve.}
\end{defn}

\begin{lemma}
{\rm (1)} Any symmetry of an antidromic curve must
leave each oval invariant. Hence all the $3$
antidromic  sextics  listed by Itenberg can only admit a mutation
preserving all the ovals.

{\rm (2)} In particular all the $3$ antidromic dividing sextics of
Itenberg (the deep nest being excluded) lack a mutation. So
Gudkov's curve, and the left-wing Rohlin curve $\frac{6}{1}2$ and
the dividing curve $\frac{4}{1}4$ are asymmetric at least under a
mutation (i.e. they cannot mutate).
\end{lemma}

\begin{proof}
(1) If there is a symmetry $\tau$ permuting somehow the ovals,
then the path in $PGL(3,\RR)$ connecting $id$ to  $\tau$ induces a
loop in the space of curves with nontrivial monodromy (namely the
ovals-permutation induced by $\tau$).

(2) As to the second assertion, our oval-preserving mutation must
necessarily invert the orientation of all ovals since it exchanges
both halves. (Recall that a mutation reverses the complex
orientation in the sense of Rohlin, since it preserves the
orientation induced by the complex structure while exchanging both
halves of the dividing curve.) Since a sense-reversing
transformation of the circle has 2 fixed points as follows from
Lefschetz's fixed-point formula. (Indeed the Lefschetz trace
number is $(+1)-(-1)=+2$ so there is a fixed point and removing it
one obtains an orientation reversing homeomorphism of the line
which has another fixed point, e.g. by applying Bolzano to the
graph of this continuous decreasing function.) (Is it true that
any sense-reversing transformation of the circle is an involution?
We do not need this anyway!)

Since Itenberg's
antidromic  curves have $r=11$ or $r=9$ ovals, we get 22 or 18
fixed points created.
Next we have:

\begin{lemma} A linear automorphism in $PGL(3, \RR)$
that fixes $4$  (or more) points of $\RR P^2$ fixes either a line
plus an isolated point or is the identity.
\end{lemma}

\begin{proof}
This follows by looking at the 4 corresponding eigen-lines in
$\RR^3$, two of which have to correspond to the same eigenvalue
(pigeonhole principle), and so there is the required fixed
projective line. The third eigenvalue left (necessarily real)
gives a third eigen-line and the announced alternative follows
depending on whether this 3rd eigenvalue differs or coincides with
the former double eigenvalue. \end{proof}

It follows that
among our 22 or 18 many fixed-points  (at least so many less one)
are aligned, but this corrupts B\'ezout.
\end{proof}

A priori it is much harder to detect an obstruction to transmute
the Gudkov curve (or its 2 antidromic cousins), and likewise hard
to show that it may be transmuted.

Assume there is a transmutation. Then since by
Kharlamov-Itenberg's lemma  the curve is antidromic the induced
permutation of ovals is trivial, but the complex orientation is
reversed. Perhaps some obstruction can be deduced from this...

[15.02.13] Two days ago, Kharlamov informed me that the Gudkov
chamber has fundamental group $\ZZ_2$ compare his letter (dated
[13.02.13]) in Sec.\,\ref{e-mail-Viro:sec}. So I presume that the
Gudkov curve can be transmuted, and that the isomorphism
$\pi_1\approx \ZZ_2$ may be realized as the monodromy acting upon
halves. Kharlamov's messages also emphasize the issue that the
$M$-curves case is somewhat easier than the other cases. Hence
while the (oval)-monodromies are completely calculated by
Itenberg, it may be the case that the determination of the
fundamental group of each chamber is somewhat harder to obtain.

Some naive questions are as follows. We presume that the deep-nest
chamber has $\pi_1=\ZZ_2$. Are the (fundamental) of chambers
always finite? This would follow (theorem of Myers, Synge, Hopf,
etc.) if there is a complete metric of positive curvature on the
hyperspace of curves (which is the case) yet the natural elliptic
metric is not complete when restricted to the chambers. Further
since Gudkov's chamber is not simply-connected, any Gudkov curve
is presumably transmutable (this amounts to say that the Kharlamov
isomorphism $\pi_1=\ZZ_2$ is realized by the monodromy of halves).
If so is the case, perhaps that even all dividing plane curves are
transmutable.

\smallskip
{\footnotesize

(Loose ideas, skip the next 2 points:)

[10.02.13, 23h49] As to the rigidity of the scheme with one oval
only we can try to investigate the degeneration $C_{2k}\to E_2
\cup C_{k-1} \cup C_{k-1}^{\sigma} $ to an ellipse via some
suitable geometric flow. Then the conjugate pairs of curves
$C_{k-1} \cup C_{k-1}^{\sigma}$ would be empty so rigid-isotopic,
etc, etc.

[11.02.13] As to Rohlin's (unproved) assertion (1978) of the total
reality of the schemes $\frac{6}{1}2$ try perhaps to look at the
cubics spanned by the hexagon through the inner ovals.

}

\subsection{Isotopic rigidity of the empty scheme:
connectedness of invisible
curves}\label{rigidity-empty-scheme-via-dyna:sec}

[08.04.13] We now come to a tortuous revelation of a basic truth,
namely the fact that the empty chamber is connected, i.e. any 2
real smooth plane curves with empty real locus can be connected by
a path of similar curves (avoiding the discriminant). This
problematic covers no less than $6$ sections up to
Sec.\,\ref{invisible-discriminant-codim-2:sec}, where it is
elucidated that the portion of the discriminant inside the locus
of empty(=invisible) real curves has codimension 2 hence cannot
effect a disconnection. Pivotal in this search was a kind letter
by Shustin explaining the linear homotopy argument between empty
curves, which shows that the empty locus (including possibly
singular curves) is connected, and actually much more like being
contractible. So we warn the reader that those six sections are
far from a geodesic toward the goal, but we had not the courage to
censure any bit of our poorly organized material as it  often
ramifies toward considerations of independent interest.


[23.01.13] In fact to be honest I realize that even the foundation
of our reasoning (in the previous section\footnote{After several
permutation in this text, it is not clear anymore what was the
``previous section'', but checking dates and contents this can be
almost surely identified as Sec.\,\ref{application-of-CCC:sec}.})
is not completely sound, namely the following fact:

\begin{theorem} {\rm (Folklore??? is it really true?
If yes where is it proved?)}
Any two empty curves are rigid-isotopic.
\end{theorem}

{\footnotesize

{\it Insertion} [08.04.13].---Folklore probably! True, certainly,
compare Shustin's letter, but also the codimension 2 lemma
(poorly) established in our
Lemma~\ref{invisible-discriminant-codim-2:lem}. Where it is
proved? We still lack a detailed reference, but apparently the
fact is so trivial that nobody took care writing down a complete
proof. It would be of interest to make a deeper historical search
of who knew first this simple result. Possible candidates:
Schl\"afli 1863 \cite{Schlaefli_1863}, Cayley, Klein 1873--1925,
C. Segre, Hilbert ca. 1891, Berzolari 1906 \cite{Berzolari_1906},
Severi e.g. 1921 \cite{Severi_1921-Vorlesungen-u-alg.-Geom-BUCH},
Brusotti 1921 or earlier, Petrovskii 1933, etc.

}

\smallskip
(Sometimes instead of saying empty curves we shall say invisible
curve, when the real locus $C_m({\Bbb R})$ is empty.)

I was sure this to be known (but completely forgot where I read
this in case I am remembering well!!!) At first sight this looks
trivial but is not. One could imagine the empty chamber to be
starlike or even convex, but this is not even evident. Somehow
invisible curves could be like the immersed half of an iceberg,
hence connected, while the visible part of the iceberg may consist
of several islands (peaks) corresponding to the menagerie of
chambers past the discriminant, well-known in Hilbert's 16th.
%
Alas even if this metaphor ought to contain some truth, it is easy
to construct an iceberg with disconnected immersed locus. Take a
letter ``E'', rotate it by $-\pi/2$ to get
the symbol ($\Pi\!\!\Pi$) considered as a tripod with 3 legs
immersed in the water. To go in 3D, just take the body of
revolution of this symbol to get an iceberg with 2 immersed
components\footnote{I should acknowledge the assistance of my
cousin \'Elias Boul\'e-Schneider for several discussions on such
topics, like topographical discussions about Hilbert's 16th.}.
After having being puzzled for while one might suspect this to be
a result \`a la Hilbert, that a form not representing zero is
something like a sum of squares...(not clear).

[24.01.13] More geometrically, one can look at the distance
between the complex locus of a real curve $C_m({\Bbb C})$ and the
real plane ${\Bbb R}P^2$. It is natural to work with the
Fubini-Study metric on ${\Bbb C}P^2$. Then  look at ``the'' point
of the real plane closest to $C_m({\Bbb C})$, which must be
generically unique. If we let $E$ be the empty chamber (a priori
not connected), we obtain so a random-map $\pi\colon E \to {\Bbb
R}P^2$ taking an invisible curve to its closest ``projection'' in
${\Bbb R}P^2$. Up to removing some subset of $E$, one could
arrange $\pi$ to be single-valued, and one would check that $\pi$
is akin to a fibration with connected fibres. The connectedness of
$E$ could follow. It is also tempting to imagine a flow driving
invisible curves to solitary nodes. This would be just the
gradient flow of the functional distance to the real locus ${\Bbb
R}P^2$. The corresponding trajectories of steepest descent could
converge
to a curve with a unique solitary node (generic case). At the
level of the Riemann surface this isotopy (given by the path of
trajectory) really amounts to the contraction of an anti-oval
toward a solitary node. So this is just a special case of Klein's
Ansatz (\ref{Klein-1876:conj-noch-entwicklungsfahig}), that a
nondividing curve can acquire a novel solitary node (by a large
deformation). (Of course recall this to be erroneous by Shustin
1985 \cite{Shustin_1985}, yet it is perhaps true for empty
curves). In fact:

\begin{lemma} Klein's Ansatz is
 trivially true for empty curves
 (just by general position and surgery
of the real locus).
\end{lemma}

\begin{proof} Indeed, given any invisible curve $C_m$, take any pencil
through it passing through a visible curve $D_m$ (with nonempty
real locus), then making this line transverse to the discriminant
we get Morse surgeries the first of which must necessarily be a
solitary node formation.
\end{proof}

{\it Insertion} [08.04.13] Maybe one can object again this proof,
by arguing that the first contact with the discriminant could be
through a pair of imaginary nodes. Paraphrasing a bit we could
imagine that travelling along the pencil spanned by $C_m,D_m$ we
hit the discriminant but then fall again in the empty chamber.
Presumably both scenarios can be avoided if we know that the
invisible discriminant has real codimension 2 (as we shall see in
the sequel, cf. Lemma~\ref{invisible-discriminant-codim-2:lem}.)

Our flow would precisely do this contraction yet in some more
organic(ized) fashion (i.e. no choices). Yet notice that we could
make the above pencil argument by choosing once for all some
visible curve $D_m$, while driving all the invisible curves $C_m$
along the line spanned by $C_m$ and $D_m$. If this does not work
look at the flow (discussed above).

Optimistically, this method may suffice to establish connectedness
of the empty locus $E \subset \vert m H\vert$. One may wonder if a
variant of the argument could not also establish Viro's open
problem on the connectedness of the pseudoline locus $P\subset
\vert m H \vert$ when $m$ is odd. (Added [08.04.13].---It seems
that this conjecture really goes back to Rohlin, if we interpreted
correctly a letter of Viro in Sec.\,\ref{e-mail-Viro:sec}, dated
[26.01.13].)

Return yet to the case of the empty locus $E$ (non void only for
curves of odd degree). By what could it be disconnected? A curve
in the discriminant $\frak D$ may well have two imaginary
conjugate singularities. But this locus has codimension 2, so it
cannot disconnect $E$. Does this suffices to prove
connectedness of $E$? Probably not as a priori it may have several
components lying ``far apart''.

Another idea is to fix an invisible conic, e.g. the ``canonical''
one $E_0\colon x_0^2+x_1^2+x_2^2=0$ (on which conj acts like an
antipodal map). To each point of the plane one can attach the
apparent contour of this ellipse (polar lines) as seen from the
given point to get a group of two points on this ellipse which is
a Riemann sphere. (This is the most synthetic way to establish the
well-known isomorphism between ${\Bbb C}P^2$ and the second
symmetric power of ${\Bbb C}P^1$.) Points of $E_0$ correspond to
groups of superposed two points, while real points in ${\Bbb
R}P^2$ maps to antipodal pair (invariant under conj when seen as a
pair). Define a function $\rho$ on ${\Bbb C}P^2$ which given a
pair measures the distance on the round sphere $S^2$ between the
corresponding 2 points. This is equal to $\pi=3.14\dots$ on ${\Bbb
R}P^2$ and vanishes on $E_0\colon x_0^2+x_1^2+x_2^2=0$.

Now given any compact sublocus of ${\Bbb C}P^2$ (in particular an
invisible curve $C_m$) one can look at the maximum of $\rho$ on
$C_m({\Bbb C})$, which is $<\pi$. This gives the functional
$$
\theta:=\max \rho\colon E \to [0,\pi[,
$$
whose ascending gradient lines should converge to solitary nodal
curves (perhaps with several such nodes). Note that the ground
(invisible) ellipse $k \cdot E_0$ (counted $k=m/2$ times) is the
unique absolute minimum of this functional. This $\theta$ is quite
likely to be a Morse function (or a slight generalization thereof
with Monkey saddles, etc.), yet the critical points (causing
annoying stagnation of the dynamics) ought to be isolated
(codimension 2 suffices), hence not affecting the connectivity of
$E$. By construction our flow tends to make an invisible curve
more ``visible'' by pushing it progressively closer to the real
plane. In the limit we expect something visible having solitary
nodes and generically just one should emerge. So upon excising
from $E$ a small set (of codimension 2, since at least two nodes
is bad) we find a (dense) subregion $E^{\ast}$ which  maps ${\Bbb
R}P^2$ by assigning the unique solitary node of the limit of $C_m
\in E^{\ast}$ under the flow at time $\theta=\pi$ (using $\theta$
as time parameter as usual for gradient flow). Now the fibre of
this map $E^{\ast}\to{\Bbb R}P^2$ is the same as the {\it bassin
d'attraction\/} of the flow which is cone-like formed by several
trajectories abutting to the same solitary node. So this cone is
connected by the end point, and connectivity of $E^{\ast}$ (hence
$E$) should follow.

Still, the main difficulty is (as usual) to show that the
$\theta$-functional lacks a local maximum preventing convergence
to a visible curve. So given any invisible curve one should
produce a small perturbation with larger $\theta$. This is
probably not too hopeless. Naively one could perturb $C_m$ inside
the pencil spanned by $k E_0$ and $C_m$. Since $k E_0$ is the most
invisible curve a deformation along it should decrease $\theta$,
while one [deformation] away [of]
it should increment $\theta$. Of course there is some objection to
this, since in a projective (real) line (a circle) it is never
clear what means ``along'' and ``away''.

All this is somewhat confuse and unconvincing. Perhaps also there
is a much more elementary argument without gradient lines. As we
said this could involve deforming all empty curves along some
fixed visible curve. But which direction of retraction should we
choose in the corresponding pencil? Since $m=2k$ is even, $\vert m
H \vert$ is of dimension
$\binom{m+2}{2}-1=(k+1)(2k+1)-1=2k^2+3k=k(2k+3)=:N$. So when $k$
is odd there must be another singular point in the foliation
induced by the {\it faisceau\/} (sheaf, bundle) of all lines
through $D_m$ (by Poincar\'e-Hopf). Can we orient this foliation?
No because when $N=2$, we have a M\"obius strip after puncturing
the basepoint of the pencil. The situation would be somewhat
simpler if we could find in the hyperspace of curves $\vert m H
\vert \approx {\Bbb R}P^N$
a hyperplane $H\ni C_m$ avoiding
the empty locus $E$. I do not know whether this is possible? Then
there would be a nice way to retract the whole complement of $H$
toward the point $C_m$ in some canonical way. In particular all
points of $E$ (invisible curve) would mark a first impact on the
real locus ${\Bbb R}P^2$. Alas it is not even obvious that
connexity of $E$ follows.

[23.01.13] Recall that the related question for odd degrees is
still an open problem.

\begin{ques} {\rm
(Viro 2008
\cite[p.\,199]{Viro_2008-From-the-16th-Hilb-to-tropical})}.---{\rm
Are all non-singular real projective curves of a given {\it odd
degree} with connected set of real points {\it rigid-isotopic} to
each other?}
\end{ques}

The emphasis is Viro, and may suggest that without odd degree the
assertion is known to be false!? If so then our (OOPS) conjecture
(\ref{OOPS:one-oval-rigid-isotopic:conj})
 would be oops!

{\it Insertion} [08.04.13].---The answer were given in Viro's
letter dated [26.01.13] in Sec.\,\ref{e-mail-Viro:sec}, and may be
summarized as follows. First in the even degree case, the problem
of rigidity of the curve with a unique oval is still open (but
probably reducible to the contraction conjecture
(\ref{Itenberg-Viro-contraction:conj})). Second, as we already
said, the question ascribed above to Viro (2008) truly goes back
to Rohlin (unpublished as far as we know).

\subsection{Rigidity of the empty scheme (Shustin's letter)}

[27.01.13] This section treats the following desideratum: given
two empty (plane) curves (hence of even degree), it is always
possible to find a path of curves linking them while avoiding the
discriminant.

In fact, I read about this fact a long time ago (ca. 2000) but
could not remember from which source. Recently (24.01.13) failing
to recover the source, we started to doubt about the truth of this
assertion. Very kindly Shustin communicated us the simple
(forgotten) argument giving a positive answer. Alas we have little
idea of who proved this first. Shustin's proof looks at first
sight extremely trivial, but on more mature thinking the story
looks a bit more tricky than expected.

\begin{lemma}\label{empty-chamber-connected-Shustin:lem}
Any two empty curves are rigid-isotopic.
\end{lemma}

\begin{proof}  (Courtesy of Eugenii Shustin [26.01.13])
The chamber of empty curves of a given (even) degree is connected.
Two such curves are defined by homogeneous polynomials $P,Q$,
supposed (w.l.o.g.) positive for all real variables not all
simultaneously zero. The linear homotopy $(1-t)P+tQ$, $0\le t\le
1$ gives then a path in the chamber of empty curves.
%
Indeed the linear path between two positive numbers
consists of positive numbers, and so all intermediate curves of
this homotopy are empty curves.

A priori one can imagine that some intermediate curve of this
homotopy (while staying empty) crosses the discriminant by
acquiring a  conjugate pair of nodes. (This eventuality was not
mentioned in Shustin's letter, but we think that is is a slight
obstacle to the argument. However it seems to be not fatal as we
shall discuss at length.)

Maybe the above argument should be supplemented by an examination
of the partial derivatives $\frac{\partial P}{\partial x_i}$.
Smoothness of a curve amounts the 3 partial derivatives (of the
defining equation) lacking a common zero. Recall Euler's relation
$m F=\sum_i \frac{\partial F}{\partial x_i} x_i$ valid for any
homogenous polynomial (form) of degree $m$. (This requires only to
be checked on monomials $x_0^i x_1^j x_2^k$ such that $i+j+k=m$.)
Combining this we infer that if there is some time $t\in [0,1]$
such that $(1-t)P+tQ$ is singular, then all its 3 partials
vanishes simultaneously at some point $(x_0,x_1,x_2)$, which by
Euler's relation would be also a zero of $(1-t)P+tQ$. This
violates however the emptiness of this curve. However for this
argument to hold good it is essential for the point
$(x_0,x_1,x_2)$ to be real, which is however not the case a
priori.

Another more qualitative argument would be to first perturb
slightly $P$ and $Q$ so that the linear pencil spanned by them is
transverse to the discriminant. In that case if some member of the
pencil acquires a singularity it will be a simple node, which
consequently must be real. This violates  emptiness of the
intermediate curves $(1-t)P+tQ$, $0\le t\le 1$.

In summary the lemma looks true, yet not in the strong sense that
any two empty smooth curves are linked by a linear homotopy
consisting only of smooth curves. This strong form amounts to
convexity of the empty chamber, and not just connectedness.
Convexity is perhaps wrong, as the curve may traverse a pair of
conjugate nodes during the linear deformation.

Can we corrupt  convexity of the empty chamber? We think yes as
follows. For the empty chamber to be nonempty, assume the degree
$m$ even. Suppose given an empty curve $C_m$ with a pair of
conjugate nodes. We shall later explain how to construct this, but
one can already  imagine in 3-space (or in the ether) a
diasymmetric Riemann surface acted upon antipodically without
fixed point, on which two symmetric vanishing cycles are
contracted. Now through the given curve $C_m$ (considered as a
point in the hyperspace of $m$-tics), trace a little rectilinear
segment transverse to the discriminant. Both extremities of the
segment will be smooth empty curves, yet the linear homotopy
connecting them  hits the discriminant (hence fails to be entirely
in the empty chamber). This argument works fine provided the
linear homotopy coincides with our little segment, instead of
being actually the ``long'' residual pieces of it in the pencil
(which is a circle).

 When $m=2$, a conic can only have one
node, when degenerating to a pair of lines. As we require $2$
conjugate nodes, let us look at quartics. Start with a pair of
empty conics with transverse complexifications intersecting in 4
points $p,p^\sigma, q, q^\sigma$ (cf.
Fig.\,\ref{ShustinEmpty:fig}a for a schematic view). By Brusotti
1921 \cite{Brusotti_1921} (or earlier workers like Pl\"ucker,
Klein, etc.) we can smooth $q,q^\sigma$ away to create an empty
quartic $\Gamma_4$
 with 2 nodes nearby $p,p^\sigma$. The corresponding
(singular) Riemann surface (complex locus $\Gamma_4({\Bbb C})$) is
visualized as a genus 3 surface (acted upon by antipody) with 2
handles shrunk to points $p, p^\sigma$ exchanged by conjugation
(Fig.\,\ref{ShustinEmpty:fig}b). We can also imagine the structure
of the discriminant near $\Gamma_4$ as being a single nappe
(Fig.\,\ref{ShustinEmpty:fig}a). Take a little rectilinear segment
transverse to this nappe. Both extremities of this segment are
smooth curves $C_4$, $D_4$ which are empty. This gives our
counterexample provided the linear homotopy $(1-t) C_4 + t D_4$,
$t\in [0,1]$ visits $\Gamma_4$.

\begin{figure}[h]
\centering
\epsfig{figure=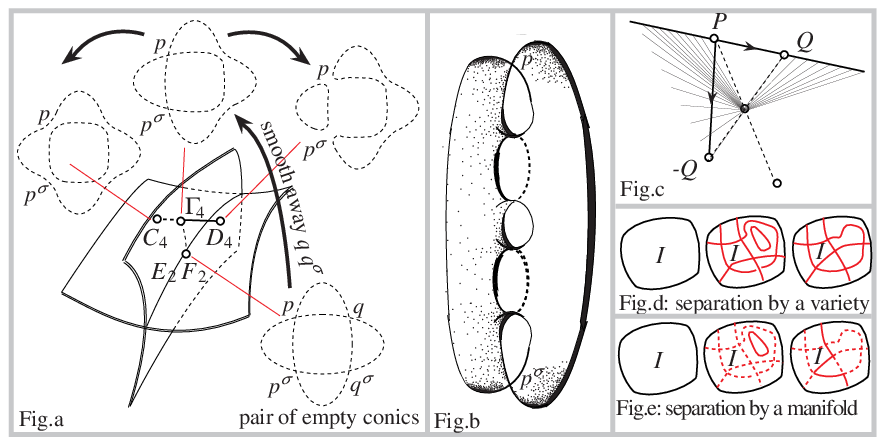,width=122mm} \vskip-5pt\penalty0
  \caption{\label{ShustinEmpty:fig}%
  A Brusotti type construction attempting to disprove
  the linear homotopy between empty curves} \vskip-5pt\penalty0
\end{figure}
\end{proof}

It is worth trying to clarify the above proof. Given two forms
$P,Q$ of degree $m$,  define the {\it linear homotopy\/} as the
path of forms $(1-t) P+t Q$ where $t\in [0,1]$. Denote it
symbolically $P\to Q$. If we look at curves (i.e. homothety
classes of forms) then between any two curves $C, D$ of degree $m$
there is a pencil of curves $\lambda C +\mu D$ which is the line
$\overline{CD}$ through both points seen in the hyperspace of
curves. A little drawing (Fig.\,\ref{ShustinEmpty:fig}c) shows
that if $P, Q$ represents $C,D$ resp., then the linear homotopy
$P\to Q$ projects to one piece of the line $\overline{PQ}$, while
the linear homotopy $P\to -Q$ describes the other road of access
in the circle $\overline{CD}$.

Shustin's argument shows that if two forms $P, Q$ not representing
zero have the same sign then the linear homotopy $P\to Q$ consists
of forms not representing zero. (Recall that the sign of even
degree forms is well-defined, because $F(\lambda x_0, \dots,
\lambda x_n)=\lambda^m F( x_0, \dots,  x_n)$.) However  it does
not say that if $P,Q$ represents nonsingular curves, then so are
all members of the linear homotopy $P\to Q$.

In our example with $C_4$ and $D_4$ we could argue that the
corresponding polynomials $P,Q$ are very near (by Brusotti's
construction) so of the same sign, and further that the linear
homotopy $P\to Q$ really passes through $R$ the defining equation
of $\Gamma_4$. In that case Shustin's argument would be
in slight jeopardy.

Does our counter-argument work? We think yes we start from $R$ a
form defining $\Gamma_4$, and perturb slightly the coefficients of
$R$ by Brusotti to get the polynomials $P$ and $Q$ so that the
linear homotopy $P\to Q$ passes through $R$. Note that $R$ has
some well-defined sign on ${\Bbb R}P^2$, and by smallness of the
perturbation $P$ and $Q$ have the same sign as $R$. Thus even when
$P$ and $Q$ do have the same sign we are not ensured that the
linear homotopy $P\to Q$ transits only through nonsingular curves.

In conclusion given two smooth\footnote{We use ``smooth'' as an
abridgement of nonsingular in the algebro-geometric sense, or if
you prefer smoothness of the complexification, but not merely of
the real locus. So a curve may look smooth while having imaginary
conjugate singularities. Such a curve is not considered as smooth
by us.} empty curves $C,D$ and choose representing forms $P,Q$
(resp.) of the same sign, then the projection of the linear
homotopy $P\to Q$ in the space of curves $\vert m H\vert$
correspond to empty curves yet not necessarily smooth.

{\footnotesize

{\it Related literature for Shustin's argument.} Maybe Wilson,
Shustin ICM, etc...

}

Of course even if the above Brusotti-type construction is correct,
it does not prove that the locus of empty curves is disconnected,
but merely that the proof via linear homotopies is insufficient.

One possible critique to our argument is that because $\Gamma_4$
has $2$ nodes it is not on a principal stratum\footnote{I borrow
this jargon from Degtyarev-Kharlamov 2000
\cite{Degtyarev-Kharlamov_2000}.} (wall) of the discriminant of
codimension 1. Yet in reality this is a pair of conjugate points
so really one point in the sense of Grothendieck's schemes (to
which we are from adhering). Perhaps our segment  not transverse
but rather tangent to the nappe of $\frak D$. However this looks
not so realist by construction. The key issue is to decide whether
our binodal curve $\Gamma_4$ is a smooth point of the
discriminant, which looks likely if we regard only real curves. Of
course in the complexified discriminant $\frak D({\Bbb C})$ there
is two nappes of passing through the binodal curves $\Gamma$.

If our reasoning is correct, we see that the problem of the
connectedness of the empty chamber
is not settled by the linear homotopy argument. Perhaps the empty
locus is even disconnected? How to approach the problem?

Let $m$ be some fixed even degree. Consider $\vert m H \vert$ the
space of all real curves, and $\frak D$ be the discriminant
parametrizing real singular curves. Denote by $I$ the {\it
invisible locus\/}, consisting of all empty curves, and let $E$ be
the {\it empty locus\/} consisting of all smooth empty curves.
Obviously $E\subset I$. In fact $E=I - \frak D$. The linear
homotopy argument shows that $I$ is connected (even convex in some
sense), but a priori the hypersurface $\frak D$ could split $I$ in
several pieces.

In general a hypersurface does not need splitting a manifold
(consider e.g. a (pseudo)line in ${\Bbb R}P^2$ or a
meridian/parallel in a torus). In our case $\partial I$ the
boundary (or frontier) of $I$ consists primarily of solitary nodal
curves (principal strata) and further subsequent lower strata.
Hence clearly,  $\partial I \subset \frak D$. But what about
$\frak D \cap \overline{I}$? A priori this does not reduce to
$\partial I$. One can imagine additional nappes of $\frak D$
moving inside $I$, or even that $\frak D$  contains spheroids (or
other closed manifolds) not directly connected to the boundary
$\partial I$ (cf. Fig.\,\ref{ShustinEmpty:fig}d).

Extending our Brusotti argument to $k=m/2$ pairs of empty conics
having transverse complexifications  shows that for any even
integer $m\ge 4$, there is an empty curve $\Gamma_m$ of degree $m$
with 2 conjugate nodes $p,p^\sigma$. Such a curve belongs to
$\frak D \cap {I}$, i.e. is both singular and invisible. Can such
a curve be  connected to $\partial I$ by a path in $\frak D$, and
hence to the ``visible world'' $V:=\vert mH \vert- I$. The answer
would be yes if $\frak D$ is connected. The latter is  a real
algebraic hypersurface,  a priori with several components.

Looking at the (singular) Riemann surface
(Fig.\,\ref{ShustinEmpty:fig}b), we can try to contract
algebraically the anti-oval winding around the middle hole toward
a solitary node. This would give a path as required. This sort of
problem was already discussed at length in another strangulation
section, yet we lack a serious procedure.

Assume now the opposite, i.e., $\frak D$  disconnected in the
sense of having a component inside $I$. The linear homotopy
argument shows that $I$ is convex in the sense that between any
two of its points the projective line joining them has one half
contained in $I$. It follows that $I$ is a contractible manifold!
(Warning since Whitehead 1936 do not draw hastily that $I$ is
homeomorphic to ${\Bbb R}^N$). Such manifolds (more generally
those which are simply-connected, or even under weaker homological
condition) are subsumed to Jordan-Brouwer separation. Under our
supposition that $\frak D$ has some component inside $I$, it would
result a separation of $I$ by $\frak D$. In fact quite
independently of this supposition even if $\frak D$ is connected
there is still a separation.

Of course we need some lemma extending Jordan separation caused by
a manifold, to a separation caused by a stratified variety (not
necessarily smooth). This looks true either by Anschauung (cf.
Fig.\,\ref{Shustin2:fig}d) or by a reduction to the manifold case
by selecting adequately strata as to manufacture first a
topological (but piecewise smooth) manifold out of the strata
(Fig.\,\ref{Shustin2:fig}e). Of course this should rest upon
Brusotti's description of the discriminant.

Whatever the method used we get a morcellation of $I$ by the
discriminant. Quite ironically the linear homotopy argument
(reminded by Shustin) seems to do exactly the opposite job than
its primary intention. More precisely it implies that $I$ is
contractible, hence subsumed to Jordan-Brouwer separation. On the
other hand our Brusotti-type construction shows that $\frak D$
appears inside $I$ (provided $m\ge 4$), hence must divide the
invisible locus. Modulo details, we believe to have proved the
following:

\begin{theorem} {\rm (Revolutionary
if true, but false!)}.---For any even integer
$m\ge 4$ the empty smooth locus (past the
discriminant) is disconnected.
\end{theorem}

If true this would wash up several misconceptions in the
literature, e.g. that the rigid-isotopy type of quartics is
unambiguously determined by the real scheme (this would be false
for the empty scheme). This rigidity is due to Klein 1876
\cite{Klein_1876_Verlauf} and well-known to Russian geometers
(e.g. Rohlin 1978, Viro 1984, 1989, 2008, etc.) Likewise it would
corrupt the same assertion in degree 6, which is included in
Nikulin's theorem (1979 \cite{Nikulin_1979/80}).

Let us again examine our argument to find our probable mistake. It
decomposes in 3 distinct steps. Remember that $I$ denotes the
invisible locus consisting of all curves having empty real locus.

(1) Linear homotopy implies that the invisible locus $I$ is a
contractible manifold.

(2) The discriminant $\frak D$ is visible inside the invisible
locus $I$ for $m\ge 4$. (This follows via simple application of
Brusotti.)

(3) Jordan separation holds true in a contractible manifold (or
more generally a simply connected one). This can be proved in
several ways, either by homology or directly by building a certain
double cover out of the hypersurface by a polarization trick going
back to Riemann (cf. e.g. Gabard 2011
\cite{Gabard_2011-Ebullition}, arXiv, ``Ebullition in foliated
surfaces vs. gravitational clumping'').

{\it Insertion} [06.04.13].---At the risk of killing some
dialectic suspense, there is a 4th issue namely the codimension of
the discriminant inside the invisible locus, as not being 1 but 2
instead!

The step which looks most fallacious is Step (1). The reason could
be the following. While there is between any two $m$-forms $P,Q$ a
path $(1-t)P+tQ$ of $m$-forms  (called the linear homotopy $P\to
Q$), and which sweeps out forms not representing zero  if $P$ and
$Q$ have the same never changing sign, it is not clear that given
(invisible) curves $C,D\in I$ there is always a consistent choice
of sign for representing forms ensuring a global retraction of $I$
to a point. Claiming this amounts finding a section of the evident
(tautological) bundle.

Let ${\cal F}_m$ be the set of all forms (=homogeneous
polynomials) of degree $m$. If we include the zero polynomial this
becomes a vector space, with the space of $m$-tics $\vert m H
\vert$ being its projectivization. Denote by
$$
\pi\colon {{\cal F}}_m \to \vert mH \vert
$$
the corresponding projection.

Choose a basepoint $D$ in $I$ (e.g. the class of the form $Q=
x_0^m+x_1^m+x_2^m$). (This is akin to Fermat's equation
$x^n+y^n=z^n$ except for lacking real points.) $Q$ has positive
sign on ${\Bbb R}P^2$. Given any point $C\in I$ choose a
representing
form $P$ which has positive sign. Then the linear homotopy
$h_P\colon P\to Q$ defined by $h_P(t)=(1-t)P+t Q$ stays in $I$,
which projected down to $\vert m H \vert$ gives a path joining $C$
to $D$. This path is actually independent of the chosen
representative $P$ of $C$ (by a variant of Thal\`es, alias
linearity). Define now
$$
H\colon I \times [0,1] \to I, \quad H(C,t)=\pi (h_P(t) ).
$$
This would be the required contraction (retraction to a point)
showing that $I$ is contractible. However the subtlety is whether
we can choose $P=s(C)$ continuously as a function of $C$. This
amounts asking  if $\pi $ (the tautological projection) admits a
continuous section above $I$. Of course $\pi$
lacks a (global) section by looking at the fundamental group
$\pi_1$ while using functoriality. Indeed the base of the
fibration has $\pi_1={\Bbb Z}_2$, while the total space has
trivial $\pi_1$.

Over the smaller subregion $I$ the situation is less obvious. Can
one compute $\pi_1(I)$? If it is non trivial then we cannot find a
section, and we are annoyed.

Can we construct a section geometrically? We can look at the
counter-image $\pi^{-1}(I)$ interpreted as the cone of forms not
representing zero (so-called {\it anisotropic\/} forms, if we
remember well some highbrow arithmetical jargon\footnote{Indeed,
we remember well, cf. e.g. EDM=Encyclopaedic Dictionary of
Mathematics 1968/87 \cite[p.\,46, Art. 13 G]{EDM_1968/87}
``[\dots] the form $f$ is anisotropic, i.e. the homogeneous
equation $f=0$ has no solution other than zero in $k$''. Or cf.
Serre's ``Cours d'arithm\'etique'', 1970--1977, who seems clever
enough to avoid the jargon, yet speaks of isotropic for quadratic
forms. I don't know who coined the term (in arithmetics), maybe
Minkowski, Hilbert, Weil? Ask a competent arithmetician.}). While
on ${\cal F}_m$ the sign of a form is well-defined at a point, on
$\pi^{-1}(I)$ it is well-defined globally. So our cone
$C:=\pi^{-1}(I)$ splits in two components $C^+, C^-$, each being
connected by the linear homotopy argument.

Now choose the hyperplane $\Pi$ through $Q$ which is orthogonal to
$Q$ seen as a vector. This could give a section.

In fact both $C^+$ and $C^-$ are contractile, being actually
starlike and even convex by the linear homotopy argument. Each of
them is fibred by rays (semi-lines=orbits under scaling by the
positive reals ${\Bbb R}_{>0}$) and the quotient of each of these
cones by the multiplicative group of positive reals ${\Bbb
R}_{>0}$ is naturally identified with $I$.

Abusing  geometric intuition we could nearly conclude that $I$ is
contractible. Yet, this is not so evident as we lack a global
cross-section, e.g. by cutting by a hyperplane selecting globally
a point in each fibres. This looks hazardous, so let us concede
some little algebraic d\'etour or rather homotopy theory
(presumably the quintessence of topology since Jordan 1865, Klein
1882 \cite{Klein_1882}, Poincar\'e 1895, Dehn-Heegard 1907 who
coined the term, and then Brouwer, H. Hopf ca. 1926--30, Hurewicz,
Borsuk ca. 1935, J.\,H.\,C. Whitehead, who else?\footnote{Joke of
Misha Gromov (yet another notorious student of V.\,A.
Ro[k]hlin).}).

The exact homotopy sequence of a fibration ${\Bbb R_{>0}} \approx
fibre \to C^+\to I$ gives
$$
0=\pi_1(fibre)\to \pi_1(C^+) \to \pi_1(I) \to \pi_0(fibre)=0,
$$
and implies that $\pi_1(I)=0$ is trivial. So there is no algebraic
obstruction to find a section over $I$ (but algebra is never
enough to ensure geometric existence!).
Pursuing in that way with the exact homotopy sequence of a
fibration of the early 1940's (Hurewicz, Hopf\footnote{Student of
E. Schmidt, himself student of D. Hilbert. So we are not to far
apart form the 16th problem.}, Stiefel, Eckmann, Steenrod, G.\,W.
Whitehead, Pontrjagin, etc.) we get
$$
0=\pi_i(fibre)\to \pi_i(C^+) \to \pi_i(I) \to \pi_{i-1}(fibre)=0,
$$
and so $\pi_i(I)=0$ for all $i=0,1,2,\dots, \infty$ (modulo
``nihil est infinito''!). Note that $I$ is connected being the
image of the connected set $C^+$. Now our space $I$ is not a bad
one (remember Viro's talk ``Compliments to bad spaces''). More
precisely, $I$ is a manifold (being an open set in the manifold
$\vert m H \vert\approx {\Bbb R}P^N$). This manifold $I$ is metric
moreover, hence it has the homotopy type of a CW-complex in the
sense of J.\,H.\,C. Whitehead\footnote{But widely anticipated by
Poincar\'e, Tietze, Brouwer,  and many others combinatorial
topologists of the early 20th century.} (compare Hanner, Borsuk,
Milnor 1959 \cite{Milnor_1959}, Palais 1962, Gabard 2006/08
\cite{Gabard-2006/08}).
By a theorem of
J.\,H.\,C Whitehead 1949, it follows that $I$ is contractile.

{\footnotesize

\smallskip
{\it Optional Remark (skip since it is not logically
required).}---From the 1940's
(Ehresmann-Feldbau=Laboureur\footnote{Laboureur means nearly
laborieux in French, and was Feldbau's pseudonym to publish
Comptes Rendus notes during  the German occupation of France
(World War II, 1939--45). Alas, it did not helped to save his life
from the Nazi persecutions. Another notorious victim of the
genocide soon afterwards was F. Hausdorff, 1944. Why so much
dramas in the human history is a puzzle to each philosopher.
Materialism, capitalism, caused by the ontological existential
fears ought to be the cause of such disasters. We can only hope
that the Riemann(=woman) surface will quickly lead us to stabler
psychological comforts (immortality, and global resurrection as to
repair such disasters).}-Hopf-Stiefel-Pontrjagin), etc., any
locally trivial fibration over a contractile base (which is
paracompact\footnote{Without this proviso it is false, e.g. the
tangent bundle to the simply-connected Pr\"ufer surface is not
trivial, for otherwise the manifold could be given a Riemann
metric tensor, hence be metricized. Compare Rad\'o 1925
(publishing a contribution of Heinz Pr\"ufer 1922),
Calabi-Rosenlicht 1953, Spivak 1970 (Vol.\,I, Appendix of
Differential Geometry), or ask Mathieu Baillif why. The latter's
e-mail is: labaffle(at)gmail.com}) is globally trivial, hence
admits a continuous section. Applying this to $C^+ \to I$ gives
the required section permitting to contract $I$ via $H$. However
all this optional remark is not really logically required.

} \smallskip

This proves the following:

\begin{theorem}
The space $I$ of invisible curves is
contractible, and so it is separated by the discriminant $\frak D$
in several components. In particular the ``chamber'' of empty
curves $E_m$ is never connected as soon as $m\ge 4$. So call it
rather the empty locus. (Insertion [06.04.13].---This last clause
is probably erroneous.)
\end{theorem}

The determination of the number of components $\iota={\rm card}
(E_m)$ is probably another pleasant game. (Let us guess that
$\iota_4=2$, and $\iota_6=3$?)

This theorem (especially its second clause) contradicts nearly
everything what has been said about the empty locus.
It shows (despite being a pure existence proof using primarily the
exact homotopy sequence of a fibration and Whitehead) that there
are obstructions to rigid-isotopy lying beyond the pure optical
level. It is of course a
marginal contribution to Hilbert's 16th problem, who primarily
asked the right opposite extreme (isotopy classification
especially of $M$-curves). Here we live in the opposite invisible
part of the mushroom (Arnold's prose) of what could be called (by
analogy with Petrovskii 1933/38 \cite{Petrowsky_1933},
\cite{Petrowsky_1938}) $m$-curves, where $m$ stands for Harnack
``minimal'' or minimalist artwork (empty locus like Mark Rothko's
monochromes\footnote{Peintre am\'ericain  d'origine russe (Dvinsk
1903--New York 1970). Il est c\'el\`ebre pour la formule
d'abstraction chromatique qu'il a \'etablie vers 1950.
(Source=Larousse Dictionnary, 1991).}). In some sense our result
of disconnectedness is reminiscent (albeit different in method) to
Marin's  disproof (1979 \cite{Marin_1979}) of the rigidity of
$M$-schemes in degree 7. In both cases the real scheme fails
 determining unambiguously a chamber of the discriminant, and this
in situations where there is no duplication by Klein's types I/II
(what Rohlin 1978 \cite{Rohlin_1978} calls schemes of indefinite
type).

\subsection{Simplifying the previous section:
disconnection of the empty locus via Jordan separation and the
exact homotopy sequence}
\label{Disconnection-of-the-empty-locus:sec}

[27.01.13]
As the former section reflects  our discovery process (as
``meandering'' as it may be) we prefer to keep its shape
unchanged. Since our conclusion   contradicts all what was
asserted about the empty locus $E$ (especially Klein 1876, and
Nikulin 1979), we shall here try to be more formal and direct,
leaving aside historical considerations, and actually simplifying
much the proof (in particular Whitehead's criterion of
contractibility via the vanishing of homotopy groups $\pi_i$ is
not needed).

[28.01.13]
Our intention is to prove the following
disconnection of the empty locus of plane curves of even order
$m\ge 4$.

\begin{theorem}\label{Gabard-anti-Klein-Rohlin:thm} {\rm
(Gabard, 27.01.13, but certainly false)}
The ``empty locus'' $E$ of all real smooth plane curves having
empty real parts of some fixed degree $m\ge 4$ is disconnected. In
other words there exists, for any even integer $m\ge 4$, two empty
smooth curves of degree $m$ which are not rigid-isotopic.
\end{theorem}

We shall emphasize that this conclusion is quite unexpected. It
seems to contradict  much that has been said about rigid-isotopy
of empty (smooth) curves. In particular, it is incompatible with
the assertion  going back to Klein 1876 \cite{Klein_1876_Verlauf}
(see also Rohlin 1978 \cite[p.\,96]{Rohlin_1978}, or Viro's
surveys 1986, 1989, 2008) that the real scheme of a quartic curve
determines uniquely its rigid-isotopy class. (More on this at the
end of this section.) It conflicts also with Nikulin's result
(1979 \cite{Nikulin_1979/80}) that for sextic real curves the real
scheme enhanced by the type data (I/II) of Klein (1876
\cite{Klein_1876}) suffices to determine the rigid-isotopy class.

Hence it is very likely that our theorem contains a serious
misconception, either at the conceptual level of definitions, and
if not so, there must be a bug in the proof below. In fact it is
well known that the empty chamber is connected. Eugenii Shustin
was kind enough to recall us the simple argument of linear
homotopy between two forms of the same sign. This is supposed to
show connectedness, yet exploiting it systematically we arrived
ironically at the opposite conclusion. In part, this discrepancy
is merely a matter of  deciding what we like to call the empty
locus. The linear homotopy argument shows connectedness of what we
call the {\it invisible locus\/} $I$ (consisting of all empty
curves), whereas by the {\it empty locus\/} $E$ we really mean the
sublocus of $I$ consisting of smooth curves. The latter space is
the more relevant one when it comes to problems of rigid-isotopy,
where the game is to travel as much as we can while avoiding the
discriminant (i.e., never strangulate the underlying Riemann
surface).

To avoid any misunderstanding, let us fix our jargon more
precisely. A {\it ternary form\/} is a homogeneous polynomial in
three variables of some degree $m$.  We shall only consider those
with real coefficients, and call them {\it real forms\/}. A {\it
real plane curve\/} is a homothety class of real forms under
scaling of the coefficients. This is nothing else that what
A.~Weil would call a plane curve defined over ${\Bbb R}$. Denote
by ${\cal F}_m$ the set of all real forms of degree $m$, and by
$\vert mH \vert$ the space of all real curves of degree $m$ (the
latter being merely the projectivization of the former). Denote by
$\pi\colon {\cal F}_m \to \vert mH \vert$ the tautological
projection which is an ${\Bbb R}^{\ast}$-bundle over $\vert mH
\vert\approx {\Bbb R}P^N$, where $N=\binom{m+2}{2}-1$.


A plane curve is {\it smooth\/} (or {\it nonsingular\/}) if the
three partial derivatives of any defining form do not vanish
simultaneously on ${\Bbb C}^3-\{0 \}$; else it is said to be {\it
singular}. The set of all real singular curves forms the {\it
discriminant\/} (hypersurface) denoted $\frak D$. Elimination
theory (or better some counting argument) shows the latter set to
be an algebraic hypersurface in the hyperspace $\vert mH \vert$ of
all $m$-tics. Note that a singular real curve may well have a
smooth real locus (in the sense of differential topology), yet it
will then have conjugate pairs of singularities exchanged by ${\rm
conj}\colon {\Bbb C} P^2 \to {\Bbb C} P^2$, $(x_0,x_1,x_2)\mapsto
(\overline{x_0},\overline{x_1},\overline{x_2})$.

A real form is {\it anisotropic\/} if it does not represents zero
(non-trivially), i.e. the sole real solution of the equation
$P(x_0,x_1,x_2)=0$ is $(x_0,x_1,x_2)=(0,0,0)$. This is tantamount
to emptiness of the real locus $C({\Bbb R})$ of the corresponding
curve. Say in this case that the real curve is {\it empty\/} or
{\it invisible}. Intersecting with any line defined over ${\Bbb
R}$, one sees that any odd degree curve
has non-void real locus. Let $I$ be the set of empty (invisible)
curves. This is nonempty iff $m$ is even, and $\pi^{-1}(I)=C$ is
the cone of anisotropic forms. Such a form has a well-defined sign
$\pm$, and accordingly the cone $C$ splits in two halves $C^+,
C^-$  invariant under ${\Bbb R}_{>0}$-scalings. (We overuse the
letter $C$, for being the cone, or the curve but no confusion
should arise.)

\smallskip
\begin{proof}
The proof of our (dubious) theorem
(\ref{Gabard-anti-Klein-Rohlin:thm}) decomposes in 3 short steps.

$\bullet$ {\it Step 1} (Simple-connectivity of the invisible locus
$I$).---We consider the fibration $\pi\colon C^+ \to I$, whose
base is the set of invisible curves, while the total space is the
space of positive-definite anisotropic form. The fibre is the
space ${\Bbb R}_{>0}$ of positive reals.

The space $C^+$ is convex. Whenever we choose 2 points in it, say
$P, Q \in C^+$, the barycentric combination $(1-t)P+tQ$ for $t\in
[0,1]$ belongs to $C^+$. Accordingly $C^+$ is certainly
contractile, and in particular {\it simply-connected\/}. (Perhaps
$C^{+}$, being starlike, is even diffeomorphic to a
genuine cell, but we do not need that presently. This   follows
perhaps from J.\,W. Alexander's lemma on isotopy, ask L.
Siebenmann or A. Marin?)

The first stage of the exact homotopy sequence of the fibering
${\Bbb R}_{>0}\approx F\to C^+\to I$ reads
$$
0=\pi_1(fibre)\to \pi_1(C^+) \to \pi_1(I) \to \pi_0(fibre)=0,
$$
and it follows that the space $I$ is also simply-connected.

$\bullet$ {\it Step 2} (Construction of invisible curves with
singularities).---By a simple application of Brusotti 1921
\cite{Brusotti_1921}, it is easy to construct invisible real
curves $C_m$ of degree $m\ge 4$ having a pair of conjugate nodes,
i.e. ordinary double points (cf. lemma below for details).
(Abstractly, from the Riemann complexification viewpoint,  imagine
a pretzel acted upon by antipody with two handles strangulated to
a pair of points $p, p^{\sigma}$ exchanged by conj.)

$\bullet$ {\it Step 3} (Jordan-Brouwer separation of the invisible
locus $I$ by the discriminant ${\frak D}$).---Paraphrasing Step~2
in our notation, this means that $\frak D \cap I$ is nonempty. The
space, we are really interested in, is the empty locus $E$
consisting of all smooth empty curves. By definition we have
$E=I-{\frak D}$. So the empty locus $E$ arises from the invisible
locus $I$ by removing a certain real algebraic hypersurface (or at
least the portion ${\frak D} \cap I $ visible in $I$). Since the
manifold $I$ is simply-connected (Step~1), it follows from the
Jordan-Brouwer separation theorem that $\frak D \cap I$
disconnects $I$. We conclude that the residual set $E=I-{\frak D}$
has at least 2 components whenever $m=2k\ge 4$.
\end{proof}

We now make more explicit the lemma required in Step~2 (for a
schematic picture in the case $m=4$, cf.
Fig.\,\ref{ShustinEmpty:fig}a.):

\begin{lemma}\label{Brusotti-binodal-invisible-curves:lem}
Given any even integer $m\ge 4$, there exists a degree $m$ real
plane curve which is invisible (empty real locus) with $2$
ordinary nodes exchanged by complex conjugation.
\end{lemma}

\begin{proof}
Take a collection of $k=m/2$ real conics $E_1, \dots, E_k$ (degree
2) each having empty real locus such that the union of their
complexifications has only normal crossings (ordinary nodes). By
Brusotti's theorem (1921 \cite{Brusotti_1921}), we may smooth away
from $E_1\cup \dots\cup E_k$ all pairs of conjugate points safe
one $p, p^{\sigma}$, where $\sigma={\rm conj}$ is complex
conjugation. (Since $k\ge 2$, there is at least 4 nodes on our
configuration of ellipses by B\'ezout.) The resulting binodal
curves is real, invisible (being manufactured by  small
perturbation of an invisible curve). The proof of the lemma is
complete.
\end{proof}

In Step~3 of the proof of
Theorem~\ref{Gabard-anti-Klein-Rohlin:thm}, we use a slightly
extended form of Jordan separation imposed by a variety (possibly
singular) and not just by a manifold. In fact we may imagine that
the structure of the discriminant permits one to deduce a sublocus
of ${\frak D} \cap I$ which is a genuine topological manifold, yet
piecewise smooth (nothing so crazy as Bing-Casson-Freedman). This
would involve aggregating suitably some principal strata of the
discriminant exploiting perhaps Brusotti's description of the
latter, or just general properties of algebraic sets.

It also conceivable that there is a direct proof by applying
directly the homological apparatus involved in Jordan separation
to the case of an algebraic hypersurface. Probably this is already
implemented somewhere (maybe by H. Kneser, Bieberbach, Whitney,
Thom, Milnor, Tognolli, Marin, Bochnak-Coste-Roy, etc.). Alas we
do not know a more precise reference. Evidently the proof of
separation within the simply-connected locus $I$ should just use
some very basic properties of real algebraic hypersurfaces.
Crudely speaking  a real algebraic hypersurface cannot ``stop''
like a manifold with boundary (via the implicit function theorem),
and so really effects a separation in the large (at least within
the simply-connected subregion $I$).
%
%
This elementary property of real algebraic variety was known for
long (e.g. when Zeuthen 1874 \cite{Zeuthen_1874} speaks of a
``branche compl\`ete'', etc.)

\subsection{Rigid isotopy of quartics:
classical sources (Schl\"afli, Zeuthen 1874, Klein 1873--76,
Rohlin 1978)} \label{Klein-rigidity-of-quartics:sec}

[28.01.13] This section  discusses in some more details some
masterpieces of classical literature  conflicting strongly with
our conclusion (Theorem~\ref{Gabard-anti-Klein-Rohlin:thm}).

A first place is Rohlin 1978 \cite[p.\,96]{Rohlin_1978}, who
ascribes the rigid-isotopy classification for $m=4$ to Klein,
while writing the following:

\smallskip

{\small

``\S 4. Isotopy.---4.1. The classical problem. By virtue of the
definition of a real plane projective algebraic curve of degree
$m$, such curves form a real projective space of dimension
$m(m+3)/2$. Singular curves, that is, curves with real or
imaginary singularities, fill out in this space a hypersurface of
degree $3(m-1)^2$, and non-singular curves fill out the
complementary open set, which splits into a finite number of
components\footnote{Read ``chamber'' if you prefer.}. It is clear
that curves that belong to one component have the same real
scheme, that is, the class of all non-singular curves with a given
real scheme consists of whole components. The investigation of
these components is a very old problem, like the investigation of
the classes themselves. It was known more than hundred years ago
that for $m\le 4$ the components coincide with the
classes\footnote{Read ``isotopy classes'', if you like.} (the
least trivial case $m=4$ was considered by Klein; see [4](Klein
1922=Klein 1876 \cite{Klein_1876_Verlauf}), p.\,112). From the
results of the previous section it follows that for $m\ge 5$ this
is not so; for the complex scheme is constant on each component,
but is capable of changing within one class for $m\ge 5$.
[\dots]''

}

\smallskip

As usual Rohlin quotes Klein's GMA=Ges.\,Math.\,Abhandl. (1922),
yet the original source is the 1876 paper ``\"Uber den Verlauf der
Abelschen Integrale bei den Kurven vierten Grades''
\cite{Klein_1876_Verlauf}. Klein's prose is as usual quite magical
(like that of Rohlin), and reads as follows:

\smallskip

{\small

`` Eine wesentliche Eigenschaft dieser Einteilung der Kurven
vierter Ordnung in sechs Arten ist in dem folgenden Satze
ausgeschprochen, der weiterhin eine fundamentale Bedeutung f\"ur
die Tragweite unserer Untersuchungen gewinnt:

{\it Von jeder allgemeinen\footnote{This seems to just mean
non-singular curve, cf. footnote in GMA.} Kurve vierter Ordnung
kann man zu jeder anderen, die derselben Art angeh\"ort, durch
allm\"ahliche reelle \"Anderung der Konstanten \"ubergehen, ohne
da{\ss} bei dem \"Ubergangsprozesse Kurven mit Doppelpunkt oder
gar allgemeine Kurven, die einer anderen Art angeh\"oren,
\"uberschritten zu werden brauchten.\/}

Ein direkter Beweis dieses Satzes hat keine
Schwierigkeit\footnote{If you are called Felix Klein!, else it may
be tricky especially if our previous theorem is right in which
case Klein-Rohlin are false! But of course, it is more likely that
Gabard missed something. [06.04.13] More seriously, it would be
interesting if a detailed account of this direct proof (alluded to
by Klein) has been worked out meanwhile. Some details are perhaps
gleanable from Degtyarev-Kharlamov 2000
\cite{Degtyarev-Kharlamov_2000}.}, aber er ist weitl\"aufig. Es
soll hier um so mehr Abstand genommen werden, als die bei ihm
n\"otig werdenden Betrachtungen mit diejenigen, die im
gegenw\"artigen Aufsatze  zu entwicklen sind, wenig
Beziehungspunkte haben. [So roughly Klein says that there is
little connections between rigid-isotopy and Abelian integrals!]
Dagegen sei angedeutet, da{\ss} man ihn verm\"oge kurzer
Zwischenbetrachtungen f\"uhren kann, wenn man auf fr\"uhere
Untersuchungen von Zeuthen und mir zur\"uckkgreift. Ich habe [in
Abh. XXXV, S.\,24, 25] gezeigt, da{\ss} ein \"ahnlicher Satz gilt
f\"ur die f\"unf Arten, welche man nach Schl\"afli bei den
allgemeinen Fl\"achen dritter Ordnung zu unterscheiden hat. Es hat
dann Zeuthen bewiesen (Math. Ann., Bd.\,7 (1874), S.\,428),
da{\ss} die Arten der Kurven vierter Ordnung den f\"unf Fl\"achen
Arten in sehr einfacher Weise entschprechen. Projiziert man die
$F_3$ von einem ihrer Punkte aus stereographisch auf eine Ebene,
so tritt als scheinbare Umh\"ullung bei den Arten I, II, III, IV
von Schl\"afli eine vierteilige, drei-, zwei-, einteilige Kurve
vierter Ordnung auf. Die Art V ergibt, bei analoger Konstruktion,
je nachdem man den Projektionspunkt auf ihrem unpaaren oder paaren
Teile annimmt, die G\"urtelkurve oder die imagin\"are Kurve.
Umgekehrt kann auch jede Kurve vierter Ordnung aus der
entschpechenden Fl\"achenart in der angegebenen Weise gewonnen
werden. Hierin liegt der vor uns gew\"unschte Beweis. Um ihn
v\"ollig zu f\"uhren, hat man nur noch die Modifikationen zu
untersuchen, welche die scheinbare Umh\"ullungskurve erf\"ahrt,
wenn der Projektionspunkt auf der fest gedachten Fl\"ache beliebig
verschoben wird. Aber auch dieses hat Zeuthen ausgef\"uhrt
[\'Etudes des propri\'et\'es de situation des surfaces cubiques;
Math Annalen, Bd.\,8, (1874/75).]''

}

Let us summarize Klein's proof of the rigid-isotopy of $C_4$'s:
first Schl\"afli in 1863 \cite{Schlaefli_1863} found five isotopy
class of real cubic surfaces $F_3$, and Klein showed them to be
rigid-isotopic (in Klein 1873
\cite{Klein_1873-Uber-Flächen-dritter-Ordn}). Then Klein exploits
the yoga of Zeuthen (which goes back to Geiser, compare Zeuthen
1874 \cite[p.\,428]{Zeuthen_1874}) of looking at the apparent
contour of the $F_3$ (projected from a point on the surface) to
get a $C_4$ (all of them arising so). In fact the $F_3$ with 2
components produces both the G\"urtelkurve or the empty $C_4$
depending on whether the center of vision is located on the
pseudo-plane or on the spherical component. This is easily
visualized if one imagine  $F_3$,  a small perturbation of a
sphere union an equatorial plane (Fig.\,\ref{Zeuthen-Klein:fig}).
Note that the empty apparent contour arises from a phenomenon of
total reality of the bundle of lines through the spherical
component (this being again reminiscent of ``Ahlfors'').

\begin{figure}[h]
\centering
\epsfig{figure=Zeuthen-Klein.eps,width=102mm} \vskip-5pt\penalty0
  \caption{\label{Zeuthen-Klein:fig}%
  Geiser-Zeuthen-Klein trick
  of the visual contour of a cubic surface}  \vskip-5pt\penalty0
\end{figure}

So a smooth $F_3$ with a marked point gives rise to a $C_4$, and
all quartics arise so. One must still study the r\^ole of the
marked point, and  there is not just the five classes of
Schl\"afli-Klein but one more due to the marking (being either on
the pseudo-plane or the spheroid). Hence Klein's argument looks
quite convincing but we should understand the details more
precisely. Maybe there is little gap in this proof when setting up
explicitly the $F_3\leftrightarrows C_4$ correspondence. Is the
contour apparent of a smooth $F_3$ always a smooth $C_4$ (and
viceversa a singular cubic to a singular quartic), so that the
Zeuthen-Klein correspondence really sets up a dictionary between
the corresponding discriminants (hence rigid-isotopy classes). In
particular to what sort of $F_3$ corresponds the empty quartic
with a pair of conjugate nodes. Those are the essential guys which
in our Theorem~\ref{Gabard-anti-Klein-Rohlin:thm} causes the
disconnection of the invisible locus $I$.

Further as the center of projection is moving in 3-space the
ZK-corres\-pondence is somewhat non-canonical, i.e. there is not a
fixed ${\Bbb P}^2$ on which to project. So the reduction proposed
by Klein is perhaps foiled somewhere, at least requires to be
modernized (and detailed) seriously. Naively  our
Theorem~\ref{Gabard-anti-Klein-Rohlin:thm} (if  correct) could be
an obstruction to completing  Klein's proof. If the
ZK-correspondence is sound, it could be
 that Klein's 1873 rigid-isotopy classification of $F_3$'s
(cubic surfaces) is foiled.

\subsection{Fixing the paradox}

[29.01.13] Can it be that the discriminant $\frak D$ while
penetrating inside the invisible locus $I$ appears there with
(real) codimension 2 hence without separating $I$?

Recall that real loci of algebraic hypersurfaces (=primals in old
British jargon, e.g. Semple-Roth) may look anomalously small. For
instance a solitary node on a plane cubic curve is merely an
isolated point of real dimension $0$ (hence real codimension $2$).

This phenomenon could foil Step~3 of the proof of
Theorem~\ref{Gabard-anti-Klein-Rohlin:thm}. However via Brusotti
description of the discriminant one could still hope that ${\frak
D} \cap I$ has real-codimension $1$. Naively the principal stratum
corresponding to a curve with a conjugate pair of nodes looks at
first of codimension 2 because there is two nodes. However by the
reality condition one of them is forced and so it is really one
closed point in the sense of Grothendieck's schemes. Without
Grothendieck such arguments of reality counting also abound in
Klein, e.g. when it comes to coverings of the Riemann sphere with
complex conjugate ramification. Here a complex ramification point
count for 2 real dimensions, while a real branch point affords one
freedom parameter. Yet under the symmetry condition both cases
actually contribute to the same. It is with this sort of argument
that Klein managed to compute the dimension of the moduli space of
real curves by aping what did Riemann over the complexes.

Alas one can argue that the stratum of the complexified
discriminant $\disc(\CC)$ with two nodes $x,y$ nearby $p,p^\sigma$
has geometric codimension 2 over the complexes, and it would
follow that the real locus $\frak D$ has at most real codimension
2 when $x,y$ lye symmetric under $\sigma=$conj. Another way  is to
start the dimension count of the discriminant from the scratch. So
we look at plane curves with a node or a higher singularity marked
on it. This gives an incidence variety $(C,p)$ consisting of
curves $C$ singular at $p$. Saying that $p$ is singular of $C$
amounts (via the Euler relation) to say that the 3 partials of a
defining equation vanish at $p$, yielding 3 linear conditions on
the coefficients. So looking at the projections $p\leftarrowtail
(C,p) \mapsto C$ we see fibres of codimension $3$, while moving
the point $p$ gives codimension $1$ (hypersurface) in the space of
curves.

Adapting this dimension count near our binodal empty curve $C$
shows that we have a pair of projections $p\leftarrowtail (C,p,
p^{\sigma}) \mapsto C$. Imposing a singularity at $p$ gives 3
linear conditions, but the point $p$ being any imaginary point of
the plane its location depend upon 2 complex parameters (4 real
parameters). When $C$ is defined over ${\Bbb R}$ the 3 equations
for the partials are again linear (in the coefficients), but
involve complex constants. Thus splitting into real and imaginary
parts yields twice so many linear equations, hence 6 of them.
Those being satisfied then $p^\sigma$ is also a node by symmetry,
and so the variety ${\frak B}$ of (binodal) curves having a
conjugate pair of nodes has real (co)dimension $4-6=-2$ in the
space of curves. In that case it seems that the principal  stratum
${\frak D}\cap I$ consisting of empty curves with a pair of
conjugate node has also real codimension 2 in $\vert mH \vert$,
and this would foil our theorem
(\ref{Gabard-anti-Klein-Rohlin:thm}).

For short let us call ${\frak D}\cap I$ the {\it invisible
discriminant}, abridged $I$-discriminant. The above argument has
to be polished by checking that the 6 linear equations are really
independent conditions. The main issue is therefore to calculate
the dimension of the invisible discriminant. If it has codimension
2 then $E=I-\frak D$ is connected, and so Klein-Rohlin were right.
(Otherwise, if of codimension $1$, then it separates and
Klein-Rohlin were wrong.)

If the $I$-discriminant has codimension 2, it is like a knot, and
it seems interesting to compute its fundamental group (of its
complement), and to look at the Picard-Lefschetz monodromic
transformation arising when winding once around a meridian of this
$I$-discriminant. Before adventuring we should first solve (more
rigorously) the dimension problem of the $I$-discriminant. This
must surely be done in Brusotti 1921 \cite{Brusotti_1921} (and
known to Gudkov, maybe in the 1974 survey \cite{Gudkov_1974/74}).
Perhaps this was already known in  the era of
Zeuthen-Klein-Harnack-Hilbert.

\subsection{The invisible discriminant has codimension $2$}
\label{invisible-discriminant-codim-2:sec}

[29.01.13] The goal of this section is to resolve our paradox,
that the locus of empty smooth curves is disconnected (violating
thereby well assessed knowledge of Klein, Rohlin-Nikulin, etc.).
It seems that our sole mistake was based on the linguistical
misconception of thinking that the discriminant-hypersurface  is a
hypersurface (throughout)! Names and terminologies are often
misleading in mathematics. In fact we shall try to convince that
inside the invisible locus (of curves with empty real loci) the
discriminant has only codimension 2, hence too small to effect any
Jordan-Brouwer separation. Once this is observed this raises some
little questions about knowing which chambers residual to the
principal strata of the discriminant (where it has really of
codimension 1) contains such smaller strata of codimension 2. In
more geometric terms, this amounts essentially  deciding which
smooth curves can acquire a pair of conjugate nodes.

Fix some even integer $m=2k$, and consider only real curves of
fixed degree $m$. Let $\vert mH \vert\approx {\Bbb R}P^N$ be the
corresponding parameter space of real $m$-tics. In this space we
pay special attention to the space $I$ of {\it invisible curves\/}
(those with empty real locus). This is clearly an open set in
$\vert mH \vert$.

Let $\frak B$ be the variety of invisible real plane curves with
at least one
singular point (hence necessarily at least a pair thereof). We
have $\frak B = {\frak D}\cap I$, the so-called {\it invisible
discriminant} (abridged $I$-discriminant).

\begin{lemma}\label{invisible-discriminant-codim-2:lem}
The $I$-discriminant has (real) codimension $2$ in the hyperspace
of all curves.
\end{lemma}

\begin{proof}
Consider the incidence relation $B=\{ (C,p) \colon C \in I, p \in
Sing C \}$. We have natural projections
$$
{\frak B}\buildrel{\pi_1}\over\longleftarrow B
\buildrel{\pi_2}\over{\longrightarrow} {\Bbb P}^2({\Bbb C}).
$$
First study the fibre of the second projection $\pi_2\colon
(C,p)\mapsto p$ . This amounts to look at all curves having a
prescribed singularity at $p$ an imaginary point. This imposes 3
linear equations (vanishing of the 3 partials, which suffices by
Euler equation for the point to be on the curve). Splitting in
real and imaginary parts gives 6 linear conditions, which looks
linearly independent. So the fibre $\pi_2^{-1}(p)\approx (I \cap
{\Bbb R}P^{N-6})\times\{p\}$. Since $\pi_2$ is surjective onto
${\Bbb P}^2({\Bbb C})-{\Bbb P}^2({\Bbb R})$ it follows that the
real dimension of $B$ is $\dim_\RR B= 4+(N-6)=N-2$. As the first
projection $\pi_1$ is generically 2-to-1 (or finite-to-one except
over special curves with multiple irreducible components), it
follows that $\dim_\RR {\frak B}=N-2$, hence of codimension 2 in
$I$ (or in $\vert m H \vert$).
\end{proof}

More generally,  imagine a visible curve (i.e. $C_m(\RR)\neq
\varnothing$) acquiring a conjugate pair of nodes. This will not
affect the real scheme (=soft isotopy class of the embedding
$C_m(\RR)\subset \RR P^2$). A priori inside each ``class'' can
penetrate a portion of the discriminant of codimension 2.

Let us be more formal.
Split the discriminant $\disc=\disc^+\sqcup \disc^{-}$ in two
parts depending on whether $Sing C$ contains a real point or not.
Precisely define $\disc^{-}$ as the set of singular curves lacking
real singularities.

The argument of the above lemma shows that $\disc^{-}$ has real
codimension two. We call it hence the {\it hypo-discriminant\/}.
The set $\disc^+$ is defined as its complement, i.e.
$\disc^+=\disc- \disc^{-}$. The latter has codimension 1 by a
variant of the above argument. (Indeed a real singularity imposes
3 linear conditions, but moving the point create 2 dimensions,
whence the defect of $-1$.) We call $\disc^+$ therefore  the {\it
hyper-discriminant}.

The assignment  $\vert mH \vert -\disc \to \frak S$ of the real
scheme to a non-singular equation is more generally defined on the
larger space $\vert m H\vert-\disc^+$ residual to the
hyper-discriminant, while being locally constant there. Crudely
put, in problems of rigid-isotopies the hypo-discriminant can be
neglected (being only of codimension 2, hence effecting no
additional separations). However when we would like not only to
study the connectivity of the chambers but also their topology
then the hypo-discriminant ought to be considered again.

A first question is whether any chamber residual to the
hyper-discriminant (abridged {\it hyper-chamber})
intersects the hypo-discriminant. (This is true for the empty
chamber residual to $\disc^+$, by our Brusotti-style lemma
\ref{Brusotti-binodal-invisible-curves:lem}.)

The general problem looks again to involve a contraction principle
of Riemann surfaces, now under a symmetric pair of vanishing
cycles. For $M$-curves, there seems to be a topological
obstruction, since strangulating two imaginary cycles $\beta,
\beta^{\sigma}$ causes a disconnection of the Riemann surface in 2
algebraic pieces $C_m\to C_k\cup C_l$ of degree $k, l$ (hence
cutting themselves in $k\cdot l$ points by B\'ezout). But $C_k,
C_l$ cuts transversally in $2$ points only, hence $k=2$, $l=1$ (up
to renumbering). Hence, this eventuality can only occur for $m=3$
(cubics), where it does occur when an $M$-cubic degenerates to
$E_2\cup L$ a conic union a disjoint line. This curve $E_2\cup L$
belongs to the hypo-discriminant since it lacks real
singularities. ([06.04.13] Further the dimension of such split
cubics is $5+2=7$, which is indeed of codimension 2 in the
hyperspace of cubics of dimension $\binom{3+2}{2}-1=\frac{5\cdot
4}{2}-1=9$.)

Given a curve in the hypo-disc $\disc^-$, it seems likely that by
genericity we may assume the latter to be a binodal curve with a
conjugate pair of nodes. By Brusotti 1921 \cite{Brusotti_1921},
smooth them away. Interpreting the process backward in time we see
a pair of imaginary conjugate vanishing cycle $\beta,
\beta^{\sigma}$ strangulating toward the nodes $p,p^{\sigma}$. If
this argument holds true we see that each hypo-discriminantal
component gives rises to a bistrangulation along imaginary cycles.
In particular:

\begin{lemma}
Safe for $m= 3$, the hyper-chamber of an $M$-curve of degree $m$
never contains the hypo-discriminant.
\end{lemma}

Perhaps this is  the sole obstruction, in the sense that any other
hyper-chamber (than those of $M$-curves) intersects the
hypo-discriminant. In fact there is perhaps still such an
obstruction for $(M-1)$-curves. Naively the latter look like an
$M$-curve safe that one ``oval'' is masked. Still if we imagine
the corresponding symmetric Riemann surface it seems that a pair
of imaginary cycles must divide.

[30.01.13] So we are led to the following general topological
question:

\begin{ques}\label{bicycle-existence:ques} (Existence of bicycles).---{\rm Given a symmetric surface in the sense of Klein 1876 (i.e. an
oriented closed surface $X$ with an orientation reversing
involution $\sigma$), when is it possible to find a pair
$\beta,\beta^{\sigma}$ of imaginary cycles(=Jordan curves) such
that $\beta \cup \beta^{\sigma}$ does not divide $X$?} Here
imaginariness means that $\beta$ has no  point fixed under
$\sigma$. It is also required that $\beta$ is disjoint from its
conjugate $\beta^{\sigma}$. We call such a pair an (imaginary)
{\it bicycle}.
\end{ques}

As well-known since Klein 1876 \cite{Klein_1876}, symmetric
surfaces are either ortho- or diasymmetric (equivalently dividing
or not, or of type~I resp. II). The type together with the number
$r$ of ``ovals'' (pointwise fixed circuit) and the genus $g$ fixes
the equivariant topology of a symmetric surface.

As to our  question it is plain that in the type~I (dividing) case
then there is a bicycle provided the surface is not maximal
$r=g+1$ (``$M$-surface''). In the maximal case such a bicycle is
not available, because the quotient $X/\sigma$ is planar hence
schlichtartig (i.e. divided by any Jordan curve). Literally
``schlichtartig'' means planar like.

Incidentally recall the implications:

\begin{lemma} \label{Gabard-five-lemma:lem}
``simply-connected$\Rightarrow$schlichtartig$\Rightarrow$
orientable'' for a topological surface.
\end{lemma}

\begin{proof} Exaggerating a bit the only rigorous proof, we are
aware of, uses the five lemma
and homology, cf. e.g. Gabard-Gauld 2011
\cite[Lemma~4.17]{Gabard-Gauld_2011-Dynamics}, Dynamics of
non-metric manifolds. It works universally without having even to
assume the surface metrizable. Hausdorffness is however crucial
(consider a branched plane and in it a Jordan curve on the upper
sheet, then you can evade via the lower sheet so as to reach the
outside of the upper sheet.)

The five lemma implies that when we have a Jordan curve $J\subset
U \subset M$ included in two nested spaces, then if it divides the
large space $M$,  it must divide the small one $U$, and viceversa
provided $H_1(U)\to H_1(M)$ is onto. This latter fact prompts the
first implication ``$\Rightarrow$'' by taking $U$ the tubular
neighborhood of $J$ (which must be trivial since otherwise there
would be an indicatrix-reversing loop violating the assumption
$\pi_1(M)=0$). The second implication ``$\Rightarrow$'' follows
from the first fact, namely a division in the large $M$ implies a
division in the small $U$, hence is particular of the tubular
neighborhood which must therefore be trivial.\end{proof}

When applied to the quotient of a symmetric surface,
Lemma~\ref{Gabard-five-lemma:lem} gives for the latter:
``ortho-sphere$\Rightarrow M$-surface$\Rightarrow$ dividing''.

Note also that if there is a bicycle on $(X,\sigma)$, then its
projection in the quotient $\bar X=X/\sigma$ is a cycle $\bar
\beta$ interiorly traced which does not divide and preserves the
indicatrix (=local orientation). (The non-division just follows
from the fact that the image of the connected set $X-(\beta \cup
\beta^{\sigma})$ has to be connected.)

It remains to answer our question (\ref{bicycle-existence:ques})
in the diasymmetric case. The lowest diasymmetric case $r=0$ is
easy since then $(X,\sigma)$ may be visualized in 3-space as a
pretzel invariant under central symmetry (antipody). Hence a
bicycle exists provided the genus $g\ge 2$. (For $g=1$ there is a
pair of cycle exchanged, but collectively they do divide. It
cannot be otherwise by Riemann's definition of the genus.)

\begin{figure}[h]
\centering
\epsfig{figure=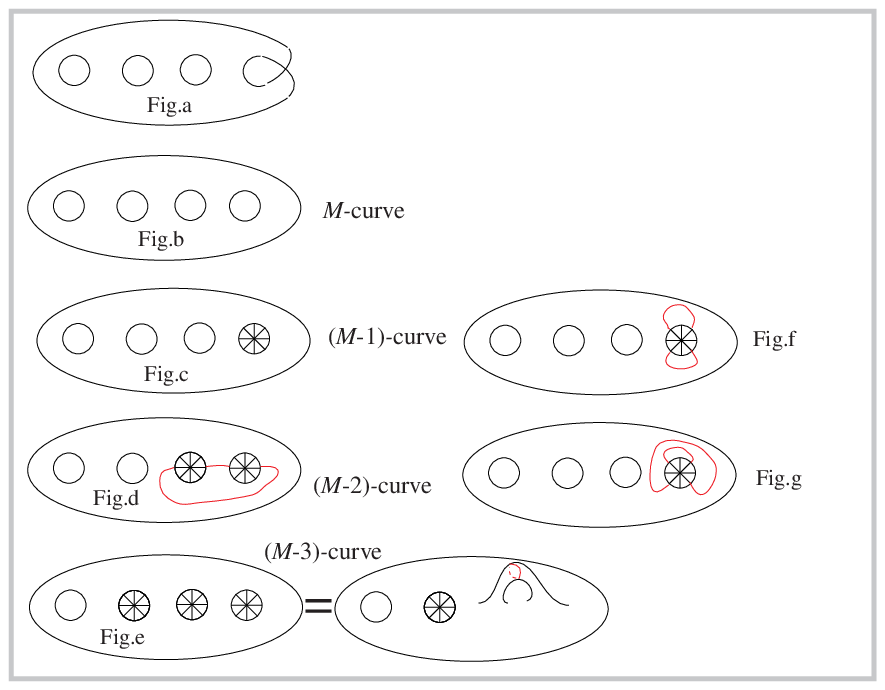,width=102mm} \vskip-5pt\penalty0
  \caption{\label{Kervaire:fig}%
  Depiction of Klein's diasymmetric surfaces
  \`a la Weichold-Kervaire vs. von Dyck}  \vskip-5pt\penalty0
\end{figure}

For the other cases there are several models. One way due to
Klein-Weichold-Kervaire (private communication of the latter in
1999) amounts to look at a M\"obius band embedded in 3-space (make
holes in it) and look at a thickening of the normal bundle of
thickness vanishing along the boundary
(Fig.\,\ref{Kervaire:fig}a).

Alternatively, we may start from the Harnack-maximal case
visualized as a planar membrane with $r=g+1$ contours
(Fig.\,\ref{Kervaire:fig}b), and kill successively contours by
cross-capping them (\`a la von Dyck 1888
\cite[p.\,479]{von-Dyck_1888}, another of Klein's student). This
operation does not alter the Euler characteristic $\chi$, and so
keep  the genus of the double unchanged as $\chi (X)=2 \chi
(X/\sigma)$. The symmetric surface is constructed abstractly via
the usual process of the double orientation cover (without
duplication of the boundary points by local orientations).

So imagine a disc with $g$ holes while cross-capping them
successively. If no cross-cap we have an $M$-curve\footnote{This
is a slight abuse of language to suit  Russian jargon (coined by
Petrovskii 1938).}, if one cross-cap an $(M-1)$-curve, if two
cross-caps an $(M-2)$-curve of type~II, etc. As soon as there is
$2$ cross-caps, connect them by a path and closing it back
(Fig.\,\ref{Kervaire:fig}d) gives a cycle $\bar \beta$ which
preserves the indicatrix (as it traverses twice the cross-caps)
and which does not divide $\bar X$ (because its apparent inside is
connected with the outside via the cross identifications). Lifting
$\bar \beta$ to $X$ gives the desired bicycle. Another, but
slightly weaker argument, is that as soon as there is 3 cross-caps
available, they can be traded against one cross-cap and one handle
(as both contribute identically to the Euler characteristic). Then
it is enough to take the meridian (or parallel) of that handle
(Fig.\,\ref{Kervaire:fig}e).

The case of where there is only one cross-cap is a bit more
tricky, and it seems that we cannot find a nondividing cycle
$\bar\beta$ which is indicatrix-preserving. (Trace a picture which
can be either like a figure ``8'' crossing the cross-cap
(Fig.\,\ref{Kervaire:fig}f) or like a figure ``$\omega$'' with
extremity linked together (Fig.\,\ref{Kervaire:fig}g). This is
abstractly just the figure 8, except that the one loop envelopes
the other one. In both cases it is seen that a division is
produced.)

Here is the obstruction:

\begin{lemma}
An $(M-1)$-surface (with $r=g$) cannot have a bicycle.
\end{lemma}

\begin{proof} Since the bicycle $\beta\cup \beta^{\sigma}$ does not divide the surface $X$,
its image $\bar \beta$ in the quotient $\bar X$ does not divide
it. The covering $X-{\rm Fix}(\sigma)\to \bar X-\partial \bar X$
 restricted to the complement of $\bar \beta$
 shows that $X-{\rm Fix}(\sigma)-(\beta \cup \beta^{\sigma})$ has
at most $2$ components (exchanged by $\sigma$). Yet it suffices to
add one oval (of ${\rm Fix}(\sigma)$) to make it connected. Hence
we have $(g-1)+2=g+1$ retrosections not disconnecting the surface,
violating Riemann's definition of the genus.
\end{proof}

In summary we have proven:

\begin{lemma}
A symmetric surface admits a bicycle iff $g\ge 2$ and  $r\le g-1$.
In other words iff it is not an $M$-surface nor an
$(M-1)$-surface.
\end{lemma}

Via Brusotti this seems to afford obstructions to the presence of
the hypo-discriminant in certain hyper-chambers. More precisely:

\begin{prop}
Inside a pre-maximal hyper-chamber (i.e. $r\ge g$) the
hypo-discriminant is vacuous. In particular a premaximal curve
(i.e. an $M$- or an $(M-1)$-curve) cannot acquire a conjugate pair
of nodes by continuous variations of its coefficients (among
smooth curves safe for the extremity). [Warning: the case $m=3$ is
the sole exception.]
\end{prop}

More
risky is the (converse) assertion that via a suitable (but very
hypothetical) contraction principle the topological presence of a
bicycle suffices to create an algebraic deformation toward a curve
$C_m$ with an imaginary pair of nodes. If optimistic about the
freedom of the joy-stick this
supports the:

\begin{conj}\label{hypo-discriminant:conj}
Any (smooth real plane) curve $C_m$ which is not premaximal can
acquire a conjugate pair of nodes (bi-node) via continuous
deformation among smooth curves safe for its extremity. In
particular the hypo-discriminant appears in all the corresponding
(``ante-maximal'') hyper-chambers.
\end{conj}

This is another  large-deformation principle a bit akin to the
Itenberg-Viro contraction conjecture
(\ref{Itenberg-Viro-contraction:conj}) for empty ovals. To get
serious prohibitions one would perhaps even require a strengthened
collective form of it. For instance if $r=2$ and if we are
dividing, we could contract several bicycles which collectively
split the Riemann surface. Browsing through increasing degrees
$m=3,4,5,6,7, \dots$ gives the genus $g=1,3,6,10,15,\dots$. By
virtue of Klein's congruence $r\equiv_2 g+1$ we look especially at
$m=4$ or $m=7$ (or $m=3k+1$). When $m=4$, we find no obstruction
to the splitting (since $g=3$ and we have 4 vanishing cycles
contracting to $4=2\cdot 2$ points in agreement with B\'ezout, see
Fig.\,\ref{Pretzel:fig}a). If $m=7$, then $g=15$ and so there is
$8$ bicycles (Fig.\,b) which strangulated toward nodes gives 16
(simple) intersections between both pieces of the degeneration
$C_7\to C_k\cup C_l$, where $k+l=7$. Testing all values
$(k,l)=(1,6),(2,5),(3,4)$ gives always the wrong number of
intersections $k\cdot l=6,10,12$ never equal to $16$. This
contradiction with B\'ezout reproves that a dividing septic cannot
have $r=2$. (Of course this is best proved as a consequence of
Rohlin-Mishachev's formula).

\begin{figure}[h]
\centering
\epsfig{figure=Pretzel.eps,width=122mm} \vskip-5pt\penalty0
  \caption{\label{Pretzel:fig}%
  Contracting bicycles (and orthocycles) on orthosymmetric surfaces} \vskip-5pt\penalty0
\end{figure}

Extending this to all degrees requires another form of
contraction. In the case of $C_5$'s: then assuming $r=1$,  there
is 3 bicycles and 1 ortho-cycle (Fig.\,c). We have a splitting
$C_5\to C_k\cup C_l$, and as $g=6$ we have 7  intersections in
$C_k\cap C_l$ after strangulation of the Riemann surface. But this
is never equal to $k\cdot l$ for $(k,l)=(1,4), (2,3)$. So this
would prove again that a quintic  with one circuit ($r=1$) cannot
be dividing. Compare (\ref{Klein-Marin-quintic:lem}) for another
proof.

{\it Insertion} [06.04.13].---Applying the same method to a
quintic with $r=3$, while imagining the underlying Riemann surface
orthosymmetric and spliced by 2 bicycles and 1 orthocycle
(Fig.\,d), the strangulation process leads to 2 algebraic pieces
$C_k\cup C_l$ intersecting in 5 points. This is again not of the
form $k\cdot \ell$ (equal as above to $4$ or $6$). This reasoning
proves that a quintic with $r=3$ cannot be dividing, which is {\it
nonsense\/} (remind the deep nest and its total reality). So the
methodology (of such imaginary contractions) appears jeopardized.
It seems at first that by ``rotating'' the pretzel of genus $6$ as
to make the five cycles into reals circuit prompts a corruption of
Itenberg-Viro (\ref{Itenberg-Viro-contraction:conj}) by the same
device. This is not so because one of the circuit is a pseudoline
(since we are in degree $m=7$ odd). But of course we could imagine
an example in even degree. In fact what protects a direct
corruption of the Itenberg-Viro conjecture
(\ref{Itenberg-Viro-contraction:conj}) is Rohlin's formula which
in case of no-nesting forces $r=k^2$, and intersecting the half of
the strangulated Riemann surface (now of degree $k$ since
exchanged by Galois) is concomitant with B\'ezout.

More generally, one can perhaps by this method get another
derivation of the Rohlin-Marin inequality $r\ge m/2$ for a plane
dividing curve of degree $m$.

One could also imagine more radical degeneration by a bicycle
$\beta, \beta^{\sigma}$ such that already $\beta$ divides
$C_m({\Bbb C})$ (fig.\,f). A such is easy to visualize in the
dividing case and would separate all imaginary handles from the
real contours. However in the case of a $C_4$ of type~I, such a
pair $\beta,\beta^\sigma$ would strangulate the surface of
genus~$3$ in three pieces of genus 1, so the degree must be at
least $3+3+3=9$, which is much greater than $4$. Hence such
contractions are unlikely to exist algebraically.

[31.03.13] When $m=4$, the above conjecture
(\ref{hypo-discriminant:conj}) looks trivial, e.g. because all
ante-maximal schemes $r\le M-2=2$ admits realization as pair of
conics (either nested or disjoint).

\subsection{Rigidity index}

[31.03.13] Another naive remark concerns the rigidity of the
``one-oval scheme'' $1$. Once the empty scheme is known to be
rigid, then via the contraction conjecture CC
(\ref{Itenberg-Viro-contraction:conj}), the one-oval scheme ought
to be rigid as well. Naively via CC one could pursue inductively
and all schemes would be rigid (which is not true as best and
first shown by Rohlin via Klein's type for $m\ge 5$). So there are
subtle obstructions coming from separation between chambers at the
next level. Despite such difficulties we call this method the {\it
rigidification procedure by reduction to the empty chamber}.

So we start from the empty chamber $0$, and then there is the
chamber $1$ (which should be still connected). Then ``attached''
to this there is the ``chamber'' $\frac{1}{1}$ and $2$, etc. Of
course here ``chamber'' should rather be ``isotopy class'' and it
should be proved that such schemes are rigid, i.e. that their
respective isotopy classes correspond to a unique chamber of the
discriminant.
%
%
%
%
\def\rig{r}
As usual in mathematics, when we are unable to prove something we
just introduce a:

\begin{defn} {\rm The {\it bifurcation index\/} $r(m)$ in degree $m$
is the smallest integer $r=\rig(m)$  such that there is a real
scheme of degree $m$ with $r$
real branches which is non-rigid, i.e. represented by 2 real
curves of degree $m$ which are not rigid-isotopic. It is set equal
to $+\infty$ if all schemes are rigid.  If finite, and diminished
by one unit it could be called the {\it rigidity index}, since
below it all schemes would be rigid. (Our terminology is a bit
awkward, because a high rigidity index truly means that the video
game is flexible.)}
\end{defn}

Basically $\rig(m)$ measure the critical level at which the above
rigidification algorithm  fails surely. Very little is known on it
as exemplified by the (still open) conjecture on the rigidity of
the one-oval scheme (\ref{OOPS:one-oval-rigid-isotopic:conj}),
which traduces into the assertion $\rig(m)\ge 2$ for all even $m$.
Even this modest estimate is pure speculation to present
knowledge.


It is trivial that $\rig(1)=+\infty$, $\rig(2)=+\infty$ (because
$PGL(3,\RR)$ acts transitively on lines or conics provided the
latter are defined by quadratic forms with the same signature).
$\rig(3)=+\infty$, i.e. all schemes of degree 3 are rigid is
already somewhat more sophisticated since there are moduli. Yet
this must follow either form Newton-Pl\"ucker or from the theory
of elliptic functions (Euler-Legendre-Abel-Weierstrass, etc.) It
is probably not completely trivial to write down the details, but
looks evident if we keep in mind a reduction to the Weierstrass
normal form $y^2=(x-a_1)(x-a_2)(x-a_3)$, where the distinct $a_i$
can either be all real or two of them imaginary conjugate.

The result of Schl\"afli-Zeuthen-Klein (cf. Klein 1876
\cite{Klein_1876_Verlauf}) implies that $\rig(4)=+\infty$. This is
already less evident, and was first proved by Klein via cubic
surfaces as discussed in
Sec.\,\ref{Klein-rigidity-of-quartics:sec}.

For $\rig(5)$ we have a scheme of indefinite type (namely $4\sqcup
J$) which is elementary to find (see
Fig.\,\ref{Gudkov-Table-quintic:fig}) and first described in
Rohlin's era (cf. e.g. 1974 \cite{Rohlin_1974/75} and 1978
\cite{Rohlin_1978}). It suffices to smooth a pair of conics plus a
line
in two different ways as indicated on the left of
Fig.\,\ref{Indefiniteodd:fig}. It follows that $\rig(5)\le 5$. By
Rohlin's inequality ($r\ge m/2$ if type~I) and Klein's congruence
holding right above, the above example is clearly minimal to
detect an obstruction to rigid-isotopy via Klein's types (see also
the argument in Sec.\ref{quintic-table-Klein-Gudkov:sec}). Thus it
is fairly clear that $\rig(5)=5$, but we know no proof without
appealing to Kharlamov's version (1981/81
\cite{Kharlamov_1981/81}) of Nikulin's classification in degree 6
(cf. next item).

Looking at the Gudkov-Rohlin table (Fig.\,\ref{Gudkov-Table3:fig})
of sextics, it is clear that $\rig(6)=5$. The proof of this rests
the deep result of Nikulin 1979 \cite{Nikulin_1979/80} that the
type enhanced real scheme suffices to encode the rigid-isotopy
class. It would be interesting to know if the CC conjecture is
able to reprove the rigidity of all sextic schemes lying below $r<
5$ the bifurcation index. Of course Nikulin tells much more.

In view of the Rohlin-Marin inequality $r\ge m/2$ for a dividing
curve,  Klein's types afford no obstruction to rigid-isotopy for
curves with few branches. A naive optimist can expect no
bifurcation below this value. By a bifurcation of a scheme we
simply mean it being stretched apart in two chambers of the
discriminant. So $m/2$ is a sort of ebullition temperature, below
which everything is frozen, i.e. only type~II schemes are
represented apart from the deep nest scheme with $r=m/2$. The
latter scheme is (pure) of type~I by the simplest form of total
reality, hence does not cause a type bifurcation, while  being
actually rigid by Nuij's theorem (1968 \cite{Nuij_1968}).

The critical temperature $r=m/2$ can be augmented by unit, because
the next $r$ is forced belonging type~II by Klein's congruence.
(All this extends to the case where $m$ is odd by taking as
critical temperature $r=(m+1)/2$ the number of branches of the
deep nest).

Hence all
schemes with $r\le [(m+1)/2]+1=[(m+3)/2]$ are necessarily of
type~II, safe for the deep nest. In particular there is no
indefinite schemes below this level, and the first such indefinite
scheme is expected to be found at height $[(m+1)/2]+2=[(m+5)/2]$.
(The height of a scheme is merely its number of components, a
jargon suggested by the diagrammatic of the Gudkov table
Fig.\,\ref{Gudkov-Table3:fig}.)

\def\indef{\beta}

\begin{defn} {\rm The {\it type bifurcation index} (or
just {\it indefiniteness\/})
$\indef=\indef(m)$ is
%
the minimal height of an indefinite scheme. It is set equal to
$+\infty$ if all schemes of degree $m$ are definite (i.e. either
of type~I or II in the sense of Rohlin 1978).}
\end{defn}

For $m\ge 5$ it seems evident that indefinite schemes always
exist, but this requires some proof. (This and more will follow
from Figs.~\ref{Indefinite:fig} and \ref{Indefiniteodd:fig}
below.)

The following is all what can be said at first sight:

\begin{lemma} We have $1 \le \rig(m)\le \indef(m)$, and
$\indef(m)\ge [(m+5)/2]$. (In fact it is a simple matter to show
that the latter is sharp for $m\ge 5$, cf.
Figs.\,\ref{Indefinite:fig} and \ref{Indefiniteodd:fig}.)
\end{lemma}

\begin{proof}
(1) The first estimate is trivial when $m$ odd, and when $m$ is
even it follows from the rigidity of the empty scheme, which is a
consequence of the fact that the invisible discriminant has real
codimension 2 (cf.
Lemma~\ref{invisible-discriminant-codim-2:lem}).

(2) The second estimate is a trivial consequence of the fact that
a rigid-isotopy induces an equivariant isotopy between the allied
symmetric Riemann surfaces. Formally the proof may require the
Ehresmann-Feldbau-Pontrjagin, etc. trivialization of a fiber
bundle over a contractible base (which is paracompact, else false
tangent bundle to the (simply-connected) Pr\"ufer surface).
Actually it requires an equivariant version thereof, but by
passing to the quotient and reconstructing the symmetric surface
as the double orientation cover we can reduce to the classical
setting.

(3) The third estimate follows from Rohlin's inequality,
conjointly with Klein's congruence, interpreted as obstructing
type~I right above the height of the  deep nest (of type~I being
totally real under a pencil of lines).
\end{proof}

The second estimate could be sharp. This dream is still much out
of reach, and so is the nightmare of refuting this via the
Fiedler-Marin method. At first the problem may look tractable, yet
it certainly does not follow from Gabard 2000 \cite{Gabard_2000}
where a sole classification of  symmetric surfaces realizable as
plane curves is given.

So below the height $[(m+3)/2]$ things are nearly pure and frozen
(no indefinite types) and naively we could expect that all schemes
are rigid below this altitude, i.e. when $r\le [(m+3)/2]$. Hence:

\begin{conj}
All schemes of degree $m$ with $r\le [(m+3)/2]$ are rigid.
\end{conj}

(Check if this  was not  disproved by Fiedler, but we do not think
so.)

We know (e.g. by the conceptual argument of Morse surgeries as
exposed in Viro 1989 \cite{Viro_1989/90-Construction}) that all
values of $r$ (number of real circuits) below Harnack's bound are
realized, by taking a generic pencil between an $M$-curve (e.g.
Harnack's) and an empty or Fermat curve with $r=1$. Likewise
either by the pedestrian argument in Gabard 2000
\cite{Gabard_2000} or perhaps a variant of Viro's conceptual
argument we know that for all intermediate values there is a
representative of type~II (safe for $M$-curves). Of course the
conceptual argument involves the Klein-Marin theorem
(\ref{Klein-Marin:lem}), since when lowering its number of
component through a Morse surgery a curve of type~II cannot become
of type~I. (Note yet that this is not enough to reprove the little
theorem of Kharlamov-Viro-Gabard exposed in Gabard 2000
\cite{Gabard_2000} in a conceptual Morse-theoretic fashion.) At
any rate it shows that the type~II is ubiquitous at all levels of
the ``pyramid'' as measured by the basic invariant $r$ (number of
real circuit), safe at the maximal $M$-level ($M=g+1$).

The behavior of the function $\rig(m)$ is highly mysterious. It
gives a measure of the flexibility of the video game allied to
Hilbert's sixteenth problem. You see on the screen $\RR P^2$
(essentially our retina) two curves of some fixed degree $m$
presenting the same isotopic topology (i.e. distribution of
ovals), can you pass continuously from one to the other by moving
the joystick in the hyperspace of all curves while avoiding the
discriminant $\disc$? If you can {\it always\/} achieve this goal
$\rig(m)=+\infty$ (you win always the game). If not $\rig(m)$
measure the smallest number of ``ovals'' (better real branches)
where you can loose the game. A priori $\rig(m)$ could be as low
as $\rig(m)=1$ when $m$ is large say $m>10^3=1000$, but some
topologist expect the Hilbert-video-game to be a more flexible
one, e.g. Rohlin-Viro-Itenberg positing rather that $\rig(m)\le 2$
for all $m$.

What is (inside each chamber of the discriminant) the curve with
largest systolic ratio, i.e. the most healthy against infarctus
(=hearth attack). Of course all this would be computed w.r.t. the
Fubini-Study metric, or maybe the uniformizing metric. It seems
likely that flows allied to such functionals ought to give some
information on the above problem, essentially because when the
systole shortens we approach the discriminant.

We learned the following from an e-mail of O.Viro (dated
[26.01.13] in Sec.\,\ref{e-mail-Viro:sec}):

\begin{lemma} If the one-oval scheme (unifolium) is rigid then the contraction
conjecture holds true for curves with one oval.
\end{lemma}

\begin{proof}
Let $C_{2k}$ be an even order smooth curve with one oval of degree
$2k$. It is enough to construct a contraction for a specific curve
$F_{2k}$ of degree $2k$. One can consider the Fermat contraction
(in affine equation) $x^{2k}+y^{2k}=\rho^{2k}$ with $\rho\to 0$.
The latter shrink to a point but alas the tangent cone is not
really that of an ordinary solitary node but rather possess $2k$
branches which form $k$ conjugate pairs corresponding to the
$2k$-th roots of $-1$. So we need a slightly different contraction
toward a solitary node.

This can probably be done explicitly or more loosely by taking a
curve of the form a circle $E_2$ union a $D_{2k-2}$ which is
invisible, while shrinking the radius of $E_2^{\rho}\colon
x^2+y^2=\rho^2$ and smoothing by Brusotti the union
$E_2^{\rho}\cup D_{2k-2}$.
\end{proof}

Now let us estimate the indefiniteness $\indef(m)$, i.e. the
lowest altitude at which there is a scheme of indefinite type.
This amounts to construct an indefinite scheme of minimum height.
This is fairly easy as shown by the following series of type~I
curves of even degree (upper row of Fig.\,\ref{Indefinite:fig}).
The bottom row of this figure exhibits type~II curves with the
same schemes. This is obtained by starting with a type~II sextic
while adding conics. Since type~II is  a genetically dominant
character, all successors are also of type~II. Hence all those
schemes are indefinite and  have the  minimal height permissible
by Rohlin's inequality (\ref{Rohlin's-inequality:cor}), namely two
units above the corresponding deep nest.

\begin{figure}[h]
\centering
\epsfig{figure=Indefinite.eps,width=122mm} \vskip-5pt\penalty0
  \caption{\label{Indefinite:fig}%
  Exhibiting indefinite schemes at minimum height} \vskip-5pt\penalty0
\end{figure}

A similar series is easy to find in odd degrees
(Fig.\,\ref{Indefiniteodd:fig}) and hardly requires any further
comments. This proves the:

\begin{figure}[h]
\centering
\epsfig{figure=Indefiniteodd.eps,width=122mm} \vskip-5pt\penalty0
  \caption{\label{Indefiniteodd:fig}%
  Exhibiting indefinite schemes at minimum height (odd degree series)} \vskip-5pt\penalty0
\end{figure}

\begin{lemma}
For $m\ge 5$, $\indef(m)=[(m+1)/2]+2=[(m+5)/2]+2$.
\end{lemma}


\begin{cor}
For $m\ge 5$, we have the estimate $\rig(m)\le [(m+5)/2]$.
\end{cor}

This implies a certain rigidity in the video game, and is merely a
consequence of Rohlin's work. Of course it would be miraculous if
this estimate is sharp, prompting a maximal flexibility of the
video game. For $m=5,6$ it is certainly sharp by the work of
Kharlamov 1981 \cite{Kharlamov_1981/81} and Nikulin 1979
\cite{Nikulin_1979/80}, respectively.
So we must concentrate on degree 7, or 8, and by Rohlin's
inequality the required obstruction must necessarily be of a
somewhat deeper nature than via Klein's types.

\subsection{Searching obstructions to rigid-isotopy below height
$DEEP+2$}

[31.01.13] Here the technology is due to Marin-Fiedler and
involves the lock allied to the subscheme $S$ of degree 7 of
symbol $\frac{3}{1}\sqcup J$ (3 ovals enveloped in one oval and a
pseudoline $J$ outside). If we trace the triangle of 3 lines
through the 3 empty ovals (Fig.\,\ref{Locks:fig}a), each line has
 7 real intersections (saturating B\'ezout).
It follows that 2 schemes of degree 7 enlarging $S$ cannot be
rigid-isotopic as soon as the distribution of the remaining ovals
past the 3 lines is different. Alas $\indef(7)=4+2=6$, hence to
beat this we must find a pair of isotopic but non rigid-isotopic
curves with $r\le 5$ circuits, which is already the height of the
Marin-Fiedler lock $S$. Hence this method seems not suited to our
goal.

\begin{figure}[h]
\centering
\epsfig{figure=Locks.eps,width=122mm} \vskip-5pt\penalty0
  \caption{\label{Locks:fig}%
  Some locks in degree 7} \vskip-5pt\penalty0
\end{figure}

We could change the lock into $\frac{2}{1} 1\sqcup J$
(Fig.\,\ref{Locks:fig}b), but this has still height $5$. Another
choice is Fig.\,\ref{Locks:fig}d but this has also height 5.
Another lock could involve a conic through 5 ovals
(Fig.\,\ref{Locks:fig}e) but this is not locked as $10< 2\cdot
7=14$.

Another idea is to use {\it pseudo-locks} like on the second row.
Alas a line is not dividing $\RR P^2$. Perhaps one can construct a
lock by aggregating the pseudoline $J$ to the lock (cf.
Fig.\,\ref{Locks:fig}f), whence our name {\it pseudo-lock}. Then
the red line union the pseudoline $J$ divides $\RR P^2$, and so
the location of a fifth oval could be an obstruction to
rigid-isotopy. Of course Fig.\,\ref{Locks:fig}g is not interesting
being saturated (maximal scheme). Fig.\,\ref{Locks:fig}h could be
employed as the former Fig.\,f. For this to work one should have
an isotopic-invariant way to distinguish both residues to the
{\it augmented-lock\/} consisting of the red line plus the
pseudoline. Alas in view of the symmetry of the lock it seems that
there is little chance to distinguish invariantly both halves (of
the augmented lock). One could imagine to move from the empty oval
to the deep oval (on Fig.\,\ref{Locks:fig}f) along the line while
choosing the route not intersecting the pseudoline $J$. W.r.t.
this oriented segment there would be a left and right hand side
residual to the lock. This concept is perhaps invariant under
isotopy, and there is some little chance to detect 2 septics with
$r=5$ which are isotopic but not rigid-isotopic.

Such a pair of septics is constructed on
Figs.\,\ref{Locks:fig}i,j, where the remaining oval lies either of
the left (Fig.\,i) or on  the right (Fig.\,j) of the oriented red
segment from {\it the} empty unnested oval  to the empty nested
oval. Does this prove  both curves being not rigid-isotopic? Maybe
not since both  are mirror images under a symmetry in
$G=PGL(3,\RR)$, which is a connected group ($\RR P^2$ being
non-orientable there is no way to reverse orientation) and so
there is a path in $G$ from the identity to the mirror
transformation. Applying this path to the first curve yields a
rigid-isotopy to the second curve. So where is our former argument
faulty? Culpability seems to be the italicized ``the'' some few
line above. Indeed there is on Fig.\,i no canonical choice for the
origin of the arrow, and if instead we had chosen it in the other
(outer) oval then the free (unlocked) oval would of course sit on
the right (instead of left) of the red arrow.

We can  try to remedy this defect by allowing only one outer oval,
but then there is another inner oval and there is no canonical way
to choose it (cf. Fig.\,\ref{Locks:fig}k). One could hope that one
of both inner ovals is distinguished, say by complex orientations
(but no chance as we are in the ``post deep-nest case'' $r=5=4+1$
hence nondividing).

[01.02.13] The situation becomes more favorable if we look at
locks in degree 9, especially the one depicted on
Fig.\,\ref{Locksdeg9:fig}l. Then there is a canonical way to trace
an arrow between the deep ovals (say from the less profound to the
more profound one as on Fig.\,\ref{Locksdeg9:fig}l). This is
invariantly defined in case the remaining oval (dashed) lies
outside the largest nonempty oval (as on
Fig.\,\ref{Locksdeg9:fig}l). Then the choice of this arrow is
canonical and it is hoped that the position of the dashed oval on
the left versus right of the arrow (augmented by the pseudoline)
affords an obstruction to rigid-isotopy. It is easy to manufacture
an algebraic curve realizing this schematic lock, cf.
Fig.\,\ref{Locksdeg9:fig}m where the free oval is righthanded. It
causes no trouble to find a similar picture with the free oval
lefthanded. This would give a nontrivial obstruction to
rigid-isotopy below Rohlin's temperature $\indef(m=9)=5+2=7$,
namely at $r=6$. In particular Fig.\,m would not be rigid-isotopic
to its mirror image. This violates however the above argument
using  connectedness of the group $PGL(3,\RR)$. Of course our
mistake is that in the nonorientable $\RR P^2$ there is no
consistent way to distinguish the left from the right. More
precisely while it is possible to orient the red line from the
less massive to the deepest oval, when the latter intercept the
pseudoline there is no way to choose a left or right sense to
bifurcate as the pseudoline itself lacks a preferred orientation.

\begin{figure}[h]
\centering
\epsfig{figure=Locksdeg9.eps,width=122mm} \vskip-5pt\penalty0
  \caption{\label{Locksdeg9:fig}%
  Some locks in degree 9} \vskip-5pt\penalty0
\end{figure}

The method becomes more effective if we permit one more ovals.
Then the two ``free'' (dashed) ovals can either be separated by
the augmented lock $L\cup J$ or not (Fig.\,\ref{Locksdeg9:fig}n).
Both cases do occur as shown by Figs.\ref{Locksdeg9:fig}o,p. Both
depicted curves are of type~II (inspect the little 3 arrows and
the negative smoothing right above it on
Fig.\,\ref{Locksdeg9:fig}o). The same local pattern appears on
Fig.\,\ref{Locksdeg9:fig}p, which is thus also of type~II. However
both curves are not rigid-isotopic, because during the
rigid-isotopy the two free ``dashed'' ovals of
Fig.\,\ref{Locksdeg9:fig}n cannot traverse the red line which is
B\'ezout-saturated nor can they traverse the pseudoline. This is a
little success of the Fiedler-Marin method, alas occurring at the
same height as the indefiniteness $\indef(m=9)=7$.

A similar example can be found already in degree 7, since we do
not actually require to orient the line, compare
Fig.\,\ref{Locks2:fig}a,b which should be self-explanatory. Note
again that both septics on Fig.\,\ref{Locks2:fig}c,d are of
type~II, yet not rigid isotopic. This would be worth stating as a
lemma since it is a little variant of the Fiedler-Marin method
(with now separation caused by the added pseudoline). However this
does not answer our puzzle of detecting obstruction to rigid
isotopy below the critical temperature $DEEP+2$.

It is then tempting to lower to degree 5, while considering the
lock Fig.\,\ref{Locks2:fig}e,f, but then alas we lack a canonical
choice for the red line. One can try other locks in degree 7, like
Fig.\,\ref{Locks2:fig}g, but then we lack again canonicalness.
Still one could  make some choice and propagate it consistently
during the isotopy. So we get Figs.\,\ref{Locks2:fig}h,i and
arrive at the fallacious conclusion that the curve is not
rigid-isotopic to itself. This nonsense helps  emphasizing the
importance of the lock being somehow God-given by the curve, and
we (human beings) making minimalist intervention upon the
creation.

\begin{figure}[h]
\centering
\epsfig{figure=Locks2.eps,width=122mm} \vskip-5pt\penalty0
  \caption{\label{Locks2:fig}%
An obstruction (Figs.\,c,d) to rigid-isotopy  in degree 7 via the
Fiedler-Marin method between 2 curves at the indefiniteness height
(the real challenge is to find such an obstruction between curves
with one less oval)} \vskip-5pt\penalty0
\end{figure}

It seems that detecting obstructions to rigid-isotopy
beyond Klein's type and below the critical temperature
(=indefiniteness $\indef(m)$) is a hard business requiring
completely new ideas, or at least some better acquaintance with
the Marin-Fiedler obstruction.

[01.02.13] Paraphrasing, the method of the lock does not seem to
obstruct rigid-isotopies below the indefiniteness $\indef(m)$,
i.e. the lowest height of an indefinite scheme. So perhaps the
first obstruction to rigid-isotopy is given by Klein's type and
occurs at height $\indef(m)=[(m+1)/2]+2=[(m+5)/2].$ In that case
the rigidity index $\rig(m)$ would be highest possible equal to
the indefiniteness $\indef(m)$.

[02.02.13] Let us summarize the discussion. For any degree $m$,
there is a deep nest with $r=[(m+1)/2]=:DEEP$ real branches. Two
units above the latter's height it is easy to construct curves
having the same real scheme yet different types (I vs. II) hence
not rigid-isotopic. Using the method of the  lock it is even
possible to exhibit at this height  curves of degree 7 or 9 having
the same real scheme and the same type~II, yet not rigid-isotopic.
Probably the method described extend to all other odd degrees.
However, it seems much more tricky and actually the locking method
seems incapable detecting obstruction below this height, starting
thus at height one unit above the height of the deep nest. Could
it be that all schemes at or below this height are rigid, i.e. any
two curves representing it are rigid-isotopic.

{\it Minor question (skip)}.---As a minor problem we suspect that
for all odd degrees $m\ge 7$ there is a non-rigid scheme at height
$DEEP+2$ containing a pair of type~II curves which are not
rigid-isotopic. This is probably easy and merely involves
extending into series the examples of Figs.\,\ref{Locks2:fig}c,d
and \ref{Locksdeg9:fig}o,p.

[03.02.13] {\it Main problem}.---So we first focus on the case
$m=7$. Let us denote by $\Delta=DEEP=[(m+1)/2]$ the height of the
deep nest. Our goal is to find obstruction to rigid-isotopy
(strictly) below height $\Delta(m)+2$. For $m=7$, we have
$\Delta=4$, and so we look at schemes with height $r=5$. Several
cases occur and are primarily the schemes
$$
\frac{3}{1}J,\quad \frac{2}{1}1J,\quad \frac{1}{1}2J,\quad 4J,
$$
where we use Gudkov's notation and $J$ denotes the pseudoline
(unique up to isotopy). The corresponding schemes are depicted on
Fig.\,\ref{Locksdeg7:fig}, where the locks are depicted as red
thick-lines which are B\'ezout saturated, while dashed-lines are
not. The philosophy of the locking method is that a free oval
cannot traverse during a rigid-isotopy the lock (without violating
B\'ezout) and therefore the distribution of additional ovals {\it
among the residual components of the lock\/} ({\it past the
lock\/} for short) has to be respected. If is {\it not}, then we
have an obstruction to rigid-isotopy. The dramaturgy in our case,
where the height is as low as $r=\Delta+1$, is that we do not have
any such additional ovals available (all having been consumed by
the lock so-to-speak). On Fig.\,a we could kill the nonempty oval,
but then we loose B\'ezout-saturation of the red-lines.

\begin{figure}[h]
\centering
\epsfig{figure=Locksdeg7.eps,width=122mm} \vskip-5pt\penalty0
  \caption{\label{Locksdeg7:fig}%
Schemes of degree 7 at height $\Delta+1= 4+1=5$}
\vskip-5pt\penalty0
\end{figure}

So we need some much deeper idea. One idea is to look how the
locked ovals themselves are separated by the lock. This seems
however to lead nowhere. Indeed examine the case of the scheme
$\frac{1}{1}2J$ (i.e. Fig.\,\ref{Locksdeg7:fig}c), where there is
a menagerie of possible disposition of the pseudoline $J$
(Fig.\,\ref{Locksdeg7:fig}e). To effect a nice separation we
include the pseudoline into the lock. The pseudoline plus the 2
thick B\'ezout-saturated lines effects a separation in 4 zones,
yet whatever the situation of $J$ the disposition of  ovals in
those zones is still the same. At least so are the number of
residual components in each of these zones,  weighted on
Fig.\,\ref{Locksdeg7:fig}e by the corresponding number of
components $2,3,3b,4$. One can even play more sophisticated games
by choosing one of this zone in some invariant manner. For
instance given the 2 points of intersections of the thick lines
with $J$, we may link them to the deep nest along the thick lines
while choosing the way avoiding the dashed line, and close this by
the piece of $J$ cutting the dashed line an even number of times
(counted by multiplicity). Since this canonical curve $J_0$ cuts
the dashed lines an even number of times, it is null-homotopic and
bounds a unique disc, which is our canonical region. Alas one
checks (experimentally) on Fig.\,\ref{Locksdeg7:fig}e that it
always contains 2 components of the scheme. There is a dual curve
constructed by taking the segments linking the points of $J \cap
L_i$ to the deep nest via the path cutting once the dashed line,
and aggregating the same portion of $J$ as above. This Jordan
curve still cuts $L_3$ (dashed line) an even number of times, and
so bounds a unique disc. The latter (alas) always contains 4
components of the curve.

Another little idea we had, is to mark for each point in $L_i\cap
J$ the vertices of the locking triangle which looks closest to the
intersection point while travelling on the given $L_i$. However as
shown by Figs.\,\ref{Locksdeg7:fig}f,g this is completely
insensitive to a variation of the position of the pseudoline.

Repeating ourselves, it seems that the {\it method of the lock\/}
fails to detect any obstruction to rigid-isotopy at height
$\Delta+1$ (or below). Accordingly one may suspect that there is
no such obstruction.

Here is an idea. Given a smooth $C_7$, there is a unique
pseudoline $J$. Let us speculate about
a large
deformation $C_7 \to C_6 \cup L_1$ toward a sextic plus a line.
This is supposed to be a path in the space of curves avoiding the
discriminant sole for its extremity. In particular the split curve
is isotopic to the original $C_7$. We call this the rectification
conjecture:

\begin{conj} {\rm (Rectification conjecture=RC)}
\label{rectif-conj:conj}
Given any (smooth, real, plane) curve of odd order $C_{2k+1}$
there is a deformation in the large toward a curve $C_{2k}\cup
L_1$ where $L_1$ is a line.
\end{conj}

{\it Objection} [07.04.13] Already for quintics, this formulation
is sloppy: take an $M$-quintic (hence with symbol $6J$), while
quartics can have at most 4 ovals.

If this large deformation is implementable (more about this soon),
then we deduce that $r(C_7)=r(C_6)+1$. If the given septic scheme
has height $r\le \Delta(7)+1=4+1=5$, then the sextic has $r\le
4=\Delta(6)+1$. But in this low range sextic schemes are rigid by
Nikulin's rigid classification enhancing the Gudkov-Rohlin table
(cf. Fig.\,\ref{Gudkov-Table3:fig}). Hence we are inclined to
think that septic schemes are rigid below  height $\Delta+1$.

Indeed given 2 septics which are (soft) isotopic, i.e. belong to
the same real scheme, we apply the rectification conjecture
(\ref{rectif-conj:conj}) twice to deduce sextics with the same
real scheme and of low height  $r\le 4$, hence rigid-isotopic. Now
using a path between the split curves of degree $6+1$ and using a
version of Brusotti's theorem with parameters (yet to be
formulated) one could argue that the 2 given septics are
rigid-isotopic. The proof would be completed.

A brief word in favor of the conjecture (\ref{rectif-conj:conj}).
Given an odd order curve there is a unique pseudoline, and one may
measure its length (w.r.t. the round elliptical geometry on the
real projective plane $\RR P^2$). Obviously the (genuine) line is
the pseudoline of minimum length, namely $\pi=3.14\dots$ if we
work on the unit sphere as preferred double cover of $\RR P^2$.
Hence for this functional (length of the pseudoline) the gradient
lines ought to converge toward curves splitting off a line. (Maybe
one can also look at the total geodesic curvature of the
pseudoline as another competing functional doing the same job.)

Having this we may dream of a grand inductive process reducing the
whole problem of rigid-isotopy (at least below the range
$\Delta+1$) to Nikulin's seminal theorem on sextics (itself
relying on deformation theory of K3 surfaces). This would lead to
a sharp estimation of the rigidity index $\rig(m)$ of the previous
section as being equal to $\Delta(m)+2$.

However even for degree 8, this looks hazardous. One could imagine
two modes of deformation of a $C_8$ to either a septic plus a line
$C_7\cup L_1$ or a $C_6\cup E_2$. The latter looks dubious for the
(8)-scheme consisting of 3 nests of depth 2 (of height
$r=6=\Delta+2$), since removing one oval one has still the line
through the two remaining nests
creating 8 intersections (too much for a $C_6$). Yet the latter is
precluded as we restrict to schemes of height $\le \Delta+1$.
Listing all of them we find in Gudkov's notation the following
list of schemes (cf. Fig.\,\ref{Locksdeg8:fig}):
\def\rau{\hskip2pt}
$$
\frac{4}{1},\rau\frac{3}{1}1,\rau\frac{2}{1}2,
\rau\frac{1}{1}3,\rau
5,\rau
(1,\frac{1}{1}2), \rau
%
(1,\frac{3}{1}), \rau (1,\frac{2}{1}1), \rau (1,\frac{2}{1})1,
\rau (1,\frac{1}{1}1)1, \rau (1,\frac{1}{1})2,\rau
\frac{1}{1}\frac{1}{1}1, \rau \frac{2}{1}\frac{1}{1}.
$$

\begin{figure}[h]
\centering
\epsfig{figure=Locksdeg8.eps,width=112mm} \vskip-5pt\penalty0
  \caption{\label{Locksdeg8:fig}%
Schemes of degree 8 at height $\Delta+1= 4+1=5$ and how they may
degenerate under conification to schemes of degree 6. Black arrows
denote the extinction of an empty oval, while red arrows denote
the liberation of an empty oval from the oval immediately
enveloping it.} \vskip-5pt\penalty0
\end{figure}

Albeit messy, our picture (Fig.\,\ref{Locksdeg8:fig}) is supposed
to
take the census of
 all possible degenerations $C_8\to C_6\cup E_2$
which are B\'ezout permissible. Of course we do not claim that all
these moves are algebraically realized, but at least B\'ezout
gives no obstruction. Alas it is far from obvious (unlike in the
odd degree case) which functional is capable effecting the large
structural deformation (LSD) of ``conification'' $C_{2k}\to
C_{2k-2}\cup E_2$ splitting off a conic $E_2$. Naively we may
expect that it is always some of the empty oval which shrinks to a
solitary node, but  that soon before getting extinct  he splits
off an infinitesimal circle (or ellipse). This could involve a
sort of isoperimetric functional measuring rotundity of ovals, and
the allied lines of steepest descent (or ascension). Note however
that for the $8$-scheme $(1,\frac{3}{1}),
(1,\frac{2}{1}1),(1,\frac{2}{1})1$ (the 3 firsts of the third row
on Fig.\,\ref{Locksdeg8:fig}) we cannot ``conify'' the empty ovals
without violating B\'ezout. Indeed removing one of the 3 possible
empty ovals leads to scheme containing the deep nest of depth 3 as
a (strict) subscheme. Of course those (8)-schemes really exist, as
depicted on Fig.\,\ref{Locksdeg8:fig}b. A priori nothing precludes
a degeneration like  Fig.\,\ref{Locksdeg8:fig}c, where a nonempty
oval would be ``conified''.

 \begin{conj} {\rm (Conical/ellipticity
conjecture=EC) [inserted 05.02.13]} Given a (smooth, real,
algebraic, plane) curve $C_{2k}$ of even degree $m=2k$  with few
ovals (i.e. $r\le \Delta(m)+1$ where $\Delta(r)=k$ is the number
of ovals of the deep nest of degree $m$) there is a deformation
(=rigid-isotopy safe its extremity) toward a curve $C_{2k-2}\cup
E_2$ where $E_2$ is an ellipse, or equivalently a circle up to
projectivity. Alas this cannot always occur by extinction of an
empty oval, but sometimes by inflation of a large oval (perhaps
via an isoperimetric gradient-flow).
\end{conj}

\def\isop{{\rm isop}}

The difficulty with this conjecture is that unlike for its odd
degree avatar (\ref{rectif-conj:conj}) we lack a canonical
functional to be minimized like the length of the pseudoline. (The
line is the shortest pseudoline, and being non-null-homotopic it
is like a systole.) In the even degree case all ovals are
null-homotopic and there is no systole in $\RR P^2$. Of course
there could be a systole on the Riemann surface of the
complexification. Alternatively one may replace the systolic
problem by an isoperimetric one taking also area into account. Let
us introduce the isoperimetric ratio ($\isop$) of an oval as its
length squared divided by the area of its bounding disc, all  in
reference to the round elliptical geometry on $\RR P^2$. In
Euclidean geometry this is minimum for a circle $(2\pi
\rho)^2/(\pi \rho^2)=4 \pi=12.566\dots$. For a large circle near
the equator this can be as close as we please to  $(2\pi)^2/ (2
\pi)=2\pi=6.28\dots$, which is smaller. This is probably the
absolute infimum if we demand the oval on the unit sphere to be
disjoint from its antipode. Now we could hope that  the minimum
isoperimetric ratio
of all ovals leads to a functional whose gradient lines tend to
inflate the most rotund oval toward an ellipse (rotundity being
measured by
the isoperimetric ratio). This could give the required
degeneration. Perhaps in the limit the most rotund oval degenerate
to a pair of lines (double line) and suppressing one of those
leads to a odd degree curve of degree one less. This would give
the other mode of degeneration:

\begin{conj}
Given any $C_{2k}$ of even degree of height $\le \Delta+1=k+1$,
there is a rigid-isotopy safe extremity toward a curve
$C_{2k-2}\cup (L_1)^2$ splitting off a double line $L_1$. More
precisely the orthogonal trajectories of the  rotundity functional
(measured by the minimum isoperimetric ratio) drives any such
curve toward such a curve in a canonical fashion.
\end{conj}
If this conjecture holds true then we would have a sharp estimate
of the rigidity index $\rig(m)$ for all degrees [end insert
05.02.13].

All what we are saying sounds very optimistical, and we are still
very far from having a decent understanding of this problem of
rigid-isotopy (strictly) below the height $\Delta+2$.

We can hope that the method of the lock is more efficient in
degree 8 than it was in degree 7 (still confining our attention to
heights $\le \Delta+1$). Fig.\,\ref{Locksdegree8:fig} depicts some
of them. The method of the lock is a jewel discovered in the late
1970's by Marin and Fiedler independently (all being inspired by
V.\,A. Rohlin's work). It involves basically the idea of attaching
in the most canonical way to a given curve a certain red
configuration acting as a separator. More precisely special
attention is paid to red thick lines which are B\'ezout-saturated,
so that the remaining ovals of the curve cannot traverse this line
during the isotopy. So basically we choose a triad of points
inside some ``deep'' ovals and link them by a triangle of lines.
Of course the choice of the points is not perfectly canonical, but
we choose them inside the disc bounding an empty oval. The
Marin-Fiedler trick is quite reminiscent of what Grothendieck
calls ``le principe des choix anodins'' (in Esquisse d'un
programme 1984
\cite{Grothendieck_1984/1997-esquisse-d'un-programme}) that
whenever we make some choices within a contractible space the
construction is nearly canonical, hence robust and fruitful. It is
also reminiscent of the moving-frame method of Darboux-Cartan
(rep\`ere mobile), since during the rigid isotopy will really move
the whole triangle.

\begin{figure}[h]
\centering
\epsfig{figure=Locksdegree8.eps,width=102mm} \vskip-5pt\penalty0
  \caption{\label{Locksdegree8:fig}%
Exploring locks for schemes of degree 8 at height $\Delta+1=
4+1=5$} \vskip-5pt\penalty0
\end{figure}

On the 3 first pictures of the 2nd row of
Fig.\,\ref{Locksdegree8:fig} we have a perfect lock by a triangle
consisting of 3 lines which are B\'ezout-saturated. Alas we have
no more ovals left to separate and the method looks inoperative.
So let us look at the next degree $9$, and list all the schemes at
height $\Delta+1=5+1=6$. It seems plain that this merely amounts
to add a pseudoline to the former configurations listed in degree
8 (cf. Fig.\,\ref{Locksdegree9:fig}).

\begin{figure}[h]
\centering
\epsfig{figure=Locksdegree9.eps,width=102mm} \vskip-5pt\penalty0
  \caption{\label{Locksdegree9:fig}%
Locks for schemes of degree 9 at height $\Delta+1= 5+1=6$}
\vskip-5pt\penalty0
\end{figure}

Again the big deception is that no elementary obstruction given by
locking appears in view. Of course one could interpret the figure
as a rectification (\ref{rectif-conj:conj}) toward octics which
are (hypothetically) rigid at height $\Delta+1$ (say via a
reduction to sextics), and so would be our curves of degree $9$.

The next real jump in complexity involves degree 10. Let us
tabulate all schemes at the critical height $\Delta+1=5+1=6$ while
avoiding any Gudkov symbolism (cf. Fig.\,\ref{Locksdegree10:fig}).
This is elaborated as follows. Start from any configuration,
especially the maximum elements sembling highly concentric and
protected medieval settlements like $(1,1,1,1,1,1)$, or
$(1,1,1,1)(1,1)$, $(1,1,1)(1,1,1)$, $(1,1)(1,1)(1,1)$, and then
apply basically two moves freeing an oval. Vertical moves
correspond to liberating a deep oval, while horizontal moves freed
a superficial oval (of small depth). Sometimes there are ovals at
3 different depths so that we have also a 3rd oblique move. The
red framed schemes are prohibited by B\'ezout, yet are useful as
generator (under the described moves) of other schemes that
otherwise are easily overlooked.

\begin{figure}[h]
\centering
\epsfig{figure=Locksdegree10.eps,width=122mm} \vskip-5pt\penalty0
  \caption{\label{Locksdegree10:fig}%
List of (all?) schemes of degree 10 at height $\Delta+1= 5+1=6$}
\vskip-5pt\penalty0
\end{figure}

[04.02.13] It seems evident at this stage that there is a
combinatorial law (which overwhelms my intelligence) impeding the
the locking method to act as an obstruction to rigid-isotopy.
Looking at all possible locks on Fig.\,\ref{Locksdegree10:fig},
 no obvious obstruction to rigid-isotopy
strikes the vision. In contrast on the basis of the same picture,
one may argue that erasing a suitable oval all our (10)-schemes
reduce to one of degree 8, and if this cancellation is geometrized
via the conification (elliptization) conjecture we could deduce
 rigidity of all the (10)-schemes at height $\Delta(10)+1=6$,
from that of the corresponding (8)-schemes, which in turn was
reduced to (6)-schemes where low-height rigidity holds true by
virtue of Nikulin's theorem (1979 \cite{Nikulin_1979/80}).

As a little experiment imagine the curve of degree $2k$ to have 3
nests of depth $d_1,d_2,d_3$. Since we are at height
$\Delta(2k)+1=k+1$, we have $k+1=r=d_1+d_2+d_3$. Let $L_1,L_2,L_3$
be 3 lines passing through the deep nests and suppose them
B\'ezout-saturated, then $d_1+d_2,d_2+d_3,d_1+d_3$ are all equal
to $2k$, and thus summing and dividing by two we infer that
$d_1+d_2+d_3=3k$, which is much greater than $r=k+1$.

If instead of 3 deep nests we have one deep nest containing 3
little ovals, then the 3 lines through them supposed
B\'ezout-saturated cut the curves in  $4+2d=2k$ real points where
$d$ is the depth of the nest. Hence $r=d+3=(d+2)+1=k+1$, so that
all ovals are exhausted by the lock (and nothing remains left to
be separated).

Such arguments seem to extend to all other schemes of
Fig.\,\ref{Locksdegree10:fig}. There sometimes we lock with only
two totally real lines like e.g. on the scheme $(1,1,1,1)2$ lying
near the center of Fig.\,\ref{Locksdegree10:fig} (right above the
2 anti-B\'ezoutian schemes). Two lines suffice to separate the
plane $\RR P^2$, but here again the construction of the lock
consumes all the ovals at disposal.
In summary it seems hopeless to find an obstruction to
rigid-isotopy at or below the height $\Delta(m)+1$ (at least via
the lock-method of Marin-Fiedler).

\begin{figure}[h]
\centering
\epsfig{figure=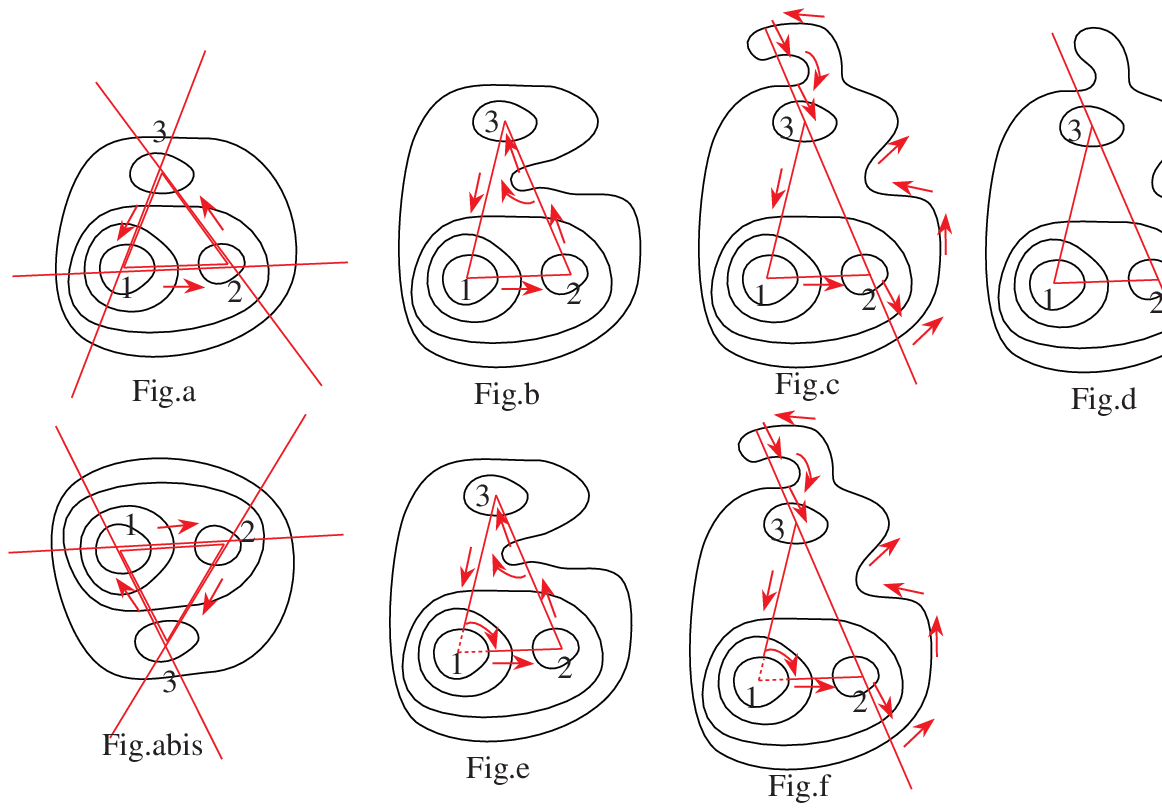,width=102mm} \vskip-5pt\penalty0
  \caption{\label{Lock10:fig}%
A distinguished scheme of degree 10 at height $\Delta+1= 5+1=6$}
\vskip-5pt\penalty0
\end{figure}

As a last chance,  consider the (10)-scheme of
Fig.\,\ref{Locksdegree10:fig} right before the ``mild'' arrow,
that is $(1,(1,\frac{1}{1}1)1)$. This is distinguished by having 3
empty ovals at different depths. So we can link them by a
$2$-simplex with boundary  oriented as
going from the deepest to the middle deep and then to the less
profound oval closed back to the deepest one. This would induce a
certain orientation on the inside of the largest oval (as usual
ovals being ordered by inclusion of their insides).
The problem however is that while the 2 lines through the deepest
oval are saturated (hence there is  preferred pathes joining them
in the inside of the maximal oval), the third
 is not and so there is no preferred way to join the middle
empty oval to the less deep one (compare Fig.\,\ref{Lock10:fig}).
However we could argue that whatsoever the way chosen we get the
same orientation (compare Figs.\,\ref{Lock10:fig}b and c). On the
latter figures we follow the line until reaching the maximal oval
$O_m$ and then follow the latter. The problem is which direction
to choose when we meet $O_m$. A priori there is no preferred sense
to bifurcate, but we may choose the path such that the circuit
$1\to 2\to 3 \to 1$ does {\it not} enclose the deep oval of depth
3 (i.e. the one containing the point $1$). This has no intrinsic
meaning unless we take the precaution of first rounding the corner
at the vertices 1 as shown on Figs.\,\ref{Lock10:fig}e,f. Note
that there is a unique way to put near $1$ an arrow circulating on
the deepest oval in such a way that we do not intercept the lines
$1,2$ and $1,3$ too frequently (i.e. only twice instead of 4
times). This as an intrinsic meaning since those lines are
saturated.

As a result any ten-ics $C_{10}$ belonging to the discussed scheme
would have a canonical orientation of the inside disc of its
maximum oval (hence of the latter as well). Of course this
(semi-)local orientation of the maximal oval propagates
continuously through a rigid-isotopy, but it seems that as $\RR
P^2$ is nonorientable no obstruction to rigid-isotopy can be
derived from this complicated trick. So even if two such curves
$C_{10}$ would have opposed canonical orientation over some region
of overlap of there maximal discs (bounding the maximal oval) this
would not impede them being rigid-isotopic.

Notwithstanding since the maximal disc (the inside of the maximal
oval) is oriented canonically, we may look at the deep line $1,2$
through the deepest empty ovals. This line does not separate $\RR
P^2$, but  certainly separates the maximal disc. Further the deep
line is oriented by going from $1$ to $2$ while staying inside the
maximum oval. Using the canonical orientation there is a left and
right hand side of this deep line inside the maximal disc. The
location of the superficial empty oval as being right- or
left-sided could give an obstruction (since the 3rd superficial
oval $O_3$ is not permitted to traverse the deep line during the
rigid-isotopy). So if like on Fig.\,\ref{Lock10:fig}a the
superficial (empty) oval $O_3$ is left-sided with respect to the
oriented line $1,2$ and the canonical orientation it will stay so
during for all curves explored by the isotopy, in particular for
the end curve. So naively it would suffices to apply a horizontal
axis $1,2$ symmetry to Fig.\,\ref{Lock10:fig}a  (and realize the
scheme geometrically which causes no difficulty via Brusotti) as
to find a curve with $O_3$ sitting on the other (right) side.
However we must really work with the canonical orientation of the
maximal disc, looking at Fig.\,\ref{Lock10:fig}abis shows that the
oval $O_3$ really sits on the left albeit sembling on the right
(where of course left has to be interpreted as the half pointed by
the canonical orientation). With all these confusing remarks, it
should be clear that there is no hope to detect an obstruction to
rigid-isotopy.

[05.02.13] We can also study the embryology of the scheme as shown
on Fig.\,\ref{Lock10bis:fig} depicting a nearly exhaustive list of
collision which an oval can acquire with the non-saturated line
$2,3$ through the 2 most superficial ovals. This represents the
possible cytoplasmic expansions of the ovals, but does not {\it
per se\/} afford obstructions to rigid-isotopy since all
configurations are linked to the initial one in some starlike
fashion.

\begin{figure}[h]
\centering
\epsfig{figure=Lock10bis.eps,width=122mm} \vskip-5pt\penalty0
  \caption{\label{Lock10bis:fig}%
Embryology of the distinguished (10)-scheme at height $\Delta+1=
5+1=6$ via cytoplasmic expansions colliding with the non-saturated
line $2,3$} \vskip-5pt\penalty0
\end{figure}

At the opposite extreme of such B\'ezout permissible moves, we
have the following 3 motions forbidden by B\'ezout where one of
the empty oval cannot traverse the saturated thick red line
(Fig.\,\ref{Lock10tris:fig}). So if we transgress the B\'ezout
obstruction by letting the oval traverse the dead-line then we get
the configuration of the second row of Fig.\,\ref{Lock10tris:fig}
which are priori could be non-rigid-isotopic to the initial one.
Alas there is still this argument of symmetry using the
connectedness of the group $PGL(3, \RR)$ which prevents one to
conclude that the configuration pre- and post-transgression are
not rigid-isotopic.

\begin{figure}[h]
\centering
\epsfig{figure=Lock10tris.eps,width=92mm} \vskip-5pt\penalty0
  \caption{\label{Lock10tris:fig}%
Cannot traverse (1st row), then transgressing that law (2nd row),
and transgressing the transgression (3rd row) via a simple
symmetry about the thick line, hence rigid-isotopic to the 1st
row} \vskip-5pt\penalty0
\end{figure}

\subsection{Trying in vain to corrupt Nikulin (via Marin-Fiedler)}
\label{Nikulin-corruption:sec}

[05.02.13] It is quite tempting (for dummies) to see if the method
of the lock (Marin-Fiedler 1979--1980) can parasite Nikulin's
rigid-isotopy classification of sextic (1979
\cite{Nikulin_1979/80}). Of course this is not to palish the glory
of Nikulin's theorem which is perhaps the deepest jewel ever
obtained along the lines of Hilbert's 16th problem, but rather an
 experimental game emphasizing the profundity of
Nikulin's result. Usually, the more a theorem looks
unbelievable, the deeper it stands.

For instance we may start with the basic scheme $\frac{3}{1}$ of
degree 6 (locked by the triad of lines through the 3 pairs of deep
ovals), and enhance it by adding 2 outer ovals to get the scheme
$\frac{3}{1}2$. We look at the distribution of outer ovals  past
the locking triangle, which a priori can be as on
Fig.\,\ref{Lock6:fig} either monopartite or bipartite. If one is
capable to exhibit two curves $C_6$ with distinct distributions
then both curves are not rigid-isotopic, for during a
rigid-isotopy the unlocked ovals cannot traverse the (moving)
triangle which is already B\'ezout-saturated. (Of course the
locking triangle works as well for $\frac{2}{1}1$, but then there
is nothing to separate, and if we add ovals then canonicalness of
the triangle is spoiled.)

On tracing explicit sextic curves $C_6$ via the small perturbation
method applied to configurations of 3 conics we always find the
same mono-partite arrangement where both ovals lies in the same
component residual to the triangle (Fig.\,\ref{Lock6:fig}).

\begin{figure}[h]
\centering
\epsfig{figure=Lock6.eps,width=122mm} \vskip-5pt\penalty0
  \caption{\label{Lock6:fig}%
Lock in degree 6 and trying to corrupt Nikulin}
\vskip-5pt\penalty0
\end{figure}

It seems impossible to corrupt Nikulin's result. As the scheme
$\frac{3}{1}2$ has height $r=6=\Delta+3$, three units above the
deep nest it is necessarily of type~II (by Klein's congruence) and
therefore Nikulin's theorem actually implies the:

\begin{lemma}
Any sextic $C_6$ belonging to the scheme $\frac{3}{1}2$ is such
that the triangle through the $3$ deep ovals does not separate the
outer ovals.
\end{lemma}

(We do not know whether this can be proved in an elementary
fashion without using the technological arsenal behind Nikulin's
theorem.) [07.04.13] Update: yes we can, cf. Le~Touz\'e in
Sec.\,\ref{LeTouze:sec}.

[06.02.13] Of course we may also add to $\frac{3}{1}$ more ovals
and examine the resulting distributions past the deep triangle.
That is we consider the schemes $\frac{3}{1}\ell$, where $2\le
\ell \le 5$ according to Gudkov's table
(Fig.\,\ref{Gudkov-Table3:fig}).

Consider first the scheme $\frac{3}{1}3$. Smoothing 3 ellipses we
can realize this scheme in two fashions either of type~I or II
(Fig.\,\ref{Lock6bis:fig}). However in both cases the locking
triangle through the deep (odd) ovals does not  separate the 3
outer ovals. This is quite surprising as both curves are not
rigid-isotopic, one could have expected that the lock-method to
detect the obstruction. We may also realize this scheme via a
variant of Hilbert's oscillation method, but again the
distribution of the 3 outer ovals is the same mono-partite one (at
least on the Walt-Disney depiction of Hilbert ``\`a la Gudkov'').

\begin{figure}[h]
\centering
\epsfig{figure=Lock6bis.eps,width=122mm} \vskip-5pt\penalty0
  \caption{\label{Lock6bis:fig}%
Locking the degree 6 scheme $\frac{3}{1}3$ and trying to corrupt
Nikulin} \vskip-5pt\penalty0
\end{figure}

Since both types I, II have representatives with the same
distribution, Nikulin's theorem implies the:

\begin{lemma}
Any sextic $C_6$ with scheme $\frac{3}{1}3$ (be it dividing or
not) is such that the triangle through the $3$ deep ovals does not
separate the $3$ outer ovals.
\end{lemma}

Next examine the scheme $\frac{3}{1}4$. Here we start with a
schematic picture \`a la Hilbert-Gudkov producing the curve
$\frac{4}{1}4$ (cf. Fig.\,\ref{Lock3-14:fig}a) which has too much
inner ovals (4 instead of the 3 desired). Such a Hilbert-vibration
is realized by Fig.\,b. A suitable smoothing gives Fig.\,c. The
latter has actually a companion generated by smoothing differently
the 3 inner nodes. In both cases however the deep triangle does
not separate the 4 outer ovals. Fig.\,d depicts  a Hilbert
vibration perturbing the union of both ellipses to the
Zeuthen-Klein G\"urtelkurve, but the quartic $C_4$ would then
intersect too frequently (at least 10 times) the conic. Such a
vibration is therefore precluded.  The dual vibration however
(Fig.\,e) is B\'ezout compatible (as the $C_4$ intersect $8$ times
the two conics). It is questionable if such a vibration exists as
the dual does not. Anyway let us (somewhat liberally) smooth
Fig.\,e to get Fig.\,f,  a somewhat funny curve belonging to the
scheme $\frac{3}{1}3$. Tracing the triangle through the deepest
ovals is somewhat challenging, but does not seem to effect a
division of the 3 outer ovals. A priori the depiction could be
like on the surrealist d\'etail (i.e., the median oval lying on
the ``left'' of the line through the other 2 inner ovals), but
this does not even seem to affect our issue about distribution of
outer ovals past the lock. Fig.\,g depicts another mode of
vibration which still overwhelms B\'ezout. Fig.\,h depicts yet
another mode of vibration essentially dual of Fig.\,b, but which
also overwhelms B\'ezout. It is a bit puzzling that not any
admissible vibration seems to admit a dual vibration.

\begin{figure}[h]
\centering
\epsfig{figure=Lock3-14.eps,width=122mm} \vskip-5pt\penalty0
  \caption{\label{Lock3-14:fig}%
Locking the degree 6 scheme $\frac{3}{1}4$ and trying to corrupt
Nikulin. For a slower depiction of the Harnack-method curve on the
right, cf. Fig.\,\ref{HarnaGudkov3-14:fig} below}
\vskip-5pt\penalty0
\end{figure}

At any rate if we believe in Nikulin's theorem (as we should since
it is Soviet mathematics of the best stock) our sole
Fig.\,\ref{Lock3-14:fig}c suffices to imply (since by Klein's
congruence our scheme $\frac{3}{1}4$ is of type~II) the following:

\begin{lemma}
Any sextic $C_6$ belonging to the scheme $\frac{3}{1}4$  is such
that the triangle through the $3$ deep ovals does not separate the
$4$ outer ovals.
\end{lemma}

Next we consider the scheme $\frac{3}{1}5$. For this we can either
look at Gudkov's construction (Fig.\,\ref{GudkovCampo-5-15:fig})
or at the easier construction via a variant of Harnack's method.
In Gudkov's setting, we must presumably consider the pull-back of
the triangle under the Cremona transformation and this a bit
tricky to depict. This should be manageable if one is in good form
but perhaps there is a more elementary direct construction via the
Harnack method.

\begin{figure}[h]
\centering
\epsfig{figure=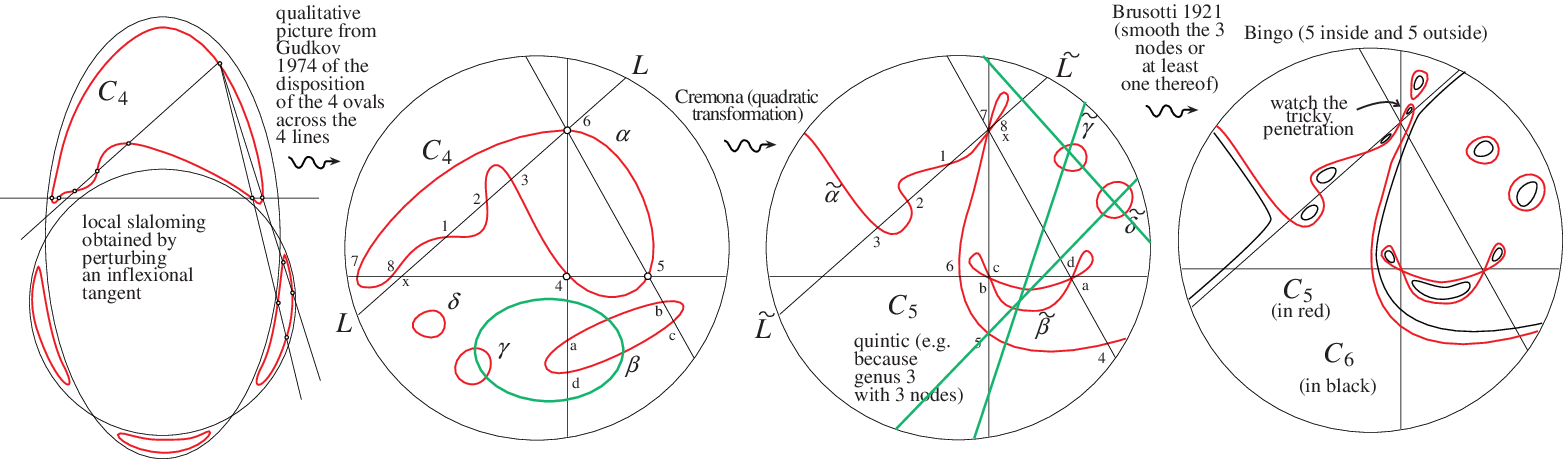,width=122mm}
\vskip-5pt\penalty0
  \caption{\label{GudkovCampoLock:fig}%
Locking the scheme $\frac{3}{1}5$ via Gudkov (but aborted)}
\vskip-5pt\penalty0
\end{figure}

So it seems fundamental to construct the over-scheme
$\frac{4}{1}5$ via the variant of Harnack's method mentioned  in
Gudkov 1974 \cite[p.\,42]{Gudkov_1974/74}, where alas no details
are supplied. The smaller scheme $\frac{3}{1}5$ we are interested
in should  then easily be deduced by taming the smoothing. Since
this has some independent interest we devote the next section to
the topic.

\subsection{Gudkov's variant
of Harnack: construction of the $(M-1)$-scheme $\frac{4}{1}5$}

[06.02.13] We now try to fix Gudkov's claim (in 1974
\cite[p.\,42]{Gudkov_1974/74}) that a suitable variant of
Harnack's method produces the $(M-1)$-scheme $\frac{4}{1}5$. {\it
Per se} this is not extremely original for we already managed (on
the shoulder of Gudkov's ``original'' construction of
$\frac{5}{1}5$, cf. Fig.\,\ref{GudkovCampo-5-15:fig}) to exhibit
this scheme, yet now a more elementary method is demanded. Despite
elementariness, if one is not so clever (like the writer) this
game can be pretty time consuming as demonstrated by the following
section. This consisted in a sequence of failing trials, and alas
TeX forced us to censure most of these instructive trials as
otherwise our text was not anymore synchronized with the images.

\begin{figure}[h]
\centering
\epsfig{figure=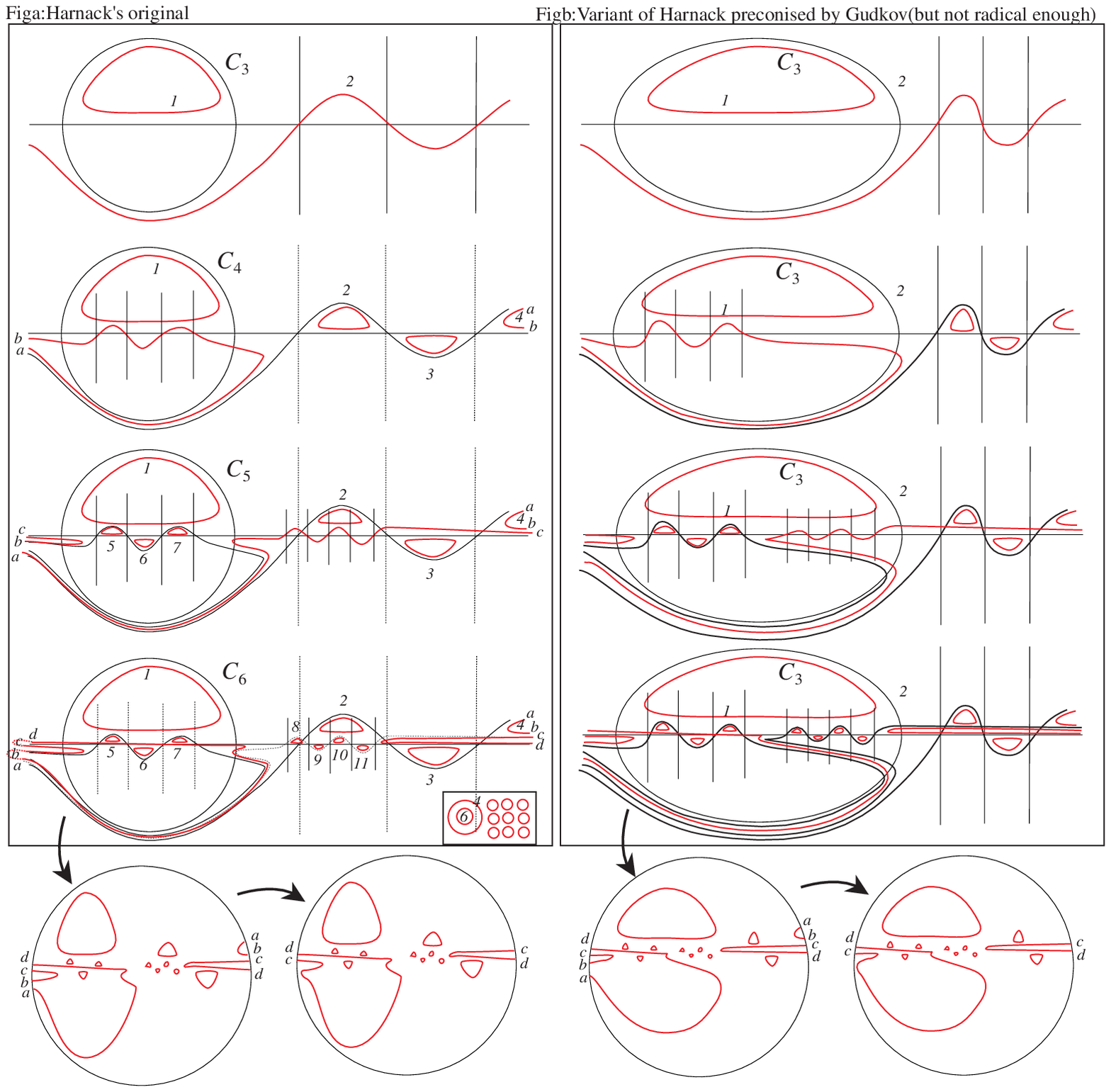,width=122mm}
\vskip-10pt\penalty0
  \caption{\label{Harnack-Gudkovvariant:fig}%
Uninspired variant of Harnack leading to Harnack's scheme again}
\vskip-5pt\penalty0
\end{figure}

For convenience the first picture
(Fig.\,\ref{Harnack-Gudkovvariant:fig}a) reminds the classical
implementation of Harnack's method of 1876 \cite{Harnack_1876}
(little warning: in the original paper the depiction is much left
to the imagination of the writer, and our picture though standard
is really inspired by nice drawings available in Viro's papers). A
first idea is to put the oscillation inside the ground circle, but
this looks too naive and we recover exactly Harnack's scheme
$\frac{1}{1}9$ (cf. Fig.\,\ref{Harnack-Gudkovvariant:fig}b).

After much efforts and trials we ultimately found (the next day
[07.02.13]) the solution as Fig.\,\ref{HarnaGudkov4-15:fig}. The
trick is to leave much room in between the vertical lines
effecting the oscillation of Harnack's method, so as to place the
subsequent vibration in between. One of the difficulty we
encountered before finding the solution is that since the desired
configuration $\frac{4}{1}5$ is an $(M-1)$-curve one is tempted to
 start with an $(M-1)$-cubic. Then one can
apparently loose much energy in the desert.

Instead we start form a Harnack-maximal cubic obtained by slight
perturbation of an ellipse $E_2\cup L$ union the horizontal line
$L$ and perturb this by a triplet of vertical lines. It results
the black depicted $C_3$ on the first row of
Fig.\,\ref{HarnaGudkov4-15:fig} intersecting thrice the horizontal
line.

The quartic curve $C_3\cup L$ is then perturbed by a quadruplets
of lines. Those could be a priori be located everywhere, but we
choose them in between as depicted on the figure. Here it seems
quite crucial that as the number of lines is even we may
concentrate  the vibration on a single oval. After this vibration
the large central oval looks like a pair of Ray-Ban eyeglasses
(viewed in perspective). We have now a $C_4$ oscillating 4 times
across the (horizontal) line $L$ and we perturb again the union.
How to do this? Always by the same method but we are free to
choose the location of the vibrator. A priori since an oval can
vibrate an even number of times across a line we may want to
choose only 4 vibrations in the ``nearby glass'' of the Ray-Ban
and one outside. This leads to something interesting namely the
scheme $\frac{4}{1}4$ (compare Fig.\,\ref{indef414:fig}, which we
transported earlier in this text). Here instead we keep the 5
vibratory lines close together (this usually maximizes their
vibratory impact), and all inside the big glass of the Ray-Ban,
cf. second row of the figure. It remains to depict the resulting
smoothing of $C_4\cup L$. The result is the red curve $C_5$
depicted but it is essential to choose this oscillation (and not
the opposite one) in which case you destroy many ovals (this will
be depicted concretely on the next
Fig.\,\ref{HarnaGudkov3-14:fig}). So there is something like a
snake visiting the nearby glass of the Ray-Ban. This gives a $C_5$
traversing 5 times the line $L$.

Smoothing the union $C_5\cup L$ produces the sextic of the 3rd row
(of Fig.\,\ref{HarnaGudkov4-15:fig}). Note that the two
(red-colored) branches nearby the horizontal line are linked
together at $\infty$ to form a single circuit, which we call the
median circuit. More generally all branches going to infinity are
connected with the diametrically opposite branch. The median
circuit of the $C_6$ is clearly the unique nonempty oval. What
appears naively in its interior is in reality a M\"obius band (due
to the diametral identification), hence its interior really
contains 4 ovals. This shows that the constructed curve realizes
the desired scheme $\frac{4}{1}5$.

\begin{figure}[h]
\centering
\epsfig{figure=HarnaGudkov4-15.eps,width=122mm}
\vskip-10pt\penalty0
  \caption{\label{HarnaGudkov4-15:fig}%
A variant of Harnack producing the scheme $\frac{4}{1}5$}
\vskip-5pt\penalty0
\end{figure}

Somewhat against our expectation this curve cannot be simplified
toward the scheme $\frac{3}{1}5$ we were interested in (in the
previous section), as the 4 inner ovals are not coming from a
vibration. Nonetheless during our exploration up to finding this
premaximal scheme  $\frac{4}{1}5$ nearly Gudkovian, we found a
variant of the exposed construction yielding the scheme
$\frac{3}{1}5$ (cf. Fig.\,\ref{HarnaGudkov3-15XXL:fig} below).

First let us choose the opposite mode of vibration as the lucky
one we first depicted. This gives Fig.\,\ref{HarnaGudkov3-14:fig}.
Now the snake oscillates around the nearby glass, wind around the
nose of the investigator, to loop around the second (distant)
glass of the Ray-Ban, etc. (As we must optically smooth the union
of both black curves $C_4\cup L$, we could a priori hope to close
up an oval with the bottom half of the first close glass, but this
forces a 6th intersection in $C_5\cap L$ violating B\'ezout.) On
smoothing $C_5\cup L$ we find a curve realizing the scheme
$\frac{3}{1}4$. Although not so exciting as $\frac{3}{1}5$, this
is already interesting for the purpose of the previous section
(namely trying to corrupt Nikulin).

So the game is to trace the triangle through the 3 inner ovals and
look upon the separation it effects upon the 4 outer ovals. Of
course our picture has poor metrical qualities as we blew it up
topologically as to see what happens within the viscera of
Harnack's method. Notwithstanding
the naive green triangle depicted seems to leave unseparated the 4
outer ovals (which to me remembered appears ``inside''). However
upon dragging the upper vertex below while staying in the outer
oval residual to the upper semi-circle, we can easily (at least on
our topological picture) effect a separation. Remind (from the
reasoning of the previous section) that if such a division occurs,
then the rigid-isotopy classification of Slava Nikulin 1979
\cite{Nikulin_1979/80} is violated. Can we infer anything serious
from such a topological picture of Harnack method? Maybe we can
via a mental contraction of some ovals restore some metrical
faithfulness in the depiction as to be sufficiently accurate to
answer the (non)separation question by the fundamental triangle
through the 3 (deep) inner ovals.

Let us start with the observation that the initial cubic $C_3$
looks on our distorted picture (Fig.\,\ref{HarnaGudkov3-14:fig})
more like a quintic (consider a line ``parallel'' to the
horizontal one passing through the unique oval of the $C_3$). So
in the real picture the central oscillating bump of the cubic is
much less pronounced. Imagine the oval of the cubic as a sun
radiating light, then there cannot be shadow lying behind the hill
formed by this bump (otherwise 4 intersections with a line too
much  for B\'ezout).

So the real picture is heuristically like the 4th row of
Fig.\,\ref{HarnaGudkov3-14:fig}. In particular the Ray-Ban glasses
(=vibrating oval of the $C_4$) is much stretched vertically. One
may argue that the Ray-Ban glass traps the oscillation, and also
the resulting 4 outer ovals created in the last step of Harnack's
iteration. Accordingly it seems sufficient to use the $C_4$ as a
sort of envelope. Since the 3 ovals of the quartic $C_4$ distinct
from the Ray-Ban are actually (modulo infinitesimal perturbations)
the 3 inner ovals of the final sextic $C_6$, and noting also that
the line through two of them cannot intersect the Ray-Ban oval
(B\'ezout), we may conclude that any triad of lines through the
inner ovals of our $C_6$ does not separate the 4 outer ovals. This
no-separation scenario is in accordance with our previous
depiction of such a curve via the more user friendly Hilbert's
method (cf. Fig.\,\ref{Lock3-14:fig}). Alas or fortunately our
reasoning does not foil Nikulin's theorem. We summarize this
trapping argument by the:

\begin{lemma}
For the sextic curve $C_6$ of scheme $\frac{3}{1}4$ realized via
Harnack's method the fundamental triangle through the deep ovals
does not separate the $4$ outer ovals.
\end{lemma}

\begin{figure}[h]
\centering
\epsfig{figure=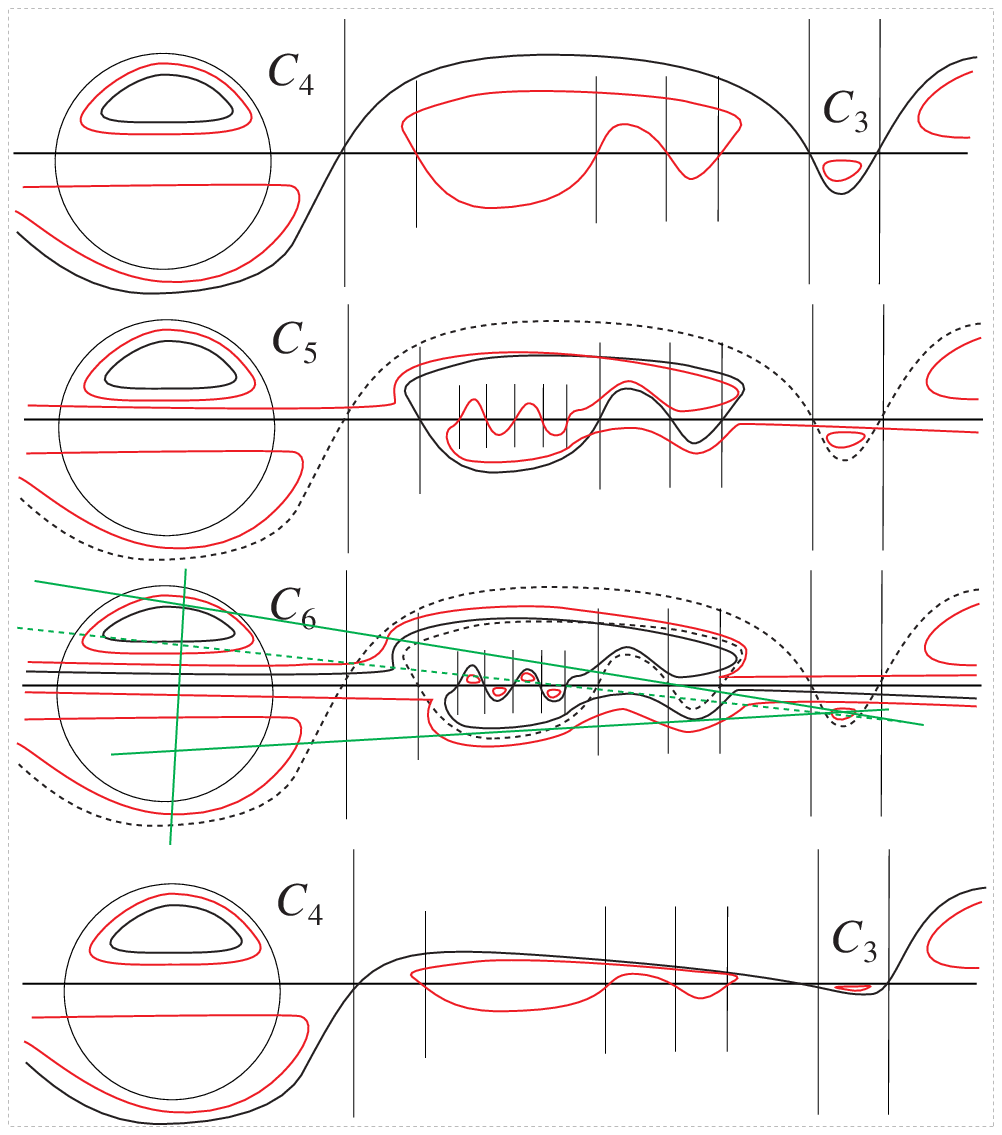,width=122mm}
\vskip-10pt\penalty0
  \caption{\label{HarnaGudkov3-14:fig}%
A variant of Harnack producing the scheme $\frac{3}{1}4$}
\vskip-5pt\penalty0
\end{figure}

At this stage, it is with a mixture  of happiness and
disappointment that Nikulin still seems to resist our naive
aggression via the Fiedler-Marin method. It remains however to
look at the scheme $\frac{3}{1}5$.

[08.02.13] How to realize it? Again several tests are required and
usually we (at least the writer) lack an understanding of the
predestination governing Harnack's method. Using the technique of
the microcosmic vibration ``in between'' we realized the schemes
$\frac{4}{1}5$, $\frac{4}{1}4$ and $\frac{3}{1}4$ all hitting
quite central positions of  Gudkov's pyramid
(Fig.\,\ref{Gudkov-Table3:fig}). But how to get $\frac{3}{1}5$
lying more ``eccentric'' on this table?
Incidentally one could dream that this Harnack method we are using
leads to the eclectic Gudkov scheme $\frac{5}{1}5$. Of course this
would corrupt experimental evidence
assembled along
centennial working
tradition (Harnack 1876, Hilbert 1891, Rohn 1888--1913, Brusotti
1910--1945, Gudkov 1954--1973, etc.). However we do not know
(personally) a theoretical obstruction impeding Harnack's method
to produce Gudkov's scheme. Arguably if well assimilated Harnack's
method reduces to a finite collection (for a fixed degree say
$m=6$) of combinatorially distinct locations for the vibratory
lines. So it suffices to explore all choices and notice that
Gudkov's scheme never appears. We do not claim to be clever enough
to complete this boring exercise, but our microfilm picture
perhaps contributes to this (cf.
Fig.\,\ref{HarnaGudkov3-15MICRO:fig}). [08.04.13] It is not clear
at this stage if this picture will be publishable in the arXiv due
to size limitations.

But let us return to our main duty of exhibiting $\frac{3}{1}5$.
Here a series of tests given in micro-film format
(Fig.\,\ref{HarnaGudkov3-15MICRO:fig} only consultable on a PC
where one can zoom and alas unreadable on the paper). Alas we
cannot give the pictures in decent format for otherwise the flow
of pictures overrun dramatically what we have to say on the topic.
We are in the realm of pure geometry were only pictures have some
weight, but alas this does not seem to please my TeX-compilator,
who accept at most two pictures per page.
Here the second column picture of this microfilm shows an
interesting variant of the scheme $\frac{3}{1}4$ where the 4 outer
ovals are not directly enveloped by the ``Ray-Ban'' oval, and so
our former argument does not readily apply here. It seems however
dubious to expect a corruption of Nikulin. Without getting
sidetracked by this issue, keep in mind our goal of realizing
$\frac{3}{1}5$.

\begin{figure}[h]
\centering
\epsfig{figure=HarnaGudkov3-15MICRO.eps,width=122mm}
\vskip-10pt\penalty0
  \caption{\label{HarnaGudkov3-15MICRO:fig}%
A microfilm cataloging several variants of Harnack trying to find
the scheme $\frac{3}{1}5$, cf. next picture for the solution}
\vskip-5pt\penalty0
\end{figure}

After several trials (cf. again the microfilm
Fig.\,\ref{HarnaGudkov3-15MICRO:fig}) we arrived at the idea of
using the same vibratory configuration of lines as for
$\frac{4}{1}5$ safe that instead of starting from an $M$-cubic we
start from an $(M-1)$-cubic. This seems to require locating one of
the vibratory line inside the circle. Our final picture is
Fig.\,\ref{HarnaGudkov3-15XXL:fig}. It hardly deserves to be
commented upon once it is found, except for saying that the
initial cubic is to be thought of as a small perturbation of the
circle $E_2$ union the line, despite sembling a large deformation
thereof. The trick in tracing  Harnack's curves is always to
exaggerate  small perturbations as to create some free room to
depict the next stage of the inductive process (vibratory
pudding). This is of course possible due to the malleability of
the continuum $\RR$ of real numbers.

On this figure (Fig.\,\ref{HarnaGudkov3-15XXL:fig}) we recognize
again our Ray-Ban oval except that it has now acquired a
``branch'' (compare 2nd row of
Fig.\,\ref{HarnaGudkov3-15XXL:fig}). Again our interest is to
apply the lock method of Fiedler-Marin. So we trace the triplet of
lines through 3 points in the deep (inner) ovals $1,2,3$, and
examine whether and how this triangle splits apart the outer ovals
{\it 1,2,3,4,5}. (Notice the importance of italicization in our
notation: italics are outer ovals while roman-arabic numbers are
the inner ovals.) In contrast to the Harnack curve realizing
$\frac{3}{1}4$ where all the 4 outer ovals were encapsulated in
the Ray-Ban oval, we notice now that the oval {\it 5} lies outside
this (Ray-Ban) oval. So our former argument does not readily
apply. Notice also that the inner oval $2$ lies inside the
Ray-Ban. All this is a bit puzzling but should not discourage us
attempting to study the division of the outer ovals by the locking
triangle for our Harnack-modified curve $C_6$. Note incidentally
that the latter is not perfectly well defined as a curve unless we
specify exactly the deformation constants involved in Harnack's
small perturbation method. Yet it seems natural to expect that the
combinatorial data involved in our Harnack-Gudkov style
description is enough to determine unambiguously a rigid-isotopy
class. Hence by the Fiedler-Marin locking argument (involving
merely B\'ezout saturation) we
infer that the distribution of outer ovals within the 4 components
past the lock is well-defined. It remains only to determine it.

\begin{figure}[h]
\centering
\epsfig{figure=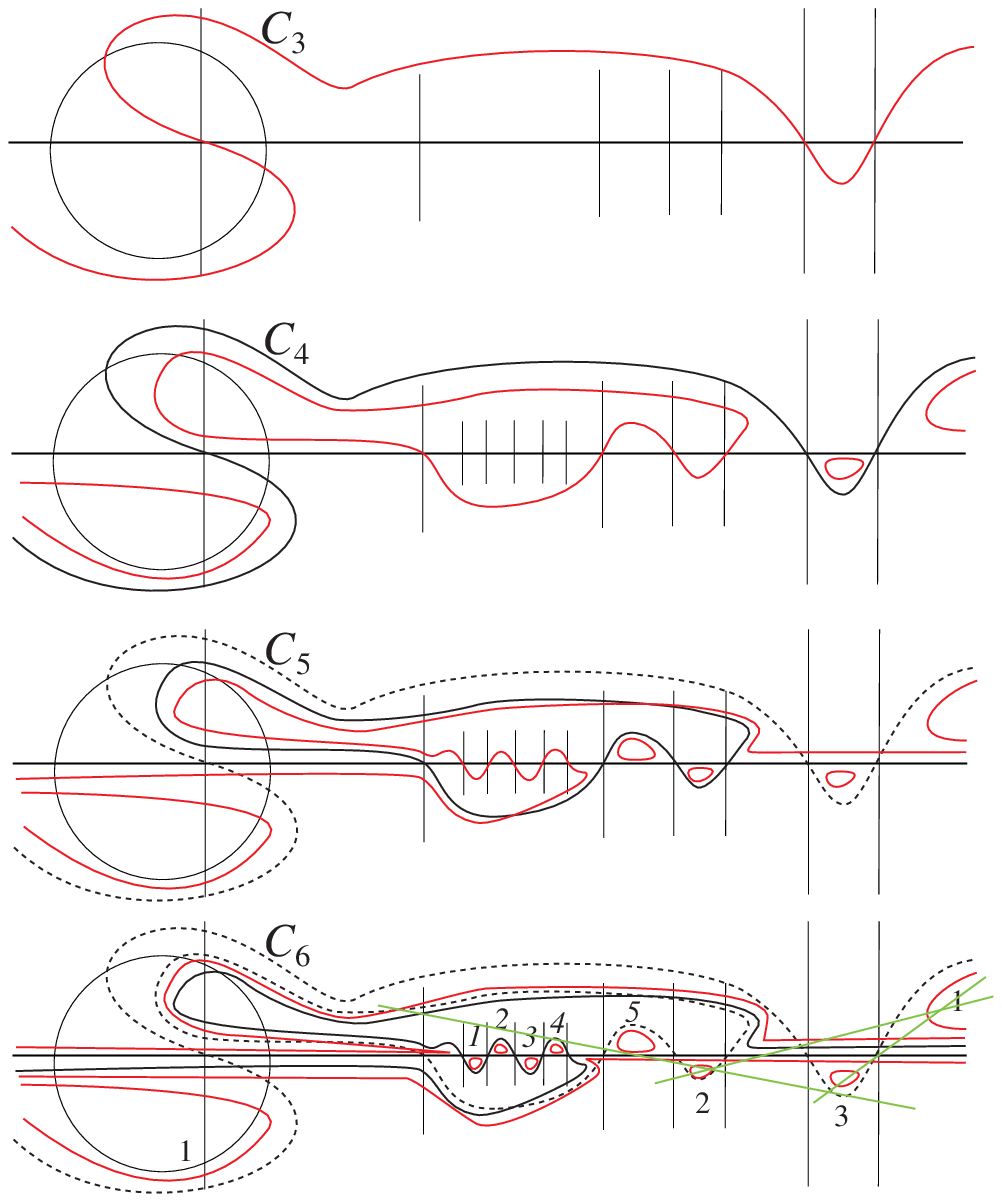,width=122mm}
\vskip-10pt\penalty0
  \caption{\label{HarnaGudkov3-15XXL:fig}%
A variant of Harnack producing the scheme $\frac{3}{1}5$ (known to
Gudkov 1954, and surely earlier: Harnack 1876?, Hilbert 1891?,
Ragsdale 1906, Brusotti 1910, etc., but tricky to find)}
\vskip-5pt\penalty0
\end{figure}

In the sequel we shall often speak of ``the line passing through
two disjoint ovals''. This is a slight abuse of language for such
a line is not uniquely defined, and is really intended to mean
choose 2 points in the insides of the 2 disjoint ovals and trace
the line joining them. Of course the phrasing ``the line'' becomes
somewhat sloppy, but when the two ovals are inner ovals then any
such line is B\'ezout saturated, and so from the viewpoint of
analysis situs there is some canonicalness.

Since ``the'' (or a) line through the ovals $1,3$ regarded on the
quartic $C_4$ cannot cut more times the $C_4$, we infer that it
does not cut the Ray-Ban oval of the $C_4$ (the one oscillating 4
times across the horizontal line $L$). Next the line through the
ovals $2,3$ interpreted on the quintic $C_5$ cannot cut more this
curve safe for a point on its pseudoline. A similar remark holds
for the line through $1,2$. All this looks a bit sterile and we
really need the geometry of the picture to understand the
distributional question. For this purpose, look at the 3
green-colored lines on Fig.\,\ref{HarnaGudkov3-15XXL:fig}, while
enlarging slightly oval~3 as to adjust the picture. It seems then
that the triangle separates oval {\it 5} from the ovals {\it
1,2,3,4}. Of course upon stretching further oval $3$, we could
arrange that the line $2,3$ passes above oval {\it 5} in which
case the locking triangle   effects no subdivision of the outer
ovals. Which of both scenarios corresponds  the reality? A priori
the first scenario looks more likely (at least in line with our
picture). Remind however the slogan (anonymous, Poincar\'e, etc.)
``La g\'eom\'etrie c'est l'art de bien raisonner sur des figures
mal dessin\'ees''.

Let us attempt a more realist depiction on the following figure
(Fig.\,\ref{HarnaGudkov3-15XXLTEST:fig}). Even the first
right-side picture (fig.\,d) is not B\'ezout permissible (the
green-line cut the cubic $C_3$ five times). Further it may be
observed that the line through $2,3$ may pass ``below'' the series
of ovals {\it 1,2,3,4}. This is a third possible scenario in which
there is no subdivision.

\begin{figure}[h]
\centering
\epsfig{figure=HarnaGudkov3-15XXL-USED.eps,width=122mm}
\vskip-10pt\penalty0
  \caption{\label{HarnaGudkov3-15XXLTEST:fig}%
Ovals distribution  past the deep triangle on Harnack's model of
 $\frac{3}{1}5$} \vskip-5pt\penalty0
\end{figure}

Admittedly this question
looks quite tricky to decide and requires some good idea or
high optical acuity. Harnack's method seems not ideally suited to
clinch the matter. Hilbert's method would be more convenient, yet
does not seem capable producing the scheme $\frac{3}{1}5$ which is
slightly more on the Harnack right-hand side of Gudkov's pyramid
(Fig.\,\ref{Gudkov-Table3:fig}). This is surely no intrinsic
reason since Harnack's scheme itself is accessible to Hilbert's
method. Inspecting carefully our former Fig.\,\ref{GudHilb8:fig}
cataloging several variants of Hilbert's method it is
pretty clear why Hilbert's method fails producing the scheme
$\frac{3}{1}5$. Indeed what comes closest to $\frac{3}{1}5$ is the
scheme $\frac{4}{1}4$ depicted near the center of
Fig.\,\ref{GudHilb8:fig}, and one may argue that the vibrating
oval has always an even number of (cytoplasmic) expansions coming
across the fundamental ellipse $E_2$ of Hilbert's construction
thereby creating an odd number of ovals. So to have 3 inner ovals
requires 2 inner expansions like on the scheme $\frac{3}{1}3$ on
the bottom of Fig.\,\ref{GudHilb8:fig} but this
dissipates too much of the oscillating energy and not enough outer
ovals are created.

As an attempt to corrupt Nikulin it would be interesting to detect
a second realization of the scheme $\frac{3}{1}5$ besides
Harnack's presented above, and study on it the distribution of
outer ovals past the fundamental triangle through the deep ovals.
Even without any scepticism about Nikulin the net effect is that
if the latter is correct
the determination of this distribution for a single curve $C_6$
belonging to the scheme $\frac{3}{1}5$ which is of type~II  would
by the Fiedler-Marin argument determines this distribution for all
curves belonging to the scheme. (That such a curve is of type~II
necessarily, follows from Arnold's congruence
(\ref{Rohlin-implies-Arnold:lem}) $(3=1-3+5=)\chi=p-n\equiv k^2
\pmod 4(=3^2=9\equiv 1)$ valid for all type~I curves.) We get so
some nice geometric theorem as a consequence of the truth of
Nikulin. (Insertion [08.04.13]: a much more elementary argument is
given by Le~Touz\'e in Sec.\,\ref{LeTouze:sec}.)

[09.02.13] We now present some tricky argument in favor of
non-separation of the outer ovals on the Harnack model constructed
above.

Referring to Fig.\,\ref{HarnaGudkov3-15XXLTEST:fig}a, let us look
at the quartic $C_4$ with $r=3$ ovals occurring as an intermediate
step of Harnack's construction. The unique oval of the $C_4$ which
oscillates 4 times across the horizontal line $L$ is referred to
as the {\it Ray-Ban oval\/}. Label the intersection points in
$C_4\cap L$ as $p_1,p_2,p_3,p_4$. Consider the 3 green-lines
through the 3 inner ovals of the sextic $C_6$
(Fig.\,\ref{HarnaGudkov3-15XXLTEST:fig}c) but imagined traced on
this $C_4$ (i.e. on plate Fig.\,a). Since the oval $2$ of the
$C_6$ is enveloped by the Ray-Ban oval (of the $C_4$) it may be
inferred that the line $2,3$ does not intercept the Ray-Ban oval
outside of the
arc $p_3,p_4$ of the $C_4$ on Fig.\,a). Hence the line $2,3$ is
actually much more horizontal than on our Fig.\,c. More precisely
we infer the following.

Since the ovals {\it 1,2,3,4} are encapsulated in the Ray-Ban oval
(of the $C_4$), the line $2,3$ passes below them. (Of course
``below'' as no absolute sense in projective geometry, but here it
has since we have another line $L$ as reference.)

Further it also passes below oval {\it 5}, for otherwise it passes
above but having to avoid the Ray-Ban it would then have to lounge
the nasal portion of the Ray-Ban while passing between oval {\it
5} and the curvilinear arc $p_2,p_3$ of $C_4$. (Recall that our
line $2,3$ is de facto B\'ezout-saturated (w.r.t. to $C_6$ or even
w.r.t. the $C_5$), hence cannot intercept any outer oval, in
particular {\it 5}.) But then the line $2,3$ would intersect twice
the horizontal line $L$, violating the simplest case of B\'ezout.

In conclusion the line $2,3$ passes below all outer ovals {\it
1,2,3,4,5}, and the real picture could be more like
Fig.\,\ref{HarnaGudkov3-15XXLTEST:fig}e. Alas this depiction does
not seem possible because the line $2,3$ already crosses 6 times
the sextic so  the oval $1$ cannot cross this line. Since it moves
above it on the right side of the picture it must resurface on the
right side below the line $2,3$, which is not the case on our
picture (Fig.\,e). Accordingly
Fig.\,\ref{HarnaGudkov3-15XXLTEST:fig}f might be more realistic,
and the conclusion would be that the fundamental triangle does not
separate the outer ovals.

Have we proved anything? Let us say ``yes'' and state the
following lemma of which we shall supply a more formal proof right
below.

\begin{lemma}
The fundamental triangle consisting of the $3$ lines passing
through the $3$ inner ovals of Harnack's curve (depicted above as
Fig.\,\ref{HarnaGudkov3-15XXL:fig} or
Fig.\,\ref{HarnaGudkov3-15XXLTEST:fig}) realizing the scheme
$\frac{3}{1}5$ does not separate the $5$ outer ovals.
\end{lemma}

\begin{proof}
The trick toward a more formal proof is to consider the
topological disc $D$ obtained from the inside of the Ray-Ban oval
by expanding it at $p_2,p_3$ linearly while smashing it inside at
$p_3,p_4$ (compare the shaded region on
Fig.\,\ref{HarnaGudkov3-15XXLTEST:fig}a). Since this region
contains the 5 outer ovals {\it 1,2,3,4,5} it suffices to check
that this disk avoids the 3 green lines through the 3 inner ovals.

This is clear for the line $1,3$ which is B\'ezout-saturated on
the $C_4$, hence can only attack our modified disc $D$ through the
arc $p_2,p_3$, but as the inside of the Ray-Ban oval is avoided it
results a second intersection with $L$ violating B\'ezout.

The same argument works for the remaining two lines $1,2$ and
$2,3$ after noticing that since those lines intercept twice the
oval $2$ of the sextic, it may be assumed that they intercept
twice the curvilinear arc $p_3,p_4$ of $C_4$. This follows merely
from the nature of the method of small perturbation. {\it
Warning}.---In fact a priori we could imagine that the line $2,3$
penetrates in the oval $2$ much more vertically than on Fig.\,c
meaning really that it intercepts the segment $p_3,p_4$ of $L$,
but in that case too, it is clear that the fundamental triangle
does not separate the outer ovals. In fact in this case our line
$2,3$ cuts the Ray-Ban oval only twice, and we may excise from $D$
the half of the trace of our line on $D$ containing $p_4$, plus a
little tubular neighborhood thereof. During this excision it is
clear that we do not loose the covering of the outer ovals, since
the line $2,3$ is B\'ezout-saturated on $C_6$.

Now each of our lines through oval $2$ is B\'ezout-saturated with
the $C_4$ hence must avoid the inside of the Ray-Ban. However our
line  cannot penetrate the arc $p_2,p_3$, for otherwise a 2nd
intersection with $L$ is created, hence has void intersection with
$D$.
\end{proof}

\subsection{Viro's construction specialized to the
scheme $\frac{3}{1}5$}

[08.04.13] This section explores other realizations than Harnack's
(especially of the scheme $\frac{3}{1}5$) which is somewhat
cumbersome. It is primarily a matter of exploring Viro's method,
but the latter turns out to be not much more suited than Harnack's
model to fix the distribution question. As we already said the
royal road is Le~Touz\'e's argument in Sec.\,\ref{LeTouze:sec}. Of
course Viro's method has supernatural appeal too, but our
exposition is far from explaining the true core of the dissipation
method which is a secret to us. Hence this section can be omitted
with loss of continuity.

\begin{figure}[h]
\centering
\epsfig{figure=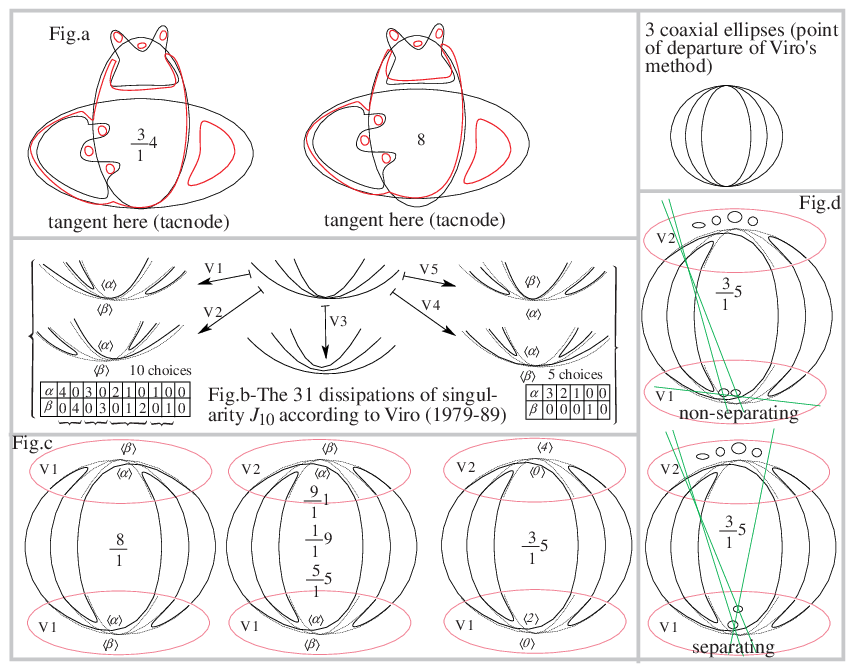,width=122mm} \vskip-10pt\penalty0
  \caption{\label{Viro3-15:fig}%
The dissipative (but all encompassing) method of Viro with special
attention to the scheme $\frac{3}{1}5$} \vskip-5pt\penalty0
\end{figure}

[08.02.13] A first idea is to use Marin's variant of Hilbert's
method but this seems only able to produce the scheme
$\frac{3}{1}4$ (cf. Fig.\,\ref{Viro3-15:fig}a). Another idea is to
use Viro's dissipation of 3 ellipses tangent at 2 points. This
being again a small perturbation method like Harnack's it is not a
priori clear that we will be in a better position to tackle the
distribution question past the deep triangle. The charming feature
of Viro's method  is its ability to create nearly all sextics as
perturbation of this configuration of 3 coaxial ellipses. To
implement this, look at Fig.\,29 in Viro 1989/90
\cite[p.\,1103]{Viro_1989/90-Construction} showing all the
possible dissipations of a germ of curve singularity of type
$J_{10}^-$ consisting of 3 real branches having a second order
tangency like on our global model of the 3 coaxial ellipses. This
Viro figure is reproduced as Fig.\,\ref{Viro3-15:fig}b below,
which includes 5 modes of dissipation denoted by us V1,V2, \dots,
V5 (V standing for Viro of course). Each of them admits an array
of permissible values for spontaneous ``champagne bubbling'' of
ovals created out of the blue. Then we can patchwork such
smoothing {\it independently\/} at both singular points of the
configuration of 3 ellipses to create a global curve with
controlled topology. (This is of course highly reminiscent of
Brusotti-Gudkov's independence of smoothing, based on Severi and
in turn upon Riemann(-Roch) via possible d\'etours through the
Plato cavern of Brill-Noether. Compare Brusotti 1921
\cite{Brusotti_1921} and Gudkov 1974 \cite{Gudkov_1974/74}.)

For instance the dissipation V1-V2 with $(\alpha,\beta)=(4,0)$ and
$(\alpha,\beta)=(4,0)$ resp. yields Hilbert's scheme
$\frac{9}{1}1$. Choosing instead for V1, $(\alpha,\beta)=(0,4)$
and for V2 $(\alpha,\beta)=(0,4)$ yields Harnack's scheme
$\frac{1}{1}9$. If we choose for V1, $(\alpha,\beta)=(4,0)$ and
for V2 $(\alpha,\beta)=(0,4)$ yields Gudkov's scheme
$\frac{5}{1}5$.
Pause a little moment at this stage, to be puzzled by the fact
that Gudkov's rare bird---which escaped the attention of all
experts during 8 decades (from Hilbert to Gudkov)---appears in the
fingers of Viro as a species not much more tropical than common
birds like Harnack and Hilbert. Perhaps this banalization of
Gudkov by Viro is against the philosophy expressed in
Sec.\,\ref{Diophantine-and-proba:sec} that Gudkov's curve(s) ought
to have some statistical and Diophantine rarity.

How to realize $\frac{3}{1}5$? It suffices to choose at V1
$(\alpha,\beta)=(2,0)$ and at V2 $(\alpha,\beta)=(0,4)$ (cf.
Fig.\,c, right).
Does Viro's method help to solve our distributional problem. A
priori not as shown by our naive depiction
Fig.\,\ref{Viro3-15:fig}d exhibiting both distribution (separating
or not) as  logically possible a priori. Recall yet that on behalf
of Nikulin's theorem both options cannot occur simultaneously.
Of course naive geometric intuition tells us that the upper part
of Fig.\,d is more likely with small ovals spread horizontally  as
a vestige of the horizontal tangent line at the singular point of
the initial configuration of 3 coaxial ellipses. But of course
even in this situation where the 2 microcosmic ovals (generated by
the dissipation of the bottom singularity ``V1'') are nearly
horizontal (and so is the line through them) this does not prevent
the ``vertical'' green-lines to separate the upper series of 4
ovals. But again on ground of some microscopical geometric
intuition it seems  realist to argue that even if the top vertex
of the green-triangle is chosen very near to the top of the banana
oval (i.e. the large inner oval of the curve $C_6$ resulting from
the fusion of the two branches $2,3$ of ellipses when labelled
from left to right), the little 4 top ovals will condense
themselves as to be non-separated by the deep triangle. Actually
if a separation would occur, then by sweeping the nearly vertical
green-line inside the pencil of lines rooted at a basepoint in the
bottom oval gives by continuity an intermediate line with 8
intersections with the sextic $C_6$ overwhelming B\'ezout (or the
smoothness of the $C_6$). Indeed if the line through one of the
bottom micro-oval has slope of ca. 135 degree it cuts twice the
micro-oval, twice the banana and twice the (largest) nonempty
oval, hence 6 times the curve. If we let this angle diminishes to
90 degree (plus $\varepsilon$) by dragging the upper (banana)
vertex up to the top of the banana while supposing that the pair
of lines effects a division of the 4 top micro-ovals for a
suitably small value of $\varepsilon $ then both lines have again
6 intersections, but in between $135$ degrees and $90+\varepsilon$
degrees there must be a line cutting 8 times the curve $C_6$
namely the one line sweeping the separated oval. Sorry for this
messy argument. Of course the key is just to observe that the
green-lines are B\'ezout-saturated, hence the distribution of
outer ovals cannot change.

By the same sort of argument, precisely by tracing the line
through two points near the bottom of the inner banana and the
outer banana we get a line which is already B\'ezout-saturated.
Pushing this line to its ultimate confinement we get the bitangent
through the most ``meridional'' points of both bananas. Following
this motion by continuity implies that the 2 bottom ovals are
pressed down below this bitangent, and so looks nearly horizontal.
Alas this does not prevent the situation of the bottom half of
Fig.\,\ref{Viro3-15:fig}d where both microscopic ovals are sitting
nearly one above the other but of course at much lesser height
than depicted, that is below the bitangent to the most meridional
portion of both bananas.

All this is quite exciting for the imagination, but does not seem
to answer our puzzle on the distribution of outer ovals past the
fundamental triangle. To analyze better the situation we should
introspect  in more detail the quantitative geometric aspects of
Viro's construction which certainly includes answers to our basic
question.

Ignoring that issue for the moment, we can introduce the concept
of the bundle spanned by two disjoint ovals. This is the
collection of all lines traced through a pair of points chosen
inside the respective ovals (``boundaries'' included). This bundle
is often supported by (or spanning)  a bordered surface
homeomorphic to a M\"obius band. If we think of both ovals as
celestial bodies (like Earth and Moon) then this bundle (or rather
its support, i.e. what is swept out by this collection of lines)
is essentially the region where one oval masks the other (at least
partially) like during an eclipse. We call thus this region the
{\it eclipsus\/} of both given ovals, or just their mutual {\it
shadow\/}.

The shadow of the 2 bottom micro-ovals in case their mutual
disposition is nearly vertical (like on the bottom half of
Fig.\,\ref{Viro3-15:fig}d) cannot intercept any further oval than
the 3 obvious one (each of the 2 protagonists plus the nonempty
oval enclosing them). In particular in that case of nearly
vertical alinement of both ``meridional'' micro-ovals their shadow
must find its way out through the little room left vacant between
the top part of the outer banana and the top 4 micro-ovals. This
is actually possible as suggested by
Fig.\,\ref{Viro3-15shadows:fig}c, provided both ovals really live
at the microscopic scale. Fig.\,\ref{Viro3-15shadows:fig}a depicts
the shadow of both bananas. Each line in this shaded region is
B\'ezout-saturated, hence no ovals can
survive in this region, hence the situation is forced to be like
on Fig.\,\ref{Viro3-15shadows:fig}b. Once we have Fig.\,c then we
must still analyze further shadows, but some thinking at the
nanoscale should convince the reader that the verticality scenario
posited by Fig.\,c is not further obstructed by Monsieur \'Etienne
B\'ezout. Philosophically algebraic plane curves are like
celestial configurations not liking to have their horizon too much
saturated by galactic nebulosity. They express a principle of
economy and purity.

\begin{figure}[h]
\centering
\epsfig{figure=Viro3-15shadows.eps,width=122mm}
\vskip-10pt\penalty0
  \caption{\label{Viro3-15shadows:fig}%
Ovals distribution past the deep triangle on Viro's model of
$\frac{3}{1}5$} \vskip-5pt\penalty0
\end{figure}

Our conclusion is that the topological aspect of Viro's method
alone does not seem sufficient to settle the distributional
question of a $C_6$ of type $\frac{3}{1}5$. However it is more
likely that for some Viro curve $C_6$ the distribution of the
bottom ovals is nearly horizontal, and therefore that the
fundamental triangle through the 3 inner ovals does not separate
the 5 outer ovals of this $C_6$. If this is true (and Nikulin's
theorem also) then we deduce first the following lemma and next
the following theorem by uniting the forces of all lemmas of the
previous Sec.\,\ref{Nikulin-corruption:sec}:

\begin{lemma}
Any sextic $C_6$ with scheme $\frac{3}{1}5$ is such that the
fundamental triangle through the $3$ inner ovals does not separate
the five outer ovals.
\end{lemma}

\begin{theorem}\label{Nikulin-Fiedler-Marin-Gabard-no-separa:thm}
Any sextic $C_6$ belonging to a scheme of the form
$\frac{3}{1}\ell$ (with $0\le \ell \le 5$ according to Gudkov's
table) is such that the fundamental triangle through the $3$ inner
ovals does not separate the outer ovals.
\end{theorem}

If this theorem is true, one may of course wonder if there is an
elementary proof circumventing the highbrow intervention of K3
surfaces, Torelli, etc. i.e. all the technology involved in
Nikulin's proof. [13.02.13] Such an elementary proof is given in
the next section (\ref{LeTouze:sec}) and was communicated by
Le~Touz\'e.

\subsection{Fiedler-Le~Touz\'e's answer (conical chromatic law)}
\label{LeTouze:sec}

[13.02.13] Three days ago, I received the following answer from
S\'everine Fiedler-Le~Touz\'e (n\'ee Le~Touz\'e, and often
abridged as a such too avoid confusion with her husband Fiedler,
who also worked in the field ``a long time ago'').

$\bullet$$\bullet$$\bullet$ [10.02.13] Bonjour Alexandre, Thomas
m'a transmis ta question. La r\'eponse est toute simple: soient
$A$, $B$, $C$ trois ovales int\'erieurs et $D$, $E$ deux ovales
exterieurs de ta sextique. Le triangle fondamental $ABC$ est
enti\`erement contenu dans l'ovale non-vide. Si $D$ et $E$ sont
dans deux triangles $ABC$ (non-fondamentaux) diff\'erents, alors
la conique passant par $A$, $B$, $C$, $D$, $E$ coupe la sextique
en $14$ points, contradiction. Avec des coniques, on montre plus
g\'en\'eralement que: Les ovales vides de la sextique sont
distribu\'es dans deux chaines (int, ext), l'ordre cyclique est
donn\'e par les pinceaux de droites bas\'es dans les ovales
interieurs. Les ovales interieurs sont dispos\'es en position
convexe dans l'ovale non-vide. Bon dimanche,  S\'everine

\smallskip
Translated in my poor English this gives the:

\begin{lemma}\label{LeTouze:lem}
Let $C_6$ be a sextic with $3$ inner ovals, and at least $2$ outer
ovals, then the latter are
in the same component
past the fundamental triangle consisting of the $3$ lines through
the deep
ovals. In particular any sextic of
type $\frac{3}{1}\ell$
has all its ovals distributed in the same component
past the fundamental triangle.
\end{lemma}

{\it Insertion} [08.04.13].---As a loose idea it could be
interesting to see if the method can be boosted as to prohibit the
scheme $\frac{3}{1}6$. Of course this follows also via total
reality of the scheme $\frac{2}{1}6$ also due to Le~Touz\'e.

\smallskip
\begin{proof}
Assume on the contrary that the 2 outer ovals are in different
subregions past the fundamental triangle. The idea is to look at
the conic passing through the $3+2=5$ ovals (3 being deep and 2
being outer ovals). This conic certainly cuts the $C_6$ in at
least $(5+1)\cdot 2=12$ points, like say on
Fig.\,\ref{LeTouze:fig}a albeit of course no depiction is required
since those intersections are so-to-speak topologically forced.
However from our supposition that the 2 outer ovals are separated
by the fundamental triangle the real picture is rather like
Fig.\,\ref{LeTouze:fig}b or Fig.\,\ref{LeTouze:fig}c yielding
$2\cdot 5+4=14$ intersections. B\'ezout is overwhelmed.

\begin{figure}[h]
\centering
\epsfig{figure=LeTouze.eps,width=122mm} \vskip-5pt\penalty0
  \caption{\label{LeTouze:fig}%
  Le~Touz\'e's argument showing that the outer ovals are not
  separated by the fundamental triangle} \vskip-5pt\penalty0
\end{figure}

It remains to find the intrinsic reason of why this holds true.
The key is to look at the order of the 5 assigned points on the
conic. On Fig.\,a the 3 inner points (in black) are not separated
by the outer ones (white-colored), while on Fig.\,b the 3 inner
points are separated by the 2 outer points. This explains the
formation of extra intersections whenever we have to
salesman-travel from the inside to the outside of the nonempty
oval. More precisely when the 3 black inner points are separated
by the 2 white outer points then we see 4 elliptical arcs with
dichromatic boundary, each of which contributing for an
intersection (with the nonempty oval separating the outer from the
inner ovals). So the whole story is reduced to the following
Hilfssatz somewhat hard to state elegantly:

\begin{lemma} {\rm (Chromatic law for conics)}.---
\label{LeTouze-Gabard-Hilfssatz:lem} Let us be given $3$
black-colored points in the plane $\RR P^2$ which are not aligned,
and let $T$ be the triangle through them. Take $2$ additional
(white colored) points outside the triangle. Then the unique conic
$E$ through the $5$  points considered is smooth. Further if both
white points are in the same component past the triangle $T$, then
the $5$ points are monochromatically distributed on the conic as
white-white-black-black-black (with only $2$ chromatic transitions
when reading cyclically). If instead both white points are
separated by the triangle, then the distribution is dichromatic as
white-black-white-black-black, which read cyclically gives $4$
chromatic transitions (from black to white or viceversa).
\end{lemma}

Applying this Hilfssatz concludes the proof of the lemma
(\ref{LeTouze:lem}).
\end{proof}

The philosophical outcome of this brilliant argument (communicated
by S\'everine
Le~Touz\'e) is that we cannot hope any corruption to Nikulin by
the Fiedler-Marin locking technique. In fact the distribution of
outer ovals past the fundamental triangle is always monopartite.
In particular we get a more conceptual and lucid proof of several
lemmas that we  tried hard to establish on the cumbersome models
of Harnack, Viro. In particular, we get an elementary proof of
Theorem~\ref{Nikulin-Fiedler-Marin-Gabard-no-separa:thm} without
the whole transcendental apparatus behind Nikulin's theorem (K3
surfaces, Torelli, etc.).


\subsection{Trying to extend Nuij-Dubrovin rigid-isotopy
of the deep nest via total
reality}\label{Nuij-Dubrovin-extended:sec}

[23.01.13] Apart from the beautiful result of Nikulin on sextics,
stating that the real scheme enhanced by Klein's types suffices to
determine unambiguously the rigid-isotopy class
(\ref{Nikulin:thm})
and the rigidity of the empty scheme,
the only positive general result available is Nuij-Dubrovin's
theorem stating that the deep nest constitutes a unique
rigid-isotopy class.

\begin{defn}\label{rigid-scheme:defn}
{\rm Let us say that an ($m$)-scheme (i.e. a scheme of order $m$)
is {\it rigid\/} if any two $m$-tics curves representing the
scheme are rigid-isotopic.}
\end{defn}

Since the deep nest is totally real (\`a la Ahlfors) under a
pencil of lines one might wonder if other totally real schemes
also enjoy rigidness. Examples of such total schemes include all
the $M$-schemes by Bieberbach-Grunsky, the $(M-2)$-schemes
$\frac{6}{1}2$ and its mirror $\frac{2}{1}6$ by an unpublished
argument of Rohlin (nobody is able to reconstruct). More easily it
includes the scheme of degree 8
$\frac{1}{1}\frac{1}{1}\frac{1}{1}\frac{1}{1}$ (4 nests of depth
2) which is total under a pencil of conics.

Conjecturally for a scheme, type~I implies maximal (Rohlin 1978).

\begin{defn}
{\rm
$\bullet$ A scheme of degree $m$ is {\it total of order $k$} if
any curve $C_m$
representing the scheme admits a total pencil of $k$-tics. (For
instance the $2k$-scheme $(1,1,\dots,1)$ consisting of
$k$ nested ovals is total of order 1.)

$\bullet$ Say that a scheme of degree $m$ is total if any curve
representing
the scheme admits a total pencil.}
\end{defn}

It seems natural to expect the:

\begin{conj} If a scheme of degree $m$ is total
then it is total of order $k$ for some universal $k$ depending
only on $m$.
\end{conj}

 At first glance
this could follow from Ahlfors 1950 \cite{Ahlfors_1950}.
At any rate we have the implications:

\begin{lemma}
Total of some order $k$ $\Rightarrow$ total $\Rightarrow$ type~I
(perhaps implying maximal).
\end{lemma}

Marin convinced me that the transition from the abstract to
embedded viewpoints might be not so easy, hence the converse of
the second arrow might
not be a trivial corollary of Ahlfors 1950, but we naively still
hope so.

One may speculate on an extension of Nuij-Dubrovin's theorem as
follows:

\begin{conj}
For a scheme, totality of order $k$ (for some fixed $k$) implies
rigidness.
\end{conj}


{\it Insertion} [08.04.13].---This conjecture (like the subsequent
one (\ref{M-schemes-rigid})) is probably also disrupted by Marin's
 obstruction (discussed in the next Sec.\,\ref{Marin:sec}).
Indeed by our extrinsic variant of Riemann-Bieberbach-Grunsky
(\ref{total-reality-of-plane-M-curves:thm})  the total reality of
plane $M$-curves of degree $m$ is exhibited by a pencil of curves
of degree $(m-2)$, hence we have totality for some universal
degree, namely $k=m-2$ (depending only on the degree $m$ of the
scheme and not on the geometry of the curve), yet no rigidity can
be observed by Marin's obstruction.


\smallskip

Perhaps it suffices to assume type~I, or maybe even maximal
implies rigidity (in increasing order of hazardousness). Here
maximality is interpreted in the sense of Rohlin (as opposed to
Harnack's more specialized sense).
Of course a priori there is very little evidence for a direct
correlation between those concepts. (Again try to look if Fiedler,
Marin give some counterexample, more on Marin soon.)

In particular we may have something like:

\begin{conj}\label{M-schemes-rigid} {\rm (Too Naive!!!,
completely false as shown by Marin, Fiedler)}.---Any $M$-scheme is
rigid.
\end{conj}

A priori if the devil of algebra does well his job this ought to
be completely false in high degrees $m\ge 8$, or say perhaps $m\ge
10^3=1000$.

\subsection{Marin's lock 1979: obstruction to rigid-isotopies}
\label{Marin:sec}

[23.01.13] Alas it seems that the above conjecture
(\ref{M-schemes-rigid}) is completely wrong for $m=7$ already,
compare Marin 1979 \cite[p.\,60--61]{Marin_1979}. Alas I was on a
bad day and could not completely understand his argument, which
looks however fairly simple involving the prose: ``La distinction
des deux courbes se fait en \'etudiant la position des ovales
ext\'erieurs par rapport aux droites joignant les ovales
impairs.'' (cf. p.\,60, of \loccit). At first this looks sloppy
justification as we are not just playing with projectivities but
with rigid-isotopies which (in marked contrast to their names) are
completely soft pathes in the residue of the discriminant. Yet I
am sure that Marin is right (as usual) but his argumentation is
for  highbrow readers?

\begin{figure}[h]
\hskip-1.2cm\penalty0 \epsfig{figure=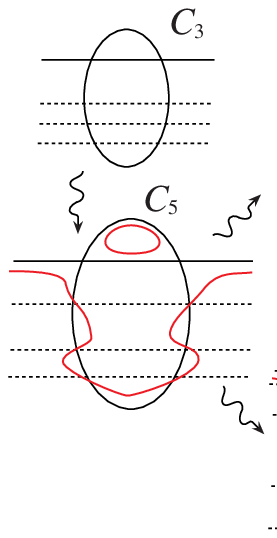,width=132mm}
\vskip-5pt\penalty0
  \caption{\label{Marin:fig}%
  Marin's counterexample to rigid-isotopy} \vskip-5pt\penalty0
\end{figure}

Ah yes the argument must be that when dragging the curve in the
parameter space (with the joystick) while choosing a triangle
through the deep ovals as on Marin's picture (reproduced as
Fig.\,\ref{Marin:fig}), then the
forced intersection with the pseudo-line gives already total
reality (or B\'ezout-saturation) of these lines with the septic
$C_7$ preventing the remaining (outer) ovals to traverse the
triangle during the motion (=rigid-isotopy).
Hence the distribution of ovals past the deep triangle is an
invariant of the rigid-isotopy class, i.e. any 2 curves exhibiting
distinct distributions of ovals past the fundamental triangle
cannot be deformed into the other.

So this is a bit like the moving frame of E.~Cartan, and  seems
indeed to corroborate Marin's clever observation!!! (Compare also
Degtyarev-Kharlamov 2000 \cite{Degtyarev-Kharlamov_2000} who call
this trick a ``lock'', while ascribing it as well to Fiedler.)
Note that even the complex orientations agree on Marin's example.
This method of the lock (or moving frame/traingle) affords
therefore an obstruction (\`a la B\'ezout) to rigid-isotopy. It
uses
the fact that a triangle in the
projective plane subdivides it in 4 pieces. One can wonder if
other (more complicated) locks are also useful. This method surely
deserves to be better explored and assimilated
(as remarked in Degtyarev-Kharlamov \loccit). For instance one can
wonder if it is enough to a lock with a pair of lines which
suffices to separate the outer ovals of the top figure. During the
isotopy we can keep track of them (at least the ovals where they
are passing through). Of course Marin's choice has the advantage
of canonicalness.  What is crucial is that the lock do not
degenerate during the isotopy, which is ensured by the fact that
the 3 inner ovals cannot become aligned without violating
B\'ezout. We have proved Marin's result:

\begin{theorem} {\rm (Marin 1979, or Fiedler)}
There is two isotopic $M$-septics, i.e. having the same real
scheme (and in fact the same complex orientations, but that
requires adding the arrows on Fig.\,\ref{Marin:fig}), yet which
are not rigid-isotopic (i.e., belong to distinct chambers past the
discriminant).
\end{theorem}

{\it Insertion} [08.04.13].---This raises of course the question
of counting the number of septics chambers corresponding to this
scheme (of degree 7). Perhaps variants of Marin's figure
(Fig.\,\ref{Marin:fig}) produce more than 2 chambers, but we are
not sure.

\subsection{Still some link between totality and rigidity?
Highbrow Nuij's principle}

[23.01.13] Is there still some link between total reality and
rigidity? A priori imagine the simplest situation of total reality
under a pencil of line
(ensured whenever we have a deep nest). Then very naively one
could imagine to contract progressively the curve to some normal
form like concentric circles and then drag it as a such toward the
other center of perspective and blow it up again along the radial
foliation toward the other curve. Of course doing so we meet
reducible curves hence the discriminant yet perturbing the path
there is some hope to avoid it completely, proving thereby the
Nuij theorem \cite{Nuij_1968}. Though  surrealist this argument is
the best we can give in favor of a connection between totality and
rigidity.

How does Nuij or Dubrovin prove their fantastic results? Can we
``do-it-yourself'' by making precise the above idea? One trick
would be to take a total pencil (vision from the innermost oval)
and perturb the nest toward concentric circles around the center
of perspective. Since both curves are (softly) isotopic there is
some chance that one path along the pencil (this being a circle
abstractly there is 2 such pathes) does not cross the discriminant
while affording a rectilinear rigid-isotopy. Then one finishes as
above.

\smallskip

$\bigstar$ {\it Long Insertion} [08.04.13].---In quintessence, the
deep nest for which rigidity holds true by Nuij 1968, is the
satellite of the conic (cf.
Sec.\,\ref{satellite-total-reality:sec}) whose rigidity can be
nearly ascribed to ancient Greeks (or Descartes, Newton,
Sylvester's law of inertia for quadratics forms, etc.). By analogy
the quadrifolium schemes of degree $m=4k$ consisting of 4 nests of
depth $k$ are total under a pencil of conics (as is fairly
trivial, and explicitly remarked in Rohlin 1978 at least for
$m=8$). This in turn is the satellite of the quartic quadrifolium
(degree $m=4$ with $r=4$ ovals) whose rigidity is known since
Klein 1876 (Sec.\,\ref{Klein-rigidity-of-quartics:sec}). Hence
this gives some evidence that totality in degree 2 implies
rigidity. Extrapolating further along a stability of rigidity
under satellites it could follow from the Rohlin-Le~Touz\'e
$(M-2)$-schemes of degree 6 (whose rigidity is ensured by
Nikulin's theorem deeper than Klein but sharing with it the r\^ole
of surfaces, viz. K3 quartics vs. cubics for Klein) that:

\begin{conj}\label{satellite-Rohlin-(6)-schemes-rigid:conj}
All satellites of Rohlin's sextic schemes $\frac{6}{1}2$ (and its
mirror $\frac{2}{1}6$) are rigid.
\end{conj}

Of course if rigidity is stable under satellites this would not
only applies to the Rohlin (totally real) schemes but to all
schemes of degree 6 which are of definite type as tabulated on the
Gudkov-Rohlin table (Fig.\,\ref{Gudkov-Table3:fig}). There are
precisely $64-2\cdot 8=48$ many such schemes. $\bigstar$ {\it End
insertion}.

\smallskip

Another method of proof (of Nuij's theorem) could be dynamical
like the one (we attempted) for CCC (\ref{CCC:conj}), yet
involving another functional a priori. Very loosely the functional
ought to have a unique attractor consisting of a series of
concentric circles (or perhaps ellipses). If the flow can be shown
to have this unique attractor (itself a certain manifold) then
every nested curve converges there and going forth and then back
we link rigidly our both curves.

Can we adapt the above synthetic method (or find an even more
synthetic method) to construct the (hypothetical) rigid-isotopy
between two (8)-schemes consisting of 4 nests of depth 2? The
naive canonical form would be the  same with circle or ellipses
yet its degree is twice too big (namely 16). ({\it Added}
[08.04.13].---Perhaps as suggested above, the canonical form is
the satellite of a quartic which has the correct degree.)

 We can first observe
that such an octic curve has, like the deep nest, all the nested
ovals forming negative pairs (this can be seen either \`a la
Fiedler on a model or \`a la Ahlfors via the total pencil which
forces the orientation to gyrate in the same sense as swept by the
pencil). As a such there is no topological obstruction to shrink
them at the microscopic scale (or apply alternatively a variant of
CCC). Once contracted at the microscopic scale our configuration
moves without resistance in the free vacuum and then may be
re-expanded at the next curve. This is very sloppy heuristic of
course,  sembling much like inter-sidereal travelling, but there
may be some truth in this. At least one sees a connection with
CCC. If this works we get a proof of the:

\begin{conj}\label{rigidity-sat-quadrifolium:conj}
The $8$-scheme $4\times\frac{1}{1}$ of $4$ nests of depth $2$ is
rigid (i.e. any two of its representatives are rigid-isotopic).
More generally the $4k$-scheme consisting of 4 nests of depth $k$
is rigid.
\end{conj}

This could be the ``degree 2'' avatar of Nuij's theorem, and ought
to be proved by an iterated variant of CCC (like C++). As yet the
total reality of the scheme was only involved to ensure that all
the ovals belonging to some nest gyrate in the same sense
according to complex orientation (what Rohlin 1978 calls negative
pairs). Those negativity may be interpreted as some depressiveness
permitting precisely the collapse to the microscopic scale where
then we can travel without friction within the ``ether''.

Given any collection of deep nests with negative pairs of ovals
(gyrating in the same sense), we may hope to contract them at the
nanoscale (via CCC or C++) and then travel and re-expand to reach
any other curve with the same complex orientation. This looks
topologically plausible yet the drawback of ignoring the total
pencil (`\`a la Ahlfors-Rohlin) is that our assumption does not
only pertain on the real scheme but also upon its complex
characteristics. So the assertion gain in generality but loose
some elegance.

As yet we have merely considered schemes which are towers (i.e.
without branching in their nesting graph). However Rohlin (1978)
claims that the $(M-2)$-schemes of degree 6 $\frac{6}{1}2$ and its
mirror $\frac{2}{1}6$ are total of order 3. Can we deduce that
those schemes are rigid by a method independent of Nikulin 1979
\cite{Nikulin_1979/80}, and analog to the one sketched above? One
should be in position to visualize the cubic pencil so as to draw
the complex orientation. Bypassing this difficult task we may
appeal to Rohlin's formula to deduce the complex orientation.
Assume the scheme to be $\frac{2}{1}6$, hence $2(\Pi^+
-\Pi^-)=r-k^2=9-9=0$ tell us that one pair is positive and the
other negative. So both inner ovals gyrate ``differently'' (in
accordance with Fig.\,\ref{GudHilb8:fig}).

At this stage one is quite puzzled, i.e. one does not see how to
bring the curve at the nanoscale via contractions. Of course
granting CCC we can shrink the empty ovals, but a priori cannot
shrink the nonempty one. So our heuristic method looks here quite
impuissant! Another quasi-paradox of our heuristic method arises
for a curve of type~I belonging to the scheme $9$, which under a
total pencil of cubics (easy to visualize, cf.
Fig.\,\ref{Fcubic:fig}) could be isotoped to the other curve(s) of
this scheme of type~II. Yet of course the second curve lacks a
total pencil to re-expand.

All this is just supposed to illustrate that we see no direct
relation between total reality and rigidity, at least via the
naive contraction approach. However this does not preclude a
deeper relationship. That would maybe involve exploiting  more the
total pencil as a tool to construct a first reduction to some
normal shape, which ought to be then easily tele-transported and
then re-expanded via the second total pencil. So the total pencil
should act as some sort of
(contracting) wormhole or as a railway guiding the curve to some
canonical shape easier to tele-transport
(at the speed of light). Alas this is much too vague to convince
us about any implication like ``total $\Rightarrow$ rigid''.

[26.01.13] After Shustin's e-mail, who remembered me the reason
why the empty scheme is rigid, one can suspect two basic scenarios
ensuring rigidity. Taking the deep nest as prototype, where
rigidity holds true by Nuij 1968 \cite{Nuij_1968} one could
suspect that rigidity is causal either of total reality \`a la
Ahlfors or by the proximity to the empty chamber (which is
connected, and actually baricentrically ``convex''). Both
phenomena could explain Nuij's rigidity of the deep nest while
affording basic intuition about guessing further rigidity results.
For instance total reality could explain the rigidity of schemes
swept out by pencil of conics (e.g. $4\times \frac{1}{1}$ in
degree 8), while the proximity to the empty locus could via CCC
prompt rigidity of the scheme with one oval (in even degree at
least).

Another naive idea I had (but which is now quite outdated) is that
while total reality could imply rigidity via Ahlfors, the avatar
of the latter for empty curves (namely  Witt 1934
\cite{Witt_1934}) could be involved in the rigidity of the empty
scheme.

\section{More total reality}


\subsection{Another attempt to prove Rohlin's total reality claim}

{\it Editorial note} [08.04.14].---The prose of this section
starts a bit abruptly, due to a permutation of section. Prior this
material came right after Le~Touz\'e's section
(Sec.\,\ref{LeTouze:sec}), explaining why the chromatic law for
conics impedes a direct corruption of Nikulin by Fiedler-Marin's
trick. This miracle of extra intersection created by dichromatism
gave me some hope to attack Rohlin's highbrow claim, but the
difficult turned out to be immense to fill.

\smallskip
[13.02.13, 12h42] Further the impact of this method of extra
intersections
gained  by dichromatism must probably also be the key behind
Rohlin's proof of the universal orthosymmetry of the sextic
schemes $\frac{6}{1}2$ and $\frac{2}{1}6$.
We call any curve having one of these schemes a {\it Rohlin curve}
as the latter Academician in his 1978 article \cite{Rohlin_1978}
was the first (and actually the unique
creature in the universe except for possible extraterrestrial
intelligences) to state the universal orthosymmetry of such
curves.

Through the 8 deep ovals (equivalently the empty ones) of such a
Rohlin curve we let pass a pencil of cubics. As above
(\ref{LeTouze-Gabard-Hilfssatz:lem}) we imagine the inner
basepoints black colored while the outer basepoints are white
colored.

Let $C_3$ be any cubic of the pencil, which we assume smooth for
simplicity. If $C_3(\RR)$ is connected, then $C_3$ intersects
$8\cdot 2=16$ plus twice the nonempty oval of $C_6$, and so the
intersection is totally real. If $C_3$ is not connected then its
splits an oval and a pseudoline. A priori the oval could visit the
6 inner points while the pseudoline the 2 outer ovals. In this
case there is no forced extra intersections.

\begin{figure}[h]
\centering
\epsfig{figure=LeTouzeRohlin.eps,width=122mm} \vskip-5pt\penalty0
  \caption{\label{LeTouzeRohlin:fig}%
  Attempts to prove Rohlin's claim} \vskip-5pt\penalty0
\end{figure}

By Le~Touz\'e's lemma~\ref{LeTouze:lem} we may infer that the
triangle through any 3 of the 6 inner ovals (assuming
$\frac{6}{1}2$) does not separate the 2 outer ovals
(Fig.\ref{LeTouzeRohlin:fig}\,c). There are $\binom{6}{3}=20$ such
triangles. But this does not seem very useful.

Let us start again. Consider a $C_6$ with scheme $\frac{6}{1}2$, a
so-called {\it left-wing Rohlin curve} (in view of its position on
the Gudkov table, Fig.\,\ref{Gudkov-Table3:fig}). Consider the
pencil of cubics through 8 deep basepoints selected inside the
$6+2=8$ empty ovals of the $C_6$. We claim (with Rohlin 1978
\cite{Rohlin_1978}) that this pencil is totally real (and
consequently the curve is of type~I).

Total reality of $C_3\cap C_6$ is clear when the cubic is
connected for then there are $2\cdot 6+2=18$ intersections, the
last two  being created while crossing the nonempty oval.

So let us assume the cubic disconnected. Two cases are possible
either it is smooth or not. If singular then the cubic may either
have a solitary point or even be a conic union a disjoint line.
However since we are free to perturb the 8 basepoints the pencil
can probably be assumed to be transverse to the discriminant so
that its singular members are uninodal curves. This rules out the
second case. In the first case of a solitary cubic then the
solitary node cannot be one of the 8 basepoints, and so the
connected pseudoline of this singular cubic visits all 8 points,
 creating thereby 2 additional intersections  with the nonempty oval.

So may assume the cubic smooth and as soon as its pseudoline
visits both inner and outer points we are finished (2 bonus
intersections are created). On the other hand if the pseudoline
visits only inner points then it must evade out of the nonempty
oval (otherwise it would be null-homotopic) and so we score again
2 extra points, and gain total reality. Hence we may assume:

\smallskip
{\bf 1st reduction:} {\it The pseudoline of our cubic $C_3$ visits
the $2$ outer points whereas its oval visits the $6$ inner points
while being traced inside the nonempty oval of the $C_6$ (cf. {\rm
Fig.\,\ref{LeTouzeRohlin:fig}d}).}
\smallskip

We have to show that this is contradictory, but are still far from
the goal. In fact at the time of writing these lines the writer
does not know if he will ever be able to complete this argument.

A first remark is that our pencil of cubic has another
(non-assigned basepoint). Where is it? We think that it must be on
the pseudoline of $C_3$ for simple vibratory reasons. Indeed look
at Fig.\,\ref{LeTouzeRohlin:fig}d and imagine a nearby cubics
$C_\epsilon$ in the real locus of the pencil. The corresponding
oval will have to oscillate about that of $C_3$ and since 6
basepoints are on the oval the oscillation closes up perfectly.
(In savant terms
the oval has a trivial tubular bundle.)

Another thing
natural to do is to cut out the inside
of the oval $O$ of the cubic $C_3$ of Fig.d in order to apply the
Poincar\'e index formula to the  foliation ${\cal F}$ induced by
the pencil. Then the situation is a bit messy but as follows. The
pencil hits the discriminant of cubics 12 times over the complexes
as $\deg \disc_3=3(m-1)^2=12$ for $m=3$. A priori not all
intersections are real, but can occur in conjugate pairs. Those
singular curves which are real are either solitary cubics or have
two real branches crossing transversally (``real bitangent'' and
we call them {\it nodal cubics}). Denote their respective number
$\sigma$ and $\beta$ (where $\sigma$ stands for ``solitary'' and
$\beta$ for ``bitangent''). We have $\sigma+\beta=12-2k$.

Using the Poincar\'e local index formula $j=1+\frac{I-E}{2}$ where
$I$, $E$ are the number of internal resp. external tangencies of
the foliation with a small circle surrounding the singularity, it
is a simple matter to compute indices. Of course a basepoint gives
a foyer of index $+1$, a solitary cubic gives a ``centre'' of
index +1 (since $I=E=2$), while a nodal cubic gives a hyperbolic
saddle of index $-1$ (as $I=0, E=4$). Applying Poincar\'e's index
formula (cf. Poincar\'e 1885 or Gabard 2011
\cite{Gabard_2011-Euler-Poincare-obst-pretzel-long-tentacles},
arXiv, ``long tentacles'') it follows
$$
9+\sigma - \beta=\chi (\RR P^2)=1.
$$
It may be deduced that $\beta \ge 8$, that is the following:

\begin{lemma}\label{nodal-cubics-8-many-in-a-pencil:lem}
Any generic pencil of cubics contains at least $8$ nodal
cubics.
\end{lemma}

Further we easily tabulate the possible value as $(\beta,
\sigma)=(8,0),(9,1), (10,2)$.

It would be however probably more interesting to apply the
Poincar\'e formula in the inside of the oval $O$ (of the cubic
$C_3$) doubled to get a sphere. (This doubling is merely a trick
yet useful to eliminate the boundary.) The difficulty in doing so
is that we do not really know a priori how the singularities of
the foliation $\cal F$ induced by the pencil are distributed
inside the oval $O$. So let us denote with subscript
``naught=0''
the corresponding quantity of singularities inside the oval $O$.
Then we have
$$
6+2 \sigma_0- 2 \beta_0=\chi(S^2)=2,
$$
where we used implicitly the fact that the 9th basepoint is not
inside our oval (nor on its periphery since it is rather located
on the pseudoline).

Alas at this stage the situation looks confuse. One idea is to
imagine a solitary node inside the oval $O$. Then there is a
unique time direction so that this oval inflates while moving
inside the pencil. Since the oval $O$ has a tube neighborhood like
Fig.\,\ref{LeTouzeRohlin:fig}e this oval cannot hit the oval $O$,
and must rather collide with the pseudoline component to form a
nodal singularity of type $\beta$. So to each solitary node is
canonically assigned a non-solitary node (all this occurring
inside $O$). It seems evident that the corresponding map is
injective, and so $\sigma_0\le \beta_0$. Alas this gives no
contradiction when injected in Poincar\'e's relation displayed
above.

\begin{Scholium}
\label{scholie:Rohlin-does-not-boils-to-Poincare} {\rm
[14.02.13]}---In fact it seems unlikely that there is a proof of
Rohlin's total reality assertion (for $\frac{6}{1}2$ and its
mirror) using Poincar\'e's formula only.
\end{Scholium}

We arrived at this conclusion after tracing a rather complicated
foliation of the plane $\RR P^2$ containing the bad cubic $C_3$ as
a leaf. (More about the discussion of the relevant pictures soon.)
The bad cubic is one which is not totally real, hence whose  oval
is necessarily enclosed in the nonempty oval of the $C_6$.

In reality such a foliation is fairly easy to construct just by
starting with the bad cubic and then merging its components
together to a nodal cubic and pursuing the depiction in a more or
less canonical fashion. At each step B\'ezout for $C_3\cap C_6$ is
respected and Poincar\'e index formula is of course verified.

What should we deduce? Could it be
that Rohlin's total reality assertion is false, while its theorem
on the type~I of his schemes is right as follows from some
highbrow topological congruence (due to himself, Kharlamov and
Marin, cf. (\ref{Kharlamov-Marin-cong:thm}))? If the latter
super-classical congruence is correct then via Ahlfors theorem the
total reality assertion is likely to hold true yet perhaps not for
{\it all\/} pencil of cubics through the 8 empty ovals (or a
priori for pencils involving curves of higher order). Perhaps the
proof should involve Abel's theorem applied on the cubics, yet
they vary so seems unlikely. Perhaps Abel has to be used on the
$C_6$?

All this is  puzzling and we frankly confess our poor
understanding which calls for a synthesis between the abstract
viewpoint of Riemann-Klein-Ahlfors and the embedded viewpoint of
Harnack-Hilbert-Gudkov-Rohlin.
This subdivision of our science is still vivid today, compare e.g.
in Russia the tradition along Natanzon vs. Kharlamov-Viro.

As we noticed earlier in this text (very optionally see
Sec.\,\ref{sec:Total-reality-Harnack-max-case}), it is also
tantalizing to trace a total pencil on the $M$-sextics, already on
those of Harnack and Hilbert. It would be of interest to know what
is the degree of curves forming a total pencil in the $M$-case,
whose existence follows in principle (modulo Marin's private
communication objection) from Ahlfors even in the simple
schlichtartig variant of Bieberbach-Grunsky. ({\it Added}
[08.04.13].---This question should  now be settled via
(\ref{total-reality-of-plane-M-curves:thm}).)

Since an $M$-sextic has 10 empty ovals it looks quite improbable
that total reality is exhibited by a pencil of cubics which has
only 9 foyer-type singularities. This must, we believe, easily
follow from Poincar\'e's index formula. Indeed each empty oval
must contain a singularity of positive index while thorn
singularities of index $+1/2$ are precluded for an algebraic
pencil.

Now let us discuss our picture leading to the announced
Scholium~\ref{scholie:Rohlin-does-not-boils-to-Poincare}. The game
is to foliate the plane by ``flexible'' cubics in the sense that
we depict only the singular nodal curves abstractly like a figure
``8''. A cubic cannot have a node plus an oval. Using this we
sometimes have the impression that there is an obstruction to
complete the foliated structure like on
Fig.\,\ref{LeTouzeRohlin:fig}. (Could it be the case that Rohlin
made such a mistake, in the sense of a too hasty inference?)

\begin{figure}[h]
\centering
\epsfig{figure=LeTouzeRohlin2.eps,width=122mm} \vskip-5pt\penalty0
  \caption{\label{LeTouzeRohlin2:fig}%
  Nearly disproving Rohlin's claim, or rather showing that
  its proof does not reduce to the analysis situs of foliations
  \`a la Poincar\'e.} \vskip-5pt\penalty0
\end{figure}

However deleting curves and starting again one finds
Fig.\,\ref{LeTouzeRohlin2:fig} which is topologically admissible.
On it each cubic looks like a cubic, B\'ezout is respected as well
as Poincar\'e's index formula (as it should). In fact our solution
shows no solitary cubics. Of course we do not claim that this
free-hand drawing foliation is algebraic (in which case Rohlin's
claim would be erroneous), but we are also not able to exclude
this eventuality. What this picture really shows is that we cannot
expect to prove Rohlin's assertion via the sole apparatus of the
combinatorial topology of foliations (\`a la Poincar\'e). So if
true Rohlin's statement has some deeper geometric significance,
and it is quite tantalizing to imagine its complexity. Further by
Ahlfors theorem it is likely that  Rohlin's claim is just the top
of the iceberg of  a
plethora
of another phenomena of total reality in higher degrees which must
all be very delightful to visualize if not quickly overburdening
any human intelligence. Again our prophecy (cf. Introd. of this
text) is that this is linked to the stability of matter at the
nano-scale, or at least that such totally real pencils describe
the dynamics of electrons about an atomic nucleus, as  ellipses
described the trajectory of Mars about the Sun in  Kepler's days
(ca. 1605).

\subsection{Total reality from the elementary viewpoint}

[15.03.13] Broadly speaking our main Leitmotiv is the question of
examining if there is any  relation between total reality \`a la
Ahlfors (or rather
Riemann-Schottky-Klein-Bieberbach-Grunsky-Ahlfors, etc.) and
Hilbert's 16th problem. Rohlin's claim about his sextic curves
posits such a deep relation. In this section we try to explore
more systematically this relation.

Before entering into the details let us pose some of the guiding
questions. Given an (abstract) dividing curve $C$ there is
according to Ahlfors 1950 \cite{Ahlfors_1950}, a totally real map
$f\colon C\to \PP^1$ of the curve to the projective line $\PP^1$
(i.e. $f^{-1}(\PP^1(\RR))=C(\RR)$).

When this curve is plane does this map extends to the projective
plane $\PP^2$ as to be induced by a pencil of curves? The question
is actually pure geometry primarily meaningful over the complexes.
Given a plane curve $C_m$ defined over $\CC$, and a holomorphic
map $C_m(\CC) \to \PP^1(\CC)$ is it true that there is a pencil of
curves so that the map induced by the pencil is the given
holomorphic map.

If this is true then it is certainly true equivariantly and
Ahlfors theorem implies that {\it any plane dividing curve admits
a total pencil of curves.}\footnote{[08.04.13] Compare with
Le~Touz\'e's article 2013
\cite{Fiedler-Le-Touzé_2013-Totally-real-pencils-Cubics}, where it
is asserted that this was  implicitly conjectured by Rohlin, in
1978. Recall that Le~Touz\'e's husband Fiedler is a direct student
of Rohlin, and so this may also be based upon some oral tradition,
in case Rohlin was too cautious to put crazy ideas on the paper.}

The next question is then to trace such pencils, and try to
control its order (i.e. the degree of the curve constituting it).
Tracing them in case of $M$-curves looks an especially hard
exercise. ([08.04.13] Okay but see
(\ref{total-reality-of-plane-M-curves:thm}).)

Without appealing to Ahlfors theorem one can also study totally
real pencils {\it per se\/}, as a tool to detect the dividing
character of
curves. Actually it is this trivial criterion (as applied to the
G\"urtelkurve, or also hyperelliptic curves) which lead the writer
to discover Ahlfors' theorem in ca. 2000--01 independently and
prior of knowing about Ahlfors' work. So whenever a curve $C_m$ is
swept out by a total pencil of curves it is dividing.

So one can examine which sort of curves are  exposed to such a
total pencil to derive in principle an infinite series of
orthosymmetry criterions. The prototype is the case of the deep
nest totally swept out by a pencil of lines through the deepest
oval. Then one would like to study pencil of conics, cubics,
quartics, etc.

Some extra difficulty arises from the distinction between
universal total reality where it is forced by the sole knowledge
of the real scheme as in the case of deep nests or Rohlin's
sextics $\frac{6}{1}2$ and it mirror $\frac{2}{1}6$, and
``versal'' total reality where the detailed geometry of the curve
is required to exhibit total reality.
This is of course much allied to what Rohlin calls schemes of
indefinite type.

For instance the octic scheme $4\times \frac{1}{1}$ consisting of
4 nest of depth 2 is universally totally real, under the pencil of
conics through 4 basepoints inside the empty ovals. This example
easily extends to schemes of degree $4k$ having 4 nests of depth
$k$. Existence of such curves for each $k$  is demonstrated by
Fig.\,\ref{Total:fig}. ([14.03.13] A somewhat more conceptual
reason is given by taking the algebraic satellites, i.e. nearby
levels of the quadrifolium quartics, so $P_4\cup
P_4+\epsilon_2\cup\dots \cup P_4+\epsilon_k$ and smoothing this
union of  reducible curve.)

\begin{figure}[h]
\centering
\epsfig{figure=Total.eps,width=122mm} \vskip-5pt\penalty0
  \caption{\label{Total:fig}%
  Curves $C_{4k}$ with 4 nests of depth $k$,  totally
  real under a pencil of conics assigned to pass
  through the deepest ovals} \vskip-5pt\penalty0
\end{figure}

{\it Insertion} [08.04.13].---It may be observed  (Fiedler's
smoothing law) that those curves are of type~I for surgical
reasons, providing another proof independent of total reality.
Further the total pencil induces the complex orientations due to
the holomorphic character of the underlying total map. It suffices
then to imagine the pencil of conics to see that the intersection
series will move along Fiedler's arrows (this we shall vaguely
refer to as dextrogyration). Understanding this properly in
general could be the source of some progresses in the field,
maybe?

\def\quadri{quadrifolium }

\begin{lemma}
Any  curve $C_{4k}$ of degree $4k$ whose real scheme consists of
$4$ nests of depth $k$ (\quadri for short) is total under the
pencil of conics through the $4$ empty ovals.
\end{lemma}

\begin{proof}
Let $C_2$ be any conic of the pencil. First, $C_2(\RR)$ is
connected, e.g. because it is a rational curve, aka as {\it
unicursal\/} in Cayley's jargon, cf. optionally the Introd. of
Harnack 1876 \cite{Harnack_1876}, which is of course not really
required on the case at hand since it suffices like in Antiquity
to sweep out the conic by lines from one of its point. It follows
that $C_2(\RR)$  has to cut our curve $C_{4k}$ in $4\cdot 2k$ real
points for topological reasons. But this is the maximum
permissible according to B\'ezout, hence the pencil is totally
real.
\end{proof}

It seems likely that the converse statement holds true, namely any
curve $C_{4k}$ totally real under a pencil of conics with 4 real
basepoints has this scheme of 4   nests of depth $k$. Note at
least that the depth of the ovals cannot be distributed otherwise
without violating B\'ezout for lines. For instance if a $C_{12}$
instead of having 4 nests of depth 3, had nests of depths say
$2,4,3,3$ then the line through the nests of depth $4$ and $3$
would have too much intersection (namely $8+6=14>12$).

It is worth noticing that the degree of such total maps are in
accordance with the bound $r+p$ announced in Gabard 2006
\cite{Gabard_2006}. Indeed if $C_{4k}$ is a quadrifolium, then
$r=4k$. Hence by the obvious Klein relation $g=(r-1)+2p$ and the
genus formula $g=\frac{(m-1)(m-2)}{2}$, where $m$ is the degree,
we find
$$
g=\textstyle\frac{(4k-1)(4k-2)}{2}=(4k-1)(2k-1), \quad \textrm{
and }
$$
$$
p=\textstyle\frac{g-(r-1)}{2}=\textstyle\frac{(4k-1)(2k-1)-(4k-1)}{2}=
\textstyle\frac{(4k-1)(2k-2)}{2}=(4k-1)(k-1).
$$
On the other hand by letting degenerate the 4 basepoints against
the deep oval we find a total morphism of degree $2\cdot
4k-4=8k-4=4(2k-1)$. This has to be compared with the $r+p$ bound
$$
r+p=4k+(4k-1)(k-1),
$$
which is indeed much greater as shown e.g. by evaluating for
$k=1,2,\dots$. We find for $k=1$, $4(2k-1)=4\le r+p=4$. For $k=2$,
$4(2k-1)=12\le r+p=8+7=15$, and so on due to quadratic growth of
$r+p$. In fact the gonality of such \quadri curves is probably
$\gamma=4(2k-1)$ (at least majored by this quantity) and so
significantly lower than the universal upper bound $r+p$ stated in
Gabard 2006 \cite{Gabard_2006}.

In fact it is worth testing the truth of this $r+p$ bound on
sextics already. Then we shall basically pass into review all the
dividing curves of the Gudkov-Rohlin table
Fig.\,\ref{Gudkov-Table3:fig}. It would be natural to start from
the top of this table but as $M$-curves are paradoxically tricky
to understand from the viewpoint of total reality, we start from
the bottom. The paradox is that the total reality phenomenon for
abstract $M$-curves is basically the schlichtartig ($p=0$) case of
Ahlfors theorem which is pretty much easier than the positive
genus case. ({\it Added} [08.04.13].---This paradox is settled via
(\ref{total-reality-of-plane-M-curves:thm}).)

So starting from the bottom we have first the deep nest $(1,1,1)$.
Then $r=3$ and $p=[g-(r-1)]/2=[10-2]/2 = 4$. The gonality
$\gamma=5 \le r+p=7$ is exhibited by the pencil of lines through a
point on the deepest oval.

For the scheme $\frac{4}{1}$ (when of type~I) total reality comes
from the conics pencil through the deep nest. It leads to a total
series of degree $\gamma \le 2\cdot 6-4=8 \le r+p=5+3=8$ and
Gabard's bound is sharply realized.
(Trick: while it is sometimes boring to compute $r+p$ it may be
remembered that this is also $\frac{r+(g+1)}{2}$, i.e. the mean
between $r$ and Harnack's bound $g+1$.)

\def\tetris{tetris }

In fact checking total reality involves here the exact geometry
and not merely knowledge of the real scheme. More specifically the
2 prototypes of curves of type $\frac{4}{1}$ are depicted on
Fig.\,\ref{R4-1:fig}. One would like to have a geometric criterion
for deciding a priori the type.

In both cases we have 4 deep (=empty) ovals. Now we may choose 4
points in them. Given such a tetrad there are 2 cases to be
distinguished (cf. Fig.\ref{Total2:fig}\,a). Crudely put, either
one of the 4 points can be inside the triangle spanned by the 3
other or not. However projectively this is a misconception as
shown by Fig.\,b. Yet as we are given a curve of type
$\frac{4}{1}$ it may look either like one of the 2 versions of
Fig.\,c. Here the deep triangles (those traced through 3 deep
ovals) have always a distinguished 2-simplex traced inside the
nonempty oval. We call any such a (fundamental) simplex and there
are 4 of them. Now two cases are possible: either one fundamental
simplex is the union of the 3 others or not. Alternatively one
oval is contained inside a fundamental simplex or not. (Check that
this is well-defined requires keeping B\'ezout in the background
memory.) Hence Fig.\,a recovers some intrinsic significance (and
amounts essentially to the 2 possible visions we may have of a
3D-tetrahedron when projected on our 2D-retina). We call the
second option a \tetris (as a short cut for stable tetrad like a
prism stably posed on the sheet of paper, in contrast to the
unstable tetrad posed on its edge hence in unstable equilibrium).

\begin{figure}[h]
\centering
\epsfig{figure=Total2.eps,width=122mm} \vskip-5pt\penalty0
  \caption{\label{Total2:fig}%
  Tetrad} \vskip-5pt\penalty0
\end{figure}

This being said we have the following recognition lemma of
Klein-Rohlin's type by pure geometry:

\begin{lemma}
If the $4$ empty ovals of a sextic curve of type $\frac{4}{1}$
form a \tetris then the pencil of conics through the deep (empty)
ovals is total, and in particular the curve is dividing (type~I).
\end{lemma}

\begin{proof}
Let us choose any conic $C_2$ of the pencil (through the 4 deep
points inside the empty ovals). Like on Fig.\,e this curve will
appear inside the largest fundamental simplex. Then the idea is to
surger the conic into 2 pseudolines $C_2=J_1+J_2$ as shown on
Fig.\,f or g. There is several way to do this but choose one. This
surgery amounts to aggregate a certain edge of the tetrahedron
which is exempt of intersection with the nonempty oval $N$
(because it is already B\'ezout-saturated).
Therefore $C_2\cap N=(J_1\cup J_2)\cap N=(J_1\cap N)\cup ( J_2\cap
N)$, but as each $J_i$ is a pseudoline
each
must intersect twice the nonempty oval $N$, and we gain 4 extra
intersections. On the other hand $C_2$ intersect twice each empty
oval and so we totalize $2\cdot 4+ 4=12$ intersections the maximum
permitted  by B\'ezout. Total reality follows, and the proof is
complete.
\end{proof}

Albeit not perfectly hygienical our proof shows
how to
%
gain extra intersections by this splitting method. (Ideally we
could hope that this method is also the key to Rohlin's total
reality claim for the curve $\frac{6}{1}2$, but this is not clear
a priori. Imagine the bad cubic whose oval is contained inside the
nonempty oval of the $C_6$, and what to do next!!!?)

Next we would like a similar optical recognition criterion of the
type for the sextic scheme $\frac{2}{1}2$. Here looking at
Fig.\,\ref{R2-12:fig} reproduced below as
Fig.\,\ref{Total2-12:fig}a,b below (which I borrowed from
Degtyarev-Kharlamov's survey 2000 \cite{Degtyarev-Kharlamov_2000})
suggests that what distinguishes both types is whether the line
through the 2 outer ovals separates or not the  two inner ovals
within the inside of the nonempty oval.

\begin{figure}[h]
\centering
\epsfig{figure=Total2-12.eps,width=122mm} \vskip-5pt\penalty0
  \caption{\label{Total2-12:fig}%
  Prismatic recognition of the type for the scheme $\frac{2}{1}2$} \vskip-5pt\penalty0
\end{figure}

Somewhat more formally, let us extract from both prototypical
curves (Fig.\,a, b resp.) some combinatorial datum. We mark in
black inner ovals by choosing a point inside, and choose also 2
white points in the outer ovals. Unlike in the previously studied
case, there is no preferred fundamental simplex inside the
nonempty oval $N$, but we have a $1$-simplex entirely traced
inside the nonempty oval $N$, which we mark by a double stroke. So
we extract Fig.\,c resp.\,d and what distinguishes both is the
issue that the line through the white vertices intercepts the line
through the black vertices along its double marking or not. Let us
call the first case (like Fig.\,c) a {\it crucifix\/} and then we
have the:

\begin{lemma}
If a sextic curve $C_6$ of real scheme $\frac{2}{1}2$ has a
crucifix, then the pencil of conics through the empty ovals is
totally real (and the curve is of type~I).
\end{lemma}

\begin{proof}
We have defacto $4\times 2=8$ real intersections coming from the
empty ovals. Take any conic of the pencil. Two cases may appear.
Either the conic is {\it dichromatic\/}, that is when we follow it
along some orientation we visit the 4 basepoints in the sequence
black-white-black-white (BWBW) in this alternating way, or it can
be monochromatic if this sequence reads BBWW (compare Fig.\,e). In
the dichromatic case 4 intersections are created, and total
reality is ensured. In the monochromatic case, we apply the
splitting method which decomposes the conic $C_2$ as an union of
two pseudolines $J_1,J_2$ traced on Fig.\,f  obtained by cutting
the conic at the two black points and adding the fundamental
$1$-simplex linking both black vertices in the inside of $N$ (the
nonempty oval). Since this $1$-simplex does not cut $N$, the
intersection with $N$ remains the same after this surgery, but
each pseudoline forces 2 intersections with $N$, and we gain the 4
required extra intersections. Total reality of the whole pencil is
proved.
\end{proof}

{\it Insertion} [09.04.13].---On a second reading of this proof,
it is not clear what prevents to apply the same argument to the
other configuration. Try to clarify this at the occasion.

To summarize the proof is the same as the previous one safe for an
intervention of chromatism (black and white reflecting the inner
and outer ovals). Note also that as $r$ is the same as in the
previous case, Gabard's bound $r+p$ is likewise verified (at least
not
quashed=invalidated)!

$\bullet$ Can we continue this game? According to Gudkov-Rohlin's
table (Fig.\,\ref{Gudkov-Table3:fig}) the next specimen to study
is $\frac{5}{1}1$. As we have now $6$ empty ovals it is evident
that a pencil of conics will not exhibit total reality. We
probably have to move to cubics with 8 basepoints assignable.

Let us use the same naive device of combinatorial extraction from
two prototypes (cf. Fig.\,\ref{Total5-11:fig}). We get some
beautiful bi-pyramid (or octahedron) and one should imagine each
face (or interface) shaded whenever the corresponding triangle is
fundamental (i.e. included in the outer oval of the sextic $C_6$).
Alas both configurations so obtained look combinatorially
equivalent, and we feel puzzled. The next idea that comes to mind
is to look at the conic through the deep black points rooted in
the inner ovals. It seems that what distinguishes both types (I
vs. II) is the location of the outer points as being resp. inside
or outside this conic. Another feature distinguishing both models
is the absence resp. presence of a line through  the outer oval
missing the nonempty oval $N$.

\begin{figure}[h]
\centering
\epsfig{figure=Total5-11.eps,width=122mm} \vskip-5pt\penalty0
  \caption{\label{Total5-11:fig}%
  Optical recognition for $\frac{5}{1}1$} \vskip-5pt\penalty0
\end{figure}

In fact even on the model it is quite difficult to guess which
pencil of cubics will exhibit total reality of the type~I
configuration (as predicted by Ahlfors' theorem). One could take
the horizontal line through the white point which cut the sextic 6
times, and take two extra basepoints on this line (perhaps in the
2 lower wings of the butterfly). Then at least the split cubic
consisting of the conic through the 5 black and 3 white points
would be totally real. Another puzzling point is that such a
pencil will have mapping-degree $3\cdot 6- 8=10$ when the 8
basepoints degenerate on the curve $C_6$, whereas Gabard predicts
one of degree $r+p=(r+g+1)/2$ the mean value of $r$ and Harnack's
bound that is $(7+11)/2=9$. So perhaps this constitutes a
(potential) counterexample to Gabard 2006 \cite{Gabard_2006}, at
least if all abstract pencils are concrete and realized by cubics
pencils. If we look at quartics pencil with
$\binom{4+2}{2}-1-1=13$ free basepoints then the degree would be
$4\cdot 6- 13=11$, still higher than Gabard's bound. For quintics
there are $\binom{5+2}{2}-1-1=19$ free basepoints and so the
degree is $5\cdot 6-19=11$, for sextics $\binom{6+2}{2}-1-1=26$,
so the degree is $6\cdot 6-26=10$, for septics
$\binom{7+2}{2}-1-1=34$, so the  degree is $7\cdot 6-34=8$. Gabard
seems rescued, yet it looks quite tantalizing to understand the
geometry of such a total pencil if it exists. (If Gabard's bound
is true and the Riemann-Hilbert transition from the abstract to
the concrete viewpoints equally holds true then such a total
pencil should exist of order at least seven!)

{\it Insertion} [09.04.13] Let us summarize this as follows:

\begin{Scholium}\label{(M-4)-sextics-corrupt-Gabard:scholium}
$(M-4)$-sextics of type~I are perhaps a good place where to
corrupt Gabard's bound $r+p$. And if not it is  at least a pi\`ece
de r\'esistance against the principle that any abstract pencil is
concrete, and therefore Ahlfors abstract theorem is unlikely to
apply without friction in Hilbert's 16th problem. In other words
Riemann's canary feels claustrophobic in the Plato cavern of
Brill-Noether-Hilbert. Perhaps the above example merely corrupts
the conception that the mapping-degree  of a total pencil is
minimized when the order of its constituting curves is. However it
could still be true that any dividing plane curve of degree $m$
has its total reality exhibited by a pencil of order $(m-2)$ (or
less), compare e.g. {\rm
(\ref{total-reality-of-plane-M-curves:thm})} for the case of
$M$-curve.
\end{Scholium}

[16.02.13] Let us leave aside this problematic concerning the
truth of Gabard's bound $r+p$ to concentrate on the existence on a
cubics pencil which is total on our sextic of symbol
$\frac{5}{1}1$. Of course the existence of the latter is more an
act of faith than a truth a priori, as it is not obviously implied
by Ahlfors' theorem. The latter probably gives the existence of a
total pencil and one may wonder what is the least possible order
of the curves in the pencil.

$\bullet$ Then we can look at the next curve $\frac{3}{1}3$ of the
Gudkov tabulation (Fig.\,\ref{Gudkov-Table3:fig}). Two models are
depicted on Fig.\,\ref{Total3-13:fig}c, and one may hope to
distinguish them by some combinatorial recipe (perhaps by looking
at the inner fundamental
 simplex and some outer simplex).

\begin{figure}[h]
\centering
\epsfig{figure=Total3-13.eps,width=122mm} \vskip-5pt\penalty0
  \caption{\label{Total3-13:fig}%
  Failing to understand total reality of the scheme $\frac{3}{1}3$} \vskip-5pt\penalty0
\end{figure}

Another idea is that since our curve has 6 empty ovals one should
look at the corresponding hexagon and at pencils of cubics spanned
by 2 triangles.
Specifically, we may choose a hexagon which visits the 3 black
inner points and the 3 white outer points in dichromatic
alternation (BWBWBW), cf. e.g. Fig.\,d. Then we may expect that
the pencil of cubics spanned by the red and green cubics is total.
Alas Fig.\,e refutes this expectation. Of course there are other
dichromatic hexagons but this is unlikely to be the right method.
For instance Fig.\,f is another dichromatic hexagon, yet the
corresponding pencil is still not total as shown by Fig.\,g.

In conclusion those $(M-4)$-schemes are a bit puzzling from the
viewpoint of total reality as there is no (obvious) canonically
defined pencils since we have 6 empty ovals, which is not the
number of basepoints of a pencil of plane curves, namely $4$ for
conics and 8 for cubics. The two extra virtual basepoints for
cubics could be chosen as high-order contacts imposed to the
pencil and this done properly could exhibit total reality. It
remains however to understand the natural geometric condition that
are so-to-speak imposed by the geometrical vision of the curve.
Remind indeed that total reality always amounts to place the
ocular system ``inside'' of the glass so that the latter has no
apparent contour (compare the baby case of the G\"urtelkurve,
Fig.\,\ref{Guertel-saturated:fig}.)

{\it $(M-2)$-curves}.---We may hope that the situation is improved
when moving to $(M-2)$-curves. The first case to study is the
scheme $\frac{8}{1}$. Fig.\,\ref{Total8-1:fig} shows  models of
both types I vs. II, but it is again quite puzzling to decide
which intrinsic criterion distinguishes both configurations. Of
course a loose answer could be that the type~I configuration is
characterized by the fact that the pencil of cubics through the 8
empty ovals is total, however one could desire a more optical
recognition algorithm. Perhaps what distinguishes the type~I is
the possibility of tracing a convex octagon through the empty
ovals. Note that convexity has some meaning since given two points
in the nonempty oval $N$ we shall always select the
half-projective line (segment) which is inside this oval $N$.
({\it Added} [09.04.13].---But the oval $N$ can be non-convex, and
so this is meaningful only when the 2 points are inside the deep
ovals of course.)

\begin{figure}[h]
\centering
\epsfig{figure=Total8-1.eps,width=122mm} \vskip-5pt\penalty0
  \caption{\label{Total8-1:fig}%
  Failing to understand total reality of the scheme $\frac{8}{1}$} \vskip-5pt\penalty0
\end{figure}

For this scheme $\frac{8}{1}$, one could nearly argue that the
pencil of cubics through the deep 8 points is always total, for we
have $8\cdot 2$ automatic intersections, plus the 2 coming from
the fact that the cubic is not null-homotopic hence must cut twice
the nonempty oval. However this would contradict Rohlin's remark
that this scheme is indefinite as shown by the above constructions
(pictures). However the sole obstruction to total reality of a
cubic in this deep pencil is that the cubic has a small oval
entirely inside some of the empty ovals. Indeed if the cubic is
connected total reality is clear as $2\cdot 8+2=18=3\cdot 6$, and
if not yet the oval of the $C_3$ visits at least two ovals of the
sextic then each oval visited contribute for 2 intersections and
total reality is evident.

{\it Insertion} [09.04.13].---This ``soleness'' looks inexact:
another obstruction occurs when the cubic splits off an oval
visiting all 8 basepoints and the pseudoline stays confined
outside the nonempty oval of the $C_6$.

One would like to show that under a suitable geometric hypothesis
(capturing the essence of the type~I scheme) this sole obstacle
cannot occur.

[17.02.13] All these questions are fairly delicate and have to be
extended to all type~I schemes listed by Rohlin (see again the
Gudkov-Rohlin Table=Fig.\,\ref{Gudkov-Table3:fig}). Precisely what
is demanded is an optical recognition procedure of the type in the
sense of Klein (orthosymmetry vs. diasymmetry) via a synthetical
device ensuring total reality of a certain class of pencils
naturally attached to the curve (or its schemes). This would
extend somehow Rohlin's claim of the absolute orthosymmetry of the
sextic schemes $\frac{6}{1}2$ and $\frac{2}{1}6$, which is the
purest manifestation of the phenomenon. Meanwhile (yesterday),
S\'everine Fiedler-Le~Touz\'e informed
us (and several other colleagues, cf. letter in
Sec.\,\ref{e-mail-Viro:sec} dated [16.02.13]) that she was able to
prove Rohlin's claim for the scheme $\frac{2}{1}6$. Probably her
argument contains crucial ideas that solve as well our slightly
generalized problematic.

{\it Insertion} [09.04.13].---This is nearly true, safe that it
turned out that Le~Touz\'e 2013
\cite{Fiedler-Le-Touzé_2013-Totally-real-pencils-Cubics} proves a
slightly weaker assertion than the full Rohlin claim, namely she
relies on the RKM-congruence ensuring type~I a priori.

Concretely we have then the scheme $\frac{4}{1}4$ for which one
need to formulate an optical recognition, and the scheme $9$. For
the latter it seems that the type~I configuration is characterized
by the fact that the pencil of cubics through 8 deep points inside
some of $8$ ovals is such that the 9th basepoint lands in the 9th
oval.
A lucky stroke!

{\it Insertion} [09.04.13].---It is nearly implicit from the above
that we posit:

\begin{conj}
Any dividing $(M-2)$-sextic has a total pencil of cubics. In a
stronger shape: any cubics-pencil assigned to pass in the $8$ deep
ovals is total, and it is permissible to let degenerate the
basepoints on the ovals themselves.
\end{conj}

Then it is also plain to see that this is much in line with
Gabard's bound $r+p=10$, since $3\cdot 6-8=10$ is the mapping
degree when the 8 basepoints are degenerated upon the ovals.

\subsection{Total reality of plane $M$-curves}

[17.02.13] Next we have the case of $M$-sextics. This puzzled me
for a while, yet it seems clear that now cubics pencils will not
exhibit total reality. The reason is that we have 10 empty ovals
but only 9 basepoints available. On the other hand it seems
evident that each empty oval must contain a singularity of
foyer-type corresponding to a basepoint of the pencil. So this
follows from Poincar\'e's index formula applied to the foliation
induced by the pencil (cf. Lemma~\ref{Poincare-lower-bound} much
earlier in this text, but restituted below in perhaps clearer
fashion).

\smallskip
{\footnotesize

{\it Insertion} [09.04.13].---To clarify a bit the sloppy logics,
the crucial point is---as vaguely explained in that lemma---that
total reality forces transversality of the foliation along the
ovals. This is only approximatively true (consider the
G\"urtelkurve swept by a pencil of lines based on  the deep oval).
Yet perturbing a bit, we can avoid this case. Then, granting
transversality of the foliation, the index formula applied to the
inside of the oval gives at least one singularity of positive
index, where half-valued integer like $+1/2$ (thorns) are
precluded by algebraicity of the foliation. So we infer, in each
deep oval, the presence of a foyer-type singularity of index $+1$
(at least), materializing in turn to a basepoint. Higher indices
are non-algebraic if semi-integral, and if $+2$ or more then they
can be considered as a dipole (coalescence of 2 foyers) hence
consuming more basepoints. So at any rate we have $D\le k^2$ by
choosing in each of the $D$ many deep ovals a foyer-type
singularity imagined of simple multiplicity (index $+1$), and
where $k$ is the order of the pencil.

}
\smallskip

This being said we shall move to pencil of quartics. The crucial
idea is to remind the synthetic proof in the abstract context of
the schlichtartig avatar of Ahlfors's theorem, i.e. the theorem
due to Riemann 1857-Schottky 1875--77-Bieberbach 1925-Grunsky
1937. More precisely we have in mind the simple argument via
Riemann-Roch rediscovered by Huisman and Gabard, yet first clearly
enunciated in Enriques-Chisini 1915. Bypassing all these
historical details, the logical argument is simply given in our
Lemma~\ref{Enriques-Chisini:lemma} (prior in this text). The idea
is merely that if one has an abstract $M$-curve (not necessarily
plane), then choosing one point on each oval (=real circuit which
is linguistically better in this abstract context) one has a group
of $g+1$ points which therefore move in its linear equivalence
class by Riemann(-Roch), or just by Abel's theorem since there are
$g$ Abelian differentials (holomorphic one-forms) imposing
magneto-hydrodynamical constraints upon the motion
of a divisor in
its linear equivalence class.

So our effective divisor of degree $g+1$ moves on the curve of
genus $g$. Since there is only one point on each ``oval'', it is
like a miniature railroad, in which there is only one train one
each track, and so there cannot be collisions and total reality is
automatic.

 Now when the $M$-curve is plane, I was frustrated to
know nothing on the degree of such total maps. However the answer
is ``{\it toute simple}''(=very simple). Indeed inspired by the
abstract proof, choose one point on each oval (and also one on the
pseudoline if there is one). Then there is a standard recipe to
construct the linear series spanned by a given group of points on
a plane curve $C_m$ (due to Brill-Noether 1873/74?,
Enriques-Chisini's book 1915, Severi's book 1921
\cite{Severi_1921-Vorlesungen-u-alg.-Geom-BUCH},  van~der
Waerden's book 1939/73 \cite{van-der-Waerden_1939/73}, Walker's
book 1950 \cite{Walker_1950/62}, who else?): just choose an
integer $k$ large enough so as to have enough free parameters to
pass a $k$-tics through the given group of points. Choose such a
curve $C_k$ and look at the residual intersection with the curve
$C_m$. Consider then all $C_k$'s passing through this residual
intersection and the latter cut on the curve  groups whose mobile
part are divisors equivalent to the given one. This method clearly
belongs to the genre of a sweeping method (balayage).

Applying this to an $M$-curve leads to the following very modest
theorem (stated as a such just because it escaped my attention for
several months, if not years):

\begin{theorem}\label{total-reality-of-plane-M-curves:thm}
Given any plane $M$-curve of degree $m$ there is a total pencil of
$k$-tics of degree $k=m-2$ (two units less than the given degree
$m$). In fact exactly like in the abstract Bieberbach-Grunsky
theorem, any equidistribution of points (i.e., one point on each
real circuit) moves in a linear system of dimension $\ge 1$ and
induces a totally real pencil by the sweeping method. In
particular each $M$-sextic is total under a pencil of quartics.
\end{theorem}

\begin{proof} As the proof involves some arithmetical nonsense
it is didactic to first handle the case of sextics. Then Harnack's
bound (in Petrovskii's notation) is
$M=g+1=\frac{(m-1)(m-2)}{2}+1=11$. The space of $k$-tics has
dimension $\dim \vert k H \vert=\binom{k+2}{2}-1$, that is $5$ for
conics, $9$ for cubics, 14 for quartics, etc. Choose an
equidistribution of $11$ points one on each oval of the $C_6$.
Then quartics have enough freedom to
visit them. Choose a $C_4$ passing through the 11 points, and the
residual group has $4\cdot 6-11=24-11=13$ points. But this is
exactly one less than the dimension of
all quartics, and so the residual series---consisting of all
curves passing through the residual group---gives the required
pencil. The total reality of the latter follows by the
non-collision principle involving the continuity argument implicit
in Enriques-Chisini's anticipation of the Bieberbach-Grunsky
theorem.
Of course the impossible-to-beat anticipation is Riemann's 1857
Nachlass \cite{Riemann_1857_Nachlass}!

The general case is merely the same numerological coincidence
worked out in general. Given any $M$-curve $C_m$ of degree $m$,
Harnack's bound is $M=g+1=\frac{(m-1)(m-2)}{2}+1$. Choose $M$
points on the real locus $C_m({\RR})$, one on each oval. Locate
the least integer $k$ such that $\dim \vert k H \vert \ge M$.
Since both the genus and the dimension of this complete linear
system are given by binomial coefficients, this traduces into
$\binom{k+2}{2}-1 \ge \binom{m-1}{2}+1$ which is first satisfied
for $k=m-2$ (but not at $k=m-3$). Indeed this amounts to
$\binom{m}{2}-1 \ge \binom{m-1}{2}+1$ which is plain as
$\binom{m}{2}=1+2+3+\dots+(m-1)$. Now the residual intersection of
a $C_k$ through the $M$ points with $C_m$ gives so many points as
the following expression, which turns out to be the dimension of
the system $\vert k H \vert $ less one unit, as shown by the
following boring calculation:
\begin{align*}
k\cdot m-M&=(m-2)m-\textstyle\frac{(m-1)(m-2)}{2}-1 \cr
&\textstyle =(m-\frac{m-1}{2})(m-2)-1
=(\frac{m+1}{2})(m-2)-1=\dots=\dim \vert k H \vert -1.
\end{align*}
%
Somewhat more elegantly,
\begin{align*}
k\cdot m-M &=(m-2) m-[1+2+\dots+(m-2)]-1 \cr
&=[(m-1)+(m-2)+\dots+2]-1=\textstyle\binom{m}{2}-2=\dim\vert kH
\vert-1,
\end{align*}
for $k=m-2$.
\end{proof}

Several questions arise as usual after discovering a trivial
truth. (Derri\`ere les montagnes encore des montagnes: Proverbe
des iles cr\'eoles, if I remember well). In the case of sextics
one can probably say therefore (modulo a more careful analysis of
the case of $(M-4)$-sextics that all dividing sextics have their
total reality exhibited by a pencil of degree $\le 4$). Perhaps it
is true in general that:

\begin{conj} Any dividing $m$-tic has its total reality exhibited by a
pencil of curves of order less than $(m-2)$.
\end{conj}

A more serious game would be to see if the above theorem
essentially due to
Riemann-Enriques-Chisini-Bieberbach-Grunsky-Wirtinger-Huisman-Gabard
does not imply when suitably complemented by foliation theory \`a
la Poincar\'e the well-known obstruction of
Hilbert-Rohn-Petrovskii-Gudkov for $M$-sextics (e.g. the highbrow
Gudkov-Rohlin congruence  mod 8 or at least the weak version
thereof mod 4 due to Arnold $\chi=p-n\equiv k^2 \pmod 4$). For
higher degrees one may even dream of new results along this
method, but all this requires more serious
work.

One could even dream that the method extends outside the realm of
$M$-curves, as to recover e.g. Rohlin's claim (meanwhile
Le~Touz\'e's theorem) but this is unlikely because the continuity
principle of no collision meets then serious difficulties, which
are precisely those making Ahlfors theorem harder than
Bieberbach-Grunsky's theorem. Gabard 2006 gives an abstract
topological algorithm
overcoming this difficulty of collisions, but it looks hard to
transplant this to the context of Hilbert's 16th problem.
%
{\it Added} [09.04.13].---The key is of course that if one has an
overpopulation of trains circulating on a track one must ensure
dextrogyration so as to avoid collisions. Then total reality is
granted.

[02.03.13] As a lovely special instance of the above theorem
(\ref{total-reality-of-plane-M-curves:thm}), Le~Touz\'e (1 March
2013 \cite{Fiedler-Le-Touzé_2013-Totally-real-pencils-Cubics})
observes the:

\begin{Scholium} {\rm (Le~Touz\'e 2013)}
\label{LeTouze-quintic:scholie} ``For an $M$-quintic $\langle J
\sqcup 6 \rangle$, one finds a suitable pencil of cubics with six
basepoints distributed on the six ovals, and two further chosen on
the odd component $J$. As this component must cut any cubic an odd
number of times, the required $15$ real intersections are
granted.''

---{\rm Gabard's addendum [09.04.13]:} As it will be observed in the
sequel, but can already be noted here, one can also avoid this
topological argument with the pseudoline, by noting that since
$2\cdot 6+2=14$ intersections are granted,  the remaining one is
forced to reality by algebra (Galois-Tartaglia\footnote{Tartaglia,
also known as Niccolo Fontana (1500?--1557) won in 1535 a
mathematical contest by solving many different cubics, and gave
his solution to Cardano (1501--1576), who published in 1539
``Artis Magnae'' alias the ``Great Art, or the Rules of Algebra'',
where complex numbers are used in Cardano's formula to express the
real roots of cubics.} involution).
\end{Scholium}

{\it Long Insertion (ca. $2\frac{1}{2}$ pages)} [09.04.13].---How
to generalize this Le~Touz\'e's Scholium to an $M$-septic? Pencil
of quintics have 19 basepoints assignable. Assign 17 of them on
ovals and 2 on the pseudoline, then there is $34+3=37$ real
intersections granted, overwhelming the $5\cdot 7=35$ of B\'ezout.
But by Harnack $M_7=g_7+1=\frac{6\cdot 5}{2}+1=16$, so we do not
have as many ovals as 17. Actually our sloppy argument reproves
Harnack's bound with one unit less. This is, by the way, not so
surprising as Enriques-Chisini 1915'a purpose was precisely to
re-derive a proof of Harnack via Riemann-Roch. Okay, but this
sounds too modest and we can surely expect more clever
generalizations of Le~Touz\'e's Scholium. Indeed, by Harnack we
have $15$ ovals on the $C_7$ (recall from M\"obius-von Staudt that
odd order curves have exactly one pseudoline). Distribute the 19
basepoints on the $15$ many ovals, plus 4 on the pseudoline. Then
$30+4=34$ real intersections are granted by topology, and the last
one is forced by algebra (Galois-Tartaglia symmetry of complex
conjugation). Hence total reality is granted. We have proven the:

\begin{lemma}
Given any $M$-septic $C_7$, the pencil of quintic assigned to pass
through the inside of all $15$ ovals (warning: it is more prudent
to assign them directly on the ovals themselves) of $C_7$ and $4$
points marked on the pseudoline is totally real.
\end{lemma}

Two questions arise. A deep one is whether this can be used to
infer something about Hilbert's 16th (distribution of ovals solved
for $m=7$ essentially by a single hero, Viro ca. 1979). Another
question is whether Le~Touz\'e's Scholium extends to all (odd)
degrees. More philosophically, it seems that our abstract argument
of Theorem~\ref{total-reality-of-plane-M-curves:thm} works in all
degrees and gives the required total pencil of degree $(m-2)$, yet
 Le~Touz\'e's scholium (\ref{LeTouze-quintic:scholie})
 looks sharper as it tells
precisely where to assign basepoints, without using the sweeping
method and the Restsatz \`a la Brill-Noether.

In degree 6, Le~Touz\'e's method suggests looking at a
Harnack-maximal $C_6$ swept by $C_4$'s. Those (quartics) have 13
basepoints assignable. Distribute $11$ of them on the ovals, but
where to place the remaining 2? A priori only 22 real
intersections are granted. Of course we still have 2 more
basepoints, but they do not force new intersections. How to ensure
total reality in this case?

Let us look at degree $m=9$. Then septics have
$\binom{7+2}{2}-2=34$ basepoints assignable (for a pencil). For
$m=9$, Harnack's bound is $M_9=\frac{8\cdot7}{2}+1=29$. So
distribute the 34 bases on the 28 ovals plus 6 on the pseudoline,
granting so $2\cdot 28+6 = 56+6=62$ intersections (just one less
than B\'ezout's $7\cdot 9=63$), but algebra forces reality of the
last man. At this stage it is evident that Le~Touz\'e's Scholium
extends to all odd degrees as follows:

\begin{theorem}\label{Le-Touzé-extended-in-odd-degree:scholium} (Extended Le~Touz\'e's Scholium).---Given
any $M$-curve of odd degree $m$, the pencil of $(m-2)$-tics
assigned to visit once all ovals and with residual collection of
basepoints assigned on the pseudoline is totally real.
Further, the  mobile part of the pencil has exactly one point
moving on each real circuit, and the allied circle map has lowest
possible mapping-degree namely the number $r=M$ of real circuits
(exactly like in the Riemann-Bieberbach-Grunsky theorem).
\end{theorem}

\begin{proof} What is first chocking is that when $m=5$
(Le~Touz\'e's Scholium) there is really an extra intersection
gained by topology (of the pseudoline), while for $m=7,9$ our
argument merely uses algebra. So a priori the argument could split
in two cases depending on some (sordid) periodicity modulo 4. The
geometer dislikes intrusion of capitalism and arithmetics in his
garden. However it should be observed that even in Le~Touz\'e's
argument one can use algebra, which is alas more capitalistic than
her geometric argument. So it is still reasonable to expect an
unified proof (without mod 4 stories) along pure arithmetical
nonsense.

This is as follows: let $C_m$ be our $M$-curve of odd degree
$m=2k+1$. We look at curves $C_{m-2}$. Those can be
assigned
 $\binom{m}{2}-2=:B$ many basepoints. But by
``M\"obius-von Staudt'' $C_m$ has $M-1$ ovals (just omit  the
pseudoline of course), and by Harnack $M-1=g=\binom{m-1}{2}$. We
 distribute the $B$ basepoints on the $g$ ovals and the $B-g$
 remaining ones on the pseudoline. Let us calculate
$B-g=[(1+2+\dots+(m-1))-2]-[1+2+\dots+(m-2)]=(m-1)-2=m-3$. So we
have $2g+(B-g)$ real intersections granted by topology (of the
ovals), and this is equal to $(m-1)(m-2)+(m-3)=m(m-2)-1$. But this
is one unity less than B\'ezout's number $m(m-2)$ of complex
intersections in $C_m\cap C_{m-2}$, so that the last intersection
is forced to reality too! Did we used the assumption that $m$ is
odd in any dramatic fashion? I would say no, but we did! Probably
the argument adapts to even degrees as well if one is a bit more
clever than we were for $m=6$.

As to the last clause, it is evident by construction. Indeed since
one basepoint is assigned on each oval, some extra (mobile)
intersection is created on the oval. Further the last intersection
granted by the Galois-Tartaglia symmetry conj, is forced to live
on the pseudoline, since each real circuit contains at least one
mobile point (by an evident sweeping principle or just the fact
that any point has a well defined image).
\end{proof}

Hence a problem of interest is to understand the even degree case,
and we hope that someone will easily tackle this question. It is
quite beautiful at this stage to feel some big harmony between
Riemann, Harnack,  Brill-Noether as well as Le~Touz\'e,
or Rohlin, at least a sort of unity between conformal and
algebraic geometry. Poincar\'e wrote something like the following:
{\it ``La pens\'ee n'est qu'un \'eclair dans la nuit, mais c'est
ce qui \'eclaire tout''}.

In fact, when $m=6$ we have 13 basepoints (for quartics) and 11
ovals on the $C_6$. Could it be useful to impose the 2 additional
basepoints as imaginary conjugate on the curve $C_6$. Those points
will not be real, yet since they are statical they do not spoil
total reality which merely involves  dynamical points. (Note: this
is not a new idea, cf. e.g. p.\,7 of Gabard's Thesis 2004
\cite{Gabard_2004}.) Note also that in the above proof
(\ref{total-reality-of-plane-M-curves:thm}) adapting
Bieberbach-Grunsky to plane $M$-curves nothing grants that
basepoints are real. In fact we start from any group of $g+1$
points equidistributed on the $M=g+1$ circuits, pass a curve of
sufficiently large degree (i.e. $(m-2)$) through them and look at
the residual intersection, which is a priori not totally real (but
just real, stable under conj). This gives perhaps some evidence
that we should for $m=6$ permit a pair of conjugate basepoints.
Doing so we really have 11 real points moving on each oval in
accordance with the train-track principle (i.e.
Bieberbach-Grunsky, or Enriques-Chisini, etc.)

So we state:

\begin{lemma}\label{Le-Touzé-scholium-deg-6:lem}
Given any $M$-sextic, the pencil of quartics assigned to visit any
$11$ points marked on the $11$ ovals and a pair of conjugate
points of $C_6$ is totally real and induces a (circle) map of
degree $11$.
\end{lemma}

\begin{proof} Quartics depend upon $\binom{4+2}{2}-1=14$
parameters and so $13$ basepoints may be assigned. Distribute them
on the 11 ovals available and fix the 2 remaining ones as a
conjugate pair of points of the $C_6$. For topological reason each
real curve of the pencil cuts once more each oval (usual closing
lemma for ovals), and so $2\cdot 11=22$ real intersections are
granted. By B\'ezout there is a total of $4\cdot 6=24$
intersections. Hence, {\it all} intersections are under control,
i.e. either the 22 real ones or 2 imaginary ones which are
statical. The latter do {\it not\/} perturb total reality, since
they are not moving ``electrons''.
\end{proof}

Probably the statement extends to all other even degrees (by
working properly the arithmetics eventually by using what we
already calculated in the odd degree case).

More geometrically (and returning to $m=6$), one may wonder if we
could not by continuity push the 2 imaginary basepoints on the
real locus so as to impose a tangential contact in the limit
(zusammenr\"ucken). Then the modest advantage is that all
basepoints would again be visible on the reals, and total reality
should be conserved by continuity. So we arrive at the:

\begin{lemma}
For any $M$-sextic, the pencil of quartics assigned to visit any
$11$ points marked on the $11$ ovals  and tangent at a $12$th
point of the $C_6$ is totally real and induces a (circle) map of
degree $11$.
\end{lemma}

Hence some contact can be imposed at any point, and the resulting
foliation will look like a dipole at this point of tangency. It
would be interesting to see if we can infer any of the deep
classical obstructions on the distributions of ovals due to
Hilbert, Rohn, Gudkov, from this method.

Another idea (to be explored better than what follows) is to
assign such imaginary basepoints on the puzzling case of
$(M-4)$-sextics (of type~I), cf.
Scho\-lium~\ref{(M-4)-sextics-corrupt-Gabard:scholium}. Then we
had a pencil of cubics, of which we distribute the 8 basepoints on
the 6 empty ovals and 2 points remain left. A priori we could
imagine that an imaginary pair counts just for one linear
condition (since after all the passage through one of them forces
passing through the conjugate). So we could impose 2 imaginary
pairs of additional basepoints, and the mapping-degree of the
pencil would be $3\cdot 6- 6 -4=8$, again in accordance with
Gabard's bound $r+p=9$ (best interpreted as the mean of $r$ and
Harnack's bound). This looks however dubious since the pencil of
cubics would then have $6+4=10$ basepoints overwhelming B\'ezout.
This can be repaired if we impose only 5 real points and 2
imaginary pairs, and then Gabard's bound is (exactly) verified,
since $3\cdot 6-5-4=9$. However it another piece of work to check
that total reality can be ensured. So let us be happy with a
conjectural (and admittedly vague) statement:

\begin{conj}
Any $(M-4)$-sextic of type~I admits a total cubics-pencil  of
degree $9$ (like Gabard) with $5$ basepoints assigned on the oval
and $2$ pairs of imaginary basepoints assigned on the curve.
\end{conj}

Here only $10$ real intersections are granted (among the $18-4=14$
which are moving) and so much work remains to be done to ensure
totality under suitable assumptions. If we impose 6 real
basepoints then 12 are granted, etc.
{\it End Long insertion.}

\smallskip

[18.02.13] Let us try the following strategy. [We now come back to
Rohlin's total reality problem.] Suppose the sextic to have the
scheme $\frac{6}{1}2$. Let us choose 8 points $p_i$ on the $8$
empty ovals, one on each oval, and consider the corresponding
pencil of cubics. We would like to show total reality of this
pencil. It is clear that any cubic cuts at least $2\cdot 8=16$
times the $C_6$. This is because the $C_3$ can be either tangent
at $p_i$ to $C_6$ or transverse. In the first case, intersection
multiplicity is two, while in the second, one side of an
infinitesimal analytic arc of the curve is inside the oval while
the outer is out, hence an extra intersection is gained by closing
the real circuit.

Take any cubic in the pencil which is {\it connected\/} (and
smooth). The latter is clearly totally real as 2 bonus
intersections are created on the nonempty oval of the $C_6$ (which
acts as a separator between the inner and outer basepoints). Note
that we use the lemma that any pencil of cubics contains a
connected cubics
(which we nearly proved via
Lemma~\ref{nodal-cubics-8-many-in-a-pencil:lem} showing that a
pencil contains in general 8 nodal cubics).

Now a simple idea ensuring total reality would be to look at a
nearby cubic $C_3'$ and look if the 2 extra intersections {\it
gyrate in the same sense\/} (dextrogyrate) on the nonempty oval
$N$. In  case of dextrogyration, total reality follows because
there would be no collision between both points on $N$. Indeed
when we move the cubic curve in the pencil then the intersection
points move continuously and none of them can suddenly change its
sense of motion, for otherwise there would be 2 curves of the
pencil (nearby hence distinct) passing through the same point.

So it suffices checking that the 2 intersections $p',q'$ of
$C_3'\cap C_6$ located on $N$ move in same sense (dextrogyrate) on
$N$, equivalently that $p',q'$ are separated by the corresponding
2 intersections $p,q$ for $C_3\cap N$.

Since the curve $C_3'$ is a small perturbation of $C_3$ it
oscillates about it (in a slaloming fashion). Now a simple picture
shows that the gyration is good (occurs in the same sense, or
dextrogyre) iff the number of basepoints inside $N$ is odd. So the
whole question reduces to knowing if the 9th (non-assigned)
basepoint of the pencil is located inside $N$ (or not). If it is
inside then we are finished and total reality follows. Alas I know
about of no argument ensuring the inside-ness of the 9th
basepoint.

The above argument (or rather strategy) relies on the existence of
a connected cubic in the pencil which must be a simple matter.
This can be bypassed if we argue differently. It is clear that the
sole obstruction to total reality
is a disconnected cubic whose oval lies inside $N$. Such a cubic
is smooth except if it has a solitary node, yet in that case
total reality is evident for the pseudoline of the solitary cubics
has to connect an inner and outer point so contribute for an extra
17th intersections. Then either by algebra or topology the 18th
intersection is real too. Given such a bad cubic $C_3$ which is
smooth and whose ovals lies inside $N$, we can again look at a
small perturbation $C_3'$ which will oscillate about $C_3$, and so
 do the corresponding ovals. Now it is clear by a slaloming
argument that the oscillation is possible iff the number of inner
basepoints (inside $N$) is even. Hence again we would have a
contradiction, if we knew that the 9th basepoint of the
cubics-pencil lies inside $N$.

Whatever the strategy adopted, Rohlin's claim seems to require
innerness of the 9th basepoint. So this gives a 2nd reduction of
Rohlin's claim.

As a metaphor it seems that such total reality proofs \`a la
Rohlin-Le~Touz\'e are akin to an Eiger-Nordwand ascension. There
are several base-camps where to rest, but as the climbing goes on
they become rarer and rarer and one is forced to follow a nearly
canonical route,
[more and more vertiginous and perilous, by the way.]
It should be noted yet that our approach is slightly weaker than
the Rohlin-Le~Touz\'e claim for we do not check total reality of
all pencils with 8 deep basepoints {\it inside\/} the empty ovals,
but merely the case when the latter 8 points are located {\it
on\/} the ovals. Our weaker variant suffices yet to detect total
reality and so the type~I of such Rohlin's schemes (e.g.
$\frac{6}{1}2$). Our tactic looks simpler, since when basepoints
are assigned in the interior of ovals, they do not create defacto
real intersections, because the cubic's oval could be
microscopically nested inside one oval of the sextic.

{\it A dubious strategy with the hexagon (skip the next 2
paragraphs)}.---Another idea was suggested by a naive look at the
Hilbert-style construction of $\frac{6}{1}2$. Here it seems that
the hexagon through the inner points is not convex. $\bigstar${\it
Inserted Objection} [09.04.13].---This is so on the naive
Walt-Disney picture Fig.\,\ref{GudHilb8:fig}, but less clear on
the more realist picture Fig.\,\ref{GudHilb6-12:fig}.$\bigstar$ On
the other hand if there is a bad cubic (one whose oval is inside
the nonempty oval of the $C_6$) then we know (since at least
Zeuthen 1874 \cite{Zeuthen_1874}) that the oval of the cubic is
convex. This is to mean that whenever we join two points inside
the oval by the rectilinear segment inside the oval it stays
entirely inside the oval. This is not well phrased and should be
formulated by saying that whenever we take 2 points inside the
cubic-oval the line through it dissected in 2 pieces by the 2
intersection points with the oval is such that the half not
meeting the pseudoline is entirely within the inside of the oval.

So Rohlin's claim would follow if it can be shown that the
fundamental hexagon of our $C_6$ of type $\frac{6}{1}2$ is
non-convex. Here the fundamental hexagon is defined as the union
of all fundamental 2-simplices (triangles). Recall that given 3
points on the inner ovals (inside $N$) the lines joining them are
B\'ezout-saturated and there is a unique full-triangle traced
inside $N$, which we call fundamental. The fundamental hexagon of
our $C_6$ (with $8$ marked points $p_i$ on the empty ovals)  is
the union of all these fundamental triangles rooted at the 6 inner
points.

\smallskip
Another strategy is as follows. Choose any 8 points on the sextic
$C_6$, one on each empty oval. To show: the pencil of cubics
through them is totally real.

{\it Step 0}.---The sole obstruction to total reality is the
presence of  a bad cubic, i.e. one with an oval entirely traced
inside the nonempty oval of $C_6$.

{\it Step 1}.---Assume that there is a bad cubic then the 6 inner
points are hexagonally distributed on the oval.

{\it Step 2}.---Imagine the 6 inner points black and the 2 outer
points white-colored. Then try to infer existence of a conic $C_2$
through 3 inner points and the 2 outer points which is
dichromatic, i.e. such that the 2 white points split the 3 black
points in 2 groups.

{\it Step 3}.---Such a conic has 4 transitions from black to
white, hence cuts the $C_6$ in $10+4=14>12=2\cdot 6$ violating
B\'ezout.

The difficult step is Step~2. To exhibit a dichromatic conic it
suffices by Le~Touz\'e's lemma (\ref{LeTouze:lem}) that some
fundamental triangle through 3 black points separates the 2 white
points.

So by contradiction assume that all black triangles does not
separate the 2 white points. But alas I do not know why this
circumstance (which is actually forced by Le~Touz\'e's lemma)
implies a contradiction with the bad cubic assumption.

\subsection{Yet another strategy via long run
evolution of the bad cubic}

[19.02.13] In this section we explore another strategy toward a
proof of the Rohlin-Le~Touz\'e claim of total reality for Rohlin's
curve $\frac{6}{1}2$. The argument perhaps adapts to its mirror
$\frac{2}{1}6$, yet we concentrate on $\frac{6}{1}2$ for
simplicity. As above, we rather attack the somewhat weaker total
reality assertion for a pencil with basepoints assigned {\it on\/}
the ovals of the curve. By letting  basepoints degenerate on the
ovals, this is probably logically implied by the
Rohlin-Le~Touz\'e's theorem (of which at the time of writing we
have not seen a proof).
It seems also that Le~Touz\'e proves rather the case of the mirror
$\frac{2}{1}6$ but probably her argument adapts to
$\frac{6}{1}2$.
Our argument is just a strategy far from a complete proof, trying
to study the dynamical evolution of a pencil lacking total reality
while hoping to detect a contradiction. So it is a dynamical
approach, but perhaps the real proof (of Rohlin and Le~Touz\'e) is
more clear-cut or based perhaps on the same idea.

Start by recalling certain trivialities, which we repeat for
convenience. In all this section, $C_6$ denotes a ``Rohlin curve''
of type $\frac{6}{1}2$, i.e. 6 ovals enveloped in a larger one
with 2 ovals outside. The 6 ovals are said to be inner ovals and
the 2 ovals outside called outer ovals.

\begin{defn}
A pencil of cubics passing through $8$ basepoints injectively
distributed on the $8$ empty ovals is said to be deep.
\end{defn}

The following statement (for me still hypothetical) is a variant
of the Rohlin 1978--Le~Touz\'e 2013 theorem, and probably weaker
than it, yet which seems to us easier to prove as there is not the
possibility of microscopic ovals passing through the assigned
basepoints yet without creating real intersections. Albeit weaker
it is sufficient for detecting the type~I of the given scheme, and
in some sense stronger as it yields circle maps (or totally real
maps) of smaller mapping-degree.

\begin{theorem}
Any deep pencil on a sextic $C_6$ of real scheme $\frac{6}{1}2$ is
totally real, i.e. any (real) curve of the pencil cuts only real
points on the $C_6$.
\end{theorem}

The sequel is an (unsuccessful) attempt of proof of this
(hypothetical) statement.

\smallskip
\begin{proof} It is divided in several steps, each justified
 in the subsequent paragraph.

$\bullet$ {\it Step 1}.---Each cubic $C_3$ of such a deep pencil
has at least 16 (real) intersections with the  $C_6$.

Indeed $C_3$ is assigned to pass through 8 points $p_i$
distributed on the 8 empty ovals of $C_6$. Two cases can occur.
Either the cubic is tangent to the sextic at $p_i$ in which case
we have intersection multiplicity $2$, or the $C_3$ is transverse
to $C_6$ in which case there is through $p_i$ a small analytic arc
of the $C_3$ with extremities both inside and outside the
corresponding oval of $C_6$. By basic properties of algebraic
projective curves, this arc of curve has to close up itself and so
a 2nd (real) intersection with $C_6$ is created. A priori $C_3$
could visit $p_i$ via just a solitary node (isolated real ordinary
double point). In that case the intersection multiplicity is still
$2$, and by the way I suspect that this case cannot occur by
elementary properties of pencils which have a foyer-type
singularity at the basepoints preventing an isolated singularity
to appear there. (All this is clumsy due to a lack of profound
algebro-geometric knowledge of the writer.)

Step~1 shows that we are quite close to total reality, where {\it
each\/} cubic curve is required to have 18 real intersections
(counted by multiplicity) with the sextic $C_6$.

$\bullet$ {\it Step 2}.---The sole obstruction to total reality is
the presence of a {\it bad cubic\/}, i.e. a smooth cubic with 2
components whose oval is contained in the nonempty oval $N$ of the
$C_6$.

If the cubic $C_3$ is connected (i.e. $C_3(\RR)$ is connected),
then as it must visit both inner and outer points
a 17th intersection is created and the 18th follows either by
algebra or topology. If $C_3$ is not connected then it is either
smooth with 2 components, or  a solitary cubic with a solitary
node. In the latter case the solitary node passes at most through
one of the eight $p_i$ (though this is improbable), yet even in
that case the pseudoline of the solitary cubics visits both inner
and outer points so has to be total. Hence the sole curve possibly
failing total reality is a smooth cubic with 2 real branches. It
has further to be monochromatic in the sense that the outer and
inner points $p_i$ have to be ``purely'' distributed on both real
circuits of the $C_3$. Else if both an inner and an outer point
among the $p_i$ land on a same circuit of $C_3$ then a 17th
intersection is created by topology, and so an 18th one by
algebra. Further if the  inner points are on the pseudoline of
$C_3$, then topology forces  a 17th intersection (else the
pseudoline would be contractible inside the bounding disc of the
nonempty oval $N$). So the inner points are on the oval of $C_3$,
and Step~2 is completed.

So from now on we shall assume  that our deep pencil contains a
bad cubic $C_3$, and try to infer a contradiction. Several basic
remarks are perhaps useful.

1. The unique oval of a cubic with 2 components is convex in some
obvious sense. (Perhaps this already implies a contradiction, but
need to be detailed.)

2. The oval of our bad cubic $C_3$  will vibrate during an
infinitesimal motion along the (deep) pencil $\Pi$. As $6$
basepoints are assigned on the oval $O$ of $C_3$, a vibratory
(slaloming) principle implies that the oval oscillates an even
number of times across itself. (Of course this may also be reduced
to homological intersection mod 2.) It follows that the 9th
basepoint of the pencil $\Pi$ is located on the pseudoline of the
bad cubic $C_3$.

Now our strategy is the naive one of studying the long-run
evolution of the bad cubic as time evolves, i.e. as  the cubic is
dragged along the pencil. Probably the real
argument of Rohlin-Le~Touz\'e is more clear-cut B\'ezout-style
obstruction without dynamical process.

So what may happen to our bad cubic as time evolves?
The discriminant of plane cubics has alas even degree
$3(m-1)^2=12$ for $m=3$ (or more generally when $m$ is odd) so
that we cannot infer presence of a singular curve in the pencil
for basic degree reasons. Yet there is surely a deeper argument
either like Klein-Marin (1876--1988 \cite{Marin_1988}) or via
Poincar\'e's index formula (1885) prompting the existence of a
connected curve in any pencil of cubics. Cf. e.g.
(\ref{nodal-cubics-8-many-in-a-pencil:lem}).

Accordingly two scenarios may occur when the bad cubic is dragged
along one of the two possible
sense along the real locus of the pencil:

SC1.---The bad cubic has its oval coalescing with its pseudoline.

SC2.---The bad cubic sees its oval shrinking to a solitary node.

Of course SC2 seems unlikely since the oval of the bad $C_3$
passes through the 6 inner points so a shrinking looks impossible
at least in the near future of $C_3$. So SC1 is the first thing to
occur when the bad cubic is propagated along the deep pencil.

A qualitative picture
(without high precision tracing instrument)  may give something
like Fig.\,\ref{Total-qualitative:fig}a showing the coalescence of
the oval of the bad (black) cubic with its pseudoline via
transition through a nodal cubic (in red). On tracing naively the
next
lilac curve one seems to get a corruption with B\'ezout as the
lilac curve seems intersecting 4 times the horizontal line. This
is fairly naive and there must be ways to avoid such a trivial
accident.

Another optical illusion is the following. On looking
Fig.\,\ref{Total-qualitative:fig} one may get the impression that
in the transition from the red curve to the lilac one along the
segment $A,B$ the cubic  must necessarily split off the line $A,B$
(and accordingly a so-called residual conic $C_2$). If so, then
the 6 remaining (assigned) basepoints have to lie on the residual
conic $C_2$ which intersects $12+2=14$ times the $C_6$, since 2
bonus intersections are forced with $N$ (by dichromatism).
(Note also to complete the argument that none of the 3 (assigned)
basepoints can be aligned as then we get $6+2=8>6$ intersection of
$C_6$ with a line.)

However on zooming (violently) the segment $A,B$ one arrives at
Fig.\,\ref{Total-qualitative:fig}b showing a transition from red
to lilac by an undulating family of (qualitative) cubics
respecting B\'ezout (at least as far as the intersection with
line $A,B$ is concerned). During this undulation no splitting off
of a line is forced.

\begin{figure}[h]
\centering
\epsfig{figure=Total-qualitative.eps,width=122mm}
\vskip-5pt\penalty0
  \caption{\label{Total-qualitative:fig}%
  Naive prelude to a dynamical approach to Rohlin-Le~Touz\'e's total
  reality assertion for sextics} \vskip-5pt\penalty0
\end{figure}

A similar depiction could settle the pseudo contradiction with
B\'ezout of Fig.\,a (involving the line $C,D$ and the lilac
curve). This is suggested on our loose picture Fig.\,c. In reality
nobody tell us that the picture is like this,  being possibly
rather like Fig.\,d or even different. It is clear at this stage
that the argument becomes much involved if possible to complete at
all. Philosophically the drawback of our strategy is that it is
indirect (by contradiction). One could dream of a direct argument,
but this surely requires different ideas.

Our indirect argument requires solid consolidations perhaps by
enumerating carefully the several Morse surgeries implied by the
evolution. By genericity those could be assumed of elementary type
(uninodal curves only). Also
during the time the oval of the bad cubic stays an oval,  its
expansion seems, by convexity,
confined within the fundamental triangle of Fig.\,a. Finally, in
the limit when we encounter the first nodal curve of the pencil,
the inside of this loop (which is also the geometric limit of the
insides of the ovals past the bad cubic) has also to be convex and
therefore contained in the fundamental triangle of Fig.\,a, plus
its companion (forming a ``David star'').
\end{proof}

\subsection{Another strategy to Rohlin-Le~Touz\'e's phenomenon}

[21.02.13] (based on hand-notes of the past 3 days).

We consider again a sextic $C_6$ of type $\frac{6}{1}2$. We
distribute $8$ points on the empty ovals of the $C_6$. The
phenomenon in question claims that the pencil of cubics through
those 8 points is totally real, i.e., each real curve of the
pencil cuts only real points on $C_6$.

First one notices that each curve of the pencil (denoted $\Pi$)
cuts at least 16 points on the empty ovals. (Here and in the
sequel, intersections are always counted by multiplicity.) Denote
by $N$ the nonempty oval of $C_6$.

\begin{lemma} If the pencil $\Pi$ is not totally real, then
it contains a bad cubic $C_3$, i.e. such that $C_3(\RR)\cap
N=\varnothing$.
\end{lemma}

\begin{proof} If all cubics of $\Pi$ cut $N$ then all have 2
extra intersections located on $N$, and so the pencil is totally
real.
\end{proof}

Such a bad cubic is necessarily smooth, because singular cubics
are either connected or have a solitary node, but in the latter
case the real pseudoline connects inner and outer points so an
interception of $N$ is forced by continuity.

Assume (by contradiction) that there is a bad cubic in $\Pi$. One
idea is  to look at the future of this bad conic along the pencil
$\Pi$. One can introduce the projection induced by the pencil as
the map
$$
\pi\colon C_6 \to \Pi
$$
taking a point of the curve to the unique curve of the pencil
passing through it. $\pi(N)=:G$ is the set of good conics, whose
complement is $B$ the set of bad conics. Under our assumption that
$B\neq \varnothing$, it is clear that $G$ is a (compact) interval
in the circle $\Pi$. In fact as the pencil is defacto nearly total
with 16 real intersections (over the 18 maximum permissible),
the map $\pi\colon N\to G$ is two-to-one. Hence given $C_3$ our
(initial) bad cubic we may let it degenerate toward one of the 2
extremities of the interval $G$ along 2 pathes consisting only of
bad cubics safe for their
extremities. As only 2 extra intersections are possible, this may
occur in 3 fashions only (by a simple continuity argument):

(I) Inner touch: the oval of $C_3$ inflates inside $N$ and
ultimately touch it from the interior.

(D) Double touch: the oval of $C_3$ inflates from inside and
collides with the pseudoline of $C_3$ on a point of $N$.

(O) Outer touch: the pseudoline of $C_3$ touches $N$ (necessarily
from outside) while the oval stays inside $N$ disjoint from it.

Further when dragging the curve along $\Pi$,  at some stage (first
touch or contact) two real points eventually appear on $N$, and
subsequently move apart along $N$ (without possible return by the
property of linear systems or holomorphic maps) to merge again on
the opposite first contact of $C_3$ with good cubics. This looks
attractive but is probably only a first step toward
a contradiction. In fact a simple picture
(Fig.\,\ref{Tot1:fig}a) shows that such a scenario is perfectly
 permissible, topologically at least.

\begin{figure}[h]
\centering
\epsfig{figure=Tot1.eps,width=122mm} \vskip-5pt\penalty0
  \caption{\label{Tot1:fig}%
  A flexible pencil of cubics and an obscure contradiction
  via the method of barrages} \vskip-5pt\penalty0
\end{figure}

This figure suggested another idea as follows. While the above
picture (Fig.\,\ref{Tot1:fig}a) is topologically
legal, the thick traced blue curve seems to violate B\'ezout upon
tracing a line through its node intercepting it 4 times. A tactic
would be to argue by Poincar\'e's index formula (applied to the
inside of the egg $E$, i.e. the unique oval of our bad cubic
$C_3$) that there is necessarily such a nodal curve in the pencil
(with node located inside $E$), and by some messy combinatorial
argument such a curve would necessarily corrupt B\'ezout,
heuristically because it has to visit too many points forcing
high-contortion like the thick blue curve above.

As to the Poincar\'e argument, look at the inside $E^{\ast}$ of
the egg $E$ with the (mildly singular) foliation induced by $\Pi$
and double it to a sphere, $2E^{\ast}\approx S^2$. We see on the
boundary (doubled!) 6 foyers of index $+1$ (locally like the
pencil of lines through a point). A priori there could be centers
(locally like concentric circles)
with index $+1$
and arising
from a solitary cubic. Finally nodal cubics (with non isolated
ordinary singular point) contributes for (hyperbolic) saddles
(locally
like the levels of $x^2-y^2$) which are of index $-1$.
Poincar\'e's index formula tells
the sum of indices
being equal to the Euler characteristic of the manifold. Hence
$6-2 S\ge \chi(S^2)=2$, where $S$ is the number of saddles, and we
deduce that there is at least two of them inside $E$. (As $E$ is
smooth they cannot be located on the boundary of $E^{\ast}$.)

The above programme sounds good (albeit requiring alienating
combinatorics!) until the moment, one realizes that a nodal  cubic
is able to salesman-travel through the 6 basepoints on $E$ without
being contorted. Remember at this stage that the 8 (assigned)
basepoints of $\Pi$ determine (by B\'ezout) a 9th one, which for
vibratory reasons has to be outside $N$ (otherwise a slight
perturbation of $E$ would intercept an odd number of times  $E$,
violating the depiction or if you prefer homological intersection
modulo 2.) Of course if the 9th unassigned basepoint lands on $N$
then total reality is evident.

Let us now depict such a nodal cubic able to visit the 6 inner
points without being contorted (Fig.\,\ref{Tot1:fig}b). [22.02.13]
On the latter all the (rational) nodal cubics occurring in the
pencil have  relatively decent looks. To formalize the lack of
contortion of such a cubic one can uses the pencil of lines
through the node which cuts a group of 3 points with 2 of them
statically monopolized by the node while the 3rd moving along the
curve. So when one looks from a nodal cubic from its node one
always see at most (an in fact exactly one) point forced to be
real. Our idea was that  at least one of the nodal cubics (ensured
via Poincar\'e's index formula) would be contorted, i.e. violating
this tightness of nodal cubics, yet our
Fig.\,\ref{Tot1:fig}b
gives little hope to complete this.

Another strategy also jeopardized by the above pictures
(Fig.\,\ref{Tot1:fig}) is that there ought to be always a cubic of
the pencil which is dichromatic in the sense that the 6 inner
points (black colored) and 2 outer points (white colored) are
lying mixed on some suitable cubic of the pencil $\Pi$ with the 2
white points separating the  collection of all 6 black points. If
so is the case,  4 extra intersections are gained on the nonempty
oval $N$, and B\'ezout is violated. Perhaps this strategy is the
right one but requires more geometrical argument \`a la
Le~Touz\'e.

Yet another idea is that by using the nodal cubics of the system
we may infer that the outer basepoints are strongly stretched
apart, while by contrast Le~Touz\'e's lemma (chromatic law for
conics, cf. \ref{LeTouze-Gabard-Hilfssatz:lem}) forces them to be
much condensed, in the sense of not being separated by any
triangle through any triplet among the 6 inner points. Remember
(from Le~Touz\'e's Sec.\,\ref{LeTouze:sec}) that if a separation
occurs then the conics through the 3 corresponding inner points
and the 2 outer points is dichromatic (with the 2 white points
separating the 3 black points) so that the corresponding conic has
$10+4$ intersections with $C_6$ (violating B\'ezout).

Another idea is that since the 9th basepoint is outside $N$ (for
the vibratory reasons already explained), all cubics of the pencil
have to oscillate about those 3 points. This is perhaps
incompatible with the  tightness of (rational) nodal cubics.

Yet another idea was that the 9th basepoint of our pencil $\Pi$
(almost canonically assigned to the $C_6$) has to land inside $N$
and this would contradict the vibratory properties of a bad cubic.
However this miraculous property looks logically much stronger
(i.e. not logically equivalent) to the Rohlin-Le~Touz\'e total
reality claim, so that this is perhaps not a realistic strategy,
at least we were not able to implement it.

Maybe what is required is an avatar of the chromatic law for
cubics instead of the version for conics
(Lemma~\ref{LeTouze-Gabard-Hilfssatz:lem}). The logics would be as
follows. Trace the ``diamond'' of all $\binom{6}{2}=15$ lines
through the 6 inner points. By the chromatic law for conics, this
diamond does not separate the 2 outer points. So by a hypothetical
chromatic law for cubics it could follow that the pencil of cubics
through the 8 points is dichromatic, i.e. contains a dichromatic
cubic. The latter would overwhelm B\'ezout. Of course all this if
it works should use the assumption of a bad cubic which implies an
hexagonal (convex) distribution of the 6 inner points on the
egg-shaped oval $E$ of the bad cubic.

The bottom foliation (Fig.\,\ref{Tot1:fig}c) extends the
right-part (Fig.\,b) of that figure (while changing slightly the
colorimetry), and shows again that there is no topological
obstruction in the large. So it seems that the contradiction (if
it exists, i.e. if Rohlin-Le~Touz\'e are right) must really
involve some deeper geometry (presumably at the level of B\'ezout,
or maybe Cayley-Bacharach, Jacobi, etc.).

Of course our global picture shows some new nodal cubics which are
highly contorted, for instance the thick-blue curves. Reminding
 tightness of nodal cubics, the inside of the loop of that cubic
must be convex. By the {\it loop\/} of a nodal cubic we mean the
unique arc joining the node to itself via the half which is
null-homotopic in $\RR P^2$. This being said, we
may from the node of the blue-thick curve $B_3$ trace a
rectilinear segment
joining the top-point of the loop of $B_3$, and lying entirely in
the inside of the loop of $B_3$ (cf. dashed line on
Fig.\,\ref{Tot1:fig}c). This segment which is linear (despite the
appearances!) cuts for topological reasons at least 4 times the
lilac-colored cubic, hence B\'ezout is corrupted, and perhaps the
Rohlin-Le~Touz\'e theorem is nearly proved.

Let us  formalize the argument. Consider the pencil $\Pi$ of
cubics through  8 basepoints (injectively) distributed on the 8
empty ovals of the $C_6$. If $\Pi$ is not totally real, there is a
bad cubic $C_3$ whose real part is disjoint from $N$, the nonempty
oval of $C_6$. Denote by $E$ the unique oval of this bad cubic
which is necessarily smooth. By applying Poincar\'e's index
formula to $E$ (or the double of its inside) we infer that there
is at least 2 saddle points inside $E$. On applying it to $\RR
P^2$ we infer that there is at least 8 saddle points on the whole
projective plane. Such saddle points correspond to nodal cubics
(and perhaps it is convenient to assume some genericity of the
pencil after dragging slightly the 8 assigned basepoints).

By B\'ezout recall that $N$ is at most intercepted twice by each
$C_3$ of the pencil, and actually exactly twice for each cubic
which is not bad (i.e. which intersects $N$). Accordingly we get
an involution with 2 fixed points on $N$ (namely the first contact
of the bad cubic with good conics). This permits to fold the
boundary of $N$ to get a topological sphere. Now depending on
whether the first contact with good cubics are inner touch, or
outer touch, or double touch (as discussed earlier) we get by the
folding different type of singularities. Specifically an inner
touch induce no singularity, and so do a double touch, while a
outer touch induces a center (do some simple local pictures to get
convinced).

So we may apply Poincar\'e inside $N$ (folded) and deduce that
there is at least 4 saddles inside $N$ (in accordance with the
picture) and perhaps at most 6 saddles (compare picture or think
hard). All this to ensure that there is at least one saddle
outside $N$ and the corresponding nodal cubic ought to have always
a loop enveloping 4 transverse arcs of another nodal cubic with
inner node (as on the picture). Remind that the existence of a
lilac-colored cubic seems to be forced by B\'ezout. If all this
works then we are finished and the general case is so-to-speak
always reducible to the one depicted.

Of course we need to be slightly formal (and clever) for instance
by defining the concept of a barred-pair of nodal cubics (or
barrage for short). This is a pair of nodal cubics such that one
of them appears 4 times inside the loop of the other. (For an
example cf. again the thickest curves of Fig.\ref{Tot1:fig}c.) It
remains then to show that such a barred pair always exists, which
requires some abstract self-confidence in combinatorics or a long
discourse. Note on the picture (at least) that if we consider
instead of the thick lilac curve the red one then there is also a
barrage consisting of 4 disjoint arcs inside the loop of the blue
curve. Hence the proof could decompose in the following 2 steps:

Step 1.---Show that there is always  a nodal cubic whose loop
visits all the 8 assigned basepoints.

Step 2.---Show that there is always another nodal cubic forming a
barrage w.r.t. the nodal cubic of step 1, i.e. which appears in
the inside of its loop as 4 pairs of transverse arcs joining the 8
basepoints in pairs.

This is perhaps a universal property of pencil of cubics (or maybe
valid only in our special situation, where the 8 basepoints are on
a $C_6$, with six of them  hexagonally distributed of the convex
egg of the bad cubic $C_3$). Universality would be better as then
the proof could be simpler, but this looks too optimistic for in
that case pencil of cubics would just not exist. Further if Step~1
looks too hard, one could imagine other types of barrages like the
one depicted on the 3rd row of Fig.\,\ref{Tot1:fig}d.

We hope that [all] this [mess] can be made clearer and perhaps
there is a simpler argument (maybe Le~Touz\'e's proof).

Albeit difficult to make formal the above proof (if it is one!)
shows the special r\^ole played by nodal cubics in the pencil
which have lowest complexity from the viewpoint of algebraic
geometry. Those are perhaps the unique ``br\`eche par laquelle on
puisse entrer dans une place r\'eput\'ee jusqu'ici imprenable''.

\subsection{Doubling, Satellites and total reality}
\label{satellite-total-reality:sec}

[23.02.13] As  discussed at length, Rohlin-Le~Touz\'e's theorem is
somewhat elusive to prove but let us assume it to be correct. Why
is it so important? Why is it fairly difficult to prove? How does
it generalize? As a last remark we note that any proof using the
bad cubic tends to be indirect, and this makes any proof a bit
frustrating. One could dream of a direct proof using maybe the
fact that any cubic of the pencil is dichromatic in the sense of
having both inner and outer points one the same component of the
cubic. This would give a direct proof but of course still much
remains to be justified.

Though quite unable to complete the proof, we may try to speculate
of what comes next, and what is the true phenomenology governing
such phenomena of total reality.

According to Ahlfors theorem (1950 \cite{Ahlfors_1950}) what is
behind total reality is basically the orthosymmetric character of
the curve. More concretely (or in the spirit of Rohlin 1978
\cite{Rohlin_1978}), total reality seems to be sometimes forced by
the sole knowledge of the real scheme. For instance, we have the
prototype of the deep nest of depth $k$ and degree $2k$ which is
totally real under a pencil of lines. The point here is that
topology forces so many intersections as algebra permits whence
total reality. Idem for a quadrifolium nest consisting of 4 nests
of depth $k$ and degree $4k$ which is total under a pencil of
conics assigned to pass through any 4 points distributed in the
deepest ovals.

Modulo technicalities, some higher intelligence should be able to
perceive total reality of Rohlin's sextics with the same ease as
in the above two examples. Of course certain aspects changes
radically, like the 9th unassigned basepoint, as well as the issue
that cubics are possibly disconnected, concomitantly with their
irrationality, or positive genus of the underlying Riemann
surfaces.
Is this a sufficient reason to mistrust the ubiquitousness of the
phenomenon of total reality, say as (partially) evidenced by
Ahlfors theorem at the abstract level?


Typical to the basic cases of total reality---sweeping of deep
nests via pencil of lines through a deep center of perspective, or
the vision of
quadrifolia through conics---is some concentric paradigm, namely
an infinite series of species totally real under the same pencil.
So one can start from a conic and imagine its unique oval
(unifolium) doubled, then tripled, etc., and so we get the series
of deep nests. The same ``satellitosis'' occurs by starting from
the quadrifolium quartic and doubling each of its ovals in a tube
neighborhood, to get an octic totally real, a twelve-tic, etc.

\begin{defn} {\rm Given a real scheme $S$ (of degree $m$) with only ovals
(=nullhomotopic curves) we may abstractly define its {\it $k$th
satellite\/} by
replicating each oval up to a certain multiplicity $k\ge 1$ and
get so the scheme $k\times S$ of degree $km$. (In Rohlin's sense,
a {\it real scheme\/} is primarily an isotopy class of embedding
of a disjoint union of circles plus some integer $m$ given in the
background memory, the so-called {\it degree\/} of the scheme.)}
\end{defn}

Note that this abstract operation can be aped algebraically just
by taking an equation of even degree  realizing the scheme $S$ (we
assume this to be possible) and then taking $k$ nearby levels
close to zero while perturbing the union to get a smooth algebraic
curve realizing the $k$th satellite schemes. This makes sense
because the sign of an even degree form is well-defined.

In particular we can take Rohlin's sextic scheme $\frac{6}{1}2$
and double it (second satellite) to get the scheme $2 \times
\frac{6}{1}2$ of degree 12, or triple it, and so on. It seems
clear that this satellite is totally real under the same pencil of
cubics as in Rohlin's (unproven) phenomenon. This is evident when
the satellite is realized as a small algebraic perturbation of
parallel levels, because in that case we have on the original
sextic curve a foliation transverse to the real locus [[10.04.13]
this is a bit sloppy but nearly true], and transversality is
topologically stable. So we get an infinite series of curves of
type~I (as forced by total reality), and it is likely that not
merely the algebraic satellites are of type~I but all the curves
belonging to the schemes. This would be the case if the schemes
were known to be rigid, i.e. each forming a unique rigid-isotopy
class (\ref{rigid-scheme:defn}). More pragmatically, the fact that
total reality is exhibited by a synthetic procedure (namely by
assigning 8 basepoints on the 8 empty ovals of the $C_6$ of
Rohlin's type or over any satellite of this scheme) makes that we
have some robust recipe ensuring total reality. So it is likely
that Rohlin-Le~Touz\'e's theorem implies the following:

\begin{lemma} {\rm (Hypothetical)}.---Any $k$th satellite of Rohlin's schemes of degree $6$ (they are 2
of them namely  $\frac{6}{1}2$ and $\frac{2}{1}6$) is again total
under a pencil of cubics and so of type~I. In particular the $2$nd
satellites of Rohlin's schemes are schemes of degree $12$ which
are of type~I.
\end{lemma}

With some good faith (or pessimism) one could
fear that this implies a corruption of Rohlin's maximality
conjecture (type~I implies maximal). The idea would be to take a
fairly complicated configuration  of 6 ellipses and smooth it \`a
la Brusotti to get a curve whose scheme enlarges $2\times(
\frac{6}{1}2)$. Fig.\,\ref{Tot2:fig} includes inconclusive
attempts along this naive tactic.

\begin{figure}[h]
\centering
\epsfig{figure=Tot2.eps,width=122mm}
\vskip-5pt\penalty0
  \caption{\label{Tot2:fig}%
  Naive pseudo-counterexample to Rohlin's maximality conjecture} \vskip-5pt\penalty0
\end{figure}

Adhering to the opposite attitude, the 2nd satellite of any one of
both Rohlin's $6$-schemes are $12$-schemes (denoted $2\times R$ or
just $2R$, cf. Fig.\,\ref{Tot2:fig}a) which are totally real in
some geometric way (pencil of cubics through the 8 empty ovals).
Hence it is likely that those schemes cannot be enlarged without
corrupting B\'ezout. More precisely assume a real $12$-scheme $S$
enlarging $2R$, then select in $S$ a replica of $2R$ and construct
the allied total pencil. Let pass a curve through one of the
deleted oval of $S$, and get a corruption with B\'ezout.

More generally if a scheme is of type~I, one may expect its
representing curves  to be totally real under a pencil of curves
in some geometrically controlled way. This posits both a
concretization of Ahlfors theorem as well as an extension of
Rohlin-Le~Touz\'e's theorem. The byproduct  would be a general
proof of  Rohlin's maximality conjecture. At this stage, we
confess to have first understood the full swing of Rohlin's
prophetical allusion when formulating his maximality conjecture:
``there is much to say in its favor'' (cf. Rohlin 1978
\cite{Rohlin_1978}).
This idea will be developed in the next section.

We can also look at the 3 possible $M$-sextics permitted by
Gudkov's classification and take their satellites to get schemes
of type~I, actually total under a pencil of quartics via
(\ref{total-reality-of-plane-M-curves:thm}). The emerging
philosophy is that the phenomenon of total reality should be
stable under satellitoses and possesses a series of minimal (or
primitive) models in each degree. Pencils of lines correspond to
deep nests. Pencils of conics (with 4 real basepoints) correspond
to the quadrifolium quartics and its satellites. Pencil of cubics
have two minimal models with  Rohlin's sextics. Pencil of quartics
have (at least) 3 minimal models given by the $M$-sextics (cf.
Theorem~\ref{total-reality-of-plane-M-curves:thm}), etc. All this
is quite vague and need perhaps strong correction, but our
intention is  to suggest the idea of a big tower of total pencils,
of which the Rohlin-Le~Touz\'e phenomenon should  just be one of
the very first cornerstones supporting a big cathedral. Admittedly
the latter may reach such altitudes, that its higher structure is
still completely dissimulated behind the clouds.

The motive behind total reality seems to be a topological
predestination forcing reality of all intersections. So the
phenomenon ought to be fairly robust. Now if we are given a scheme
of type~I, then any curve representing it is totally real by
Ahlfors theorem.
[$\star$ Not even obvious!] Yet to attack Rohlin's maximality
conjecture (RMC) we need more namely total reality forced by
topological reasons. This amounts essentially to a synthetic
knowledge a priori of the location of the basepoints.
In this case let us say that the scheme is photovoltaic, more
precisely:

\begin{defn}
A real scheme is photovoltaic (PV) if there is a canonical
recipe($\approx$algorithm$\approx$Turing machine) exhibiting a
total pencil of curves on it. When the recipe is as simple as
saying ``by assigning basepoints on the empty ovals'' of any
representing curve of the scheme, the scheme is said to be
photographic.
\end{defn}

We have ``photovoltaic'' implies ``type~I'', and even
``photovoltaic'' implies ``maximal''.
Of course the problem is that our ``canonical recipe'' is poorly
defined, but one may just understand some algorithm. For instance
Rohlin's $6$-schemes are photographic by the Rohlin-Le~Touz\'e's
theorem, while the $M$-schemes of degree 6 are photovoltaic since
there is an algorithm to construct a total pencil via some
residual series (cf.
Theorem~\ref{total-reality-of-plane-M-curves:thm}). [$\star$
But compare also (\ref{Le-Touzé-extended-in-odd-degree:scholium})
showing that, in odd degree at least, there is a more concrete
recipe for the total reality of $M$-schemes.] One chance to go
around the conceptual difficulty of the ill-posedness of our
definition would be the following miracle:

\begin{conj} All schemes of type~I which are
not $M$-schemes are actually
photographic.
\end{conj}

[$\star$ Again in view of
(\ref{Le-Touzé-extended-in-odd-degree:scholium}) it is likely that
$M$-schemes have not to be excluded.]

This is true for sextics (granting the Rohlin-Le~Touz\'e theorem),
and deserves to be investigated in general.

If the conjecture is true in general,  type~I implies photovoltaic
(by virtue of Theorem~\ref{total-reality-of-plane-M-curves:thm}),
and hence maximal, and Rohlin's conjecture would be settled.

Of course we are using the implication ``PV'' implies maximal. It
looks hard to prove it because ``PV'' is ill-defined, but we
really may avoid this concept since $M$-schemes are automatically
maximal (Harnack 1876), while the other are photographic (by the
conjecture) so that Rohlin's maximality conjecture follow form
the:

\begin{lemma} (Hypothetical!!!)
If a scheme is photographic then it is maximal.
\end{lemma}

\begin{proof}~[$\star$ too vague!]
Suppose by contradiction that $S\subset E$ is an enlargement of
the photographic scheme $S$. Choose $E_m$ a real curve
representing the scheme $E$, and select a sublocus $\Sigma_m$ of
$E_m$ realizing the scheme $S$. Alas we loose algebraicity doing
so. However this sublocus $\Sigma_m$ is total under a pencil of
curves with basepoints assigned (say) on the empty ovals of
$\Sigma_m$ for ``robust'' topological reasons. This is to mean
that any curve of the pencil of $k$-tics cuts $k\cdot m$ points on
$\Sigma_m$ for topological reasons (e.g., like for the deep
nests). Then we could conclude to a contradiction with B\'ezout by
letting pass a curve of the pencil through the extra oval of
$E_m$.
%
%
%
%
\end{proof}

The above clumsy proof imposes a refinement of the definition
making the above lemma true with ``photogenic'' instead of the
``photographic'' assumption.

\begin{defn}
An $m$-scheme is photogenic if any (differentiable, or real
analytic) curve $\Sigma$ representing it admits a ``total'' pencil
of $k$-tics such that each curve of the pencil cuts at least
$k\cdot m$ points on $\Sigma$. It may even be assumed that the
basepoints of such a pencil are assigned in the insides of the
empty ovals.
\end{defn}

\begin{figure}[h]
\centering
\epsfig{figure=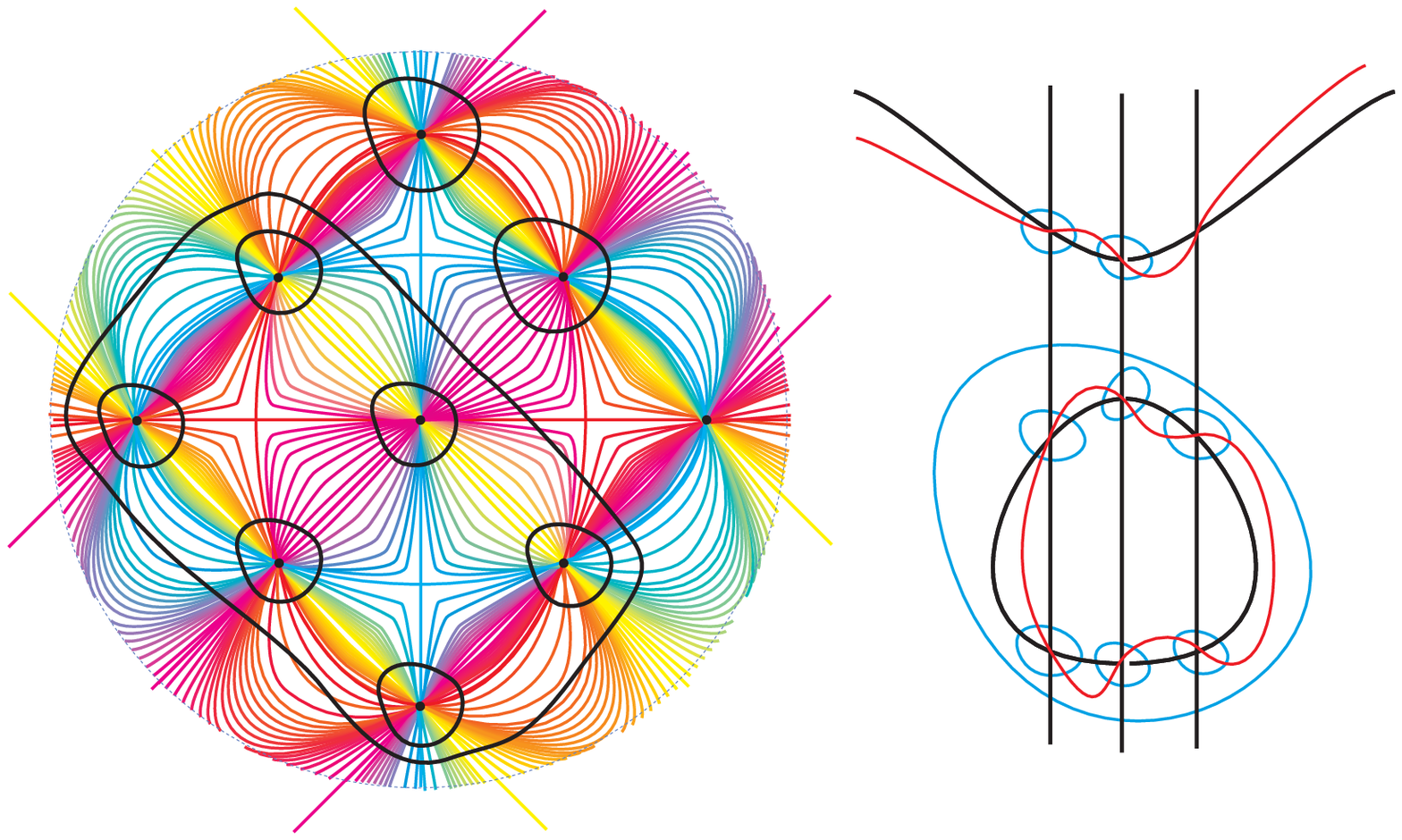,width=122mm} \vskip-5pt\penalty0
  \caption{\label{Tot3:fig}%
  Against a purely topological proof of
  Rohlin-Le~Touz\'e's theorem} \vskip-5pt\penalty0
\end{figure}

This ``photogeny'' is a violent evasion outside the
algebro-geometric realm, yet the deep nest as well as the
quadrifolium of depth $k$ are photogenic schemes in this sense. It
would be interesting to know if the Rohlin's $6$-schemes are
photogenic amounting to say that Rohlin-Le~Touz\'e's theorem
accepts a purely topological proof. This seems already quite
unlikely, and the right part of Fig.\,\ref{Tot3:fig}
 supplies a simple counterexample.
Here we consider a smooth cubic curve $C_3$ and triad of lines
$D_3$ (both black colored). We trace (in the smooth category) the
blue curve $C_6$ realizing the $6$-scheme $\frac{6}{1}2$ with 8
small ovals about 8 of the 9 intersections of the cubics, plus one
large oval enveloping the oval of the smooth cubic. Now the pencil
of cubics spanned by $C_3$ and $D_3$ may be interpreted as the
pencil of cubics assigned to pass through the insides of the 8
empty ovals of the flexible curve $C_6$, yet it fails to be total
as $C_3\cap C_6$ has only 16 points (and not $18$ the product of
their degrees).

So we cannot expect to be so naive as to be photogenic. This
relates to the fact that pencil of cubics (or higher order curves)
generally contains disconnected curves. One crude way to ensure
total reality could be to use degenerate pencils lying entirely in
the discriminant and more than that consisting only of rational
curves (forcing via Harnack-Klein or less, like L\"uroth-Clebsch,
or Cayley)
the curve to be connected. However this looks
overspecialized and probably not even suited to detect the
universal orthosymmetry of Rohlin's $6$-schemes.

It remains to clarify several aspects. Is total reality stable
under satellites? In particular is there an infinite series of
examples above  Rohlin-Le~Touz\'e's phenomenon of total reality.
What are the higher order avatars of  Rohlin-Le~Touz\'e's theorem,
and how frequent is the phenomenon? More precisely which schemes
are photographic? This looks of course extremely hard requiring a
highbrow extension of the Rohlin-Le~Touz\'e's theorem. Are
photographic schemes stable under satellites? If yes this is the
trivial part  of an iterative
propagation of each total reality phenomenon. A priori one can
speculate that photographic schemes are quite rare and essentially
exhausted by pencil of lines, conics and cubics. In contrast one
may expect the phenomenon to be
ubiquitous and so frequent that  all schemes of type~I (safe
perhaps some $M$-schemes)
[$\star$ this proviso looks not justified anymore, cf.
(\ref{Le-Touzé-extended-in-odd-degree:scholium})] are
photographic. In that case there is some little chance to tackle
Rohlin's maximality conjecture (the part thereof post-Shustin's
disproof). Alas even that looks difficult. One may also wonder how
frequent are schemes of type~I, again rarity versus abundance is
quite puzzling.

\subsection{Stability of type~I under satellites}

[24.02.13] Are schemes of type~I stable under satellites? The
first case to test is $2\times R$ the 2nd satellite of Rohlin's
$6$-scheme $R:=\frac{6}{1}2$. Of course taking a perturbation of
the double of Hilbert's realization of $\frac{6}{1}2$
(Figs.\,\ref{GudHilb8:fig} or \ref{GudHilb6-12:fig}) it is likely
that we find a dividing curve, and perhaps  Rohlin-Le~Touz\'e's
theorem is sufficiently robust as to imply universally the type~I
of this $12$-scheme. If not it may be  that the $12$-scheme
$2\times R$ is indefinite. A priori curves of degree 12 could be
sufficiently messy as to allow a type~II realization of the
12-scheme $2R$, or in contrast Rohlin-Le~Touz\'e's phenomenon
could be sufficiently robust as to propagate to satellites.

The data of a curve plus a totally real pencil of ``adjoint''
curves is called a flash, and we say that the curve is flashed by
the pencil. If a curve of even degree is flashed by a pencil then
the doubled curve (and more generally its $k$th satellite)
obtained by small perturbation of $k$ concentric levels is flashed
by the same pencil. Note that for an algebraic satellite to be
defined it is convenient to take an affine chart in which the
whole curve is visible. This works certainly for  Hilbert's
realization of Rohlin's $6$-schemes. Hence it is clear that the
12-scheme $2R$ contains a representatives of type~I (hence is not
a scheme of type~II). The question is to decide if this scheme is
of type~I or indefinite.

One idea could be to realize the 8
nests of depth 2 by an octic and then add two ellipses to get
$2R$. However,  passing a (connected) cubic through the 8 deep
nests
creates $4\cdot 8=32>3\cdot 8=24$ many intersections, and B\'ezout
is much overwhelmed. Replacing the octic by a curve of degree 10
is still insufficient $(32>3\cdot 10=30)$.

\begin{figure}[h]
\centering
\epsfig{figure=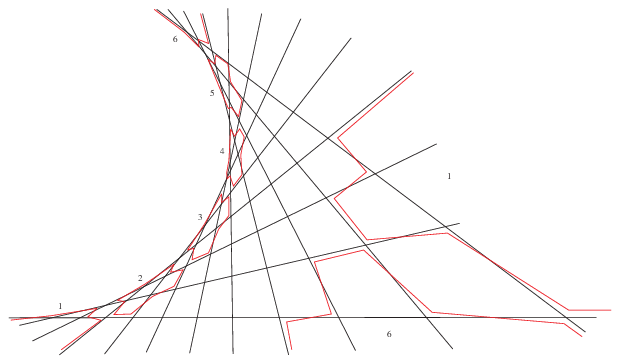,width=122mm} \vskip-5pt\penalty0
  \caption{\label{Tot4:fig}%
  Attempt to construct the doubled Rohlin's scheme via lines} \vskip-5pt\penalty0
\end{figure}

A priori one could hope to find a type~II realization of the
$12$-scheme $2R$ by perturbing an arrangement of 12 lines. This is
a bit messy to depict. The most approaching object we could trace
is shown on Fig.\,\ref{Tot4:fig}. This is rather akin to the
double of the (other) Rohlin scheme $\frac{2}{1}6$, but alas there
is not enough free room left to build the prescribed
configuration.

\begin{figure}[h]
\centering
\epsfig{figure=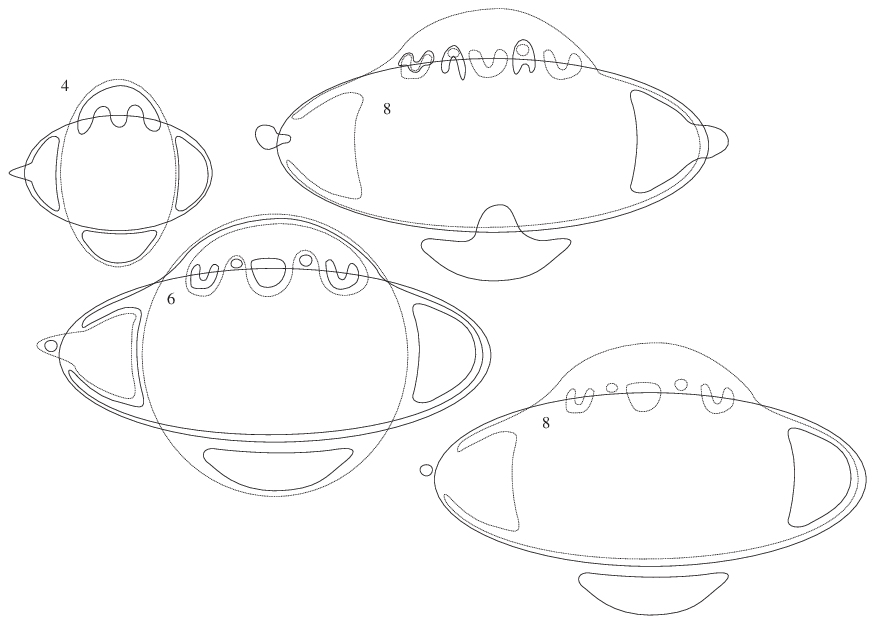,width=122mm} \vskip-5pt\penalty0
  \caption{\label{Tot5:fig}%
  Attempt to construct the doubled Rohlin's scheme via Hilbert} \vskip-5pt\penalty0
\end{figure}

Another idea is to use Hilbert's method, but the latter does not
seem ideally suited for the generation of nest of depth 2
(Fig.\,\ref{Tot5:fig} of very poor quality). Let us shamefully
leave this delicate question, as we sincerely hope that total
reality is ubiquitous (in particular stable under satellites).

\subsection{Satellites of curves of odd degrees}\label{Satellite-odd-degree:sec}

{\it Inserted} [16.03.13].---It seems evident that the
construction of satellites  extends to curve of odd degrees.
Of course there is a slight complication coming from the fact that
the pseudo-line lacks a trivial tube-neighborhood, and so we
cannot replicate so canonically as in the even degree case.
As a simple example consider a cubic with 2 circuits (one oval and
a pseudoline). Doubling its oval and ``doubling'' its pseudoline
will lead to a curve of degree 6 which (for a suitable smoothing)
will be a nest of depth 3, hence again totally real under a pencil
of lines.

By analogy if we look at the next odd degree, namely 5, we have
examples of total reality given by the $M$-quintics (cf.
Le~Touz\'e's Scholie \ref{LeTouze-quintic:scholie}). So when
taking its satellite we are supposed to find a nice example of
total reality in degree 10 for a scheme of the form $(1,6\times
1)$, i.e. 6 nests of depth 2 enveloped in a larger oval. So it is
natural to conjecture that this scheme of degree 10 is of type~I.

If this is possible to prove this is quite interesting because the
scheme in question has $r=13$ ovals which is fairly low in
comparison to Harnack's bound $M=37$, when $m=10$ as
$g=\frac{(m-1)(m-2)}{2}=\frac{9\cdot 8}{2}=9\cdot 4=36$. Of
course, the type~I of this scheme is not covered by the
RKM-congruence for $(M-2)$-curves. So this gives a certain
addendum to Rohlin's (somehow denigrating) remark that the method
of total reality is somehow subsumed to the RKM-congruence (cf.
his remark in 1978 which reads
%
%
``However, all the schemes that we have so far succeeded in coping
with by means of these devices are covered by Theorem~3.4 and 3.5.
[i.e., the congruences]'', compare
(\ref{Rohlin1978-total-reality:quote}) for the integral citation.

Of course  total reality (hence type~I) is also observed for
satellites of the unifolium or quadrifolium having the same
property of being at lesser altitude than $(M-2)$-schemes. Thus
the phenomenon under examination is formally not new, but those
examples being so trivial they were probably not taken seriously
enough. So:

\begin{conj}\label{satellite-of-M-quintic-total:conj}
The scheme of degree $10$ of symbol $(1, 6\times \frac{1}{1})$
(cf. Fig.\,\ref{satellite-of-Harnack's-quintic:fig}a) arising as
the 2nd satellite of Harnack's $M$-quintic with symbol $6 \sqcup
J$ (a unique rigid-isotopy class by Kharlamov-Nikulin) is totally
real under a pencil of cubics, hence in particular of type~I.
\end{conj}

To prove this we use the method of Le~Touz\'e's scholie, namely to
assign the 8 basepoints of a cubics-pencil on the 6 ovals of the
quintic plus 2 on the pseudoline. Then we have 14 intersections
and the last one is forced to reality either by algebra
(Galois-Tartaglia) or topology (M\"obius-von Staudt).

{\it Optional side remark}.---It may be observed that the scheme
in question is not prohibited by Rohlin's formula. Hint:
decomposes the Hilbert tree of the scheme in $x$ and $y$ many
branches of length 2 (so $x+y=6$) which are resp. positively or
negatively charged. By Rohlin's formula we have
$2(\pi-\eta)=r-k^2=13-25=-12$. By the signs-law (cf.
Fig.\,\ref{Signs-law-dyad:fig}) we find $\pi-\eta=x-3y$, and thus
$x-3y=-6$, $x+y=6$. Eliminating $x$ gives $-4y=-12$, so $y=3$ and
$x=3$. Rohlin's equation with signs is thus soluble.

A priori Le~Touz\'e's total reality should adapt to the double.
Formally we assign 6 basepoints on the deep ovals of the 6 nests
and 2 on the maximal oval. Crudely speaking we await $6\cdot
4+3\cdot 2=24+6=30=3\cdot 10$ and total reality would be granted.
However in reality we get less than that on basic topological
grounds. However it is clear that for a small deformation of the
doubled quintic we can expect total reality and the hope is that
this propagates to the full scheme (chamber which a priori is not
even known to be connected!) Maybe there is some deep reason
ensuring total reality like in Le~Touz\'e's argument.

{\it Insertion} [10.04.13].---It also interesting to compute the
mapping-degree of the allied circle map. The number of mobile
points of Le~Touz\'e's series will be $6+12+4=22$, with 1 point
circulating on each of the $6$ deep ovals, 2 on the 6 ovals
immediately surrounding them and 4 moving on the doubled
pseudoline (maximal oval). This degree of $22$ can be compared
with Gabard's bound $(r+M)/2=(13+37)/2=25$ and turns out to be
compatible with it.

A validation of the conjecture
(\ref{satellite-of-M-quintic-total:conj}) could be of interest for
the following reason related to Rohlin's maximality
conjecture(=RMC). If we think globally at the satellite operation
and the arithmetics of small integers
factorized into primes ($1,2,3,4=2\cdot 2,5,6=2\cdot
3,7,8=2^3,9=3^2,10=2\cdot 5,11,12=2^2 \cdot 3=2\cdot 6$, etc.), we
remark that the first nontrivial satellite {\it not\/} totally
real under a pencil of lines or conics truly arises in degree 10.
Degree 9 involves the prime $3$, yet the 3rd satellite of the
cubic with 2 components is merely the deep nest $4\sqcup J$
(totally real under lines). In degree 6 we have indeed the 2 total
realities of Rohlin-Le~Touz\'e yet they are
primitive manifestations (not satellites), and of course in
adequation with Gudkov's classification and Rohlin's maximality
conjecture. Hence degree 10 is the first case where (granting our
conjecture) we get a type~I scheme which possibly is not maximal
(in case we are skeptical about the truth of RMC). Of course this
would be against our own philosophy that Ahlfors has much to say
about Hilbert's 16th. Yet we must keep in mind this eventuality.


Thus  we can ask for a curve of degree 10 enlarging the doubled
$M$-quintic (with 6 bifolia, plus one large oval,
Fig.\,\ref{satellite-of-Harnack's-quintic:fig}a). The naive
construction of Fig.\,\ref{satellite-of-Harnack's-quintic:fig}b
does not even reach more than 4 bifolia nested in a larger oval.
So other techniques of construction are demanded, perhaps Harnack,
Hilbert or Viro.

\begin{figure}[h]
\centering
\epsfig{figure=Biquintic.eps,width=96mm} \vskip-5pt\penalty0
  \caption{\label{satellite-of-Harnack's-quintic:fig}%
  A possible type~I scheme of degree 10 (2nd satellite), and a
  naive construction} \vskip-5pt\penalty0
\end{figure}

If the posited phenomenon
 of total reality holds
true,
and is sufficiently explicit to imply maximality (as nearly
evident by B\'ezout) then our scheme of degree 10
(Fig.\,\ref{satellite-of-Harnack's-quintic:fig}a called say the
closed sextibifolia) would be maximal in the hierarchy of all
schemes of degree 10. Hence this scheme could not be enlarged and
it results, in one stroke, a myriad of prohibitions upon Hilbert's
16th problem in degree 10 (still wide open), basically by virtue
of the sole idea of total reality which goes back virtually to
Riemann's thesis 1851/57, then Schottky, Klein, Bieberbach,
Teichm\"uller, Ahlfors just to name the heros.

So the question looks nephralgic. As a philosophical detail, while
in Hilbert's 16th there is some traditional focus upon curves
maximizing the number of ovals (so-called $M$-curves since
Petrovskii 1933/38), we see here in contrast that the lower the
number of ovals is (for a curve subsumed to total reality) the
stronger will be its prohibitive impact upon the higher stages of
the pyramid. Of course the prototype is the deep nest, but this is
merely the trivial case.

In degree 5 there are also $(M-2)$-curves which are totally real,
typically under a pencil of conics, yet the corresponding scheme
is not of type~I. Its double will be of degree 10 and have a
bi-quadrifolium (double couche) nested in a larger oval.

\subsection{Toward a census of all type~I (totally real) schemes
or at least extension of the Rohlin-Le~Touz\'e's phenomenon
prompted by the Rohlin-Kharlamov-Marin congruence for
$(M-2)$-schemes}
\label{census-and-extension-of-Rohlin-Le-Touzé:sec}

[24.02.13] Another
problem is to list all schemes of type~I, or at least those
totally real under a pencil of curves. We  restrict attention to
even degrees schemes ($m=2k$) and call  {\it order\/} the degree
$d$ of curves involved in the total pencil.

In degree $m=2$ we have a single oval (unifolium, denoted $1$ in
Gudkov's notation) which is total under a pencil of lines $d=1$
with center of perspective chosen inside the oval.

In degree $m=4$, we have the nest of depth 2 (denoted
$\frac{1}{1}=(1,1)$ in Gudkov's symbolism) total under a pencil of
lines, and the quadrifolium $4$ total under a pencil of conics.

For $m=6$, we have the nest of depth 3 (Gudkov symbol $(1,1,1)$)
total for $d=1$, and for $d=3$ the 2 Rohlin's schemes
$\frac{6}{1}2$, $\frac{2}{1}6$ (Rohlin-Le~Touz\'e's theorem), as
well as the three $M$-schemes $\frac{9}{1}1$, $\frac{5}{1}5$,
$\frac{1}{1}9$ of Hilbert, Gudkov, Harnack respectively. The
latter are total for $d=4$ (cf.
Theorem~\ref{total-reality-of-plane-M-curves:thm}).

Before attacking the case $m=8$, some remarks are in order. What
is expected is that sometimes total reality is ensured by
B\'ezout. This is the case when $d=1$, or $d=2$ where we have the
extra knowledge of the connectivity  of all members of the pencil
due to the rationality (unicursality) of all genus $0$ curves.
This property is lost when passing to high-orders $d\ge 3$
pencils. Still we can hope a priori that some other reasons prompt
total reality in some favorable cases, which would be more
elementary than the Rohlin-Le~Touz\'e's theorem for $(m,d)=(6,3)$.
It is with this naive hope that we feel encouraged to adventure in
the jungle of $m=8$. Another vague motivation is that schemes of
type~I are conjecturally maximal, so we are really attacking a
sort of simplified Hilbert's 16th problem
(or rather surfing on its upper envelope).

A last remark is that we
(sentimentally) expect via Ahlfors theorem that the dividing
character of curves is always exhibited by a {\it linear\/} pencil
inducing a map to the projective line (whose complexification is
the Riemann sphere with its standard ``equatorial'' real structure
like our planet Earth).

{\it Optional (non-linear pencils).}---However
any totally real map to more complicated dividing curve suffices
to exhibit the dividing character of the covering curve. So
perhaps we should keep in mind to explore also such nonlinear
pencil. How do they arise concretely is another question. A naive
guess of mine was via a plane cubic $E_3$ with 2 circuits and its
dual curve (variety of tangents). We could hope that any point of
some curve $C_m$ determines a unique tangent to $E_3$ but alas
there are (generally)
6 of them passing through a given point
%
(intersect with the polar curve, a conic here). Perhaps a suitable
adaptation of this idea
leads somewhere.
However it is fairly standard, and we briefly discussed this (in
Part~I devoted to the abstract theory of Riemann-Klein-Ahlfors)
that generally speaking curves mapping to irrational curves of
positive genus have specialized moduli. Hence it is quite unlikely
that we shall gain a general methodology, though plane curves
themselves are modularly confined.

$\bigstar$ For $m=8$, we have the deep nest of depth 4, denoted
$(1,1,1,1)$ total for $d=1$, and the doubled quadrifolium
$\frac{1}{1}\frac{1}{1}\frac{1}{1}\frac{1}{1}=2\times 4$ which is
total for $d=2$.

$\bullet$ We examine next $d=3$ (cubics-pencils). We have then 8
basepoints available assumed all real, and so we look at curves
with this number of empty ovals. Imposing the 8 basepoints on the
empty ovals, we are ensured for twice so many intersections (i.e.
$16$) but this is still much less than $3\cdot 8=24$ (B\'ezout's
upper bound). Total reality looks hard to ensure. Of course we may
envelop our empty ovals by some nonempty ovals. Remember that
there is at most 4 nonempty ovals (as the doubled quadrifolium
$2\times 4$ is total under a pencil of conics hence maximal). So
we may range our 8 empty ovals in 4 groups of ovals and consider
the scheme $\frac{k}{1} \frac{\ell}{1} \frac{m}{1} n$, where
$k+\ell+m+n=8$. If optimistic each of the 3 nonempty ovals
contributes for 2 intersections, and we arrive at $16+6=22$ which
is still less than $24$. Hence it seems nearly impossible to find
a (naive) phenomenon of total reality for $(m,d)=(8,3)$, but this
does not of course exclude the possibility of such a phenomenon
prompted by deep geometric reasons \`a la Rohlin-Le~Touz\'e. Note
also that the case $m=6$ showed that there is no (absolute) total
reality for $d=2$, and so we cannot expect a priori to observe the
phenomenon for each preassigned order and therefore let us skip
the present value $d=3$.

$\bullet$ Assume next $d=4$. We have then $13$ basepoints
assignable. Remember indeed that the space $\vert 4 H \vert$ of
all quartics has dimension $\binom{4+2}{2}-1=\frac{6\cdot
5}{2}-1=14$, so that 13 conditions leave the mobility of a pencil.
So we are again directed toward curves with 13 empty ovals (over
which we shall as usual distribute our 13 basepoints), ensuring so
26 intersections, which is less than the $d m=4\cdot 8=32$
required. Enveloping our ovals in the at most 3 possible nonempty
ovals (if 4 of them we reduce to the doubled quadrifolium) we get
the schemes $\frac{k}{1} \frac{\ell}{1} \frac{m}{1} n$, where
$k+\ell+m+n=13$. So the number of ovals is $r=13+3=16$ and we have
an $(M-6)$-scheme (as $M=g+1=22$ for $m=8$). Perhaps like in the
case $(m,d)=(6,3)$ some of them are total for deep geometrical
reasons.
If lucky, the 3 nonempty ovals creates 6 additional intersections,
so reaching $26+6=32$ B\'ezout's bound, and total reality would be
granted.
 Remember yet that for $(m,d)=(6,2)$ there is no
phenomenon of total reality, at least of the purest form where
only knowledge of the real scheme is required. Repeating
ourselves,  we cannot expect a priori that total reality
prevails for each value of $d$ given a fixed $m$.

$\bullet$ Examine next $d=5$. Then $\dim \vert 5 H \vert =
\frac{7\cdot 6}{2}-1=20$ so that 19 basepoints are assignable. By
the same token, we look at schemes $\frac{k}{1} \frac{\ell}{1}
\frac{m}{1} n$, where $k+\ell+m+n=19$, ensuring $2\cdot 19=38<40=d
m=5 \cdot 8$ intersections. So in fact outside from the 19 empty
ovals it is enough to have one nonempty oval being intercepted to
gain total reality. We list the following candidates:

$\bullet$ $\frac{k}{1}\ell$, with $k+\ell=19$ and $r=20=M-2$.

$\bullet$ $\frac{k}{1}\frac{\ell}{1}m$, with $k+\ell+m=19$ and
$r=21=M-1$ so cannot be total by Klein's congruence $r\equiv g+1
\pmod 2$.

$\bullet$ $\frac{k}{1} \frac{\ell}{1} \frac{m}{1} n$, with
$k+\ell+m+n=19$ and $r=19+3=22=M$, which are $M$-curves.

In the first class of $(M-2)$-schemes we may appeal to the
Rohlin-Kharlamov-Marin congruence (\ref{Kharlamov-Marin-cong:thm})
to
select the serious candidates. (It seems fairly plausible that
this mode of reasoning was the true motivation
behind Rohlin's  Ansatz of total reality for the $6$-scheme
$\frac{6}{1}2$ and its mirror, and that he found his (unpublished)
proof a posteriori of this deeper knowledge.) This RKM-congruence
states that an $(M-2)$-curve of degree $m=2k$ and type~II
satisfies the congruence $\chi=p-n\equiv k^2$ or $k^2\pm 2 \pmod
8$.
This looks a priori undigest, but as merely to be interpreted as a
deviation from Gudkov's hypothesis(=congruence proved by Rohlin),
compare, e.g., the diagrammatic in degree 6
(Fig.\,\ref{Gudkov-Table3:fig}).
 So
when the congruence is violated a scheme of type~I is granted.
Since $\chi=1-k+\ell$ and $k^2=4^2=16\equiv_8 0$ (beware the
overuse of the letter $k$ but no risk of confusion). Recall the
Swiss cheese algorithm for the Euler characteristic $\chi$ (that
whenever we make a hole in the sense of removing a disc, $\chi$
drops
 by one unit, cf. Listing-Klein-von Dyck 1888, etc.).
Starting with the scheme $\frac{19}{1}$, we find
$\chi=1-19=-18\equiv_8 -2$. Then we have $\frac{18}{1}1$ for which
$\chi=1-18+1$ equal to the former plus 2 units, and there is
always an increment of 2. Running through the full list of such
schemes we find that the condition $\chi\equiv_8 4$ ensuring
type~I occurs with periodicity 4 for the following schemes:
\begin{equation}\label{octics-five-examples-RKM:eq}
\frac{16}{1}3,\quad \frac{12}{1}7,\quad \frac{8}{1}11,\quad
\frac{4}{1}15,\quad 20.
\end{equation}
The latter case (of the scheme $20$) seems to disprove our
collective contraction conjecture (\ref{CCC:conj}).

\begin{theorem} \label{CCC-conj-disproof:thm}
{\rm (ERRONEOUS---cf. $\bigstar$ below)}.---The collective
contraction conjecture (of Gabard positing a wild extension of the
contraction conjecture of Itenberg-Viro) is false in degree $8$
already.
\end{theorem}

$\bigstar$ [05.03.13] {\it Corrigendum}.---It follows easily from
the so-called Thom conjecture that the scheme $20$ is not realized
by an algebraic curve of degree 8 (necessarily of type~I if it
existed by the RKM-congruence), compare
Theorem~\ref{Thom-Ragsdale:thm}. So the given argument is not a
disproof of CCC. [07.03.13] In fact a simpler obstruction of this
scheme comes from Rohlin's formula (\ref{Rohlin-formula:thm}), as
$2(\pi-\eta)=r-k^2$ but the left-side is zero, so $r=k^2=16$,
which is no the case. [10.04.13] Further this scheme is also
prohibited (and this was historically the first proof available)
by Petrovskii's inequality (\ref{Petrovskii's-inequalities:thm}),
which reads $\chi\le \frac{3}{2}k(k-1)+1=18+1=19$. This  in
contrast to the proofs via Thom or Rohlin does not use the
dividing character of the curve prompted by RKM. Further it should
be noted that our Thom-style theorem cited above is erroneous in
the generality stated, yet sufficient to imply the present
application as in the case at hand the filled surface is
orientable (since we only glue disc to the half, and so there is
no risk to create a twisted handle like in Klein's bottle). For
definiteness, let us briefly work out the argument. Since $20$ is
an $(M-2)$-curve (of type~I by RKM), we may split the Riemann
surface and fill one half by the $20$ discs bounding the ovals. It
will result a surface of genus $1$, whose fundamental class has
degree $8/2=4$. However Thom conjectured (and Kronheimer-Mrowka,
and others, proved) that the genus is minimized by algebraic
(smooth) curves, hence at least $3$ in degree 4. Since our surface
beats this bound, the real curve $20$ is prohibited.

\smallskip

\begin{proof}~[Outdated,
%
%
but keep in mind the 2nd part of the proof (strangulation
argument), which under CCC would provide another obstruction of
the scheme $20$, of a fairly intuitive character, though hard to
implement with present technology.]
---It seems clear that this 8-scheme $20$
is realized by (a variant of) Hilbert's method.\footnote{Again
this claim is a mistake: an obstruction follows from Thom's
conjecture, meanwhile the theorem of Kronheimer-Mrowka 1994
\cite{Kronheimer-Mrowka_1994}.} (I should still work out this in
some more detail.) The resulting curve is of type~I by the just
cited congruence of Rohlin-Kharlamov-Marin, which is essentially
based either on the deep Rohlin's signature theorem for spin
4-manifolds, or perhaps on K\"ahler geometry in the presentation
of Kharlamov (unpublished?). Alternatively on the model at hand
(via Hilbert's method) the type~I of this curve realizing $20$ may
be checked more elementarily via Fiedler's sense-preserving
smoothing law (elementary surgeries). However the resulting curve
cannot be contracted collectively by shrinking simultaneously all
its ovals to points, for if it could, then $C_8\to C_4 \cup
C_4^\sigma$ would degenerate to a pair of conjugate quartics
obtained by strangulating the Riemann surface $C_8(\CC)$ along all
the separating ovals, and so $C_4\cap C_4^\sigma$ would consist of
20 solitary nodes.  B\'ezout is overwhelmed.
\end{proof}

Moreover the above 5 schemes are avatars of the total reality
claim of Rohlin-Le~Touz\'e for $(m,d)=(6,3)$, i.e. sextics flashed
by cubics, while now octics are flashed by quintics. In both cases
we note the r\^ole of curves of order 3 units less than the given
degree $m$, and one seems being sidetracked to the theory of
adjoint curves \`a la Brill-Noether, etc. Recall indeed that
adjoints of order $(m-3)$  cut out the so-called {\it canonical
series\/} on the given plane curve, and thus there is perhaps some
conceptual reason ensuring total reality of all these linear
systems. This is perhaps the royal road to attack the
Rohlin-Le~Touz\'e's assertion/theorem (and extension thereof
prompted by the RKM-congruence).

In both cases $m=6$ or $8$ we have $(M-2)$-curves which are of
type~I, and swept out by a pencil of order $d=m-3$. The latter
cuts the canonical series of the curve $C_m$ of degree $2g-2$ and
dimension $(2g-2)-g=g-2$ or rather $g-1$?? In view of Gabard 2006
\cite{Gabard_2006} we may expect to find a total morphism of
degree the mean value of $r=M-2$ and Harnack's bound $M$, hence of
degree $M-1$. So we could choose so many points on the
ovals of $C_m$ while putting two of them on the nonempty oval.
It is then hoped that the 2 points situated on the same oval will
{\it dextrogyrate} (i.e. move along one orientation of the oval
without entering in collision) and then total reality is ensured.

While any collection of $M=g+1$ points on a curve of genus $g$
moves in its linear equivalence class, only special collections of
$M-1=g$ points will move but perhaps this is enough to ensure
total reality hence recover the type~I of the above list of
schemes predicted by the RKM-congruence
(\ref{Kharlamov-Marin-cong:thm}).

It is not entirely clear if the Rohlin-Le~Touz\'e's phenomenon is
true in full generality or only for special groups of points (at
least this is the naive intuition coming from the abstract Ahlfors
and Gabard viewpoint). Note in this respect that Rohlin's claim is
a priori less strongly formulated than  Le~Touz\'e's assertion, in
claiming only that a pencil of cubic exhibit total reality and not
that all of them with deeply assigned 8 basepoints are total
(compare Rohlin 1978 \cite[p.\,94]{Rohlin_1978} with Le~Touz\'e's
2013 announcement in Sec.\,\ref{e-mail-Viro:sec}). ($\star$
[10.04.13]---Meanwhile see also the article
\cite{Fiedler-Le-Touzé_2013-Totally-real-pencils-Cubics}.)
It is likely (say by analogy with the trivial case $m=4$) that
Rohlin had in mind the strong assertion of Le~Touz\'e, and that it
is only the extreme compression of Rohlin's exposition that forced
him to his somewhat looser version of the statement.

For $m=8$ we may assign 19 basepoints on the empty ovals of any of
the $(M-2)$-scheme listed above (the scheme $20$ is exceptional in
a sense that remains to be clarified\footnote{[08.03.13]
Fortunately this schemes is not realized.}). Then we can pass a
quintic $C_5$ through these points and one extra point. Choose the
latter on the nonempty oval $N$ of the $C_8$. Another intersection
is created by topology. So we see 21 real points. The residual
intersection of $C_5\cap C_8$ consist of $40-21=19$ points, etc...

Of course the difficulty looks immense but let us postulate the
following avatar of Le~Touz\'e's theorem (i.e. strong form of
Rohlin's total reality claim):

\begin{conj}\label{tot-real-for-octics:conj}
For any octic representing one of the $8$-schemes listed above
(Eq.~\ref{octics-five-examples-RKM:eq}) the pencil of quintics
assigned to pass through the $19$ empty ovals (or even their
insides) is totally real. (Of course the scheme $20$ deserves a
modified statement of which we do not know yet the exact
shape.)\footnote{[06.03.13] This special treatment can be
dispensed as this scheme is prohibited by Thom conjecture, cf.
Theorem~\ref{Thom-Ragsdale:thm}, or more elementarily by Rohlin's
formula.}
\end{conj}

If so is the case we get a total morphism $C_8\to \PP^1$ of degree
$40-19=21$, which is the mean value of $r=20$ and $M=g+1=22$, in
accordance with Gabard 2006 \cite{Gabard_2006}. Of course the
latter affords only weak evidence as its result is subsumed to
high suspicion.

It could be expected (granting Gabard's result as correct) that a
suitable interpretation thereof (at the level of extrinsic
algebraic geometry \`a la Brill-Noether) could supply a proof of
the conjecture at least in the weak form of a special
configuration of $19+2=21$ points two of them being distributed on
the same oval (while dextrogyrating). However ideally we would
like a purely synthetical proof say as elementary as B\'ezout
without incursion of such transcendental philosophers like
Abel-Riemann. This is perhaps possible as some highbrow variant of
the Rohlin-Le~Touz\'e's theorem but remains to be explored and is
quite likely to be extremely elusive. On the other hand it could
be of vital interest that the Abel-Riemann abstract viewpoint may
help to see clearer what happens in the Plato cavern of Hilbert's
16th problem as twisted by Rohlin's synthesis with Klein's
Riemannian viewpoint. In this optimistic scenario we may hope to
get when enlarging further $m$ above $m=8$ an infinite series of
schemes of type~I totally flashed by pencil of $(m-3)$-curves. So
the phenomenon of total reality \`a la Rohlin-Le~Touz\'e  would be
fairly frequent (and in some sense an extrinsic reflection of
Ahlfors abstract theorem).

Perhaps Le~Touz\'e's proof (2013) adapts to give the above the
conjecture, but alas as I do not know yet the details. It can also
be the case that additional difficulties occur while the
combinatorics viz. geometry becomes more involved and the argument
more tedious. The argument is likely to start as follows. As we
have $2\cdot 19=38$ intersections granted, only two are missing to
reach total reality at $40=5\cdot 8$. The sole obstruction is
therefore a bad quintic $C_5$ in the pencil disjoint from the
nonempty oval $N$ of the $C_8$. But this mean that there is an
Abelian differential without zero on this oval, and try to derive
a contradiction.

Remember that the differential has $2g-2$ zeros
(Riemann-Poincar\'e index formula) and we may hope to infer
something. When we look at the trajectories of a (generic)
holomorphic $1$-form we see only hyperbolic saddles of index $-1$
explaining the degree of the canonical class as
boiling down to Poincar\'e's index formula (1881/85). In the case
$m=6$ for the scheme $\frac{6}{1}2$ we have 8 zeros assigned each
creating a companion on the same oval (so 16) and a total of
$3\cdot 6=18$ zeros, in accordance with $2g-2$ for $g=10$ the
genus of sextics. The above looks
a numerical miracle alike, but is not for $g=\frac{(m-1)(m-2)}{2}$
so that $\deg K=2g-2=(m-1)(m-2)-2=m(m-3)$, showing that adjoint
curves of degree $m-3$ are indeed involved in the canonical class.
There is of course a more intrinsic reason allied to the
adjunction formula.

The dream would be that there is some metaphysical/toplogical
principle ensuring total reality on the basis of holomorphic
1-forms which can often be interpreted in terms of incompressible
fluids. Note yet that the argument cannot be too abstract
(else  could imply that all $(M-2)$-curves are of type~I
regardless of the isotopy class, a nonsense already for $m=4$),
 yet ideally it would be as simple as in the case $m=4$ so
$d=m-3=1$ where we really see total reality on the G\"urtelkurve
$C_4$ with 2 nested ovals.  Here we see some potential place of
action for old stuff \`a la Abel-Riemann-Klein ([10.04.13] or
maybe also Thurston's argument in Gross-Harris 1981), but alas
this escaped much from my memory or my curriculum. And still there
is always this objection that the argument should really use the
assumption on the real scheme, and so ought to be more in the
Arnold-Rohlin spirit.

Have we listed all schemes of degree 8 of type~I? Probably not  as
our search was far from exhaustive. It remains of course to list
$M$-schemes (and this is a classical still open problem for some
few exceptional cases). Compare works by Viro, Korchagin, etc. It
is however likely that our list is exhaustive for $(M-2)$-schemes.
([08.03.13] Not even true, as we shall soon see!)

$\bigstar$ What is next? Degree $m=10$ of course. Here we have the
deep nest of depth 5, totally flashed by a pencil of lines. The
quadrifolium $4$ does not give nothing by taking its satellite (as
our $m=10$ is not a multiple of 4). Next we move directly to
adjoint curves of order $d=m-3=7$ (septics). One has $\dim \vert
7H \vert= \binom{7+2}{2}-1=\frac{9\cdot 8}{2}-1=35$. So $34$
basepoints may be assigned freely, and we can force $34 \cdot
2=68<70=7\cdot 10$ nearly all points to be real by assigning the
basepoints to be located on distinct ovals. Again it is
a reasonably folly (by analogy with Rohlin-Le~Touz\'e's assertion)
to expect that under adding an extra nonempty oval enveloping some
of the ovals and if furthermore the RKM-congruence is satisfied
that the resulting scheme (being of type~I) is totally real under
the described  pencil. Precisely we look at the schemes
$\frac{k}{1} \ell $, with $k+\ell =34 $. So we have $r=35$ ovals
and $M=g+1=\frac{(m-1)(m-2)}{2}+1=\frac{9\cdot8}{2}+1=37$ is 3
units above the number of basepoints (no surprise as  the
dimension of the space of curves and the genus both involve the
same binomial coefficient). So we are again in presence of
$(M-2)$-curves. The RKM-congruence says that type~II forces $\chi$
to be either $k^2=25\equiv_8 1, k^2\pm 2\equiv 3,-1\equiv 7 \pmod
8$ so that $\chi\equiv_8 5$ forces type~I. Applying the Swiss
cheese recipe to
 $\frac{34}{1}$, we find $\chi=
1-34=-33\equiv_8 -1 $ and then running through all subsequent
schemes $\frac{33}{1}1$, etc.,  $\chi$ always increases by two
units. So we first met $\chi=5$ for $\frac{31}{1}3$, and find
using fourfold periodicity  the following list of schemes
(potentially totally real):
\begin{equation}\label{RKM-schemes-deg-10:planar:eq}
\frac{31}{1}3, \quad\frac{27}{1}7, \quad\frac{23}{1}11,
\quad\frac{19}{1}15, \quad\frac{15}{1}19,
\quad\frac{11}{1}23,\quad\frac{7}{1}27,\quad\frac{3}{1}31.
\end{equation}
(Like for $m=6$ (but unlike the case $m=8$) the list
is symmetrical under the evident mirror of partnership in the
jargon of Kharlamov-Finashin.) Again we expect the following total
reality:

\begin{conj}
All curves $C_{10}$ of degree $10$ representing any one of the
schemes of the previous display
\eqref{RKM-schemes-deg-10:planar:eq} are totally real under the
pencil of septics assigned to
visit $34$ basepoints injectively distributed among the $34$ empty
ovals of the ten-ics $C_{10}$. Actually the last item
$\frac{3}{1}31$ of the list is prohibited by either  Thom {\rm
(\ref{Thom-Ragsdale:thm})} or by Rohlin's formula
$2(\Pi^+-\Pi^-)=r-k^2=35-25=10$, since $\Pi^+-\Pi^-\le \Pi:=\Pi^+
+\Pi^-=3$.
\end{conj}

{\it Insertion} [10.04.13].---The last scheme of the series is
(alas) {\it not\/} prohibited by Petrovskii's inequality
(\ref{Petrovskii's-inequalities:thm}), but it is by  the strong
Petrovskii estimate of Arnold (1971), cf.
(\ref{Strong-Petrovskii-Arnold-ineq:thm}). This states $p-n^{-}\le
\frac{3}{2}k(k-1)+1$, where $n^{-}=0$ here (negative hyperbolic
ovals), so $p\le 31$ while our scheme as $p=32$. (It may be
useful---if you are better in geography than in arithmetics---to
visualize all this on Fig.\,\ref{Degree10:fig}. Crudely put,
Arnold is as strong as the Ragsdale conjecture.) Further our
assertion regarding the prohibition by Thom is certainly foiled as
there is no reason ensuring orientability of the Arnold surface in
the present case, as we are really attaching a 3-holed disc to the
half of the complexification. (More explanations in
Sec.\,\ref{Thom:sec}.)

\smallskip

It is evident (due to the little arithmetical coincidence between
the genus and the dimension of the curve-hyperspace, plus the
universal validity of the RKM-congruence) that this series of
$(M-2)$-schemes propagates in each degree $m\ge 10$, and so we get
an infinite (nearly tautological) repetition  of the above
conjectures for each even integer $m$. {\it Can all these
conjectures
 be proven in a single stroke?\/} This would be a highbrow
extension of the Rohlin-Le~Touz\'e's theorem. This would give an
infinity and certain abundance to the phenomenon of total reality
as one could have suspected from Ahlfors' theorem. Of course we do
not claim that this will supply an exhaustive list of the
phenomenon, but perhaps it is modulo the operation of satellites.
More precisely:

\begin{conj}\label{primitive-manifestation-of-tot-real:conj}
Any primitive manifestation of the phenomenon of total reality  on
a curve $C_m$ of (even) degree $m$ arises either as an
$(M-2)$-scheme totally real under a pencil of adjoint curves of
order $m-3$ assigned to pass through the empty ovals of $C_m$, or
as a pencil of curves of order $m-2$ if $C_m$ is an $M$-curve (cf.
{\rm Theorem~\ref{total-reality-of-plane-M-curves:thm}\/}).

If not primitive then the scheme is a satellite of either:

\noindent $\bullet$ the unifolium scheme $1$ of degree $2$  total
under a pencil of lines (this gives the series of deep nests which
exist in all degrees $m$), or

\noindent $\bullet$ the quadrifolium $4$ of degree $4$ total under
a pencil of conics, with satellites in all degrees multiples of
$4$, or

\noindent $\bullet$ the other $(M-2)$-schemes of Rohlin
$\frac{6}{1}2$ or $\frac{2}{1}6$ with satellites in all degrees
multiple of $6$, and so on inductively as the satellites of
$(M-2)$-schemes predicted by the Rohlin-Kharlamov-Marin
congruence, or finally

\noindent $\bullet$ as satellites of $M$-schemes of lower degrees
always dividing the given one $m$.
\end{conj}

If this conjecture is true we would have a complete classification
of the phenomenon of total reality for plane curves. This is
surely somewhat premature and probably requires some
slight adjustments to reach more respectableness.

{\it Insertion} [10.04.13].---In particular, one must probably
also takes into account satellites of curves of odd orders, cf.
Sec.\,\ref{Satellite-odd-degree:sec}. For instance in degree
$m=10$, there is probably the 2nd satellite of Harnack's
$M$-quintic playing a r\^ole.

At any rate note that our initial expectation that some phenomenon
of total reality is purely prompted by B\'ezout in a very
primitive way is apparently never borne out. It seems rather that
apart from the satellites of the elementary schemes (unifolium and
quadrifolium) flashed resp. by the trivial pencil of lines and
conics the phenomenon of total reality is at least as hard as
Rohlin-Le~Touz\'e's theorem, but perhaps not much harder. At least
both ought to be connected by the geometry of the canonical
series.

\subsection{More $(M-2)$-schemes in degree 8 of type~I}
\label{RKM-schemes-deg-8-MORE:sec}

[26.02.13] In fact it is clear that even for $m=8$ we have not
listed all $(M-2)$-schemes of type~I for we have only considered
those with one nonempty oval, but we must also consider those with
2, or 3 nonempty ovals. Tabulating a complete list is merely an
exercise of combinatorics.

Geometrically, it may not be essential to
assign
basepoints on empty ovals but some can be located on nonempty
ovals, and we may expect total reality provided the RKM-congruence
is fulfilled. The sole problem is that we then lack some recipe to
assign basepoints, and so the game becomes somewhat obscure
[but quite challenging].

First the RKM-congruence (\ref{Kharlamov-Marin-cong:thm}) can be
more conveniently paraphrased  as:

\begin{theorem}\label{RKM-congruence-reformulated:thm}
{\rm (Rohlin 1978-Kharlamov 197?-Marin 1979)} Any $(M-2)$-scheme
of degree $m=2k$ such that $\chi\equiv k^2+4 \pmod 8$ is of
type~I.
\end{theorem}

[08.03.13] {\it Little Warning}.---There is a minor metaphysical
trouble with this statement. Indeed when $m=8$ (or for larger $m$)
we have the scheme $20$  which satisfies the RKM-congruence, but
which is not realized algebraically as follows either from Thom
(\ref{Thom-Ragsdale:thm}) or from Rohlin's formula
$2(\Pi^+-\Pi^-)=r-k^2=20-16=4$, since $\Pi^+-\Pi^-\le \Pi:=\Pi^+
+\Pi^-=0$ due to the absence of nesting in $20$. There is two ways
to go around this trouble, either add the assumption that the
scheme is algebraic, or interpret Rohlin's definition of the types
of schemes by declaring (usual logical nonsense allied to the
empty~set) that a non-realized scheme
is simultaneously of type~I and type~II (but not of indefinite
type which needs algebraic representatives in both types). Of
course when tracing pyramids, e.g. the Gudkov-Rohlin table
(Fig.\,\ref{Gudkov-Table3:fig}) we ascribe the type~I label
(red-circle) only to those schemes which are of type~I in the
concrete sense that the scheme is {\it algebraic\/} (and all its
realizations are of type~I).

\smallskip

\begin{proof}
The RKM-congruence (Theorem~\ref{Kharlamov-Marin-cong:thm}) says
that an $(M-2)$-curve of degree $m=2k$ and type~II satisfies the
congruence $\chi\equiv k^2 $ or $k^2 \pm 2 \pmod 8$.

So the theorem follows after checking the following basic fact:
{\it an $(M-2)$-curve of order $m=2k$ verifies universally $\chi
\equiv k^2 \pmod 2$.}
%
(From the diagrammatic of pyramids, e.g.
Fig.\,\ref{Gudkov-Table3:fig}, this is the fairly trivial matter
that the rhombic equilateral lattice underlying the pyramid is
adjusted at $k^2$, and one may infer that the claimed congruence
holds more generally on all $(M-2i)$-levels.)

This is easy to prove using the relations $\chi=p-n$, $r=p+n=M-2$
$$
\chi=p-n=(p+n)-2n\equiv_2 (p+n) = r=M-2,
$$
and by Harnack's bound and the genus formula
$g=\frac{(m-1)(m-2)}{2}$ we have
$$
M=g+1=\textstyle\frac{(2k-1)(2k-2)}{2}+1=(2k-1)(k-1)+1=2k^2-3k+2,
$$
whence
$$
\chi\equiv_2 M-2=2k^2-3k \equiv_2 k \equiv_2 k^2.
$$
\end{proof}

Is the converse statement true? Remember that if a scheme of
degree $m=2k$ is of type~I, then it satisfies  Arnold's congruence
$\chi\equiv k^2 \pmod 4$. Hence $\chi\equiv k^2, k^2+4 \pmod 8$,
and the second option leads to type~I, but I do not know if the
first option necessarily implies type~II or  indefinite type. We
know only that this converse holds true for $m\le 6$ by the
Gudkov-Rohlin table (=our Fig.\,\ref{Gudkov-Table3:fig}) which
involve explicit constructions. So the RKM-congruence detects many
$(M-2)$-schemes of type~I, but it is not clear (to me) if it
detects all of them.

\begin{conj}\label{RKM-converse:conj}
An $(M-2)$-scheme of degree $m=2k$ which is of type~I necessarily
satisfies the RKM-congruence $\chi\equiv k^2+4 \pmod 8$.
\end{conj}

{\it Insertion} [11.04.13].---Perhaps an answer can be found in
Rohlin 1978 \cite[p.\,93--94]{Rohlin_1978}, esp. Art. 3.5 and the
end of 3.6, where it seems that an extremal property of the strong
Arnold's inequalities observed by Zvonilov-Wilson prompts type~I
in situation apparently not covered by the congruence. Alas, I had
not yet the time to assimilate this properly, but look forward
with great excitement to do so in the future (after some long
editorial duty).

[08.03.13] Again there is little worry about definitions. For
instance when $m=6$ we can consider the (non-algebraic) scheme
$(1,1,1)6$ which is thus of type~I (in the logical sense but of
course also of type~II), yet with $\chi=(1-1+1)+6=7\neq
3^2+4\equiv_8 5$. So we tacitely assume the scheme of type~I in
the strong sense that it is algebraically realized.

Of course the conjecture \ref{RKM-converse:conj} is true for $m=6$
(look at the Gudkov-Rohlin Table=Fig.\,\ref{Gudkov-Table3:fig}),
which depends upon explicit construction of curves of type~II for
all schemes which are not RKM. Already for $m=8$, the conjecture
seems to demand a menagerie of construction. One could hope that
there is some theoretical argument.

Let us leave this question aside, as we merely want to list
schemes of type~I potentially subsumed to the phenomenon of total
reality.

Let us now tackle the combinatorial aspect of dressing the list of
all $(M-2)$-schemes of degree $m=8$ satisfying the RKM-congruence
(hence of type~I).

We may start with schemes with zero or one empty ovals and list
all of them using the fourfold periodicity as we already did.  Yet
to be more systematic we start with $\frac{16}{1}3$ expand its
Gudkov's symbol as $\frac{16}{1}\frac{0}{1}\frac{0}{1}1$ to be of
the shape $\frac{x}{1}\frac{k}{1}\frac{\ell}{1}m$ and then we
trace a cubical lattice (Fig.\,\ref{RKM-schemes-deg-8:fig}) in
3-space encoding all variations of this symbol for varied values
of $(k,\ell,m)$. To aid visualization it turned useful to ascribe
 colors to the different levels: the ground floor is
orange, the 1st floor is lilac, the 2nd floor blue, the 3rd floor
is cyan, the 4th floor is yellow-green.  As we are interested in
$(M-2)$-schemes we have the relation $x+k+\ell+m+3=M-2=20$. The
crucial point is that when $k$ or $\ell$ increases by one unity,
then one hole is traded against another hole, so that $\chi$ is
left unchanged. In contrast an increment of $m$ reduces $x$ by
one, and so a hole in ``$x$'' is traded against a disc outside, so
that $\chi$ increases by two. Hence as the RKM-congruence is
modulo 8, we have a 4-fold periodicity until reaching the same
value for $\chi$. (This explains the vertical motion along the
cubical lattice.) Those symbols surrounded by dashed lines are
doubloons (non-normalized Gudkov's symbol), yet useful to stop the
combinatorial proliferation. Underbraced symbols are those whose
Gudkov's symbol admits a shorter expression given below the brace
(when enough room is left available).

\begin{figure}[h]
\centering
\epsfig{figure=TotTable.eps,width=122mm} \vskip-5pt\penalty0
  \caption{\label{RKM-schemes-deg-8:fig}%
  (Nearly) exhaustive list of all $(M-2)$-schemes of degree 8
  satisfying the Rohlin-Kharlamov-Marin congruence
  (hence of type~I); for a more complete list cf. Lemma~\ref{RKM-schemes-deg-8-census:lem}} \vskip-5pt\penalty0
\end{figure}

All those schemes  are avatars of the 2 Rohlin's $(M-2)$-schemes
of degree 6 (subsumed to total reality). It is a simple matter to
count them. First collect on the nearby face of
Fig.\,\ref{RKM-schemes-deg-8:fig} all schemes lying in perspective
beyond the  red/thick numbers. Adding them vertically gives the
blue/big numbers on the bottom row, yielding a total of $5 \cdot
11+9= 55+9=64$ schemes. (That this is a power of 2, incidentally
the same as Rohlin's count of all schemes of degree 6 decorated by
types, is probably a mere coincidence.)

This is somewhat amazing combinatorics, and the geometrical
conjecture would be that all these schemes are total under a
pencil of quintics with suitably assigned 19 basepoints. The case
which looks most appealing is when there are exactly 19 empty
ovals. Those corresponds to the 4 schemes forming the vertical
left 1-simplex of Fig.\,\ref{RKM-schemes-deg-8:fig}, which we call
the monolith. The monolith has some obvious structure of a
3-simplex, stratified as follows into sub-simplices:

$\bullet$ 0-simplex corresponding to the scheme $20$ with zero
nonempty oval,

$\bullet$  $1$-simplex corresponding to the 4 schemes with 1
nonempty oval and so admitting a Gudkov's symbol
$\frac{k}{1}\ell$,

$\bullet$  $2$-simplex corresponding to the 20 schemes with 2
nonempty ovals  admitting a Gudkov writing
$\frac{k}{1}\frac{\ell}{1}m$,

$\bullet$  $3$-simplex corresponding to the 39 schemes with 3
nonempty ovals  admitting a Gudkov writing
$\frac{k}{1}\frac{\ell}{1}\frac{m}{1}n$.

Only  the schemes forming the one simplex  have exactly 19 empty
ovals. For the other categories (with resp. 20, 19, 18, 17 empty
ovals) one may assign the 19 basepoints among the nonempty ovals
(by choosing say the most massive nonempty oval, i.e. containing
the largest number of empty ovals). Of course this is pure
speculation and maybe the exact opposite has to be done.

There is probably here work for several generations of computing
machines, unless one is able to crack all total reality phenomenon
in a single stroke.
Somewhat brutally in comparison to our low understanding of where
to assign basepoints, we posit that whenever the RKM-congruence is
verified then there is a phenomenon of total reality:

\begin{conj}
Suppose given an $(M-2)$-curve of degree $m=2k$ verifying the
RKM-congruence $\chi\equiv k^2+4 \pmod 8$. Then the pencil of
adjoint curves of order $(m-3)$ ascribed to pass through the empty
ovals plus some other points distributed on the nonempty ovals is
totally real.
\end{conj}

As
to the crude arithmetics, remember that the number $B$ of
basepoints assignable to adjoints of order $(m-3)$ is given by the
binomial coefficient
$$
B=\textstyle\binom{(m-3)+2}{2}-1-1,
$$
while the pre-Harnack bound
$$
M-2=(g+1)-2=\textstyle\binom{m-1}{2}+1-2,
$$
so that
$$
B=M-3.
$$
This means that we have one basepoint less than the number of
ovals, and so we may canonically distribute them when there is one
nonempty oval.

\begin{rem}\label{M-2-curve-degree-like-Gabard:rem}
{\rm [03.03.13] Assuming  we are capable to ensure total reality
of the pencil, it may be observed that the degree of the induced
total map to $\PP^1$ would be in accordance with Gabard's bound
$r+p$, which is also the mean value of $r$ and $M=g+1$. This
follows from a simple calculation. First
$$
2B=2(M-3)=2(\textstyle\frac{(m-1)(m-2)}{2}+1-3)=(m-1)(m-2)-4=m(m-3)-2,
$$
and so the degree of the map is
$$
m(m-3)-B=2B+2-B=B+2=M-1,
$$
which is Gabard's bound i.e. the mean of $r=M-2$ and $M$. }
\end{rem}

Of course if one as some self-confidence in Gabard 2006
\cite{Gabard_2006}, then there is a total map of that degree on
each dividing $(M-2)$-curve (in particular those satisfying the
RKM-congruence which are universally of type~I). By some
concretization yoga this map would be induced by a total pencil,
and by a dubious reverse engineering of the above arithmetics this
would be a pencil of $(m-3)$-tics. This gives some very weak
evidence for the:

\begin{conj}
Any dividing $(M-2)$-curve of degree $m$ is totally real under a
pencil of curves of order $(m-3)$.
\end{conj}

[27.02.13] In fact we are not even sure that the above cubical
lattice (Fig.\,\ref{RKM-schemes-deg-8:fig}) gives an  exhaustive
list of RKM-schemes in  degree 8, where we use the jargon:

\begin{defn}
A scheme of degree $2k$ is an RKM-scheme if it is an
$(M-2)$-scheme satisfying the Rohlin-Kharlamov-Marin congruence
$\chi\equiv k^2+4 \pmod 8$ which forces type~I (alias
orthosymmetry) of the scheme, i.e. that the real locus disconnects
the complex one.
\end{defn}

While any RKM-scheme is of type~I, we do not know whether the
converse holds true (for $(M-2)$-schemes). It is true for $m=6$ as
follows from the Gudkov-Rohlin classification
(Fig.\,\ref{Gudkov-Table3:fig}).
It seems that degree 8 is a
perfidious iceberg killing any
naive conjecture arising
from contemplation of low order curves (say $\deg \le 6$).
Specific illustration of this vague principle are Shustin's
disproof of the one-half of Rohlin's maximality conjecture, and
concomitantly the disproof of Klein's Ansatz that nondividing
curves may always acquire a solitary node. Another remark along
the same line is the disproof (using the RKM-congruence) of the
naive CCC-conjecture, cf. Theorem~\ref{CCC-conj-disproof:thm}.
([06.03.13] Alas this disproof of CCC is disproved by Thom's
conjecture, as remarked there!)

Now back to our classification of RKM-schemes of degree 8 we may
wonder if there is one containing $(1,1,1)$ the nest of depth 3.
As we focus on $(M-2)$-schemes and since $M=22$ when $m=8$, we may
start with this configuration plus 17 outer ovals. In Gudkov's
notation this is the scheme $(1,1,1)17$. Any scheme $S$ of even
degree is bounded by the Ragsdale
orientable membrane
$S^{\ast}$ with $\chi(S^{\ast})=p-n$. In our case
$\chi((1,1,1)17^{\ast})=1-1+1+17=18\equiv 2 \pmod 8$ and not 4 as
posited by the RKM-congruence. If we trade outer ovals against
inner ovals lying deepest then $\chi$ is left unchanged. However
if the trading is made for ovals at intermediate depth then the
outer discs of the
Ragsdale membrane becomes holes and $\chi$ diminishes by two. So
the RKM-congruence is first arranged for the scheme
$(1,\frac{1}{1}3)14$ (with $\chi=12$), and then using 4-fold
periodicity the list is augmented as:
\begin{equation}\label{RKM-scheme-deg-8-four-primitive-type:eq}
(1,\frac{1}{1}3)14, \quad (1,\frac{1}{1}7)10, \quad
(1,\frac{1}{1}11)6, \quad (1,\frac{1}{1}15)2,
\end{equation}
which are all RKM-schemes (containing the nest of depth 3). Once
the Euler characteristic is adjusted to satisfy the
RKM-congruence,  we may trade outer ovals with innermost oval at
depth 3 without altering $\chi$. So each of these schemes produces
a list of derived schemes also RKM. Namely the 15 schemes
\begin{align*}
(1,\frac{1}{1}3)14,& (1,\frac{2}{1}3)13, (1,\frac{3}{1}3)12,
(1,\frac{4}{1}3)11, (1,\frac{5}{1}3)10, (1,\frac{6}{1}3)9,
(1,\frac{7}{1}3)8, (1,\frac{8}{1}3)7, \cr (1,\frac{9}{1}3)6,&
(1,\frac{10}{1}3)5, (1,\frac{11}{1}3)4, (1,\frac{12}{1}3)3,
(1,\frac{13}{1}3)2, (1,\frac{14}{1}3)1, (1,\frac{15}{1}3),
\end{align*}
and then the 11 schemes
\begin{align*}
(1,\frac{1}{1}7)10,& (1,\frac{2}{1}7)9, (1,\frac{3}{1}7)8,
(1,\frac{4}{1}7)7,  (1,\frac{5}{1}7)6,\cr  (1,\frac{6}{1}7)5,&
(1,\frac{7}{1}7)4, (1,\frac{8}{1}7)3, (1,\frac{9}{1}7)2,
(1,\frac{10}{1}7)1, (1,\frac{11}{1}7),
\end{align*}
and likewise the 7 schemes
\begin{align*}
(1,\frac{1}{1}11)6,\quad (1,\frac{2}{1}11)5,\quad
(1,\frac{3}{1}11)4,\quad (1,\frac{4}{1}11)3,\quad
(1,\frac{5}{1}11)2,\quad (1,\frac{6}{1}11)1,\quad
(1,\frac{7}{1}11),
\end{align*}
and finally, the 3 schemes
\begin{align*}
(1,\frac{1}{1}15)2,\quad (1,\frac{2}{1}15)1,\quad
(1,\frac{3}{1}15).
\end{align*}

(All together this gives $15+11+7+3=36$ additional schemes to be
added to the 64 tabulated on Fig.\,\ref{RKM-schemes-deg-8:fig},
hence a total of $64+36=100$.)

If we neglect the largest nonempty oval we have 19 ovals and we
may expect total reality of the quintic pencil ascribed to pass
trough any (injective) distribution of the 19 basepoints on those
19 ovals (which albeit not all empty are the deepest items of the
combinatorial scheme).

At this stage we hope to have exhausted the RKM-schemes of degree
8:

\begin{lemma}\label{RKM-schemes-deg-8-census:lem}
Any RKM-scheme of degree $8$ has either depth $\le 2$ in which
case it is catalogued as one of the $64$ schemes  of
Fig.\,\ref{RKM-schemes-deg-8:fig}, or it has depth $3$ in which
case it is one of the $4$ displayed schemes or one of the $36$
derived products where an outer oval is traded against an
innermost oval (cf. the last 4 display formulae tabulating the
corresponding $36$ Gudkov's symbols). In particular there are
exactly $64+36=100$ schemes of degree $8$ which are RKM, and hence
of type~I (and therefore potentially subsumed to the phenomenon of
total reality). [06.03.13] Addendum: One of them (at least),
namely $20$ is not realized as it violates the Thom conjecture
(cf. {\rm \ref{Thom-Ragsdale:thm}}), or better the Rohlin formula.
[13.03.13]
 {\it Warning}.---The list of 100 schemes is far from exhaustive,
cf. remarks right after that Lemma~\ref{RKM-scheme-ruled-out:lem}.
\end{lemma}


\begin{proof} (pseudo-proof)
Alas we are not even sure that this list is now exhaustive albeit
it might be likely by using the concept of depth of a scheme (the
longest chain of ovals totally ordered by inclusions of their
insides, i.e. the unique bounding disc given by the Schoenflies
theorem in its smooth variant implicit in M\"obius 1863
\cite{Moebius_1863}, Hilbert (tacit), Dehn ca. 1899 (unpublished),
Osgood 1902, Schoenflies 1906, etc., cf. e.g. Siebenmann 2005
\cite{Siebenmann_2005} for some historical background and the
literature cited therein).

Given any scheme of degree $8$, its depth is at most 4. If equal
to 4 it contains the deep nest and so the scheme is saturated
(i.e. it cannot be enlarged without corrupting B\'ezout). If the
depth is 3 then its contains $(1,1,1)$ the nest of depth 3, and if
we were not too bad in combinatorics  our recipe of 36 schemes
above was exhaustive. For the same vague reason, when the depth is
$\le 2$ then the catalogue of 64 schemes is exhaustive.

So we have one RKM-scheme of depth 1, $63$ such schemes of depth
$2$, and 36 RKM-schemes of depth 3, while the unique scheme of
depth 4 (deep nest (1,1,1,1)) is not an $(M-2)$-scheme hence not
an RKM-scheme.
\end{proof}

Further, it could be that sophisticated B\'ezout-style
obstructions \`a la Fiedler-Viro
(\ref{Viro-Fiedler-prohibition:thm}) prohibit the realizability of
some of those schemes in the algebraic realm. So perhaps several
items albeit schemes in the abstract sense of Rohlin are not
algebraically realized. (Improvising terminology and to conflict
even more with the Grothendieck-Rohlin collapse of jargon we could
speak of a Hilbert-scheme (H-scheme) when the scheme is realized
algebraically.) So I do not know if the 100 RKM-schemes listed
above are H-schemes. ({\it Update} [06.03.13] At least one of them
$20$ is not realized as follows from Thom's conjecture, cf.
Theorem~\ref{Thom-Ragsdale:thm}, or better just apply Rohlin's
formula.) Taking another naive look at the Gudkov-Rohlin table for
degree $m=6$ (Fig.\,\ref{Gudkov-Table3:fig}) shows that
$(M-2)$-schemes are subjected to no restriction and so we can
speculate the same for $m=8$, in which case all our 100 schemes
would be $H$-schemes. At any rate we note that enlarging $m=6$ by
just 2 units, involve a de-multiplication by the factor $50$ of
all RKM-schemes.

(Skip this paragraph.) Have we really listed every RKM-schemes? We
could start from another elementary configuration like 2 nests of
depth 2 (Gudkov symbol $\frac{1}{1}\frac{1}{1}=(1,1)(1,1)$) and
then add 16 outer ovals to get $(1,1)(1,1) 16$. Then $\chi=
1-1+1-1+16\equiv 0 \pmod 8$, while the good RKM-value is 4. So we
trade outer ovals for inner ovals at depth 1, and so $\chi$
diminishes by 2 units. Hence we find first
$(1,3)(1,1)14=\frac{3}{1}\frac{1}{1}14$ or
$(1,2)(1,2)14=\frac{2}{1}\frac{2}{1}14$. Those are already
catalogued on Fig.\,\ref{RKM-schemes-deg-8:fig}. Then as the Euler
characteristic is adjusted we may apply the same trick of trading
outer ovals for innermost ovals without changing $\chi$, yet doing
so we create schemes with 2 nests, one of depth 3 and one of depth
2, so that B\'ezout is violated by tracing the line through their
``centers''. So it seems that no new candidates for total reality
occurs along this way.

\subsection{How to assign basepoints?}\label{How-to-assign-basept:sec}

[01.03.13] As a matter of extending the Rohlin-Le~Touz\'e theorem
(still unpublished and abridged RLT) to degree 8, we would like to
know where to assign basepoints on each item of our list of 100
RKM-schemes. A priori not all of them are totally real in some
uniform way despite the presence of Ahlfors theorem. Recall that
our census of 100 RKM-schemes may be interpreted as five families:

(1).---the scheme $20$, (not realized by Thom
\ref{Thom-Ragsdale:thm}, or Rohlin's formula).

(2).---4 schemes of the form $\frac{k}{1}\ell$,

(3).---20 schemes of the form $\frac{k}{1}\frac{\ell}1 m$,

(4).---39 schemes of the form $\frac{k}{1}\frac{\ell}1 \frac{m}{1}
n$,

(5).---36 schemes enlarging the nest of depth 3.

The class (2) consist precisely of those elements having 19 empty
ovals. And those are the most direct candidates for an avatar of
the RLT-theorem. However in the class (5) there is also 19
preferred deep ovals, namely all those which are either empty or
if nonempty which are not maximal (for the usual order on ovals
given by inclusion of their bounding discs). So there is a family
of 40 schemes where a direct extension of the RLT-theorem is
straightforward (at least to state, but maybe not to prove).

On the other hand it could be the case that there is an extended
formulation including all those 100 schemes. At least one idea
would be to consider the notion of dextrogyre oval (abridged
dextro-oval).

For a dividing curve, we say that an oval $O$ is a dextro-oval if
its porous-inside $O^{\star}$, that is the inside minus the
insides of all ovals directly
inside it,  has
complex orientation matching  $\partial
O^{\star}$ that arising as boundary of the porous-inside.

Of course any empty oval is dextro. As an example consider the
G\"urtelkurve $C_4$ of degree 4 with 2 nested ovals. Then either
by using the total pencil of lines through a center in the
innermost of the nest or by Fiedler's law of positive smoothings
the complex orientation consist of 2 concentric circles with the
``same'' orientation. So the nonempty oval is not dextro, while of
course the inner oval is (being empty). Here we see that the
pencil of lines is total precisely when its basepoint is assigned
on the dextro-oval.

We may therefore expect that the pencil of quintics on our curves
of degree 8 is total whenever the $19$ basepoints are distributed
on  $19$ dextro-ovals supposed available.

Consider e.g. the scheme $\frac{3}{1}\frac{1}{1} 14$. Then by
Rohlin's formula $2(\Pi^{+}-\Pi^{-})=r-k^2=20-16=4$, the
difference $\Pi^{+}-\Pi^{-}$ is 2, while $\Pi^{+}+\Pi^{-}=4$, so
that $\Pi^{+}=3$ and $\Pi^{-}=1$. From this one infers (picture)
that there is at least one dextro-oval which is not empty. And so
we have here precisely 19 dextro-ovals.

Obviously one can extend this to some other schemes, and running
through the catalogue one could dress an exhaustive list of all
schemes with 19 dextro-ovals (of course 20 will not belong to it),
and expect the phenomenon of total reality for the latter.

Of course this method is somewhat ad hoc as it uses the dividing
character of the curve while in its purest form (say as a way of
taking independence of the RKM-congruence) one would like to avoid
this knowledge.

\subsection{Back to degree 6: weak form of RLT}\label{RLT:sec}

[01.03.13] Obviously we were moving too fast by looking at degree
8 and need to return to degree 6, to get rid off of those
combinatorial difficulties while concentrating on the geometry of
total reality.

Our point now is that we may be   pessimistic about
Rohlin-Le~Touz\'e's theorem as being  false in the generality
announced by Le~Touz\'e 2013\footnote{[02.03.13] This is a
misconception of mine and Le~Touz\'e's statement is finer and so
caricatural (or strong) as I misinterpreted her announcement, more
discussion about this soon. However I still do not know whether
the strong caricatural statement is wrong or not.} (or at least
difficult to prove). Recall moreover that Rohlin's cryptical
statement is not as strong as Le~Touz\'e (at least leaves some
free room for interpretation). Even if Le~Touz\'e's claim of total
reality is correct, it could be that total reality is easier to
prove for special 8 assigned basepoints. This weaker statement
would still be sufficient to detect the dividing character of
curves having an RKM-scheme, i.e. $\frac{6}{1}2$ and
$\frac{2}{1}6$.

\begin{lemma}\label{9th-basepointI:lem}
Suppose given a sextic $C_6$ of type $\frac{6}{1}2$. Assign $8$
basepoints on the empty ovals, and look at the corresponding
pencil $\Pi$ of cubics. Then there is a 9th basepoint $p_9$ of
$\Pi$. To ensure total reality of $\Pi$ it is enough that $p_9$ is
either on $N$ the nonempty oval of $C_6$ or more generally in its
inside $N^{\ast}$.
\end{lemma}

\begin{proof}
As 16 intersections are forced by topology, it remains only to
gain 2 extra intersections for totality. If $p_9\in N$ then this
is clear, and if $p_9\in N^{\ast}$ then by taking a cubic $C_3$ of
$\Pi$ through a point $p\in N$ which is connected then it is easy
to see for vibratory reasons that the other point $q$ of $C_3\cap
N$ (whose existence is ensured either by topology or by algebra)
will dextrogyrate on $N$. This is to mean that when $C_3$ is
slightly perturbed both points $p,q$ will move on $N$ along the
same orientation. This suffices to ensure total reality, as both
points cannot then enter in collision to disappear in the
imaginary locus.

This argument uses existence of a connected cubic in any pencil of
cubics. Another slight variant is to argue by contradiction by
appealing to a bad cubic $C_3$, i.e. disjoint from $N$. If $p_9$
is inside $N$ then the oval of $C_3$ cannot vibrate properly, and
we reach a contradiction.
\end{proof}

So the whole problem of Rohlin-Le~Touz\'e in weakened  form can be
reduced to the:

\begin{conj}
Given any $C_6$ of type  $\frac{6}{1}2$, there exists an
(injective) distribution of $8$ points $p_1,\dots, p_8$ on the $8$
empty ovals such that the pencil of cubics $\Pi$ interpolating
them has its $9$th basepoint located in the (sealed) inside
$N^{\ast}$ (i.e. $N$ included) of the nonempty oval. (In
particular such a calibrated pencil is totally real.)
\end{conj}

One could hope to prove this weaker assertion by pure topology
without having to enter in fine B\'ezout-style considerations
(upon which Le~Touz\'e's proof is likely to rest).

How to prove this conjecture? One very naive way would be just to
use that $N$ divides the plane and so it would suffice to find two
octuplets $p_i$ such that  $p_9$ is resp. inside and outside $N$.
By continuity of the 9th basepoint as a function $\beta$ of the 8
assigned one, we would be ensured of another intermediate octuplet
such that $p_9\in N$.

A less naive way would be to assume by contradiction that $p_9$
always misses $N^{\ast}$ (the sealed interior) and then retract on
the core of the residual M\"obius band. This seems to give an
essential (non null-homotopic) map. On the other hand as $\beta$
extends to a (16-dimensional) cell (given by allowing the $p_i$ to
explore the insides of the 8 empty ovals) it must be
null-homotopic. Alas the first step of argument is not easy to
complete.

At any rate, designating by $E_i$ the 8 empty ovals of the $C_6$,
we can define the Rohlin body of the $C_6$ as the image of the
Rohlin map
$$
\beta \colon E_1\times\dots\times E_8 \to \RR P^2
$$
taking the 8 assigned basepoints of $\Pi$ to the 9th unassigned
basepoint $p_9$.

This map is well-defined by virtue of  the following easy lemma.

\begin{lemma}\label{independent-cond:lem} Our $8$ basepoints impose independent conditions on
the space of cubics, and this holds true more generally when the
basepoints are allowed to vary in the insides $E_i^{\ast}$ of the
empty ovals $E_i$.
\end{lemma}

\begin{proof}
Let $\Pi_i$ be the linear system of curves passing through the
first $i$ points $p_1,\dots,p_i$, $i=1,\dots, 8$. We have a
filtration
$$
\vert 3H \vert \supset \Pi_1 \supset \Pi_2 \supset \Pi_3 \supset
\Pi_4\supset \Pi_5\supset \Pi_6\supset \Pi_7\supset \Pi_8=\Pi,
$$
and one checks that all inclusions are strict by exhibiting an
appropriate curve. Strictness of the first inclusions is trivial
by taking appropriate configuration of 3 lines
(Fig.\,\ref{TotFiltra:fig}). For  strictness of $\Pi_6\supset
\Pi_7$,  consider a conic $C_2$ passing through 5 inner points
plus a line through the 6th inner point, yet missing the 2 outer
points $7,8$. We have to check that this $C_3$ does not pass
through $7$. If it does then $C_2\cap C_6$ would consist of
$5\cdot 2+2+2=14>2\cdot 6$ points as two extra intersections are
created with $N$, so that B\'ezout is violated. For the last
strictness $\Pi_7\supset \Pi_8$ one takes the same $C_2$ and
aggregate the line $L$ through the points $7,8$.

\begin{figure}[h]
\centering
\epsfig{figure=TotFiltra.eps,width=122mm} \vskip-5pt\penalty0
  \caption{\label{TotFiltra:fig}%
  Checking independence of conditions} \vskip-5pt\penalty0
\end{figure}

\end{proof}

Denoting by $E:=E_1\times \dots \times E_8$ the Cartesian product
of the empty ovals and by $E^\ast:=E_1^{\ast}\times \dots \times
E_8^{\ast}$ that of their insides $E_i^{\ast}$, we have a
factorization of $\beta$ as
$$
\beta\colon E \to E^{\ast}\to \RR P^2,
$$
where the first map is the inclusion and the second $\beta^{\ast}$
is given by the same recipe as $\beta$. It follows that $\beta$ is
null-homotopic (since $E^{\ast}$ is a 16-cell).

On the other hand we may hope to show that if $\beta$ avoids
$N^{\ast}$ the inside of the empty oval that  the induced
co-restriction map $\bar\beta\colon E\to \RR P^2-N^\ast$ whose
target is homotopically a circle is {\it essential\/} (i.e. not
null-homotopic). For this it would be enough to show that the
induced morphism $\pi_1(\bar\beta)$ hits a odd multiple of the
generator of the $\pi_1$ of $\RR P^2-N^\ast$ which is a M\"obius
band.  This contradiction would prove the conjecture and so total
reality of a suitable pencil. However this strategy demands some
geometric understanding that presently elude us.

Another more naive strategy (using less topology) would be that
the map $\beta\colon E \to \RR P^2$ is open (say as a vestige of
the holomorphic character of the  underlying complexification).
Then it would be plain that $\beta(E)$ is compact (hence closed)
and open, hence equal to all $\RR P^2$ by connectedness of the
latter space. (As to be soon discussed this surjectivity of
$\beta$ would however conflict with the Le~Touz\'e theorem)

At any rate we see that the problem of total reality \`a la
Rohlin-Le~Touz\'e in its weakened form due to us (or perhaps
Rohlin depending on the interpretation of his cryptical statement)
is fairly basic at first sight. We are given a curve $C_6$ of type
$\frac{6}{1}2$. Given any 8 points on the empty ovals $E_i$ we
have some predestination mapping $\beta$ assigning to them the 9th
basepoint (by a somewhat elusive recipe) and total reality of the
pencil $\Pi$ of curves passing through the 8 points is ensured if
$p_9$ lands in the inside $N^{\ast}$ of the nonempty oval $N$ of
$C_6$. So crudely speaking we have one chance over two that total
reality holds true, for a given configuration of point. We would
like to show that it is always possible to have total reality for
a clever configuration of 8 points, while the stronger
Le~Touz\'e's announcement claims it for all choices of 8 points.

Our hope is that independently of whether this stronger statement
is right or false there ought to be a simpler proof of the weaker
assertion by say essentially topological methods. By using the
dextrogyration argument it is clear that we have the:

\begin{lemma}\label{9th-basepoint-totalII:lem} The
pencil $\Pi$ is total iff $p_9$ its non-assigned basepoint belongs
to $N^{\ast}$ (the sealed inside of the nonempty oval $N$).
\end{lemma}

\begin{proof} If $p_9\in N$ total reality is clear. If $p_9$ is in
the open inside (interior) of $N$, then we have an odd number of
basepoints insides. Looking at the oscillation of a connected
member $C_3$ of the pencil about its basepoints it results that
both points of $C_3\cap N$ (there cannot be more than two by
B\'ezout) will dextrogyre, i.e. move along the same orientation of
$N$ (compare Fig.\,\ref{Dextrogyre:fig}a). Hence no collision can
occur in the long run, since by the holomorphic character of the
map a point cannot reverse spontaneously its sense of  motion as
the curve $C_3$ is dragged along the real locus of the pencil.

{\it Insertion} [11.04.13].---We can dispense  connectedness
 of the $C_3$ (though quite easy to arrange) as follows.
Choose any point on $N$ and consider the cubic $C_3$ through this
point. If $C_3\cap C_6$ is not totally real, the oval of $C_3$
 is necessarily inside $N$. For
vibratory reasons this oval must visit an even number of
basepoints, and actually must visit all 6 assigned inner
basepoints (otherwise total reality of $C_3\cap C_6$ is granted).
The 9th basepoint cannot be located on the oval of the $C_3$ (else
it cannot vibrate properly), hence it is situated on its
pseudoline and we may again conclude dextrogyration  by the
slaloming argument across an odd number of basepoints (compare
Fig.\,\ref{Dextrogyre:fig}b).

\begin{figure}[h]
\centering
\epsfig{figure=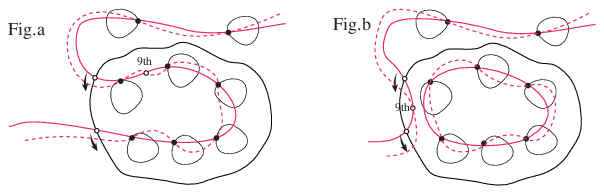,width=122mm} \vskip-5pt\penalty0
  \caption{\label{Dextrogyre:fig}%
A dextrogyration argument forced by an odd number of slaloming}
\vskip-5pt\penalty0
\end{figure}

Conversely if the 9th basepoint $p_9$ lies outside of $N$ then
there is an even number of basepoints inside $N$ and the 2 points
of a connected $C_3$ located on $N$ will anti-dextrogyre. In that
case there will be a collision in the long run and total reality
is foiled.
\end{proof}

Hence the option $\beta$ surjective would contradict the
Le~Touz\'e's theorem. Of course this is not a serious objection
against her theorem because usually holomorphic maps restricted to
real loci fails blatantly to be surjective (a key prototype of
this phenomenon opposite to total reality, is when one projects on
a line an ellipse from an outer point).

{\it Optional paragraph.}---Another slight variant: we may assign
7 basepoints on all safe one empty ovals and the 8th one $p_8$ on
$N$. Then a little advantage is that  whenever $p_8$ is collinear
with 2 points of the first seven $p_i$, we get a special split
cubic in the pencil $\Pi$, namely the line plus a residual conic.
It is not very clear if this little advantage is really useful,
and leave open this discussion.

\subsection{New meditation after
reception of Le~Touz\'e's article}

[02.03.13] Yesterday (01.03.13), we conjointly received (with
several other Russian colleagues Kharlamov, Viro, Nikulin, etc.) a
copy of Le~Touz\'e's article vindicating Rohlin's cryptical
assertion of total reality of the RKM-sextics of type
$\frac{6}{1}2$ and $\frac{2}{1}6$ under a pencil of cubics.

This fascinating paper helped me to rectify several misconceptions
of mine about the content of her earlier announcement. In
particular, she writes the following illuminating remark
(Le~Touz\'e\footnote{The official name of this  author is
Fiedler-Le~Touz\'e, yet as well-known n\'ee Le~Touz\'e and
abridged as a such in the sequel.} 2013
\cite[p.\,3]{Fiedler-Le-Touzé_2013-Totally-real-pencils-Cubics}):

``By a congruence due to Kharlamov [8](=Kharlamov-Viro 1988/91
\cite{Kharlamov-Viro_1988/91}), the real schemes $\la 2 \vc 1 \la
6\ra\ra$ or $\la 6 \vc 1 \la 2\ra\ra$ are both of type~I. We use
this fact in the proof, so we confirm that the sextics with this
two schemes do not contradict Rokhlin's conjecture. Rokhlin
claimed that he could prove the very same statement {\it
without\/} using the fact that the sextics are dividing. It's a
stronger result. Unfortunately, his proof was never published and
is now lost\footnote{This seems alas to be the bitter state of
affairs as follows from a recent consultation (January 2013) of
the leading experts (Viro, Marin, Kharlamov, Fiedler, Le~Touz\'e,
etc.); compare e-mails gathered in Sec.\,\ref{e-mail-Viro:sec}}.''

The Rokhlin conjecture alluded to by Le~Touz\'e is made explicit
in (p.\,2--3 of \loccit), which is again worth quoting:

``The $M$-curves are clearly dividing. The so-called hyperbolic
curves are also dividing. A hyperbolic curve of degree $m=2k$ or
$2k+1$, consists in $k$ nested ovals, plus one pseudo-line if $m$
is odd. A pencil of lines whose base point is chosen in the
innermost oval sweeps out the curve in such a way that the $m$
intersections are always real. One says that this pencil of lines
is totally real with respect to the hyperbolic curve. Starting
from this observation, Rokhlin [10](=1978 \cite{Rohlin_1978})
presents a beautiful argument\footnote{I agree, but the argument
is nearly trivial in the sense that it just uses the fact that the
image of a connected set is connected (Marin's oral remark during
my talk in Grenoble ca. 1999). Personally, I knew this argument
since 1999 (arguing with pathes prior to Marin's oral
simplification of it). We cannot record if we rediscovered it
independently of Rohlin 1978 (but do record that we may have found
some indirect inspiration from Gross-Harris 1981, who treat  the
case of hyperelliptic curves $y^2=f(x)$, with $f(x)>0$
throughout). At any rate it is evident that Rohlin's argument can
be drastically simplified. Rohlin uses a lot a certain fibering
while it is plain that it suffices to use the map to the
equatorial (orthosymmetric) sphere, cf. e.g. Gabard 2006
\cite{Gabard_2006}.} proving that if an algebraic curve is swept
out by a totally real pencil of lines, then this curve is
dividing. The argument generalizes to pencils of curves of higher
degrees. Can conversely any dividing curve be endowed with some
totally real pencil?\footnote{My opinion was always that a
positive answer should be a trivial consequence of Ahlfors theorem
(cf. e.g. Gabard's Thesis 2004 \cite[p.\,7]{Gabard_2004}). However
since Marin warned me in January 2013 (cf.
Sec.\,\ref{e-mail-Viro:sec}) it may be the case that the
transition from the abstract to the embedded viewpoints is not so
easy. Yet I am still confident that it holds true. The point is
primarily a matter of projective algebraic geometry, namely the
question if any abstract morphism on a concrete plane curve is
induced by a pencil of ambient curves. This is either trivially
true or trivially wrong, but alas I do not know the answer due to
failing memory about the foundations of algebraic geometry.
[12.04.13] Additionally, it can also be that sometimes imaginary
basepoints have to be allowed, and so total reality really the
mobile part of the pencil. We hope to be capable presenting this
more clearly in the future, but see perhaps already
(\ref{Le-Touzé-scholium-deg-6:lem}).} A weaker conjecture
suggested implicitly in [10](=Rohlin 1978 \cite{Rohlin_1978}) is
that any curve whose real scheme is of type~I may be endowed with
a totally real pencil\footnote{Whether this is implicit or not is
an interpretation-matter, unless of course some direct contact
with Rohlin (e.g. via the husband)
testifies such a conjecture of Rohlin. Again it could be the case
that Ahlfors theorem nearly trivially implies this (novel)
conjecture of Rohlin. Even if so in abstracto then the game is
probably far from finished as one would like to get synthetic
descriptions of the total pencils. It seems quite likely that this
game can keep busy several generation of workers. [12]}. It turns
out that the $M$-curves may indeed be endowed with suitable
pencils of degree $(m-2)$, see [6], page
348(=Theorem~\ref{total-reality-of-plane-M-curves:thm} in this
work=Gabard 2012/13 \cite{Gabard_2012/13},  pagination
may
have fluctuated meanwhile).''

Le~Touz\'e's article clarified several misinterpretation of mine
about her announced result (dated 16 f\'evrier 2013, cf.
Sec.\,\ref{e-mail-Viro:sec}) but which we reproduce now as we
misunderstood it:

$\bullet\bullet\bullet$ samedi 16 f\'evrier 2013 17:54:55

Dear Alexandre, dear other colleagues,

I have managed to prove that a pencil of cubics with eight
base~points distributed in the eight empty ovals of a sextic
$2 \cup 1(6)$ is necessarily totally real. Details will follow
soon in a paper. Yours,
S\'everine

In fact my misconception was to think that Le~Touz\'e claims total
reality for {\it any} such pencil. So shame on me for not having
read her message more carefully as she expressly writes {\it a
pencil\/}. For me it is still unclear if the stronger claim of
total reality for all such pencils holds true. Such a strong form
of total reality holds in the basic cases (pencil of rational
curves of degree $\le 2$) but perhaps the case of cubics is
completely different for such curves need not being connected and
also there is a 9th predestined basepoint which cannot be freely
assigned (magneto repulsion well-known at least since Euler 1748
as reported e.g. in Griffiths-Harris 1978
\cite[p.\,673]{Griffiths-Harris_1978/94}).

Taking (deliberately) the risk of being too cavalier let us put
forward the strongest form of total reality (\`a la
Rohlin-Le~Touz\'e but perhaps too coarsely interpreted
than what those authors ever wrote):

{\it Added in proof} [08.03.13].---Meanwhile Le~Touz\'e validated
the following conjecture, for explanations, cf. her message (dated
5 March 2013 in Sec.\,\ref{e-mail-Viro:sec}).

\begin{conj}\label{SRLT:conj} {\rm (SRLT=Strong
Rohlin-Le~Touz\'e total reality, as misinterpreted by Gabard)}
Given any (smooth) sextic $C_6$ of type $\frac{6}{1}2$ or
$\frac{2}{1}6$ and any (injective) distribution of $8$ points on
the $8$ empty oval of $C_6$, the pencil of cubics through the $8$
assigned basepoints is totally real.
\end{conj}

First, the basic Lemma~\ref{independent-cond:lem} (on independency
of conditions) seems to imply that the given linear system is
really a pencil. The next step would be to prove (or disprove)
this strong conjecture of total reality. The (more
cautious) statement of Le~Touz\'e
does not imply the conjecture, yet it would be interesting to know
if Le~Touz\'e is aware of an obstruction
refuting the conjecture.
Let us reproduce Le~Touz\'e's theorem (cf. Le~Touz\'e 2013
\cite{Fiedler-Le-Touzé_2013-Totally-real-pencils-Cubics}):

\begin{theorem}\label{LeTouze-big-thm-01-March-13:thm} {\rm (Le~Touz\'e 1st March 2013,
announced 17 February 2013)} Any $(M-2)$-sextic with real scheme
$\frac{2}{1}6$ or $\frac{6}{1}2$ may be endowed with a totally
real pencil of cubics with $8$ basepoints distributed in the $8$
empty ovals.
\end{theorem}

Also very interesting (though somewhat disappointing) is the issue
that Le~Touz\'e's proof uses the RKM-congruence (as brilliantly
nuanced by herself), and so does not recover the {\it synthetic a
priori\/} character of Rohlin's claim (which so becomes even more
cryptical than it ever was). This dependence of Le~Touz\'e's proof
on the RKM-congruence (which she ascribes to Kharlamov alone,
like in the original text Rohlin 1978)
has some overlap  with thoughts I had yesterday (especially
Sec.\,\ref{How-to-assign-basept:sec}).

Le~Touz\'e's work makes very acute the desire to find the
so-called {\it lost proof of  Rohlin\/}. We thought on the
question from various angles yet our present understanding is very
confused due to overwork and nervous collapse
(i.e., plethora of strategies and panoply  of statements of
varying strength).

The naive strategy (toward conjecture SRLT=(\ref{SRLT:conj}) or
just a weakened form thereof) at which we arrived yesterday
(Sec.\,\ref{RLT:sec}), was to use the map $\beta$ assigning to the
8 assigned basepoints the 9th (unassigned) basepoint. With this
the lost proof of Rohlin could be harpooned by an essentially
topological argument, yet we failed to overcome the last
difficulty.

For convenience let us repeat some of our ideas along this
topological tactic. Given a $C_6$ say of type $\frac{6}{1}2 $ for
simplicity. Denote by $E_i$ the 8 empty ovals. For any injective
distribution of points on those 8 ovals we have a pencil $\Pi$ of
cubics through them (Lemma~\ref{independent-cond:lem}), and the
latter is totally real iff the 9th basepoint of $\Pi$ lands in the
(sealed) inside $N^{\ast}$ of $N$ the nonempty oval of $C_6$ (cf.
Lemmas \ref{9th-basepointI:lem} and
\ref{9th-basepoint-totalII:lem}).

Formally we can so introduce the {\it $9$th basepoint map\/}
$$
\beta\colon E_1\times \dots \times E_8 \to \RR P^2
$$
taking the octuplet $(p_1, \dots , p_8) $ to the 9th basepoint  of
the pencil $\Pi$ of cubics passing through $p_1, \dots, p_8$.

A priori we could hope  this mapping to be onto for reasons of
Brouwer's topological degree of a map (in homology modulo 2).
Precisely if the induced morphism $H_2(\beta, \ZZ_2)$ is
non-trivial, the mapping $\beta$ would be surjective (for
otherwise a point is missed and the map factorizes through a
punctured (hence open) Hausdorff manifold whose top-dimensional
homology vanishes). Of course all this general theory can  be
dispensed as we  just have a punctured projective plane
homotopically equivalent to a circle. In this surjectivity
scenario for $\beta$, it would suffice to take an inner point
$p_9\in N^{\ast}$ and any lift of it via $\beta$ yields an
octuplet
inducing a total pencil on the $C_6$.
In contrast,
taking a point outside $N$ would imply that the
strong form of total reality (SRLT)=(\ref{SRLT:conj}) for all
octuplets fails.

However it is unlikely that this Brouwer-style surjectivity
criterion works because the map $\beta$ extends (still by virtue
of Lemma~\ref{independent-cond:lem}) to the (sealed) insides
$E_i^{\ast}$ (bounding discs) of the $E_i$ as to give the map
$$
\beta^{\ast}\colon E_1^{\ast}\times \dots \times E_8^{\ast} \to
\RR P^2,
$$
also defined by taking the 8 assigned basepoints to the 9th one.
Since the source of the map $\beta^{\ast}$ is a contractible
$16$-cell (hypercube) the induced map $H_2(\beta, \ZZ_2)=0$ is
trivial.

Still there is some hope to do something good. If the
range(=image) of the map $\beta\colon E:=\times_{i=1}^8 E_i \to
\RR P^2$ meets $N^{\ast}$ we have total reality at least in the
weak form (RLT). (Asserting the strong form SRLT amounts knowing
that $\beta(E)\subset N^{\ast}$.) Hence to prove RLT it suffices
to show that the option $\beta(E)$ disjoint form $N^{\ast}$ leads
to a contradiction. Our naive idea is then that the map $\beta$
would have its range confined to the M\"obius band $M:=\RR
P^2-N^{\ast}$ residual to $N^\ast$. If one is able to show that
$\beta\colon E \to M$ is not null-homotopic (e.g. by showing that
it hits an odd multiple of the generator of $\pi_1(M)=\ZZ$), then
it would follow that $\pi_1 (\beta)$ is nontrivial, violating the
above factorization $\beta^{\ast}$ through a contractible 16-cell.
This contradiction would prove the weak form of RLT, i.e.
existence of an octuplet inducing a total pencil.

Of course the above tactic to be completed requires a proof of the
hypothetical fact that there is a loop in the 8-torus $E$ taken by
$\beta$ to the nontrivial element of $\pi_1(\RR P^2)$. This seems
hard to prove and requires at any rate some geometric
understanding of pencil of cubics, especially of the
predestination process creating the 9th basepoint as a function of
the 8 assigned ones.

More pragmatically we could define the {\it Rohlin body\/} of a
given $C_6$ (of RKM-type) as the image $B:=\beta (E)$. This
compactum and especially its location w.r.t. the nonempty oval $N$
will govern much of the total reality question of the $C_6$.
Essentially we have a trichotomy of alternatives:

$\bullet$ Either $B \subset N^{\ast}$ in which case the $C_6$ is
strongly totally real, in the sense that any octuplet (in $E$)
induces a totally real pencil, or

$\bullet$   $B$ overlaps $N^{\ast}$ without being contained in it,
in which case some octuplets induce a totally real pencil and some
other do not, or finally

$\bullet$  $B$ is disjoint of $N^{\ast}$ in which case all
octuplets fail inducing a total pencil. (This is nearly
incompatible with the RLT total reality phenomenon).

In the first scenario we could say that the sextic $C_6$ is {\it
strongly totally real\/}, in the second that  is {\it (weakly)
totally real\/}, and in the 3rd case that is ``anti-real''.

We believe that a continuity/degeneration argument applied to
Le~Touz\'e's theorem (2013
\cite{Fiedler-Le-Touzé_2013-Totally-real-pencils-Cubics}) prevents
the 3rd option, by letting degenerate the 8 basepoints to the
empty ovals. Alas it does not seem that Le~Touz\'e's theorem
(\ref{LeTouze-big-thm-01-March-13:thm}) gives sufficiently many
inner permissible basepoints so as to ensure
a degeneration to the ovals.

As to the  first two scenarios, we do not know if both of them do
occur, and if not which one is ubiquitous. Paraphrasing, we do not
know if there is a single sextic of type $\frac{6}{1}2$ (or it
mirror) such that any octuplet chosen on the empty ovals induces a
totally real pencil of cubics,
nor can we preclude the option that all sextics verify this
property.

{\it Insertion} [12.04.13].---If we interpreted correctly the last
news of Le~Touz\'e, it seems that total reality holds in the
strongest possible sense, i.e. for all octuplets and even when the
latter are inside of the ovals and not on themselves directly.
Compare her article (2013
\cite{Fiedler-Le-Touzé_2013-Totally-real-pencils-Cubics}, plus
eventually  her letter in Sec.\,\ref{e-mail-Viro:sec}.

\section{Esquisse
d'un programme (d\'ej\`a esquiss\'e): the Ahlfors-Rohlin
Verschmelzung} \label{Esquisse-dun-prog-deja-esquiss:sec}

[07.03.13] The programme in question was probably nearly implicit
in Rohlin 1978 \cite{Rohlin_1978}, safe that he seems to have
missed the possible connection with Ahlfors theorem. It is only in
this respect that our programme bears some originality, yet
presently we are unable to substantiate it in any serious fashion.

\subsection{Large scale structure of total reality as it pertains
to Hilbert's 16th problem}

[04.03.13] Let us brush a sloppy summary of the situation. Rohlin
1978 \cite{Rohlin_1978} in somewhat cryptical fashion asserted
existence of a totally real pencil of cubics on all sextics curves
of type $\frac{6}{1}2$ or its mirror, yielding therefore a
geometrization of the RKM-congruence asserting that $\chi\equiv
k^2+4 \pmod 8$ forces the curve of degree $2k$ being of type~I.
After consulting several specialists (Viro, Marin, Kharlamov,
Fiedler, Le~Touz\'e), it seems that this proof
is now lost forever
(or dormant in some celestial Eden). It seems extremely
challenging to rediscover it if it ever existed, i.e. if Rohlin's
argument was sound and complete,
as opposed to a Fermat-style cryptical allusion destined to
challenge
geometers. On reading the survey Degtyarev-Kharlamov 2000
\cite{Degtyarev-Kharlamov_2000}, it seems that this is not the
sole prophetical allusion of Rohlin in his fantastic 1978 survey.
Also the question of estimating the number of points through which
one can pass a rational connected curve is also sloppily stated by
Rohlin without proof, and Degtyarev-Kharlamov consider the problem
as still open. Perhaps the status of this problem evolved
meanwhile.

Le~Touz\'e 2013
\cite{Fiedler-Le-Touzé_2013-Totally-real-pencils-Cubics} supplied
the first (and unique) written proof of a weak form of Rohlin's
total reality phenomenon, yet assuming the RKM-congruence, and so
the curve to be of type~I.
As suggested in Le~Touz\'e's article (2013 \loccit), Rohlin seems
to have
cultivated
a large
expansion of the
phenomenon of total reality, as follows:

\begin{conj} {\rm TOR=Total reality (Le~Touz\'e 2013,
Gabard 2004 \cite{Gabard_2004}, who ascribed this as implicit in
Ahlfors 1950, and Teichm\"uller 1941 \cite{Teichmueller_1941}, who
loosely claims this to be found in Klein's works).}---Any dividing
(plane) curve admits a totally real pencil of curves.
\end{conj}

\begin{conj} {\rm ROTOR=Rohlin's total reality (implicit
in Rohlin 1978 according to Le~Touz\'e 2013)}.---Any dividing
curve representing a scheme of type~I admits a totally real
pencil.
\end{conj}

\begin{conj} {\rm RMC=Radio
Monte Carlo=Rohlin's maximality conjecture}.---Any real scheme of
type~I is maximal in the hierarchy of all schemes (of some fixed
degree).
\end{conj}

As noted by Le~Touz\'e (always \loccit), TOR implies ROTOR, and it
is always tempting to believe that the latter implies RMC.

At the very source of that string of implications, we should have
the Ahlfors theorem (ATR=abstract total reality) (which
Teichm\"uller 1941 \cite{Teichmueller_1941} ascribes to Klein)

\begin{theorem} {\rm ATR=Abstract total reality (Ahlfors 1950,
Gabard 2004--06).}---Any (abstract) dividing curve $C$ (or what
Klein calls an orthosymmetric Riemann surface) admits a totally
real map $C\to \PP^1$ to the projective line, i.e. such that
$f^{-1}(\PP^1 (\RR))=C(\RR)$. Further the degree of such a map can
be arranged $\le g+1$, where $g$ is the genus of $C$. According to
Gabard 2006, but still deserves to be better analyzed, there
should even be such a total map of degree $\le \frac{r+(g+1)}{2}$
where $r$ is the number of real circuits of the curve $C$.
\end{theorem}

ATR should imply TOR modulo some little difficulties. A first
difficulty is merely due my incompetence, of knowing if an
abstract morphism from a plane curve to the line is necessarily
induced by a pencil of curves. Another little difficulty is that a
priori some basepoints (of a total pencil given by Ahlfors) may be
imaginary conjugate, so that total reality has perhaps to include
a slightly broader definition (admittedly more cumbersome) than
the one used in Le~Touz\'e. So probably the correct definition of
a totally real pencil has to involve total reality of the moving
points of the linear series, while the (statical) basepoints
themselves being permitted to be non real(=imaginary, as we say
since Tartaglia-Cardano (1535/39), essentially).

At the very end of the string
ATR$\Rightarrow$TOR$\Rightarrow$ROTOR$\Rightarrow$RMC, we seem to
have (accepting the elusive RMC) a sort of subordination of all
highbrow congruences modulo 8, that is the Gudkov-Rohlin
congruence for $M$-curves, the Gudkov-Krakhnov-Kharlamov
congruence for $(M-1)$-curves
(\ref{Gudkov-Krakhnov-Kharlamov-cong:thm}), to the
(Rohlin)-Kharlamov-Marin congruence
(\ref{Kharlamov-Marin-cong:thm}). To appreciate this fact in
degree 6, contemplate once more the Gudkov-Rohlin
Table(=Fig.\,\ref{Gudkov-Table3:fig}). Then we see that
(virtually) all congruences are subordinated to that of
Kharlamov-Marin modulo the truth of RMC, safe the prohibition of
the schemes $11$ (Hilbert) and $\frac{10}1$ (Rohn). The latter may
be  prohibited either via Arnold's congruence or Rohlin's formula.

So in crude approximation (and modulo RMC) ``all'' the
prohibitions of Hilbert's 16th problem (say perhaps apart from
refined B\'ezout-style obstructions \`a la Fiedler-Viro virulent
in degree $8$ or higher) are subsumed to the Kharlamov-Marin
congruence ensuring type~I whenever $\chi\equiv k^2+4 \pmod 8$
(abridged RKM, where R stands for Rohlin, albeit I am not sure
about his exact contribution, while most writers credit Kharlamov
and Marin only, cf. e.g. Kharlamov-Viro 1988/91
\cite{Kharlamov-Viro_1988/91}).
({\it Insertion} [12.04.13].---In Rohlin's 1978
\cite[p.\,93]{Rohlin_1978} the result is credited to Kharlamov
alone, but it is remarked that the proof of this theorem was still
unpublished.)

Accepting the above string of implications, we arrive essentially
at the conclusion that all of Hilbert's 16th problem could be
governed by the phenomenon of total reality. This would be
especially true if there is a geometrization of the RKM-congruence
via total reality. This requires a highbrow extension of the
Rohlin-Le~Touz\'e theorem to all other $(M-2)$-schemes verifying
the RKM-congruence. This looks hard but according to the previous
sections it is likely that the total pencil is always of order
$(m-3)$ involving the so-called adjoint curves of Brill-Noether
incarnating the canonical series. The resulting degree of the map
would then also be in accordance with Gabard's version of Ahlfors
theorem (compare Remark~\ref{M-2-curve-degree-like-Gabard:rem}).

To caricature a bit, as to emphasize our philosophy (especially as
it pertains to the title of the present text devoted to Ahlfors),
we formulate the:

\begin{Scholium} All (or most) of Hilbert's 16th problem
can be reduced to the phenomenon of total reality, due to various
authors. In in its most primitive schlichtartig form, this
involves primarily Riemann 1857 \cite{Riemann_1857_Nachlass},
Schottky 1875--77 \cite{Schottky_1877}, Bieberbach 1925
\cite{Bieberbach_1925}, Grunsky, etc. and in general to Klein
according to Teichm\"uller 1941 \cite{Teichmueller_1941}, Ahlfors
1950 \cite{Ahlfors_1950}).
\end{Scholium}

Basically the phenomenon of total reality involves a linearization
in the large of the curve via the concept of branched coverings,
where any dividing curve is reduced to its most primitive
incarnation, namely the line $\PP^1$. Any extraterrestrial planet,
possibly  with handles and several equators ({\it exoplanet\/}),
yet not so exotic as to share with our planet Earth the character
of orthosymmetry(=type~I=dividing) involving 2 distinguishable
hemispheres, can be
conformally shrunk to our equatorial sphere so that fibers above
the equator $\PP^1(\RR)$ are totally real, i.e. on the
exo-equators of the exoplanet\footnote{[12.04.13] (Please skip
this footnote, if you believe in capitalism).---We invented this
exoplanet metaphor in 2004, as to sell our postdoc-research
programme to an FNS-administrator (FNS=SNF=Schweizerische National
Fond), specialized in astronomy (at some Geneva observatory). The
success was very limited, no funding were ever obtained and much
energy and time wasted for nothing. Some few weeks later another
Swiss cooperative stole me 15'000 Euros of economies. Life then
started to require environmental punch (nutrition in the
containers, and other pleasant duties like bicycling the heavy
nutriments over steep mountains). Can we develop a more tolerant
science enrolling more people on less restrictive financial
constraints, especially more modest retribution of the workers?
Ahlfors is far from a hero in this respect (elitist attitude than
looks much overdone in view of the  little originality of his
contributions to science, compare what he borrowed from Courant,
Hurwitz, Riemann, Gr\"otzsch, Teichm\"uller, etc.). The real
question is of course: can we get rid off of capitalism, granting
the fact that a sufficient motor of life is to reach immortality
(for free and for all), as it was ever encoded in our genes since
the amoebic morphogenesis.}.

This is an abstract theorem yet it should imply a concrete result
of total reality like the above TOR. This process could be termed
 {\it descent\/} from the Riemannian universe (Gromov's prose)  to the terrestrial
 Plato cavern (of Hilbert's 16th problem). Especially important
would be a quantitative control on the order of those total
pencils obtained by descent of the abstract
Riemann-Schottky-Klein-Teichm\"uller-Ahlfors circle maps
(equivalently total maps to $\PP^1$). For a simple implementation
in the case of $M$-curves, cf. our
Theorem~\ref{total-reality-of-plane-M-curves:thm}, which is a
trivial adaptation of the no collision principle of
Riemann-Enriques-Chisini-Bieberbach-Wirtinger, etc.

Next when we move down to $(M-2)$-curves the first concrete (and
nontrivial) phenomenon of total reality is the Rohlin-Le~Touz\'e
theorem (abridged RLT after their author's names or for (total)
ReaLiTy). Historiographically, it is noteworthy that Rohlin does
not seem to have ever been aware of Ahlfors theorem, and it seems
that the gap between
both traditions (Riemann vs. Hilbert) as not yet been fully
bridged. It is also notorious (either from the viewpoint of
geometric function theory \`a la Riemann,
or the algebro-geometric perspective) that total reality is much
easier for $M$-curves as it involves a schlichtartig semi-Riemann
surface (planar orthosymmetric half). This is also evidenced by
the fact that there is no collision between a group of $g+1$
distributed on the $M=g+1$ ovals. In the non-Harnack maximal case,
any group of $g+1$ points moves (Riemann-Roch) but then several of
them being distributed on the same oval a risky collision can
occur in the long run, foiling
total reality. A subtle condition of dextrogyration must be
ensured to gain total reality.

So what remains to be done?

Project 1.---Try to clarify the above logical implications between
ATR, TOR, ROTOR, RMT. This is basically a Riemann-ification of
Hilbert's 16th problem.

Project 2.---Try to understand better the lost proof of Rohlin,
and how the statement extends to curves of higher order.

Two routes toward (Project 2) seems a priori possible making the
exploration somehow elusive or at least time consuming. Either
work from the scratch in the Plato cavern of  plane curves
grooving the plane $\PP^2$ where Hilbert's problem is formulated,
or attempt a descent of the abstract result \`a la
Riemann-Schottky-Klein-Enriques-Chisini-Bieberbach-Grunsky-Wirtinger-Teichm\"uller-Ahlfors-Ga\-bard.
The method of descent looks a priori delicate but  worked for
(plane) $M$-curves, cf. again
Theorem~\ref{total-reality-of-plane-M-curves:thm}. The drawback of
this method of descent is that it would (like Le~Touz\'e's proof)
depend upon a knowledge a priori of the dividing character of the
curve. If  optimistic it could be  that Gabard's theorem
(2004-2006 \cite{Gabard_2006}) could imply total reality for
dividing $(M-2)$-curves  of degree $m$ via a pencil of (adjoint)
curves of degree $(m-3)$ (cf. again
Remark~\ref{M-2-curve-degree-like-Gabard:rem} for some weak
evidence). As a special case it could be the case that Gabard's
theorem implies a weak-form of the RLT-theorem, yet much work
along ATR$\Rightarrow$TOR is required to bridge the gap.

\subsection{Rohlin's intuition vindicable via Ahlfors?
Algorithmic r\^ole of RMC for plotting machines}

[05.03.13] All what we have to say is now essentially well-known,
but alas still in embryonal state. Let us try to make the
philosophy behind the grand dessin imagined by Rohlin more
palatable, than in the previous section.

What is Hilbert's 16th problem at all about? Answer: Topology of
real plane algebraic curves. But in reality we are geometers and
topology is merely a weakening of what wanted to study earlier
geometers (say  Diophante, Euclid, Archimedes,
the algebro-geometrization of  Fermat, Descartes, Newton, etc.)
What is called upon is an understanding of the big video game
where given an equation one traces the corresponding curve in real
time. One should imagine a powerful enough machine showing us in
real time the curve evolving when dragging with a joystick the
coefficients of the equation.

Any machine able to do this presumably
request at an algorithm telling
when to stop the tracing procedure as the real locus has been
represented within sufficient accuracy as to
infer the exact topology
of the curve.

It is at this stage that Rohlin's maximality conjecture may enter
into the scene. Indeed for a given degree there exists certain
distinguished schemes (in the sense of Rohlin) representing
so-to-speak fully crystallized extremal shapes not susceptible of
any further apparition of ovals (as small as they may be). As
popularized by Rohlin 1978, this intuition of saturated schemes
truly goes back to Klein 1876, yet in some primitive sense
distinct from Rohlin's interpretation.

Consider for instance a curve of degree 6 whose real locus
contains a deep nest of depth 3, then it is already saturated and
there cannot be any further oval. More generally we have certain
schemes of type~I, which according to Rohlin's intuition ought to
be maximal, hence incarnating the maximum topological complexity
permissible for the given degree.

The problem of understanding this ontological truth splits in two
parts: Why is it true, and why is it useful?

First we try answering the second question. As already pointed
out, Rohlin's maximality conjecture (RMC) incarnates a sort of
stopping process for the plotting machine (realizing the dynamical
 video game or just its statical variant). From the
viewpoint of Hilbert's 16th problem, the fact that the deep nest
of depth 3 is a maximal scheme of degree 6 forbids a menagerie of
other schemes enlarging it which a priori could exist, but do not
essentially by virtue of B\'ezout's theorem. Hence Rohlin's type~I
schemes (granting their maximality) are like advanced sentinels
prohibiting schemes of higher topological complexity. Without such
prohibitions Hilbert's 16th problem would be even more intractable
than it already is. Hence Rohlin's maximality conjecture is a sort
of upper bound for the complexity of Hilbert's 16th problem. This
 should be sufficient reasons for answering the
utilitarian aspect.

As to the first question, it is somehow ironical that Klein seems
to have been much in touch with both aspects of our question.
First, as we said he is regarded (by Rohlin himself) as
a precursor of Rohlin's maximality conjecture. Second (but this is
more elusive to
testify
with high accuracy), Klein is credited by Teichm\"uller 1941
\cite{Teichmueller_1941} as the true forerunner of Ahlfors
theorem. At some broader scale, all goes back to Riemann
(especially if his life would not have been so short).

Rohlin's maximality conjecture should according to our intuition
(Gabard ca. 1st January 2013) reduces to Ahlfors theorem, via what
Le~Touz\'e calls the Rohlin total reality conjecture, cf. (ROTOR)
of the previous section. So it is quite interesting to see that
the extremal shapes (maximal schemes) of Rohlin are induced by
schemes of type~I, and what makes this possible is the phenomenon
of total reality. Behind the latter there is of course Ahlfors
circle maps, and so basically an extension of the Riemann mapping
theorem. This in turn is governed by potential theory, itself
concomitant of the calculus of variation of Euler-Lagrange as
applied to Laplace's equation. All this d\'etour to make apparent
that the algebro-geometric extremal principle posited by Rohlin's
maximality conjecture seems governed by another extremal
principle, namely those ensuring solvability of Dirichlet's
principle. So be it via Abel (and what some like to call Hodge
theory) or directly via
Euler-Lagrange-Laplace-Dirichlet-Riemann-etc-Ahlfors we have the
phenomenon of total reality for curves of type~I, and when the
scheme itself is of type~I this phenomenon acquires an extra punch
of universality, making the scheme maximal in the
Hilbert-Gudkov-Rohlin hierarchy.

Coarsely, our thesis could be that what missed to Rohlin to
complete his programme (amounting essentially to bound the
complexity of Hilbert's 16th problem) is merely a rather simple
theorem of function theory (due basically to
Riemann-Schottky-Klein-Teichm\"uller-Ahlfors). The latter in turn
being not much more than a bordered avatar of Riemann existence
theorem exhibiting any closed Riemann surface  as a branched cover
of the projective line. In other words any abstract Riemann
surface (\`a la Riemann-Prym-Klein-Weyl-Rad\'o) becomes a concrete
one (\`a la Abel-Galois-Cauchy-Puiseux-Riemann-Weierstrass, etc.)

Next it comes to make this Ahlfors-Rohlin fusion a reality. Here
starts alas some little difficulties which we hope to able to
overcome in the future. The problem breaks in two steps:

$\bullet$ Step~1. Make the Ahlfors theorem concrete by
specializing it to dividing plane curves, and conclude the
existence of a total pencil of ``adjoint'' curves exhibiting total
reality.

$\bullet$ Step~2. Prove that a totally real  pencil (often
abridged total pencil) implies maximality of the scheme.

Start with a scheme of type~I of degree $m$. This means by
definition that any curve $C_m$ representing it is of type~I. By
ROTOR (i.e. a concretized version of Ahlfors theorem), there is a
total pencil of curves $\Pi$ of order say $k$.

This means that each curve $C_k$ of the pencil cuts $C_m$ only
along real points [as soon as they are mobile]. [{\it Brackets
Added} [12.04.13].---Probably imaginary basepoint have to be
allowed, yet the mobile part of the pencil should be totally
real.] Such a balayage seems actually to supply a fast algorithm
to trace the real locus by reduction to a problem in one variable
(consider e.g. the G\"urtelkurve swept out by a pencil of lines),
which in turn could involve the Newton root-finding algorithm via
linearization. Of course in general,  curves of the pencil are not
rational and so this asks for a tricky extension of Newton.

Once  existence of a total pencil is granted we seems nearly
finished, because if there are more ovals then passing a curve of
the pencil through the additional point would corrupt B\'ezout.
The notorious bug of that argument is that we assume implicitly
the over-scheme to be realized by an augmentation of the given
algebraic curve, which is priori not the case. [12.04.13] Further
a total pencil exist as well on some type~I curves belonging to
indefinite (yet non-maximal) schemes, compare the case of sextics.
This is of course another obstruction to completing crudely the
just sketched argument.

However what makes that the argument works for the deep nest or
the satellites of the quadrifolium? Is it the fact that the pencil
curves are rational? or is it a sort of geometric intuition of the
pencil, or perhaps some canonicalness of it? Maybe to get the RMC
we need not only existence of a total pencil but some sort of
uniqueness (up to anodyne choices like center of perspectives
chosen in the ``deepest'' ovals).

It seems that a scheme of type~I incarnates a family of curves
which are so-to-speak totally real in some canonical way, and the
total pencil is virtually God-given. (Beware yet that the family
in question is not necessarily connected in the hyperspace of
curves, cf. Marin 1979 \cite{Marin_1979} or
Fig.\,\ref{Marin:fig}.)

At this stage it seems important to remember a metaphor allied to
total reality. Total reality means that all intersections are
visible on the reals. Using a pencil means essentially that we
choose a mode of vision of the curve. Basepoints are eyes of some
insect having several eyes and curves of the pencil are optical
rays enhancing how the animal perceives the curve (Gebilde). What
is strange is that total reality amounts saying that the vision is
purely transverse and so the object is in reality invisible (no
apparent contour). To make  this concrete consider the example of
the G\"urtelkurve $C_4$ (2 nested ovals) projected from an
innermost point inside the deepest oval. Paraphrasing in a real
life metaphor, looking at a glass of wine from outside you see its
apparent contour, but when placed inside of it, it suddenly
becomes invisible. This is total reality.

It is
tantalizing that total reality (via Ahlfors theorem) seems so
close to prove RMC(=Rohlin's maximality conjecture) but apparently
fails. As we (or better Rohlin) suspected the assumption that the
scheme is of type~I must impose the corresponding curves being
strongly harpooned by total reality. But how to make this idea
precise?

We can imagine the space of all curves representing the scheme,
and think about this a universal curve of type~I. There should
then be a version of Ahlfors theorem for family of curves or (as
 Teichm\"uller, Ahlfors, Bers, liked to say)
 for a {\it variable\/} Riemann surface.

The net effect would be that total reality is
genetically imbued in the curve(s) itself in such a strong fashion
that the scheme is maximal. So any scheme of type~I has a
canonical vision making it totally real, amounting essentially to
look at the world form inside the glass (or bottle) of wine. This
is akin to the photoelectric effect. (Compare with the known
examples of the unifolium and its satellites, alias deep nests in
the jargon of Hilbert and the Russian school, or the quadrifolium,
and its satellite total under a pencil of conics, plus the
(elusive) Rohlin-Le~Touz\'e phenomenon for sextics).

Once this photoelectric vision of the curve is given then nothing
more can appear in the blue sky and so the scheme is maximal. This
is the intuition of why  RMC holds true, but how to convert this
in a mathematical proof. What seems to be in
demand is a mechanism which from the shape alone of the scheme
identifies the total vision of the curve. This we call the
photorealism or photogenism.

If a scheme is of type~I then it is photogenic, and then it must
be maximal.

This seems to
request for a general mechanism of where to assign basepoints
which would extend the total reality of unifolium, quadrifolium,
and 9-folium of Rohlin-Le~Touz\'e flashed resp. by
 by pencil of lines, conics and cubics. As we
discussed in a earlier section the case of degree 8 schemes looks
a bit puzzling, where by the RKM-congruence we have plenty of
$(M-2)$-schemes of type~I (ca. 100 if we were not too bad in
counting). The center of vision (basepoints) are then quite hard
to predict. In general there are $B=M-3$ of them (where $M$ is
Harnack's bound), so $B=19$ for $m=8$, and alas it is presently
not very clear where to assign them in full generality. Making all
this explicit could solve the question of giving a precise sense
to our notion of photogeny, and as a by product crack the RMC.

Is this a realist strategy? Is there a more abstract argument? If
not, we really need some highbrow extensions of the
Rohlin-Le~Touz\'e theorem dictating us for all schemes of type~I
where to assign basepoints. This seems to call first for a
classification of schemes of type~I.

{\it Long (paragraph) Insertion} [12.04.13].---To tell the truth
it should be remarked that even in the case of Rohlin-Le~Touz\'e
(degree $m=6$) we lack presently a proof of the desideratum that
the vision of total reality (via the pencil of cubics) is strong
enough as to ensure maximality of the scheme. (This conclusion is
of course true via the Gudkov census but we lack a direct proof
along the philosophy of total reality.) Perhaps Rohlin knew a
proof, but as far as we know it was not published too.) Note two
things. First, the more naive principle of maximality of
Klein-Marin 1876/1988, when combined with Itenberg's contraction
affords another approach to the problem of RMC, which is perhaps
easier to implement (though in general only based on the
conjectural principle of contraction). Second, it seems that in
Rohlin's approach we lack some flexible medium to carry the
enlarged scheme to the original one. This could involve trying to
approximated a diffeomorphism of $\RR P^2$ by something more
algebraic (maybe a Cremona transformation), but then it looks hard
to finish the job. So maybe $A+B+C=Rmc^2$, i.e. Ahlfors, plus
B\'ezout, plus Cremona implies Rohlin's maximality conjecture. One
may also wonder if there is not a much more flexible proof of RMC
say akin to Rohlin's formula where merely soft topology is used
(while avoiding any contraction principle). We  have then
basically 2 curves, one enlarging the other, and one of which
universally of type~I. So one could fill the half \`a la Rohlin,
i.e. by all discs (like in the proof of Rohlin's formula) and
inspect the intersection of the homology class of degree $k$ with
the homology class of the enlarging curve. The sequel is certainly
hard to complete. (It is at this stage that we had the idea of
using Mangler to isotope the enlarged curve back to the original,
cf. Sec.\,\ref{RMC-via-Mangler:sec} for more detail on this
strategy to attack RMC via Ahlfors, plus Mangler.)

Here we know quite little, but as said earlier in this text, it
could be the case that the RKM-congruence is a universal detector
of $(M-2)$-schemes of type~I, while all other type~I schemes arise
as satellites of schemes of type~I of lower orders (dividing the
given degree $m$). So when $m=2p$ is twice a prime number there
should be no such satellite (except that of the unifolium) and all
new type~I schemes would be concentrated at the $(M-2)$-level. Of
course we can always make abstraction of the $M$-schemes where RMC
holds trivially true. So we see some sort of higher arithmetic
structure emerging in Hilbert's 16th problem as boosted by
Rohlin's conceptions, namely a sort of inductive process that
could progressively step-by-step enumerate all schemes of type~I,
merely as $(M-2)$-schemes of type~I satisfying the RKM-congruence
mod 8, or as satellite of earlier such schemes, and for all of
them expect a synthetical revelation  of the type by a canonical
pencil \`a la Rohlin-Le~Touz\'e incarnating primitive forms of the
phenomenon of total reality. This deserves nearly the name of
Rohlin's divination.

As the list of photogenic schemes increases at each step $m$, we
may conclude  RMC by having exhibited in some ad hoc fashion the
total reality of all type~I schemes, and the RMC would follow
step-by-steps. Needless to say this requires an immense effort,
and the induction required to validate RMC in all degrees looks a
priori extremely tricky.

Furthermore one could imagine that all this ascension effected in
autarchy from Ahlfors theorem by using rather ad hoc optical
recognition procedure via total reality. This would be parallel to
the evident total reality of the satellites of the unifolium
(alias deep nests) and idem for the quadrifolium, or  Rohlin's
schemes of degree 6 (modulo the lost proof of Rohlin). In that
case the theory would be purely Rohlinian and this is probably
essentially what Rohlin envisioned.

In contradistinction, when  attacking RMC, we know a priori the
scheme being of type~I so there could be some inference of Ahlfors
theorem permitting to shortcut the (pure) total reality vision of
Rohlin. This inference could increase  the (ascensional) speed
conceding some abstractness in the verification of RMC. Yet, as
observed, even this looks hard unless we get a better grip upon
the abstract total reality of Ahlfors.

A first modest (but nontrivial) exercise is to write down a clear
version of Rohlin-Le~Touz\'e's total reality claim, and using it
deduce the maximality of those 2 schemes. Here again notice that
exploiting the type~I assumption as do Le~Touz\'e is not a
concession since we are interested in RMC. So here total reality
seems sufficiently strong (canonical) to ensure maximality of the
schemes and we rederive so from Le~Touz\'e's result the
prohibition of Gudkov, etc. (compare
Table~\ref{Gudkov-Table3:fig}). In particular RMC holds true in
degree 6 for some intrinsic reason allied to total reality, as
opposed to being a byproduct of the full classification of Gudkov.

A more highbrow project would be to inject the function theory \`a
la Ahlfors in the problem. Assume given a scheme of type~I, we can
for each representing curve $C_m$ choose a total pencil $\Pi$
which is a line in the space $\vert kH\vert$ of $k$-tics curves.
It seems plausible that  the dependence can be made continuous.
Then we have a universal family of photoelectric effects on $C_m$
and its deformations (possibly in different chambers of the
discriminant) in which case the line  $\Pi$ may jump, a priori
even in different hyperspaces indexed by different $k$.

On applying the $k$-tuple Veronese embedding $v_k$---i.e. the
holomorphic map $\PP^2\to \PP^N$ induced by the linear system of
all k-tics---the total reality of $v_k(C_m)$ would appear under a
pencil of hyperplane, hence the curve would be total real under a
pencil and therefore located as several spires gyrating around the
base locus (plane of codimension 2 in $\PP^N$). Now it may be
expected that the phenomenon of total reality is as evident as it
was for the deep nest (i.e. reduction to the case of a linear
pencil) and that we may conclude maximality from B\'ezout (applied
of course now in the Veronese hyperspace).

After this little psychoanalysis of  Rohlin's
secret garden, we see that ``la r\'ealit\'e totale nous colle \`a
la peau.'' In some sense the phenomenon ought to be so inherent to
a curve belonging to a scheme of type~I that maximality of the
scheme should follow via the photoelectric effect. By the latter
we really mean that the total pencil being saturated nothing more
is allowed to appear in the blue sky without corrupting B\'ezout.
(Prototype: a deep nest with a pencil of lines through the deepest
oval, or satellites of the quadrifolium in degree $4k$.)

As a foundational detail, I always thought that possibly imaginary
basepoints have to be permitted in the definition of a totally
real pencil, so that merely moving points of the series are real
(cf. Gabard 2004, p.\,7). Now I am not sure that this is really
required.

In all basic examples of total reality (i.e. the deep nests
interpretable as satellites of the unifolium or the quadrifolium
and its satellites) the permissible basepoints ensuring total
reality are always varying through a contractible union of cells
as they are located inside the deepest ovals.  This is probably
also true for the sextics of Rohlin-Le~Touz\'e. If this is a
general phenomenon then this is a bit in line with our desideratum
that the total pencil ought to be almost canonically associated to
the dividing scheme of type~I. If instead we are interested in
total maps of lowest possible degree then we are inclined to let
degenerate the basepoints on the ovals themselves and so the total
pencils of this sort are parametrized rather by tori.

Another idea [developed in Sec.\,\ref{Thom:sec}] is to fill the
plane curve by the (orientable) membrane of $\RR P^2$, to get a
certain smooth surface in $\CC P^2$ whose fundamental class is
$kH\in H_2(\CC P^2, \ZZ)$. Smoothing its corner and applying
Thom's conjecture (=Kronheimer-Mrowka's theorem) could lead to
some interesting consequence. (More about this soon, cf.
Theorem~\ref{Thom-Ragsdale:thm}.) Of course this is basically
related to the ideas of Arnold and Rohlin.

Following our main theme, the idea would be that there is always
for a scheme of type~I some preferred (up to the ambiguity of a
contractible space of parameters) total pencil, which we call a
photon. This would naively speaking be obtained by assigning
basepoints among the deepest ovals. All this works good for
degrees $\le 6$. In degree $8$, the RKM-scheme 20 already affords
a little problem as quintics have 19 basepoints assignable and it
is not clear which ovals have to be used as ``anchor'' basepoints.
[But this scheme is prohibited by Thom, cf. again
Theorem~\ref{Thom-Ragsdale:thm}, or argue via Rohlin's formula.]

Once we have a photon (i.e. a canonical total pencil) then we
would like to argue that its satisfies the photoelectric effect,
and RMC would follow.

In the case of $(M-2)$-curves of type~I we have $B=M-3$ basepoints
for a pencil of $(m-3)$-tics that are freely assignable. This is
one unit less than the number of ovals and it is not clear which
one can be dispensed of being marked by a basepoint. We could
imagine that we could always dispense the oval whose porous inside
has the most negative Euler characteristic. To make this serious
compare Fig.\,\ref{RKM-schemes-deg-8:fig}, where we find however
schemes where such a dispensed oval is not uniquely defined, e.g.
$\frac{6}{1}\frac{6}{1}\frac{4}{1}1$. So our recipe is certainly
dubious.

We hope to have made the nature of the question clear enough. It
seems first that there is no direct reduction of RMC to Ahlfors
theorem, except perhaps if one as some deeper grasp upon the
geometry of a total pencil (photon) so as to ensure via the
photoelectric effect the RMC. As we said at the beginning of the
section, the net impact would be a sort of upper bound upon the
complexity of Hilbert's 16th problem. In fact it would really be
the clef de vo\^ute yielding some insights upon the architecture
of the pyramid of all schemes of some fixed degree $m$, namely
type~I schemes ought to be maximal element (yet not the sole ones
cf. Shustin 1985
\cite{Shustin_1985/85-ctrexpls-to-a-conj-of-Rohlin}). This
conjectural maximality explains nearly all prohibitions (at least
in the case $m=6$).

Despite our obsession to harpoon RMC via Ahlfors' total reality
(what Viro was sceptical about) and we have only the vague
suggestions of the present section to propose [see also maybe the
next Sec.\,\ref{RMC-via-Mangler:sec}], it may be argued that even
if it worked by a clever trick, it would perhaps not be as
satisfactory as the full and slow progression along the menagerie
of all schemes that must be tabulated along the higher order cases
of Hilbert's 16th problem.
In other words we would like to know the whole pyramid and not
just its maximal element. From the viewpoint of the French
revolution we would like to know the whole folk and not just the
aristocrats. More seriously we want to ``see'' the exact geometry
of the phenomenon of total reality, and not just its
capitalistical/hierchical impact via the photoelectric effect upon
a validation of the RMC. [12.04.13] Also interesting is the
question of the density of schemes below the aristocrats. Of
course it seems that the pyramids as dense below the maximal
elements: philosophically because an aristocrat seems unable to
provide for its wants without the force of all its servitors.

\subsection{An isotopic attack on  RMC via Mangler 1939}
\label{RMC-via-Mangler:sec}

[12.04.13] As we often experimented RMC would follow from Ahlfors
if we knew that the enlarged scheme lies in a tube neighborhood of
the given curve of type~I. The difficulty is  that a priori the
enlarging curve is very distant and hard to compare to the
original one. Crudely put one could expect to reduce always to the
easy case by using an isotopy.

So one could try to isotope the diffeomorphism taking the small
curve to its enlargement to the identity. Recall this to be
possible since the mapping class group of $\RR P^2$ is trivial)
[Mangler  1939 \cite{Mangler_1939}, often used by Teichm\"uller
1939]. One could then perhaps try to extend this isotopy to $\CC
P^2$ so as to get reduced to the case where the enlargement is a
small perturbation (actually the identity). In this case Ahlfors
suffice to imply RMC, since Ahlfors's pencil affords something
like a transverse structure, and one gets an easy corruption with
B\'ezout by letting pass a curve of the pencil through the new
oval (not within the tube neighborhood).

The difficulty is of course to check that B\'ezout (which is some
something algebraic rigid) is conserved during the very plastical
deformation of isotopy. Yet perhaps we may reinterpret
intersections homologically as to gain more flexibility. Further
the isotopy could be compatible with reality (equivariant and
respecting $\RR P^2$ and its complement of imaginary points).
Finally due to the geometric interpretation of intersection
numbers (in homology) their values will be clearly conserved by
the isotopy. So we arrive at the:

\begin{Scholium}
There is perhaps a trivial proof of RMC via isotopy of $\RR P^2$
equivariantly extended to $\CC P^2$, so that RMC reduces trully to
Ahlfors.
\end{Scholium}

\begin{proof}
Suppose $C_m$ to be a curve of degree $m$ belonging to a scheme of
type~I. Let $D$ be a curve (of the same degree) whose scheme
enlarges that of $C_m$. Fix a diffeomorphism $f$ of pair $(\RR
P^2, C_m)\to (\RR P^2, D_{\ast})$ where $D_{\ast}$ is $D$ less one
oval (w.l.o.g. or more ovals in general).

By Mangler 1939, we can isotope $f$ to the identity of $\RR P^2$.
Now it seems reasonable to expect that there is a natural way to
extend an isotopy of $\RR P^2$ to one of $\CC P^2$. This does not
need to be strongly unique but merely to exist in some sense that
it preserves real parts and maybe can be chosen equivariant w.r.t.
conj. I.e. carrying a point along the isotopy up to time
$t\in[0,1]$ commutes with the symmetry conj. Knowing that the
quotient $\CC P^2/ conj$ is $S^4$ could be of
valuable assistance to construct the extended-isotopy.

So we have $f_t$ an isotopy of $\RR P^2$ say with $f_1=f$ and
$f_0=id$, and $F_t$ and extension thereof to $\CC P^2$.  The map
$f$ pushes injectively the ovals of the first good curve $C_m$
into those of the hypothetical enlargement $D$. So operating
backward in time along the isotopy $F_t$ we may retract the
complexified curve $D(\CC)$ so that its real part becomes close to
that of $C_m$ (and even identic to it). Denote $D_0$ this
``temporal retraction'', which is  a ``flexible'' Riemannian
surface, with fundamental class still of degree $m$, by
homotopy-invariance of homology.

Now by total reality the first curve being of type~I it admits a
total pencil, all of whose members have B\'ezout saturated
intersections with the curve. Taking a curve $P_k$ of the pencil
passing through the additional oval of $D$ isotoped backward in
time ($t=0$), create one extra intersection (that will count
positively because the extended isotopy is orientation
preserving). All other other intersections also counts positively
if we are capable arranging the large isotopy $F_t$ to respect
somehow the normal bundle of $\RR P^2$. Hence the pull-back $D_0$
will have excessive number of intersections with $P_k$. The
homological B\'ezout (i.e. Poincar\'e, Lefschetz, etc.) is
therefore corrupted. Rohlin's maximality conjecture would be
proved by ``soft topology'' plus some Ahlfors.

The critique to this argument however is that a priori it applies
to any dividing curve supporting a total pencil and those can be
of indefinite type (yet not maximal), cf. the case of degree $m=6$
where there is plenty of such examples
(Fig.\,\ref{Gudkov-Table3:fig}). So a serious gap requests to be
filled. Maybe this will be an easy game for Alexis Marin?

To be optimistic, our argument looks so close to prove the big
desideratum (=RMC) that it is certainly worth exploring further.
In particular since the argument is spoiled by the objection of
indefinite schemes there must be (for instance in degree $6$ where
the Ahlfors total pencil are very easy to describe explicitly,
e.g. for the scheme $9$ where we have the simple
Fig.\,\ref{Fcubic:fig}) there must be some obstruction to extend
the Mangler isotopy to $\CC P^2$ (at least in a fashion that
positivity of intersections are conserved).

Understanding this obstruction, and assuming that one capable to
show that it vanishes if the scheme is of type~I could afford a
proof of Rohlin's maximality conjecture. The proof is likely to
involve some 4D-topology (say \`a la Marin-Siebenmann-Alexander).
\end{proof}

\subsection{Additional remarks on Rohlin-Le~Touz\'e
total reality for $(M-2)$-sextics of RKM-type}

[04.03.13] We concentrate again on the case of a sextic $C_6$ of
type $\frac{6}{1}2$. Before entering into the elusive technical
details we recall the basic problematic yet quite elusive for the
moment. We would like to show a phenomenon of total reality, yet
its exact shape is still obscure to us.

Either we can use the RKM-congruence to infer a priori that the
curve is of type~I. Then we could apply the (abstract) result of
either Ahlfors 1950, or Gabard 2006 and hope to effect a descent
in the plane, to get a total pencil (hopefully of cubics). This
descent probably requires some theory \`a la Brill-Noether as a
Plato-cavern-style reflection of Riemann's work. Alternatively, we
can try to follow the (direct concrete) route proposed by
Le~Touz\'e 2013
\cite{Fiedler-Le-Touzé_2013-Totally-real-pencils-Cubics} depending
upon a detailed analysis of pencils of cubics.

Finally we could dream to recover the lost proof of Rohlin, i.e.
without assuming the dividing character  while trying directly to
ensure total reality of such a curve under a pencil of cubics with
suitably assigned $8$ basepoints. Here the problem is that we know
presently very little of how general the phenomenon of total
reality is, i.e. which octuplets induce total reality. More about
this soon.

A very first (extremely basic remark) that we all use
subconsciously is a sort of closing lemma for algebraic curves:

\begin{lemma}\label{Closing-lemma:lem} {\rm (Closing lemma)}.---Given any real plane projective curve, its real locus (if non
empty) consists of closed circuits (Jordan curves in $\RR P^2$) or
possible poly-cycles like figure 8, etc, or eventually an isolated
singularity. Topologically it is always locally a multi-node
consisting of a certain number of branches crossing transversally,
or an isolated point. When the point of the real curve is
non-isolated then there is at least one Jordan curve based on the
given point.
\end{lemma}

\begin{proof}
This is a mixture of algebraic-geometry, implicit function
theorem, and topological compactness of $\RR P^2$, and abstract
classification of  compact (Hausdorff) $1$-manifolds (just the
circle), and some singular ``graph'' avatars.
\end{proof}

This pertains to the very first step of our total reality story as
follows:

\begin{lemma}
Assume given 8 basepoints distributed on the empty ovals $E_i$ of
our $C_6$. Let $\Pi$ be the pencil of cubics passing through the 8
points. Then all curves of $\Pi$ cut the curve $C_6$ in at least
16 real points (i.e. just 2 units less than the maximum
permissible by B\'ezout). Say in that case that the pencil is
quasi-total.
\end{lemma}

\begin{proof}
Our 8 basepoints forces 8 real intersections, but by the closing
lemma each intersection has at least one companion (possibly the
same yet then with a tangency and so counted with multiplicity 2).
More precisely we look at one of the basepoint and choose any
curve of the pencil. First note that the basepoint cannot be an
isolated real point of the curve $C_3$, and so there is by the
closing lemma a topological Jordan curve in the real locus of
$C_3$ which has thus at least 2 intersections with the given oval.
\end{proof}

So to reach Rohlin's (lost) theorem just 2 real intersections are
missing (quasi-total) but the gap toward total reality is still
immense.

Let us first observe the following extension where basepoints are
located inside the ovals, as opposed to the former case where they
were directly imposed on the ovals themselves:

\begin{lemma} (Interior distribution quasi-total).---
If the $8$ basepoints are assigned in the insides of the empty
ovals $E_i$, then the pencil is also quasi-total, i.e. $C_3\cap
C_6$ has $16$ real intersection for all $C_3\in \Pi$.
\end{lemma}

\begin{proof} Suppose given such a $C_3$, hence visiting
the basepoints $p_i$ labelled as to be in the insides of the
$E_i$. Each $p_i$ forces 2 intersections in $C_3\cap C_6$, except
if the circuit of $C_3$ through $p_i$ is a small oval inside
$E_i$. But then the residual pseudoline of $C_3$ has to visit all
7 remaining points, and so is forced to intercept $N$ the nonempty
oval of $C_6$. We count then $14+2=16$ intersections.

{\it Addendum}.---Further in the above situation of a small oval
of $C_3$ inside $E_i$, then for vibratory reasons the 9th
basepoint of $\Pi$ has to be on it.

Hence if $p_i$ is an inner point of $N$,  then (even) total
reality is fulfilled (Lemmas \ref{9th-basepointI:lem} and
\ref{9th-basepoint-totalII:lem}).

If $p_i$ is an outer point then the residual pseudoline of $C_3$
will intercept (twice) $N$, and we have again $14+2=16$ (real)
intersections.
\end{proof}

The above {\it Addendum} is not formally required for the proof of
the lemma but we include it as it nearly give a hope to attack
Rohlin's total reality claim (abridged RTR in the sequel). Beware
that Rohlin's statement is very loose in the original paper
(Rohlin 1978 \cite{Rohlin_1978}) and so RTR should not be given a
too strong connotation from the scratch. Part of the problem is to
decide with which level of generality  Rohlin's claim is correct.

With our zero knowledge, we can distinguish several layers of
interpretation for RTR. In its strongest form this would be the
assertion:

\begin{conj} {\rm (Inside total reality)=ITR}
Denote by $E_i^{\ast}$ the (sealed) insides of the empty ovals
$E_i$ of the $C_6$, then the pencil $\Pi$ of cubics through any
points $(p_1,\dots, p_8)\in E_1^{\ast}\times \dots\times
E_8^{\ast}$ is totally real.
\end{conj}

Then there are several weaker variants, namely the same conclusion
under the assumption that the $p_i$ belong to the ovals $E_i$
themselves. The corresponding statement is called OTR, for oval
total reality. Another weakening is to relax the conclusion by
claiming only total reality of $\Pi$  for a suitable octuplet, as
opposed to claiming it for all of them. This relaxed form induces
statements  called WITR resp. WOTR, where the ``W'' stands for
weak total reality. Though being weak this would be enough to
geometrize the degree $m=6$ case of the RKM-congruence.
Le~Touz\'e's theorem is essentially WITR, i.e. weak inside total
reality modulo the fact that Le~Touz\'e assumes (or infers from
the RKM-congruence) the dividing character of the curve.

Of course we have formal implications like the following
commutative square
$$
\begin{array}{rcl}
\textrm{ITR} & \Rightarrow & \textrm{OTR}\\
\Downarrow \quad &      &\quad \Downarrow\\
\textrm{WITR} & \Leftarrow & \textrm{WOTR}.\\
\end{array}
$$
Alas we know very little about those statements. We do not know if
the strong versions (upper row) are true, and if foiled it could a
priori still be the case that they hold true for special sextics
$C_6$. To be factual at the time of writing (and modulo an
understanding of Le~Touz\'e's proof) the only available knowledge
is that the weakest form WITR holds true, and even as we said
under the assumption that the curve is of type~I as may be
inferred from the RKM-congruence. Hence of course the whole square
can be extended to a cube with another square face of statements
assuming the dividing character of the curves. Modulo RKM-both
squares are actually formally equivalent, but a very purist could
prefer eliminating this dependency. As asserted (but never proved)
by Rohlin 1978, one could hope to do more and prove one of the
above statement {\it ex nihilo\/} (without reliance upon RKM).

As usual in mathematics (or in the world of bird of preys) one
should always start  attacking the weakest prey, namely WITR. This
is a bit strange because quasi-total reality is slightly easier to
establish when the basepoints are located on the ovals. Further
keep in the subconscious part of the brain, that Rohlin's hints
are so vague that it is not even clear that our 4-fold strategy
covers all what is permissible (for instance it could be useful to
assign the basepoint not on the empty ovals but one also on $N$.
This looks exotic, but perhaps  useful in extreme case of
desperation).

So what is a reasonable strategy toward WITR?

We may start from the observation that the pencil $\Pi$ is totally
real iff the 9th basepoint $p_9$ of $\Pi$ is in the sealed inside
$N^{\ast}$ of the $C_6$ (cf.
Lemma~\ref{9th-basepoint-totalII:lem}). In reality this lemma
holds true for basepoints assigned on the ovals, but probably
extends to the broader setting. This leads to the following:

\begin{conj} (Hypothetical lemma)
Assume $p_9$ to be in $N^{\ast}$ and the $p_i\in E_i^\ast$ in the
(sealed) insides. Then $\Pi$ is totally real (abridged total).
\end{conj}

Let us show where the naive proof breaks down. Assume given any
curve $C_3$ of the pencil. A priori $C_3$ may pass through an
inner point $p_i$ (i.e. inside $N$) via a microscopic oval $E$ of
$C_3$ entirely inside $E_i$, thereby creating no real
intersections. Of course then the residual pseudoline of $C_3$
(i.e. $J=C_3(\RR)-E$) intercepts (twice) the nonempty oval $N$,
but this affords altogether only $14+2=16$ real intersections.
Note that $p_9$ has for vibratory reasons necessarily to be
located on
$E$, yet this is no contradiction. Perhaps I missed something and
there is a more clever argument establishing this modest technical
conjecture. (Le~Touz\'e probably has some idea.)

Let us skip this conjecture, while attacking rather the stronger
looking WOTR proposition, as in the latter case total reality is
easier to ensure. Of course doing so we loose some freedom for the
parameters as the large  16-dimensional cell
$E^\ast=\times_{i=1}^8 E_i^\ast$ is traded against the
8-dimensional torus $E=\times_{i=1}^8 E_i$ but perhaps this
suffices to conclude. Further the advantage would be to get a
total map of lower degree, namely one corroborating Gabard's
bound.

So it is a delicate matter to decide which strategy ``insides of
the oval versus the ovals themselves'' is more likely to give a
proof of RTR (=Rohlin's total reality claim).

For the moment we have no better idea than the topological
approach sketched in one of the previous section, i.e. to ensure
that the 9th basepoint lands in $N^{\ast}$, and so abort this
delicate question.

\section{Thom's conjecture vs. Hilbert's 16th}
\label{Thom:sec}

[21.03.13] {\it Warning.}---All this Sec.\,\ref{Thom:sec} is
poorly organized for reasons to be soon explained. In particular
it contains several mistakes, but also such fundamental results as
Petrovskii inequalities,  the strong-Petrovskii-Arnold
inequalities. Some higher Gudkov tables of periodic elements (e.g.
Fig.\,\ref{Degree10:fig}) show the geographical impact of
Petrovskii-Arnold as compared to Ragsdale's conjecture (briefly
discussed in Sec.\,\ref{Ragsdale-conj:sec}). The importance of the
Petrovskii-Arnold results was pointed out to me by Thomas Fiedler,
who corrected several benign mistakes and one much more fatal bug
of mine. This section should thus be read with extreme
discernment, as it mixes both the best (Petrovskii-Arnold and even
the marvellous construction of Itenberg-Viro) and the worst
(Gabard). Several footnotes and WARNINGS should aid the reader to
avoid going into the same pitfall as I did. All those WARNINGS are
due to kind letters of Fiedler who fixed all my misconceptions and
posed me challenging problems. We hope in the future to be able to
reorganize the text in a more decent fashion after exploring in
more depth a possible fascinating interplay between Hilbert,
Ragsdale, Thom=Kronheimer-Mrowka 1994
\cite{Kronheimer-Mrowka_1994} (independently Morgan-Szab\'o-Taubes
1995/96 \cite{Morgan-Szabo-Taubes_1996}), and the work of
Petrovskii-Arnold (1938--1971). Of course the interested reader is
invited to consult more professional sources,
notably Mikhalkin 1994
\cite{Mikhalkin_1994-adjunction-Thom}.

[21.03.13] This section was built around the fundamental result
$\chi\le k^2$ for any dividing curve of even degree $2k$
(Theorem~\ref{Thom-Ragsdale:thm}) directly inferred from the
so-called {\it Thom conjecture\/} (which he humbly considered
himself as rather belonging to the folklore, compare footnote in
Lee Rudolph 1984 \cite{Rudolph_1984}). This should have implied a
clear-cut impact of Thom upon Hilbert's 16th problem. The summit
of our fictional ``Gabard-Thom'' theory went so far as to
establish one-half of Ragsdale's conjecture, (still open for
$M$-curves) (cf.
Lemma~\ref{Thom-implies-one-half-of-Ragsdale:lem}) and to show
that Ren\'e Thom was on his 31, i.e. can be stronger than the
conjunction of all Russian estimates, congruences and formulas
(due primarily to Petrovskii 1938, Gudkov 1969, Arnold 1971 and
Rohlin 1972--74--78), cf.
Theorem~\ref{Alsatian-scheme-Thom-strong-Petrov-Arnold:thm}. This
would have refuted a belief of Th. Fiedler (cf. his letter ca. 13
March in Sec.\,\ref{e-mail-Viro:sec}).

$\bigstar\bigstar\bigstar$ Fortunately, Fiedler brought us back to
reality by showing that our reasoning is wrong as it overlooks the
issue that despite being constructed by pasting two orientable
pieces---namely Klein's orthosymmetric half married with  Miss
Ragsdale's membrane bounding the curve from inside---the so-called
{\it Arnold surface\/} (1971) does {\it not\/} need to be
orientable. My mistake is thus nearly as basic as having
overlooked that one  can create (like in Klein's bottle)
non-orientable objects merely by pasting a handle to itself in a
twisted fashion (this reminds me some lovely pictures in the
Fuks-Rohlin ``Beginner's course on topology'').

This being confessed, most of this section is foiled and  we are
much indebted to Th. Fiedler for having catched our mistake at the
right moment and stimulated our investigations.  Albeit much of
the sequel is foiled we have decided to keep it for didactic
reasons.
In our case it was so pathetic to write ca. 40 pages based upon a
misconception without noticing anything (prior to Fiedler's
correction) that we would by no mean that somebody else do the
same mistake.

More positively many questions arises through Fiedler's
correction.

(1) Where (in particular in which degree) lives the first
counterexample to the Gabard-Thom bound $\chi\le k^2$ (no false
modesty in calling it so since it is false) for all dividing
curves?
[The sole counterexample I know is the Itenberg-Viro curve
corrupting Ragsdale. This is a beautiful picture, see
Fig.\,\ref{Itenberg:fig}.]

(2) If ``Gabard-Thom'' is obviously false for theoretical reasons,
why does it look nearly true as implying one-half of the vestiges
of Ragsdale's conjecture?

(3) Under which condition is the estimate $\chi\le k^2$ still
true? As noted by Fiedler, the answer seems rather clear namely
iff the Arnold surface is orientable. This is in turn the case iff
all primitive pairs are positive. The last condition is perhaps
not an ``iff''. In crude approximation one could say that the
Arnold surface is orientable iff the Rohlin tree is positively
charged throughout. At least if this is the case (a rather
stringent condition) then Arnold's surface is orientable, and the
estimate  $\chi\le k^2$ holds true (cf. the limpid proof of the
erroneous Theorem~\ref{Thom-Ragsdale:thm}). Further it seems
evident that when the tree is positively charged throughout then
the Rohlin mass $\pi-\eta$ is maximized (We always set
$\pi:=\Pi^+$, $\eta:=\Pi^-$ to abridge Rohlin's notation). Alas
even that is false (more details soon).

So Thom's conjecture has still something to say on Hilbert's 16th
yet its impact is much more subtle than expected when doing the
fundamental mistake.

Finally, all this section on ``Thom'' contains above all
comparative study with the strength of Rohlin's formula
$2(\pi-\eta)=r-k^2$ (see (\ref{Rohlin-formula:thm})), hence a
cuneiform formalism that really pertains to Rohlin's complex
orientations as it electrifies the Hilbert tree by putting charges
(distributions of signs on the edge). Albeit our exploration of
this topic was completely
random (and biased by our erroneous Thom estimate), it should have
some independent interest. Reorganizing all this material, without
loosing any bit of information will take us several weeks, and
cannot be done on the present edition. We hope in the future being
able to give a more structured exposition of this cuneiform
formalism (Hilbert's tree with signs, alias Rohlin's trees) and
our messy account can motivate others to clarify this.

\subsection{What can be salvaged after
Fiedler's earthquake?}

[22.03.13] As spotted by Fiedler,  we overlooked  that  Arnold's
surface (arising by pasting Klein's half with Ragsdale's membrane)
is not necessarily orientable though both its constituents are.
Hence one cannot apply Thom so straightforwardly. Incidentally we
hope that Thom applies without trouble coming from the necessity
of rounding  corners. This is folklore but it would be nice to
find adequate reference (Thom?, Cerf?, Hirzebruch? Milnor?, Wall?,
etc.)

Though Thom was a heuristic way to discover the inequality
$\chi\le k^2$ it could be  that it holds true for more elementary
reasons directly rooted in Rohlin's formula. Let us briefly
explain how. [$\bigstar$ Non-sense (!) by the Itenberg-Viro
counterexample in degree $10$, Fig.\,\ref{Itenberg:fig}.]

In the sequel we use the jargon of the {\it Rohlin tree} which is
simply Hilbert's tree (encoding the distribution of ovals via a
POSET whose order relation comes from the inclusion of the insides
of the ovals), plus a decoration of its edges by signs coming from
Rohlin's complex orientations of a curve of type~I (in the sense
of Klein, also called latter by him orthosymmetric curves). So if
any plane real curve has a Hilbert tree, dividing curves have an
extra distribution of signs on the primitive edges (those of
length one) which by the signs-law (\ref{Signs-law:lem})
propagates consistently along the whole injective pairs of the
tree.

First if Rohlin's tree is positively charged (i.e. if all
primitive pairs of ovals are positive in the sense of Rohlin) then
Arnold's surface $A=C^+\cup R$ (=Klein's half $C^+$ pasted with
Ragsdale's (orientable) membrane $R$) is orientable too! In that
case via Thom we have $\chi\le k^2$ (cf.
Theorem~\ref{Thom-Ragsdale:thm} and its simple proof). [$\bigstar$
Okay, but as we shall soon see this may also be inferred from
Rohlin's formula!]

(1) Moreover for a positively charged tree we have $\pi-\eta=n$,
i.e. the Rohlin mass $\pi-\eta$ is equal to the number $n$ of
negative=odd ovals (compare
Lemma~\ref{Rohlin-mass-of-a-positively-charged-tree:lem} below).

(2) Rohlin's formula reads $2(\pi-\eta)=r-k^2$, hence fixes the
Rohlin mass $\mu:=\pi-\eta$ which is something coming from the
``complexification'' (i.e. the Riemannian) in terms of real
characteristics ($r$ being the number of ovals and $k=m/2$ the
semi-degree). This formula can be used to express the
Euler-Ragsdale characteristic $\chi=\chi(R)$ as follows
\begin{align}\label{Rohlin-to-Arnold-tris:eq}
\chi=p-n=(p+n)-2n&=r-2n\cr
&=[2(\pi-\eta)+k^2]-2n\cr
&=k^2+2[(\pi -\eta) -n].
\end{align}

Combining (1) and (2) gives the following lemma (incidentally
remarked in Fiedler's letter):

\begin{lemma}
If all primitive pairs are positive (equivalently if the Rohlin
tree is positively charged) then $\chi=k^2$.
\end{lemma}

Now it seemed to us realist to posit a sort of ``positive mass
conjecture'' (POSMASS) stipulating that the Rohlin mass
$\mu:=\pi-\eta$ (of a signed tree) is maximized when all primitive
edges/pairs (of the tree) are positively charged. (We formulated
this some 7 days ago, cf. optionally
\ref{positive-mass-conjecture:conj} dated [15.03.13].)

This looks quite appealing, albeit we have little evidence for the
truth of this principle which if optimistic could be pure
combinatorics (i.e. valid for any signed tree with charges
respecting the signs-law) or be merely valid for such Rohlin
signed trees arising via  dividing curves. Of course it could be
also false in this restricted case. Our evidence for POSMASS is
presently only derived from the case of chains of length $\le 4$
(or so), compare the signs-law for dyads
(Fig.\,\ref{Signs-law-dyad:fig}), that for triads
(Fig.\,\ref{Signs-law-triad:fig}), and that for tetrads
(Fig.\,\ref{Signs-law-tetrad:fig}).

If POSMASS is true, then $\pi-\eta\le n$ and so the ``Euler-Rohlin
formula''  \eqref{Rohlin-to-Arnold-tris:eq} implies  $\chi\le
k^2$, i.e. the so-called Gabard-Thom (dubious) estimate. The
striking issue is that the Gabard-Thom theorem would be still true
but merely as a logical consequence of Rohlin's formula (plus some
combinatorics required to validate POSMASS). In particular the
intervention of Thom could be completely dispensed.

This scenario looks  risky, since Fiedler claims our Gabard-Thom
theorem to be wrong. However, it is[=was] not clear to me if there
is really a counter-example to the theorem, or if Fiedler just
stated wrongness of our reasoning.

More ironically I forgot to remember that even 5 days before
formulating POSMASS, I had a simple counterexample to it in the
combinatorial setting (cf.
Theorem~\ref{Garidi-mass-conj-is-FALSE:thm}). Hence there is no
chance to prove the Gabard-Thom estimate via pure combinatorics.

Of course, it is also likely that the POSMASS conjecture is false
also for Rohlin's trees arising as dividing curves but that
deserves an explicit example. Of course the method should be to
ape algebraically the combinatorial structure of a bat\^onnet that
foils the mass conjecture (cf.
Fig.\,\ref{Garidi-mass-false:fig}a). A b\^atonnet is merely a
usual tree with a trunk that ramifies strongly into several
branches at depth 2 (look at that picture
Fig.\,\ref{Garidi-mass-false:fig}a.)

Accordingly, my first idea was to look back in Gabard 2000
\cite[Fig.\,13, p.\,154]{Gabard_2000} where is traced a classical
(variant of) Hilbert's construction of an $M$-curve of degree 8,
which has nearly the required b\^atonnet structure. We shall soon
reproduce this and related pictures. Moreover why degree 8? Simply
because in degree 6 the Thom-Gabard estimate $\chi\le k^2$ (for
dividing curves) is trivially true as follows by glancing at
Gudkov's table (=Fig.\,\ref{Gudkov-Table3:fig}) of which we merely
use Hilbert's intuition/theorem that the unnested $M$-scheme
(symbol~$11$) is not algebraic, plus the fact that the unnested
$(M-1)$-scheme (symbol~$10$) has no dividing realization (as
follows from Klein's congruence $r\equiv_2 g+1$).

This being said let us do an iterated Hilbert construction
(Fig.\,\ref{HilbGab1:fig}). This gives first the well-known
$M$-sextic $C_6$ of Hilbert (symbol $\frac{9}{1}1$), and then an
$M$-octic $C_8$ with $\chi=16=k^2$, and then an $M$-curve of
degree 10, $C_{10}$ with only $\chi=9$, yet still $\le k^2=25$
(and congruent to it mod 8 as it should by virtue of Gudkov
hypothesis). Of course the Gabard-Thom estimate $\chi\le k^2$ has
little chance to be corrupted so, since it formally implies
one-half of Ragsdale's conjecture (which is still open for
$M$-curves, cf.
Lemma~\ref{Thom-implies-one-half-of-Ragsdale:lem}). A more
intrinsic reason is of course that Ragsdale conjectures were
calibrated along a deep contemplation of the Harnack-Hilbert
method. Hence historical continuity is fighting against our
attempt to corrupt Gabard-Thom via Hilbert's construct. As the
(maximal) $M$-Ragsdale conjecture is still open (not succumbing
even to the Viro-Itenberg method, cf. Itenberg-Viro 1996
\cite{Itenberg-Viro_1996-disproves-Ragsdale}), it is much more
likely  to corrupt Gabard-Thom by using non-maximal curves. Keep
this idea in mind for later.

\begin{figure}[h]
\hskip-2.8cm\penalty0 \epsfig{figure=HilbGab1.eps,width=182mm}
\vskip-5pt\penalty0
  \caption{\label{HilbGab1:fig}%
Hilbert construction} \vskip-5pt\penalty0
\end{figure}

Before adventuring outside the realm of $M$-curves, let us do also
a plate for Harnack $M$-curves constructed \`a la Hilbert
(Fig.\,\ref{HilbGab2:fig}). Those Harnack-style curves  have a
priori a larger $\chi$, so better suited to corrupt Gabard-Thom.
Precisely Harnack's $M$-sextic has $\chi=9=k^2$, then we get a
$C_8$ with $\chi=16=k^2$, and then a $C_{10}$ with $\chi=25=k^2$.
Our depicted $C_{12}$ (right) is not the most natural choice as we
switched to an ``internal vibration'', while in the first steps
$C_6\rightsquigarrow C_{8}\rightsquigarrow C_{10}$ we consistently
opted for an external oscillation (typical of Harnack's curves
reckoned \`a la Hilbert).   It should be noted than even the
natural choice of an external vibration does not lead to a
$C_{12}$ with $\chi=36$, but we were only able to get one with
$\chi=28$ (in accordance with Gudkov hypothesis). [{\it Added in
proof}.---THIS IS A MISTAKE, due to the fact that I reported ovals
at the wrong place] However Miss Ragsdale in 1906 was surely much
more clever than we are [THIS IS NOT PERTINENT ANY MORE], and I
presume she was able to reach always the ``Gabard-Thom
upper-bound'' $\chi\le k^2$, since this really amounts to
(one-half of) her conjecture (cf. again
Lemma~\ref{Thom-implies-one-half-of-Ragsdale:lem}).

\begin{figure}[h]
\hskip-2.8cm\penalty0 \epsfig{figure=HilbGab2.eps,width=182mm}
\vskip-5pt\penalty0
  \caption{\label{HilbGab2:fig}%
Harnack curves constructed \`a la Hilbert} \vskip-5pt\penalty0
\end{figure}

[23.03.13] Correcting my mistake, an extension of this figure
(Fig.\,\ref{HilbGab2:fig}) shows the following result (compare
Hilbert 1891, and Ragsdale 1906):

\begin{lemma}
There is an infinite series of $M$-curves of degree $2k$ with
$\chi=k^2$. Hence if Ragsdale's conjecture $\chi\le k^2$ (for
$M$-curves) is true, then it is sharp.
\end{lemma}

\begin{proof} Look at the first steps of Fig.\,\ref{HilbGab2:fig}
and do not commit the mistake of making an inner vibration, but
choose always external vibrations as we did for $C_6, C_8$. On
looking at the Hilbert trees of $C_6, C_8, C_{10}$ one easily
derive the general evolution of the Hilbert tree of $C_{2k}$.
Namely the number of outer ovals $9,17,27$ augments along the
increment $+8, +10, +12, etc.$, while the  tree itself is always
pushed down one step deeper while acquiring new branches on its
top, compare the windows on the figure for $C_6, C_8, C_{10}$. In
view of the extreme regularity of the construction it is easy to
extrapolate the nested structure of $C_{2k}$. Writing down a
general formula looks not even necessary, and the Gudkov symbol
will be something like
$$
(1,(2k-6) (1, 2k-8 (1, 2k-10) \dots )) (9+8+10+\dots+2k).
$$
In more geometric terms this means that the Hilbert tree of
$C_{2k}$ has $9+8+10+\dots+(2k-2)+(2k)$ outer ovals, and a trunk
of length $2k-7$ with branches hanging on. It should be an easy
matter to compute directly the Euler characteristic of $C_{2k}$ by
the evolution rule for the tree.
\end{proof}

[23.03.13] Perhaps instead of using Hilbert's method one must
really uses the more time-consuming Harnack original method which
amounts to oscillate around a ground-line instead of the ellipse
used in Hilbert's method. (Incidentally, I wonder if one wants a
fast-Hilbert method, if it is possible to vibrate across a split
quartic, union of 2 ellipses.)

All this is good but will perhaps only confirm the intuition (of
Ragsdale) that her estimate $\chi\le k^2$ is sharp for $M$-curves
(if true at all). Our object is somewhat different namely to
refute the Gabard-Thom estimate $\chi\le k^2$ for all dividing
curves. Of course it could be possible to disprove the Ragsdale
$M$-estimate $\chi\le k^2$ yet this deserves highbrow methods as
even the powerful Viro-Itenberg method apparently failed as yet
 in that game (compare Itenberg-Viro 1996
\cite{Itenberg-Viro_1996-disproves-Ragsdale}).

So how to construct non-Harnack-maximal dividing curves? As we
know from degree 6 (Fig.\,\ref{GudHilb8:fig}), the trick is just
to disperse the vibratory energy on several ovals, as opposed to
the monopole of Hilbert's vibration where a single oval is
oscillating. By Klein's congruence we look at $(M-2)$-curves, and
in degree 6 the specimen with largest $\chi$ is the scheme with
symbol $9$, which can be obtained by (a variant of) Hilbert's
method where the vibration is dispatched on 2 ovals (cf. the
earlier Fig.\,\ref{GudHilb8:fig} or right below
Fig.\,\ref{HilbGab3:fig}). It is now hoped that when iterating the
construction to higher degrees we get a refutation of Gabard-Thom.

Applying this idea  we get the following series of $(M-2)$-curves
(Fig.\,\ref{HilbGab3:fig}), all of type~I by Fiedler's law of
smoothing
dictated by (and dictating) complex orientations. Of course at the
higher steps we choose again a monopolized vibration as otherwise
we descend further the energetic level and reach curves with less
ovals that $(M-2)$. We use now a quicker depiction mode where
vibration and smoothing are depicted on the same plate at each
step of Hilbert's inductive process. Further on the diagram
standing right below each curve, we depict Hilbert's nested tree
encoding the distribution of ovals, and some easy calculation of
topological characteristics of the curves so constructed. The
conclusion is that we get an infinite series of $(M-2)$-curves
with $\chi=k^2$ as shown by the figure up to $k=7$, and the
regularity of the procedure is so evident that this property
easily follows for all $k$. In particular it may be observed that
the number of outer ovals evolves along the progression $9, 17,
27, 39, 53,\dots$,  incrementing along the progression
$+8,+10,+12, +14, \dots$, while the nested portion of the tree is
simply pushed down at each step, with $(2k-6)$ new branches
arising on the top.

\begin{figure}[h]
\hskip-2.8cm\penalty0 \epsfig{figure=HilbGab3.eps,width=182mm}
\vskip-5pt\penalty0
  \caption{\label{HilbGab3:fig}%
A series of dividing $(M-2)$-curves with $\chi=k^2$}
\vskip-5pt\penalty0
\end{figure}

Hence,
modulo some arithmetical nonsense, we have proved the:

\begin{lemma}
There exists an infinite series of dividing $(M-2)$-curves
$C_{2k}$ of degree $2k$ with $\chi=k^2$.
\end{lemma}

Alas, this does not refute the Gabard-Thom estimate, but rather
show its sharpness in case it would be correct. (Another proof of
the sharpness can be derived by perturbing ellipses, cf. Remark
right after Thm~\ref{Thom-Ragsdale:thm}).

(The contrast between the present regularity and the lack thereof
for $M$-curves was so striking that it permitted us to correct the
earlier mistake that we did above.)

But where to find a counter-example to the Gabard-Thom-bound (as
promised by Fiedler's claim of erroneousness)?

Lacking some imagination let us redo the Hilbert vibration for
$M$-curves more systematically by always vibrating from
``outside''. This gives Fig.\,\ref{HilbGab4:fig}. The same
regularity is observed while the general pattern becomes evident
after some few
iterations. Hence one can reduce to the depiction of the trees.
The latter get always more profound by one unit as $k$ increments,
while the number of deepest ovals belongs to the series $9, 17,
27, 39, \dots$ which regularly increments by $+8,+10, +12, \dots$
so that one can predict the future evolution of the tree, compare
the very bottom row. This shows that one times over two we will
attain the Gabard-Thom bound $\chi=k^2$, while of course the sign
of $\chi$ oscillates between negative and positive values. It is
evident that we will not get a counterexample to Gabard-Thom in
this fashion.

\begin{figure}[h]
\hskip-2.8cm\penalty0 \epsfig{figure=HilbGab4.eps,width=182mm}
\vskip-5pt\penalty0
  \caption{\label{HilbGab4:fig}%
A series of Hilbert $M$-curves with $\chi=k^2$ once time over 2
(oscillating Euler-Ragsdale's $\chi$)} \vskip-5pt\penalty0
\end{figure}

Then we can still vary more the constructions, e.g. by starting
with a Harnack-like outer vibration of the $C_6$ like on
Fig.\,\ref{HilbGab2:fig} and perform the dissipation (leading to
$(M-2)$-curves) at the next step on the octic $C_8$. Alas we still
found $\chi=k^2$ for the $C_8$ and even the $C_{10}$ (details of
the picture on p.\,AR-114 of my hand-notes).

At this stage one gets a bit depressed. It seems hard to corrupt
Gabard-Thom by construction \`a la Hilbert-Harnack. Maybe I missed
something, or perhaps one should appeal to more sophisticated
constructions like Viro-Itenberg. [This turned out to be the good
idea, more soon.] At this stage I cannot therefore preclude the
option that the Gabard-Thom estimate is true (of course for
another reason than the gapped proof given in
Theorem~\ref{Thom-Ragsdale:thm}).

So a priori 4 scenarios may happen:

(1).---Either there is an elementary counterexample to
Gabard-Thom(=GT) via elementary constructions \`a la
Harnack-Hilbert (and we missed it due to lack of cleverness).

(2).---There is a counterexample to Gabard-Thom via highbrow
constructions \`a la Viro-Itenberg. In particular it is not
impossible that the counterexample described by them to Ragsdale
conjecture also supplies a counterexample to Gabard-Thom. This
requires controlling the type in their construction (which is an
issue known to them). More on this in the next
Sec.\,\ref{Itenberg-Viro-disprove-Gab-Thom:sec}.

(3).---There is a counterexample to GT via another construction
not covered by Viro-Itenberg. (To my knowledge there is no theorem
stating that any algebraic curve is constructible via their
method.)

(4).---GT is true for another reason  than the one exposed in
Theorem~\ref{Thom-Ragsdale:thm}. This would however conflict with
Fiedler's assertion that our theorem is wrong. (Again, it is not
clear if Fiedler merely stated wrongness in the proof or of the
statement.)

The next section gives a clear-cut answer to this puzzle.

\subsection{A formal disproof of Gabard-Thom via Itenberg-Viro's
patchwork and Kharlamov-Marin}
\label{Itenberg-Viro-disprove-Gab-Thom:sec}

[23.03.13] Here we say ``formal'' just because we are not yet
familiar with the patchwork method due to Viro, and elaborated by
Itenberg later into the so-called $T$-curves context. Of course
the method has nothing formal: it is crystallography of the best
stock as we shall soon see.

As we failed along strategy (1) (cf. previous section), let us
look at (2) which invites to take a better look on  Itenberg-Viro
seminal paper (1996 \cite{Itenberg-Viro_1996-disproves-Ragsdale}).
There on Fig.\,2 (p.\,20) we find a remarkable picture reproduced
below as Fig.\,\ref{Itenberg:fig}. (Our sole change is to have
traced the curve with less thickness to see better what happens
than on the downloaded black-and-white pdf.) Note that the
underlying triangulation (by triangles all of area one-half of the
unit square) is symmetrical under the dihedral group (Vierergruppe
$D_4\approx \ZZ_2\times \ZZ_2$, cf. Fig.\,a). I did not noticed
this for a while when reproducing this figure! Then there is a
signs-distribution that gives the red curve via a bisection
procedure of all triangles. Interpret this merely as a
piecewise-linear random walk if you have zero-knowledge like the
writer. Tracing this curve is a miracle, very enjoying if one is
computer assisted. Representing the curve alone gives Fig.\,b (a
simple cut-and-past operation for the computer). As usual in
projective geometry the boundary of the rhombs must be identified
antipodically. Hence the 5 semi-ovals on the bottom left-side of
the rhombs are really just capping off (closing) the long
contorted oval occupying the oriental (Siberian) part of the
rhombs. We count on Fig.\,c precisely 29 outer ovals. Hence the
curve in question (which admits an algebraic realization by a deep
theorem of Viro-Itenberg) has the scheme depicted on Fig.\,d,
hence Gudkov symbol $(1, (1,2)2)29$. The total number of oval is
$r=29+6=35$, i.e. 2 unit less than $M=37$ (temperature of the
human body), so its an  $(M-2)$-curve. Hence there no obstruction
(via Klein's congruence) for the curve being dividing. Looking
optionally at the corresponding Hilbert tree (Fig.\,e) gives
quickly $\chi=29+1-3+2=29$ (variant look at the scheme Fig.\,d and
apply a Swiss cheese recipe \`a la Euler-Listing, etc.). Arnold's
congruence mod 4 ($\chi\equiv_4 k^2$ if type~I) is verified and so
there still no obstruction for the curve being dividing. In the
fact, the wonderful RKM-congruence $\chi\equiv k^2+4 \pmod 8$
(Rohlin-Kharlamov-Marin, see
(\ref{RKM-congruence-reformulated:thm})) implies  this scheme
being of type~I (and so is in particular any curve representing
it). Hence the Itenberg-Viro curve is of type~I, and it gives the
long searched counter-example to our Gabard-Thom pseudo-theorem.
Probably, there is a more elementary (organical) way to deduce the
dividing character of the curve by an avatar of Fiedler's
signs-law in the realm of the Viro method. (This goes back to the
early 1980's, and perhaps in the present $T$-curve context is due
to Itenberg. Parenti's thesis 199X \cite{Parenti_199X} looks also
involved in this topic.)

\begin{figure}[h]
\hskip-0.4cm\penalty0 \epsfig{figure=Itenberg.eps,width=132mm}
\vskip-5pt\penalty0
  \caption{\label{Itenberg:fig}%
The Itenberg-Viro patchwork construction of a (dividing)
$(M-2)$-curve of degree 10 killing Ragsdale's conjecture, as well
as the (dubious) Gabard-Thom conjecture ($\chi\le k^2$ for all
dividing curves)!} \vskip-5pt\penalty0
\end{figure}

So we are at this stage nearly in the paradise! To add some
suspense note that this curve being of type~I (even in the very
strong sense that its scheme is) a vague conjecture  of us
(founded on Ahlfors) posits that there should be a total pencil of
$(m-3)$-adjoint curves exhibiting the dividing character of the
curve (as in the Rohlin-Le~Touz\'e 1978--2013 theorem for sextics,
cf. Le~Touz\'e 2013
\cite{Fiedler-Le-Touzé_2013-Totally-real-pencils-Cubics}).
According to our (hypothetical) extension of the Le~Touz\'e
theorem, there should be a pencil of septics cutting only real
points on the Itenberg-Viro curve. In general, the degree of such
a pencil should be in accordance with Gabard's bound $r+p$ on the
degree of circle maps (compare
Sec.~\ref{census-and-extension-of-Rohlin-Le-Touzé:sec}, the
numerology therein and especially
Remark~\ref{M-2-curve-degree-like-Gabard:rem}). Pencil of septics
have 34 basepoints freely assignable (in general this is $M-3$),
whereas the Itenberg-Viro curve has 33 empty ovals. So where to
assign the remaining basepoint? This is again a bit puzzling.
However it is likely by the philosophy of dextrogyration that
Rohlin's complex orientation plays some r\^ole there. Rohlin's
formula applied to the Itenberg-Viro curve gives
$2(\pi-\eta)=r-k^2=35-25=10$, hence $\pi-\eta=5$, while
$\pi+\eta=7$ as is apparent from the Hilbert tree (Fig.\,e, where
one counts 2 additional pairs of length 2). Hence $\pi=6$ and
$\eta=1$. The presence of a negative pair is no surprise, as
Fiedler noticed this being the gap in our (erroneous) proof of the
Gabard-Thom estimate (Theorem~\ref{Thom-Ragsdale:thm}). By the
signs-law (\ref{Signs-law:lem}) it is clear that the (unique)
negative pair must have length one (primitive pair), otherwise all
primitive edges are positively charged but then there are 2
negative pairs (those of length 2). Since there is always a
bijective correspondence between primitive edges of Hilbert's tree
and (non-maximal) ovals by taking the bottom of the edge which is
always uniquely defined. So there is one negative oval (i.e. whose
edge) and thinking more with the signs-law it seems that the
negative pair is forced to be the trunk (i.e. the edge that
ramifies at depth 2, cf. Fig.\,e and lemma below). So this ovals
(negatively charged right above it) is perhaps the good candidate
of where to assign the remaining basepoint. All this deserves of
course to be better understood, and is digressed upon in the next
section (\ref{Galton-brett:sec}).

Let us verify the simple:

\begin{lemma}
The Rohlin tree of the Itenberg-Viro curve has a unique negative
charge that is forced to be on the trunk.
\end{lemma}

\begin{proof} As noted above Rohlin's formula implies that there
is a unique negative pair ($\eta=1$). Since in Rohlin's
arithmetics $+\times +=-$ (i.e. consanguinity is bad) the deepest
subtree of length 2 (``Y''-shaped modulo horizontal mirror) cannot
be positively charged else there would be 2 negative pairs of
length 2. So a primitive minus-charge  must be located on the
``Y''-subtree. If it is on one of the 2 branches (as opposed to
being on the trunk at depth $0-1$) then by uniqueness, the trunk
and the (other) branch are positively charged, so that it results
a minus-charge on their concatenation (of length 2), violating
$\eta=1$. So the negative charge must be located on the trunk as
asserted.
\end{proof}

\subsection{Where to assign basepoints to ensure total reality:
toward a Galton-Brett algorithm?} \label{Galton-brett:sec}

[24.03.13] The above fantastic example of Itenberg-Viro
(Fig.\,\ref{Itenberg:fig}) raises again the general problem of
deciding  where to assign basepoints as to ensure total reality of
an adjoint pencil. This means as usual a pencil of curves cutting
only real points on a dividing curve. By general topology (the
image of a connected set is connected) a dividing curve presents
no obstruction to the existence of such a pencil. More than that,
Ahlfors theorem says that there is no conformal obstruction to do
this (being after all just an extended Riemann mapping theorem for
bordered surfaces of higher topological structure that the disc).
If one has a total map to the line $\PP^1$ then we have a branched
cover taking boundary to boundary and interior to interior, and
the map restricted to the real locus is an (unbranched) covering.
Accordingly there is a phenomenon of dextrogyration, i.e. when the
image-point circulates once around the fundamental circle $\PP^1
(\RR)$ the counter-images (fibre) circulate along the complex
orientation of the abstract curve, i.e. as the boundary of the
half. (This follows of course from the holomorphic, hence sense
preserving, character of  Ahlfors maps.)

If the curve is plane ($C\subset \PP^2$) we would like a procedure
predicting where to assign basepoints. By the above dextrogyration
principle, there should be some relation with Rohlin's complex
orientations, which measure merely the (abstract) complex
orientations as compared with those of rings (annuli) bounding
(injective) pair of ovals in the plane $\RR P^2$. By the signs law
it suffices to determine Rohlin's signs for primitive pairs of
ovals.

A vague idea is as follows. Given a plane curve we have the
Hilbert tree and Rohlin's complex orientations decorate its edges
with signs (pluses, minuses). We can imagine the resulting Rohlin
tree as a ``Galton Brett'', i.e. Galton's table where
billiard-balls fall downwards along an inclined table interspersed
with a distribution of nails. Whenever meeting one of those nails
the ball is deflected left or right with probability one-half. For
the usual equilateral distribution we recover so the
Chinese-Pascal binomial distribution.

Our naive idea is to interpret the Rohlin tree as a Galton-Brett,
while putting balls at the top of  Hilbert's tree and looking
where they stabilize to an equilibrium. It is imagined that
negative pairs are inclined so that  balls fall gravitionally
along them. Consider the example of the deep nest. We know then
either from Rohlin's formula or via the dextrogyration argument
applied to the obvious pencil of lines through the deep nest that
all primitive pairs are negative. Here the Galton-Brett reduces to
a simple track (without branching) always negatively charged, and
the ball descends right up to its bottom. This is in agreement
with the fact that total reality of a lines-pencil is ensured when
assigning the basepoint in the deepest oval.

Similar considerations hold for the quadrifolium and its
satellites, i.e. curves of degree divisible by four and totally
really under a pencil of conics.

On the next example of Rohlin-Le~Touz\'e's sextics, e.g. that of
type $\frac{6}{1}2$, the Rohlin tree has 6 branches emanating from
an oval and Rohlin's formula $2(\pi-\eta)=r-k^2=9-9=0$ shows that
$\pi=\eta=3$ since we have a total $\pi+\eta=6$ of six pairs. So
our tree has 3 negative and 3 positive edges. Our metaphor of the
Galton-Brett already looks dubious on that example, since it would
prescribe imposing basepoints only on the $3$ ovals surmounted by
a negative charge (and of course the 2 outer ovals). So it remains
to understand if an improved Galton-Brett principle permits to
understand where  to assign basepoints in function of a knowledge
of Rohlin's complex orientations.

Maybe an improved rule is to let balls fall-down regardless of
signs along the Rohlin tree, and some few of them could stay in
levitation (unstable equilibrium being blocked by a needle) with
special signs-property, like being surmounted by a negative charge
while branching down below (so-called hyperbolic ovals). This
complicated condition comes to mind when looking at the deep nest
plus the above Itenberg-Viro curve, where the trunk (of the tree)
is negatively charged, which is the only reasonable
signs-distribution compatible with Rohlin's formula (cf. lemma
above). So on the tree of Fig.\,\ref{Itenberg:fig}e balls would
fall along the Hilbert tree and stabilize of course in each
``deep'' ovals (aka empty ovals), but some nontrivial equilibrium
arises at the vertex at depth 2 which branches further (alias
hyperbolic oval). Of course hyperbolicity alone is not enough as
shown by Rohlin-Le~Touz\'e's theorem. However hyperbolicity plus a
negative charge above it could give an equilibrium, i.e. a place
where to assign a basepoint. Though a bit complicated this looks
even reasonable from the viewpoint of Galton's Brett. Namely
hyperbolic ovals are those where there is an indetermination
(bifurcation) when falling down, while negativity of the edge
above is a sort of kinetic impulse giving the particle some
momentum forcing it to move against the bifurcation, whence an
``unstable'' equilibrium (crystallizing thereby in the formation
of a basepoint). Of course all this need to be further explored,
and to be related to more intrinsic properties of dividing plane
curves.

For the Viro-Itenberg curve this algorithm would assign basepoints
of the septics-pencil
on the 33 empty ovals (stable equilibrium) plus one unstable
equilibrium materialized by the unique hyperbolic oval at depth 1.
This would give the 34 basepoints required in a pencil of septics,
and total reality could follow (assuming that our Galton algorithm
is somehow compatible with dextrogyration or perhaps indexes
formulae \`a la Gauss-Kronecker-Poincar\'e-von Dyck).

Let us look at more examples. For octics we have 4 basic schemes
listed in Eq.~\ref{octics-five-examples-RKM:eq} which satisfy the
RKM-congruence (cf. also the Gudkov table in degree 8,
Fig.\,\ref{Degree8:fig}). Those were
\begin{equation*}
\frac{16}{1}3,\quad \frac{12}{1}7,\quad \frac{8}{1}11,\quad
\frac{4}{1}15,\quad 20.
\end{equation*}
The last of which is precluded as it violates either Petrovskii's
inequality (\ref{Petrovskii's-inequalities:thm}) or the Thom
estimate $\chi\le k^2$ which is valid when there is no nesting (or
even the more elementary Rohlin's formula). However for all other
schemes it may be reasonable to expect total reality for a pencil
of quintics which has 19 basepoints (recall that $B=M-3$ for the
number of basepoints in terms of Harnack's bound $M$) and all of
them are ascribed on the empty ovals (in accordance with our
Galton principle).

All those 4 RKM-schemes are just the top of the  iceberg depicted
on Fig.\,\ref{RKM-schemes-deg-8:fig}. On that tabulation we find
e.g. the scheme $\frac{3}{1}\frac{1}{1}14$. This has 18 empty
ovals, and we need a 19th basepoint. Our algorithm of the negative
hyperbolic oval fails to give it since looking at the Hilbert tree
of the scheme we see a unique hyperbolic oval, and this has no
edge above it! So our method fails and deserves to be further
improved. Less likely, our method could be right and then it could
preclude existence of those schemes in type~I. (Note that our
Galton method could have killed the scheme $20$, since we expect
19 basepoints but there are 20 stable equilibriums.) Of course all
this must be further explored. Summarizing, a fundamental question
seems to be:

\begin{ques}
Is there a general algorithm telling one where to ascribe
basepoints in terms of the combinatorics of Rohlin's tree encoding
his complex orientations? If so then we get a mechanical device
extending the Rohlin-Le~Touz\'e phenomenon of total reality for
$(M-2)$-sextics satisfying the RKM-congruence.
\end{ques}

Viceversa, suppose zero-knowledge on the complex orientations one
could argue that the principle of total reality via dextrogyration
is a good recipe to infer a knowledge of them (e.g. as it is
flagrant in the trivial deep nest case). Presently very little
seems to be known in general, and this is of course much
reminiscent of the lost proof of Rohlin's (last) theorem.

\subsection{When is Arnold's surface orientable and
Ragsdale via Bieberbach-Grunsky?}
\label{Ragsdale-via-Riemann-Bieberbach-Grunsky:sec}

[27.03.13] Here we propose a (naive) attack upon one half of the
Ragsdale conjecture for $M$-curves. This may be translated as the
condition $\chi\le k^2$ (cf.
Lemma~\ref{Thom-implies-one-half-of-Ragsdale:lem}). As shown by
the proof of the erroneous Theorem~\ref{Thom-Ragsdale:thm}, under
the additional assumption that the Arnold surface is orientable,
Thom applies and gives promptly the (upper) Ragsdale estimate
$\chi\le k^2$. (Note: The full Ragsdale amounts to the pinching
$-k^2\le \chi \le k^2$, equivalently $\vert \chi \vert \le k^2$.)

\begin{figure}[h]
\centering
\epsfig{figure=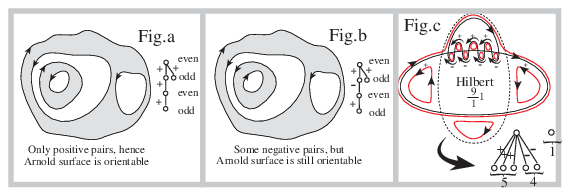,width=122mm} \vskip-5pt\penalty0
  \caption{\label{Ragsdale:fig}%
When is Arnold's surface orientable? It is, if only positive
pairs, but more generally when Rohlin's tree is positively charged
on odd edges.} \vskip-5pt\penalty0
\end{figure}

So the core of the question is to know when Arnold's surface is
orientable. We shall discuss this soon. The net impact could be as
follows:

\begin{conj}\label{Arnold-surface-M-curves-orient:conj}
The Arnold surface of a (plane) $M$-curve is always orientable. If
this is true then the upper-Ragsdale estimate $\chi\le k^2$
follows from Thom (cf. proof of (\ref{Thom-Ragsdale:thm})).
$\bigstar$ Alas, the sequel shows that this naive conjecture fails
already for Hilbert's $M$-sextic, cf. Fig.\,\ref{Ragsdale:fig}c.
\end{conj}

Recall that the {\it Arnold surface\/} of a dividing plane curve
of even degree $2k$, is simply Klein's half of the curve filled by
Ragsdale's membrane in $\RR P^2$ bounding the curve from inside.
It will be orientable iff the orientation coming from the
complexification and the real Ragsdale membrane match together in
some sense made precise below.

First it is plain that if Rohlin's tree is positively charged (on
all its primitive edges) then  Arnold's surface is orientable (cf.
Fig.\,\ref{Ragsdale:fig}a). This positive-charge assumption is
very stringent and implies actually much more, namely that the
Rohlin mass $\pi-\eta$ equals $n$ (cf.
Lemma~\ref{Rohlin-mass-of-a-positively-charged-tree:lem}). Via
Rohlin's formula rewritten as $\chi=k^2+2[(\pi-\eta)-n]$, this
implies in turn that $\chi=k^2$ exactly.

However to derive the (elusive) upper-Ragsdale-estimate $\chi\le
k^2$ from Thom, it suffices that the Arnold surface is orientable.
A small picture (Fig.\,\ref{Ragsdale:fig}b) convinces one that
this holds more generally whenever Rohlin's tree is positively
charged on odd edges. Here we always define the depth of an edge
in reference to that of its bottom vertex (the latter being
uniquely defined by---and defining uniquely---the given edge).

So we get the:

\begin{lemma}\label{Arnold-surface-orientable-iff-oddly-charged:lem}
The Arnold surface is orientable iff the Rohlin tree is positively
charged on odd edges (say then that it is oddly charged).
\end{lemma}

This condition is much weaker, but still implies Ragsdale via
Thom. So the upper Ragsdale $M$-conjecture (URMC) reduces to the:


\begin{conj}\label{oddly-charged-M-curves:conj}
Any $M$-curve is oddly charged.
\end{conj}

A guess could be to use our translation
(Theorem~\ref{total-reality-of-plane-M-curves:thm}) of the
Bieberbach-Grunsky theorem (truly due to Riemann 1857). The idea
is that for $M$-curves we have a fairly explicit way to construct
a total pencil via curves of degree $(m-2)$ assigned to pass
through any distribution of $g+1$ points (one on each oval) and
then by looking at the residual group of points, while assigning
them as basepoints. For more details cf. proof of
Theorem~\ref{total-reality-of-plane-M-curves:thm}, but we repeat
the general recipe of the construction of a total series on a
plane $M$-curve $C_m$ of degree $m$:

1.---Choose any distribution $D$ of $g+1$ points one on each oval.

2.---Let pass a curve $\Gamma_{m-2}=:\Gamma$ of degree $(m-2)$
through $D$.

3.---Consider $R$ the residual intersection $\Gamma\cap C$ less
the points of $D$.

4.---Assign $R$ as the basepoints to the system of curves of
degree $(m-2)$, and get (or choose) a pencil $\Pi$ putting the
initial group $D$ into motion.

5.---By continuity the pencil $\Pi$ is total since there is only
one point one each circuit hence no risk of collision. Total
reality follows.

The dream would be that this procedure is sufficiently explicit as
to control complex orientations, especially the issue that the
Rohlin tree is oddly charged. If this is possible we get a proof
of the upper-half $\chi\le k^2$ of Ragsdale's conjecture.

Concretely, once the distribution $D$ is fixed we are assured that
the curve $\Gamma$ will cut $C$ once more along each ovals (by the
closing lemma for algebraic circuits \ref{Closing-lemma:lem}). So
we have $2(g+1)$ real intersections in $\Gamma\cap C$, i.e.
$2(g+1)=2(\frac{(m-1)(m-2)}{2}+1)=(m-1)(m-2)+2=m^2-3m+4$. This is
less than the $m(m-2)$ expected.

\begin{lemma}\label{disproof-orientability-Arnold-M-curve:lem} Alas, the oddly-charged conjecture
(\ref{oddly-charged-M-curves:conj}) or equivalently the
orientability of the Arnold surface of an $M$-curve
(\ref{Arnold-surface-M-curves-orient:conj}) fails already in
degree 6, e.g. for Hilbert's $M$-sextic as shown on
Fig.\,\ref{Ragsdale:fig}c.
\end{lemma}

\begin{proof} Indeed reporting complex orientation via Fiedler's
transmission-law we get Fig.\,\ref{Ragsdale:fig}c. Here we report
first the orientation induced on the quadrifolium quartic $C_4$
from the dashed ellipse, and then smooth $C_4\cup E_2$ (where
$E_2$ is the thick-ellipse) along positive orientation and receive
so the complex orientations of Hilbert's sextic. We see that among
the 9 nested ovals, 5 are dominated by a positive pair, while 4
are by a negative pair. This is in accordance with Rohlin's
formula, $2(\pi-\eta)=r-k^2=11-9=2$, i.e. $\pi-\eta=1$ while
$\pi+\eta=9$. Hence $2\pi=10$, i.e. $\pi=5$ and $\eta=4$. However
this refutes our very naive conjecture
(\ref{oddly-charged-M-curves:conj}), which diagrammatically
amounts saying that the Rohlin tree is positively charged on edges
at odd depths. Of course the equivalent formulation in terms of
the orientability of the Arnold surface
(\ref{Arnold-surface-M-curves-orient:conj}) is killed in the same
stroke.
\end{proof}

So our naive strategy fails severely but of course Hilbert's
sextic has $\chi$ very negative ($\chi=2-9=-8$). Hence there is
perhaps a refined argument, that can establish Ragsdale. Alas it
seems that what we just did kill definitively an approach via Thom
which requires orientability of the Arnold surface. Of course one
could expect a tricky case distinction along the sign of $\chi$,
and a strengthened conjecture stating orientability of the Arnold
surface (of an $M$-curve)  provided $\chi>0$. Even this is  easily
disproved, e.g. by looking at the $M$-curve $C_{10}$ with $\chi=9$
of Fig.\,\ref{HilbGab1:fig}, and reporting the complex
orientations via Fiedler's law. This is a bit tedious but
straightforward and gives Fig.\,\ref{HilbGab1bis:fig}. We find
that the 11 edges at depth 3 splits into 7 positive pairs and 4
negative ones. Rohlin's formula can be verified via the signs-law.
However Rohlin's tree is not positively charged at the odd depth
3.

\begin{figure}[h]
\hskip-3.5cm\penalty0 \epsfig{figure=HilbGab1bis.eps,width=188mm}
\vskip-5pt\penalty0
  \caption{\label{HilbGab1bis:fig}%
An non-orientable Arnold surface in degree 10, with $\chi>0$}
\vskip-5pt\penalty0
\end{figure}

Of course one could still expect that Rohlin's tree is
oddly-charged when $\chi>k^2$, and this would suffices via Thom to
prove the upper Ragsdale conjecture, but we are obviously playing
a sterile arithmetical game, without much geometrical penetration.
Alternatively, if not via Thom one could hope to use directly
Rohlin's formula, but again some external information on complex
orientations must be gained via some deep geometric procedure. As
we said one can dream that a version of the Bieberbach-Grunsky
theorem do the job, but that deserves investigating with much more
care and patience than we are presently able to do. Good luck to
anybody who still feel optimistic. Of course it may also be
 that Ragsdale's upper bound is just false by a highbrow
variant of the Itenberg-Viro construction (though since
Itenberg-Viro 1996 \cite{Itenberg-Viro_1996-disproves-Ragsdale}
nobody apparently ever succeeded), and it is evident that we
(personally) lack experimental data to feel really secure in
claiming the Ragsdale bound. Hence we abort this problem for the
moment.

\subsection{A sporadic (?) obstruction via Thom
(Kronheimer-Mrowka 1994)}

\def\TKM{Thom-Kronheimer-Mrowka}

[22.03.13] It should be noted that the first examples of Thom's
conjecture as applied to Hilbert's 16th where we fill an unnested
curve by discs are not affected by Fiedler's correction. So in
particular the ``elementary'' degree 3 case of Thom due to
Kervaire-Milnor 1961 (yet relying massively upon Rohlin's early
work ca. 1951--52) really implies e.g. a purely topological proof
of Hilbert's intuition of nesting for $M$-sextics. This is
detailed below. More generally, Thom's conjecture forbids in all
(even) degrees $m=2k\ge 6$ the possibility of an unnested
$M$-curve (symbol $M$). This was first proved by Petrovskii 1938,
and can  also be deduced from Rohlin's formula
$0=2(\pi-\eta)=r-k^2$, since $M$ is strictly larger than $k^2$
when $k\ge 3$.

[05.03.13] Now just a little remark along the Thom conjecture
(=the Kron\-heimer-Mrowka theorem 1994
\cite{Kronheimer-Mrowka_1994}, abridged as ``Thom'' in the
sequel). If we look at the $(M-2)$-scheme $20$ of degree 8, and
fill one half by the canonical orientable membrane we get a
surface of genus $p=1$ whose homology class is $4H$ (where $H$ is
the natural generator of $H_2(\CC P^2,\ZZ)=\ZZ$, the so-called
hyperplane-section, here a line). By Thom the genus should be at
least as big as that of a (smooth) quartic, hence 3. So we get
the:

\begin{lemma}
The scheme $20$ is not realized algebraically by a curve of degree
$8$ (necessarily of type~I by the RKM-congruence
\ref{Kharlamov-Marin-cong:thm}).
\end{lemma}

If we take the scheme $\frac{4}{1}15$ (cf.
Fig.\,\ref{RKM-schemes-deg-8:fig}), then the genus of the filled
membrane will be $1+4=5$ and so Thom's principle is not violated.
No other scheme of that table are prohibited by Thom.

If we look at $m=6$, and the Gudkov table
Fig.\,\ref{Gudkov-Table3:fig}, especially the seminal intuition of
Hilbert ca. 1891--1900, that the unnested scheme $11$ does not
exist algebraically, then again we see that this may be inferred
from Thom's conjecture (meanwhile a theorem). Indeed making the
canonical filling of the half Riemann surface by the canonical
(Ragsdale) membrane (often denoted $\RR P^2_+$) we get a surface
of genus of 0 realizing the homology class $3H$, hence violating
Thom's conjecture.

{\it Insertion} [21.03.13].---As a matter of fact, this special
degree 3 case of Thom's conjecture was first established by
Kervaire-Milnor 1961 \cite[Cor.\,2,
p.\,1652]{Kervaire-Milnor_1961}, basing themselves much upon
Rohlin's work of 1951--52. Thom is alas not mentioned there
(KM61). So historiographically, it is worth emphasizing that the
Kervaire-Milnor paper afforded so-to-speak the first purely
topological proof of Hilbert's 1891 semi-intuition/semi-theorem
that a Harnack-maximal sextic cannot have all its 11 ovals
unnested! Prior to this we had only available: (1) the
algebro-geometric (stratificational) proofs of Hilbert (1891
unpublished) and Rohn (1911--13), plus technical refinements of
the same method by Gudkov ca. 1954 and (2) the proof of Petrovskii
1933/38 half analytical (Euler-Jacobi analytical interpolation)
and half topological (Morse theory). Nowadays we have of course
the proof via Rohlin's formula (1974 \cite{Rohlin_1974/75}), which
is more elementary and purely topological (or the related one via
Arnold's congruence mod 4).

Several ideas arise.

(1).---Of course this Thom argument is not the most elementary
prohibition of the scheme $10$ of degree 6, but maybe it is a good
way to prohibit the scheme $20$ in degree $8$ (at least I know no
other method for the moment).

{\it Update} [07.03.13]: the prohibition of this scheme $20_{8}$
follows more elementarily from Rohlin's formula
$2(\Pi^+-\Pi^-)=r-k^2$ (\ref{Rohlin-formula:thm}), since the
absence of nesting implies vanishing of the left-hand side, hence
$r=k^2$ has to be a square (even $16$ as $k=4$). Private anecdote,
I missed this consequence of Rohlin, and noticed it while
completing the Gudkov Table in degree 8 (cf.
Fig.\,~\ref{Degree8:fig}). [21.03.13] Another way to prohibit this
scheme $20$, I presume the first historically found, involves the
Petrovskii inequalities, cf.
(\ref{Petrovskii's-inequalities:thm}).

More generally what schemes can be prohibited by Thom, and did it
affords new obstructions (not known before Kronheimer-Mrowka)? The
first question is answered by Theorem~\ref{Thom-Ragsdale:thm}
below [alas erroneous!], while the second is perhaps answered via
Theorem~\ref{Alsatian-scheme-Thom-strong-Petrov-Arnold:thm}.

(2).---Looking at this canonical membrane filling  might be a good
device toward understanding the RKM-congruence. (But this is
merely a matter of reading once carefully the Kharlamov or Marin
arguments.) At least it is tempting to calculate the
self-intersection of this filled membrane with itself (or its
companion) to get some numerical relation. Doing so we probably
obtain nothing new but what exactly? (guess the Arnold
congruence).

(3).---Despite having corners this filled Riemann surface is
perhaps a good object to do conformal geometry with (compare
especially works by H.\,A. Schwarz ca. 1870 and his student Koebe
ca. 1906--07, also that of Hilb ca. 1907, NB: Hilb is not Hilbert
misprinted but a less well-known conformal geometer of that
period).

(4).---As this Thom argument prohibits the scheme $20$ [true
despite Fiedler's correction], which was an obstacle toward
assigning the 19 basepoint, try to pursue the game of
understanding the order 8 avatars of the Rohlin-Le~Touz\'e
theorem.

About (2), let $C_m$ be a dividing curve of even degree $m=2k$.
Denote by $F$ the closed surface $C_m^{+}\cup R$, where $R:=\RR
P^2_+$ is the canonical ``Ragsdale'' membrane (my own jargon but
historically justified I think after reading Viro's admiration for
Miss Ragsdale, as US-Studentin of Klein-Hilbert). On the one hand
$F$ is homologous to $kH$. To compute the self-intersection $F^2$
we use a vector field on the canonical membrane $R$, which is
either transverse or tangent along the boundary and with finitely
many zeros inside $R$. We have by Kronecker-Poincar\'e's index
formula $\sum indices= \chi$, where $\chi$ is the Euler
characteristic of the membrane $R$. Multiplying this vector field
by $i=\sqrt{-1}$ permits to push one replica of $F$ in general
position. Naively it seems to follow that $k^2=F^2=\chi$, but one
needs to count better indices....

Another remark is to write down the Thom conjecture inequality for
the filled surface, and this gives the following (which is
certainly not new [$\bigstar$ but alas false!!], cf. maybe
Degtyarev-Kharlamov 2000 \cite{Degtyarev-Kharlamov_2000}, but
after a rapid check it does not seem to be explicitly stated
there):

\begin{theorem}\label{Thom-Ragsdale:thm}---{\rm ERRONEOUS AT LEAST IN THE
WEAK SENSE THAT THERE IS A BASIC BUG, YET NO CONCRETE EXAMPLE
KNOWN TO ME\footnote{[27.03.13] Meanwhile the simplest
counterexample, I was able to find is the Itenberg-Viro curve
constructed on Fig.\,\ref{Itenberg:fig}.}} {\rm (Thom applied to
Klein-Hilbert-Ragsdale)}.---Let $C_m$ be a dividing plane curve of
degree $m=2k$. Then $\chi \le k^2$, where $\chi$ denotes (as
usual) the Euler characteristic of the Ragsdale membrane.
\end{theorem}

{\it Insertion.}[14.03.13].---For each $k$, it is a simple matter
to convince that the estimate is sharp, compare
Figs.\,\ref{CCCRoses:fig} and \ref{CCCRoses2:fig}.

{\it Inserted (optional reading).} [16.03.13].---At first sight
this result looks so limpid [$\bigstar$ outdated now!] that one
may wonder if it extends to higher dimensions, e.g. to algebraic
surfaces in $\PP^3$. First one requires an extension of the Thom
conjecture for surfaces in $\CC P^3$. This is perhaps quite
straightforward, by replacing the genus by the Euler
characteristic and arguing experimentally that surgeries (aka
spherical modifications increase the topological complexity, yet
without changing the homology class). However the second step of
the proof fails blatantly as we lack a natural extension of the
concept of dividing curves to surfaces though several peoples
(notably Viro) proposed extensions requiring e.g. that the
homology class of the real part mod 2 vanishes in the
complexification. Alas, the real locus of an algebraic surface,
having (real) codimension 2 in its complexification,  never
divides. Even if it would the Ragsdale membrane has only (real)
dimension 3, hence not ideally suited to cap off the 4D-half (if
it existed). This could be remedied by looking at surfaces in
$\PP^4$ instead, yet we still lack a way to split the
complexification by the real locus. Of course all this failure is
somewhat akin to the lack of a good extension of Rohlin's formula
to surfaces as deplored upon, e.g. in Degtyarev-Kharlamov 2000
\cite{Degtyarev-Kharlamov_2000}. Basically, the ideas behind
Rohlin and the Thom estimate above are very similar, namely to
fill the ``imaginary'' half by a ``real'' membrane coming from the
real locus (either the Ragsdale membrane or bounding discs for
ovals).

\smallskip

\begin{proof} (of (\ref{Thom-Ragsdale:thm})).---We paste to the (bordered) half $C_m^+$
of the dividing curve $C_m$ the Ragsdale membrane $R$ which is the
orientable surface bounding $C_m(\RR)$. The resulting closed
surface $F$ is orientable [HERE IS THE MISTAKE (22.03.13)!!!] and
realizes the homology class $kH$ (of halved degree) in the group
$H_2(\CC P^2, \ZZ)\approx \ZZ$. It is plain\footnote{[21.03.13]
Find accurate references, by Thom, Cerf, Hirzebruch, Milnor, Wall,
etc. I confess that I lack a precise reference.} that we can
smooth the ``corners'' arising along the ``cut-and-paste-locus''
to get a nearby smooth surface still denoted $F$. By Thom's
conjecture (=meanwile the Kronheimer-Mrowka theorem 1994
\cite{Kronheimer-Mrowka_1994}) we infer that the genus of $F$, say
$f:=g(F)$, is at least as big as that of a smooth curve of the
same degree, i.e.,
$$
f\ge g(k)=\textstyle\frac{(k-1)(k-2)}{2}.
$$
On the other hand we have by additivity of the characteristic
$$
\chi(F)=\chi (C_m^+)+\chi (R).
$$
For the same reason $2\chi(C_m^+)=\chi(C_m)=2-2g(m)$, and so
\begin{align*}
\chi&:=\chi(R)=\chi(F)-\chi (C_m^+)=(2-2f)-1+g(m)=1-2f+g(m) \cr
&\le 1-(k-1)(k-2)+\textstyle\frac{(2k-1)(2k-2)}{2}\cr
&=1-(k-1)(k-2)+(2k-1)(k-1) =1+(k-1)(k+1)=k^2.
\end{align*}
\end{proof}

As already discussed, this has some interesting applications, e.g.
to the prohibition of Hilbert's (unnested) scheme $11$ of degree
6, and to the scheme $20$ in degree 8. (However all this can be
more elementarily  deduced from Rohlin's formula.) [21.03.13] Yet
compare
Theorem~\ref{Alsatian-scheme-Thom-strong-Petrov-Arnold:thm} below
for an example showing that Thom's estimate is sometimes stronger
than the conjunction of several powerful prohibitions of the
Russian school (Petrovskii 1938, Gudkov 1969, Arnold 1971, Rohlin
1972/74). Also pleasant is the direct link of this estimate with
those conjectured decades prior to Thom by Virginia Ragsdale in
1906 (cf. Sec.\,\ref{Ragsdale-conj:sec}). As I was informed by Th.
Fiedler, it seems that it is Mikhalkin who first investigated
systematically the repercussion of Thom-Kronheimer-Mrowka upon
Hilbert's 16th.

[07.03.13] When we look back at Gudkov's Table
(Fig.\,\ref{Gudkov-Table3:fig}) we see that we get a nearly
complete system of prohibition by using total reality and the
Rohlin maximality conjecture (RMC), while combining it with the
Thom obstruction. Remember that RMC ought to be a reliable
principle whenever total reality is exhibited in some concrete
fashion as in the Rohlin-Le~Touz\'e theorem. Hence what misses is
a prohibition of the scheme $\frac{10}{1}$ in degree 6. One may
thus wonder if there is an avatar of Thom's conjecture for
non-orientable surfaces in $\CC P^2$, able to prohibit the sextic
scheme $\frac{10}{1}$ (of Rohn). Cavalier, one could  put forward
something like the:

\begin{conj}
Every prohibition of Hilbert's 16th problem, is either
interpretable via total reality and the allied Rohlin maximality
principle to the effect that a scheme flashed by a total pencil is
maximal, or is a consequence of Thom's conjecture plus an avatar
thereof including non-orientable membranes.
\end{conj}

Of course for the sextic scheme $\frac{10}{1}$, the idea would be
to fill by the non-orientable membrane (residual to the Ragsdale
membrane).
Further our conjecture is certainly much premature unless it takes
into account advanced B\'ezout-style prohibitions \`a la
Fiedler-Viro (cf. Theorem~\ref{Viro-Fiedler-prohibition:thm}), and
Petrovskii-Arnold style prohibitions (cf.
Theorems~\ref{Petrovskii's-inequalities:thm} and
\ref{Strong-Petrovskii-Arnold-ineq:thm}).

\subsection{Ragsdale's conjecture (Ragsdale 1906, Petrovskii 1938,
Viro 1979/80, Itenberg 1993, Thom 19XX-Kronheimer-Mrowka 1995, and
still open, Fiedler)}\label{Ragsdale-conj:sec}

[18.03.13] As I was made (personally) aware by Fiedler (cf. his 9
March 2013 letter reproduced in Sec.\,\ref{e-mail-Viro:sec}) there
ought to be some connection between Thom's and Ragsdale's
conjecture, which is still open for $M$-curves, despite the
disproofs due to Viro 1979 (=Viro 1980/80
\cite{Viro_1980-degree-7-8-and-Ragsdale}) and Itenberg 1993
\cite{Itenberg_1993-ctrex-a-Ragsdale} (cf. also Itenberg-Viro 1996
\cite{Itenberg-Viro_1996-disproves-Ragsdale}). I can also remember
an oral discussion with Thomas Fiedler (Geneva ca. 2011), where
Thomas alluded to all the effort he invested on the Ragsdale
problem (for $M$-curves). At that time (and arguably still today),
I could not appreciate the full swing of this investment.

This section makes no pretence of any breakthrough in the field.
It is rather  a humble  attempt to get familiarized with the
topic.
%
Despite our incompetence, let us make some remarks. From our
viewpoint of total reality much allied to Ahlfors theorem, which
quite paradoxically seems  more familiar to complex/conformal
geometers than purely real ones (having in mind the antagonism
between Riemann-Schottky-Klein-Bieberbach-Teichm\"uller-Ahlfors
versus
Harnack-Hilbert-Ragsdale-Rohn-Petrovskii-Gudkov-Arnold-Rohlin,
etc.) we could expect a connection of Ragsdale's conjecture to our
paradigm of total reality, e.g. via
Theorem~\ref{total-reality-of-plane-M-curves:thm} as a first basic
step. This vague suggestion  should probably not be taken too
seriously. Another vague idea is to wonder if there is some
connection of Ragsdale with the contraction conjectures of
Itenberg-Viro, or perhaps our version thereof called CCC, cf.
(\ref{CCC:conj}).

After those abrupt remarks, let us be more pedestrian. First what
is Ragsdale's conjecture at all about? What is known on it and
what is not? Especially does it connect to 4D-dimensional topology
as the whole Hilbert problem was realized to be since Arnold's
breakthrough 1971 \cite{Arnold_1971/72} and the deeper
investigations of Rohlin (e.g., the validation of Gudkov's
hypothesis $\chi\equiv_8 k^2$). In particular how does Ragsdale
connect with Thom's conjecture which is basically a problem of
embedded differential topology  of smooth surfaces in the complex
projective plane $\CC P^2$ (arguably the  4-manifold simplest to
visualize as the configuration space of all unordered pairs
grooving on the 2-sphere). As we shall see, the link Thom-Ragsdale
is very clear-cut, at least for one half of the Ragsdale
conjecture (cf.
Lemma~\ref{Thom-implies-one-half-of-Ragsdale:lem}). [$\bigstar$
Okay, but alas this is based on our erroneous estimate $\chi\le
k^2$!]

First, Virginia Ragsdale, coming from the U.\,S. was a student of
both Klein and Hilbert in G\"ottingen ca. 1906. Building upon a
careful inspection of the features of Harnack's and Hilbert's
constructions (of small vibratory perturbations), she posited a
conjecture on the numbers $p,n$ of even resp. odd ovals of real
plane algebraic curves. Although the chance of deriving any
transcendental truth from such a specific mode of generation looks
a priori very meagre, the conjecture in question turned out to be
extremely robust requiring at least ca. 7 decades up to being
disproved [Viro, Itenberg]. Yet some respectable vestiges remains
open, and deserves further efforts.

Via some naive acquaintance with Gudkov's Table in degree 6
(Fig.\,\ref{Gudkov-Table3:fig}) (the little that we personally
have at disposal) and the allied geometry of pyramids, it seems
 that  the Harnack and Hilbert constructions explore
only the superficies of the pyramid, while the profound part of
the puzzle is cracked in Gudkov's revolution (ca. 1969--72)
constructing the scheme $\frac{5}{1}5$
(which is so-to-speak the Pharaoh chamber). The Ragsdale
conjecture (in modernized shape) may be stated as the estimate
$\vert \chi \vert \le k^2$, hence it is perhaps not too surprising
that the particular methods of Harnack and Hilbert lead to sharp
estimates at least for $M$-curves. What appears historically first
is likely to be the most superficial objects, hence extremalizing
the functional $\vert \chi\vert$ which roughly measures the level
of superficiality in the pyramid.

We  state now precisely  Ragsdale's original statement (compare
Ragsdale 1906 \cite{Ragsdale_1906}).

\def\Petrov{\textstyle\frac{3}{2}k(k-1)}

\begin{conj}\label{Ragsdale-conj:conj} {\rm (Ragsdale 1906,
disproved for the number $n$ of odd ovals by Viro 1979, and in
general by Itenberg 1993)}---For any curve of even degree $m=2k$,
we have
$$
p\le \textstyle\frac{3}{2}k(k-1)+1, \qquad n\le
\textstyle\frac{3}{2}k(k-1)
$$
\end{conj}

As explained in Itenberg-Viro 1996
\cite{Itenberg-Viro_1996-disproves-Ragsdale} (especially p.\,24)
this was
refuted by Viro in 1979 for the number $n$ of odd ovals, and by
Itenberg in general. Moreover it is explained (in \loccit) that
Petrovskii made similar conjectures, being apparently unaware of
Ragsdale's paper. In particular the so-called {\it Petrovskii
inequality} (cf. (\ref{Petrovskii's-inequalities:thm})) is
considered there as having been formerly conjectured by Ragsdale
as a weak form of her conjecture. Finally it is remarked that
Petrovskii himself (1938) formulated a version of Ragsdale's
conjecture (\ref{Ragsdale-conj:conj}), yet more cautious by one
unit than Ragsdale's on the number $n$, so that both bounds are
identic equal to $\textstyle\frac{3}{2}k(k-1)+1$.

Despite the disproof (by Viro-Itenberg) of both the Ragsdale and
the weaker Petrovskii conjectures, the interesting quick is that
the case of $M$-curves is still open (at least in the weaker
formulation of Petrovskii). Precisely

\begin{conj} {\rm (Ragsdale's conjecture on $M$-curves 1906,
still open)}---For any $M$-curve of even degree $m=2k$,  the Euler
characteristic $\chi=p-n$ of the Ragsdale membrane is bounded by
the square of the semi-degree $k$, i.e.
$$
\vert \chi \vert \le k^2.
$$
\end{conj}

We borrowed this from Itenberg-Viro 1996
\cite[p.\,24]{Itenberg-Viro_1996-disproves-Ragsdale} (cf. also
Kharlamov-Viro (undated)
\cite[p.\,15]{Kharlamov-Viro_XXXX-UNDATED}). It should be remarked
that one half of this conjecture (namely the estimate $\chi\le
k^2$) follows directly (cf. Theorem~\ref{Thom-Ragsdale:thm}) from
Thom's conjecture proved by Kronheimer-Mrowka in 1994
\cite{Kronheimer-Mrowka_1994}. Curiously, this is not pointed out
in the Itenberg-Viro 1996 article (presumably due to backlog
reasons). Actually as noted in Theorem~\ref{Thom-Ragsdale:thm},
the estimate $\chi\le k^2$ holds more generally for dividing
curves.

$\bigstar\bigstar$ {\it Insertion.} [23.03.13] This historical
puzzle is now completely fixed by Fiedler's correction of my
mistake of overlooking that the Arnold surface is not necessarily
orientable.

Some few words are required to understand why the above conjecture
$\vert \chi \vert \le k^2$ is termed Ragsdale's conjecture. (It
seems to me that Itenberg-Viro 1996
\cite[p.\,24]{Itenberg-Viro_1996-disproves-Ragsdale} contains a
serious misprint at this place, specifically on p.\,24 in the
statement of the Ragsdale conjecture on $M$-curves the equivalent
conditions $p\ge \frac{(k-1)(k-2)}{2}$ and
$n\ge\frac{(k-1)(k-2)}{2}$ looks to me erroneous; and the same
misprint appears in Kharlamov-Viro
\cite[p.\,15]{Kharlamov-Viro_XXXX-UNDATED}) Let us clarify this as
follows:

\begin{lemma}\label{Thom-implies-one-half-of-Ragsdale:lem}
For $M$-curves of degree $2k$, the condition $\vert \chi\vert\le
k^2$ is equivalent to Petrovskii's cautious version of the
Ragsdale conjecture, i.e.
$$
p\le \textstyle\frac{3}{2}k(k-1)+1, \qquad n\le
\textstyle\frac{3}{2}k(k-1)+1.
$$
More precisely the upper estimate on $\chi$ (i.e. $\chi\le k^2$)
is equivalent to the bound on $p$, while the lower estimate
$-k^2\le \chi$ is equivalent to the bound on $n$. Further the
first upper bound $\chi\le k^2$ follows from Thom's bound (cf.
{\rm Theorem~\ref{Thom-Ragsdale:thm}} valid more generally for any
dividing curve), while the other is perhaps still open, though one
could dream reducing it to Thom too, after taking maybe an
orientable cover (but looks dubious), or maybe by reducing it via
differential geometry to Gauss-Bonnet and Wirtinger (as discussed
below).
\end{lemma}

\begin{proof} Start from the condition $-k^2\le\chi\le k^2$. By
definition $\chi=p-n$, and $M=r=p+n$. As usual Harnack's bound  is
$M=g+1=\frac{(2k-1)(2k-2)}{2}+1=(2k-1)(k-1)+1=2k^2-3k+2$. Adding
$\chi=p-n\le k^2$ to $p+n=M=2k^2-3k+2$ gives
$$
2p\le k^2+2k^2-3k+2=3k^2-3k+2=3k(k-1)+2,
$$
whence Ragsdale's bound on $p$. On the other hand, the lower
estimate on $\chi$, i.e. $-k^2\le \chi=p-n$,  rewritten as $k^2\ge
n-p$, gives when added to $p+n=M$ the 2nd Ragsdale estimate on
$n$.
\end{proof}

Again the text of Itenberg-Viro 1996
\cite[p.\,24]{Itenberg-Viro_1996-disproves-Ragsdale} which reads
as follows, seems not perfectly up-to-date [SORRY MY MISTAKE!] in
view of Kronheimer-Mrowka's validation of Thom's conjecture:

\begin{quota}[Itenberg-Viro 1996]{\rm ``Which of Ragsdale's questions are still
open now? The inequalities\footnote{That is the cautious
Petrovskii version of Ragsdale's estimates.}
$$
p\le \textstyle\frac{3}{2}k(k-1)+1, \qquad n\le
\textstyle\frac{3}{2}k(k-1)+1.
$$
have been neither proved\footnote{It seems to me that the
estimates on $p$ follows from Thom's conjecture, as explained in
Lemma~\ref{Thom-implies-one-half-of-Ragsdale:lem}} nor disproved
for $M$-curves.''}
\end{quota}

As shown by Lemma~\ref{Thom-implies-one-half-of-Ragsdale:lem}
(implicit in Itenberg-Viro's article modulo the misprint), Thom's
conjecture implies one half of Ragsdale conjecture [ALAS NOT
TRUE], namely the ``positive'' half concerning $p$ where it
matches exactly with Petrovskii's subsequent rediscovery of the
conjecture. However the second half looks much out of reach, as it
requires the estimate $-k^2\le\chi=p-n$ which seems to take care
of the non-orientable (``anti-Ragsdale'') membrane (not ideally
suited to Thom).

First without any idea the lower bound on $\chi=\chi(B^+)$ the
characteristic of the (orientable) Ragsdale membrane $B^+$, i.e.
$-k^2\le \chi$, can be transmuted using $B^+\cup B^-=\RR P^2$ into
$k^2\ge\chi(B^-)-1$.

So one seems forced to study this non-orientable membrane $B^-$.
One idea to explore is to arrange an orientable membrane via the
usual trick of the double orienting cover (essentially due to
Gauss, M\"obius, Klein, Teichm\"uller 1939, etc.) but which have
now to be implemented in some embedded fashion. This looks dubious
as we do not know what to do along the boundary of $B^-$. At least
the surface we get (granting that there is some natural way to
construct an oriented Verdoppelung=double) would have tripodal
singularities along the boundary. This is common in soap film
experiment, yet a priori outside the
tolerance
permitted in Thom's conjecture.

As a completely different strategy there could be a result (dual
to Thom's) stating that for smooth surfaces the genus cannot be
too big when attention is confined to smooth surfaces arising by
rounding
corners of the half of a dividing curve capped off by the Ragsdale
membrane. Of course in general the genus of a smooth surface of
prescribed (homological) degree can be made as large as we please
(just attach small handles), yet perhaps the surfaces that arise
by the ``Ragsdale filling procedure'' of Klein's orthosymmetric
half (a method truly inaugurated by Arnold 1971, and Rohlin 1974,
etc.) are of a special type  subsumed to an upper-bound upon the
genus in terms of the degree.

Recall the Wirtinger's inequalities stating that complex
projective varieties (or even K\"ahler manifolds) minimize the
volume among differential-geometric submanifolds in a given
homology class. Perhaps this combined with Gauss-Bonnet can supply
the required upper-estimate upon the genus dual to Thom's
estimate, hence validating the remaining half of Ragsdale's
conjecture (as modified by Petrovskii 1938).

[19.03.13] In fact I do not know if the lower estimate of the
pinching $-k^2\le \chi\le k^2$ holds true more generally for
dividing curves as do the upper estimate $\chi \le k^2$ by virtue
of  Thom's (genus) bound (Theorem~\ref{Thom-Ragsdale:thm}) [FALSE,
cf. Itenberg-Viro's curve on Fig.\,\ref{Itenberg:fig}]. It is
worth first noting that the Petrovskii jargon (1938
\cite{Petrowsky_1938}) of $p$ and $n$ as positive and negative
ovals is quite good as they contribute positively resp. negatively
to $\chi$ of the Ragsdale membrane.
[Of course the sign of $\chi$ is a matter of convention that
varied through the ages, but at least now we seem to all agree
about its sign ($\chi(pt)=+1$).]
Then it is also useful to keep in mind the geography of the
generic Gudkov pyramid (cf. Fig.\,\ref{Pyramidragsdale:fig})
though this is a coarse simplification of the real one which is a
multidimensional (non-planar) object as soon as $m\ge 8$. Now the
point is that as shown by
Lemma~\ref{Thom-implies-one-half-of-Ragsdale:lem} the Thom
estimate $\chi\le k^2$ and its dual $-k^2\le\chi$ formally implies
for $M$-curves the Ragsdale-Petrovskii estimates $p\le P$, and
$n\le P$ respectively, where we set $P:=\frac{3}{2}k(k-1)+1$ (for
Petrovskii's bound). Diagrammatically, this amounts saying that
{\it Ragsdale's zone\/} ($p,n\le P$) arises from Thom's vertical
strip $\vert\chi\vert\le k^2$ by reflecting vertical rays at angle
of 60 degrees (cf. Fig.\,\ref{Pyramidragsdale:fig}). We get so a
pentagonal diamond (the Ragsdale diamond) that was supposed to
contain all algebraic schemes by virtue of the Ragsdale-Petrovskii
conjecture, which alas turned out wrong by Itenberg-Viro
[Fig.\,\ref{Itenberg:fig}]. However on the top face of the diamond
($M$-curves) the conjecture is still robust. Further, Thom's strip
$-k^2\le \chi \le k^2$ is perhaps a container for all dividing
curves [FALSE, cf. again the Itenberg-Viro
Fig.\,\ref{Itenberg:fig}]. This holds true on the right positive
side (still Theorem~\ref{Thom-Ragsdale:thm}) [FALSE!] but the
lower estimate is presently more dubious. A look on Gudkov's
table(=Fig.\,\ref{Gudkov-Table3:fig}) shows that $-k^2\le\chi$
holds true in degree $6$ (actually for all curves regardless of
being dividing), yet this low-degree case is probably atypical.

\begin{figure}[h]
\centering
\epsfig{figure=Pyramidragsdale.eps,width=102mm}
\vskip-5pt\penalty0
  \caption{\label{Pyramidragsdale:fig}%
Generic Gudkov table in degree $2k$ showing the geography of the
Ragsdale conjecture as a diamond $p,n\le P:=\frac{3}{2}k(k-1)+1$
expanding Thom's house $-k^2\le \chi \le k^2$} \vskip-5pt\penalty0
\end{figure}

By analogy with Thom's conjecture, let us put forward the:

\begin{conj}\label{anti-Thom-conj:conj} {\rm (Anti-Thom conjecture, due to Gabard hence
probably not serious at all)} Any plane dividing curve of even
degree $m=2k$ respects a Thom-style lower
 bound
$$
-k^2\le \chi.
$$
\end{conj}

If true this would imply the remaining half of the Ragsdale
conjecture for $M$-curves (via
Lemma~\ref{Thom-implies-one-half-of-Ragsdale:lem}). It would be
interesting to know if the Viro-Itenberg method (viz.
counter-examples) already disproves this naive conjecture. [QUITE
PROBABLE!] If not, one would like to imagine a proof along the
above sketched line (Gauss-Bonnet-Wirtinger), or via an oriented
double cover of the non-orientable Ragsdale membrane $B^-$.
Another third strategy would be to use an eversion (cf.
Sec.~\ref{Eversion:sec}) to reduce to the case of Thom, yet this
looks hazardous as it requires a large deformation (for which very
little is known apart vague speculation of us that
differential-geometric flows could do such jobs).

At this stage it is wise to contemplate the higher Gudkov's
pyramids in degree 8 and 10 (cf. resp. Figs.\,\ref{Degree8:fig}
and \ref{Degree10:fig}). In degree 8 (Fig.\,\ref{Degree8:fig}) we
see that the anti-Thom line $-k^2\le \chi$ is well adjusted to the
blue-rhombs materializing Arnold's congruence mod 4 for dividing
curves. Hence there is little chance to corrupt our anti-Thom
conjecture (\ref{anti-Thom-conj:conj}). In contrast in the degree
10 table (Fig.\,\ref{Degree10:fig}) there is a myriad of 7 schemes
adventuring outside the anti-Thom line. Those are given by the
symbols $\frac{34}{1}2$, $\frac{33}{1}1$, $\frac{32}{1}$ and
$\frac{31}{1}3$, $\frac{30}{1}2$, $\frac{29}{1}1$, $\frac{28}{1}$.
If any one of those schemes admits a type~I(=orthosymmetric)
realization our conjecture is faulty. This problem can either be
approached by Harnack or Hilbert's method of vibrations or by the
Viro-Itenberg patchworking. Note that it is unlikely that those
schemes (in type~I) are prohibited by Rohlin's formula. [21.03.13]
However the first 3 listed (with $\chi=-31$) are prohibited by
Petrovskii's inequality (\ref{Petrovskii's-inequalities:thm}).

For the 4 remaining schemes (with $\chi=-27$) one can expect to do
naive Hilbert constructions like on Fig.\,13 of Gabard 2000
\cite{Gabard_2000}, and look what happens. That requires some
concentration and is differed to latter.

[21.03.13] It may be noted that the strong Petrovskii inequality
(\ref{Strong-Petrovskii-Arnold-ineq:thm}) $n-p^- \le
\frac{3}{2}k(k-1)=30$ specialized to the range of our diagram
(Fig.\,\ref{Degree10:fig}) involving  simple symbols of the form
$\frac{x}{1}y$ where $p^-=1$ (one hyperbolic oval, i.e. which
ramifies in Hilbert's tree) implies that $n\le 31$ and so Ragsdale
line is corroborated. This does not kill any of the 4 schemes in
candidature above. Note further that it must be a general issue
that the strong Petrovskii-Arnold estimate gives Ragsdale in the
``planar'' range of the pyramid involving symbols $\frac{x}{1}y$.
On the right-side of Fig.\,\ref{Degree10:fig}, the (other dual)
strong Petrovskii-Arnold estimate $p-n^-\le
\frac{3}{2}k(k-1)=30+1$ (with now $n^-=0$) also kills all schemes
lying ``above'' Ragsdale's line. This implies a severe crumbling
in the corners of the pyramid (cf. Fig.\,\ref{Degree10:fig}) on
both the right and left side of it, yet note that a priori the
scheme $\frac{31}{1}3$ (if it exists) seems to imply some
mysterious asymmetry in the architecture. Further this scheme,
being of type~I, could be an interesting place to look for a
counterexample to Rohlin's maximality conjecture.

\subsection{Thom versus Rohlin}

[07.03.13] As noticed above the Thom obstruction
(\ref{Thom-Ragsdale:thm}) is at least for degree 6 (and to some
extend in degree 8) subsumed to Rohlin's formula
(\ref{Rohlin-formula:thm}). The latter also prohibits the scheme
$\frac{10}{1}$ (cf. Fig.\,\ref{Gudkov-Table3:fig}) of degree 6
without that we have to worry about a dubious non-orientable
extension of Rohlin's formula. One can wonder if in general the
information derived from Thom (\ref{Thom-Ragsdale:thm}) is always
subsumed to Rohlin's formula.

Let us notice the following:

\begin{lemma}
The scheme $20$ cannot be realized in degree $8$.
\end{lemma}

\begin{proof}
Let $C_8$ be a (hypothetical)  octic of type $20$. By the
RKM-congruence (\ref{Kharlamov-Marin-cong:thm}) or better its
reformulation as (\ref{RKM-congruence-reformulated:thm}) the curve
has to be of type~I, but then its existence is
ruled out by Rohlin's formula (\ref{Rohlin-formula:thm}).

[21.03.13] Another proof of the lemma (historically sharper) is to
appeal to Petrovskii's inequality (1933/38), cf.
(\ref{Petrovskii's-inequalities:thm}).
\end{proof}

It is worth then comparing the Gudkov table in degree 8
(Fig.\,\ref{Degree8:fig}), which shows  that several obstructions
are not readily interpreted via total reality and Rohlin's allied
principle of maximality (look especially at the upper-right corner
of that figure).

Now we turn to the question of deciding if Thom is subsumed to
Rohlin (at least in the realm of Hilbert's 16th problem). Glancing
at  Fig.\,\ref{Degree8:fig} the answers seems to be yes for degree
8 (at least for schemes of the form $\frac{x}{1}y$ as those
represented on that figure). If true in general this should follow
from a simple combinatorial argument.

If we look at the degree $m=10$ table, we find the following
structure (Fig.\,\ref{Degree10:fig}). This picture is built as
usual. First one compute Harnack's bound
$M=g+1=\frac{(m-1)(m-2)}{2}+1=\frac{9\cdot 8}{2}+1=36+1=37$. So
one extends the previous pyramid in degree 8, up to that level
$37$. Then there is the sawtooth broken line \`a la Gudkov. Its
upper undulations have to be adjusted at the Gudkov-Rohlin
congruence $\chi\equiv k^2 \pmod 8$, here $k^2=25$. Remind that on
such pictures Euler-Ragsdale's $\chi$ may (always) be interpreted
as the abscissa (``$x$-axis'', i.e. horizontal axis). So we have
the $M$-schemes lying at the top of the sawtooth broken line,
while in their depressions (``creux'') we have the RKM-schemes
with $\chi\equiv k^2+4 \pmod 8$, that are forced being of type~I.
All this can be extended into the lattice of blue-rhombs where
type~I curves are forced to live (Arnold's congruence mod 4).
[{\it Warning.\/}---On our Fig.\,\ref{Degree10:fig} this lattice
of blue rhombs is correct on the upper half, but need to be
adjusted (mentally) on the lower part, where we just copied the
pyramid in degree 8. However we thought it would be more
instructive to see this lower object intact as to appreciate
better the growing mode of pyramids.]

By Thom (\ref{Thom-Ragsdale:thm}) we  have $\chi\le k^2=25$ for
curves  of type~I. It seems (at first sight) that several schemes
permissible for Rohlin are prohibited by Thom, yet the ultimate
answer will be nearly the opposite one.

\begin{figure}[h]
\hskip-3.2cm\penalty0 \epsfig{figure=Degree10.eps,width=172mm}
\vskip-5pt\penalty0
  \caption{\label{Degree10:fig}%
Partial Gudkov's table in degree 10: trying to enhance that the
Thom obstruction is not subsumed to Rohlin's formula in general}
\vskip-5pt\penalty0
\end{figure}

First the $M$-scheme $\frac{2}{1}34$ is prohibited by Thom
[DUBIOUS INFERENCE OF THOM!] though it is not by the Gudkov-Rohlin
congruence. (Incidentally this scheme is prohibited by Petrovskii
(\ref{Petrovskii's-inequalities:thm}).) Is this scheme prohibited
by Rohlin's formula $2(\Pi^+-\Pi^-)=r-k^2$? Here we have two
nested pairs hence $\Pi=\Pi^+ +\Pi^-=2$, and $r-k^2=37-25=12$. So
$\Pi^+-\Pi^-\le \Pi=2$ and Rohlin's formula cannot be fulfilled.
So Thom brings nothing new.

Then there is the scheme $\frac{1}{1}33$, here $\Pi=1$, so
$\Pi^+-\Pi^-\le \Pi=1$, but $r-k^2=35- 25=10$ and Rohlin's formula
cannot hold.

For the scheme $33$, Rohlin's formula cannot hold as well (here
$\Pi=0$ so $\Pi^+-\Pi^-=0$) and $r$ has to be a square.

Next we have the RKM-scheme $\frac{3}{1}31$. Here $\Pi=3$ and so
$\Pi^+-\Pi^-\le 3$, while $r-k^2=35-25=10$ and again Rohlin's
formula cannot be verified.

Below the former we have $\frac{2}{1}30$. Here $\Pi=2$ and so
$\Pi^+-\Pi^-\le 2$, while $r-k^2=33-25=8$ and Rohlin cannot be
fulfilled.

Below, we have $\frac{1}{1}29$, where  $\Pi=1$ and so
$\Pi^+-\Pi^-\le 1$, while $r-k^2=31-25=6$ and Rohlin cannot be
fulfilled.

Below, we have $29$ which is not realized in type~I, as $r$ is not
a square.

Conclusion: {\it all schemes prohibited by Thom are actually also
prohibited by Rohlin.} (at least within the range of
Fig.\,\ref{Degree10:fig}).

Perhaps Rohlin even prohibits more than Thom. The next boy is the
scheme $\frac{6}{1}30$. Here $\Pi=6$, and so $\Pi^+-\Pi^-\le 6$,
while $r-k^2=37-25=12$ and Rohlin can by now be fulfilled. So no
obstruction.

All this little experiments points out to a subsumation of Thom to
Rohlin (at least for schemes of the form $\frac{x}{1}y$).

\begin{lemma}
For schemes of the type $\frac{x}{1}y$, Thom's inequality is
subsumed to Rohlin's formula.
\end{lemma}

\begin{proof}
For a scheme of this form the total number of pairs (denoted
$\Pi$) is $\Pi= x$, hence $\Pi^+-\Pi^-\le x$. By Rohlin's formula
(\ref{Rohlin-formula:thm}) we infer
$$
2x\ge 2(\Pi^+ -\Pi^-)=r-k^2=(1+x+y)-k^2.
$$
Calculating $\chi$ gives
$$
\chi=1-x+y\le k^2,
$$
by the above estimate.
\end{proof}

Perhaps this is even true in general, but this deserves another
argument. (Update the sequel, will show that the contrary is true,
cf e.g.
Theorem~\ref{Alsatian-scheme-Thom-strong-Petrov-Arnold:thm}.)

First the argument extends to schemes of the form
$\frac{x}{1}\frac{y}{1}z$. Then $\Pi=x+y$, and so $\Pi^+-\Pi^-\le
\Pi=x+y$. Then writing down Rohlin's formula
$$
2(x+y)\ge 2(\Pi^+-\Pi^-)=r-k^2=(2+x+y+z)-k^2,
$$
from which it is inferred that
$$
\chi=1-x+1-y+z=2-(x+y)+z\le k^2,
$$
i.e. Thom's estimate.

And so on, it seems that the passage from Rohlin to Thom will
always succeed as long as the depth is at most one. So for a
counterexample we shall investigate deeper schemes.

For instance in degree 8, we can look at an extension of the
$3$-nest (of depth 3). For instance the $M$-scheme
$(1,1,1)19=(3\times 1) 19 $ (in our satellite notation). Here
$\Pi=3$, and Rohlin's formula $2(\Pi^+ -\Pi^-)=r-k^2=22-16=6$ is
verified for $\Pi^+=3$, and $\Pi^-=0$. But Thom's estimate
$\chi\le k^2=16$ is not (as $\chi=1-1+1+19=20$). So here we get an
example where Thom is not subsumed to Rohlin. However
our example is
artificial
being prohibited by the Gudkov-Rohlin congruence $20=\chi\equiv
k^2=16 \pmod 8$. Yet our example makes unlikely a general
subordination of Thom to Rohlin's formula alone.

A basic idea is to adjust $\chi=20$ at $\chi=16$  to make it
Gudkov-Rohlin compatible. Starting from the above scheme, we may
trade an outer oval for one at depth 1 in the $3$-nest. Each such
trading diminishes $\chi$ by 2, and so two trades are required to
adjust $\chi=16$. Doing this we get the scheme
$(1,\frac{1}{1}2)17$ but alas now Thom's inequality is verified.
(It is also easy to check that Rohlin's formula is satisfied.)

So our game becomes: find a { \it French
scheme\/}, i.e. one prohibited by Thom yet not
succumbing under the armada of Russian prohibitions (Gudkov,
Arnold, and above all Rohlin, and its companions especially
Kharlamov-Marin).

Let us look at degree 10,  and to a deep scheme, say extending the
$4$-nest. As the Harnack bound is $M=g+1=9\cdot 4+1=37$, we look
at the scheme $(1,1,1,1)33=(4\times 1)33$. Now
$\chi=(1-1+1-1)+33=33\equiv k^2=25 \pmod 8$, i.e. the
Gudkov-Rohlin congruence is fulfilled. However the scheme is
prohibited by Thom's estimate $\chi \le k^2$. Further Rohlin's
formula $2(\Pi^+ -\Pi^-)=r-k^2=37-25=12$, and $\Pi=\binom{4}{2}=6$
(count all pair in the $4$-nest) so that $\Pi^+=6$ and $\Pi^-=0$.
So Rohlin's formula affords no prohibition. We have proven the:

\begin{theorem}\label{French-scheme:thm} {\rm (ERRONEOUS, cf.
Corrigendum right below, and for a corrected version cf.
Theorem~\ref{French-scheme-corrected:thm} below)} There exists a
French scheme, i.e. where Thom is not subsumed to Rohlin's formula
nor to the Gudkov-Rohlin congruence mod $8$. The scheme in
question is even an $M$-scheme of degree $10$, namely
$(1,1,1,1)33$. However for schemes of depth $\le 2$, Rohlin's
formula is as strong (and of course stronger) than the Thom
obstruction.
\end{theorem}

[10.03.13] {\it Corrigendum}.---Th. Fiedler objected as follows to
the above theorem (compare his letter dated [09.03.13] in
Sec.\,\ref{e-mail-Viro:sec}):

``The $M$-curve of degree 10 mentioned in your Thm
28.11[=\ref{French-scheme:thm}] is in fact ruled out by Rokhlin's
formula. I think that you have mixed $\Pi^+$ with $\Pi^-$. In a
positive couple the orientations are just opposite. So, four
nested ovals can contribute at most $+2$ to Rokhlin's formula.''

\begin{figure}[h]
\centering
\epsfig{figure=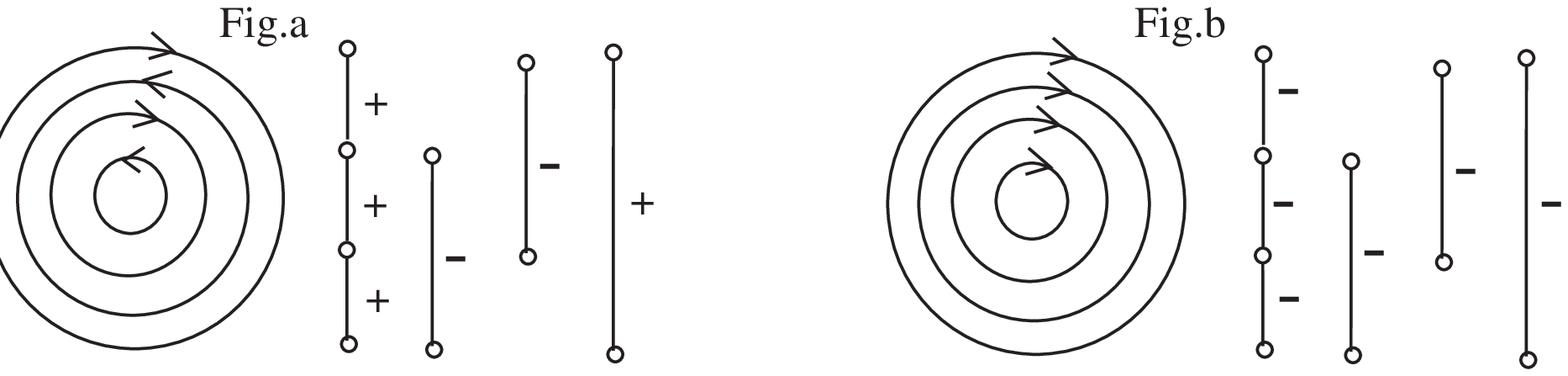,width=122mm}
\vskip-5pt\penalty0
  \caption{\label{Fiedler-correction:fig}%
  Fiedler's correction} \vskip-5pt\penalty0
\end{figure}

After some hesitation, Gabard realized of course that Fiedler is
perfectly right. Let me paraphrase his explanation differently.
This is of course allied to the signs-law of
Fig.\,\ref{Signs-law-dyad:fig}, but let us be more specific. If we
consider a 4-nest and choose on it complex orientations forming
positive pairs at each immediately successive nested ovals then we
get  Fig.\,\ref{Fiedler-correction:fig}a. Pairs of length 2
becomes negative, while the unique pair of length 3 is positive
again. (This is seen either by looking at the picture or if one
like extracting an arithmetical law one finds the twisted
signs-law of (\ref{Signs-law:lem})
akin to usual arithmetics modulo a twisted sign, so $+\times
+=-$,\; $+\times -=+$,\; $-\times +=+$,\; $-\times -=-$. This is
best memorized by saying that ``mixing the genes is good, while
consanguinity is bad''!) So the contribution to $\Delta
\Pi:=\Pi^{+}-\Pi^{-}$ is at most 2, and certainly never equal to
$6$ (though it can be $-6$ as on
Fig.\,\ref{Fiedler-correction:fig}b).

It remains now to see how Fiedler's remark generalizes as to see
if Thom's estimate $\chi\le k^2$ is a formal consequence of
Rohlin's formula. If not then it remains to find another (more
serious!) French scheme.

The naive scenario would be that Rohlin always implies Thom (at
least within the realm of Hilbert's 16th). The only chance to
prove this seems to involve an estimation of the corrector term in
the ``Rohlin-to-Arnold formula'' \eqref{Rohlin-to-Arnold:eq},
which we reproduce for convenience
\begin{align}\label{Rohlin-to-Arnold-bis:eq}
\chi=p-n=(p+n)-2n&=r-2n\cr
&=[2(\Pi^+-\Pi^-)+k^2]-2n\cr
&=k^2+2(\Pi^+ -\Pi^- -n).
\end{align}
So setting $\Delta \Pi:=\Pi^+ -\Pi^-$ we would like to show the:

\begin{conj} {\rm (Garidi\footnote{Tarik Garidi (aus der
Nordseek\"uste) is a well-known scientist in Geneva (student of
Piron),
 specialized in anti-de-Sitter and notorious for having
introduced a mass concept which can be negative-valued, like the
signed difference $\Delta \Pi:=\Pi^+-\Pi^-$ of Rohlin.} mass
conjecture.)}---It holds universally $\Delta \Pi \le n$.
\end{conj}

Of course the conjecture is true (at least if one believes in the
Kronheimer-Mrowka validation of the Thom conjecture) [THIS
REASONING IS DUBIOUS BEING BASED ON OUR FALSE THOM ESTIMATE], so
that the true meaning of our conjecture is an independent
derivation of the estimate via pure combinatorics. Let us be more
precise.

\begin{defn}
{\rm A  signed or {\it Rohlin  tree} is a combinatorial object
consisting of a (finite) directed set plus a distribution of signs
$\pm$ on its edges such that the signs-law (of
Lemma~\ref{Signs-law:lem}) is verified. By a directed set we mean
a finite POSET such that each element as at most one superior,
i.e. an element larger and minimal with this property. Recall also
that the signs-law can be easily remembered by saying that
consanguinity is bad, i.e. $+\times +=-$, $-\times -=-$, while
mixing the genes is good $+\times -= +$ and $-\times += +$ (this
exotic signs-law is the exact opposite of the usual convention).}
\end{defn}

Of course this concept arises naturally when taking a smooth
dividing plane (algebraic) curve $C_m$ of even degree $m=2k$ and
assigning to it its Hilbert tree (encoding the distribution of
ovals), while decorating the edges with signs coming from the
complex orientations as in Rohlin's formula
(\ref{Rohlin-formula:thm}).

More precisely, given an oriented real scheme (i.e. an isotopy
class of embedding of a disjoint union of circles in $\RR P^2$
supplied with an orientation), we can assign to it a Rohlin tree.
Conversely it is clear that any Rohlin tree arises in this
fashion. To a real plane dividing curve is assigned a complex
orientation (uniquely defined up to reversal of all orientations),
yet this leaves invariant the concept of positive or negative
pairs as defined by Rohlin. Hence to be very formal, we have first
the map taking a dividing real curve to its real scheme with
complex orientation, which is a ``projectively'' oriented real
scheme (weel-defined up to reversing all orientations), which in
turn defines unambiguously a Rohlin tree. Diagrammatically,
$$
\textrm{dividing plane}  \to \textrm{oriented real schemes } \to
\textrm{Rohlin trees}.
$$
\vskip-23pt
$$
\hskip-80pt\textrm{curves ($m=2k$) \hskip0.5cm (mod reversion) }
$$

So the precise meaning of the above conjecture is that any Rohlin
tree (not necessarily induced by a real algebraic curve) satisfies
the above estimate $\Delta \Pi \le n$. After several hours of
attempting to prove this ``Garidi mass conjecture'', one finds a
simple counterexample as follows.

\begin{figure}[h]
\centering
\epsfig{figure=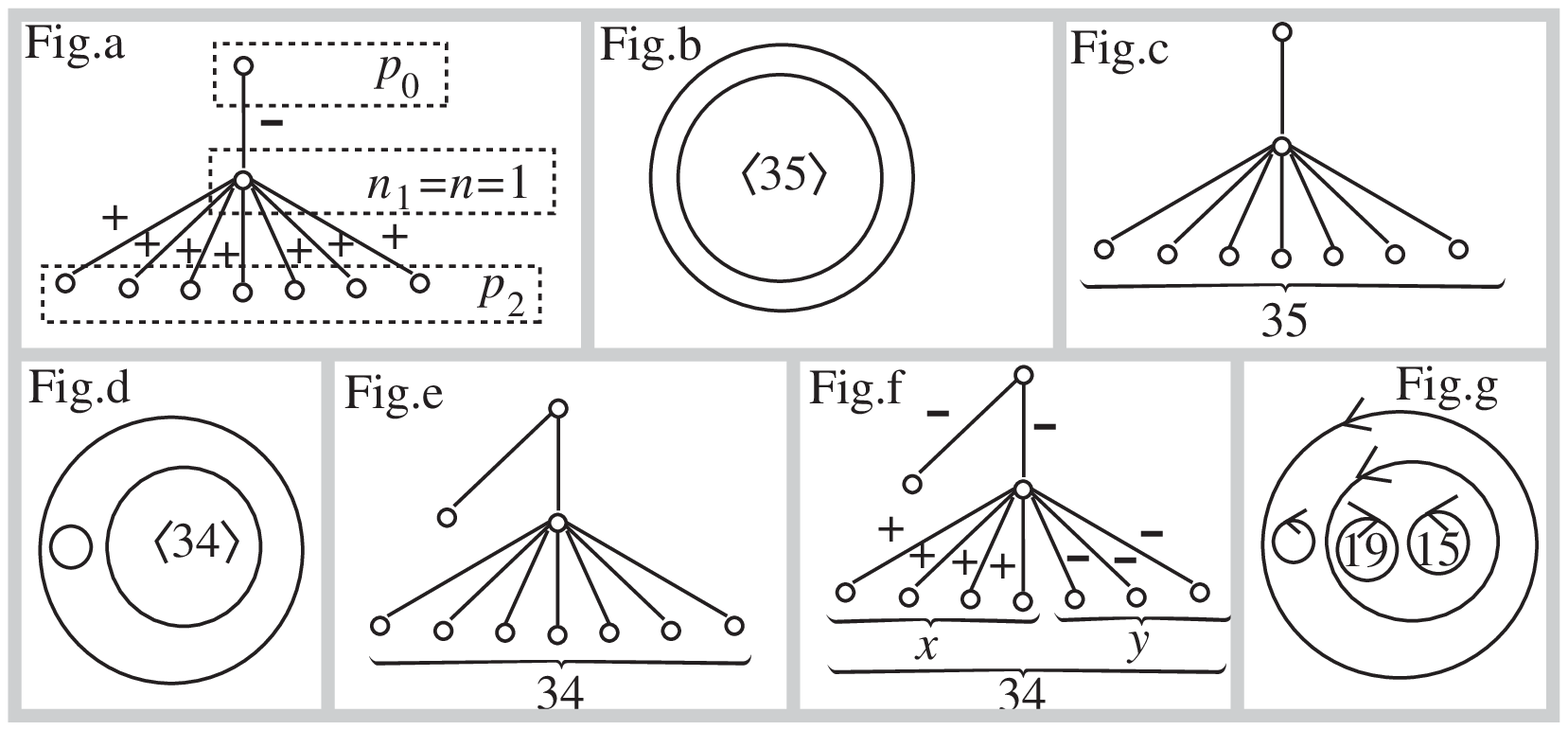,width=122mm} \vskip-5pt\penalty0
  \caption{\label{Garidi-mass-false:fig}%
  A (basic) counterexample to the Garidi mass conjecture} \vskip-5pt\penalty0
\end{figure}

To keep $n$ small (say $n=1$), we consider a tree with only one
vertex at depth 1, but ramifying (violently) at depth 2 (cf
Fig.\,\ref{Garidi-mass-false:fig}a).

Let us introduce the sign distribution of Fig.\,a, and denote by
$p_2$ the number of (even) vertices at depth 2. We find by using
the twisted signs-law:
$$
\Delta \Pi= (-1)+p_2+p_2=2p_2-1.
$$
(Here the first $-1$ comes from the top edge (visible on Fig.\,a),
the second term $p_2$ is the contribution of the $+$-signs visible
on Fig.\,a), while the 3rd term $p_2$ comes from the $p_2$ pairs
of length 2 obtained by concatenation of elementary edges. The
signs-law in question (based on Fig.\,\ref{Signs-law-dyad:fig}) is
the same as the usual one modulo a twist by $-1$. So here $+\times
-$ gives $+$ (the opposite of the usual sign rule!). The displayed
formula is justified.)

Now as soon as $p_2\ge 2$, the above $\Delta \Pi $ will be $\ge
3$, foiling thereby the mass conjecture.

{\it Insertion} [22.03.13] For later reference, let us state this
as a:

\begin{theorem} \label{Garidi-mass-conj-is-FALSE:thm}
The Garidi mass conjecture is false, and therefore the positive
mass conjecture (\ref{positive-mass-conjecture:conj}) is erroneous
too. For a simple counterexample cf.
Fig.\,\ref{Garidi-mass-false:fig}a right above.
\end{theorem}

This basic corruption aids us to detect a more serious French
scheme (where Thom is not subsumed to Rohlin). As above we
consider curves of degree $m=10$. Harnack's bound is
$M=g+1=\frac{9\cdot 8}{2}+1=9\cdot 4 +1=37$ (temperature of the
human body). Our tree converts then to the scheme of
Fig.\ref{Garidi-mass-false:fig}\,b, where we see 2 nested ovals
containing 35 unnested ones. The characteristic of the ``Ragsdale
membrane'' is $\chi=(1-1)+35=35 \nleqslant 25=k^2$, so that the
scheme is prohibited by Thom. Is this scheme (with Gudkov symbol
$(1,1,35)$) prohibited by Rohlin's formula?

Remember that Rohlin's formula implies Arnold's congruence (cf.
(\ref{Rohlin-implies-Arnold:lem})), so the answer is an (indirect)
yes since $\chi\equiv k^2=25\neq 35 \pmod 4$.

As we have an $M$-scheme, let us even adjust to the (stronger)
Gudkov congruence mod 8: $\chi\equiv k^2=25=33 \pmod 8$. So in
order to diminish $\chi$ by 2 (from 35 to 33), let us trade a deep
oval (at depth 2) against one at depth 1 (cf. transition from
Fig.\,b to Fig.\,d). We have now $\chi=33$, and so Thom is still
violated. Is this new scheme (symbol $(1,1(1,34))$)  prohibited by
Rohlin's formula?

We count (e.g. via Fig.\,e) that the total number of pair is
$\Pi=2+34+34=70$ (this can be viewed as an application of the
formula $\Pi=n_1+2p_2+3n_3+4p_4+ etc$, cf. proof of
(\ref{Stalin:lemma})). To abridge Rohlin's (heavy ``Cyrillic'')
notation let us set $\pi:=\Pi^+$, and $\eta:=\Pi^-$, and Rohlin's
formula then reads $2(\pi-\eta)=r-k^2=37-25=12$. So we get the
pair of equations $\pi-\eta=6$ and $\pi+\eta=\Pi=70$. Adding them
gives $2\pi=76$, whence $\pi=38$, and $\eta=32$. So Rohlin's is
(formally) soluble but is there a distribution of signs compatible
with the signs law?

To answer this let us consider a ``variable'' distribution of
signs like on Fig.\,f with $x$ many $+$ and $y$ many $-$ for the
edges rooted ``at depth 1'', while both edges rooted at depth 0
have signs $-$. Of course we assume $x+y=34$. Counting the number
of positive pairs $\pi$ we find $\pi=x+x=2x$ (were the second $x$
term comes from the sign-law $+\times -= +$ the opposite of the
usual convention!). For the number $\eta$ of negative pairs we
find $\eta=2+y+y=2+2y$, where the 3rd $y$ term comes again from
the exotic sign-law (``of Rohlin''). Combining with the previous
paragraph, gives $x=19$ and $y=15$. All equations are then
verified! Conclusion there is a distribution of signs on the tree
which satisfies Rohlin's formula, which therefore does nor
prohibit the scheme under examination, i.e. $(1,1(1,34))$. The
latter is therefore a French scheme. So we hope to have this time
proven the:

\begin{theorem}\label{French-scheme-corrected:thm}
There exists a ``French'' $M$-scheme of degree $10$, namely that
with Gudkov symbol $(1,1(1,34))$ (cf. {\rm
Fig.\,\ref{Garidi-mass-false:fig}d\/} above), i.e. which is
prohibited by Thom but not by the armada of Russian congruences
(especially Gudkov's) nor by Rohlin's formula.
\end{theorem}

{\it Insertion}---[17.03.13] Thomas Fiedler kindly reacted as
follows to this statement, cf. his [12.03.13]-letter in
Sec.\,\ref{e-mail-Viro:sec} reproduced below for convenience (our
brackets are just automatized updates of labels):

``sorry, but all your $M$-schemes of degree 10 in
Thm~30.14[=\ref{French-scheme-corrected:thm}] and
30.15[=\ref{French-scheme-corrected-bis:thm}] have $n=2$ and are
ruled out simply by Petrovskis inequality. I don't think that
genus bounds give anything new for real schemes
alone\footnote{[18.03.13] I think that (modest)
Theorem~\ref{Alsatian-schemes:thm} below corrupts this belief of
Th. Fiedler (who left the subject a long time ago), yet this does
not jeopardize at all his invaluable help (and incredible memory!)
in view  of all the crucial corrections he took care to make on
the present text.} but they definitely do so for configurations of
several real curves. Just take a look on Mikhalkin's paper.''

Though a pertinent remark (since Petrovskii is not a French guy),
Fiedler's remark does not affect the modest truth of our statement
but points to the Petrovskii inequality as another sharp Russian
weapon. (Shamefully, I confess to have not properly appreciated
this fundamental statement prior to Fiedler's comment.) The latter
states:

\begin{theorem} {\rm (Petrovskii 1933/38 \cite{Petrowsky_1938})}
\label{Petrovskii's-inequalities:thm} For any real plane smooth
curve of even degree $m=2k$ we have the (so-called) Petrovskii
inequalities (which are pure jewels nearly coming out of the blue
safe for having been apparently  anticipated conjecturally by Miss
Ragsdale in 1906)
$$
- \textstyle\frac{3}{2} k(k-1)\le \chi \le \textstyle\frac{3}{2}
k(k-1)+1.
$$
(An analogous but more complicated statement holds for curves of
odd degrees.)
\end{theorem}

\begin{proof} We make just some few remarks.

{\it Historical substance.}---Petrovskii 1938
\cite[p.\,191]{Petrowsky_1938} comments that his method of proof
is based on two ingredients:

(1) a formula of Jacobi-Euler 1768--70 (and also cite en passant
Kronecker, which as we know is also one of the forerunner of
Poincar\'e's index theory via Hermite's transmissive r\^ole), and,

(2) on the consideration of the deformations of lines $F(x,y)=C$
when $C$ crosses the critical values of $F(x,y)$. These last
investigations being identified as  analogous to those of Morse
(1925 \cite{Morse_1925?}) on the critical points of a function.

{\it Neo-expressionist proofs}.---Another proof of Petrovskii's
inequalities namely the ``Preuve d'Arnol'd dans une pr\'esentation
de A. Marin'' is given in A'Campo 1979
\cite[p.\,537--17]{A'Campo_1979}, where the result is stated as:
$$
\vert 2\chi-1 \vert \le \textstyle\frac{3m^2}{4}-\frac{3m}{2}+1,
$$
which (after setting $m=2k$) is readily seen to be equivalent to
the above formulation. Indeed $\vert 2\chi-1 \vert \le
\textstyle\frac{3(2k)^2}{4}-\frac{3(2k)}{2}+1=3k^2-3k+1=3k(k-1)+1$.
Hence $3k(k-1) \le 2\chi\le 3k(k-1)+2$, and the equivalence is now
obvious.

{\it The original statement differs slightly.}---On adapting to
our (fairly standard modern) notations, Petrovskii's original
result is stated as follows (cf. p.\,190 of Petrovskii 1938
\cite{Petrowsky_1938})
$$
\vert p-n \vert \le \textstyle\frac{3m^2-6m}{8}+1.
$$
As $m=2k$ and $\chi=p-n$, this gives indeed
$$
\vert \chi \vert\le 3
\textstyle\frac{(2k)^2-2(2k)}{8}+1=\textstyle\frac{3}{2}(k^2-k)+1
=\textstyle\frac{3}{2}k(k-1)+1,
$$
which is essentially the announced bound modulo a discrepancy on
the lower-bound by one unit. In fact we copied the stated lower
bound from Rohlin 1978 and hope that there is no misprint there.

Let us look at the example $k=3$ of sextics. Then $\frac{3}{2}
k(k-1)=\frac{3}{2} 3\cdot 2=9$, and so $-9\le \chi $, hence even
the stronger version written down by Rohlin 1978 (also in Wilson
1978, p.\,55) does not prohibit ``Rohn's scheme'' $\frac{10}{1}$
(cf. Gudkov's Table=Fig.\,\ref{Gudkov-Table3:fig}).
\end{proof}

Applied to our situation $k=5$, Petrovskii's theorem shows that
$\chi \le \frac{3}{2} 5\cdot 4+1=31$ and so our scheme with
$\chi=33$ is prohibited by Petrovskii.

In general, for an even degree $m=2k$ we have Harnack's bound
$M=g+1=(2k-1)(k-1)+1=2k^2-3k+2\approx 2k^2$, the universal
Petrovskii's bounds $-P\le \chi\le P+1$, where $P=\frac{3}{2}
k(k-1)\approx \frac{3}{2} k^2$, and finally Thom's bound $\chi \le
k^2 $ for dividing curves (only). So when $k$ is large the
``Hilbert-Petrovskii-Gudkov'' pyramid looks as follows
(Fig.\,\ref{PyramidPetrov:fig}), and of course Thom will
asymptotically be stronger than Petrovskii (at least for dividing
curves and on the right-wing of the pyramid where $\chi$ is
positive).

\begin{figure}[h]
\centering
\epsfig{figure=PyramidPetrov.eps,width=122mm} \vskip-5pt\penalty0
  \caption{\label{PyramidPetrov:fig}%
  Petrovskii versus Thom} \vskip-5pt\penalty0
\end{figure}

So I am not sure not adhere completely with Fiedler's illuminating
comment, because if we take a larger $m=2k$ than $10$ then it will
be possible to arrange Petrovskii's bound yet not Thom's one,
while  further taking care of respecting Gudkov's congruence and
Rohlin's formula.

Let us first take $m=12$ (so $k=6$) then Petrovskii's upper-bound
is $P+1=\frac{3}{2} 6\cdot 5+1=46$, while Thom's is the sharper
$k^2=36$. Arranging Gudkov's hypothesis $\chi\equiv k^2 \pmod 8$
permits to take $36+8=44$ which is still lower than Petrovskii's
upper-bound. Now $M=g+1=\frac{11\cdot 10}{2}+1=56$. Hence to
arrange $\chi=44$, we transplant 6  ovals of the unnested
configuration at depth 1 (each such move drops $\chi$ by 2 units)
to get the (12)-scheme\footnote{Our notational trick is to denote
with parenthesis the degree of the scheme, since without
parenthesis e.g. in an $M$-scheme the magnitude in front is
traditionally not the degree but the number of ovals of the
scheme.} $(1,6)49$ (with $\chi=44$). By Rohlin's formula
$2(\pi-\eta)=r-k^2=56-36=20$, hence $\pi-\eta=10$. As there are no
deep nesting the signs-law is negligible and we merely have
$\pi+\eta=6$, so that $2\pi=16$, whence $\pi=8$ and we get an
obstruction.

\smallskip

[17.03.13]{\it Optional reading (skip if you do not want to loose
the main-flow and move to $\clubsuit\clubsuit$)}.---This argument
extends to the following formal consequence of Rohlin (probably
subsumed to Thom, yet much more elementary).

\begin{lemma}\label{Rohlin-consequence-for-M-curves:lem}
An $M$-curve of degree $2k$ and of type $\frac{x}{1}y$ with few
nested ovals in the sense that $x<\frac{(k-1)(k-2)}{2}$ is
prohibited by Rohlin's formula. In particular there is no unnested
$M$-curve provided $k\ge 3$.
\end{lemma}

\begin{proof} First recall that Harnack's bound is $r=M=g+1=
(2k-1)(k-1)+1=2k^2-3k+2$. By Rohlin's formula
$2(\pi-\eta)=r-k^2=k^2-3k+2=(k-1)(k-2)$. So
$\pi-\eta=\frac{(k-1)(k-2)}{2}=:\binom{k-1}{2}$. But $\pi+\eta=x$,
so that $2\pi=x+\binom{k-1}{2}$, and hence
$\pi=(x+\binom{k-1}{2})/2$. Yet the equation $\pi+\eta=x$ is
impossible whenever $\pi>x$, that is when
$\pi=(x+\binom{k-1}{2})/2>x$, i.e. as $(x+\binom{k-1}{2})>2x$, so
when $x<\binom{k-1}{2}$, which is the asserted condition.
\end{proof}

\begin{figure}[h]
\centering
\epsfig{figure=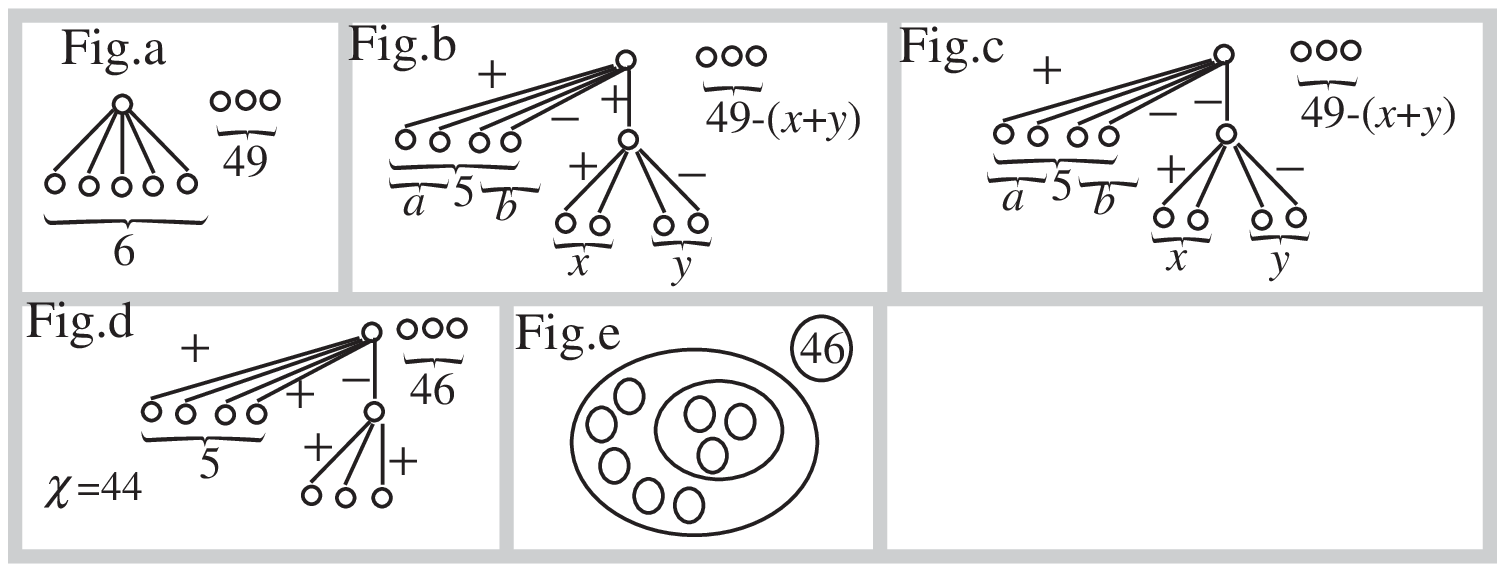,width=122mm} \vskip-5pt\penalty0
  \caption{\label{Fied3:fig}%
  Looking for a  French scheme in degree
  12 compatible with Petrovskii} \vskip-5pt\penalty0
\end{figure}

$\clubsuit\clubsuit$ [18.03.13] Back to our main object of the
(12)-scheme $(1,6)49$, our idea is to remove this Rohlin
obstruction by injecting more freedom gained by transferring some
ovals at depth 2 (leaving thus $\chi$ unchanged). So starting from
the (12)-scheme $(1,6)49$, whose tree is depicted as
Fig.\,\ref{Fied3:fig}a, we transplant ovals at depth 2 to get
Fig.\,b with a certain quantity $x+y$ of ovals at depth 2. By
Rohlin's formula $\pi-\eta=10$, and from Fig.\,b we have
$\pi+\eta=6+2(x+y)$. Adding the last equations gives
$2\pi=16+2(x+y)$, whence $\pi=8+(x+y)$. The condition $\pi\le
\pi+\eta$ becomes so $8+(x+y)\le 6+2(x+y)$, i.e. $x+y\ge 2$, which
is necessary to solve Rohlin's equation.

Now introduce signs as on Fig.\,b by putting $a$ many pluses on
the 5 edges at depth 1 and $b$ many minuses (so $a+b=5$), and a
$+$-sign  on the trunk ({\it Warning.}---This choice is fatally
bad as we shall see, and a minus sign is more fruitful as we shall
experiment soon). Applying the signs-law (cf.
Fig.\,\ref{Signs-law-dyad:fig}) to Fig.\,b we find by sorting out
contribution according to their edge-length (which appears
underbraced)
$$
\pi-\eta=\underbrace{a-b+1+x-y}_1+\underbrace{(-x+y)}_2=a-b+1.
$$
But as $\pi-\eta=10$, we deduce the system $a-b=9$, $a+b=5$,
whence $2a=14$ and $a=7$ which is incompatible with $a+b=5$.

(Then we repeated such calculation in higher degrees $14, 16$
finding always the same obstruction, though for $m=16$ there are
even two possible values of $\chi$ permissible under Gudkov and
Petrovskii).

However if we take a $-$-sign on the trunk (Fig.\,c), then the
signs-law gives
$$
\pi-\eta=\underbrace{a-b-1+x-y}_1+\underbrace{(x-y)}_2=
a-b-1+2(x-y).
$$
As by Rohlin we still have $\pi-\eta=10$, this gives the system
$a-b=11-2(x-y)$, $a+b=5$, whence $2a=16-2(x-y)$, i.e. $a=8-(x-y)$.
Let us fix $x-y=3$, so that $a=5$, $b=0$. Since $x+y\ge 2$ (cf.
necessary condition discussed few lines above), we may choose
$x=3$ and $y=0$. It is worth at his stage checking that the
resulting signs distribution indeed solves Rohlin's equation, and
we have proven:

\begin{theorem}\label{Alsatian-schemes:thm} There exists
an Alsatian scheme,
i.e. a  French scheme which furthermore respects Petrovskii's
inequalities. More precisely, there is an $M$-scheme of degree
$12$ for instance $(1, 5(1,3))46$ (cf. {\rm
Figs.\,\ref{Fied3:fig}d,e\/}) which respects both the Petrovskii
bound, the Gudkov congruence and Rohlin's formula, yet which is
prohibited by Thom's bound.
\end{theorem}

{\it Insertion}.---[20.03.13] Thomas Fiedler was kind enough to
object (once more) to this result as follows (cf. his message
dated 19 March 2013 in Sec.\,\ref{e-mail-Viro:sec}, yet reproduced
here for convenience):

``sorry again, but your curve has $p=50$ and is ruled out by
Arnold's inequality : $p\le 3/2 k(k-1) + 1 + n_-$\footnote{The
usual notation is $n^-$ at least to be conform with Rohlin 1978
\cite[p.\,86--87]{Rohlin_1978}, but I keep it so to stay faithful
to the message of Fiedler.}, which is $47$ in this case. In fact
Arnold's inequalities are by fare the strongest result in the
whole field.''

\smallskip
We will try to react to this objection right below the proof, cf.
$\bigstar\bigstar$ below.

\smallskip

\begin{proof}[Proof of (\ref{Alsatian-schemes:thm})] The assertion is clear by our search, but as it is
easy to make mistakes, let us do an ad hoc self-contained
verification. The scheme in question is depicted on
Fig.\,\ref{Fied3:fig}e. It has $\chi=(1-6+3)+46=44$. Petrovskii's
inequalities says $-P\le\chi\le P+1$, where
$P:=\frac{3}{2}k(k-1)=\frac{3}{2}6\cdot 5=45$, and this is
satisfied by our scheme. Gudkov's congruence $\chi\equiv_8
k^2=6^2=36$ is also verified. Finally Rohlin's formula is
fulfilled when like on Fig.\,\ref{Fied3:fig}d all signs are
positive safe that on the trunk at depth 1 which has further
ramifications at depth 2. (The depth of an edge is that of the
unique vertex below it.) In that case the signs-law gives
$$
\pi-\eta=\underbrace{5-1+3}_1+\underbrace{3}_2,
$$
where contributions are underbraced along the length of the pairs.
Hence $\pi-\eta=10$ in accordance with Rohlin's formula
$2(\pi-\eta)=r-k^2=56-6^2=20$. However the scheme in question is
prohibited by Thom's estimate $\chi\le k^2=36$.
\end{proof}

It should be easy  to extend the result to other schemes but it
looks artificial to strive toward maximum generality as our
purpose was merely to find an example where Thom affords valuable
information.

\smallskip

$\bigstar\bigstar$ {\it Trying to fix Fiedler's new objection
based on Arnold's strong Petrovskii inequalities}.---The result
mentioned by Fiedler is the following sometimes called the {\it
strong Petrovskii inequalities}. Those are really due to Arnold
1971 \cite{Arnold_1971/72}, and sharper than Petrovskii's original
inequalities of 1933/38. Apparently (cf. Rohlin 1978
\cite{Rohlin_1978}, p.\,87, footnote), Arnold's original statement
contained further unnecessary restrictions that were relaxed in
Rohlin 1974 \cite{Rohlin_1974/75}. The final shape of the result
is as follows:

\begin{theorem}\label{Strong-Petrovskii-Arnold-ineq:thm} {\rm (Strong Petrovskii inequalities,
aka Arnold inequalities).---(Arnold 1971, Rohlin 1974)}.---For any
curve of even degree $m=2k$,
$$
n-p^- \le \textstyle\frac{3}{2} k(k-1), \qquad p-n^- \le
\textstyle\frac{3}{2} k(k-1)+1,
$$
where $p^-, n^-$ are  the number of positive=even resp.
negative=odd ovals which are hyperbolic (cf. {\rm
Definition~\ref{hyperbol-ovals:def}\/} right below).
\end{theorem}

\begin{defn} \label{hyperbol-ovals:def}
$\bullet$ {\rm An oval of a plane real algebraic curve (or a
scheme=distri\-bution of ovals) is {\it elliptic}, {\it parabolic}
or {\it hyperbolic}\footnote{As far as I am informed the general
coinage of this trichotomy is due to Felix Klein (in geometry) and
Dubois-Reymond (in PDE's).} depending on whether its
poros\footnote{From the Greek ``poros''=``hole'' (aping a bit
Grothendieck's ``topos'' or ``topoi'').} (cf. below) has positive,
zero or negative Euler characteristic\footnote{It is crucial here
to adopt the modern convention regarding the sign of Euler's
$\chi$. This is courtesy of Michel Kervaire, that turned out  to
be correct when looking at old texts, like perhaps Listing, von
Dyck 1988, Poincar\'e 1885--1895, etc., where the opposite sign
convention was used!}.

$\bullet$ The {\it poros} of an oval of a plane curve $C_m(\RR)$
is
the inside of the oval
minus the insides of all ovals immediately nested in the given one
(equivalently  remove the insides of all subordinated ovals).}
\end{defn}

An oval is hyperbolic iff the Hilbert tree of the scheme ramifies
at the corresponding vertex. This basic remark is the key to the
little problem suggested by Fiedler's objection.

Now our Alsatian question is again: is there a scheme where Thom
is stronger than the conjunction of strong-Petrovskki=Arnold,
Gudkov's hypothesis and Rohlin's formula.
Our basic algorithm to do this is always same:
\smallskip

(1) Start from the unnested configuration of $M$-ovals.

\smallskip
(2) Then adjust $\chi$ to $k^2$ as to verify Gudkov's hypothesis.
(This can be done by transplanting outer ovals at depth 1,
dragging them say inside a fixed oval.)

\smallskip
(3) Then Thom's estimate $\chi\le k^2$ is verified, but corrupt it
by incrementing $\chi$ by 8. Do this as many times as Petrovskii's
bound $\chi\le \frac{3}{2}k(k-1)+1$ permits.

\smallskip
(4) Next without changing $\chi$ transplants outer ovals at depth
2 (as many as you want), while introducing a branched structure
making Rohlin's equation soluble. (This being essentially inspired
by our counter-example to the Garidi mass conjecture discussed
above.)
\smallskip

\begin{figure}[h]
\centering
\epsfig{figure=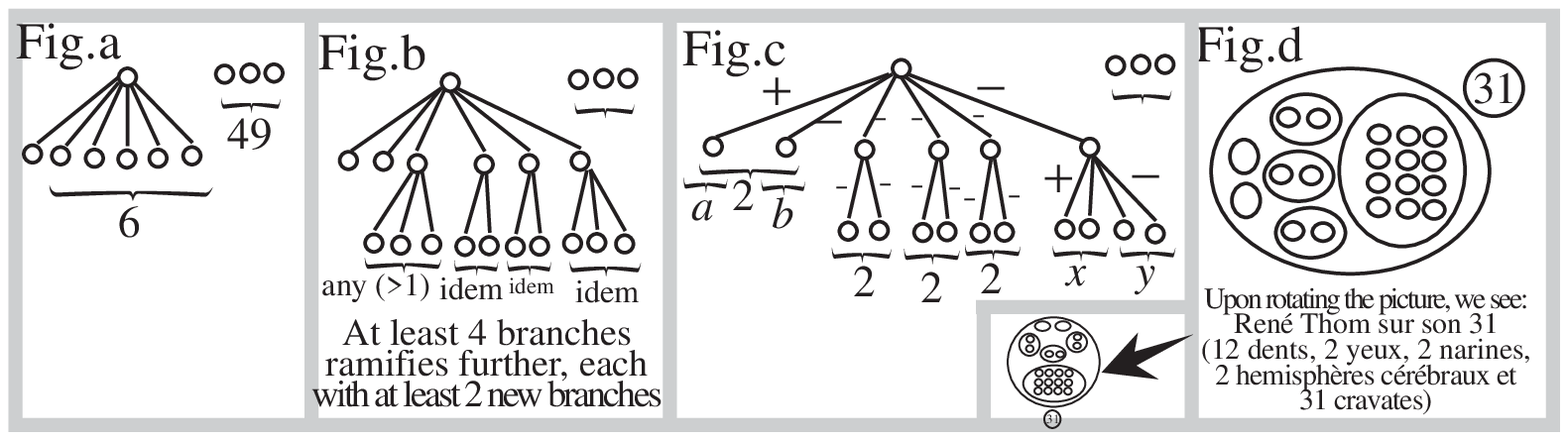,width=122mm} \vskip-5pt\penalty0
  \caption{\label{Fied4:fig}%
  Constructing an Alsatian scheme of degree 12 where Thom is stronger
  than all the Russian estimates, congruences and formulas:
  alias ``Ren\'e Thom sur son 31''} \vskip-5pt\penalty0
\end{figure}

Let us be more explicit. Suppose $k=6$ (so $m=2k=12$). Harnack's
bound is $M=g+1=\frac{11\cdot 10}{2}+1=55+1=56$. Adjust to
$\chi=k^2=36$ (Gudkov) and increment by 8 to get $44,52$.
Petrovskii 1938 says $\chi \le \frac{3}{2}k(k-1)+1=3\cdot 3\cdot
5+1=46$. So we consider $\chi=44$. Hence we transplant from the
unnested scheme $M=56$, precisely $(M-\chi)/2=(56-44)/2=6$ ovals
at depth 1 to get adjusted at $\chi=44$. This gives the scheme on
Fig.\,\ref{Fied4:fig}a with symbol $(1,6)49$. This scheme is
prohibited by Rohlin's formula $2(\pi-\eta)=r-k^2=56-36=20$, hence
$\pi-\eta=10$, while $\pi+\eta=6$. This is of course soluble as
$\pi=8$ and $\eta=-2$ (ruled out as $\eta\ge 0$ is a cardinal!).
Variant: imagine signs on the edges of the tree (of Fig.\,a) we
can have at most 6 pluses and so $\pi-\eta\le 6$ hence cannot be
$10$. To arrange Rohlin the idea is to transplant ovals at depth 2
creating thereby more freedom to solve Rohlin's equation.

The extra challenge is to take care of the (strong)
Petrovskii-Arnold estimate. Note first that our configuration
$(1,6)49$ has $p=50$ and this will not change under transfers at
depth 2. Hence to respect Arnold's estimate
$$
p\le \textstyle\frac{3}{2}k(k-1)+1+n^-=46+n^-,
$$
it suffices to arrange $n^-=4$. Remember that $n^-$ counts the
number of hyperbolic negative(=odd) ovals, so  we just  have to
transplant ovals at 4 different places of Fig.\,a to get a tree
 like Fig.\,b with at least two extra branches growing at
each 4 places. To ensure hyperbolicity it is sufficient to have
branchings of ``order 2'', so we look to Fig.\,c where the 3 first
branches ramifies by 2, while the fourth by a magnitude $x+y$ (yet
undetermined) safe for $x+y\ge 2$. On the left of the tree (still
Fig.\,c), we put $+$-signs at $a$ many places (w.l.o.g. on the
``left'' though this has no intrinsic meaning here), and $-$-signs
at $b$ many places. Hence $a+b=2$. Likewise on the right of the
tree we introduce $x$ many $+$ resp. $y$ many $-$ as depicted on
Fig.\,c. On the ``center'' of the tree where we have 3 branches
like ``Y''-letters inverted, we plug everywhere $-$-signs as those
depressive guys are simplest to calculate by the signs-law
($-\times -=-$). This rigidification looks a reasonable Ansatz for
there should be already enough free parameters available with
$x,y$ and $a,b$ (subjected to $a+b=2$). Once the combinatorics of
the tree is fixed and the signs-distribution too (modulo the free
parameters) we can compute the Rohlin mass $\pi-\eta$ according to
the signs-law. First note that each 3 central subtrees
(``Y-shaped'') contributes for 5 pairs (3 visible, plus 2
concatenation) each being negatively charged. Hence the
contribution of each such subtree is $-5$. Globally on the whole
tree,  we find therefore (upon remembering the signs-law
(\ref{Signs-law:lem}) saying that mixing the genes is good so
$+\times -=+$, while consanguinity is bad, e.g. $-\times -=-$)
$$
\pi-\eta=a-b-15-1+2x-2y=a-b-16+2(x-y).
$$
By Rohlin's formula $\pi-\eta=10$, and we get the system
$a-b=26-2(x-y)$, $a+b=2$. Hence $2a=28-2(x-y)$, i.e. $a=14-(x-y)$.
Choose (freewill vs. predestination!) $x-y=12$, so that $a=2$,
$b=0$, and choose again $x=12$, and $y=0$ as a special solution.
This proves the:

\begin{theorem}\label{Alsatian-scheme-Thom-strong-Petrov-Arnold:thm} There exists an Alsatian $M$-scheme where Thom is
stronger than the conjunction of (strong) Petrovskii-Arnold 1971,
Gudkov hypothesis 1969--72 (proved by Rohlin-Marin), and Rohlin's
formula. Specifically, there is such a scheme in degree $12$,
namely the one allied to {\rm Fig.\,\ref{Fied4:fig}\,c\/} for
$x+y=12$, whose scheme is depicted on {\rm Fig.\,d\/} called
``Ren\'e Thom sur son $31$'' since the nested portion of the tree
involves (counting along increasing depths)
$p_0+n_1+p_2=1+6+18=25$ ovals so that it remains left $56-25=31$
outer ovals. The $M$-scheme in question has Gudkov symbol
$(1,2(1,2)(1,2)(1,2)(1,12))31$.
\end{theorem}

\begin{proof} Let us do again an ad hoc verification. From the
scheme (Fig.\,c with $x+y=12$ or Fig.\,d) we have
$\chi=1-6+18+31=44$. So Gudkov $\chi\equiv_8 k^2$ is happy.
Petrovskii 1938, i.e. $\chi \le \frac{3}{2}k(k-1)+1=45+1=46$ is
also satisfied. Now the strong version of Petrovskii due to
Arnold, reads
$$
p-n^-\le \textstyle\frac{3}{2}k(k-1)+1=46.
$$
But our scheme has (cf. Fig.\,c or d) altogether $r^-=5$
hyperbolic ovals (i.e. containing immediately at least 2 other
ovals in their insides), yet only 4 of them are at odd depth 1, so
$n^-=4$. Though equivalent, working with the tree (as opposed to
the scheme) looks more convenient to see this (at least when one
is tired). On the other hand our scheme has either by its
construction (cf. Fig.\,a) or by counting $p=50$, since
$p=p_0+p_2=(1+31)+(6+12)=32+18=50$. Hence the Petrovskii-Arnold
estimate is verified.

Note also that the other strong Petrovskii inequality $n-p^-\le
\frac{3}{2}k(k-1)=45$ is verified, as $n=6$ (either from Fig.\,c
or from $p+n=r=M=56$, where $r$ is as usual in our notation the
number of ``reellen Z\"uge'', denoted $l$ in Russia). Although not
needed it may be noted that $p^-=1$ either via Fig.\,c or via the
relation $p^- + n^-=r^-$ (splitting hyperbolic ovals according to
their parities).

Finally Rohlin's formula is verified for the signs-distribution of
Fig.\,c where $x=12$, $y=0$ and $a=2$. Indeed the signs-law gives
$\pi-\eta=2-3\cdot 5+(-1)+12+12=2-16+24=10$, in accordance with
Rohlin's formula $2(\pi-\eta)=r-k^2=56-36=20$. The verification is
complete.
\end{proof}

\subsection{Old material (to skip or reorganize)}

{\it Sequel of my text (prior to Fiedler's objection(s) via
Petrovskii, and then via Arnold)}.---One can paraphrase the
statement (\ref{French-scheme-corrected:thm}) by saying that there
is a complex scheme (i.e. with orientation) which satisfies
Rohlin's formula, but which is not realized algebraically (being
ruled out by Thom's $\chi\le k^2$). Our example is of degree 10,
and is simply the complex scheme associated to the Rohlin tree of
Fig.\,\ref{Garidi-mass-false:fig}f for $(x,y)=(19,15)$, hence is
representable as the complex scheme of
Fig.\,\ref{Garidi-mass-false:fig}\,g. In fact such schemes already
exist in degree 6 (cf. optionally
Theorem~\ref{no-chance-to-reduce-Gudkov-to-Rohlin:thm}), where it
is just a matter of solving Rohlin's equation  for the two
$M$-schemes of degree 6 which are prohibited by Gudkov's
hypothesis.

Of course the example proposed probably belongs to a larger list
of such French schemes. More about this soon. On the other hand it
could be nice to know if there is a French scheme in degree 8
already.

{\it Degree $m=8$}.---Then the Harnack bound is
$M=g+1=\frac{7\cdot 6}{2}+1=22$. Applying the same method, we
start with the $M$-scheme $(1,1,20)$ with $\chi=20$. This has to
be adjusted to the Gudkov hypothesis $\chi\equiv_8 k^2=16$, but
then Thom's inequality $\chi\le k^2$ is satisfied, except if we
could move up to $\chi=24$ but this violates the basic estimate
$\chi\le M$. [Proof: $\chi=p-n\le p+n=r\le M$ by Harnack.] So it
seems that there is no French in degree 8, but a more systematic
study is required. As a loose evidence for the absence of French
scheme in degree 8, we note that the RKM-congruence for
$(M-2)$-schemes $\chi\equiv k^2+4 \pmod 8$ (ensuring type~I)
forces under Thom's inequality ($16\le \chi \le M=22$ and
Harnack's bound) to have $\chi=20$. But then $p-n=20$ and
$p+n=r=20$, so that $2p=40$, hence $p=20$ and $n=0$. So our scheme
is forced to be unnested and is $20$, which is prohibited by
Rohlin's formula.

\begin{figure}[h]
\centering
\epsfig{figure=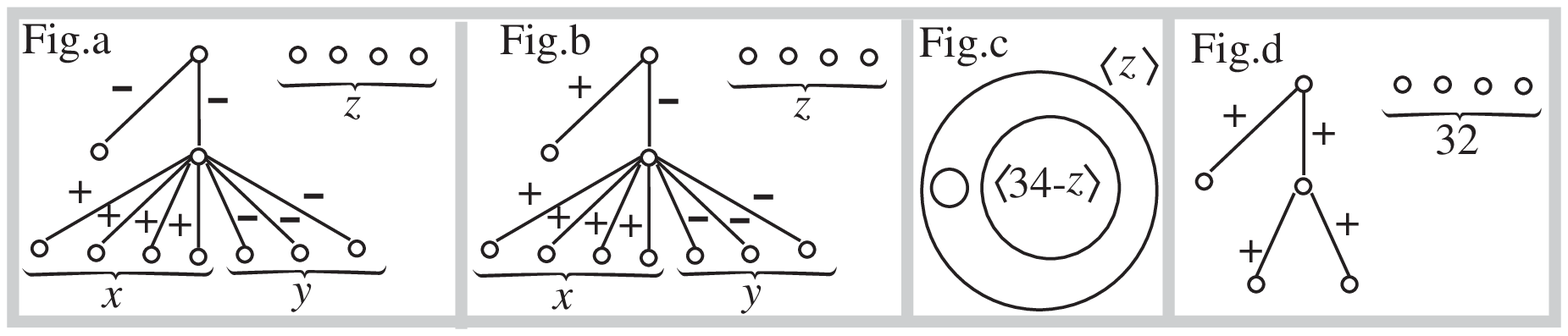,width=122mm} \vskip-5pt\penalty0
  \caption{\label{Fiedler3:fig}%
  More French schemes of degree 10} \vskip-5pt\penalty0
\end{figure}

{\it More French schemes in degree $10$}.---By the above method we
now proceed to find more French schemes. The idea is merely that
starting from the scheme on Fig.\,d we may move innermost ovals at
depth 2 outside at depth $0$ without changing $\chi=33$ (as forced
by the Gudkov hypothesis). So we consider a scheme whose tree is
like Fig.\,\ref{Fiedler3:fig}a where there are $z$ outer ovals
which are empty. We introduce $x$ and $y$ many free signs plus and
minus resp. as on Fig.\ref{Fiedler3:fig}\,a, with both top signs
negative. To get an $M$-scheme we impose $x+y+z=34$. As before we
seek to solve the Rohlin's equation. Recall that we abridge
$\pi:=\Pi^+$, $\eta:=\Pi^-$. By Rohlin's formula
$2(\pi-\eta)=r-k^2=37-25=12$ so $\pi-\eta=6$. But on the other
hand by the signs-law, we have (by looking at Fig.\,a)
$$
\pi=2x \qquad \textrm{and} \qquad \eta=2+2y.
$$
If $z$ is given (a priori in the range $0\le z \le 34$, but we
shall soon see that some more constraint are required), we solve
in $x,y$ as to satisfy Rohlin's formula. This gives $\pi-\eta=6$
and $\pi+\eta=2(1+x+y)$, so adding $2\pi=8+2(x+y)$, hence
$\pi=4+x+y=38-z$. So $\eta=32-z$. Finally we find
$x=\frac{\pi}{2}=19-\frac{z}{2}$. This requires so the assumption
$z$ even ($\bigstar$!!!), which looks anomalous but more about
this soon. And finally $y=\frac{\eta-2}{2}=\frac{30-z}{2}$, so
that we must assume $z\le 30$. ($\bigstar$ HYPOTHESIS to add!)

Now when $z$ is odd we proceed similarly, but the trick is to
change one of the top sign as on Fig.\,b into a plus. We still
have $x+y+z=34$, but now by the signs-law applied to the new
diagram (Fig.\,b):
$$
\pi=1+2x \qquad \textrm{and} \qquad \eta=1+2y.
$$
By Rohlin's formula we still have $\pi-\eta=6$, and now
$\pi+\eta=2(1+x+y)$ (actually like above!) and so repeating the
above $\pi=38-z$, $\eta=32-z$. Solving gives
$x=\frac{\pi-1}{2}=\frac{37-z}{2}$ and
$y=\frac{\eta-1}{2}=\frac{31-z}{2}$. So we assume $z\le 31$
($\bigstar$ Hypothesis!).

All this should prove the admittedly insignificant following
result:

\begin{theorem}\label{French-scheme-corrected-bis:thm}
For any integer $0\le z \le 31$, the $M$-scheme $(1,1(1,34-z))z$
(compare Fig.\,\ref{Fiedler3:fig}c) is a French scheme, i.e. it is
prohibited by Thom but not by Gudkov's hypothesis nor by Rohlin's
formula.
\end{theorem}

From the proof $z=31$ is sharp with this property, but it is
perhaps tranquilizing to check this more experimentally. Then the
tree reduces to Fig.\,d, and the total number of pair is
$\Pi=4+2=6$, but as $\pi-\eta=6$ (by Rohlin's formula) we are
forced to have only positive pairs, but this is impossible by the
signs-law (since in Rohlin's arithmetics $+$ times $+$ is minus!)

What to do next?
A
naive game would be to classify all French scheme in degree 10.
Another more serious problem would be to detect some universal
rule as to understand better the prohibitions given by Rohlin's
formula, or the lack thereof when the Hilbert tree can be given a
sign distribution so that Rohlin's formula is satisfied. All this
looks a bit unappealing combinatorics, yet the lack of
conceptualization in our above account surely ask for  a better
understanding of the Rohlin tree. One would like to understand all
Rohlin trees of dividing curves. Any such must satisfy the Rohlin
formula, but as we saw this is not the sole obstruction (as
sometimes Thom imposes additional constraints).

As a more specific goal we could try to find a French
$(M-2)$-scheme in degree 10. As noted earlier this is impossible
to do on the ``planar'' face of the Gudkov pyramid (cf. our
Fig.\,\ref{Degree10:fig}). There Rohlin's formula is always as
strong as Thom. (This follows also from the truth of the Garidi
mass principle for such simple schemes. We leave as a loose-end
exercise to exhibit classes of schemes for which the mass
conjecture holds true, albeit disproved in general for
``bat\^onnet'' like schemes, cf.
Fig.\,\ref{Garidi-mass-false:fig}.)

We now proceed to find $(M-2)$-schemes of degree 10 where Thom is
stronger than Rohlin's formula. The method is similar as above,
but we repeat the detail by unifying somewhat the proof. The
bat\^onnet structure of schemes violating the mass conjecture
suggests looking at the scheme $(1,1,33)$ (cf.
Fig.\,\ref{Fiedler4:fig}a). Here $\chi=33$, but in order to apply
Thom we have to ensure type~I, and the best known recipe to do
this is to adjust to the RKM-congruence $\chi\equiv k^2+4=25+4=29
\pmod 8$. So as $\chi=p-n\le p+n=r=35$, we cannot move up to 37,
but instead lower down $\chi=33$ to $29$. This is achieved by
delocalizing two ovals at depth 2 toward ovals at depth 1, and we
get the scheme $(1,2(1,31))$ (cf. Fig.\,\ref{Fiedler4:fig}b).

\begin{lemma}
The scheme $(1,2(1,31))$ and more generally its companions
$(1,2(1,31-z))z$ (cf. Fig.\,\ref{Fiedler4:fig}c) are by RKM of
type~I, but prohibited by Thom $\chi\le k^2=25$. However, provided
$z\le 29$, all these schemes are not prohibited by Rohlin's
formula.
\end{lemma}

\begin{proof}
As above we prove that there is a distribution of signs on the
Hilbert tree of the scheme compatible with the signs-law and with
Rohlin's formula. For this we consider the diagram of
Fig.\,\ref{Fiedler4:fig}d, where we have free parameters $x,y$
counting the number of positive resp. negative edges at depth 2.
We introduce also $\epsilon, \delta$ counting positive resp.
negative signs on the only 2 available  edges at depth 1 (lacking
prolongation). So $\epsilon+\delta=2$. The edge at depth 1
prolonging to depth 2 is given the sign $-1$. Further we have $z$
isolated vertices at depth 0 corresponding to the outer empty
ovals of the scheme of Fig.\,c. We have thus $x+y+z=31$. .

\begin{figure}[h]
\centering
\epsfig{figure=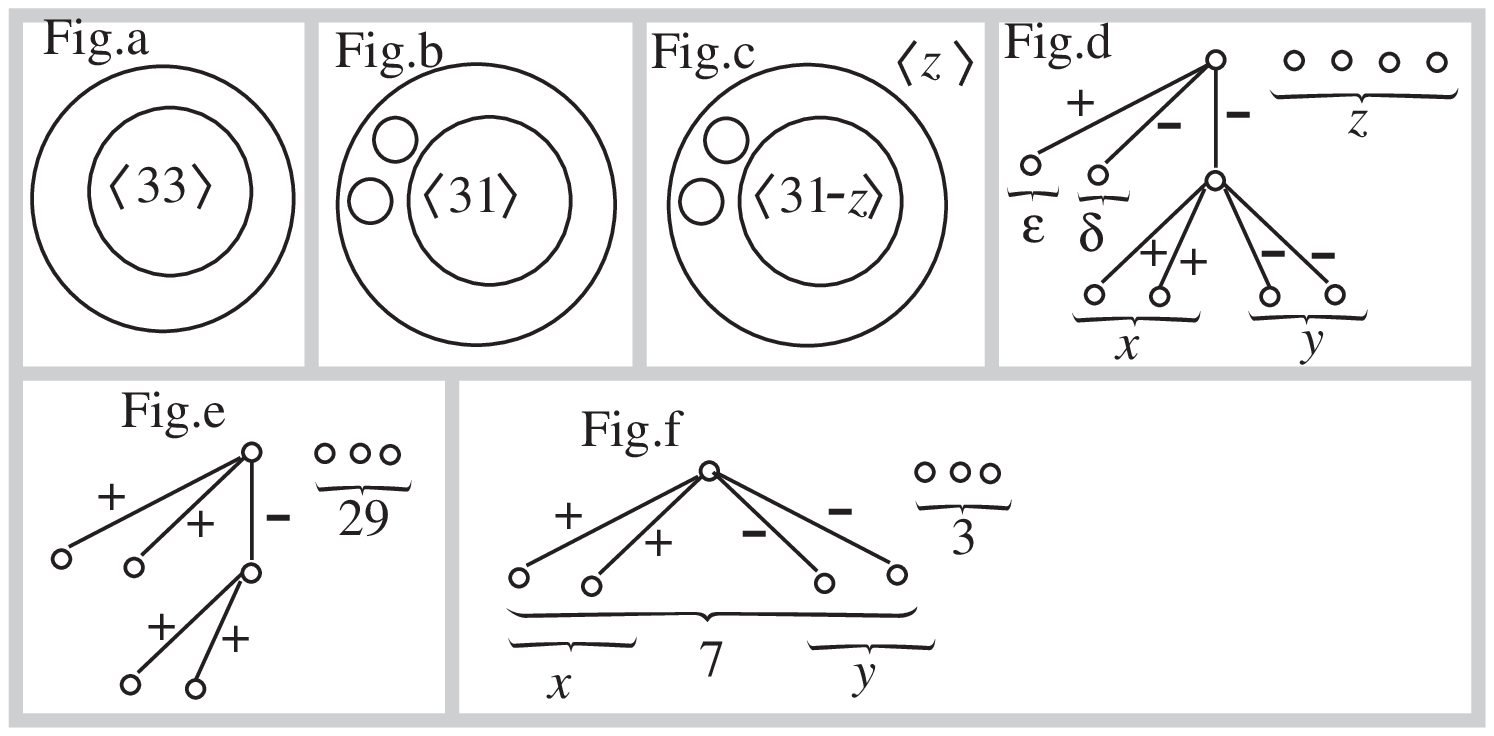,width=122mm} \vskip-5pt\penalty0
  \caption{\label{Fiedler4:fig}%
  French $(M-2)$-schemes of degree 10} \vskip-5pt\penalty0
\end{figure}

By Rohlin's formula $2(\pi-\eta)=r-k^2=35-25=10$, we find
$\pi-\eta=5$. By applying the signs-law to Fig.\,d, we find
$$
\pi=2x+\epsilon \qquad \textrm{and} \qquad \eta=2y+\delta+1.
$$
So $\pi+\eta=2(x+y)+3$, and thus $2\pi=8+2(x+y)$, whence
$\pi=4+(x+y)=35-z$ and $\eta=\pi-5=30-z$. So we solve
$$
x=\frac{\pi-\epsilon}{2}=\frac{35-z-\epsilon}{2} \quad
\textrm{and} \quad
y=\frac{\eta-\delta-1}{2}=\frac{29-z-\delta}{2}.
$$
It remains just to ensure integrality by choosing appropriately
the free parameters. Specifically,

$\bullet$ if $z$ is odd, then for $x$ to be integral choose
$\epsilon$ even (and then $\delta$ is even and so $y$ as well, but
sorry this parenthetical stuff is automatic!);

$\bullet$ if $z$ is even, choose $\epsilon$ odd.

Note that $z\le 29$ (to ensure $y\ge 0$), while for $z=29$ the
system of equations is soluble for $\delta=0$, $y=0$ and
$\epsilon=2$, whence $x=\frac{35-29-2}{2}=2$. (Optionally, one can
check diagrammatically that Rohlin's formula is verified on this
extremal example, cf. Fig.\,e. where we have $\pi=4+2=6$ and
$\eta=1$, so $\pi-\eta=5$ as it should by Rohlin.)
\end{proof}

What is the moral of all these messy calculations with this exotic
signs-law of Rohlin? Is there some respectable way to extract a
general result. A vague moral is that it is always quite boring to
apply Rohlin's formula and it is quite easy to make mistakes, as
was pointed out by Fiedler. Perhaps some higher intelligence than
me is able to discern some order in this chaos of Rohlin trees and
proves valuable corollaries of Rohlin's formula. However it is
also clear that the latter has some
limitation in failing to prohibit schemes
ruled out by Thom.

In the same spirit we wondered (earlier in this text,
Sec.\,\ref{Gudkov-hyp-via-Rohlin's-formula:sec}) if Rohlin's
formula could imply the Gudkov hypothesis. Perhaps a variant of
the method used right above can detect a counterexample, i.e. a
scheme prohibited by Gudkov but not by Rohlin. Actually we have
candidates in degree 6 already, namely the schemes $\frac{7}{1}3$
and its mirror $\frac{3}{1}7$. Taking the first its diagram-tree
is like on Fig.\,\ref{Fiedler4:fig}f, where we have $x$ many
positive edges and $y$ many negative ones. Of course $x+y=7$.
Since here we have no adjacent edges to compose, and we simply
have $\pi=x$ and $\pi=y$ and Rohlin's equation
$2(\pi-\eta)=r-k^2=11-9=2$ is trivial to solve. Indeed, $x-y=1$
and $x+y=7$ gives $2x=8$, whence $x=4$ and $y=3$. This proves the:

\begin{theorem}\label{no-chance-to-reduce-Gudkov-to-Rohlin:thm}
There is no chance to reduce Gudkov's hypothesis to Rohlin's
formula, unless one is able to put more stringent restriction on
the complex orientations (via geometric procedures).
\end{theorem}

$\bigstar$---{\it Old stuff, pre-Fiedler's correction, hence to be
read with
discernment.}---Several questions arise:

1.---On working out more carefully the combinatorics, try to
understand if a French scheme already exists in degree 8.

2.---It seems that as the degree increases French schemes will be
more and more frequent and so the r\^ole of Thom increases more
and more and becomes a valuable complement to the Russian
congruences, and Rohlin's formula.

3.---Glancing at the Gudkov table in degree 10
(Fig.\,\ref{Degree10:fig}), we see that the famous Gudkov broken
line is subjected to a severe deformation. More precisely, the
$M$-schemes $\frac{2}{1}34$ is prohibited by either Rohlin's
formula or by Thom [or even by Petrovskii's inequality
(\ref{Petrovskii's-inequalities:thm})], and so is the
$(M-2)$-scheme $\frac{3}{1}31$ for the same reasons plus (the full
punch of) the RKM-congruence
(\ref{RKM-congruence-reformulated:thm}). The diagrammatic
consequence is a distortion of the Gudkov front-line which is
penetrated from above to become the lilac line on
Fig.\,\ref{Degree10:fig}.

Hence directly below $\frac{3}{1}31$ we have 2 schemes (whose
symbols are not even worth writing down), and which according to
the diagrammatic of Fig.\,\ref{Degree10:fig} are not maximal
schemes. In contrast below the scheme $\frac{2}{1}34$ we have the
two schemes $\frac{2}{1}33=:S_1$ and $\frac{1}{1}34=:S_2$. Those
guys are interesting birds because they are maximal (at least in
the planar model of the Hilbert-Gudkov pyramid, as follows from
Rohlin's formula), yet they are certainly not of type~I (by
Klein's congruence $r\equiv_2 g+1$ of 1876). So there is some
little opportunity here to corrupt Rohlin's maximality conjecture
(RMC). Of course as all what we are using here was well-known to
Rohlin, our expectation will certainly quickly turn to
disillusion.

{\it Insertion.} [19.03.13]---Apart from the fact that (as
discussed in the sequel) it is hard to prohibit all extensions of
any one of those schemes $S_1,S_2$, a sharper look at the Russian
architecture of Fig.\,\ref{Degree10:fig} shows that both those
schemes are simply ruled out by Petrovskii's inequality
(\ref{Petrovskii's-inequalities:thm}) telling us that $\chi\le
\frac{3}{2}k(k-1)+1=\frac{3}{2}5(5-1)+1=3\cdot 5\cdot 2+1=31$.

{\footnotesize A first question is whether those two
$(M-1)$-schemes $S_1,S_2$ are realized\footnote{[19.03.13] The
answer is no and follows from Petrovskii's inequalities.}. If so,
then it remains to check that they are maximal. So we have to list
all their possible enlargements (extensions).

Let us do this for $S_2=\frac{1}{1}34$ (which looks more appealing
as it is ``less nested''), we have (up to isotopy) 4 possible
extensions depending upon the additional oval is added:

(1) outside the nest, which option leads to $\frac{1}{1}35$, which
is prohibited by the Gudkov-Rohlin congruence;

(2) inside the nest at depth 1, which option leads to
$\frac{2}{1}34$, which is prohibited by Rohlin's formula (or
Thom);

(3) inside the nest at depth 2, which option leads to $(1,1,1)34$.
Then $\chi=(1-1+1)+34=35\neq k^2=25 \pmod 8 $ so that the scheme
is prohibited by the Gudkov-Rohlin congruence, while one can even
note that its forerunner namely the Arnold congruence mod 4
suffices;

(4) inside an outer oval, which option leads to
$\frac{1}{1}\frac{1}{1}33$. Then $\chi=(1-1)+(1-1)+33=k^2=25\pmod
8$ so that Gudkov-Rohlin is not violated. However Thom's
inequation is violated, and so this scheme is not realized. As to
Rohlin's formula we have $\Pi=2$, and so
$2(\Pi^+-\Pi^-)=r-k^2=37-25=12$ cannot be verified as
$(\Pi^+-\Pi^-)\le(\Pi^++\Pi^-)=\Pi$. Intuitively we see that when
$\chi$ is large Rohlin's formula forces some nesting.

At this stage we have clearly explored all possible extensions of
$S_2$ and prohibited all of them via the classical congruence
(Gudkov-Arnold-Rohlin) or the Rohlin formula.

Yet even that is wrong as the additional oval needs not to be a
small one injected as above, but it can also surround other ovals.
Hence we have at least the following two families of schemes
enlarging $S_2$:

(5) $(1,\frac{1}{1}x)y$ where $x+y=34$, and

(6) $\frac{1}{1} \frac{x}{1} y $, where $x+y=34$.

This can of course be much diminished by the Gudkov-Rohlin
congruence, as follows. For (5) and $x=0$, we have the primitive
scheme $(1,1,1)34$ with $\chi=(1-1+1)+34=35$ which 2 unit above
$33=25 \pmod 8$. Hence by trading outer ``$y$'' ovals against
``$x$'' ovals at depth 1 we decrease $\chi$ by 2, so that the
first scheme with correct $\chi$ is $(1,\frac{1}{1}1)33$, and then
the list extends by using 4-fold periodicity as
$$
(1,\frac{1}{1}1)33, \quad, (1,\frac{1}{1}5)29, etc.
$$
with resp. $\chi=33, 25, etc.$ (descent by 8). Hence the first is
still prohibited by Thom, yet the subsequent schemes are unlikely
to be prohibited by Thom nor by Rohlin's formula. At this stage
our naive project to corrupt Rohlin's maximality conjecture breaks
down. For instance the second listed Thom compatible scheme
$(1,\frac{1}{1}5)29$ has $\Pi=\binom{3}{2}+5=8$ and so Rohlin's
formula $2(\Pi^+-\Pi^-)=r-k^2=37-25=12$ is soluble for the pair
$\Pi^{\pm}=(7,1)$.

At this stage the only vestiges of expectance to foil RMC via our
naive strategy would be that all those extended schemes are
prohibited by a degree-10-extension of the Fiedler-Viro theorem
(cf. Theorem~\ref{Viro-Fiedler-prohibition:thm} in degree 8), but
that looks a dubious expectation.

Of course it is much more likely that either Harnack, Hilbert or
Viro's method of construction realize one of those schemes, and so
that our scheme $S_2$ (and likewise $S_1$) are not maximal.

}

\subsection{Some intuition behind Hilbert, Thom, Rohlin?}

[12.03.13] The theory of abstract real algebraic curves is
essentially a closed chapter of mathematics. More precisely since
Riemann 1857 \cite{Riemann_1857} paper on Abelian integrals, we
have a clear-cut vision of all complex curves and their moduli.
Avatars thereof in the bordered setting or even nonorientable
realm were worked out by Klein 1882 \cite{Klein_1882} (upon the
heritage of Gauss, M\"obius, Listing the discoverers of ``il
nastro di M\"obius''). Teichm\"uller 1939 \cite{Teichmueller_1939}
put the nearly final touch by approaching the moduli problem via
quasiconformal maps along the philosophy developed by Gr\"otzsch
1928, Lavrentieff 1929, Ahlfors ca. 1930--35. Notwithstanding our
thesis is that from the viewpoint of total reality (or what
amounts to the same branched covers of the disc), some
quantitative aspects have not yet attained their ultimate
perfection and sharpness, though it may be only a matter of
assimilating the heritage of previous generations. This is
discussed at length in the Introduction of this text, but briefly
we may recall that in the closed case the ultimate perfection
regarding Riemann surfaces expressed as branched cover of the
line(=Riemann sphere) is achieved in Riemann 1857, Brill-Noether
1874 \cite{Brill-Noether_1874}, and especially Meis 1960
\cite{Meis_1960} (which I had never the occasion to consult). When
it comes to bordered surfaces (equivalently orthosymmetric real
curves), the first steps belongs to Riemann 1857 (Nachlass)
\cite{Riemann_1857}, Schottky 1875--1877, Bieberbach 1925, Grunsky
1937, Ahlfors 1950, etc., up to perhaps Gabard 2004
\cite{Gabard_2004}, 2006 \cite{Gabard_2006}, whose result still
deserves better introspection, but whose sharpness is
adhered upon (perhaps too hastily) in several works (e.g.
Fraser-Schoen 2011 \cite{Fraser-Schoen_2011}, Coppens 2011
\cite{Coppens_2011}. Actually Coppens' result although adhering to
the truth of Gabard's  is so logically independent that its truth
may not be  jeopardized by a disproof of Gabard's bound $\gamma\le
r+p$.

Next we may move in the Plato cavern of plane curves and
contemplate the so-called {\it Hilbert's problem\/} on the
topology of real plane algebraic curves. Here most of the
difficulty is allied to a certain combinatorial mess arising from
the nested structures of ovals and their distribution. The whole
point is to look at a certain hierarchical structure (POSET)
arising from the inclusion between the insides of ovals. We call
it the {\it Hilbert tree}.

As early as 1891, Hilbert was the first to formulate the intuition
that this Hilbert tree has a certain verticality, in the sense
that an algebraic curve cannot reduce to an unnested collection of
ovals. He formulates this for $M$-curves of degree 6, and of
course the assertion holds only for maximal curves (cf. the Gudkov
table=Fig.\,\ref{Gudkov-Table3:fig} or for a specific example
Fig.\,\ref{KleinRo-sextic:fig}). The Hilbert tree of such an
unnested scheme (Gudkov's symbol $M$, where $M=g+1$ is Harnack's
bound) just reduce to a collection of vertices at depth 0, without
nesting (hence without edges).

At some stage Thom had the following intuition.

\smallskip
{\footnotesize {\it Historical notice (Thom 1982, Rudolph
1984)}.---[21.03.13]---Alas the whole designation ``Thom
conjecture'' looks poorly documented in print. In Donaldson 1989
\cite{Donaldson_1989},  we read a loose ``usually ascribed to
R.~Thom''. Maybe the oldest source where the term ``Thom
conjecture'' appears is Kirby's problem list (1970
\cite{Kirby_1970--95} regularly updated by its author). Of course
as early as 1961 Kervaire-Milnor 1961 \cite{Kervaire-Milnor_1961}
were able to tackle the first non-trivial case of degree 3
(building upon Rohlin's works) of the so-called ``Thom
conjecture'', but do not use this nomenclature. The best
clarification I could found is given by Lee Rudolph 1984
\cite{Rudolph_1984}, whose  footnote reads: ``Professor Thom has
remarked (personal communication, November 19, 1982) that the
conjecture perhaps more properly belongs to folklore.''

}

\begin{theorem} {\rm (Thom conjecture, Kronheimer-Mrowka theorem 1994, and
independently Morgan-Szab\'o-Taubes 1995/96
\cite{Morgan-Szabo-Taubes_1996})}.---The genus of a smooth
embedded oriented surface in $\CC P^2$ is at least as large than
that of ``the'' algebraic smooth curve realizing the same homology
class, alias the degree.
\end{theorem}

This prompts a wide extension of Hilbert's verticality principle
[or nesting Ansatz] that the Hilbert tree cannot be to flat. This
is alas a bit like in human feodal systems where a certain
concentration of power, slavery and subjection is observed.
Precisely this is given by Theorem~\ref{Thom-Ragsdale:thm} stating
that any dividing curve of degree $m=2k$ satisfies $\chi \le k^2$.
As the degree increases to 8, or 10, etc. this prohibits more and
more schemes on the right wing of the Gudkov pyramid (cf. e.g.
Fig.\,\ref{Degree10:fig}), yet to some extend those schemes can
also be more elementarily ruled out by Rohlin's formula.

Since the Harnack bound is
$M=g+1=\frac{(2k-1)(2k-2)}{2}+1=(2k-1)(k-1)+1=2k^2-3k+2$, the
asymptotic location as $m=2k\to \infty$ of  Thom's bound $k^2$
[ALAS INCORRECT, MY MISTAKE CORRECTED BY FIEDLER] is the half
value of Harnack's bound $M\approx 2k^2$, and so nearly one
quarter of the $M$-curves are prohibited by Thom (cf. e.g.
Fig.\,\ref{Degree10:fig}). However this figure only represents
schemes with one nonempty oval so that the real pyramid is much
less amputated than its planar sheet leads one to suspect.

On the other hand the Hilbert hierarchies cannot be too deep.
Indeed one cannot observe in degree $m$ a nested chain of ovals of
depth larger than $k=m/2$ as follows directly from B\'ezout for
lines. More generally using conics through 5 points there cannot
be schemes with 5 chains of total length larger that $2m$. This is
the well-known topics of Hilbert's bound on the depth of nests.

Further when those deepest configurations are attained then the
curve is necessarily of type~I via the phenomenon of total
reality.

Summarizing we see that the Hilbert tree cannot be too deep (as
follows from reality consideration and B\'ezout) nor can it be too
superficial (as follows from the complexification, e.g. Thom's
principle or the Petrovskii inequalities
\ref{Petrovskii's-inequalities:thm}). Further the total reality
phenomenon ought to play some big r\^ole as envisioned by Rohlin.
Of course all this needs to be made more precise (further
explored).

\subsection{On an arithmetical problem valorizing Thom
in the detriment of Rohlin (and sometimes Gudkov): yet another
numerical coincidence regarding Hilbert's 16th}

[13.04.13] This section (and its title!) is somewhat naive and
misleading (as I failed to exploit the full punch of Rohlin's
formula), but diverges to a lovely basic arithmetic problem (on
which I have little grasp, but must be very classical, surely
Gauss? or probably much older Diophante?). Hence it was not
censured but can safely be omitted.

[16.03.13] As we said, but repeat it once more, Hilbert arrived
ca. 1891 (vgl. Hilbert 1891
\cite{Hilbert_1891_U-die-rellen-Zuege}) at the intuition (at least
for sextics) that ovals of algebraic curves necessarily exhibit
some feodal structure of
nesting
impeding all the ovals  lying outside each others. Is this really
true in general? [OF COURSE, E.G., VIA ROHLIN'S FORMULA, AND FIRST
KNOWN TO PETROVSKII 1933/38.] Of course Hilbert posited this for
$M$-curves (of order $m\ge 6$) and nearly proved it only for $m=
6$, yet presumably we are allowed to extrapolate his thoughts. In
fact, the original text proceeds more carefully and reads as
follows (cf. Hilbert 1891, \loccit, Fussnote, S.\,418, in Ges.
Abhandl., Bd.\,II, Algebra, Invariantentheorie, Geometrie):

``Diesen Fall $n=6$ habe ich einer weiteren eingehenden
Untersuchung unterworfen, wobei ich\,---\,freilich auf einem
au{\ss}erordentlich umst\"andlichen Wege\,--- fand, da{\ss} die
elf Z\"uge einer Kurve 6-ter Ordnung keinesfalls s\"amtlich
au{\ss}erhalb und voneinander getrennt verlaufen k\"onnen. Dieses
Resultat erscheint mir deshalb von Interesse, weil er zeigt,
da{\ss} f\"ur Kurven mit  der Maximalzahl von Z\"ugen der
topologisch einfachste Fall nicht immer\footnote{As we shall soon
recall, since Thom/Kronheimer-Mrowka 1995 we may replace ``nicht
immer'' by ``never''!!! (provided $m\ge 6$) and $m$ even.
[17.03.13] In fact this was known much earlier since Petrovskii
1938, cf. Lemma~\ref{Hilbert's-nesting-intuition:lem} below. An
elementary proof also follows from Rohlin's formula as reminded in
the same lemma. It would be interesting to say more on the odd
degree case (again Petrovskii 1938 \cite{Petrowsky_1938} should
suffice to corroboration ``Hilbert thesis'' for $m\ge 7$).}
m\"oglich ist.''

Hilbert addressed this again, and more generally the question of
elucidating the isotopic classification of sextic curves (or even
quartics surfaces), into his well-known 16th problem at the Paris
Congress of 1900. Circa 7 decades came the deep semi-experimental
work of Academician Dimitrii Andreevich Gudkov solving Hilbert's
isotopic problem for sextics curves, though according to Arnold
(e.g. 1997/00
\cite{Arnold_1997/200X-Symplectization-Complexification}) this
left the supervising teacher Petrovskii quite dubitative, not to
say skeptical). Later V.\,A. Rohlin made a series of jokes by
noticing that his 1952 theorem on the divisibility by 16 of the
signature (coined by Weyl 1923, in Analisis Situs Combinatorio
written in Spanish with the assistance of his linguist wife) of a
spinorial smooth $4$-manifold turns out to imply the Gudkov
hypothesis $\chi\equiv k^2 \pmod 8$ for $M$-curves (extrapolating
widely the phenomenology observed in degree $6$), where $\chi$
denotes as usual the Euler characteristic of the ``Ragsdale''
orientable membrane bounding the ovals from ``inside''.

In the same
elan,  Rohlin 1974 (and 1978) wrote down the {\it Rohlin
formula\/} $2(\Pi^+-\Pi^-)=r-k^2$  (cf. \ref{Rohlin-formula:thm})
which is so fundamental that it seems worth (to save ink) to
abridge notation by letting $\pi:=\Pi^+$ and $\eta:=\Pi^-$. When
the curve has no nesting this formula implies that $r=k^2$ is a
square (whenever the curve is dividing, as it is automatically the
case for $M$-curves).

{\it Warning.}---[17.03.13]  The sequel is much ill-posed as I
(very shamefully) missed to notice that $r$ is not any square, but
that of the semi-degree $k$. Yet since it seems to involve
pleasant arithmetics, hence we did not censured it!

Restricting attention to $M$-curves of even degree $m=2k$, we are
therefore invited to study the following arithmetical problem as a
way to corrupt Hilbert's feudalistic intuition (forced presence of
nesting for $M$-curves of high-degrees $m\ge 6$).

\begin{prob} {\rm (Quadrature of the Harnack bound).}---Given any integer $k\ge$, set $m=2k$ (interpreted as the order of
the curve) and let
$M=g+1=\frac{(2k-1)(2k-2)}{2}+1=(2k-1)(k-1)+1=2k^2-3k+2$ be the
corresponding Harnack bound. For which values of $k$ is  $M$
predestined to be a square, so that Rohlin does not prohibit the
unnested $M$-scheme (THIS IS FALSE), and additionally try to
arrange the Gudkov congruence $M=\chi\equiv k^2 \pmod 8$.
\end{prob}

Assume that there is such an integer $k$ then the unnested
$M$-scheme (Gudkov symbol $M$ also!) is not prohibited by Gudkov
nor by Rohlin and so constitutes a potential violation of
Hilbert's principle. To kill the suspense, as far as I know this
scenario could not have been precluded until the Kronheimer-Mrowka
validation of Thom's conjecture. Recall (from
\ref{Thom-Ragsdale:thm}) that Thom implies $\chi\le k^2$ and so
$\chi$ cannot be as large as $M$ which is asymptotically twice so
big (except of course for low degrees $m\le 4$).

So curiously I would say (personal feeling probably foiled due to
a lack of Russian wisdoms) that before Thom-Kronheimer-Mrowka 1994
\cite{Kronheimer-Mrowka_1994} Hilbert's intuition could have been
completely wrong depending upon a resolution of the above
arithmetical problem.

Of course I arrived at the problem by  modest acquaintance with
the geometry of Gudkov's pyramid which look basically as follows
(Fig.\,\ref{Pyramid:fig}):

\begin{figure}[h]
\centering
\epsfig{figure=Pyramid.eps,width=122mm} \vskip-5pt\penalty0
  \caption{\label{Pyramid:fig}%
  Obstruction upon curves of type~I: Thom versus Rohlin seen in
  simplified 2D-diagrammatic not reflecting the full combinatorial
  mess allied to Hilbert's 16th, yet sufficient for our purpose} \vskip-5pt\penalty0
\end{figure}

Now being as bad in arithmetics than in combinatorics, we propose
to tackle the arithmetical problem naively. We want to know first
of all  when Harnack's bound is a square. We know that Harnack's
bound evolves like the genus of plane curves as a quadratic
function (Gauss freshman calculus at 5 years old) namely
$g=g(m)=1+2+3+\dots+(m-2)=\frac{(m-1)(m-2)}{2}$. Here we
concentrate upon even degrees $m=2k$ and we may tabulate the
values of $M$ by adding a linear progression of increment 4 (cf.
Fig.\,\ref{Pyramid:fig} right-side table). Then we compare this
with the list of all squares $1,4,9, 16\dots$. Apart from the
trivial values $k=1,2$, the Harnack bound is never a square until
we reach the value $k=17$ where $M=529$ is the square of $23$.
Since $17=16+1$ this is yet another numerological coincidence in
Hilbert's 16th (yet surely not as deep of those of Gudkov-Rohlin).
In view of this I naively expected that the next quadrature of
Harnack's bound occurs at $k=2\cdot 16+1=33$. But a simple
calculation shows this to be a wrong intuition. What about $k=4
\cdot 16+1=65$, also not good! At this stage we stopped  guessing
and continued the tabulation by hand, and found the next square
Harnack bound at $k=46$, namely $M=\frac{91\cdot 90}{2}+1=91\cdot
45+1=4096$ which is $64^2$.

It is now time looking at the Gudkov congruence. For $k=17$, the
Harnack bound $M=529$ (which is also $\chi$ the characteristic of
the Ragsdale membrane as we suppose no nesting) reduces modulo 8
to $529\equiv_8 49 $ (after removing $480$) and then to $1$ (after
removing 48). On the other hand $k^2\equiv_8 17^2\equiv_8 1^2
\equiv_8 1$. Hence Gudkov's congruence is satisfied. So the
$M$-scheme (with Gudkov symbol $529$) could be a potential
counter-example to Hilbert's intuition, if Thom would not
salvage
it!

For the next value $k=46$, $\chi=M=4096\equiv_8 96 \equiv_8 16
\equiv_8 0$, whereas $k^2=(46)^2\equiv_8 (-2)^2=4$. So Gudkov
hypothesis is not verified and here Hilbert's intuition is already
vindicated by Gudkov (and Rohlin who proved it). Notice that
Arnold's congruence would be not enough not rescue Hilbert.

This is essentially all what we have to say on this problem.
Though the arithmetical problem looks attractive, I do not know
how to solve it in general. Of course it looks an easy Diophantine
equation
$$
x^2=2k^2-3k+2,
$$
of degree 2 but this does not boils down to  studying the rational
points on a conic (via the usual method of sweeping like for
Pythagorean triplets), as we are here really concerned with
integral solutions of an equation laking homogeneity. Of course I
presume it is a well-known and easy problem of arithmetics, yet
actually its solution has little impact upon Hilbert's sixteenth
problem, since the corresponding (unnested) schemes are defacto
prohibited by Thom. Still the particular solution exhibited above
$(k,x)=(17, 23)$ offers another instance of French scheme where
Thom is stronger that the conjunction of Gudkov and Rohlin. I do
not if prior to Kronheimer-Mrowka one was able to prohibit the
existence of the corresponding $M$-scheme of degree $34$ with
$M=529$ unnested ovals.

[17.03.13] In fact the answer is a trivial consequence of Rohlin's
formula, since $0=2(\pi-\eta)=r-k^2$, but $r=529$ is not the
square of $17$ (as $17^2=289$). More generally, all the
problematic of this section is spoiled since I omitted at the
beginning of the discussion to notice that Rohlin's formula does
not merely implies in the unnested case that $r$ is the square,
but is indeed the square of the semi-degree $k$. In fact as shown,
e.g. by Lemma~\ref{Rohlin-consequence-for-M-curves:lem} it is a
trivial consequence of Rohlin's formula that Hilbert's nesting
intuition occurs for $m\ge 6$. The simplest way to argue is as
follows.

\begin{lemma}\label{Hilbert's-nesting-intuition:lem}
{\rm (Hilbert's nesting intuition,
validated via Rohlin 1974 or via Petrovskii 1938)} Any $M$-curve
of even degree $m=2k$ with $m\ge 6$ cannot be unnested, i.e. have
all its ovals outside another. Actually, a dividing unnested curve
is forced to have $r=k^2$ ovals.
\end{lemma}

\begin{proof} A simple way to argue is via
Rohlin's formula 1974 (\ref{Rohlin-formula:thm}). We
have $0=2(\pi-\eta)=r-k^2$, hence $r=k^2$. But by the $M$-curve
assumption $r=M=g+1=(2k-1)(k-1)+1=2k^2-3k+2$, which is strictly
larger than $r=k^2$ as soon as $k\ge 3$.

The first clause also follows by specializing Petrovskii's
upper-bound $\chi\le \frac{3}{2}k(k-1)+1=:P$ valid for all curves.
It is plain indeed to check that  $P<M$ (``Petrovskii is sharper
than Harnack'') for $k\ge 3$. Indeed the difference
$M-P=(2k^2-3k+2)-(\frac{3}{2}k(k-1)+1)=\frac{1}{2}k^2
-\frac{3}{2}k+\frac{1}{2}=\frac{1}{2}\underbrace{(k^2-3k+1)}$. The
roots of the underbraced quadratics are $k_{1,2}=\frac{3\pm
\sqrt{9-4}}{2}=\frac{3\pm \sqrt{5}}{2}\le \frac{3+
\sqrt{9}}{2}=3$. Variant: $\sqrt{5}\approx 2.24$, so $k_1\approx
2.62$.
\end{proof}

[17.03.13] Misha Gabard (my father) and his skills in Excel
calculated the next values of $k$ making $M=2k^2-3k+2$ into a
square. The complete list for $M\le 1'000'000=10^6$ is
$$
k=1,2,17,46,553,1538,18761,52222,637297,
$$
where $M$ is resp. the square of
$$
x=1,2,23,64,781,2174,26531,73852,901273.
$$
On looking at the successive ratio of $k$, one finds the list
$$
R:=\frac{k_{n+1}}{k_n}=2, 8.5, 2.71, 12.02, 2.78, 12.20, 2.7835,
12.2036.
$$
So naively the progression oscillates between a factor of about
2.78 and one about 12.20. So there seems to some extreme
regularity, and one can predict in advance the size of the
solutions.

\subsection{The long quest of a Caucasian scheme,
i.e. where Rohlin is stronger than Thom and Gudkov (united)}

[12.03.13] The following section has been written in realm time
out of paper notes, and we have reproduced our long march toward
the trivial truth. The pressed reader can directly jump to
Theorem~\ref{Caucasian-scheme:thm}. (WHICH IS ACTUALLY FALSE) So
one must really read the first few sections which collects loosely
organized thoughts on the problem of finding a scheme prohibited
by Rohlin but not by Thom. Of course to avoid trivialities one
must add more assumptions like the curve being an $M$-curves
satisfying the Arnold or even Gudkov congruence.

[12.03.13] As a consequence of Thom's estimate $\chi\le k^2$
(\ref{Thom-Ragsdale:thm}), it seems that the Gudkov pyramids are
much amputated on their right wings causing thereby a certain
asymmetry of them. It is at this stage of some interest to wonder
what happens on the opposite left-wings where Thom becomes useless
but perhaps Rohlin's formula has still some prohibitive impact.

In view of Fig.\,\ref{Degree10:fig} let us look at the scheme
$\frac{34}{1}2$  (not prohibited by Gudkov's hypothesis). Here
$\chi=1-34+2=-31$ and we have no prohibition by Thom. Is this
scheme prohibited by Rohlin's formula? Certainly not because as
the Rohlin tree is simple there is no signs-law and the Rohlin
equation $2(\pi-\eta)=r-k^2=37-25=12$ is soluble under the obvious
additional relation $\pi+\eta=34$. Indeed $2\pi=40$, so $\pi=20$
and $\eta=14$. (Optional question: Is this scheme realized
algebraically?) ([21.03.13] Answer: negative as follows from
Petrovskii 1933/38 (\ref{Petrovskii's-inequalities:thm})!)

\begin{figure}[h]
\centering
\epsfig{figure=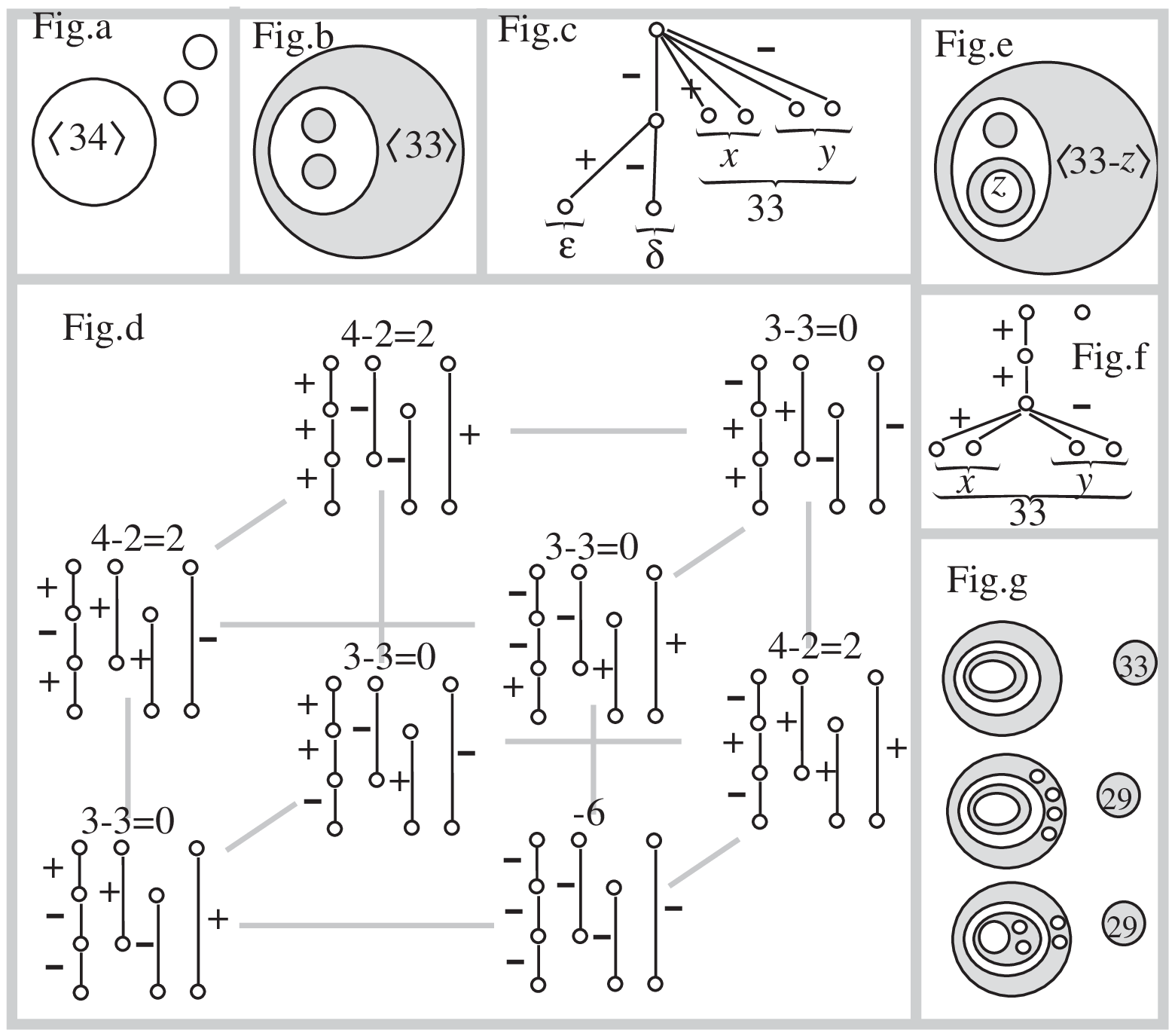,width=122mm}
\vskip-5pt\penalty0
  \caption{\label{Signs-law-triad:fig}%
  Seeking a Caucasian scheme of degree 10} \vskip-5pt\penalty0
\end{figure}

Without changing $\chi$ we may delocalize the 2 outer islands to
make them islands in a lake. We find so the scheme $(1,(1,2)33)$
(cf. Fig.\,\ref{Signs-law-triad:fig}b). Is this schemes prohibited
by Rohlin? We think the answer is no and the proof proceeds along
the usual algorithm of solving the Rohlin equation under the
signs-law. Recall that the latter can be easily remembered by
saying that consanguinity is bad, i.e. $+\times +=-$, $-\times
-=+$, while mixing the genes is good, i.e. $+\times -= +$ and
$-\times += -$. This exotic signs law is the exact opposite of the
usual convention.

Fig.\,c depicts the Hilbert tree of the scheme of Fig.\,b, and we
decorate it with a sign distribution as depicted, i.e. with
$x$-many plus, and $y$-many minus, and likewise $\epsilon$ plus
and $\delta$ minus at the indicated place. We have thus $x+y=33$,
and $\epsilon+\delta=2$. Using the signs-law we find for the
number of positive $\pi:=\Pi^+$ resp. negative pair $\eta:=\Pi^-$
the following expressions:
$$
\pi=x+2\epsilon, \quad \textrm{and} \quad \eta =y+2\delta+1.
$$
Adding gives $\pi+\eta=(x+y)+4+1=38$. By Rohlin's formula
$2(\pi-\eta)=r-k^2=37-25=12$, so $\pi-\eta=6$. Hence $2\pi=44$, so
$\pi=22$ and $\eta=\pi-6=16$. Solving finally in $x$ gives
$x=\pi-2\epsilon=22-2\epsilon$. We are free to choose say
$\epsilon=2$, then $x=22-4=18$ and $y=15$. So Rohlin's equation is
soluble.

At this stage the question becomes under which condition is
Rohlin's equation not soluble? A priori it may be remarked that
Rohlin's formula $2(\pi-\eta)=r-k^2$ is coupled with the formula
$\pi+\eta=\Pi$ the total number of pair. Both right-hand sides of
this pair of equations are entirely determined by the real scheme
without having to worry about orientations. So under reasonable
hypothesis, like $r-k^2$ even which is a disguised (Russian)
version of Klein's congruence $r\equiv g+1 \pmod 2$---as follows
from the boring calculation
$$
r\equiv_2 g+1=(2k-1)(k-1)+1=2k^2-3k+2\equiv_2 3k \equiv_2 k
\equiv_2 k^2,
$$
---we could expect that Rohlin's equation is always soluble,
especially if the (concomitant) Arnold congruence is satisfied.

However this is not the case as shown by Fiedler's corrigendum.
There we considered the $M$-scheme $(1,1,1,1)33$ (of degree 10)
which is prohibited by Rohlin's formula since $\pi-\eta\le 2$ by
the signs-law for triad (compare
Fig.\,\ref{Signs-law-triad:fig}d). But on the other hand Rohlin's
formula forces $2(\pi-\eta)=r-k^2=37-25=12$, hence $\pi-\eta=6$.

Hence the philosophy behind Fiedler's corrigendum seems to be that
Rohlin's formula imposes restriction, when there is much
predestination forced by the signs-law, and this is naively
speaking the case when there is much nesting like on the example
just given. It seems of course to be of some interest to
generalize the estimate $\Delta \Pi=\pi-\eta\le 2$ for triads to
tetrads, etc.

So the scheme $(1,1,1,1)33$ is prohibited by Rohlin (in its strong
form of the signs-law), but it is also by Thom as $\chi=33
\nleqslant k^2=25$. Now however we would rather be interested in
schemes prohibited by Rohlin but not by Thom.

A first idea we had was to start from the scheme on Fig.\,b, and
move $z$-many ovals inside to get Fig.\,e. As to maximize the
``predestination'' it seems wise to move the empty oval outside,
and to assume $z=33$. But still under these circumstances it
turned out that the Rohlin equation is soluble for a suitable
distribution of sign, cf. Fig.\,f, from which we infer by the
signs-law the following (where underbraced is the length of the
corresponding edges=pairs)
$$
\pi-\eta=\underbrace{x-y+2}_{1}+\underbrace{(-x+y-1)}_{2}
+\underbrace{(+x-y)}_{3}=x-y+1.
$$
(In this calculation it is useful to remind the consanguinity law,
and the Fig.\,d for a triad.) Hence as $\pi-\eta=6$ (by Rohlin),
we $x-y=5$, and $x+y=33$. So $2x=38$, whence $x=19$, and $y=14$.
So Rohlin's equations is soluble.

Philosophically it seems to be that whenever we have the free
parameters $x,y$ then we can solve despite the predestination
forced by the signs law.

Keep in mind that our question is whether Rohlin implies
prohibitions on the left wing of the pyramid (i.e. for $\chi \ll
0$ much negative) where Thom tells nothing. In fact $\chi \le k^2$
is enough for Thom to be non-prohibitive. Hence our next idea was
to start from Fiedler's example $(1,1,1,1)33$ with $\chi=33$, and
lower down to $\chi=25$. This lowering may be achieved by trading
the (oceanic) outer islands against lakes, i.e. ovals at odd
depths, and this requires to be done four times. (This can be done
in several ways, cf. Fig.\,g.) However on applying the usual
algorithm of sign distributions, we were always able to solve the
Rohlin equation in a way compatible with the signs-law. Details on
p.\,AR91 of my hand-notes, but the philosophy seems to be
basically that Rohlin's formula gives one equation and the
signs-law another but as soon as there free-parameter $x,y$
counting the number of signs on branches there is enough freedom
to solve all equations consistently.

\begin{figure}[h]
\centering
\epsfig{figure=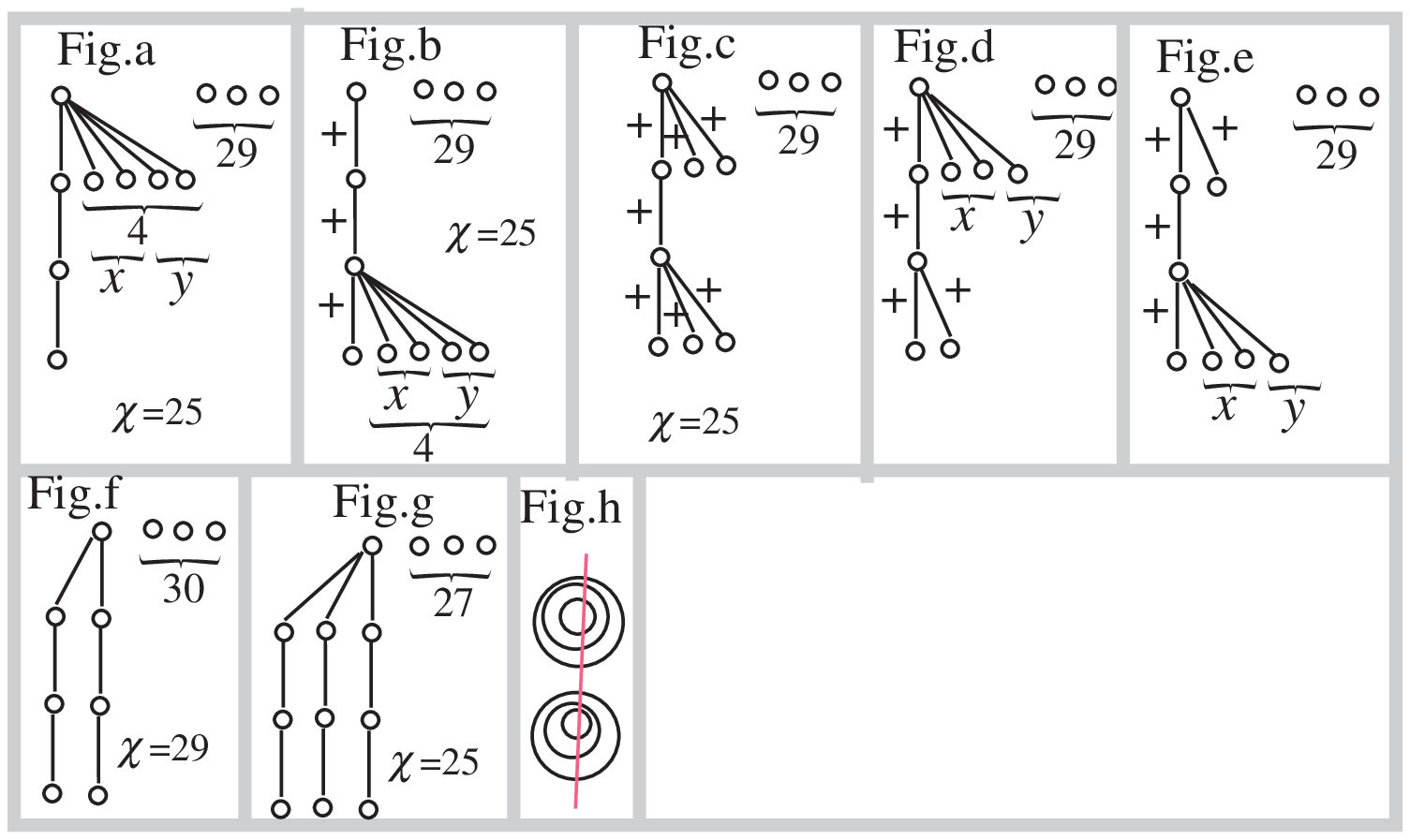,width=122mm} \vskip-5pt\penalty0
  \caption{\label{Fied9:fig}%
  Seeking a Caucasian scheme of degree 10} \vskip-5pt\penalty0
\end{figure}

In fact let us write down the argument. Starting from Fiedler's
example $(1,1,1,1)33$ with $\chi=33$, we lower down to $\chi=25$
by dragging 4 outer ovals at odd depth. This can be done in
several fashions as we said, for instance like on
Fig.\,\ref{Fied9:fig}a by putting the 4 ovals at depth 1. However
this new $M$-scheme is not prohibited by Rohlin's formula. The
latter says that $2(\pi-\eta)=r-k^2=37-25=12$, $\pi-\eta=6$. On
the other hand by using the signs-law, and splitting the 4 empty
ovals at depth 1 into $4=x+y$, where $x$ are positive, and $y$ are
negative (w.r.t. the sign of the edge right above it), and
assuming further for simplicity that the trunk is everywhere
positive, we can still solve the equation. Indeed by the signs-law
we have $\pi-\eta=2+x-y$, and so $x-y=4$, $x+y=4$ whence $x=4$,
$y=0$.

If instead we drag the 4 ovals at depth 3 we get
Fig.\,\ref{Fied9:fig}b. Now the signs-law becomes more involved,
but we find splitting according to the length of the pairs
(underbraced index) the following expression (assuming for
simplicity $+$-signs already fixed on the trunk and as above there
are $x$ many $+$ and $y$ many $-$):
$$
\pi-\eta=\underbrace{3+x-y}_{1}+ \underbrace{(-2-x+y)}_{2}+
\underbrace{(+1+x-y)}_{3}=2+x-y,
$$
where we used the signs-laws for triad
(Fig.\,\ref{Signs-law-triad:fig}d) As $\pi-\eta=6$, this gives
$x-y=4$, $x+y=4$, which is soluble.

Then we can also disperse the 4 ovals at different heights like on
Fig.\,\ref{Fied9:fig}c where 2 are at depth 1 and the 2 others at
depth 3. Instead of applying the method of indeterminate signs, we
content to give a solving sign distribution as depicted with only
$+$-signs. One checks that
$\pi-\eta=\underbrace{7}_{1}+\underbrace{(-4)}_{2}+
\underbrace{(+3)}_{3}=6$, in agreement with Rohlin's formula.

Finally it remains to analyze the case of Fig.\,d, where we
rigidify already some signs. Here $x+y=3$. Calculating via the
signs-law gives
$$
\pi-\eta=\underbrace{4+x-y}_{1}+\underbrace{(-3)}_{2}+
\underbrace{(+2)}_{3}=3+x-y,
$$
whence, as $\pi-\eta=6$, the system $x-y=3$, $x+y=3$, which is
soluble (integrally) as $x=3$, $y=0$. (One checks mentally that
this everywhere positive distribution works, as one do quickly
mistakes in such calculation!)

Last the case of Fig.\,\ref{Fied9:fig}e, is also handled by the
same method. Here
$$
\pi-\eta=\underbrace{4+x-y}_{1}+\underbrace{(-2-x+y)}_{2}+
\underbrace{(+1+x-y)}_{3}=3+x-y,
$$
like above (!) hence soluble.

Some further idea would be to increase the ``predestination'' by
adding one or more triads as on Fig.\,f,g, yet such schemes are
already prohibited by B\'ezout. Actually to lower $\chi$ down to
Thom's range $\chi\le k^2$ we look at Fig.\,g, but the latter is
not even prohibited by Rohlin as $\pi-\eta=T_1+T_2+T_3$ is
contributed by 3 trunks each $T_i\in\{2,0,-6\}$ (by
Fig.\,\ref{Signs-law-triad:fig}d) hence soluble as by Rohlin
$\pi-\eta=6$. On Fig.\,f instead we have only 2 trunks so the
scheme is prohibited by Rohlin (of course more elementarily by
B\'ezout), yet it is also by Thom. So it does not solve our
problem of finding a scheme where Rohlin is stronger than Thom.

Hence Fiedler's example looks a typical case of predestination of
the Rohlin mass $\Delta \Pi=\pi-\eta$. Of course it may be
generalized in higher degrees than $10$, by looking at $M$-schemes
of the form $(1,1,1,1)M-4$. Then by Rohlin
$2(\pi-\eta)=r-k^2=M-k^2=(g+1)-k^2=(2k-1)(k-1)+1-k^2
=k^2-3k+2=(k-1)(k-2)$. But on the other hand by the signs-law for
triad (Fig.\,\ref{Signs-law-triad:fig}d) we have
$\Delta\Pi:=\pi-\eta\le 2$. It follows that the scheme considered
is prohibited as soon as $\pi-\eta=(k-1)(k-2)/2$ is $\ge 3$ that
is for $k\ge 4$. (The case $k=4$ being quite stupid for B\'ezout
would have sufficed.)

Okay, but such schemes are also prohibited by Thom. One could also
try to deepen the nest as the degree increase. Yet our goal is
really to find a ``Caucasian'' scheme, i.e. obstructed by Rohlin,
but not by Thom (nor by Gudkov or Arnold).

[14.03.13] So let us approach this problem more systematically.
First when $m=6$ it is clear that there is no Caucasian $M$-scheme
as follows from Gudkov's table (=Fig.\,\ref{Gudkov-Table3:fig}).
Indeed all the Thom permissible $M$-schemes are prohibited by
Gudkov.

\begin{figure}[h]
\centering
\epsfig{figure=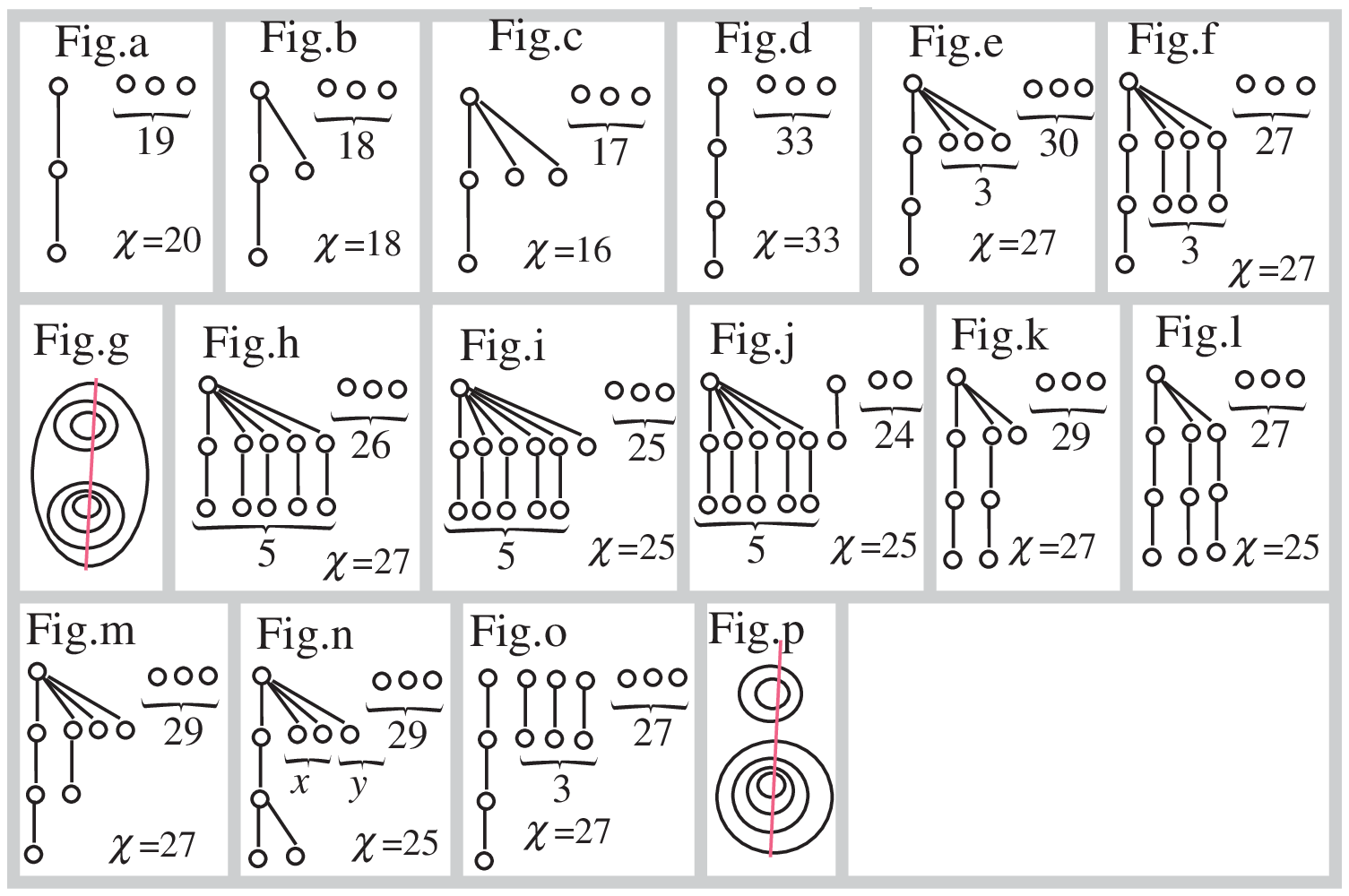,width=122mm} \vskip-5pt\penalty0
  \caption{\label{Fied1:fig}%
  Searching Caucasian schemes} \vskip-5pt\penalty0
\end{figure}

So we move to $m=8$, so $M=22$. Here we start with the $M$-scheme
$(1,1,1)19$ with $\chi=20$ (cf. Fig.\,\ref{Fied1:fig}a). As the
3-nest corresponds to a dyad (2 pairs of length 1) they contribute
by the signs-law to at most 1 to $\pi-\eta$. Hence Rohlin's
formula $2(\pi-\eta)=r-k^2=22-16=6$ (i.e. $\pi-\eta=3$) cannot be
solved. Of course the scheme in question is also prohibited for
deeper reasons (at via deeper results) like Gudkov hypothesis
$\chi\equiv k^2=16\pmod 8$, or Thom's inequality $\chi \le k^2$.
Our goal is to find an $M$-scheme where the ``elementary'' Rohlin
formula becomes stronger than the conjunction of 2 deep results
(Gudkov hypothesis proved by Rohlin-Rohlin/Atiyah-Singer-Marin and
Thom proved by Kronheimer-Mrowka).

Can we adjust the invariants as to
neutralize Gudkov and Thom? A first idea is to drag one outer oval
at depth 1. Then the contribution to the Rohlin (signed) mass
$\pi-\eta$ is still $\le 2$, and so his formula cannot be solved.
Yet doing so we have $\chi=18$ and the scheme is also prohibited
by Gudkov or Thom.

If we delocalize one more outer oval at depth 1 we get Fig.\,c,
where now Rohlin's formula can be satisfied, and as $\chi=16$ both
Gudkov and Thom are happy. So we see some annoying (as far as our
Caucasian policy is concerned) concomitance between the 3 forces
involved Rohlin, Gudkov, Thom.

This phenomenon is not specific to degree 8 and repeats itself in
degree $m=10$, where $M=37$ (temperature of the body of a
primate). Let us experiment this concretely. As B\'ezout now
permits, we start now with a deep nest of profundity 4 (4 nested
ovals) and add 33 outer ovals to reach Harnack's bound $M=37$ (cf.
Fig.\,d, which involve a ``triad'' chain with 3 consecutive
edges). By Rohlin's formula $2(\pi-\eta)=r-k^2=37-25=12$, so
$\pi-\eta=6$. However by the signs-law for triad the contribution
of the triad is at most $2$ (cf. Fig.\,\ref{Signs-law-triad:fig}),
and thus the scheme $(1,1,1,1)33$ is prohibited by Rohlin. As
$\chi=33\equiv_8 k^2=25$, the scheme is not prohibited by Gudkov,
but it is by Thom. So our Caucasian goal is not achieved.

To improve the situation we have to lower $\chi$ down to $k^2=25$.
As far as the signs-law is involved we can transfer at most 3 oval
at depth 1 (like on Fig.\,e), so as to have a contribution to
$\pi-\eta$ still $\le 2+3=5<6$. Doing so $\chi=27$, and the scheme
is still prohibited by Thom (and anew by Gudkov). Another idea is
to use 3 dyads as on Fig.\,f as the latter also contribute to at
most $1$ to  Rohlin's mass $\pi-\eta$. Alas doing so does not
diminish $\chi$ in the Thom range, and actually violates B\'ezout
(cf. Fig.\,g). This can be remedied if we abort the triad, and
look at a configuration with 5 dyads (the maximum possible while
still taking care to making Rohlin's equation $\pi-\eta=6$
impossible). This gives Fig.\,h with alas still $\chi=27$. So this
schemes is prohibited by both Rohlin, Gudkov, and Thom (but as far
as I see not by B\'ezout even for conics).

If we nest one more outer oval we may get Fig.\,i with $\chi=25$,
but suddenly Rohlin's equation is now soluble (choose e.g.
$+$-signs throughout). Likewise we may consider Fig.\,j, but
Rohlin is likewise soluble. Considering Fig.\,k still calibrated
as to make Rohlin's equation impossible (as each triad contribute
for $\le 2$), we only reach $\chi=27$ (of course this
configuration is anti-B\'ezout). We can push it further to
Fig.\,l, which despite being prohibited by B\'ezout it is not by
Gudkov nor by Thom, yet alas not by Rohlin since $2+2+2=6$ and so
Rohlin's equation is soluble taking $+$-signs on all edges.

Our naive cuneiform construction can still be more varied, yet it
seems unlikely that we will ever find a Caucasian scheme by this
method. We still consider Fig.\,m (not interesting). Next look at
Fig.\,n with $\chi=25$. Denoting by $x$, $y$ the number of $+$'s
 resp. $-$'s on the edge right above the corresponding letter, the
signs-law gives (after fixing $+$ on the trunks)
$$
\pi-\eta=\underbrace{4+x-y}_{1}+\underbrace{(-3)}_{2}
+\underbrace{(+2)}_{3}=3+x-y,
$$
whence the system $x-y=3$, $x+y=3$ soluble as $(x,y)=(3,0)$. So
the scheme is not obstructed by Rohlin.

We can also transmute Fig.\,d into Fig.\,o, where again 3 branches
are added so as to keep Rohlin's  formula impossible, but again
$\chi$ drops only to $27$ (and not $25$). Of course such a scheme
is defacto prohibited by the B\'ezout-Hilbert bound on the depth
of nests (cf. Fig.\,p).

If it is not possible to find a Caucasian scheme in degree $10$,
what about degree 12. First we need to extend the signs-law to
tetrads. While the latter involves for dyads a square (4 possible
products of two signs), and for triads a cube (with 8 possible
signs combinations), we have now a 4D-hypercube with $2^4=16$
combinations. The signs-law for tetrads is depicted below
(Fig.\,\ref{Signs-law-tetrad:fig}). It corroborates the  a-priori
reasonable expectation that the maximum contribution arises when
all 4 signs are $+$, in which case the contribution to Rohlin's
mass $\pi-\eta$ is 2. Still a priori we may expect that our
(Caucasian) game will not become easier since Harnack's bound
$M=g+1=(2k-1)(k-1)+1=2k^2-3k+2$ increases much faster than Thom's
bound $k^2$ and so we will have more pain to lower down $\chi$ in
Thom's range. Despite these objections {\it a priori\/} let us
track down our prey more slowly.


\begin{figure}[h]
\centering
\epsfig{figure=Signs-law-tetrad.eps,width=122mm}
\vskip-5pt\penalty0
  \caption{\label{Signs-law-tetrad:fig}%
  Signs-law for tetrads (hypercube in 4D-space)} \vskip-5pt\penalty0
\end{figure}

[15.03.13] Hence:

\begin{lemma}
The contribution to Rohlin's mass $\pi-\eta$ of a deep nest is:

$\bullet$ for a dyad either $+1$ or $-3$,

$\bullet$ for a triad either $+2,0$ or $-6$,

$\bullet$ for a tetrad either $+2,0,-2$ or $-10$.
\end{lemma}

\begin{figure}[h]
\centering
\epsfig{figure=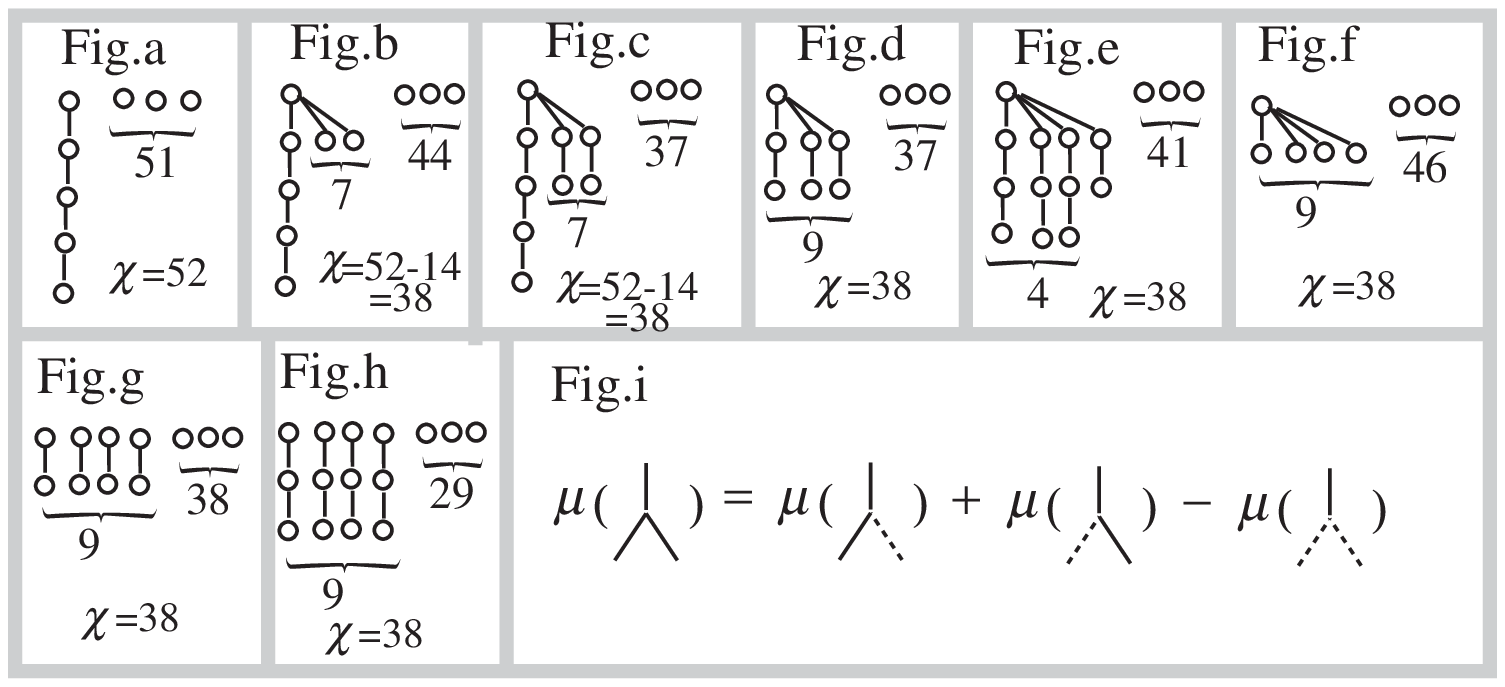,width=122mm} \vskip-5pt\penalty0
  \caption{\label{Fied2:fig}%
  Signs-law for tetrad (hypercube in 4D-space)} \vskip-5pt\penalty0
\end{figure}

With some combinatorial ingeniousness it should be easy to extend
to the general case. However let us first tackle our Caucasian
problem in degree 12. Again we resort to the cuneiform formalism
used above. When $m=12$, $M=g+1=\frac{11\cdot 10}{2}+1=55+1=56$.
We consider first Fig.\,\ref{Fied2:fig}a. Here $\chi=52$, while
Gudkov says $\chi\equiv_8 k^2=36=44=52$, which is verified. By
Rohlin's formula $2(\pi-\eta)=r-k^2=56-36=20$, so $\pi-\eta=10$,
but the contribution of the tetrad is at most 2 by the signs-law,
and thus the scheme is prohibited by Rohlin. Of course it is also
by Thom $\chi \le k^2$. As above the game is to lower $\chi$, by
transferring outer ovals at depth 1. As to keep Rohlin in
defeat, we may add at most 7 branches as on Fig.\,b, and obtain a
scheme with $\chi=52-2\cdot 7=52-14=38$, again 2 unit above Thom's
bound. The next idea is to let branches of dyads (dyadic branches)
hang on. Each contributes at most $1$ to Rohlin's mass $\pi-\eta$,
and so keeping Rohlin's formula in check
we may add 7 of them, but of course $\chi$ remains unchanged to
$\chi=38$ (as we merely traded outer ovals at depth 0 for ones at
depth 2). Next Fig.\,d involves only dyads contributing for at
most 1, so keeping Rohlin check-mate we can plug 9 of them, and
the resulting $\chi$ is still $38$. Using instead triadic branches
as on Fig.\,e which contributes for 2 (at most) we may plug 4 of
them and the remaining unit is consumed by inserting a dyadic
branch, and we find of course again $\chi=38$. Considering only
monadic branches as on Fig.\,f (und zwar 9 of them to Rohlin in
check) yields again $\chi=38$. If like on Fig.\,g the monadic
branches are not subsumed to a single dominator we find again
introducing 9 of them, $\chi=38$.  Using instead 9 dyadic branches
of the same sort, we get Fig.\,h, where still $\chi=38$.

So it looks again hard to find a Caucasian scheme where Rohlin is
stronger than Thom (and Gudkov united). Either we are looking at
the wrong place or there is some subsumation of Rohlin to Thom,
for some trivial arithmetical reasons. That is assume you have an
$M$-scheme with $\chi\equiv_8 k^2$ (Gudkov) and $\chi\le k^2$ then
Rohlin's equation is always soluble. Indeed write formally
$2(\pi-\eta)=r-k^2$. By the $M$-curve assumption $r-k^2$ is
$r-k^2=(2k-1)(k-1)+1-k^2=k^2-3k+2$, etc.

One could also ask if for $M$-curves the Thom and Gudkov
obstructions are the sole one, but this is probably corrupted by
B\'ezout-style obstruction \`a la Fiedler-Viro (in degree 8
already).

\subsection{Some weak evidence against Caucasian schemes}

[15.03.13] In the previous section we tried (hard) to find a
``Caucasian'' scheme, i.e. prohibited by Rohlin's formula
$2(\pi-\eta)=r-k^2$ but not by Thom $\chi\le k^2$, but failed. Of
course as Rohlin's formula formally implies Arnold congruence it
is actually a simple matter to find such a scheme, e.g. any
$M$-scheme violating Arnold's congruence but not Thom do the job.
For instance in degree 6, the $M$-scheme $\frac{2}{1}8$ is
prohibited by Rohlin, but not by Thom. Likewise for the
$(M-2)$-scheme $\frac{1}{1}7$ in type~I.

However if we add the Gudkov hypothesis as a side condition (or
perhaps just the Arnold congruence) then it seemed difficult to
find a Caucasian scheme where Rohlin is stronger than Thom. The
sequel tries to give some evidence that it is impossible to find a
Caucasian scheme, but our argument will be somewhat loose. The
main difficulty is the mess arising with the signs-law and so the
difficulty looks merely combinatorial. Of course it is not
impossible that we missed a trivial counter-example that impedes
the completion of the present programme.

First we described in the previous section the mess arising form
the signs-law applied on dyad, triads, up to tetrads, which are
totally ordered chains. In general the situation is more tricky as
the Hilbert tree of the scheme may be highly branched.

Recall that to a dividing curve is assigned complex orientations
(up to global reversion of all of them), which in turn decorates
the Hilbert tree (of the curve) with a signs-distribution
(abridged charge) making it into what we call the Rohlin tree. (Of
course the tree can be a ``forest'', i.e. have several components,
and it really ``branches downwards'' as to look more like Arnold's
paradigm of the mushrooms.)

Each such (Rohlin) signed tree is completely determined by the
signs ascribed to the edges of length 1 as it then extends to
longer edges by the signs-law. Further each such tree has a Rohlin
mass $\mu=\pi-\eta$ which is the difference between $\pi:=\Pi^+$
the number of positive pairs and $\eta:=\Pi^-$ the number of
negative pairs.

When the Rohlin tree is the one induced by a curve of degree
$m=2k$, (the marvellous) Rohlin's formula $2(\pi-\eta)=r-k^2$
fixes the mass in function of $r$ the number of ovals and $k$ the
semi-degree. Further $\pi+\eta=\Pi$ is the total number of pair,
and so one gets the wrong impression that Rohlin's equation is
always soluble but there is some hidden r\^ole played by the mass
$\mu:=\pi-\eta$ which for certain configuration can only be very
low (especially on chains), cf. Fiedler's example with
$(1,1,1,1)33$ in degree 10.

So the core of the problem is to understand  the behavior of the
Rohlin mass $\mu$. In the previous section we nearly understood
this for vertical chains. Especially easy, is the case where all
signs are positive, in which case it is a simple matter to show
the:

\begin{lemma}
Given an $n$-chain of $n$ consecutive edges all positively charged
then the Rohlin mass of the chain is the integer part
$[\frac{n+1}{2}]$
\end{lemma}

\begin{proof} Look and see (i.e. make pictures). Indeed looking at
the Signs-laws for triad ($n=3$) or tetrad ($n=4$), we get resp.
$\mu=3-2+1=1+1=2$ and $\mu=(4-3)+(2-1)=1+1=2$. For a 5-chain this
extends as $\mu=(5-4)+(3-2)+1=3$, for a 6-chain as
$\mu=(6-5)+(4-3)+(2-1)=3$, and so on.
\end{proof}

In general for a signed tree there ought to be a sort of skein
relation permitting an iterated evaluation of the Rohlin mass
$\mu$, based on the formula that $\mu$ of an inverted ``Y''
looking like a $\Lambda$ surmounted by a chain is equal to $\mu$
of the left maximal chain in the inverted ``Y'' plus the right
chain, minus $\mu$ of the common trunk. If this is not clear,
please compare Fig.\,\ref{Fied2:fig}i. This formula is of course a
formal consequence of the inclusion-exclusion principle in
combinatorics or measure theory.

Now in sloppy fashion first, the idea  is that under certain
assumptions (like the conjunction of Gudkov and Thom's $\chi\le
k^2$) one could show the existence of a charge (=distribution of
signs) solving Rohlin's equation. Perhaps this can be done via a
sort of linear algebra modulo 2.

First the charge in question is merely a $\pm 1$-valued function
on the set $E$ of all edges of the tree. The set of all such
distributions denoted $\cal E$ can be turned into a vector space
over the field ${\Bbb F}_2$ with 2 elements. Define indeed the sum
of two charges $\epsilon, \delta$ as
$(\epsilon+\delta)(e)=\epsilon(e)\times \delta(e)$, where $\times$
is the Rohlin product given by the signs-law (i.e. the opposite of
the usual sign convention for products). The neutral element is
$\epsilon_0$ the minus distribution, as
$(\epsilon+\epsilon_0)(e)=\epsilon(e)\times (-1)=\epsilon(e)$.
Also each element has order 2, as it should. (More boring details
in p.\,AR99=hand-notes, especially the inversed charge where all
signs are switched is not the inverse charge!) The multiplication
by a scalar is naturally defined. The Rohlin mass is the function
$\mu\colon \cal E \to {\Bbb Z}$.

Now that  we have a good vector space, we could hope  our problem
reducible to linear algebra! Intuitively if $\chi\le k^2$ then
there must be enough edges as to solve Rohlin's equation. More
precisely Rohlin's formula fixes the mass via $2\mu =r-k^2$, and
via the skein relation the mass of the tree reduces to that of
chains which in turn can be reduced to that of edges via the
signs-law, i.e. the knowledge of the charge $\epsilon$ itself. If
one is good in combinatorics there is a little hope to show that
each skein relation induces a linear equation and count that there
is enough free parameter as to solve the equation. {\it
Warning.}---As remarked in more details latter, already in degree
6 we have the scheme $5$ whose type~I realization is prohibited by
Rohlin but not by Thom, so there is no chance to complete this
programme, unless extra assumptions are added, e.g. that of being
an $M$-curves (which further must satisfy the Arnold congruence,
else prohibited by Rohlin but not by Thom).

As a trivial example we may extend the observations of the
previous section. Consider an $M$-scheme without nesting. Then
$\chi=M$. Now to diminish $\chi$ we introduce edges (i.e. nested
pairs), cf. Fig.\,\ref{Fied2:fig}g for an example. We have
$r=M=g+1=(2k-1)(k-1)+1=2k^2-3k+2$, so $r-k^2=k^2-3k+2=(k-1)(k-2)$.
By Rohlin's formula $\mu:=\pi-\eta=(k-1)(k-2)/2$. Hence to keep
Rohlin's formula in default, we introduce only $\mu-1$ edges. Then
we compute $\chi$, and find $\chi=r-2(\mu-1)$ and a boring
calculation shows this to be $k^2+2$. So if obstructed by Rohlin
then also by Thom. (More details in p.AR.99 of the hand-notes.) Of
course this is not the general case as the scheme has a very
specific structure akin to Fig.\,\ref{Fied2:fig}g. Can we
generalize, perhaps but requires to work out some messy
combinatorics.

A somewhat more appealing idea is that if the scheme is obstructed
by Rohlin then it is because its mass $\mu$ is strictly less than
$\pi-\eta=(r-k^2)/2$ even when the tree is positively charged, and
conjecturally this should maximize the mass. All this is vague but
points to the right direction.
Namely it gives the idea of computing the Rohlin mass of a tree
with positive charges only. The answer turns out to be simple and
elegant:

\begin{lemma}\label{Rohlin-mass-of-a-positively-charged-tree:lem}
The Rohlin mass $\mu$ of a positively charged (signed) tree $T$ is
equal to
$$
\mu(T)=n_1+n_3+n_5+\dots=n,
$$
the number of ovals at odd depths, where $n_1$ counts those at
depth $1$, $n_3$ at depth $3$, and so on.
\end{lemma}

\begin{proof} Make a picture of a tree with possibly several
components. Put plus signs everywhere as stipulated. By additivity
we may focus attention on a single component of the tree. Each
vertex at depth $\ge 1$ has exactly one edge above it. All pairs
are enumerated by starting from vertices at depth 1 and looking at
edges above them gives the contribution $n_1$. Next we look at the
$p_2$ many vertices at depth 2, each inducing two pairs above it
(of length 1 and 2 resp.). The first contribute for $+1$, while
the other has sign $-1$ by the signs-law $+\times +=-$ (recall
that consanguinity is bad). So  ovals at depth 2 contributes for
$p_2-p_2=0$. Continuing in this fashion we find:
$$
\mu=n_1+(p_2-p_2)+(n_3-n_3+n_3)+(p_4-p_4+p_4-p_4)+\dots,
$$
which implies the announced formula.
\end{proof}

As implicit above we posit the:

\begin{conj}\label{positive-mass-conjecture:conj}
The Rohlin mass of a signed (Rohlin) tree is maximized when the
tree is positively charged throughout.
\end{conj}

Some evidence comes from the case of chains (as tabulated on the
signs-law tables, e.g. Fig.\,\ref{Signs-law-tetrad:fig}). There is
perhaps a simple argument.

But do we really need this? Let us make another observation based
on the lemma.

Assume that Thom holds, i.e. $\chi\le k^2$. In general, we have:
$$
\chi=p_0-n_1+p_2-n_3+p_4-\dots,
$$
where each symbol $p_i, n_i$ counts the number of ovals at depth
$i$, where $p,n$ are just ``residue'' of Petrovskii notations for
positive and negative but best interpreted in terms of even or odd
depth resp. (The notation are nearly consistent in
French-Swiss-German, where ``even=pair'' and ``odd=uNgerade''.)

Besides, the total number of ovals, denoted $r$, is expressible as
$$
r=p_0+n_1+p_2+n_3+p_4+\dots,
$$
so that subtracting the double of the  Rohlin mass $\mu$ of the
positively charged Rohlin tree (as calculated in the lemma) gives
the relation:
$$
r-2\mu=\chi,
$$
which holds universally when the tree is positively charged.

So if Thom is verified, i.e. $\chi\le k^2$, we find
$2\mu=r-\chi\ge r-k^2$. This means that there is no obstruction
{\it a priori\/} to solve Rohlin's equation, since Rohlin's mass
is as large as it should by virtue of Rohlin's formula
$2\mu=r-k^2$. Paraphrasing, Rohlin's equation is virtually
soluble.

\begin{lemma}
If Thom's estimate $\chi\le k^2$ is fulfilled, then there is no
``quantitative'' obstruction to solve Rohlin's equation. Yet
beware that there may of course be finer arithmetical reasons
impeding solubility as with the scheme $5$ of degree $6$ which has
no type~I realization.
\end{lemma}

This prompts some evidence toward the:

\begin{conj}
There is no Caucasian scheme, where Rohlin is stronger than Thom
(at least modulo adding some suitable hypotheses, e.g. that of an
$M$-scheme).
\end{conj}

Can we find a formal proof? A crudely idea is to notice that if we
charge the tree negatively throughout then by the signs-law all
pairs are negative as $-\times -=-$. So the Rohlin mass of the
negatively charged tree is $-\Pi$, where $\Pi$ is the total number
of pairs. Hence by a dubious mean-value theorem (in the
discontinuous realm of the arithmetics of quanta) we would like to
infer  existence of a charge fulfilling Rohlin's formula.

Another idea is to introduce indeterminate signs and try to solve
a system of linear equations. We did this frequently formerly, but
we had some grasp on the geometry of the tree. Whether this can be
done in abstracto is not
clear to me, and may of course converge to the first strategy
using linear algebra on the spaces of all charges (plus the
skein-relation).

Of course recall that we have the Rohlin-Marin inequality that a
dividing curve $C_{m=2k}$ has $r\ge m/2=k$, i.e. at least as many
ovals as the semi-degree. This is a formal consequence of Rohlin's
formula and precludes in degree $6$ a type~I incarnation of the
scheme $1$ (unifolium), which is not prohibited by Arnold's
congruence mod 4 (cf. the Gudkov-Rohlin
table=Fig.\,\ref{Gudkov-Table3:fig}). This example (or also the
scheme $5$ in degree 6) are of course trivial counterexamples to
the above conjecture (freed from the parenthetical proviso). The
latter gains however some more credibility when the curve is
assumed to be an $M$-curve (or perhaps even an $(M-2)$-curve).

It is evident that our whole problem is somewhat ill-posed, yet we
hope to have demonstrated that some
complicity between Rohlin and Thom requires to be elucidated.

The simple example of the scheme $5$ in degree 6, where Thom's
estimate $\chi\le k^2$ as well as Arnold's congruence are
fulfilled, but whose realizability in type~I is precluded by
Rohlin's formula (as $r$ is not a square and there is no nesting
hence $\Pi=0$, and so a fortiori $\pi=\eta=0$) shows that our
above desideratum of solving Rohlin's equation under the sole
assumption of being in Thom's range is not realistic. So one must
really add some extra assumptions, typically that of being an
$M$-curve, which can perhaps be somewhat relaxed.

Let us close the discussion via precise conjectures:

\begin{conj}
$\bullet$ An $M$-scheme verifying the Gudkov congruence
$\chi\equiv k^2 \pmod 8$ and the Thom estimate $\chi\le k^2$ is
never prohibited by Rohlin's formula.

$\bullet$ An $(M-2)$-scheme verifying the
RKM=Rohlin-Kharlamov-Marin congruence $\chi\equiv k^2 +4 \pmod 8$
(ensuring the scheme to be of type~I) and the Thom estimate
$\chi\le k^2$ is never prohibited by Rohlin's formula.

\end{conj}

\subsection{A basic mistake in the search of a Caucasian scheme}

DO NOT READ THE SEQUEL IT IS FALSE!

\begin{figure}[h]
\centering
\epsfig{figure=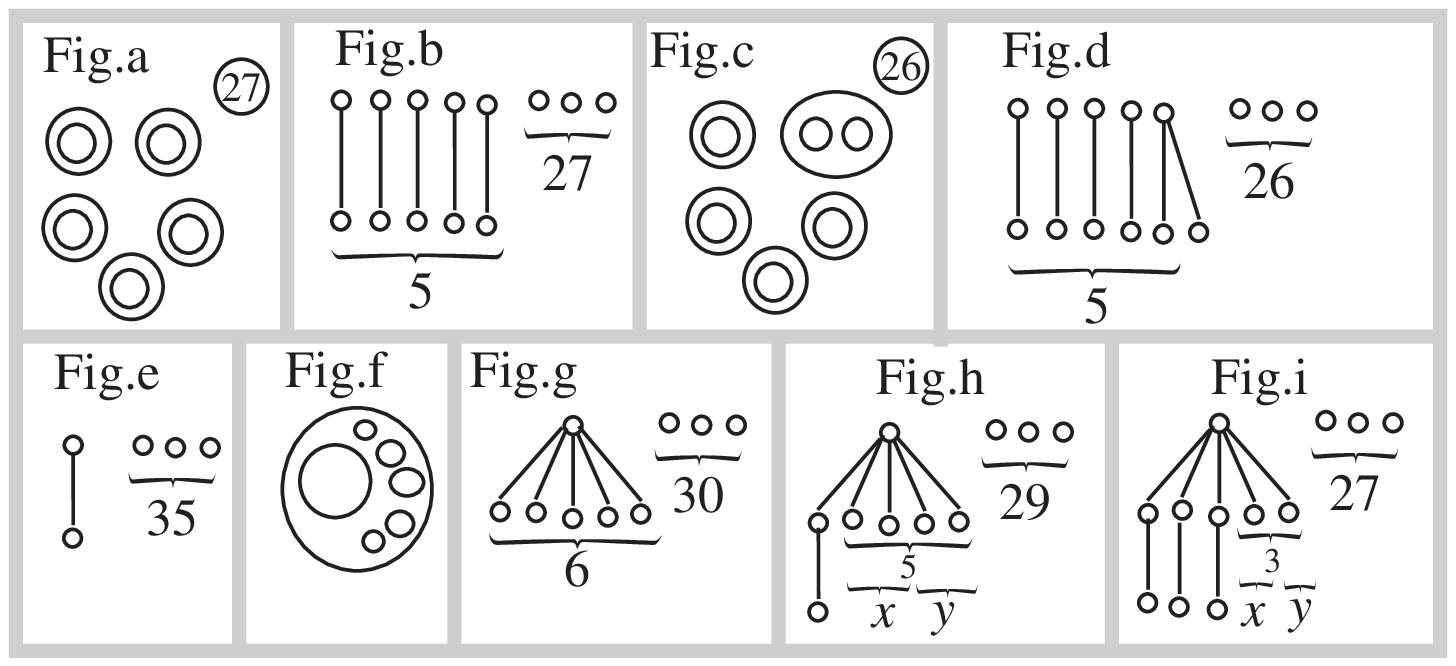,width=122mm} \vskip-5pt\penalty0
  \caption{\label{Fiedler6:fig}%
  Finding a Caucasian scheme of degree 10} \vskip-5pt\penalty0
\end{figure}

Our next idea was to look at an $M$-scheme extending 5 nests of
depth 2 (cf. Fig.\,\ref{Fiedler6:fig}a,b). By Rohlin's formula we
still have $\pi-\eta=6$, but $\pi-\eta\le \pi+\eta=\Pi=5$ so
Rohlin is violated, but alas $\chi=27$ hence the scheme is also
prohibited by Thom. Can we diminish $\chi$? Yes as usual by
trading an outer oval at depth 0 against one at depth 0, cf.
Fig.\,c. But then the corresponding tree (Fig.\,d) has $\Pi=6$
pairs and so Rohlin's formula is soluble with $\pi=6$ (i.e. all
edges positive).

Finally, we started from the scheme $(1,1)35$ (cf. Fig.\,e) with
$\chi=35$ and to lower to Thom's range $\chi=25$, we make 5 moves
to get Fig.\,f and its allied tree on Fig.\,g. But again we have
$\Pi=6$ and so can solve Rohlin's equation by putting only
positive signs.

From this last configuration, we decided to drag one of the ovals
at depth 2 to get Fig.\,h, which has still $\chi=25$. By Rohlin
$\pi-\eta=6$. Now introducing free variables $x,y$ counting
positive resp. negative signs on the 5 edges and choosing any
distribution of signs on the trunk of length 2, we know that the
latter will contribute for at most $\le 2$ to
$\Delta\Pi:=\pi-\eta$ (by Fig.\,\ref{Signs-law-triad:fig}d) and
its contribution $T$ is either $2,0,-6$, which is at any rate
even. WARNING: HERE I MADE A BASIC CONFUSION IN THE LENGTH OF THIS
CHAIN!!! On calculating $\pi-\eta$ by the signs-law we get
$$
\pi-\eta=T+x-y,
$$
so $x-y=6-T$ and $x+y=5$ so that $2x=11-T$ which is impossible
modulo 2! So we found our first scheme prohibited by Rohlin but
not by Thom, an therefore:

\begin{theorem}\label{Caucasian-scheme:thm} (WARNING=ERRONEOUS)
There exists a ``Caucasian'' scheme\footnote{So-called because
V.\.A. Rohlin was born in Bakou, from parents themselves coming
from Odessa, and if I am not wrong in Geography Bakou belongs to
Caucasus.}, where Rohlin is stronger than Gudkov and Thom. More
precisely the $M$-scheme of degree 10 of Fig.\,h that is
$(1,(1,1)5)29$---in Gudkov's notation---is prohibited by Rohlin's
formula, but not by Gudkov's hypothesis $\chi=k^2 \pmod 8$ nor by
Thom's inequality $\chi\le k^2$. The example has $\chi=25$ (by
construction).
\end{theorem}

Of course the above argument extends to the case where we drop 3
ovals at depth 2 (Fig.\,i) so that we have $x+y=3$ still odd.
(This could even be $x+y=1$.) Indeed denoting by $T_1,T_2,T_3$ the
contributions of the 3 trunks of length 2 on Fig.\,i and by $T$
their sum (which is even by signs-law for triad), we find by the
signs-law
$$
\pi-\eta=T+x-y,
$$
and therefore as $\pi-\eta=6$ (by Rohlin's formula) we have
$x-y=6-T$ and $x+y=3$, so that summing $2x=9-T$, which is
impossible modulo 2.

The phenomenon just discovered is probably not new and perhaps
related to Slepyan's law (also a formal consequence of Rohlin's
formula), cf. perhaps Rohlin 1978 or Degtyarev-Kharlamov 2000
\cite{Degtyarev-Kharlamov_2000}.

\begin{figure}[h]
\centering
\epsfig{figure=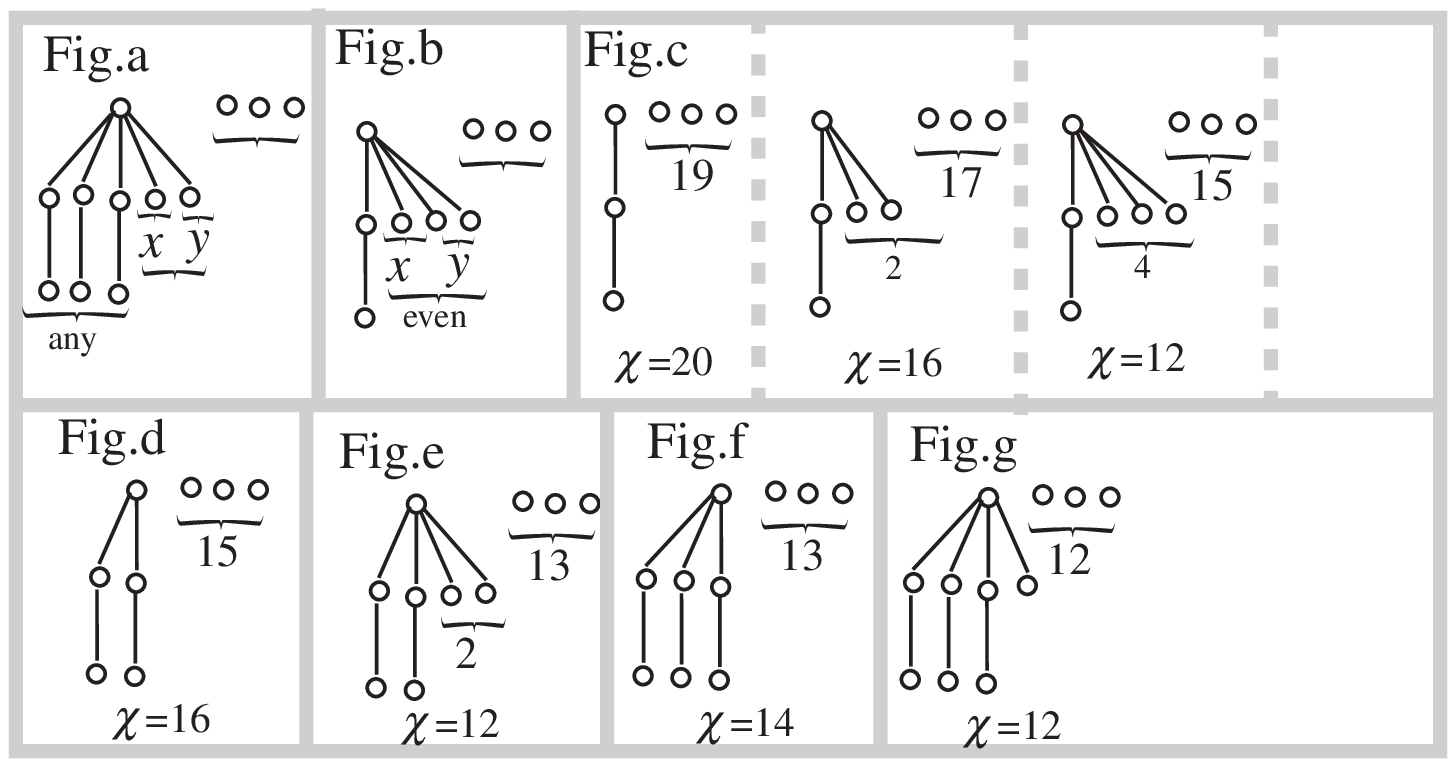,width=122mm} \vskip-5pt\penalty0
  \caption{\label{Fiedler7:fig}%
  More Caucasian schemes of degree 10 and some in degree 8} \vskip-5pt\penalty0
\end{figure}

[13.03.13] It is quite evident that we may generalize somewhat the
result. What seems essential to the argument is that the Hilbert
tree of the scheme as the structure of Figs.\,h.i with an odd
number $x+y$ of empty ovals at depth 1, so to be like
Fig.\,\ref{Fiedler7:fig}a with $x+y$ odd.

In fact let us be more general and leave degree 10. So suppose to
have an $M$-scheme of (arbitrary) degree $2k$ (or more generally a
scheme of type~I, but we reserve this for latter) so that the
dividing character is granted and therefore Rohlin's formula
applies.

We suppose additionally that the scheme is like Fig.\,a with $x+y$
denoting the number of empty ovals at depth 1, partitioned into
$x$ many positive pairs and $y$ many negative pairs (when looking
at the unique edge right above those ovals).

The argument is then to compute $\pi-\eta$ in two fashions. One
way involves Rohlin's formula $2(\pi-\eta)=r-k^2$, while the other
route involves the signs-law for triad (cf.
Fig.\,\ref{Signs-law-triad:fig}d) and gives
$$
\pi-\eta=T+x-y,
$$
where $T$ is the contribution of the trunks of length 2 which is
necessarily even (again by Fig.\,\ref{Signs-law-triad:fig}d). Now
if $r-k^2\equiv 0\pmod 4$, by Rohlin $\pi-\eta$ is even and the
signs-law equation is corrupted if $x-y\equiv_2 x+y$ is odd.
Viceversa if $r-k^2$ is not divisible by 4, then Rohlin says that
$\pi-\eta$ is odd, but the signs-law that it is even provided
$x+y$ is even.

Finally for $M$-curves it is a simple matter to check that
$r-k^2\equiv 0\pmod 4$ iff $k\equiv 1,2 \pmod 4$. Indeed
$$
r-k^2=(g+1)-k^2=(2k-1)(k-1)+1-k^2=k^2-3k+2,
$$
which is mod 4 for $k=1$, $1-3+2=0$ and for $k=2$, $4-6+2=0$,
while for $k=3$, it is $9-9+2=2$ and for $k=0=4$, it is $2$.

Hence we have proved the:

\begin{theorem}
Define a dendritic scheme as one like depicted on
Fig.\,\ref{Fiedler7:fig}a, i.e. with Gudkov symbol of the form
$(1,(1,1)\dots(1,1) x+y ) z$.

If $k\equiv 1,2 \pmod 4$, there is no $M$-curve of degree $2k$
with dentritic scheme having an odd number of empty ovals at depth
1.

If instead $k\equiv 1,2 \pmod 4$, there is no such curve having an
even number of empty ovals at depth 1.
\end{theorem}

Some additional remarks are in order, which we detail right after.

1.---First it is a simple matter to see that Caucasian schemes
exists already in degree 8.

2.---Second it seems clear that we may formulate of the theorem
for $(M-2)$-schemes satisfying the RKM-congruence ensuring type~I.

3.---Third, we could expect to extend the result to the case where
there are several dendrites, or deeper nests.

1.---Indeed for definiteness we may assume that there is a single
trunk of length 2 (like on Fig.\,b). For $m=8$, $M=22$ and so for
an $M$-curve we have $2(\pi-\eta)=r-k^2=22-16=6$, so that
$\pi-\eta$ is odd. Yet calculating via the signs-law prompts that
$\pi-\eta=T+x-y$ which is even provided the number $x+y$ of
branches of length 1 is even. (THIS IS AGAIN A MISTAKE CAUSED BY
CONFUSION IN THE LENGTH OF THE CHAIN: BASICALLY IT INVOLVES 3
OVALS BUT THE LENGTH IS 2!!!) So we get schemes prohibited by
Rohlin along the series depicted on Fig.\,c, which traduced in
Gudkov's symbols gives the list
$$
(1,1,1)19, (1,(1,1)2)17, (1,(1,1)4)15, (1,(1,1)6)13, (1,(1,1)8)11,
\dots, (1,(1,1)18)1,
$$
where $\chi$ is first $\chi=(1-1+1)+19=20$ (hence the scheme is
also prohibited by Gudkov or Thom) and then successively drops by
4 units, so that the second listed scheme $(1,(1,1)2)17$ has
$\chi=16$ (hence not prohibited by Gudkov nor by Thom, but
prohibited by Rohlin), and so on. At $\chi=8$ we find another
Caucasian scheme $(1,(1,1)6)13$.

2.---For $(M-2)$-schemes satisfying the RKM-congruence $\chi\equiv
k^2+4 \pmod 8$ we are ensured to be of type~I. We consider again a
dentritic scheme like on Fig.\,\ref{Fiedler7:fig}a and repeat the
above argument. By Rohlin $2(\pi-\eta)=r-k^2$, while on the other
hand $\pi-\eta=T+x-y$ where $T$ is the contribution from the
trunks which is even. As $x+y\equiv_2 x-y$, we get a contradiction
as soon as $x+y$ and $(r-k^2)/2$ have opposite parities. When
$r=M-2$, $r-k^2$ is congruent to 0 mod 4 precisely when $k\equiv
0,3 \pmod 4$ (just shift by 2 the previous calculation).

To get an example, consider $m=8$ and take just one trunk like on
Fig.\,c, but now consider the $(M-2)$-scheme $(1,1,1)17$. It has
$\chi=18$. But we impose the RKM-congruence, hence adjust $\chi$
to $12$. Hence we consider $(1,(1,1)3)14$, and this is prohibited
by Rohlin. Indeed his formula becomes $2(\pi-\eta)=r-k^2=20-16=4$,
which is divisible by 4 (as predicted above), while by the
signs-law $\pi-\eta=T+x-y$ is odd since there are $x+y=3$ little
branches hanging on. [21.03.13] THIS IS ERRONEOUS BUT MAYBE CAN BE
CORRECTED BY CHANGING OF PARITY, as the trunks contributes for an
odd number.

So the real outcome of this method of prohibition based on
Rohlin's formula and the signs-law seems to be a powerful tool for
prohibition. It remains of course to examine its exact
significance, and how it generalize to nest of deeper structure.

It may be noted that earlier in this text we attempted a complete
classification of RKM-schemes of degree 8. This was rather a
census, i.e. a weak form of combinatorially possible schemes yet
without any claim of realizability. Now with the present method we
see that some of them are prohibited. It remains to understand
which of them are prohibited by Rohlin enhanced by the signs-law.
Clearly the argument given above extends and implies the following
lemma.

\begin{lemma}\label{RKM-scheme-ruled-out:lem}
(ERRONEOUS) The four primitive types of RKM-schemes, i.e.
$$
(1,(1,1)3)14, (1,(1,1)7)10, (1,(1,1)11)6, (1,(1,1)15)2,
$$
already listed in
Equation~\ref{RKM-scheme-deg-8-four-primitive-type:eq}, are
prohibited by Rohlin's formula.
\end{lemma}

However their derived products looks harder to prohibit as they
ramify and do not anymore belong to the dendrite type.

For instance for the scheme $(1,(1,2)3)13$ one can easily solve
Rohlin's equation with a distribution of sign, and so probably for
all other derived products. That remains to be checked.

On the other hand, it is clear that our census of 100 schemes was
far from exhaustive. Indeed we may consider a dendrite with 2
trunks like on Fig.\,d. This has $\chi=16$, and to adjust the
RKM-congruence $\chi\equiv_8 k^2+4=20\equiv_8 12$ we move down to
$\chi=12$ by transferring 2 outer ovals at depth 1 to get Fig.\,e.
This scheme is not prohibited by Rohlin. Indeed
$2(\pi-\eta)=r-k^2=20-16=4$ (so $\pi-\eta=2$),  and by the
signs-law $\pi-\eta=T_1+T_2+x-y$ where each trunks contributes to
$T_i\in \{2,0,-6\}$ by Fig.\,\ref{Signs-law-triad:fig}d. Even if
we impose $T_1=T_2=2$, the system is still soluble, being $x-y=2$
and $x+y=2$, hence $(x,y)=(2,0)$.

As usual, the RKM-scheme of Fig.\,e whose (Gudkov) symbol is
$(1,\frac{1}{1}\frac{1}{1}2)13$ has a myriad of companions. One
can either:

$\bullet$ without changing $\chi$ transfer outer ovals at depth 2
(in various ways);

$\bullet$ drop $\chi$ by 8 units by transfer quanta of 4 ovals at
depth 1.

All this is a bit messy to write down, and this will still not be
exhaustive as we can start from the configuration with 3 trunks
(Fig.\,f) which has $\chi=14$, and to adjust to RKM we make one
transfer at depth 1 (Fig.\,g) and get so another RKM-scheme. This
times the scheme is prohibited by Rohlin as
$2(\pi-\eta)=r-k^2=20-16=4$ (so $\pi-\eta=2$), but by the
signs-law $\pi-\eta=T+x-y$ where each 3 trunks contributes evenly
while $x-y\equiv_2 x+y=1$ is odd. Likewise if from Fig.\,g we
transfer quanta of 4 outer ovals at depth 1 the number of branches
$x+y$ is still odd, and so those schemes are prohibited too.
However those schemes derived from Fig.\,g by transferring outer
ovals at depth 2 are probably not prohibited.

Then there is still the cases of 4, 5, etc. trunks and the
classification looks quite messy to obtain. Fortunately the story
the story as soon as the tree contains 4 disjoint edges, since
this correspond to 4 nest of depth 2 (a configuration which is
saturated by B\'ezout or better Rohlin's maximality principle).
Note incidentally that this principle also prohibits the schemes
like Fig.\,g with 3 trunks since there are extra branches, so that
the above prohibition via Rohlin's formula can be subsumed to
total reality. [Warning this last sentence looks dubious!] Albeit
messy, it could be of primary importance to get a good view of
what happens along the way to extend the Rohlin-Le~Touz\'e
phenomenon of total reality from degree 6 to degree 8.

So a pivotal question is whether there is any reasonable way to
list all RKM-schemes of degree 8? If so then make some cleaning by
ruling out those prohibited by Rohlin's formula (enhanced by the
signs-law), and finally try to understand which are realized
algebraically (simplified form of Hilbert's 16th nearly solved by
the experts, but by far by myself, as one requires certainly the
Viro method). Once this s achieved try to understand if all those
schemes are subjected to the phenomenon of total reality (probably
under pencil of quintics, as we discussed earlier). If so then
there is some chance that Rohlin's maximality conjecture holds
true in degree $8$.

All this requires either Herculean forces or some good idea.

\subsection{Toward a complete census of RKM-schemes of degree~8}

[13.03.13] Our goal is to list all RKM-schemes of degree 8. Those
are $(M-2)$-scheme satisfying the RKM-congruence $\chi\equiv
k^2+4=20 \pmod 8$. Of course there is a menagerie of them, but we
have also upper bound given by the saturation principle of Rohlin
allied to total reality, once the depth is 4 the configuration is
saturated and cannot develop further. Actually the nest of depth
4, is not an $(M-2)$-scheme and so the depth is at most 3 (or 2
depending on the way you count). Likewise the pencil of conics
shows that that there can be at most 4 nests of depth 2, or when
translated in the cuneiform language of Hilbert's tree there is at
most 4 edges which are disjoint.

With this upper bound in mind, there is some little hope to make a
complete classification. Further as the number of ovals $r$ is
fixed to $M-2=20$, we may kill all empty ovals lacking a superior
(so-called outer ovals) and thus condense a bit notation. So to
each scheme is assigned a ``skeleton'' (kill the outer ovals) from
which we may recover the scheme unambiguously.

We abort this project as it is quite overwhelming.

\subsection{Trying to corrupt Rohlin's maximality conjecture (RMC)}

WARNING DO NOT READ: FULL OF MISTAKES. But try to correct at the
occasion. Compare p.AR95--96 for the original, and keep in MIND
the example of p.AR96. This is the $M$-scheme of degree 8 with
symbol $\frac{1}{3}\frac{1}{17}$ which is prohibited by Rohlin,
since $\pi-\eta=3$ (via Rohlin), but $\pi+\eta=20$, so $2\pi=23$
which is not soluble integrally. BUT yet another MISTAKE, this
scheme is already prohibited by Gudkov (or even Arnold!) as
$\chi=(1-3)+(1-17)=-18$

[14.03.13] Another basic question is whether from all schemes (in
particular $M$-schemes) that we prohibited via Rohlin (and the
signs-law), if it is not possible to corrupt RMC by finding an
$(M-1)$-scheme right below which is maximal, but of course not of
type~I by Klein's congruence $r\equiv g+1 \pmod 2$ for dividing
curves.

One such scheme in degree 8 was given  by $(1,\frac{1}{1}2)17$
(cf. Fig.\,\ref{Fiedler7:fig}c) with $\chi=16$. WARNING THIS IS A
MISTAKE, as Rohlin's equation can be solved with all signs
positive!!!
 Of course there is
a myriad of other such $M$-schemes verifying the Gudkov congruence
$\chi\equiv k^2 \pmod 8$ but prohibited for different reasons
(Rohlin with signs-law, or by Thom (\ref{Thom-Ragsdale:thm})). So
the vague idea is that if the scheme is French or Caucasian (i.e.
prohibited by Thom resp. Rohlin but not by Gudkov) then the scheme
is nearly realized in the sense that killing one of its oval then
the GKK-congruence $\chi\equiv k^2\pm 1 \pmod 8$ (cf.
\ref{Gudkov-Krakhnov-Kharlamov-cong:thm}) is satisfied and so
there is some hope to construct some $(M-1)$-curve. Further, and
this is the most dubious part of the game, we would like that the
resulting $(M-1)$-curve cannot be enlarged, which requires to
inspect a menagerie of schemes.

Of course we play this game {\it \`a contre coeur} as it is
against our philosophy that the phenomenon of total reality is
ubiquitous, and as posited by Rohlin 1978, that it governs the
saturation principle (alias Rohlin's maximality conjecture) saying
that a scheme of type~I is maximal in the hierarchy of all
schemes.

Note actually our logical MISTAKE, namely our project only
disproof the half of RMC already disproved by Shustin, as it will
exhibit a maximal $(M-1)$-scheme which is not of type~I. However
the harder game is to find a scheme of type~I which is not
maximal.

Let us however work out an example of this disproof strategy for
RMC as a potential application of the Rohlin formula and the
signs-law obstruction.

\begin{figure}[h]
\centering
\epsfig{figure=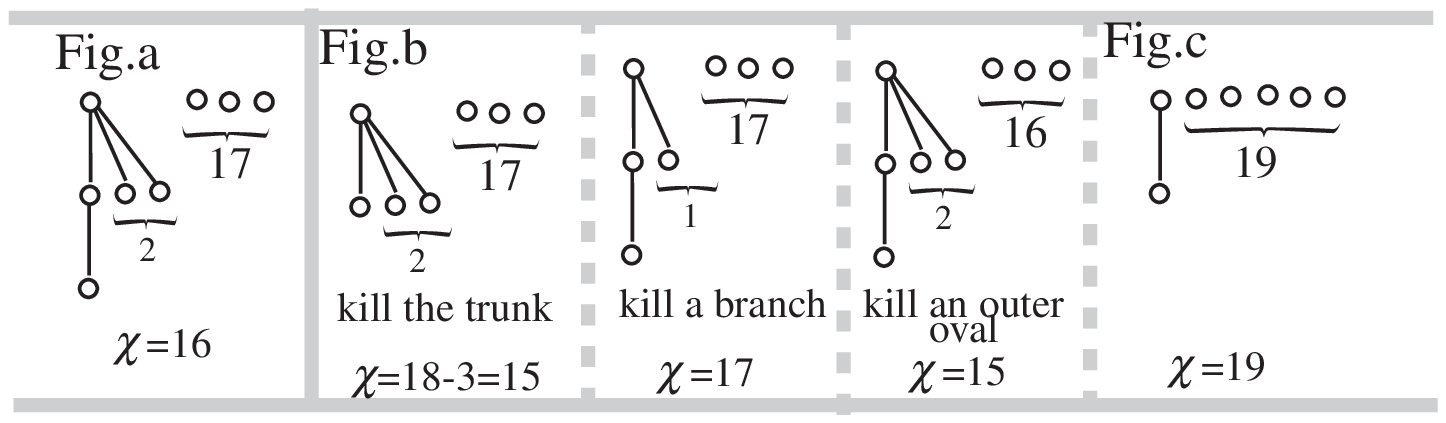,width=122mm} \vskip-5pt\penalty0
  \caption{\label{Fiedler8:fig}%
  Trying to corrupt Rohlin's maximality conjecture,
  but merely the sense already disproved by Shustin
  1985} \vskip-5pt\penalty0
\end{figure}

We start with an $M$-scheme which is prohibited by Rohlin though
not by Gudkov. Our (fairly random) candidate is, as said above,
the $M$-scheme in degree 8 given  by $(1,\frac{1}{1})$ (cf.
Fig.\,\ref{Fiedler7:fig}c or Fig.\,\ref{Fiedler8:fig}) with
$\chi=16$. The latter can be diminished to an $(M-1)$-schemes, 3
typical ways being depicted on Fig.\,b, where either the trunk is
killed, or a branch  or an outer oval. Of course one could also
kill the maximal nonempty oval like on Fig.\,c, but then the
GKK-congruence is not verified as $\chi=19$.

The 3rd specimen of Fig.\,b...

WARNING A THIS STAGE I HAD TO STOP AS I NOTICED AN EARLIER
MISTAKE!!!!

\section{E-mail correspondence}\label{e-mail-Viro:sec}

[09.01.13] This section gathers responses given by experts (Viro,
Marin, Orevkov, Kharlamov, Shustin, Le~Touz\'e, Fiedler, etc.) to
some naive questions of mine about the work of Rohlin. Here are
the original messages in chronological order (inserted with the
tacit approval of their authors).
I acknowledge most sincerely their authors for the stimulating
atmosphere it created and their generous answers.
Messages are left in their original shapes safe for adding
brackets [ ] supplying electronically-updated label-links to the
present text.

$\bullet$ On Wed, 9 Jan 2013 13:33:23 +0100
(alexandregabard@hotmail.com) wrote to Kharlamov, Marin, Viro,
Fiedler, Orevkov, and Mikhalkin a collective e-mail titled ``Two
naive questions on Rohlin 1978'':


Dear Viatcheslav, Alexis, Oleg, Thomas, Stepan and Grisha,

Sorry for disturbing so many experts among yours with some little
aspect of the work of academician Vladimir Abramovich. (I should
have written this message in French, yet cannot remember exactly
about Oleg's progresses over the last 6 years in that language.)

I was those last days quite fascinated by reading Rohlin's 1978
survey on complex topological characteristics of real curves in
some more detail. As you all know, he gave a quite spectacular
enhancement of Gudkov's pyramid for all schemes of sextics by
enriching it with the data of Klein's type~I/II (1876). (Compare
optionally Fig.71[=\ref{Gudkov-Table3:fig}] on page 208 of the
attached pdf file giving a graphical snapshot view of Rohlin's
achievement.)

My two questions are as follows.

(1) First Rohlin (1978) claims to have a certain synthetic
argument (via pencils of cubics) able to show the type~I of the
schemes 6/1 2 and 2/1 6. He confesses however his argument to be a
complicated one. Let me cite Rohlin exactly:

``...when we apply it to curves of degree 3, we can establish (in
a rather complicated way) that the schemes $\frac{6}{1}2$ and
$\frac{2}{1} 6$ of degree 6 belong to type~I. However, all the
schemes that we have so far succeeded in coping with by means of
these devices are covered by Theorem~3.4 and 3.5.''

My first question is whether Rohlin's synthetic argument has ever
been published (assuming its truth of course)? I suspect the proof
to be quite beautiful, but I am myself not quite able to write it
down for the moment. Did one of you ever worked out the argument
in detail, or remember about some exposition during Rohlin's
seminar? Is it of the same order of difficulty as the Hilbert-Rohn
method, requiring ``roughness'' \`a la Andronov-Pontrjagin to turn
round? Would it be didactically useful to publish (on the arXiv)
an account of Rohlin's argument if one is able to reconstruct it?
Many thanks if you have some ideas (or recent references) on those
or related questions...

(2) The second question is of course the general Rohlin's
maximality conjecture (a scheme is of type~I iff it is maximal in
the hierarchy of all real schemes of some fixed degree). As
reported in Viro's survey (1986 Progresses over the last 6 years)
it seems that one implication was disproved by Polotovskii and
Shustin (combined efforts ca. 1982, 1985). Yet one implication
looks still possible, namely ``type~I implies maximal'' (if I am
not wrong). It seems to me that this (last vestige of the) Rohlin
conjecture could be proved (somewhat eclectically) in two lines
via Ahlfors theorem (1950) on the total reality of orthosymmetric
curves (alias type~I). Namely if the curve is of type~I, then
there is a pencil of curves cutting only real points on the curve,
so its real scheme cannot be enlarged without violating B\'ezout.
q.e.d. Alas, some more thinking shows of course this argument to
be insufficient but maybe there is a (clever) way to complete it.
Qu'en pensez-vous?

Many thanks for your attention, and also for all your fantastic
papers (I am presently trying to digest). So do not take the pain
to answer me properly if my questions sound too naive. I apologize
again for this collective message, but as the material is quite
old, most of you probably forgot some details. So I hoped to
maximize some chance of getting an answer from a collective chat
room.

Best regards, Alex (Gabard)

PS: The attachment\footnote{I.e. the present text as it was on the
date of the 09.01.13, meanwhile pagination may have changed.} is a
copy of an informal text of mine on the Ahlfors map. Section 24
(pp. 205--229) is more specifically devoted to Rohlin's
conjecture, yet contains nothing original (except being poorly
organized).


$\bullet$ [Viro's answer the same day (09.01.13) ca. 20h00,
additional footnotes are mine (Gabard)]


Dear Alexandre,

Thank you for your message and manuscript. I was not aware about
the Ahlfors theorem\footnote{This is not perfectly true, as I
exposed Ahlfors theorem to Oleg Viro during its last visit in
Geneva (ca. 2010--11, his talk on fields of char 1, Connes,
tropical geometry, etc.), yet probably my explanations where so
obscure that Oleg immediately forgot about it.}. It seems to be
very interesting.

I doubt though if it can be used for proving the half of Rokhlin
conjecture. It gives a proof for impossibility of raising the
number of components of a type~I curve by a single algebraic Morse
modification (what I called Klein's thesis).

I do not remember if I even ever heard about Rokhlin's proof that
you ask about, but the fact follows from the congruence.
Slava\footnote{=Vladimir Abramovich Rohlin, of course.} did not
mention it when he proved the corresponding congruence (at the
moment the type was not yet considered). I learned this theorem
from Slava in September 1977 and wrote down Slava's proof to my
notebook then. I guess the first proofs was\footnote{Sic, singular
or plural? Not so important of course.} published by Slava Nikulin
(among many other statements) and Alexis Marin. Marin's proof
looks simpler, but requires Pin- structures.

Best regards, Oleg


$\bullet$ Gabard's reply [Same day (09.01.13) ca. 21h00]


Dear Oleg,

Many thanks for your rapid and illuminating responses, plus all
the historical details. If you see no objection, I would be very
happy to cut-and-paste them in my survey. I still need to
assimilate some congruences of the early phase (Rohlin,
Gudkov-Krakhnov-Kharlamov, etc.) Hence you cannot imagine how your
hints are illuminating my modest understanding of that golden
period. Regarding Ahlfors, as you say, there is little hope to
crack the big fish, yet of course I shall keep you informed if I
get not too depressed by the immense difficulty.

All the best, and so many thanks again, Alex


$\bullet$ 10 Jan 2013 (Marin's answer)

Cher Gabard

En plein d\'em\'enagement, je met un peu plus de temps \`a vous
r\'epondre que Viro.

Comme Viro, je ne connais pas la preuve de Rohlin pour votre
premi\`ere question (c'est pourquoi j'avais imagin\'e la preuve
dont parle Viro qui est dans "Quelques remarques sur les courbes
alg\'ebriques planes r\'eelle", votre r\'ef\'erence
742[=\cite{Marin_1979}]) Cependant ce s\'eminaire de Paris VII est
dans un carton et y restera tant que je n'aurai p\^u trouver un
nouvel appartement assez grand pour contenir ma biblioth\`eque et,
n'ayant le temps d'aller \`a la biblioth\`eque, ma m\'emoire ne me
permet pas de vous en dire plus que Viro.

Pour la seconde question par contre je peux vous r\'epondre, c'est
\`a dire lever votre aveux d'incompr\'ehension en fin (p. 226) de
preuve du Lemme 24.20\footnote{i.e. Lemma~\ref{Klein-Marin:lem},
warning meanwhile the numbering may changes, but the one in this
footnote is automatic (hence the right one)}.

Soit une courbe s\'eparante gagnant un ovale de plus apr\`es
franchiment d'un point quadratique ordinaire. Un argument de
congruence (utilisant $d > 2$ dans le cas plan ou une hypoth\`ese
dans le cas g\'en\'eral donnant que la d\'esingularis\'ee de cette
courbe de franchiment est irr\'eductible : l'ensemble de ses
points complexe est connexe) donne que cette d\'esingularis\'ee de
la courbe de franchiment est non s\'eparante.

Ainsi deux points non r\'eels conjugu\'es de la courbe de
franchiment sont li\'e par un arc \'evitant la partie r\'eelle, en
particulier le point singulier, et par extension des
isotipie\footnote{``isotopie'' of course.} un tel arc subsite dans
toute d\'eformation vers l'un des des deux c\^ot\'es du
discriminant, en particulier avant le franchiment la courbe est
non s\'eparante ce qui contredit l'hypoth\`ese.

Par contre si le franchiment du discriminant se fait en un point
singulier plus compliqu\'e il me semble que l'on peut augmenter le
nombre de composantes connexes d'une courbe s\'eparante. Je crois
me souvenir que selon les constructions de Viro (ou peut \^etre
seulement apr\`es avec la pr\'esentation Itenbergienne de cette
m\'ethode de Viro) il y a une courbe singuli\`ere de degr\'e 6
dont tout voisinage contient tous les types.

N'\'etant plus familier du sujet depuis plus de 20 ans je ne peux
vous en dire plus, par contre pour les surfaces de degr\'e 3 vous
trouverez dans le second tome des oeuvres de Klein un magnifique
article illustr\'e de non moins magnifiques figures o\`u il
\'etabli que tous les types de surface cubique s'obtiennent par
d\'eformation de la (unique \`a changement projectif de
coordonn\'ees) surface cubique qui a 4 points quadratiques
ordinaire.

Merci de votre long article que j'essayerai de lire quand
d\'em\'enagement, vente,.... seront termin\'es.

Bien cordialement et bonne ann\'ee.

Alexis Marin

PS 1 Je trouve Viro un peu "oublieux" d'\'ecrire " (at the moment
the type was not yet considered)": en parcourant le second tome
des oeuvres de Klein vous vous appercevrez qu'un sci\`ecle avant
Viro "tout" \'etait chez Klein!

2 Vous trouverez un article historique, beaucoup plus court* et
sur un autre sujet en mettant dans la boite de recherche d'Arxiv
le mot clef ``troupeau''

*il fait 6 pages table des mati\`eres comprise et tout est dit (de
fa\c{c}on ``autocontenue'') dans le r\'esum\'e en fran\c{c}ais de
la premi\`ere page, mais si vous remontez \`a toutes les
r\'ef\'erences** des commentaires bibliographiques cel\`a peut
vous prendre un peu de temps.

**accesibles \`a travers la "biblioth\`eque des sophomores
http://alexis.marin.free.fr/BIB/

$\bullet$ Gabard's answer [12.01.13 ca. 23h00]

Cher Alexis, C'est avec une immense joie que j'ai re\c{c}u votre
message. N'ayant pas d'internet \`a la maison, je l'ai seulement
d\'ecouvert ce soir en visitant mon p\`ere, qui lui est
connect\'e. Je vais donc tenter d'assimiler toutes vos remarques
savantes, et si vous le permettez, de les int\'egrer dans mon
survey, en sp\'ecifiant bien s\^ur qu'il s'agit de vos
contributions. De mon c\^ot\'e, je me demande si une courbe
non-s\'eparante peut toujours acqu\'erir un point double ordinaire
solitaire. (C'est semi-implicite dans Klein 1876 qui \'ecrivait
``noch entwicklungsf\"ahig'', mais il me semble que \c{c}a
contredit le r\'esultat de Shustin 1985 (contre-exemple \`a la
conjecture de Rohlin), dont la logique m'\'echappe quelque peu,
mais j'ai s\^urement rat\'e une subtilit\'e).

Gr\^ace \`a vos commentaires je devrais pouvoir produire
prochainement une version plus solide et limpide de la section
correspondante du survey, que je vous enverrai d\`es que possible.
L'interaction avec Ahlfors me semble aussi prometteuse...

Amiti\'es, et merci infiniment pour votre message, Alex

PS J'esp\`ere que le d\'em\'enagement se passe bien. Restez-vous
\`a Grenoble, ou bien s'agit-il d'une op\'eration plus
cons\'equente?

PPS: J'ai bien \`a la maison votre article de Paris VII, qui a
toujours \'et\'e mon meilleur compagnon (en 1999-2000), et je suis
content de le retrouver pour ce point encore plus profond.

PPPS: je me suis procur\'e une copie de l'article sur ``il capo'',
qui me semble fabuleux. Merci beaucoup. C'est exactement l'analyse
que l'on rencontre \`a proximit\'e de Dirichlet, etc jusqu'\`a
Ahlfors, et Rogosinski, et que je dois essayer \`a l'occasion
d'apprivoiser...

PS 1 Je trouve Viro un peu "oublieux" d'\'ecrire " (at the moment
the type was not yet considered)": en parcourant le second tome
des oeuvres de Klein vous vous appercevrez qu'un sci\`ecle avant
Viro "tout" \'etait chez Klein!

Vous avez parfaitement raison, et je suis moi m\^eme tr\`es
``sp\'ecialis\'e '' dans l'oeuvre de Klein. Cependant le gros
quiz, c'est l'assertion de Teichm\"uller 1941, qui pr\'etend que
Klein 18XX? anticipe Ahlfors 1950, de 70 ans environ. Toute courbe
s\'eparante (ou surface de Riemann orthosymm\'etrique, pour
reprendre le jargon klein\'een) admet un morphisme r\'eel vers la
droite dont les fibres au dessus des points r\'eels sont toutes
exclusivement form\'ees de points r\'eels. C'est cet \'enonc\'e
fondamental qui me semble \^etre sous-exploit\'e!
Evidemment\footnote{[26.03.13] This ``evidently'' sounds to me a
bit sloppy, as the full credit for this remark goes to Viro.}
comme la not\'e Viro, il implique la partie facile de l'assertion
de Klein (1876): une courbe s\'eparante ne peut gagner un ovale
spontan\'ement comme une bulle de champagne surgit du n\'eant.

$\bullet$ R\'eponse de Marin (le lendemain 13 Jan 2013 ca. 09h00)

de les int\'egrer dans mon survey, en sp\'ecifiant bien s\^ur
qu'il s'agit de vos contributions.

A part l'explication de votre doute (o\`u relativement \`a
l'article que vous citez il n'y a que les mots ``extension des
isotopies'' en plus) ce ne sont que de tr\`es vagues souvenirs que
je vous conseille de v\'erifier (\'eventuellement aupr\`es de plus
comp\'etent : Viro, Itenberg,... avant de les int\'egrer)

De mon c\^ot\'e, je me demande si une courbe non-s\'eparante peut
toujours acqu\'erir un point double ordinaire solitaire.

voulez-vous dire dont les deux directions tangentes sont complexes
conjugu\'ee? cel\`a me parait tr\`es tr\`es optimiste.

(C'est semi-implicite dans Klein 1876 qui \'ecrivait ``noch
entwicklungsf\"ahig'',

\^Etes vous s\^ur que c'est ce que pensait Klein, ou incluait-il
dans ce terme les modification par franchiment d'une courbe ayant
un unique point double qui est quadratique ordinaire \`a tangentes
r\'eelles ``apparu en rapprochant deux points d'un m\^eme ovale''?

mais il me semble que \c{c}a contredit le r\'esultat de Shustin
1985 (contre-exemple \`a la conjecture de Rohlin), dont la logique
m'\'echappe quelque peu, mais j'ai s\^urement rat\'e une
subtilit\'e).

PS J'esp\`ere que le d\'em\'enagement se passe bien.

oui mais c'est long, \`a ce propos, vous trouverez sur

http://alexis.marin.free.fr/BIB/papier/

la liste des livres que j'ai en plusieurs exemplaires et (sauf
ceux dont la colonne "h\'eritier" est remplie (par Vinel et/ou
Guillou)) qui sont \`a la disposition de qui (en particulier vous)
les demande. Restez-vous \`a Grenoble, ou bien s'agit-il d'une
op\'eration plus cons\'equente?

Je reste pr\`es de Grenoble (mon adresse est dans la signature
\'electronique ci-dessous

PPPS: je me suis procur\'e une copie de l'article sur ``il capo'',

Voulez vous dire "Le capo"?

Cependant le gros quiz, c'est l'assertion de Teichm\"uller 1941,
qui pr\'etend que Klein 18XX? anticipe Ahlfors 1950, de 70 ans
environ. Toute courbe s\'eparante (ou surface de Riemann
orthosymm\'etrique, pour reprendre le jargon klein\'een) admet un
morphisme r\'eel vers la droite dont les fibres au dessus des
points r\'eels sont toutes exclusivement form\'ees de points
r\'eels.

Voulez-vous dire rev\^etement d'espace total l'ensemble des
ovales? Il y a-t-il quelque chose de plus pr\'ecis sur le degr\'e
et sa r\'epartition parmis les ovales? Les r\'ef\'erences
sont-elles dans votre article?

C'est cet \'enonc\'e fondamental qui me semble \^etre
sous-exploit\'e! Evidemment comme la not\'e Viro, il implique la
partie facile de l'assertion de Klein (1876): une courbe
s\'eparante ne peut gagner un ovale spontan\'ement comme une bulle
de champagne surgit du n\'eant.

Soyez plus pr\'ecis pourquoi un tel morphisme admettrait-il une
d\'ef\-ormation le long de la modification d'adjonction d'un
ovale?

Amiti\'es.

Alexis

-- http://le-tonneau-de-thales.tumblr.com/

Alexis Marin, chez Danielle Bozonat 6 All\'ee de la roseraie,
38240 Meylan fixe : 04 76 00 96 54 port. : 06 38 29 33 99,
00351925 271 040

$\bullet$ Gabard 13 Jan 2013 ca. 13h30

Cher Alexis, Merci pour votre message. Je vais en effet essayer
d'int\'egrer vos commentaires de mani\`ere cibl\'ee et prudente.
De toute mani\`ere avant d'arXiver une nouvelle version d'ici six
mois environ, j'aurai l'occasion de vous montrer pr\'ecisement la
prose que je vous aurez emprunt\'e. J'essaye maintenant de
r\'epondre \`a vos questions:

De mon c\^ot\'e, je me demande si une courbe non-s\'eparante peut
toujours acqu\'erir un point double ordinaire solitaire.

voulez-vous dire dont les deux directions tangentes sont complexes
conjugu\'ee? cel\`a me parait tr\`es tr\`es optimiste.

REPONSE: Oui, exactement \`a tangentes imaginaires conjugu\'ees.
Cela me parait aussi tr\`es optimiste. Klein semble le pr\'etendre
semi-implicitement (du moins qu'il n' y a a priori pas
d'obstruction topologique \`a la formation de telles bulles de
champagne). Cependant si ce truc fou (``Klein-vache'') est vrai
alors un des sens de la conjecture de Rohlin 1978 (type~I iff
maximal real scheme) est v\'erifi\'e. Malheureusement, ce que
donne ``Klein-vache'' est le sens de Rohlin d\'etruit par Shustin
1985 (dont je n'ai cependant pas compris l'argument). Mais vous
avez surement raison ``Klein-vache'' est probablement beaucoup
trop optimiste...

\^Etes vous s\^ur que c'est ce que pensait Klein, ou incluait-il
dans ce terme les modification par franchiment d'une courbe ayant
un unique point double qui est quadratique ordinaire \`a tangentes
r\'eelles ``apparu en rapprochant deux points d'un m\^eme ovale''?

REPONSE: je pense que oui, car Klein pr\'ecise ``isolierte reelle
Doppeltangente'', comparez ma Quote 24.2\footnote{Meanwhile this
numbering may have changed into
Quote~\ref{Klein_1876-niemals-isolierte:quote}.} page 205 de mon
survey (si vous n'avez pas le volume 2 de Klein sous la main).
Ainsi il me semble que votre interpr\'etation moderne (Marin 1988)
diff\`ere un peu de l'original Klein\'een, en \'etant toutefois
plus puissant que l'assertion d'origine.

PS J'esp\`ere que le d\'em\'enagement se passe bien.

oui mais c'est long, \`a ce propos, vous trouverez sur

http://alexis.marin.free.fr/BIB/papier/

la liste des livres que j'ai en plusieurs exemplaires et (sauf
ceux dont la colonne "h\'eritier" est remplie (par Vinel et/ou
Guillou)) qui sont \`a la disposition de qui (en particulier vous)
les demande.

C'est une magnifique liste de tr\'esor. Je voudrais bien les
acqu\'erir, mais je me demande si mon hygi\`ene de vie (overwork)
rend une telle acquisition raisonable...(Il faudrait que je passe
\`a Grenoble avec la camionnette de mon oncle pour r\'ecup\'erer
les ``invendus''. Il est pr\'ef\'erable en effet de trouver des
preneurs plus comp\'etents que moi. Si en dernier recours, vous ne
trouvez pas de preneurs je pourrais r\'ecup\'erer les volumes
restants en vrac...Merci infiniment pour cette g\'en\'ereuse
proposition. Moi m\^eme je suis tr\`es marginal financi\`erement
et spatialement, petit appartement \`a Gen\`eve partag\'e avec ma
m\`ere (avec environ 8 tonnes de litt\'erature math\'ematique),
mais dans le futur je pourrai peut \^etre m'installer dans une
ferme fribourgoise, o\`u il reste de l'espace pour expandre la
biblioth\`eque...)

PPPS: je me suis procur\'e une copie de l'article sur ``il capo'',

Voulez vous dire "Le capo"? Oui, j'essayais d'improviser en
italien, mais c'est une langue plus subtil que vous utilisez...

Cependant le gros quiz, c'est l'assertion de Teichm\"uller 1941,
qui pr\'etend que Klein 18XX? anticipe Ahlfors 1950, de 70 ans
environ. Toute courbe s\'eparante (ou surface de Riemann
orthosymm\'etrique, pour reprendre le jargon klein\'een) admet un
morphisme r\'eel vers la droite dont les fibres au dessus des
points r\'eels sont toutes exclusivement form\'ees de points
r\'eels.

Voulez- vous dire rev\^etement d'espace total l'ensemble des
ovales? Il y a-t-il quelque chose de plus pr\'ecis sur le degr\'e
et sa r\'epartition parmis les ovales? Les r\'ef\'erences
sont-elles dans votre article?

OUI, toute surface de Riemann \`a bord (=membrane compacte)
s'exprime comme rev\^etement holomorphe ramifi\'e du disque. C'est
juste une version relative (\`a bord) du th\'eor\`eme d'existence
de Riemann qui concr\`etise toute surface de Riemann close comme
revetement conforme de la sph\`ere (ronde). Il a fallut toutefois
attendre la contribution d'Ahlfors 1950 qui donne en plus un
contr\^ole sur le degr\'e d'un tel rev\^etement conforme, \`a
savoir r+2p, o\`u r est le nombre d'``ovales'' (mieux le nombre de
contours de la membrane), et p son genre. La Th\`ese de moi-m\^eme
(Gabard 2004, et l'article de 2006 au Commentarii Math. Helv.)
donne un meilleur contr\^ole, \`a savoir $r+p$, en \'economisant
donc une cartouche pour chaque anse. Les r\'ef\'erences pr\'ecises
sont dans le survey. L'\'enonc\'e d'Ahlfors \'etait vachement
anticip\'e dans le cas $p=0$ (membrane planaire ou schlichtartig
pour reprendre la terminologie de Paul Koebe) par la grande
lign\'ee Riemann 1857 (Nachlass), Schottky 1875-77, Bieberbach
1925 et Grunsky 1937. Lorsqu'on passe au double de Schottky-Klein
de la surface \`a bord on obtient (via Ahlfors) une courbe
s\'eparante avec un morphisme totalement r\'eel vers la droite
projective. Inversement toute courbe s\'eparante est totalement
r\'eelle, puisqu'il suffit d'appliquer Ahlfors \`a une des
moiti\'es orthosym\'etrique de Klein.

[Gabard] C'est cet \'enonc\'e fondamental qui me semble \^etre
sous-exploit\'e! Evidemment comme la not\'e Viro, il implique la
partie facile de l'assertion de Klein (1876): une courbe
s\'eparante ne peut gagner un ovale spontan\'ement comme une bulle
de champagne surgit du n\'eant.

[Marin] Soyez plus pr\'ecis pourquoi un tel morphisme
admettrait-il une d\'eform\-ation le long de la modification
d'adjonction d'un ovale?

[Gabard] Je pense que \c{c}a marche car lorsque la courbe est
plong\'ee dans le plan, le morphisme total d'Ahlfors admet une
r\'ealisation projective comme un pinceau de courbes planes dont
tous les membres d\'ecoupent seulement des points r\'eels sur la
courbe orthosymm\'etrique (=s\'eparante). Par cons\'equent, en
tra\c{c}ant la courbe du pinceau total qui passe par un point de
l'oval spontan\'ement cr\'e\'e, on obtient une contradiction avec
B\'ezout. Donc Ahlfors 1950 implique Klein 1876, mais votre
d\'emonstration de 1988$-\varepsilon$ (votre preuve est d\'ej\`a
mentionn\'ee dans Viro 1986) est surement plus intrins\`eque et
voisine de l'argument d'origine de Klein (s'il en avait un au
del\`a de la pure contemplation empirique des quartiques
notamment...)

Merci infiniment pour vos messages, et d'ici tout bient\^ot (3-4
jours) je vous enverrai une version mise-\`a-jour du survey qui
clarifiera peut-\^etre les assertions pr\'ec\'edentes. Toutefois
les grands probl\`emes et plein de d\'etails m'\'echappent encore
dans la pyramide Gudkovo-Rohlinienne. Quelle splendide pyramide
qui joint \`a la perfection Klein et Hilbert! Un d\'etail qui
m'\'echappe, c'est le fait que le discriminant est de degr\'e
$3(m-1)^2=75$ pour $m=6$, tandis que que du point  de vue des
chirurgies ``de Morse'' il y a des cycles de longueur 4 dans la
pyramide de Gudkov. Donc il y a un probl\`eme de parit\'e si on
d\'eforme le long d'un pinceau g\'en\'erique (transverse au
discriminant)...D\'esol\'e, de vous emb\^eter avec ces d\'etails
que j'ai honte de ne pas r\'eussir \`a clarifier depuis quelques
jours.

Amiti\'es, et bon courage pour la suite du d\'em\'enagement, Alex

$\bullet$ [16h40 15.01.13] Cher Alexis, Merci encore pour vos
messages et vos remarques fascinantes que je dois encore bien
dig\'erer. De mon c\^ot\'e, j'ai fait de minimes progr\`es, et
vous envoie malgr\'e votre d\'em\'enagement une version ajourn\'ee
de mon survey. Il me semble que le truc fou dont nous parlions il
y a quelques jours, que j'appele depuis ``Klein-vache'', i.e. la
possiblilit\'e de faire naitre un noeud solitaire (\`a tangentes
conjugu\'ees) depuis n'importe quelle courbe diasym\'etrique est
vrai pour les sextiques. Pour cela j'utilise un argument qui
combine Rohlin 1978, Klein-Marin 1988, et Nikulin 1979
(classification isotopique) et un r\'esultat reli\'e de Itenberg
1994 (possibilit\'e de contracter n'importe quel ovale vide, i.e.
sans autre ovales dans son int\'erieur, sur un tel noeud isol\'e).
Les d\'etails de la preuve sont expos\'es dans la
Prop.24.24[meanwhile this is \ref{Klein-vache-deg-6:prop}], page
235 du fichier ci-joint. ? Evidemment, en principe ``Klein-vache''
n'a aucune chance d'\^etre vrai en degr\'e sup\'erieur. Cependant
la seule obstruction que je connaisse est ce r\'esultat de Shustin
1985, dont je ne comprends toujours pas la logique de base (sans
m\^eme parler du fait que c'est fond\'e sur la m\'ethode de Viro,
dissipation de singularit\'es tacnodales..., une technologie que
je n'ai jamais maitris\'ee). Mes objections naives \`a l'argument
de Shustin se trouvent en page 248 (dans le paragraphe qui
pr\'ecède la Figure 94[=meanwile Fig.\,\ref{Shustin:fig}]). Dans
cette figure, je ne sais pas comment prohiber le $(M-1)$-sch\'emas
encadr\'e par le carr\'e vert (\`a mi-hauteur de la figure), et
dans son article de 1985 Shustin n'est pas tr\'es explicite. Mais
bon, il s'agit la d'une question assez ennuyeuse et en fait je
vais peut-\^etre prendre l'initiative d'\'ecrire un nouveau
message collectif pour clarifier ce point d'ici quelques heures.
Merci infiniment encore pour vos messages, et meilleurs voeux de
courage pour la suite du d\'em\'enagement, Amiti\'es, Alex PS:
Pour l'instant j'ai in\'egr\'e en vrac tous nos \'echanges e-mail
dans le survey (p.219 et suivante), mais bien entendu d\`es que
possible je censurerai les remarques plus confidentielles..., et
masquerai les r\'ep\'et\'etitions, voire l'int\'egralit\'e de la
discussion si je parviens bien \`a r\'esumer votre apport malgr\'e
mon anglais catastrophique. Cependant en relisant vos remarques,
elles apportent une prose substantielle que je ne saurais jamais
reproduire en anglais, donc je trouverais tr\`es dommage de
censurer vos souvenirs en vracs!!! Evidemment rien ne presse et je
suis d\'esol\'e de vous avoir d\'erang\'e durant cette d\'elicate
op\'eration du d\'em\'enagement inter-grenoblois. Amiti\'es,
encore, et je vous tiens au courant d'\'eventuelles progr\`es...Je
suis surtout curieux des r\'eponses de Shustin (et Viro) s'ils
parviennent \`a \'eclairer ma lanterne. PPS: Je joins une copie de
la note de Shustin, si jamais, mais je ne veux pas vous distraire
de votre t\^ache prioritaire...

$\bullet$ [15.01.13--18h30] Dear Evgenii, Ilia, Oleg and Alexis
(and Felix Klein),

I was much fascinated those last days by Evgenii's counterexample
to (one part of) Rohlin's maximality conjecture to the effect that
a real scheme is of type~I iff it is maximal in the hierarchy of
all schemes. Quite interestingly this work of you (Shustin) also
destroys an old (semi-)conjecture of Klein (1876) positing that
any nondividing plane curve can acquire a solitary node by
crossing only once the discriminant (the resulting Morse surgery
then sembling like the formation a champagne bubble arising like a
blue sky catastrophe of little green men's coming with flying
saucers).

Alas from Shustin's note of 1985 (in its English translation), I
was not quite able to understand your proof (compare optionally
the attached file, on page 248, in the paragraph right before
Figure 94[=Fig.\,\ref{Shustin:fig}]). In fact I do not know how to
prohibit the $(M-1)$-scheme $4/1 2/1 1/1 11$ enlarging Shustin's
(M-2)-scheme. Alas I am not an expert in the field and I feel
quite shameful disturbing you with such a detail. Despite having
myself full Leningradian origins (through my father), I do not
master the Russian language so that it may well be the case that
the original Russian text is more detailed than its translation.
Of course it is much more likely that I missed something
well-known, that you perhaps may not have made completely explicit
in the note? (Incidentally I send you a copy of Shustin's note for
convenience!)

I apologize for this question of detail, yet it seems quite
important to me for your result of 1985 is the only obstruction (I
am aware of) to the naive desideratum of truth about Klein's
conjecture. Klein himself is extremely clever and quite ambiguous
about stating this as a conjecture or as a result (compare
optionally Klein's original quote reproduced on page 206 of the
attachment). Today I managed as a simple exercise to check the
truth of Klein's hypothesis in degree 6, via an armada of Russian
results (especially Itenberg 1994 contraction principle for empty
ovals), plus the Klein-Marin theorem (for the details of this
exercise cf. optionally Prop.24.24[=\ref{Klein-vache-deg-6:prop}]
on page 235 of the attached text).

You, Oleg Viro, in the preface of that volume presenting
Itenberg's article (1994) advanced the (crazy?) conjecture that
one might always be able to contract empty ovals!!! Do you know if
there is meanwhile some counterexample (in high degrees)? Of
course there is some vague parallelism between Itenberg's
contraction and the one required to implement Klein's hypothesis
(which must amount shrinking an anti-oval, i.e. an invariant
circle acted upon antipodically by conj).

Sorry again for disturbing you with all these naive questions, and
do not take the pain answering me properly if you are overwhelmed
by other more important duties. Many thanks for all your
attention. Sincerely yours, Alex (Gabard)

$\bullet$$\bullet$$\bullet$ [16.01.13--02h57: Oleg Viro]

Dear Alexandre,

I do not mind to pose crazy conjectures. I do not mind if my crazy
conjecture would be disproved.

However, I suspect that my conjecture is not as crazy as
possibility of shrinking of an anti-oval. The difference between
the oval and an anti-oval is that the oval is assumed to exist and
be empty, i.e., not linked with the complex curve in whatever
sense, while the anti-oval apparently has none of these
properties.

I am not aware about any counter-examples that you ask about. I do
not bet that they do not exist, but find the question stimulating,
and better motivated than the conjecture that was proven to be
wrong.

Best, Oleg

$\bullet$$\bullet$$\bullet$ [16.01.13--14h56: Stepan Orevkov]

A small remark:

It is wrong that $11 U 1<1> U 1<2> U 1<4>$ is not a part of an
$(M-1)$-scheme. It is\footnote{Not clear how to interpret this?
Does it mean that Shustin's claim is wrong, or simply that this
scheme is an $(M-1)$-scheme. My question  was whether this
$(M_1)$-scheme is realized algebraically, of course. Yet, I admit
that my question was a bit ill posed.}. Moreover, there is no
known example of $(M-2)$-curve of type~II which cannot be obtained
from an $(M-1)$-curve by removing an empty oval.

In contrary, there are $(M-1)$-curves of degree $8$ (which are
necessarily of type~II) which do not come from any $M$-curve.
These are:

$3<6>$

$4 U 1<2> U 2<6>$

$8 U 2<2> U 1<6>$

$12 U 3<2>$

Constru[r]ction (inspired by Shustin's construction of $4 U
3<5>$):

Consire[der] a tricuspidal quartic $Q_{sing}$ symmetric by a
rotation $R$ by $120$ degree and perturb[e] is[=it] so that each
cusp gives an oval (we assume that this perturbation is very
small). Let $Q$ be the perturbed curve. Two flex points appear on
$Q$ near each cusp of $Q_{sing}$. We chose flex points $p_0, p_1,
p_2$ (one flex point near each cusp) so that $R(p_0)=p_1,
R(p_1)=p_2, R(p2)=p_0$. We choose homogeneous coordinates $(x_0 :
x_1 : x_2)$ so that the line $x_i = 0$ is tangent to $Q$ at $p_i$
$(i = 0,1,2)$.

Let $C$ be the image of $Q$ under the Cremona transformation $(x_0
: x_1 : x_2) \mapsto (x_1x_2 : x_2x_0 : x_0x_1)$. Then $C$ has 3
singular points, each singular point has two irreducible local
branches: a branch with $E6$ and a smooth branch which cuts it
``transversally''. By a perturbation of $C$ we obtain all the four
curves mentioned above.

The fact that these curves cannot be obtained from $M$-curves
immediately follows from the fact that, for any $M$-curve of
degree 8 of the form $b U 1<a_1> U 1<a_2> U 1<a_3>$, all the
numbers $a_1$, $a_2$, $a_3$ are odd\footnote{[24.01.13] On reading
Viro' survey 1989/90 \cite[p.\,1085,
2.5.H]{Viro_1989/90-Construction}, one should easily locate the
source for this assertion. References seems to be Viro 1983
\cite{Viro_1983/84-new-prohibitions}, and the survey Viro 1986
\cite{Viro_1986/86-Progress}. As explained there (Viro \loccit
1989/90, p.\,1085) this is a prohibition not coming from topology,
but from B\'ezout. In fact this result is mentioned again in Viro
1989/90 \cite[p.\,1126, 5.3.E]{Viro_1989/90-Construction}, with
the exact cross-reference as being Viro 1983
\cite{Viro_1983/84-new-prohibitions}}.

Best regards Stepa O

$\bullet$ [17.01.13 ca. 23h00]

Dear Oleg and Stepa,

Many thanks for all your fascinating remarks and detailed answers.
I look forward digesting them carefully tomorrow.

Sorry for my late reply as I have no internet at home and was
quite busy trying to understand some basic facts, notably that one
may have some ``eversion'' of a real scheme when the oval explodes
at infinity undergoing a Morse surgery not affecting its
connectedness. This implies that there is some hidden passages in
the Gudkov-Rohlin pyramid of all sextics changing a Gudkov symbol
$k/l \ell$ to its mirror $\ell/1 k$. The resulting combinatorics
of this graph looks quite formidable and I wonder if it is known
whether each of those secret edges corresponding to eversions
(except those linking $M$-curves) can be explored algebraically.
Perhaps the problem is related to Ilia's shrinking process for
empty ovals, but seems to involve yet another species of
``anti-ovals'', namely those with two fixed points under conj, yet
located on the same oval. All what I am saying is for sure
well-known to you since time immemorial, yet I was very happy to
understand this point which solved several paradoxes of mine,
notably those related to the degree of the discriminant and the
contiguity graph between chambers residual to the discriminant
under elementary algebraic Morse surgeries, as Oleg says. Of
course, I shall send you an updated version of my file, when I
manage to reorganize slightly the exposition.

Many many thanks for all your excellent answers! All the best,
Alex

$\bullet$ [18.01.13 ca. 10h00, Viatcheslav Kharlamov]

Dear Alex,

I followed rather attentively the discussions, but kept silence
since had no much to add to the reaction of the others.

This ``eversion'', as you call it, played some important r\^ole in
the prehistory of the Gudkov conjecture. As you probably know, the
first classification declared by Gudkov was wrong, and it is one
of his "thesis referees", Prof. Morosov, who had objected the
first classification exactly because of a small irregularity with
respect to ``eversion'' of the answer. Repairing this asymmetry
Gudkov came to his final result, and, if my memory is correct,  in
particular, at this stage discovered the missing $M$-curve.

If honestly, I don't remember did somebody ever before discussed
seriously any conceptual explanation to this ``eversion''.
However, it was implicitly present in all results obtained through
$K3$ and their lattices. Recently, studying the shadows of cubic
surfaces with Sergey Finashin and having proven, to our own
surprise, for them a very similar ``symmetry'', which we have
called ``partners relation'', we have formalized it as follows.

First level of explanation is coming from lattices of double
coverings: the partner relation consists indeed in transferring an
$U$-summand (unimodular even lattice of rang $2$ and signature
$0$) from one eighenlattice to another.  Second level of
explanation is coming from moduli in terms of periods: each
partner in the partner pair can be deformed to a triple conic,
near the triple conic the family looks as $Q^3+tbQ^2+t^2cQ+d=0$,
and switching of the sign of $t$ (passing through the triple
conic) replace curves of one deformation class by curves from the
partners class: moreover, such degenerations are deformationally
unique. Literally the same  explanation (and with much easier
proofs at the both levels) works for nonsingular sextics (the
shadows are sextic curves with $6$ cusps on a conic; remarkably,
in many respects they behave in a way more similar to that of
nonsingular sextics, than other sextics with singularities).

Yours, Viatcheslav Kharlamov

$\bullet$ [18.01.13, Kharlamov, title of message=Correction]
Writing the message a bit in a hurry I did not describe fully and
appropriately the partner relation at the lattice language. The
summand $U$ does play a crucial r\^ole, and it should be moved
from one eighenlattice to another, but then additionally one
should exchange the eighenlattices. In fact this $U$ contains
indeed the $2$-polarization vector, $h¨2=2$, and thus the
eighenlattice containing this distinguished U is aways
$(-1)$-eighenlattice.

The existing exception to the partner relation (as I remember, in
the nonsingular case, there is only one) is the case when the
$(-1)$-eighenlattice does not contain such a pair $(U,h)$.

Sorry, for being in a hurry, but I should stop at this point. Hope
that now it is more clear.

$\bullet$ [18.01.13, Gabard, ca. 21h00] Dear Viatcheslav, Oleg and
Stepa (and all the others),

So many thanks for all the excellent comments, especially on
Morosov. There was some allusion to this issue in Viro's survey
from 2006, in Japanese Journal of Math, as to the lack of symmetry
in Gudkov's initial answer. Yet Morosov was not mentioned if I
remember well...

On my side I was quite stimulated by the last letter from Oleg,
about the contraction conjecture, as looking indeed much more
realist than Klein's Ansatz on the champagne bubbling in any
nondividing curve. I attempted today to imagine what sort of proof
one could expect to find for this fascinating Itenberg-Viro
contraction conjecture of empty ovals.

After some trials with orthogonal trajectories to the functional
computing the area of the empty oval, I arrived at some sort of
strategy (probably completely fantasist) consisting in using the
Riemann mapping theorem as applied to the interior of the empty
oval. Naively as the contour is algebraic so is the Riemann map
and hence the concentric sublevels of it ought to be algebraic
curves of the same degree!!??? This would give the shrinking.

I am sure that tomorrow while checking more carefully the details
all this argument will crash down. Hence sorry for this premature
message. Some more details about this and my naive understanding
of ``eversions'' are in the attached file, especially Section
24.15 (p.251) and p.242 (Section 24.12 for eversions). Regarding
eversions I wonder which edges in the Gudkov pyramid are actually
realized algebro-geometrically? All, except those connecting the
$M$-schemes is my naive guess, yet it is probably too
optimistic...

Many thanks again for sharing all your knowledge on that
fascinating topic, and all your exciting letters. All the best,
Alex

[21.01.13, ca. 20h00]

Dear real geometers,

Thank you again, Oleg, Alexis, Stepa, and Viatcheslav, for all
your messages which I have carefully integrated in my TeX-notes,
and to which I frequently refer for citation in my text. Your
messages suggested me several ideas I would never have explored
without your precious hints.

On my side, I noticed of course that the cavalier Riemann mapping
strategy toward the (Itenberg-Viro 1994) contraction conjecture
(CC) of empty ovals fails blatantly (cf. Section
25.7(=\ref{CC-via-Riemann:sec}), pages 255-258, roughly even if
the Riemann map of an algebraic oval would be algebraic then its
degree seems to be twice as big as it should, or better the
polynomials arising as norms of algebraic Riemann maps are not the
most general representatives of their degree!!!). Perhaps the
Riemann method works for special ovals, but of course they are
unlikely to be interspersed in all chambers of the discriminant!
This failure drifted me toward another formulation of the
contraction conjecture which I call CCC, for collective
contraction conjecture. This posits that all empty ovals of a real
algebraic curve can be contracted simultaneously toward solitary
nodes (by a path having solely its end-point in the discriminant).
This looks even more ``crazy'' than CC, but I found no
counterexamples (in my pockets). I would much appreciate if you
already thought about this natural variant, especially if you
detected some counterexample (perhaps arising from the
Viro-Itenberg patchworking method or the dissipation of higher
singularities, with which I am alas still unfamiliar with, like in
Shustin's counterexample to Rohlin's maximality conjecture).

Here are the trivialities I managed to prove. Via Brusotti 1921,
it is plain that CCC implies the usual contraction conjecture (CC)
(cf. details in Lemma 25.22(=\ref{CCviaCCC-Brusotti:lem}) on page
263 of the attached file). On the other hand CCC implies (as a
large deformation principle) several well-known prohibitions. E.g.
a two-seconds proof of the Hilbert-Rohn-Petrovskii prohibition of
the sextic $M$-scheme $11$ (eleven ovals without nesting), as well
as Rohlin's prohibition of the sextic scheme $5$ of type~I (by the
way causing the unique asymmetry in the Gudkov-Rohlin table of
sextics). Under CCC, all these facts appear as trivial
consequences of B\'ezout (compare Section 25.8(=\ref{CCC:sec}) on
pages 258-259). I found this simplicity quite exciting (even
though it leads to nothing new as compared to Arnold-Rohlin). One
can wonder if Hilbert already used this, at least as a heuristic
tool??? Philosophically, I found also interesting that such large
deformation conjectures produce prohibitions, in contradistinction
to small perturbations as being primarily a method of construction
(Harnack-Hilbert-Brusotti, etc.). There is accordingly some nice
duality between Luigi Brusotti and Ilia-Oleg's contraction
conjecture. Of course you surely noted this issue a long time ago,
yet for me it was a happy discovery (yesterday).

Perhaps CCC and CC are actually equivalent, yet this looks more
hazardous but maybe not completely improbable... (One would just
have to synchronize the death of all ovals posited by CC.) This is
all the modest news I have collected during the week-end. Of
course I still have some naive hope that CCC (hence CC) could be
attacked via some gradient flow, but it looks quite difficult to
locate the right functional (or Morse function). Looking at the
area (or length) of all empty ovals is probably too
naive...Perhaps some ``degree of roughness'' \`a la Gudkov could
be projectively more  intrinsic and useful...

Thank you so much for your attention and all your brilliant
letters and answers, while apologizing me for sending you only
easy doodlings. Best regards, Alex

[26.01.13, ca 20h00]

Dear Oleg, Stepa, Viatcheslav, Evgenii, Alexis, Thomas, etc.

I continued my naive investigations of real plane curves. What a
beautiful story! I finally ``understood'' and studied in detail
the marvellous construction of Gudkov $\frac{5}{1} 5$ (ca.
1971-73), as to understand the more tricky (but related)
construction proposed by Stepa, which I attempted to depict on
Fig. 111(=\ref{Orevkov2:fig}) of the attached file. (I did not as
yet assimilated the full details but feel on the good way. In fact
I tried to use the dissipation of $Z_15$ in Viro's survey from
1989/90 in Leningrad Math. J., which I hope is the same as the
$E_6$ advocated by Stepa. Sorry for being very ignorant about
singularities...)

Yesterday, I also finally understood the correctedness of
Evgenii's argument. (As helped by Stepa's e-mail, the point which
I missed is this obstruction of Viro extending that of Fiedler)
for $M$-schemes of degree 8 as having necessarily ``odd content''.

On the other hand, I was scared (since three days) by the fact
that something which I subconsciously thought as evident (or
rather which I was sure to have read somewhere) is perhaps not
true. My (naive) question is whether two empty curves are
necessarily rigid-isotopic? This looks at first between
metaphysical nonsense and ``triviality''? Maybe it is unknown,
when $m$ is large enough. (m=6 follows from Nikulin 1979, and as
far as I know there is not a simpler proof, say valid for all
(even) degrees). So I am quite shameful asking you about this
point: Is the empty chamber always connected?

I tried a dynamical approach (to this problem) in Section
25.12(\ref{rigidity-empty-scheme-via-dyna:sec}), but it is not
very convincing. On the other hand, if the empty room is
connected, then maybe the space of all curves with one component
is also connected? (Naively one would apply the Itenberg-Viro
contraction conjecture, to reduce to the empty case, move there
for a while to resurface at the other curve (the contraction
thereof). Perturbing this path in the ``visible world'' would
conclude the proof modulo some difficulties...) Again, you Oleg,
in your wonderful survey of 2008 (in Japanese J. Math) lists as an
open problem the question of deciding the rigid-isotopy of curves
of odd degree having a unique real circuit. As you emphasize the
word ``odd degree'', I wondered if the case of even degree (again
with only one oval) is already settled?

In Section 25.10(=\ref{CCCviaDynamics:sec}), I have attempted a
naive dynamical approach to the collective contraction
conjecture(CCC). This states that we can shrink simultaneously all
the empty ovals toward solitary nodes. This is a bit like a
perfect landing in flight simulator where all wheels touch the
ground simultaneously. My naive strategy is just to study the
gradient lines of the functional measuring the total area of all
empty ovals, but it  is surely not serious. It would be exciting,
in my opinion, to describe a counterexample to CCC if there is
one.

Many thanks for the attention, all your patience about my naive
reasonings, and above all for the brilliant answers you already
gave me. Best regards, Alex

$\bullet$$\bullet$$\bullet$ samedi 26 janvier 2013 20:15:54, the
prompt response of Eugenii Shustin:

Dear Alex,

The chamber of empty curves of a given (even) degree is indeed
connected: two such curves can be defined by homogeneous
polynomials, positive for any real not all zero variables, and
their linear homotopy $(1-t)P+tQ$, $0\le t\le 1$, gives a path in
the chamber of empty curves.

By the way, another (well) known connected chamber consists of
hyperbolic curves (i.e. those which have totally real intersection
with lines of certain pencil) - this is a consequence of Nuij W. A
note on hyperbolic polynomials. Math. Scandinavica 23 (1968), no.
1, 69--72.

With best wishes, Eugenii

$\bullet$$\bullet$$\bullet$ samedi 26 janvier 2013 21:08:27, Oleg
Viro:

Dear Alex,

The counterpart of the Rokhlin conjecture\footnote{[30.01.13] This
is another very interesting (inedited?) piece of information, as
to my knowledge this never appeared under the (printed) pen of
V.\,A. Rohlin. So here Viro tell us something very interesting not
yet available in print (as far as I know).} about rigid-isotopy of
any two curves of odd degree with one component is the obvious
observation described by Evgenii, about empty curves of even
degree. The question about curves of even degree with a single
oval is equivalent to the question about removing this single oval
by an algebraic Morse modification. I don't think it was ever
discussed, but I could miss it.

$Z_{15}$ is not $E_6$. The easiest way to construct the Gudkov
$M$-curve is by perturbing two $J_{10}$ singularities of the union
of 3 non-singular conics tangent to each other at 2
points.\footnote{[30.01.13] This is probably true yet this
requires a more highbrow dissipation theory, than in Gudkov's
second existence proof which apart from the trick with Cremona
uses only dissipation of ordinary nodes \`a la Brusotti 1921
(so-called $A_1$-singularities).}

Best, Oleg

$\bullet$ [28.01.13, lundi 28 janvier 2013 20:03:58] Gabard wrote
euphorically\footnote{[30.01.13] This alas turned out to be quite
foiled as the invisible part of the discriminant as real
codimension 2.} an e-mail titled ``Some more metaphysical
non-sense about the rigid-isotopy of empty curves?'':

Dear Eugenii, Oleg, Viatcheslav, Alexis, Stepa, etc.

So many thanks, Eugenii, for putting me again on the right track,
and recalling me the argument which I shamefully forgot about.
Yesterday, I was quite excited by trying to digest your argument
(albeit it
seems so simple). In fact the little detail that worried me is
that I do not know why during the linear homotopy $(1-t)P+tQ$ the
variable curve could not acquire (while staying of course empty if
$P,Q$ have both the same sign) a pair of conjugate nodal
singularities. This puzzled me for a while, and then using
systematically your argument, I arrived at the somewhat opposite
conclusion that the empty (smooth) chamber must be disconnected
(for all even degrees $m\ge 4$)!!! This violates all what we know
since Rohlin 1978 (and surely Gudkov as well??), while the former
refers directly back to the argument of Klein 1876 based on
Schl\"afli cubics surfaces $F_3$'s and Zeuthen correspondence
between cubic surfaces and quartics (via the apparent contour).
Klein's proof is a bit tricky and uses as well his rigidification
(Klein 1873) of Schl\"afli's isotopic classification. Needless to
say I could not follow Klein's reasoning completely, as I just
studied it today for ca. 2 hours. Marin informed me recently that
he, in contrast, was able to digest all of those Kleinian works!

So using your method of linear homotopy, one sees quickly that the
(cone) space $C^+$ of positive anisotropic (=not representing
zero) forms is contractile (convex actually) hence
simply-connected. Its projection in the space of curves is the
invisible locus $I$ consisting of all empty curves. Since the
latter is merely a quotient of $C^+$ (by positive homotheties) it
follows that it is also simply-connected (via the exact homotopy
sequence of a fibering). But the discriminant is visible inside
this invisible locus $I$, since it is a simple matter via Brusotti
(1921) to construct empty curves with a pair of conjugate nodes.
Thus we see inside the simply-connected manifold $I$ a certain
hypersurface (namely a portion of the discriminant), which by
Jordan-Brouwer (or a slight extension thereof) should separate
this manifold $I$ in pieces (at least so is my naive intuition).
It follows that our empty chamber (consisting of smooth curves) is
disconnected!!!! This is my proof in its broad lines (for more
details, compare Section 25.13, page 282-283 of the attachement,
Theorem 25.29 and its proof on page 283). This is just one page
long...

Since this conclusion contradicts violently what is asserted by
Klein 1876 (and approved by Rohlin 1978), it is of course very
likely that my proof contains a serious flaw, or at least that I
am confusing somehow the basic conceptions. However presently I do
not see where is my mistake! Of course, my pseudo-theorem also
violates the part of Nikulin 1979 concerned with the rigid-isotopy
of the empty chamber of sextics.

Many thanks again for your attention, and sorry for overflowing
your mail boxes with my naive questions (and dubious reasonings).
Thank you again so much for all your excellent and detailed
responses (especially on $E_6$ and $Z_15$).

Best wishes, Alex PS: I send you a copy of my TeX-file in case
someone would like to work out a specific passage. At the occasion
I would also be happy to send you my figures in zipped format so
that one of you can continue the project in case I make a fatal
bicycle accident (like Academician V.I. Arnold?)

$\bullet$ [30.01.13, 18h10]

Dear Oleg, Eugenii, and the other experts,

I think that I found the mistake in my ``proof'' of the
disconnection of the empty locus (that you certainly noticed
meanwhile in case my explanation is the correct one). The reason
seems to be simply that the discriminant inside the invisible
locus has only real codimension 2, hence cannot separate anything.
I have attempted to explain this in Section 25.14 on page 288. If
this is not wrong it seems that the next natural question is to
decide which chambers residual to the principal stratum of the
discriminant contains such smaller pieces of the discriminant
shrunk to codimension 2. I think to have found a topological
obstacle for $M$ and $(M-1)$-curves, and conjecture (very naively)
this to be the sole obstruction. In more geometric terms, this
amounts essentially to decide which smooth curves can acquire a
pair of imaginary conjugate nodes.

Many thanks, Eugenii and Oleg, for your detailed answers. As you
said, it seems that the (Itenberg-Viro) contraction conjecture of
empty ovals implies the rigidity conjecture for even order curves
with a unique oval. However I should probably still try to
understand this implication in some more details. Perhaps it is
somehow related to the previous codimension 2 phenomenon inside
the ``invisible'' chamber.

Sorry for all my confusing messages, and many thanks again for all
your kind efforts in trying to educate myself. All the best, Alex

$\bullet$ [01.02.13, ca. 20h00]
Obstruction to rigid-isotopy (strictly) below height DEEP+2?

Dear Oleg, Eugenii, Alexis, Thomas, Stepa, etc.

Many thanks for all your brilliant messages and articles I am
still slowly trying to assimilate properly. I hope not taking too
much of your precious time. Albeit I met all of you only rarely, I
remind very accurately your brilliant talks (in Geneva or Rennes),
and so it is a special pleasure to remind each of yours while
trying to explore this fantastic topic.

On my side I was those last two days fascinated by the conjecture
that the one-oval scheme ought to be rigid, as Oleg or Rokhlin
conjectures. (Let me say that a scheme is rigid, if all the curves
representing it are rigid-isotopic.) Given a degree $m$, one may
wonder what is the smallest height $r(m)$ at which there is a
non-rigid scheme. (For me the height of a scheme just means its
number of components.)

For any degree $m$, there is of course the deep nest with
$r=[(m+1)/2]=:DEEP$ real branches. Two units above the latter's
height, it is easy to construct (for each $m$) curves having the
same real scheme yet different types (I vs. II) hence not
rigid-isotopic. (This is a simple iteration of Rohlin's
construction in degree $5$, cf. Figs. 102, 103 in my file). Using
the Marin-Fiedler method of the  lock it is even possible to
exhibit at this height $DEEP+2$  curves of degree 7 or 9 having
the same real scheme and the same type~II, yet not rigid-isotopic
(Figs. 105, 104). (Probably the method extends to all other odd
degrees.) However, it seems much more tricky (and the lock-method
seems ineffective) to detect obstruction below this height
$DEEP+1$ (i.e. one unit above the height of the deep nest). Could
it be the case that all schemes at or below this height are rigid?
Of course this looks super-optimistic as we do not even know
rigidity at height one, but I was unable to find a counterexample.
I would be very happy if you know one? If there is a simple
candidate, I hope to detect it alone during the next few days\dots
. So do not take care answering me if my question is trivial. (As
I just work on this since two days, I probably missed something
accessible.)

Paraphrasing slightly, I found quite puzzling, that the very
explicit function $r(m)$ measuring the smallest height of a
non-rigid scheme is only subsumed to the large pinching $1 \le
r(m)\le [(m+1)/2]+2=DEEP+2$. Of course a better lower bound seems
out of reach, but perhaps you know better upper bounds.

I also wondered if there is an extension of the Nuij-Dubrovin
rigidity of the deep nest to, say, the totally real scheme of
degree 8 consisting of 4 nests of depth 2. I should think more
seriously on this at the occasion.

Many thanks for your attention, and sorry again for all my
enthusiastic and naive e-mails. All the best, Alex

$\bullet$ [written 08.02.13 and sent 09.01.13]

Dear Oleg, Eugenii, Stepa, Viatcheslav, etc.

I still continued my trip through real plane curves and cannot say
that my curiosity is starting to fade out. I tried for several
days to find a counter-example to the conjecture (of mine so
probably quite wrong) that all schemes below height $DEEP+2$ are
rigid, where $DEEP=[(m+1)/2]$ is the number of branches of the
deep nest of degree $m$. At least the method of the lock
(Fiedler-Marin) seems quite inoperant to detect an obstruction to
rigid-isotopy at such low altitudes. If true, the proof probably
involves a geometric flow collapsing either the pseudoline to a
line (by shortening its length like a systole) or improving the
rotundity of some oval to a circle (via an isoperimetric
functional?). If all this works, it would reduce the low-altitude
rigidity conjecture to Nikulin's theorem (or maybe even Klein's on
$C_4$) as the starting step of a big recursive process. Of course
this seems still quite out reach (canary music) unless one feels
very motivated!

Next I tried to corrupt the truth of Slava's remarkable
rigid-isotopic classification (Nikulin 1979) of sextics via the
Marin-Fiedler locking argument using B\'ezout saturation. Of
course I have nothing against Slava, but this was rather intended
to test experimentally the power of Nikulin's result. Specifically
I looked at sextic schemes of the form $3/1 \ell$, and wondered if
for some specific curves the distribution of the $\ell$ outer
ovals away the fundamental triangle traced through the 3 inner
ovals (those enveloped by the unique nonempty oval) could be
different for different curves. On all examples I tested it seems
that the outer ovals are never separated by the ``deep'' triangle.
So we find no violation of Nikulin's theorem, and the latter
rather implies that as soon as we are able to visualize the
distribution for a single curve it will be the same for all curves
belonging to this scheme. The case most tricky to understand is
the maximal permissible, namely $3/1 5$. I managed to construct it
\`a la Harnack (as preconized in Gudkov 1974 or 1954). But being
quite unable to decide from this model the distributional question
of the outer ovals past the fundamental triangle, I decided to
switch to Oleg's method of construction via dissipation of the
singularities of a  triplet of coaxial ellipses. I played this
game yesterday but could not decide the distributional question
for this Viro curve (cf. especially Fig.126 on page 320 and the
hypothetical Theorem 26.29 on page 322). In fact today I tried
again to inspect directly Harnack construction and found Lemma
26.27 on page 319 whose proof seemed to me very transparent until
I found the little warning, which I think is not fatal. In
conclusion I believe now that there is no separation by the
fundamental triangle!!!??

Of course I imagine that, if I am not completely wrong, what I am
investigating must be quite familiar to you. I would much
appreciate if you know if this hypothetical theorem (26.29 page
322) is true. It amounts essentially to check whether in Viro's
construction of $3/1 5$ the triangle through the 3 deep inner
ovals does not separate the 5 outer ovals. I find this question
quite attractive as it seems to require some understanding of the
geometric location of the microscopic ovals arising in Viro's
method (optionally compare Fig. 127 (page 322) which shows a
scenario with the two bottom micro-ovals aligned vertically in
which case the fundamental triangle would separate the outer
ovals). This scenario seems to me quite unlikely but it does not
seem to be impeded by naive B\'ezout obstructions.

Many thanks for your attention, and sorry again for all my naive
and confuse questions. Thank you very much again for your precious
guidance and answers.

All the best, Alex

$\bullet$$\bullet$$\bullet$ (10.02.13) Bonjour Alexandre, Thomas
m'a transmis ta question. La r\'eponse est toute simple: soient
$A$, $B$, $C$ trois ovales int\'erieurs et $D$, $E$ deux ovales
exterieurs de ta sextique. Le triangle fondamental $ABC$ est
enti\'erement contenu dans l'ovale non-vide. Si $D$ et $E$ sont
dans deux triangles $ABC$ (non-fondamentaux) diff\'erents, alors
la conique passant par $A$, $B$, $C$, $D$, $E$ coupe la sextique
en $14$ points, contradiction. Avec des coniques, on montre plus
g\'en\'eralement que: Les ovales vides de la sextique sont
distribu\'es dans deux chaines (int, ext), l'ordre cyclique est
donn\'e par les pinceaux de droites bas\'es dans les ovales
interieurs. Les ovales interieurs sont dispos\'es en position
convexe dans l'ovale non-vide. Bon dimanche,   S\'everine

$\bullet$ [12.02.13] Is the Gudkov chamber simply-connected?

Dear S\'everine, Viatcheslav, Ilia, Oleg, and all the other
experts,

First many thanks, S\'everine, for your excellent answer on my
distribution question of ovals of sextic, and sorry for my late
reply on it as I lack an Internet connection at home.

I tried today to understand when a dividing (plane) curve admits a
transmutation, i.e. a rigid-isotopy permuting both halves of the
curve. I also studied the weaker notion of mutation of when there
is a linear automorphism of the plane permuting both halves. Using
the Kharlamov-Itenberg calculation of the monodromy of sextics I
think that I managed to get some obstruction to mutability,
especially for the 3 dividing curves which have trivial
monodromies (compare Lemma 26.6, page 288, which is hopefully
correct). However I don't know if the Gudkov chamber (or the 2
other related ``antidromic'' chambers, i.e. having trivial
monodromies) is simply-connected. I hoped to detect some non
simple-connectivity by looking at the monodromy induced on the
halves instead of the ovals. At least this works of course for the
deep-nest chamber which is not simply-connected since there is a
symmetric model which can be mutated. So my (hopefully not too
naive) question is the following: is it known whether or not the
Gudkov chamber is simply-connected? (equivalently is the Gudkov
curve transmutable?) The same question looks attractive for the
other 2 antidromic curves, i.e. the left wing ``Rohlin curve''
$6/1 2$ and $4/1 4$ in type~I.

Thank you so much for all your attention and patience, and in
advance for your answer if it is known. Best regards, Alex

$\bullet$$\bullet$ $\bullet$  mercredi 13 f\'evrier 2013 04:34:12

Dear Alexandre, If I understand correctly the question then the
answer is not, if I state the question appropriately then the
answer is yes. I mean the following precise statements. Let
consider the part of the projective space of real sextics that is
represented by maximal sextics of Gudkov's type. Then the
fundamental group of this part is $Z/2$. It becomes simply
connected after taking quotient by the natural action of
$SL(3,R)$. In fact, before factorization it is a fibration over
contractible base with the fiber $SL(3,R)$. These results (and
there analogs for other maximal sextics and certain curves of
lower degree) are contained in my talk On monodromies of real
plane algebraic curves at one of Petrovsky seminars in 80th, I
guess (short summary should be found in Russian Surveys). The
proof (in the case of sextics) is rather straightforward as soon
as based on the $K3$ surfaces periods uniformization. As it
happens rather often with this approach, to treat the maximal
curves is extremely easy, since the corresponding eighenlattices
become unimodular. In general the period domain, which is the
product of two polyhedra in the real case, represents the studied
sextics (or associated K3 surfaces) only up to codimension 2.
Which makes laborious to treat the fundamental group. But,
surprise, in the case of maximal curves there are no codimension 2
phenomena, since such holes appear only as traces of $(-2)$-cycles
having nontrivial components in the both eighenspaces, which is
impossible since in the maximal case the components are integral
and the eighenlattices are even. I don't remember by heart the
final result for other maximal sextics. It should be pointed in
the same summary and by the way easy to get following the same
approach I have pointed. The key is that even if it is no more a
pure fibration - it has special fibers which are quotients of
$SL(3,R)$ by the corresponding monodromy group (which indeed
coincides with the maximal possible group of symmetries for the
given type of sextics) - its fundamental group is exactly the
fundamental group of this special quotient. Yours VK

$\bullet$ mercredi 13 février 2013 11:46:20

Dear Colleagues, Am I alone who did not receive a copy of
Severine's letter? I would be happy to know its content :) Yours
VK

$\bullet$ mercredi 13 f\'evrier 2013 13:43:10

Dear colleagues, I had written only to Alexandre, sorry! My answer
was this: let $A, B, C$ be three inner ovals, and $D, E$ be two
outer ovals of the sextic. The fundamental triangle $ABC$ is
entirely contained in the nonempty oval. If $D$ and $E$ are in two
different (non-fundamental) triangles $ABC$, then the conic
through $A, B, C, D, E$ cuts the sextic at 14 points,
contradiction. Using conics, one proves more generally that there
is a natural cyclic ordering of the empty ovals, given by the
pencils of lines based at the inner ovals. The empty ovals are
distributed in two consecutive chains (inner, outer). The inner
ovals lie in convex position in the nonempty oval. Best regards,
S\'everine

$\bullet$ [14.02.13]

Dear Viatcheslav, S\'everine and all the other colleagues,

Thank you very much for this beautiful answer on the Gudkov
chamber. I look forward to digest properly all that incredible
technology that you and Nikulin developed. Again many thanks also
to S\'everine for the clever argument which I digested yesterday
with great pleasure, and integrated in my notes in Section
26.10(=\ref{LeTouze:sec}) pages 332--334. This gave me yesterday
some motivation again to attack the very first question of all our
chat room, namely Rohlin's claim that the pencil of cubics through
the 8 deep basepoints located inside the 8 empty ovals of any
sextic curve of type $6/1 2$ or its mirror $2/1 6$ is totally
real, hence of type~I (also called orthosymmetry by Klein ca.
1881-82 and his student Weichold 1883).

In fact I (naively) hoped to prove this Rohlin claim via
Poincar\'e's index theorem, yet the qualitative picture (Fig. 133
on page 337) rather inclined me to believe that the proof cannot
reduce to mere combinatorial topology of foliations (i.e.
Poincar\'e's index formula of 1885). So I am still puzzled, but
perhaps an argument like S\'everine's one do the job. At any rate
I would be very excited if someone manages to reconstruct this
proof asserted by Rohlin (1978) if it is not too tantalizing for
the brain.

Otherwise I am also much frustrated by failing to visualize
totally real pencil on the three $M$-sextics, whose existence is I
think predicted by Ahlfors theorem of 1950 (or better the special
zero-genus case thereof known to Riemann 1857, and reworked by
Schottky 1875-77, or even Bieberbach 1925 and his more respectable
student Grunsky 1937). Marin warned me recently that the
transition from the abstract Riemann surface viewpoint to the
planar context ``of Hilbert's 16th problem'' may be not so easy as
I always assumed subconsciously. (If necessary, all the
correspondence I received from all the colleagues is gathered in
Section 24.6, p.221). Overpassing this difficulty (which I hope is
not fatal) there should be on all $M$-curves (more generally
dividing curves) auxiliary pencils which are totally real. Alas
for $M$-sextics (even $M$-quintics), I am completely unable to
trace them and know nothing about the degree of the curves
involved (in the pencil). I hope to be able to tackle such
questions in the future, but perhaps you have better ideas (or
motivations) than I do have.

Thank you very much again for all your brilliant answers, and kind
messages. All the best, Alex

$\bullet\bullet\bullet$ samedi 16 f\'evrier 2013 17:54:55

Dear Alexandre, dear other colleagues,

I have managed to prove that a pencil of cubics with eight
base~points distributed in the eight empty ovals of a sextic
$2 \cup 1(6)$ is necessarily totally real. Details will follow
soon in a paper. Yours,

S\'everine

$\bullet$ [16.02.13,19h41]

Dear S\'everine and colleagues,

Congratulations for this fantastic achievement. I am sure the
proof must be very beautiful. On my side I tried to work out for
all sextics of type~I an optical recognition procedure of the type
by some synthetical procedure akin to Rohlin's claim, yet this is
still much in embryo. In particular the case of $(M-4)$-sextics is
quite puzzling as it seems to contradict the version of Ahlfors
theorem due to myself (existence of a totally real map of degree
the mean value the number of ovals and Harnack's bound). I hope to
send you more palatable material soon, but confess that the
questions look quite hard and I seem much less efficient than
S\'everine. So I suppose that Rohlin's claim is one among several
other (less pure) total reality result. So I look forward with
great interest to see S\'everine's article.

All the best, Alex

[19.02.13] Dear colleagues,

Many congratulations again to S\'everine for your fantastic
achievement. Sorry to have been brief in my last letter, as I
wrote (lacking an internet connection at home) from a friend of
mine who had a romantic party with his girlfriend, and I do not
wanted to disturb too long his romantic evening. Meanwhile I also
tried hard to concentrate on a proof of the Rohlin-le Touz\'e's
theorem, which still overwhelms my intelligence. The last things
that I have written are on pages 336--352 (Sections 27.1, 27.2,
and 27.3), but this is poorly organized and supplies no serious
proof of the Rohlin-Le~Touz\'e's theorem.

Some few days ago, I got Theorem 27.5 (on page 346), which (if it
is correct) answers one of the question I asked in my penultimate
e-mail (as well as desideratum of Alexis), namely the question of
estimating the order of curves involved  in a total pencil on an
$M$-curve. It seems that there is always such a pencil of order
$(m-2)$, i.e. two units less than the given degree $m$ of the
$M$-curve. In fact, the proof is a nearly trivial adaptation of
the abstract argument going back to several peoples (in
chronological order Riemann 1857, Schottky 1875, Enriques-Chisini
1915, Bieberbach 1925, Grunsky 1937, Courant 1939, Wirtinger 1942,
Ahlfors 1947, 1950, a myriad of Japaneses, a myriad of Russians
including Golusin 1953/57, etc....., up to Huisman 2000, and
Gabard 2001/2006, who else?). The point is that total reality is
trivial in the case of $M$-curves since we have one point
circulating on each oval (such a group moves by Riemann-Roch!!!)
and so we have like a train-track with only one train on each
track, hence no collision can occur and total reality is
automatic. If we work with plane curves we only need to take
curves of order $(m-2)$ which have enough free parameters to pass
through any given distribution of $M$ points (one on each oval),
and this works by looking at the residual group of points (details
in the proof on page 346). So this is quite interesting but
probably only a first step toward deeper things. (One could dream
to recover all the Gudkov-Rohlin/Arnold congruence via this method
but that looks hard work...)

After this little discovery I focused again on the
Rohlin-S\'everine theorem, yet without any success. So I have not
more to report for the moment.

Thanks a lot for the attention, and all my congratulations again
to S\'everine for your deep advance. Best wishes, Alex

$\bullet$ 19.02.13 Dear Alex,

Let me ask you a question from your previous field of interest. Do
you know any example of a non-Hausdorff 1-manifold which does not
admit a differential structure? I heard about existence and could
easily construct examples of exotic, i.e., homeomorphic but not
diffeomorphic non-Hausdorff 1-manifolds. See

http://www.map.mpim-bonn.mpg.de/1-manifolds

Sincerely, Oleg

20.02.13

Dear Oleg, David and Mathieu,

Many thanks, Oleg, for your lovely question, and best greetings to
the other friends. Alas my memory is failing quite dramatically,
so my answer will be of poor quality. If I remember well I asked
myself the same question some 3-4 years ago, but I cannot record
to have ever found an answer. Thus I forward your question to
David and Mathieu, the leading experts of non-metric surfaces who
perhaps will supply a better answer. On my side I hope to think
more seriously to your question when I see clearer with
Rohlin-Le~Touz\'e's sextics.

Maybe a first idea is that there ought to be a (non-canonical)
``twistor construction'' assigning to each non-Hausdorff curve a
Hausdorff surface fibered by (real) lines. This construction
should go back to Haefliger's very first note in the colloque de
topologie de Strasbourg ca. 1955-1956  (yet it is not very
detailed). In substance it is like a train-track construction \`a
la Penner-Thurston…(some intuition about this is given in my
article `Ebullition and gravitational clumping, arXiv, 2011). Do
not worry if you don't understand me, as I myself remember only
vague souvenirs and are not so convinced by what I am saying!!! In
fact Haefliger (ca. 1956) claims this construction only for second
countable curve (even with a proviso on the fundamental group),
but when I was in touch with the subject I was fairly convinced
that it must work universally.

OPTIONAL REMARK: Haefliger, and Haefliger-Reeb 1957 use this
construction to prove that any simply-connected curve (second
countable) arises as the leaf space of a foliation of the plane.
(Sketch of proof: take the twistor of the given curve which is by
the exact sequence of a fibering 1-connected and (by
Poincaré-Volterra) second countable, hence it is the plane, q.e.d)

So the idea would be to descend a smooth structure on the surface
to get one on the curve. Alas, it is a well-known open problem
whether any (non-metric but Hausdorff) surface admits a smooth
structure (Spivak 1971, Nyikos, etc.) However quite puzzlingly
Siebenmann 2005 (Russian Math Surveys) claims (and even prove in
some details) that a PL structure exists universally on all such
surfaces, merely as a consequence of Schoenflies theorem. So
perhaps Siebenmann argument work as well for DIFF structures, and
the metaphysical problem of Spivak-Nyikos is cracked. If this
works (ask maybe Siebenmann, or an Indian in the
States(=Ramachandran) who albeit not an expert was fairly
convinced that there should be no asymmetry between PL and DIFF in
dimension 2), then there is perhaps some chance to get a smooth
structure on all non-Hausdorff curves. Of course there is perhaps
a more direct strategy without transiting through surfaces.

Otherwise, regarding exotic smooth structures on curves the
original reference is Haefliger-Reeb 1957 article in
L'Enseignement Math. Perhaps you could quote this in your
brilliant web-page.

Sorry for this vague answer, but at the moment my brain is much
concentrated on this Rohlin-Le Fiedler total reality claim which
still puzzles me a lot!!!

Best greetings to all, as well as to Rachel and Chiara. All the
best, Alex

$\bullet$ [22.02.13] Dear colleagues (especially S\'everine),

I worked hard (but without success) on the Le~Touz\'e's theorem,
at least for 8 basepoints assigned on the nonempty ovals of a
sextic of type  $6/1 2$. If I understood well S\'everine's
announcement, you rather handle the case of $2/1 6$ and assign
more generally the points in the insides of the empty ovals (but
of course I suppose that your argument adapts to $6/1 2$). Even in
my weaker form I am not really able to conclude but send you my
last thinking on the question (Section 27.4, p.352--356, esp.
Fig.141). Ultimately I found a method which I call ``barrages''. A
special r\^ole is played by nodal cubics of the pencil, and I try
to get a corruption with B\'ezout by looking at nodal curves with
a barrage, i.e. such that 4 arcs of some other cubic joins
pairwise the 8 basepoints distributed on the loop of the original
cubic. (By the loop of a nodal cubic, I mean the unique path from
the node to itself which is null-homotopic in the plane $RP^2$.)
Of course I am not sure that details can be decently completed,
but for the moment it is the only reasonable strategy I could
imagine. I am sure that S\'everine's argument is much more elegant
and convincing. My reasoning is completely conditioned by Fig.141,
and I am probably too naive in believing that it reflects the
general situation.

Sorry for sending you this very coarse material, and of course do
not take the pain to react to this message. Many thanks again a
lot to all for sharing so generously your knowledge and for all
your answers. Best regards, Alex

[25.02.13] Dear real geometers,

I was still much fascinated by the Rohlin-Le Touzé theorem (RLT)
albeit still not able to prove it. Being frustrated by my failing
attempts (probably due to a lack of stubbornness and competence in
algebraic geometry) I decided to speculate a bit of why it is so
important or at least to explore how the statement could
generalize.

In its most elementary incarnation involving pencil of lines and
conics, the phenomenon of total reality occurs along infinite
series stable under the operation of satellite of a real scheme
(of even order). Satellite just amounts to trace each oval with a
certain multiplicity $k$ (jargon obviously borrowed from knot
theory). So the unifolium scheme of degree 2 (allied to a conic)
gives rise to the deep nests, and the quadrifolium scheme of
degree 4 gives rise by taking its satellites to an infinite series
of schemes of order multiples of 4 which are totally real under a
pencil of conics (assigned to pass through the deepest ovals). It
seems therefore natural to ask if the satellites (e.g. the second
satellite) of the Rohlin's scheme $6/1 2$ (or its partner $2/1 6$)
are also totally real (and hence of type~I) under the ``same''
pencil of cubics as posited by the Rohlin-Le~Touz\'e theorem. Alas
I was not even able to settle this question. (Of course this seems
evident (granting RLT) for a small perturbation of the algebraic
double (essentially $F \cup F+\epsilon$), since total reality
forces transversality of the foliation induced by the pencil with
the curve.)

Next, I tried to understand what are the higher order avatars of
the RLT-theorem (in the hope that it is not an isolated phenomenon
as vaguely suggested by Ahlfors theorem). I found using the
Rohlin-Kharlamov-Marin congruence ensuring the type~I-ness
(=orthosymmetry) of some $(M-2)$-schemes an (obvious) infinite
series of avatars of the Rohlin's $(M-2)$-schemes of degree 6 .
Those are also $(M-2)$-schemes and total reality seems to be
possible for a pencil of curves of order $(m-3)$, exactly like for
the G\"urtelkurve of Zeuthen-Klein (bifolium quartic with 2 nested
ovals totally flashed by a pencil of line through the deep nest)
or for the Rohlin's sextic (flashed by a pencil of cubics). So it
seems that the theory of adjoint curves of order $(m-3)$ plays
some special r\^ole in this question of Rohlin-S\'everine. I would
be very happy if one of you knows if it is reasonable to expect an
extension the RLT total reality theorem to all this schemes whose
type~I ness is ensured by Rohlin-Kharlamov-Marin congruence (sorry
if I am not hundred percent right in crediting as I could not
extract the exact history of this subliminal result). Specifically
I have Conjectures 27.17 and 27.18 (page 365 and 367 resp.) which
list some candidate-schemes for total reality in degree 8 and 10.
If the conjectures are right, it would be of great interest to
know if S\'everine's proof adapts to them. Sorry if I am too naive
about the real difficulty of such problems, but I found exciting
to wonder if there is something more general behind the cryptical
allusion of Rohlin. Of course I presume that he derived the
synthetic result a posteriori from highbrow topology (or K\"ahler
geometry in Kharlamov's case?), but perhaps there is a simple
explanation with (``basic'') algebraic geometry and total reality
as S\'everine was able to do? As Oleg knows my problem is that I
wasted too much time with non-metric manifolds and so forgot all
the little I ever knew about algebraic geometry.

During the way, I think to have found a counterexample to the
conjecture of mine (inspired by the Itenberg-Viro contraction
conjecture of empty ovals), and according to which all empty ovals
could be contracted simultaneously to solitary nodes. This
counter-example is Thm 27.16 on page 364 (which I hope is correct
and sharp as far as the degree is concerned).

Thanks a lot for the attention, and sorry for all the modest news
(you surely thought about in sharper form already). All the best,
Alex PS: The material summarized in this message occupies page
357-367 (Sections 27.5, 27.6, 27.7), as usual I had not much time
to polish, but I hope it is still readable.

[27.02.13] A census of 100 octic $(M-2)$-schemes of type~I
satisfying the RKM-congruence, plus a little addendum for Oleg's
non-Hausdorff curves

Dear colleagues,

I have pursued some preliminary study toward the total reality
phenomenon, yet merely in its combinatorial aspect prompted by the
modulo 8 RKM-congruence (for Rohlin-Kharlamov-Marin) ensuring the
type~I of $(M-2)$-schemes of degree $2k$ with $\chi = k^2+4 \pmod
8$. Accordingly, I call an RKM-scheme any  $(M-2)$-scheme
satisfying this congruence. While any RKM-scheme is of type~I, I
do not know alas whether the converse statement is true. If it is
known I would be extremely grateful if someone can tell me (and
our collective chat room) the answer. Further I noticed that the
list given in my previous e-mail of RKM-schemes of degree 8 can be
much enlarged. If I am not too bad in combinatorics, there are
precisely 100 such schemes in degree 8, all of them being
potentially subsumed to the phenomenon of total reality under a
pencil of quintics akin to the Rohlin-Le~Touz\'e theorem (for
sextics flashed by cubics). This modest material is to be found in
Section 27.8, p.368-373 (especially Fig. 146 page 370 and Lemma
27.24, p.372, plus all the 36 Gudkov symbols on page 372). I hope
of course that I missed nobody in this catalogue. Extrapolating a
bit using the (hypothetical) converse statement to RKM, I would
say that there are precisely 100 schemes of type~I which are
$(M-2)$-schemes. Is this well-known and correct?

Actually, I do not really know if all these 100 schemes are
realized algebraically, but presume that most of them (all?) are.
Possibly I am much too naive. Of course it is quite amazing to see
that the only two RKM-schemes of degree 6 (namely $6/1 2$ and $2/1
6$) demographically explodes to a menagerie of 100 such schemes in
degree $8$, but that should be no surprise for you much acquainted
with the higher cases of Hilbert's 16th problem. It would be even
more crazy if all those 100 schemes (or at least a good portion
thereof) are subsumed to the phenomenon of total reality. If you
have some ideas on those circle of ideas, I would be extremely
thankful.

Many thanks again for the patience and attention, and I hope that
what I am telling is nearly correct (not too surrealist). Very
best regards, Alex

PS: For Oleg, regarding my loose answer on smooth structures on
non-Hausdorff 1-manifolds, I would like to add another
philosophical remark related to the method of Haefliger's
``twistor''. This is of course like a thickening along a normal
bundle except that there is no ambient manifold (safe the ether)
and so the construction must be intrinsic. To my knowledge it was
never exposed in details (albeit Haefliger's 1st article ca.
1955-56 in Colloque de Topologie de Strasbourg uses implicitly
this construction). Now my point is that albeit the twistor method
looks somewhat indirect, I think that it is fairly useful. For
instance, I was since 2006-07 puzzled by the naive question if the
fundamental group of a one-manifold is always a free group. (Of
course such non-Hausdorff curves resemble somehow graphs, whence
some intuition). For instance the line with 2 origins has
$\pi_1=Z$ as follows quickly from Seifert-van Kampen (and if 3
origins or 2 doubled origins then $\pi_1=F_2$ is free of rank 2).
Ultimately in 2011 I found a general answer to this ``freeness''
puzzle by using the Haefliger twistor construction, while showing
first that all open (non-metric) surfaces have free fundamental
groups. (This is actually a very modest extension of the metric
case, which to my knowledge is first treated in Ahlfors-Sario book
of 1960, albeit it may have belonged to the folklore much earlier,
say Kerekjarto, H. Kneser, Rado, in the 1920's, Papakyriakopoulos
in the 1940's???). This material is exposed in some details in my
arXiv note of ca. 2011 (Ebullition in Foliated surfaces versus
gravitational clumping). I hope that those results are nearly
correct but they certainly require more professional treatments
and exposition than I was able to do. I hope this little remark
makes perhaps more plausible that the approach via the
(Hausdorffizing) Haefliger twistor is also reasonable for your
problem of DIFF structures.

$\bullet\bullet\bullet$ vendredi 1 mars 2013 19:07:12

Dear Alexandre, dear other colleagues,   here is the note I had
promised to send you.  There are still many open questions, as
Alexandre wrote. It would be also interesting to know whether one
could find a totally real pencil with respect to the dividing
$M-2$-sextics with real scheme of indefinite type. I will think
about it when I have more time. Best regards,

S\'everine

(01.03.13, 22h15)

Dear S\'everine and the other colleagues,

So many thanks S\'everine for sending us your splendid article. I
am much excited to read the details tomorrow, as myself started
today to doubt about the whole result (at least in the strong form
that any points 8 points distributed on the empty ovals ensures
total reality). (If I am not wrong the whole phenomenon depends
upon the location of the 9th base point, namely the pencil is
totally real iff the 9th base point lands in the inside of the
nonempty oval.) So I was much depressed and lost in my poorly
organized thoughts. So your sending arrives as a true deliverance
for my brain.

Many congratulations again to S\'everine for this fantastic work.
Very best regards, Alex

[02.03.13] Can total reality fail for a distribution of 8 points
on the empty ovals?

Dear S\'everine and the other geometers (especially Professor
Nikulin),

I enjoyed much a detailed look at your splendid article full of
illuminating remarks. I will probably need much more time to
digest the impressive technology you use, and need to print the
material to make a deeper reading (especially of the former works
upon which your argument seems to depend). So many thanks again
for sending us your work in so rapid delay. I wrote some naive
reactions in Section 27.11, where I mostly copied your sayings,
and tried to add hopefully pertinent footnotes.

Regarding your question ``Can conversely any dividing curve be
endowed with some totally real pencil?'', I still wonder if a
positive answer is not a trivial consequence of Ahlfors theorem
(compare very optionally Gabard's Thesis 2004, page 7). However
since Marin warned me in January 2013 (cf. Section of e-mails) it
may be the case that the transition from the abstract conception
of Riemann-Schottky-Klein to the embedded viewpoints of
Hilbert-Gudkov-Rohlin is not so easy. Yet I am still confident (or
naive enough) to believe that it holds true. The point seems to be
primarily a matter of projective algebraic geometry, namely the
question if any abstract morphism on a concrete plane curve to the
line $\PP^1$ is induced by a (linear) pencil of ambient curves.
This is either trivially true or trivially wrong, but alas I do
not know the answer due to my failing memory about the foundations
of algebraic geometry.

Your article already helped much as I suffered under the
misconception that your result states that any distribution of 8
points on the empty ovals induces a totally real pencil. Your
statement is much more subtle, yet personally I do not know if
this stronger (universal) form of total reality is wrong! If you
know a counterexample foiling universal total reality I would be
very happy. It could then still be the case that there is some
special sextics for which universal total reality holds true, i.e.
for all octuplets distributed on the empty ovals. (Perhaps reading
more carefully your article, especially the aspect related to
Nikulin-Kharlamov's rigid-isotopic classification already answers
those questions?)

(The newest material of mine (as usual confusing and poorly
organized) occupies Section 27.9--27.11 on pages 373--378. Here I
attempted a topological approach to the existence of octuplets
inducing a totally real pencil, but alas was not able to conclude,
presumably because I know too little on the predestination process
creating the 9th basepoint as a function of the 8 assigned ones.)

Many congratulations again to S\'everine for this breakthrough.
Best regards, Alex

dimanche 3 mars 2013 18:07:57 a new version with small
corrections?

Dear Alexandre, dear other colleagues,  I owe you some apologies:
the Theorem was slightly incorrect, as Alexandre pointed out. I
let you discover this new version, where I have reformulated the
Theorem, and added a few words in the end of the proof. Best
regards,  S\'everine

(04.03.13) Dear S\'everine and the other colleagues,

Many thanks for the new version. In fact, it seems that the main
change is that you now assign the 8 basepoints ON the empty ovals
instead of IN their insides. Rereading my previous message, I
realize that I misstated your original statement and so it is pure
chance that assignation on the ovals turned out to be ``more
correct''.

Your fascinating article gave me new forces to think about the
problem, but alas still without success. For instance, I still do
not know if there exist octuplets (on the empty ovals) failing to
induce a totally real pencil. Of course assigning them in the
insides gives more freedom, but presently it looks to me harder to
ensure total reality. So despite your correction, it could still
be the case (in my modest understanding) that the pencil is total
for all octuplets chosen in the insides  of the empty ovals.
Perhaps you know a counterexample to this strongest form of the
statement?

Many thanks again for the article, which guided much my thinkings.
I hope to send you more exciting news soon, but the whole problem
which you call ``the lost proof of Rohlin'' seems to me still much
out of reach. All the best, Alex

$\bullet$$\bullet$$\bullet$ answer to Alexandre's questions (mardi
5 mars 2013 13:30:42)

Dear Alexandre, dear other colleagues,  let me try to answer the
question with a new formulation.

Assume first that the base~points are distributed {\it inside} of
the empty ovals. Applying your nice ``dextrogyration argument'' to
all nine ovals gives the following lemma:

{\it The pencil is totally real iff 9 lies inside of the non-empty
oval $O$ and outside of the empty ovals.}

If 9 is outside of $O$, the bad cubics are as shown in Figure 2 of
the paper. If 9 is inside of an empty oval $X$, the bad cubics
have an oval passing through the two base~points 9 and X only, and
this oval is entirely contained in the empty oval denoted also
$X$. To get rid of this latter possibility, it suffices to take
the base~points {\it on} the empty ovals.

In ii), I give an explicit description of the pencil, valuable for
any generic choice of the eight base~points {\it inside} of the
eight empty ovals. (It turns out that the only possible
non-generic situation is that of a pencil with a double base point
$9=2$, this means that the points 1, ..8 lie on a nodal cubic with
node at 2.) Recall that 2 is the base point chosen in the extreme
inner oval forming a positive pair with $O$. For this pencil, the
only possibly bad cubics are those with an oval passing through 9
and 2 only. To grant total reality, it suffices to choose the base
point 2 {\it on} the corresponding empty oval, the other
base~points lie arbitrarily in the inside discs of the other empty
ovals. Thus, your conjecture 27.29[=\ref{SRLT:conj}] is true, and
an even stronger result holds for the sextic with six inner ovals.
Best regards, S\'everine

[07.03.13] Little news from Alex, and so many thanks to S\'everine
for the answer

Dear Colleagues,

First many thanks to S\'everine for your very detailed answer
(which I will study in detail tomorrow). Sorry for being always a
bit differed in time due to my lack of internet at home.

I added some material in my loose notes. In Section 28.1--28.2
(pp.384--392), I tried once more to explore the grand programme
that Rohlin might have had in mind, namely total reality and its
connection with his maximality conjecture. As I often said it
seems to me that the missing link could be played by Ahlfors
theorem, or perhaps Rohlin had a grand vision that he could
arrange total reality by purely synthetical processes extending in
all degree the already tricky theorem of Rohlin-Le~Touz\'e in
degree $m=6$. This idea when explored in full looks to me
extremely vertiginous, but its net impact would be a sort of upper
bound upon the complexity of Hilbert's 16th problem, and in some
sense subsume all prohibitions (\`a la Gudkov et cie.) to the
paradigm of total reality. All this necessitates to be made much
more precise, but I [have] attempted to make a psychoanalysis of
what Rohlin may have had in the brain, without that he himself
ventured to put it on the paper due to his own modesty and
pragmatism.

Next I discovered the little Theorem
28.7[=\ref{Thom-Ragsdale:thm}] (p.393), which is just a matter of
making explicit the consequence of Thom's conjecture
(=Kronheimer-Mrowka theorem) as it pertains to Hilbert's 16th
problem. The result is the lovely estimate\footnote{[26.03.13] It
turned that this is wrong, cf. Fiedler's letter dated [21.03.13].}
$\chi \le k^2$ for a curve of type~I and degree $2k$. With this I
realized that my former counterexample (with the scheme $20$ in
degree 8) to CCC(=collective contraction conjecture) is actually
killed by Thom, and realized (later only!!) that it is also killed
by Rohlin's formula. So CCC is again resuscitated but probably not
for long!?

Then I tried to make a comparative study of Rohlin's formula
versus the Thom obstruction. It seems that the latter is often
implied by Rohlin's formula, but not always. More in Section 28.4
(p.393). It seems however that at least for degree $m \ge 10$
there is some cases where Thom really affords new information not
covered by Russian congruences or Rohlin's formula (cf. Thm 28.11,
p.396). Finally using the Gudkov table in degree 10 (=Fig.148 on
page 395), I got some naive hope to disprove the Rohlin maximality
conjecture, but this quickly turned into disillusion (cf. Point 3
on p.396--397).

Sorry for all these messy remarks, yet I found the r\^ole of Thom
quite pleasant. I am sure that this is not new, and that I read it
somewhere, but again could not recover where precisely. (I thought
it was in Degtyarev-Kharlamov 2000's survey but apparently not,
though Kronheimer-Mrowka is alluded to.) If you remember some
anecdotes about the r\^ole of Thom's conjecture in Hilbert's 16th
problem, and who puts it first into action as a such, I would be
extremely happy to insert your remarks in my (messy) survey.

Thanks a lot for the attention, Best regards, Alex

$\bullet\bullet\bullet$ Thomas Fiedler wrote (samedi 9 mars 2013
17:32:59)

Dear Alexandre,

I am no longer in business in this field, but let me just make
some remarks which could be perhaps helpful.

The $M$-curve of degree 10 mentioned in your Thm
28.11[=\ref{French-scheme:thm}] is in fact ruled out by Rokhlin's
formula. I think that you have mixed $Pi^+$ with $Pi^-$. In a
positive couple the orientations are just opposite. So, four
nested ovals can contribute at most $+2$ to Rokhlin's formula.

It is an interesting idea to apply the Thom conjecture to real
algebraic curves. To my knowledge the only new result obtained
this way is contained in G. Mikhalkin ``Adjunction inequality for
real algebraic curves''.

Let me formulate the problem (which exists certainly already
somewhere).

GENERALIZED THOM PROBLEM.

Let $X$ be a simply connected smooth closed $4$-manifold and let
$h$ be a non trivial integer $2$-dimensional homology class. Let
$F$ be a smoothly embedded oriented surface which represents $h$
and such that the components of $F$ represent classes which are
linearly independent over $Z/2Z$. What is the maximal Euler
characteristic of $F$?

In the complex projective plane this boils down to the Thom
problem, because evidently the surface F has to be connected in
this case. However, it becomes interesting in a more general
complex surface. It is an easy matter to make a non connected
surface $F$ connected but the opposite is quite hard. It is
equivalent to finding an embedded ``membrane'' with trivial normal
Euler number. I don't know wether Seiberg-Witten theory nowadays
can give a sufficient criterium to ensure the existence of such a
membrane (as stretching the neck of the surface $F$). But it seems
to me that this is the place to look at.

If one considers as $X$ the double cover of the projective plane
ramified in the complexification of an $M$-curve of even degree
then one can consider as $F$ the fix point set of one of the two
induced anti-holomorphic involutions on $X$. We know
$2b_0(F)+b_1(F)$ from Harnacks equality. Hence the Euler
characteristic $2b_0(F)-b_1(F)$ is maximal when $b_1(F)$ is
minimal, i.e. the numbers $p$ or $n$ of the real curve are
maximal. So this problem is closely related to the still open
Ragsdale conjecture for $M$-curves.

Best regards, Thomas

$\bullet$ [Gabard, 09.03.13, ca. 21h00] Does the Gudkov hypothesis
(mod 8) reduces to Rohlin's (complex orientation) formula

Dear Colleagues,

On reading recently the Degtyarev-Kharlamov 2000 survey, I learned
the (simple) fact that the Arnold congruence mod 4 (weak Gudkov
hypothesis) can be reduced to Rohlin's formula. I wrote down a
proof of this simple issue in Lemma
24.16[=\ref{Rohlin-implies-Arnold:lem}], on page 219 of the messy
survey.

Of course then I wondered if a sharpened combinatorial argument
taking into account the signs distribution on the edges of the
``Hilbert tree'' (encoding the distribution of ovals) prompted by
Rohlin's complex orientations could likewise subsume the Gudkov
hypothesis to Rohlin's formula. The combinatorics becomes much
more messy and I lack a good idea on how to exploit the $M$-curve
assumption (apart from the dubious idea of using the
dextrogyration allied to a totally real pencil, but this looks
somewhat ad hoc!!!?). If feasible (so or otherwise), the net
impact would be that the highbrow topology used in the
Rohlin-Rohlin/Atiyah-Singer-Marin proofs could be replaced perhaps
by basic algebraic geometry. This issue is merely didactic of
course, yet it looks perhaps technically challenging to implement
this modest ``dream''. My (unsuccessful) attempt to tackle this
reduction is given in Section 24.5, pages 220--223.

Of course, I am sure that you already tried hard along this way
and that this is a pot-pourri naive problem. So I write you the
letter, only in the hope that someone already worked this
reduction successfully, though I doubt (as otherwise it would
certainly have been mentioned in the Degtyarev-Kharlamov 2000
survey).

Many thanks again for all your attention and all your precious
hints, Alex

PS: I apologize much to S\'everine for not having yet found the
time to study carefully enough the last brilliant explanations,
but look forward doing so in the best delay.

PPS: Many thanks also to Thomas for the brilliant answer on Thom
and Mikhalkin, Seiberg-Witten, etc. I will include your letter in
(the next version of) my notes so that anybody can contemplate it.

$\bullet$ [11.03.13]

Dear Thomas and the other Colleagues,

So many thanks to Thomas for having catched my mistake regarding
Thom versus Rohlin (existence of schemes prohibited by Thom but
not by Rohlin's formula). I have written down a corrected version
of the Theorem (whose first clause is I think still true despite
Thomas's corrigendum) as Thm 30.14 on page 382. I hope this time
it is correct! The subsequent Thm 30.15 and Lemma 30.16 describe
larger family of such schemes. Alas, what I have written is not
extremely appealing, but I hope still readable for such experts as
you are (though I confess that it is not extremely exciting!). The
philosophy is just of course that there is no subsumation of
Thom's estimate $\chi\le k^2$ to Rohlin's formula.

Many thanks again to all for all your kind answers, and especially
to Thomas for taking care to bring me on the right track!

Best regards, Alex

PS: I have drifted the Section of your e-mails at the end of the
text, in particular the last letter of Thomas is to be found on
page 435. Alas, I had not yet the time to digest its full swing,
but the idea looks very promising!

$\bullet\bullet\bullet$ mardi 12 mars 2013 13:56:45

Dear Alexander,

sorry, but all your $M$-schemes of degree $10$ in Thm 30.14 and
30.15 have $n=2$ and are ruled out simply by Petrovskis
inequality. I don't think that genus bounds give anything new for
real schemes alone but they definitely do so for configurations of
several real curves. Just take a look on Mikhalkin's paper.

Best regards Thomas

[16.03.13]

Dear Geometers,

Many thanks again to Thomas for his former correction. I got so
drifted in a sort of cuneiform formalism of trees with signs but
mostly lost myself into dubious combinatorics of (what I call) the
Rohlin tree (=Hilbert's nested tree with signs materializing of
course Rohlin's complex orientations). Ultimately after several
basic combinatorial mistakes, I did NOT even succeeded in finding
a Caucasian $M$-scheme where Rohlin's formula is stronger than
Thom $\chi \le k^2$. This problem is discussed in Sections 30.7
and 30.8 (pages 393--402). It seems to me that there is (for
$M$-curves at least) a certain concomitance between Rohlin and
Thom, i.e. you cannot corrupt Rohlin's formula without corrupting
simultaneously Thom's bound. So maybe it is reasonable to
conjecture that if Thom's estimate is fulfilled then Rohlin's
equation is always soluble for a suitable distribution of signs on
the edges of the tree. Sorry that my summary is vague as I do not
myself understand properly what happens.

Today, I switched on a somewhat more pleasant arithmetical problem
exposed in Section 30.6, p.391--392. This is hopefully more
readable, and you surely studied this a long time ago. Here the
question is to find an $M$-scheme without nesting of high-degree
(say $m\ge 6$ assuming zero-knowledge). As we know, Hilbert
posited the intuition that $M$-curves are forced to exhibit
nesting.  (In fact on reading this afternoon more carefully
Hilbert's text, he is not so categoric but let us assume so, to
add some suspense to our story!) In view of Rohlin's formula
(which forces the number of ovals of an unnested dividing curve to
be a square) and the Gudkov congruence mod 8, one is invited to
ask when the Harnack bound $M=g+1$ of a degree $2k$ curve is a
(perfect) SQUARE? This leads to a little arithmetical problem
which admits a nontrivial solution at $k=17$ for $M=529$ which is
by a lucky stroke equal to $23$ squared. (Of course I found this
just by an ``exhaustive'' tabulation.)  (Note at this stage that
$17=16+1$ yet another subconscious coincidence in Hilbert's
numbering!!??) It turns out moreover that the Gudkov congruence is
then fulfilled! So this leads me to ask if, you in Russia, knew a
way (prior to Thom-Kronheimer-Mrowka 1994) to prohibit this
$M$-scheme of degree $2k=34$ with 529 ovals (all lying outside
another, i.e. no nesting). Personally, I would be very happy if
someone can tell me an answer (in case I did not foiled the
arithmetics!!!) Incidentally, I am so ignorant in that field that
I do not know how to solve in general the Diophantine problem of
the quadrature of Harnack's bound, i.e. for which $k$ is Harnack's
bound $M=g(2k)+1$ a square (=Problem 30.19, on page 391).

Finally, this morning, building upon S\'everine's (elementary)
remark on $M$-quintics totally real under a pencil of cubics
(nothing so hard as the Rohlin-Le~Touz\'e theorem for sextics), I
had the idea to extend the construction of satellites to curves of
odd degrees, and so found a scheme of degree 10 (the 2nd satellite
of Harnack $M$-quintic). This schemes, according to the philosophy
of total reality, should be of type~I (stability under
satellites). This material is exposed in Section 29.8 page
352--353. The scheme in question has 13 ovals one of them
enclosing 6 nests of depth 2 (cf. Fig 149 p.353). Assuming that
this scheme is of type~I, it could be a counterexample to Rohlin's
maximality conjecture, in case someone ever saw a curve of degree
10 enlarging it. Otherwise, more in line with Rohlin's philosophy,
this scheme being totally real in some explicit way it should be
maximal and it results a myriad of prohibitions on all schemes
enlarging it! If you know some experimental construction \`a la
Viro-Itenberg that may help to see clearer, I would be extremely
thankful. Of course, it would also be very interesting to know if
S\'everine thinks that there is some good chance that total
reality holds true for this satellite.

Many thanks for your attention, and I hope my questions are not
too ill-posed for such experts as you are! All the best, Alex

PPS: Many thanks for the message meanwhile received from Thomas
which I hope make my message not to obsolete. Best regards.

PPPS: Dear Viatcheslav, I try twice to send you my message as I
received delivery failure notification! I hope I am not
overloading your mail-box. I try now with a zipped file hopefully
toujours lisible pour vous. Amiti\'es. Alexandre

$\bullet\bullet\bullet$ mardi 19 mars 2013 17:54:39

Dear colleagues, my note "totally real pencils. . ." is now
available on the archiv, here is the reference:
http://arxiv.org/abs/1303.4341  Best regards,  S\'everine

$\bullet$ [19.03.13, ca. 19h30] An Alsatian scheme and planning to
put a version on arXiv in ca. 10 days (please confirm me if you
accept the insertion of your kind letters in my messy survey)

Dear Colleagues,

Many thanks to Thomas again for pointing out the marvellous
Petrovskii's inequalities (which I confess I had not assimilated
properly before, despite all the brilliant surveys available).
Shame on me! So of course most of the questions of my previous
message were rather stupid. In particular that relating to the
arithmetical problem was completely ill-posed as I missed to use
the full punch of Rohlin's formula. I apologize much for all these
inconsistencies.

Meanwhile I have attempted to find what I call an Alsatian scheme
where Thom is stronger than the conjunction of Petrovskii 1933/38,
Gudkov hypothesis (1969--72=Rohlin's semi-proof) and Rohlin's
formula (1974--78). I think that such an Alsatian scheme do exist
(cf. Thm 30.22 on page 392). This answers (I hope correctly) a
question raised by Thomas in the message reproduced right below,
where if I understand it well Thomas expected that Thom says
nothing new already known earlier in Russia.

Then in view of Thomas' stimulating messages and also the
marvellous Itenberg-Viro 1996 Math. Intelligencer article, I
adventured slightly in the province of Miss Ragsdale. Here first I
learned that the Ragsdale conjecture for $M$-curves really boils
down to Thom (at least one-half thereof)\footnote{[26.03.13] This
is again a misconception of mine, since the estimate $\chi\le k^2$
of Theorem~\ref{Thom-Ragsdale:thm} was based on the erroneous
supposition that the Arnold surface is always orientable, cf.
again Fiedler's letter dated [21.03.13].}. I found strange that
Kronheimer-Mrowka 1994, was not cited in Itenberg-Viro 1996
probably due to backlog reasons (by Intelligencer)\footnote{For
the same reason as the previous footnote this comment is not
pertinent anymore.}? Further it seems to me that there is a
misprint in the statement of the RAGSDALE CONJECTURE ON $M$-CURVES
(on p.24 of Itenberg-Viro 1996). At any rate, I tried to write
down my own naive account on Ragsdale in Sect.
30.4[=\ref{Ragsdale-conj:sec}] (p.380--384) which details more
slowly what I understood (hopefully correctly). In fact, I wonder
at my (premature) stage if the full Ragsdale conjecture could not
follow \`a la Thom via the lower estimate $-k^2 \le \chi$ which I
very naively conjecture to be true for all dividing curves? (I
confess that I made no experiments even with the methods of
constructions I am aware of, i.e. Hilbert-Harnack).

So as you see, life really starts becoming exciting. Alas I fear
that I will be strongly interrupted due to editorial reasons of
our Journal L'Enseign. Math. in Geneva of which I am the
TeX-editor. So perhaps I should stop thinking and try to polish a
bit the big mess I produced (during ca. 10 days), before putting
it on the arXiv (prior to my long editorial job). Alas, I had
already great difficulties to submit the previous version of my
file (due to size limitation policy of arXiv), but received an
exceptional derogation to do so. Now my text on Ahlfors + Rohlin
is twice as large as it was before (600 vs. 300 pages) so there is
little chance that I get accepted. I will try to see if an arXiv
administrator looks optimistic. Meanwhile, I would be very happy
if you all confirm me that you accept the integration of your
marvellous letters in my modest text. Of course answer me only in
case of objection. I plan to submit the new version ca. the 1
March 2013. Of course my text is so messy that your letters alas
are not properly pushed into evidence, but I expect to produce a
better text in the next months or years.

Many thanks again to all for your letters, indulgence, patience
and kind answers!

Best regards, Alex

$\bullet$$\bullet$$\bullet$ {\it Reproduction of Fiedler's former
letter (already above but the other colleagues did not saw it)}
mardi 12 mars 2013 13:56:45

Dear Alexander,

sorry, but all your $M$-schemes of degree $10$ in Thm 30.14 and
30.15 have $n=2$ and are ruled out simply by Petrovskis
inequality. I don't think that genus bounds give anything new for
real schemes alone but they definitely do so for configurations of
several real curves. Just take a look on Mikhalkin's paper.

Best regards Thomas

$\bullet\bullet\bullet$ mardi 19 mars 2013 22:17:27

Dear Alexandre,

sorry again, but your curve has $p=50$ and is ruled out by
Arnold's inequality:
$p\le 3/2 k(k-1) + 1 + n_-$, which is $47$ in this case. In fact
Arnold's inequalities are by fare the strongest result in the
whole field.

Best regards Thomas

$\bullet$ [20.03.13] ``Ren\'e Thom sur son 31?''

Lieber Thomas (and the other colleagues),

Es wird jadoch immer spannender, oder wat? Many thanks for the new
challenge raised by Thomas (which I received as a special gift on
my birthday date, today 20st of March). As Thomas demonstrated
yesterday evening (cf. message reproduced below), my example of
Alsatian scheme (where Thom is stronger than the Soviet Red army)
(cf. Thm 30.22[=\ref{Alsatian-schemes:thm}], p.392) collapses
under the strong Petrovskii inequality of Academician Vladimir
Igorevich Arnold. It took me just some few hours (as I was quite
tired) to find a stronger candidate of Alsatian scheme where Thom
looks stronger that the conjunction of (strong) Petrovskii-Arnold
1971, Gudkov 1969 (proved by Rohlin 1972-Marin ca. 1977), and
Rohlin's formula (1974--78). For the exact statement cf. Thm
30.25[=\ref{Alsatian-scheme-Thom-strong-Petrov-Arnold:thm}] on
page 395. I hope I made no mistake (I checked the details twice).
By (3-fold) experience, I am quite confident that Thomas will find
a new obstruction in his pockets killing this new example.

Hence many thanks again to Thomas for all this precious guidance
that oriented much my modest working. As I said I am now under
strong temporal constraints, and will not be able to pursue any
thinking for a long period of circa 1 month. So please feel free
to elaborate more upon the direction indicated by Thomas if it
looks hard stopping the inertia. Myself find the Alsatian topic
quite pleasant, yet it would (I presume) be interesting to
penetrate deeper in the geographical question (by really
understanding the diagrammatic impact of Petrovskii-Arnold,
Ragsdale, etc. upon the higher Gudkov's tables of periodic
elements) I made several pictures of big pyramids (e.g Fig.
153[=\ref{Degree10:fig}] p. 386) and that could be a first step
toward understanding better what happens. Of course all this must
be familiar to you, yet personally I still lack a good algorithm
to make good (color-plates) maps evidencing ``all'' obstructions.

My main worry now, is how to publish (and polish) the 600 pages
long article that I have produced, especially in view of stringent
size restrictions imposed by the arXiv administrators.

Many thanks that you seem to all approve the insertion of your
letters in my survey, which looks to me essential as you
influenced much my (chaotical) trajectory. Of course in the future
we will have the occasion to clean better in case some bugs are
detected, which is quite unlikely (apart of course in my own
letters full of inconsistencies).

All the best, Alex

$\bullet\bullet\bullet$ jeudi 21 mars 2013 07:51:45

Dear Alexandre,

unfortunately your Theorem 30.7[=\ref{Thom-Ragsdale:thm}] is
wrong. Half of the curve together with R is usually called
Arnold's surface. It is an orientable surface iff R has an
orientation which induces the complex orientation on its boundary,
i.e. the real curve. Hence if there is a negative pair of ovals in
in the boundary of one component of R then Arnold's surface is not
orientable.

By the way it is well known and easy to prove that if Arnold's
surface is orientable then
$p-n=k^2$.

Best regards, Thomas

$\bullet$ [21.03.2013, ca. 22h00] Dear Thomas and the other
Colleagues,

Many thanks to Thomas for having spotted out my fundamental
mistake. I apologize much hence for all the dubious letters that I
sent you the former days. I hope I will still be able to repair a
bit the situation in the next weeks, though this will require
dramatic changes in my messy text! So many thanks again to Thomas
for having detected this great Harnaque!!!

All the best, Alex

$\bullet\bullet\bullet$ vendredi 22 mars 2013 07:54:47

Dear Alexandre,

don't worry. Just read Mikhalkins paper and find other interesting
applications of his method. Having a reducible dividing curve you
can switch the canonical orientation of exactly one component.
This allows to construct immersed surfaces which are far from
being complex curves. Then indeed genus bounds start to work.

Best regards, Thomas

$\bullet$ [24.03.13] Is there a simpler counter-example to
``Gabard-Thom'' than via Itenberg-Viro mirabilis $(M-2)$-curve of
degree 10 disproving Ragsdale

Dear Colleagues,

As pointed out by Thomas (cf. my last message) my proof of the
``Gabard-Thom'' estimate $\chi \le k^2$ for all dividing curves of
degree $2k$ was highly fraudulent as it was based on the
(erroneous) assumption that the Arnold surface (=Klein's half
married with Miss Ragsdale) is always orientable. Of course Thom
has still something to say on Hilbert's 16th e.g. in the very
special case of no nesting. As a historical curiosity one can
notice that the elementary case due to Kervaire-Milnor (1961) of
Thom in degree 3 gives in my opinion the first purely topological
proof of Hilbert's Ansatz of nesting for $M$-sextics. Prior to
that we had only Hilbert-Rohn (stratificational as explained by
Eugenii), and Petrovskii 1933/38 which involves the
Euler-Jacobi-Kronecker stuff (interpolation formula). Of course
Kervaire-Milnor is merely Rohlin's early work 1951 disguised!

Next I was a bit puzzled by Thomas's claim that my ``Gabard-Thom''
theorem is wrong, since I was not able to find an explicit
counter-example. I tried a while with elementary constructions \`a
la Harnack-Hilbert but could not find a single counterexample. So
I found the pleasant plates Fig.
152-153-154-155[=\ref{HilbGab1:fig}--\ref{HilbGab4:fig}] on pages
382 and ff. Here we get nice infinite families of curves with
$\chi=k^2$ that were surely known to Hilbert and Ragsdale. Alas
this inclined more toward thinking that Gabard-Thom is sharp,
rather than disproving it. Then of course I had the idea to take a
closer look to Itenberg-Viro 1996's article (disproof of
Ragsdale), which came like a deliverance and killed in the same
stroke the Gabard-Thom dubious estimate. I used the
Kharlamov-Marin congruence to check type~I of the Itenberg-Viro
curve (which I reproduced on Fig. 156[=\ref{Itenberg:fig}],
p.387), though there is surely a more elementary argument \`a la
Fiedler-Viro-Itenberg-Parenti. Many congratulations by the way to
Ilia and Oleg for this geometric paradise, which I contemplated
with much pleasure and extreme admiration yesterday evening!!! A
naive question of mine, is whether we really need to resort to the
patchwork method to disprove Gabard-Thom, in the sense that
perhaps I missed (a non-maximal) counterexample via
Harnack-Hilbert.

Finally this marvellous curve (Itenberg-Viro's) poses again the
question of total reality, this time under a pencil of septics (in
general I conjecture total reality for adjoint curves of order
$m-3$ when it comes to $(M-2)$-curves). The question is where to
assign exactly the anchor basepoints! I tried some guesses in Sec.
31.3[\ref{Galton-brett:sec}] (p. 388). This deserve perhaps much
deeper investigations than what I am presently able to do. Perhaps
S\'everine and Thomas already have some good ideas. My own dream
is that there should be a combinatorial recipe telling one from
the sole knowledge of the Rohlin tree (=Hilbert's one with signs
given by complex orientations) where to assign basepoints. One
should perhaps imagine the Rohlin tree as a Galton-Brett=table
with billiards balls falling downwards to the empty ovals and
perhaps stabilizing at some other (unstable) equilibriums when
they meet a ``nail'', which in first approximation could be a
hyperbolic oval negatively charged on the edge right above it.

All my thanks again to Thomas for having catched my fundamental
mistake at the right moment. I hope to polish a bit the text
during the next days before submitting it to the arXives. Then I
must move toward doing more boring editorial duties for the Swiss
journal L' Ens. Math. (close to collapse by the way).

Good Sunday to all, and best regards, Alex

PS: Many thanks also for Thomas's last letter which I just
discover now. It is an invitation to read Mikhalkin's article,
which alas I had not yet the time to do properly...


\section{Synoptic tabulations}

This is an attempt to gather information scattered through the
literature. The first synoptic project  compiles a list of
nomenclatures. A second tabulation reflects how Ahlfors work
(existence of circle maps) has been appreciated by subsequent
workers of a slightly dissident nature in the sense that they cite
conjointly other sources.

\subsection{Nomenclature project}

This section tries to get sharp lower bounds
on the basic nomenclature of our topic. As Poincar\'e tried to
convince Felix Klein ``{\it Name ist Schall und Rauch\/}'' (cf.
e.g. Klein 1923 \cite[p.\,611]{Klein-Werke-III_1923}), but it is
somehow pleasant to investigate the historical background of some
jargons to
use them hopefully appropriately.

$\bullet$ (Gauss 1825/1844, F.\,T. Schubert earlier?) {\bf
Conformal mapping=konf\-orme Abbildung}, maybe the first
non-trivial result is to be found in Gauss 1825 \cite{Gauss_1825},
yet the word ``conformal'' itself appears in Gauss 1844 in the
first paper on higher geodesy: ``{\it ich werde daher dieselben
conforme Abbildungen oder \"Ubertragungen nennen, indem ich diesem
sonst vagen Beiworte eine mathematisch scharf bestimmte Bedeutung
beilege\/}'' [Werke IV, p.\,262]. As noted in Struik 1933
\cite[p.\,164]{Struik_1933} (via Cantor), the word ``conformal''
is already used prior to Gauss by F.\,T. Schubert in ``De
projectione sphaeroidis ellipticae geographica'', {\it Nova Acta
Petr.}, p.\,130--146.

$\bullet$ (1865) {\bf Riemann surface}, maybe first coined by
C. Neumann  1865 \cite{Neumann_1865}, followed by L\"uroth
1871 \cite{Lueroth_1871}, Clebsch 1872 \cite{Clebsch_1872},
Klein 1874--76 \cite{Klein_1876}, \cite{Klein_1876}, Clifford
1877 \cite{Clifford_1877} and then too many to record.

$\bullet$ {\bf Berandete  (Riemannsche) Fl\"achen, Compact
bordered Riemann surfaces, finite Riemann surface, membranes}. The
first appellation appears often in Klein 1882 \cite{Klein_1882}
(reprint in Klein 1923 \cite[p.\,569,\S 23]{Klein-Werke-III_1923})
and others. The second appellation is coined and popularized in
Ahlfors-Sario's 1960 book \cite{Ahlfors-Sario_1960}, whereas the
third competing name is used in Schiffer-Spencer's book of 1954
\cite{Schiffer-Spencer_1954}. The term membrane also occurs (in
this context) by Klein in his lecture notes.

$\bullet$ (1907?) {\bf Uniformization} probably a coinage of
Poincar\'e. In 1883, just the word ``fonction uniforme'' appears
and the word ``uniformization'' as a such, came in vogue ca. two
decades latter in Poincar\'e 1907 \cite{Poincare_1907} and Koebe
1907 \cite{Koebe_1907_UbaK1}.

$\bullet$ (1908) {\bf Kreisnormierungsprinzip} coined and
proved (in fairly general special cases: finite connectivity
and symmetric under complex conjugation) by Koebe in 1908
\cite{Koebe_1908_UbaK3}.

$\bullet$ (1912) {\bf Schwarz's lemma}. The coinage as a such
appears first in Carath\'eo\-dory 1912
\cite{Caratheodory_1912}, but already published in the modern
fashion in 1907 by the same writer \cite{Caratheodory_1907},
acknowledging the argument of E. Schmidt.

$\bullet$ (1916) {\bf Extremal problems=Extremalprobleme} used
in function theory by Bieberbach 1916 \cite{Bieberbach_1916}.

$\bullet$ (ca. 1914) {\bf Circle mapping=Kreisabbildung.} This is
used (at least) since Bieberbach 1914
\cite[p.\,100]{Bieberbach_1914}, Koebe 1915 \cite{Koebe_1915},
Bergman[n] 1922 \cite[p.\,238]{Bergman_1922}, Bochner 1922
\cite[p.\,184]{Bochner_1922}, with the English translation
appearing first in Garabedian-Schiffer 1950
\cite{Garabedian-Schiffer_1950}.

$\bullet$ (ca. 1975--1977) {\bf Total reality, ``total reell'',
etc.} The adjective ``total reell'' (totally real) is first used
(in the generality) in Geyer-Martens 1977 \cite[p.\,101,
p.\,103]{Geyer-Martens_1977}, where the connection with Ahlfors
theorem is made explicit along the line already suggested in
Alling-Greenleaf 1969 \cite{Alling-Greenleaf_1969}. Geyer-Martens
ascribe (cf. p.\,101) the concept of total reality (when
paraphrased in the language of field extensions) to J.\,T. Knight
1969 \cite{Knight_1969}. Somewhat earlier in lesser generality of
Galois coverings ``total reell'' appears already in Martens 1975
\cite{Martens_1975}. Meanwhile I think that ``totally real'' is
quite widespread, especially in the growing field of real
enumerative algebraic geometry (e.g. works by Sottile).

$\bullet$ {\bf Erster Art=Type I, orthosymmetric and dividing
curves} All this jargon is due to Klein. More precisely, Erster
Art appears in Klein 1876 \cite{Klein_1876} and is much used in
Russian literature (meanwhile diffusing in the west), cf. e.g.
Rohlin 1978 \cite{Rohlin_1978} (and Gudkov 1974
\cite{Gudkov_1974/74}???). I remember some irony of Grisha
Mikhalkin during a talk by Orevkov, where he found Klein's
subsequent jargon ``orthosymmetric'' (first in print in Weichold
1883 \cite{Weichold_1883}) quite awkward. Yet, Klein himself
turned to be quite proud of this more intrinsic coinage. This
turned to be quite influential, adhered by eminent workers like
Koebe, J. Douglas, etc. albeit quite in desuetude today. The
reason is mostly due to synonyms like dividing (or separating)
curves.

Summarizing, the following words are used resp by:

$\bullet$ ``Erster Art=Type I'' first coined in Klein 1876
\cite{Klein_1876}, and adhered to by Rohlin 1978
\cite[p.\,90]{Rohlin_1978}, and then much of the subsequent
Russian literature,

$\bullet$ orthosymetrisch Klein ca. 1882 (lectures), adhered to
Weichold 1883 \cite{Weichold_1883} (first occurrence in print),
then Klein 1891/92 \cite{Klein_1891--92_Vorlesung-Goettingen}
(Vorles. G\"ottingen), and followed by Koebe 1907
\cite{Koebe_1907_UrAK} (etc.), Fatou 1930 (in Appel-Goursat 1930
\cite{Appell-Goursat-Fatou_1930}), Julia 1932 \cite{Julia_1932},
Douglas 1936 \cite{Douglas_1936-Some-new-results}--1939
\cite{Douglas_1939-The-most-general}, etc.

$\bullet$ ``zerteilend vs. nichtzerteilend'' in Fiedler 1981
\cite[p.\,7]{Fiedler_1981}

$\bullet$ ``divide'' alone is briefly mentioned in Arnold 1971
\cite{Arnold_1971/72} (yet only as a property of $M$-curves)

$\bullet$ ``dividing curves'' is used by Wilson 1978
\cite[p.\,66]{Wilson_1978}, Viro 1986/86
\cite[p.\,58]{Viro_1986/86-Progress} Kharlamov-Viro 1988/91
\cite[p.\,359]{Kharlamov-Viro_1988/91}, Gilmer 1991
\cite{Gilmer_1991}, Degtyarev-Kharlamov 2000
\cite[p.\,736,\,737]{Degtyarev-Kharlamov_2000}.

$\bullet$ ``separating curves'', occurs in Fiedler 1982/83
\cite[p.\,162]{Fiedler_1982/83-Pencil}, Dubrovin 1983/85
\cite{Dubrovin_1983/85}, Nikulin 1983/84 \cite{Nikulin_1983/84},
Benedetti-Risler 1990 \cite{Benedetti-Risler_1990}, Natanzon 1990
\cite{Natanzon_1990/90} or 1999
\cite{Natanzon_1999-Moduli-real-alg-surf.superanal-differ-spinors},
Coppens 2011 \cite{Coppens_2011},

$\bullet$ ``courbes s\'eparantes'' in Marin 1979
\cite{Marin_1979}, ``courbe qui s\'epare sa complexifi\'ee in
Marin 1988 \cite{Marin_1988}, Gabard 2006 \cite{Gabard_2006} (alas
in French it sounds strange to say ``courbe divisante'').

 The following concept is a priori foreign to our survey,
albeit it would be interesting to see if the methods of
Gr\"otzsch-Teichm\"uller are of some relevance to the Ahlfors
mapping of 1950. This is another mathematical question, but here
we content ourselves with a point of terminology:

$\bullet$ (1928/1935) {\bf Quasiconformal mappings=quasikonforme
Abbildungen}. This nomenclature is usually ascribed to Ahlfors
1935, who however could not remember precisely from where he
borrowed the jargon, according to K\"uhnau 1997
\cite[p.\,133]{Kuehnau_1997}), which is worth quoting:

\begin{quota}[K\"uhnau 1997] {\small \rm
Der Name Gr\"otzsch ist wohl bei vielen vor allem mit der Theorie
der quasikonformen Abbildungen verbunden, die er ab 1928
begr\"undete. Die Bezeichnung ``Quasikonforme Abbildungen'' wurde
allerdings erst sp\"ater von L.\,V. Ahlfors eingef\"uhrt.
(Freilich sagte mir Ahlfors Februar 1992 in Oberwolfach, da{\ss}
er diese Bezeichnung bei jemandem ``gestohlen'' habe, er wisse nur
nicht mehr bei wem.)

}
\end{quota}

Maybe it contributes to the question to remember that the jargon
``{\it quasikonform}'' appears already in 1914, und zwar bei
Carath\'eodory 1914 \cite[\S 16]{Caratheodory_1914}(=page 294 in
the pagination of the Ges.\,Math.\,Schriften,\,Bd.\,3).

\subsection{Dissidence from Ahlfors}\label{dissident:sec}

[31.08.12] Sec.\,\ref{Ahlfors-proof:sec}  attempted to present
Ahlfors' proof in full details, but failed to digest  the details.
%
This deplorable issue motivated us to tabulate a list of
``dissident'' authors, who instead of quoting the original source
Ahlfors 1950 \cite{Ahlfors_1950} adhered
to subsequent
treatments. Two accounts emerge with high rating, namely:

$\bullet$ Heins 1950 \cite{Heins_1950}

$\bullet$ Royden 1962 \cite{Royden_1962}

Of course, our ``dissident'' writers (quoting beside Ahlfors some
derived product) never (as far as I know) criticizes directly the
1950 work of Ahlfors. At least there dissidence may suggest that
themselves were not completely happy with (resp. convinced by) the
original text  finding more convenient another implementation.
Albeit nobody ever expressed  frontal objections against Ahlfors
1950 \cite{Ahlfors_1950}, it is not to be excluded (yet of very
low probability ca. $10^{-14}$) that somebody once detected some
little bug, explaining perhaps the numerous initiatives to reprove
Ahlfors' result from different viewpoints. (We mention again the
articles by Mizumoto 1960 \cite{Mizumoto_1960} and Kuramochi 1952
\cite{Kuramochi_1952} (undigest?), and refers for a extensive
tabulation of such initiatives to the circled item of
Fig.\,\ref{Map:fig}).

Here is
a sample of dissident authors (grouped according to their
preferred source) with relevant extracts in ``\dots'':

\medskip
{\sc Voting for Heins 1950:}

$\bullet$ Stout 1972 \cite[p.\,345]{Stout_1972}: ``\dots a
theorem of Ahlfors [2](=Ahlfors 1950 \cite{Ahlfors_1950})
shows that $\cal H(R)$ contains many inner functions. (See
also the elegrant [sic!] construction of Heins [15](=Heins
1950~\cite{Heins_1950}) as well as the earlier paper of
Bieberbach [3](=Bieberbach 1925 \cite{Bieberbach_1925}) which
deals with the case of planar domains.)''

$\bullet$ Khavinson 1984 \cite[p.\,377]{Khavinson-Dimitri_1984}:
``The following theorem is a classical result of Bieberbach and
Grunsky (see [6](=Golusin 1952/57~\cite{Golusin_1952/57}),
[8](=Grunsky 1978~\cite{Grunsky_1978})). For a different approach
due to L. Ahlfors, see [1](=Ahlfors 1950~\cite{Ahlfors_1950}). Our
proof, although discovered independently, is almost the same as
that due to M.~Heins in [11](=preprint=now published as Heins
1985~\cite{Heins_1985-Extreme-normalized-LIKE-AHLF}) or H.~Grunsky
in [8](=Grunsky 1978 \cite{Grunsky_1978}). {\sc Theorem 3.} {\it
Let $\zeta_1, \dots, \zeta_n$ be arbitrary fixed points on
$\gamma_1, \dots,\gamma_n$ respectively. Then, for each $j$,
$\phi(z)$ is the unique function giving a conformal mapping of $G$
onto an $n$-sheeted right half-plane such that
$\phi(\zeta_j)=\infty$, for all $j$, $\phi(z_0)=1$.}

$\bigstar$ admittedly, this Khavinson' extract in not hundred
percent pertinent to our present purpose inasmuch as the
Bieberbach-Grunsky theorem is confined to the planar case.

\medskip
{\sc Voting for Royden 1962:}

$\bullet$ Stout 1965 \cite{Stout_1965}: ``In order to
establish our result, we shall need to make use of a result of
Ahlfors [1](=Ahlfors 1950 \cite{Ahlfors_1950}). (For an
alternative proof, one may consult Royden [15](=Royden
1962~\cite{Royden_1962}.)

Theorem~3.1 {\it There exists a function $P$ holomorphic on a
neighborhood of $\bar R$ which maps $R$ onto the open unit
disc in an one-to-one manner for some $n$ and which satisfies
$\vert P \vert =1$ on $\partial R$.}''

$\bigstar$ Of course the above ``one-to-one'' is a typo to be read
as ``$n$-to-one''.

$\bullet$ Alling 1966 \cite[p.\,346]{Alling_1966}: ``Finally,
I am indebted to Professor Royden for his excellent paper,
{\it The boundary values of analytic and harmonic functions},
[24](=Royden 1962 \cite{Royden_1962}), which not only gave a
new proof of the existence of the Ahlfors' map, but also gave
generalizations of the classical boundary value theorems over
the disc. \dots''

$\bullet$ Stout 1966/67 \cite[p.\,366]{Stout_1966/67}: ``Let
$R$ be a finite open Riemann surface whose boundary $\Gamma$
consists of $N$ analytic, pairwise disjoint, simple closed
curves. Let $\eta$ be an analytic mapping from $R$ onto $U$,
the open unit disc which is holomorphic on a neighborhood of
$\overline R$ and which is of modulus one on $\Gamma$. That
such functions exists was first established by Ahlfors
[1](=Ahlfors 1950 \cite{Ahlfors_1950}); another proof of their
existence is in the paper [12](=Royden 1962
\cite{Royden_1962}).''

$\bullet$ Stout 1967 \cite{Stout_1967-Interpolation}: ``It is
convenient to make use of an {\it Ahlfors map\/} for $R$,
i.e., a  function continuous on $\overline R$ and holomorphic
in $R$ which is constantly of modulus one on $\Gamma$. The
existence of such function was established by Ahlfors in
[1](=Ahlfors 1950 \cite{Ahlfors_1950}); an alternative proof
of their existence is in [4](=Royden 1962
\cite{Royden_1962}).''

$\bullet$ O'Neill-Wermer 1968 \cite{O'Neill-Wermer_1968}:
``Let $W$ be a region on some Riemann surface whose boundary
is the union of a finite number of analytic simple closed
curves and with $W$ having compact closure. In ``Open Riemann
surfaces and extremal problems on compact subregions'',
(1950), L. Ahlfors considers the following extremal problem:

{\it Problem\/} I. {\it Let $a,b$ be points of $W$. among the
functions $F$ analytic on $W$ with $\vert F(z)\vert\le 1$ on
$W$ and $F(a)=0$, it is required to find the one which makes
$\vert F(b)\vert$ a maximum.}

He shows that this problem has a unique
solution\footnote{Presumably, the authors omit the rotational
ambiguity.} $f$ which maps $W$ in an $n$-to-$1$
fashion\footnote{Of course Ahlfors' statement is somewhat stronger
giving $r \le n\le r+2p$, where $r$ is the number of contours and
$p$ the genus.} onto the unit disk, for some $n$. His method of
proof depends on a certain associated extremal problem introduced
by P.\,R. Garabedian in his Thesis. (See Garabedian 1949
\cite{Garabedian_1949}). Another proof is given by H. Royden,
``The boundary values of analytic and harmonic functions,'' Math.
Z. 78 (1962), 1--24.''

$\bullet$ Stanton 1971 \cite[p.\,293]{Stanton_1971}: ``Our
argument rests on the following theorem of Ahlfors [1](=1950).
{\sc Theorem.} {\it There is a function $f$ which is analytic
on $W\cup \Gamma$ and which maps [the interior] $W$ onto $U$
and $\Gamma$ onto $T$.\/} This theorem is also proved in
Royden [7](=1962). A function $f$ of the kind described in
this theorem is called an {\it Ahlfors mapping}.''

$\bigstar$ Upon recalling, that Stanton is a Royden student
this may eventually be counted as a self-voting.

$\bullet$ Hejhal 1972 \cite[p.\,119]{Hejhal_1972}: ``Suppose
first of all that $W$ is the interior of a compact bordered
surface $\overline W$. L.~Ahlfors [2](=1950) and H. Royden
[24](=1962) have studied the present linear extremal problem
on such $W$ at least for the case $\chi \equiv {\rm constant}$
and $\frak L [f]=f(b)$ with $b\in W$. \dots''

$\bullet$ Gamelin 1973 \cite[p.\,3]{Gamelin_1973-Extremal-I}:
``\dots the paper of H.\,L. Royden deals with finite bordered
Riemann surfaces.''

$\bullet$ Gamelin 1973 \cite[p.\,1105]{Gamelin_1973-BAMS}:
``For dual extremal problems on Riemann surfaces, see
[2](=Ahlfors 1950 \cite{Ahlfors_1950}) and [36](=Royden 1962
\cite{Royden_1962}).''

$\bullet$ Fisher 1973 \cite[p.\,1183]{Fisher_1973}: ``A
similar problem [\dots] has been investigated by L.~Ahlfors
[A1], H. Royden [R], and others. In that case, the class of
competing function is convex, the solution is unique, is
analytic across the boundary $\Gamma$, and has modulus one on
$\Gamma$.'' And further on page 1187: ``Let $F$ be the
solution to the Ahlfors-Royden extremal problem described in
the introduction. \dots''

$\bullet$ Lund 1974 \cite[p.\,495]{Lund_1974}: ``Let $U$ be
the open unit disk in ${\Bbb C}$. We call $F$ an unimodular
function if $F$ is analytic in a neighborhood of
$\overline{R}$ and maps $\overline R$ onto $\overline U$ so
that $F$ is $n$-to-one if we count the multiplicity of $F$
where $dF$ vanishes. If $T$ is the unit circle, then $F$ maps
$\Gamma$ onto $T$. The existence of such a function was first
proved by Ahlfors [1](=1950). Later, Royden [4](=1962) gave
another proof of this result.''

$\bullet$ Kirsch 2005 \cite{Kirsch_2005}: ``Ahlfors
generalized Garabedian's result to regions on Riemann surfaces
[2](=Ahlfors 1950 \cite{Ahlfors_1950}); see Royden's paper
[159](=Royden 1962 \cite{Royden_1962}) for another treatment
as well as further references to the literature.''

\medskip
{\sc Other votes:}

$\bullet$ Alpay-Vinnikov 2000
\cite[p.\,240]{Alpay-Vinnikov_2000}: ``It has been shown by
Ahlfors [4](=Ahlfors 1950 \cite{Ahlfors_1950}) that such a
function [=ramified $n$-sheeted covering of the unit disk]
always exists, and it may be chosen to have the minimal
possible degree $g+1$; see also [5](=Alling-Greenleaf 1971
\cite{Alling-Greenleaf_1971}), [19](=Fay 1973
\cite{Fay_1973}), and [21](=Fedorov 1991
\cite{Fedorov_1991}).''

Apart from the fact that the writer (Gabard) does not adhere with
Alpay-Vinnikov's claim about $g+1$ being the minimal possible
degree for such a mapping ($g$ is the genus of the double, cf.
op.\,cit. p.\,230), the three proposed references are in our
opinion not perfectly adequate as substitute to Ahlfors 1950
\cite{Ahlfors_1950}. Alling-Greenleaf \cite[p.\,16,
Thm\,1.3.6]{Alling-Greenleaf_1971} only states Ahlfors' result yet
without reproving it, whereas both Fay and Fedorov recover the
result in the planar case only.

\subsection{Acknowledgements}

The author wishes to thank the following long list of
geometers (in chronological order of their interaction with
the writer in connection to the present text)

$\bullet$ Felice Ronga (ca. 1997/98 for his explanation of
Brusotti's theorem),

$\bullet$ Claude Weber, Michel Kervaire (for their explanations on
how to classify Klein's symmetric surfaces by looking at the
quotient bordered surface)

$\bullet$ Fr\'ed\'eric Bihan for pleasant discussions about real
algebraic geometry,

$\bullet$ Lee Rudolph (ca. 1999 for explaining to us what is the
natural topological model for a real elliptic curve with only one
``oval'', namely just a torus $S^1\times S^1$ acted upon by factor
permutation $(x,y)\mapsto (y,x)$ fixing thereby the diagonal
circle),

$\bullet$ Alexis Marin, Viatcheslav Kharlamov, Oleg Viro,
Jean-Jacques Risler, Thierry Vust, Michel Kervaire, Pierre de la
Harpe, John Steinig (for their comments and corrections improving
the shape of the article Gabard 2000 \cite{Gabard_2000})

$\bullet$ Ragahavan Narasimhan, Jacek Bochnack (ca. 1999 for {\it
not\/} having been in touch with Ahlfors' result of 1950
\cite{Ahlfors_1950} enabling me some free gestation about thinking
on the problem)

$\bullet$  Manfred Knebusch for his kind interest in the
modest work Gabard 2000 \cite{Gabard_2000},

$\bullet$ Johannes Huisman  for his constant interest
(2001--04--06), and his care about correcting  bugs in both my
Thesis and the article Gabard 2006 \cite{Gabard_2006},

$\bullet$ Sergei Finashin for an exciting discussion in Rennes
2001,

$\bullet$ Jean-Claude Hausmann (ca. 2000/01) for telling me about
the standard surjectivity criterion via the Brouwer degree, which
was decisive to complete Gabard 2006 \cite{Gabard_2006},

$\bullet$ Antonio Costa,  for  his fascinating talks in Geneva,

$\bullet$ Bujalance for his surely over-enthusiastic
Zentralblatt review of my article (Gabard
2006~\cite{Gabard_2006}),

$\bullet$ Fraser-Schoen, whose brilliant work revived my
interest in the theory of the Ahlfors' mapping (ca. the 13
March 2011) at a stage were I was mostly
sidetracked by ``non-metric manifolds'' thanks to efforts of
Mathieu Baillif and David Gauld.

$\bullet$
Stepan Orevkov, Oleg Viro (2011) for
their talks and pleasant discussions,

$\bullet$ Marc Coppens (2011--12) for e-mails, and  his
work on the separating gonality (2011 \cite{Coppens_2011})
adumbrating sharper insights on the degree of the Ahlfors function
(or rather the more general allied circle maps). His turning-point
result
appeals to a better conciliation of the analytic theory of Ahlfors
with the algebro-geometric viewpoint.

$\bullet$ (2011/12) Hugo Parlier, Peter Buser, Alexandre Girouard,
Gerhard Wanner und Martin Gander are acknowledged for their recent
e-mail exchanges.

$\bullet$ (Oct. 2012) Daniel Coray for enlarging the capacity of
my TeX-compilator, and for his lovely (Cambridge-style) teaching
about the geometric Galois action in ca. 1999.

$\bullet$ (16 Nov.\;2012) Franc Forstneri\v{c}
for pointing out his text with Wold (2012
\cite{Forstneric-Wold_2012}) showing the  state of the art on the
proper holomorphic embedding problem.

$\bullet$ (09 Jan. 2013) Oleg Viro for his excellent answers on
some naive questions on Rohlin's paper (cf.
Sec.\,\ref{e-mail-Viro:sec}).

$\bullet$ (10 Jan. 2013) Alexis Marin for his invaluable insights
on Klein's intuitions and
more (cf. Sec.\,\ref{e-mail-Viro:sec}).

$\bullet$ (January-February-March 2013) Viatcheslav Kharlamov,
Stepa Orevkov, Eugenii Shustin, S\'everine Le~Touz\'e, Thomas
Fiedler for all their letters reproduced in
Sec.\,\ref{e-mail-Viro:sec}.

\section{Bibliographic comments}
The writer does not pretend that the following bibliography is
complete (nor that he absorbed all those fantastic contributions
in full details). More extensive bibliographies (overlapping
ours), but covering more material include those of:

$\bullet$ Ahlfors-Sario 1960 \cite{Ahlfors-Sario_1960} (ca. 40
pages times 25 items per pages=1000 entries covering such topics
as the Dirichlet problem, extremal problems, the type problem, the
allied classification theory, etc.);

$\bullet$ Grunsky 1978 \cite{Grunsky_1978} (=562 refs,
including 48 Books).

Most entries of our bibliography are followed by some comment
explaining briefly the connection to our primary topic of the
Ahlfors map. The following symbolism is used:

$\clubsuit$ serves to point out a special connection to Ahlfors
1950 (especially alternative proofs).

$\spadesuit$ gives other comments (attempting to summarize the
paper contents or to explicit the connection in which we cite
it).

$\bigstar$ marks sources, I could not as yet procure a copy.

\def\Ah47{\textsf{A47}}

\def\A50{\textsf{A50}}

\def\AS60{\textsf{AS60}}
\def\G78{\textsf{G78}}

$\bullet$ the stickers/sigles \AS60, \G78 are assigned when the
source has already been cited in Ahlfors-Sario 1960
\cite{Ahlfors-Sario_1960} resp. Grunsky 1978 \cite{Grunsky_1978}.

$\bullet$ \A50 designates those references citing the paper
Ahlfors 1950 \cite{Ahlfors_1950} (there represents circa 106
articles on ``Google''), and occasionally \Ah47 those quoting
Ahlfors 1947 \cite{Ahlfors_1947}.

{\bf $\heartsuit$n} is something like the indicator of the US
rating agency (to be read ``liked by $n$''). It indicates the
cardinal number {\bf n} of citations of the paper as measured by
``Google Scholar''. The latter machine often misses
cross-citations, especially those in old books, or old articles
with references given in footnotes format. Many sources cited in
Grunsky's book (1978 \cite{Grunsky_1978}) are never cited
electronically. Accordingly, those rating numbers only supply a
statistical idea of the literature ramifications lying beyond a
given entry. Also
 low-citation articles are sometimes the most
 polished product ripe for museum entrance. Forelli 1979
\cite{Forelli_1979}  is typical: self-contained, elegant and
polished proof of Ahlfors result, yet  only rated by 3.

Our bibliography is somewhat conservative with comparatively few
modern references. Our excuse is two-fold: modern expressionism is
sometimes harder to grasp, and recent references are usually well
detected through computer search.

(Papers are listed in alphabetical, and then chronological order,
regardless of shared co-authorship.)

The primary focus is on  the Ahlfors map and the weaker (but more
general) circle maps. As a such the topic overflow slightly over
the territory of real algebraic geometry. Ahlfors-Sario's book
\AS60 address Riemann surfaces, whereas Grunsky's book \G78
focuses to the case of planar domains. Hence both bibliographies
\AS60, \G78 are quite complementary, and ours is essentially a
fusion of both, but we gradually included more and more recent
contribution. Still additional references are welcome.

For conformal maps, it is helpful when browsing the vast
literature to keep in mind the basic question: {\it what result
through which method?}

{\bf Results.} Objects traditionally range along increasing order
of generality through: simply-connected regions,
multiply-connected ones and finally Riemann surfaces. We often add
a humble compactness proviso, as the passage to open objects is
traditionally achieved through the exhaustion trick (going back at
least to Poincar\'e 1883 \cite{Poincare_1883}, and see also Koebe
1907 \cite{Koebe_1907_UbaK1}), and active in recent time (e.g.
Garabedian-Schiffer 1950 \cite{Garabedian-Schiffer_1950}.)

As to the mappings, they may all be interpreted in some way or
another as ramification of RMT (Riemann's mapping to a
circle=disc). We distinguish primarily:

$\bullet$ CM=circle maps (usually not univalent, but multi-sheeted
disc with branch, or winding points=Windungspunkte)

$\bullet$ KNP=Kreisnormierung(sprinzip) (univalent map to a
circular domain)

$\bullet$ SM=slit mappings for various types of them (parallel,
circular, radial, logarithmic spiral, etc.). Those are all allied
to certain natural foliation of the sphere, and some extreme
generality in this respect is achieved in Schramm's Thesis where
any foliation is permitted as support for the slits.

{\bf Methods.} They may be classified
in two broad classes quantitative vs. qualitative (each having
some branchings):

\noindent$\bigstar$ (Quantitative) variational methods,
including:

$\bullet$ DP=Dirichlet principle (or more broadly speaking,
potential theory=PT, centering around such concepts as the
Green's function, harmonic measures (i.e. harmonic function
with special null/one boundary prescription of the various
contour), etc. Of course, there is a standard yoga between
Dirichlet and Green, so all this is essentially one and the
same method.

$\bullet$ IM=Iterative methods (originators: Koebe and
Carath\'eodory), and by extension this may proliferate up to
including the circle packings.

$\bullet$ EP=extremal problems (e.g. the one of maximizing the
derivative amongst the class of function bounded-by-one) and
leads to the Ahlfors map.

$\bullet$ BK=Bergman kernel (or Szeg\"o kernel), here the
fundamental ideas rest upon Hilbert's space methods, and the
idea of
orthogonal system. Initially, the method is also inspired by Ritz,
and Bieberbach extremal problem (1914 \cite{Bieberbach_1914}) for
the area swept out by the function. Since the middle 1940's, there
were found several conformal identities among so-called domain
functions (Green's, Neumann's, etc.) and the kernel functions so
that virtually this is now highly connected to DP$\approx$PT. Also
the Ahlfors map is  expressible in term of the Bergman kernel (cf.
e.g., Nehari 1950 \cite{Nehari_1950}) so that this heading is
strongly connected to EP.

$\bullet$ PP=Plateau problem style methods (for RMT, this starts
with the observation of Douglas 1931
\cite{Douglas_1931-Solution}). This strongly allied to DP, albeit
some distinction is useful to keep in mind just for cataloguing
purposes.

\noindent$\bigstar$ (Qualitative) topological methods:

$\bullet$ the {\it continuity method}, as old as Schl\"afli,
(as Koebe notices somewhere) is involved in the accessory
parameters of Schwarz-Christoffel, in Klein-Poincar\'e's
uniformization through automorphic functions, Brouwer
(invariance of the domain), Koebe, etc., e.g. Golusin 1952/57
\cite{Golusin_1952/57})

$\bullet$ Brouwer topological degree and the allied
surjectivity criterion (cf. e.g., Mizumoto 1960
\cite{Mizumoto_1960}, Gabard 2006 \cite{Gabard_2006}). Here
the idea is that there is some topological stability of the
embedding of a curve into its Jacobian via the Abel mapping in
the sense that its homological feature are unsensitive to
variation of the complex (analytic) structure (moduli), and
this enables one to draw universal statement by purely
topological considerations.

Finally we have attempted to manufacture a genealogy map  showing
the affiliation between the authors. The picture turned out to be
so large that TeX prefers reject it at the very end of the file.

\begin{figure}[h]
\vskip-1.2cm\penalty0
    \hskip-4.2cm\penalty0\epsfig{figure=Geneal.eps,width=205mm}
\vskip-5pt\penalty0
  \caption{\label{Geneal:fig}%
  The authors involved in this bibliography (and genealogy)}
\vskip-5pt\penalty0
\end{figure}

\def\fakeref{\hskip-18pt [$\bigstar$]\;\;}

{\small

[15.10.12]  When I reached 884 references, I unfortunatel met the
so-called ``TeX capacity exceeded, sorry.'' obstruction (cf.
Knuth's ``The TeX Book'', p.\,300 for more details). Thus I had to
deactivate some references which are not used for cross-citation,
albeit they clearly belong to our topic.
%
[16.10.12] This problem was ultimately solved by
my advisor Daniel Coray, to whom I express my deepest gratitude
for enlarging the TeX capacity of my compilator.

}


{\small

}

\bigskip
\begin{minipage}[b]{0.6\linewidth} Alexandre
Gabard

Universit\'e de Gen\`eve

Section de Math\'ematiques

2-4 rue du Li\`evre, CP 64

CH-1211 Gen\`eve 4

Switzerland

alexandregabard@hotmail.com
\end{minipage}



\begin{thebibliography}{9}
\bibitem{Gabard-Gauld_2011-Dynamics}
A. Gabard, D. Gauld, {\it Dynamics of non-metric manifolds}, arXiv
2011. [$\spadesuit$ this is just cited for the proof of the
implication: simply-connected $\Rightarrow$ schlichtartig
$\Rightarrow$ orientable, which may be reduced to the five
lemma]~\hfill {\bf $\heartsuit$2}
\end{thebibliography}
\end{document}